\documentclass{article}

\usepackage{rotating}
\usepackage{tikz-cd}
\usepackage[all]{xy}
\usepackage{amssymb}
\usepackage{amsmath}
\usepackage{amsthm}
\newtheorem{thm}{Theorem}
\newtheorem{prop}{Proposition}
\newtheorem{lem}{Lemma}
\newtheorem{rem}{Remark}
\newtheorem{cor}{Corollary}
\newtheorem{defi}{Definition}
\newtheorem{defiprop}{Definition-Proposition}

\def\Set{\mathop{\rm Set}\nolimits}
\def\Top{\mathop{\rm Top}\nolimits}
\def\RTop{\mathop{\rm RTop}\nolimits}
\def\Sch{\mathop{\rm Sch}\nolimits}
\def\CH{\mathop{\rm CH}\nolimits}

\def\hom{\mathop{\rm hom}\nolimits}

\def\codim{\mathop{\rm codim}\nolimits}
\def\dim{\mathop{\rm dim}\nolimits}

\def\Tr{\mathop{\rm Tr}\nolimits}
\def\Cor{\mathop{\rm Cor}\nolimits}

\def\SmVar{\mathop{\rm SmVar}\nolimits}
\def\Gr{\mathop{\rm Gr}\nolimits}

\def\PSmVar{\mathop{\rm PSmVar}\nolimits}
\def\Var{\mathop{\rm Var}\nolimits}
\def\Hom{\mathop{\rm Hom}\nolimits}
\def\Spec{\mathop{\rm Spec}\nolimits}
\def\Coh{\mathop{\rm Coh}\nolimits}
\def\supp{\mathop{\rm supp}\nolimits}
\def\QPVar{\mathop{\rm QPVar}\nolimits}
\def\AnSp{\mathop{\rm AnSp}\nolimits}
\def\CW{\mathop{\rm CW}\nolimits}
\def\Cw{\mathop{\rm Cw}\nolimits}
\def\PVar{\mathop{\rm PVar}\nolimits}

\def\Tot{\mathop{\rm Tot}\nolimits}
\def\sing{\mathop{\rm sing}\nolimits}

\def\Im{\mathop{\rm Im}\nolimits}

\def\Cone{\mathop{\rm Cone}\nolimits}
\def\ad{\mathop{\rm ad}\nolimits}

\def\nul{\mathop{\rm nul}\nolimits}
\def\log{\mathop{\rm log}\nolimits}
\def\Diff{\mathop{\rm Diff}\nolimits}
\def\An{\mathop{\rm An}\nolimits}
\def\Ho{\mathop{\rm Ho}\nolimits}
\def\PSh{\mathop{\rm PSh}\nolimits}
\def\AnSm{\mathop{\rm AnSm}\nolimits}
\def\Tr{\mathop{\rm Tr}\nolimits}

\def\AnDM{\mathop{\rm AnDM}\nolimits}
\def\pt{\mathop{\rm pt}\nolimits}
\def\DA{\mathop{\rm DA}\nolimits}
\def\AnDA{\mathop{\rm AnDA}\nolimits}

\def\Fun{\mathop{\rm Fun}\nolimits}
\def\Sh{\mathop{\rm Sh}\nolimits}
\def\Ab{\mathop{\rm Ab}\nolimits}
\def\Bti{\mathop{\rm Bti}\nolimits}
\def\TM{\mathop{\rm TM}\nolimits}
\def\Ouv{\mathop{\rm Ouv}\nolimits}
\def\coker{\mathop{\rm coker}\nolimits}
\def\Ext{\mathop{\rm Ext}\nolimits}
\def\CS{\mathop{\rm CS}\nolimits}

\def\Cat{\mathop{\rm Cat}\nolimits}
\def\RCat{\mathop{\rm RCat}\nolimits}
\def\TriCat{\mathop{\rm TriCat}\nolimits}
\def\an{\mathop{\rm an}\nolimits}
\def\holim{\mathop{\rm holim}\nolimits}
\def\ho{\mathop{\rm ho}\nolimits}

\def\Ring{\mathop{\rm Ring}\nolimits}
\def\cRing{\mathop{\rm cRing}\nolimits}
\def\Mod{\mathop{\rm Mod}\nolimits}
\def\Der{\mathop{\rm Der}\nolimits}
\def\Shv{\mathop{\rm Shv}\nolimits}
\def\Proj{\mathop{\rm Proj}\nolimits}
\def\Inj{\mathop{\rm Inj}\nolimits}
\def\Cyl{\mathop{\rm Cyl}\nolimits}
\def\nul{\mathop{\rm nul}\nolimits}
\def\Dec{\mathop{\rm Dec}\nolimits}

\def\dar[#1]{\ar@<2pt>[#1]\ar@<-2pt>[#1]}

\title{The Hodge realization functor on the derived category of relative motives}

\author{Johann Bouali}

\begin{document}

\maketitle

\begin{abstract}
We give, for a complex algebraic variety $S$, a Hodge realization functor $\mathcal F_S^{Hdg}$
from the (unbounded) derived category of constructible motives $\DA_c(S)$ over $S$ 
to the (undounded) derived category $D(MHM(S))$ of algebraic mixed Hodge modules over $S$.
Moreover, for $f:T\to S$ a morphism of complex quasi-projective algebraic varieties,
$\mathcal F_{-}^{Hdg}$ commutes with the four operations $f^*,f_*,f_!,f^!$ on $\DA_c(-)$ and $D(MHM(-))$,
making in particular the Hodge realization functor a morphism of 2-functor 
on the category of complex quasi-projective algebraic varieties which for a given $S$ 
sends $\DA_c(S)$ to $D(MHM(S))$, moreover $\mathcal F_S^{Hdg}$ commutes with tensor product.
We also give an algebraic and analytic Gauss-Manin realization functor from which we
obtain a base change theorem for algebraic De Rham cohomology and for all smooth morphisms a relative 
version of the comparaison theorem of Grothendieck between the algebraic De Rahm cohomology and the analytic De Rahm cohomology.
\end{abstract}

\tableofcontents

\section{Introduction}

Saito's theory of mixed Hodge modules associate to each complex algebraic variety $S$ a category $MHM(S)$
which is a full subcategory of $\PSh_{\mathcal D(1,0)fil,rh}(S/(\tilde S_I))\times_IC_{fil}(S)$
which extend variations of mixed Hodge structure and admits a canonical monoidal structure
given by tensor product, and associate to each morphism of complex algebraic varieties
$f:X\to S$, four functor $Rf_{Hdg!},Rf_{Hdg*},f^{\hat*Hdg},f^{*Hdg}$. 
In the case of a smooth proper morphism $f:X\to S$ with $S$ and $X$ smooth,
$H^nRf_{Hdg*}\mathbb Z_X^{Hdg}$ is the variation of Hodge structure
given by the Gauss-Manin connexion and the local system $H^nRf_*\mathbb Z_X$. 
Moreover, these functors induce the six functor formalism of Grothendieck.
We thus have, for a complex algebraic variety $S$ a canonical functor 
\begin{equation*}
MH(/S):\Var(\mathbb C)/S\to D(MHM(S)), \; (f:X\to S)\mapsto Rf_{!Hdg}\mathbb Z_X^{Hdg}
\end{equation*}
and
\begin{equation*}
MH(/-):\Var(\mathbb C)\to\TriCat, \; S\mapsto (MH(S):\Var(\mathbb C)/S\to D(MHM(S))),
\end{equation*}
is a morphism of 2-functor. 
In this work, we extend $MH(/-)$ to motives by constructing, for each complex algebraic variety $S$,
a canonical functor $\mathcal F_S^{Hdg}:\DA(S)\to D(MHM(S))$ which is monoidal, that is commutes with tensor product,
together with, for each morphism of complex algebraic varieties $g:T\to S$ 
a canonical transformation map $T(g,\mathcal F^{Hdg})$, which make
\begin{equation*}
\mathcal F_{-}^{Hdg}:\Var(\mathbb C)\to\TriCat, \; S\mapsto (\mathcal F_S^{Hdg}:\DA(S)\to D(MHM(S))),
\end{equation*}
is a morphism of 2-functor : this is the contain of theorem \ref{main}. 
A partial result in this direction has been obtained by Ivorra in \cite{Ivorra} using a different approach.
We already have a Betti realization functor 
\begin{equation*}
\Bti_{-}:\Var(\mathbb C)\to\TriCat, \; S\mapsto (\Bti_S^*:\DA(S)\to D(S)),
\end{equation*}
which extend the Betti realization.
The functor $\mathcal F_{-}^{Hdg}:=(\mathcal F_{-}^{FDR},\Bti_{-})$ is obtained by constructing the De Rham part 
\begin{equation*}
\mathcal F_{-}^{FDR}:\Var(\mathbb C)\to\TriCat, \; 
S\mapsto (\mathcal F_S^{FDR}:\DA(S)\to D_{\mathcal D(1,0)fil,\infty}(S/(\tilde S_I))),
\end{equation*}
which takes values in the derived category of filtered algebraic $D$-modules 
$D_{\mathcal D(1,0)fil,\infty}(S/(\tilde S_I)):=K_{\mathcal D(1,0)fil,\infty}(S/(\tilde S_I)))([E_1]^{-1})$
obtained by inverting the classes filtered Zariski local equivalences $[E_1]$ modulo $\infty$-filtered homotopy
and then using the following key theorem (theorem \ref{Be})
\begin{thm}\label{BeDEF}
\begin{itemize}
\item[(i)]Let $S\in\Var(\mathbb C)$. Let $S=\cup_{i\in I}S_i$ an open cover such that there exists
closed embedding $i_i:S_i\hookrightarrow\tilde S_i$ with $\tilde S_i\in\SmVar(\mathbb C)$. 
Then the full embedding
\begin{eqnarray*}
\iota_S:MHM(S)\hookrightarrow C_{\mathcal D(1,0)fil,rh}(S/(\tilde S_I))\times_I D_{fil}(S^{an}) 
\end{eqnarray*}
induces a full embedding
\begin{equation*}
\iota_S:D(MHM(S))\hookrightarrow D_{\mathcal D(1,0)fil,rh}(S/(\tilde S_I))\times_I D_{fil}(S^{an}) 
\end{equation*}
whose image consists of 
$(((M_I,F,W),u_{IJ}),(K,W),\alpha)\in D_{\mathcal D(1,0)fil,rh}(S/(\tilde S_I))\times_I D_{fil}(S^{an})$ such that 
\begin{equation*}
((H^n(M_I,F,W),H^n(u_{IJ})),H^n(K,W),H^n\alpha)\in MHM(S) 
\end{equation*}
for all $n\in\mathbb Z$ and such that for all $p\in\mathbb Z$, 
the differential of $(\Gr_W^pM_I,F)$ are strict for the filtration $F$ 
(in particular, the differentials of $(M_I,F,W)$ are strict for the filtration $F$).
\item[(ii)]Let $S\in\Var(\mathbb C)$. Let $S=\cup_{i\in I}S_i$ an open cover such that there exists
closed embedding $i_i:S_i\hookrightarrow\tilde S_i$ with $\tilde S_i\in\SmVar(\mathbb C)$. 
Then the full embedding
\begin{eqnarray*}
\iota_S:MHM(S)\hookrightarrow C_{\mathcal D(1,0)fil,rh}(S/(\tilde S_I))\times_I D_{fil}(S^{an}) 
\end{eqnarray*}
induces a full embedding
\begin{equation*}
\iota_S:D(MHM(S))\hookrightarrow D_{\mathcal D(1,0)fil,\infty,rh}(S/(\tilde S_I))\times_I D_{fil}(S^{an}) 
\end{equation*}
whose image consists of 
$(((M_I,F,W),u_{IJ}),(K,W),\alpha)\in D_{\mathcal D(1,0)fil,\infty,rh}(S/(\tilde S_I))\times_I D_{fil}(S^{an})$ such that 
\begin{equation*}
((H^n(M_I,F,W),H^n(u_{IJ})),H^n(K,W),H^n\alpha)\in MHM(S) 
\end{equation*}
for all $n\in\mathbb Z$ and such that there exist $r\in\mathbb N$ and an $r$-homotopy equivalence
$((M_I,F,W),u_{IJ})\to((M'_I,F,W),u_{IJ})$ such that for all $p\in\mathbb Z$, 
the differential of $(\Gr_W^pM'_I,F)$ are strict for the filtration $F$ 
(in particular, the differentials of $(M'_I,F,W)$ are strict for the filtration $F$).
\end{itemize}
\end{thm}
Note that the category $D_{\mathcal D(1,0)fil,\infty,rh}(S/(\tilde S_I))$ is NOT triangulated.
More precisely the canonical triangles of $D_{\mathcal D(1,0)fil,\infty,rh}(S/(\tilde S_I))$ does NOT satisfy
the 2 of 3 axiom of a triangulated category. 
Moreover there exist canonical triangles of $D_{\mathcal D(1,0)fil,\infty,rh}(S/(\tilde S_I))$ which are NOT
the image of distinguish triangles of $\pi_S(D(MHM(S)))$.
This method can be seen as a relative version of the construction of F.Lecomte and N.Wach in \cite{LW}.  

In section 6.1.1 and 6.2.1, we construct an algebraic and analytic Gauss-Manin realization functor, 
but this functor does NOT give a complex of filtered $D$-module, BUT a complex of filtered $O$-modules whose cohomology 
sheaves have a structure of filtered $D_S$ modules. Hence, it does NOT get to the desired category. 
Moreover the Hodge filtration is NOT the right one : see proposition \ref{TGMFDRprop} and proposition \ref{FFXD} 
However this functor gives some interesting results.
Let $S\in\Var(\mathbb C)$ and $S=\cup_{i=1}^l S_i$ an open cover such that there exist closed embeddings
$i_i:S_i\hookrightarrow\tilde S_i$  with $\tilde S_i\in\SmVar(\mathbb C)$ smooth. 
For $I\subset\left[1,\cdots l\right]$, denote by $S_I=\cap_{i\in I} S_i$ and $j_I:S_I\hookrightarrow S$ the open embedding.
We then have closed embeddings $i_I:S_I\hookrightarrow\tilde S_I:=\Pi_{i\in I}\tilde S_i$.
We define the filtered algebraic Gauss-Manin realization functor defined as
\begin{eqnarray*}
\mathcal F_S^{GM}:C(\Var(\mathbb C)^{sm}/S)\to C_{Ofil,\mathcal D}(S/(\tilde S_I))^{\vee}, \; M\mapsto \\
\mathcal F_S^{GM}(F):=((e(\tilde S_I)_*\mathcal Hom^{\bullet}(L(i_{I*}j_I^*F),
E_{et}(\Omega^{\bullet}_{/\tilde S_I}),F_b))[-d_{\tilde S_I}],u^q_{IJ}(F)),
\end{eqnarray*} 
see definition \ref{DRalgdefFunctGM}.
Note that the canonical triangles of $D_{Ofil,\mathcal D,\infty}(S/(\tilde S_I))$ 
does NOT satisfy the 2 of 3 axiom of a triangulated category.
The filtered algebraic Gauss-Manin realization functor induces by proposition \ref{projwachGM}
\begin{eqnarray*}
\mathcal F_S^{GM}:\DA_c(S)^{op}\to D_{Ofil,\mathcal D,\infty}(S/(\tilde S_I)), \; M\mapsto \\
\mathcal F_S^{GM}(M):=((e(\tilde S_I)_*\mathcal Hom^{\bullet}(L(i_{I*}j_I^*F),
E_{et}(\Omega^{\bullet}_{/\tilde S_I}),F_b))[-d_{\tilde S_I}],u^q_{IJ}(F))
\end{eqnarray*}
where $F\in C(\Var(\mathbb C)^{sm}/S)$ is such that $M=D(\mathbb A^1,et)(F)$. 
We then prove (theorem \ref{mainthmGM}):
\begin{thm}\label{mainthmGMDEF}
\begin{itemize}
\item[(i)]Let $g:T\to S$ is a morphism with $T,S\in\Var(\mathbb C)$. 
Assume there exist a factorization $g:T\xrightarrow{l}Y\times S\xrightarrow{p_S}$
with $Y\in\SmVar(\mathbb C)$, $l$ a closed embedding and $p_S$ the projection. 
Let $S=\cup_{i=1}^lS_i$ be an open cover such that 
there exists closed embeddings $i_i:S_i\hookrightarrow\tilde S_i$ with $\tilde S_i\in\SmVar(\mathbb C)$.
Then, for $M\in\DA_c(S)$
\begin{eqnarray*}
T(g,\mathcal F^{GM})(M):Rg^{*mod[-],\Gamma}\mathcal F_S^{GM}(M)\to\mathcal F_T^{GM}(g^*M)
\end{eqnarray*}
is an isomorphism in $D_{O_Tfil,\mathcal D,\infty}(T/(Y\times\tilde S_I))$.
\item[(ii)]Let $g:T\to S$ is a morphism with $T,S\in\SmVar(\mathbb C)$. Then, for $M\in\DA_c(S)$
\begin{eqnarray*}
T^O(g,\mathcal F^{GM})(M):Lg^{*mod}\mathcal F_S^{GM}(M)\to\mathcal F_T^{GM}(g^*M)
\end{eqnarray*}
is an isomorphism in $D_{O_T}(T)$.
\item[(iii)] A base change theorem for algebraic De Rham cohomology :
Let $g:T\to S$ is a morphism with $T,S\in\SmVar(\mathbb C)$. 
Let $h:U\to S$ a smooth morphism with $U\in\Var(\mathbb C)$. Then the map (see definition \ref{TDw})
\begin{equation*}
T^O_w(g,h):Lg^{*mod}Rh_*(\Omega^{\bullet}_{U/S},F_b)\to Rh'_*(\Omega^{\bullet}_{U_T/T},F_b)
\end{equation*}
is an isomorphism in $D_{O_T}(T)$.
\end{itemize}
\end{thm}
Let $S\in\Var(\mathbb C)$ and $S=\cup_{i=1}^l S_i$ an open cover such that there exist closed embeddings
$i_i:S_i\hookrightarrow\tilde S_i$  with $\tilde S_i\in\SmVar(\mathbb C)$. 
For $I\subset\left[1,\cdots l\right]$, denote by $S_I=\cap_{i\in I} S_i$ and $j_I:S_I\hookrightarrow S$ the open embedding.
We then have closed embeddings $i_I:S_I\hookrightarrow\tilde S_I:=\Pi_{i\in I}\tilde S_i$.
We define the filtered analytic Gauss-Manin realization functor defined as
\begin{eqnarray*}
\mathcal F_S^{GM}:\DA_c(S)^{op}\to D_{Ofil,\mathcal D,\infty}(S/(\tilde S_I))^{\vee}, \; M\mapsto \\
\mathcal F_S^{GM}(M):=((e(\tilde S_I)_*\mathcal Hom^{\bullet}(\An_{\tilde S_I}^*L(i_{I*}j_I^*F),
E_{et}(\Omega^{\bullet}_{/\tilde S_I}),F_b))[-d_{\tilde S_I}],u^q_{IJ}(F))
\end{eqnarray*}
where $F\in C(\Var(\mathbb C)^{sm}/S)$ is such that $M=D(\mathbb A^1,et)(F)$, 
see definition \ref{DRalgdefFunctGMan}. We then prove (theorem \ref{mainthmGMan}):
\begin{thm}\label{mainthmGManDEF}
\begin{itemize}
\item[(i)]Let $g:T\to S$ is a morphism with $T,S\in\Var(\mathbb C)$. 
Assume there exist a factorization $g:T\xrightarrow{l}Y\times S\xrightarrow{p_S}$
with $Y\in\SmVar(\mathbb C)$, $l$ a closed embedding and $p_S$ the projection. 
Let $S=\cup_{i=1}^lS_i$ be an open cover such that 
there exists closed embeddings $i_i:S_i\hookrightarrow\tilde S_i$ with $\tilde S_i\in\SmVar(\mathbb C)$.
Then, for $M\in\DA_c(S)$
\begin{eqnarray*}
T(g,\mathcal F_{an}^{GM})(M):Rg^{*mod[-],\Gamma}\mathcal F_{S,an}^{GM}(M)\to\mathcal F_{T,an}^{GM}(g^*M)
\end{eqnarray*}
is an isomorphism in $D_{O_Tfil,\mathcal D^{\infty},\infty}(T/(Y\times\tilde S_I))$.
\item[(ii)]Let $g:T\to S$ is a morphism with $T,S\in\SmVar(\mathbb C)$. Then, for $M\in\DA_c(S)$
\begin{eqnarray*}
T(g,\mathcal F_{an}^{GM})(M):Lg^{*mod[-]}\mathcal F_{S,an}^{GM}(M)\to\mathcal F_{T,an}^{GM}(g^*M)
\end{eqnarray*}
is an isomorphism in $D_{O_T}(T)$.
\end{itemize}
\end{thm}
A consequence of the construction of the transformation map between the algebraic and analytic Gauss-Manin realization functor
is the following (theorem \ref{mainthmGAGAGM})
\begin{thm}\label{mainthmGAGAGMDEF}
\begin{itemize}
\item[(i)] Let $S\in\Var(\mathbb C)$. Then, for $M\in\DA_c(S)$
\begin{eqnarray*}
\mathcal J_S(-)\circ H^nT(\An,\mathcal F_{an}^{GM})(M):
J_S(H^n(\mathcal F_S^{GM}(M))^{an})\xrightarrow{\sim}H^n\mathcal F_{S,an}^{GM}(M)
\end{eqnarray*}
is an isomorphism in $\PSh_{\mathcal D}(S^{an}/(\tilde S_I^{an}))$.
\item[(ii)] A relative version of Grothendieck GAGA theorem for De Rham cohomology
Let $h:U\to S$ a smooth morphism with $S,U\in\SmVar(\mathbb C)$. Then,
\begin{equation*}
\mathcal J_S(-)\circ J_ST^O_{\omega}(an,h):
J_S((R^nh_*\Omega^{\bullet}_{U/S})^{an})\xrightarrow{\sim} R^nh_*\Omega^{\bullet}_{U^{an}/S^{an}}
\end{equation*}
is an isomorphism in $\PSh_{\mathcal D}(S^{an})$.
\end{itemize}
\end{thm}
In section 6.1.2, using results of sections 2, 4 and 5,
we construct the algebraic filtered De Rham realization functor $\mathcal F_{-}^{FDR}$.
We construct it via a larger category and use theorem \ref{Be}: 
Let $S\in\Var(\mathbb C)$ and $S=\cup_{i=1}^l S_i$ an open cover such that there exist closed embeddings
$i_i:S_i\hookrightarrow\tilde S_i$  with $\tilde S_i\in\SmVar(\mathbb C)$. 
For $I\subset\left[1,\cdots l\right]$, denote by $S_I=\cap_{i\in I} S_i$ and $j_I:S_I\hookrightarrow S$ the open embedding.
We then have closed embeddings $i_I:S_I\hookrightarrow\tilde S_I:=\Pi_{i\in I}\tilde S_i$.
we define in definition \ref{DRalgdefFunct}(ii) which use definition \ref{wtildew} and definition \ref{RCHhatdef}, 
the filtered algebraic De Rahm realization functor defined as
\begin{eqnarray*}
\mathcal F_S^{FDR}:C(\Var(\mathbb C)^{sm}/S)\to C_{\mathcal Dfil}(S/(\tilde S_I)), 
F\mapsto\mathcal F_S^{FDR}(F):= \\
(e'(\tilde S_I)_*\mathcal Hom^{\bullet}(\hat R^{CH}(\rho_{\tilde S_I}^*L(i_{I*}j_I^*F)),
E_{zar}(\Omega^{\bullet,\Gamma,pr}_{/\tilde S_I},F_{DR}))[-d_{\tilde S_I}],u^q_{IJ}(F)).
\end{eqnarray*}
By proposition \ref{projwach}(ii), it induces
\begin{eqnarray*}
\mathcal F_S^{FDR}:\DA_c(S)\to D_{\mathcal D(1,0)fil}(S/(\tilde S_I)), 
M\mapsto\mathcal F_S^{FDR}(M):= \\
(e'(\tilde S_I)_*\mathcal Hom^{\bullet}(\hat R^{CH}(\rho_{\tilde S_I}^*L(i_{I*}j_I^*(F,W))),
E_{zar}(\Omega^{\bullet,\Gamma,pr}_{/\tilde S_I},F_{DR}))[-d_{\tilde S_I}],u^q_{IJ}(F))
\end{eqnarray*}
where $F\in C(\Var(\mathbb C)^{sm}/S)$ is such that $M=D(\mathbb A^1,et)(F)$.
We compute this functor for an homological motive in proposition \ref{keyalgsing1} and we get by corollary \ref{FDRMHM},
for $S\in\Var(\mathbb C)$ and $M\in\DA_c(S)$, $\mathcal F_S^{FDR}(M)\in\pi_S(D(MHM(S))$,
and the following (theorem \ref{mainthm}):
\begin{thm}\label{mainthmDEF}
\begin{itemize}
\item[(i)]Let $g:T\to S$ a morphism, with $S,T\in\Var(\mathbb C)$. 
Assume we have a factorization $g:T\xrightarrow{l}Y\times S\xrightarrow{p_S}S$
with $Y\in\SmVar(\mathbb C)$, $l$ a closed embedding and $p_S$ the projection. Let $M\in\DA_c(S)$. 
Then map in $\pi_T(D(MHM(T)))$ 
\begin{eqnarray*}
T(g,\mathcal F^{FDR})(M):g^{\hat*mod}_{Hdg}\mathcal F_S^{FDR}(M)\xrightarrow{\sim}\mathcal F_T^{FDR}(g^*M)
\end{eqnarray*} 
given in definition \ref{TgDRdefsing} is an isomorphism.
\item[(ii)] Let $f:X\to S$ a morphism with $X,S\in\Var(\mathbb C)$. Assume there exist a factorization
$f:X\xrightarrow{l}Y\times S\xrightarrow{p_S}S$ with $Y\in\SmVar(\mathbb C)$, $l$ a closed embedding and $p_S$ the projection.
Then, for $M\in\DA_c(X)$, the map given in definition \ref{SixTalg}
\begin{equation*}
T_!(f,\mathcal F^{FDR})(M):Rf^{Hdg}_!\mathcal F_X^{FDR}(M)\xrightarrow{\sim}\mathcal F_S^{FDR}(Rf_!M)
\end{equation*}
is an isomorphism in $\pi_S(D(MHM(S))$.
\item[(iii)] Let $f:X\to S$ a morphism with $X,S\in\Var(\mathbb C)$, $S$ quasi-projective. Assume there exist a factorization
$f:X\xrightarrow{l}Y\times S\xrightarrow{p_S}S$ with $Y\in\SmVar(\mathbb C)$, $l$ a closed embedding and $p_S$ the projection.
We have, for $M\in\DA_c(X)$, the map given in definition \ref{SixTalg}
\begin{equation*}
T_*(f,\mathcal F^{FDR})(M):\mathcal F_S^{FDR}(Rf_*M)\xrightarrow{\sim}Rf^{Hdg}_*\mathcal F_X^{FDR}(M)
\end{equation*}
is an isomorphism in $\pi_S(D(MHM(S))$.
\item[(iv)] Let $f:X\to S$ a morphism with $X,S\in\Var(\mathbb C)$, $S$ quasi-projective. Assume there exist a factorization
$f:X\xrightarrow{l}Y\times S\xrightarrow{p_S}S$ with $Y\in\SmVar(\mathbb C)$, $l$ a closed embedding and $p_S$ the projection.
Then, for $M\in\DA_c(S)$, the map given in definition \ref{SixTalg}
\begin{equation*}
T^!(f,\mathcal F^{FDR})(M):\mathcal F_X^{FDR}(f^!M)\xrightarrow{\sim}f^{*mod}_{Hdg}\mathcal F_S^{FDR}(M)
\end{equation*}
is an isomorphism in $\pi_X(D(MHM(X))$. 
\item[(v)]Let $S\in\Var(\mathbb C)$ and $S=\cup_{i=1}^l S_i$ an open affine covering and denote, 
for $I\subset\left[1,\cdots l\right]$, $S_I=\cap_{i\in I} S_i$ and $j_I:S_I\hookrightarrow S$ the open embedding.
Let $i_i:S_i\hookrightarrow\tilde S_i$ closed embeddings, with $\tilde S_i\in\SmVar(\mathbb C)$. 
Then, for $M,N\in\DA_c(S)$, the map in $\pi_S(D(MHM(S)))$ 
\begin{eqnarray*}
T(\mathcal F_S^{FDR},\otimes)(M,N): 
\mathcal F_S^{FDR}(M)\otimes^{Hdg}_{O_S}\mathcal F_S^{FDR}(N)\xrightarrow{\sim}\mathcal F_S^{FDR}(M\otimes N)
\end{eqnarray*}
given in definition \ref{SixTalg} is an isomorphism.
\end{itemize}
\end{thm}
We also have a canonical transformation map between the Gauss-Manin and the De Rham functor given
in definition \ref{GMFDRdef} wich satisfy (see proposition \ref{TGMFDRprop}) :
\begin{prop}\label{TGMFDRpropDEF}
Let $S\in\Var(\mathbb C)$ and $S=\cup_{i=1}^l S_i$ an open cover such that there exist closed embeddings
$i_i:S_i\hookrightarrow\tilde S_i$ with $\tilde S_i\in\SmVar(\mathbb C)$. 
\begin{itemize}
\item[(i)] For $M\in\DA_c(S)$ the map in $D_{O_S,\mathcal D}(S/(\tilde S_I))=D_{O_S,\mathcal D}(S)$
\begin{equation*}
o_{fil}T(\mathcal F^{GM}_S,\mathcal F^{FDR}_S)(M):
o_{fil}\mathcal F_S^{GM}(L\mathbb D_SM)\xrightarrow{\sim} o_{fil}\mathcal F^{FDR}_S(M)
\end{equation*}
given in definition \ref{GMFDRdef} is an isomorphism if we forgot the Hodge filtration $F$.
\item[(ii)]For $M\in\DA_c(S)$ and all $n,p\in\mathbb Z$, the map in $\PSh_{O_S,\mathcal D}(S/(\tilde S_I))$
\begin{equation*}
F^pH^nT(\mathcal F^{GM}_S,\mathcal F^{FDR}_S)(M):F^pH^n\mathcal F_S^{GM}(L\mathbb D_SM)\hookrightarrow F^pH^n\mathcal F^{FDR}_S(M)
\end{equation*}
given in definition \ref{GMFDRdef} is a monomorphism.
Note that $F^pH^nT(\mathcal F^{GM}_S,\mathcal F^{FDR}_S)(M)$ is NOT an isomorphism in general :
take for example $M(S^o/S)^{\vee}=D(\mathbb A^1,et)(j_*E_{et}(\mathbb Z(S^o/S)))$ 
for an open embedding $j:S^o\hookrightarrow S$, then 
\begin{equation*}
\mathcal F_S^{GM}(L\mathbb D_SM(S^o/S)^{\vee})=\mathcal F_S^{GM}(\mathbb Z(S^o/S))=j_*E(O_{S^o},F_b)\notin\pi_S(MHM(S)) 
\end{equation*}
and hence NOT isomorphic to $\mathcal F_S^{GM}(L\mathbb D_SM(S^o/S)^{\vee})\in\pi_S(MHM(S))$, 
see remark \ref{remHdgkey}. 
It is an isomorphism in the very particular cases where $M=D(\mathbb A^1,et)(\mathbb Z(X/S))$
or $M=D(\mathbb A^1,et)(\mathbb Z(X^o/S))$ for $f:X\to S$ is a smooth proper morphism and $n:X^o\hookrightarrow X$ is an open subset
such that $X\backslash X^o=\cup D_i$ is a normal crossing divisor 
and such that $f_{|D_i}=f\circ i_i:D_i\to X$ are SMOOTH morphism with $i_i:D_i\hookrightarrow X$ the closed embedding and
considering $f_{|X^o}=f\circ n:X^o\to S$ (see proposition \ref{FFXD}).
\end{itemize}
\end{prop}

Let $S\in\Var(\mathbb C)$. Let $S=\cup_iS_i$ an open cover such that there exists
closed embedding $i_i:S\hookrightarrow\tilde S_i$ with $\tilde S_i\in\SmVar(\mathbb C)$.
We define the Hodge realization functor as, using definition \ref{DRalgdefFunct},
\begin{eqnarray*}
\mathcal F_S^{Hdg}:=(\mathcal F_S^{FDR},\Bti_S^*):
C(\Var(\mathbb C)^{sm}/S)\to D_{\mathcal D(1,0)fil}(S/(\tilde S_I))\times_I D_{fil}(S^{an}) 
\end{eqnarray*}
which is on objects given by, for $F\in C(\Var(\mathbb C)^{sm}/S)$, taking $(F,W)\in C_{fil}(\Var(\mathbb C)^{sm}/S)$
such that $D(\mathbb A^1,et)(F,W)$ gives the weight structure on $D(\mathbb A^1,et)(F)$,
\begin{eqnarray*}
\mathcal F_S^{Hdg}(F):=(\mathcal F_S^{FDR}(F),\Bti_S^*F,\alpha(F)):= \\
(e(S)_*\mathcal Hom((\hat R_{\tilde S_I}^{CH}(\rho_{\tilde S_I}^*Li_{I*}j_I^*(F,W)),
\hat R^{CH}(T^q(D_{IJ})(-))),(E_{zar}(\Omega^{\bullet,\Gamma,pr}_{/\tilde S_I},F_{DR}),T_{IJ})), \\
e(S)_*\underline{\sing}_{\mathbb D^*}\An_S^*L(F,W),\alpha(F)) 
\in D_{\mathcal D(1,0)fil}(S/(\tilde S_I))\times_I D_{fil}(S^{an}),
\end{eqnarray*}
where $\alpha(F)$ is given in definition \ref{HodgeRealDAsing}
and on morphism, for $m:F_1\to F_2$ with $F_1,F_2\in C(\Var(\mathbb C)^{sm}/S)$,
\begin{eqnarray*}
\mathcal F_S^{Hdg}(m):=(\mathcal F_S^{FDR}(m),\Bti_S^*(m),\theta(m)),
\end{eqnarray*} 
where $\theta(m)$ is given in definition \ref{HodgeRealDAsing}.
This functor induces by proposition \ref{projwach} the functor
\begin{eqnarray*}
\mathcal F_S^{Hdg}:=(\mathcal F_S^{FDR},\Bti_S^*):
\DA(S)\to D_{\mathcal D(1,0)fil}(S/(\tilde S_I))\times_I D_{fil}(S^{an}), \\ 
M=D(\mathbb A^1,et)(F)\mapsto
\mathcal F_S^{Hdg}(M):=\mathcal F_S^{Hdg}(F)=(\mathcal F_S^{FDR}(M),\Bti_S^*M,\alpha(M)),
\end{eqnarray*}
with $\alpha(M)=\alpha(F)$. 

The main result of this article is the following (theorem \ref{main}) :

\begin{thm}\label{mainDEF}
Let $k\subset\mathbb C$ be a subfield.
\begin{itemize}
\item[(i)] For $S\in\Var(\mathbb C)$, we have $\mathcal F_S^{Hdg}(\DA_c(S))\subset D(MHM(S))$,  
\begin{equation*}
\iota_S:D(MHM(S))\hookrightarrow D_{\mathcal D(1,0)fil}(S/(\tilde S_I))\times_I D_{fil}(S^{an}) 
\end{equation*}
being a full embedding by theorem \ref{Be}.
\item[(ii)] The Hodge realization functor $\mathcal F_{Hdg}(-)$ define a morphism of 2-functor on $\Var(\mathbb C)$
\begin{equation*}
\mathcal F^{Hdg}_{-}:\Var(\mathbb C)\to(\DA_c(-)\to D(MHM(-)))
\end{equation*}
whose restriction to $\QPVar(\mathbb C)$ is an homotopic 2-functor in sense of Ayoub. More precisely,
\begin{itemize}
\item[(ii0)] for $g:T\to S$ a morphism, with $T,S\in\QPVar(\mathbb C)$, and $M\in\DA_c(S)$, the
the maps of definition \ref{TgDRdefsing} and of definition \ref{TgBti} induce an isomorphism in $D(MHM(T))$
\begin{eqnarray*}
T(g,\mathcal F^{Hdg})(M):=(T(g,\mathcal F^{FDR})(M),T(g,bti)(M),0): \\
g^{\hat*Hdg}\mathcal F_S^{Hdg}(M):=\iota_T^{-1}(g^{\hat{*}mod}_{Hdg}\mathcal F_S^{FDR}(M),g^{*w}\Bti_S(M),g^*(\alpha(M))) \\
\xrightarrow{\sim}\iota_T^{-1}(\mathcal F_T^{FDR}(g^*M),\Bti_T^*(g^*M),\alpha(g^*M))=:\mathcal F_T^{Hdg}(g^*M),
\end{eqnarray*} 
\item[(ii1)] for $f:T\to S$ a morphism, with $T,S\in\QPVar(\mathbb C)$, and $M\in\DA_c(T)$,  
the maps of definition \ref{SixTalg} and of definition \ref{TBtiSix} induce an isomorphism in $D(MHM(S))$
\begin{eqnarray*}
T_*(f,\mathcal F^{Hdg})(M):=(T_*(f,\mathcal F^{FDR})(M),T_*(f,bti)(M),0): \\
Rf_{Hdg*}\mathcal F_T^{Hdg}(M):=\iota_S^{-1}(Rf^{Hdg}_*\mathcal F_T^{FDR}(M),Rf_{*w}\Bti_T(M),f_*(\alpha(M))) \\ 
\xrightarrow{\sim}\iota_S^{-1}(\mathcal F_S^{FDR}(Rf_*M),\Bti_S^*(Rf_*M),\alpha(Rf_*M))=:\mathcal F_S^{Hdg}(Rf_*M),
\end{eqnarray*}  
\item[(ii2)] for $f:T\to S$ a morphism, with $T,S\in\QPVar(\mathbb C)$, and $M\in\DA_c(T)$, 
the maps of definition \ref{SixTalg} and of definition \ref{TBtiSix} induce an isomorphism in $D(MHM(S))$
\begin{eqnarray*}
T_!(f,\mathcal F^{Hdg})(M):=(T_!(f,\mathcal F^{FDR})(M),T_!(f,bti)(M),0): \\
Rf_{!Hdg}\mathcal F_T^{Hdg}(M):=\iota_S^{-1}(Rf^{Hdg}_!\mathcal F_T^{FDR}(M),Rf_{!w}\Bti_T^*(M),f_!(\alpha(M))) \\ 
\xrightarrow{\sim}\iota_S^{-1}(\mathcal F_S^{FDR}(Rf_!M),\Bti_S^*(Rf_!M),\alpha(f_!M))=:\mathcal F_S^{Hdg}(Rf_!M),
\end{eqnarray*} 
\item[(ii3)] for $f:T\to S$ a morphism, with $T,S\in\QPVar(\mathbb C)$, and $M\in\DA_c(S)$,
the maps of definition \ref{SixTalg} and of definition \ref{TBtiSix} induce an isomorphism in $D(MHM(T))$
\begin{eqnarray*}
T^!(f,\mathcal F^{Hdg})(M):=(T^!(f,\mathcal F^{FDR})(M),T^!(f,bti)(M),0): \\
f^{*Hdg}\mathcal F_S^{Hdg}(M):=\iota_T^{-1}(f^{*mod}_{Hdg}\mathcal F_S^{FDR}(M),f^{!w}\Bti_S(M),f^!(\alpha(M))) \\ 
\xrightarrow{\sim}\iota_T^{-1}(\mathcal F_T^{FDR}(f^!M),\Bti_T^*(f^!M),\alpha(f^!M))=:\mathcal F_T^{Hdg}(f^!M),
\end{eqnarray*}
\item[(ii4)] for $S\in\Var(\mathbb C)$, and $M,N\in\DA_c(S)$,
the maps of definition \ref{SixTalg} and of definition \ref{TBtiSix} 
induce an isomorphism in $D(MHM(S))$
\begin{eqnarray*}
T(\otimes,\mathcal F^{Hdg})(M,N):=(T(\otimes,\mathcal F_S^{FDR})(M,N),T(\otimes,bti)(M,N),0): \\
\iota_S^{-1}(\mathcal F_S^{FDR}(M)\otimes^{Hdg}_{O_S}\mathcal F_S^{FDR}(N),
\Bti_S(M)\otimes\Bti_S(N),\alpha(M)\otimes\alpha(N)) \\
\xrightarrow{\sim}\mathcal F_S^{Hdg}(M\otimes N):=
\iota_S^{-1}(\mathcal F_S^{FDR}(M\otimes N),\Bti_S(M\otimes N),\alpha(M\otimes N)).  
\end{eqnarray*}
\end{itemize}
\item[(iii)] For $S\in\Var(\mathbb C)$, the following diagram commutes :
\begin{equation*}  
\xymatrix{\Var(\mathbb C)/S\ar[rrr]^{MH(/S)}\ar[d]_{M(/S)} & \, & \, & D(MHM(S))\ar[d]^{\iota^S} \\
\DA(S)\ar[rrr]^{\mathcal F_S^{Hdg}} & \, & \, & D_{\mathcal D(1,0)fil}(S/(\tilde S_I))\times_I D_{fil}(S^{an})}
\end{equation*}
\end{itemize}
\end{thm}
We obtain theorem \ref{mainDEF} from theorem \ref{mainthmDEF} and from the result on the Betti factor
after checking the compatibility of these transformation maps with the isomorphisms $\alpha(M)$.

I am grateful to F.Mokrane for his help and support during the preparation of this work 
as well as J.Wildeshaus for the interest and remarks that he made on a first version of this text.
I also thank J.Ayoub, C.Sabbah and M.Saito for the interest they have brought to this work.

\section{Generalities and Notations}

\subsection{Notations}

\begin{itemize}

\item After fixing a universe, we denote by
\begin{itemize}
\item $\Set$ the category of sets, 
\item $\Top$ the category of topological spaces, 
\item $\Ring$ the category of rings and $\cRing\subset\Ring$ the full subscategory of commutative rings,
\item $\RTop$ the category of ringed spaces, 
\begin{itemize}
\item whose set of objects is
$\RTop:=\left\{(X,O_X), \; X\in\Top, \, O_X\in\PSh(X,\Ring)\right\}$
\item whose set of morphism is 
$\Hom((T,O_T),(S,O_S)):=\left\{((f:T\to S),(a_f:f^*O_S\to O_T))\right\}$
\end{itemize}
and by $ts:\RTop\to\Top$ the forgetfull functor.
\item $\Cat$ the category of small categories which comes with the forgetful functor $o:\Cat\to\Fun(\Delta^1,\Set)$, 
where $\Fun(\Delta^1,\Set)$ is the category of simplicial sets, 
\item $\RCat$ the category of ringed topos
\begin{itemize}
\item whose set of objects is 
$\RCat:=\left\{(\mathcal X,O_X), \; \mathcal X\in\Cat, \, O_X\in\PSh(\mathcal X,\Ring)\right\}$,
\item whose set of morphism is 
$\Hom((\mathcal T,O_T),(\mathcal S,O_S)):=\left\{((f:\mathcal T\to \mathcal S),(a_f:f^*O_S\to O_T)),\right\}$
\end{itemize}
and by $tc:\RCat\to\Cat$ the forgetfull functor.
\end{itemize}

\item Let $F:\mathcal C\to\mathcal C'$ be a functor with $\mathcal C,\mathcal C'\in\Cat$. 
For $X\in\mathcal C$, we denote by $F(X)\in\mathcal C'$ the image of $X$, 
and for $X,Y\in\mathcal C$,  we denote by $F^{X,Y}:\Hom(X,Y)\to\Hom(F(X),F(Y))$ the corresponding map.

\item For $\mathcal C\in\Cat$, we denote by $\mathcal C^{op}\in\Cat$ the opposite category whose set of object
is the one of $\mathcal C$ : $(\mathcal C^{op})^0=\mathcal C^0$, and whose morphisms are the morphisms of $\mathcal C$
with reversed arrows.

\item Let $\mathcal C\in\Cat$. For $S\in\mathcal C$, we denote by $\mathcal C/S$ the category 
\begin{itemize}
\item whose set of objects $(\mathcal C/S)^0=\left\{X/S=(X,h)\right\}$ consist of the morphisms $h:X\to S$ with $X\in\mathcal C$,
\item whose set of morphism $\Hom(X'/S,X/S)$ between $X'/S=(X',h'),X/S=(X,h)\in\mathcal C/S$ 
consits of the morphisms $(g:X'\to X)\in\Hom(X',X)$ such that $h\circ g=h'$. 
\end{itemize}
We have then, for $S\in\mathcal C$, the canonical forgetful functor 
\begin{eqnarray*}
r(S):\mathcal C/S\to\mathcal C, \; \; X/S\mapsto r(S)(X/S)=X, \; (g:X'/S\to X/S)\mapsto r(S)(g)=g
\end{eqnarray*}
and we denote again $r(S):\mathcal C\to\mathcal C/S$ the corresponding morphism of (pre)sites.
\begin{itemize}
\item Let $F:\mathcal C\to\mathcal C'$ be a functor with $\mathcal C,\mathcal C'\in\Cat$. 
Then for $S\in\mathcal C$, we have the canonical functor 
\begin{eqnarray*}
F_S:\mathcal C/S\to\mathcal C'/F(S), \; \; X/S\mapsto F(X/S)=F(X)/F(S), \\ (g:X'/S\to X/S)\mapsto (F(g):F(X')/F(S)\to F(X)/F(S))
\end{eqnarray*}
\item Let $\mathcal S\in\Cat$. Then, for a morphism $f:X'\to X$ with $X,X'\in\mathcal S$
we have the functor  
\begin{eqnarray*}
C(f):\mathcal S/X'\to\mathcal S/X, \; \; Y/X'=(Y,f_1)\mapsto C(f)(Y/X'):=(Y,f\circ f_1)\in\mathcal S/X, \\ 
(g:Y_1/X'\to Y_2/X')\mapsto (C(f)(g):=g:Y_1/X\to Y_2/X)
\end{eqnarray*}
\item Let $\mathcal S\in\Cat$ a category which admits fiber products. 
Then, for a morphism $f:X'\to X$ with $X,X'\in\mathcal S$, we have the pullback functor
\begin{eqnarray*}
P(f):\mathcal S/X\to\mathcal S/X', \; \; Y/X\mapsto P(f)(Y/X):=Y\times_X X'/X'\in\mathcal S/X', \\ 
(g:Y_1/X\to Y_2/X)\mapsto (P(f)(g):=(g\times I):Y_1\times_X X'\to Y_2\times_X X')
\end{eqnarray*} 
which is right adjoint to $C(f):\mathcal S/X'\to\mathcal S/X$,
and we denote again $P(f):\mathcal S/X'\to\mathcal S/X$ the corresponding morphism of (pre)sites.
\end{itemize}

\item Let $\mathcal C,\mathcal I\in\Cat$. Assume that $\mathcal C$ admits fiber products.
For $(S_{\bullet})\in\Fun(\mathcal I^{op},\mathcal C)$, 
we denote by $\mathcal C/(S_{\bullet})\in\Fun(\mathcal I,\Cat)$ the diagram of category given by
\begin{itemize}
\item for $I\in\mathcal I$, $\mathcal C/(S_{\bullet})(I):=\mathcal C/S_I$,
\item for $r_{IJ}:I\to J$, $\mathcal C/(S_{\bullet})(r_{IJ}):=P(r_{IJ}):\mathcal C/S_I\to\mathcal C/S_J$,
where we denoted again $r_{IJ}:S_J\to S_I$ the associated morphism in $\mathcal C$. 
\end{itemize}

\item Let $(F,G):\mathcal C\leftrightarrows\mathcal C'$ an adjonction between two categories.
\begin{itemize}
\item For $X\in C$ and $Y\in C'$, we consider the adjonction isomorphisms
\begin{itemize}
\item $I(F,G)(X,Y):\Hom(F(X),Y)\to\Hom(X,G(Y)), \; (u:F(X)\to Y)\mapsto (I(F,G)(X,Y)(u):X\to G(Y))$
\item $I(F,G)(X,Y):\Hom(X,G(Y))\to\Hom(F(X),Y), \; (v:X\to G(Y))\mapsto (I(F,G)(X,Y)(v):F(X)\to Y)$. 
\end{itemize}
\item For $X\in\mathcal C$, we denote by 
$\ad(F,G)(X):=I(F,G)(X,F(X))(I_{F(X)}):X\to G\circ F(X)$. 
\item For $Y\in\mathcal C'$ we denote also by 
$\ad(F,G)(Y):=I(F,G)(G(Y),Y)(I_{G(Y)}):F\circ G(Y)\to Y$.
\end{itemize}
Hence, 
\begin{itemize}
\item for $u:F(X)\to Y$ a morphism with $X\in C$ and $Y\in C'$, we have 
$I(F,G)(X,Y)(u)=G(u)\circ\ad(F,G)(X)$,
\item for $v:X\to G(Y)$ a morphism with $X\in C$ and $Y\in C'$, we have 
$I(F,G)(X,Y)(v)=\ad(F,G)(Y)\circ F(v)$.
\end{itemize}

\item Let $\mathcal C$ a category. 
\begin{itemize}
\item We denote by $(\mathcal C,F)$ the category of filtered objects :
$(X,F)\in(\mathcal C,F)$ is a sequence $(F^{\bullet}X)_{\bullet\in\mathbb Z}$ indexed by $\mathbb Z$
with value in $\mathcal C$ together with monomorphisms $a_p:F^pX\hookrightarrow F^{p-1}X\hookrightarrow X$.
\item We denote by $(\mathcal C,F,W)$ the category of bifiltered objects : $(X,F,W)\in(\mathcal C,F,W)$ 
is a sequence $(W^{\bullet}F^{\bullet}X)_{\bullet,\bullet}\in\mathbb Z^2$ indexed by $\mathbb Z^2$ with value in $\mathcal C$ 
together with monomorphisms $W^qF^pX\hookrightarrow F^{p-1}X$, $W^qF^pX\hookrightarrow W^{q-1}F^pX$.
\end{itemize}

\item For $\mathcal C$ a category and $\Sigma:\mathcal C\to\mathcal C$ an endofunctor, we denote
by $(\mathcal C,\Sigma)$ the corresponding category of spectra, whose objects are sequence of objects of $\mathcal C$
$(T_i)_{i\in\mathbb Z}\in\Fun(\mathbb Z,\mathcal C)$ together with morphisms $s_i:T_i\to\Sigma T_{i+1}$, and whose
morphism from $(T_i)$ to $(T'_i)$ are sequence of morphisms $T_i\to T'_i$ which commutes with the $s_i$.

\item Let $\mathcal A$ an additive category.
\begin{itemize}
\item We denote by $C(\mathcal A):=\Fun(\mathbb Z,\mathcal A)$ 
the category of (unbounded) complexes with value in $\mathcal A$, 
where we have denoted $\mathbb Z$ the category whose set of objects is $\mathbb Z$,
and whose set of morphism between $m,n\in\mathbb Z$ consists of one element (identity) if $n=m$,
of one elemement if $n=m+1$ and is $\emptyset$ in the other cases. 
\item We have the full subcategories $C^b(\mathcal A)$, $C^-(\mathcal A)$, $C^+(\mathcal A)$ of 
$C(\mathcal A)$ consisting of bounded, resp. bounded above, resp. bounded below complexes.
\item We denote by $K(\mathcal A):=\Ho(C(\mathcal A))$ the homotopy category of $C(\mathcal A)$
whose morphisms are equivalent homotopic classes of morphism and by
$Ho:C(\mathcal A)\to K(\mathcal A)$ the full homotopy functor. 
The category $K(\mathcal A)$ is in the standard way a triangulated category.
\end{itemize}

\item Let $\mathcal A$ an additive category.
\begin{itemize}
\item We denote by $C_{fil}(\mathcal A)\subset(C(\mathcal A),F)=C(\mathcal A,F)$ the full
additive subcategory of filtered complexes of $\mathcal A$ such that the filtration is biregular : 
for $(A^{\bullet},F)\in(C(\mathcal A),F)$, we say that $F$ is biregular if $F^{\bullet}A^r$ is finite for all $r\in\mathbb Z$.
\item We denote by $C_{2fil}(\mathcal A)\subset(C(\mathcal A),F,W)=C(\mathcal A,F,W)$ the full
subcategory of bifiltered complexes of $\mathcal A$ such that the filtration is biregular.
\item For $A^{\bullet}\in C(\mathcal A)$, we denote by $(A^{\bullet},F_b)\in(C(\mathcal A),F)$ the complex endowed
with the trivial filtration (filtration bete) : $F^pA^n=0$ if $p\geq n+1$ and $F^pA^n=A^n$ if $p\leq n$. 
Obviously, a morphism $\phi:A^{\bullet}\to B^{\bullet}$, with $A^{\bullet},B^{\bullet}\in C(\mathcal A)$
induces a morphism $\phi:(A^{\bullet},F_b)\to (B^{\bullet},F_b)$.
\item For $(A^{\bullet},F)\in C(\mathcal A,F)$, we denote by $(A^{\bullet},F(r))\in C(\mathcal A,F)$ the filtered
complex where the filtration is given by $F^p(A^{\bullet},F(r)):=F^{p+r}(A^{\bullet},F)$.
\item Two morphisms $\phi_1,\phi_2:(M,F)\to(N,F)$ with $(M,F),(N,F)\in C(\mathcal A,F)$ are said to be $r$-filtered homotopic
if there exist a morphism in $\Fun(\mathbb Z,(\mathcal A,F))$
\begin{equation*}
h:(M,F(r-1))[1]\to(N,F), h:=(h^n:(M^{n+1},F(r-1))\to (N^n,F))_{n\in\mathbb Z}, 
\end{equation*}
where $\mathbb Z$ have only trivial morphism (i.e. $h$ is a graded morphism  but not a morphism of complexes) 
such that $d'h+hd=\phi_1-\phi_2$, where $d$ is the differential of $M$ and $d'$ is the differential of $N$, 
and we have $h(F^pM^{n+1})\subset F^{p-r+1}N^n$, note that by definition $r$ does NOT depend on $p$ and $n$ ; we then say that 
\begin{equation*}
(h,\phi_1,\phi_2):(M,F)[1]\to(N,F) 
\end{equation*}
is an $r$-filtered homotopy.
By definition, an $r$-filtered homotopy $(h,\phi_1,\phi_2):(M,F)[1]\to(N,F)$ is an $r'$-filtered homotopy for all $r'\geq r$,
and a $1$-filtered homotopy is an homotopy of $C(\mathcal A,F)$. 
By definition, an $r$-filtered homotopy $(h,\phi_1,\phi_2):(M,F)[1]\to(N,F)$ gives
if we forgot filtration an homotopy $(h,\phi_1,\phi_2):M[1]\to N$ in $C(\mathcal A)$.
\item Two morphisms $\phi_1,\phi_2:(M,F)\to(N,F)$ with $(M,F),(N,F)\in C(\mathcal A,F)$ are said to be $\infty$-filtered homotopic
if there exist $r\in\mathbb N$ such that $\phi_1,\phi_2:(M,F)\to(N,F)$ are $r$-filtered homotopic.
Hence, if $\phi_1,\phi_2:(M,F)\to(N,F)$ are $\infty$-filtered homotopic, then $\phi_1,\phi_2:M\to N$ are homotopic ;
of course the converse is NOT true since $r$ does NOT depend on $p,n\in\mathbb Z$.
\item We will use the fact that by definition if $\phi:M\to N$ with $M,N\in C(\mathcal A)$ is an homotopy equivalence,
then $\phi:(M,F_b)\to\phi(N,F_b)$ is a $2$-filtered homotopy equivalence.
\item A morphism $\phi:(M,F)\to(N,F)$ with $(M,F),(N,F)\in C(\mathcal A,F)$ is said to be an $r$-filtered homotopy
equivalence if there exist a morphism $\phi':(N,F)\to(M,F)$ such that
\begin{itemize}
\item $\phi'\circ\phi:(M,F)\to(M,F)$ is $r$-filtered homotopic to $I_M$ and 
\item $\phi\circ\phi':(N,F)\to(N,F)$ is $r$-filtered homotopic to $I_N$.
\end{itemize}
If $\phi:(M,F)\to(N,F)$ is an $r$-filtered homotopy equivalence, then it is an $s$-filtered homotopy equivalence for $s\geq r$.
If $\phi:(M,F)\to(N,F)$ is an $r$-filtered homotopy equivalence, $\phi:M\to N$ is an homotopy equivalence
\item A morphism $\phi:(M,F)\to(N,F)$ with $(M,F),(N,F)\in C(\mathcal A,F)$ is said to be an $\infty$-filtered homotopy
equivalence if there exist $r\in\mathbb Z$ such that
$\phi:(M,F)\to(N,F)$ with $(M,F),(N,F)\in C(\mathcal A,F)$ is an $r$-filtered homotopy equivalence.
If $\phi:(M,F)\to(N,F)$ is an $\infty$-filtered homotopy equivalence, $\phi:M\to N$ is an homotopy equivalence ;
the converse is NOT true since $r$ does NOT depend on $p,n\in\mathbb Z$.
\item We denote by $K_r(\mathcal A,F):=\Ho_r(C(\mathcal A,F))$ the homotopy category of $C(\mathcal A,F)$
whose objects are those of $C(\mathcal A,F)$ 
and whose morphisms are the morphisms of $C(\mathcal A,F)$ modulo $r$-filtered homotopies,
and by $\Ho_r:C(\mathcal A,F)\to K_r(\mathcal A,F)$ the full homotopy functor. 
However, for $r>1$ the category $K_r(\mathcal A,F)$ with the canonical triangles the standard ones does NOT satify the 2 of 3
axiom of a triangulated category. Indeed, for $r>1$, a commutative diagram in $K_r(\mathcal A,F)$
\begin{equation*}
\xymatrix{(A^{\bullet},F)\ar[r]^{m}\ar[d]_{\phi} & (B^{\bullet},F)\ar[r]^{i_2}\ar[d]_{\psi} & 
\Cone(m)=((A^{\bullet},F)[1]\oplus(B^{\bullet},F),d,d'-m)\ar[r]^{p_1} & 
(A^{\bullet},F)[1]\ar[d]^{\phi[1]} \\
(A^{'\bullet},F)\ar[r]^{m'} & (B^{'\bullet},F)\ar[r]^{i_2} & 
\Cone(m')=((A^{'\bullet},F)[1]\oplus(B^{'\bullet},F),d,d'-m')\ar[r]^{p_1} & (A^{'\bullet},F)[1]},
\end{equation*}
i.e. which commutes modulo $r$-filtered homotopy can NOT be completed with a third vertical arrow.
In particular, for $r>1$ and commutative diagram in $C(\mathcal A,F)$
\begin{equation*}
\xymatrix{(A^{\bullet},F)\ar[r]^{m}\ar[d]_{\phi} & (B^{\bullet},F)\ar[r]^{i_2}\ar[d]_{\psi} & 
\Cone(m)=((A^{\bullet},F)[1]\oplus(B^{\bullet},F),d,d'-m)\ar[r]^{p_1}\ar[d]^{(\phi[1],\psi)} & 
(A^{\bullet},F)[1]\ar[d]^{\phi[1]} \\
(A^{'\bullet},F)\ar[r]^{m'} & (B^{'\bullet},F)\ar[r]^{i_2} & 
\Cone(m')=((A^{'\bullet},F)[1]\oplus(B^{'\bullet},F),d,d'-m')\ar[r]^{p_1} & (A^{'\bullet},F)[1]},
\end{equation*}
if $\psi$ and $\phi$ are $r$-filtered homotopy equivalence, then $(\phi[1],\psi)$ is NOT necessary
an $r$-filtered homotopy equivalence.
\item We denote by $K_{fil,r}(\mathcal A):=\Ho_r(C_{fil}(\mathcal A))$ the homotopy category of $C_{fil}(\mathcal A)$
whose objects are those of $C_{fil}(\mathcal A,)$ 
and whose morphisms are the morphisms of $C_{fil}(\mathcal A)$ modulo $r$-filtered homotopies,
and by $\Ho_r:C_{fil}(\mathcal A)\to K_{fil,r}(\mathcal A)$ the full homotopy functor. 
However, for $r>1$ the category $K_{fil,r}(\mathcal A)$ with the canonical triangles the standard ones does NOT satify the 2 of 3
axiom of a triangulated category.
\item We denote by $K_{fil,\infty}(\mathcal A):=\Ho_{\infty}(C_{fil}(\mathcal A))$ the homotopy category of $C_{fil}(\mathcal A)$
whose objects are those of $C_{fil}(\mathcal A,)$ 
and whose morphisms are the morphisms of $C_{fil}(\mathcal A)$ modulo $\infty$-filtered homotopies,
and by $\Ho_{\infty}:C_{fil}(\mathcal A)\to K_{fil,\infty}(\mathcal A)$ the full homotopy functor. 
However, the category $K_{fil,\infty}(\mathcal A)$ with the canonical triangles the standard ones does NOT satify the 2 of 3
axiom of a triangulated category.
\item We have the Deligne decalage functor
\begin{eqnarray*}
\Dec: C(\mathcal A,F)\to C(\mathcal A,F), \; (M,F)\mapsto \Dec(M,F):=(M,\Dec F), \\
\Dec F^pM^n:=F^{p+n}M^n\cap d^{-1}(F^{p+n+1}M^{n+1})
\end{eqnarray*}
It is the right adjoint of the shift functor
\begin{eqnarray*}
S:C(\mathcal A,F)\to C(\mathcal A,F), \; (M,F)\mapsto S(M,F):=(M,SF), \; SF^pM^n:=F^{p-n}M^n
\end{eqnarray*}
The dual decalage functor
\begin{eqnarray*}
\Dec^{\vee}: C(\mathcal A,F)\to C(\mathcal A,F), \; (M,F)\mapsto \Dec^{\vee}(M,F):=(M,\Dec^{\vee} F), \\
\Dec^{\vee} F^pM^n:=F^{p+n}M^n+d(F^{p+n-1}M^{n+1})
\end{eqnarray*}
is the left adjoint of the shift functor.
Note that $\Dec((M,F)[1])\neq(\Dec(M,F))[1]$, $\Dec^{\vee}((M,F)[1])\neq(\Dec^{\vee}(M,F))[1]$
and $S((M,F)[1])\neq(S(M,F))[1]$.
\end{itemize}

\item Let $\mathcal A$ be an abelian category. 
Then the additive category $(\mathcal A,F)$ is an exact category which admits kernel and cokernel 
(but is NOT an abelian category). 
A morphism $\phi:(M,F)\to(N,F)$ with $(M,F)\in(\mathcal A,F)$ is strict if 
the inclusion $\phi(F^nM)\subset F^nN\cap\Im(\phi)$ is an equality, i.e. if $\phi(F^nM)=F^nN\cap\Im(\phi)$.

\item Let $\mathcal A$ be an abelian category. 
\begin{itemize}
\item For $(A^{\bullet},F)\in C(\mathcal A,F)$, considering $a_p:F^pA^{\bullet}\hookrightarrow A^{\bullet}$
the structural monomorphism of of the filtration, we denote by, for $n\in\mathbb N$, 
\begin{equation*}
H^n(A^{\bullet},F)\in(\mathcal A,F), \; 
F^pH^n(A^{\bullet},F):=\Im(H^n(a_p):H^n(F^pA^{\bullet})\to H^n(A^{\bullet}))\subset H^n(A^{\bullet})
\end{equation*}
the filtration induced on the cohomology objects of the complex.
In the case $(A^{\bullet},F)\in C_{fil}(\mathcal A)$, 
the spectral sequence $E^{p,q}_r(A^{\bullet},F)$ associated to $(A^{\bullet},F)$ 
converge to $\Gr_F^pH^{p+q}(A^{\bullet},F)$, that is for all $p,q\in\mathbb Z$, there exist $r_{p+q}\in\mathbb N$,
such that $E^{p,q}_s(A^{\bullet},F)=\Gr_F^pH^{p+q}(A^{\bullet},F)$ for all $s\leq r_{p+q}$.
\item A morphism $m:(A^{\bullet},F)\to(B^{\bullet},F)$ with $(A^{\bullet},F),(B^{\bullet},F)\in C(\mathcal A,F)$
is said to be a filtered quasi-isomorphism if for all $n,p\in\mathbb Z$, 
\begin{equation*}
H^n\Gr^p_F(m):H^n(\Gr^p_FA^{\bullet})\xrightarrow{\sim}H^n(\Gr^p_FB^{\bullet})
\end{equation*}
is an isomorphism in $\mathcal A$. 
Consider a commutative diagram in $C(\mathcal A,F)$
\begin{equation*}
\xymatrix{(A^{\bullet},F)\ar[r]^{m}\ar[d]_{\phi} & (B^{\bullet},F)\ar[r]^{i_2}\ar[d]_{\psi} & 
\Cone(m)=((A^{\bullet},F)[1]\oplus(B^{\bullet},F),d,d'-m)\ar[r]^{p_1}\ar[d]^{(\phi[1],\psi)} & 
(A^{\bullet},F)[1]\ar[d]^{\phi[1]} \\
(A^{'\bullet},F)\ar[r]^{m'} & (B^{'\bullet},F)\ar[r]^{i_2} & 
\Cone(m')=((A^{'\bullet},F)[1]\oplus(B^{'\bullet},F),d,d'-m')\ar[r]^{p_1} & (A^{'\bullet},F)[1]}
\end{equation*}
If $\phi$ and $\psi$ are filtered quasi-isomorphisms, then $(\phi[1],\psi)$ is an filtered quasi-isomorphism.
That is, the filtered quasi-isomorphism satisfy the 2 of 3 property for canonical triangles.
\item If two morphisms $\phi_1,\phi_2:(M,F)\to(N,F)$ with $(M,F),(N,F)\in C(\mathcal A,F)$ are $r$-filtered homotopic, 
then for all $p,q\in\mathbb Z$ and $s\geq r$.
\begin{equation*}
E^{p,q}_s(\phi_1)=E^{p,q}_s(\phi_2):E^{p,q}_s(M,F)\to E^{p,q}_s(M,F).
\end{equation*}
Hence if $\phi:(M,F)\to(N,F)$ with $(M,F),(N,F)\in C(\mathcal A,F)$ is an $r$-filtered homotopy equivalence
then for all $p,q\in\mathbb Z$ and $s\geq r$.
\begin{equation*}
E^{p,q}_r(\phi):E^{p,q}_r(M,F)\xrightarrow{\sim} E^{p,q}_r(N,F)
\end{equation*}
is an isomorphism in $\mathcal A$.
\item Let $r\in\mathbb N$.
A morphism $m:(A^{\bullet},F)\to(B^{\bullet},F)$ with $(A^{\bullet},F),(B^{\bullet},F)\in C(\mathcal A,F)$
is said to be an $r$-filtered quasi-isomorphism if it belongs to the
submonoid of arrows generated by filtered quasi-isomorphism and $r$-filtered homotopy equivalence, that is if
there exists $m_i:(C_i^{\bullet},F)\to(C_{i+1}^{\bullet},F)$, $0\leq i\leq s$, with $(C_i^{\bullet},F)\in C(\mathcal A,F)$
$(C_0^{\bullet},F)=(A^{\bullet},F)$ and $(C_s^{\bullet},F)=(B^{\bullet},F)$, such that
\begin{equation*}
m=m_s\circ\cdots\circ m_i\circ\cdots\circ m_0:(A^{\bullet},F)\to(B^{\bullet},F)
\end{equation*}
and $m_i:(C_i^{\bullet},F)\to(C_{i+1}^{\bullet},F)$ either a filtered quasi-isomorphism or an $r$-filtered homotopy equivalence.
Note that our definition is stronger then the one given in \cite{CG} in order to get a multiplicative system.
Indeed, if $m:(A^{\bullet},F)\to(B^{\bullet},F)$ with $(A^{\bullet},F),(B^{\bullet},F)\in C(\mathcal A,F)$
is an $r$-filtered quasi-isomorphism then for all $p,q\in\mathbb Z$ and $s\geq r$,
\begin{equation*}
E^{p,q}_r(m):E^{p,q}_r(A^{\bullet},F)\xrightarrow{\sim} E^{p,q}_r(B^{\bullet},F)
\end{equation*}
is an isomorphism in $\mathcal A$, but the converse is NOT true. 
If a morphism $m:(A^{\bullet},F)\to(B^{\bullet},F)$, with $(A^{\bullet},F),(B^{\bullet},F)\in C_{fil}(\mathcal A)$
is an $r$-filtered quasi-isomorphism, then for all $n\in\mathbb Z$
\begin{equation*}
H^n(m):H^n(A^{\bullet},F)\xrightarrow{\sim} H^n(B^{\bullet},F)
\end{equation*}
is a filtered isomorphism, i.e. an isomorphism in $(\mathcal A,F)$. 
The converse is true if there exist $N_1,N_2\in\mathbb Z$ such that 
$H^n(A^{\bullet})=H^n(B^{\bullet})=0$ for $n\leq N_1$ or $n\geq N_2$.
A filtered quasi-isomorphism is obviously a $1$-filtered quasi-isomorphism. 
However for $r>1$, the $r$-filtered quasi-isomorphism does NOT satisfy 
the 2 of 3 property for morphisms of canonical triangles.
\item A morphism $m:(A^{\bullet},F)\to(B^{\bullet},F)$ with $(A^{\bullet},F),(B^{\bullet},F)\in C(\mathcal A,F)$
is said to be an $\infty$-filtered quasi-isomorphism if there exist $r\in\mathbb N$ such that
$m:(A^{\bullet},F)\to(B^{\bullet},F)$ an $r$-filtered quasi-isomorphism.
If a morphism $m:(A^{\bullet},F)\to(B^{\bullet},F)$, with $(A^{\bullet},F),(B^{\bullet},F)\in C_{fil}(\mathcal A)$
is an $\infty$-filtered quasi-isomorphism, then for all $n\in\mathbb Z$
\begin{equation*}
H^n(m):H^n(A^{\bullet},F)\xrightarrow{\sim} H^n(B^{\bullet},F)
\end{equation*}
is a filtered isomorphism.
The converse is true if there exist $N_1,N_2\in\mathbb Z$ such that 
$H^n(A^{\bullet})=H^n(B^{\bullet})=0$ for $n\leq N_1$ or $n\geq N_2$.
The $\infty$-filtered quasi-isomorphism does NOT satisfy the 2 of 3 property for morphisms of canonical triangles.
If $m:(A^{\bullet},F)\to(B^{\bullet},F)$ is such that $m:A^{\bullet}\to B^{\bullet}$ is a quasi-isomorphism
but $m:(A^{\bullet},F)\to(B^{\bullet},F)$ is not an $\infty$-filtered quasi-isomorphism, then
it induces an isomorphisms
$H^n(m):H^n(A^{\bullet})\xrightarrow{\sim} H^n(B^{\bullet})$, hence injective maps
\begin{equation*}
H^n(m):F^pH^n(A^{\bullet},F)\hookrightarrow F^pH^n(B^{\bullet},F)
\end{equation*}
which are not isomorphism (the non strict case). 
If $m:(A^{\bullet},F)\to(B^{\bullet},F)$ is such that $m:A^{\bullet}\to B^{\bullet}$ is a quasi-isomorphism
but $m:(A^{\bullet},F)\to(B^{\bullet},F)$ is not an $\infty$-filtered quasi-isomorphism (the non strict case), then
$H^n\Cone(m)=0$ for all $n\in\mathbb N$, hence $\Cone(m)\to 0$ is an $\infty$-filtered quasi-isomorphism ;
this shows that the $\infty$-filtered quasi-isomorphism does NOT satisfy the 2 of 3 property for morphism of canonical triangles.
\end{itemize}

\item Let $\mathcal A$ be an abelian category.
\begin{itemize}
\item We denote by $D(\mathcal A)$ the localization of $K(\mathcal A)$
with respect to the quasi-isomorphisms and by $D:K(\mathcal A)\to D(\mathcal A)$ the localization functor. 
The category $D(\mathcal A)$ is a triangulated category in the unique way such that $D$ a triangulated functor.
\item We denote by $D_{fil}(\mathcal A)$ the localization of $K_{fil}(\mathcal A)$
with respect to the filtered quasi-isomorphisms and by 
$D:K_{fil}(\mathcal A)\to D_{fil}(\mathcal A)$ the localization functor.
\end{itemize}

\item Let $\mathcal A$ be an abelian category.
We denote by $\Inj(A)\subset A$ the full subcategory of injective objects,
and by $\Proj(A)\subset A$ the full subcategory of projective objects.

\item For $\mathcal S\in\Cat$ a small category, we denote by 
\begin{itemize}
\item $\PSh(\mathcal S):=\PSh(\mathcal S,\Ab):=\Fun(\mathcal S,\Ab)$ the category of presheaves on $\mathcal S$,
i.e. the category of presheaves of abelian groups on $\mathcal S$,
\item $\PSh(\mathcal S,\Ring):=\Fun(\mathcal S,\Ring)$ the category of presheaves of ring on $\mathcal S$,
and $\PSh(\mathcal S,\cRing)\subset\PSh(\mathcal S,\Ring)$ the full subcategory of presheaves of commutative ring.
\item for $F\in\PSh(\mathcal S)$ and $X\in\mathcal S$, $F(X)=\Gamma(X,F)$ the sections on $X$ and for $h:X'\to X$ a morphism
with $X,X'\in\mathcal S$, $F(h):=F^{X,Y}(h):F(X)\to F(X')$ the morphism of abelian groups,
\item $C(\mathcal S)=\PSh(\mathcal S,C(\mathbb Z))=C(\PSh(\mathcal S))=\PSh(\mathcal S\times\mathbb Z)$
the big abelian category of complexes of presheaves on $\mathcal S$ with value in abelian groups,
\item $K(\mathcal S):=K(\PSh(\mathcal S))=\Ho(C(\mathcal S))$
In particular, we have the full homotopy functor $Ho:C(\mathcal S)\to K(\mathcal S)$,
\item $C_{(2)fil}(\mathcal S):=C_{(2)fil}(\PSh(\mathcal S))\subset C(\PSh(\mathcal S),F,W)$
the big abelian category of (bi)filtered complexes of presheaves on $\mathcal S$ with value in abelian groups
such that the filtration is biregular, and $\PSh_{(2)fil}(\mathcal S)=(\PSh(\mathcal S),F,W)$,
\item $K_{fil}(\mathcal S):=K_{fil}(\PSh(\mathcal S))=\Ho(C_{fil}(\mathcal S))$,
\item $K_{fil,r}(\mathcal S):=K_{fil,r}(\PSh(\mathcal S))=\Ho_r(C_{fil}(\mathcal S))$,
$K_{fil,\infty}(\mathcal S):=K_{fil,\infty}(\PSh(\mathcal S))=\Ho_{\infty}(C_{fil}(\mathcal S))$.
\end{itemize}
For $f:\mathcal T\to\mathcal S$ a morphism a presite with $\mathcal T,\mathcal S\in\Cat$,
given by the functor $P(f):\mathcal S\to\mathcal T$,
we will consider the adjonctions given by the direct and inverse image functors : 
\begin{itemize}
\item $(f^*,f_*)=(f^{-1},f_*):\PSh(\mathcal S)\leftrightarrows\PSh(\mathcal T)$,
which induces $(f^*,f_*):C(\mathcal S)\leftrightarrows C(\mathcal T)$, we denote, 
for $F\in C(\mathcal S)$ and $G\in C(\mathcal T)$ by
\begin{equation*}
\ad(f^*,f_*)(F):F\to f_*f^*F \; , \; \ad(f^*,f_*)(G):f^*f_*G\to G
\end{equation*}
the adjonction maps,
\item $(f_*,f^{\bot}):\PSh(\mathcal T)\leftrightarrows\PSh(\mathcal S)$,
which induces $(f_*,f^{\bot}):C(\mathcal T)\leftrightarrows C(\mathcal S)$, we denote
for $F\in C(\mathcal S)$ and $G\in C(\mathcal T)$ by
\begin{equation*}
\ad(f_*,f^{\bot})(F):G\to f^{\bot}f_* G \; , \; \ad(f_*,f^{\bot})(G):f_*f^{\bot}F\to F
\end{equation*}
the adjonction maps. 
\end{itemize}

\item For $(\mathcal S,O_S)\in\RCat$ a ringed topos, we denote by 
\begin{itemize}
\item $\PSh_{O_S}(\mathcal S)$ the category of presheaves of $O_S$ modules on $\mathcal S$, 
whose objects are $\PSh_{O_S}(\mathcal S)^0:=\left\{(M,m),M\in\PSh(\mathcal S),m:M\otimes O_S\to M\right\}$,
together with the forgetful functor $o:\PSh(\mathcal S)\to \PSh_{O_S}(\mathcal S)$,
\item $C_{O_S}(\mathcal S)=C(\PSh_{O_S}(\mathcal S))$ 
the big abelian category of complexes of presheaves of $O_S$ modules on $\mathcal S$,
\item $K_{O_S}(\mathcal S):=K(\PSh_{O_S}(\mathcal S))=\Ho(C_{O_S}(\mathcal S))$,
in particular, we have the full homotopy functor $Ho:C_{O_S}(\mathcal S)\to K_{O_S}(\mathcal S)$,
\item $C_{O_S(2)fil}(\mathcal S):=C_{(2)fil}(\PSh_{O_S}(\mathcal S))\subset C(\PSh_{O_S}(\mathcal S),F,W)$,
the big abelian category of (bi)filtered complexes of presheaves of $O_S$ modules on $\mathcal S$ such that the filtration is biregular
and $\PSh_{O_S(2)fil}(\mathcal S)=(\PSh_{O_S}(\mathcal S),F,W)$,
\item $K_{O_Sfil}(\mathcal S):=K_{fil}(\PSh_{O_S}(\mathcal S))=\Ho(C_{O_Sfil}(\mathcal S))$,
\item $K_{O_Sfil,r}(\mathcal S):=K_{fil,r}(\PSh_{O_S}(\mathcal S))=\Ho_r(C_{O_Sfil}(\mathcal S))$,
$K_{O_Sfil,\infty}(\mathcal S):=K_{fil,\infty}(\PSh_{O_S}(\mathcal S))=\Ho_{\infty}(C_{O_Sfil}(\mathcal S))$.
\end{itemize}

\item For $\mathcal S_{\bullet}\in\Fun(\mathcal I,\Cat)$ a diagram of (pre)sites, with $\mathcal I\in\Cat$ a small category, 
we denote by 
\begin{itemize}
\item $\Gamma\mathcal S_{\bullet}\in\Cat$ the associated diagram category
\begin{itemize}
\item whose objects are $\Gamma\mathcal S_{\bullet}^0:=\left\{(X_I,u_{IJ})_{I\in\mathcal I}\right\}$, 
with $X_I\in\mathcal S_I$, and for $r_{IJ}:I\to J$ with $I,J\in\mathcal I$,  
$u_{IJ}:X_J\to r_{IJ}(X_I)$ are morphism in $\mathcal S_J$
noting again $r_{IJ}:\mathcal S_I\to\mathcal S_J$ the associated functor,
\item whose morphism are $m=(m_I):(X_I,u_{IJ})\to(X'_I,v_{IJ})$ satisfying 
$v_{IJ}\circ m_I =r_{IJ}(m_J)\circ u_{IJ}$ in $\mathcal S_J$, 
\end{itemize}
\item $\PSh(\mathcal S_{\bullet}):=\PSh(\Gamma\mathcal S_{\bullet},\Ab)$ the category of presheaves on $\mathcal S_{\bullet}$,
\begin{itemize}
\item whose objects are $\PSh(\mathcal S_{\bullet})^0:=\left\{(F_I,u_{IJ})_{I\in\mathcal I}\right\}$, 
with $F_I\in\PSh(\mathcal S_I)$, and for $r_{IJ}:I\to J$ with $I,J\in\mathcal I$, 
$u_{IJ}:F_I\to r_{IJ*}F_J$ are morphism in $\PSh(\mathcal S_I)$, 
noting again $r_{IJ}:\mathcal S_J\to\mathcal S_I$ the associated morphism of presite,
\item whose morphism are $m=(m_I):(F_I,u_{IJ})\to(G_I,v_{IJ})$ satisfying 
$v_{IJ}\circ m_I =r_{IJ*}m_J\circ u_{IJ}$ in $\PSh(\mathcal S_I)$, 
\end{itemize}
\item $\PSh(\mathcal S_{\bullet},\Ring):=\PSh(\Gamma\mathcal S_{\bullet},\Ring)$ 
the category of presheaves of ring on $\mathcal S_{\bullet}$ given in the same way,
and $\PSh(\mathcal S_{\bullet},\cRing)\subset\PSh(\mathcal S_{\bullet},\Ring)$ 
the full subcategory of presheaves of commutative ring.
\item $C(\mathcal S_{\bullet}):=C(\PSh(\mathcal S_{\bullet}))$
the big abelian category of complexes of presheaves on $\mathcal S_{\bullet}$ with value in abelian groups,
\item $K(\mathcal S_{\bullet}):=K(\PSh(\mathcal S_{\bullet}))=\Ho(C(\mathcal S_{\bullet}))$,
in particular, we have the full homotopy functor $Ho:C(\mathcal S_{\bullet})\to K(\mathcal S_{\bullet})$,
\item $C_{(2)fil}(\mathcal S_{\bullet}):=C_{(2)fil}(\PSh(\mathcal S_{\bullet}))\subset C(\PSh(\mathcal S_{\bullet}),F,W)$
the big abelian category of (bi)filtered complexes of presheaves on $\mathcal S_{\bullet}$ with value in abelian groups
such that the filtration is biregular, and $\PSh_{(2)fil}(\mathcal S_{\bullet})=(\PSh(\mathcal S_{\bullet}),F,W)$,
by definition $C_{(2)fil}(\mathcal S_{\bullet})$ is the category
\begin{itemize}
\item whose objects are $C_{(2)fil}(\mathcal S_{\bullet})^0:=\left\{((F_I,F,W),u_{IJ})_{I\in\mathcal I}\right\}$, 
with $(F_I,F,W)\in C_{(2)fil}(\mathcal S_I)$, and for $r_{IJ}:I\to J$ with $I,J\in\mathcal I$, 
$u_{IJ}:(F_I,F,W)\to r_{IJ*}(F_J,F,W)$ are morphism in $C_{(2)fil}(\mathcal S_I)$, 
noting again $r_{IJ}:\mathcal S_J\to\mathcal S_I$ the associated morphism of presite,
\item whose morphism are $m=(m_I):((F_I,F,W),u_{IJ})\to((G_I,F,W),v_{IJ})$ satisfying 
$v_{IJ}\circ m_I =r_{IJ*}m_J\circ u_{IJ}$ in $C_{(2)fil}(\mathcal S_I)$, 
\end{itemize}
\item $K_{fil}(\mathcal S_{\bullet}):=K_{fil}(\PSh(\mathcal S_{\bullet}))=\Ho(C_{fil}(\mathcal S_{\bullet}))$,
\item $K_{fil,r}(\mathcal S_{\bullet}):=K_{fil,r}(\PSh(\mathcal S_{\bullet}))=\Ho_r(C_{fil}(\mathcal S_{\bullet}))$,
$K_{fil,\infty}(\mathcal S_{\bullet}):=K_{fil,\infty}(\PSh(\mathcal S_{\bullet}))=\Ho_{\infty}(C_{fil}(\mathcal S_{\bullet}))$.
\end{itemize}
Let $\mathcal I,\mathcal I'\in\Cat$ be small categories.
Let $(f_{\bullet},s):\mathcal T_{\bullet}\to\mathcal S_{\bullet}$ a morphism a diagrams of (pre)site with 
$\mathcal T_{\bullet}\in\Fun(\mathcal I,\Cat),\mathcal S_{\bullet}\in\Fun(\mathcal I',\Cat)$,
which is by definition given by a functor $s:\mathcal I\to \mathcal I'$ and morphism of functor 
$P(f_{\bullet}):\mathcal S_{s(\bullet)}:=\mathcal S_{\bullet}\circ s\to\mathcal T_{\bullet}$.
Here, we denote for short, $\mathcal S_{s(\bullet)}:=\mathcal S_{\bullet}\circ s\in\Fun(\mathcal I,\Cat)$.
We have then, for $r_{IJ}:I\to J$ a morphism, with $I,J\in\mathcal I$, a commutative diagram in $\Cat$ 
\begin{equation*}
D_{fIJ}:=\xymatrix{\mathcal S_{s(J)}\ar[r]^{r^s_{IJ}} & \mathcal S_{s(I)} \\
\mathcal T_J\ar[r]^{r^t_{IJ}}\ar[u]^{f_J} & \mathcal T_I\ar[u]^{f_I}}.
\end{equation*}
We will consider the the adjonction given by the direct and inverse image functors :
\begin{eqnarray*}
((f_{\bullet},s)^*,(f_{\bullet},s)^*)=((f_{\bullet},s)^{-1},(f_{\bullet},s)_*):
\PSh(\mathcal S_{s(\bullet)})\leftrightarrows\PSh(\mathcal T_{\bullet}), \\
F=(F_I,u_{IJ})\mapsto (f_{\bullet},s)^*F:=(f_I^*F_I,T(D_{fIJ})(F_J)\circ f_I^*u_{IJ}), \\ 
G=(G_I,v_{IJ})\mapsto (f_{\bullet},s)_*G:=(f_{I*}G_I,f_{I*}v_{IJ}).
\end{eqnarray*} 
It induces the adjonction
$((f_{\bullet},s)^*,(f_{\bullet},s)_*):C(\mathcal S_{s(\bullet)})\leftrightarrows C(\mathcal T_{\bullet})$.
We denote, for $(F_I,u_{IJ})\in C(\mathcal S_{s(\bullet)})$ and $(G_I,v_{IJ})\in C(\mathcal T_{\bullet})$ by
\begin{eqnarray*}
\ad((f_{\bullet},s)^*,(f_{\bullet},s)_*)((F_I,u_{IJ})):(F_I,u_{IJ})\to (f_{\bullet},s)_*(f_{\bullet},s)^*(F_I,u_{IJ}), \\ 
\ad((f_{\bullet},s)^*,(f_{\bullet},s)_*)((G_I,v_{IJ})):(f_{\bullet},s)^*(f_{\bullet},s)_*(G_I,v_{IJ})\to (G_I,v_{IJ})
\end{eqnarray*}
the adjonction maps.

\item Let $\mathcal I\in\Cat$ a small category.
For $(\mathcal S_{\bullet},O_{S_{\bullet}})\in\Fun(\mathcal I,\RCat)$ a diagram of ringed topos, we denote by 
\begin{itemize}
\item $\PSh_{O_{S_{\bullet}}}(\mathcal S_{\bullet}):=\PSh_{O_{\Gamma\mathcal S_{\bullet}}}(\Gamma\mathcal S_{\bullet})$ 
the category of presheaves of modules on $(\mathcal S_{\bullet},O_{S_{\bullet}})$,
\begin{itemize}
\item whose objects are $\PSh_{O_{S_{\bullet}}}(\mathcal S_{\bullet})^0:=\left\{(F_I,u_{IJ})_{I\in\mathcal I}\right\}$, 
with $F_I\in\PSh_{O_{S_I}}(\mathcal S_I)$, and for $r_{IJ}:I\to J$ with $I,J\in\mathcal I$, 
$u_{IJ}:F_I\to r_{IJ*}F_J$ are morphism in $\PSh_{O_{S_I}}(\mathcal S_I)$, 
noting again $r_{IJ}:\mathcal S_J\to\mathcal S_I$ the associated morphism of presite,
\item whose morphism are $m=(m_I):(F_I,u_{IJ})\to(G_I,v_{IJ})$ satisfying 
$v_{IJ}\circ m_I =r_{IJ*}m_J\circ u_{IJ}$ in $\PSh_{O_{S_I}}(\mathcal S_I)$, 
\end{itemize}
\item $C_{O_{S_{\bullet}}}(\mathcal S_{\bullet}):=C(\PSh_{O_{S_{\bullet}}}(\mathcal S_{\bullet}))$,
\item $K_{O_{S_{\bullet}}}(\mathcal S_{\bullet}):=K(\PSh_{O_{S_{\bullet}}}(\mathcal S_{\bullet}))
=\Ho(C_{O_{S_{\bullet}}}(\mathcal S_{\bullet}))$,
in particular, we have the full homotopy functor $Ho:C(\mathcal S_{\bullet})\to K(\mathcal S_{\bullet})$,
\item $C_{O_{S_{\bullet}}(2)fil}(\mathcal S_{\bullet}):=C_{O_{S_{\bullet}}(2)fil}(\PSh(\mathcal S_{\bullet}))
\subset C(\PSh_{O_{S_{\bullet}}}(\mathcal S_{\bullet}),F,W)$
the big abelian category of (bi)filtered complexes of presheaves of modules on $(\mathcal S_{\bullet},O_{S_{\bullet}})$
such that the filtration is biregular, and 
$\PSh_{O_{S_{\bullet}}(2)fil}(\mathcal S_{\bullet})=(\PSh_{O_{S_{\bullet}}}(\mathcal S_{\bullet}),F,W)$,
by definition $C_{O_{S_{\bullet}}(2)fil}(\mathcal S_{\bullet})$ is the category
\begin{itemize}
\item whose objects are 
$C_{O_{S_{\bullet}}(2)fil}(\mathcal S_{\bullet})^0:=\left\{((F_I,F,W),u_{IJ})_{I\in\mathcal I}\right\}$, 
with $(F_I,F,W)\in C_{O_{S_I}(2)fil}(\mathcal S_I)$, and for $r_{IJ}:I\to J$ with $I,J\in\mathcal I$, 
$u_{IJ}:(F_I,F,W)\to r_{IJ*}(F_J,F,W)$ are morphism in $C_{O_{S_I}(2)fil}(\mathcal S_I)$, 
noting again $r_{IJ}:\mathcal S_J\to\mathcal S_I$ the associated morphism of presite,
\item whose morphism are $m=(m_I):((F_I,F,W),u_{IJ})\to((G_I,F,W),v_{IJ})$ satisfying 
$v_{IJ}\circ m_I =r_{IJ*}m_J\circ u_{IJ}$ in $C_{O_{S_I}(2)fil}(\mathcal S_I)$, 
\end{itemize}
\item $K_{O_{S_{\bullet}}(2)fil}(\mathcal S_{\bullet}):=K_{fil}(\PSh_{O_{S_{\bullet}}(2)}(\mathcal S_{\bullet}))
=\Ho(C_{O_{S_{\bullet}}(2)fil}(\mathcal S_{\bullet}))$,
\item $K_{O_{S_{\bullet}}fil,r}(\mathcal S_{\bullet}):=K_{fil,r}(\PSh_{O_{S_{\bullet}}}(\mathcal S_{\bullet}))
=\Ho_r(C_{O_{S_{\bullet}}fil,r}(\mathcal S_{\bullet}))$,
$K_{O_{S_{\bullet}}fil,\infty}(\mathcal S_{\bullet}):=K_{\infty,r}(\PSh_{O_{S_{\bullet}}}(\mathcal S_{\bullet}))
=\Ho_{\infty}(C_{O_{S_{\bullet}}fil,\infty}(\mathcal S_{\bullet}))$.
\end{itemize}

\item Let $\mathcal S\in\Cat$. For $\Sigma:C(\mathcal S)\to C(\mathcal S)$ an endofunctor, we denote by
$C_\Sigma(\mathcal S)=(C(\mathcal S),\Sigma)$ the corresponding category of spectra.

\item Denote by $\Sch\subset\RTop$ the full subcategory of schemes. 
For a field $k$, we consider $\Sch/k:=\Sch/\Spec k$ the category of schemes over $\Spec k$. We then denote by
\begin{itemize}
\item $\Var(k)\subset\Sch/k$ the full subcategory consisting of algebraic varieties over $k$, 
i.e. schemes of finite type over $k$,
\item $\PVar(k)\subset\QPVar(k)\subset\Var(k)$ the full subcategories consisting of 
quasi-projective varieties and projective varieties respectively, 
\item $\PSmVar(k)\subset\SmVar(k)\subset\Var(k)$ the full subcategories consisting of 
smooth varieties and smooth projective varieties respectively.
\end{itemize}
A morphism $h:U\to S$ with $U,S\in\Var(k)$ is said to be smooth if it is flat with smooth fibers.
A morphism $r:U\to X$ with $U,X\in\Var(k)$ is said to be etale if it is non ramified and flat.
In particular an etale morphism $r:U\to X$ with $U,X\in\Var(k)$ 
is smooth and quasi-finite (i.e. the fibers are either the empty set or a finite subset of $X$)

\item Denote by $\Top^2$ the category whose set of objects is 
\begin{equation*}
(\Top^2)^0:=\left\{(X,Z), \; Z\subset X \; \mbox{closed}\right\}\subset\Top\times\Top
\end{equation*}
and whose set of morphism between $(X_1,Z_1),(X_2,Z_2)\in\Top^2$ is
\begin{eqnarray*}
\Hom_{\Top^2}((X_1,Z_1),(X_2,Z_2)):= 
\left\{(f:X_1\to X_2), \; \mbox{s.t.} \; Z_1\subset f^{-1}(Z_2)\right\}\subset\Hom_{\Top}(X_1,X_2)
\end{eqnarray*}
For $S\in\Top$, $\Top^2/S:=\Top^2/(S,S)$ is then by definition the category whose set of objects is 
\begin{eqnarray*}
(\Top^2/S)^0:=\left\{((X,Z),h), h:X\to S, \; Z\subset X \; \mbox{closed} \;\right\}\subset\Top/S\times\Top
\end{eqnarray*}
and whose set of morphisms between $(X_1,Z_1)/S=((X_1,Z_1),h_1),(X_2,Z_2)/S=((X_2,Z_2),h_2)\in\Top^2/S$
is the subset
\begin{eqnarray*}
\Hom_{\Top^2/S}((X_1,Z_1)/S,(X_2,Z_2)/S):= \\
\left\{(f:X_1\to X_2), \; \mbox{s.t.} \; h_1\circ f=h_2 \; \mbox{and} \; Z_1\subset f^{-1}(Z_2)\right\}
\subset\Hom_{\RTop}(X_1,X_2)
\end{eqnarray*}
We denote by
\begin{eqnarray*}
\mu_S:\Top^{2,pr}/S:=\left\{((Y\times S,Z),p), p:Y\times S\to S, \; Z\subset Y\times S \; \mbox{closed} \;\right\}
\hookrightarrow\Top^2/S
\end{eqnarray*}
the full subcategory whose objects are those with $p:Y\times S\to S$ a projection,
and again $\mu_S:\Top^2/S\to\Top^{2,pr}/S$ the corresponding morphism of sites. 
We denote by 
\begin{eqnarray*}
\Gr_S^{12}:\Top/S\to\Top^{2,pr}/S, \; X/S\mapsto\Gr_S^{12}(X/S):=(X\times S,\bar X)/S, \\
(g:X/S\to X'/S)\mapsto\Gr_S^{12}(g):=(g\times I_S:(X\times S,\bar X)\to(X'\times S,\bar X'))
\end{eqnarray*}
the graph functor, $X\hookrightarrow X\times S$ being the graph embedding (which is a closed embedding if $X$ is separated),
and again $\Gr_S^{12}:\Top^{2,pr}/S\to\Top/S$ the corresponding morphism of sites.

\item Denote by $\RTop^2$ the category whose set of objects is 
\begin{equation*}
(\RTop^2)^0:=\left\{((X,O_X),Z), \; Z\subset X \; \mbox{closed}\right\}\subset\RTop\times\Top
\end{equation*}
and whose set of morphism between $((X_1,O_{X_1}),Z_1),((X_2,O_{X_2}),Z_2)\in\RTop^2$ is
\begin{eqnarray*}
\Hom_{\RTop^2}(((X_1,O_{X_1}),Z_1),((X_2,O_{X_2}),Z_2)):= \\
\left\{(f:(X_1,O_{X_1})\to (X_2,O_{X_2})), \; \mbox{s.t.} \; Z_1\subset f^{-1}(Z_2)\right\}
\subset\Hom_{\RTop}((X_1,O_{X_1}),(X_2,O_{X_2}))
\end{eqnarray*}
For $(S,O_S)\in\RTop$, $\RTop^2/(S,O_S):=\RTop^2/((S,O_S),S)$ is then by definition the category whose set of objects is 
\begin{eqnarray*}
(\RTop^2/(S,O_S))^0:= \\
\left\{(((X,O_X),Z),h), h:(X,O_X)\to(S,O_S), \; Z\subset X \; \mbox{closed} \;\right\}\subset\RTop/(S,O_S)\times\Top
\end{eqnarray*}
and whose set of morphisms between $(((X_1,O_{X_1}),Z_1),h_1),(((X_2,O_{X_2}),Z_2),h_2)\in\RTop^2/(S,O_S)$
is the subset
\begin{eqnarray*}
\Hom_{\RTop^2/(S,O_S)}(((X_1,O_{X_1}),Z_1)/(S,O_S),((X_2,O_{X_2}),Z_2)/(S,O_S)):= \\
\left\{(f:(X_1,O_{X_1})\to (X_2,O_{X_2})), \; \mbox{s.t.} \; h_1\circ f=h_2 \; \mbox{and} \; Z_1\subset f^{-1}(Z_2)\right\} \\ 
\subset\Hom_{\RTop}((X_1,O_{X_1}),(X_2,O_{X_2}))
\end{eqnarray*}
We denote by 
\begin{eqnarray*}
\mu_S:\RTop^{2,pr}/S:=\left\{(((Y\times S,q^*O_Y\otimes p^*O_S),Z),p), p:Y\times S\to S, \; 
Z\subset Y\times S \; \mbox{closed} \;\right\}\hookrightarrow\RTop^2/S
\end{eqnarray*}
the full subcategory whose objects are those with $p:Y\times S\to S$ is a projection,
and again $\mu_S:\RTop^2/S\to\RTop^{2,pr}/S$ the corresponding morphism of sites. 
We denote by 
\begin{eqnarray*}
\Gr_S^{12}:\RTop/S\to\RTop^{2,pr}/S, \\ 
(X,O_X)/(S,O_S)\mapsto\Gr_S^{12}((X,O_X)/(S,O_S)):=((X\times S,q^*O_X\otimes p^*O_S),\bar X)/(S,O_S), \\
(g:(X,O_X)/(S,O_S)\to (X',O_{X'})/(S,O_S))\mapsto \\
\Gr_S^{12}(g):=(g\times I_S:((X\times S,q^*O_X\otimes p^*O_S),\bar X)\to((X'\times S,q^*O_X\otimes p^*O_S),\bar X'))
\end{eqnarray*}
the graph functor, $X\hookrightarrow X\times S$ being the graph embedding (which is a closed embedding if $X$ is separated),
$p:X\times S\to S$, $q:X\times S\to X$ the projections,
and again $\Gr_S^{12}:\RTop^{2,pr}/S\to\RTop/S$ the corresponding morphism of sites.

\item We denote by $\Sch^2\subset\RTop^2$ the full subcategory such that the first factors are schemes. 
For a field $k$, we denote by $\Sch^2/k:=\Sch^2/(\Spec k,\left\{\pt\right\})$ and by
\begin{itemize}
\item $\Var(k)^2\subset\Sch^2/k$ the full subcategory such that the first factors are algebraic varieties over $k$, 
i.e. schemes of finite type over $k$,
\item $\PVar(k)^2\subset\QPVar(k)^2\subset\Var(k)^2$ the full subcategories such that the first factors are
quasi-projective varieties and projective varieties respectively, 
\item $\PSmVar(k)^2\subset\SmVar(k)^2\subset\Var(k)^2$ the full subcategories such that the first factors are
smooth varieties and smooth projective varieties respectively.
\end{itemize}
In particular we have, for $S\in\Var(k)$, the graph functor 
\begin{eqnarray*}
\Gr_S^{12}:\Var(k)/S\to\Var(k)^{2,pr}/S, \; X/S\mapsto\Gr_S^{12}(X/S):=(X\times S,X)/S, \\
(g:X/S\to X'/S)\mapsto\Gr_S^{12}(g):=(g\times I_S:(X\times S,X)\to(X'\times S,X'))
\end{eqnarray*}
the graph embedding $X\hookrightarrow X\times S$ is a closed embedding 
since $X$ is separated in the subcategory of schemes $\Sch\subset\RTop$,
and again $\Gr_S^{12}:\Var(k)^{2,pr}/S\to\Var(k)/S$ the corresponding morphism of sites.

\item Denote by $\CW\subset\Top$ the full subcategory of $CW$ complexes, by $\CS\subset\CW$ the full subcategory of $\Delta$ complexes,
by $\TM(\mathbb R)\subset\CW$ the full subcategory of topological (real) manifolds
which admits a CW structure (a topological manifold admits a CW structure if it admits a differential structure)
and by $\Diff(\mathbb R)\subset\RTop$ the full subcategory of differentiable (real) manifold. 
We denote by $\CW^2\subset\Top^2$ the full subcategory such that the first factors are $CW$ complexes,
by $\TM(\mathbb R)^2\subset\CW^2$ the full subcategory such that the first factors are topological (real) manifolds
and by $\Diff(\mathbb R)^2\subset\RTop^2$ the full subcategory such that the first factors are differentiable (real) manifold. 

\item Denote by $\AnSp(\mathbb C)\subset\RTop$ the full subcategory of analytic spaces over $\mathbb C$,
and by $\AnSm(\mathbb C)\subset\AnSp(\mathbb C)$ the full subcategory of smooth analytic spaces (i.e. complex analytic manifold).
A morphism $h:U\to S$ with $U,S\in\AnSp(\mathbb C)$ is said to be smooth if it is flat with smooth fibers.
A morphism $r:U\to X$ with $U,X\in\AnSp(\mathbb C)$ is said to be etale if it is non ramified and flat.
By the Weirstrass preparation theorem (or the implicit function theorem if $U$ and $X$ are smooth),
a morphism $r:U\to X$ with $U,X\in\AnSp(\mathbb C)$ is etale if and only if it is an isomorphism local.

We denote by $\AnSp(\mathbb C)^2\subset\RTop^2$ the full subcategory such that the first factors are analytic spaces over $\mathbb C$,
and by $\AnSm(\mathbb C)^2\subset\AnSp(\mathbb C)^2$ the full subcategory such that the first factors are
smooth analytic spaces (i.e. complex analytic manifold).
In particular we have, for $S\in\AnSp(\mathbb C)$, the graph functor 
\begin{eqnarray*}
\Gr_S^{12}:\AnSp(\mathbb C)/S\to\AnSp(\mathbb C)^{2,pr}/S, \; X/S\mapsto\Gr_S^{12}(X/S):=(X\times S,X)/S, \\
(g:X/S\to X'/S)\mapsto\Gr_S^{12}(g):=(g\times I_S:(X\times S,X)\to(X'\times S,X'))
\end{eqnarray*}
the graph embedding $X\hookrightarrow X\times S$ is a closed embedding since $X$ is separated in $\RTop$,
and again $\Gr_S^{12}:\AnSp(\mathbb C)^{2,pr}/S\to\AnSp(\mathbb C)/S$ the corresponding morphism of sites.

\item For $V\in\Var(\mathbb C)$, we denote by $V^{an}\in\AnSp(\mathbb C)$   
the complex analytic space associated to $V$ with the usual topology induced by the usual topology of $\mathbb C^N$. 
For $W\in\AnSp(\mathbb C)$, we denote by $W^{cw}\in\AnSp(\mathbb C)$ the topological space given by $W$ which is a $CW$ complex.
For simplicity, for $V\in\Var(\mathbb C)$, we denote by $V^{cw}:=(V^{an})^{cw}\in\CW$. We have then  
\begin{itemize}
\item the analytical functor $\An:\Var(\mathbb C)\to\AnSp(\mathbb C)$, $\An(V)=V^{an}$,
\item the forgetful functor $\Cw=tp:\AnSp(\mathbb C)\to\CW$, $\Cw(W)=W^{cw}$,
\item the composite of these two functors $\widetilde\Cw=\Cw\circ\An:\Var(\mathbb C)\to\CW$, $\widetilde\Cw(V)=V^{cw}$. 
\end{itemize}
We have then  
\begin{itemize}
\item the analytical functor $\An:\Var(\mathbb C)^2\to\AnSp(\mathbb C)^2$, $\An((V,Z))=(V^{an},Z^{an})$,
\item the forgetful functor $\Cw=tp:\AnSp(\mathbb C)^2\to\CW^2$, $\Cw((W,Z))=(W^{cw},Z^{cw})$,
\item the composite of these two functors $\widetilde\Cw=\Cw\circ\An:\Var(\mathbb C)^2\to\CW^2$, 
$\widetilde\Cw((V,Z))=(V^{cw},Z^{cw})$. 
\end{itemize}
\end{itemize}

\subsection{Additive categories, abelian categories and tensor triangulated categories}

Let $\mathcal A$ an additive category. 
\begin{itemize}
\item For $\phi:F^{\bullet}\to G^{\bullet}$ a morphism with $F^{\bullet},G^{\bullet}\in C(\mathcal A)$,
we have the mapping cylinder 
$\Cyl(\phi):=((F^n\oplus F^{n+1}\oplus G^{n+1},(\partial^n_F,\partial^{n+1}_F,\phi^{n+1}+\partial^n G)\in C(\mathcal A)$.
and the mapping cone $\Cone(\phi):=((F^n\oplus G^{n+1},(\partial^n_F,\phi^{n+1}+\partial^n G)\in C(\mathcal A)$.
\item The category $K(\mathcal A):=\Ho(C(\mathcal A))$ is a triangulated category
with distinguish triangles $F^{\bullet}\xrightarrow{i_F}\Cyl(\phi)\xrightarrow{q_F}\Cone(\phi)\xrightarrow{r_F} F^{\bullet}[1]$.
\item The category $(\mathcal A,F)$ is obviously again an additive category.
\item Let $\phi:F^{\bullet}\to G^{\bullet}$ a morphism with $F^{\bullet},G^{\bullet}\in C(\mathcal A)$.
Then it is obviously a morphism of filtered complex $\phi:(F^{\bullet},F_b)\to (G^{\bullet},F_b)$, 
where we recall that $F_b$ is the trivial filtration $(F^{\bullet},F_b),(G^{\bullet},F_b)\in C_{fil}(\mathcal A)$.
\end{itemize}

We recall the following property of the internal hom functor if it exists of a tensor triangulated category 
and the definition of compact and cocompact object.

\begin{prop}\label{inthom}
Let $(\mathcal T,\otimes)$ a tensor triangulated category admitting countable direct sum and product
compatible with the triangulation.
Assume that $\mathcal T$ has an internal hom (bi)functor $R\mathcal Hom(.,.):\mathcal T^2\to\mathcal T$
which is by definition the right adjoint to $(\cdot\otimes\cdot):\mathcal T^2\to\mathcal T$.
Then,
\begin{itemize}
\item for $N\in\mathcal T$, the functor 
$R\mathcal Hom(\cdot,N):\mathcal T\to\mathcal T$
commutes with homotopy colimits : 
for $M=\ho\lim_{\to i\in I} M_i$, where $I$ is a countable category, we have
\begin{equation*}
R\mathcal Hom(M,N)\xrightarrow{\sim}\ho\lim_{\leftarrow i\in I}R\mathcal Hom(M_i,N).
\end{equation*}
\item dually, for $M\in\mathcal T$, 
the functor $R\mathcal Hom(M,\cdot):\mathcal T\to\mathcal T$
commutes with homotopy limits : 
for $N={\ho\lim_{\leftarrow}}_{i\in I} N_i$, where $I$ is a countable category, we have
\begin{equation*}
R\mathcal Hom(M,N)\xrightarrow{\sim}\ho\lim_{\leftarrow}R\mathcal Hom(M,N_i).
\end{equation*}
\end{itemize}
\end{prop}

\begin{proof}
Standard.
\end{proof}

Let $(\mathcal T,\otimes)$ a tensor triangulated category admitting countable direct sum and product
compatible with the triangulation. 
Assume that $\mathcal T$ has an internal hom functor $R\mathcal Hom(.,.):\mathcal T\to\mathcal T$.
\begin{itemize}
\item For $N\in\mathcal T$, the functor 
$R\mathcal Hom(\cdot,N):\mathcal T\to\mathcal T$
does not commutes in general with homotopy limits : 
for $M={\ho\lim_{\leftarrow}}_{i\in I} M_i$, where $I$ is a countable category, the canonical map
\begin{equation*}
{\ho\lim_{\to}}_{i\in I}R\mathcal Hom(M_i,N)\to R\mathcal Hom(M,N) 
\end{equation*}
is not an isomorphism in general if $I$ is infinite. It commutes if and only if $N$ is compact. 
\item Dually, for $M\in\mathcal T$, the functor 
$R\mathcal Hom(M,\bullet):\mathcal T\to\mathcal T$
does not commutes in general with infinite homotopy colimits. It commutes if and only if $M$ is cocompact. 
\end{itemize}

Most triangulated category comes from the localization of the category of complexes of an abelian category 
with respect to quasi-isomorphisms.
In the case where the abelian category have enough injective or projective object, 
the triangulated category is the homotopy category of the complexes of injective, resp. projective, objects.

\begin{prop}\label{projinj1}
Let $\mathcal A$ an abelian category with enough injective and projective.
\begin{itemize}
\item A quasi-isomorphism $\phi:Q^{\bullet}\to F^{\bullet}$, 
with $F^{\bullet},Q^{\bullet}\in C^-(\mathcal A)$ such that the $Q^n$ are projective
is an homotopy equivalence.
\item Dually,a quasi-isomorphism $\phi:F^{\bullet}\to I^{\bullet}$, 
with $F^{\bullet},I^{\bullet}\in C^+(\mathcal A)$ such that the $I^n$ are projective
is an homotopy equivalence.
\end{itemize}
\end{prop}

\begin{proof}
Standard.
\end{proof}

\begin{prop}\label{projinj2}
Let $\mathcal A$ an abelian category with enough injective and projective satisfying AB3 
(i.e. countable direct sum  of exact sequences are exact sequence). 
\begin{itemize}
\item Let $K(P)\subset K(\mathcal A)$ be the thick subcategory generated by (unbounded) complexes of projective objects. Then,
$K(P)\hookrightarrow K(\mathcal A)\xrightarrow{D} D(\mathcal A)$
is an equivalence of triangulated categories.
\item Similarly, let $K(I)\subset K(\mathcal A)$ be the thick subcategory generated by (unbounded) complexes of injective objects. Then
$K(I)\hookrightarrow K(\mathcal A)\xrightarrow{D} D(\mathcal A)$
is an equivalence of triangulated categories.
\end{itemize}
\end{prop}

\begin{proof}
It follows from proposition \ref{projinj1} : see \cite{Neeman}.
\end{proof}

Let $\mathcal A\subset\Cat$ an abelian category.
Let $\phi:(M,F)\to(N,F)$ a morphism with $(M,F),(N,F)\in C_{fil}(\mathcal A)$. 
Then the distinguish triangle
\begin{equation*}
(M,F)\xrightarrow{\phi}(N,F)\xrightarrow{i_1}\Cone(\phi)=((M,F)[1]\oplus(N,F),(d,d'-\phi)\xrightarrow{p_1}(M,F)[1]
\end{equation*}
gives a sequence
\begin{equation*}
\cdots\to H^n(M,F)\xrightarrow{H^n(\phi)}H^n(N,F)\xrightarrow{H^n(i_1)}H^n(\Cone(\phi))\xrightarrow{H^n(p_1)}H^{n+1}(M,F)\to\cdots
\end{equation*}
which, if we forgot filtration is a long exact sequence in $\mathcal A$ ; however the morphism are NOT strict in general.

\subsection{Presheaves on a site and on a ringed topos}

\subsubsection{Functorialities}

Let $\mathcal S\in\Cat$ a small category.
For $X\in\mathcal S$ we denote by $\mathbb Z(X)\in\PSh(\mathcal S)$ the presheaf represented by $X$. 
By Yoneda lemma, a representable presheaf $\mathbb Z(X)$ is projective.

\begin{prop}
\begin{itemize}
\item Let $\mathcal S\in\Cat$ a small category.
The projective presheaves $\Proj(\PSh(\mathcal S))\subset\PSh(\mathcal S)$ are the direct summand of the representable presheaves
$\mathbb Z(X)$ with $X\in\mathcal S$. 
\item More generally let $(\mathcal S,O_S)\in\RCat$ a ringed topos.
The projective presheaves $\Proj(\PSh_{O_S}(\mathcal S))\subset\PSh_{O_S}(\mathcal S)$ of $O_S$ modules
are the direct summand of the representable presheaves $\mathbb Z(X)\otimes O_S$ with $X\in\mathcal S$. 
\end{itemize}
\end{prop}

\begin{proof}
Standard.
\end{proof}

Let $f:\mathcal T\to\mathcal S$ a morphism of presite with $\mathcal T,\mathcal S\in\Cat$.
For $h:U\to S$ a morphism with $U,S\in\mathcal S$, we have $f^*\mathbb Z(U/S)=\mathbb Z(P(f)(U/S))$.

We will consider in this article filtered complexes of presheaves on a site.
Let $f:\mathcal T\to\mathcal S$ a morphism of presite with $\mathcal T,\mathcal S\in\Cat$.
\begin{itemize}
\item The functor $f_*:C(\mathcal T)\to C(\mathcal S)$ gives, by functoriality, the functor
\begin{eqnarray*}
f_*:C_{(2)fil}(\mathcal T)\to C_{(2)fil}(\mathcal S), \; (G,F)\mapsto f_*(G,F):=(f_*G,f_*F,),
\end{eqnarray*}
since $f_*$  preserves monomorphisms.
\item The functor $f^*:C(\mathcal S)\to C(\mathcal T)$ gives, by functoriality, the functor
\begin{eqnarray*}
f^*:C_{(2)fil}(\mathcal S)\to C_{(2)fil}(\mathcal T), \; (G,F)\mapsto f^*(G,F), \; 
F^p(f^*(G,F)):=\Im(f^*F^pG\to f^*G).
\end{eqnarray*}
In the particular case where $f^*:\PSh(\mathcal S)\to\PSh(\mathcal T)$ preserves monomorphisms, 
we have $f^*(G,F)=(f^*G,f^*F)$.
\item The functor $f^{\bot}:C(\mathcal S)\to C(\mathcal T)$ gives, by functoriality, the functor
\begin{eqnarray*}
f^{\bot}:C_{(2)fil}(\mathcal T)\to C_{(2)fil}(\mathcal S), \; (G,F)\mapsto f^{\bot}(G,F):=(f^{\bot}G,f^{\bot}F),
\end{eqnarray*}
since $f^{\bot}:C(\mathcal S)\to C(\mathcal T)$  preserves monomorphisms.
\end{itemize}

Let $f:\mathcal T\to\mathcal S$ a morphism of presite with $\mathcal T,\mathcal S\in\Cat$.
\begin{itemize}
\item The adjonction
$(f^*,f_*)=(f^{-1},f_*):C(\mathcal S)\leftrightarrows C(\mathcal T)$,
gives an adjonction
\begin{eqnarray*}
(f^*,f_*):C_{(2)fil}(\mathcal S)\leftrightarrows C_{(2)fil}(\mathcal T), \;
(G,F)\mapsto f^*(G,F) \;, \; (G,F)\mapsto f_*(G,F), 
\end{eqnarray*}
with adjonction maps, for $(G_1,F)\in C_{(2)fil}(\mathcal S)$ and $(G_2,F)\in C_{(2)fil}(\mathcal T)$
\begin{eqnarray*}
\ad(f^*,f_*)(G_1,F):(G_1,F)\to f_*f^*(G_1,F) \; , \; \ad(f^*,f_*)(G_2,F):f^*f_*(G_2,F)\to (G_2,F).
\end{eqnarray*}

\item The adjonction
$(f_*,f^{\bot}):C(\mathcal S)\leftrightarrows C(\mathcal T)$,
gives an adjonction
\begin{eqnarray*}
(f_*,f^{\bot}):C_{(2)fil}(\mathcal T)\leftrightarrows C_{(2)fil}(\mathcal S), \;
(G,F)\mapsto f_*(G,F) \; , \; (G,F)\mapsto f^{\bot}(G,F),
\end{eqnarray*}
with adjonction maps, for $(G_1,F)\in C_{(2)fil}(\mathcal S)$ and $(G_2,F)\in C_{(2)fil}(\mathcal S)$
\begin{eqnarray*}
\ad(f^*,f_*)(G_2,F):(G_2,F)\to f^{\bot}f_*(G_2,F) \; , \; \ad(f^*,f_*)(G_1,F):f_*f^{\bot}(G_1,F)\to (G_1,F).
\end{eqnarray*}

\end{itemize}

\begin{rem}
Let $\mathcal T,\mathcal S\in\Cat$ small categories and $f:\mathcal T\to\mathcal S$ a morphism of presite.
Then the functor $f^*:\PSh(\mathcal S)\to\PSh(\mathcal T)$ preserve epimorphism but does NOT preserve monomorphism in general
(the colimits involved are NOT filetered colimits). However it preserve monomorphism between projective presheaves by Yoneda and we thus
set for $(Q,F)\in C_{fil}(\Proj(\PSh(\mathcal S)))$, that is $F^pQ^n\in\Proj(\PSh(\mathcal S))$ for all $p,n\in\mathbb Z$,
$f^*(Q,F):=(f^*Q,f^*F)$.
\end{rem}

For a commutative diagram of presite : 
\begin{equation*}
D=\xymatrix{ 
\mathcal T'\ar[r]^{g_2}\ar[d]^{f_2} & \mathcal S'\ar[d]^{f_1} \\
\mathcal T\ar[r]^{g_1} & \mathcal S},
\end{equation*}
with $\mathcal T,\mathcal T'\mathcal S,\mathcal S'\in\Cat$, we denote by, for $F\in C(\mathcal S')$, 
\begin{eqnarray*}
T(D)(F): g_1^*f_{1*} F\xrightarrow{g_1^*f_{1*}\ad(g_2^*,g_{2*})(F)}g_1^*f_{1*}g_{2*}g_2^*F=g_1^*g_{1*}f_{2*}g_2^*F
\xrightarrow{\ad(g_1^*g_{1*})(f_{2*}g_2^*F)} f_{2*}g_2^*F 
\end{eqnarray*}
the canonical transformation map in $C(\mathcal T)$, and for $(G,F)\in C_{fil}(\mathcal S')$, 
\begin{eqnarray*}
T(D)(G,F): g_1^*f_{1*}(G,F)\xrightarrow{g_1^*f_{1*}\ad(g_2^*,g_{2*})(G,F)}
g_1^*f_{1*}g_{2*}g_2^*(G,F)=g_1^*g_{1*}f_{2*}g_2^*(G,F) \\
\xrightarrow{\ad(g_1^*g_{1*})(f_{2*}g_2^*(G,F))} f_{2*}g_2^*(G,F). 
\end{eqnarray*}
the canonical transformation map in $C_{fil}(\mathcal T)$ given by the adjonction maps.

We will use the internal hom functor and the tensor product for presheaves on a site or for presheaves of modules on a ringed topos. 
We recall the definition in the filtrered case.
\begin{itemize}

\item Let $(\mathcal S,O_S)\in\RCat$. 
We have the tensor product bifunctor 
\begin{eqnarray*}
(\cdot)\otimes(\cdot):\PSh(\mathcal S)^2\to\PSh(\mathcal S),  
(F,G)\longmapsto(X\in\mathcal S\mapsto (F\otimes G)(X):=F(X)\otimes G(X)
\end{eqnarray*}
It induces a bifunctor :
\begin{eqnarray*}
(\cdot)\otimes(\cdot):C(\mathcal S)\times C(\mathcal S)\to C(\mathcal S),  
(F,G)\mapsto F\otimes G:=\Tot(F^{\bullet}\otimes G^{\bullet}), (F\otimes G)^n=\oplus_{r\in\mathbb Z}F^r\otimes G^{n-r}
\end{eqnarray*}
and a bifunctor 
\begin{eqnarray*}
(\cdot)\otimes(\cdot):C(\mathcal S)\times C_{O_S}(\mathcal S)\to C_{O_S}(\mathcal S), \; 
\alpha.(F\otimes G):=F\otimes(\alpha.G)
\end{eqnarray*}
For $(G_1,F),(G_2,F)\in C_{fil}(\mathcal S)$, $G_3\in C(\mathcal S)$, we define 
(note that tensor product preserve monomorphism only after tensoring with $\mathbb Q_S\in\PSh(\mathcal S)$)
\begin{itemize}
\item $F^p((G_1,F)\otimes G_3):=\Im(F^pG_1\otimes G_3\to G_1\otimes G_3)$ and
$F^p(G_3\otimes (G_1,F)):=\Im(G_3\otimes F^pG_3\to G_3\otimes G_1)$,
\item $F^pF^q((G_1,F)\otimes (G_2,F)):=\Im(F^pG_1\otimes F^qG_2\to G_1\otimes G_2)$ and
\begin{equation*}
F^k((G_1,F)\otimes (G_2,F)):=F^k\Tot_{FF}((G_1,F)\otimes (G_2,F)):=
\oplus_{p\in\mathbb Z}\Im(F^pG_1\otimes F^{k-q}G_2\to G_1\otimes G_2)
\end{equation*}
\end{itemize}
Note that in the case where $G_1^n=0$ for $n<0$, we have $(G_1,F_b)\otimes(G_2,F)=G_1\otimes(G_2,F)$.
We get the bifunctors 
\begin{eqnarray*}
(-)\otimes(-):C_{fil}(\mathcal S)^2\to C_{fil}(\mathcal S), \; \; 
(-)\otimes(-):C_{fil}(\mathcal S)\times C_{O_Sfil}(\mathcal S)\to C_{O_Sfil}(\mathcal S).
\end{eqnarray*}
We have the tensor product bifunctor 
\begin{eqnarray*}
(\cdot)\otimes_{O_S}(\cdot):\PSh_{O_S}(\mathcal S)^2\to\PSh(\mathcal S),  
(F,G)\longmapsto(X\in\mathcal S\mapsto (F\otimes_{O_S}G)(X):=F(X)\otimes_{O_S(X)}G(X)
\end{eqnarray*}
It induces a bifunctor :
\begin{eqnarray*}
(\cdot)\otimes_{O_S}(\cdot):C_{O_S}(\mathcal S)\times C_{O_S}(\mathcal S)\to C(\mathcal S),  
(F,G)\mapsto F\otimes_{O_S} G:=\Tot(F^{\bullet}\otimes_{O_S} G^{\bullet})
\end{eqnarray*}
For $(G_1,F),(G_2,F)\in C_{O_Sfil}(\mathcal S)$, $G_3\in C_{O_S}(\mathcal S)$, we define similarly
$(G_1,F)\otimes_{O_S} G_3$, $G_3\otimes_{O_S} (G_1,F)$, and
\begin{equation*}
F^k((G_1,F)\otimes_{O_S} (G_2,F)):=F^k\Tot_{FF}((G_1,F)\otimes_{O_S} (G_2,F)):=
\oplus_{p\in\mathbb Z}\Im(F^pG_1\otimes_{O_S} F^{k-q}G_2\to G_1\otimes_{O_S} G_2)
\end{equation*}
Note that in the case where $G_1^n=0$ for $n<0$, we have $(G_1,F_b)\otimes_{O_S}(G_2,F)=G_1\otimes_{O_S}(G_2,F)$. This gives
\begin{itemize}
\item in all case it gives the bifunctor 
$(-)\otimes_{O_S}(-):C_{O_S^{op}fil}(\mathcal S)\otimes C_{O_Sfil}(\mathcal S)\to C_{fil}(\mathcal S)$. 
\item in the case $O_S$ is commutative, it gives the bifunctor 
$(-)\otimes_{O_S}(-):C_{O_Sfil}(\mathcal S)^2\to C_{O_Sfil}(\mathcal S)$. 
\end{itemize}

\item Let $(\mathcal S,O_S)\in\RCat$. We have the internal hom bifunctor 
\begin{eqnarray*}
\mathcal Hom(\cdot,\cdot):\PSh(\mathcal S)^2\to\PSh(\mathcal S), \\ 
(F,G)\longmapsto(X\in\mathcal S\mapsto\mathcal Hom(F,G)(X):=\Hom(r(X)_*F,r(X)_*G)
\end{eqnarray*}
with $r(X):\mathcal S\to\mathcal S/X$ (see subsection 2.1). It induces a bifunctors :
\begin{eqnarray*}
\mathcal Hom(\cdot,\cdot):C(\mathcal S)\times C(\mathcal S)\to C(\mathcal S), \; (F,G)\mapsto\mathcal Hom^{\bullet}(F,G) 
\end{eqnarray*}
and a bifunctor 
\begin{eqnarray*}
\mathcal Hom(\cdot,\cdot):C(\mathcal S)\times C_{O_S}(\mathcal S)\to C_{O_S}(\mathcal S), \; 
\alpha.\mathcal Hom(F,G):=\mathcal Hom(F,\alpha.G)
\end{eqnarray*}
For $(G_1,F),(G_2,F)\in C_{fil}(\mathcal S)$, $G_3\in C(\mathcal S)$, we define 
\begin{itemize}
\item $F^p\mathcal\Hom(G_3,(G_1,F)):=\mathcal Hom(G_3,F^pG_1)\hookrightarrow\mathcal Hom(G_3,G_1))$,
note that the functor $G\mapsto\mathcal Hom(F,G)$ preserve monomorphism,
\item the dual filtration $F^{-p}\mathcal Hom((G_1,F),G_3):=\ker(\mathcal Hom(G_1,G_3)\to\mathcal\Hom(F^pG_1,G_3))$
\item $F^pF^q\mathcal Hom((G_1,F),(G_2,F)):=\ker(\mathcal Hom(G_1,F^pG_2)\to\mathcal\Hom(F^qG_1,F^pG_2))$, and 
\begin{eqnarray*}
F^k\mathcal Hom^{\bullet}((G_1,F),(G_2,F)):=\Tot_{FF}\mathcal Hom((G_1,F),(G_2,F)):= \\
\oplus_{p\in\mathbb Z}\ker(\mathcal Hom(G_1,F^{k+p}G_2)\to\mathcal Hom(F^pG_1,F^{k+p}G_2))
\end{eqnarray*}
\end{itemize}
We get the bifunctors 
\begin{eqnarray*}
\mathcal Hom(\cdot,\cdot):C_{fil}(\mathcal S)\times C_{fil}(\mathcal S)\to C_{fil}(\mathcal S), \; \;
\mathcal Hom(\cdot,\cdot):C_{fil}(\mathcal S)\times C_{O_Sfil}(\mathcal S)\to C_{O_Sfil}(\mathcal S). 
\end{eqnarray*}
We have the internal hom bifunctor 
\begin{eqnarray*}
\mathcal Hom_{O_S}(\cdot,\cdot):\PSh_{O_S}(\mathcal S)\times\PSh_{O_S}(\mathcal S)\to\PSh(\mathcal S) \\
(F,G)\longmapsto(X\in\mathcal S\mapsto\mathcal Hom_{O_S}(F,G)(X):=\Hom_{O_S}(r(X)_*F,r(X)_*G).
\end{eqnarray*}
It gives similarly
\begin{itemize}
\item in all case a bifunctor 
$\mathcal Hom_{O_S}(\cdot,\cdot):C_{filO_S}(\mathcal S)\times C_{filO_S}(\mathcal S)\to C_{fil}(\mathcal S)$, 
\item the case $O_S$ is commutative, a bifunctor 
$\mathcal Hom_{O_S}(\cdot,\cdot):C_{filO_S}(\mathcal S)\times C_{filO_S}(\mathcal S)\to C_{O_Sfil}(\mathcal S)$. 
\end{itemize}

\end{itemize}

Let $\phi:A\to B$ of rings.
\begin{itemize}
\item Let $M$ a $A$ module. We say that $M$ admits a $B$ module structure if there
exits a structure of $B$ module on the abelian group $M$ which is compatible with $\phi$ together with the $A$ module structure on $M$.
\item For $N_1$ a $A$-module and $N_2$ a $B$ module. 
$I(A/B)(N_1,N_2):\Hom_A(N_1,N_2)\xrightarrow{\sim}\Hom_B(N_1\otimes_A B,N_2)$ is the adjonction between 
the restriction of scalars and the extension of scalars.
\item For $N',N''$ a $A$-modules, 
$ev_A(\hom,\otimes)(N',N'',B):\Hom_A(N',N'')\otimes_A B\to\Hom_A(N',N''\otimes_A B)$. 
is the evaluation classical map.
\end{itemize}
Let $\phi:(\mathcal S,O_{1})\to(\mathcal S,O_{2})$ a morphism of presheaves of ring on $\mathcal S\in\Cat$.
\begin{itemize}
\item Let $M\in\PSh_{O_1}(\mathcal S)$. We say that $M$ admits an $O_2$ module structure if there
exits a structure of $O_2$ module on $M\in\PSh(\mathcal S)$ which is compatible with $\phi$ together with the $O_1$ module structure on $M$.
\item For $N_1\in C_{O_1}(\mathcal S)$ and $N_2\in C_{O_2}(\mathcal S)$, 
\begin{equation*}
I(O_1/O_2)(N_1,N_2):\Hom_{O_1}(N_1,N_2)\xrightarrow{\sim}\Hom_{O_2}(N_1\otimes_{O_1} B,N_2)
\end{equation*}
is the adjonction between the restriction of scalars and the extension of scalars.
\item For $N',N''\in C_{O_1}(\mathcal S)$, 
\begin{equation*}
ev_{O_1}(\hom,\otimes)(N',N'',O_2):\mathcal Hom_{O_1}(N',N'')\otimes_{O_1} {O_2}\to\mathcal Hom_{O_1}(N',N''\otimes_{O_1} O_2). 
\end{equation*}
is the classical evaluation map.
\end{itemize}

Let $(\mathcal S,O_S)\in\RCat$. 
\begin{itemize}
\item For $F_1,F_2,G_1,G_2\in C(\mathcal S)$, we denote by 
\begin{equation*}
T(\otimes,\mathcal Hom)(F_1,F_2,G_1,G_2):
\mathcal Hom(F_1,G_1)\otimes\mathcal Hom(F_2,G_2)\to\mathcal Hom(F_1\otimes F_2,G_1\otimes G_2)
\end{equation*}
the canonical map.
\item For $G_3\in C(\mathcal S)$ and $G_1,G_2\in C_{O_S}(\mathcal S)$, we denote by
\begin{eqnarray*}
ev(hom,\otimes)(G_3,G_1,G_2):\mathcal Hom(G_3,G_1)\otimes_{O_S}G_2\to\mathcal Hom(G_3,G_1\otimes_{O_S}G_2) \\
\phi\otimes s\mapsto (s'\mapsto\phi(s')\otimes s)
\end{eqnarray*}
\item Let $\mathcal S\in\Cat$ a small category. Let $(H_X:C(\mathcal S/X)\to C(\mathcal S/X))_{X\in\mathcal S}$ a familly of functors
which is functorial in $X$. We have by definition, for $F_1,F_2\in C(\mathcal S)$, the canonical transformation map
\begin{eqnarray}\label{THhom0}
T(H,hom)(F_1,F_2):H(\mathcal Hom^{\bullet}(F_1,F_2))\to\mathcal Hom^{\bullet}(H(F_1),H(F_2)), \\  
\phi\in\Hom(F_{1|X},F_{2|X})\mapsto H^{F_{1|X},F_{2|X}}(\phi)\in\Hom(H(F_{1|X}),H(F_{2|X}))
\end{eqnarray}
in $C(\mathcal S)$.
\end{itemize}

Let $\mathcal T,\mathcal S\in\Cat$ small categories and $f:\mathcal T\to\mathcal S$ a morphism of presite.
\begin{itemize}
\item For $F_1,F_2\in C(\mathcal T)$ we have by definition $f_*(F_1\otimes F_2)=f_*F_1\otimes f_*F_2$. 
For $G_1,G_2\in C(\mathcal S)$,
we have a canonical isomorphism $f^*G_1\otimes f^*G_2\xrightarrow{\sim}f^*(G_1\otimes G_2)$ 
since the tensor product is a right exact functor,
and a canonical map $f^{\bot}G_1\otimes f^{\bot}G_2\to f^{\bot}(G_1\otimes G_2)$.
\item We have for $F\in C(\mathcal S)$ and $G\in C(\mathcal T)$ the adjonction isomorphim,   
\begin{equation}\label{Ifhom}
I(f^*,f_*)(F,G):f_*\mathcal Hom^{\bullet}(f^*F,G)\xrightarrow{\sim}\mathcal Hom^{\bullet}(F,f_*G).  
\end{equation}
\item Let $O_S\in\PSh(\mathcal S,\Ring)$ by a presheaf of ring so that $(\mathcal S,O_S),(\mathcal T,f^*O_S)\in\RCat$.
We have for $F\in C_{O_S}(\mathcal S)$ and $G\in C_{f^*O_S}(\mathcal T)$ the adjonction isomorphim,   
\begin{equation}\label{IfhomO}
I(f^*,f_*)(F,G):f_*\mathcal Hom_{f^*O_S}^{\bullet}(f^*F,G)\xrightarrow{\sim}\mathcal Hom_{O_S}^{\bullet}(F,f_*G),  
\end{equation}
and 
\begin{itemize}
\item the map $\ad(f^*,f_*)(F):F\to f_*f^*F$ in $C(\mathcal S)$ is $O_S$ linear, that is is a map in $C_{O_S}(\mathcal S)$,
\item the map $\ad(f^*,f_*)(G):f^*f_*G\to G$ in $C(\mathcal T)$ is $f^*O_S$ linear, that is is a map in $C_{f^*O_S}(\mathcal T)$.
\end{itemize}
\item For $F_1,F_2\in C(\mathcal T)$, we have the canonical map
\begin{eqnarray}\label{Tfhom0}
T_*(f,hom)(F_1,F_2):=T(f_*,hom):f_*\mathcal Hom^{\bullet}(F_1,F_2)\to\mathcal Hom^{\bullet}(f_*F_1,f_*F_2), \\ 
\mbox{for} \, X\in\mathcal S, \, \phi\in\Hom(F_{1|f^*(X)},F_{2|f^*(X)})\mapsto 
{f_*}^{F_{1|f^*(X)},F_{2|f^*(X)}}(\phi)\in\Hom(f_*F_{1|f^*(X)},f_*F_{2|f^*(X)})
\end{eqnarray}
given by evaluation.
\item For $G_1,G_2\in C(\mathcal S)$, we have the following canonical transformation in $C(\mathcal T)$
\begin{eqnarray}\label{Tfhom}
T(f,hom)(G_1,G_2):=T(f^*,hom)(G_1,G_2): \\
f^*\mathcal Hom^{\bullet}(G_1,G_2)\xrightarrow{f^*\mathcal Hom(G_1,\ad(f^*,f_*)(G_2))}
f^*\mathcal Hom^{\bullet}(G_1,f_*f^*G_2)\xrightarrow{f^*I(f^*,f_*)(G_1,G_2)} \\
f^*f_*\mathcal Hom^{\bullet}(f^*G_1,f^*G_2)\xrightarrow{\ad(f^*,f_*)(\mathcal Hom(f^*G_1,f^*G_2))}
\mathcal Hom^{\bullet}(f^*G_1,f^*G_2),  
\end{eqnarray}
\item Let $O_S\in\PSh(\mathcal S,\Ring)$ by a presheaf of ring so that $(\mathcal S,O_S),(\mathcal T,f^*O_S)\in\RCat$.
For $G_1,G_2\in C_{O_S}(\mathcal S)$, we have the following canonical transformation in $C_{f^*O_S}(\mathcal T)$
\begin{eqnarray}\label{TfhomO}
T(f,hom)(G_1,G_2):=T(f^*,hom)(G_1,G_2): \\
f^*\mathcal Hom_{O_S}^{\bullet}(G_1,G_2)\xrightarrow{f^*\mathcal Hom_{O_S}(G_1,\ad(f^*,f_*)(G_2))}
f^*\mathcal Hom_{O_S}^{\bullet}(G_1,f_*f^*G_2)\xrightarrow{f^*I(f^*,f_*)(G_1,G_2)} \\
f^*f_*\mathcal Hom_{f^*O_S}^{\bullet}(f^*G_1,f^*G_2)\xrightarrow{\ad(f^*,f_*)(\mathcal Hom_{f^*O_S}(f^*G_1,f^*G_2))}
\mathcal Hom^{\bullet}_{f^*O_S}(f^*G_1,f^*G_2),  
\end{eqnarray}
\item Let $O_S\in\PSh(\mathcal S,\Ring)$ by a presheaf of ring so that $(\mathcal S,O_S),(\mathcal T,f^*O_S)\in\RCat$.
For $M\in C_{O_S}(\mathcal S)$ and $N\in C_{f^*O_S}(\mathcal T)$, we denote by
\begin{eqnarray}\label{Tfotimes}
T(f,\otimes)(M,N): M\otimes_{O_S} f_*N\xrightarrow{\ad(f^*,f_*)(M\otimes_{O_S} f_*N)} \\
f_*f^*(M\otimes_{O_S} f_*N)=f_*(f^*M\otimes_{f^*O_S} f^*f_*N)
\xrightarrow{\ad(f^*,f_*)(N)}f_*(f^*M\otimes_{f^*O_S}N)
\end{eqnarray}
the canonical transformation map.
\end{itemize}

Let $f:(\mathcal T,O_T)\to(\mathcal S,O_S)$ a morphism with $(\mathcal S,O_S),(\mathcal T,O_T)\in\RCat$.
We have the adjonction
\begin{equation*}
(f^{*mod},f_*):C_{O_S}(\mathcal S)\leftrightarrows C_{O_T}(\mathcal T)
\end{equation*}
with $f^{*mod}G:=f^*G\otimes_{f^*O_S}O_T$. If $f^*:C(\mathcal S)\to C(\mathcal T)$
preserve monomorphisms, it induces the adjonction
\begin{equation*}
(f^{*mod},f_*):C_{O_Sfil}(\mathcal S)\leftrightarrows C_{O_Tfil}(\mathcal T)
\end{equation*}
with $f^{*mod}(G,F):=f^*(G,F)\otimes_{f^*O_S}O_T$.

For a commutative diagram in $\RCat$ : 
\begin{equation*}
D=\xymatrix{ 
(\mathcal T',O'_2)\ar[r]^{g_2}\ar[d]^{f_2} & (\mathcal S',O'_1)\ar[d]^{f_1} \\
(\mathcal T,O_2)\ar[r]^{g_1} & (\mathcal S,O_1)},
\end{equation*}
we denote by, for $F\in C_{O'_1}(\mathcal S')$, 
\begin{eqnarray*}
T^{mod}(D)(F): g_1^{*mod}f_{1*} F\xrightarrow{g_1^{*mod}f_{1*}\ad(g_2^{*mod},g_{2*})(F)} 
g_1^{*mod}f_{1*}g_{2*}g_2^{*mod}F=g_1^{*mod}g_{1*}f_{2*}g_2^{*mod}F \\
\xrightarrow{\ad(g_1^{*mod}g_{1*})(f_{2*}g_2^{*mod}F)} f_{2*}g_2^{*mod}F 
\end{eqnarray*}
the canonical transformation map in $C_{O_2}(\mathcal T)$ and, for $(G,F)\in C_{O'_1fil}(\mathcal S')$, 
\begin{eqnarray*}
T^{mod}(D)(G,F): g_1^{*mod}f_{1*}(G,F)\xrightarrow{g_1^{*mod}f_{1*}\ad(g_2^{*mod},g_{2*})(G,F)} 
g_1^{*mod}f_{1*}g_{2*}g_2^{*mod}(G,F)=g_1^{*mod}g_{1*}f_{2*}g_2^{*mod}(G,F) \\
\xrightarrow{\ad(g_1^{*mod}g_{1*})(f_{2*}g_2^{*mod}(G,F))} f_{2*}g_2^{*mod}(G,F) 
\end{eqnarray*}
the canonical transformation map in $C_{O_2fil}(\mathcal T)$ given by the adjonction maps.

Let $f:(\mathcal T,O_T)\to(\mathcal S,O_S)$ a morphism with $(\mathcal S,O_S),(\mathcal T,O_T)\in\RCat$.
\begin{itemize}
\item We have, for $M,N\in C_{O_S}(\mathcal S)$ the canonical transformation map in $C_{O_T}(\mathcal T)$
\begin{eqnarray*}
T^{mod}(f,hom)(M,N):f^{*mod}\mathcal\Hom_{O_1}(M,N)\xrightarrow{T(f,hom)(M,N)\otimes_{f^*O_1} O_2}
\mathcal Hom_{f^*O_1}(f^*M,f^*N)\otimes_{f^*O_1}O_2 \\ \xrightarrow{e(\hom,\otimes)(f^*M,f^*N)}\mathcal Hom_{f^*O_1}(f^*M,f^{*mod}N)
\xrightarrow{I(f^*O_1/O_2)(f^*M,f^{*mod}N)}\mathcal Hom_{O_2}(f^{*mod}M,f^{*mod}N)
\end{eqnarray*}

\item We have, for $M\in C_{O_S}(\mathcal S)$ and $N\in C_{O_T}(\mathcal T)$, the canonical transformation map in $C_{O_T}(\mathcal T)$
\begin{eqnarray}\label{Tfotimesmod}
T^{mod}(f,\otimes)(M,N): M\otimes_{O_{\mathcal S}} f_*N\xrightarrow{\ad(f^{*mod},f_*)(M\otimes_{O_S} f_*N)} \\
f_*f^{*mod}(M\otimes_{O_S} f_*N)=f_*(f^{*mod}M\otimes_{O_T} f^{*mod}f_*N)
\xrightarrow{\ad(f^{*mod},f_*)(N)}f_*(f^{*mod}M\otimes_{O_T}N)
\end{eqnarray}
the canonical transformation map.
\end{itemize}

We now give some properties of the tensor product functor and hom functor given above

\begin{prop}\label{otimesderQ}
Let $(\mathcal S,O_S)\in\RCat$. Then, the functors
\begin{itemize}
\item $(-)\otimes(-):C(\mathcal S)^2\to C(\mathcal S), \; C(\mathcal S)\times C_{O_S}(\mathcal S)\to C_{O_S}(\mathcal S)$ 
\item $(-)\otimes_{O_S}(-):C_{O_S^{op}}(\mathcal S)\times C_{O_S}(\mathcal S)\to C(\mathcal S)$ and in case $O_S$ is commutative 
$(-)\otimes_{O_S}(-):C_{O_S}(\mathcal S)^2\to C_{O_S}(\mathcal S)$
\end{itemize}
are left Quillen functor for the projective model structure. In particular, 
\begin{itemize}
\item for $L\in C(\mathcal S)$ is such that $L^n\in\PSh(\mathcal S)$ are projective for all $n\in\mathbb Z$, 
and $\phi:F\to G$ is a quasi-isomorphism with $F,G\in C(\mathcal S)$, then
$\phi\otimes I:F\otimes L\to G\otimes L$ is a quasi-isomorphism, 
\item for $L\in C_{O_S}(\mathcal S)$ is such that $L^n\in\PSh_{O_S}(\mathcal S)$ are projective for all $n\in\mathbb Z$, 
and $\phi:F\to G$ is a quasi-isomorphism with $F,G\in C_{O_S}(\mathcal S)$, then
$\phi\otimes I:F\otimes_{O_S} L\to G\otimes_{O_S}L$ is a quasi-isomorphism. 
\end{itemize}
\end{prop}

\begin{proof}
Standard.
\end{proof}

\begin{prop}\label{homderQ}
Let $(\mathcal S,O_S)\in\RCat$. Then, the functors
\begin{itemize}
\item $\mathcal Hom(\cdot,\cdot):
C(\mathcal S)\times C_{O_S}(\mathcal S)\to C_{O_S}(\mathcal S), \; C(\mathcal S)\times C_{O_S}(\mathcal S)\to C_{O_S}(\mathcal S)$,
\item $\mathcal Hom_{O_S}(\cdot,\cdot):
C_{O_S^{op}}(\mathcal S)\times C_{O_S}(\mathcal S)\to C(\mathcal S)$ and in the case $O_S$ is commutative
$\mathcal Hom_{O_S}(\cdot,\cdot):C_{O_S}(\mathcal S)\times C_{O_S}(\mathcal S)\to C_{O_S}(\mathcal S)$,
\end{itemize}
are on the left hand side left Quillen functor for the projective model structure. In particular, 
\begin{itemize}
\item for $L\in C(\mathcal S)$ is such that $L^n\in\PSh(\mathcal S)$ are projective for all $n\in\mathbb Z$, 
and $\phi:F\to G$ is a quasi-isomorphism with $F,G\in C(\mathcal S)$, then
$\mathcal Hom(L,\phi):\mathcal Hom^{\bullet}(L,F)\to \mathcal Hom^{\bullet}(L,G)$
is a quasi-isomorphism, 
\item for $L\in C_{O_S}(\mathcal S)$ is such that $L^n\in\PSh_{O_S}(\mathcal S)$ are projective for all $n\in\mathbb Z$, 
and $\phi:F\to G$ is a quasi-isomorphism with $F,G\in C_{O_S}(\mathcal S)$, then
$\mathcal Hom_{O_S}(L,\phi):\mathcal Hom^{\bullet}_{O_S}(L,F)\to \mathcal Hom^{\bullet}_{O_S}(L,G)$
is a quasi-isomorphism. 
\end{itemize}
\end{prop}

\begin{proof}
Standard.
\end{proof}

Let $\mathcal S\in\Cat$ a site endowed with topology $\tau$.
Denote by $a_{\tau}:\PSh(\mathcal S)\to\Sh(\mathcal S)$ the sheaftification functor
A morphism $\phi:F^{\bullet}\to G^{\bullet}$ with $F^{\bullet},G^{\bullet}\in C(\mathcal S))$
is said to be a $\tau$ local equivalence if 
\begin{equation*}
a_{\tau}H^n(\phi):a_{\tau}H^n(F^{\bullet})\to a_{\tau}H^n(G^{\bullet})
\end{equation*}
is an isomorphism for all $n\in\mathbb Z$, where $a_{\tau}$ is the sheaftification functor.
Recall that $C_{fil}(\mathcal S)\subset(C(\mathcal S),F)=C(\PSh(\mathcal S),F)$ denotes
the category of filtered complexes of abelian presheaves on $\mathcal S$ whose filtration is biregular.
\begin{itemize}
\item A morphism $\phi:(F^{\bullet},F)\to (G^{\bullet},F)$ with $(F^{\bullet},F),(G^{\bullet},F)\in C_{fil}(\mathcal S)$
is said to be a filtered $\tau$ local equivalence or an $1$-filtered $\tau$ local equivalence if 
\begin{equation*}
a_{\tau}H^n(\phi):a_{\tau}H^n(\Gr_F^pF^{\bullet})\xrightarrow{\sim} a_{\tau}H^n(\Gr_F^pG^{\bullet})
\end{equation*}
is an isomorphism for all $n,p\in\mathbb Z$.
\item Let $r\in\mathbb N$. A morphism
$\phi:(F^{\bullet},F)\to (G^{\bullet},F)$ with $(F^{\bullet},F),(G^{\bullet},F)\in C_{fil}(\mathcal S)$
is said to be an $r$-filtered $\tau$ local equivalence if it belongs to the
submonoid of arrows generated by filtered $\tau$ local equivalences and $r$-filtered homotopy equivalences,
that is if there exists $\phi_i:(C_i^{\bullet},F)\to(C_{i+1}^{\bullet},F)$, $0\leq i\leq s$, 
with $(C_i^{\bullet},F)\in C_{fil}(\mathcal S)$,
$(C_0^{\bullet},F)=(F^{\bullet},F)$ and $(C_s^{\bullet},F)=(G^{\bullet},F)$, such that
\begin{equation*}
\phi=\phi_s\circ\cdots\circ\phi_i\circ\cdots\circ\phi_0:(F^{\bullet},F)\to(G^{\bullet},F)
\end{equation*}
and $\phi_i:(C_i^{\bullet},F)\to(C_{i+1}^{\bullet},F)$ either a filtered $\tau$ local equivalence
or an $r$-filtered homotopy equivalence.
If $\phi:(F^{\bullet},F)\to (G^{\bullet},F)$ with $(F^{\bullet},F),(G^{\bullet},F)\in C_{fil}(\mathcal S)$
is an $r$-filtered $\tau$ local equivalence, then for all $p,q\in\mathbb Z$,
\begin{equation*}
a_{\tau}E^{p,q}_r(\phi):a_{\tau}E^{p,q}_r(F^{\bullet},F)\xrightarrow{\sim}a_{\tau}E^{p,q}_r(G^{\bullet},F)
\end{equation*}
is an isomorphism but the converse is NOT true. 
Note that if $\phi$ is an $r$-filtered $\tau$ local equivalence, that it is an $s$-filtered $\tau$ local
equivalence for all $s\geq r$.
\item A morphism $\phi:(F^{\bullet},F)\to (G^{\bullet},F)$ with $(F^{\bullet},F),(G^{\bullet},F)\in C_{fil}(\mathcal S)$
is said to be a $\infty$-filtered $\tau$ local equivalence if there exists $r\in\mathbb N$
such that $\phi$ is an $r$-filtered $\tau$ local equivalence. 
If a morphism $\phi:(F^{\bullet},F)\to (G^{\bullet},F)$ with $(F^{\bullet},F),(G^{\bullet},F)\in C_{fil}(\mathcal S)$
is an $\infty$-filtered $\tau$ local equivalence then, for all $n\in\mathbb Z$,
\begin{equation*}
a_{\tau}H^n(\phi):a_{\tau}H^n(F^{\bullet},F)\to a_{\tau}H^n(G^{\bullet},F)
\end{equation*}
is an isomorphism of filtered sheaves on $\mathcal S$. Recall the converse is true in the case there exist $N_1,N_2\in\mathbb Z$, 
such that $H^n(F^{\bullet},F)=H^n(G^{\bullet},F)=0$ for $n\leq N_1$ or $n\geq N_2$.
\end{itemize}
Let $(\mathcal S,O)$ a ringed topos where $\mathcal S\in\Cat$ is a site endowed with topology $\tau$.
A morphism $\phi:(F^{\bullet},F)\to (G^{\bullet},F)$ with $(F^{\bullet},F),(G^{\bullet},F)\in C_{O_Sfil}(\mathcal S)$
is said to be a filtered $\tau$ local equivalence or an $1$-filtered $\tau$ local equivalence if 
$o\phi:(F^{\bullet},F)\to (G^{\bullet},F)$ is one in $C_{fil}(\mathcal S)$, i.e.
\begin{equation*}
a_{\tau}H^n(\phi):a_{\tau}H^n(\Gr_F^pF^{\bullet})\xrightarrow{\sim} a_{\tau}H^n(\Gr_F^pG^{\bullet})
\end{equation*}
is an isomorphism for all $n,p\in\mathbb Z$.
Let $r\in\mathbb N$.
A morphism $\phi:(F^{\bullet},F)\to(G^{\bullet},F)$ with $(F^{\bullet},F),(G^{\bullet},F)\in C_{O_Sfil}(\mathcal S)$
is said to be an $r$-filtered $\tau$ local equivalence if 
there exists $\phi_i:(C_i^{\bullet},F)\to(C_{i+1}^{\bullet},F)$, $0\leq i\leq s$, with $(C_i^{\bullet},F)\in C_{O_Sfil}(\mathcal S)$
$(C_0^{\bullet},F)=(F^{\bullet},F)$ and $(C_s^{\bullet},F)=(G^{\bullet},F)$, such that
\begin{equation*}
\phi=\phi_s\circ\cdots\circ\phi_i\circ\cdots\circ\phi_0:(F^{\bullet},F)\to(G^{\bullet},F)
\end{equation*}
and $\phi_i:(C_i^{\bullet},F)\to(C_{i+1}^{\bullet},F)$ either a filtered $\tau$ local equivalence
or an $r$-filtered homotopy equivalence.

Let $\mathcal S\in\Cat$ a site which admits fiber product, endowed with topology $\tau$.
A complex of presheaves $F^{\bullet}\in C(\mathcal S)$ is said to be $\tau$ fibrant if it satisfy descent for covers in $\mathcal S$,
i.e. if for all $X\in\mathcal S$ and all $\tau$ covers $(c_i:U_i\to X)_{i\in I}$ of $X$,
denoting $U_J:=(U_{i_0}\times_SU_{i_1}\times_S\cdots U_{i_r})_{i_k\in J}$ and for $I\subset J$, $p_{IJ}:U_J\to U_I$ is the projection,
\begin{equation*}
F^{\bullet}(c_i):F^{\bullet}(X)\to\Tot(\oplus_{card I=\bullet} F^{\bullet}(U_I))
\end{equation*}
is a quasi-isomorphism of complexes of abelian groups.
\begin{itemize}
\item A complex of filtered presheaves $(F^{\bullet},F)\in C_{fil}(\mathcal S)$ is said to be filtered $\tau$ fibrant
or $1$-filtered $\tau$ fibrant if it satisfy descent for covers in $\mathcal S$,
i.e. if for all $X\in\mathcal S$ and all $\tau$ covers $(c_i:U_i\to X)_{i\in I}$ of $X$,
\begin{equation*}
(F^{\bullet},F)(c_i):(F^{\bullet},F)(X)\to\Tot(\oplus_{card I=\bullet} (F^{\bullet},F)(U_I))
\end{equation*}
is a filtered quasi-isomorphism of filtered complexes of abelian groups.
\item Let $r\in\mathbb N$. A complex of filtered presheaves $(F^{\bullet},F)\in C_{fil}(\mathcal S)$ 
is said to be $r$-filtered $\tau$ fibrant if there exist an $r$-filtered homotopy equivalence 
$m:(F^{\bullet},F)\to (F^{'\bullet},F)$ with $(F^{'\bullet},F)\in C_{fil}(\mathcal S)$ filtered $\tau$ fibrant.
If $(F^{\bullet},F)\in C_{fil}(\mathcal S)$ is $r$-filtered $\tau$ fibrant,
then for all $X\in\mathcal S$ and all $\tau$ covers $(c_i:U_i\to X)_{i\in I}$ of $X$,
\begin{equation*}
E_r^{p,q}(F^{\bullet},F)(c_i):E_r^{p,q}(F^{\bullet},F)(X)\to E_r^{p,q}(\Tot(\oplus_{card I=\bullet}(F^{\bullet},F)(U_I)))
\end{equation*}
is an isomorphism for all $n,p\in\mathbb Z$, but the converse is NOT true.
Note that if $(F^{\bullet},F)$ is $r$-filtered $\tau$ fibrant, then it is $s$-filtered $\tau$ fibrant for all $s\geq r$. 
\item A complex of filtered presheaves $(F^{\bullet},F)\in C_{fil}(\mathcal S)$ is said to be $\infty$-filtered $\tau$ fibrant
if there exist $r\in\mathbb N$ such that $(F^{\bullet},F)$ is $r$-filtered $\tau$ fibrant.
If a complex of filtered presheaves $(F^{\bullet},F)\in C_{fil}(\mathcal S)$ 
is $\infty$-filtered $\tau$ fibrant, then for all $X\in\mathcal S$ and all $\tau$ covers $(c_i:U_i\to X)_{i\in I}$ of $X$,
\begin{equation*}
H^n(F^{\bullet},F)(c_i):H^n(F^{\bullet},F)(X)\to H^n\Tot(\oplus_{card I=\bullet}(F^{\bullet},F)(U_I))
\end{equation*}
is a filtered isomorphism for all $n\in\mathbb Z$.
\end{itemize}
Let $(\mathcal S,O)$ a ringed topos where $\mathcal S\in\Cat$ is a site endowed with topology $\tau$. Let $r\in\mathbb N$.
\begin{itemize}
\item A complex of presheaves $F^{\bullet}\in C_{O_S}(\mathcal S)$ is said to be $\tau$ fibrant 
if $F^{\bullet}=oF^{\bullet}\in C(\mathcal S)$ is $\tau$ fibrant.
\item A complex of presheaves $(F^{\bullet},F)\in C_{O_Sfil}(\mathcal S)$ is said to be filtered $\tau$ fibrant 
if $(F^{\bullet},F)=(oF^{\bullet},F)\in C_{fil}(\mathcal S)$ is filtered $\tau$ fibrant.
\item A complex of presheaves $(F^{\bullet},F)\in C_{O_Sfil}(\mathcal S)$ is said to be $r$-filtered $\tau$ fibrant 
if there exist an $r$-filtered homotopy equivalence
$m:(F^{\bullet},F)\to (F^{'\bullet},F)$ with $(F^{'\bullet},F)\in C_{O_Sfil}(\mathcal S)$ filtered $\tau$ fibrant.
\end{itemize}

\subsubsection{Canonical flasque resolution of a presheaf on a site or a presheaf of module on a ringed topos}

Let $\mathcal S\in\Cat$ a site with topology $\tau$. 
Denote $a_{\tau}:\PSh(\mathcal S)\to\Shv(\mathcal S)$ the sheaftification functor.
There is for $F\in C(\mathcal S)$ an explicit $\tau$ fibrant replacement : 
\begin{itemize}
\item $k:F^{\bullet}\to E_{\tau}^{\bullet}(F^{\bullet}):=\Tot (E_{\tau}^{\bullet}(F^{\bullet}))$,
if $F^{\bullet}\in C^{+}(\mathcal S)$,  
\item $k:F^{\bullet}\to E_{\tau}^{\bullet}(F^{\bullet}):=\holim\Tot (E_{\tau}^{\bullet}(F^{\bullet\geq n}))$, 
if $F^{\bullet}\in C(\mathcal S)$ is not bounded below.  
\end{itemize}
The bicomplex $E^{\bullet}(F^{\bullet}):=E_{\tau}^{\bullet}(F^{\bullet})$ 
together with the map $k:F^{\bullet}\to E^{\bullet}(F^{\bullet})$ 
is given inductively by
\begin{itemize}
\item considering $p_{\mathcal S}:\mathcal S^{\delta}\to\mathcal S$ the morphism of site from the discrete
category $\mathcal S^{\tau}$ whose objects are the points of the topos $\mathcal S$ and we take
\begin{eqnarray*}
k_0:=\ad(p_S^*,p_{S*})(F^{\bullet})\to E^0(F^{\bullet}):=p_{S*}p_S^*F^{\bullet}:=
\bigoplus_{s\in\mathcal S^{\tau}}\lim_{X\in\mathcal S, s\in X}F^{\bullet}(X),
\end{eqnarray*}  
then $a_{\tau}k_0:a_{\tau}F^{\bullet}\to E^0(F^{\bullet})$ is injective and $E^0(F^{\bullet})$ is $\tau$ fibrant,
\item denote $Q^0(F^{\bullet}):=a_{\tau}\coker(k_0:F^{\bullet}\to E^0(F^{\bullet}))$ and take the composite 
\begin{equation*}
E^0(F^{\bullet})\to Q^0(F^{\bullet})\to E^1(F^{\bullet}):=E^0(Q^0(F^{\bullet})).
\end{equation*}
\end{itemize}
Note that $k:F^{\bullet}\to E^{\bullet}(F^{\bullet})$ is a $\tau$ local equivalence 
and that $a_{\tau}k:a_{\tau}F^{\bullet}\to E^{\bullet}(F^{\bullet})$ is injective by construction. 

Since $E^0$ is functorial, $E$ is functorial:
for $m:F^{\bullet}\to G^{\bullet}$ a morphism, with $F^{\bullet},G^{\bullet}\in C(\mathcal S)$,
we have a canonical morphism $E(m):E(F)\to E(G)$ such that $E(m)\circ k=k'\circ m$, with $k:F\to E(F)$ and $k':G\to E(G)$. 
Note that $E^0$, hence $E$ preserve monomorphisms.
This gives, for $(F^{\bullet},F)\in C_{fil}(\mathcal S)$, 
a filtered $\tau$ local equivalence $k:(F^{\bullet},F)\to E^{\bullet}(F^{\bullet},F)$ 
with $E^{\bullet}(F^{\bullet},F)$ filtered $\tau$ fibrant.

Moreover, we have a canonical morphism $E(F)\otimes E(G)\to E(F\otimes G)$.

There is, for $g:\mathcal T\to\mathcal S$ a morpism of presite with $\mathcal T,\mathcal S\in\Cat$ two site,
and $F^{\bullet}\in C(\mathcal S)$, a canonical transformation
\begin{equation}\label{gE}
T(g,E)(F^{\bullet}):g^*E(F^{\bullet})\to E(g^*F^{\bullet})
\end{equation}
given inductively by
\begin{itemize}
\item $T(g,E^0)(F):=T(g,p_S)(p_S^*F):g^*E^0(F)=g^*p_{S*}p_S^*F\to p_{T*}g^*p_S^*F=p_{T*}p_T^*g^*F=E^0(g^*F)$,
$T(g,Q^0)(F):=\overline{T(g,E^0)(F)}:g^*Q^0(F)=\coker(g^*F\to g^*E^0(F))\to Q^0(g^*F)=\coker(g^*F\to E^0(g^*F)$
\item $T(g,Q^1)(F):g^*E^1(F)=g^*E^0(Q^0(F))\xrightarrow{T(g,E^0)(Q^0(F))}E^0(g^*Q^0(F))
\xrightarrow{E^0(T(g,Q^0)(F))}E^0(Q^0(g^*F))=E^1(g^*F)$.
\end{itemize}

Let $(\mathcal S,O_S)\in\RCat$ with topology $\tau$. 
Then, for $F^{\bullet}\in C_{O_S}(\mathcal S)$, $E_{\tau}(F^{\bullet})$ is naturally a complex of $O_S$ modules
such that $k:F^{\bullet}\to E_{\tau}(F^{\bullet})$ is $O_S$ linear, that is is a morphism in $C_{O_S}(\mathcal S)$.

\subsubsection{Canonical projective resolution of a presheaf of module on a ringed topos}

Let $(\mathcal S,O_S)\in\RCat$. 
We recall that we denote by, for $U\in\mathcal S$,
$\mathbb Z(U)\in\PSh(\mathcal S)$ the presheaf represented by $U$ : for $V\in\mathcal S$
$\mathbb Z(U)(V)=\mathbb Z\Hom(V,U)$, and for $h:V_1\to V_2$ a morphism in $\mathcal S$, and $h_1:V_1\to U$
$\mathbb Z(U)(h):h_1\to h\circ h_1$, and $s$ is the morphism of presheaf given by $s(V_1)(h_1)=F(h_1)(s)\in F(V_1)$.
There is for $F\in C_{O_S}(\mathcal S)$ a complex of $O_S$ module an explicit projective replacement : 
\begin{itemize}
\item $q:L_O^{\bullet}(F^{\bullet}):=\Tot (L_O^{\bullet}(F^{\bullet}))\to F^{\bullet}$, if $F^{\bullet}\in C^{-}(\mathcal S)$,  
\item $q:L_O^{\bullet}(F^{\bullet}):=\holim\Tot (L_O^{\bullet}(F^{\bullet\leq n}))$
if $F^{\bullet}\in C(\mathcal S)$ is not bounded above.  
\end{itemize}
For $O_S=\mathbb Z_{\mathcal S}$, we denote $L^{\bullet}_{\mathbb Z_{\mathcal S}}(F^{\bullet})=:L(F^{\bullet})$.
The bicomplex $L_O^{\bullet}(F^{\bullet})$ together with the map $q:L_O^{\bullet}(F^{\bullet})\to F^{\bullet}$ is given inductively by
\begin{itemize}
\item considering the pairs $\left\{U\in\mathcal S, s\in F(U)\right\}$,
where $U$ is an object of $\mathcal S$ and $s$ a section of $F$ over $U$ we take
\begin{equation*}
q_0:L_O^0(F):=\bigoplus_{(U\in\mathcal S,s\in F(U))}\mathbb Z(U)\otimes O_S\xrightarrow{s} F,
\end{equation*}  
then $q_0$ is surjective and $L_O^0(F)$ is projective,
this construction is functorial : for $m:F\to G$ a morphism in $\PSh(\mathcal S)$ the following diagram commutes
\begin{equation*}
\xymatrix{\bigoplus_{(U\in\mathcal S,s\in F(U))}\mathbb Z(U)\otimes O_S\ar[rr]^{q_0}\ar[d]_{L_O(m)} & \, & F\ar[d]^{m} \\
\bigoplus_{(U\in\mathcal S,s'\in G(U))}\mathbb Z(U)\otimes O_S\ar[rr]^{q_0} & \, & G}
\end{equation*}
where $(L_O(m)_{|(U,s)})_{(U,m(U)(s))}=I_{\mathbb Z(U)}$ and $(L_O(m)_{|(U,s)})_{(U,s')}=0$ if $s'\neq m(U)(s)$, 
\item denote $K_O^0(F):=\ker(q_0:L_O^0(F)\to F))$ and take the composite 
\begin{equation*}
q_1:L_O^1(F^{\bullet}):=L_O^0(K_O^0(F^{\bullet}))\xrightarrow{q_0(K_O^0(F))} K_O^0(F^{\bullet})\hookrightarrow L_O^0(F^{\bullet}).
\end{equation*}
\end{itemize}
Note that $q=q(F):L(F^{\bullet})\to F^{\bullet}$ is a surjective quasi-isomorphism by construction. 
Since $L_O^0$ is functorial, $L_O$ is functorial :
for $m:F^{\bullet}\to G^{\bullet}$ a morphism, with $F^{\bullet},G^{\bullet}\in C(\mathcal S)$,
we have a canonical morphism $L_O(m):L_O(F)\to L_O(G)$ such that $q'\circ L_O(m)=m\circ q'$, with $q:L_O(F)\to F$ and $q':L_O(G)\to G$. 
Note that $L_O^0$ and hence $L_O$ preserve monomorphisms.
In particular, it gives for $(F^{\bullet},F)\in C_{O_S}(\mathcal S)$, a filtered quasi-isomorphism
$q:L_O(F^{\bullet},F)\to(F^{\bullet},F)$.
Moreover, we have a canonical morphism $L_O(F)\otimes L_O(G)\to L_O(F\otimes G)$.

Let $g:\mathcal T\to\mathcal S$ a morphism of presite with $\mathcal T,\mathcal S\in\Cat$ two sites.
\begin{itemize}
\item Let $F^{\bullet}\in C(\mathcal S)$. Since $g^*L(F^{\bullet})$ is projective and $q(g^*F):L(g^*F^{\bullet})\to g^*F^{\bullet}$ is
a surjective quasi-isomorphism, there is a canonical transformation 
\begin{equation}\label{gL}
T(g,L)(F^{\bullet}):g^*L(F^{\bullet})\to L(g^*F^{\bullet})
\end{equation}
unique up to homotopy such that $q(g^*F)\circ T(g,L)(F^{\bullet})=g^*q(F)$.
\item Let $F^{\bullet}\in C(\mathcal S)$. Since $L(g^*F^{\bullet})$ is projective and $g^*q(F):g^*L(F^{\bullet})\to g^*F^{\bullet}$ is
a surjective quasi-isomorphism, there is a canonical transformation 
\begin{equation}\label{gL3}
T(g,L)(F^{\bullet}):L(g^*F^{\bullet})\to g^*L(F^{\bullet})
\end{equation}
unique up to homotopy such that $g^*q(F)\circ T(g,L)(F^{\bullet})=q(g^*F)$.
\item Let $F^{\bullet}\in C(\mathcal T)$. Since $L(g_*F^{\bullet})$ is projective and $g_*q(F):g_*L(F^{\bullet})\to g_*F^{\bullet}$ is
a surjective quasi-isomorphism, there is a canonical transformation 
\begin{equation}\label{gL2}
T_*(g,L)(F^{\bullet}):L(g_*F^{\bullet})\to g_*L(F^{\bullet})
\end{equation}
unique up to homotopy such that $g_*q(F)\circ T_*(g,L)(F^{\bullet})=q(g_*F)$.
\end{itemize}

Let $g:(\mathcal T,O_T)\to(\mathcal S,O_S)$ a morphism with $(\mathcal T,O_T),(\mathcal S,O_S)\in\RCat$.
Let $F^{\bullet}\in C_{O_S}(\mathcal S)$. 
Since $g^{*mod}L_O(F^{\bullet})$ is projective and $q(g^{*mod}F):L_O(g^{*mod}F^{\bullet})\to g^{*mod}F^{\bullet}$ is
a surjective quasi-isomorphism, there is a canonical transformation 
\begin{equation}\label{gLO}
T(g,L_O)(F^{\bullet}):g^{*mod}L_O(F^{\bullet})\to L_O(g^{*mod}F^{\bullet})
\end{equation}
unique up to homotopy such that $q(g^{*mod}F)\circ T(g,L_O)(F^{\bullet})=g^{*mod}q(F)$.

Let $p:(\mathcal S_{12},O_{S_{12}})\to(\mathcal S_1,O_{S_1})$ a morphism 
with $(\mathcal S_{12},O_{S_{12}}),(\mathcal S_1,O_{S_1})\in\RCat$,
such that the structural morphism $p^*O_{S_1}\to O_{S_{12}}$ is flat.
Let $F^{\bullet}\in C_{O_S}(\mathcal S)$. 
Since $L_O(p^{*mod}F^{\bullet})$ is projective and $p^{*mod}q(F):p^{*mod}L_O(F^{\bullet})\to p^{*mod}F^{\bullet}$ is
a surjective quasi-isomorphism, there is also in this case a canonical transformation 
\begin{equation}\label{pL2}
T(p,L_O)(F^{\bullet}):L_O(p^{*mod}F^{\bullet})\to p^{*mod}L_O(F^{\bullet})
\end{equation}
unique up to homotopy such that $p^{*mod}q(F)\circ T(p,L_O)(F^{\bullet})=q(p^{*mod}F)$.

\subsubsection{The De Rham complex of a ringed topos and functorialities}

Let $A\in\cRing$ a commutative ring. For $M\in\Mod(A)$, we denote by 
\begin{equation*}
\Der_A(A,M)\subset\Hom(A,M)=\Hom_{\Ab}(A,M)
\end{equation*}
the abelian subgroup of derivation. 
Denote by $I_A=\ker(s_A:A\otimes A\to A)\subset A\otimes A$ the diagonal ideal with $s_A(a_1,a_2)=a_1-a_2$. 
Let $\Omega_A:=I_A/I_A^2\in\Mod(A)$ together with its derivation map $d=d_A:A\to\Omega_A$. 
Then, for $M\in\Mod(A)$ the canonical map
\begin{equation*}
w(M):\Hom_A(\Omega_A,M)\xrightarrow{\sim}\Der_A(A,M), \; \psi\mapsto \phi\circ d
\end{equation*}
is an isomorphism, that is $\Omega_A$ is the universal derivation.
In particular, its dual $T_A:=\mathbb D^A(\Omega_A)=D^A(I_A/I_A^2)$ is isomorphic to the derivations group :
$w(A):T_A\xrightarrow{\sim}\Der_A(A,A)$. Also note that $\Der_A(A,A)\subset\Hom(A,A)$ is a Lie subalgebra.
If $\phi:A\to B$ is a morphism of commutative ring, we have a canonical morphism of abelian group 
$\Omega_{(B/A)}\phi:\Omega_A\to\Omega_B$.

Let $(\mathcal S,O_S)\in\RCat$, with $O_S\in\PSh(\mathcal S,\cRing)$ commutative. For $G\in\PSh_{O_S}(\mathcal S)$, we denote by 
\begin{equation*}
\Der_{O_S}(O_S,G)\subset\mathcal Hom(O_S,G)=\mathcal Hom_{\Ab}(O_S,G)
\end{equation*}
the abelian subpresheaf of derivation. 
Denote by $\mathcal I_S=\ker(s_S:O_S\otimes O_S\to O_S)\in\PSh_{O_S\times O_S}(\mathcal S)$ the diagonal ideal 
with $s_S(X)=s_{O_S(X)}$ for $X\in\mathcal S$. 
Then $\Omega_{O_S}:=\mathcal I_S/\mathcal I_S^2\in\PSh_{O_S}(\mathcal S)$ 
together with its derivation map $d:O_S\to\Omega_{O_S}$ is the universal derivation $O_S$-module : 
the canonical map
\begin{equation*}
w(G):\mathcal Hom_{O_S}(\Omega_{O_S},G)\xrightarrow{\sim}\Der_{O_S}(O_S,G), \; \phi\mapsto \phi\circ d
\end{equation*}
is an isomorphism. In particular, its dual 
$T_{O_S}:=\mathbb D^O_{O_S}(\Omega_{O_S})=\mathbb D_{\mathcal S}^O(\mathcal I_S/\mathcal I_S^2)$ 
is isomorphic to the presheaf of derivations :
$w(O_S):T_{O_S}\xrightarrow{\sim}\Der_{O_S}(O_S,O_S)$ and $\Der_{O_S}(O_S,O_S)\subset\mathcal\Hom(O_S,O_S)$ is a Lie subalgebra.
The universal derivation $d=d_{O_S}:O_S\to\Omega_{O_S}$ induces the De Rham complex 
\begin{equation*}
DR(O_S):\Omega^{\bullet}_{\mathcal S}:=\wedge^{\bullet}\Omega_{O_S}\in C(\mathcal S) 
\end{equation*}
A morphism $\phi:O'_S\to O_S$ with $O_S,O'_S\PSh(\mathcal S,\cRing)$ induces by the universal property canonical morphisms
\begin{equation*}
\Omega_{O'_S/O_S}:\Omega_{O'_S}\to\Omega_{O_S} \; , \; \mathbb D^O_{O_S}\Omega_{O'_S/O_S}:T_{O_S}\to T_{O'_S}
\end{equation*}
in $\PSh_{O_S}(\mathcal S)$.
\begin{itemize}
\item In the particular cases where $S=(S,O_S)\in\Var(\mathbb C)$ or $S=(S,O_S)\in\AnSp(\mathbb C)$, 
we denote as usual $\Omega_S:=\Omega_{O_S/\mathbb C_S}$, $T_S:=T_{O_S/\mathbb C_S}$ 
and $DR(S):=DR(O_S/\mathbb C_S):\Omega_{S}^{\bullet}\in C(S)$.
\item In the particular cases where $S=(S,O_S)\in\Diff(\mathbb R)$ is a differential manifold, 
we denote as usual $\mathcal A_S:=\Omega_{O_S/\mathbb R_S}$, $T_S:=T_{O_S/\mathbb R_S}$ and 
$DR(S):=DR(O_S/\mathbb R_S):\mathcal A_{S}^{\bullet}\in C(S)$.
\end{itemize}

For $f:(\mathcal X,O_X)\to(\mathcal S,O_S)$ with $(\mathcal X,O_X),(\mathcal S,O_S)\in\RCat$ such that $O_X$ and $O_S$ are commutative, 
we denote by
\begin{equation*}
\Omega_{O_X/f^*O_S}:=
\coker(\Omega_{O_X/f^*O_S}:\Omega_{f^*O_S}\to\Omega_{O_X})\in\PSh_{f^*O_S}(\mathcal X)
\end{equation*}
the relative cotangent sheaf. The surjection $q=q_{O_X/f}:\Omega_{O_X}\to\Omega_{O_X/f^*O_S}$ gives the
derivation $w(\Omega_{O_X/f^*O_S})(q)=d_{O_X/f}:O_X\to\Omega_{O_X/f^*O_S}$.
It induces the surjections $q^p:=\wedge^pq:\Omega^p_{O_X}\to\Omega^p_{O_X/f^*O_S}$. We then have the realtive De Rham complex
\begin{equation*}
DR(O_X/f^*O_S):=\Omega^{\bullet}_{O_X/f^*O_S}:=\wedge^{\bullet}\Omega_{O_X/f^*O_S}\in C_{f^*O_S}(\mathcal X). 
\end{equation*}
whose differnetials are given by
\begin{equation*}
\mbox{for}\; X^o\in\mathcal X \;\mbox{and} \; \omega\in\Gamma(X^o,\Omega^p_{O_X})\; d(q^p(\omega)):=q^{p+1}(d(\omega))
\end{equation*}
Note that $\Omega^{\bullet}_{O_X/f^*O_S}\in C_{f^*O_S}(\mathcal S)$ is a complex of $f^*O_S$ modules, 
but is NOT a complex of $O_X$ module since the differential is a derivation hence NOT $O_X$ linear.
On the other hand, the canonical map in $\PSh_{f^*O_S}(\mathcal S)$
\begin{equation*}
T(f,hom)(O_S,O_S):f^*\mathcal Hom(O_S,O_S)\to\mathcal Hom(f^*O_S,f^*O_S)
\end{equation*}
induces morphisms 
\begin{equation*}
T(f,hom)(O_S,O_S):f^*T_{O_S}\to T_{f^*O_S} \mbox{ and } \mathbb D^O_{f^*O_S}T(f,hom)(O_S,O_S):\Omega_{f^*O_S}\to f^*\Omega_{O_S}.
\end{equation*}
In this article, We will be interested in the following particular cases : 
\begin{itemize}
\item In the particular case where $O_S\PSh(\mathcal S,\cRing)$ is a sheaf, $\Omega_{O_S},T_{O_S}\in\PSh_{O_S}(\mathcal S)$ are sheaves.
Hence, 
\begin{equation*}
T(f,hom)(O_S,O_S):a_{\tau}f^*T_{O_S}\xrightarrow{\sim} T_{a_{\tau}f^*O_S} \mbox{\; and \; } 
\mathbb D^O_{f^*O_S}T(f,hom)(O_S,O_S):\Omega_{a_{\tau}f^*O_S}\xrightarrow{\sim} a_{\tau}f^*\Omega_{O_S}
\end{equation*}
are isomorphisms,where $a_{\tau}:\PSh(\mathcal S)\to\Shv(\mathcal S)$ is the sheaftification functor.
We will note again in this case by abuse (as usual)
$f^*O_S:=a_{\tau}f^*O_S$, $f^*\Omega_{O_S}:=a_{\tau}f^*\Omega_{O_S}$ and $f^*T_{O_S}:=a_{\tau}f^*T_{O_S}$, so that 
\begin{equation*}
\Omega_{f^*O_S}=f^*\Omega_{O_S} \mbox{\; and \;} f^*T_{O_S}=T_{f^*O_S}
\end{equation*}
\item In the particular cases where $S=(S,O_S),X=(X,O_X)\in\Var(\mathbb C)$ or $S=(S,O_S),X=(X,O_X)\in\AnSp(\mathbb C)$,
we denote as usual $\Omega_{X/S}:=\Omega_{O_X/f^*O_S}$, $q_{X/S}:=q_{O_X/S}:\Omega_X\to\Omega_{X/S}$ 
and $DR(X/S):=DR(O_X/f^*O_S):\Omega_{X/S}^{\bullet}\in C_{f^*O_S}(S)$.
\item In the particular cases where $S=(S,O_S),X=(X,O_X)\in\Diff(\mathbb R)$, 
we denote as usual $\mathcal A_{X/S}:=\Omega_{O_X/f^*O_S}$, $q_{X/S}:=q_{O_X/S}:\mathcal A_X\to\mathcal A_{X/S}$ and
 $DR(X/S):=DR(O_X/f^*O_S):\mathcal A_{X/S}^{\bullet}\in C_{f^*O_S}(S)$.
\end{itemize} 

\begin{defi}\label{TDw}
For a commutative diagram in $\RCat$ 
\begin{equation*}
D=\xymatrix{(\mathcal X,O_X)\ar[r]^{f} & (\mathcal S,O_S) \\
(\mathcal X',O_{X'})\ar[u]^{g'}\ar[r]^{f'} & (\mathcal T,O_T)\ar[u]^{g}},
\end{equation*}
whose structural presheaves are commutative sheaves, the map in $C_{g^{'*}O_Xfil}(\mathcal X')$
\begin{equation*}
\Omega_{O_{X'}/g^{'*}O_X}:g^{'*}(\Omega^{\bullet}_{O_X},F_b)=(\Omega^{\bullet}_{g^{'*}O_X},F_b)\to(\Omega^{\bullet}_{O_{X'}},F_b)  
\end{equation*}
pass to quotient to give the map in $C_{g^{'*}O_Xfil}(\mathcal X')$
\begin{eqnarray*}
\Omega_{(O_Y/g^{'*}O_X)/(O_T/g^*O_S)}:=(\Omega_{O_Y/g^{'*}O_X})^q: \\
g^{'*}(\Omega^{\bullet}_{O_X/f^*O_S},F_b)=(\Omega^{\bullet}_{g^{'*}O_X/g^{'*}f^*O_S},F_b)\to(\Omega^{\bullet}_{O_Y/f^{'*}O_T},F_b) 
\end{eqnarray*}
It is in particular given for $X^{'o}\in\mathcal X'$, $g'^*(X^o)\leftarrow X^{'o}$ and $\hat\omega\in\Gamma(X^o,\Omega^p_{O_X/f^*O_S})$, 
\begin{eqnarray*}
\Omega_{(O_{X'}/g^{'*O_X})/(O_T/g^*O_S)}(X^{'o})(\omega):=
q_{O_{X'}/f'}(\Omega_{O_{X'}/g^{'*}O_X}(\omega))\in\Gamma(X^{'o},\Omega^p_{O_{X'}/f^{'*}O_T}).
\end{eqnarray*} 
where $\omega\in\Gamma(X^o,\Omega^p_{O_X})$ such that $q_{O_X/f}(\omega)=\hat\omega$.
We then have the following canonical transformation map in $C_{O_Tfil}(\mathcal T)$
\begin{eqnarray*}
T^O_{\omega}(D): 
g^{*mod}L_Of_*E(\Omega^{\bullet}_{O_X/f^*O_S},F_b)
\xrightarrow{q}
g^*f_*E(\Omega^{\bullet}_{O_X/f^*O_S},F_b)\otimes_{g^*O_S}O_T \\
\xrightarrow{T(g',E)(-)\circ T(D)(E(\Omega^{\bullet}_{O_X/f^*O_S}))}   
(f'_*E(g^{'*}(\Omega^{\bullet}_{O_X/f^*O_S},F_b)))\otimes_{g^*O_S}O_T \\
\xrightarrow{E(\Omega_{(O_{X'}/g^{'*}O_X)/(O_T/g^*O_S}))} 
f'_*E(\Omega^{\bullet}_{O_{X'}/f^{'*}O_T},F_b)\otimes_{g^*O_S}O_T 
\xrightarrow{m}
f'_*E(\Omega^{\bullet}_{O_{X'}/f^{'*}O_T},F_b),
\end{eqnarray*}
with $m(n\otimes s)=s.n$.
\end{defi}

\subsection{Presheaves on diagrams of sites or on diagrams of ringed topos}

Let $\mathcal I,\mathcal I'\in\Cat$ and
$(f_{\bullet},s):\mathcal T_{\bullet}\to\mathcal S_{\bullet}$ a morphism of diagrams of presites 
with $\mathcal T_{\bullet}\in\Fun(\mathcal I,\Cat),\mathcal S_{\bullet}\in\Fun(\mathcal I',\Cat)$.
Recall it is by definition given by a functor $s:\mathcal I\to I'$ and morphism of functor 
$P(f_{\bullet}):\mathcal S_{s(\bullet)}:=\mathcal S_{\bullet}\circ s\to\mathcal T_{\bullet}$.
and that we denote for short, $\mathcal S_{s(\bullet)}:=\mathcal S_{\bullet}\circ s\in\Fun(\mathcal I,\Cat)$.
Recall that, for $r_{IJ}:I\to J$ a morphism, with $I,J\in\mathcal I$, $D_{fIJ}$ is the commutative diagram in $\Cat$ 
\begin{equation*}
D_{fIJ}:=\xymatrix{\mathcal S_{s(J)}\ar[r]^{r^s_{IJ}} & \mathcal S_{s(I)} \\
\mathcal T_J\ar[r]^{r^t_{IJ}}\ar[u]^{f_J} & \mathcal T_I\ar[u]^{f_I}}.
\end{equation*}
The adjonction
\begin{eqnarray*}
((f_{\bullet},s)^*,(f_{\bullet},s)_*)=((f_{\bullet},s)^{-1},(f_{\bullet},s)_*):
C(\mathcal S_{s(\bullet)})\leftrightarrows C(\mathcal T_{\bullet}), \\
G=(G_I,u_{IJ})\mapsto (f_{\bullet},s)^*(G):=(f_I^*(G_I),T(D_{fIJ})(G_J)\circ f_I^*u_{IJ})  \\
G=(G_I,u_{IJ})\mapsto (f_{\bullet},s)_*(G):=(f_{I*}(G_I),f_{I*}u_{IJ})
\end{eqnarray*}
gives an adjonction
\begin{eqnarray*}
((f_{\bullet},s)^*,(f_{\bullet},s)_*):C_{(2)fil}(\mathcal S_{s(\bullet)})\leftrightarrows C_{(2)fil}(\mathcal T_{\bullet}), \\
(G,F)=((G_I,F),u_{IJ})\mapsto (f_{\bullet},s)^*(G,F):=(f_I^*(G_I,F),T(D_{fIJ})(G_J,F)\circ f_I^*u_{IJ})  \\
(G,F)=((G_I,F),u_{IJ})\mapsto (f_{\bullet},s)_*(G,F):=(f_{I*}(G_I,F),f_{I*}u_{IJ}). 
\end{eqnarray*}

For a commutative diagram of diagrams of presite : 
\begin{equation*}
D=\xymatrix{ 
\mathcal T'_{\bullet}\ar[r]^{(g_2,s'_2)}\ar[d]^{(f_2,s_2)} & \mathcal S'_{s'_2(\bullet)}\ar[d]^{(f_1,s_1)} \\
\mathcal T_{s_2(\bullet)}\ar[r]^{(g_1,s'_1)} & \mathcal S_{s(\bullet)}},
\end{equation*}
with $\mathcal I,\mathcal I',\mathcal J,\mathcal J'\in\Cat$ and
$\mathcal T_{\bullet}\in\Fun(\mathcal I,\Cat),\mathcal T'_{\bullet}\in\Fun(\mathcal I',\Cat),
\mathcal S_{\bullet}\in\Fun(\mathcal J,\Cat),\mathcal S'_{\bullet}\in\Fun(\mathcal J',\Cat)$, 
and $s=s_1\circ s'_2=s_2\circ s'_1:\mathcal I'\to J$,
we denote by, for $F=(F_I,u_{IJ})\in C(\mathcal S'_{s'_2(\bullet)})$, 
\begin{eqnarray*}
T(D)(F): g_1^*f_{1*} F\xrightarrow{g_1^*f_{1*}\ad(g_2^*,g_{2*})(F)}g_1^*f_{1*}g_{2*}g_2^*F=g_1^*g_{1*}f_{2*}g_2^*F
\xrightarrow{\ad(g_1^*g_{1*})(f_{2*}g_2^*F)} f_{2*}g_2^*F 
\end{eqnarray*}
the canonical transformation map in $C(\mathcal T_{s_2(\bullet)})$, and
for $(G,F)=((G_I,F),u_{IJ})\in C_{fil}(\mathcal S'_{s'_2(\bullet)})$, 
\begin{eqnarray*}
T(D)(G,F): g_1^*f_{1*}(G,F)\xrightarrow{g_1^*f_{1*}\ad(g_2^*,g_{2*})(G,F)}
g_1^*f_{1*}g_{2*}g_2^*(G,F)=g_1^*g_{1*}f_{2*}g_2^*(G,F) \\
\xrightarrow{\ad(g_1^*g_{1*})(f_{2*}g_2^*(G,F))} f_{2*}g_2^*(G,F) 
\end{eqnarray*}
the canonical transformation map in $C_{fil}(\mathcal T_{s_2(\bullet)})$ given by the adjonction maps.

Let $\mathcal S_{\bullet}\in\Fun(\mathcal I,\RCat)$ a diagram of ringed topos with $\mathcal I\in\Cat$.
We have the tensor product bifunctor 
\begin{eqnarray*}
(\cdot)\otimes(\cdot):\PSh(\mathcal S_{\bullet})^2\to\PSh(\mathcal S_{\bullet}), \\ 
((F_I,u_{IJ},(G_I,u_{IJ}))\mapsto (F_I,u_{IJ})\otimes (G_I,v_{IJ}):=(F_I\otimes G_I,u_{IJ}\otimes v_{IJ})
\end{eqnarray*}
We get the bifunctors 
\begin{eqnarray*}
(-)\otimes(-):C_{fil}(\mathcal S_{\bullet})^2\to C_{fil}(\mathcal S_{\bullet}), \; \; 
(-)\otimes(-):C_{fil}(\mathcal S_{\bullet})\times C_{O_Sfil}(\mathcal S_{\bullet})\to C_{O_Sfil}(\mathcal S_{\bullet}).
\end{eqnarray*}
We have the tensor product bifunctor 
\begin{eqnarray*}
(\cdot)\otimes_{O_S}(\cdot):\PSh_{O_S}(\mathcal S_{\bullet})^2\to\PSh(\mathcal S_{\bullet}), \\  
((F_I,u_{IJ},(G_I,u_{IJ}))\mapsto (F_I,u_{IJ})\otimes_{O_{S_I}} (G_I,v_{IJ}):=
(F_I\otimes_{O_{S_I}} G_I,u_{IJ}\otimes v_{IJ})
\end{eqnarray*}
which gives,
\begin{itemize}
\item in all case it gives the bifunctor 
$(-)\otimes_{O_S}(-):C_{O_S^{op}fil}(\mathcal S_{\bullet})\otimes C_{O_Sfil}(\mathcal S_{\bullet})\to 
C_{fil}(\mathcal S_{\bullet})$. 
\item in the case $O_S$ is commutative, it gives the bifunctor 
$(-)\otimes_{O_S}(-):C_{O_Sfil}(\mathcal S_{\bullet})^2\to C_{O_Sfil}(\mathcal S_{\bullet})$. 
\end{itemize}

Let $(f_{\bullet},s):(\mathcal T_{\bullet},O_T)\to(\mathcal S_{\bullet},O_S)$ a morphism with 
$(\mathcal S_{\bullet},O_S)\in\Fun(\mathcal I',\RCat),(\mathcal T_{\bullet},O_T)\in\Fun(\mathcal I,\RCat)$
and $\mathcal I,\mathcal I'\in\Cat$. 
which is by definition given by a functor $s:\mathcal I\to \mathcal I'$ and morphism of ringed topos 
$f_{\bullet}:(\mathcal T_{\bullet},O_T)\to(\mathcal S_{s(\bullet)},O_S)$.
As before, we denote for short, 
$(\mathcal S_{s(\bullet)},O_S):=(\mathcal S_{\bullet},O_S)\circ s\in\Fun(\mathcal I,\RCat)$.
Denote as before, for $r_{IJ}:I\to J$ a morphism, with $I,J\in\mathcal I$, $D_{fIJ}$ the commutative diagram in $\RCat$ 
\begin{equation*}
D_{fIJ}:=\xymatrix{\mathcal S_{s(J)}\ar[r]^{r^s_{IJ}} & \mathcal S_{s(I)} \\
\mathcal T_J\ar[r]^{r^t_{IJ}}\ar[u]^{f_J} & \mathcal T_I\ar[u]^{f_I}}.
\end{equation*}
We have then the adjonction
\begin{eqnarray*}
((f_{\bullet},s)^{*mod},(f_{\bullet},s)_*):C_{O_S}(\mathcal S_{s(\bullet)})\leftrightarrows C_{O_T}(\mathcal T_{\bullet}), \\
(G_I,u_{IJ})\mapsto(f_{\bullet},s)^{*mod}(G_I,u_{IJ}):=(f_I^{*mod}G_I,T^{mod}(D_{fIJ})(G_J)\circ f_I^{*mod}u_{IJ}), \\
(G_I,u_{IJ})\mapsto(f_{\bullet},s)_*(G_I,u_{IJ}):=(f_{I*}G_I,f_{I*}u_{IJ}). 
\end{eqnarray*}
which induces the adjonction
\begin{eqnarray*}
((f_{\bullet},s)^{*mod},(f_{\bullet},s)*):
C_{O_Sfil}(\mathcal S_{s(\bullet)})\leftrightarrows C_{O_Tfil}(\mathcal T_{\bullet}), \\
((G_I,F),u_{IJ})\mapsto(f_{\bullet},s)^{*mod}((G_I,F),u_{IJ}):=(f_I^{*mod}(G_I,F),T^{mod}(D_{fIJ})(G_J)\circ f_I^{*mod}u_{IJ}), \\
((G_I,F),u_{IJ})\mapsto (f_{\bullet},s)_*((G_I,F),u_{IJ}):=(f_{I*}(G_I,F),f_{I*}u_{IJ}).
\end{eqnarray*}

For a commutative diagram of diagrams of ringed topos, : 
\begin{equation*}
D=\xymatrix{ 
(\mathcal T'_{\bullet},O'_2)\ar[r]^{(g_2,s'_2)}\ar[d]^{(f_2,s_2)} & (\mathcal S'_{s'_2(\bullet)},O'_1)\ar[d]^{(f_1,s_1)} \\
(\mathcal T_{\bullet},O_2)\ar[r]^{(g_1,s'_1)} & (\mathcal S_{s(\bullet)},O_1)},
\end{equation*}
with $\mathcal I,\mathcal I',\mathcal J,\mathcal J'\in\Cat$ and
$\mathcal T_{\bullet}\in\Fun(\mathcal I,\Cat),\mathcal T'_{\bullet}\in\Fun(\mathcal I',\Cat),
\mathcal S_{\bullet}\in\Fun(\mathcal J,\Cat),\mathcal S'_{\bullet}\in\Fun(\mathcal J',\Cat)$, 
and $s=s_1\circ s'_2=s'_1\circ s_2:\mathcal I'\to\mathcal J$,
we denote by, for $F=(F_I,u_{IJ})\in C_{O'_1}(\mathcal S'_{s'_2(\bullet)})$, 
\begin{eqnarray*}
T^{mod}(D)(F): g_1^{*mod}f_{1*} F\xrightarrow{g_1^{*mod}f_{1*}\ad(g_2^{*mod},g_{2*})(F)} 
g_1^{*mod}f_{1*}g_{2*}g_2^{*mod}F=g_1^{*mod}g_{1*}f_{2*}g_2^{*mod}F \\
\xrightarrow{\ad(g_1^{*mod}g_{1*})(f_{2*}g_2^{*mod}F)} f_{2*}g_2^{*mod}F 
\end{eqnarray*}
the canonical transformation map in $C_{O_2}(\mathcal T_{s_2(\bullet)})$, and
for $G=((G_I,F),u_{IJ})\in C_{O'_1fil}(\mathcal S'_{s'_2(\bullet)})$, 
\begin{eqnarray*}
T^{mod}(D)(G,F): g_1^{*mod}f_{1*}(G,F)\xrightarrow{g_1^{*mod}f_{1*}\ad(g_2^{*mod},g_{2*})(G,F)} 
g_1^{*mod}f_{1*}g_{2*}g_2^{*mod}(G,F)=g_1^{*mod}g_{1*}f_{2*}g_2^{*mod}(G,F) \\
\xrightarrow{\ad(g_1^{*mod}g_{1*})(f_{2*}g_2^{*mod}(G,F))} f_{2*}g_2^{*mod}(G,F) 
\end{eqnarray*}
the canonical transformation map in $C_{O_2fil}(\mathcal T_{s_2(\bullet)})$ given by the adjonction maps.

Let $(\mathcal S_{\bullet},O_S)\in\Fun(\mathcal I,\RCat)$ a diagram of ringed topos with 
$\mathcal I\in\Cat$ and, for $I\in\mathcal I$, $\mathcal S_I$ is endowed with topology $\tau_I$
and for $r:I\to J$ a morphism with $I,J\in\mathcal I$, $r_{IJ}:\tilde S_J\to\tilde S_I$ is continous. 
Then the diagram category $(\Gamma\mathcal S_{\bullet},O_S)\in\RCat$ 
is endowed with the associated canonical topology $\tau$, and then
\begin{itemize}
\item A morphism $\phi=(\phi_I):((F_I,F),u_{IJ})\to((G_I,F),u_{IJ})$ with 
$((F_I,F),u_{IJ}),((G_I,F),u_{IJ})\in C_{O_Sfil}(\mathcal S_{\bullet})$
is a filtered $\tau$ local equivalence if and only if the $\phi_I$ are filtered $\tau$ local equivalences
for all $I\in\mathcal I$.
\item Let $r\in\mathbb N$. A morphism $\phi=(\phi_I):((F_I,F),u_{IJ})\to((G_I,F),v_{IJ})$ with 
$((F_I,F),u_{IJ}),((G_I,F),v_{IJ})\in C_{O_Sfil}(\mathcal S_{\bullet})$
is an $r$-filtered $\tau$ local equivalence if and only if 
there exists $\phi_i:((C_{iI},F),u_{iIJ})\to((C_{(i+1)I},F),u_{(i+1)IJ})$, $0\leq i\leq s$, 
with $((C_{iI},F),u_{iIJ})\in C_{O_Sfil}(\mathcal S_{\bullet})$,
$((C_{0I},F),u_{0IJ})=((F_I,F),u_{IJ})$ and $((C_{sI},F),u_{sIJ})=((G_I,F),v_{IJ})$, such that
\begin{equation*}
\phi=\phi_s\circ\cdots\circ\phi_i\circ\cdots\circ\phi_0:((F_I,F),u_{IJ})\to((G_I,F),v_{IJ})
\end{equation*}
and $\phi_i:(C_i^{\bullet},F)\to(C_{i+1}^{\bullet},F)$ either a filtered $\tau$ local equivalence
or an $r$-filtered homotopy equivalence.
the $\phi_I$ are $r$-filtered $\tau$ local equivalences
for all $I\in\mathcal I$.
\item A complex of presheaves $((G_I,F),u_{IJ})\in C_{O_Sfil}(\mathcal S_{\bullet})$ is filtered $\tau$ fibrant
if and only if the $(G_I,F)\in C_{O_Sfil}(\mathcal S_I)$ are filtered $\tau$ fibrant for all $I\in\mathcal I$.
\item Let $r\in\mathbb N$. A complex of presheaves
$((G_I,F),u_{IJ})\in C_{O_Sfil}(\mathcal S_{\bullet})$ is $r$-filtered $\tau$ fibrant
if there exist an $r$-filtered homotopy equivalence $m:((G_I,F),u_{IJ})\to ((G'_I,F),u_{IJ})$ with
$(G'_I,F)\in C_{O_Sfil}(\mathcal S_{\bullet})$ filtered $\tau$ fibrant.
\end{itemize}

\subsection{Presheaves on topological spaces and presheaves of modules on a ringed spaces}

In this subsection, we will consider the particular case of presheaves on topological spaces.

Let $f:T\to S$ a continous map with $S,T\in\Top$. We denote as usual the adjonction
\begin{equation*}
(f^*,f_*):=(P(f)^*,P(f)_*):\PSh(S)\leftrightarrows\PSh(T)
\end{equation*}
induced by the morphism of site given by the pullback functor 
\begin{equation*}
P(f):\Ouv(S)\to\Ouv(T), \; (S^o\subset S)\mapsto P(f)(S^o):=S^o\times_S T\xrightarrow{\sim} f^{-1}(S^o)\subset T
\end{equation*}
Since the colimits involved in the definition of $f^*=P(f)^*$ are filtered, $f^*$ also preserve monomorphism. Hence, we get an adjonction
\begin{equation*}
(f^*,f_*):\PSh_{fil}(S)\leftrightarrows\PSh_{fil}(T), \; f^*(G,F):=(f^*G,f^*F)
\end{equation*}

Let $f:(T,O_T)\to (S,O_S)$ a morphism with $(S,O_S),(T,O_T)\in\Top$. We have then the adjonction
\begin{equation*}
(f^{*mod},f_*):=(P(f)^{*mod},P(f)_*):\PSh_{O_Sfil}(S)\leftrightarrows\PSh_{O_Tfil}(T), \; f^{*mod}(G,F):=f^*(G,F)\otimes_{f^*O_S}O_T
\end{equation*}

Recall $\CW\subset\Top$ is the full subcategory whose objects consists of $CW$ complexes.
Denote, for $n\in\mathbb N$, $\mathbb I^n:=[0,1]^n,S^n:=\mathbb I^n/\partial\mathbb I^n\in\CW$ and 
$\Delta^n\subset I^n$ the $n$ dimensional simplex. We get $\mathbb I^*,\Delta^*\in\Fun(\Delta,\CW)$ 
Denote for $S\in\Top$, $\Sigma_1S:=S\times\mathbb I^1/((\left\{0\right\}\times S)\cup(\left\{1\right\}\times S))\in\Top$.
\begin{itemize}
\item Let $f:T\to S$ a morphism with $T,S\in\Top$. 
We have the mapping cylinder $\Cyl(f):=(T\times I^1)\sqcup_f S\in\Top$ 
and the mapping cone $\Cone(f):=(T\times\mathbb I^1)\sqcup_f S\in\Top$.
We have then the quotient map $q_f:\Cyl(f)\to\Cone(f)$ and a canonical retraction $r_f:\Cone(f)\to\Sigma^1T$
\item Recall two morphisms $f,g:T\to S$ with $T,S\in\Top$ are homotopic if there
exist $H:T\times I^1\to S$ continous such that $H\circ (I\times i_0)=f$ and $H\circ (I\times i_1)=g$.
Then $K(\Top):=\Ho_{I^1}(\Top)$ is a triangulated category with distinguish triangle
\begin{equation*}
T\xrightarrow{i_T}\Cyl(f)\xrightarrow{q_f}\Cone(f)\xrightarrow{r_f}\Sigma^1T.
\end{equation*}
\item For $X\in\Top$, denote for $n\in\mathbb N$, $\pi_n(X):\Hom_{K(\Top)}(S^n,X)$ the homotopy groups.
For $f:T\to S$ a morphism with $T,S\in\Top$, we have for $n\in\mathbb N$ the morphisms of abelian groups
\begin{equation*}
f_*:\pi_n(T)\to\pi_n(S), h\mapsto f\circ h
\end{equation*}
Recall two morphisms $f,g:T\to S$ with $T,S\in\Top$ are weakly homotopic if 
$f_*=g_*:\pi_n(T)\to\pi_n(S)$ for all $n\in\mathbb N$.
\item For $X\in\Top$, denote by
$C^{\sing}_*(X):=\mathbb Z\Hom(\Delta^*,X)\in C^-(\mathbb Z)$ the complex of singular chains and by 
$C_{\sing}^*(X):=\mathbb D^{\mathbb Z}C^{\sing}_*(X):=\mathbb D\mathbb Z\Hom(\Delta^*,X)\in C^-(\mathbb Z)$ 
the complex of singular cochains.
For $f:T\to S$ a morphism with $T,S\in\Top$, we have 
\begin{itemize}
\item the morphism of complexes of abelian groups
\begin{equation*}
f_*:C^{\sing}_*(T)\to C^{\sing}_*(S), \sigma\mapsto f\circ\sigma, \\ 
\end{equation*}
\item the morphism of complexes of abelian groups
\begin{eqnarray*}
f^*:=\mathbb D^{\mathbb Z}f_*:C_{\sing}^*(T)\to C_{\sing}^*(S), \; 
\alpha\longmapsto f^*\alpha:(\sigma\mapsto f^*\alpha(\sigma):=\alpha(f\circ\sigma))
\end{eqnarray*}
\end{itemize}
We denote by $C^*_{X,\sing}\in C^+(X)$ the complex of presheaves of singular cochains given by,
\begin{eqnarray*}
(U\subset X)\mapsto C^*_{X,\sing}(U):=C^*_{X,\sing}(U):=C_{\sing}^*(U):=\mathbb D^{\mathbb Z}\mathbb Z\Hom(\Delta^*,U), \\
(j:U_2\hookrightarrow U_1)\mapsto (j^*:C_{\sing}^*(U_1)\to C_{\sing}^*(U_2)
\end{eqnarray*}
and by $c_X:\mathbb Z_X\to C^*_{X,\sing}$ the inclusion map.
For $f:T\to S$ a morphism with $T,S\in\Top$, we have the morphism of complexes of presheaves 
\begin{equation*}
f^*:C^*_{S,\sing}\to f_*C^*_{T,\sing}
\end{equation*}
in $C(S)$.
\end{itemize}

\begin{thm}
\begin{itemize}
\item[(i)] If two morphisms $f,g:T\to S$ with $T,S\in\Top$ are weakly homotopic, then 
\begin{equation*}
H^n(f_*)=H^n(g_*):H_{n,\sing}(T,\mathbb Z):=H^nC^{\sing}_*(T)\to H_{n,\sing}(S,\mathbb Z):=H^nC^{\sing}_*(S).
\end{equation*}
\item[(ii)] For $S\in\Top$ there exists $CW(S)\in\CW$ together with a morphism $L_S:CW(S)\to S$
which is a weakly homotopic equivalence, that is $L_{S*}:\pi_n(CW(S))\xrightarrow{\sim}\pi_n(S)$ are isomorphisms
of abelian groups for all $n\in\mathbb N$.
\item[(ii)'] For $f:T\to S$ a morphism, with $T,S\in\Top$, and $L_S:CW(S)\to S$, $L_S:CW(T)\to T$
weakly homotopy equivalence with $CW(S),CW(T)\in\CW$ there exist a morphism
$L(f):CW(T)\to CW(S)$ unique up to homotopy such that the following diagram in $\Top$ commutes
\begin{equation*}
\xymatrix{CW(S)\ar[r]^{L_S} & S \\
CW(T)\ar[r]^{L_T}\ar[u]^{L(f)} & T\ar[u]^f}.
\end{equation*}
In particular, for $S\in\Top$, $CW(S)$ is unique up to homotopy.
\end{itemize}
\end{thm}

\begin{proof}
See \cite{Hatcher}.
\end{proof}

We have Kunneth formula for singular cohomology :

\begin{prop}\label{KunnethTop}
Let $X_1,X_2\in\Top$. Denote by $p_1:X_1\times X_2\to X_1$ and $p_2:X_1\times X_2\to X_2$ the projections. Then
\begin{equation*}
p_1^*\otimes p_2^*:C^*_{\sing}(X_1)\otimes C^*_{\sing}(X_2)\to C^*_{\sing}(X_1\times X_2)
\end{equation*}
is a quasi-isomorphism.
\end{prop}

\begin{proof}
Standard (see \cite{Hatcher} for example): 
follows from the fact that for all $p\in\mathbb N$, $H^nC_{\sing}^*(\Delta^p)=0$ for all $n\in\mathbb Z$.
\end{proof}

\begin{rem}
By definition, $X\in\Top$ is locally contractile if an only if 
the inclusion map $c_X:\mathbb Z_X\to C^*_{X,\sing}$ is an equivalence top local.
In this case it induce, by taking injective resolutions, for $n\in\mathbb Z$ isomorphisms
\begin{equation*}
H^nc_X^k: H^n(X,\mathbb Z_X)\xrightarrow{\sim}\mathbb H^n(X,C^*_{X,\sing})=H^nC^*_{\sing}(X)=:H^n_{\sing}(X,\mathbb Z).
\end{equation*}
\end{rem}

We will use the following easy propositions :

\begin{prop}\label{projformula}
\begin{itemize}
\item[(i)]Let $(S,O_S)\in\RTop$. 
Then, if $K^{\bullet}\in C_{O_S}^-(S)$ is a bounded above complex such that $K^n\in\PSh_{O_S}(S)$ are locally free for all $n\in\mathbb Z$, 
and $\phi:F^{\bullet}\to G^{\bullet}$ is a top local equivalence with $F,G\in C_{O_S}(\mathcal S)$, then
$\phi\otimes I:F^{\bullet}\otimes_{O_S} L^{\bullet}\to G^{\bullet}\otimes_{O_S} L^{\bullet}$ is an equivalence top local. 
\item[(ii)] Let $f:(T,O_T)\to (S,O_S)$ a morphism with $(T,O_T),(S,O_S)\in\RTop$. 
Then, if $K\in C^b_{O_S}(S)$ is a bounded complex such that $K^n\in\PSh_{O_S}(S)$ are locally free 
for all $n\in\mathbb Z$, and $N\in C_{O_T}(T)$
\begin{equation*}
k\circ T^{mod}(f,\otimes)(M,E(N)):K\otimes_{O_S}f_*E(N)\to f_*((f^{*mod}K)\otimes_{O_T}E(N))\to f_*E((f^{*mod}K)\otimes_{O_T}E(N)) 
\end{equation*}
is an equivalence top local.
\end{itemize}
\end{prop}

\begin{proof}
Standard.
\end{proof}

\begin{prop}\label{Tiotimes}
Let $i:(Z,O_Z)\hookrightarrow (S,O_S)$ a closed embedding of ringed spaces, with $Z,S\in\Top$.
Then for $M\in C_{O_S}(S)$ and $M\in C_{i^*O_S}(Z)$,
\begin{equation*}
T(i,\otimes)(M,N):M\otimes_{O_S} i_*N\to i_*(i^*M\otimes_{i^*O_S} N)
\end{equation*}
is an equivalence top local.
\end{prop}

\begin{proof}
Standard. Follows form the fact that $j^*i_*N=0$. 
\end{proof}

We note the following :

\begin{prop}\label{NakInj}
Let $(S,O_S)\in\Sch$ such that $O_{S,s}$ are reduced local rings for all $s\in S$. 
For $s\in S$ consider $q:L_{O_{S,s}}(k(s))\to k(s)$ the canonical projective resolution of 
the $O_{S,s}$ module $k(s):=O_{S,s}/m_s$ (the residual field) of $s\in S$.
For $s\in S$ denote by $i_s:\left\{s\right\}\hookrightarrow S$ the embedding. 
Let $\phi:F\to G$ a morphism with $F,G\in C_{O_S,c}(S)$ i.e. such that $a_{zar}H^nF,a_{zar}H^nG\in\Coh(S)$ . If 
\begin{equation*}
i_s^*\phi\otimes_{i_s^*O_S}L_{i_s^*O_S}(k(s)):
i_s^*F\otimes_{i_s^*O_S}L_{i_s^*O_S}(k(s))\to i_s^*G\otimes_{i_s^*O_S}L_{i_s^*O_S}(k(s))
\end{equation*}
is a quasi-isomorphism for all $s\in S$, then $\phi:F\to G$ is an equivalence top local.
\end{prop}

\begin{proof}
Let $s\in S$. Since tensorizing with $L_{i_s^*O_S}(k(s))$ is an exact functor, 
we have canonical isomorphism $\alpha(F)$,$\alpha(G)$ fiting in a commutative diagram
\begin{equation*}\label{NakInjdia}
\xymatrix{
H^n(i_s^*F\otimes_{i_s^*O_{S}}L_{i_s^*O_S}(k(s)))\ar[rrr]^{H^n(i_s^*\phi\otimes_{i_s^*O_{S}}L_{i_s^*O_S}(k(s)))}\ar[d]^{\alpha(F)}
& \, & \, & H^n(i_s^*G\otimes_{i_s^*O_{S}}L_{i_s^*O_S}(k(s)))\ar[d]^{\alpha(G)} \\
i_s^*(H^nF)\otimes_{i_s^*O_{S}}L_{i_s^*O_S}(k(s))\ar[rrr]^{i_s^*(H^n\phi)\otimes_{i_s^*O_{S}}L_{i_s^*O_S}(k(s))} & \, & \, & 
i_s^*(H^nG)\otimes_{i_s^*O_{S}}L_{i_s^*O_S}(k(s))}
\end{equation*}
Let $n\in\mathbb Z$. By hypothesis
\begin{equation*}
H^n(i_s^*\phi\otimes_{i_s^*O_{S}}L_{i_s^*O_S}(k(s))):
H^n(i_s^*F\otimes_{i_s^*O_{S}}L_{i_s^*O_S}(k(s)))\xrightarrow{\sim} H^n(i_s^*G\otimes_{i_s^*O_{S}}L_{i_s^*O_S}(k(s)))
\end{equation*}
is an isomorphism. Hence, the diagram \ref{NakInjdia} implies that
\begin{equation*}
i_s^*(H^n\phi)\otimes_{i_s^*O_{S}}L_{i_s^*O_S}(k(s)):
i_s^*(H^nF)\otimes_{i_s^*O_{S}}L_{i_s^*O_S}(k(s))\xrightarrow{\sim} i_s^*(H^nG)\otimes_{i_s^*O_{S}}L_{i_s^*O_S}(k(s))
\end{equation*}
is an isomorphism. We conclude on the one hand that $i_s^*H^n\phi:i_s^*H^nF\to i_s^*H^nG$ is surjective by Nakayama lemma
since $i_s^*H^nF$, $i_s^*H^nG$ are $O_{S,s}$ modules of finite type as $F,G\in C_{O_S,c}(S)$ has coherent cohomology sheaves,
and on the other hand that the rows of the following commutative diagram are isomorphism
\begin{equation*}
\xymatrix{
H^0(i_s^*(H^nF)\otimes_{i_s^*O_{S}}L_{i_s^*O_S}(k(s)))
\ar[rrr]^{H^0(i_s^*(H^n\phi)\otimes_{i_s^*O_{S}}L_{i_s^*O_S}(k(s)))\sim}\ar[d]^{=} & \, & \, & 
H^0(i_s^*(H^nG)\otimes_{i_s^*O_{S}}L_{i_s^*O_S}(k(s)))\ar[d]^{=} \\
i_s^*(H^nF)\otimes_{i_s^*O_{S}}k(s)\ar[rrr]^{i_s^*(H^n\phi)\otimes_{i_s^*O_{S}}k(s)\sim} & \, & \, & 
i_s^*(H^nG)\otimes_{i_s^*O_{S}}k(s)}.
\end{equation*}
Since 
\begin{equation*}
i_s^*(H^n\phi)\otimes_{i_s^*O_{S}}k(s):i_s^*(H^nF)\otimes_{i_s^*O_{S}}k(s)\xrightarrow{\sim}i_s^*(H^nF)\otimes_{i_s^*O_{S}}k(s)
\end{equation*}
is an isomorphism for all $s\in S$, $O_{S,s}=:i_s^*O_S$ are reduced, and $a_{zar}H^nF,a_{zar}H^nG$ are coherent, 
$i_s^*H^n\phi:i_s^*H^nF\to i_s^*H^nG$ are injective.
\end{proof}

Let $i:Z\hookrightarrow S$ a closed embedding, with $S,Z\in\Top$. 
Denote by $j:S\backslash Z\hookrightarrow S$ the open embedding of the complementary subset.
We have the adjonction
\begin{equation*}
(i_*,i^!):=(i_*,i^{\bot}):C(Z)\to C(S), \; \mbox{with \, in \, this \, case} \; i^!F:=\ker(F\to j_*j^*F). 
\end{equation*}
It induces the adjonction $(i_*,i^!):C_{(2)fil}(Z)\to C_{(2)fil}(S)$ (we recall that $i^!:=i^{\bot}$ preserve monomorphisms).

Let $i:Z\hookrightarrow S$ a closed embedding, with $S,Z\in\Top$. 
Denote by $j:S\backslash Z\hookrightarrow S$ the open embedding of the complementary subset.
We have the support section functors :
\begin{itemize}
\item We have the functor
\begin{equation*}
\Gamma_Z:C(S)\to C(S), \; F\mapsto\Gamma_Z(F):=\Cone(\ad(j^*,j_*)(F):F\to j_*j^*F)[-1], 
\end{equation*}
together with the canonical map $\gamma_Z(F):\Gamma_ZF\to F$.
We have the factorization 
\begin{equation*}
\ad(i_*,i^!)(F):i_*i^!F\xrightarrow{\ad(i_*,i^!)(F)^{\gamma}}\Gamma_ZF\xrightarrow{\gamma_Z(F)}F,
\end{equation*}
and $\ad(i_*,i^!)(F)^{\gamma}:i_*i^!F\to\Gamma_ZF$ is an homotopy equivalence.
Since $\Gamma_Z$ preserve monomorphisms, it induce the functor
\begin{equation*}
\Gamma_Z:C_{fil}(S)\to C_{fil}(S), \; (G,F)\mapsto\Gamma_Z(G,F):=(\Gamma_ZG,\Gamma_ZF), 
\end{equation*}
together with the canonical map $\gamma_Z((G,F):\Gamma_Z(G,F)\to (G,F)$.
\item We have also the functor
\begin{equation*}
\Gamma^{\vee}_Z:C(S)\to C(S), \; 
F\mapsto\Gamma^{\vee}_ZF:=\Cone(\ad(j_!,j^*)(F):j_!j^*F\to F), 
\end{equation*}
together with the canonical map $\gamma^{\vee}_Z(F):F\to\Gamma^{\vee}_ZF$.
We have the factorization 
\begin{equation*}
\ad(i^*,i_*)(F):F\xrightarrow{\gamma_Z^{\vee}(F)}\Gamma_Z^{\vee}F\xrightarrow{\ad(i^*,i_*)(F)^{\gamma}}i_*i^*F,
\end{equation*}
and $\ad(i^*,i_*)(F)^{\gamma}:\Gamma^{\vee}_ZF\to i_*i^*F$ is an homotopy equivalence.
Since $\Gamma^{\vee}_Z$ preserve monomorphisms, it induce the functor
\begin{equation*}
\Gamma_Z:C_{fil}(S)\to C_{fil}(S), \; (G,F)\mapsto\Gamma^{\vee}_Z(G,F):=(\Gamma^{\vee}_ZG,\Gamma^{\vee}_ZF), 
\end{equation*}
together with the canonical map $\gamma^{\vee}_Z(G,F):(G,F)\to\Gamma_Z^{\vee}(G,F)$.
\end{itemize}

\begin{defiprop}\label{gamma1sect2}
\begin{itemize}
\item[(i)] Let $g:S'\to S$ a morphism and $i:Z\hookrightarrow S$ a closed embedding with $S',S,Z\in\Top$. 
Then, for $(G,F)\in C_{fil}(S)$, there is a canonical map in $C_{fil}(S')$
\begin{equation*}
T(g,\gamma)(G,F):g^*\Gamma_{Z}(G,F)\to\Gamma_{Z\times_S S'}g^*(G,F)
\end{equation*}
unique up to homotopy such that $\gamma_{Z\times_S S'}(g^*(G,F))\circ T(g,\gamma)(G,F)=g^*\gamma_{Z}(G,F)$.
\item[(ii)] Let $i_1:Z_1\hookrightarrow S$, $i_2:Z_2\hookrightarrow Z_1$ be closed embeddings with $S,Z_1,Z_2\in\Top$.
Then, for $(G,F)\in C_{fil}(S)$, 
\begin{itemize}
\item there is a canonical map $T(Z_2/Z_1,\gamma)(G,F):\Gamma_{Z_2}(G,F)\to\Gamma_{Z_1}(G,F)$ in $C_{fil}(S)$
unique up to homotopy such that $\gamma_{Z_1}(G,F)\circ T(Z_2/Z_1,\gamma)(G,F)=\gamma_{Z_2}(G,F)$ 
together with a distinguish triangle in $K_{fil}(S):=K(\PSh_{fil}(S))$
\begin{equation*}
\Gamma_{Z_2}(G,F)\xrightarrow{T(Z_2/Z_1,\gamma)(G,F)}\Gamma_{Z_1}(G,F)\xrightarrow{\ad(j_2^*,j_{2*})(\Gamma_{Z_1}(G,F))}
\Gamma_{Z_1/\backslash Z_2}(G,F)\to\Gamma_{Z_2}(G,F)[1]
\end{equation*} 
\item there is a canonical map $T(Z_2/Z_1,\gamma^{\vee})(G,F):\Gamma_{Z_1}^{\vee}(G,F)\to\Gamma_{Z_2}^{\vee}(G,F)$ 
in $C_{fil}(S)$
unique up to homotopy such that $\gamma^{\vee}_{Z_2}(G,F)=T(Z_2/Z_1,\gamma^{\vee})(G,F)\circ\gamma^{\vee}_{Z_1}(G,F)$. 
together with a distinguish triangle in $K_{fil}(S)$
\begin{equation*}
\Gamma_{Z_1\backslash Z_2}^{\vee}(G,F)\xrightarrow{\ad(j_{2!},j_2^*)(G,F)}
\Gamma_{Z_1}^{\vee}(G,F)\xrightarrow{T(Z_2/Z_1,\gamma^{\vee})(G,F))}\Gamma^{\vee}_{Z_2}(G,F)\to\Gamma_{Z_2\backslash Z_1}^{\vee}(G,F)[1]
\end{equation*} 
\end{itemize}
\item[(iii)] Consider a morphism $g:(S',Z')\to(S,Z)$ with $(S',Z'),(S,Z)\in\Top^2$. We denote, for $G\in C(S)$ the composite 
\begin{equation*}
T(D,\gamma^{\vee})(G):g^*\Gamma_Z^{\vee}G\xrightarrow{\sim}\Gamma^{\vee}_{Z\times_S S'}g^*G
\xrightarrow{T(Z'/Z\times_S S',\gamma^{\vee})(G)}\Gamma_{Z'}^{\vee}g^*G
\end{equation*}
and we have then the factorization
$\gamma_{Z'}^{\vee}(g^*G):g^*G\xrightarrow{g^*\gamma_Z^{\vee}(G)} g^*\Gamma_Z^{\vee}G
\xrightarrow{T(D,\gamma^{\vee})(G)}\Gamma_{Z'}^{\vee}g^*G$.
\end{itemize}
\end{defiprop}

\begin{proof}
\noindent(i): We have the cartesian square
\begin{equation*}
\xymatrix{S\backslash Z\ar[r]^{j} & S \\
S'\backslash Z\times_S S'\ar[r]^{j'}\ar[u]^{g'} & S'\ar[u]^{g}}
\end{equation*} 
and the map is given by 
\begin{equation*}
(I,T(g,j)(j^*G)):\Cone(g^*G\to g^*j_*j^*G)\to\Cone(g^*G\to j'_*j^{'*}g^*G=j'_*g^{'*}j^*G).
\end{equation*}

\noindent(ii): Follows from the fact that $j_1^*\Gamma_{Z_2}G=0$ and $j_1^*\Gamma^{\vee}_{Z_2}G=0$, 
with $j_1:S\backslash Z_1\hookrightarrow S$ the closed embedding.

\noindent(iii): Obvious.
\end{proof}

Let $(S,O_S)\in\RTop$. Let $Z\subset S$ a closed subset.
Denote by $j:S\backslash Z\hookrightarrow S$ the open complementary embedding, 
\begin{itemize}
\item For $G\in C_{O_S}(S)$, $\Gamma_ZG:=\Cone(\ad(j^*,j_*)(G):F\to j_*j^*G)[-1]$ has a (unique) structure of $O_S$ module
such that $\gamma_Z(G):\Gamma_ZG\to G$ is a map in $C_{O_S}(S)$. This gives the functor
\begin{equation*}
\Gamma_Z:C_{O_Sfil}(S)\to C_{filO_S}(S), \; (G,F)\mapsto\Gamma_Z(G,F):=(\Gamma_ZG,\Gamma_ZF), 
\end{equation*}
together with the canonical map $\gamma_Z((G,F):\Gamma_Z(G,F)\to (G,F)$.
Let $Z_2\subset Z$ a closed subset.
Then, for $G\in C_{O_S}(S)$, $T(Z_2/Z,\gamma)(G):\Gamma_{Z_2}G\to\Gamma_ZG$ is a map in $C_{O_S}(S)$ (i.e. is $O_S$ linear).

\item For $G\in C_{O_S}(S)$, $\Gamma^{\vee}_ZG:=\Cone(\ad(j_!,j^*)(G):j_!j^*G\to G)$ has a unique structure of $O_S$ module,
such that $\gamma^{\vee}_Z(G):G\to\Gamma_Z^{\vee}G$ is a map in $C_{O_S}(S)$. This gives the functor
\begin{equation*}
\Gamma^{\vee}_Z:C_{O_Sfil}(S)\to C_{filO_S}(S), \; (G,F)\mapsto\Gamma^{\vee}_Z(G,F):=(\Gamma^{\vee}_ZG,\Gamma^{\vee}_ZF), 
\end{equation*}
together with the canonical map $\gamma^{\vee}_Z((G,F):(G,F)\to\Gamma^{\vee}_Z(G,F)$.
Let $Z_2\subset Z$ a closed subset.
Then, for $G\in C_{O_S}(S)$, $T(Z_2/Z,\gamma^{\vee})(G):\Gamma_Z^{\vee}G\to\Gamma_{Z_2}^{\vee}G$ 
is a map in $C_{O_S}(S)$ (i.e. is $O_S$ linear).

\item For $G\in C_{O_S}(S)$, we will use
\begin{eqnarray*}
\Gamma_Z^{\vee,h} G:&=&\mathbb D_S^OL_O\Gamma_ZE(\mathbb D^O_SG) \\
:&=&\Cone(\mathbb D^O_SL_O\ad(j_*,j^*)(E(\mathbb D^O_SG)):\mathbb D^O_SL_Oj_*j^*E(\mathbb D^O_SG)\to\mathbb D^O_SL_OE(\mathbb D^O_SG))
\end{eqnarray*}
and we have the canonical map $\gamma_Z^{\vee,h}(G):M\to\Gamma_Z^{\vee,h}G$ of $O_S$ module. The factorization  
\begin{eqnarray*}
\ad(j_!,j^*)(L_OM):j_!j^*L_OG\xrightarrow{(k\circ\mathbb D^OI(j_!,j^*)(\mathbb D^Oj^*L_OG)\circ d(j_!j^*L_OG))^q} \\
\mathbb D^O_SL_Oj_*j^*E(\mathbb D^O_SL_OG)\xrightarrow{\ad(j_*,j^*)(E(\mathbb D^O_SL_OG))}
\mathbb D^O_SL_OE(\mathbb D^O_SL_OG)
\end{eqnarray*}
gives the factorization 
$\gamma_Z^{\vee,h}(L_OG):L_OG\xrightarrow{\gamma_Z^{\vee}(L_OG)}\Gamma_Z^{\vee}L_OG
\xrightarrow{(k\circ\mathbb D^OI(j_!,j^*)(\mathbb D^Oj^*L_OG)\circ d(j_!j^*L_OG))^q}\Gamma_Z^{\vee,h}L_OG$.
We get the functor
\begin{eqnarray*}
\Gamma_Z^{\vee,h}:C_{O_Sfil}(S)\to C_{O_Sfil}(S), \;  
(G,F)\mapsto\Gamma_Z^{\vee,h}(G,F):=\mathbb D_S^OL_O\Gamma_ZE(\mathbb D^O_S(G,F)),
\end{eqnarray*}
together with the factorization 
\begin{eqnarray*}
\gamma_Z^{\vee,h}(L_O(G,F)):L_O(G,F)\xrightarrow{\gamma_Z^{\vee}(L_O(G,F))}\Gamma_Z^{\vee}L_O(G,F) \\
\xrightarrow{(k\circ\mathbb D^OI(j_!,j^*)(\mathbb D_S^Oj^*L_O(G,F))\circ d(j_!j^*L_O(G,F)))^q}\Gamma_Z^{\vee,h}L_O(G,F),
\end{eqnarray*}

\item Consider $\mathcal I\subset O_S$ a right ideal of $O_S$ such that $\mathcal I^o_Z\subset\mathcal I$, 
where $\mathcal I^o_Z\subset O_S$ is the left and right ideal consisting of section which vanish on $Z$.
\begin{itemize}
\item For $G\in\PSh_{O_S}(S)$, we consider, $S^o\subset S$ being an open subset, 
\begin{equation*}
\mathcal IG(S^o)=<\left\{f.m, m\in G(S^o),f\in\mathcal I(S^o)\right\}>\subset G(S^o) 
\end{equation*}
since $\mathcal I$ is a right ideal,
and we denote by $b_I(G):\mathcal IG\to G$ the injective morphism of $O_S$ modules and by $c_Z(G):G\to G/\mathcal IG$ the quotient map. 
The adjonction map $\ad(j_!,j^*)(G):j_!j^*G\to G$ factors trough $b_I(G)$ : 
\begin{equation*}
ad(j_!,j^*)(G):j_!j^*G\xrightarrow{b^I_{Z/S}(G)}\mathcal IG\xrightarrow{b_I(G)}G
\end{equation*}
We have then the support section functor,
\begin{eqnarray*}
\Gamma^{\vee,O,I}_Z:C_{O_S}(S)\to C_{O_S}(S), \; G\mapsto\Gamma_Z^{\vee,O,I}G:=\Cone(b_I(G):\mathcal IG\to G)
\end{eqnarray*}
together with the canonical map $\gamma_Z^{\vee,O}(G):G\to\Gamma_Z^{\vee,O}G$ which factors through
\begin{equation*}
\gamma_Z^{\vee,O,I}(G):G\xrightarrow{\gamma_Z^{\vee}(G)}\Gamma^{\vee}_ZG\xrightarrow{b^I_{S/Z}(G)}\Gamma^{\vee,O}_ZG.
\end{equation*}
By the exact sequence $0\to\mathcal IG\xrightarrow{b_I(G)}G\xrightarrow{c_I(G)}G/\mathcal IG\to 0$, we have an homotopy equivalence 
$c_I(G):\Gamma_Z^{\vee,O,I}G\to G/\mathcal IG$.
\item For $G\in\PSh_{O_S}(S)$, we consider 
\begin{equation*}
b'_I(G):G\to G\otimes_{O_S}\mathbb D^O_S(\mathcal I):=G\otimes_{O_S}\mathcal Hom(\mathcal I,O_S)
\end{equation*} 
The adjonction map $\ad(j^*,j_*)(G):G\to j_*j^*G$ factors trough $b'_I(G)$ : 
\begin{equation*}
ad(j^*,j_*)(G):G\xrightarrow{b'_I(G)}G\otimes_{O_S}\mathbb D^O_S(\mathcal I)\xrightarrow{b^{'I}_{Z/S}(G)}j_*j^*G
\end{equation*}
We have then the support section functor,
\begin{eqnarray*}
\Gamma^{O,I}_Z:C_{O_S}(S)\to C_{O_S}(S), \; G\mapsto\Gamma_Z^{O,I}G:=\Cone(b'_I(G):G\to G\otimes_{O_S}\mathbb D^O_S(\mathcal I))[-1]
\end{eqnarray*}
together with the canonical map $\gamma_Z^{O}(G):\Gamma_Z^OG\to G$ which factors through
\begin{equation*}
\gamma_Z^{O,I}(G):\Gamma_Z^OG\xrightarrow{b^I_{S/Z}(G)}\Gamma_ZG\xrightarrow{\gamma_Z(G)}G.
\end{equation*}
\item By definition, we have for a canonical isomorphism 
\begin{equation*}
I(D,\gamma^O)(G):\mathbb D^O_S\Gamma^{\vee,O,I}G\xrightarrow{\sim}\Gamma_Z^{O,I}\mathbb D^O_S G
\end{equation*}
which gives the transformation map in $C_{O_S}(S)$
\begin{eqnarray*}
T(D,\gamma^O)(G):\Gamma^{\vee,O,I}\mathbb D^O_SG\xrightarrow{d(-)}\mathbb D^{O,2}_S\Gamma^{\vee,O,I}\mathbb D^O_SG 
\xrightarrow{\mathbb D^O_SI(D,\gamma^O)(\mathbb D^O_SG)^{-1}} \\
\mathbb D^O_S\Gamma_Z^{O,I}\mathbb D^{O,2}_S G\xrightarrow{\mathbb D^O_S\Gamma_Z^{O,I}d(G)}\mathbb D^O_S\Gamma_Z^{O,I}G
\end{eqnarray*}
\end{itemize}
\end{itemize}

\begin{defiprop}\label{gamma1sect2mod}
\begin{itemize}
\item[(i)]Let $g:(S',O_{S'})\to (S,O_S)$ a morphism and $i:Z\hookrightarrow S$ a closed embedding with $(S',O_{S'},(S,O_S)\in\RTop$. 
Then, for $(G,F)\in C_{O_Sfil}(S)$, there is a canonical map in $C_{O_{S'}fil}(S')$
\begin{equation*}
T^{mod}(g,\gamma)(G,F):g^{*mod}\Gamma_{Z}(G,F)\to\Gamma_{Z\times_S S'}g^{*mod}(G,F)
\end{equation*}
unique up to homotopy, such that $\gamma_{Z\times_S S'}(g^{*mod}G)\circ T^{mod}(g,\gamma)(G)=g^{*mod}\gamma_{Z}G$.
\item[(ii)]Let $i_1:(Z_1,O_{Z_1})\hookrightarrow (S,O_S)$, $i_2:(Z_2,O_{Z_2})\hookrightarrow (Z_1,O_{Z_1})$ 
be closed embeddings with $S,Z_1,Z_2\in\Top$. Then, for $(G,F)\in C_{O_Sfil}(S)$, there is a canonical map in $C_{O_Sfil}(S)$
\begin{equation*}
T(Z_2/Z_1,\gamma^{\vee,O})(G,F):\Gamma_{Z_1}^{\vee,O}(G,F)\to\Gamma_{Z_2}^{\vee,O}(G,F) 
\end{equation*}
unique up to homotopy such that $\gamma^{\vee,O}_{Z_2}(G,F)=T(Z_2/Z_1,\gamma^{\vee,O})(G,F)\circ\gamma^{\vee,O}_{Z_1}(G,F)$. 
\item[(iii)] Consider a morphism $g:((S',O_{S'}),Z')\to((S,O_S),Z)$ with $((S',O_{S'}),Z')\to((S,O_S),Z)\in\RTop^2$.
We denote, for $M\in C_{O_S}(S)$ the composite 
\begin{equation*}
T^{mod}(D,\gamma^{\vee,O})(G):g^{*mod}\Gamma_Z^{\vee,O}G\xrightarrow{\sim}\Gamma^{\vee,O}_{Z\times_S S'}g^{*mod}G
\xrightarrow{T(Z'/Z\times_S S',\gamma^{\vee,O})(G)}\Gamma_{Z'}^{\vee,O}g^{*mod}G
\end{equation*}
and we have then the factorization
\begin{equation*}
\gamma_{Z'}^{\vee,O}(g^{*mod}M):g^{*mod}G\xrightarrow{g^{*mod}\gamma_Z^{\vee,O}(G)} g^{*mod}\Gamma_Z^{\vee,O}G
\xrightarrow{T^{mod}(D,\gamma^{\vee,O})(G)}\Gamma_{Z'}^{\vee,O}g^{*mod}G
\end{equation*}
\end{itemize}
\end{defiprop}

\begin{proof}
\noindent(i): We have the cartesian square
\begin{equation*}
\xymatrix{S\backslash Z\ar[r]^{j} & S \\
S'\backslash Z\times_S S'\ar[r]^{j'}\ar[u]^{g'} & S'\ar[u]^{g}}
\end{equation*} 
and the map is given by 
\begin{equation*}
(I,T^{mod}(g,j)(j^*G)):\Cone(g^{*mod}G\to g^{*mod}j_*j^*G)\to\Cone(g^{*mod}G\to j'_*j^{'*}g^{*mod}G=j'_*g^{'*mod}j^*G).
\end{equation*}

\noindent(ii):Obvious.

\noindent(iii):Obvious.
\end{proof}

\begin{defiprop}\label{TDwgamma}
Consider a commutative diagram in $\RTop$
\begin{equation*}
D_0=\xymatrix{f: (X,O_X)\ar[r]^{i} & (Y,O_Y)\ar[r]^{p} & (S,O_S) \\
f':(X',O_{X'}\ar[r]^{i'}\ar[u]^{g'} & (Y',O_{Y'})\ar[u]^{g''}\ar[r]^{p'} & (T,O_T)\ar[u]^{g} }
\end{equation*}
with $i$, $i'$ being closed embeddings. Denote by $D$ the right square of $D_0$. 
The closed embedding $i':X'\hookrightarrow Y'$ factors through
$i':X'\xrightarrow{i'_1} X\times_Y Y'\xrightarrow{i'_0} Y'$ where $i'_1,i'_0$ are closed embeddings. 
\begin{itemize}
\item[(i)] We have the canonical map,
\begin{eqnarray*}
E(\Omega_{((O_{Y'}/g^{''*}O_Y))/(O_T/g^*O_S)})\circ T(g'',E)(-)\circ T(g'',\gamma)(-): \\
g^{''*}\Gamma_{X}E(\Omega^{\bullet}_{O_Y/p^*O_S},F_b)\to\Gamma_{X\times_Y Y'}E(\Omega^{\bullet}_{O_{Y'}/p^{'*}O_T},F_b)
\end{eqnarray*}
unique up to homotopy such that the following diagram in $C_{g^{''*}p^*O_Sfil}(Y')=C_{p^{'*}g^*O_Sfil}(Y')$ commutes
\begin{equation*}
\xymatrix{g^{''*}\Gamma_{X}E(\Omega^{\bullet}_{O_Y/p^*O_S},F_b)
\ar[rrrr]^{E(\Omega_{((O_{Y'}/Y))/(O_T/S)})\circ T(g'',E)(-)\circ T(g'',\gamma)(-)}\ar[d]_{\gamma_X(-)} & \, & \, & \, &
\Gamma_{X\times_Y Y'}E(\Omega^{\bullet}_{O_{Y'}/p^{'*}O_T},F_b)\ar[d]^{\gamma_{X\times_Y Y'}(-)} \\
g^{''*}E(\Omega^{\bullet}_{O_Y/p^*O_S},F_b)\ar[rrrr]^{E(\Omega_{((O_{Y'}/g^{''*}O_Y)/(O_T/S))}\circ T(g'',E)(-)} 
& \, & \, & \, &  E(\Omega^{\bullet}_{O_Y'/p^{'*}O_T},F_b)}.
\end{equation*}
\item[(ii)] There is a canonical map,
\begin{equation*}
T^O_{\omega}(D)^{\gamma}:g^{*mod}L_Op_*\Gamma_{X}E(\Omega^{\bullet}_{O_Y/p^*O_S},F_b)\to
p'_*\Gamma_{X\times_Y Y'}E(\Omega^{\bullet}_{O_{Y'}/p^{'*}O_T},F_b)
\end{equation*}
unique up to homotopy such that the following diagram in $C_{O_{T}fil}(T)$ commutes
\begin{equation*}
\xymatrix{g^{*mod}L_Op_*\Gamma_{X}E(\Omega^{\bullet}_{O_Y/p^*O_S},F_b)
\ar[rrr]^{T_{\omega}^O(D)^{\gamma}}\ar[d]_{\gamma_X(-)} & \, & \, & 
p'_*\Gamma_{X\times_Y Y'}E(\Omega^{\bullet}_{O_{Y'}/p^{'*}O_T},F_b)\ar[d]^{\gamma_{X\times_Y Y'}(-)} \\
g^{*mod}L_Op_*E(\Omega^{\bullet}_{O_Y/p^*O_S},F_b)\ar[rrr]^{T^O_{\omega}(D)} & \, & \, &  
p'_*E(\Omega^{\bullet}_{O_{Y'}/p^{'*}O_T},F_b)}.
\end{equation*}
\item[(iii)] We have the canonical map in $C_{f^{'*}O_T}(Y')$
\begin{equation*}
T(X'/X\times_Y Y',\gamma)(E(\Omega^{\bullet}_{O_{Y'}/p^{'*}O_T},F_b)):
\Gamma_{X'}E(\Omega^{\bullet}_{O_{Y'}/p^{'*}O_T},F_b)\to\Gamma_{X\times_Y Y'}E(\Omega^{\bullet}_{O_{Y'}/p^{'*}O_T},F_b)
\end{equation*}
unique up to homotopy such that $\gamma_{X\times_Y Y'}(-)\circ T(X'/X\times_Y Y',\gamma)(-)=\gamma_{X'}(-)$.
\end{itemize}
\end{defiprop}

\begin{proof}
Immediate from definition. We take for the map of point (ii) the composite
\begin{eqnarray*}
T^O_{\omega}(D)^{\gamma}: 
g^{*mod}L_Op_*\Gamma_{X}E(\Omega^{\bullet}_{O_Y/p^*O_S},F_b)
\xrightarrow{q} 
g^*p_*\Gamma_{X}E(\Omega^{\bullet}_{O_Y/p^*O_S},F_b)\otimes_{g^*O_S}O_T \\ 
\xrightarrow{T(g'',E)(-)\circ T(g'',\gamma)(-)\circ T(D)(E(\Omega^{\bullet}_{O_X/p^*O_S}))}   
(p'_*\Gamma_{X\times_Y Y'}E(g^{''}\Omega^{\bullet}_{O_Y/p^*O_S},F_b))\otimes_{g^*O_S}O_T \\ 
\xrightarrow{E(\Omega_{(O_{Y'}/g^{''*}O_Y)/(O_T/g^*O_S}))} 
p'_*\Gamma_{X\times_Y Y'}E(\Omega^{\bullet}_{O_{Y'}/p^{'*}O_T},F_b)\otimes_{g^*O_S}O_T 
\xrightarrow{m}
p'_*\Gamma_{X\times_Y Y'}E(\Omega^{\bullet}_{O_{Y'}/p^{'*}O_T},F_b),
\end{eqnarray*}
with $m(n\otimes s)=s.n$.
\end{proof}

\begin{defi}\label{ZS}
\begin{itemize}
\item[(i)] Let $S\in\Top$. For $Z\subset S$ a closed subset, we denote by $C_{Z}(S)\subset C(S)$
the full subcategory consisting of complexes of presheaves $F\in C(S)$ such that $a_{top}H^n(j^*F)=0$ for all $n\in\mathbb Z$,
where $j:S\backslash Z\hookrightarrow S$ is the complementary open embedding and $a_{top}$ is the sheaftification functor.
\item[(i)'] More generally, let $(S,O_S)\in\RTop$. For $Z\subset S$ a closed subset, we denote by 
\begin{equation*}
C_{O_S,Z}(S)\subset C_{O_S}(S) \; ,\;  \mathcal QCoh_Z(S)\subset\mathcal QCoh(S) 
\end{equation*}
the full subcategories consisting of complexes of presheaves $G\in C_{O_S}(S)$ such that $a_{top}H^n(j^*F)=0$ for all $n\in\mathbb Z$,
resp. quasi-coherent sheaves $G\in\mathcal QCoh(S)$ such that $j^*F=0$. 
\item[(ii)] Let $S\in\Top$. For $Z\subset S$ a closed subset, we denote by $C_{fil,Z}(S)\subset C_{fil}(S)$ 
the full subcategory consisting of filtered complexes of presheaves $(G,F)\in C_{fil}(S)$ 
such that there exist $r\in\mathbb N$ 
and an $r$-filtered homotopy equivalence $\phi:(G,F)\to(G',F)$ with $(G',F)\in C_{fil}(S)$
such that $a_{top}j^*H^n\Gr_F^p(G',F)=0$ for all $n,p\in\mathbb Z$,
where $j:S\backslash Z\hookrightarrow S$ is the complementary open embedding and $a_{top}$ is the sheaftification functor.
Note that this definition say that this $r$ does NOT depend on $n$ and $p$.
\item[(ii)'] More generally, let $(S,O_S)\in\RTop$. For $Z\subset S$ a closed subset, we denote by 
\begin{equation*}
C_{O_S,fil,Z}(S)\subset C_{O_S,fil}(S), \mathcal QCoh_{fil,Z}(S)\subset\mathcal QCoh(S) 
\end{equation*}
the full subcategories consisting of filtered complexes of presheaves $(G,F)\in C_{O_Sfil}(S)$ 
such that there exist $r\in\mathbb N$ 
and an $r$-filtered homotopy equivalence $\phi:(G,F)\to(G',F)$ with $(G',F)\in C_{fil}(S)$
such that $a_{top}j^*H^n\Gr_F^p(G',F)=0$ for all $n,p\in\mathbb Z$,
resp. filtered quasi-coherent sheaves $(G,F)\in\mathcal QCoh(S)$ 
and an $r$-filtered homotopy equivalence $\phi:(G,F)\to(G',F)$ with $(G',F)\in C_{fil}(S)$
such that there exist $r\in\mathbb N$ such that $j^*H^n\Gr_F^p(G,F)=0$ for all $n,p\in\mathbb Z$.
Note that this definition say that this $r$ does NOT depend on $n$ and $p$.
\end{itemize}
\end{defi}

Let $(S,O_S)\in\RTop$ and $Z\subset S$ a closed subset. 
\begin{itemize}
\item For $(G,F)\in C_{fil}(S)$, we have 
$\Gamma_Z(G,F),\Gamma_Z^{\vee}(G,F)\in C_{fil,Z}(S)$. 
\item For $(G,F)\in C_{O_Sfil}(S)$, we have 
$\Gamma_Z(G,F),\Gamma_Z^{\vee}(G,F),\Gamma_Z^{\vee,h}(G,F),\Gamma_Z^{\vee,O}(G,F)\in C_{O_Sfil,Z}(S)$.
\end{itemize}

\begin{prop}\label{KS}
Let $S\in\Top$ and $Z\subset S$ a closed subspace. Denote by $i:Z\hookrightarrow S$ the closed embedding.
\begin{itemize}
\item[(i)] The functor $i^*:\Shv_{Z}(S)\to\Shv(Z)$ is an equivalence of category whose inverse is
$i_*:\Shv(Z)\to\Shv_Z(S)$. 
More precisely $\ad(i_*,i^*)(H):i^*i_*H\to H$ is an isomorphism if $H\in\Shv(Z)$ and
$\ad(i_*,i^*)(G):G\to i_*i^*G$ is an isomorphism if $G\in\Shv_Z(S)$. 
\item[(ii)]: The functor $i^*:\Shv_{fil,Z}(S)\to\Shv_{fil}(Z)$ is an equivalence of category whose inverse is
$i_*:\Shv_{fil}(Z)\to\Shv_{fil,Z}(S)$. 
More precisely $\ad(i_*,i^*)(H,F):i^*i_*(H,F)\to (H,F)$ is an isomorphism if $(H,F)\in\Shv(Z)$ and
$\ad(i_*,i^*)(G,F):(G,F)\to i_*i^*(G,F)$ is an isomorphism if $(G,F)\in\Shv_Z(S)$.
\item[(iii)]: The functor $i^*:D_{\tau,fil,Z}(S)\to D_{\tau,fil}(Z)$ is an equivalence of category whose inverse is
$i_*:D_{\tau,fil}(Z)\to D_{\tau,fil,Z}(S)$. 
More precisely $\ad(i_*,i^*)(H,F):i^*i_*(H,F)\to (H,F)$ is an equivalence top local if $(H,F)\in C_{fil}(Z)$ 
and $\ad(i_*,i^*)(G,F):(G,F)\to i_*i^*(G,F)$ is an equivalence top local if $(G,F)\in C_{fil,Z}(S)$.
\end{itemize}
\end{prop}

\begin{proof}
\noindent(i):Standard.

\noindent(ii): Follows from (i).

\noindent(iii): Follows from (ii).
\end{proof}

Let $S\in\Top$ and $Z\subset S$ a closed subspace.
By proposition \ref{KS}, if $G\in C(S)$, $\ad(i_*,i^*)(\Gamma_ZG):\Gamma_ZG\to i_*i^*\Gamma_ZG$ 
is an equivalence top local since $\Gamma_ZG\in C_Z(S)$.

Let $(S,O_S)\in\RTop$. 
Let $D=\cup_i D_i\subset X$ a normal crossing divisor, denote by $j:S\backslash D\hookrightarrow S$ the open embedding, 
and consider $\mathcal I_D\subset O_S$ the ideal of vanishing function on $D$ which is invertible. We set, for $M\in C_{O_S}(S)$,
\begin{equation*}
M(*D):=\lim_n\mathcal Hom_{O_S}(\mathcal I^n,M),
\end{equation*}
and we denote by $a_D(F):F\to F(*D)$ the surjective morphism of presheaves. 
The adjonction map $\ad(j^*,j_*)(F):F\to j_*j^*F$ factors trough $a_D(F)$ : 
\begin{equation*}
ad(j^*,j_*)(F):F\xrightarrow{a_D(F)}F(*D)\xrightarrow{a_{S/D}(F)} j_*j^*F
\end{equation*}

\begin{rem}
\begin{itemize}
\item Let $j:U\hookrightarrow X$ an open embedding, with $(X,O_X)\in\RTop$. 
Then if $F\in\Coh_{O_U}(U)$ is a coherent sheaf of $O_U$ module, $j_*F$ is quasi-coherent but NOT coherent in general. 
In particular for $F\in C_{O_U}(U)$ whose cohomology sheaves $a_{tau}H^nF$ are coherent for all $n\in\mathbb Z$, 
the cohomology sheaves $R^nj_*F:=a_{\tau}H^nj_*E(F)$ of $Rj_*F=j_*E(F)$ are quasi-coherent but NOT coherent.   
\item Let $j:U\hookrightarrow X$ an open embedding, with $X\in\Sch$. Then if $F\in\Coh(U)$ is a coherent sheaf of $O_U$ module,
$j_*F$ is quasi-coherent but NOT coherent. However, there exist an $O_X$ submodule $\tilde F\subset j_*F$ such that $j^*\tilde F=F$ and
$\tilde F\in\Coh(X)$.
\end{itemize}
\end{rem}

The following propositions are true for schemes but NOT for arbitrary ringed spaces like analytic spaces :

\begin{prop}\label{star}
\begin{itemize}
\item[(i)]Let $X=(X,O_X)\in\Sch$ a noetherien scheme and $D\subset X$ a closed subset. 
Denote by $j:U=X\backslash D\hookrightarrow X$ an open embedding.
Then for $F\in\mathcal QCoh_{O_U}(U)$ a quasi coherent sheaf, 
$j_*F\in\mathcal QCoh_{O_X}(X)$ is quasi-coherent and is the direct limit of its coherent subsheaves.
\item[(ii)]Let $X=(X,O_X)$ a noetherien scheme and $D=\cup D_i\subset X$ a normal crossing divisor. 
Denote by $j:U=X\backslash D\hookrightarrow X$ an open embedding.
Then for $F\in\mathcal QCoh_{O_U}(U)$ a quasi coherent sheaf, the canonical map 
$a_{X/D}(F):F(*D)\xrightarrow{\sim}j_*F$ is an isomorphism. 
\end{itemize}
\end{prop}

\begin{proof}
Standard.
\end{proof}

\begin{prop}\label{KSO}
Let $S=(S,O_S)\in\Sch$ and $Z\subset S$ a closed subscheme. Denote by $i:Z\hookrightarrow S$ the closed embedding.
\begin{itemize}
\item[(i)] For $G\in\mathcal QCoh_{Z}(S)$, $i^*G$ has a canonical structure of $O_Z$ module.
Moreover, the functor $i^*:\mathcal QCoh_{Z}(S)\to\mathcal QCoh(Z)$ is an equivalence of category whose inverse is
$i_*:\mathcal QCoh(Z)\to\mathcal QCoh_Z(S)$. 
\item[(ii)]: The functor $i^*:\mathcal QCoh_{fil,Z}(S)\to\mathcal QCoh_{fil}(Z)$ is an equivalence of category whose inverse is
$i_*:\mathcal QCoh_{fil}(Z)\to\mathcal QCoh_{fil,Z}(S)$.
\item[(iii)]: The functor $i^*:D_{O_Sfil,Z,qc}(S)\to D_{O_Zfil,qc}(Z)$ is an equivalence of category whose inverse is
$i_*:D_{O_Zfil,qc}(Z)\to D_{O_Sfil,Z,qc}(S)$.
\end{itemize}
\end{prop}

\begin{proof}
\noindent(i):Standard.

\noindent(ii): Follows from (i).

\noindent(iii): Follows from (ii) since $i^*$ and $i_*$ are exact functors.
\end{proof}

\begin{defi}
Let $(S,O_S)\in\RTop$ a locally ringed space with $O_S$ commutative. Consider an $\kappa_S\in C_{O_S}(S)$. 
Let $\mathcal I\subset O_S$ an ideal subsheaf and $Z=V(\mathcal I)\subset S$ the associated closed subset.
For $G\in\PSh_{O_S}(S)$, we denote by $\hat G_Z:=\hat G_I:=\lim_kG/\mathcal I^kG$ the completion with respect to the ideal $\mathcal I$
and by $c^{\infty}_Z(G):G\to\hat G_Z$ the quotient map. Then, the canonical map 
\begin{eqnarray*}
d_{\kappa_S,Z}(G):G\xrightarrow{d(G)}\mathbb D^{O,2}_SG\xrightarrow{T^{mod}(\otimes\kappa_S,hom)(\mathbb D^O_SG,O_S)} \\
\mathcal Hom_{O_S}(\mathbb D^O_SG\otimes_{O_S}\kappa_S,\kappa_S)\xrightarrow{T^{mod}(\Gamma_ZE,hom)(-,-)}
\mathcal Hom_{O_S}(\Gamma_ZE(\mathbb D^O_SG\otimes_{O_S}\kappa_S),\Gamma_ZE(\kappa_S)) 
\end{eqnarray*}
factors through
\begin{equation*}
d_{\kappa_S,Z}(G):G\xrightarrow{c^{\infty}_Z(G)}\hat G_Z\xrightarrow{d_{\kappa_S,Z}(G)}
\mathcal Hom_{O_S}(\Gamma_ZE(\mathbb D^O_SG\otimes_{O_S}\kappa_S),\Gamma_ZE(\kappa_S)) 
\end{equation*}
Clearly if $G\in C_{O_S}(S)$ then $d_{\kappa_S,Z}(G)$ is a map in $C_{O_S}(S)$.
On the other hand, we have a commutative diagram
\begin{equation*}
\xymatrix{
\Omega^p_{O_S}\ar[rrr]^{d_{\kappa_S,Z}(\Omega_{O_S}^p)}\ar[d]^{d} & \, & \, &
\mathcal Hom_{O_S}(\Gamma_ZE(\mathbb D^O_S\Omega_{O_S}^p\otimes_{O_S}\kappa_S),\Gamma_ZE(\kappa_S))\ar[d]^{d^{\kappa_S,\gamma}} \\ 
\Omega^{p+1}_S\ar[rrr]^{d_{\kappa_S,Z}(\Omega_{O_S}^{p+1})} & \, & \, & 
\mathcal Hom_{O_S}(\Gamma_ZE(\mathbb D^O_S\Omega_{O_S}^{p+1}\otimes_{O_S}\kappa_S),\Gamma_ZE(\kappa_S))}
\end{equation*}
so that $d_{\kappa_S,Z}(\Omega_S^{\bullet})\in C(S)$.
\end{defi}

The following theorem is the from \cite{Harshorne}

\begin{thm}
Let $S\in\Var(\mathbb C)$. Let $Z=V(\mathcal I)\subset S$ a closed subset. Denote by $K_S\in\PSh_{O_S}(S)$ the canonical bundle.
Then, for $G\in C_{O_S,c}(S)$,
\begin{equation*}
d_{K_S,Z}(G):\hat G_Z\to\mathcal Hom_{O_S}(\Gamma_ZE(\mathbb D^O_SG\otimes_{O_S} K_S),\Gamma_ZE(K_S)) 
\end{equation*}
is an equivalence Zariski local.
\end{thm}

Let $f:(X,O_X)\to(S,O_S)$ a morphism with $(S,O_S)\in\RTop$. In the particular case where $O_S$ is a commutative sheaf of ring, 
$T_{O_S}\in\PSh_{O_S}(S)$ and $\Omega_{O_S}=\mathbb D_{O_S}T_{O_S}\in\PSh_{O_S}(S)$ are sheaves and the morphism in $\PSh(X)$
\begin{equation*}
T(f,\hom)(O_S,O_S):f^*\mathcal Hom(O_S,O_S)\to\mathcal Hom(f^*O_S,f^*O_S)
\end{equation*}
induces isomorphisms
$T(f,\hom)(O_S,O_S):f^*T_{O_S}\xrightarrow{\sim}T_{f^*O_S}$ and 
$\mathbb D_{f^*O_S}T(f,\hom)(O_S,O_S):\Omega_{f^*O_S}\to f^*\Omega_{O_S}$
where for $F\in\Shv(S)$, we denote again (as usual) by abuse $f^*F:=a_{\tau}f^*F\in\Shv(S)$,
$a_{tau}:\PSh(S)\to\Shv(S)$ being the sheaftification functor.

\begin{defi}
\begin{itemize}
\item[(i)]Let $(X,O_X)\in\RTop$. 
A foliation $(X,O_X)/\mathcal F$ on $(X,O_X)$ is an $O_X$ module $\Omega_{O_X/\mathcal F}\in\PSh_{O_X}(X)$ 
together with a derivation map $d:=d_{\mathcal F}:O_X\to\Omega_{O_X/\mathcal F}$ such that 
\begin{itemize}
\item the associated map $q:=q_{\mathcal F}:=\omega_{X}(d):\Omega_{O_X}\to\Omega_{O_X/\mathcal F}$ is surjective  
\item satisfy the integrability condition $d(\ker q)\subset\ker q$ which implies that the map 
$d:\Omega_{O_X}^p\to\Omega_{O_X}^{p+1}$ induce factors trough 
\begin{equation*}
\xymatrix{\Omega^p_{O_X}\ar[r]^{d}\ar[d]_{q^p:=\wedge^pq} & \Omega^{p+1}_{O_X}\ar[d]^{q^{p+1}:=\wedge^pq} \\
\Omega^p_{O_X/\mathcal F}\ar[r]^{d} & \Omega^{p+1}_{O_X/\mathcal F}}
\end{equation*}
and $d:\Omega^p_{O_X/\mathcal F}\to\Omega^{p+1}_{O_X/\mathcal F}$ is neccessary unique by the surjectivity of 
$q^p:\Omega_{O_X}^p\to\Omega^p_{O_X/\mathcal F}$. 
\end{itemize}
In the particular case where $\Omega_{O_X/\mathcal F}\in\PSh_{O_X}(X)$ is a locally free sheaf of $O_X$ module,
$\mathbb D_{O_X}q:T_{O_X/\mathcal F}:=\mathbb D_{O_X}\Omega_{O_X/\mathcal F}\to T_{O_X}$ is injective and the second condition is then equivalent
to the fact that the sub $O_X$ module $T_{O_X/\mathcal F}\subset T_{O_X}$ is a Lie subalgebra, that is
$[T_{O_X/\mathcal F},T_{O_X/\mathcal F}]\subset T_{O_X/\mathcal F}$.
\item[(ii)] A piece of leaf a foliation $(X,O_X)/\mathcal F$ with $(X,O_X)\in\RTop$ such that $O_X$ is a commutative sheaf of ring 
is an injective morphism of ringed spaces $l:(Z,O_Z)\hookrightarrow(X,O_X)$ such that
$\Omega_{i^*O_X/O_Z}::\Omega_{i^*O_X}\to\Omega_{O_Z}$ factors trough an isomorphism 
\begin{eqnarray*}
\Omega_{i^*O_X/O_Z}:\Omega_{i^*O_X}\xrightarrow{\mathbb D_{i^*O_X}T(i,hom)(O_X,O_X)}i^*\Omega_{O_X}
\xrightarrow{i^*q}i^*\Omega_{O_X/\mathcal F}\to\Omega_{O_Z}.
\end{eqnarray*}
\item[(iii)]If $f:(X,O_X)\to(S,O_S)$ is a morphism with $(X,O_X),(S,O_S)\in\RTop$, we have the foliation $(X,O_X)/(S,O_S):=((X,O_X),f)$ on 
$(X,O_X)$ given by the surjection 
\begin{equation*}
q:\Omega_{O_X}\to\Omega_{O_X/f^*O_S}:=\coker(\Omega_{O_X/f^*O_S}:\Omega_{f^*O_S}\to\Omega_{O_X}). 
\end{equation*}
The fibers $i_{X_s}:(X_s,O_{X_s})\hookrightarrow (X,O_X)$ for each $s\in S$ are the leaves of the foliation.
\item[(iv)] We have the category $FolRTop$
\begin{itemize}
\item whose objects are foliated ringed spaces $(X,O_X)/\mathcal F$ with $O_X$ a commutatif sheaf of ring and 
\item whose morphisms $f:(X,O_X)/\mathcal F\to (S,O_S)/\mathcal G$ are morphisms of ringed spaces $f:(X,O_X)\to (S,O_S)$ such that
$\Omega_{O_X/f^*O_S}:\Omega_{f^*O_S}\to\Omega_{O_X}$ factors through 
\begin{equation*}
\xymatrix{ 
f^*\Omega_{O_S}\ar[d]_{f^*q_{\mathcal G}}\ar[rr]^{\mathbb D_{f^*O_S}T(f,\hom)(O_S,O_S)^{-1}} & \, &
 \Omega_{f^*O_S}\ar[rr]^{\Omega_{O_X/f^*O_S}} & \, & \Omega_{O_X}\ar[d]^{q_{\mathcal F}} \\
f^*\Omega_{O_S/\mathcal G}\ar[rrrr]^{\Omega_{O_X/f^*O_S}^q} & \, & \, & \, & \Omega_{O_X/\mathcal F}}.
\end{equation*}
\end{itemize}
This category admits inverse limits with
$(X,O_X)/\mathcal F\times(Y,O_Y)/\mathcal G=(X\times Y,p_X^*O_X\otimes p_Y^*O_Y)/p_X^*\mathcal F\otimes p_Y^*\mathcal G$ and 
\begin{eqnarray*}
(X,O_X)/\mathcal F\times_{(S,O_S)/\mathcal H}(Y,O_Y)/\mathcal G=
(X\times_S Y,\delta_S^*(p_X^*O_X\otimes p_Y^*O_Y))/p_X^*\mathcal F\otimes p_Y^*\mathcal G
\end{eqnarray*}
with $\delta_S:X\times_SY\hookrightarrow X\times Y$ the embedding given by the diagonal $\delta_S:S\hookrightarrow S\times S$.
\end{itemize}
\end{defi}

Let $S\in\Top$. Let $S=\cup_{i=1}^l S_i$ an open cover and denote by $S_I=\cap_{i\in I} S_i$.
Let $i_i:S_i\hookrightarrow\tilde S_i$ closed embeddings, with $\tilde S_i\in\Top$. 
For $I\subset\left[1,\cdots l\right]$, denote by $\tilde S_I=\Pi_{i\in I}\tilde S_i$.
We then have closed embeddings $i_I:S_I\hookrightarrow\tilde S_I$ and, for $J\subset I$, the following commutative diagram
\begin{equation*}
D_{IJ}=\xymatrix{ S_I\ar[r]^{i_I} & \tilde S_I \\
S_J\ar[u]^{j_{IJ}}\ar[r]^{i_J} & \tilde S_J\ar[u]^{p_{IJ}}}  
\end{equation*}
where $j_{IJ}:S_J\hookrightarrow S_I$ is the open embedding so that $j_I\circ j_{IJ}=j_J$.
and $p_{IJ}:\tilde S_J\to\tilde S_I$ the projection. 
This gives the diagram of topological spaces $(\tilde S_I)\in\Fun(\mathcal P(\mathbb N),\Top)$ which
which gives the diagram $(\tilde S_I)\in(\Ouv(\tilde S_I))\in\Fun(\mathcal P(\mathbb N),\Cat)$
Denote $m:\tilde S_I\backslash(S_I\backslash S_J)\hookrightarrow\tilde S_I$ the open embedding.

\begin{defi}
Let $S\in\Top$. Let $S=\cup_{i=1}^l S_i$ an open cover and denote by $S_I=\cap_{i\in I} S_i$.
Let $i_i:S_i\hookrightarrow\tilde S_i$ closed embeddings, with $\tilde S_i\in\Top$. 
We denote by $C_{fil}(S/(\tilde S_I))\subset C_{fil}((\tilde S_I))$ the full subcategory  
\begin{itemize}
\item whose objects are $(G,F)=((G_I,F)_{I\subset\left[1,\cdots l\right]},u_{IJ})$,
with $(G_I,F)\in C_{fil,S_I}(\tilde S_I)$, and  
$u_{IJ}:m^*(G_I,F)\to m^*p_{IJ*}(G_J,F)$ are $\infty$-filtered top local equivalences satisfying for $I\subset J\subset K$, 
$p_{IJ*}u_{JK}\circ u_{IJ}=u_{IK}$ in $C_{fil}(\tilde S_I)$,
\item the morphisms $m:((G_I,F),u_{IJ})\to((H_I,F),v_{IJ})$ being (see section 2.1) a family of morphism of complexes,  
\begin{equation*}
m=(m_I:(G_I,F)\to (H_I,F))_{I\subset\left[1,\cdots l\right]}
\end{equation*}
such that $v_{IJ}\circ m_I=p_{IJ*}m_J\circ u_{IJ}$ in $C_{fil}(\tilde S_I)$.
\end{itemize}  
A morphism $m:((G_I,F),u_{IJ})\to((H_I,F),v_{IJ})$ is an $r$-filtered top local equivalence if
there exists $\phi_i:((C_{iI},F),u_{i,IJ})\to((C_{(i+1)I},F),u_{(i+1)IJ})$, $0\leq i\leq s$, 
with $((C_{iI},F),u_{iIJ})\in C_{fil}(S/(\tilde S_I))$
$((C_{0I},F),u_{0IJ})=((G_I,F),u_{IJ})$ and $((C_{sI},F),u_{sIJ})=((H_I,F),v_{IJ})$, such that
\begin{equation*}
\phi=\phi_s\circ\cdots\circ\phi_i\circ\cdots\circ\phi_0:((G_I,F),u_{IJ})\to((H_I,F),v_{IJ})
\end{equation*}
and $\phi_i:((C_{iI},F),u_{iIJ})\to((C_{(i+1)I},F),u_{(i+1)IJ})$ either a filtered $\tau$ local equivalence
or an $r$-filtered homotopy equivalence. 
\end{defi}

Denote $L=[1,\ldots,l]$ and for $I\subset L$,
$p_{0(0I)}:S\times\tilde S_I\to S$, $p_{I(0I)}:S\times\tilde S_I\to S_I$ the projections.By definition, we have functors 
\begin{itemize}
\item $T(S/(\tilde S_I)):C_{fil}(S)\to C_{fil}(S/(\tilde S_I))$, $(G,F)\mapsto (i_{I*}j_I^*(G,F),I)$
\item $T((\tilde S_I)/S):C_{fil}(S/(\tilde S_I))\to C_{fil}(S)$,
$((G_I,F),u_{IJ})\mapsto \ho\lim_{I\subset L}p_{0(0I)*}\Gamma^{\vee}_{S_I}p_{I(0I)}^*(G_I,F)$.
\end{itemize}
Note that the functors $T(S/(\tilde S_I)$ are embedding, since
\begin{equation*}
\ad(i_I^*,i_{I*})(j_I^*F):i_I^*i_{I*}j_I^*F\to j_I^*F
\end{equation*}
are top local equivalence.

Let $f:X\to S$ a morphism, with $X,S\in\Top$.
Let $S=\cup_{i=1}^{l} S_i$ and $X=\cup_{i=1}^l X_i$ be open covers
and $i_i:S_i\hookrightarrow\tilde S_i$, $i'_i:X_i\hookrightarrow\tilde X_i$ be closed embeddings,
such that, for each $i\in[1,l]$, $f_i:=f_{|X_i}:X_i\to S_i$ lift to a morphism $\tilde f_i:\tilde X_i\to\tilde S_i$. 
Then, $f_I=f_{|X_I}:X_I=\cap_{i\in I} X_i\to S_I=\cap_{i\in I} S_i$ lift to the morphism 
\begin{equation*}
\tilde f_I=\Pi_{i\in I}\tilde f_i:\tilde X_I=\Pi_{i\in I}\tilde X_i\to \tilde S_I=\Pi_{i\in I}\tilde S_i
\end{equation*}
Denote by $p_{IJ}:\tilde S_J\to\tilde S_I$ and $p'_{IJ}:\tilde X_J\to\tilde X_I$ the projections. 
Consider for $J\subset I$ the following commutative diagrams
\begin{equation*}
D_{IJ}=\xymatrix{ S_I\ar[r]^{i_I} & \tilde S_I \\
S_J\ar[u]^{j_{IJ}}\ar[r]^{i_J} & \tilde S_J\ar[u]^{p_{IJ}}} \;, \;  
D'_{IJ}=\xymatrix{ X_I\ar[r]^{i'_I} & \tilde X_I \\
X_J\ar[u]^{j'_{IJ}}\ar[r]^{i'_J} & \tilde X_J\ar[u]^{p'_{IJ}}} \; , \;  
D_{fI}=\xymatrix{ S_I\ar[r]^{i_I} & \tilde S_I \\
X_I\ar[u]^{f_I}\ar[r]^{i'_I} & \tilde X_I\ar[u]^{\tilde f_I}}  
\end{equation*}
We have then following commutative diagram
\begin{equation*}
\xymatrix{ \, & X_I\ar[r]^{n'_I} & \tilde X_I & \tilde X_I\backslash X_I\ar[l]^{n'_I} \\
i'_J:X_J\ar[r]^{l_{IJ}}\ar[ru]^{j'_{IJ}}& X_I\times X_I\times\tilde X_{J\backslash I}\ar[u]^{p'_{IJ}}\ar[r]^{n'_I\times I} & 
\tilde X_J\ar[u]^{p'_{IJ}} & \tilde X_J\backslash X_J\ar[l]^{n'_J}\ar[u]^{p'_{IJ}}}.
\end{equation*}
whose square are cartesian.
We then have the pullback functor
\begin{eqnarray*}
f^*:C_{(2)fil}(S/(\tilde S_I))\to C_{(2)fil}(X/(\tilde X_I)), \\
((G_I,F),u_{IJ})\mapsto f^*((G_I,F),u_{IJ}):=(\Gamma^{\vee}_{X_I}\tilde f_I^*(G_I,F),\tilde f_J^*u_{IJ})
\end{eqnarray*}
with 
\begin{eqnarray*}
\tilde f_J^*u_{IJ}:\Gamma^{\vee}_{X_I}\tilde f_I^*(G_I,F)
\xrightarrow{\ad(p_{IJ}^{'*},p'_{IJ*})(-)}p'_{IJ*}p_{IJ}^{'*}\Gamma^{\vee}_{X_I}\tilde f_I^*(G_I,F) 
\xrightarrow{T_{\sharp}(p_{IJ},n'_I)(-)^{-1})}
p'_{IJ*}\Gamma^{\vee}_{X_I\times\tilde X_{J\backslash I}}p_{IJ}^{'*}\tilde f_I^*(G_I,F) \\
\xrightarrow{p'_{IJ*}\gamma^{\vee}_{X_J}(-)}
p'_{IJ*}\Gamma^{\vee}_{X_J}p_{IJ}^{'*}\tilde f_I^*(G_I,F)=p'_{IJ*}\Gamma^{\vee}_{X_J}\tilde f_J^*p_{IJ}^*(G_I,F)
\xrightarrow{\Gamma^{\vee}_{X_J}\tilde f_J^*I(p_{IJ}^*,p_{IJ*})(-,-)(u_{IJ})}
\Gamma^{\vee}_{X_J}\tilde f_J^*(G_J,F)
\end{eqnarray*}
Let $(G,F)\in C_{fil}(S)$. Since, $j_I^{'*}i'_{I*}j_I^{'*}f^*(G,F)=0$,
the morphism $T(D_{fI})(j_I^*(G,F)):\tilde f_I^*i_{I*}j_I^*(G,F)\to i'_{I*}j_I^{'*}f^*(G,F)$ factors trough
\begin{equation*}
T(D_{fI})(j_I^*(G,F)):\tilde f_I^*i_{I*}j_I^*(G,F)\xrightarrow{\gamma_{X_I}^{\vee}(-)}
\Gamma_{X_I}^{\vee}\tilde f_I^*i_{I*}j_I^*(G,F)\xrightarrow{T^{\gamma}(D_{fI})(j_I^*(G,F))} i'_{I*}j_I^{'*}f^*(G,F)
\end{equation*}
We have then, for $(G,F)\in C_{fil}(S)$, the canonical transformation map
\begin{equation*}
\xymatrix{
f^*T(S/(\tilde S_I))(G,F)\ar[rrr]^{T(f,T(0/I))(G,F)}\ar[d]_{=}\ & \, & \, & 
T(X/(\tilde X_I))(f^*(G,F))\ar[d]^{=} \\
(\Gamma^{\vee}_{X_I}\tilde f_I^*i_{I*}j_I^*(G,F),\tilde f_J^*I)\ar[rrr]^{T^{\gamma}(D_{fI})(j_I^*(G,F))} 
& \, & \, & (i'_{I*}j_I^{'*}f^*(G,F),I)}
\end{equation*}

\begin{prop}\label{TSSItop}
Let $S\in\Top$. Let $S=\cup_{i=1}^l S_i$ an open cover and denote by $S_I=\cap_{i\in I} S_i$.
Let $i_i:S_i\hookrightarrow\tilde S_i$ closed embeddings, with $\tilde S_i\in\Top$. 
Denote by $D_{(2)fil,\infty}((S/(\tilde S_I)))=\Ho_{top,\infty}(C_{(2)fil}((S/(\tilde S_I))))$ 
the localization of $C_{(2)fil}((S/(\tilde S_I)))$ with respect to top local equivalences.
The functor $T(S/(\tilde S_I))$ induces an equivalence of category 
\begin{equation*}
T(S/(\tilde S_I)):D_{(2)fil,\infty}(S)\xrightarrow{\sim} D_{(2)fil,\infty}((S/(\tilde S_I)))
\end{equation*}
with inverse $T((\tilde S_I)/S)$
\end{prop}

\begin{proof}
Follows from the fact that for $(G,F)\in C_{fil}(S)$, 
\begin{equation*}
\ho\lim_{I\subset L}p_{0(0I)*}\Gamma^{\vee}_{S_I}p_{I(0I)}^*(i_{I*}j_I^*(G,F))\to p_{0(0I)*}\Gamma^{\vee}_{S_I}j_I^*(G,F)
\end{equation*}
is an equivalence top local.
\end{proof}

For $f:T\to S$ a morphism with $T,S\in\Top$ locally compact (in particular Hausdorf), e.g. $T,S\in\CW$,
there is also a functor $f_!:C(T)\to C(S)$ given by the section which have compact support over $f$, and,
for $K_1,K_2\in C(T)$, we have a canonical map
\begin{equation*}
T_!(f,hom):f_*\mathcal Hom(K_1,K_2)\to \mathcal Hom(f_!K_1,f_!K_2)
\end{equation*}

The main result on presheaves on locally compact spaces is the following :
\begin{thm}\label{Vdual}
Let $f:T\to S$ a morphism with $T,S\in\Top$ locally compact.
\begin{itemize}
\item[(i)] The derived functor $Rf_!:D(T)\to D(S)$ has a right adjoint $f^!$ (Verdier duality) 
and, for $K_1,K_2\in D(T)$ and $K_3,K_4\in D(S)$, we have canonical isomorphisms 
\begin{itemize}
\item $Rf_*R\mathcal Hom^{\bullet}(Rf_!K_1,K_3)\xrightarrow{\sim}R\mathcal Hom^{\bullet}(K_1,f^!K_3)$
\item $f^!R\mathcal Hom^{\bullet}(K_3,K_4)\xrightarrow{\sim}R\mathcal Hom^{\bullet}(f^*K_3,f^!K_4)$
\end{itemize}
\item[(ii)] Denote, for $K\in D(S)$, $D(K)=R\mathcal Hom^{\bullet}(K,a_S^!\mathbb Z)\in D(S)$ the Verdier dual of $K$.
Then, if $K\in D_c(S)$, the evaluation map $ev^c(S)(K):K\to D(D(K)$ is an isomorphism.
\item[(iii)] Assume we have a factorization $f:T\xrightarrow{l}Y\xrightarrow{p}S$ of $f$ with
$l$ a closed embedding and $p$ a smooth morphism of relative dimension $d$. Then $f^!K=i^!p^*K[d]$
\end{itemize}
\end{thm}

\begin{proof}
\noindent(i):Standard, the proof is formal (see \cite{VerdierC}).

\noindent(ii): See \cite{VerdierC}.

\noindent(iii): The fact that $p^!K=p^*K[d]$ follows by Poincare duality for topological manifold.
\end{proof}

We have by theorem \ref{Vdual} a pair of adjoint functor
\begin{equation*}
(Rf_!,f^!):D(T)\leftrightarrows D(S) 
\end{equation*}
\begin{itemize}
\item with $f_!=f_*$ if $f$ is proper,
\item with $f^!=f^*[d]$ if $f$ is smooth of relative dimension $d$. 
\end{itemize}

\subsection{Presheaves on the big Zariski site or on the big etale site}

For $S\in\Var(\mathbb C)$, we denote by $\rho_S:\Var(\mathbb C)^{sm}/S\hookrightarrow\Var(\mathbb C)/S$ be the full subcategory 
consisting of the objects $U/S=(U,h)\in\Var(\mathbb C)/S$ such that the morphism $h:U\to S$ is smooth. 
That is, $\Var(\mathbb C)^{sm}/S$ is the category  
\begin{itemize}
\item whose objects are smooth morphisms 
$U/S=(U,h)$, $h:U\to S$ with $U\in\Var(\mathbb C)$, 
\item whose morphisms $g:U/S=(U,h_1)\to V/S=(V,h_2)$ 
is a morphism $g:U\to V$ of complex algebraic varieties such that $h_2\circ g=h_1$. 
\end{itemize}
We denote again $\rho_S:\Var(\mathbb C)/S\to\Var(\mathbb C)^{sm}/S$ the associated morphism of site.
We will consider 
\begin{equation*}
r^s(S):\Var(\mathbb C)\xrightarrow{r(S)}\Var(\mathbb C)/S\xrightarrow{\rho_S}\Var(\mathbb C)^{sm}/S
\end{equation*}
the composite morphism of site.
For $S\in\Var(\mathbb C)$, we denote by $\mathbb Z_S:=\mathbb Z(S/S)\in \PSh(\Var(\mathbb C)^{sm}/S)$ the constant presheaf
By Yoneda lemma, we have for $F\in C(\Var(\mathbb C)^{sm}/S)$, $\mathcal Hom(\mathbb Z_S,F)=F$.
For $f:T\to S$ a morphism, with $T,S\in\Var(\mathbb C)$, we have the following commutative diagram of sites
\begin{equation}\label{pf0}
\xymatrix{\Var(\mathbb C)/T\ar[d]^{P(f)}\ar[r]^{\rho_T} & \Var(\mathbb C)^{sm}/T\ar[d]^{P(f)} \\ 
\Var(\mathbb C)/S\ar[r]^{\rho_S} & \Var(\mathbb C)^{sm}/S} 
\end{equation}
We denote, for $S\in\Var(\mathbb C)$, the obvious morphism of sites 
\begin{equation*}
\tilde e(S):\Var(\mathbb C)/S\xrightarrow{\rho_S}\Var(\mathbb C)^{sm}/S\xrightarrow{e(S)}\Ouv(S)  
\end{equation*}
where $\Ouv(S)$ is the set of the Zariski open subsets of $S$, given by the inclusion functors
$\tilde e(S):\Ouv(S)\hookrightarrow\Var(\mathbb C)^{sm}/S\hookrightarrow\Var(\mathbb C)/S$.
By definition, for $f:T\to S$ a morphism with $S,T\in\Var(\mathbb C)$, the commutative diagram of sites (\ref{pf0})
extend a commutative diagram of sites :
\begin{equation}\label{empf0}
\xymatrix{
\tilde e(T):\Var(\mathbb C)/T\ar[d]^{P(f)}\ar[rr]^{\rho_T} & \, & \Var(\mathbb C)^{sm}/T\ar[d]^{P(f)}\ar[rr]^{e(T)} & \, & 
\Ouv(T)\ar[d]^{P(f)} \\
\tilde e(S):\Var(\mathbb C)/S\ar[rr]^{\rho_S} & \, & \Var(\mathbb C)^{sm}/S\ar[rr]^{e(S)} & \, & \Ouv(S)}
\end{equation}

\begin{itemize}
\item As usual, we denote by
\begin{equation*}
(f^*,f_*):=(P(f)^*,P(f)_*):C(\Var(\mathbb C)^{sm}/S)\to C(\Var(\mathbb C)^{sm}/T)
\end{equation*}
the adjonction induced by $P(f):\Var(\mathbb C)^{sm}/T\to \Var(\mathbb C)^{sm}/S$.
Since the colimits involved in the definition of $f^*=P(f)^*$ are filtered, $f^*$ also preserve monomorphism. 
Hence, we get an adjonction
\begin{equation*}
(f^*,f_*):C_{fil}(\Var(\mathbb C)^{sm}/S)\leftrightarrows C_{fil}(\Var(\mathbb C)^{sm}/T), \; f^*(G,F):=(f^*G,f^*F)
\end{equation*}
\item As usual, we denote by
\begin{equation*}
(f^*,f_*):=(P(f)^*,P(f)_*):C(\Var(\mathbb C)/S)\to C(\Var(\mathbb C)/T)
\end{equation*}
the adjonction induced by $P(f):\Var(\mathbb C)/T\to \Var(\mathbb C)/S$.
Since the colimits involved in the definition of $f^*=P(f)^*$ are filtered, $f^*$ also preserve monomorphism. 
Hence, we get an adjonction
\begin{equation*}
(f^*,f_*):C_{fil}(\Var(\mathbb C)/S)\leftrightarrows C_{fil}(\Var(\mathbb C)/T), \; f^*(G,F):=(f^*G,f^*F)
\end{equation*}
\end{itemize}

\begin{itemize}
\item For $h:U\to S$ a smooth morphism with $U,S\in\Var(\mathbb C)$, 
the pullback functor $P(h):\Var(\mathbb C)^{sm}/S\to \Var(\mathbb C)^{sm}/U$ 
admits a left adjoint $C(h)(X\to U)=(X\to U\to S)$.
Hence, $h^*:C(\Var(\mathbb C)^{sm}/S)\to C(\Var(\mathbb C)^{sm}/U)$ admits a left adjoint
\begin{equation*}
h_{\sharp}:C(\Var(\mathbb C)^{sm}/U)\to C(\Var(\mathbb C)^{sm}/S), \; 
F\mapsto((V,h_0)\mapsto\lim_{(V',h\circ h')\to(V,h_0)}F(V',h'))
\end{equation*}
Note that we have for $V/U=(V,h')$ with $h':V\to U$ a smooth morphism 
we have $h_{\sharp}(\mathbb Z(V/U))=\mathbb Z(V'/S)$ with $V'/S=(V',h\circ h')$.
Hence, since projective presheaves are the direct summands of the representable presheaves,
$h_{\sharp}$ sends projective presheaves to projective presheaves.
For $F^{\bullet}\in C(\Var(\mathbb C)^{sm}/S)$ and $G^{\bullet}\in C(\Var(\mathbb C)^{sm}/U)$,
we have the adjonction maps
\begin{equation*}
\ad(h_{\sharp},h^*)(G^{\bullet}):G^{\bullet}\to h^*h_{\sharp}G^{\bullet} \; , \;
\ad(h_{\sharp},h^*)(F^{\bullet}):h_{\sharp}h^*F^{\bullet}\to F^{\bullet}.
\end{equation*}
For a smooth morphism $h:U\to S$, with $U,S\in\Var(\mathbb C)$, we have the adjonction isomorphism, 
for $F\in C(\Var(\mathbb C)^{sm}/U)$ and $G\in C(\Var(\mathbb C)^{sm}/S)$,   
\begin{equation}\label{Ihhom}
I(h_{\sharp},h^*)(F,G):\mathcal Hom^{\bullet}(h_{\sharp}F,G)\xrightarrow{\sim}h_*\mathcal Hom^{\bullet}(F,h^*G).  
\end{equation}

\item For $f:T\to S$ any morphism with $T,S\in\Var(\mathbb C)$, 
the pullback functor $P(f):\Var(\mathbb C)/T\to \Var(\mathbb C)/S$ 
admits a left adjoint $C(f)(X\to T)=(X\to T\to S)$.
Hence, $f^*:C(\Var(\mathbb C)/S)\to C(\Var(\mathbb C)/T)$ admits a left adjoint
\begin{equation*}
f_{\sharp}:C(\Var(\mathbb C)/T)\to C(\Var(\mathbb C)/S), \; 
F\mapsto((V,h_0)\mapsto\lim_{(V',f\circ h')\to(V,h_0)}F(V',h'))
\end{equation*}
Note that we have, for $(V/T)=(V,h)$, $f_{\sharp}\mathbb Z(V/T)=\mathbb Z(V/S)$ with $V/S=(V,f\circ h)$.
Hence, since projective presheaves are the direct summands of the representable presheaves,
$h_{\sharp}$ sends projective presheaves to projective presheaves.
For $F^{\bullet}\in C(\Var(\mathbb C)/S)$ and $G^{\bullet}\in C(\Var(\mathbb C)/T)$,
we have the adjonction maps
\begin{equation*}
\ad(f_{\sharp},f^*)(G^{\bullet}):G^{\bullet}\to f^*f_{\sharp}G^{\bullet} \; , \;
\ad(f_{\sharp},f^*)(F^{\bullet}):f_{\sharp}f^*F^{\bullet}\to F^{\bullet}.
\end{equation*}
For a morphism $f:T\to S$, with $T,S\in\Var(\mathbb C)$, we have the adjonction isomorphism, 
for $F\in C(\Var(\mathbb C)/T)$ and $G\in C(\Var(\mathbb C)/S)$,   
\begin{equation}\label{Ifsharphom}
I(f_{\sharp},f^*)(F,G):\mathcal Hom^{\bullet}(f_{\sharp}F,G)\xrightarrow{\sim}f_*\mathcal Hom^{\bullet}(F,f^*G).  
\end{equation}
\end{itemize}

\begin{itemize}
\item For a commutative diagram in $\Var(\mathbb C)$ : 
\begin{equation*}
D=\xymatrix{ 
V\ar[r]^{g_2}\ar[d]^{h_2} & U\ar[d]^{h_1} \\
T\ar[r]^{g_1} & S},
\end{equation*}
where $h_1$ and $h_2$ are smooth,
we denote by, for $F^{\bullet}\in C(\Var(\mathbb C)^{sm}/U)$, 
\begin{equation*}
T_{\sharp}(D)(F^{\bullet}): h_{2\sharp}g_2^*F^{\bullet}\to g_1^*h_{1\sharp}F^{\bullet}
\end{equation*}
the canonical map in $C(\Var(\mathbb C)^{sm}/T)$ given by adjonction. 
If $D$ is cartesian with $h_1=h$, $g_1=g$ $f_2=h':U_T\to T$, $g':U_T\to U$,
\begin{equation*}
T_{\sharp}(D)(F^{\bullet})=:T_{\sharp}(g,h)(F^{\bullet}):
h'_{\sharp}g^{'*}F^{\bullet}\xrightarrow{\sim} g^*h_{\sharp}F^{\bullet}
\end{equation*}
is an isomorphism and for $G^{\bullet}\in C(\Var(\mathbb C)^{sm}/T)$
\begin{equation*}
T(D)(G^{\bullet})=:T(g,h)(G^{\bullet}):g^*h_*G^{\bullet}\xrightarrow{\sim} h'_*g^{'*}G^{\bullet}
\end{equation*}
is an isomorphism.

\item For a commutative diagram in $\Var(\mathbb C)$ : 
\begin{equation*}
D=\xymatrix{ 
V\ar[r]^{g_2}\ar[d]^{f_2} & X\ar[d]^{f_1} \\
T\ar[r]^{g_1} & S},
\end{equation*}
we denote by, for $F^{\bullet}\in C(\Var(\mathbb C)/X)$, 
\begin{equation*}
T_{\sharp}(D)(F^{\bullet}): f_{2\sharp}g_2^*F^{\bullet}\to g_1^*f_{1\sharp}F^{\bullet}
\end{equation*}
the canonical map in $C(\Var(\mathbb C)/T)$ given by adjonction. 
If $D$ is cartesian with $h_1=h$, $g_1=g$ $f_2=h':X_T\to T$, $g':X_T\to X$,
\begin{equation*}
T_{\sharp}(D)(F^{\bullet})=:T_{\sharp}(g,f)(F^{\bullet}):
f'_{\sharp}g^{'*}F^{\bullet}\xrightarrow{\sim} g^*f_{\sharp}F^{\bullet}
\end{equation*}
is an isomorphism and for $G^{\bullet}\in C(\Var(\mathbb C)/T)$
\begin{equation*}
T(D)(G^{\bullet})=:T(g,h)(G^{\bullet}):f^*g_*G^{\bullet}\xrightarrow{\sim} g'_*f^{'*}G^{\bullet}
\end{equation*}
is an isomorphism.
\end{itemize}

For $f:T\to S$ a morphism with $S,T\in\Var(\mathbb C)$, 
\begin{itemize}
\item we get for $F\in C(\Var(\mathbb C)^{sm}/S)$ from the a commutative diagram of sites (\ref{empf0}) 
the following canonical transformation 
\begin{equation*}
T(e,f)(F^{\bullet}):f^*e(S)_*F^{\bullet}\to e(T)_*f^*F^{\bullet},
\end{equation*}
 which is NOT a quasi-isomorphism in general. 
However, for $h:U\to S$ a smooth morphism with $S,U\in\Var(\mathbb C)$, 
$T(e,h)(F^{\bullet}):h^*e(S)_*F^{\bullet}\xrightarrow{\sim} e(T)_*h^*F^{\bullet}$ is an isomorphism.
\item we get for $F\in C(\Var(\mathbb C)/S)$ from the a commutative diagram of sites (\ref{empf0}) 
the following canonical transformation 
\begin{equation*}
T(e,f)(F^{\bullet}):f^*e(S)_*F^{\bullet}\to e(T)_*f^*F^{\bullet},
\end{equation*}
 which is NOT a quasi-isomorphism in general. 
However, for $h:U\to S$ a smooth morphism with $S,U\in\Var(\mathbb C)$, 
$T(e,h)(F^{\bullet}):h^*e(S)_*F^{\bullet}\xrightarrow{\sim} e(T)_*h^*F^{\bullet}$ is an isomorphism.
\end{itemize}

Let $S\in\Var(\mathbb C)$, 
\begin{itemize}
\item We have for $F,G\in C(\Var(\mathbb C)^{sm}/S)$, 
\begin{itemize}
\item $e(S)_*(F\otimes G)=(e(S)_*F)\otimes (e(S)_*G)$ by definition 
\item the canonical forgetfull map 
\begin{equation*}
T(S,hom)(F,G):e(S)_*\mathcal Hom^{\bullet}(F,G)\to\mathcal Hom^{\bullet}(e(S)_*F,e(S)_*G).
\end{equation*}
which is NOT a quasi-isomorphism in general.
\end{itemize}
By definition, we have for $F\in C(\Var(\mathbb C)^{sm}/S)$, $e(S)_*E_{zar}(F)=E_{zar}(e(S)_*F)$.
\item We have for $F,G\in C(\Var(\mathbb C)/S)$, 
\begin{itemize}
\item $e(S)_*(F\otimes G)=(e(S)_*F)\otimes (e(S)_*G)$ by definition 
\item the canonical forgetfull map 
\begin{equation*}
T(S,hom)(F,G):e(S)_*\mathcal Hom^{\bullet}(F,G)\to\mathcal Hom^{\bullet}(e(S)_*F,e(S)_*G).
\end{equation*}
which is NOT a quasi-isomorphism in general.
\end{itemize}
By definition, we have for $F\in C(\Var(\mathbb C)/S)$, $e(S)_*E_{zar}(F)=E_{zar}(e(S)_*F)$.
\end{itemize}

Let $S\in\Var(\mathbb C)$. 
\begin{itemize}
\item We have the dual functor
\begin{eqnarray*}
\mathbb D_S:C(\Var(\mathbb C)^{sm}/S)\to C(\Var(\mathbb C)^{sm}/S), \; 
F\mapsto\mathbb D_S(F):=\mathcal Hom(F,E_{et}(\mathbb Z_S))
\end{eqnarray*}
It induces the functor
\begin{eqnarray*}
L\mathbb D_S:C(\Var(\mathbb C)^{sm}/S)\to C(\Var(\mathbb C)^{sm}/S), \; 
F\mapsto L\mathbb D_S(F):=\mathbb D_S(LF):=\mathcal Hom(LF,E_{et}(\mathbb Z_S))
\end{eqnarray*}
\item We have the dual functor
\begin{eqnarray*}
\mathbb D_S:C(\Var(\mathbb C)/S)\to C(\Var(\mathbb C)/S), \; 
F\mapsto\mathbb D_S(F):=\mathcal Hom(F,E_{et}(\mathbb Z_S))
\end{eqnarray*}
It induces the functor
\begin{eqnarray*}
L\mathbb D_S:C(\Var(\mathbb C)/S)\to C(\Var(\mathbb C)/S), \; 
F\mapsto L\mathbb D_S(F):=\mathbb D_S(LF):=\mathcal Hom(LF,E_{et}(\mathbb Z_S))
\end{eqnarray*}
\end{itemize}

The adjonctions 
\begin{equation*}
(\tilde e(S)^*,\tilde e(S)_*):C(\Var(\mathbb C)/S)\leftrightarrows C(S) \; , \; 
(e(S)^*,e(S)_*):C(\Var(\mathbb C)^{sm}/S)\leftrightarrows C(S)
\end{equation*}
induce adjonctions
\begin{equation*}
(\tilde e(S)^*,\tilde e(S)_*):C_{fil}(\Var(\mathbb C)/S)\leftrightarrows C_{fil}(S) \; , \;
(e(S)^*,e(S)_*):C_{fil}(\Var(\mathbb C)^{sm}/S)\leftrightarrows C_{fil}(S)
\end{equation*}
given by $e(S)_*(G,F):=(e(S)_*G,e(S)_*F)$, since $e(S)_*$ and $e(S)^*$ preserve monomorphisms. Note that
\begin{itemize}
\item for $F\in\PSh(\Var(\mathbb C)^{sm}/S)$, $e(S)_*F$ is simply the restriction of $F$ to the small Zariski site of $X$, 
\item for $F\in\PSh(\Var(\mathbb C)/S)$, $\tilde e(S)_*F=e(S)_*\rho_{S*}F$ 
is simply the restriction of $F$ to the small Zariski site of $X$,
$\rho_{S*}F$ being the restriction of $F$ to $\Var(\mathbb C)^{sm}/S$.  
\end{itemize}
Together with the internal hom functor, we get the bifunctor,
\begin{eqnarray}
e(S)_*\mathcal Hom(\cdot,\cdot):
C_{fil}(\Var(\mathbb C)^{sm}/S)\times C_{fil}(\Var(\mathbb C)^{sm}/S)\to C_{2fil}(S), \\
((F,W),(G,F))\mapsto
e(S)_*\mathcal Hom^{\bullet}((F^{\bullet},W),(G^{\bullet},F)).
\end{eqnarray}

For $i:Z\hookrightarrow S$ a closed embedding, with $Z,S\in\Var(\mathbb C)$, we denote by
\begin{equation*}
(i_*,i^!):=(P(i)_*,P(i)^{\bot}):C(\Var(\mathbb C)^{sm}/Z)\leftrightarrows C(\Var(\mathbb C)^{sm}/S)
\end{equation*}
the adjonction induced by the morphism of site $P(i):\Var(\mathbb C)^{sm}/Z\to\Var(\mathbb C)^{sm}/S$
For $i:Z\hookrightarrow S$ a closed embedding, $Z,S\in\Var(\mathbb C)$, we denote 
\begin{equation*}
\mathbb Z_{Z,S}:=\Cone(\ad(i^*,i_*)(\mathbb Z_S):\mathbb Z_S\to i_*\mathbb Z_Z)
\end{equation*}

We have the support section functors of a closed embedding $i:Z\hookrightarrow S$ for presheaves on the big Zariski site.
\begin{defi}\label{gamma}
Let $i:Z\hookrightarrow S$ be a closed embedding with $S,Z\in\Var(\mathbb C)$ 
and $j:S\backslash Z\hookrightarrow S$ be the open complementary subset.
\begin{itemize}
\item[(i)] We define the functor
\begin{equation*}
\Gamma_Z:C(\Var(\mathbb C)^{sm}/S)\to C(\Var(\mathbb C)^{sm}/S), \;
G^{\bullet}\mapsto\Gamma_Z G^{\bullet}:=\Cone(\ad(j^*,j_*)(G^{\bullet}):G^{\bullet}\to j_*j^*G^{\bullet})[-1],
\end{equation*}
so that there is then a canonical map $\gamma_Z(G^{\bullet}):\Gamma_ZG^{\bullet}\to G^{\bullet}$.
\item[(ii)] We have the dual functor of (i) :
\begin{equation*}
\Gamma^{\vee}_Z:C(\Var(\mathbb C)^{sm}/S)\to C(\Var(\mathbb C)^{sm}/S), \; 
F\mapsto\Gamma^{\vee}_Z(F^{\bullet}):=\Cone(\ad(j_{\sharp},j^*)(G^{\bullet}):j_{\sharp}j^*G^{\bullet}\to G^{\bullet}), 
\end{equation*}
together with the canonical map $\gamma^{\vee}_Z(G):F\to\Gamma^{\vee}_Z(G)$.
\item[(iii)] For $F,G\in C(\Var(\mathbb C)^{sm}/S)$, we denote by 
\begin{equation*}
I(\gamma,hom)(F,G):=(I,I(j_{\sharp},j^*)(F,G)^{-1}):\Gamma_Z\mathcal Hom(F,G)\xrightarrow{\sim}\mathcal Hom(\Gamma^{\vee}_ZF,G)
\end{equation*}
the canonical isomorphism given by adjonction.
\end{itemize}
\end{defi}

Let $i:Z\hookrightarrow S$ be a closed embedding with $S,Z\in\Var(\mathbb C)$ 
and $j:S\backslash Z\hookrightarrow S$ be the open complementary.
\begin{itemize}
\item For $G\in C(\Var(\mathbb C)^{sm}/S)$, the adjonction map $\ad(i_*,i^!)(G):i_*i^!G\to G$
factor through $\gamma_Z(G)$ :
\begin{equation*}
\ad(i_*,i^!)(G):i_*i^!G\xrightarrow{\ad(i_*,i^!)(G)^{\gamma}} \Gamma_Z(G)\xrightarrow{\gamma_Z(G)} G. 
\end{equation*}
However, note that when dealing with the big sites $P(i):\Var(\mathbb C)^{sm}/Z\to\Var(\mathbb C)^{sm}/S$,
if $G\in C(\Var(\mathbb C)^{sm}/S)$ is not $\mathbb A_S^1$ local and Zariski fibrant,
\begin{equation*}
\ad(i_*,i^!)(G)^{\gamma}:i_*i^!G\to\Gamma_Z(G)
\end{equation*}
is NOT and homotopy equivalence, and $\Gamma_ZG\in C(\Var(\mathbb C)^{sm}/S)$ is NOT in general in the image of the functor 
$i_*:C(\Var(\mathbb C)^{sm}/Z)\to C(\Var(\mathbb C)^{sm}/S)$. 
\item For $G\in C(\Var(\mathbb C)^{sm}/S)$, the adjonction map $\ad(i^*,i_*)(G):G\to i_*i^*G$
factor through $\gamma^{\vee}_Z(G)$ : 
\begin{equation*}
\ad(i^*,i_*)(G):G\xrightarrow{\gamma_Z^{\vee}(G)}\Gamma_Z^{\vee}G\xrightarrow{\ad(i^*,i_*)(G)^{\gamma}}i_*i^*G,
\end{equation*}
and as in (i), $\ad(i^*,i_*)(G)^{\gamma}:\Gamma^{\vee}_Z(G)\to i_*i^*G$ is NOT an homotopy equivalence but
\end{itemize}

Let $i:Z\hookrightarrow S$ be a closed embedding with $S,Z\in\Var(\mathbb C)$.
\begin{itemize}
\item Since $\Gamma_Z:C(\Var(\mathbb C)^{sm}/S)\to C(\Var(\mathbb C)^{sm}/S)$ preserve monomorphism, it induces a functor  
\begin{equation*}
\Gamma_Z:C_{fil}(\Var(\mathbb C)^{sm}/S)\to C_{fil}(\Var(\mathbb C)^{sm}/S), \; 
(G,F)\mapsto\Gamma_Z(G,F):=(\Gamma_ZG,\Gamma_ZF)
\end{equation*}
\item Since $\Gamma^{\vee}_Z:C(\Var(\mathbb C)^{sm}/S)\to C(\Var(\mathbb C)^{sm}/S)$ preserve monomorphism, it induces a functor  
\begin{equation*}
\Gamma^{\vee}_Z:C_{fil}(\Var(\mathbb C)^{sm}/S)\to C_{fil}(\Var(\mathbb C)^{sm}/S), \; 
(G,F)\mapsto\Gamma^{\vee}_Z(G,F):=(\Gamma^{\vee}_ZG,\Gamma^{\vee}_ZF)
\end{equation*}
\end{itemize}

\begin{defiprop}\label{gamma2sect2}
\begin{itemize}
\item[(i)] Let $g:S'\to S$ a morphism and $i:Z\hookrightarrow S$ a closed embedding with $S',S,Z\in\Var(\mathbb C)$. 
Then, for $(G,F)\in C_{fil}(\Var(\mathbb C)^{sm}/S)$, there exist a map in $C_{fil}(\Var(\mathbb C)^{sm}/S')$
\begin{equation*}
T(g,\gamma)(G,F):g^*\Gamma_{Z}(G,F)\to\Gamma_{Z\times_S S'}g^*(G,F)
\end{equation*}
unique up to homotopy such that $\gamma_{Z\times_S S'}(g^*(G,F))\circ T(g,\gamma)(G,F)=g^*\gamma_{Z}(G,F)$.
\item[(ii)] Let $i_1:Z_1\hookrightarrow S$, $i_2:Z_2\hookrightarrow Z_1$ be closed embeddings with $S,Z_1,Z_2\in\Var(\mathbb C)$.
Then, for $(G,F)\in C_{fil}(\Var(\mathbb C)^{sm}/S)$, 
\begin{itemize}
\item there exist a canonical map $T(Z_2/Z_1,\gamma)(G,F):\Gamma_{Z_2}(G,F)\to\Gamma_{Z_1}(G,F)$ in $C_{fil}(\Var(\mathbb C)^{sm}/S)$ 
unique up to homotopy such that $\gamma_{Z_1}(G,F)\circ T(Z_2/Z_1,\gamma)(G,F)=\gamma_{Z_2}(G,F)$, 
together with a distinguish triangle 
\begin{equation*}
\Gamma_{Z_2}(G,F)\xrightarrow{T(Z_2/Z_1,\gamma)(G,F)}\Gamma_{Z_1}(G,F)\xrightarrow{\ad(j_2^*,j_{2*})(\Gamma_{Z_1}(G,F))}
\Gamma_{Z_1\backslash Z_2}(G,F)\to\Gamma_{Z_2}(G,F)[1]
\end{equation*}
in $K_{fil}(\Var(\mathbb C)^{sm}/S):=K(\PSh_{fil}(\Var(\mathbb C)^{sm}/S))$,
\item there exist a map $T(Z_2/Z_1,\gamma^{\vee})(G,F):\Gamma_{Z_1}^{\vee}(G,F)\to\Gamma_{Z_2}^{\vee}(G,F)$ 
in $C_{fil}(\Var(\mathbb C)^{sm}/S)$
unique up to homotopy such that $\gamma^{\vee}_{Z_2}(G,F)=T(Z_2/Z_1,\gamma^{\vee})(G,F)\circ\gamma^{\vee}_{Z_1}(G,F)$, 
together with a distinguish triangle 
\begin{equation*}
\Gamma^{\vee}_{Z_1\backslash Z_2}(G,F)\xrightarrow{\ad(j_{2\sharp},j_2^*)(\Gamma^{\vee}_{Z_1}(G,F))}\Gamma^{\vee}_{Z_1}(G,F)
\xrightarrow{T(Z_2/Z_1,\gamma^{\vee})(G,F)}\Gamma^{\vee}_{Z_2}(G,F)\to\Gamma^{\vee}_{Z_1\backslash Z_2}(G,F)[1]
\end{equation*}
in $K_{fil}(\Var(\mathbb C)^{sm}/S)$.
\end{itemize}
\item[(iii)] Consider a morphism $g:(S',Z')\to(S,Z)$ with $(S',Z'),(S,Z)\in\Var^2(\mathbb C)$
We denote, for $G\in C(\Var(\mathbb C)^{sm}/S)$ the composite 
\begin{equation*}
T(D,\gamma^{\vee})(G):g^*\Gamma_Z^{\vee}G\xrightarrow{\sim}\Gamma^{\vee}_{Z\times_S S'}g^*G
\xrightarrow{T(Z'/Z\times_S S',\gamma^{\vee})(G)}\Gamma_{Z'}^{\vee}g^*G
\end{equation*}
and we have then the factorization
$\gamma_{Z'}^{\vee}(g^*G):g^*G\xrightarrow{g^*\gamma_Z^{\vee}(G)}g^*\Gamma_Z^{\vee}G
\xrightarrow{T(D,\gamma^{\vee})(G)}\Gamma_{Z'}^{\vee}g^*G$.
\end{itemize}
\end{defiprop}

\begin{proof}
\noindent(i): We have the cartesian square
\begin{equation*}
\xymatrix{S\backslash Z\ar[r]^{j} & S \\
S'\backslash Z\times_S S'\ar[r]^{j'}\ar[u]^{g'} & S'\ar[u]^{g}}
\end{equation*} 
and the map is given by 
\begin{equation*}
(I,T(g,j)(j^*G)):\Cone(g^*G\to g^*j_*j^*G)\to\Cone(g^*G\to j'_*j^{'*}g^*G=j'_*g^{'*}j^*G).
\end{equation*}

\noindent(ii): Follows from the fact that $j_1^*\Gamma_{Z_2}G=0$ and $j_1^*\Gamma^{\vee}_{Z_2}G=0$, 
with $j_1:S\backslash Z_1\hookrightarrow S$ the closed embedding.

\noindent(iii): Obvious.
\end{proof}

The following easy proposition concern the restriction from the big Zariski site to the small site Zariski site :

\begin{prop}\label{keylem2}
For $f:T\to S$ a morphism and $i:Z\hookrightarrow S$ a closed embedding, with $Z,S,T\in\Var(\mathbb C)$, we have
\begin{itemize}
\item[(i)] $e(S)_*f_*=f_*e(T)_*$ and $e(S)^*f_*=f_*e(T)^*$ 
\item[(ii)] $e(S)_*\Gamma_Z=\Gamma_Ze(S)_*$.
\end{itemize}
\end{prop}

\begin{proof}
\noindent(i):The first equality $e(S)_*f_*=f_*e(T)_*$ is given by the diagram (\ref{empf0}). 
The second equality is immediate from definition after a direct computation.

\noindent(ii) For $G^{\bullet}\in C(\Var(\mathbb C)^{sm}/S)$, we have the canonical equality
\begin{eqnarray*}
e(S)_*\Gamma_Z(G^{\bullet})=e(S)_*\Cone (G\to j_*j^*G^{\bullet})[-1]
&=&\Cone (e(S)_*G^{\bullet}\to e(S)_*j_*j^*G^{\bullet})[-1] \\
&=&\Cone(e(S)_*G^{\bullet}\to j_*j^*e(S)_*G^{\bullet})[-1] \\
&=&\Gamma_Ze(S)_*G^{\bullet}
\end{eqnarray*}
by (i) and since $j:S\backslash Z\hookrightarrow S$ is a smooth morphism.
\end{proof}

\begin{defi}
For $S\in\Var(\mathbb C)$, we denote by 
\begin{equation*}
C_{O_S}(\Var(\mathbb C)^{sm}/S):=C_{e(S)^*O_S}(\Var(\mathbb C)^{sm}/S) 
\end{equation*}
the category of complexes of presheaves on $\Var(\mathbb C)^{sm}/S$ endowed with a structure of $e(S)^*O_S$ module, and by 
\begin{equation*}
C_{O_Sfil}(\Var(\mathbb C)^{sm}/S):=C_{e(S)^*O_Sfil}(\Var(\mathbb C)^{sm}/S) 
\end{equation*}
the category of filtered complexes of presheaves on $\Var(\mathbb C)^{sm}/S$endowed with a structure of $e(S)^*O_S$ module.
\end{defi}

Let $S\in\Var(\mathbb C)$. Let $Z\subset S$ a closed subset.
Denote by $j:S\backslash Z\hookrightarrow S$ the open complementary embedding, 
\begin{itemize}
\item For $G\in C_{O_S}(\Var(\mathbb C)^{sm}/S)$, $\Gamma_ZG:=\Cone(\ad(j^*,j_*)(G):F\to j_*j^*G)[-1]$ 
has a (unique) structure of $e(S)^*O_S$ module such that $\gamma_Z(G):\Gamma_ZG\to G$ is a map in $C_{O_S}(\Var(\mathbb C)^{sm}/S)$. 
This gives the functor
\begin{equation*}
\Gamma_Z:C_{O_Sfil}(\Var(\mathbb C)^{sm}/S)\to C_{filO_S}(\Var(\mathbb C)^{sm}/S), \; 
(G,F)\mapsto\Gamma_Z(G,F):=(\Gamma_ZG,\Gamma_ZF), 
\end{equation*}
together with the canonical map $\gamma_Z((G,F):\Gamma_Z(G,F)\to (G,F)$.
Let $Z_2\subset Z$ a closed subset. Then, for $G\in C_{O_S}(\Var(\mathbb C)^{sm}/S)$, 
$T(Z_2/Z,\gamma)(G):\Gamma_{Z_2}G\to\Gamma_ZG$ is a map in $C_{O_S}(\Var(\mathbb C)^{sm}/S)$ (i.e. is $e(S)^*O_S$ linear).

\item For $G\in C_{O_S}(\Var(\mathbb C)^{sm}/S)$, $\Gamma^{\vee}_ZG:=\Cone(\ad(j_{\sharp},j^*)(G):j_{\sharp}j^*G\to G)$ 
has a unique structure of $e(S)^*O_S$ module, such that $\gamma^{\vee}_Z(G):G\to\Gamma_Z^{\vee}G$ 
is a map in $C_{O_S}(\Var(\mathbb C)^{sm}/S)$. This gives the the functor
\begin{equation*}
\Gamma^{\vee}_Z:C_{O_Sfil}(S)\to C_{filO_S}(S), \; (G,F)\mapsto\Gamma^{\vee}_Z(G,F):=(\Gamma^{\vee}_ZG,\Gamma^{\vee}_ZF), 
\end{equation*}
together with the canonical map $\gamma^{\vee}_Z((G,F):(G,F)\to\Gamma^{\vee}_Z(G,F)$.
Let $Z_2\subset Z$ a closed subset. Then, for $G\in C_{O_S}(\Var(\mathbb C)^{sm}/S)$, 
$T(Z_2/Z,\gamma^{\vee})(G):\Gamma_Z^{\vee}G\to\Gamma_{Z_2}^{\vee}G$ is a map in $C_{O_S}(\Var(\mathbb C)^{sm}/S)$ 
(i.e. is $e(S)^*O_S$ linear).
\end{itemize}

\begin{defi}\label{ZSvar}
Let $S\in\Var(\mathbb C)$. Let $Z\subset S$ a closed subset.
\begin{itemize}
\item[(i)] We denote by 
\begin{equation*}
C_Z(\Var(\mathbb C)^{sm}/S)\subset C(\Var(\mathbb C)^{sm}/S)
\end{equation*}
the full subcategory consisting of complexes of presheaves 
$F^{\bullet}\in C(\Var(\mathbb C)^{sm}/S)$ such that $a_{et}H^n(j^*F^{\bullet})=0$ for all $n\in\mathbb Z$,
where $j:S\backslash Z\hookrightarrow S$ is the complementary open embedding and $a_{et}$ is the sheaftification functor.
\item[(i)'] We denote by 
\begin{equation*}
C_{O_S,Z}(\Var(\mathbb C)^{sm}/S)\subset C_{O_S}(\Var(\mathbb C)^{sm}/S), 
\end{equation*}
the full subcategory consisting of complexes of presheaves 
$F^{\bullet}\in C(\Var(\mathbb C)^{sm}/S)$ such that $a_{et}H^n(j^*F^{\bullet})=0$ for all $n\in\mathbb Z$,
where $j:S\backslash Z\hookrightarrow S$ is the complementary open embedding and $a_{et}$ is the sheaftification functor.
\item[(ii)] We denote by 
\begin{equation*}
C_{fil,Z}(\Var(\mathbb C)^{sm}/S)\subset C_{fil}(\Var(\mathbb C)^{sm}/S) 
\end{equation*}
the full subcategory consisting of filtered complexes of presheaves $(F^{\bullet},F)\in C_{fil}(\Var(\mathbb C)^{sm}/S)$ 
such that there exist $r\in\mathbb N$ and an $r$-filtered homotopy equivalence 
$\phi:(F^{\bullet},F)\to(F^{'\bullet},F)$ with $(F^{'\bullet},F)\in C_{fil}(\Var(\mathbb C)^{sm}/S)$
such that $a_{et}j^*H^n\Gr_F^p(F^{'\bullet},F)=0$ for all $n,p\in\mathbb Z$,
where $j:S\backslash Z\hookrightarrow S$ is the complementary open embedding and $a_{et}$ is the sheaftification functor.
Note that by definition this $r$ does NOT depend on $n$ and $p$.
\item[(ii)'] We denote by 
\begin{equation*}
C_{O_Sfil,Z}(\Var(\mathbb C)^{sm}/S)\subset C_{O_Sfil}(\Var(\mathbb C)^{sm}/S) 
\end{equation*}
the full subcategory consisting of filtered complexes of presheaves $(F^{\bullet},F)\in C_{O_Sfil}(\Var(\mathbb C)^{sm}/S)$ 
such that there exist $r\in\mathbb N$ and an $r$-filtered homotopy equivalence 
$\phi:(F^{\bullet},F)\to(F^{'\bullet},F)$ with $(F^{'\bullet},F)\in C_{O_Sfil}(\Var(\mathbb C)^{sm}/S)$
such that $a_{et}j^*H^n\Gr_F^p(F^{'\bullet},F)=0$ for all $p,q\in\mathbb Z$,
where $j:S\backslash Z\hookrightarrow S$ is the complementary open embedding and $a_{et}$ is the sheaftification functor.
Note that by definition this $r$ does NOT depend on $n$ and $p$.
\end{itemize}
\end{defi}

Let $S\in\Var(\mathbb C)$ and $Z\subset S$ a closed subset. 
\begin{itemize}
\item For $(G,F)\in C_{fil}(\Var(\mathbb C)^{sm}/S)$, we have 
$\Gamma_Z(G,F),\Gamma_Z^{\vee}(G,F)\in C_{fil,Z}(\Var(\mathbb C)^{sm}/S)$.
\item For $(G,F)\in C_{O_Sfil}(\Var(\mathbb C)^{sm}/S)$, we have 
$\Gamma_Z(G,F),\Gamma_Z^{\vee}(G,F)\in C_{O_Sfil,Z}(\Var(\mathbb C)^{sm}/S)$.
\end{itemize}

Let $S_{\bullet}\in\Fun(\mathcal I,\Var(\mathbb C))$ with $\mathcal I\in\Cat$, a diagram of algebraic varieties.
It gives the diagram of sites $\Var(\mathbb C)^2/S_{\bullet}\in\Fun(\mathcal I,\Cat)$.  
\begin{itemize}
\item Then $C_{fil}(\Var(\mathbb C)/S_{\bullet})$ is the category  
\begin{itemize}
\item whose objects $(G,F)=((G_I,F)_{I\in\mathcal I},u_{IJ})$,
with $(G_I,F)\in C_{fil}(\Var(\mathbb C)/S_I)$,
and $u_{IJ}:(G_I,F)\to r_{IJ*}(G_J,F)$ for $r_{IJ}:I\to J$, denoting again $r_{IJ}:S_I\to S_J$, are morphisms
satisfying for $I\to J\to K$, $r_{IJ*}u_{JK}\circ u_{IJ}=u_{IK}$ in $C_{fil}(\Var(\mathbb C)/S_I)$,
\item the morphisms $m:((G,F),u_{IJ})\to((H,F),v_{IJ})$ being (see section 2.1) a family of morphisms of complexes,  
\begin{equation*}
m=(m_I:(G_I,F)\to (H_I,F))_{I\in\mathcal I}
\end{equation*}
such that $v_{IJ}\circ m_I=p_{IJ*}m_J\circ u_{IJ}$ in $C_{fil}(\Var(\mathbb C)/S_I)$.
\end{itemize}
\item Then $C_{fil}(\Var(\mathbb C)^{sm}/S_{\bullet})$ is the category  
\begin{itemize}
\item whose objects $(G,F)=((G_I,F)_{I\in\mathcal I},u_{IJ})$,
with $(G_I,F)\in C_{fil}(\Var(\mathbb C)^{sm}/S_I)$,
and $u_{IJ}:(G_I,F)\to r_{IJ*}(G_J,F)$ for $r_{IJ}:I\to J$, denoting again $r_{IJ}:S_I\to S_J$, are morphisms
satisfying for $I\to J\to K$, $r_{IJ*}u_{JK}\circ u_{IJ}=u_{IK}$ in $C_{fil}(\Var(\mathbb C)^{sm}/S_I)$,
\item the morphisms $m:((G,F),u_{IJ})\to((H,F),v_{IJ})$ being (see section 2.1) a family of morphisms of complexes,  
\begin{equation*}
m=(m_I:(G_I,F)\to (H_I,F))_{I\in\mathcal I}
\end{equation*}
such that $v_{IJ}\circ m_I=p_{IJ*}m_J\circ u_{IJ}$ in $C_{fil}(\Var(\mathbb C)^{sm}/S_I)$.
\end{itemize}
\end{itemize}
As usual, we denote by
\begin{eqnarray*}
(f_{\bullet}^*,f_{\bullet*}):=(P(f_{\bullet})^*,P(f_{\bullet})_*): 
C(\Var(\mathbb C)^{(sm)}/S_{\bullet})\to C(\Var(\mathbb C)^{(sm)}/T_{\bullet})
\end{eqnarray*}
the adjonction induced by 
$P(f_{\bullet}):\Var(\mathbb C)^{(sm)}/T_{\bullet}\to \Var(\mathbb C)^{(sm)}/S_{\bullet}$.
Since the colimits involved in the definition of $f_{\bullet}^*=P(f_{\bullet})^*$ are filtered, 
$f_{\bullet}^*$ also preserve monomorphism. Hence, we get an adjonction
\begin{eqnarray*}
(f_{\bullet}^*,f_{\bullet*}):
C_{fil}(\Var(\mathbb C)^{(sm)}/S_{\bullet})\leftrightarrows C_{fil}(\Var(\mathbb C)^{(sm)}/T_{\bullet}), \\
f_{\bullet}^*((G_I,F),u_{IJ}):=((f_I^*G_I,f_I^*F),T(f_I,r_{IJ})(-)\circ f_I^*u_{IJ}).
\end{eqnarray*}

Let $S\in\Var(\mathbb C)$. Let $S=\cup_{i=1}^l S_i$ an open affine cover and denote by $S_I=\cap_{i\in I} S_i$.
Let $i_i:S_i\hookrightarrow\tilde S_i$ closed embeddings, with $\tilde S_i\in\Var(\mathbb C)$. 
For $I\subset\left[1,\cdots l\right]$, denote by $\tilde S_I=\Pi_{i\in I}\tilde S_i$.
We then have closed embeddings $i_I:S_I\hookrightarrow\tilde S_I$ and for $J\subset I$ the following commutative diagram
\begin{equation*}
D_{IJ}=\xymatrix{ S_I\ar[r]^{i_I} & \tilde S_I \\
S_J\ar[u]^{j_{IJ}}\ar[r]^{i_J} & \tilde S_J\ar[u]^{p_{IJ}}}  
\end{equation*}
where $p_{IJ}:\tilde S_J\to\tilde S_I$ is the projection
and $j_{IJ}:S_J\hookrightarrow S_I$ is the open embedding so that $j_I\circ j_{IJ}=j_J$.
This gives the diagram of algebraic varieties $(\tilde S_I)\in\Fun(\mathcal P(\mathbb N),\Var(\mathbb C))$ which
the diagram of sites $\Var(\mathbb C)^{sm}/(\tilde S_I)\in\Fun(\mathcal P(\mathbb N),\Cat)$. 
Denote by $m:\tilde S_I\backslash(S_I\backslash S_J)\hookrightarrow\tilde S_I$ the open embedding.
Then $C_{fil}(\Var(\mathbb C)^{sm}/(\tilde S_I)$ is the category  
\begin{itemize}
\item whose objects $(G,F)=((G_I,F),u_{IJ})$ with $(G_I,F)\in C_{fil}(\Var(\mathbb C)^{sm}/\tilde S_I)$,
and $u_{IJ}:(G_I,F)\to p_{IJ*}(G_J,F)$ are morphisms
satisfying for $I\subset J\subset K$, $p_{IJ*}u_{JK}\circ u_{IJ}=u_{IK}$ in $C_{fil}(\Var(\mathbb C)^{sm}/\tilde S_I)$,
\item the morphisms $m:((G,F),u_{IJ})\to((H,F),v_{IJ})$ being a family of morphisms of complexes,  
\begin{equation*}
m=(m_I:(G_I,F)\to (H_I,F))_{I\in\mathcal I}
\end{equation*}
such that $v_{IJ}\circ m_I=p_{IJ*}m_J\circ u_{IJ}$ in $C_{fil}(\Var(\mathbb C)^{sm}/\tilde S_I)$.
\end{itemize}

\begin{defi}
Let $S\in\Var(\mathbb C)$. Let $S=\cup_{i=1}^l S_i$ an open cover and denote by $S_I=\cap_{i\in I} S_i$.
Let $i_i:S_i\hookrightarrow\tilde S_i$ closed embeddings, with $\tilde S_i\in\Var(\mathbb C)$. 
We will denote by 
$C_{fil}(\Var(\mathbb C)^{sm}/(S/(\tilde S_I)))\subset C_{fil}(\Var(\mathbb C)^{sm}/(\tilde S_I))$ 
the full subcategory  
\begin{itemize}
\item whose objects $(G,F)=((G_I,F)_{I\subset\left[1,\cdots l\right]},u_{IJ})$,
with $(G_I,F)\in C_{fil,S_I}(\Var(\mathbb C)^{sm}/\tilde S_I)$,
and $u_{IJ}:m^*(G_I,F)\to m^*p_{IJ*}(G_J,F)$ for $I\subset J$, are $\infty$-filtered Zariski local equivalence,
satisfying for $I\subset J\subset K$, $p_{IJ*}u_{JK}\circ u_{IJ}=u_{IK}$ in $C_{fil}(\Var(\mathbb C)^{sm}/\tilde S_I)$,
\item the morphisms $m:((G,F),u_{IJ})\to((H,F),v_{IJ})$ being (see section 2.1) a family of morphisms of complexes,  
\begin{equation*}
m=(m_I:(G_I,F)\to (H_I,F))_{I\subset\left[1,\cdots l\right]}
\end{equation*}
such that $v_{IJ}\circ m_I=p_{IJ*}m_J\circ u_{IJ}$ in $C_{fil}(\Var(\mathbb C)^{sm}/\tilde S_I)$.
\end{itemize}
A morphism $m:((G_I,F),u_{IJ})\to((H_I,F),v_{IJ})$ is an $r$-filtered Zariski, resp. etale local, equivalence,
if there exists $\phi_i:((C_{iI},F),u_{iIJ})\to((C_{(i+1)I},F),u_{(i+1)IJ})$, $0\leq i\leq s$, 
with $((C_{iI},F),u_{iIJ})\in C_{fil}(\Var(\mathbb C)^{sm}/(S/(\tilde S_I)))$
$((C_{0I},F),u_{0IJ})=((G_I,F),u_{IJ})$ and $((C_{sI},F),u_{sIJ})=((H_I,F),v_{IJ})$, such that
\begin{equation*}
\phi=\phi_s\circ\cdots\circ\phi_i\circ\cdots\circ\phi_0:((G_I,F),u_{IJ})\to((H_I,F),v_{IJ})
\end{equation*}
and $\phi_i:((C_{iI},F),u_{iIJ})\to((C_{(i+1)I},F),u_{(i+1)IJ})$ either a filtered Zariski, resp. etale, local equivalence
or an $r$-filtered homotopy equivalence.  
\end{defi}

Denote $L=[1,\ldots,l]$ and for $I\subset L$,
$p_{0(0I)}:S\times\tilde S_I\to S$, $p_{I(0I)}:S\times\tilde S_I\to S_I$ the projections.
By definition, we have functors 
\begin{itemize}
\item $T(S/(\tilde S_I)):C_{fil}(\Var(\mathbb C)^{sm}/S)\to C_{fil}(\Var(\mathbb C)^{sm}/(S/(\tilde S_I)))$,
$(G,F)\mapsto (i_{I*}j_I^*(G,F),T(D_{IJ})(j_I^*(G,F)))$,
\item $T((\tilde S_I)/S):C_{fil}(\Var(\mathbb C)^{sm}/(S/(\tilde S_I)))\to C_{fil}(\Var(\mathbb C)^{sm}/S)$,
$((G_I,F),u_{IJ})\mapsto \ho\lim_{I\subset L}p_{0(0I)*}\Gamma^{\vee}_{S_I}p_{I(0I)}^*(G_I,F)$. 
\end{itemize}
Note that the functors $T(S/(\tilde S_I)$ are NOT embedding, since
\begin{equation*}
\ad(i_I^*,i_{I*})(j_I^*F):i_I^*i_{I*}j_I^*F\to j_I^*F
\end{equation*}
are Zariski local equivalence but NOT isomorphism since we are dealing
with the morphism of big sites $P(i_I):\Var(\mathbb C)^{sm}/S_I\to\Var(\mathbb C)^{sm}/\tilde S_I$.
However, by theorem \ref{2functDM}, these functors induce full embeddings 
\begin{equation*}
T(S/(\tilde S_I)):D_{fil}(\Var(\mathbb C)^{sm}/S)\to D_{fil}(\Var(\mathbb C)^{sm}/(S/(\tilde S_I)))
\end{equation*}
since for $F\in C(\Var(\mathbb C)^{sm}/S)$, 
\begin{equation*}
\ho\lim_{I\subset L}p_{0(0I)*}\Gamma^{\vee}_{S_I}p_{I(0I)}^*(i_{I*}j_I^*F)\to p_{0(0I)*}\Gamma^{\vee}_{S_I}j_I^*F
\end{equation*}
is an equivalence Zariski local.

Let $f:X\to S$ a morphism, with $X,S\in\Var(\mathbb C)$.
Let $S=\cup_{i=1}^{l} S_i$ and $X=\cup_{i=1}^l X_i$ be affine open covers
and $i_i:S_i\hookrightarrow\tilde S_i$, $i'_i:X_i\hookrightarrow\tilde X_i$ be closed embeddings. 
Let $\tilde f_i:\tilde X_i\to\tilde S_i$ be a lift of the morphism $f_i=f_{|X_i}:X_i\to S_i$.
Then, $f_I=f_{|X_I}:X_I=\cap_{i\in I} X_i\to S_I=\cap_{i\in I} S_i$ lift to the morphism 
\begin{equation*}
\tilde f_I=\Pi_{i\in I}\tilde f_i:\tilde X_I=\Pi_{i\in I}\tilde X_i\to \tilde S_I=\Pi_{i\in I}\tilde S_i
\end{equation*}
Denote by $p_{IJ}:\tilde S_J\to\tilde S_I$ and $p'_{IJ}:\tilde X_J\to\tilde X_I$ the projections.  
Consider for $J\subset I$ the following commutative diagrams
\begin{equation*}
D_{IJ}=\xymatrix{ S_I\ar[r]^{i_I} & \tilde S_I \\
S_J\ar[u]^{j_{IJ}}\ar[r]^{i_J} & \tilde S_J\ar[u]^{p_{IJ}}} \;, \;  
D'_{IJ}=\xymatrix{ X_I\ar[r]^{i'_I} & \tilde X_I \\
X_J\ar[u]^{j'_{IJ}}\ar[r]^{i'_J} & \tilde X_J\ar[u]^{p'_{IJ}}} \; , \;  
D_{fI}=\xymatrix{ S_I\ar[r]^{i_I} & \tilde S_I \\
X_I\ar[u]^{f_I}\ar[r]^{i'_I} & \tilde X_I\ar[u]^{\tilde f_I}}  
\end{equation*}
We have then following commutative diagram
\begin{equation*}
\xymatrix{ \, & X_I\ar[r]^{n'_I} & \tilde X_I & \tilde X_I\backslash X_I\ar[l]^{n'_I} \\
i'_J:X_J\ar[r]^{l_{IJ}}\ar[ru]^{j'_{IJ}}& X_I\times X_I\times\tilde X_{J\backslash I}\ar[u]^{p'_{IJ}}\ar[r]^{n'_I\times I} & 
\tilde X_J\ar[u]^{p'_{IJ}} & \tilde X_J\backslash X_J\ar[l]^{n'_J}\ar[u]^{p'_{IJ}}}.
\end{equation*}
whose square are cartesian.
We then have the pullback functor
\begin{eqnarray*}
f^*:C_{(2)fil}(\Var(\mathbb C)^{sm}/S/(\tilde S_I))\to C_{(2)fil}(\Var(\mathbb C)^{sm}/X/(\tilde X_I)), \\
((G_I,F),u_{IJ})\mapsto f^*((G_I,F),u_{IJ}):=(\Gamma^{\vee}_{X_I}\tilde f_I^*(G_I,F),\tilde f_J^*u_{IJ})
\end{eqnarray*}
with 
\begin{eqnarray*}
\tilde f_J^*u_{IJ}:\Gamma^{\vee}_{X_I}\tilde f_I^*(G_I,F)
\xrightarrow{\ad(p_{IJ}^{'*},p'_{IJ*})(-)}p'_{IJ*}p_{IJ}^{'*}\Gamma^{\vee}_{X_I}\tilde f_I^*(G_I,F) 
\xrightarrow{T_{\sharp}(p_{IJ},n'_I)(-)^{-1})}
p'_{IJ*}\Gamma^{\vee}_{X_I\times\tilde X_{J\backslash I}}p_{IJ}^{'*}\tilde f_I^*(G_I,F) \\
\xrightarrow{p'_{IJ*}\gamma^{\vee}_{X_J}(-)}
p'_{IJ*}\Gamma^{\vee}_{X_J}p_{IJ}^{'*}\tilde f_I^*(G_I,F)=p'_{IJ*}\Gamma^{\vee}_{X_J}\tilde f_J^*p_{IJ}^*(G_I,F)
\xrightarrow{\Gamma^{\vee}_{X_J}\tilde f_J^*I(p_{IJ}^*,p_{IJ*})(-,-)(u_{IJ})}
\Gamma^{\vee}_{X_J}\tilde f_J^*(G_J,F)
\end{eqnarray*}
Let $(G,F)\in C_{fil}(\Var(\mathbb C)^{sm}/S)$. Since, $j_I^{'*}i'_{I*}j_I^{'*}f^*(G,F)=0$,
the morphism $T(D_{fI})(j_I^*(G,F)):\tilde f_I^*i_{I*}j_I^*(G,F)\to i'_{I*}j_I^{'*}f^*(G,F)$ factors trough
\begin{equation*}
T(D_{fI})(j_I^*(G,F)):\tilde f_I^*i_{I*}j_I^*(G,F)\xrightarrow{\gamma_{X_I}^{\vee}(-)}
\Gamma_{X_I}^{\vee}\tilde f_I^*i_{I*}j_I^*(G,F)\xrightarrow{T^{\gamma}(D_{fI})(j_I^*(G,F))} i'_{I*}j_I^{'*}f^*(G,F)
\end{equation*}
We have then, for $(G,F)\in C_{fil}(S)$, the canonical transformation map
\begin{equation*}
\xymatrix{
f^*T(S/(\tilde S_I))(G,F)\ar[rrr]^{T(f,T(0/I))(G,F)}\ar[d]_{=}\ & \, & \, & 
T(X/(\tilde X_I))(f^*(G,F))\ar[d]^{=} \\
(\Gamma^{\vee}_{X_I}\tilde f_I^*i_{I*}j_I^*(G,F),\tilde f_J^*I)\ar[rrr]^{T^{\gamma}(D_{fI})(j_I^*(G,F))} 
& \, & \, & (i'_{I*}j_I^{'*}f^*(G,F),I)}
\end{equation*}

To show that the cohomology sheaves of the filtered De Rham realization functor of constructible motives are mixed hodge modules, 
we will need to take presheaves of the following form

\begin{defi}\label{projBMmotdef}
\begin{itemize}
\item[(i)]Let $f:X\to S$ a morphism with $X,S\in\Var(\mathbb C)$. Assume that there exist a factorization
$f:X\xrightarrow{i} Y\times S\xrightarrow{p} S$, with $Y\in\SmVar(\mathbb C)$, 
$i:X\hookrightarrow Y$ is a closed embedding and $p$ the projection. We then consider 
\begin{equation*}
Q(X/S):=p_{\sharp}\Gamma_X^{\vee}\mathbb Z_{Y\times S}\in C(\Var(\mathbb C)^{sm}/S).
\end{equation*}
By definition $Q(X/S)$ is projective.
\item[(ii)]Let $f:X\to S$ and $g:T\to S$ two morphism with $X,S,T\in\Var(\mathbb C)$. Assume that there exist a factorization
$f:X\xrightarrow{i} Y\times S\xrightarrow{p} S$, with $Y\in\SmVar(\mathbb C)$, 
$i:X\hookrightarrow Y$ is a closed embedding and $p$ the projection. 
We then have the following commutative diagram whose squares are cartesian
\begin{equation*}
\xymatrix{
f:X\ar[r]^{i} & Y\times S\ar[r]^{p} & S \\
f':X_T\ar[r]^{i'}\ar[u]^{g'} & Y\times T\ar[r]^{p'}\ar[u]^{g'':=(I\times g)} & T\ar[u]^{g}}
\end{equation*}
We then have the canonical isomorphism in $C(\Var(\mathbb C)^{sm}/T)$
\begin{eqnarray*}
T(f,g,Q):g^*Q(X/S):=g^*p_{\sharp}\Gamma_X^{\vee}\mathbb Z_{Y\times S}
\xrightarrow{T_{\sharp}(g,p)(-)^{-1}}p'_{\sharp}g^{''*}\Gamma_X^{\vee}\mathbb Z_{Y\times S} \\
\xrightarrow{p'_{\sharp}T(g'',\gamma^{\vee})(-)^{-1}}p'_{\sharp}\Gamma_{X_T}^{\vee}\mathbb Z_{Y\times T}=:Q(X_T/T).
\end{eqnarray*}
\item[(iii)]Let $f:X\to S$ a morphism with $X,S\in\Var(\mathbb C)$. Assume that there exist a factorization
$f:X\xrightarrow{i} Y\times S\xrightarrow{p} S$, with $Y\in\SmVar(\mathbb C)$, 
$i:X\hookrightarrow Y$ is a closed embedding and $p$ the projection. We then consider 
\begin{equation*}
Q^h(X/S):=p_*\Gamma_XE_{et}(\mathbb Z_{Y\times S})\in C(\Var(\mathbb C)^{sm}/S).
\end{equation*}
\item[(iv)]Let $f:X\to S$ and $g:T\to S$ two morphism with $X,S,T\in\Var(\mathbb C)$. Assume that there exist a factorization
$f:X\xrightarrow{i} Y\times S\xrightarrow{p} S$, with $Y\in\SmVar(\mathbb C)$, 
$i:X\hookrightarrow Y$ is a closed embedding and $p$ the projection. 
We then have the following commutative diagram whose squares are cartesian
\begin{equation*}
\xymatrix{
f:X\ar[r]^{i} & Y\times S\ar[r]^{p} & S \\
f':X_T\ar[r]^{i'}\ar[u]^{g'} & Y\times T\ar[r]^{p'}\ar[u]^{g'':=(I\times g)} & T\ar[u]^{g}}
\end{equation*}
We then have the canonical morphism in $C(\Var(\mathbb C)^{sm}/T)$
\begin{eqnarray*}
T(f,g,Q^h):g^*Q^h(X/S):=g^*p_*\Gamma_XE_{et}(\mathbb Z_{Y\times S})
\xrightarrow{T(g,p)(-)}p'_*g^{''*}\Gamma_XE_{et}(\mathbb Z_{Y\times S}) \\
\xrightarrow{p'_*T(g'',\gamma)(-)}p'_*\Gamma_{X_T}E_{et}(\mathbb Z_{Y\times T})=:Q^h(X_T/T).
\end{eqnarray*}
\end{itemize}
\end{defi}

We now give the definition of the $\mathbb A^1$ local property :

Denote by
\begin{eqnarray*}
p_a:\Var(\mathbb C)^{(sm)}/S\to\Var(\mathbb C)^{(sm)}/S, \;  
X/S=(X,h)\mapsto (X\times\mathbb A^1)/S=(X\times\mathbb A^1,h\circ p_X), \\ 
(g:X/S\to X'/S)\mapsto ((g\times I_{\mathbb A^1}):X\times\mathbb A^1/S\to X'\times\mathbb A^1/S)
\end{eqnarray*}
the projection functor and again by $p_a:\Var(\mathbb C)^{(sm)}/S\to\Var(\mathbb C)^{(sm)}/S$
the corresponding morphism of site.

\begin{defi}\label{a1locdef}
Let $S\in\Var(\mathbb C)$. Denote for short $\Var(\mathbb C)^{(sm)}/S$ 
either the category $\Var(\mathbb C)/S$ or the category $\Var(\mathbb C)^{sm}/S$.
\begin{itemize}
\item[(i0)] A complex $F\in C(\Var(\mathbb C)^{(sm)}/S)$ is said to be $\mathbb A^1$ homotopic if
$\ad(p_a^*,p_{a*})(F):F\to p_{a*}p_a^*F$ is an homotopy equivalence.
\item[(i)] A complex $F\in C(\Var(\mathbb C)^{(sm)}/S)$ is said to be $\mathbb A^1$ invariant 
if for all $U/S\in\Var(\mathbb C)^{(sm)}/S$,
\begin{equation*}
F(p_U):F(U/S)\to F(U\times\mathbb A^1/S) 
\end{equation*}
is a quasi-isomorphism, where $p_U:U\times\mathbb A^1\to U$ is the projection.
Obviously, if a complex $F\in C(\Var(\mathbb C)^{(sm)}/S)$ is $\mathbb A^1$ homotopic then it is $\mathbb A^1$ invariant.
\item[(ii)] Let $\tau$ a topology on $\Var(\mathbb C)$. 
A complex $F\in C(\Var(\mathbb C)^{(sm)}/S)$ is said to be $\mathbb A^1$ local for
the topology $\tau$, if for a (hence every) $\tau$ local equivalence $k:F\to G$ with $k$ injective and
$G\in C(\Var(\mathbb C)^{(sm)}/S)$ $\tau$ fibrant, e.g. $k:F\to E_{\tau}(F)$, 
$G$ is $\mathbb A^1$ invariant for all $n\in\mathbb Z$.
\item[(iii)] A morphism $m:F\to G$ with $F,G\in C(\Var(\mathbb C)^{(sm)}/S)$ is said to an $(\mathbb A^1,et)$ local equivalence 
if for all $H\in C(\Var(\mathbb C)^{(sm)}/S)$ which is $\mathbb A^1$ local for the etale topology
\begin{equation*}
\Hom(L(m),E_{et}(H)):\Hom(L(G),E_{et}(H))\to\Hom(L(F),E_{et}(H)) 
\end{equation*}
is a quasi-isomorphism.
\end{itemize}
\end{defi}

Denote $\square^*:=\mathbb P^*\backslash\left\{1\right\}$
\begin{itemize}
\item Let $S\in\Var(\mathbb C)$. For $U/S=(U,h)\in\Var(\mathbb C)^{sm}/S$, we consider 
\begin{equation*}
\square^*\times U/S=(\square^*\times U,h\circ p)\in\Fun(\Delta,\Var(\mathbb C)^{sm}/S).
\end{equation*}
For $F\in C^-(\Var(\mathbb C)^{sm}/S)$, it gives the complex
\begin{equation*}
C_*F\in C^-(\Var(\mathbb C)^{sm}/S), U/S=(U,h)\mapsto C_*F(U/S):=\Tot F(\square^*\times U/S)
\end{equation*}
together with the canonical map $c_F:=(0,I_F):F\to C_*F$.
For $F\in C(\Var(\mathbb C)^{sm}/S)$, we get
\begin{equation*}
C_*F:=\holim_n C_*F^{\leq n}\in C(\Var(\mathbb C)^{sm}/S),
\end{equation*}
together with the canonical map $c_F:=(0,I_F):F\to C_*F$.
For $m:F\to G$ a morphism, with $F,G\in C(\Var(\mathbb C)^{sm}/S)$,
we get by functoriality the morphism $C_*m:C_*F\to C_*G$.
\item Let $S\in\Var(\mathbb C)$. For $U/S=(U,h)\in\Var(\mathbb C)/S$, we consider 
\begin{equation*}
\square^*\times U/S=(\mathbb A^*\times U,h\circ p)\in\Fun(\Delta,\Var(\mathbb C)/S).
\end{equation*}
For $F\in C^-(\Var(\mathbb C)/S)$, it gives the complex
\begin{equation*}
C_*F\in C^-(\Var(\mathbb C)/S), U/S=(U,h)\mapsto C_*F(U/S):=\Tot F(\square^*\times U/S)
\end{equation*}
together with the canonical map $c=c(F):=(0,I_F):F\to C_*F$.
For $F\in C(\Var(\mathbb C)/S)$, we get
\begin{equation*}
C_*F:=\holim_n C_*F^{\leq n}\in C(\Var(\mathbb C)/S),
\end{equation*}
together with the canonical map $c_F:=(0,I_F):F\to C_*F$.
For $m:F\to G$ a morphism, with $F,G\in C(\Var(\mathbb C)/S)$,
we get by functoriality the morphism $C_*m:C_*F\to C_*G$.
\end{itemize}

\begin{prop}\label{ca1Var}
\begin{itemize}
\item[(i)]Let $S\in\Var(\mathbb C)$.
Then for $F\in C(\Var(\mathbb C)^{sm}/S)$, $C_*F$ is $\mathbb A^1$ local for the etale topology
and  $c(F):F\to C_*F$ is an equivalence $(\mathbb A^1,et)$ local.
\item[(ii)]A morphism $m:F\to G$ with $F,G\in C(\Var(\mathbb C)^{(sm)}/S)$ is an $(\mathbb A^1,et)$ local equivalence
if and only if there exists 
\begin{equation*}
\left\{X_{1,\alpha}/S,\alpha\in\Lambda_1\right\},\ldots,\left\{X_{r,\alpha}/S,\alpha\in\Lambda_r\right\}
\subset\Var(\mathbb C)^{(sm)}/S 
\end{equation*}
such that we have in $\Ho_{et}(C(\Var(\mathbb C)^{(sm)}/S))$
\begin{eqnarray*}
\Cone(m)\xrightarrow{\sim}
\Cone(\oplus_{\alpha\in\Lambda_1}\Cone(\mathbb Z(X_{1,\alpha}\times\mathbb A^1/S)\to\mathbb Z(X_{1,\alpha}/S)) \\
\to\cdots\to\oplus_{\alpha\in\Lambda_r}\Cone(\mathbb Z(X_{r,\alpha}\times\mathbb A^1/S)\to\mathbb Z(X_{r,\alpha}/S)))
\end{eqnarray*}
\end{itemize}
\end{prop}

\begin{proof}
Standard : see Ayoub's thesis for example.
\end{proof}

\begin{defiprop}\label{projmodstr}
Let $S\in\Var(\mathbb C)$.
\begin{itemize}
\item[(i)]With the weak equivalence the $(\mathbb A^1,et)$ local equivalence and 
the fibration the epimorphism with $\mathbb A^1_S$ local and etale fibrant kernels gives
a model structure on  $C(\Var(\mathbb C)^{sm}/S)$ : the left bousfield localization
of the projective model structure of $C(\Var(\mathbb C)^{sm}/S)$. 
We call it the projective $(\mathbb A^1,et)$ model structure.
\item[(ii)]With the weak equivalence the $(\mathbb A^1,et)$ local equivalence and 
the fibration the epimorphism with $\mathbb A^1_S$ local and etale fibrant kernels gives
a model structure on  $C(\Var(\mathbb C)/S)$ : the left bousfield localization
of the projective model structure of $C(\Var(\mathbb C)/S)$. 
We call it the projective $(\mathbb A^1,et)$ model structure.
\end{itemize}
\end{defiprop}

\begin{proof}
See \cite{C.D}.
\end{proof}

\begin{prop}\label{g1}
Let $g:T\to S$ a morphism with $T,S\in\Var(\mathbb C)$.
\begin{itemize}
\item[(i)] The adjonction $(g^*,g_*):C(\Var(\mathbb C)^{sm}/S)\leftrightarrows C(\Var(\mathbb C)^{sm}/T)$
is a Quillen adjonction for the $(\mathbb A^1,et)$ projective model structure (see definition-proposition \ref{projmodstr}).
\item[(i)'] Let $h:U\to S$ a smooth morphism with $U,S\in\Var(\mathbb C)$.
The adjonction $(h_{\sharp},h^*):C(\Var(\mathbb C)^{sm}/U)\leftrightarrows C(\Var(\mathbb C)^{sm}/S)$
is a Quillen adjonction for the $(\mathbb A^1,et)$ projective model structure.
\item[(i)''] The functor $g^*:C(\Var(\mathbb C)^{sm}/S)\to C(\Var(\mathbb C)^{sm}/T)$
sends quasi-isomorphism to quasi-isomorphism,
sends equivalence Zariski local to equivalence Zariski local, and equivalence etale local to equivalence etale local,
sends $(\mathbb A^1,et)$ local equivalence to $(\mathbb A^1,et)$ local equivalence.
\item[(ii)] The adjonction $(g^*,g_*):C(\Var(\mathbb C)/S)\leftrightarrows C(\Var(\mathbb C)/T)$
is a Quillen adjonction for the $(\mathbb A^1,et)$ projective model structure (see definition-proposition \ref{projmodstr}).
\item[(ii)'] The adjonction $(g_{\sharp},g^*):C(\Var(\mathbb C)/T)\leftrightarrows C(\Var(\mathbb C)/S)$
is a Quillen adjonction for the $(\mathbb A^1,et)$ projective model structure.
\item[(ii)''] The functor $g^*:C(\Var(\mathbb C)/S)\to C(\Var(\mathbb C)/T)$
sends quasi-isomorphism to quasi-isomorphism,
sends equivalence Zariski local to equivalence Zariski local, and equivalence etale local to equivalence etale local,
sends $(\mathbb A^1,et)$ local equivalence to $(\mathbb A^1,et)$ local equivalence.
\end{itemize}
\end{prop}

\begin{proof}
Standard : see \cite{C.D} for example.
\end{proof}

\begin{prop}\label{rho1}
Let $S\in\Var(\mathbb C)$. 
\begin{itemize}
\item[(i)] The adjonction $(\rho_S^*,\rho_{S*}):C(\Var(\mathbb C)^{sm}/S)\leftrightarrows C(\Var(\mathbb C)/S)$
is a Quillen adjonction for the $(\mathbb A^1,et)$ projective model structure.
\item[(ii)]The functor $\rho_{S*}:C(\Var(\mathbb C)/S)\to C(\Var(\mathbb C)^{sm}/S)$
sends quasi-isomorphism to quasi-isomorphism,
sends equivalence Zariski local to equivalence Zariski local, and equivalence etale local to equivalence etale local,
sends $(\mathbb A^1,et)$ local equivalence to $(\mathbb A^1,et)$ local equivalence.
\end{itemize}
\end{prop}

\begin{proof}
Standard : see \cite{C.D} for example.
\end{proof}

Let $S\in\Var(\mathbb C)$. Let $S=\cup_{i=1}^l S_i$ an open affine cover and denote by $S_I=\cap_{i\in I} S_i$.
Let $i_i:S_i\hookrightarrow\tilde S_i$ closed embeddings, with $\tilde S_i\in\Var(\mathbb C)$.
For $(G_I,K_{IJ})\in C(\Var(\mathbb C)^{(sm)}/(\tilde S_I)^{op})$ and 
$(H_I,T_{IJ})\in C(\Var(\mathbb C)^{(sm)}/(\tilde S_I))$, we denote
\begin{eqnarray*}
\mathcal Hom((G_I,K_{IJ}),(H_I,T_{IJ})):=(\mathcal Hom(G_I,H_I),u_{IJ}((G_I,K_{IJ}),(H_I,T_{IJ})))
\in C(\Var(\mathbb C)^{(sm)}/(\tilde S_I))
\end{eqnarray*}
with
\begin{eqnarray*}
u_{IJ}((G_I,K_{IJ})(H_I,T_{IJ})):\mathcal Hom(G_I,H_I) \\
\xrightarrow{\ad(p_{IJ}^*,p_{IJ*})(-)}p_{IJ*}p_{IJ}^*\mathcal Hom(G_I,H_I)
\xrightarrow{T(p_{IJ},hom)(-,-)}p_{IJ*}\mathcal Hom(p_{IJ}^*G_I,p_{IJ}^*H_I) \\
\xrightarrow{\mathcal Hom(p_{IJ}^*G_I,T_{IJ})}p_{IJ*}\mathcal Hom(p_{IJ}^*G_I,H_J)
\xrightarrow{\mathcal Hom(K_{IJ},H_J)}p_{IJ*}\mathcal Hom(G_J,H_J).
\end{eqnarray*}
This gives in particular the functor
\begin{eqnarray*}
C(\Var(\mathbb C)^{(sm)}/(\tilde S_I))\to C(\Var(\mathbb C)^{(sm)}/(\tilde S_I)^{op}),
(H_I,T_{IJ})\mapsto(H_I,T_{IJ}).
\end{eqnarray*}
Let $S\in\Var(\mathbb C)$. Let $S=\cup_{i=1}^l S_i$ an open affine cover and denote by $S_I=\cap_{i\in I} S_i$.
Let $i_i:S_i\hookrightarrow\tilde S_i$ closed embeddings, with $\tilde S_i\in\SmVar(\mathbb C)$.
The functor $p_a$ extend to a functor
\begin{eqnarray*}
p_a:\Var(\mathbb C)^{(sm)}/(\tilde S_I)\to\Var(\mathbb C)^{(sm)}/(\tilde S_I), \;  
(X_I/\tilde S_I,s_{IJ})\mapsto (X\times\mathbb A^1/\tilde S_I,s_{IJ}\times I), \\ 
(g=(g_I):(X_I/\tilde S_I,s_{IJ})\to (X'_I/\tilde S_I,s'_{IJ}))\mapsto 
((g\times I_{\mathbb A^1}):(X_I\times\mathbb A^1/\tilde S_I,s_{IJ})\to(X'_I\times\mathbb A^1/\tilde S_I,s'_{IJ}))
\end{eqnarray*}
the projection functor and again by $p_a:\Var(\mathbb C)^{(sm)}/(\tilde S_I)\to\Var(\mathbb C)^{(sm)}/(\tilde S_I)$
the corresponding morphism of site.

\begin{defi}\label{a1locdefIJ}
Let $S\in\Var(\mathbb C)$. Let $S=\cup_{i=1}^l S_i$ an open affine cover and denote by $S_I=\cap_{i\in I} S_i$.
Let $i_i:S_i\hookrightarrow\tilde S_i$ closed embeddings, with $\tilde S_i\in\Var(\mathbb C)$.
\begin{itemize}
\item[(i0)]A complex $(F_I,u_{IJ})\in C(\Var(\mathbb C)^{(sm)}/(\tilde S_I))$ is said to be $\mathbb A^1$ homotopic if 
$\ad(p_a^*,p_{a*})((F_I,u_{IJ})):(F_I,u_{IJ})\to p_{a*}p_a^*(F_I,u_{IJ})$ is an homotopy equivalence.
\item[(i)] A complex  $(F_I,u_{IJ})\in C(\Var(\mathbb C)^{(sm)}/(\tilde S_I))$ is said to be $\mathbb A^1$ invariant 
if for all $(X_I/\tilde S_I,s_{IJ})\in\Var(\mathbb C)^{(sm)}/(\tilde S_I)$ 
\begin{equation*}
(F_I(p_{X_I})):(F_I(X_I/\tilde S_I),F_J(s_{IJ}))\circ u_{IJ}(-)\to 
(F_I(X_I\times\mathbb A^1/\tilde S_I),F_J(s_{IJ}\times I)\circ u_{IJ}(-)) 
\end{equation*}
is a quasi-isomorphism, where $p_{X_I}:X_I\times\mathbb A^1\to X_I$ are the projection,
and $s_{IJ}:X_I\times\tilde S_{J\backslash I}/\tilde S_J\to X_J/\tilde S_J$.
Obviously a complex $(F_I,u_{IJ})\in C(\Var(\mathbb C)^{(sm)}/(\tilde S_I))$ is $\mathbb A^1$ invariant
if and only if all the $F_I$ are $\mathbb A^1$ invariant.
\item[(ii)]Let $\tau$ a topology on $\Var(\mathbb C)$. 
A complex $F=(F_I,u_{IJ})\in C(\Var(\mathbb C)^{(sm)}/(\tilde S_I))$ is said to be $\mathbb A^1$ local 
for the $\tau$ topology induced on $\Var(\mathbb C)/(\tilde S_I)$, 
if for an (hence every) $\tau$ local equivalence $k:F\to G$ with $k$ injective and 
$G=(G_I,v_{IJ})\in C(\Var(\mathbb C)^{(sm)}/(\tilde S_I))$ $\tau$ fibrant,
e.g. $k:(F_I,u_{IJ})\to (E_{\tau}(F_I),E(u_{IJ}))$, $G$ is $\mathbb A^1$ invariant.
\item[(iii)] A morphism $m=(m_I):(F_I,u_{IJ})\to (G_I,v_{IJ})$ with 
$(F_I,u_{IJ}),(G_I,v_{IJ})\in C(\Var(\mathbb C)^{(sm)}/(\tilde S_I))$ 
is said to be an $(\mathbb A^1,et)$ local equivalence 
if for all $(H_I,w_{IJ})\in C(\Var(\mathbb C)^{(sm)}/(\tilde S_I))$ which is $\mathbb A^1$ local for the etale topology
\begin{eqnarray*}
(\Hom(L(m_I),E_{et}(H_I))):\Hom(L(G_I,v_{IJ}),E_{et}(H_I,w_{IJ}))\to\Hom(L(F_I,u_{IJ}),E_{et}(H_I,w_{IJ})) 
\end{eqnarray*}
is a quasi-isomorphism.
Obviously, if a morphism $m=(m_I):(F_I,u_{IJ})\to (G_I,v_{IJ})$ with 
$(F_I,u_{IJ}),(G_I,u_{IJ})\in C(\Var(\mathbb C)^{(sm)}/(\tilde S_I)^{op})$ 
is an $(\mathbb A^1,et)$ local equivalence, 
then all the $m_I:F_I\to G_I$ are $(\mathbb A^1,et)$ local equivalence.
\item[(iv)] A morphism $m=(m_I):(F_I,u_{IJ})\to (G_I,v_{IJ})$ with 
$(F_I,u_{IJ}),(G_I,v_{IJ})\in C(\Var(\mathbb C)^{(sm)}/(\tilde S_I)^{op})$ 
is said to be an $(\mathbb A^1,et)$ local equivalence 
if for all $(H_I,w_{IJ})\in C(\Var(\mathbb C)^{(sm)}/(\tilde S_I))$ which is $\mathbb A^1$ local for the etale topology
\begin{eqnarray*}
(\Hom(L(m_I),E_{et}(H_I))):\Hom(L(G_I,v_{IJ}),E_{et}(H_I,w_{IJ}))\to\Hom(L(F_I,u_{IJ}),E_{et}(H_I,w_{IJ})) 
\end{eqnarray*}
is a quasi-isomorphism.
Obviously, if a morphism $m=(m_I):(F_I,u_{IJ})\to (G_I,v_{IJ})$ with 
$(F_I,u_{IJ}),(G_I,u_{IJ})\in C(\Var(\mathbb C)^{(sm)}/(\tilde S_I)^{op})$ 
is an $(\mathbb A^1,et)$ local equivalence, 
then all the $m_I:F_I\to G_I$ are $(\mathbb A^1,et)$ local equivalence.
\end{itemize}
\end{defi}

\begin{prop}\label{ca1VarIJ}
Let $S\in\Var(\mathbb C)$. Let $S=\cup_{i=1}^l S_i$ an open affine cover and denote by $S_I=\cap_{i\in I} S_i$.
Let $i_i:S_i\hookrightarrow\tilde S_i$ closed embeddings, with $\tilde S_i\in\Var(\mathbb C)$.
\begin{itemize}
\item[(i)]A morphism $m:F\to G$ with $F,G\in C(\Var(\mathbb C)^{(sm)}/(\tilde S_I))$ 
is an $(\mathbb A^1,et)$ local equivalence if and only if there exists 
\begin{eqnarray*}
\left\{(X_{1,\alpha,I}/\tilde S_I,u^1_{IJ}),\alpha\in\Lambda_1\right\},\ldots,
\left\{(X_{r,\alpha,I}/\tilde S_I,u^r_{IJ}),\alpha\in\Lambda_r\right\}
\subset\Var(\mathbb C)^{(sm)}/(\tilde S_I)
\end{eqnarray*}
with $u^l_{IJ}:X_{l,\alpha,I}\times\tilde S_{J\backslash I}/\tilde S_J\to X_{l,\alpha,J}/\tilde S_J$,
such that we have in $\Ho_{et}(C(\Var(\mathbb C)^{(sm)}/(\tilde S_I)))$
\begin{eqnarray*}
\Cone(m)\xrightarrow{\sim}\Cone( \oplus_{\alpha\in\Lambda_1} 
\Cone((\mathbb Z(X_{1,\alpha,I}\times\mathbb A^1/\tilde S_I),\mathbb Z(u_{IJ}^1\times I)) 
\to(\mathbb Z(X_{1,\alpha,I}/\tilde S_I),\mathbb Z(u_{IJ}^1))) \\
\to\cdots\to\oplus_{\alpha\in\Lambda_r}
\Cone((\mathbb Z(X_{r,\alpha,I}\times\mathbb A^1/\tilde S_I),\mathbb Z(u_{IJ}^r\times I)) 
\to(\mathbb Z(X_{r,\alpha,I}/\tilde S_I),\mathbb Z(u^r_{IJ}))))
\end{eqnarray*}
\item[(ii)]A morphism $m:F\to G$ with $F,G\in C(\Var(\mathbb C)^{(sm)}/(\tilde S_I)^{op})$ 
is an $(\mathbb A^1,et)$ local equivalence if and only if there exists 
\begin{eqnarray*}
\left\{(X_{1,\alpha,I}/\tilde S_I,u^1_{IJ}),\alpha\in\Lambda_1\right\},\ldots,
\left\{(X_{r,\alpha,I}/\tilde S_I,u^r_{IJ}),\alpha\in\Lambda_r\right\}
\subset\Var(\mathbb C)^{(sm)}/(\tilde S_I)^{op}
\end{eqnarray*}
with $u^l_{IJ}:X_{l,\alpha,J}/\tilde S_J\to X_{l,\alpha,I}\times\tilde S_{J\backslash I}/\tilde S_J$,
such that we have in $\Ho_{et}(C(\Var(\mathbb C)^{(sm)}/(\tilde S_I)^{op}))$
\begin{eqnarray*}
\Cone(m)\xrightarrow{\sim}\Cone( \oplus_{\alpha\in\Lambda_1} 
\Cone((\mathbb Z(X_{1,\alpha,I}\times\mathbb A^1/\tilde S_I),\mathbb Z(u_{IJ}^1\times I)) 
\to(\mathbb Z(X_{1,\alpha,I}/\tilde S_I),\mathbb Z(u_{IJ}^1))) \\
\to\cdots\to\oplus_{\alpha\in\Lambda_r}
\Cone((\mathbb Z(X_{r,\alpha,I}\times\mathbb A^1/\tilde S_I),\mathbb Z(u_{IJ}^r\times I)) 
\to(\mathbb Z(X_{r,\alpha,I}/\tilde S_I),\mathbb Z(u^r_{IJ}))))
\end{eqnarray*}
\end{itemize}
\end{prop}

\begin{proof}
Standard. See Ayoub's thesis for example.
\end{proof}

\begin{itemize}
\item For $f:X\to S$ a morphism with $X,S\in\Var(\mathbb C)$, we denote as usual (see \cite{C.D} for example),
$\mathbb Z^{tr}(X/S)\in\PSh(\Var(\mathbb C)/S)$ the presheaf given by
\begin{itemize}
\item for $X'/S\in\Var(\mathbb C)/S$, with $X'$ irreducible, 
$\mathbb Z^{tr}(X/S)(X'/S):=\mathcal Z^{fs/X}(X'\times_S X)\subset\mathcal Z_{d_{X'}}(X'\times_S X)$
which consist of algebraic cycles $\alpha=\sum_in_i\alpha_i\in\mathcal Z_{d_{X'}}(X'\times_S X)$ such that,
denoting $\supp(\alpha)=\cup_i\alpha_i\subset X'\times_S X$ its support and $f':X'\times_S X\to X'$ the projection,
$f'_{|\supp(\alpha)}:\supp(\alpha)\to X'$ is finite surjective, 
\item for $g:X_2/S\to X_1/S$ a morphism, with $X_1/S,X_2/S\in\Var(\mathbb C)/S$,
\begin{equation*}
\mathbb Z^{tr}(X/S)(g):\mathbb Z^{tr}(X/S)(X_1/S)\to\mathbb Z^{tr}(X/S)(X_2/S), \;
\alpha\mapsto (g\times I)^{-1}(\alpha)
\end{equation*}
with $g\times I:X_2\times_S X\to X_1\times_S X$, noting that, by base change,
$f_{2|\supp((g\times I)^{-1}(\alpha))}:\supp((g\times I)^{-1}(\alpha))\to X_2$ is finite surjective,
$f_2:X_2\times_S X\to X_2$ being the projection.
\end{itemize}
\item For $f:X\to S$ a morphism with $X,S\in\Var(\mathbb C)$ and $r\in\mathbb N$, 
we denote as usual (see \cite{C.D} for example),
$\mathbb Z^{equir}(X/S)\in\PSh(\Var(\mathbb C)/S)$ the presheaf given by
\begin{itemize}
\item for $X'/S\in\Var(\mathbb C)/S$, with $X'$ irreducible, 
$\mathbb Z^{equir}(X/S)(X'/S):=\mathcal Z^{equir/X}(X'\times_S X)\subset\mathcal Z_{d_{X'}}(X'\times_S X)$
which consist of algebraic cycles $\alpha=\sum_in_i\alpha_i\in\mathcal Z_{d_{X'}}(X'\times_S X)$ such that,
denoting $\supp(\alpha)=\cup_i\alpha_i\subset X'\times_S X$ its support and $f':X'\times_S X\to X'$ the projection,
$f'_{|\supp(\alpha)}:\supp(\alpha)\to X'$ is dominant, with fibers either empty or of dimension $r$, 
\item for $g:X_2/S\to X_1/S$ a morphism, with $X_1/S,X_2/S\in\Var(\mathbb C)/S$,
\begin{equation*}
\mathbb Z^{equir}(X/S)(g):\mathbb Z^{equir}(X/S)(X_1/S)\to\mathbb Z^{equir}(X/S)(X_2/S), \;
\alpha\mapsto (g\times I)^{-1}(\alpha)
\end{equation*}
with $g\times I:X_2\times_S X\to X_1\times_S X$, noting that, by base change,
$f_{2|\supp((g\times I)^{-1}(\alpha))}:\supp((g\times I)^{-1}(\alpha))\to X_2$ is obviously dominant,
with fibers either empty or of dimension $r$, $f_2:X_2\times_S X\to X_2$ being the projection.
\end{itemize}
\item Let $S\in\Var(\mathbb C)$. We denote by 
$\mathbb Z_S(d):=\mathbb Z^{equi0}(S\times\mathbb A^d/S)[-2d]$
the Tate twist. For $F\in C(\Var(\mathbb C)/S)$, we denote by $F(d):=F\otimes\mathbb Z_S(d)$.
\end{itemize}

For $S\in\Var(\mathbb C)$, let $\Cor(\Var(\mathbb C)^{sm}/S)$ be the category 
\begin{itemize}
\item whose objects are smooth morphisms 
$U/S=(U,h)$, $h:U\to S$ with $U\in\Var(\mathbb C)$, 
\item whose morphisms $\alpha:U/S=(U,h_1)\to V/S=(V,h_2)$ 
is finite correspondence that is $\alpha\in\oplus_i\mathcal Z^{fs}(U_i\times_S V)$, 
where $U=\sqcup_i U_i$, with $U_i$ connected (hence irreducible by smoothness), and $\mathcal Z^{fs}(U_i\times_S V)$ 
is the abelian group of cycle finite and surjective over $U_i$. 
\end{itemize}
We denote by 
$\Tr(S):\Cor(\Var(\mathbb C)^{sm}/S)\to\Var(\mathbb C)^{sm}/S$ 
the morphism of site
given by the inclusion functor
$\Tr(S):\Var(\mathbb C)^{sm}/S\hookrightarrow\Cor(\Var(\mathbb C)^{sm}/S)$
It induces an adjonction
\begin{equation*}
(\Tr(S)^*\Tr(S)_*):C(\Var(\mathbb C)^{sm}/S)\leftrightarrows C(\Cor(\Var(\mathbb C)^{sm}/S))
\end{equation*}
A complex of preheaves $G\in C(\Var(\mathbb C)^{sm}/S)$ is said to admit transferts
if it is in the image of the embedding
\begin{equation*}
\Tr(S)_*:C(\Cor(\Var(\mathbb C)^{sm}/S)\hookrightarrow C(\Var(\mathbb C)^{sm}/S),
\end{equation*}
that is $G=\Tr(S)_*\Tr(S)^*G$.

We will use to compute the algebraic Gauss-Manin realization functor the following

\begin{thm}\label{DDADM}
Let $\phi:F^{\bullet}\to G^{\bullet}$ an etale local equivalence with $F^{\bullet},G^{\bullet}\in C(\Var(\mathbb C)^{sm}/S)$.
If $F^{\bullet}$ and $G^{\bullet}$ are $\mathbb A^1$ local and admit tranferts 
then $\phi:F^{\bullet}\to G^{\bullet}$ is a Zariski local equivalence.
Hence if $F\in C(\Var(\mathbb C)^{sm}/S)$ is $\mathbb A^1$ local and admits transfert 
\begin{equation*}
k:E_{zar}(F)\to E_{et}(E_{zar}(F))=E_{et}(F) 
\end{equation*}
is a Zariski local equivalence.
\end{thm}

\begin{proof}
See \cite{C.D}.
\end{proof}

\subsection{Presheaves on the big Zariski site or the big etale site of pairs}

We recall the definition given in subsection 5.1 :
For $S\in\Var(\mathbb C)$, $\Var(\mathbb C)^2/S:=\Var(\mathbb C)^2/(S,S)$ is by definition (see subsection 2.1)
the category whose set of objects is 
\begin{eqnarray*}
(\Var(\mathbb C)^2/S)^0:=
\left\{((X,Z),h), h:X\to S, \; Z\subset X \; \mbox{closed} \;\right\}\subset\Var(\mathbb C)/S\times\Top
\end{eqnarray*}
and whose set of morphisms between $(X_1,Z_1)/S=((X_1,Z_1),h_1),(X_1,Z_1)/S=((X_2,Z_2),h_2)\in\Var(\mathbb C)^2/S$
is the subset
\begin{eqnarray*}
\Hom_{\Var(\mathbb C)^2/S}((X_1,Z_1)/S,(X_2,Z_2)/S):= \\
\left\{(f:X_2\to X_2), \; \mbox{s.t.} \; h_1\circ f=h_2 \; \mbox{and} \; Z_1\subset f^{-1}(Z_2)\right\} 
\subset\Hom_{\Var(\mathbb C)}(X_1,X_2)
\end{eqnarray*}
The category $\Var(\mathbb C)^2$ admits fiber products : $(X_1,Z_1)\times_{(S,Z)}(X_2,Z_2)=(X_1\times_S X_2,Z_1\times_Z Z_2)$.
In particular, for $f:T\to S$ a morphism with $S,T\in\Var(\mathbb C)$, we have the pullback functor
\begin{equation*}
P(f):\Var(\mathbb C)^2/S\to\Var(\mathbb C)^2/T, P(f)((X,Z)/S):=(X_T,Z_T)/T, P(f)(g):=(g\times_S f)
\end{equation*}
and we note again $P(f):\Var(\mathbb C)^2/T\to\Var(\mathbb C)^2/S$ the corresponding morphism of sites.

We will consider in the construction of the filtered De Rham realization functor the
full subcategory $\Var(\mathbb C)^{2,sm}/S\subset\Var(\mathbb C)^2/S$ such that the first factor is a smooth morphism :
We will also consider, in order to obtain a complex of D modules in the construction of the filtered De Rham realization functor,
the restriction to the full subcategory $\Var(\mathbb C)^{2,pr}/S\subset\Var(\mathbb C)^2/S$ 
such that the first factor is a projection :

\begin{defi}\label{PVar12S}
\begin{itemize}
\item[(i)]Let $S\in\Var(\mathbb C)$. We denote by
\begin{equation*}
\rho_S:\Var(\mathbb C)^{2,sm}/S\hookrightarrow\Var(\mathbb C)^2/S 
\end{equation*}
the full subcategory
consisting of the objects $(U,Z)/S=((U,Z),h)\in\Var(\mathbb C)^2/S$ such that the morphism $h:U\to S$ is smooth.
That is, $\Var(\mathbb C)^{2,sm}/S$ is the category  
\begin{itemize}
\item whose objects are $(U,Z)/S=((U,Z),h)$, with $U\in\Var(\mathbb C)$, $Z\subset U$ a closed subset, 
and $h:U\to S$ a smooth morphism,
\item whose morphisms $g:(U,Z)/S=((U,Z),h_1)\to (U',Z')/S=((U',Z'),h_2)$ 
is a morphism $g:U\to U'$ of complex algebraic varieties such that $Z\subset g^{-1}(Z')$ and  $h_2\circ g=h_1$. 
\end{itemize}
We denote again $\rho_S:\Var(\mathbb C)^2/S\to\Var(\mathbb C)^{2,sm}/S$ the associated morphism of site.We have 
\begin{equation*}
r^s(S):\Var(\mathbb C)^2\xrightarrow{r(S):=r(S,S)}\Var(\mathbb C)^2/S\xrightarrow{\rho_S}\Var(\mathbb C)^{2,sm}/S
\end{equation*}
the composite morphism of site.
\item[(ii)]Let $S\in\Var(\mathbb C)$. We will consider the full subcategory 
\begin{equation*}
\mu_S:\Var(\mathbb C)^{2,pr}/S\hookrightarrow\Var(\mathbb C)^2/S
\end{equation*}
whose subset of object consist of those whose morphism is a projection to $S$ : 
\begin{eqnarray*}
(\Var(\mathbb C)^{2,pr}/S)^0:=\left\{((Y\times S,X),p), \; Y\in\Var(\mathbb C), \;
 p:Y\times S\to S \; \mbox{the projection}\right\}\subset(\Var(\mathbb C)^2/S)^0.
\end{eqnarray*}
\item[(iii)]We will consider the full subcategory 
\begin{equation*}
\mu_S:(\Var(\mathbb C)^{2,smpr}/S)\hookrightarrow\Var(\mathbb C)^{2,sm}/S
\end{equation*}
whose subset of object consist of those whose morphism is a smooth projection to $S$ : 
\begin{eqnarray*}
(\Var(\mathbb C)^{2,smpr}/S)^0:=\left\{((Y\times S,X),p), \; Y\in\SmVar(\mathbb C), \;
 p:Y\times S\to S \; \mbox{the projection}\right\}\subset(\Var(\mathbb C)^2/S)^0
\end{eqnarray*}
\end{itemize}
\end{defi}
For $f:T\to S$ a morphism with $T,S\in\Var(\mathbb C)$, we have by definition, the following commutative diagram of sites
\begin{equation}\label{muf}
\xymatrix{\Var(\mathbb C)^2/T\ar[rr]^{\mu_T}\ar[dd]_{P(f)}\ar[rd]^{\rho_T} & \, & 
\Var(\mathbb C)^{2,pr}/T\ar[dd]^{P(f)}\ar[rd]^{\rho_T} & \, \\
\, & \Var(\mathbb C)^{2,sm}/T\ar[rr]^{\mu_T}\ar[dd]_{P(f)} & \, & \Var(\mathbb C)^{2,smpr}/T\ar[dd]^{P(f)} \\
\Var(\mathbb C)^2/S\ar[rr]^{\mu_S}\ar[rd]^{\rho_S} & \, & \Var(\mathbb C)^{2,pr}/S\ar[rd]^{\rho_S} & \, \\
\, & \Var(\mathbb C)^{2,sm}/S\ar[rr]^{\mu_S} & \, & \Var(\mathbb C)^{2,smpr}/S}.
\end{equation}

Recall we have (see subsection 2.1), for $S\in\Var(\mathbb C)$, the graph functor 
\begin{eqnarray*}
\Gr_S^{12}:\Var(\mathbb C)/S\to\Var(\mathbb C)^{2,pr}/S, \; X/S\mapsto\Gr_S^{12}(X/S):=(X\times S,X)/S, \\
(g:X/S\to X'/S)\mapsto\Gr_S^{12}(g):=(g\times I_S:(X\times S,X)\to(X'\times S,X'))
\end{eqnarray*}
Note that $\Gr_S^{12}$ is fully faithfull.
For $f:T\to S$ a morphism with $T,S\in\Var(\mathbb C)$, we have by definition, the following commutative diagram of sites
\begin{equation}\label{Grf}
\xymatrix{\Var(\mathbb C)^{2,pr}/T\ar[rr]^{\Gr_T^{12}}\ar[dd]_{P(f)}\ar[rd]^{\rho_T} & \, & 
\Var(\mathbb C)/T\ar[dd]^{P(f)}\ar[rd]^{\rho_T} & \, \\
\, & \Var(\mathbb C)^{2,smpr}/T\ar[rr]^{\Gr_T^{12}}\ar[dd]_{P(f)} & \, & \Var(\mathbb C)^{sm}/T\ar[dd]^{P(f)} \\
\Var(\mathbb C)^{2,pr}/S\ar[rr]^{\Gr_S^{12}}\ar[rd]^{\rho_S} & \, & \Var(\mathbb C)/S\ar[rd]^{\rho_S} & \, \\
\, & \Var(\mathbb C)^{2,sm}/S\ar[rr]^{\Gr_S^{12}} & \, & \Var(\mathbb C)^{sm}/S}.
\end{equation}
where we recall that $P(f)((X,Z)/S):=((X_T,Z_T)/T)$, since smooth morphisms are preserved by base change.

\begin{itemize}
\item As usual, we denote by
\begin{equation*}
(f^*,f_*):=(P(f)^*,P(f)_*):C(\Var(\mathbb C)^{2,(sm)}/S)\to C(\Var(\mathbb C)^{2,(sm)}/T)
\end{equation*}
the adjonction induced by $P(f):\Var(\mathbb C)^{2,(sm)}/T\to \Var(\mathbb C)^{2,(sm)}/S$.
Since the colimits involved in the definition of $f^*=P(f)^*$ are filtered, $f^*$ also preserve monomorphism. 
Hence, we get an adjonction
\begin{equation*}
(f^*,f_*):C_{fil}(\Var(\mathbb C)^{2,(sm)}/S)\leftrightarrows C_{fil}(\Var(\mathbb C)^{2,(sm)}/T), \; f^*(G,F):=(f^*G,f^*F)
\end{equation*}
For $S\in\Var(\mathbb C)$, 
we denote by $\mathbb Z_S:=\mathbb Z((S,S)/(S,S))\in\PSh(\Var(\mathbb C)^{2,(sm)}/S)$ the constant presheaf.
By Yoneda lemma, we have for $F\in C(\Var(\mathbb C)^{2,(sm)}/S)$, $\mathcal Hom(\mathbb Z_S,F)=F$.
\item As usual, we denote by
\begin{equation*}
(f^*,f_*):=(P(f)^*,P(f)_*):C(\Var(\mathbb C)^{2,(sm)pr}/S)\to C(\Var(\mathbb C)^{2,(sm)pr}/T)
\end{equation*}
the adjonction induced by $P(f):\Var(\mathbb C)^{2,(sm)pr}/T\to \Var(\mathbb C)^{2,(sm)pr}/S$.
Since the colimits involved in the definition of $f^*=P(f)^*$ are filtered, $f^*$ also preserve monomorphism. 
Hence, we get an adjonction
\begin{equation*}
(f^*,f_*):C_{fil}(\Var(\mathbb C)^{2,(sm)pr}/S)\leftrightarrows C_{fil}(\Var(\mathbb C)^{2,(sm)pr}/T), \; f^*(G,F):=(f^*G,f^*F)
\end{equation*}
For $S\in\Var(\mathbb C)$, 
we denote by $\mathbb Z_S:=\mathbb Z((S,S)/(S,S))\in\PSh(\Var(\mathbb C)^{2,sm}/S)$ the constant presheaf.
By Yoneda lemma, we have for $F\in C(\Var(\mathbb C)^{2,sm}/S)$, $\mathcal Hom(\mathbb Z_S,F)=F$.
\end{itemize}

\begin{itemize}
\item For $h:U\to S$ a smooth morphism with $U,S\in\Var(\mathbb C)$, 
$P(h):\Var(\mathbb C)^{2,sm}/S\to\Var(\mathbb C)^{2,sm}/U$ admits a left adjoint
\begin{equation*}
C(h):\Var(\mathbb C)^{2,sm}/U\to\Var(\mathbb C)^{2,sm}/S, \; C(h)((U',Z'),h')=((U',Z'),h\circ h').
\end{equation*}
Hence $h^*:C(\Var(\mathbb C)^{2,sm}/S)\to C(\Var(\mathbb C)^{2,sm}/U)$ admits a left adjoint
\begin{eqnarray*}
h_{\sharp}:C(\Var(\mathbb C)^{2,sm}/U)\to C(\Var(\mathbb C)^{2,sm}/S), \\
F\mapsto (h_{\sharp}F:((U,Z),h_0)\mapsto\lim_{((U',Z'),h\circ h')\to ((U,Z),h_0)} F((U',Z')/U))
\end{eqnarray*}
\item For $h:X\to S$ a morphism with $X,S\in\Var(\mathbb C)$, 
$P(h):\Var(\mathbb C)^2/S\to\Var(\mathbb C)^2/X$ admits a left adjoint
\begin{equation*}
C(h):\Var(\mathbb C)^2/X\to\Var(\mathbb C)^2/S, \; C(h)((X',Z'),h')=((X',Z'),h\circ h').
\end{equation*}
Hence $h^*:C(\Var(\mathbb C)^2/S)\to C(\Var(\mathbb C)^2/X)$ admits a left adjoint
\begin{eqnarray*}
h_{\sharp}:C(\Var(\mathbb C)^2/X)\to C(\Var(\mathbb C)^{2,sm}/S), \\
F\mapsto (h_{\sharp}F:((X,Z),h_0)\mapsto\lim_{((X',Z'),h\circ h')\to ((X,Z),h_0)} F((X',Z')/X))
\end{eqnarray*}
\item For $p:Y\times S\to S$ a projection with $Y,S\in\Var(\mathbb C)$ with $Y$ smooth, 
$P(p):\Var(\mathbb C)^{2,smpr}/S\to\Var(\mathbb C)^{2,smpr}/Y\times S$ admits a left adjoint
\begin{eqnarray*}
C(p):\Var(\mathbb C)^{2,smpr}/Y\times S\to\Var(\mathbb C)^{2,smpr}/S, \\ 
C(p)((Y'\times S,Z'),p')=((Y'\times S,Z'),p\circ p').
\end{eqnarray*}
Hence $p^*:C(\Var(\mathbb C)^{2,smpr}/S)\to C(\Var(\mathbb C)^{2,smpr}/Y\times S)$ admits a left adjoint
\begin{eqnarray*}
p_{\sharp}:C(\Var(\mathbb C)^{2,smpr}/Y\times S)\to C(\Var(\mathbb C)^{2,smpr}/S), \\
F\mapsto (p_{\sharp}F:((Y_0\times S,Z),p_0)\mapsto
\lim_{((Y'\times Y\times S,Z'),p\circ p')\to ((Y_0\times S,Z),p_0)} F((Y'\times Y\times S,Z')/Y\times S))
\end{eqnarray*}
\item For $p:Y\times S\to S$ a projection with $Y,S\in\Var(\mathbb C)$, 
$P(p):\Var(\mathbb C)^{2,pr}/S\to\Var(\mathbb C)^{2,pr}/Y\times S$ admits a left adjoint
\begin{equation*}
C(p):\Var(\mathbb C)^{2,pr}/Y\times S\to\Var(\mathbb C)^{2,pr}/S, \; C(p)((Y'\times S,Z'),p')=((Y'\times S,Z'),p\circ p').
\end{equation*}
Hence $p^*:C(\Var(\mathbb C)^{2,pr}/S)\to C(\Var(\mathbb C)^{2,pr}/Y\times S)$ admits a left adjoint
\begin{eqnarray*}
p_{\sharp}:C(\Var(\mathbb C)^{2,pr}/Y\times S)\to C(\Var(\mathbb C)^{2,pr}/S), \\
F\mapsto (p_{\sharp}F:((Y_0\times S,Z),p_0)\mapsto
\lim_{((Y'\times Y\times S,Z'),p\circ p')\to ((Y_0\times S,Z),p_0)} F((Y'\times Y\times S,Z')/Y\times S))
\end{eqnarray*}
\end{itemize}

Let $S\in\Var(\mathbb C)$. 
\begin{itemize}
\item We have the dual functor
\begin{eqnarray*}
\mathbb D_S:C(\Var(\mathbb C)^{2,(sm)}/S)\to C(\Var(\mathbb C)^{2,(sm)}/S), \; 
F\mapsto\mathbb D_S(F):=\mathcal Hom(F,E_{et}(\mathbb Z((S,S)/S)))
\end{eqnarray*}
It induces the functor
\begin{eqnarray*}
L\mathbb D_S:C(\Var(\mathbb C)^{2,(sm)}/S)\to C(\Var(\mathbb C)^{2,(sm)}/S), \; 
F\mapsto L\mathbb D_S(F):=\mathbb D_S(LF):=\mathcal Hom(LF,E_{et}(\mathbb Z((S,S)/S)))
\end{eqnarray*}
\item We have the dual functor
\begin{eqnarray*}
\mathbb D_S:C(\Var(\mathbb C)^{2,(sm)pr}/S)\to C(\Var(\mathbb C)^{2,(sm)pr}/S), \; 
F\mapsto\mathbb D_S(F):=\mathcal Hom(F,E_{et}(\mathbb Z((S,S)/S)))
\end{eqnarray*}
It induces the functor
\begin{eqnarray*}
L\mathbb D_S:C(\Var(\mathbb C)^{2,(sm)pr}/S)\to C(\Var(\mathbb C)^{2,(sm)pr}/S), \; 
F\mapsto L\mathbb D_S(F):=\mathbb D_S(LF):=\mathcal Hom(LF,E_{et}(\mathbb Z((S,S)/S)))
\end{eqnarray*}
\end{itemize}

\begin{prop}\label{Grhom}
\begin{itemize}
\item[(i)]Let $S\in\Var(\mathbb C)$. Let $h:U\to S$ a smooth morphism with $U\in\Var(\mathbb C)$.
Then for $F\in C(\Var(\mathbb C)^{sm}/S)$, the canonical map in $C(\Var(\mathbb C)^{2,smpr}/S)$
\begin{eqnarray*}
T(\Gr_S^{12},hom)(\mathbb Z(U/S),F):\Gr_S^{12*}\mathcal Hom(\mathbb Z(U/S),F)\xrightarrow{\sim}
\mathcal Hom(\Gr_S^{12*}\mathbb Z(U/S),\Gr_S^{12*}F)
\end{eqnarray*}
is an isomorphism.
\item[(ii)]Let $S\in\Var(\mathbb C)$. Let $h:U\to S$ a morphism with $U\in\Var(\mathbb C)$.
Then for $F\in C(\Var(\mathbb C)/S)$, the canonical map in $C(\Var(\mathbb C)^{2,pr}/S)$
\begin{eqnarray*}
T(\Gr_S^{12},hom)(\mathbb Z(U/S),F):\Gr_S^{12*}\mathcal Hom(\mathbb Z(U/S),F)\xrightarrow{\sim}
\mathcal Hom(\Gr_S^{12*}\mathbb Z(U/S),\Gr_S^{12*}F)
\end{eqnarray*}
is an isomorphism.
\end{itemize}
\end{prop}

\begin{proof}
\noindent(i): We have, for $(X\times S,Z)/S\in\Var(\mathbb C)^{2,smpr}/S$ the following commutative diagram
\begin{equation*}
\xymatrix{\Gr_S^{12*}\mathcal Hom(\mathbb Z(U/S),F)((X\times S,Z)/S)
\ar[rr]^{T(\Gr_S^{12},hom)(\mathbb Z(U/S),F)((X\times S,Z)/S)}\ar[d]^{=} & \, &
\mathcal Hom(\Gr_S^{12*}\mathbb Z(U/S),\Gr_S^{12*}F)((X\times S,Z)/S)\ar[d]^{=} \\
\lim_{((X\times S,Z)/S)\to\Gr_S^{12}(V/S)}F(U\times_S V)\ar[rr] & \, &
\lim_{((X\times U,Z)/Z\times_SU)\to\Gr_S^{12}(W/S)F(W)}}
\end{equation*}
We then note that the map 
$\left\{(((X\times S,Z)/S)\to\Gr_S^{12}(V/S))\right\}\to\left\{((X\times U,Z)/Z\times_SU)\to\Gr_S^{12}(W/S)\right\}$
obviously admits an inverse since a map $(X\times U,Z\times_SU)/S\to(W\times S,W)/S)$ is uniquely determined by
a map $g:X\to W$ such that $(g\times I_S)(Z)\subset W$.
\noindent(ii):Similar to (i).
\end{proof}

We have the support section functors of a closed embedding $i:Z\hookrightarrow S$ for presheaves on the big Zariski site of pairs.
\begin{defi}\label{gamma12}
Let $i:Z\hookrightarrow S$ be a closed embedding with $S,Z\in\Var(\mathbb C)$ 
and $j:S\backslash Z\hookrightarrow S$ be the open complementary subset.
\begin{itemize}
\item[(i)] We define the functor
\begin{equation*}
\Gamma_Z:C(\Var(\mathbb C)^{2,sm}/S)\to C(\Var(\mathbb C)^{2,sm}/S), \;
G^{\bullet}\mapsto\Gamma_Z G^{\bullet}:=\Cone(\ad(j^*,j_*)(G^{\bullet}):G^{\bullet}\to j_*j^*G^{\bullet})[-1],
\end{equation*}
so that there is then a canonical map $\gamma_Z(G^{\bullet}):\Gamma_ZG^{\bullet}\to G^{\bullet}$.
\item[(ii)] We have the dual functor of (i) :
\begin{equation*}
\Gamma^{\vee}_Z:C(\Var(\mathbb C)^{2,sm}/S)\to C(\Var(\mathbb C)^{2,sm}/S), \; 
F\mapsto\Gamma^{\vee}_Z(F^{\bullet}):=\Cone(\ad(j_{\sharp},j^*)(G^{\bullet}):j_{\sharp}j^*G^{\bullet}\to G^{\bullet}), 
\end{equation*}
together with the canonical map $\gamma^{\vee}_Z(G):F\to\Gamma^{\vee}_Z(G)$.
\item[(iii)] For $F,G\in C(\Var(\mathbb C)^{2,sm}/S)$, we denote by 
\begin{equation*}
I(\gamma,hom)(F,G):=(I,I(j_{\sharp},j^*)(F,G)^{-1}):\Gamma_Z\mathcal Hom(F,G)\xrightarrow{\sim}\mathcal Hom(\Gamma^{\vee}_ZF,G)
\end{equation*}
the canonical isomorphism given by adjonction.
\end{itemize}
\end{defi}

Note that we have similarly for $i:Z\hookrightarrow S$, $i':Z'\hookrightarrow Z$ closed embeddings,
$g:T\to S$ a morphism with $T,S,Z\in\Var(\mathbb C)$ and $F\in C(\Var(\mathbb C)^{2,sm}/S)$, 
the canonical maps in $C(\Var(\mathbb C)^{2,sm}/S)$
\begin{itemize}
\item $T(g,\gamma)(F):g^*\Gamma_ZF\xrightarrow{\sim}\Gamma_{Z\times_S T}g^*F$, 
$T(g,\gamma^{\vee})(F):\Gamma^{\vee}_{Z\times_S T}g^*F\xrightarrow{\sim}g^*\Gamma_ZF$
\item $T(Z'/Z,\gamma)(F):\Gamma_{Z'}F\to\Gamma_Z F$, $T(Z'/Z,\gamma^{\vee})(F):\Gamma^{\vee}_ZF\to\Gamma_{Z'}^{\vee}F$
\end{itemize}
but we will not use them in this article.

Let $S_{\bullet}\in\Fun(\mathcal I,\Var(\mathbb C))$ with $\mathcal I\in\Cat$, a diagram of algebraic varieties.
It gives the diagram of sites $\Var(\mathbb C)^2/S_{\bullet}\in\Fun(\mathcal I,\Cat)$.  
\begin{itemize}
\item Then $C_{fil}(\Var(\mathbb C)^{2,(sm)}/S_{\bullet})$ is the category  
\begin{itemize}
\item whose objects $(G,F)=((G_I,F)_{I\in\mathcal I},u_{IJ})$,
with $(G_I,F)\in C_{fil}(\Var(\mathbb C)^{2,(sm)}/S_I)$,
and $u_{IJ}:(G_I,F)\to r_{IJ*}(G_J,F)$ for $r_{IJ}:I\to J$, denoting again $r_{IJ}:S_I\to S_J$, are morphisms
satisfying for $I\to J\to K$, $r_{IJ*}u_{JK}\circ u_{IJ}=u_{IK}$ in $C_{fil}(\Var(\mathbb C)^{2,(sm)}/S_I)$,
\item the morphisms $m:((G,F),u_{IJ})\to((H,F),v_{IJ})$ being (see section 2.1) a family of morphisms of complexes,  
\begin{equation*}
m=(m_I:(G_I,F)\to (H_I,F))_{I\in\mathcal I}
\end{equation*}
such that $v_{IJ}\circ m_I=p_{IJ*}m_J\circ u_{IJ}$ in $C_{fil}(\Var(\mathbb C)^{2,(sm)}/S_I)$.
\end{itemize}
\item Then $C_{fil}(\Var(\mathbb C)^{2,(sm)pr}/S_{\bullet})$ is the category  
\begin{itemize}
\item whose objects $(G,F)=((G_I,F)_{I\in\mathcal I},u_{IJ})$,
with $(G_I,F)\in C_{fil}(\Var(\mathbb C)^{2,(sm)pr}/S_I)$,
and $u_{IJ}:(G_I,F)\to r_{IJ*}(G_J,F)$ for $r_{IJ}:I\to J$, denoting again $r_{IJ}:S_I\to S_J$, are morphisms
satisfying for $I\to J\to K$, $r_{IJ*}u_{JK}\circ u_{IJ}=u_{IK}$ in $C_{fil}(\Var(\mathbb C)^{2,(sm)}/S_I)$,
\item the morphisms $m:((G,F),u_{IJ})\to((H,F),v_{IJ})$ being (see section 2.1) a family of morphisms of complexes,  
\begin{equation*}
m=(m_I:(G_I,F)\to (H_I,F))_{I\in\mathcal I}
\end{equation*}
such that $v_{IJ}\circ m_I=p_{IJ*}m_J\circ u_{IJ}$ in $C_{fil}(\Var(\mathbb C)^{2,(sm)pr}/S_I)$.
\end{itemize}
\end{itemize}
For $s:\mathcal I\to\mathcal J$ a functor, with $\mathcal I,\mathcal J\in\Cat$, and
$f_{\bullet}:T_{\bullet}\to S_{s(\bullet)}$ a morphism with 
$T_{\bullet}\in\Fun(\mathcal J,\Var(\mathbb C))$ and $S_{\bullet}\in\Fun(\mathcal I,\Var(\mathbb C))$, 
we have by definition, the following commutative diagrams of sites
\begin{equation}\label{mufIJ}
\xymatrix{\Var(\mathbb C)^2/T_{\bullet}\ar[rr]^{\mu_{T_{\bullet}}}\ar[dd]_{P(f_{\bullet})}\ar[rd]^{\rho_{T_{\bullet}}} & \, & 
\Var(\mathbb C)^{2,pr}/T_{\bullet}\ar[dd]^{P(f_{\bullet})}\ar[rd]^{\rho_{T_{\bullet}}} & \, \\
\, & \Var(\mathbb C)^{2,sm}/T_{\bullet}\ar[rr]^{\mu_{T_{\bullet}}}\ar[dd]_{P(f_{\bullet})} & \, & 
\Var(\mathbb C)^{2,smpr}/T_{\bullet}\ar[dd]^{P(f_{\bullet})} \\
\Var(\mathbb C)^2/S_{\bullet}\ar[rr]^{\mu_{S_{\bullet}}}\ar[rd]^{\rho_{S_{\bullet}}} & \, & 
\Var(\mathbb C)^{2,pr}/S_{\bullet}\ar[rd]^{\rho_{S_{\bullet}}} & \, \\
\, & \Var(\mathbb C)^{2,sm}/S_{\bullet}\ar[rr]^{\mu_{S_{\bullet}}} & \, & \Var(\mathbb C)^{2,smpr}/S_{\bullet}}.
\end{equation}
and
\begin{equation}\label{GrfIJ}
\xymatrix{\Var(\mathbb C)^{2,pr}/T_{\bullet}
\ar[rr]^{\Gr_{T_{\bullet}}^{12}}\ar[dd]_{P(f_{\bullet})}\ar[rd]^{\rho_{T_{\bullet}}} & \, & 
\Var(\mathbb C)/T\ar[dd]^{P(f_{\bullet})}\ar[rd]^{\rho_{T_{\bullet}}} & \, \\
\, & \Var(\mathbb C)^{2,smpr}/T_{\bullet}\ar[rr]^{\Gr_{T}^{12}}\ar[dd]_{P(f_{\bullet})} & \, & 
\Var(\mathbb C)^{sm}/T_{\bullet}\ar[dd]^{P(f_{\bullet})} \\
\Var(\mathbb C)^{2,pr}/S_{\bullet}\ar[rr]^{\Gr_{S_{\bullet}}^{12}}\ar[rd]^{\rho_{S_{\bullet}}} & \, & 
\Var(\mathbb C)/S_{\bullet}\ar[rd]^{\rho_{S_{\bullet}}} & \, \\
\, & \Var(\mathbb C)^{2,sm}/S_{\bullet}\ar[rr]^{\Gr_{S_{\bullet}}^{12}} & \, & 
\Var(\mathbb C)^{sm}/S_{\bullet}}.
\end{equation}
Let $s:\mathcal I\to\mathcal J$ a functor, with $\mathcal I,\mathcal J\in\Cat$, and
$f_{\bullet}:T_{\bullet}\to S_{s(\bullet)}$ a morphism with 
$T_{\bullet}\in\Fun(\mathcal J,\Var(\mathbb C))$ and $S_{\bullet}\in\Fun(\mathcal I,\Var(\mathbb C))$.
\begin{itemize}
\item As usual, we denote by
\begin{eqnarray*}
(f_{\bullet}^*,f_{\bullet*}):=(P(f_{\bullet})^*,P(f_{\bullet})_*): 
C(\Var(\mathbb C)^{2,(sm)}/S_{\bullet})\to C(\Var(\mathbb C)^{2,(sm)}/T_{\bullet})
\end{eqnarray*}
the adjonction induced by 
$P(f_{\bullet}):\Var(\mathbb C)^{2,(sm)}/T_{\bullet}\to \Var(\mathbb C)^{2,(sm)}/S_{\bullet}$.
Since the colimits involved in the definition of $f_{\bullet}^*=P(f_{\bullet})^*$ are filtered, 
$f_{\bullet}^*$ also preserve monomorphism. Hence, we get an adjonction
\begin{eqnarray*}
(f_{\bullet}^*,f_{\bullet*}):
C_{fil}(\Var(\mathbb C)^{2,(sm)}/S_{\bullet})\leftrightarrows C_{fil}(\Var(\mathbb C)^{2,(sm)}/T_{\bullet}), \\
f_{\bullet}^*((G_I,F),u_{IJ}):=((f_I^*G_I,f_I^*F),T(f_I,r_{IJ})(-)\circ f_I^*u_{IJ})
\end{eqnarray*}
\item As usual, we denote by
\begin{eqnarray*}
(f_{\bullet}^*,f_{\bullet*}):=(P(f_{\bullet})^*,P(f_{\bullet})_*):
C(\Var(\mathbb C)^{2,(sm)pr}/S_{\bullet})\to C(\Var(\mathbb C)^{2,(sm)pr}/T_{\bullet})
\end{eqnarray*}
the adjonction induced by 
$P(f_{\bullet}):\Var(\mathbb C)^{2,(sm)pr}/T_{\bullet}\to \Var(\mathbb C)^{2,(sm)pr}/S_{\bullet}$.
Since the colimits involved in the definition of $f_{\bullet}^*=P(f_{\bullet})^*$ are filtered, 
$f_{\bullet}^*$ also preserve monomorphism. Hence, we get an adjonction
\begin{eqnarray*}
(f_{\bullet}^*,f_{\bullet*}):
C_{fil}(\Var(\mathbb C)^{2,(sm)pr}/S_{\bullet})\leftrightarrows C_{fil}(\Var(\mathbb C)^{2,(sm)pr}/T_{\bullet}), \\
f_{\bullet}^*((G_I,F),u_{IJ}):=((f_I^*G_I,f_I^*F),T(f_I,r_{IJ})(-)\circ f_I^*u_{IJ})
\end{eqnarray*}
\end{itemize}

Let $S\in\Var(\mathbb C)$. Let $S=\cup_{i=1}^l S_i$ an open affine cover and denote by $S_I=\cap_{i\in I} S_i$.
Let $i_i:S_i\hookrightarrow\tilde S_i$ closed embeddings, with $\tilde S_i\in\Var(\mathbb C)$. 
For $I\subset\left[1,\cdots l\right]$, denote by $\tilde S_I=\Pi_{i\in I}\tilde S_i$.
We then have closed embeddings $i_I:S_I\hookrightarrow\tilde S_I$ and for $J\subset I$ the following commutative diagram
\begin{equation*}
D_{IJ}=\xymatrix{ S_I\ar[r]^{i_I} & \tilde S_I \\
S_J\ar[u]^{j_{IJ}}\ar[r]^{i_J} & \tilde S_J\ar[u]^{p_{IJ}}}  
\end{equation*}
where $p_{IJ}:\tilde S_J\to\tilde S_I$ is the projection
and $j_{IJ}:S_J\hookrightarrow S_I$ is the open embedding so that $j_I\circ j_{IJ}=j_J$.
This gives the diagram of algebraic varieties $(\tilde S_I)\in\Fun(\mathcal P(\mathbb N),\Var(\mathbb C))$ 
which gives the diagram of sites $\Var(\mathbb C)^2/(\tilde S_I)\in\Fun(\mathcal P(\mathbb N),\Cat)$.  
This gives also the diagram of algebraic varieties $(\tilde S_I)^{op}\in\Fun(\mathcal P(\mathbb N)^{op},\Var(\mathbb C))$
which gives the diagram of sites $\Var(\mathbb C)^2/(\tilde S_I)^{op}\in\Fun(\mathcal P(\mathbb N)^{op},\Cat)$.  
\begin{itemize}
\item Then $C_{fil}(\Var(\mathbb C)^{2,(sm)}/(\tilde S_I))$ is the category  
\begin{itemize}
\item whose objects $(G,F)=((G_I,F)_{I\subset\left[1,\cdots l\right]},u_{IJ})$,
with $(G_I,F)\in C_{fil}(\Var(\mathbb C)^{2,(sm)}/\tilde S_I)$,
and $u_{IJ}:(G_I,F)\to p_{IJ*}(G_J,F)$ for $I\subset J$, are morphisms
satisfying for $I\subset J\subset K$, $p_{IJ*}u_{JK}\circ u_{IJ}=u_{IK}$ in $C_{fil}(\Var(\mathbb C)^{2,(sm)}/\tilde S_I)$,
\item the morphisms $m:((G,F),u_{IJ})\to((H,F),v_{IJ})$ being (see section 2.1) a family of morphisms of complexes,  
\begin{equation*}
m=(m_I:(G_I,F)\to (H_I,F))_{I\subset\left[1,\cdots l\right]}
\end{equation*}
such that $v_{IJ}\circ m_I=p_{IJ*}m_J\circ u_{IJ}$ in $C_{fil}(\Var(\mathbb C)^{2,(sm)}/\tilde S_I)$.
\end{itemize}
\item Then $C_{fil}(\Var(\mathbb C)^{2,(sm)pr}/(\tilde S_I))$ is the category  
\begin{itemize}
\item whose objects $(G,F)=((G_I,F)_{I\subset\left[1,\cdots l\right]},u_{IJ})$,
with $(G_I,F)\in C_{fil}(\Var(\mathbb C)^{2,(sm)pr}/\tilde S_I)$,
and $u_{IJ}:(G_I,F)\to p_{IJ*}(G_J,F)$ for $I\subset J$, are morphisms
satisfying for $I\subset J\subset K$, $p_{IJ*}u_{JK}\circ u_{IJ}=u_{IK}$ in $C_{fil}(\Var(\mathbb C)^{2,(sm)pr}/\tilde S_I)$,
\item the morphisms $m:((G,F),u_{IJ})\to((H,F),v_{IJ})$ being (see section 2.1) a family of morphisms of complexes,  
\begin{equation*}
m=(m_I:(G_I,F)\to (H_I,F))_{I\subset\left[1,\cdots l\right]}
\end{equation*}
such that $v_{IJ}\circ m_I=p_{IJ*}m_J\circ u_{IJ}$ in $C_{fil}(\Var(\mathbb C)^{2,(sm)pr}/\tilde S_I)$.
\end{itemize}
\item Then $C_{fil}(\Var(\mathbb C)^{2,(sm)}/(\tilde S_I)^{op})$ 
is the category  
\begin{itemize}
\item whose objects $(G,F)=((G_I,F)_{I\subset\left[1,\cdots l\right]},u_{IJ})$,
with $(G_I,F)\in C_{fil}(\Var(\mathbb C)^{2,(sm)}/\tilde S_I)$,
and $u_{IJ}:(G_J,F)\to p_{IJ}^*(G_I,F)$ for $I\subset J$, are morphisms
satisfying for $I\subset J\subset K$, $p_{JK}^*u_{IJ}\circ u_{JK}=u_{IK}$ in $C_{fil}(\Var(\mathbb C)^{2,(sm)}/\tilde S_K)$,
\item the morphisms $m:((G,F),u_{IJ})\to((H,F),v_{IJ})$ being (see section 2.1) a family of morphisms of complexes,  
\begin{equation*}
m=(m_I:(G_I,F)\to (H_I,F))_{I\subset\left[1,\cdots l\right]}
\end{equation*}
such that $v_{IJ}\circ m_J=p_{IJ}^*m_I\circ u_{IJ}$ in $C_{fil}(\Var(\mathbb C)^{2,(sm)}/\tilde S_J)$.
\end{itemize}
\item Then $C_{fil}(\Var(\mathbb C)^{2,(sm)pr}/(\tilde S_I)^{op})$ 
is the category  
\begin{itemize}
\item whose objects $(G,F)=((G_I,F)_{I\subset\left[1,\cdots l\right]},u_{IJ})$,
with $(G_I,F)\in C_{fil}(\Var(\mathbb C)^{2,(sm)pr}/\tilde S_I)$,
and $u_{IJ}:(G_J,F)\to p_{IJ}^*(G_I,F)$ for $I\subset J$, are morphisms
satisfying for $I\subset J\subset K$, $p_{JK}^*u_{IJ}\circ u_{JK}=u_{IK}$ in $C_{fil}(\Var(\mathbb C)^{2,(sm)pr}/\tilde S_K)$,
\item the morphisms $m:((G,F),u_{IJ})\to((H,F),v_{IJ})$ being (see section 2.1) a family of morphisms of complexes,  
\begin{equation*}
m=(m_I:(G_I,F)\to (H_I,F))_{I\subset\left[1,\cdots l\right]}
\end{equation*}
such that $v_{IJ}\circ m_J=p_{IJ}^*m_I\circ u_{IJ}$ in $C_{fil}(\Var(\mathbb C)^{2,(sm)pr}/\tilde S_J)$.
\end{itemize}
\end{itemize}

We now define the Zariski and the etale topology on $\Var(\mathbb C)^2/S$.

\begin{defi}\label{tau12}
Let $S\in\Var(\mathbb C)$. 
\begin{itemize}
\item[(i)]Denote by $\tau$ a topology on $\Var(\mathbb C)$, e.g. the Zariski or the etale topology. 
The $\tau$ covers in $\Var(\mathbb C)^2/S$ of $(X,Z)/S$ are the families of morphisms 
\begin{eqnarray*}
\left\{(c_i:(U_i,Z\times_X U_i)/S\to(X,Z)/S)_{i\in I}, \; 
\mbox{with} \; (c_i:U_i\to X)_{i\in I} \, \tau \, \mbox{cover of} \, X \, \mbox{in} \, \Var(\mathbb C)\right\}
\end{eqnarray*}
\item[(ii)]Denote by $\tau$ the Zariski or the etale topology on $\Var(\mathbb C)$. 
The $\tau$ covers in $\Var(\mathbb C)^{2,sm}/S$ of $(U,Z)/S$ are the families of morphisms 
\begin{eqnarray*}
\left\{(c_i:(U_i,Z\times_U U_i)/S\to(U,Z)/S)_{i\in I}, \; 
\mbox{with} \; (c_i:U_i\to U)_{i\in I} \, \tau \, \mbox{cover of} \, U \, \mbox{in} \, \Var(\mathbb C)\right\}
\end{eqnarray*}
\item[(iii)]Denote by $\tau$ the Zariski or the etale topology on $\Var(\mathbb C)$. 
The $\tau$ covers in $\Var(\mathbb C)^{2,(sm)pr}/S$ of $(Y\times S,Z)/S$ are the families of morphisms 
\begin{eqnarray*}
\left\{(c_i\times I_S:(U_i\times S,Z\times_{Y\times S} U_i\times S)/S\to(Y\times S,Z)/S)_{i\in I}, \; 
\mbox{with} \; (c_i:U_i\to Y)_{i\in I} \, \tau \, \mbox{cover of} \, Y \, \mbox{in} \, \Var(\mathbb C)\right\}
\end{eqnarray*}
\end{itemize}
\end{defi}

Let $S\in\Var(\mathbb C)$. Denote by $\tau$ the Zariski or the etale topology on $\Var(\mathbb C)$.
In particular, denoting 
$a_{\tau}:\PSh(\Var(\mathbb C)^{2,(sm)}/S)\to\Shv(\Var(\mathbb C)^{2,(sm)}/S)$ and
$a_{\tau}:\PSh(\Var(\mathbb C)^{2,(sm)pr}/S)\to\Shv(\Var(\mathbb C)^{2,(sm)pr}/S)$
the sheaftification functors,
\begin{itemize}
\item a morphism $\phi:F\to G$, with $F,G\in C(\Var(\mathbb C)^{2,(sm)}/S)$,
is a $\tau$ local equivalence if $a_{\tau}H^n\phi:a_{\tau}H^nF\to a_{\tau}H^nG$ is an isomorphism,
a morphism $\phi:F\to G$, with $F,G\in C(\Var(\mathbb C)^{2,(sm)pr}/S)$,
is a $\tau$ local equivalence if $a_{\tau}H^n\phi:a_{\tau}H^nF\to a_{\tau}H^nG$ is an isomorphism ;
\item $F^{\bullet}\in C(\Var(\mathbb C)^{2,(sm)}/S)$ is $\tau$ fibrant 
if for all $(U,Z)/S\in\Var(\mathbb C)^{2,(sm)}/S$ and all $\tau$ covers 
$(c_i:(U_i,Z\times_U U_i)/S\to(U,Z)/S)_{i\in I}$ of $(U,Z)/S$,
\begin{equation*}
F^{\bullet}(c_i):F^{\bullet}((U,Z)/S)\to\Tot(\oplus_{card I=\bullet} F^{\bullet}((U_I,Z\times_UU_I)/S))
\end{equation*}
is a quasi-isomorphism of complexes of abelian groups,
$F^{\bullet}\in C(\Var(\mathbb C)^{2,(sm)pr}/S)$ is $\tau$ fibrant 
if for all $(Y\times S,Z)/S\in\Var(\mathbb C)^{2,(sm)pr}/S$ and all $\tau$ covers 
$(c_i\times I_S:(U_i\times S,Z\times_{Y\times S}U_i\times S)/S\to(Y\times S,Z)/S)_{i\in I}$ of $(Y\times S,Z)/S$,
\begin{equation*}
F^{\bullet}(c_i\times I_S):F^{\bullet}((Y\times S,Z)/S)\to
\Tot(\oplus_{card I=\bullet} F^{\bullet}((U_I\times S,Z_I\times_YU_J)/S))
\end{equation*}
is a quasi-isomorphism of complexes of abelian groups ;
\item a morphism $\phi:(G_1,F)\to (G_2,F)$, with $(G_1,F),(G_2,F)\in C_{fil}(\Var(\mathbb C)^{2,(sm)}/S)$, 
is an filtered $\tau$ local equivalence if for all $n,p\in\mathbb Z$,
\begin{equation*}
a_{\tau}H^n\Gr_F^p(\phi):a_{\tau}H^n\Gr_F^p(G_1,F)\to a_{\tau}H^n\Gr_F^p(G_2,F) 
\end{equation*}
is an isomorphism of sheaves on $\Var(\mathbb C)^{2,(sm)}/S$, 
a morphism $\phi:(G_1,F)\to (G_2,F)$, with $(G_1,F),(G_2,F)\in C_{fil}(\Var(\mathbb C)^{2,(sm)pr}/S)$, 
is an filtered $\tau$ local equivalence if for all $n,p\in\mathbb Z$, 
\begin{equation*}
a_{\tau}H^n\Gr_F^p(\phi):a_{\tau}H^n\Gr_F^p(G_1,F)\to a_{\tau}H^n\Gr_F^p(G_2,F) 
\end{equation*}
is an isomorphism of sheaves on $\Var(\mathbb C)^{2,(sm)pr}/S$ ;
\item let $r\in\mathbb N$,
a morphism $\phi:(G_1,F)\to (G_2,F)$, with $(G_1,F),(G_2,F)\in C_{fil}(\Var(\mathbb C)^{2,(sm)}/S)$, 
is an $r$-filtered $\tau$ local equivalence if there exists $\phi_i:(C_i,F)\to(C_{i+1},F)$, $0\leq i\leq s$, 
with $(C_i,F)\in C_{fil}(\Var(\mathbb C)^{2,(sm)}/S)$, $(C_0,F)=(G_1,F)$ and $(C_s,F)=(G_2,F)$, such that
\begin{equation*}
\phi=\phi_s\circ\cdots\circ\phi_i\circ\cdots\circ\phi_0:(G_1,F)\to(G_2,F)
\end{equation*}
and $\phi_i:(C_i,F)\to(C_{i+1},F)$ either a filtered $\tau$ local equivalence
or an $r$-filtered homotopy equivalence, 
a morphism $\phi:(G_1,F)\to (G_2,F)$, with $(G_1,F),(G_2,F)\in C_{fil}(\Var(\mathbb C)^{2,(sm)pr}/S)$, 
is an $r$-filtered $\tau$ local equivalence if there exists $\phi_i:(C_i,F)\to(C_{i+1},F)$, $0\leq i\leq s$, 
with $(C_i,F)\in C_{fil}(\Var(\mathbb C)^{2,(sm)pr}/S)$, $(C_0,F)=(G_1,F)$ and $(C_s,F)=(G_2,F)$, such that
\begin{equation*}
\phi=\phi_s\circ\cdots\circ\phi_i\circ\cdots\circ\phi_0:(G_1,F)\to(G_2,F)
\end{equation*}
and $\phi_i:(C_i,F)\to(C_{i+1},F)$ either a filtered $\tau$ local equivalence
or an $r$-filtered homotopy equivalence ;
\item $(F^{\bullet},F)\in C_{fil}(\Var(\mathbb C)^{2,(sm)}/S)$ is filtered $\tau$ fibrant if
for all $(U,Z)/S\in\Var(\mathbb C)^{2,(sm)}/S$ and all $\tau$ covers 
$(c_i:(U_i,Z\times_UU_i)/S\to(U,Z)/S)_{i\in I}$ of $(U,Z)/S$,
\begin{eqnarray*}
H^n\Gr_F^p(F^{\bullet},F)(c_i):(F^{\bullet},F)((U,Z)/S)\to \\
H^n\Gr_F^p(\Tot(\oplus_{card I=\bullet}(F^{\bullet},F)((U_I,Z\times_UU_I)/S)))
\end{eqnarray*}
is an isomorphism of of abelian groups for all $n,p\in\mathbb Z$ ;
$(F^{\bullet},F)\in C_{fil}(\Var(\mathbb C)^{2,(sm)pr}/S)$ is filtered $\tau$ fibrant 
if for all $(Y\times S,Z)/S\in\Var(\mathbb C)^{2,(sm)pr}/S$ and all $\tau$ covers 
$(c_i\times I_S:(U_i\times S,Z\times_{Y\times S}U_i\times S)/S\to(Y\times S,Z)/S)_{i\in I}$ of $(Y\times S,Z)/S$,
\begin{eqnarray*}
H^n\Gr_F^p(F^{\bullet},F)(c_i\times I_S):H^n\Gr_F^p(F^{\bullet},F)((Y\times S,Z)/S)\xrightarrow{\sim} \\
H^n\Gr_F^p(\Tot(\oplus_{card I=\bullet} (F^{\bullet},F)((U_I\times S,Z\times_YU_I)/S)))
\end{eqnarray*}
is an isomorphism of abelian groups for all $n,p\in\mathbb Z$.
\end{itemize}
Let $S_{\bullet}\in\Fun(\mathcal I,\Var(\mathbb C))$ with $\mathcal I\in\Cat$.
\begin{itemize}
\item A morphism $m:((G_I,F),u_{IJ})\to((H_I,F),v_{IJ})$ in $C_{fil}(\Var(\mathbb C)^{2,(sm)}/S_{\bullet})$
is an $r$-filtered Zariski, resp. etale local, equivalence,
if there exists $\phi_i:((C_{iI},F),u_{iIJ})\to((C_{(i+1)I},F),u_{(i+1)IJ})$, $0\leq i\leq s$, 
with $((C_{iI},F),u_{iIJ})\in C_{fil}(\Var(\mathbb C)^{2,(sm)}/S_{\bullet})$
$((C_{0I},F),u_{0IJ})=((G_I,F),u_{IJ})$ and $((C_{sI},F),u_{sIJ})=((H_I,F),v_{IJ})$, such that
\begin{equation*}
\phi=\phi_s\circ\cdots\circ\phi_i\circ\cdots\circ\phi_0:((G_I,F),u_{IJ})\to((H_I,F),v_{IJ})
\end{equation*}
and $\phi_i:((C_{iI},F),u_{iIJ})\to((C_{(i+1)I},F),u_{(i+1)IJ})$ either a filtered Zariski, resp. etale, local equivalence
or an $r$-filtered homotopy equivalence.  
\item A morphism $m:((G_I,F),u_{IJ})\to((H_I,F),v_{IJ})$ in $C_{fil}(\Var(\mathbb C)^{2,(sm)pr}/S_{\bullet})$
is an $r$-filtered Zariski, resp. etale local, equivalence,
if there exists $\phi_i:((C_{iI},F),u_{iIJ})\to((C_{(i+1)I},F),u_{(i+1)IJ})$, $0\leq i\leq s$, 
with $((C_{iI},F),u_{iIJ})\in C_{fil}(\Var(\mathbb C)^{2,(sm)pr}/S_{\bullet})$
$((C_{0I},F),u_{0IJ})=((G_I,F),u_{IJ})$ and $((C_{sI},F),u_{sIJ})=((H_I,F),v_{IJ})$, such that
\begin{equation*}
\phi=\phi_s\circ\cdots\circ\phi_i\circ\cdots\circ\phi_0:((G_I,F),u_{IJ})\to((H_I,F),v_{IJ})
\end{equation*}
and $\phi_i:((C_{iI},F),u_{iIJ})\to((C_{(i+1)I},F),u_{(i+1)IJ})$ either a filtered Zariski, resp. etale, local equivalence
or an $r$-filtered homotopy equivalence. 
\end{itemize}

Will now define the $\mathbb A^1$ local property on $\Var(\mathbb C)^2/S$.

Denote $\square^*:=\mathbb P_{\mathbb C}^*\backslash\left\{1\right\}$
\begin{itemize}
\item Let $S\in\Var(\mathbb C)$. For $(X,Z)/S=((X,Z),h)\in\Var(\mathbb C)^{2,(sm)}/S$, we consider 
\begin{equation*}
(\square^*\times X,\square^*\times Z)/S=((\square^*\times X,\square^*\times Z,h\circ p)\in\Fun(\Delta,\Var(\mathbb C)^{2,(sm)}/S).
\end{equation*}
For $F\in C^-(\Var(\mathbb C)^{2,(sm)}/S)$, it gives the complex
\begin{equation*}
C_*F\in C^-(\Var(\mathbb C)^{2,(sm)}/S), (X,Z)/S=((X,Z),h)\mapsto C_*F((X,Z)/S):=\Tot F((\square^*\times X,\square^*\times Z/S)
\end{equation*}
together with the canonical map $c_F:=(0,I_F):F\to C_*F$.
For $F\in C(\Var(\mathbb C)^{2,(sm)}/S)$, we get
\begin{equation*}
C_*F:=\holim_n C_*F^{\leq n}\in C(\Var(\mathbb C)^{2,(sm)}/S),
\end{equation*}
together with the canonical map $c_F:=(0,I_F):F\to C_*F$.
For $m:F\to G$ a morphism, with $F,G\in C(\Var(\mathbb C)^{2,(sm)}/S)$,
we get by functoriality the morphism $C_*m:C_*F\to C_*G$.
\item Let $S\in\Var(\mathbb C)$. For $(Y\times S,Z)/S=((Y\times S,Z),h)\in\Var(\mathbb C)^{2,(sm)pr}/S$, we consider 
\begin{equation*}
(\square^*\times Y\times S,\square^*\times Z)/S=(\square^*\times Y\times S,\square^*\times Z,h\circ p)
\in\Fun(\Delta,\Var(\mathbb C)/S).
\end{equation*}
For $F\in C^-(\Var(\mathbb C)^{2,(sm)pr}/S)$, it gives the complex
\begin{eqnarray*}
C_*F\in C^-(\Var(\mathbb C)^{2,(sm)pr}/S), \\ 
(Y\times S,Z)/S=((Y\times S,Z),h)\mapsto C_*F((Y\times S,Z)/S):=\Tot F(\square^*\times Y\times S,\square^*\times Z)/S)
\end{eqnarray*}
together with the canonical map $c=c(F):=(0,I_F):F\to C_*F$.
For $F\in C(\Var(\mathbb C)^{2,(sm)pr}/S)$, we get
\begin{equation*}
C_*F:=\holim_n C_*F^{\leq n}\in C(\Var(\mathbb C)^{2,(sm)pr}/S),
\end{equation*}
together with the canonical map $c=c(F):=(0,I_F):F\to C_*F$.
For $m:F\to G$ a morphism, with $F,G\in C(\Var(\mathbb C)^{2,(sm)pr}/S)$,
we get by functoriality the morphism $C_*m:C_*F\to C_*G$.
\item Let $S\in\Var(\mathbb C)$. Let $S=\cup_{i=1}^l S_i$ an open affine cover and denote by $S_I=\cap_{i\in I} S_i$.
Let $i_i:S_i\hookrightarrow\tilde S_i$ closed embeddings, with $\tilde S_i\in\Var(\mathbb C)$.
For $F=(F_I,u_{IJ})\in C(\Var(\mathbb C)^{2,(sm)}/(\tilde S_I))$, it gives the complex
\begin{equation*}
C_*F=(C_*F_I,C_*u_{IJ})\in C(\Var(\mathbb C)^{2,(sm)}/(\tilde S_I)), 
\end{equation*}
together with the canonical map $c_F:=(0,I_F):F\to C_*F$.
\item Let $S\in\Var(\mathbb C)$. Let $S=\cup_{i=1}^l S_i$ an open affine cover and denote by $S_I=\cap_{i\in I} S_i$.
Let $i_i:S_i\hookrightarrow\tilde S_i$ closed embeddings, with $\tilde S_i\in\Var(\mathbb C)$.
For $F=(F_I,u_{IJ})\in C(\Var(\mathbb C)^{2,(sm)}/(\tilde S_I))$, it gives the complex
\begin{equation*}
C_*F=(C_*F_I,C_*u_{IJ})\in C(\Var(\mathbb C)^{2,(sm)}/(\tilde S_I)^{op}), 
\end{equation*}
together with the canonical map $c_F:=(0,I_F):F\to C_*F$.
\item Let $S\in\Var(\mathbb C)$. Let $S=\cup_{i=1}^l S_i$ an open affine cover and denote by $S_I=\cap_{i\in I} S_i$.
Let $i_i:S_i\hookrightarrow\tilde S_i$ closed embeddings, with $\tilde S_i\in\Var(\mathbb C)$.
For $F=(F_I,u_{IJ})\in C(\Var(\mathbb C)^{2,(sm)pr}/(\tilde S_I))$, it gives the complex
\begin{equation*}
C_*F=(C_*F_I,C_*u_{IJ})\in C(\Var(\mathbb C)^{2,(sm)pr}/(\tilde S_I)), 
\end{equation*}
together with the canonical map $c_F:=(0,I_F):F\to C_*F$.
\item Let $S\in\Var(\mathbb C)$. Let $S=\cup_{i=1}^l S_i$ an open affine cover and denote by $S_I=\cap_{i\in I} S_i$.
Let $i_i:S_i\hookrightarrow\tilde S_i$ closed embeddings, with $\tilde S_i\in\Var(\mathbb C)$.
For $F=(F_I,u_{IJ})\in C(\Var(\mathbb C)^{2,(sm)pr}/(\tilde S_I))$, it gives the complex
\begin{equation*}
C_*F=(C_*F_I,C_*u_{IJ})\in C(\Var(\mathbb C)^{2,(sm)pr}/(\tilde S_I)^{op}), 
\end{equation*}
together with the canonical map $c_F:=(0,I_F):F\to C_*F$.
\end{itemize}

Let $S\in\Var(\mathbb C)$. Denote for short $\Var(\mathbb C)^{2,(sm)}/S$ 
either the category $\Var(\mathbb C)^2/S$ or the category $\Var(\mathbb C)^{2,sm}/S$. Denote by
\begin{eqnarray*}
p_a:\Var(\mathbb C)^{2,(sm)}/S\to\Var(\mathbb C)^{2,(sm)}/S, \\ 
(X,Z)/S=((X,Z),h)\mapsto (X\times\mathbb A^1,Z\times\mathbb A^1)/S=((X\times\mathbb A^1,Z\times\mathbb A^1,h\circ p_X), \\ 
(g:(X,Z)/S\to (X',Z')/S)\mapsto 
((g\times I_{\mathbb A^1}):(X\times\mathbb A^1,Z\times\mathbb A^1)/S\to (X'\times\mathbb A^1,Z'\times\mathbb A^1)/S)
\end{eqnarray*}
the projection functor and again by $p_a:\Var(\mathbb C)^{2,(sm)}/S\to\Var(\mathbb C)^{2,(sm)}/S$
the corresponding morphism of site.
Let $S\in\Var(\mathbb C)$.Denote for short $\Var(\mathbb C)^{2,(sm)}/S$ 
either the category $\Var(\mathbb C)^2/S$ or the category $\Var(\mathbb C)^{2,sm}/S$. 
Denote for short $\Var(\mathbb C)^{2,(sm)pr}/S$ 
either the category $\Var(\mathbb C)^{2,pr}/S$ or the category $\Var(\mathbb C)^{2,smpr}/S$. Denote by
\begin{eqnarray*}
p_a:\Var(\mathbb C)^{2,(sm)pr}/S\to\Var(\mathbb C)^{2,(sm)pr}/S, \\ 
(Y\times S,Z)/S=((Y\times S,Z),p_S)\mapsto 
(Y\times S\times\mathbb A^1,Z\times\mathbb A^1)/S=((Y\times S\times\mathbb A^1,Z\times\mathbb A^1,p_S\circ p_{Y\times S}), \\ 
(g:(Y\times S,Z)/S\to (Y'\times S,Z')/S)\mapsto 
((g\times I_{\mathbb A^1}):(Y\times S\times\mathbb A^1,Z\times\mathbb A^1)/S\to 
(Y'\times S\times\mathbb A^1,Z'\times\mathbb A^1)/S)
\end{eqnarray*}
the projection functor and again by $p_a:\Var(\mathbb C)^{2,(sm)pr}/S\to\Var(\mathbb C)^{2,(sm)pr}/S$
the corresponding morphism of site.

\begin{defi}\label{a1loc12def}
\begin{itemize}
\item[(i0)]A complex $F\in C(\Var(\mathbb C)^{2,(sm)}/S)$ is said to be $\mathbb A^1$ homotopic if 
$\ad(p_a^*,p_{a*})(F):F\to p_{a*}p_a^*F$ is an homotopy equivalence.
\item[(i0)']A complex $F\in C(\Var(\mathbb C)^{2,(sm)pr}/S)$ is said to be $\mathbb A^1$ homotopic if 
$\ad(p_a^*,p_{a*})(F):F\to p_{a*}p_a^*F$ is an homotopy equivalence.
\item[(i)] A complex  $F\in C(\Var(\mathbb C)^{2,(sm)}/S)$, is said to be $\mathbb A^1$ invariant 
if for all $(X,Z)/S\in\Var(\mathbb C)^{2,(sm)}/S$ 
\begin{equation*}
F(p_X):F((X,Z)/S)\to F((X\times\mathbb A^1,(Z\times\mathbb A^1))/S) 
\end{equation*}
is a quasi-isomorphism, where $p_X:(X\times\mathbb A^1,(Z\times\mathbb A^1))\to (X,Z)$ is the projection.
Obviously, if a complex  $F\in C(\Var(\mathbb C)^{2,(sm)}/S)$ is $\mathbb A^1$ homotopic, then it is $\mathbb A^1$ invariant.
\item[(i)'] A complex $G\in C(\Var(\mathbb C)^{2,(sm)pr}/S)$, 
is said to be $\mathbb A^1$ invariant if for all $(Y\times S,Z)/S\in\Var(\mathbb C)^{2,(sm)pr}/S$ 
\begin{equation*}
G(p_{Y\times S}):G((Y\times S,Z)/S)\to G((Y\times\mathbb A^1\times S,(Z\times\mathbb A^1))/S) 
\end{equation*}
is a quasi-isomorphism of abelian group.
Obviously, if a complex  $F\in C(\Var(\mathbb C)^{2,(sm)pr}/S)$ is $\mathbb A^1$ homotopic, then it is $\mathbb A^1$ invariant.
\item[(ii)]Let $\tau$ a topology on $\Var(\mathbb C)$. 
A complex $F\in C(\Var(\mathbb C)^{2,(sm)}/S)$ is said to be $\mathbb A^1$ local 
for the $\tau$ topology induced on $\Var(\mathbb C)^2/S$, 
if for an (hence every) $\tau$ local equivalence $k:F\to G$ with $k$ injective and $G\in C(\Var(\mathbb C)^{2,(sm)}/S)$ $\tau$ fibrant,
e.g. $k:F\to E_{\tau}(F)$, $G$ is $\mathbb A^1$ invariant.
\item[(ii)']Let $\tau$ a topology on $\Var(\mathbb C)$. 
A complex $F\in C(\Var(\mathbb C)^{2,(sm)pr}/S)$ is said to be $\mathbb A^1$ local 
for the $\tau$ topology induced on $\Var(\mathbb C)^{2,pr}/S$, 
if for an (hence every) $\tau$ local equivalence $k:F\to G$ with $k$ injective and $G\in C(\Var(\mathbb C)^{2,(sm)pr}/S)$ $\tau$ fibrant,
e.g. $k:F\to E_{\tau}(F)$, $G$ is $\mathbb A^1$ invariant.
\item[(iii)] A morphism $m:F\to G$ with $F,G\in C(\Var(\mathbb C)^{2,(sm)}/S)$ is said to an $(\mathbb A^1,et)$ local equivalence 
if for all $H\in C(\Var(\mathbb C)^{2,(sm)}/S)$ which is $\mathbb A^1$ local for the etale topology
\begin{equation*}
\Hom(L(m),E_{et}(H)):\Hom(L(G),E_{et}(H))\to\Hom(L(F),E_{et}(H)) 
\end{equation*}
is a quasi-isomorphism.
\item[(iii)'] A morphism $m:F\to G$ with $F,G\in C(\Var(\mathbb C)^{2,(sm)pr}/S)$ is said to an $(\mathbb A^1,et)$ local equivalence 
if for all $H\in C(\Var(\mathbb C)^{2,(sm)pr}/S)$ which is $\mathbb A^1$ local for the etale topology
\begin{equation*}
\Hom(L(m),E_{et}(H)):\Hom(L(G),E_{et}(H))\to\Hom(L(F),E_{et}(H)) 
\end{equation*}
is a quasi-isomorphism.
\end{itemize}
\end{defi}

\begin{prop}\label{ca1Var12}
\begin{itemize}
\item[(i)]Let $S\in\Var(\mathbb C)$.
Then for $F\in C(\Var(\mathbb C)^{2,(sm)}/S)$, $C_*F$ is $\mathbb A^1$ local for the etale topology
and  $c(F):F\to C_*F$ is an equivalence $(\mathbb A^1,et)$ local.
\item[(i)']Let $S\in\Var(\mathbb C)$.
Then for $F\in C(\Var(\mathbb C)^{2,(sm)pr}/S)$, $C_*F$ is $\mathbb A^1$ local for the etale topology
and  $c(F):F\to C_*F$ is an equivalence $(\mathbb A^1,et)$ local.
\item[(ii)]A morphism $m:F\to G$ with $F,G\in C(\Var(\mathbb C)^{2,(sm)}/S)$ is an $(\mathbb A^1,et)$ local equivalence
if and only if $a_{et}H^nC_*\Cone(m)=0$ for all $n\in\mathbb Z$.
\item[(ii)']A morphism $m:F\to G$ with $F,G\in C(\Var(\mathbb C)^{2,(sm)pr}/S)$ is an $(\mathbb A^1,et)$ local equivalence
if and only if $a_{et}H^nC_*\Cone(m)=0$ for all $n\in\mathbb Z$.
\item[(iii)]A morphism $m:F\to G$ with $F,G\in C(\Var(\mathbb C)^{2,(sm)}/S)$ is an $(\mathbb A^1,et)$ local equivalence
if and only if there exists 
\begin{eqnarray*}
\left\{(X_{1,\alpha},Z_{1,\alpha})/S,\alpha\in\Lambda_1\right\},\ldots,
\left\{(X_{r,\alpha},Z_{r,\alpha})/S,\alpha\in\Lambda_r\right\}\subset\Var(\mathbb C)^{2,(sm)}/S
\end{eqnarray*}
such that we have in $\Ho_{et}(C(\Var(\mathbb C)^{2,(sm)}/S))$
\begin{eqnarray*}
\Cone(m)\xrightarrow{\sim}\Cone(\oplus_{\alpha\in\Lambda_1}
\Cone(\mathbb Z((X_{1,\alpha}\times\mathbb A^1,Z_{1,\alpha}\times\mathbb A^1)/S)\to\mathbb Z((X_{1,\alpha},Z_{1,\alpha})/S)) \\
\to\cdots\to\oplus_{\alpha\in\Lambda_r}
\Cone(\mathbb Z((X_{r,\alpha}\times\mathbb A^1,Z_{r,\alpha}\times\mathbb A^1)/S)\to\mathbb Z((X_{r,\alpha},Z_{r,\alpha})/S)))
\end{eqnarray*}
\item[(iii)']A morphism $m:F\to G$ with $F,G\in C(\Var(\mathbb C)^{2,(sm)pr}/S)$ is an $(\mathbb A^1,et)$ local equivalence
if and only if there exists 
\begin{eqnarray*}
\left\{(Y_{1,\alpha}\times S,Z_{1,\alpha})/S,\alpha\in\Lambda_1\right\},\ldots,
\left\{(Y_{r,\alpha}\times S,Z_{r,\alpha})/S,\alpha\in\Lambda_r\right\}\subset\Var(\mathbb C)^{2,(sm)pr}/S
\end{eqnarray*}
such that we have in $\Ho_{et}(C(\Var(\mathbb C)^{2,(sm)}/S))$
\begin{eqnarray*}
\Cone(m)\xrightarrow{\sim}\Cone(\oplus_{\alpha\in\Lambda_1}
\Cone(\mathbb Z((Y_{1,\alpha}\times\mathbb A^1\times S,Z_{1,\alpha}\times\mathbb A^1)/S)\to
\mathbb Z((Y_{1,\alpha}\times S,Z_{1,\alpha})/S)) \\
\to\cdots\to\oplus_{\alpha\in\Lambda_r}
\Cone(\mathbb Z((Y_{r,\alpha}\times\mathbb A^1\times S,Z_{r,\alpha}\times\mathbb A^1)/S)\to
\mathbb Z((Y_{r,\alpha}\times S,Z_{r,\alpha})/S)))
\end{eqnarray*}
\end{itemize}
\end{prop}

\begin{proof}
Standard : see Ayoub's thesis section 4 for example. 
Indeed, for (iii), by definition, if $\Cone(m)$ is of the given form, 
then it is an equivalence $(\mathbb A^1,et)$ local, on the other hand if $m$ is an equivalence $(\mathbb A^1,et)$ local,
we consider the commutative diagram
\begin{equation*}
\xymatrix{F\ar[r]^{c(F)}\ar[d]^{m} & C_*F\ar[d]^{c_*m} \\ G\ar[r]^{c(G)} & C_*G}
\end{equation*}
to deduce that $\Cone(m)$ is of the given form.
\end{proof}

\begin{defiprop}\label{projmodstr12}
Let $S\in\Var(\mathbb C)$.
\begin{itemize}
\item[(i)]With the weak equivalence the $(\mathbb A^1,et)$ local equivalence and 
the fibration the epimorphism with $\mathbb A^1_S$ local and etale fibrant kernels gives
a model structure on  $C(\Var(\mathbb C)^{2,(sm)}/S)$ : the left bousfield localization
of the projective model structure of $C(\Var(\mathbb C)^{2,(sm)}/S)$. 
We call it the projective $(\mathbb A^1,et)$ model structure.
\item[(ii)]With the weak equivalence the $(\mathbb A^1,et)$ local equivalence and 
the fibration the epimorphism with $\mathbb A^1_S$ local and etale fibrant kernels gives
a model structure on  $C(\Var(\mathbb C)^{2,(sm)pr}/S)$ : the left bousfield localization
of the projective model structure of $C(\Var(\mathbb C)^{2,(sm)pr}/S)$. 
We call it the projective $(\mathbb A^1,et)$ model structure.
\end{itemize}
\end{defiprop}

\begin{proof}
Similar to the proof of proposition \ref{projmodstr}.
\end{proof}

We have, similarly to the case of single varieties the following :

\begin{prop}\label{g12}
Let $g:T\to S$ a morphism with $T,S\in\Var(\mathbb C)$.
\begin{itemize}
\item[(i)] The adjonction $(g^*,g_*):C(\Var(\mathbb C)^{2,(sm)}/S)\leftrightarrows C(\Var(\mathbb C)^{2,(sm)}/T)$
is a Quillen adjonction for the projective $(\mathbb A^1,et)$ model structure (see definition-proposition \ref{projmodstr12})
\item[(i)'] The functor $g^*:C(\Var(\mathbb C)^{2,(sm)}/S)\to C(\Var(\mathbb C)^{2,(sm)}/T)$
sends quasi-isomorphism to quasi-isomorphism,
sends equivalence Zariski local to equivalence Zariski local, and equivalence etale local to equivalence etale local,
sends $(\mathbb A^1,et)$ local equivalence to $(\mathbb A^1,et)$ local equivalence.
\item[(ii)] The adjonction $(g^*,g_*):C(\Var(\mathbb C)^{2,(sm)pr}/S)\leftrightarrows C(\Var(\mathbb C)^{2,(sm)pr}/T)$
is a Quillen adjonction for the projective $(\mathbb A^1,et)$ model structure (see definition-proposition \ref{projmodstr12})
\item[(ii)'] The functor $g^*:C(\Var(\mathbb C)^{2,(sm)pr}/S)\to C(\Var(\mathbb C)^{2,(sm)pr}/T)$
sends quasi-isomorphism to quasi-isomorphism,
sends equivalence Zariski local to equivalence Zariski local, and equivalence etale local to equivalence etale local,
sends $(\mathbb A^1,et)$ local equivalence to $(\mathbb A^1,et)$ local equivalence.
\end{itemize}
\end{prop}

\begin{proof}
\noindent(i):Follows immediately from definition.
\noindent(i)': Since the functor $g^*$ preserve epimorphism and also monomorphism (the colimits involved being filetered),
$g^*$ sends quasi-isomorphism to quasi-isomorphism. Hence it preserve Zariski and etale local equivalence.
The fact that it preserve $(\mathbb A^1,et)$ local equivalence then follows similarly to the single case by the fact
that $g_*$ preserve by definition $\mathbb A^1$ equivariant presheaves.

\noindent(ii) and (ii)': Similar to (i) and (i)'.
\end{proof}

\begin{prop}\label{rho12}
Let $S\in\Var(\mathbb C)$. 
\begin{itemize}
\item[(i)] The adjonction $(\rho_S^*,\rho_{S*}):C(\Var(\mathbb C)^{2,sm}/S)\leftrightarrows C(\Var(\mathbb C)^2/S)$
is a Quillen adjonction for the $(\mathbb A^1,et)$ projective model structure.
\item[(i)']The functor $\rho_{S*}:C(\Var(\mathbb C)^2/S)\to C(\Var(\mathbb C)^{2,sm}/S)$
sends quasi-isomorphism to quasi-isomorphism,
sends equivalence Zariski local to equivalence Zariski local, and equivalence etale local to equivalence etale local,
sends $(\mathbb A^1,et)$ local equivalence to $(\mathbb A^1,et)$ local equivalence.
\item[(ii)] The adjonction $(\rho_S^*,\rho_{S*}):C(\Var(\mathbb C)^{2,smpr}/S)\leftrightarrows C(\Var(\mathbb C)^{2,pr}/S)$
is a Quillen adjonction for the $(\mathbb A^1,et)$ projective model structure.
\item[(ii)']The functor $\rho_{S*}:C(\Var(\mathbb C)^{2,pr}/S)\to C(\Var(\mathbb C)^{2,smpr}/S)$
sends quasi-isomorphism to quasi-isomorphism,
sends equivalence Zariski local to equivalence Zariski local, and equivalence etale local to equivalence etale local,
sends $(\mathbb A^1,et)$ local equivalence to $(\mathbb A^1,et)$ local equivalence.
\end{itemize}
\end{prop}

\begin{proof}
Similar to the proof of proposition \ref{rho1}.
\end{proof}

\begin{prop}\label{mu12}
Let $S\in\Var(\mathbb C)$. 
\begin{itemize}
\item[(i)] The adjonction $(\mu_S^*,\mu_{S*}):C(\Var(\mathbb C)^{2,pr}/S)\leftrightarrows C(\Var(\mathbb C)^2/S)$
is a Quillen adjonction for the $(\mathbb A^1,et)$ projective model structure.
\item[(i)']The functor $\mu_{S*}:C(\Var(\mathbb C)^2/S)\to C(\Var(\mathbb C)^{2,pr}/S)$
sends quasi-isomorphism to quasi-isomorphism,
sends equivalence Zariski local to equivalence Zariski local, and equivalence etale local to equivalence etale local,
sends $(\mathbb A^1,et)$ local equivalence to $(\mathbb A^1,et)$ local equivalence.
\item[(ii)] The adjonction $(\mu_S^*,\mu_{S*}):C(\Var(\mathbb C)^{2,smpr}/S)\leftrightarrows C(\Var(\mathbb C)^{2,pr}/S)$
is a Quillen adjonction for the $(\mathbb A^1,et)$ projective model structure.
\item[(ii)']The functor $\mu_{S*}:C(\Var(\mathbb C)^{2,sm}/S)\to C(\Var(\mathbb C)^{2,smpr}/S)$
sends quasi-isomorphism to quasi-isomorphism,
sends equivalence Zariski local to equivalence Zariski local, and equivalence etale local to equivalence etale local,
sends $(\mathbb A^1,et)$ local equivalence to $(\mathbb A^1,et)$ local equivalence.
\end{itemize}
\end{prop}

\begin{proof}
Similar to the proof of proposition \ref{rho1}. 
Indeed, for (i)' or (ii)', if $m:F\to G$ with $F,G\in C(\Var(\mathbb C)^{2,(sm)})$ is an equivalence $(\mathbb A^1,et)$
local then (see proposition \ref{ca1Var12}), there exists 
\begin{equation*}
\left\{(X_{1,\alpha},Z_{1,\alpha})/S,\alpha\in\Lambda_1\right\},\ldots,
\left\{(X_{r,\alpha},Z_{r,\alpha})/S,\alpha\in\Lambda_r\right\}
\subset\Var(\mathbb C)^{2,(sm)}/S 
\end{equation*}
such that we have in $\Ho_{et}(C(\Var(\mathbb C)^{2,(sm)}/S))$
\begin{eqnarray*}
\Cone(m)\xrightarrow{\sim}\Cone(\oplus_{\alpha\in\Lambda_1}
\Cone(\mathbb Z((X_{1,\alpha}\times\mathbb A^1,Z_{1,\alpha}\times\mathbb A^1)/S)
\to\mathbb Z((X_{1,\alpha},Z_{1,\alpha})/S)) \\
\to\cdots\to\oplus_{\alpha\in\Lambda_r}
\Cone(\mathbb Z((X_{r,\alpha}\times\mathbb A^1,Z_{r,\alpha}\times\mathbb A^1)/S)\to
\mathbb Z((X_{r,\alpha},Z_{r,\alpha})/S))) \\
\xrightarrow{\sim}\Cone(
\Cone(\oplus_{\alpha\in\Lambda_1}\mathbb Z((X_{1,\alpha},Z_{1,\alpha})/S)\otimes\mathbb Z(\mathbb A^1,\mathbb A^1)/S\to
\oplus_{\alpha\in\Lambda_1}\mathbb Z((X_{1,\alpha},Z_{1,\alpha})/S)) \\
\to\cdots\to
\Cone(\oplus_{\alpha\in\Lambda_r}\mathbb Z((X_{r,\alpha},Z_{r,\alpha})/S)\otimes\mathbb Z((\mathbb A^1,\mathbb A^1)/S)\to
\oplus_{\alpha\in\Lambda_r}\mathbb Z((X_{r,\alpha},Z_{r,\alpha})/S)),
\end{eqnarray*}
this gives in $\Ho_{et}(C(\Var(\mathbb C)^{2,(sm)pr}/S))$
\begin{eqnarray*}
\Cone(\mu_{S*}m)\xrightarrow{\sim}\Cone( \\
\Cone((L\mu_{S*}\oplus_{\alpha\in\Lambda_1}\mathbb Z((X_{1,\alpha},Z_{1,\alpha})/S))\otimes\mathbb Z((\mathbb A^1,\mathbb A^1)/S)
\to(L\mu_{S*}\oplus_{\alpha\in\Lambda_1}\mathbb Z((X_{1,\alpha},Z_{1,\alpha})/S)) \\
\to\cdots\to
\Cone((L\mu_{S*}\oplus_{\alpha\in\Lambda_r}\mathbb Z((X_{r,\alpha},Z_{r,\alpha})/S))\otimes\mathbb Z((\mathbb A^1,\mathbb A^1)/S)
\to(L\mu_{S*}\oplus_{\alpha\in\Lambda_r}\mathbb Z((X_{1,\alpha},Z_{1,\alpha})/S))))
\end{eqnarray*}
hence $\mu_{S*}m:\mu_{S*}F\to\mu_{S*}G$ is an equivalence $(\mathbb A^1,et)$ local.
\end{proof}

We also have

\begin{prop}\label{Gra1}
Let $S\in\Var(\mathbb C)$. 
\begin{itemize}
\item[(i)] The adjonction $(\Gr_S^{12*},\Gr_{S*}^{12}):C(\Var(\mathbb C)/S)\leftrightarrows C(\Var(\mathbb C)^{2,pr}/S)$
is a Quillen adjonction for the $(\mathbb A^1,et)$ projective model structure.
\item[(ii)] The adjonction $(\Gr_S^{12*}\Gr_{S*}^{12}:C(\Var(\mathbb C)^{sm}/S)\leftrightarrows C(\Var(\mathbb C)^{2,smpr}/S)$
is a Quillen adjonction for the $(\mathbb A^1,et)$ projective model structure.
\end{itemize}
\end{prop}

\begin{proof}
Immediate from definition.
\end{proof}

\begin{itemize}
\item For $f:X\to S$ a morphism with $X,S\in\Var(\mathbb C)$ and $Z\subset X$ a closed subset, we denote
$\mathbb Z^{tr}((X,Z)/S)\in\PSh(\Var(\mathbb C)^2/S)$ the presheaf given by
\begin{itemize}
\item for $(X',Z')/S\in\Var(\mathbb C)^2/S$, with $X'$ irreducible, 
\begin{equation*}
\mathbb Z^{tr}((X,Z)/S)((X',Z')/S):=\left\{\alpha\in\mathcal Z^{fs/X}(X'\times_S X),s.t. p_X(p_{X'}^{-1}(Z'))\subset Z\right\}
\subset\mathcal Z_{d_{X'}}(X'\times_S X) 
\end{equation*}
\item for $g:(X_2,Z_2)/S\to (X_1,Z_1)/S$ a morphism, with $(X_1,Z_1)/S,(X_2,Z_2)/S\in\Var(\mathbb C)^2/S$,
\begin{equation*}
\mathbb Z^{tr}((X,Z)/S)(g):\mathbb Z^{tr}((X,Z)/S)((X_1,Z_1)/S)\to\mathbb Z^{tr}((X,Z)/S)((X_2,Z_2)/S), \;
\alpha\mapsto (g\times I)^{-1}(\alpha)
\end{equation*}
with $g\times I:X_2\times_S X\to X_1\times_S X$.
\end{itemize}
\item For $f:X\to S$ a morphism with $X,S\in\Var(\mathbb C)$, $Z\subset X$ a closed subset and $r\in\mathbb N$, we denote
$\mathbb Z^{equir}((X,Z)/S)\in\PSh(\Var(\mathbb C)^2/S)$ the presheaf given by
\begin{itemize}
\item for $(X',Z')/S\in\Var(\mathbb C)^2/S$, with $X'$ irreducible,
\begin{equation*}
\mathbb Z^{equir}((X,Z)/S)((X',Z')/S):=\left\{\alpha\in\mathcal Z^{equir/X}(X'\times_S X),s.t. p_X(p_{X'}^{-1}(Z'))\right\}
\subset\mathcal Z_{d_{X'}}(X'\times_S X) 
\end{equation*}
\item for $g:(X_2,Z_2)/S\to (X_1,Z_1)/S$ a morphism, with $(X_1,Z_1)/S,(X_2,Z_2)/S\in\Var(\mathbb C)^2/S$,
\begin{equation*}
\mathbb Z^{equir}((X,Z)/S)(g):\mathbb Z^{equir}((X,Z)/S)((X_1,Z_1)/S)\to\mathbb Z^{equir}((X,Z)/S)((X_2,Z_2)/S), \;
\alpha\mapsto (g\times I)^{-1}(\alpha)
\end{equation*}
with $g\times I:X_2\times_S X\to X_1\times_S X$.
\end{itemize}
\item Let $S\in\Var(\mathbb C)$. We denote by 
$\mathbb Z_S(d):=\mathbb Z^{equi0}((S\times\mathbb A^d,S\times\mathbb A^d)/S)[-2d]$
the Tate twist. For $F\in C(\Var(\mathbb C)^2/S)$, we denote by $F(d):=F\otimes\mathbb Z_S(d)$.
\end{itemize}

For $S\in\Var(\mathbb C)$, let $\Cor(\Var(\mathbb C)^{2,(sm)}/S)$ be the category 
\begin{itemize}
\item whose objects are those of $\Var(\mathbb C)^{2,(sm)}/S$, i.e. 
$(X,Z)/S=((X,Z),h)$, $h:X\to S$ with $X\in\Var(\mathbb C)$, $Z\subset X$ a closed subset, 
\item whose morphisms $\alpha:(X',Z)/S=((X',Z),h_1)\to (X,Z)/S=((X,Z),h_2)$ 
is finite correspondence that is $\alpha\in\oplus_i\mathbb Z^{tr}((X_i,Z)/S)((X',Z')/S)$, 
where $X'=\sqcup_i X'_i$, with $X'_i$ connected, 
the composition being defined in the same way as the morphism $\Cor(\Var(\mathbb C)^{(sm)}/S)$. 
\end{itemize}
We denote by 
$\Tr(S):\Cor(\Var(\mathbb C)^{2,(sm)}/S)\to\Var(\mathbb C)^{2,(sm)}/S$ 
the morphism of site
given by the inclusion functor
$\Tr(S):\Var(\mathbb C)^{2,(sm)}/S\hookrightarrow\Cor(\Var(\mathbb C)^{2,(sm)}/S)$
It induces an adjonction
\begin{equation*}
(\Tr(S)^*\Tr(S)_*):C(\Var(\mathbb C)^{2,(sm)}/S)\leftrightarrows C(\Cor(\Var(\mathbb C)^{2,(sm)}/S))
\end{equation*}
A complex of preheaves $G\in C(\Var(\mathbb C)^{2,(sm)}/S)$ is said to admit transferts
if it is in the image of the embedding
\begin{equation*}
\Tr(S)_*:C(\Cor(\Var(\mathbb C)^{2,(sm)}/S)\hookrightarrow C(\Var(\mathbb C)^{2,(sm)}/S),
\end{equation*}
that is $G=\Tr(S)_*\Tr(S)^*G$.
We then have the full subcategory $\Cor(\Var(\mathbb C)^{2,(sm)pr}/S)\subset\Cor(\Var(\mathbb C)^{2,(sm)}/S)$
consisting of the objects of $\Var(\mathbb C)^{2,(sm)pr}/S)$.
We have the adjonction
\begin{equation*}
(\Tr(S)^*\Tr(S)_*):C(\Var(\mathbb C)^{2,(sm)pr}/S)\leftrightarrows C(\Cor(\Var(\mathbb C)^{2,(sm)pr}/S))
\end{equation*}
A complex of preheaves $G\in C(\Var(\mathbb C)^{2,(sm)pr}/S)$ is said to admit transferts
if it is in the image of the embedding
\begin{equation*}
\Tr(S)_*:C(\Cor(\Var(\mathbb C)^{2,(sm)pr}/S)\hookrightarrow C(\Var(\mathbb C)^{2,(sm)pr}/S),
\end{equation*}
that is $G=\Tr(S)_*\Tr(S)^*G$.

Let $S\in\Var(\mathbb C)$. Let $S=\cup_{i=1}^l S_i$ an open affine cover and denote by $S_I=\cap_{i\in I} S_i$.
Let $i_i:S_i\hookrightarrow\tilde S_i$ closed embeddings, with $\tilde S_i\in\Var(\mathbb C)$.
\begin{itemize}
\item For $(G_I,K_{IJ})\in C(\Var(\mathbb C)^{2,(sm)}/(\tilde S_I)^{op})$ and 
$(H_I,T_{IJ})\in C(\Var(\mathbb C)^{2,(sm)}/(\tilde S_I))$, we denote
\begin{eqnarray*}
\mathcal Hom((G_I,K_{IJ}),(H_I,T_{IJ})):=(\mathcal Hom(G_I,H_I),u_{IJ}((G_I,K_{IJ}),(H_I,T_{IJ})))
\in C(\Var(\mathbb C)^{2,(sm)}/(\tilde S_I))
\end{eqnarray*}
with
\begin{eqnarray*}
u_{IJ}((G_I,K_{IJ})(H_I,T_{IJ})):\mathcal Hom(G_I,H_I) \\
\xrightarrow{\ad(p_{IJ}^*,p_{IJ*})(-)}p_{IJ*}p_{IJ}^*\mathcal Hom(G_I,H_I)
\xrightarrow{T(p_{IJ},hom)(-,-)}p_{IJ*}\mathcal Hom(p_{IJ}^*G_I,p_{IJ}^*H_I) \\
\xrightarrow{\mathcal Hom(p_{IJ}^*G_I,T_{IJ})}p_{IJ*}\mathcal Hom(p_{IJ}^*G_I,H_J)
\xrightarrow{\mathcal Hom(K_{IJ},H_J)}p_{IJ*}\mathcal Hom(G_J,H_J).
\end{eqnarray*}
This gives in particular the functor
\begin{eqnarray*}
C(\Var(\mathbb C)^{2,(sm)}/(\tilde S_I))\to C(\Var(\mathbb C)^{2,(sm)}/(\tilde S_I)^{op}),
(H_I,T_{IJ})\mapsto(H_I,T_{IJ}).
\end{eqnarray*}
\item For $(G_I,K_{IJ})\in C(\Var(\mathbb C)^{2,(sm)pr}/(\tilde S_I)^{op})$ and 
$(H_I,T_{IJ})\in C(\Var(\mathbb C)^{2,(sm)pr}/(\tilde S_I))$, we denote
\begin{eqnarray*}
\mathcal Hom((G_I,K_{IJ}),(H_I,T_{IJ})):=(\mathcal Hom(G_I,H_I),u_{IJ}((G_I,K_{IJ}),(H_I,T_{IJ})))
\in C(\Var(\mathbb C)^{2,(sm)pr}/(\tilde S_I))
\end{eqnarray*}
with
\begin{eqnarray*}
u_{IJ}((G_I,K_{IJ})(H_I,T_{IJ})):\mathcal Hom(G_I,H_I) \\
\xrightarrow{\ad(p_{IJ}^*,p_{IJ*})(-)}p_{IJ*}p_{IJ}^*\mathcal Hom(G_I,H_I)
\xrightarrow{T(p_{IJ},hom)(-,-)}p_{IJ*}\mathcal Hom(p_{IJ}^*G_I,p_{IJ}^*H_I) \\
\xrightarrow{\mathcal Hom(p_{IJ}^*G_I,T_{IJ})}p_{IJ*}\mathcal Hom(p_{IJ}^*G_I,H_J)
\xrightarrow{\mathcal Hom(K_{IJ},H_J)}p_{IJ*}\mathcal Hom(G_J,H_J).
\end{eqnarray*}
This gives in particular the functor
\begin{eqnarray*}
C(\Var(\mathbb C)^{2,(sm)pr}/(\tilde S_I))\to C(\Var(\mathbb C)^{2,(sm)pr}/(\tilde S_I)^{op}),
(H_I,T_{IJ})\mapsto(H_I,T_{IJ}).
\end{eqnarray*}
\end{itemize}
Let $S\in\Var(\mathbb C)$. Let $S=\cup_{i=1}^l S_i$ an open affine cover and denote by $S_I=\cap_{i\in I} S_i$.
Let $i_i:S_i\hookrightarrow\tilde S_i$ closed embeddings, with $\tilde S_i\in\Var(\mathbb C)$.
The functors $p_a$ naturally extend to functors
\begin{eqnarray*}
p_a:\Var(\mathbb C)^{2,(sm)}/(\tilde S_I)\to\Var(\mathbb C)^{2,(sm)}/(\tilde S_I), \\ 
((X,Z)/\tilde S_I,u_{IJ})\mapsto ((X\times\mathbb A^1,Z\times\mathbb A^1)/\tilde S_I,u_{IJ}\times I), \\
(g:((X,Z)/\tilde S_I,u_{IJ})\to ((X',Z')/\tilde S_I,u_{IJ}))\mapsto \\
((g\times I_{\mathbb A^1}):((X\times\mathbb A^1,Z\times\mathbb A^1)/\tilde S_I,u_{IJ}\times I)\to 
((X'\times\mathbb A^1,Z'\times\mathbb A^1)/\tilde S_I,u_{IJ}\times I))
\end{eqnarray*}
the projection functor and again by $p_a:\Var(\mathbb C)^{2,(sm)}/(\tilde S_I)\to\Var(\mathbb C)^{2,(sm)}/(\tilde S_I)$
the corresponding morphism of site, and
\begin{eqnarray*}
p_a:\Var(\mathbb C)^{2,(sm)pr}/(\tilde S_I)\to\Var(\mathbb C)^{2,(sm)pr}/(\tilde S_I), \\ 
((Y\times\tilde S_I,Z)/\tilde S_I,u_{IJ})\mapsto 
((Y\times\tilde S_I\times\mathbb A^1,Z\times\mathbb A^1)/\tilde S_I,u_{IJ}\times I), \\ 
(g:((Y\times\tilde S_I,Z)/\tilde S_I,u_{IJ})\to ((Y'\times\tilde S_I,Z')/\tilde S_I,u_{IJ}))\mapsto \\
((g\times I_{\mathbb A^1}):((Y\times\tilde S_I\times\mathbb A^1,Z\times\mathbb A^1)/\tilde S_I,u_{IJ}\times I),  
((Y'\times\tilde S_I\times\mathbb A^1,Z'\times\mathbb A^1)/\tilde S_I,u_{IJ}\times I)), 
\end{eqnarray*}
the projection functor and again by $p_a:\Var(\mathbb C)^{2,(sm)pr}/(\tilde S_I)\to\Var(\mathbb C)^{2,(sm)pr}/(\tilde S_I)$
the corresponding morphism of site.
These functors also gives the morphisms of sites
$p_a:\Var(\mathbb C)^{2,(sm)}/(\tilde S_I)^{op}\to\Var(\mathbb C)^{2,(sm)}/(\tilde S_I)^{op}$ and
$p_a:\Var(\mathbb C)^{2,(sm)pr}/(\tilde S_I)^{op}\to\Var(\mathbb C)^{2,(sm)pr}/(\tilde S_I)^{op}$.

\begin{defi}\label{a1loc12defIJ}
Let $S\in\Var(\mathbb C)$. Let $S=\cup_{i=1}^l S_i$ an open affine cover and denote by $S_I=\cap_{i\in I} S_i$.
Let $i_i:S_i\hookrightarrow\tilde S_i$ closed embeddings, with $\tilde S_i\in\Var(\mathbb C)$.
\begin{itemize}
\item[(i0)]A complex $(F_I,u_{IJ})\in C(\Var(\mathbb C)^{2,(sm)}/(\tilde S_I))$ is said to be $\mathbb A^1$ homotopic if 
$\ad(p_a^*,p_{a*})((F_I,u_{IJ})):(F_I,u_{IJ})\to p_{a*}p_a^*(F_I,u_{IJ})$ is an homotopy equivalence.
\item[(i0)']A complex $(F_I,u_{IJ})\in C(\Var(\mathbb C)^{2,(sm)pr}/(\tilde S_I))$ is said to be $\mathbb A^1$ homotopic if 
$\ad(p_a^*,p_{a*})((F_I,u_{IJ})):(F_I,u_{IJ})\to p_{a*}p_a^*(F_I,u_{IJ})$ is an homotopy equivalence.
\item[(i)] A complex  $(F_I,u_{IJ})\in C(\Var(\mathbb C)^{2,(sm)}/(\tilde S_I))$ is said to be $\mathbb A^1$ invariant 
if for all $((X_I,Z_I)/\tilde S_I,s_{IJ})\in\Var(\mathbb C)^{2,(sm)}/(\tilde S_I)$ 
\begin{equation*}
(F_I(p_{X_I})):(F_I((X_I,Z_I)/\tilde S_I),F_J(s_{IJ})\circ u_{IJ}(-)\to 
(F_I((X_I\times\mathbb A^1,(Z_I\times\mathbb A^1))/\tilde S_I),F_J(s_{IJ}\times I)\circ u_{IJ}(-)) 
\end{equation*}
is a quasi-isomorphism, where $p_{X_I}:(X_I\times\mathbb A^1,(Z_I\times\mathbb A^1))\to (X_I,Z_I)$ are the projection,
and $s_{IJ}:(X_I\times\tilde S_{J\backslash I},Z_I)/\tilde S_J\to(X_J,Z_J)/\tilde S_J$.
Obviously a complex $(F_I,u_{IJ})\in C(\Var(\mathbb C)^{2,(sm)}/(\tilde S_I))$ is $\mathbb A^1$ invariant
if and only if all the $F_I$ are $\mathbb A^1$ invariant.
\item[(i)'] A complex  $(G_I,u_{IJ})\in C(\Var(\mathbb C)^{2,(sm)pr}/(\tilde S_I))$ is said to be $\mathbb A^1$ invariant 
if for all $((Y\times\tilde S_I,Z_I)/\tilde S_I,s_{IJ})\in\Var(\mathbb C)^{2,(sm)pr}/(\tilde S_I)$ 
\begin{eqnarray*}
(G_I(p_{Y\times\tilde S_I})):
(G_I((Y\times\tilde S_I,Z_I)/\tilde S_I),G_J(s_{IJ})\circ u_{IJ}(-))\to \\
(G_I((Y\times\tilde S_I\times\mathbb A^1,(Z_I\times\mathbb A^1))/\tilde S_I),G_J(s_{IJ}\times I)\circ u_{IJ}(-)) 
\end{eqnarray*}
is a quasi-isomorphism. 
Obviously a complex  $(G_I,u_{IJ})\in C(\Var(\mathbb C)^{2,(sm)pr}/(\tilde S_I))$ is $\mathbb A^1$ invariant
if and only if all the $G_I$ are $\mathbb A^1$ invariant.
\item[(ii)]Let $\tau$ a topology on $\Var(\mathbb C)$. 
A complex $F=(F_I,u_{IJ})\in C(\Var(\mathbb C)^{2,(sm)}/(\tilde S_I))$ is said to be $\mathbb A^1$ local 
for the $\tau$ topology induced on $\Var(\mathbb C)^2/(\tilde S_I)$, 
if for an (hence every) $\tau$ local equivalence $k:F\to G$ with $k$ injective and 
$G=(G_I,v_{IJ})\in C(\Var(\mathbb C)^{2,(sm)}/(\tilde S_I))$ $\tau$ fibrant,
e.g. $k:(F_I,u_{IJ})\to (E_{\tau}(F_I),E(u_{IJ}))$, $G$ is $\mathbb A^1$ invariant.
\item[(ii)']Let $\tau$ a topology on $\Var(\mathbb C)$. 
A complex $F=(F_I,u_{IJ})\in C(\Var(\mathbb C)^{2,(sm)pr}/(\tilde S_I))$ is said to be $\mathbb A^1$ local 
for the $\tau$ topology induced on $\Var(\mathbb C)^2/(\tilde S_I)$, 
if for an (hence every) $\tau$ local equivalence $k:F\to G$ with $k$ injective and 
$G=(G_I,u_{IJ})\in C(\Var(\mathbb C)^{2,(sm)pr}/(\tilde S_I))$ $\tau$ fibrant,
e.g. $k:(F_I,u_{IJ})\to (E_{\tau}(F_I),E(u_{IJ}))$, $G$ is $\mathbb A^1$ invariant.
\item[(iii)] A morphism $m=(m_I):(F_I,u_{IJ})\to (G_I,v_{IJ})$ with 
$(F_I,u_{IJ}),(G_I,v_{IJ})\in C(\Var(\mathbb C)^{2,(sm)}/(\tilde S_I))$ 
is said to be an $(\mathbb A^1,et)$ local equivalence 
if for all $H=(H_I,w_{IJ})\in C(\Var(\mathbb C)^{2,(sm)}/(\tilde S_I))$ 
which is $\mathbb A^1$ local for the etale topology 
\begin{eqnarray*}
(\Hom(L(m_I),E_{et}(H_I))):\Hom(L(G_I,v_{IJ}),E_{et}(H_I,w_{IJ}))\to\Hom(L(F_I,u_{IJ}),E_{et}(H_I,w_{IJ})) 
\end{eqnarray*}
is a quasi-isomorphism.
Obviously, if a morphism $m=(m_I):(F_I,u_{IJ})\to (G_I,v_{IJ})$ with 
$(F_I,u_{IJ}),(G_I,u_{IJ})\in C(\Var(\mathbb C)^{2,(sm)}/(\tilde S_I))$ 
is an $(\mathbb A^1,et)$ local equivalence, 
then all the $m_I:F_I\to G_I$ are $(\mathbb A^1,et)$ local equivalence.
\item[(iii)'] A morphism $m=(m_I):(F_I,u_{IJ})\to (G_I,v_{IJ})$ with 
$(F_I,u_{IJ}),(G_I,v_{IJ})\in C(\Var(\mathbb C)^{2,(sm)pr}/(\tilde S_I))$ 
is said to be an $(\mathbb A^1,et)$ local equivalence 
if for all $(H_I,w_{IJ})\in C(\Var(\mathbb C)^{2,(sm)pr}/(\tilde S_I))$ 
which is $\mathbb A^1$ local for the etale topology
\begin{eqnarray*}
(\Hom(L(m_I),E_{et}(H_I))):\Hom(L(G_I,v_{IJ}),E_{et}(H_I,w_{IJ}))\to\Hom(L(F_I,u_{IJ}),E_{et}(H_I,w_{IJ})) 
\end{eqnarray*}
is a quasi-isomorphism.
Obviously, if a morphism $m=(m_I):(F_I,u_{IJ})\to (G_I,v_{IJ})$ with 
$(F_I,u_{IJ}),(G_I,u_{IJ})\in C(\Var(\mathbb C)^{2,(sm)pr}/(\tilde S_I))$ 
is an $(\mathbb A^1,et)$ local equivalence, 
then all the $m_I:F_I\to G_I$ are $(\mathbb A^1,et)$ local equivalence.
\item[(iv)] A morphism $m=(m_I):(F_I,u_{IJ})\to (G_I,v_{IJ})$ with 
$(F_I,u_{IJ}),(G_I,v_{IJ})\in C(\Var(\mathbb C)^{2,(sm)}/(\tilde S_I)^{op})$ 
is said to be an $(\mathbb A^1,et)$ local equivalence 
if for all $H=(H_I,w_{IJ})\in C(\Var(\mathbb C)^{2,(sm)}/(\tilde S_I))$ 
which is $\mathbb A^1$ local for the etale topology 
\begin{eqnarray*}
(\Hom(L(m_I),E_{et}(H_I))):\Hom(L(G_I,v_{IJ}),E_{et}(H_I,w_{IJ}))\to\Hom(L(F_I,u_{IJ}),E_{et}(H_I,w_{IJ})) 
\end{eqnarray*}
is a quasi-isomorphism.
Obviously, if a morphism $m=(m_I):(F_I,u_{IJ})\to (G_I,v_{IJ})$ with 
$(F_I,u_{IJ}),(G_I,u_{IJ})\in C(\Var(\mathbb C)^{2,(sm)}/(\tilde S_I)^{op})$ 
is an $(\mathbb A^1,et)$ local equivalence, 
then all the $m_I:F_I\to G_I$ are $(\mathbb A^1,et)$ local equivalence.
\item[(iv)'] A morphism $m=(m_I):(F_I,u_{IJ})\to (G_I,v_{IJ})$ with 
$(F_I,u_{IJ}),(G_I,v_{IJ})\in C(\Var(\mathbb C)^{2,(sm)pr}/(\tilde S_I)^{op})$ 
is said to be an $(\mathbb A^1,et)$ local equivalence 
if for all $(H_I,w_{IJ})\in C(\Var(\mathbb C)^{2,(sm)pr}/(\tilde S_I))$ 
which is $\mathbb A^1$ local for the etale topology
\begin{eqnarray*}
(\Hom(L(m_I),E_{et}(H_I))):\Hom(L(G_I,v_{IJ}),E_{et}(H_I,w_{IJ}))\to\Hom(L(F_I,u_{IJ}),E_{et}(H_I,w_{IJ})) 
\end{eqnarray*}
is a quasi-isomorphism.
Obviously, if a morphism $m=(m_I):(F_I,u_{IJ})\to (G_I,v_{IJ})$ with 
$(F_I,u_{IJ}),(G_I,u_{IJ})\in C(\Var(\mathbb C)^{2,(sm)pr}/(\tilde S_I)^{op})$ 
is an $(\mathbb A^1,et)$ local equivalence, 
then all the $m_I:F_I\to G_I$ are $(\mathbb A^1,et)$ local equivalence.
\end{itemize}
\end{defi}

\begin{prop}\label{ca1Var12IJ}
Let $S\in\Var(\mathbb C)$. Let $S=\cup_{i=1}^l S_i$ an open affine cover and denote by $S_I=\cap_{i\in I} S_i$.
Let $i_i:S_i\hookrightarrow\tilde S_i$ closed embeddings, with $\tilde S_i\in\Var(\mathbb C)$.
\begin{itemize}
\item[(i)]Then for $F\in C(\Var(\mathbb C)^{2,(sm)}/(\tilde S_I)^{op})$, 
$C_*F$ is $\mathbb A^1$ local for the etale topology
and  $c(F):F\to C_*F$ is an equivalence $(\mathbb A^1,et)$ local.
\item[(i)']Then for $F\in C(\Var(\mathbb C)^{2,(sm)pr}/(\tilde S_I)^{op})$, 
$C_*F$ is $\mathbb A^1$ local for the etale topology
and  $c(F):F\to C_*F$ is an equivalence $(\mathbb A^1,et)$ local.
\item[(ii)]A morphism $m:F\to G$ with $F,G\in C(\Var(\mathbb C)^{2,(sm)}/(\tilde S_I)^{op})$ 
is an $(\mathbb A^1,et)$ local equivalence
if and only if $a_{et}H^nC_*\Cone(m)=0$ for all $n\in\mathbb Z$.
\item[(ii)']A morphism $m:F\to G$ with $F,G\in C(\Var(\mathbb C)^{2,(sm)pr}/(\tilde S_I)^{op})$ 
is an $(\mathbb A^1,et)$ local equivalence
if and only if $a_{et}H^nC_*\Cone(m)=0$ for all $n\in\mathbb Z$.
\item[(iii)]A morphism $m:F\to G$ with $F,G\in C(\Var(\mathbb C)^{2,(sm)}/(\tilde S_I)^{op})$ 
is an $(\mathbb A^1,et)$ local equivalence if and only if there exists 
\begin{eqnarray*}
\left\{((X_{1,\alpha,I},Z_{1,\alpha,I})/\tilde S_I,u^1_{IJ}),\alpha\in\Lambda_1\right\},\ldots,
\left\{((X_{r,\alpha,I},Z_{r,\alpha,I})/\tilde S_I,u^r_{IJ}),\alpha\in\Lambda_r\right\}
\subset\Var(\mathbb C)^{2,(sm)}/(\tilde S_I)^{op}
\end{eqnarray*}
with 
\begin{equation*}
u^l_{IJ}:(X_{l,\alpha,J},Z_{l,\alpha,J})/\tilde S_J\to 
(X_{l,\alpha,I}\times\tilde S_{J\backslash I},Z_{l,\alpha,I}\times\tilde S_{J\backslash I})/\tilde S_J
\end{equation*}
such that we have in $\Ho_{et}(C(\Var(\mathbb C)^{2,(sm)}/(\tilde S_I)^{op}))$
\begin{eqnarray*}
\Cone(m)\xrightarrow{\sim}\Cone( \\ \oplus_{\alpha\in\Lambda_1} 
\Cone((\mathbb Z((X_{1,\alpha,I}\times\mathbb A^1,Z_{1,\alpha,I}\times\mathbb A^1)/\tilde S_I),\mathbb Z(u_{IJ}^1\times I)) 
\to(\mathbb Z((X_{1,\alpha,I},Z_{1,\alpha,I})/\tilde S_I),\mathbb Z(u_{IJ}^1))) \\
\to\cdots\to \\ \oplus_{\alpha\in\Lambda_r}
\Cone((\mathbb Z((X_{r,\alpha,I}\times\mathbb A^1,Z_{r,\alpha,I}\times\mathbb A^1)/\tilde S_I),\mathbb Z(u_{IJ}^r\times I)) 
\to(\mathbb Z((X_{r,\alpha,I},Z_{r,\alpha,I})/\tilde S_I),\mathbb Z(u^r_{IJ}))))
\end{eqnarray*}
\item[(iii)']A morphism $m:F\to G$ with $F,G\in C(\Var(\mathbb C)^{2,(sm)pr}/(\tilde S_I)^{op})$ 
is an $(\mathbb A^1,et)$ local equivalence if and only if there exists 
\begin{eqnarray*}
\left\{((Y_{1,\alpha,I}\times\tilde S_I,Z_{1,\alpha,I})/\tilde S_I,u^1_{IJ}),\alpha\in\Lambda_1\right\},\ldots,
\left\{((Y_{r,\alpha,I}\times\tilde S_I,Z_{r,\alpha,I})/\tilde S_I,u^r_{IJ}),\alpha\in\Lambda_r\right\} \\
\subset\Var(\mathbb C)^{2,(sm)pr}/(\tilde S_I)
\end{eqnarray*}
with 
\begin{equation*}
u^l_{IJ}:(Y_{l,\alpha,J}\times\tilde S_J,Z_{l,\alpha,J})/\tilde S_J\to 
(Y_{l,\alpha,I}\times\tilde S_J,Z_{l,\alpha,I}\times\tilde S_{J\backslash I})/\tilde S_J
\end{equation*}
such that we have in $\Ho_{et}(C(\Var(\mathbb C)^{2,(sm)}/(\tilde S_I)^{op}))$
\begin{eqnarray*}
\Cone(m)\xrightarrow{\sim}\Cone(\oplus_{\alpha\in\Lambda_1} \\
\Cone((\mathbb Z((Y_{1,\alpha,I}\times\mathbb A^1\times\tilde S_I,Z_{1,\alpha,I}\times\mathbb A^1)/\tilde S_I),
\mathbb Z(u_{IJ}^1\times I)) 
\to(\mathbb Z((Y_{1,\alpha,I}\times S,Z_{1,\alpha,I})/\tilde S_I),\mathbb Z(u_{IJ}))) \\
\to\cdots\to\oplus_{\alpha\in\Lambda_r} \\
\Cone((\mathbb Z((Y_{r,\alpha,I}\times\mathbb A^1\times\tilde S_I,Z_{r,\alpha,I}\times\mathbb A^1)/\tilde S_I),
\mathbb Z(u_{IJ}^r\times I)) 
\to(\mathbb Z((Y_{r,\alpha,I}\times\tilde S_I,Z_{r,\alpha})/\tilde S_I),\mathbb Z(u_{IJ}^r)))
\end{eqnarray*}
\item[(iv)] A similar statement then (iii) holds for equivalence $(\mathbb A^1,et)$ local
$m:F\to G$ with $F,G\in C(\Var(\mathbb C)^{2,(sm)}/(\tilde S_I))$
\item[(iv)'] A similar statement then (iii) holds for equivalence $(\mathbb A^1,et)$ local
$m:F\to G$ with $F,G\in C(\Var(\mathbb C)^{2,(sm)pr}/(\tilde S_I))$
\end{itemize}
\end{prop}

\begin{proof}
Similar to the proof of proposition \ref{ca1Var12}. See Ayoub's thesis for example.
\end{proof}

In the filtered case we also consider :

\begin{defi}
Let $S\in\Var(\mathbb C)$. Let $S=\cup_{i=1}^l S_i$ an open affine cover and denote by $S_I=\cap_{i\in I} S_i$.
Let $i_i:S_i\hookrightarrow\tilde S_i$ closed embeddings, with $\tilde S_i\in\SmVar(\mathbb C)$. 
\begin{itemize}
\item[(i)]A filtered complex $(G,F)\in C_{fil}(\Var(\mathbb C)^{2,(sm)}/S)$ 
is said to be $r$-filtered $\mathbb A^1$ homotopic if 
$\ad(p_a^*,p_{a*})(G,F):(G,F)\to p_{a*}p_a^*(G,F)$ is an $r$-filtered homotopy equivalence.
\item[(i)']A filtered complex $(G,F)\in C_{fil}(\Var(\mathbb C)^{2,(sm)}/(\tilde S_I))$ 
is said to be $r$-filtered $\mathbb A^1$ homotopic if 
$\ad(p_a^*,p_{a*})(G,F):(G,F)\to p_{a*}p_a^*(G,F)$ is an $r$-filtered homotopy equivalence.
\item[(ii)]A filtered complex $(G,F)\in C_{fil}(\Var(\mathbb C)^{2,(sm)pr}/S)$ 
is said to be $r$-filtered $\mathbb A^1$ homotopic if 
$\ad(p_a^*,p_{a*})(G,F):(G,F)\to p_{a*}p_a^*(G,F)$ is an $r$-filtered homotopy equivalence.
\item[(ii)']A filtered complex $(G,F)\in C_{fil}(\Var(\mathbb C)^{2,(sm)pr}/(\tilde S_I))$ 
is said to be $r$-filtered $\mathbb A^1$ homotopic if 
$\ad(p_a^*,p_{a*})(G,F):(G,F)\to p_{a*}p_a^*(G,F)$ is an $r$-filtered homotopy equivalence.
\end{itemize}
\end{defi}

We will use to compute the algebraic De Rahm realization functor the followings

\begin{thm}\label{DDADM12}
Let $S\in\Var(\mathbb C)$. 
\begin{itemize}
\item[(i)]Let $\phi:F^{\bullet}\to G^{\bullet}$ an etale local equivalence 
with $F^{\bullet},G^{\bullet}\in C(\Var(\mathbb C)^{2,sm}/S)$.
If $F^{\bullet}$ and $G^{\bullet}$ are $\mathbb A^1$ local and admit tranferts 
then $\phi:F^{\bullet}\to G^{\bullet}$ is a Zariski local equivalence.
Hence if $F\in C(\Var(\mathbb C)^{2,sm}/S)$ is $\mathbb A^1$ local and admits transfert 
\begin{equation*}
k:E_{zar}(F)\to E_{et}(E_{zar}(F))=E_{et}(F) 
\end{equation*}
is a Zariski local equivalence.
\item[(ii)]Let $\phi:F^{\bullet}\to G^{\bullet}$ an etale local equivalence
with $F^{\bullet},G^{\bullet}\in C(\Var(\mathbb C)^{2,smpr}/S)$.
If $F^{\bullet}$ and $G^{\bullet}$ are $\mathbb A^1$ local and admit tranferts 
then $\phi:F^{\bullet}\to G^{\bullet}$ is a Zariski local equivalence.
Hence if $F\in C(\Var(\mathbb C)^{2,smpr}/S)$ is $\mathbb A^1$ local and admits transfert 
\begin{equation*}
k:E_{zar}(F)\to E_{et}(E_{zar}(F))=E_{et}(F) 
\end{equation*}
is a Zariski local equivalence.
\end{itemize}
\end{thm}

\begin{proof}
Similar to the proof of theorem \ref{DDADM}.
\end{proof}

\begin{thm}\label{DDADM12fil}
Let $S\in\Var(\mathbb C)$. Let $S=\cup_{i=1}^l S_i$ an open affine cover and denote by $S_I=\cap_{i\in I} S_i$.
Let $i_i:S_i\hookrightarrow\tilde S_i$ closed embeddings, with $\tilde S_i\in\Var(\mathbb C)$. 
\begin{itemize}
\item[(i)]Let $\phi:(F^{\bullet},F)\to (G^{\bullet},F)$ a filtered etale local equivalence 
with $(F^{\bullet},F),(G^{\bullet},F)\in C_{fil}(\Var(\mathbb C)^{2,sm}/S)$.
If $(F^{\bullet},F)$ and $(G^{\bullet},F)$ are $r$-filtered $\mathbb A^1$ homotopic and admit tranferts 
then $\phi:(F^{\bullet},F)\to (G^{\bullet},F)$ is an $r$-filtered Zariski local equivalence.
Hence if $(G,F)\in C(\Var(\mathbb C)^{2,sm}/S)$ is $r$-filtered $\mathbb A^1$ homotopic and admits transfert 
\begin{equation*}
k:E_{zar}(G,F)\to E_{et}(E_{zar}(G,F))=E_{et}(G,F) 
\end{equation*}
is an $r$-filtered Zariski local equivalence.
\item[(i)']Let $\phi:(F^{\bullet},F)\to (G^{\bullet},F)$ a filtered etale local equivalence 
with $((F_I^{\bullet},F),u_{IJ}),((G_I^{\bullet},F),v_{IJ})\in C_{fil}(\Var(\mathbb C)^{2,sm}/(\tilde S_I))$.
If $((F^{\bullet},F),u_{IJ})$ and $((G^{\bullet},F),v_{IJ})$ are $r$-filtered $\mathbb A^1$ homotopic and admit tranferts 
then $\phi:((F^{\bullet},F),u_{IJ})\to ((G^{\bullet},F),v_{IJ})$ is an $r$-filtered Zariski local equivalence.
Hence if $((G_I,F),u_{IJ})\in C(\Var(\mathbb C)^{2,sm}/S)$ is $r$-filtered $\mathbb A^1$ homotopic and admits transfert 
\begin{equation*}
k:(E_{zar}(G_I,F),u_{IJ})\to (E_{et}(E_{zar}(G_I,F)),u_{IJ})=(E_{et}(G,F),u_{IJ}) 
\end{equation*}
is an $r$-filtered Zariski local equivalence.
\item[(ii)]Let $\phi:(F^{\bullet},F)\to (G^{\bullet},F)$ a filtered etale local equivalence
with $(F^{\bullet},F),(G^{\bullet},F)\in C_{fil}(\Var(\mathbb C)^{2,smpr}/S)$.
If $F^{\bullet}$ and $G^{\bullet}$ are $r$-filtered $\mathbb A^1$ homotopic and admit tranferts 
then $\phi:(F^{\bullet},F)\to (G^{\bullet},F)$ is an $r$-filtered Zariski local equivalence.
Hence if $(G,F)\in C(\Var(\mathbb C)^{2,smpr}/S)$ is $r$-filtered $\mathbb A^1$ homotopic and admits transfert 
\begin{equation*}
k:E_{zar}(F)\to E_{et}(E_{zar}(F))=E_{et}(F) 
\end{equation*}
is an $r$-filtered Zariski local equivalence.
\item[(ii)']Let $\phi:(F^{\bullet},F)\to (G^{\bullet},F)$ a filtered etale local equivalence 
with $((F_I^{\bullet},F),u_{IJ}),((G_I^{\bullet},F),v_{IJ})\in C_{fil}(\Var(\mathbb C)^{2,smpr}/(\tilde S_I))$.
If $((F^{\bullet},F),u_{IJ})$ and $((G^{\bullet},F),v_{IJ})$ are $r$-filtered $\mathbb A^1$ homotopic and admit tranferts 
then $\phi:((F^{\bullet},F),u_{IJ})\to ((G^{\bullet},F),v_{IJ})$ is an $r$-filtered Zariski local equivalence.
Hence if $((G_I,F),u_{IJ})\in C(\Var(\mathbb C)^{2,smpr}/S)$ is $r$-filtered $\mathbb A^1$ homotopic and admits transfert 
\begin{equation*}
k:(E_{zar}(G_I,F),u_{IJ})\to (E_{et}(E_{zar}(G_I,F)),u_{IJ})=(E_{et}(G,F),u_{IJ}) 
\end{equation*}
is an $r$-filtered Zariski local equivalence.
\end{itemize}
\end{thm}

\begin{proof}
Similar to the proof of theorem \ref{DDADM12}.
\end{proof}

We have the following canonical functor :

\begin{defi}\label{eta12def}
\begin{itemize}
\item[(i)]For $S\in\Var(\mathbb C)$, we have the functor 
\begin{eqnarray*}
(-)^{\Gamma}:C(\Var(\mathbb C)^{sm}/S)\to C(\Var(\mathbb C)^{2,sm}/S), \\ 
F\longmapsto F^{\Gamma}:(((U,Z)/S)=((U,Z),h)\mapsto F^{\Gamma}((U,Z)/S):=(\Gamma^{\vee}_Zh^*LF)(U/U), \\ 
(g:((U',Z'),h')\to((U,Z),h))\mapsto \\
(F^{\Gamma}(g):(\Gamma^{\vee}_Zh^*LF)(U/U)\xrightarrow{i_{(\Gamma^{\vee}_Zh^*LF)(U/U)}}(g^*(\Gamma^{\vee}_Zh^*LF))(U'/U') \\ 
\xrightarrow{T(g,\gamma^{\vee})(h^*LF)(U'/U')}(\Gamma^{\vee}_{Z\times_U U'}g^*h^*LF)(U'/U') \\ 
\xrightarrow{T(Z'/Z\times_U U',\gamma^{\vee})(g^*h^*LF)(U'/U')}(\Gamma^{\vee}_{Z'}g^*h^*LF)(U'/U')))
\end{eqnarray*}
where $i_{(\Gamma^{\vee}_Zh^*LF)(U/U)}$ is the canonical arrow of the inductive limit.
Similarly, we have, for $S\in\Var(\mathbb C)$, the functor
\begin{eqnarray*}
(-)^{\Gamma}:C(\Var(\mathbb C)/S)\to C(\Var(\mathbb C)^2/S), \\ 
F\longmapsto F^{\Gamma}:(((X,Z)/S)=((X,Z),h)\mapsto F^{\Gamma}((X,Z)/S):=(\Gamma^{\vee}_Zh^*F)(X/X), \\ 
(g:((X',Z'),h')\to((X,Z),h))\mapsto 
(F^{\Gamma}(g):(\Gamma^{\vee}_Zh^*LF)(X/X)\to (\Gamma^{\vee}_{Z'}h^{'*}LF)(X'/X')))
\end{eqnarray*}
Note that for $S\in\Var(\mathbb C)$, $I(S/S):\mathbb Z((S,S)/S)\to\mathbb Z(S/S)^{\Gamma}$ given by
\begin{eqnarray*}
I(S/S)((U,Z),h):\mathbb Z((S,S)/S)(((U,Z),h))\xrightarrow{\gamma^{\vee}_Z(\mathbb Z(U/U))(U/U)} 
\mathbb Z(S/S)^{\Gamma}((U,Z),h):=(\Gamma_Z^{\vee}\mathbb Z(U/U))(U/U)
\end{eqnarray*} 
is an isomorphism.
\item[(ii)]Let $f:T\to S$ a morphism with $T,S\in\Var(\mathbb C)$.
For $F\in C(\Var(\mathbb C)^{sm}/S)$, we have the canonical morphism in $C(\Var(\mathbb C)^{2,sm}/T)$
\begin{eqnarray*}
T(f,\Gamma)(F):=T^*(f,\Gamma)(F):f^*(F^{\Gamma})\to(f^*F)^{\Gamma},\\
T(f,\Gamma)(F)((U',Z')/T=((U',Z'),h')): \\
f^*(F^{\Gamma})((U',Z'),h'):=
\lim_{((U',Z'),h')\xrightarrow{l}((U_T,Z_T),h_T)\xrightarrow{f_U}((U,Z),h)}(\Gamma^{\vee}_Zh^*LF)(U/U) \\
\xrightarrow{F^{\Gamma}(f_U\circ l)}(\Gamma^{\vee}_{Z'}l^*f_U^*h^*LF)(U'/U')=(\Gamma^{\vee}_{Z'}h^{'*}f^*LF)(U'/U') \\
\xrightarrow{(\Gamma^{\vee}_{Z'}h^{'*}T(f,L)(F))(U'/U')}(\Gamma^{\vee}_{Z'}h^{'*}Lf^*F)(U'/U')=:(f^*F)^{\Gamma}((U',Z'),h') 
\end{eqnarray*}
where $f_U:U_T:U\times_S T\to U$ and $h_T:U_T:=U\times_S T\to T$ are the base change maps,
the equality following from the fact that $h\circ f_U\circ l=f\circ h_T\circ l=f\circ h'$.
For $F\in C(\Var(\mathbb C)/S)$, we have similarly the canonical morphism in $C(\Var(\mathbb C)^2/T)$
\begin{equation*}
T(f,\Gamma)(F):f^*(F^{\Gamma})\to (f^*F)^{\Gamma}.
\end{equation*}
\item[(iii)]Let $h:U\to S$ a smooth morphism with $U,S\in\Var(\mathbb C)$. 
We have, for $F\in C(\Var(\mathbb C)^{sm}/U)$, the canonical morphism in $C(\Var(\mathbb C)^{2,sm}/S)$
\begin{eqnarray*}
T_{\sharp}(h,\Gamma)(F):h_{\sharp}(F^{\Gamma})\to(h_{\sharp}LF)^{\Gamma}, \\
T_{\sharp}(h,\Gamma)(F)((U',Z'),h'):h_{\sharp}(F^{\Gamma})((U',Z'),h'):=
\lim_{((U',Z'),h')\xrightarrow{l}((U,U),h)}(\Gamma^{\vee}_{Z'}l^*LF)(U'/U') \\
\xrightarrow{(\Gamma^{\vee}_{Z'}l^*\ad(h_{\sharp},h^*)(LF))(U'/U')}
(\Gamma^{\vee}_{Z'}l^*h^*h_{\sharp}LF)(U'/U')=:(h_{\sharp}LF)^{\Gamma}((U',Z')/h')
\end{eqnarray*}
\item[(iv)]Let $i:Z_0\hookrightarrow S$ a closed embedding with $Z_0,S\in\Var(\mathbb C)$. 
We have the canonical morphism in $C(\Var(\mathbb C)^{2,sm}/S)$
\begin{eqnarray*}
T_*(i,\Gamma)(\mathbb Z(Z_0/Z_0)):i_*((\mathbb Z(Z_0/Z_0))^{\Gamma}\to(i_*\mathbb Z(Z/Z))^{\Gamma}, \\
T_*(i,\Gamma)(\mathbb Z(Z_0/Z_0))((U,Z),h):
i_*((\mathbb Z(Z_0/Z_0))^{\Gamma}((U,Z),h):=(\Gamma^{\vee}_{Z\times_SZ_0}\mathbb Z(Z_0/Z_0))(U\times_S Z_0) \\
\xrightarrow{T(i_*,\gamma^{\vee})(\mathbb Z(Z_0/Z_0))(U\times_S Z_0)}
(\Gamma^{\vee}_Zi_*\mathbb Z(Z_0/Z_0))(U\times_S Z_0)=:(i_*\mathbb Z(Z/Z))^{\Gamma}((U,Z),h)
\end{eqnarray*}
\end{itemize}
\end{defi}

\begin{defi}\label{GrGamma}
Let $S\in\Var(\mathbb C)$. We have for $F\in C(\Var(\mathbb C)^{sm}/S)$ the canonical map in $C(\Var(\mathbb C)^{sm}/S)$
\begin{eqnarray*}
\Gr(F):\Gr^{12}_{S*}\mu_{S*}F^{\Gamma}\to F, \\
\Gr(F)(U/S):\Gamma_U^{\vee}p^*F(U\times S/U\times S)\xrightarrow{\ad(l^*,l_*)(p^*F)(U\times S/U\times S)}h^*F(U/U)=F(U/S)
\end{eqnarray*}
where $h:U\to S$ is a smooth morphism with $U\in\Var(\mathbb C)$ and $h:U\xrightarrow{l}U\times S\xrightarrow{p}S$
is the graph factorization with $l$ the graph embedding and $p$ the projection.
\end{defi}

\begin{prop}\label{eta12prop}
Let $S\in\Var(\mathbb C)$. 
\begin{itemize}
\item[(i)] Then,
\begin{itemize}
\item if $m:F\to G$ with $F,G\in C(\Var(\mathbb C)^{sm}/S)$ is a quasi-isomorphism,
$m^{\Gamma}:F^{\Gamma}\to G^{\Gamma}$ is a quasi-isomorphism in $C(\Var(\mathbb C)^{2,sm}/S)$,
\item if $m:F\to G$ with $F,G\in C(\Var(\mathbb C)^{sm}/S)$ is a Zariski local equivalence,
$m^{\Gamma}:F^{\Gamma}\to G^{\Gamma}$ is a Zariski local equivalence in $C(\Var(\mathbb C)^{2,sm}/S)$,
if $m:F\to G$ with $F,G\in C(\Var(\mathbb C)^{sm}/S)$ is an etale local equivalence,
$m^{\Gamma}:F^{\Gamma}\to G^{\Gamma}$ is an etale local equivalence in $C(\Var(\mathbb C)^{2,sm}/S)$,
\item if $m:F\to G$ with $F,G\in C(\Var(\mathbb C)^{sm}/S)$ is an $(\mathbb A^1,et)$ local equivalence,
$m^{\Gamma}:F^{\Gamma}\to G^{\Gamma}$ is an $(\mathbb A^1,et)$ local equivalence in $C(\Var(\mathbb C)^{2,sm}/S)$.
\end{itemize}
\item[(ii)] Then,
\begin{itemize}
\item if $m:F\to G$ with $F,G\in C(\Var(\mathbb C)/S)$ is a quasi-isomorphism,
$m^{\Gamma}:F^{\Gamma}\to G^{\Gamma}$ is a quasi-isomorphism in $C(\Var(\mathbb C)^2/S)$,
\item if $m:F\to G$ with $F,G\in C(\Var(\mathbb C)^{sm}/S)$ is a Zariski local equivalence,
$m^{\Gamma}:F^{\Gamma}\to G^{\Gamma}$ is a Zariski local equivalence in $C(\Var(\mathbb C)^{2,sm}/S)$,
if $m:F\to G$ with $F,G\in C(\Var(\mathbb C)/S)$ is an etale local equivalence,
$m^{\Gamma}:F^{\Gamma}\to G^{\Gamma}$ is an etale local equivalence in $C(\Var(\mathbb C)^2/S)$,
\item if $m:F\to G$ with $F,G\in C(\Var(\mathbb C)^{sm}/S)$ is an $(\mathbb A^1,et)$ local equivalence,
$m^{\Gamma}:F^{\Gamma}\to G^{\Gamma}$ is an $(\mathbb A^1,et)$ local equivalence in $C(\Var(\mathbb C)^2/S)$.
\end{itemize}
\end{itemize}
\end{prop}

\begin{proof}
\noindent(i): Follows immediately from the fact that for $((U,Z),h)\in\Var(\mathbb C)^{2,sm}/S$,
\begin{itemize}
\item if $m:F\to G$ with $F,G\in C(\Var(\mathbb C)^{sm}/S)$ is a quasi-isomorphism,
$\Gamma_Z^{\vee}h^*LF(m):\Gamma_Z^{\vee}h^*LF\to\Gamma_Z^{\vee}h^*LG$ is a quasi-isomorphism
\item if $m:F\to G$ with $F,G\in C(\Var(\mathbb C)^{sm}/S)$ is a is a Zariski or etale local equivalence,
$\Gamma_Z^{\vee}h^*LF(m):\Gamma_Z^{\vee}h^*LF\to\Gamma_Z^{\vee}h^*LG$ is a Zariski, resp. etale, local equivalence,
\item if $m:F\to G$ with $F,G\in C(\Var(\mathbb C)^{sm}/S)$ is an $(\mathbb A^1,et)$ local equivalence,
$\Gamma_Z^{\vee}h^*LF(m):\Gamma_Z^{\vee}h^*LF\to\Gamma_Z^{\vee}h^*LG$ is an $(\mathbb A^1,et)$ local equivalence.
\end{itemize}

\noindent(ii): Similar to (i).
\end{proof}

\subsection{Presheaves on the big analytical site}

For $S\in\AnSp(\mathbb C)$, we denote by $\rho_S:\AnSp(\mathbb C)^{sm}/S\hookrightarrow\AnSp(\mathbb C)/S$ be the full subcategory 
consisting of the objects $U/S=(U,h)\in\AnSp(\mathbb C)/S$ such that the morphism $h:U\to S$ is smooth. 
That is, $\AnSp(\mathbb C)^{sm}/S$ is the category  
\begin{itemize}
\item whose objects are smooth morphisms 
$U/S=(U,h)$, $h:U\to S$ with $U\in\AnSp(\mathbb C)$, 
\item whose morphisms $g:U/S=(U,h_1)\to V/S=(V,h_2)$ 
is a morphism $g:U\to V$ of complex algebraic varieties such that $h_2\circ g=h_1$. 
\end{itemize}
We denote again $\rho_S:\AnSp(\mathbb C)/S\to\AnSp(\mathbb C)^{sm}/S$ the associated morphism of site.
We will consider 
\begin{equation*}
r^s(S):\AnSp(\mathbb C)\xrightarrow{r(S)}\AnSp(\mathbb C)/S\xrightarrow{\rho_S}\AnSp(\mathbb C)^{sm}/S
\end{equation*}
the composite morphism of site.
For $S\in\AnSp(\mathbb C)$, we denote by $\mathbb Z_S:=\mathbb Z(S/S)\in \PSh(\AnSp(\mathbb C)^{sm}/S)$ the constant presheaf
By Yoneda lemma, we have for $F\in C(\AnSp(\mathbb C)^{sm}/S)$, $\mathcal Hom(\mathbb Z_S,F)=F$.
For $f:T\to S$ a morphism, with $T,S\in\AnSp(\mathbb C)$, we have the following commutative diagram of sites
\begin{equation}\label{pf0an}
\xymatrix{\AnSp(\mathbb C)/T\ar[d]^{P(f)}\ar[r]^{\rho_T} & \AnSp(\mathbb C)^{sm}/T\ar[d]^{P(f)} \\ 
\AnSp(\mathbb C)/S\ar[r]^{\rho_S} & \AnSp(\mathbb C)^{sm}/S} 
\end{equation}
We denote, for $S\in\AnSp(\mathbb C)$, the obvious morphism of sites 
\begin{equation*}
\tilde e(S):\AnSp(\mathbb C)/S\xrightarrow{\rho_S}\AnSp(\mathbb C)^{sm}/S\xrightarrow{e(S)}\Ouv(S)  
\end{equation*}
where $\Ouv(S)$ is the set of the open subsets of $S$, given by the inclusion functors
$\tilde e(S):\Ouv(S)\hookrightarrow\AnSp(\mathbb C)^{sm}/S\hookrightarrow\AnSp(\mathbb C)/S$.
By definition, for $f:T\to S$ a morphism with $S,T\in\AnSp(\mathbb C)$, the commutative diagram of sites (\ref{pf0an})
extend a commutative diagram of sites :
\begin{equation}\label{empf0an}
\xymatrix{
\tilde e(T):\AnSp(\mathbb C)/T\ar[d]^{P(f)}\ar[rr]^{\rho_T} & \, & \AnSp(\mathbb C)^{sm}/T\ar[d]^{P(f)}\ar[rr]^{e(T)} 
& \, & \Ouv(T)\ar[d]^{P(f)} \\
\tilde e(S):\AnSp(\mathbb C)/S\ar[rr]^{\rho_S} & \, & \AnSp(\mathbb C)^{sm}/S\ar[rr]^{e(S)} & \, & \Ouv(S)}
\end{equation}

\begin{itemize}
\item As usual, we denote by
\begin{equation*}
(f^*,f_*):=(P(f)^*,P(f)_*):C(\AnSp(\mathbb C)^{sm}/S)\to C(\AnSp(\mathbb C)^{sm}/T)
\end{equation*}
the adjonction induced by $P(f):\AnSp(\mathbb C)^{sm}/T\to \AnSp(\mathbb C)^{sm}/S$.
Since the colimits involved in the definition of $f^*=P(f)^*$ are filtered, $f^*$ also preserve monomorphism. 
Hence, we get an adjonction
\begin{equation*}
(f^*,f_*):C_{fil}(\AnSp(\mathbb C)^{sm}/S)\leftrightarrows C_{fil}(\AnSp(\mathbb C)^{sm}/T), \; f^*(G,F):=(f^*G,f^*F)
\end{equation*}
\item As usual, we denote by
\begin{equation*}
(f^*,f_*):=(P(f)^*,P(f)_*):C(\AnSp(\mathbb C)/S)\to C(\AnSp(\mathbb C)/T)
\end{equation*}
the adjonction induced by $P(f):\AnSp(\mathbb C)/T\to\AnSp(\mathbb C)/S$.
Since the colimits involved in the definition of $f^*=P(f)^*$ are filtered, $f^*$ also preserve monomorphism. 
Hence, we get an adjonction
\begin{equation*}
(f^*,f_*):C_{fil}(\AnSp(\mathbb C)/S)\leftrightarrows C_{fil}(\AnSp(\mathbb C)/T), \; f^*(G,F):=(f^*G,f^*F)
\end{equation*}
\end{itemize}

\begin{itemize}
\item For $h:U\to S$ a smooth morphism with $U,S\in\AnSp(\mathbb C)$, 
the pullback functor $P(h):\AnSp(\mathbb C)^{sm}/S\to \AnSp(\mathbb C)^{sm}/U$ 
admits a left adjoint $C(h)(X\to U)=(X\to U\to S)$.
Hence, $h^*:C(\AnSp(\mathbb C)^{sm}/S)\to C(\AnSp(\mathbb C)^{sm}/U)$ admits a left adjoint
\begin{equation*}
h_{\sharp}:C(\AnSp(\mathbb C)^{sm}/U)\to C(\AnSp(\mathbb C)^{sm}/S), \; 
F\mapsto((V,h_0)\mapsto\lim_{(V',h\circ h')\to(V,h_0)}F(V',h'))
\end{equation*}
Note that for $h':V'\to V$ a smooth morphism, $V',V\in\AnSp(\mathbb C)$, 
we have $h_{\sharp}(\mathbb Z(V'/V))=\mathbb Z(V'/S)$ with $V'/S=(V',h\circ h')$.
Hence, since projective presheaves are the direct summands of the representable presheaves,
$h_{\sharp}$ sends projective presheaves to projective presheaves.
For $F^{\bullet}\in C(\AnSp(\mathbb C)^{sm}/S)$ and $G^{\bullet}\in C(\AnSp(\mathbb C)^{sm}/U)$,
we have the adjonction maps
\begin{equation*}
\ad(h_{\sharp},h^*)(G^{\bullet}):G^{\bullet}\to h^*h_{\sharp}G^{\bullet} \; , \;
\ad(h_{\sharp},h^*)(F^{\bullet}):h_{\sharp}h^*F^{\bullet}\to F^{\bullet}.
\end{equation*}
For a smooth morphism $h:U\to S$, with $U,S\in\AnSp(\mathbb C)$, we have the adjonction isomorphism, 
for $F\in C(\AnSp(\mathbb C)^{sm}/U)$ and $G\in C(\AnSp(\mathbb C)^{sm}/S)$,   
\begin{equation}\label{Ihhoman}
I(h_{\sharp},h^*)(F,G):\mathcal Hom^{\bullet}(h_{\sharp}F,G)\xrightarrow{\sim}h_*\mathcal Hom^{\bullet}(F,h^*G).  
\end{equation}

\item For $f:T\to S$ any morphism with $T,S\in\AnSp(\mathbb C)$, 
the pullback functor $P(f):\AnSp(\mathbb C)/T\to\AnSp(\mathbb C)/S$ 
admits a left adjoint $C(f)(X\to T)=(X\to T\to S)$.
Hence, $f^*:C(\AnSp(\mathbb C)/S)\to C(\AnSp(\mathbb C)/T)$ admits a left adjoint
\begin{equation*}
f_{\sharp}:C(\AnSp(\mathbb C)/T)\to C(\AnSp(\mathbb C)/S), \; 
F\mapsto((V,h_0)\mapsto\lim_{(V',h\circ h')\to(V,h_0)}F(V',h'))
\end{equation*}
Note that we have for $h':V'\to V$ a morphism, $V',V\in\AnSp(\mathbb C)$, 
$h_{\sharp}(\mathbb Z(V'/V))=\mathbb Z(V'/S)$ with $V'/S=(V',h\circ h')$.
Hence, since projective presheaves are the direct summands of the representable presheaves,
$h_{\sharp}$ sends projective presheaves to projective presheaves.
For $F^{\bullet}\in C(\AnSp(\mathbb C)/S)$ and $G^{\bullet}\in C(\AnSp(\mathbb C)/T)$,
we have the adjonction maps
\begin{equation*}
\ad(f_{\sharp},f^*)(G^{\bullet}):G^{\bullet}\to f^*f_{\sharp}G^{\bullet} \; , \;
\ad(f_{\sharp},f^*)(F^{\bullet}):f_{\sharp}f^*F^{\bullet}\to F^{\bullet}.
\end{equation*}
For a morphism $f:T\to S$, with $T,S\in\AnSp(\mathbb C)$, we have the adjonction isomorphism, 
for $F\in C(\AnSp(\mathbb C)/T)$ and $G\in C(\AnSp(\mathbb C)/S)$,   
\begin{equation}\label{Ifsharphoman}
I(f_{\sharp},f^*)(F,G):\mathcal Hom^{\bullet}(f_{\sharp}F,G)\xrightarrow{\sim}f_*\mathcal Hom^{\bullet}(F,f^*G).  
\end{equation}
\end{itemize}

\begin{itemize}
\item For a commutative diagram in $\AnSp(\mathbb C)$ : 
\begin{equation*}
D=\xymatrix{ 
V\ar[r]^{g_2}\ar[d]^{h_2} & U\ar[d]^{h_1} \\
T\ar[r]^{g_1} & S},
\end{equation*}
where $h_1$ and $h_2$ are smooth,
we denote by, for $F^{\bullet}\in C(\AnSp(\mathbb C)^{sm}/U)$, 
\begin{equation*}
T_{\sharp}(D)(F^{\bullet}): h_{2\sharp}g_2^*F^{\bullet}\to g_1^*h_{1\sharp}F^{\bullet}
\end{equation*}
the canonical map in $C(\AnSp(\mathbb C)^{sm}/T)$ given by adjonction. 
If $D$ is cartesian with $h_1=h$, $g_1=g$ $f_2=h':U_T\to T$, $g':U_T\to U$,
\begin{equation*}
T_{\sharp}(D)(F^{\bullet})=:T_{\sharp}(g,h)(F^{\bullet}):
h'_{\sharp}g^{'*}F^{\bullet}\xrightarrow{\sim} g^*h_{\sharp}F^{\bullet}
\end{equation*}
is an isomorphism and for $G^{\bullet}\in C(\AnSp(\mathbb C)^{sm}/T)$
\begin{equation*}
T(D)(G^{\bullet})=:T(g,h)(G^{\bullet}):g^*h_*G^{\bullet}\xrightarrow{\sim} h'_*g^{'*}G^{\bullet}
\end{equation*}
is an isomorphism.

\item For a commutative diagram in $\AnSp(\mathbb C)$ : 
\begin{equation*}
D=\xymatrix{ 
V\ar[r]^{g_2}\ar[d]^{f_2} & X\ar[d]^{f_1} \\
T\ar[r]^{g_1} & S},
\end{equation*}
we denote by, for $F^{\bullet}\in C(\AnSp(\mathbb C)/X)$, 
\begin{equation*}
T_{\sharp}(D)(F^{\bullet}): f_{2\sharp}g_2^*F^{\bullet}\to g_1^*f_{1\sharp}F^{\bullet}
\end{equation*}
the canonical map in $C(\AnSp(\mathbb C)/T)$ given by adjonction. 
If $D$ is cartesian with $h_1=h$, $g_1=g$ $f_2=h':X_T\to T$, $g':X_T\to X$,
\begin{equation*}
T_{\sharp}(D)(F^{\bullet})=:T_{\sharp}(g,f)(F^{\bullet}):
f'_{\sharp}g^{'*}F^{\bullet}\xrightarrow{\sim} g^*f_{\sharp}F^{\bullet}
\end{equation*}
is an isomorphism and for $G^{\bullet}\in C(\AnSp(\mathbb C)/T)$
\begin{equation*}
T(D)(G^{\bullet})=:T(g,h)(G^{\bullet}):f^*g_*G^{\bullet}\xrightarrow{\sim} g'_*f^{'*}G^{\bullet}
\end{equation*}
is an isomorphism.
\end{itemize}

For $f:T\to S$ a morphism with $S,T\in\AnSp(\mathbb C)$, 
\begin{itemize}
\item we get for $F\in C(\AnSp(\mathbb C)^{sm}/S)$ from the a commutative diagram of sites (\ref{empf0an}) 
the following canonical transformation 
\begin{equation*}
T(e,f)(F^{\bullet}):f^*e(S)_*F^{\bullet}\to e(T)_*f^*F^{\bullet},
\end{equation*}
 which is NOT a quasi-isomorphism in general. 
However, for $h:U\to S$ a smooth morphism with $S,U\in\AnSp(\mathbb C)$, 
$T(e,h)(F^{\bullet}):h^*e(S)_*F^{\bullet}\xrightarrow{\sim} e(T)_*h^*F^{\bullet}$ is an isomorphism.
\item we get for $F\in C(\AnSp(\mathbb C)/S)$ from the a commutative diagram of sites (\ref{empf0an}) 
the following canonical transformation 
\begin{equation*}
T(e,f)(F^{\bullet}):f^*e(S)_*F^{\bullet}\to e(T)_*f^*F^{\bullet},
\end{equation*}
 which is NOT a quasi-isomorphism in general. 
However, for $h:U\to S$ a smooth morphism with $S,U\in\AnSp(\mathbb C)$, 
$T(e,h)(F^{\bullet}):h^*e(S)_*F^{\bullet}\xrightarrow{\sim} e(T)_*h^*F^{\bullet}$ is an isomorphism.
\end{itemize}

Let $S\in\AnSp(\mathbb C)$, 
\begin{itemize}
\item We have for $F,G\in C(\AnSp(\mathbb C)^{sm}/S)$, 
\begin{itemize}
\item $e(S)_*(F\otimes G)=(e(S)_*F)\otimes (e(S)_*G)$ by definition 
\item the canonical forgetfull map 
\begin{equation*}
T(S,hom)(F,G):e(S)_*\mathcal Hom^{\bullet}(F,G)\to\mathcal Hom^{\bullet}(e(S)_*F,e(S)_*G).
\end{equation*}
which is NOT a quasi-isomorphism in general.
\end{itemize}
By definition, we have for $F\in C(\AnSp(\mathbb C)^{sm}/S)$, $e(S)_*E_{usu}(F)=E_{usu}(e(S)_*F)$.
\item We have for $F,G\in C(\AnSp(\mathbb C)/S)$, 
\begin{itemize}
\item $e(S)_*(F\otimes G)=(e(S)_*F)\otimes (e(S)_*G)$ by definition 
\item the canonical forgetfull map 
\begin{equation*}
T(S,hom)(F,G):e(S)_*\mathcal Hom^{\bullet}(F,G)\to\mathcal Hom^{\bullet}(e(S)_*F,e(S)_*G).
\end{equation*}
which is NOT a quasi-isomorphism in general.
\end{itemize}
By definition, we have for $F\in C(\AnSp(\mathbb C)/S)$, $e(S)_*E_{usu}(F)=E_{usu}(e(S)_*F)$.
\end{itemize}

Let $S\in\AnSp(\mathbb C)$. 
We have the support section functor of a closed subset $Z\subset S$ for presheaves on the big analytical site.
\begin{defi}\label{gammaan}
Let $S\in\AnSp(\mathbb C)$. 
Let $Z\subset S$ a closed subset. Denote by $j:S\backslash Z\hookrightarrow S$ be the open complementary subset.
\begin{itemize}
\item[(i)] We define the functor
\begin{equation*}
\Gamma_Z:C(\AnSp(\mathbb C)^{sm}/S)\to C(\AnSp(\mathbb C)^{sm}/S), \;
G^{\bullet}\mapsto\Gamma_Z G^{\bullet}:=\Cone(\ad(j^*,j_*)(G^{\bullet}):G^{\bullet}\to j_*j^*G^{\bullet})[-1],
\end{equation*}
so that there is then a canonical map $\gamma_Z(G^{\bullet}):\Gamma_ZG^{\bullet}\to G^{\bullet}$.
\item[(ii)] We have the dual functor of (i) :
\begin{equation*}
\Gamma^{\vee}_Z:C(\AnSp(\mathbb C)^{sm}/S)\to C(\AnSp(\mathbb C)^{sm}/S), \; 
F\mapsto\Gamma^{\vee}_Z(F^{\bullet}):=\Cone(\ad(j_{\sharp},j^*)(G^{\bullet}):j_{\sharp}j^*G^{\bullet}\to G^{\bullet}), 
\end{equation*}
together with the canonical map $\gamma^{\vee}_Z(G):F\to\Gamma^{\vee}_Z(G)$.
\item[(iii)] For $F,G\in C(\AnSp(\mathbb C)^{sm}/S)$, we denote by 
\begin{equation*}
I(\gamma,hom)(F,G):=(I,I(j_{\sharp},j^*)(F,G)):\Gamma_Z\mathcal Hom(F,G)\xrightarrow{\sim}\mathcal Hom(\Gamma^{\vee}_ZF,G)
\end{equation*}
the canonical isomorphism given by adjonction.
\end{itemize}
\end{defi}

Let $S\in\AnSp(\mathbb C)$ and $Z\subset S$ a closed subset.
\begin{itemize}
\item Since $\Gamma_Z:C(\AnSp(\mathbb C)^{sm}/S)\to C(\AnSp(\mathbb C)^{sm}/S)$ preserve monomorphism, it induces a functor  
\begin{equation*}
\Gamma_Z:C_{fil}(\AnSp(\mathbb C)^{sm}/S)\to C_{fil}(\AnSp(\mathbb C)^{sm}/S), \; 
(G,F)\mapsto\Gamma_Z(G,F):=(\Gamma_ZG,\Gamma_ZF)
\end{equation*}
\item Since $\Gamma^{\vee}_Z:C(\AnSp(\mathbb C)^{sm}/S)\to C(\AnSp(\mathbb C)^{sm}/S)$ preserve monomorphism, it induces a functor  
\begin{equation*}
\Gamma^{\vee}_Z:C_{fil}(\AnSp(\mathbb C)^{sm}/S)\to C_{fil}(\AnSp(\mathbb C)^{sm}/S), \; 
(G,F)\mapsto\Gamma^{\vee}_Z(G,F):=(\Gamma^{\vee}_ZG,\Gamma^{\vee}_ZF)
\end{equation*}
\end{itemize}

\begin{defiprop}\label{gamma3sect2}
\begin{itemize}
\item[(i)] Let $g:S'\to S$ a morphism and $i:Z\hookrightarrow S$ a closed embedding with $S',S,Z\in\AnSp(\mathbb C)$. 
Then, for $(G,F)\in C_{fil}(\AnSp(\mathbb C)^{sm}/S)$, there exist a map in $C_{fil}(\AnSp(\mathbb C)^{sm}/S')$
\begin{equation*}
T(g,\gamma)(G,F):g^*\Gamma_{Z}(G,F)\to\Gamma_{Z\times_S S'}g^*(G,F)
\end{equation*}
unique up to homotopy, such that $\gamma_{Z\times_S S'}(g^*(G,F))\circ T(g,\gamma)(G,F)=g^*\gamma_{Z}(G,F)$.
\item[(ii)] Let $i_1:Z_1\hookrightarrow S$, $i_2:Z_2\hookrightarrow Z_1$ be closed embeddings with $S,Z_1,Z_2\in\AnSp(\mathbb C)$.
Then, for $(G,F)\in C_{fil}(\AnSp(\mathbb C)^{sm}/S)$, 
\begin{itemize}
\item there exist a canonical map $T(Z_2/Z_1,\gamma)(G,F):\Gamma_{Z_2}(G,F)\to\Gamma_{Z_1}(G,F)$ 
in $C_{fil}(\AnSp(\mathbb C)^{sm}/S)$ unique up to homotopy such that 
$\gamma_{Z_1}(G,F)\circ T(Z_2/Z_1,\gamma)(G,F)=\gamma_{Z_2}(G,F)$, together with a distinguish triangle 
\begin{equation*}
\Gamma_{Z_2}(G,F)\xrightarrow{T(Z_2/Z_1,\gamma)(G,F)}\Gamma_{Z_1}(G,F)\xrightarrow{\ad(j_2^*,j_{2*})(\Gamma_{Z_1}(G,F))}
\Gamma_{Z_1\backslash Z_2}(G,F)\to\Gamma_{Z_2}(G,F)[1]
\end{equation*}
in $K_{fil}(\AnSp(\mathbb C)^{sm}/S)$,
\item there exist a map $T(Z_2/Z_1,\gamma^{\vee})(G,F):\Gamma_{Z_1}^{\vee}(G,F)\to\Gamma_{Z_2}^{\vee}(G,F)$ 
in $C_{fil}(\AnSp(\mathbb C)^{sm}/S)$ unique up to homotopy such that 
$\gamma^{\vee}_{Z_2}(G,F)=T(Z_2/Z_1,\gamma^{\vee})(G,F)\circ\gamma^{\vee}_{Z_1}(G,F)$, 
together with a distinguish triangle 
\begin{equation*}
\Gamma^{\vee}_{Z_1\backslash Z_2}(G,F)\xrightarrow{\ad(j_{2\sharp},j_2^*)(\Gamma^{\vee}_{Z_1}G)}\Gamma^{\vee}_{Z_1}(G,F)
\xrightarrow{T(Z_2/Z_1,\gamma^{\vee})(G,F)}\Gamma^{\vee}_{Z_2}(G,F)\to\Gamma^{\vee}_{Z_1\backslash Z_2}(G,F)[1]
\end{equation*}
in $K_{fil}(\AnSp(\mathbb C)^{sm}/S)$. 
\end{itemize}
\item[(iii)] Consider a morphism $g:(S',Z)\to (S,Z)$ with $(S',Z)\to (S,Z)\in\AnSp(\mathbb C)$.
We denote, for $G\in C(\AnSp(\mathbb C)^{sm}/S)$ the composite 
\begin{equation*}
T(D,\gamma^{\vee})(G):g^*\Gamma_Z^{\vee}G\xrightarrow{\sim}\Gamma^{\vee}_{Z\times_S S'}g^*G
\xrightarrow{T(Z'/Z\times_S S',\gamma^{\vee})(G)}\Gamma_{Z'}^{\vee}g^*G
\end{equation*}
and we have then the factorization
$\gamma_{Z'}^{\vee}(g^*G):g^*G\xrightarrow{g^*\gamma_Z^{\vee}(G)}g^*\Gamma_Z^{\vee}G
\xrightarrow{T(D,\gamma^{\vee})(G)}\Gamma_{Z'}^{\vee}g^*G$.
\end{itemize}
\end{defiprop}

\begin{proof}
Similar to definition-proposition \ref{gamma1sect2} or definition-proposition \ref{gamma2sect2}.
\end{proof}

\begin{defi}
For $S\in\AnSp(\mathbb C)$, we denote by 
\begin{equation*}
C_{O_S}(\AnSp(\mathbb C)^{sm}/S):=C_{e(S)^*O_S}(\AnSp(\mathbb C)^{sm}/S) 
\end{equation*}
the category of complexes of presheaves on $\AnSp(\mathbb C)^{sm}/S$ endowed with a structure of $e(S)^*O_S$ module, and by 
\begin{equation*}
C_{O_Sfil}(\AnSp(\mathbb C)^{sm}/S):=C_{e(S)^*O_Sfil}(\AnSp(\mathbb C)^{sm}/S) 
\end{equation*}
the category of filtered complexes of presheaves on $\Var(\mathbb C)^{sm}/S$endowed with a structure of $e(S)^*O_S$ module.
\end{defi}

Let $S\in\AnSp(\mathbb C)$. Let $Z\subset S$ a closed subset.
Denote by $j:S\backslash Z\hookrightarrow S$ the open complementary embedding, 
\begin{itemize}
\item For $G\in C_{O_S}(\AnSp(\mathbb C)^{sm}/S)$, $\Gamma_ZG:=\Cone(\ad(j^*,j_*)(G):F\to j_*j^*G)[-1]$ 
has a (unique) structure of $e(S)^*O_S$ module such that $\gamma_Z(G):\Gamma_ZG\to G$ is a map in $C_{O_S}(\AnSp(\mathbb C)^{sm}/S)$. 
This gives the functor
\begin{equation*}
\Gamma_Z:C_{O_Sfil}(\AnSp(\mathbb C)^{sm}/S)\to C_{filO_S}(\AnSp(\mathbb C)^{sm}/S), \; 
(G,F)\mapsto\Gamma_Z(G,F):=(\Gamma_ZG,\Gamma_ZF), 
\end{equation*}
together with the canonical map $\gamma_Z((G,F):\Gamma_Z(G,F)\to (G,F)$.
Let $Z_2\subset Z$ a closed subset. Then, for $G\in C_{O_S}(\AnSp(\mathbb C)^{sm}/S)$, 
$T(Z_2/Z,\gamma)(G):\Gamma_{Z_2}G\to\Gamma_ZG$ is a map in $C_{O_S}(\AnSp(\mathbb C)^{sm}/S)$ (i.e. is $e(S)^*O_S$ linear).

\item For $G\in C_{O_S}(\Var(\mathbb C)^{sm}/S)$, $\Gamma^{\vee}_ZG:=\Cone(\ad(j_{\sharp},j^*)(G):j_{\sharp}j^*G\to G)$ 
has a unique structure of $e(S)^*O_S$ module, such that 
$\gamma^{\vee}_Z(G):G\to\Gamma_Z^{\vee}G$ is a map in $C_{O_S}(\AnSp(\mathbb C)^{sm}/S)$. 
This gives the the functor
\begin{equation*}
\Gamma^{\vee}_Z:C_{O_Sfil}(S)\to C_{filO_S}(S), \; (G,F)\mapsto\Gamma^{\vee}_Z(G,F):=(\Gamma^{\vee}_ZG,\Gamma^{\vee}_ZF), 
\end{equation*}
together with the canonical map $\gamma^{\vee}_Z((G,F):(G,F)\to\Gamma^{\vee}_Z(G,F)$.
Let $Z_2\subset Z$ a closed subset. Then, for $G\in C_{O_S}(\AnSp(\mathbb C)^{sm}/S)$, 
$T(Z_2/Z,\gamma^{\vee})(G):\Gamma_Z^{\vee}G\to\Gamma_{Z_2}^{\vee}G$ is a map in $C_{O_S}(\AnSp(\mathbb C)^{sm}/S)$ 
(i.e. is $e(S)^*O_S$ linear).
\end{itemize}

\begin{defi}\label{ZSan}
Let $S\in\AnSp(\mathbb C)$. Let $Z\subset S$ a closed subset.
\begin{itemize}
\item[(i)] We denote by 
\begin{equation*}
C_{Z}(\AnSp(\mathbb C)^{sm}/S)\subset C(\AnSp(\mathbb C)^{sm}/S), 
\end{equation*}
the full subcategory consisting of complexes of presheaves 
$F^{\bullet}\in C(\AnSp(\mathbb C)^{sm}/S)$ such that $a_{usu}H^n(j^*F^{\bullet})=0$ for all $n\in\mathbb Z$,
where $j:S\backslash Z\hookrightarrow S$ is the complementary open embedding and $a_{usu}$ is the sheaftification functor.
\item[(i)'] We denote by 
\begin{equation*}
C_{O_S,Z}(\AnSp(\mathbb C)^{sm}/S)\subset C_{O_S}(\AnSp(\mathbb C)^{sm}/S), 
\end{equation*}
the full subcategory consisting of complexes of presheaves 
$F^{\bullet}\in C(\AnSp(\mathbb C)^{sm}/S)$ such that $a_{usu}H^n(j^*F^{\bullet})=0$ for all $n\in\mathbb Z$,
where $j:S\backslash Z\hookrightarrow S$ is the complementary open embedding and $a_{usu}$ is the sheaftification functor.
\item[(ii)] We denote by 
\begin{equation*}
C_{fil,Z}(\AnSp(\mathbb C)^{sm}/S)\subset C_{fil}(\AnSp(\mathbb C)^{sm}/S) 
\end{equation*}
the full subcategory consisting of filtered complexes of presheaves $(F^{\bullet},F)\in C_{fil}(\AnSp(\mathbb C)^{sm}/S)$ 
such that there exist $r\in\mathbb N$ and an $r$-filtered homotopy equivalence $\phi:(F^{\bullet},F)\to(F^{'\bullet},F)$
with $(F^{'\bullet},F)\in C_{fil}(\AnSp(\mathbb C)^{sm}/S)$
such that $a_{usu}j^*H^n\Gr_F^p(F^{'\bullet},F)=0$ for all $n,p\in\mathbb Z$,
where $j:S\backslash Z\hookrightarrow S$ is the complementary open embedding and $a_{usu}$ is the sheaftification functor.
\item[(ii)'] We denote by 
\begin{equation*}
C_{O_Sfil,Z}(\AnSp(\mathbb C)^{sm}/S)\subset C_{O_Sfil}(\AnSp(\mathbb C)^{sm}/S) 
\end{equation*}
the full subcategory consisting of filtered complexes of presheaves $(F^{\bullet},F)\in C_{O_Sfil}(\AnSp(\mathbb C)^{sm}/S)$ 
such that there exist $r\in\mathbb N$ and an $r$-filtered homotopy equivalence $\phi:(F^{\bullet},F)\to(F^{'\bullet},F)$
with $(F^{'\bullet},F)\in C_{O_Sfil}(\AnSp(\mathbb C)^{sm}/S)$
such that $a_{usu}j^*H^n\Gr_F^p(F^{\bullet},F)=0$ for all $n,p\in\mathbb Z$,
where $j:S\backslash Z\hookrightarrow S$ is the complementary open embedding and $a_{usu}$ is the sheaftification functor.
\end{itemize}
\end{defi}

Let $S\in\AnSp(\mathbb C)$ and $Z\subset S$ a closed subset. 
\begin{itemize}
\item For $(G,F)\in C_{fil}(\AnSp(\mathbb C)^{sm}/S)$, we have 
$\Gamma_Z(G,F),\Gamma^{\vee}_Z(G,F)\in C_{fil,Z}(\AnSp(\mathbb C)^{sm}/S)$.
\item For $(G,F)\in C_{O_Sfil}(\AnSp(\mathbb C)^{sm}/S)$, we have 
$\Gamma_Z(G,F),\Gamma^{\vee}_Z(G,F)\in C_{O_Sfil,Z}(\AnSp(\mathbb C)^{sm}/S)$.
\end{itemize}

Let $S_{\bullet}\in\Fun(\mathcal I,\AnSp(\mathbb C))$ with $\mathcal I\in\Cat$, a diagram of algebraic varieties.
It gives the diagram of sites $\AnSp(\mathbb C)^2/S_{\bullet}\in\Fun(\mathcal I,\Cat)$.  
\begin{itemize}
\item Then $C_{fil}(\AnSp(\mathbb C)/S_{\bullet})$ is the category  
\begin{itemize}
\item whose objects $(G,F)=((G_I,F)_{I\in\mathcal I},u_{IJ})$,
with $(G_I,F)\in C_{fil}(\AnSp(\mathbb C)/S_I)$,
and $u_{IJ}:(G_I,F)\to r_{IJ*}(G_J,F)$ for $r_{IJ}:I\to J$, denoting again $r_{IJ}:S_I\to S_J$, are morphisms
satisfying for $I\to J\to K$, $r_{IJ*}u_{JK}\circ u_{IJ}=u_{IK}$ in $C_{fil}(\AnSp(\mathbb C)/S_I)$,
\item the morphisms $m:((G,F),u_{IJ})\to((H,F),v_{IJ})$ being (see section 2.1) a family of morphisms of complexes,  
\begin{equation*}
m=(m_I:(G_I,F)\to (H_I,F))_{I\in\mathcal I}
\end{equation*}
such that $v_{IJ}\circ m_I=p_{IJ*}m_J\circ u_{IJ}$ in $C_{fil}(\Var(\mathbb C)/S_I)$.
\end{itemize}
\item Then $C_{fil}(\AnSp(\mathbb C)^{sm}/S_{\bullet})$ is the category  
\begin{itemize}
\item whose objects $(G,F)=((G_I,F)_{I\in\mathcal I},u_{IJ})$,
with $(G_I,F)\in C_{fil}(\AnSp(\mathbb C)^{sm}/S_I)$,
and $u_{IJ}:(G_I,F)\to r_{IJ*}(G_J,F)$ for $r_{IJ}:I\to J$, denoting again $r_{IJ}:S_I\to S_J$, are morphisms
satisfying for $I\to J\to K$, $r_{IJ*}u_{JK}\circ u_{IJ}=u_{IK}$ in $C_{fil}(\AnSp(\mathbb C)^{sm}/S_I)$,
\item the morphisms $m:((G,F),u_{IJ})\to((H,F),v_{IJ})$ being (see section 2.1) a family of morphisms of complexes,  
\begin{equation*}
m=(m_I:(G_I,F)\to (H_I,F))_{I\in\mathcal I}
\end{equation*}
such that $v_{IJ}\circ m_I=p_{IJ*}m_J\circ u_{IJ}$ in $C_{fil}(\AnSp(\mathbb C)^{sm}/S_I)$.
\end{itemize}
\end{itemize}
As usual, we denote by
\begin{eqnarray*}
(f_{\bullet}^*,f_{\bullet*}):=(P(f_{\bullet})^*,P(f_{\bullet})_*): 
C(\AnSp(\mathbb C)^{(sm)}/S_{\bullet})\to C(\AnSp(\mathbb C)^{(sm)}/T_{\bullet})
\end{eqnarray*}
the adjonction induced by 
$P(f_{\bullet}):\AnSp(\mathbb C)^{(sm)}/T_{\bullet}\to \AnSp(\mathbb C)^{(sm)}/S_{\bullet}$.
Since the colimits involved in the definition of $f_{\bullet}^*=P(f_{\bullet})^*$ are filtered, 
$f_{\bullet}^*$ also preserve monomorphism. Hence, we get an adjonction
\begin{eqnarray*}
(f_{\bullet}^*,f_{\bullet*}):
C_{fil}(\AnSp(\mathbb C)^{(sm)}/S_{\bullet})\leftrightarrows C_{fil}(\AnSp(\mathbb C)^{(sm)}/T_{\bullet}), \\
f_{\bullet}^*((G_I,F),u_{IJ}):=((f_I^*G_I,f_I^*F),T(f_I,r_{IJ})(-)\circ f_I^*u_{IJ}).
\end{eqnarray*}

Let $S\in\AnSp(\mathbb C)$. Let $S=\cup_{i=1}^l S_i$ an open cover and denote by $S_I=\cap_{i\in I} S_i$.
Let $i_i:S_i\hookrightarrow\tilde S_i$ closed embeddings, with $\tilde S_i\in\AnSp(\mathbb C)$. 
For $I\subset\left[1,\cdots l\right]$, denote by $\tilde S_I=\Pi_{i\in I}\tilde S_i$.
We then have closed embeddings $i_I:S_I\hookrightarrow\tilde S_I$, and for $J\subset I$ the following commutative diagram
\begin{equation*}
D_{IJ}=\xymatrix{ S_I\ar[r]^{i_I} & \tilde S_I \\
S_J\ar[u]^{j_{IJ}}\ar[r]^{i_J} & \tilde S_J\ar[u]^{p_{IJ}}}  
\end{equation*}
where $p_{IJ}:\tilde S_J\to\tilde S_I$ is the projection
and $j_{IJ}:S_J\hookrightarrow S_I$ is the open embedding so that $j_I\circ j_{IJ}=j_J$.
This gives the diagram of analytic spaces $(\tilde S_I)\in\Fun(\mathcal P(\mathbb N),\AnSp(\mathbb C))$ which
which gives the diagram of sites $\AnSp(\mathbb C)^{sm}/(\tilde S_I)\in\Fun(\mathcal P(\mathbb N),\Cat)$.
Denote by $m:\tilde S_I\backslash(S_I\backslash S_J)\hookrightarrow\tilde S_I$ the open embedding.
Then $C_{fil}(\AnSp(\mathbb C)^{sm}/(\tilde S_I)$ is the category  
\begin{itemize}
\item whose objects $(G,F)=((G_I,F),u_{IJ})$ with $(G_I,F)\in C_{fil}(\AnSp(\mathbb C)^{sm}/\tilde S_I)$,
and $u_{IJ}:(G_I,F)\to p_{IJ*}(G_J,F)$ are morphisms
satisfying for $I\subset J\subset K$, $p_{IJ*}u_{JK}\circ u_{IJ}=u_{IK}$ in $C_{fil}(\AnSp(\mathbb C)^{sm}/\tilde S_I)$,
\item the morphisms $m:((G,F),u_{IJ})\to((H,F),v_{IJ})$ being a family of morphisms of complexes,  
\begin{equation*}
m=(m_I:(G_I,F)\to (H_I,F))_{I\in\mathcal I}
\end{equation*}
such that $v_{IJ}\circ m_I=p_{IJ*}m_J\circ u_{IJ}$ in $C_{fil}(\AnSp(\mathbb C)^{sm}/\tilde S_I)$.
\end{itemize}

\begin{defi}
Let $S\in\AnSp(\mathbb C)$. Let $S=\cup_{i=1}^l S_i$ an open cover and denote by $S_I=\cap_{i\in I} S_i$.
Let $i_i:S_i\hookrightarrow\tilde S_i$ closed embeddings, with $\tilde S_i\in\AnSp(\mathbb C)$. 
We denote by the full subcategory 
$C_{fil}(\AnSp(\mathbb C)^{sm}/(S/(\tilde S_I)))\subset C_{fil}(\AnSp(\mathbb C)^{sm}/(\tilde S_I))$  
the full subcategory 
\begin{itemize}
\item whose objects $(G,F)=((G_I,F)_{I\subset\left[1,\cdots l\right]},u_{IJ})$,
with $(G_I,F)\in C_{fil,S_I}(\AnSp(\mathbb C)^{sm}/\tilde S_I)$,
and $u_{IJ}:m^*(G_I,F)\to m^*p_{IJ*}(G_J,F)$ for $I\subset J$, are $\infty$-filtered usu local equivalence,
satisfying for $I\subset J\subset K$, $p_{IJ*}u_{JK}\circ u_{IJ}=u_{IK}$ in $C_{fil}(\AnSp(\mathbb C)^{sm}/\tilde S_I)$,
\item the morphisms $m:((G,F),u_{IJ})\to((H,F),v_{IJ})$ being (see section 2.1) a family of morphisms of complexes,  
\begin{equation*}
m=(m_I:(G_I,F)\to (H_I,F))_{I\subset\left[1,\cdots l\right]}
\end{equation*}
such that $v_{IJ}\circ m_I=p_{IJ*}m_J\circ u_{IJ}$ in $C_{fil}(\AnSp(\mathbb C)^{sm}/\tilde S_I)$.
\end{itemize}
A morphism $m:((G_I,F),u_{IJ})\to((H_I,F),v_{IJ})$ is an $r$-filtered usu local equivalence,
if there exists $\phi_i:((C_{iI},F),u_{iIJ}\to((C_{(i+1)I},F),u_{(i+1)IJ})$, $0\leq i\leq s$, 
with $((C_{iI},F),u_{iIJ}\in C_{fil}(\AnSp(\mathbb C)^{2,(sm)}/(S/(\tilde S_I)))$
$((C_{0I},F),u_{0IJ})=((G_I,F),u_{IJ})$ and $((C_{sI},F),u_{sIJ}=((H_I,F),v_{IJ})$, such that
\begin{equation*}
\phi=\phi_s\circ\cdots\circ\phi_i\circ\cdots\circ\phi_0:((G_I,F),u_{IJ})\to((H_I,F),v_{IJ})
\end{equation*}
and $\phi_i:((C_{iI},F),u_{iIJ})\to((C_{(i+1)I},F),u_{(i+1)IJ})$ either a filtered usu local equivalence
or an $r$-filtered homotopy equivalence.  
\end{defi}

Denote $L=[1,\ldots,l]$ and for $I\subset L$,
$p_{0(0I)}:S\times\tilde S_I\to S$, $p_{I(0I)}:S\times\tilde S_I\to S_I$ the projections.
By definition, we have functors 
\begin{itemize}
\item $T(S/(\tilde S_I)):C_{fil}(\AnSp(\mathbb C)^{sm}/S)\to C_{fil}(\AnSp(\mathbb C)^{sm}/(S/(\tilde S_I)))$,
$(G,F)\mapsto (i_{I*}j_I^*F,T(D_{IJ})(j_I^*(G,F)))$,
\item $T((\tilde S_I)/S):C_{fil}(\AnSp(\mathbb C)^{sm}/(S/(\tilde S_I)))\to C_{fil}(\AnSp(\mathbb C)^{sm}/S)$,
$((G_I,F),u_{IJ})\mapsto \ho\lim_{I\subset L}p_{0(0I)*}\Gamma^{\vee}_{S_I}p_{I(0I)}^*(G_I,F)$. 
\end{itemize}
Note that the functors $T(S/(\tilde S_I)$ are NOT embedding, since
\begin{equation*}
\ad(i_I^*,i_{I*})(j_I^*F):i_I^*i_{I*}j_I^*F\to j_I^*F
\end{equation*}
are usu local equivalence but NOT isomorphism since we are dealing
with the morphism of big sites $P(i_I):\AnSp(\mathbb C)^{sm}/S_I\to\AnSp(\mathbb C)^{sm}/\tilde S_I$.
However, these functors induces full embeddings 
\begin{equation*}
T(S/(\tilde S_I)):D_{fil}(\AnSp(\mathbb C)^{sm}/S)\to D_{fil}(\AnSp(\mathbb C)^{sm}/(S/(\tilde S_I)))
\end{equation*}
since for $F\in C(\AnSp(\mathbb C)^{sm}/S)$, 
\begin{equation*}
\ho\lim_{I\subset L}p_{0(0I)*}\Gamma_{S_I}p_{I(0I)}^*(i_{I*}j_I^*F)\to p_{0(0I)*}\Gamma_{S_I}j_I^*F
\end{equation*}
is an equivalence usu local.

Let $f:X\to S$ a morphism, with $X,S\in\AnSp(\mathbb C)$.
Let $S=\cup_{i=1}^{l} S_i$ and $X=\cup_{i=1}^l X_i$ be affine open covers
and $i_i:S_i\hookrightarrow\tilde S_i$, $i'_i:X_i\hookrightarrow\tilde X_i$ be closed embeddings. 
Let $\tilde f_i:\tilde X_i\to\tilde S_i$ be a lift of the morphism $f_i=f_{|X_i}:X_i\to S_i$.
Then, $f_I=f_{|X_I}:X_I=\cap_{i\in I} X_i\to S_I=\cap_{i\in I} S_i$ lift to the morphism 
\begin{equation*}
\tilde f_I=\Pi_{i\in I}\tilde f_i:\tilde X_I=\Pi_{i\in I}\tilde X_i\to \tilde S_I=\Pi_{i\in I}\tilde S_i
\end{equation*}
Denote by $p_{IJ}:\tilde S_J\to\tilde S_I$ and $p'_{IJ}:\tilde X_J\to\tilde X_I$ the projections.  
Consider for $J\subset I$ the following commutative diagrams
\begin{equation*}
D_{IJ}=\xymatrix{ S_I\ar[r]^{i_I} & \tilde S_I \\
S_J\ar[u]^{j_{IJ}}\ar[r]^{i_J} & \tilde S_J\ar[u]^{p_{IJ}}} \;, \;  
D'_{IJ}=\xymatrix{ X_I\ar[r]^{i'_I} & \tilde X_I \\
X_J\ar[u]^{j'_{IJ}}\ar[r]^{i'_J} & \tilde X_J\ar[u]^{p'_{IJ}}} \; , \;  
D_{fI}=\xymatrix{ S_I\ar[r]^{i_I} & \tilde S_I \\
X_I\ar[u]^{f_I}\ar[r]^{i'_I} & \tilde X_I\ar[u]^{\tilde f_I}}  
\end{equation*}
We have then following commutative diagram
\begin{equation*}
\xymatrix{ \, & X_I\ar[r]^{n'_I} & \tilde X_I & \tilde X_I\backslash X_I\ar[l]^{n'_I} \\
i'_J:X_J\ar[r]^{l_{IJ}}\ar[ru]^{j'_{IJ}}& X_I\times X_I\times\tilde X_{J\backslash I}\ar[u]^{p'_{IJ}}\ar[r]^{n'_I\times I} & 
\tilde X_J\ar[u]^{p'_{IJ}} & \tilde X_J\backslash X_J\ar[l]^{n'_J}\ar[u]^{p'_{IJ}}}.
\end{equation*}
whose square are cartesian.
We then have the pullback functor
\begin{eqnarray*}
f^*:C_{(2)fil}(\AnSp(\mathbb C)^{sm}/S/(\tilde S_I))\to C_{(2)fil}(\AnSp(\mathbb C)^{sm}/X/(\tilde X_I)), \\
((G_I,F),u_{IJ})\mapsto f^*((G_I,F),u_{IJ}):=(\Gamma^{\vee}_{X_I}\tilde f_I^*(G_I,F),\tilde f_J^*u_{IJ})
\end{eqnarray*}
with 
\begin{eqnarray*}
\tilde f_J^*u_{IJ}:\Gamma^{\vee}_{X_I}\tilde f_I^*(G_I,F)
\xrightarrow{\ad(p_{IJ}^{'*},p'_{IJ*})(-)}p'_{IJ*}p_{IJ}^{'*}\Gamma^{\vee}_{X_I}\tilde f_I^*(G_I,F) 
\xrightarrow{T_{\sharp}(p_{IJ},n'_I)(-)^{-1})}
p'_{IJ*}\Gamma^{\vee}_{X_I\times\tilde X_{J\backslash I}}p_{IJ}^{'*}\tilde f_I^*(G_I,F) \\
\xrightarrow{p'_{IJ*}\gamma^{\vee}_{X_J}(-)}
p'_{IJ*}\Gamma^{\vee}_{X_J}p_{IJ}^{'*}\tilde f_I^*(G_I,F)=p'_{IJ*}\Gamma^{\vee}_{X_J}\tilde f_J^*p_{IJ}^*(G_I,F)
\xrightarrow{\Gamma^{\vee}_{X_J}\tilde f_J^*I(p_{IJ}^*,p_{IJ*})(-,-)(u_{IJ})}
\Gamma^{\vee}_{X_J}\tilde f_J^*(G_J,F)
\end{eqnarray*}
Let $(G,F)\in C_{fil}(\AnSp(\mathbb C)^{sm}/S)$. Since, $j_I^{'*}i'_{I*}j_I^{'*}f^*(G,F)=0$,
the morphism $T(D_{fI})(j_I^*(G,F)):\tilde f_I^*i_{I*}j_I^*(G,F)\to i'_{I*}j_I^{'*}f^*(G,F)$ factors trough
\begin{equation*}
T(D_{fI})(j_I^*(G,F)):\tilde f_I^*i_{I*}j_I^*(G,F)\xrightarrow{\gamma_{X_I}^{\vee}(-)}
\Gamma_{X_I}^{\vee}\tilde f_I^*i_{I*}j_I^*(G,F)\xrightarrow{T^{\gamma}(D_{fI})(j_I^*(G,F))} i'_{I*}j_I^{'*}f^*(G,F)
\end{equation*}
We have then, for $(G,F)\in C_{fil}(S)$, the canonical transformation map
\begin{equation*}
\xymatrix{
f^*T(S/(\tilde S_I))(G,F)\ar[rrr]^{T(f,T(0/I))(G,F)}\ar[d]_{=}\ & \, & \, & 
T(X/(\tilde X_I))(f^*(G,F))\ar[d]^{=} \\
(\Gamma^{\vee}_{X_I}\tilde f_I^*i_{I*}j_I^*(G,F),\tilde f_J^*I)\ar[rrr]^{T^{\gamma}(D_{fI})(j_I^*(G,F))} 
& \, & \, & (i'_{I*}j_I^{'*}f^*(G,F),I)}
\end{equation*}

We have similarly to the algebraic case, we have:

\begin{defi}\label{projBMmotandef}
\begin{itemize}
\item[(i)]Let $f:X\to S$ a morphism with $X,S\in\AnSp(\mathbb C)$. Assume that there exist a factorization
$f:X\xrightarrow{i} Y\times S\xrightarrow{p} S$, with $Y\in\AnSm(\mathbb C)$, 
$i:X\hookrightarrow Y$ is a closed embedding and $p$ the projection. We then consider 
\begin{equation*}
Q(X/S):=p_{\sharp}\Gamma^{\vee}_X\mathbb Z_{Y\times S}:=
\Cone(\mathbb Z((Y\times S)\backslash X/S)\to\mathbb Z(Y\times S/S))\in C(\AnSp(\mathbb C)^{sm}/S).
\end{equation*}
By definition $Q(X/S)$ is projective since it is a complex of two representative presheaves.
\item[(ii)]Let $f:X\to S$ and $g:T\to S$ two morphism with $X,S,T\in\AnSp(\mathbb C)$. Assume that there exist a factorization
$f:X\xrightarrow{i} Y\times S\xrightarrow{p} S$, with $Y\in\AnSm(\mathbb C)$, 
$i:X\hookrightarrow Y$ is a closed embedding and $p$ the projection. 
We then have the following commutative diagram whose squares are cartesian
\begin{equation*}
\xymatrix{
f:X\ar[r]^{i} & Y\times S\ar[r]^{p} & S \\
f':X_T\ar[r]^{i'}\ar[u]^{g'} & Y\times T\ar[r]^{p'}\ar[u]^{g'':=(I\times g)} & T\ar[u]^{g}}
\end{equation*}
We then have the canonical isomorphism in $C(\AnSp(\mathbb C)^{sm}/T)$
\begin{eqnarray*}
T(f,g,Q):=T_{\sharp}(g,p)(-)^{-1}\circ T_{\sharp}(g'',j)(-)^{-1}: \\
g^*Q(X/S):=g^*p_{\sharp}\Gamma^{\vee}_X\mathbb Z_{Y\times S}
\xrightarrow{\sim}p'_{\sharp}\Gamma^{\vee}_{X_T}\mathbb Z_{Y\times T}=:Q(X_T/T)
\end{eqnarray*}
with $j:Y\times S\backslash X\hookrightarrow Y\times S$ the closed embedding.
\end{itemize}
\end{defi}

Let $S\in\AnSp(\mathbb C)$. Denote for short $\AnSp(\mathbb C)^{(sm)}/S$ 
either the category $\AnSp(\mathbb C)/S$ or the category $\AnSp(\mathbb C)^{sm}/S$. Denote by
\begin{eqnarray*}
p_a:\AnSp(\mathbb C)^{(sm)}/S\to\AnSp(\mathbb C)^{(sm)}/S, \;  
X/S=(X,h)\mapsto (X\times\mathbb D^1)/S=(X\times\mathbb D^1,h\circ p_X), \\ 
(g:X/S\to X'/S)\mapsto ((g\times I_{\mathbb D^1}):X\times\mathbb D^1/S\to X'\times\mathbb D^1/S)
\end{eqnarray*}
the projection functor and again by $p_a:\AnSp(\mathbb C)^{(sm)}/S\to\AnSp(\mathbb C)^{(sm)}/S$
the corresponding morphism of site.

We now define the $\mathbb D^1$ localization property :
\begin{defi}\label{d1locdef}
Let $S\in\Var(\mathbb C)$. 
\begin{itemize}
\item[(i0)] A complex $F\in C(\AnSp(\mathbb C)^{(sm)}/S)$ is said to be $\mathbb D^1$ homotopic if 
$\ad(p_a^*,p_{a*})(F):F\to p_{a*}p_a^*F$ is an homotopy equivalence.
\item[(i)] A complex $F\in C(\AnSp(\mathbb C)^{(sm)}/S)$ is said to be $\mathbb D^1$ invariant if for all 
$U/S\in\AnSp(\mathbb C)^{(sm)}/S$,
\begin{equation*}
F(p_U):F(U/S)\to F(U\times\mathbb D^1/S) 
\end{equation*}
is a quasi-isomorphism, where $p_U:U\times\mathbb D^1\to U$ is the projection.
Obviously, a $\mathbb D^1$ homotopic complex is $\mathbb D^1$ invariant.
\item[(ii)] A complex $F\in C(\AnSp(\mathbb C)^{(sm)}/S)$ is said to be $\mathbb D^1$ local for
the usual topology, if for a (hence every) usu local equivalence $k:F\to G$ with $k$ injective and
$G\in C(\AnSp(\mathbb C)^{(sm)}/S)$ usu fibrant, e.g. $k:F\to E_{usu}(F)$, 
$G$ is $\mathbb D^1$ invariant for all $n\in\mathbb Z$.
\item[(iii)] A morphism $m:F\to G$ with $F,G\in C(\AnSp(\mathbb C)^{(sm)}/S)$ is said to an $(\mathbb D^1,usu)$ local equivalence 
if for all $H\in C(\AnSp(\mathbb C)^{(sm)}/S)$ which is $\mathbb A^1$ local for the etale topology
\begin{equation*}
\Hom(L(m),E_{usu}(H)):\Hom(L(G),E_{usu}(H))\to\Hom(L(F),E_{usu}(H)) 
\end{equation*}
is a quasi-isomorphism.
\end{itemize}
\end{defi}

\begin{prop}\label{d1loceqprop}
A morphism $m:F\to G$ with $F,G\in C(\AnSp(\mathbb C)^{(sm)}/S)$ is an $(\mathbb D^1,et)$ local equivalence
if and only if there exists 
\begin{equation*}
\left\{X_{1,\alpha}/S,\alpha\in\Lambda_1\right\},\ldots,\left\{X_{r,\alpha}/S,\alpha\in\Lambda_r\right\}
\subset\AnSp(\mathbb C)^{(sm)}/S
\end{equation*}
such that we have in $\Ho_{et}(C(\Var(\mathbb C)^{(sm)}/S))$
\begin{equation*}
\Cone(m)\xrightarrow{\sim}\Cone(\oplus_{\alpha\in\Lambda_1}
\Cone(\mathbb Z(X_{1,\alpha}\times\mathbb D^1/S)\to\mathbb Z(X_{1,\alpha}/S))
\to\cdots\to\oplus_{\alpha\in\Lambda_r}\Cone(\mathbb Z(X_{r,\alpha}\times\mathbb D^1/S)\to\mathbb Z(X_{r,\alpha}/S)))
\end{equation*}
\end{prop}

\begin{proof}
Standard.
\end{proof}

\begin{defiprop}
Let $S\in\AnSp(\mathbb C)$
\begin{itemize}
\item[(i)] With the weak equivalence the $(\mathbb D^1,usu)$ local equivalence and 
the fibration the epimorphism with $\mathbb D^1_S$ local and usu fibrant kernels gives
a model structure on  $C(\AnSp(\mathbb C)^{sm}/S)$ : the left bousfield localization
of the projective model structure of $C(\AnSp(\mathbb C)^{sm}/S)$.
We call it the $(\mathbb D^1,usu)$ projective model structure.
\item[(ii)] With the weak equivalence the $(\mathbb D^1,usu)$ local equivalence and 
the fibration the epimorphism with $\mathbb D^1_S$ local and usu fibrant kernels gives
a model structure on  $C(\AnSp(\mathbb C)/S)$ : the left bousfield localization
of the projective model structure of $C(\AnSp(\mathbb C)/S)$.
We call it the $(\mathbb D^1,usu)$ projective model structure.
\end{itemize}
\end{defiprop}

\begin{proof}
Similar to the proof of definition-proposition \ref{projmodstr}.
\end{proof}

\begin{prop}\label{g1an}
Let $g:T\to S$ a morphism with $T,S\in\AnSp(\mathbb C)$.
\begin{itemize}
\item[(i)] The adjonction $(g^*,g_*):C(\AnSp(\mathbb C)^{sm}/S)\leftrightarrows C(\AnSp(\mathbb C)^{sm}/T)$
is a Quillen adjonction for the $(\mathbb D^1,et)$ projective model structure.
\item[(i)'] Let $h:U\to S$ a smooth morphism with $U,S\in\AnSp(\mathbb C)$.
The adjonction $(h_{\sharp},h^*):C(\AnSp(\mathbb C)^{sm}/U)\leftrightarrows C(\AnSp(\mathbb C)^{sm}/S)$
is a Quillen adjonction for the $(\mathbb D^1,et)$ projective model structure.
\item[(i)''] The functor $g^*:C(\AnSp(\mathbb C)^{sm}/S)\to C(\AnSp(\mathbb C)^{sm}/T)$
sends quasi-isomorphism to quasi-isomorphism, sends equivalence usu local to equivalence usu local,
sends $(\mathbb D^1,et)$ local equivalence to $(\mathbb D^1,et)$ local equivalence.
\item[(ii)] The adjonction $(g^*,g_*):C(\AnSp(\mathbb C)/S)\leftrightarrows C(\AnSp(\mathbb C)/T)$
is a Quillen adjonction for the $(\mathbb D^1,et)$ projective model structure (see definition \ref{projmodstr}).
\item[(ii)'] The adjonction $(g_{\sharp},g^*):C(\AnSp(\mathbb C)/T)\leftrightarrows C(\AnSp(\mathbb C)/S)$
is a Quillen adjonction for the $(\mathbb D^1,et)$ projective model structure (see definition \ref{projmodstr}).
\item[(ii)''] The functor $g^*:C(\AnSp(\mathbb C)/S)\to C(\AnSp(\mathbb C)/T)$
sends quasi-isomorphism to quasi-isomorphism, sends equivalence usu local to equivalence usu local,
sends $(\mathbb D^1,et)$ local equivalence to $(\mathbb D^1,et)$ local equivalence.
\end{itemize}
\end{prop}

\begin{proof}
Similar to the proof of proposition \ref{g1}.
\end{proof}

\begin{prop}\label{rho1an}
Let $S\in\AnSp(\mathbb C)$. 
\begin{itemize}
\item[(i)] The adjonction $(\rho_S^*,\rho_{S*}):C(\AnSp(\mathbb C)^{sm}/S)\leftrightarrows C(\AnSp(\mathbb C)/S)$
is a Quillen adjonction for the $(\mathbb D^1,et)$ projective model structure.
\item[(ii)]The functor $\rho_{S*}:C(\AnSp(\mathbb C)/S)\to C(\AnSp(\mathbb C)^{sm}/S)$
sends quasi-isomorphism to quasi-isomorphism, sends equivalence usu local to equivalence usu local,
sends $(\mathbb D^1,usu)$ local equivalence to $(\mathbb D^1,usu)$ local equivalence.
\end{itemize}
\end{prop}

\begin{proof}
Similar to the proof of proposition \ref{rho1}.
\end{proof}

Let $S\in\AnSp(\mathbb C)$. Let $S=\cup_{i=1}^l S_i$ an open affinoid cover and denote by $S_I=\cap_{i\in I} S_i$.
Let $i_i:S_i\hookrightarrow\tilde S_i$ closed embeddings, with $\tilde S_i\in\AnSp(\mathbb C)$.
For $(G_I,K_{IJ})\in C(\AnSp(\mathbb C)^{(sm)}/(\tilde S_I)^{op})$ and 
$(H_I,T_{IJ})\in C(\AnSp(\mathbb C)^{(sm)}/(\tilde S_I))$, we denote
\begin{eqnarray*}
\mathcal Hom((G_I,K_{IJ}),(H_I,T_{IJ})):=(\mathcal Hom(G_I,H_I),u_{IJ}((G_I,K_{IJ}),(H_I,T_{IJ})))
\in C(\AnSp(\mathbb C)^{(sm)}/(\tilde S_I))
\end{eqnarray*}
with
\begin{eqnarray*}
u_{IJ}((G_I,K_{IJ})(H_I,T_{IJ})):\mathcal Hom(G_I,H_I) \\
\xrightarrow{\ad(p_{IJ}^*,p_{IJ*})(-)}p_{IJ*}p_{IJ}^*\mathcal Hom(G_I,H_I)
\xrightarrow{T(p_{IJ},hom)(-,-)}p_{IJ*}\mathcal Hom(p_{IJ}^*G_I,p_{IJ}^*H_I) \\
\xrightarrow{\mathcal Hom(p_{IJ}^*G_I,T_{IJ})}p_{IJ*}\mathcal Hom(p_{IJ}^*G_I,H_J)
\xrightarrow{\mathcal Hom(K_{IJ},H_J)}p_{IJ*}\mathcal Hom(G_J,H_J).
\end{eqnarray*}
This gives in particular the functor
\begin{eqnarray*}
C(\AnSp(\mathbb C)^{(sm)}/(\tilde S_I))\to C(\AnSp(\mathbb C)^{(sm)}/(\tilde S_I)^{op}),
(H_I,T_{IJ})\mapsto(H_I,T_{IJ}).
\end{eqnarray*}
Let $S\in\AnSp(\mathbb C)$. Let $S=\cup_{i=1}^l S_i$ an open affinoid cover and denote by $S_I=\cap_{i\in I} S_i$.
Let $i_i:S_i\hookrightarrow\tilde S_i$ closed embeddings, with $\tilde S_i\in\AnSm(\mathbb C)$.
The functor $p_a$ extend to a functor
\begin{eqnarray*}
p_a:\AnSp(\mathbb C)^{(sm)}/(\tilde S_I)\to\AnSp(\mathbb C)^{(sm)}/(\tilde S_I), \;  
(X_I/\tilde S_I,s_{IJ})\mapsto (X\times\mathbb D^1/\tilde S_I,s_{IJ}\times I), \\ 
(g=(g_I):(X_I/\tilde S_I,s_{IJ})\to (X'_I/\tilde S_I,s'_{IJ}))\mapsto 
((g\times I_{\mathbb D^1}):(X_I\times\mathbb D^1/\tilde S_I,s_{IJ})\to(X'_I\times\mathbb D^1/\tilde S_I,s'_{IJ}))
\end{eqnarray*}
the projection functor and again by $p_a:\AnSp(\mathbb C)^{(sm)}/(\tilde S_I)\to\AnSp(\mathbb C)^{(sm)}/(\tilde S_I)$
the corresponding morphism of site.

\begin{defi}\label{d1locdefIJ}
Let $S\in\AnSp(\mathbb C)$. Let $S=\cup_{i=1}^l S_i$ an open affinoid cover and denote by $S_I=\cap_{i\in I} S_i$.
Let $i_i:S_i\hookrightarrow\tilde S_i$ closed embeddings, with $\tilde S_i\in\AnSp(\mathbb C)$.
\begin{itemize}
\item[(i0)]A complex $(F_I,u_{IJ})\in C(\AnSp(\mathbb C)^{(sm)}/(\tilde S_I))$ is said to be $\mathbb D^1$ homotopic if 
$\ad(p_a^*,p_{a*})((F_I,u_{IJ})):(F_I,u_{IJ})\to p_{a*}p_a^*(F_I,u_{IJ})$ is an homotopy equivalence.
\item[(i)] A complex  $(F_I,u_{IJ})\in C(\AnSp(\mathbb C)^{(sm)}/(\tilde S_I))$ is said to be $\mathbb D^1$ invariant 
if for all $(X_I/\tilde S_I,s_{IJ})\in\AnSp(\mathbb C)^{(sm)}/(\tilde S_I)$ 
\begin{equation*}
(F_I(p_{X_I})):(F_I(X_I/\tilde S_I),F_J(s_{IJ}))\circ u_{IJ}(-)\to 
(F_I(X_I\times\mathbb A^1/\tilde S_I),F_J(s_{IJ}\times I)\circ u_{IJ}(-)) 
\end{equation*}
is a quasi-isomorphism, where $p_{X_I}:X_I\times\mathbb A^1\to X_I$ are the projection,
and $s_{IJ}:X_I\times\tilde S_{J\backslash I}/\tilde S_J\to X_J/\tilde S_J$.
Obviously a complex $(F_I,u_{IJ})\in C(\AnSp(\mathbb C)^{(sm)}/(\tilde S_I))$ is $\mathbb D^1$ invariant
if and only if all the $F_I$ are $\mathbb A^1$ invariant.
\item[(ii)]Let $\tau$ a topology on $\AnSp(\mathbb C)$. 
A complex $F=(F_I,u_{IJ})\in C(\AnSp(\mathbb C)^{(sm)}/(\tilde S_I))$ is said to be $\mathbb D^1$ local 
for the $\tau$ topology induced on $\AnSp(\mathbb C)/(\tilde S_I)$, 
if for an (hence every) $\tau$ local equivalence $k:F\to G$ with $k$ injective and 
$G=(G_I,v_{IJ})\in C(\AnSp(\mathbb C)^{(sm)}/(\tilde S_I))$ $\tau$ fibrant,
e.g. $k:(F_I,u_{IJ})\to (E_{\tau}(F_I),E(u_{IJ}))$, $G$ is $\mathbb D^1$ invariant.
\item[(iii)] A morphism $m=(m_I):(F_I,u_{IJ})\to (G_I,v_{IJ})$ with 
$(F_I,u_{IJ}),(G_I,v_{IJ})\in C(\AnSp(\mathbb C)^{(sm)}/(\tilde S_I))$ 
is said to be an $(\mathbb D^1,usu)$ local equivalence 
if for all $(H_I,w_{IJ})\in C(\AnSp(\mathbb C)^{(sm)}/(\tilde S_I))$ which is $\mathbb D^1$ local for the usual topology
\begin{eqnarray*}
(\Hom(L(m_I),E_{usu}(H_I))):\Hom(L(G_I,v_{IJ}),E_{et}(H_I,w_{IJ}))\to\Hom(L(F_I,u_{IJ}),E_{usu}(H_I,w_{IJ})) 
\end{eqnarray*}
is a quasi-isomorphism.
Obviously, if a morphism $m=(m_I):(F_I,u_{IJ})\to (G_I,v_{IJ})$ with 
$(F_I,u_{IJ}),(G_I,u_{IJ})\in C(\AnSp(\mathbb C)^{(sm)}/(\tilde S_I)^{op})$ 
is an $(\mathbb D^1,usu)$ local equivalence, 
then all the $m_I:F_I\to G_I$ are $(\mathbb D^1,usu)$ local equivalence.
\item[(iv)] A morphism $m=(m_I):(F_I,u_{IJ})\to (G_I,v_{IJ})$ with 
$(F_I,u_{IJ}),(G_I,v_{IJ})\in C(\AnSp(\mathbb C)^{(sm)}/(\tilde S_I)^{op})$ 
is said to be an $(\mathbb D^1,usu)$ local equivalence 
if for all $(H_I,w_{IJ})\in C(\AnSp(\mathbb C)^{(sm)}/(\tilde S_I))$ which is $\mathbb D^1$ local for the etale topology
\begin{eqnarray*}
(\Hom(L(m_I),E_{et}(H_I))):\Hom(L(G_I,v_{IJ}),E_{usu}(H_I,w_{IJ}))\to\Hom(L(F_I,u_{IJ}),E_{usu}(H_I,w_{IJ})) 
\end{eqnarray*}
is a quasi-isomorphism.
Obviously, if a morphism $m=(m_I):(F_I,u_{IJ})\to (G_I,v_{IJ})$ with 
$(F_I,u_{IJ}),(G_I,u_{IJ})\in C(\AnSp(\mathbb C)^{(sm)}/(\tilde S_I)^{op})$ 
is an $(\mathbb D^1,usu)$ local equivalence, 
then all the $m_I:F_I\to G_I$ are $(\mathbb D^1,usu)$ local equivalence.
\end{itemize}
\end{defi}

\begin{prop}\label{cd1VarIJ}
Let $S\in\AnSp(\mathbb C)$. Let $S=\cup_{i=1}^l S_i$ an open affinoid cover and denote by $S_I=\cap_{i\in I} S_i$.
Let $i_i:S_i\hookrightarrow\tilde S_i$ closed embeddings, with $\tilde S_i\in\AnSp(\mathbb C)$.
\begin{itemize}
\item[(i)]A morphism $m:F\to G$ with $F,G\in C(\AnSp(\mathbb C)^{(sm)}/(\tilde S_I))$ 
is an $(\mathbb D^1,usu)$ local equivalence if and only if there exists 
\begin{eqnarray*}
\left\{(X_{1,\alpha,I}/\tilde S_I,u^1_{IJ}),\alpha\in\Lambda_1\right\},\ldots,
\left\{(X_{r,\alpha,I}/\tilde S_I,u^r_{IJ}),\alpha\in\Lambda_r\right\}
\subset\AnSp(\mathbb C)^{(sm)}/(\tilde S_I)
\end{eqnarray*}
with $u^l_{IJ}:X_{l,\alpha,I}\times\tilde S_{J\backslash I}/\tilde S_J\to X_{l,\alpha,J}/\tilde S_J$,
such that we have in $\Ho_{usu}(C(\AnSp(\mathbb C)^{(sm)}/(\tilde S_I)))$
\begin{eqnarray*}
\Cone(m)\xrightarrow{\sim}\Cone( \oplus_{\alpha\in\Lambda_1} 
\Cone((\mathbb Z(X_{1,\alpha,I}\times\mathbb D^1/\tilde S_I),\mathbb Z(u_{IJ}^1\times I)) 
\to(\mathbb Z(X_{1,\alpha,I}/\tilde S_I),\mathbb Z(u_{IJ}^1))) \\
\to\cdots\to\oplus_{\alpha\in\Lambda_r}
\Cone((\mathbb Z(X_{r,\alpha,I}\times\mathbb D^1/\tilde S_I),\mathbb Z(u_{IJ}^r\times I)) 
\to(\mathbb Z(X_{r,\alpha,I}/\tilde S_I),\mathbb Z(u^r_{IJ}))))
\end{eqnarray*}
\item[(ii)]A morphism $m:F\to G$ with $F,G\in C(\AnSp(\mathbb C)^{(sm)}/(\tilde S_I)^{op})$ 
is an $(\mathbb D^1,usu)$ local equivalence if and only if there exists 
\begin{eqnarray*}
\left\{(X_{1,\alpha,I}/\tilde S_I,u^1_{IJ}),\alpha\in\Lambda_1\right\},\ldots,
\left\{(X_{r,\alpha,I}/\tilde S_I,u^r_{IJ}),\alpha\in\Lambda_r\right\}
\subset\AnSp(\mathbb C)^{(sm)}/(\tilde S_I)^{op}
\end{eqnarray*}
with $u^l_{IJ}:X_{l,\alpha,J}/\tilde S_J\to X_{l,\alpha,I}\times\tilde S_{J\backslash I}/\tilde S_J$,
such that we have in $\Ho_{et}(C(\AnSp(\mathbb C)^{(sm)}/(\tilde S_I)^{op}))$
\begin{eqnarray*}
\Cone(m)\xrightarrow{\sim}\Cone(\oplus_{\alpha\in\Lambda_1} 
\Cone((\mathbb Z(X_{1,\alpha,I}\times\mathbb D^1/\tilde S_I),\mathbb Z(u_{IJ}^1\times I)) 
\to(\mathbb Z(X_{1,\alpha,I}/\tilde S_I),\mathbb Z(u_{IJ}^1))) \\
\to\cdots\to\oplus_{\alpha\in\Lambda_r}
\Cone((\mathbb Z(X_{r,\alpha,I}\times\mathbb D^1/\tilde S_I),\mathbb Z(u_{IJ}^r\times I)) 
\to(\mathbb Z(X_{r,\alpha,I}/\tilde S_I),\mathbb Z(u^r_{IJ}))))
\end{eqnarray*}
\end{itemize}
\end{prop}

\begin{proof}
Standard. See Ayoub's Thesis for example.
\end{proof}

\subsection{Presheaves on the big analytical site of pairs}

We recall the definition given in subsection 5.1 :
For $S\in\AnSp(\mathbb C)$, $\AnSp(\mathbb C)^2/S:=\AnSp(\mathbb C)^2/(S,S)$ is by definition (see subsection 2.1)
the category whose set of objects is 
\begin{eqnarray*}
(\AnSp(\mathbb C)^2/S)^0:=
\left\{((X,Z),h), h:X\to S, \; Z\subset X \; \mbox{closed} \;\right\}\subset\AnSp(\mathbb C)/S\times\Top
\end{eqnarray*}
and whose set of morphisms between $(X_1,Z_1)/S=((X_1,Z_1),h_1),(X_1,Z_1)/S=((X_2,Z_2),h_2)\in\AnSp(\mathbb C)^2/S$
is the subset
\begin{eqnarray*}
\Hom_{\AnSp(\mathbb C)^2/S}((X_1,Z_1)/S,(X_2,Z_2)/S):= \\
\left\{(f:X_2\to X_2), \; \mbox{s.t.} \; h_1\circ f=h_2 \; \mbox{and} \; Z_1\subset f^{-1}(Z_2)\right\} 
\subset\Hom_{\AnSp(\mathbb C)}(X_1,X_2)
\end{eqnarray*}
The category $\AnSp(\mathbb C)^2$ admits fiber products : $(X_1,Z_1)\times_{(S,Z)}(X_2,Z_2)=(X_1\times_S X_2,Z_1\times_Z Z_2)$.
In particular, for $f:T\to S$ a morphism with $S,T\in\AnSp(\mathbb C)$, we have the pullback functor
\begin{equation*}
P(f):\AnSp(\mathbb C)^2/S\to\AnSp(\mathbb C)^2/T, P(f)((X,Z)/S):=(X_T,Z_T)/T, P(f)(g):=(g\times_S f)
\end{equation*}
and we note again $P(f):\AnSp(\mathbb C)^2/T\to\AnSp(\mathbb C)^2/S$ the corresponding morphism of sites.

We will consider in the construction of the filtered De Rham realization functor the
full subcategory $\AnSp(\mathbb C)^{2,sm}/S\subset\AnSp(\mathbb C)^2/S$ such that the first factor is a smooth morphism :
We will also consider, in order to obtain a complex of D modules in the construction of the filtered De Rham realization functor,
the restriction to the full subcategory $\AnSp(\mathbb C)^{2,pr}/S\subset\AnSp(\mathbb C)^2/S$ 
such that the first factor is a projection :

\begin{defi}\label{PAnSp12S}
\begin{itemize}
\item[(i)]Let $S\in\AnSp(\mathbb C)$. We denote by
\begin{equation*}
\rho_S:\AnSp(\mathbb C)^{2,sm}/S\hookrightarrow\AnSp(\mathbb C)^2/S 
\end{equation*}
the full subcategory
consisting of the objects $(U,Z)/S=((U,Z),h)\in\AnSp(\mathbb C)^2/S$ such that the morphism $h:U\to S$ is smooth.
That is, $\AnSp(\mathbb C)^{2,sm}/S$ is the category  
\begin{itemize}
\item whose objects are $(U,Z)/S=((U,Z),h)$, with $U\in\AnSp(\mathbb C)$, $Z\subset U$ a closed subset, 
and $h:U\to S$ a smooth morphism,
\item whose morphisms $g:(U,Z)/S=((U,Z),h_1)\to (U',Z')/S=((U',Z'),h_2)$ 
is a morphism $g:U\to U'$ of complex algebraic varieties such that $Z\subset g^{-1}(Z')$ and  $h_2\circ g=h_1$. 
\end{itemize}
We denote again $\rho_S:\AnSp(\mathbb C)^2/S\to\AnSp(\mathbb C)^{2,sm}/S$ the associated morphism of site. We have 
\begin{equation*}
r^s(S):\AnSp(\mathbb C)^2\xrightarrow{r(S):=r(S,S)}\AnSp(\mathbb C)^2/S\xrightarrow{\rho_S}\AnSp(\mathbb C)^{2,sm}/S
\end{equation*}
the composite morphism of site.
\item[(ii)]Let $S\in\AnSp(\mathbb C)$. We will consider the full subcategory 
\begin{equation*}
\mu_S:\AnSp(\mathbb C)^{2,pr}/S\hookrightarrow\AnSp(\mathbb C)^2/S
\end{equation*}
whose subset of object consist of those whose morphism is a projection to $S$ : 
\begin{eqnarray*}
(\AnSp(\mathbb C)^{2,pr}/S)^0:=\left\{((Y\times S,X),p), \; Y\in\AnSp(\mathbb C), \;
 p:Y\times S\to S \; \mbox{the projection}\right\}\subset(\AnSp(\mathbb C)^2/S)^0.
\end{eqnarray*}
\item[(iii)]We will consider the full subcategory 
\begin{equation*}
\mu_S:(\AnSp(\mathbb C)^{2,smpr}/S)\hookrightarrow\AnSp(\mathbb C)^{2,sm}/S
\end{equation*}
whose subset of object consist of those whose morphism is a smooth projection to $S$ : 
\begin{eqnarray*}
(\AnSp(\mathbb C)^{2,smpr}/S)^0:=\left\{((Y\times S,X),p), \; Y\in\SmVar(\mathbb C), \;
 p:Y\times S\to S \; \mbox{the projection}\right\}\subset(\AnSp(\mathbb C)^2/S)^0
\end{eqnarray*}
\end{itemize}
\end{defi}
For $f:T\to S$ a morphism with $T,S\in\AnSp(\mathbb C)$, we have by definition, the following commutative diagram of sites
\begin{equation}\label{mufan}
\xymatrix{\AnSp(\mathbb C)^2/T\ar[rr]^{\mu_T}\ar[dd]_{P(f)}\ar[rd]^{\rho_T} & \, & 
\AnSp(\mathbb C)^{2,pr}/T\ar[dd]^{P(f)}\ar[rd]^{\rho_T} & \, \\
\, & \AnSp(\mathbb C)^{2,sm}/T\ar[rr]^{\mu_T}\ar[dd]_{P(f)} & \, & \AnSp(\mathbb C)^{2,smpr}/T\ar[dd]^{P(f)} \\
\AnSp(\mathbb C)^2/S\ar[rr]^{\mu_S}\ar[rd]^{\rho_S} & \, & \AnSp(\mathbb C)^{2,pr}/S\ar[rd]^{\rho_S} & \, \\
\, & \AnSp(\mathbb C)^{2,sm}/S\ar[rr]^{\mu_S} & \, & \AnSp(\mathbb C)^{2,smpr}/S}.
\end{equation}

Recall we have (see subsection 2.1), for $S\in\Var(\mathbb C)$, the graph functor 
\begin{eqnarray*}
\Gr_S^{12}:\AnSp(\mathbb C)/S\to\AnSp(\mathbb C)^{2,pr}/S, \; X/S\mapsto\Gr_S^{12}(X/S):=(X\times S,X)/S, \\
(g:X/S\to X'/S)\mapsto\Gr_S^{12}(g):=(g\times I_S:(X\times S,X)\to(X'\times S,X'))
\end{eqnarray*}
For $f:T\to S$ a morphism with $T,S\in\AnSp(\mathbb C)$, we have by definition, the following commutative diagram of sites
\begin{equation}\label{Grfan}
\xymatrix{\AnSp(\mathbb C)^{2,pr}/T\ar[rr]^{\Gr_T^{12}}\ar[dd]_{P(f)}\ar[rd]^{\rho_T} & \, & 
\AnSp(\mathbb C)/T\ar[dd]^{P(f)}\ar[rd]^{\rho_T} & \, \\
\, & \AnSp(\mathbb C)^{2,smpr}/T\ar[rr]^{\Gr_T^{12}}\ar[dd]_{P(f)} & \, & \AnSp(\mathbb C)^{sm}/T\ar[dd]^{P(f)} \\
\AnSp(\mathbb C)^{2,pr}/S\ar[rr]^{\Gr_S^{12}}\ar[rd]^{\rho_S} & \, & \AnSp(\mathbb C)/S\ar[rd]^{\rho_S} & \, \\
\, & \AnSp(\mathbb C)^{2,sm}/S\ar[rr]^{\Gr_S^{12}} & \, & \AnSp(\mathbb C)^{sm}/S}.
\end{equation}
where we recall that $P(f)((X,Z)/S):=((X_T,Z_T)/T)$, since smooth morphisms are preserved by base change.

As usual, we denote by
\begin{equation*}
(f^*,f_*):=(P(f)^*,P(f)_*):C(\AnSp(\mathbb C)^{2,sm}/S)\to C(\AnSp(\mathbb C)^{2,sm}/T)
\end{equation*}
the adjonction induced by $P(f):\AnSp(\mathbb C)^{2,sm}/T\to \AnSp(\mathbb C)^{2,sm}/S$.
Since the colimits involved in the definition of $f^*=P(f)^*$ are filtered, $f^*$ also preserve monomorphism. Hence, we get an adjonction
\begin{equation*}
(f^*,f_*):C_{fil}(\AnSp(\mathbb C)^{2,sm}/S)\leftrightarrows C_{fil}(\AnSp(\mathbb C)^{2,sm}/T), \; f^*(G,F):=(f^*G,f^*F)
\end{equation*}
For $S\in\AnSp(\mathbb C)$, we denote by 
$\mathbb Z_S:=\mathbb Z((S,S)/(S,S))\in\PSh(\AnSp(\mathbb C)^{2,sm}/S)$ the constant presheaf
By Yoneda lemma, we have for $F\in C(\AnSp(\mathbb C)^{2,sm}/S)$, $\mathcal Hom(\mathbb Z_S,F)=F$.

For $h:U\to S$ a smooth morphism with $U,S\in\AnSp(\mathbb C)$, 
$P(h):\AnSp(\mathbb C)^{2,sm}/S\to\AnSp(\mathbb C)^{2,sm}/U$ admits a left adjoint
\begin{equation*}
C(h):\AnSp(\mathbb C)^{2,sm}/U\to\AnSp(\mathbb C)^{2,sm}/S, \; C(h)((U',Z'),h')=((U',Z'),h\circ h').
\end{equation*}
Hence $h^*:C(\AnSp(\mathbb C)^{2,sm}/S)\to C(\AnSp(\mathbb C)^{2,sm}/U)$ admits a left adjoint
\begin{equation*}
h_{\sharp}:C(\AnSp(\mathbb C)^{2,sm}/U)\to C(\AnSp(\mathbb C)^{2,sm}/S), \;
F\mapsto (h_{\sharp}F:((U,Z),h_0)\mapsto\lim_{((U',Z'),h\circ h')\to ((U,Z),h_0)} F((U',Z')/U))
\end{equation*}
For $F^{\bullet}\in C(\AnSp(\mathbb C)^{2,sm}/S)$ and $G^{\bullet}\in C(\AnSp(\mathbb C)^{2,sm}/U)$,
we have the adjonction maps
\begin{equation*}
\ad(h_{\sharp},h^*)(G^{\bullet}):G^{\bullet}\to h^*h_{\sharp}G^{\bullet} \; , \;
\ad(h_{\sharp},h^*)(F^{\bullet}):h_{\sharp}h^*F^{\bullet}\to F^{\bullet}.
\end{equation*}
For a smooth morphism $h:U\to S$, with $U,S\in\AnSp(\mathbb C)$, we have the adjonction isomorphism, 
for $F\in C(\AnSp(\mathbb C)^{2,sm}/U)$ and $G\in C(\AnSp(\mathbb C)^{2,sm}/S)$,   
\begin{equation}\label{Ihhom12an}
I(h_{\sharp},h^*)(F,G):\mathcal Hom^{\bullet}(h_{\sharp}F,G)\xrightarrow{\sim}h_*\mathcal Hom^{\bullet}(F,h^*G).  
\end{equation}

For a commutative diagram in $\AnSp(\mathbb C)$ : 
\begin{equation*}
D=\xymatrix{ 
V\ar[r]^{g_2}\ar[d]^{h_2} & U\ar[d]^{h_1} \\
T\ar[r]^{g_1} & S},
\end{equation*}
where $h_1$ and $h_2$ are smooth,
we denote by, for $F^{\bullet}\in C(\AnSp(\mathbb C)^{2,sm}/U)$, 
\begin{equation*}
T_{\sharp}(D)(F^{\bullet}): h_{2\sharp}g_2^*F^{\bullet}\to g_1^*h_{1\sharp}F^{\bullet}
\end{equation*}
the canonical map given by adjonction. 
If $D$ is cartesian with $h_1=h$, $g_1=g$ $f_2=h':U_T\to T$, $g':U_T\to U$,
\begin{equation*}
T_{\sharp}(D)(F^{\bullet})=:T_{\sharp}(g,h)(F):h'_{\sharp}g^{'*}F^{\bullet}\xrightarrow{\sim} g^*h_{\sharp}F^{\bullet}
\end{equation*}
is an isomorphism.

We have the support section functors of a closed embedding $i:Z\hookrightarrow S$ for presheaves on the big analytical site of pairs.
\begin{defi}\label{gamma12an}
Let $i:Z\hookrightarrow S$ be a closed embedding with $S,Z\in\AnSp(\mathbb C)$ 
and $j:S\backslash Z\hookrightarrow S$ be the open complementary subset.
\begin{itemize}
\item[(i)] We define the functor
\begin{equation*}
\Gamma_Z:C(\AnSp(\mathbb C)^{2,sm}/S)\to C(\AnSp(\mathbb C)^{2,sm}/S), \;
G^{\bullet}\mapsto\Gamma_Z G^{\bullet}:=\Cone(\ad(j^*,j_*)(G^{\bullet}):G^{\bullet}\to j_*j^*G^{\bullet})[-1],
\end{equation*}
so that there is then a canonical map $\gamma_Z(G^{\bullet}):\Gamma_ZG^{\bullet}\to G^{\bullet}$.
\item[(ii)] We have the dual functor of (i) :
\begin{equation*}
\Gamma^{\vee}_Z:C(\AnSp(\mathbb C)^{2,sm}/S)\to C(\AnSp(\mathbb C)^{2,sm}/S), \; 
F\mapsto\Gamma^{\vee}_Z(F^{\bullet}):=\Cone(\ad(j_{\sharp},j^*)(G^{\bullet}):j_{\sharp}j^*G^{\bullet}\to G^{\bullet}), 
\end{equation*}
together with the canonical map $\gamma^{\vee}_Z(G):F\to\Gamma^{\vee}_Z(G)$.
\item[(iii)] For $F,G\in C(\AnSp(\mathbb C)^{2,sm}/S)$, we denote by 
\begin{equation*}
I(\gamma,hom)(F,G):=(I,I(j_{\sharp},j^*)(F,G)^{-1}):\Gamma_Z\mathcal Hom(F,G)\xrightarrow{\sim}\mathcal Hom(\Gamma^{\vee}_ZF,G)
\end{equation*}
the canonical isomorphism given by adjonction.
\end{itemize}
\end{defi}

Note that we have similarly for $i:Z\hookrightarrow S$, $i':Z'\hookrightarrow Z$ closed embeddings,
$g:T\to S$ a morphism with $T,S,Z\in\AnSp(\mathbb C)$ and $F\in C(\AnSp(\mathbb C)^{2,sm}/S)$, 
the canonical maps in $C(\AnSp(\mathbb C)^{2,sm}/S)$
\begin{itemize}
\item $T(g,\gamma)(F):g^*\Gamma_ZF\xrightarrow{\sim}\Gamma_{Z\times_S T}g^*F$, 
$T(g,\gamma^{\vee})(F):\Gamma^{\vee}_{Z\times_S T}g^*F\xrightarrow{\sim}g^*\Gamma_ZF$
\item $T(Z'/Z,\gamma)(F):\Gamma_{Z'}F\to\Gamma_Z F$, $T(Z'/Z,\gamma^{\vee})(F):\Gamma^{\vee}_ZF\to\Gamma_{Z'}^{\vee}F$
\end{itemize}
but we will not use them in this article.

Let $S_{\bullet}\in\Fun(\mathcal I,\AnSp(\mathbb C))$ with $\mathcal I\in\Cat$, a diagram of algebraic varieties.
It gives the diagram of sites $\AnSp(\mathbb C)^2/S_{\bullet}\in\Fun(\mathcal I,\Cat)$.  
\begin{itemize}
\item Then $C_{fil}(\AnSp(\mathbb C)^{2,(sm)}/S_{\bullet})$ is the category  
\begin{itemize}
\item whose objects $(G,F)=((G_I,F)_{I\in\mathcal I},u_{IJ})$,
with $(G_I,F)\in C_{fil}(\AnSp(\mathbb C)^{2,(sm)}/S_I)$,
and $u_{IJ}:(G_I,F)\to r_{IJ*}(G_J,F)$ for $r_{IJ}:I\to J$, denoting again $r_{IJ}:S_I\to S_J$, are morphisms
satisfying for $I\to J\to K$, $r_{IJ*}u_{JK}\circ u_{IJ}=u_{IK}$ in $C_{fil}(\AnSp(\mathbb C)^{2,(sm)}/S_I)$,
\item the morphisms $m:((G,F),u_{IJ})\to((H,F),v_{IJ})$ being (see section 2.1) a family of morphisms of complexes,  
\begin{equation*}
m=(m_I:(G_I,F)\to (H_I,F))_{I\in\mathcal I}
\end{equation*}
such that $v_{IJ}\circ m_I=p_{IJ*}m_J\circ u_{IJ}$ in $C_{fil}(\Var(\mathbb C)^{2,(sm)}/S_I)$.
\end{itemize}
\item Then $C_{fil}(\AnSp(\mathbb C)^{2,(sm)pr}/S_{\bullet})$ is the category  
\begin{itemize}
\item whose objects $(G,F)=((G_I,F)_{I\in\mathcal I},u_{IJ})$,
with $(G_I,F)\in C_{fil}(\AnSp(\mathbb C)^{2,(sm)pr}/S_I)$,
and $u_{IJ}:(G_I,F)\to r_{IJ*}(G_J,F)$ for $r_{IJ}:I\to J$, denoting again $r_{IJ}:S_I\to S_J$, are morphisms
satisfying for $I\to J\to K$, $r_{IJ*}u_{JK}\circ u_{IJ}=u_{IK}$ in $C_{fil}(\AnSp(\mathbb C)^{2,(sm)}/S_I)$,
\item the morphisms $m:((G,F),u_{IJ})\to((H,F),v_{IJ})$ being (see section 2.1) a family of morphisms of complexes,  
\begin{equation*}
m=(m_I:(G_I,F)\to (H_I,F))_{I\in\mathcal I}
\end{equation*}
such that $v_{IJ}\circ m_I=p_{IJ*}m_J\circ u_{IJ}$ in $C_{fil}(\AnSp(\mathbb C)^{2,(sm)pr}/S_I)$.
\end{itemize}
\end{itemize}
For $s:\mathcal I\to\mathcal J$ a functor, with $\mathcal I,\mathcal J\in\Cat$, and
$f_{\bullet}:T_{\bullet}\to S_{s(\bullet)}$ a morphism with 
$T_{\bullet}\in\Fun(\mathcal J,\AnSp(\mathbb C))$ and $S_{\bullet}\in\Fun(\mathcal I,\AnSp(\mathbb C))$, 
we have by definition, the following commutative diagrams of sites
\begin{equation}\label{mufIJan}
\xymatrix{\AnSp(\mathbb C)^2/T_{\bullet}\ar[rr]^{\mu_{T_{\bullet}}}\ar[dd]_{P(f_{\bullet})}\ar[rd]^{\rho_{T_{\bullet}}} & \, & 
\AnSp(\mathbb C)^{2,pr}/T_{\bullet}\ar[dd]^{P(f_{\bullet})}\ar[rd]^{\rho_{T_{\bullet}}} & \, \\
\, & \AnSp(\mathbb C)^{2,sm}/T_{\bullet}\ar[rr]^{\mu_{T_{\bullet}}}\ar[dd]_{P(f_{\bullet})} & \, & 
\AnSp(\mathbb C)^{2,smpr}/T_{\bullet}\ar[dd]^{P(f_{\bullet})} \\
\AnSp(\mathbb C)^2/S_{\bullet}\ar[rr]^{\mu_{S_{\bullet}}}\ar[rd]^{\rho_{S_{\bullet}}} & \, & 
\AnSp(\mathbb C)^{2,pr}/S_{\bullet}\ar[rd]^{\rho_{S_{\bullet}}} & \, \\
\, & \AnSp(\mathbb C)^{2,sm}/S_{\bullet}\ar[rr]^{\mu_{S_{\bullet}}} & \, & \AnSp(\mathbb C)^{2,smpr}/S_{\bullet}}.
\end{equation}
and
\begin{equation}\label{GrfIJan}
\xymatrix{\AnSp(\mathbb C)^{2,pr}/T_{\bullet}
\ar[rr]^{\Gr_{T_{\bullet}}^{12}}\ar[dd]_{P(f_{\bullet})}\ar[rd]^{\rho_{T_{\bullet}}} & \, & 
\AnSp(\mathbb C)/T\ar[dd]^{P(f_{\bullet})}\ar[rd]^{\rho_{T_{\bullet}}} & \, \\
\, & \AnSp(\mathbb C)^{2,smpr}/T_{\bullet}\ar[rr]^{\Gr_{T}^{12}}\ar[dd]_{P(f_{\bullet})} & \, & 
\AnSp(\mathbb C)^{sm}/T_{\bullet}\ar[dd]^{P(f_{\bullet})} \\
\AnSp(\mathbb C)^{2,pr}/S_{\bullet}\ar[rr]^{\Gr_{S_{\bullet}}^{12}}\ar[rd]^{\rho_{S_{\bullet}}} & \, & 
\AnSp(\mathbb C)/S_{\bullet}\ar[rd]^{\rho_{S_{\bullet}}} & \, \\
\, & \AnSp(\mathbb C)^{2,sm}/S_{\bullet}\ar[rr]^{\Gr_{S_{\bullet}}^{12}} & \, & 
\AnSp(\mathbb C)^{sm}/S_{\bullet}}.
\end{equation}
Let $s:\mathcal I\to\mathcal J$ a functor, with $\mathcal I,\mathcal J\in\Cat$, and
$f_{\bullet}:T_{\bullet}\to S_{s(\bullet)}$ a morphism with 
$T_{\bullet}\in\Fun(\mathcal J,\AnSp(\mathbb C))$ and $S_{\bullet}\in\Fun(\mathcal I,\AnSp(\mathbb C))$.
\begin{itemize}
\item As usual, we denote by
\begin{equation*}
(f_{\bullet}^*,f_{\bullet*}):=(P(f_{\bullet})^*,P(f_{\bullet})_*):
C(\AnSp(\mathbb C)^{2,(sm)}/S_{\bullet})\to C(\AnSp(\mathbb C)^{2,(sm)}/T_{\bullet})
\end{equation*}
the adjonction induced by 
$P(f_{\bullet}):\AnSp(\mathbb C)^{2,(sm)}/T_{\bullet}\to\AnSp(\mathbb C)^{2,(sm)}/S_{\bullet}$.
Since the colimits involved in the definition of $f_{\bullet}^*=P(f_{\bullet})^*$ are filtered, 
$f_{\bullet}^*$ also preserve monomorphism. Hence, we get an adjonction
\begin{eqnarray*}
(f_{\bullet}^*,f_{\bullet*}):
C_{fil}(\AnSp(\mathbb C)^{2,(sm)}/S_{\bullet})\leftrightarrows C_{fil}(\AnSp(\mathbb C)^{2,(sm)}/T_{\bullet}), \\
f_{\bullet}^*((G_I,F),u_{IJ}):=((f_I^*G_I,f_I^*F),T(f_I,r_{IJ})(-)\circ f_I^*u_{IJ})
\end{eqnarray*}
\item As usual, we denote by
\begin{equation*}
(f_{\bullet}^*,f_{\bullet*}):=(P(f_{\bullet})^*,P(f_{\bullet})_*):
C(\AnSp(\mathbb C)^{2,(sm)pr}/S_{\bullet})\to C(\AnSp(\mathbb C)^{2,(sm)pr}/T_{\bullet})
\end{equation*}
the adjonction induced by 
$P(f_{\bullet}):\AnSp(\mathbb C)^{2,(sm)pr}/T_{\bullet}\to\AnSp(\mathbb C)^{2,(sm)pr}/S_{\bullet}$.
Since the colimits involved in the definition of $f_{\bullet}^*=P(f_{\bullet})^*$ are filtered, 
$f_{\bullet}^*$ also preserve monomorphism. Hence, we get an adjonction
\begin{eqnarray*}
(f_{\bullet}^*,f_{\bullet*}):
C_{fil}(\AnSp(\mathbb C)^{2,(sm)pr}/S_{\bullet})\leftrightarrows C_{fil}(\AnSp(\mathbb C)^{2,(sm)pr}/T_{\bullet}), \\
f_{\bullet}^*((G_I,F),u_{IJ}):=((f_I^*G_I,f_I^*F),T(f_I,r_{IJ})(-)\circ f_I^*u_{IJ})
\end{eqnarray*}
\end{itemize}

Let $S\in\AnSp(\mathbb C)$. Let $S=\cup_{i=1}^l S_i$ an open affinoid cover and denote by $S_I=\cap_{i\in I} S_i$.
Let $i_i:S_i\hookrightarrow\tilde S_i$ closed embeddings, with $\tilde S_i\in\AnSp(\mathbb C)$. 
For $I\subset\left[1,\cdots l\right]$, denote by $\tilde S_I=\Pi_{i\in I}\tilde S_i$.
We then have closed embeddings $i_I:S_I\hookrightarrow\tilde S_I$ and for $J\subset I$ the following commutative diagram
\begin{equation*}
D_{IJ}=\xymatrix{ S_I\ar[r]^{i_I} & \tilde S_I \\
S_J\ar[u]^{j_{IJ}}\ar[r]^{i_J} & \tilde S_J\ar[u]^{p_{IJ}}}  
\end{equation*}
where $p_{IJ}:\tilde S_J\to\tilde S_I$ is the projection
and $j_{IJ}:S_J\hookrightarrow S_I$ is the open embedding so that $j_I\circ j_{IJ}=j_J$.
This gives the diagram of analytic spaces $(\tilde S_I)\in\Fun(\mathcal P(\mathbb N),\AnSp(\mathbb C))$
which gives the diagram of sites $\AnSp(\mathbb C)^2/(\tilde S_I)\in\Fun(\mathcal P(\mathbb N),\Cat)$.  
This gives in the same way the diagram of analytic spaces 
$(\tilde S_I)^{op}\in\Fun(\mathcal P(\mathbb N)^{op},\AnSp(\mathbb C))$
which gives the diagram of sites $\AnSp(\mathbb C)^2/(\tilde S_I)^{op}\in\Fun(\mathcal P(\mathbb N)^{op},\Cat)$.  

\begin{itemize}
\item Then $C_{fil}(\AnSp(\mathbb C)^{2,(sm)}/(\tilde S_I))$ is the category  
\begin{itemize}
\item whose objects $(G,F)=((G_I,F)_{I\subset\left[1,\cdots l\right]},u_{IJ})$,
with $(G_I,F)\in C_{fil}(\AnSp(\mathbb C)^{2,(sm)}/\tilde S_I)$,
and $u_{IJ}:(G_I,F)\to p_{IJ*}(G_J,F)$ for $I\subset J$, are morphisms
satisfying for $I\subset J\subset K$, $p_{IJ*}u_{JK}\circ u_{IJ}=u_{IK}$ in $C_{fil}(\AnSp(\mathbb C)^{2,(sm)}/\tilde S_I)$,
\item the morphisms $m:((G,F),u_{IJ})\to((H,F),v_{IJ})$ being (see section 2.1) a family of morphisms of complexes,  
\begin{equation*}
m=(m_I:(G_I,F)\to (H_I,F))_{I\subset\left[1,\cdots l\right]}
\end{equation*}
such that $v_{IJ}\circ m_I=p_{IJ*}m_J\circ u_{IJ}$ in $C_{fil}(\AnSp(\mathbb C)^{2,(sm)}/\tilde S_I)$.
\end{itemize}
\item Then $C_{fil}(\AnSp(\mathbb C)^{2,(sm)pr}/(\tilde S_I))$ is the category  
\begin{itemize}
\item whose objects $(G,F)=((G_I,F)_{I\subset\left[1,\cdots l\right]},u_{IJ})$,
with $(G_I,F)\in C_{fil}(\AnSp(\mathbb C)^{2,(sm)pr}/\tilde S_I)$,
and $u_{IJ}:(G_I,F)\to p_{IJ*}(G_J,F)$ for $I\subset J$, are morphisms
satisfying for $I\subset J\subset K$, $p_{IJ*}u_{JK}\circ u_{IJ}=u_{IK}$ in $C_{fil}(\AnSp(\mathbb C)^{2,(sm)pr}/\tilde S_I)$,
\item the morphisms $m:((G,F),u_{IJ})\to((H,F),v_{IJ})$ being (see section 2.1) a family of morphisms of complexes,  
\begin{equation*}
m=(m_I:(G_I,F)\to (H_I,F))_{I\subset\left[1,\cdots l\right]}
\end{equation*}
such that $v_{IJ}\circ m_I=p_{IJ*}m_J\circ u_{IJ}$ in $C_{fil}(\AnSp(\mathbb C)^{2,(sm)pr}/\tilde S_I)$.
\end{itemize}
\item Then $C_{fil}(\AnSp(\mathbb C)^{2,(sm)}/(\tilde S_I)^{op})$ 
is the category  
\begin{itemize}
\item whose objects $(G,F)=((G_I,F)_{I\subset\left[1,\cdots l\right]},u_{IJ})$,
with $(G_I,F)\in C_{fil}(\AnSp(\mathbb C)^{2,(sm)}/\tilde S_I)$,
and $u_{IJ}:(G_J,F)\to p_{IJ}^*(G_I,F)$ for $I\subset J$, are morphisms
satisfying for $I\subset J\subset K$, $p_{JK}^*u_{IJ}\circ u_{JK}=u_{IK}$ in $C_{fil}(\AnSp(\mathbb C)^{2,(sm)}/\tilde S_K)$,
\item the morphisms $m:((G,F),u_{IJ})\to((H,F),v_{IJ})$ being (see section 2.1) a family of morphisms of complexes,  
\begin{equation*}
m=(m_I:(G_I,F)\to (H_I,F))_{I\subset\left[1,\cdots l\right]}
\end{equation*}
such that $v_{IJ}\circ m_J=p_{IJ}^*m_I\circ u_{IJ}$ in $C_{fil}(\AnSp(\mathbb C)^{2,(sm)}/\tilde S_J)$.
\end{itemize}
\item Then $C_{fil}(\AnSp(\mathbb C)^{2,(sm)pr}/(\tilde S_I)^{op})$ 
is the category  
\begin{itemize}
\item whose objects $(G,F)=((G_I,F)_{I\subset\left[1,\cdots l\right]},u_{IJ})$,
with $(G_I,F)\in C_{fil}(\AnSp(\mathbb C)^{2,(sm)pr}/\tilde S_I)$,
and $u_{IJ}:(G_J,F)\to p_{IJ}^*(G_I,F)$ for $I\subset J$, are morphisms
satisfying for $I\subset J\subset K$, $p_{JK}^*u_{IJ}\circ u_{JK}=u_{IK}$ in $C_{fil}(\AnSp(\mathbb C)^{2,(sm)pr}/\tilde S_K)$,
\item the morphisms $m:((G,F),u_{IJ})\to((H,F),v_{IJ})$ being (see section 2.1) a family of morphisms of complexes,  
\begin{equation*}
m=(m_I:(G_I,F)\to (H_I,F))_{I\subset\left[1,\cdots l\right]}
\end{equation*}
such that $v_{IJ}\circ m_J=p_{IJ}^*m_I\circ u_{IJ}$ in $C_{fil}(\AnSp(\mathbb C)^{2,(sm)pr}/\tilde S_J)$.
\end{itemize}
\end{itemize}

We now define the usual topology on $\AnSp(\mathbb C)^2/S$.

\begin{defi}\label{tau12an}
Let $S\in\AnSp(\mathbb C)$. 
\begin{itemize}
\item[(i)]Denote by $\tau$ a topology on $\AnSp(\mathbb C)$, e.g. the usual topology. 
The $\tau$ covers in $\AnSp(\mathbb C)^2/S$ of $(X,Z)/S$ are the families of morphisms 
\begin{eqnarray*}
\left\{(c_i:(U_i,Z\times_X U_i)/S\to(X,Z)/S)_{i\in I}, \; 
\mbox{with} \; (c_i:U_i\to X)_{i\in I} \, \tau \, \mbox{cover of} \, X \, \mbox{in} \, \AnSp(\mathbb C)\right\}
\end{eqnarray*}
\item[(ii)]Denote by $\tau$ the usual or the etale topology on $\AnSp(\mathbb C)$. 
The $\tau$ covers in $\AnSp(\mathbb C)^{2,sm}/S$ of $(U,Z)/S$ are the families of morphisms 
\begin{eqnarray*}
\left\{(c_i:(U_i,Z\times_U U_i)/S\to(U,Z)/S)_{i\in I}, \; 
\mbox{with} \; (c_i:U_i\to U)_{i\in I} \, \tau \, \mbox{cover of} \, U \, \mbox{in} \, \AnSp(\mathbb C)\right\}
\end{eqnarray*}
\item[(iii)]Denote by $\tau$ the usual or the etale topology on $\AnSp(\mathbb C)$. 
The $\tau$ covers in $\AnSp(\mathbb C)^{2,(sm)pr}/S$ of $(Y\times S,Z)/S$ are the families of morphisms 
\begin{eqnarray*}
\left\{(c_i\times I_S:(U_i\times S,Z\times_{Y\times S} U_i\times S)/S\to(Y\times S,Z)/S)_{i\in I}, \; 
\mbox{with} \; (c_i:U_i\to Y)_{i\in I} \, \tau \, \mbox{cover of} \, Y \, \mbox{in} \, \AnSp(\mathbb C)\right\}
\end{eqnarray*}
\end{itemize}
\end{defi}

Let $S\in\AnSp(\mathbb C)$. Denote by $\tau$ the usual topology on $\AnSp(\mathbb C)$.
In particular, denoting 
$a_{\tau}:\PSh(\AnSp(\mathbb C)^{2,(sm)}/S)\to\Shv(\AnSp(\mathbb C)^{2,(sm)}/S)$ and
$a_{\tau}:\PSh(\AnSp(\mathbb C)^{2,(sm)pr}/S)\to\Shv(\AnSp(\mathbb C)^{2,(sm)pr}/S)$
the sheaftification functors,
\begin{itemize}
\item a morphism $\phi:F\to G$, with $F,G\in C(\AnSp(\mathbb C)^{2(sm)}/S)$,
is a $\tau$ local equivalence if $a_{\tau}H^n\phi:a_{\tau}H^nF\to a_{\tau}H^nG$ is an isomorphism,
a morphism $\phi:F\to G$, with $F,G\in C(\AnSp(\mathbb C)^{2,(sm)pr}/S)$,
is a $\tau$ local equivalence if $a_{\tau}H^n\phi:a_{\tau}H^nF\to a_{\tau}H^nG$ is an isomorphism,
\item $F^{\bullet}\in C(\AnSp(\mathbb C)^{2,(sm)}/S)$ is $\tau$ fibrant 
if for all $(U,Z)/S\in\AnSp(\mathbb C)^{2,(sm)}/S$ and all $\tau$ covers 
$(c_i:(U_i,Z\times_U U_i)/S\to(U,Z)/S)_{i\in I}$ of $(U,Z)/S$,
\begin{equation*}
F^{\bullet}(c_i):F^{\bullet}((U,Z)/S)\to\Tot(\oplus_{card I=\bullet} F^{\bullet}((U_I,Z\times_UU_I)/S))
\end{equation*}
is a quasi-isomorphism of complexes of abelian groups ;
$F^{\bullet}\in C(\AnSp(\mathbb C)^{2,(sm)pr}/S)$ is $\tau$ fibrant 
if for all $(Y\times S,Z)/S\in\AnSp(\mathbb C)^{2,smpr}/S$ and all $\tau$ covers 
$(c_i\times I_S:(U_i\times S,Z\times_{Y\times S}U_i\times S)/S\to(Y\times S,Z)/S)_{i\in I}$ of $(Y\times S,Z)/S$,
\begin{equation*}
F^{\bullet}(c_i\times I_S):F^{\bullet}((Y\times S,Z)/S)\to
\Tot(\oplus_{card I=\bullet} F^{\bullet}((U_I\times S,Z_I\times_YU_I)/S))
\end{equation*}
is a quasi-isomorphism of complexes of abelian groups,
\item a morphism $\phi:(G_1,F)\to (G_2,F)$, with $(G_1,F),(G_2,F)\in C_{fil}(\AnSp(\mathbb C)^{2,(sm)}/S)$, 
is a filtered $\tau$ local equivalence if for all $n,p\in\mathbb Z$, 
\begin{equation*}
a_{\tau}H^n\Gr_F^p\phi:a_{\tau}H^n\Gr_F^p(G_1,F)\to a_{\tau}H^n\Gr_F^p(G_2,F) 
\end{equation*}
is an isomorphism of sheaves on $\AnSp(\mathbb C)^{2,(sm)}/S$ ;
a morphism $\phi:(G_1,F)\to (G_2,F)$, with $(G_1,F),(G_2,F)\in C_{fil}(\AnSp(\mathbb C)^{2,(sm)pr}/S)$, 
is an filtered $\tau$ local equivalence if for all $n,p\in\mathbb Z$
\begin{equation*}
a_{\tau}H^n\Gr_F^p\phi:a_{\tau}H^n\Gr_F^p(G_1,F)\to a_{\tau}H^n\Gr_F^p(G_2,F) 
\end{equation*}
is an isomorphism of sheaves on  $\AnSp(\mathbb C)^{2,(sm)pr}/S$ ;
\item a morphism $\phi:(G_1,F)\to (G_2,F)$, with $(G_1,F),(G_2,F)\in C_{fil}(\AnSp(\mathbb C)^{2,(sm)}/S)$, 
is an $r$-filtered $\tau$ local equivalence if 
there exists $\phi_i:(C_i,F)\to(C_{i+1},F)$, $0\leq i\leq s$, 
with $(C_i,F)\in C_{fil}(\AnSp(\mathbb C)^{2,(sm)}/S)$, $(C_0,F)=(G_1,F)$ and $(C_s,F)=(G_2,F)$, such that
\begin{equation*}
\phi=\phi_s\circ\cdots\circ\phi_i\circ\cdots\circ\phi_0:(G_1,F)\to(G_2,F)
\end{equation*}
and $\phi_i:(C_i,F)\to(C_{i+1},F)$ either a filtered $\tau$ local equivalence
or an $r$-filtered homotopy equivalence,
a morphism $\phi:(G_1,F)\to (G_2,F)$, with $(G_1,F),(G_2,F)\in C_{fil}(\AnSp(\mathbb C)^{2,(sm)pr}/S)$, 
is an $r$-filtered $\tau$ local equivalence if
there exists $\phi_i:(C_i,F)\to(C_{i+1},F)$, $0\leq i\leq s$, 
with $(C_i,F)\in C_{fil}(\AnSp(\mathbb C)^{2,(sm)pr}/S)$, $(C_0,F)=(G_1,F)$ and $(C_s,F)=(G_2,F)$, such that
\begin{equation*}
\phi=\phi_s\circ\cdots\circ\phi_i\circ\cdots\circ\phi_0:(G_1,F)\to(G_2,F)
\end{equation*}
and $\phi_i:(C_i,F)\to(C_{i+1},F)$ either a filtered $\tau$ local equivalence
or an $r$-filtered homotopy equivalence ;
\item $(F^{\bullet},F)\in C_{fil}(\AnSp(\mathbb C)^{2,(sm)}/S)$ is filtered $\tau$ fibrant 
for all $(U,Z)/S\in\AnSp(\mathbb C)^{2,(sm)}/S$ and all $\tau$ covers 
$(c_i:(U_i,Z\times_UU_i)/S\to(U,Z)/S)_{i\in I}$ of $(U,Z)/S$,
\begin{eqnarray*}
H^n\Gr_F^p(F^{\bullet},F)(c_i):H^n\Gr_F^p(F^{\bullet},F)((U,Z)/S)\to \\
H^n\Gr_F^p(\Tot(\oplus_{card I=\bullet} (F^{\bullet},F)((U_I,Z\times_UU_I)/S)))
\end{eqnarray*}
is an isomorphism of abelian groups for all $n,p\in\mathbb Z$ ;
$(F^{\bullet},F)\in C_{fil}(\AnSp(\mathbb C)^{2,(sm)pr}/S)$ is filtered $\tau$ fibrant 
for all $(Y\times S,Z)/S\in\AnSp(\mathbb C)^{2,(sm)pr}/S$ and all $\tau$ covers 
$(c_i\times I_S:(U_i\times S,Z\times_{Y\times S}U_i\times S)/S\to(Y\times S,Z)/S)_{i\in I}$ of $(Y\times S,Z)/S$,
\begin{eqnarray*}
H^n\Gr_F^p(F^{\bullet},F)(c_i\times I_S):H^n\Gr_F^p(F^{\bullet},F)((Y\times S,Z)/S)\to \\
H^n\Gr_F^p(\Tot(\oplus_{card I=\bullet}(F^{\bullet},F)((U_I\times S,Z\times_YU_I)/S)))
\end{eqnarray*}
is an isomorphism of abelian groups for all $n,p\in\mathbb Z$.
\end{itemize}
Let $S_{\bullet}\in\Fun(\mathcal I,\AnSp(\mathbb C))$ with $\mathcal I\in\Cat$.
\begin{itemize}
\item A morphism $m:((G_I,F),u_{IJ})\to((H_I,F),v_{IJ})$ in $C_{fil}(\AnSp(\mathbb C)^{2,(sm)}/S_{\bullet})$
is an $r$-filtered Zariski, resp. etale local, equivalence,
if there exists $\phi_i:((C_{iI},F),u_{iIJ})\to((C_{(i+1)I},F),u_{(i+1)IJ})$, $0\leq i\leq s$, 
with $((C_{iI},F),u_{iIJ})\in C_{fil}(\AnSp(\mathbb C)^{2,(sm)}/S_{\bullet})$
$((C_{0I},F),u_{0IJ})=((G_I,F),u_{IJ})$ and $((C_{sI},F),u_{sIJ})=((H_I,F),v_{IJ})$, such that
\begin{equation*}
\phi=\phi_s\circ\cdots\circ\phi_i\circ\cdots\circ\phi_0:((G_I,F),u_{IJ})\to((H_I,F),v_{IJ})
\end{equation*}
and $\phi_i:((C_{iI},F),u_{iIJ})\to((C_{(i+1)I},F),u_{(i+1)IJ})$ either a filtered Zariski, resp. etale, local equivalence
or an $r$-filtered homotopy equivalence.  
\item A morphism $m:((G_I,F),u_{IJ})\to((H_I,F),v_{IJ})$ in $C_{fil}(\AnSp(\mathbb C)^{2,(sm)pr}/S_{\bullet})$
is an $r$-filtered Zariski, resp. etale local, equivalence,
if there exists $\phi_i:((C_{iI},F),u_{iIJ})\to((C_{(i+1)I},F),u_{(i+1)IJ})$, $0\leq i\leq s$, 
with $((C_{iI},F),u_{iIJ})\in C_{fil}(\AnSp(\mathbb C)^{2,(sm)pr}/S_{\bullet})$
$((C_{0I},F),u_{0IJ})=((G_I,F),u_{IJ})$ and $((C_{sI},F),u_{sIJ})=((H_I,F),v_{IJ})$, such that
\begin{equation*}
\phi=\phi_s\circ\cdots\circ\phi_i\circ\cdots\circ\phi_0:((G_I,F),u_{IJ})\to((H_I,F),v_{IJ})
\end{equation*}
and $\phi_i:((C_{iI},F),u_{iIJ})\to((C_{(i+1)I},F),u_{(i+1)IJ})$ either a filtered Zariski, resp. etale, local equivalence
or an $r$-filtered homotopy equivalence. 
\end{itemize}

Will now define the $\mathbb D^1$ local property on $\AnSp(\mathbb C)^2/S$.
Let $S\in\AnSp(\mathbb C)$. Denote for short $\AnSp(\mathbb C)^{2,(sm)}/S$ 
either the category $\AnSp(\mathbb C)^2/S$ or the category $\AnSp(\mathbb C)^{2,sm}/S$. Denote by
\begin{eqnarray*}
p_a:\AnSp(\mathbb C)^{2,(sm)}/S\to\AnSp(\mathbb C)^{2,(sm)}/S, \\ 
(X,Z)/S=((X,Z),h)\mapsto (X\times\mathbb D^1,Z\times\mathbb D^1)/S=((X\times\mathbb D^1,Z\times\mathbb D^1,h\circ p_X), \\ 
(g:(X,Z)/S\to (X',Z')/S)\mapsto 
((g\times I_{\mathbb D^1}):(X\times\mathbb D^1,Z\times\mathbb D^1)/S\to (X'\times\mathbb D^1,Z'\times\mathbb D^1)/S)
\end{eqnarray*}
the projection functor and again by $p_a:\AnSp(\mathbb C)^{2,(sm)}/S\to\AnSp(\mathbb C)^{2,(sm)}/S$
the corresponding morphism of site.
Denote for short $\AnSp(\mathbb C)^{2,(sm)pr}/S$ 
either the category $\AnSp(\mathbb C)^{2,pr}/S$ or the category $\AnSp(\mathbb C)^{2,smpr}/S$. Denote by
\begin{eqnarray*}
p_a:\AnSp(\mathbb C)^{2,(sm)pr}/S\to\AnSp(\mathbb C)^{2,(sm)pr}/S, \\ 
(Y\times S,Z)/S=((Y\times S,Z),p_S)\mapsto 
(Y\times S\times\mathbb D^1,Z\times\mathbb D^1)/S=((Y\times S\times\mathbb D^1,Z\times\mathbb D^1,p_S\circ p_{Y\times S}), \\ 
(g:(Y\times S,Z)/S\to (Y'\times S,Z')/S)\mapsto 
((g\times I_{\mathbb D^1}):(Y\times S\times\mathbb D^1,Z\times\mathbb D^1)/S\to 
(Y'\times S\times\mathbb D^1,Z'\times\mathbb D^1)/S)
\end{eqnarray*}
the projection functor and again by $p_a:\AnSp(\mathbb C)^{2,(sm)pr}/S\to\AnSp(\mathbb C)^{2,(sm)pr}/S$
the corresponding morphism of site.

\begin{defi}\label{d1loc12def}
Let $S\in\AnSp(\mathbb C)$. 
\begin{itemize}
\item[(i0)] A complex  $F\in C(\AnSp(\mathbb C)^{2,(sm)}/S)$, is said to be $\mathbb D^1$ homotopic 
if $\ad(p_a^*,p_{a*})(F):F\to F$ is an homotopy equivalence.
\item[(i)] A complex  $F\in C(\AnSp(\mathbb C)^{2,(sm)}/S)$, is said to be $\mathbb D^1$ invariant 
if for all $(X,Z)/S\in\AnSp(\mathbb C)^{2,(sm)}/S$ 
\begin{equation*}
F(p_X):F((X,Z)/S)\to F((X\times\mathbb D^1,(Z\times\mathbb D^1))/S) 
\end{equation*}
is a quasi-isomorphism, where $p_X:(X\times\mathbb D^1,(Z\times\mathbb D^1))\to (X,Z)$ is the projection.
\item[(i0)'] Similarly, a complex  $F\in C(\AnSp(\mathbb C)^{2,(sm)pr}/S)$, is said to be $\mathbb D^1$ homotopic 
if $\ad(p_a^*,p_{a*})(F):F\to F$ is an homotopy equivalence.
\item[(i)']Similarly, a complex $F\in C(\AnSp(\mathbb C)^{2,(sm)pr}/S)$ is said to be $\mathbb D^1$ invariant 
if for all $(Y\times S,Z)/S\in\AnSp(\mathbb C)^{2,(sm)pr}/S$ 
\begin{equation*}
F(p_{Y\times S}):F((Y\times S,Z)/S)\to F((Y\times S\times\mathbb D^1,(Z\times\mathbb D^1))/S) 
\end{equation*}
is a quasi-isomorphism
\item[(ii)] A complex $F\in C(\AnSp(\mathbb C)^{2,(sm)}/S)$ is said to be $\mathbb D^1$ local 
for the $\tau$ topology induced on $\AnSp(\mathbb C)^2/S$, 
if for an (hence every) $\tau$ local equivalence $k:F\to G$ with $k$ injective and $G\in C(\AnSp(\mathbb C)^{2,(sm)}/S)$ $\tau$ fibrant,
(e.g. $k:F\to E_{\tau}(F)$), $G$ is $\mathbb D^1$ invariant.
\item[(ii)'] Similarly, a complex $F\in C(\AnSp(\mathbb C)^{2,(sm)pr}/S)$ is said to be $\mathbb D^1$ local 
for the $\tau$ topology induced on $\AnSp(\mathbb C)^2/S$, 
if for an (hence every) $\tau$ local equivalence $k:F\to G$ with $k$ injective 
and $G\in C(\AnSp(\mathbb C)^{2,(sm)pr}/S)$ $\tau$ fibrant, e.g. $k:F\to E_{\tau}(F)$, $G$ is $\mathbb D^1$ invariant.
\item[(iii)] A morphism $m:F\to G$ with $F,G\in C(\AnSp(\mathbb C)^{2,(sm)}/S)$ is said to an $(\mathbb D^1,usu)$ local equivalence 
if for all $H\in C(\AnSp(\mathbb C)^{2,(sm)}/S)$ which is $\mathbb D^1$ local for the usual topology
\begin{equation*}
\Hom(L(m),E_{usu}(H)):\Hom(L(G),E_{usu}(H))\to\Hom(L(F),E_{usu}(H)) 
\end{equation*}
is a quasi-isomorphism.
\item[(iii)']Similarly, a morphism $m:F\to G$ with $F,G\in C(\AnSp(\mathbb C)^{2,(sm)pr}/S)$ 
is said to be an $(\mathbb D^1,usu)$ local equivalence 
if for all $H\in C(\AnSp(\mathbb C)^{2,(sm)pr}/S)$ which is $\mathbb D^1$ local for the usual topology
\begin{equation*}
\Hom(L(m),E_{usu}(H)):\Hom(L(G),E_{usu}(H))\to\Hom(L(F),E_{usu}(H)) 
\end{equation*}
is a quasi-isomorphism.
\end{itemize}
\end{defi}

\begin{prop}\label{d1loceqprop12}
\begin{itemize}
\item[(i)]A morphism $m:F\to G$ with $F,G\in C(\AnSp(\mathbb C)^{2,(sm)}/S)$ is an $(\mathbb D^1,usu)$ local equivalence
if and only if there exists 
\begin{eqnarray*}
\left\{(X_{1,\alpha},Z_{1,\alpha})/S,\alpha\in\Lambda_1\right\},\ldots,
\left\{(X_{r,\alpha},Z_{r,\alpha})/S,\alpha\in\Lambda_r\right\}\subset\AnSp(\mathbb C)^{2,(sm)}/S
\end{eqnarray*}
such that we have in $\Ho_{et}(C(\AnSp(\mathbb C)^{2,(sm)}/S))$
\begin{eqnarray*}
\Cone(m)\xrightarrow{\sim}\Cone(\oplus_{\alpha\in\Lambda_1}
\Cone(\mathbb Z((X_{1,\alpha}\times\mathbb D^1,Z_{1,\alpha}\times\mathbb D^1)/S)\to\mathbb Z((X_{1,\alpha},Z_{1,\alpha})/S)) \\
\to\cdots\to\oplus_{\alpha\in\Lambda_1}
\Cone(\mathbb Z((X_{1,\alpha}\times\mathbb D^1,Z_{1,\alpha}\times\mathbb D^1)/S)\to\mathbb Z((X_{1,\alpha},Z_{1,\alpha})/S)))
\end{eqnarray*}
\item[(ii)]A morphism $m:F\to G$ with $F,G\in C(\AnSp(\mathbb C)^{2,(sm)pr}/S)$ is an $(\mathbb D^1,usu)$ local equivalence
if and only if there exists 
\begin{eqnarray*}
\left\{(Y_{1,\alpha}\times S,Z_{1,\alpha})/S,\alpha\in\Lambda_1\right\},\ldots,
\left\{(Y_{r,\alpha}\times S,Z_{r,\alpha})/S,\alpha\in\Lambda_r\right\}\subset\AnSp(\mathbb C)^{2,(sm)}/S
\end{eqnarray*}
such that we have in $\Ho_{et}(C(\AnSp(\mathbb C)^{2,(sm)pr}/S))$
\begin{eqnarray*}
\Cone(m)\xrightarrow{\sim}\Cone(\oplus_{\alpha\in\Lambda_1}
\Cone(\mathbb Z((Y_{1,\alpha}\times S\times\mathbb D^1,Z_{1,\alpha}\times\mathbb D^1)/S)\to
\mathbb Z((Y_{1,\alpha}\times S,Z_{1,\alpha})/S)) \\
\to\cdots\to\oplus_{\alpha\in\Lambda_1}
\Cone(\mathbb Z((Y_{1,\alpha}\times S\times\mathbb D^1,Z_{1,\alpha}\times\mathbb D^1)/S)\to
\mathbb Z((Y_{1,\alpha}\times S,Z_{1,\alpha})/S)))
\end{eqnarray*}
\end{itemize}
\end{prop}

\begin{proof}
Standard.
\end{proof}

\begin{defiprop}
Let $S\in\AnSp(\mathbb C)$.
\begin{itemize}
\item[(i)]With the weak equivalence the $(\mathbb D^1,et)$ equivalence and 
the fibration the epimorphism with $\mathbb D^1$ local and etale fibrant kernels gives
a model structure on  $C(\AnSp(\mathbb C)^{2,(sm)}/S)$ : the left bousfield localization
of the projective model structure of $C(\AnSp(\mathbb C)^{2,(sm)}/S)$. 
We call it the projective $(\mathbb D^1,et)$ model structure.
\item[(ii)]With the weak equivalence the $(\mathbb D^1,et)$ equivalence and 
the fibration the epimorphism with $\mathbb D^1$ local and etale fibrant kernels gives
a model structure on  $C(\AnSp(\mathbb C)^{2,(sm)pr}/S)$ : the left bousfield localization
of the projective model structure of $C(\AnSp(\mathbb C)^{2,(sm)pr}/S)$. 
We call it the projective $(\mathbb D^1,et)$ model structure.
\end{itemize}
\end{defiprop}

\begin{proof}
Similar to the proof of proposition \ref{projmodstr}.
\end{proof}

We have, similarly to the case of single varieties the following :

\begin{prop}\label{g12an}
Let $g:T\to S$ a morphism with $T,S\in\AnSp(\mathbb C)$.
\begin{itemize}
\item[(i)] The adjonction $(g^*,g_*):C(\AnSp(\mathbb C)^{2,(sm)}/S)\leftrightarrows C(\AnSp(\mathbb C)^{2,(sm)}/T)$
is a Quillen adjonction for the $(\mathbb D^1,usu)$ model structure.
\item[(i)'] The functor $g^*:C(\AnSp(\mathbb C)^{2,(sm)}/S)\to C(\AnSp(\mathbb C)^{2,(sm)}/T)$
sends quasi-isomorphism to quasi-isomorphism and equivalence usu local to equivalence usu local,
sends $(\mathbb D^1,usu)$ local equivalence to $(\mathbb D^1,usu)$ local equivalence.
\item[(ii)] The adjonction $(g^*,g_*):C(\AnSp(\mathbb C)^{2,(sm)pr}/S)\leftrightarrows C(\AnSp(\mathbb C)^{2,(sm)pr}/T)$
is a Quillen adjonction for the $(\mathbb D^1,usu)$ model structure.
\item[(ii)'] The functor $g^*:C(\AnSp(\mathbb C)^{2,(sm)pr}/S)\to C(\AnSp(\mathbb C)^{2,(sm)pr}/T)$
sends quasi-isomorphism to quasi-isomorphism and equivalence usu local to equivalence usu local,
sends $(\mathbb D^1,usu)$ local equivalence to $(\mathbb D^1,usu)$ local equivalence.
\end{itemize}
\end{prop}

\begin{proof}
Similar to the proof of proposition \ref{g12}.
\end{proof}

\begin{prop}\label{rho12an}
Let $S\in\Var(\mathbb C)$. 
\begin{itemize}
\item[(i)] The adjonction $(\rho_S^*,\rho_{S*}):C(\AnSp(\mathbb C)^{2,sm}/S)\leftrightarrows C(\AnSp(\mathbb C)^2/S)$
is a Quillen adjonction for the $(\mathbb A^1,et)$ projective model structure.
\item[(i)']The functor $\rho_{S*}:C(\AnSp(\mathbb C)^2/S)\to C(\AnSp(\mathbb C)^{2,sm}/S)$
sends quasi-isomorphism to quasi-isomorphism, sends equivalence usu local to equivalence usu local,
sends $(\mathbb D^1,usu)$ local equivalence to $(\mathbb D^1,usu)$ local equivalence.
\item[(ii)] The adjonction $(\rho_S^*,\rho_{S*}):C(\AnSp(\mathbb C)^{2,smpr}/S)\leftrightarrows C(\AnSp(\mathbb C)^{2,pr}/S)$
is a Quillen adjonction for the $(\mathbb A^1,et)$ projective model structure.
\item[(ii)']The functor $\rho_{S*}:C(\AnSp(\mathbb C)^{2,pr}/S)\to C(\AnSp(\mathbb C)^{2,smpr}/S)$
sends quasi-isomorphism to quasi-isomorphism, sends equivalence usu local to equivalence usu local,
sends $(\mathbb D^1,usu)$ local equivalence to $(\mathbb D^1,usu)$ local equivalence.
\end{itemize}
\end{prop}

\begin{proof}
Similar to the proof of proposition \ref{rho1}.
\end{proof}

We also have

\begin{prop}\label{Gra1an}
Let $S\in\AnSp(\mathbb C)$. 
\begin{itemize}
\item[(i)] The adjonction $(\Gr_S^{12*},\Gr_{S*}^{12}):C(\AnSp(\mathbb C)/S)\leftrightarrows C(\AnSp(\mathbb C)^{2,pr}/S)$
is a Quillen adjonction for the $(\mathbb D^1,usu)$ projective model structure.
\item[(ii)] The adjonction 
$(\Gr_S^{12*}\Gr_{S*}^{12}:C(\AnSp(\mathbb C)^{sm}/S)\leftrightarrows C(\AnSp(\mathbb C)^{2,smpr}/S)$
is a Quillen adjonction for the $(\mathbb D^1,usu)$ projective model structure.
\end{itemize}
\end{prop}

\begin{proof}
Immediate from definition.
\end{proof}

Let $S\in\AnSp(\mathbb C)$. Let $S=\cup_{i=1}^l S_i$ an open affine cover and denote by $S_I=\cap_{i\in I} S_i$.
Let $i_i:S_i\hookrightarrow\tilde S_i$ closed embeddings, with $\tilde S_i\in\AnSp(\mathbb C)$.
\begin{itemize}
\item For $(G_I,K_{IJ})\in C(\AnSp(\mathbb C)^{2,(sm)}/(\tilde S_I)^{op})$ and 
$(H_I,T_{IJ})\in C(\AnSp(\mathbb C)^{2,(sm)}/(\tilde S_I))$, we denote
\begin{eqnarray*}
\mathcal Hom((G_I,K_{IJ}),(H_I,T_{IJ})):=(\mathcal Hom(G_I,H_I),u_{IJ}((G_I,K_{IJ}),(H_I,T_{IJ})))
\in C(\AnSp(\mathbb C)^{2,(sm)}/(\tilde S_I))
\end{eqnarray*}
with
\begin{eqnarray*}
u_{IJ}((G_I,K_{IJ})(H_I,T_{IJ})):\mathcal Hom(G_I,H_I) \\
\xrightarrow{\ad(p_{IJ}^*,p_{IJ*})(-)}p_{IJ*}p_{IJ}^*\mathcal Hom(G_I,H_I)
\xrightarrow{T(p_{IJ},hom)(-,-)}p_{IJ*}\mathcal Hom(p_{IJ}^*G_I,p_{IJ}^*H_I) \\
\xrightarrow{\mathcal Hom(p_{IJ}^*G_I,T_{IJ})}p_{IJ*}\mathcal Hom(p_{IJ}^*G_I,H_J)
\xrightarrow{\mathcal Hom(K_{IJ},H_J)}p_{IJ*}\mathcal Hom(G_J,H_J).
\end{eqnarray*}
This gives in particular the functor
\begin{eqnarray*}
C(\AnSp(\mathbb C)^{2,(sm)}/(\tilde S_I))\to C(\AnSp(\mathbb C)^{2,(sm)}/(\tilde S_I)^{op}),
(H_I,T_{IJ})\mapsto(H_I,T_{IJ}).
\end{eqnarray*}
\item For $(G_I,K_{IJ})\in C(\AnSp(\mathbb C)^{2,(sm)pr}/(\tilde S_I)^{op})$ and 
$(H_I,T_{IJ})\in C(\AnSp(\mathbb C)^{2,(sm)pr}/(\tilde S_I))$, we denote
\begin{eqnarray*}
\mathcal Hom((G_I,K_{IJ}),(H_I,T_{IJ})):=(\mathcal Hom(G_I,H_I),u_{IJ}((G_I,K_{IJ}),(H_I,T_{IJ})))
\in C(\AnSp(\mathbb C)^{2,(sm)pr}/(\tilde S_I))
\end{eqnarray*}
with
\begin{eqnarray*}
u_{IJ}((G_I,K_{IJ})(H_I,T_{IJ})):\mathcal Hom(G_I,H_I) \\
\xrightarrow{\ad(p_{IJ}^*,p_{IJ*})(-)}p_{IJ*}p_{IJ}^*\mathcal Hom(G_I,H_I)
\xrightarrow{T(p_{IJ},hom)(-,-)}p_{IJ*}\mathcal Hom(p_{IJ}^*G_I,p_{IJ}^*H_I) \\
\xrightarrow{\mathcal Hom(p_{IJ}^*G_I,T_{IJ})}p_{IJ*}\mathcal Hom(p_{IJ}^*G_I,H_J)
\xrightarrow{\mathcal Hom(K_{IJ},H_J)}p_{IJ*}\mathcal Hom(G_J,H_J).
\end{eqnarray*}
This gives in particular the functor
\begin{eqnarray*}
C(\AnSp(\mathbb C)^{2,(sm)pr}/(\tilde S_I))\to C(\AnSp(\mathbb C)^{2,(sm)pr}/(\tilde S_I)^{op}),
(H_I,T_{IJ})\mapsto(H_I,T_{IJ}).
\end{eqnarray*}
\end{itemize}
Let $S\in\AnSp(\mathbb C)$. Let $S=\cup_{i=1}^l S_i$ an open affinoid cover and denote by $S_I=\cap_{i\in I} S_i$.
Let $i_i:S_i\hookrightarrow\tilde S_i$ closed embeddings, with $\tilde S_i\in\AnSp(\mathbb C)$.
We have the functors
\begin{eqnarray*}
p_a:\AnSp(\mathbb C)^{2,(sm)}/(\tilde S_I)\to\AnSp(\mathbb C)^{2,(sm)}/(\tilde S_I), \\ 
((X,Z)/\tilde S_I,u_{IJ})\mapsto ((X\times\mathbb D^1,Z\times\mathbb D^1)/\tilde S_I,u_{IJ}\times I), \\ 
(g:((X,Z)/\tilde S_I,u_{IJ})\to ((X',Z')/\tilde S_I,u_{IJ}))\mapsto 
((g\times I_{\mathbb D^1}):((X\times\mathbb D^1,Z\times\mathbb D^1)/\tilde S_I,u_{IJ}\times I)\to 
((X'\times\mathbb D^1,Z'\times\mathbb D^1)/\tilde S_I,u_{IJ}\times I))
\end{eqnarray*}
the projection functor and again by $p_a:\AnSp(\mathbb C)^{2,(sm)}/(\tilde S_I)\to\Var(\mathbb C)^{2,(sm)}/(\tilde S_I)$
the corresponding morphism of site, and
\begin{eqnarray*}
p_a:\AnSp(\mathbb C)^{2,(sm)pr}/(\tilde S_I)\to\AnSp(\mathbb C)^{2,(sm)pr}/(\tilde S_I), \\ 
((Y\times\tilde S_I,Z)/\tilde S_I,u_{IJ})\mapsto 
((Y\times\tilde S_I\times\mathbb D^1,Z\times\mathbb D^1)/\tilde S_I,u_{IJ}\times I), \\ 
(g:((Y\times\tilde S_I,Z)/\tilde S_I,u_{IJ})\to ((Y'\times\tilde S_I,Z')/\tilde S_I,u_{IJ}))\mapsto 
((g\times I_{\mathbb D^1}):((Y\times\tilde S_I\times\mathbb D^1,Z\times\mathbb D^1)/\tilde S_I,u_{IJ}\times I), \\ 
((Y'\times\tilde S_I\times\mathbb D^1,Z'\times\mathbb D^1)/\tilde S_I,u_{IJ}\times I)), 
\end{eqnarray*}
the projection functor and again by $p_a:\AnSp(\mathbb C)^{2,(sm)pr}/(\tilde S_I)\to\AnSp(\mathbb C)^{2,(sm)pr}/(\tilde S_I)$
the corresponding morphism of site.
These functors induces also morphism of sites
$p_a:\AnSp(\mathbb C)^{2,(sm)}/(\tilde S_I)^{op}\to\AnSp(\mathbb C)^{2,(sm)}/(\tilde S_I)^{op}$
and $p_a:\AnSp(\mathbb C)^{2,(sm)pr}/(\tilde S_I)^{op}\to\AnSp(\mathbb C)^{2,(sm)pr}/(\tilde S_I)^{op}$.

\begin{defi}\label{d1loc12defIJ}
Let $S\in\AnSp(\mathbb C)$. Let $S=\cup_{i=1}^l S_i$ an open affine cover and denote by $S_I=\cap_{i\in I} S_i$.
Let $i_i:S_i\hookrightarrow\tilde S_i$ closed embeddings, with $\tilde S_i\in\AnSp(\mathbb C)$.
\begin{itemize}
\item[(i0)]A complex $(F_I,u_{IJ})\in C(\AnSp(\mathbb C)^{2,(sm)}/(\tilde S_I))$ is said to be $\mathbb D^1$ homotopic if 
$\ad(p_a^*,p_{a*})((F_I,u_{IJ})):(F_I,u_{IJ})\to p_{a*}p_a^*(F_I,u_{IJ})$ is an homotopy equivalence.
\item[(i0)']A complex $(F_I,u_{IJ})\in C(\AnSp(\mathbb C)^{2,(sm)pr}/(\tilde S_I))$ is said to be $\mathbb D^1$ homotopic if 
$\ad(p_a^*,p_{a*})((F_I,u_{IJ})):(F_I,u_{IJ})\to p_{a*}p_a^*(F_I,u_{IJ})$ is an homotopy equivalence.
\item[(i)] A complex  $(F_I,u_{IJ})\in C(\AnSp(\mathbb C)^{2,(sm)}/(\tilde S_I))$ is said to be $\mathbb D^1$ invariant 
if for all $((X_I,Z_I)/\tilde S_I,s_{IJ})\in\AnSp(\mathbb C)^{2,(sm)}/(\tilde S_I)$ 
\begin{equation*}
(F_I(p_{X_I})):(F_I((X_I,Z_I)/\tilde S_I),F_J(s_{IJ})\circ u_{IJ}(-)\to 
(F_I((X_I\times\mathbb D^1,(Z_I\times\mathbb D^1))/\tilde S_I),F_J(s_{IJ}\times I)\circ u_{IJ}(-)) 
\end{equation*}
is a quasi-isomorphism, where $p_{X_I}:(X_I\times\mathbb D^1,(Z_I\times\mathbb D^1))\to (X_I,Z_I)$ are the projection,
and $s_{IJ}:(X_I\times\tilde S_{J\backslash I},Z_I)/\tilde S_J\to(X_J,Z_J)\tilde S_J$.
\item[(i)'] A complex  $(G_I,u_{IJ})\in C(\AnSp(\mathbb C)^{2,(sm)pr}/(\tilde S_I))$ is said to be $\mathbb D^1$ invariant 
if for all $((Y\times\tilde S_I,Z_I)/\tilde S_I,s_{IJ})\in\AnSp(\mathbb C)^{2,(sm)pr}/(\tilde S_I)$ 
\begin{eqnarray*}
(G_I(p_{Y\times\tilde S_I})):
(G_I((Y\times\tilde S_I,Z_I)/\tilde S_I),G_J(s_{IJ})\circ u_{IJ}(-)\to \\
(G_I((Y\times\tilde S_I\times\mathbb D^1,(Z_I\times\mathbb D^1))/\tilde S_I),G_J(s_{IJ}\times I)\circ u_{IJ}(-)) 
\end{eqnarray*}
is a quasi-isomorphism. 
\item[(ii)]Let $\tau$ a topology on $\AnSp(\mathbb C)$. 
A complex $F=(F_I,u_{IJ})\in C(\AnSp(\mathbb C)^{2,(sm)}/(\tilde S_I))$ is said to be $\mathbb D^1$ local 
for the $\tau$ topology induced on $\AnSp(\mathbb C)^2/(\tilde S_I)$, 
if for an (hence every) $\tau$ local equivalence $k:F\to G$ with $k$ injective and 
$G=(G_I,v_{IJ})\in C(\AnSp(\mathbb C)^{2,(sm)}/(\tilde S_I))$ $\tau$ fibrant,
e.g. $k:(F_I,u_{IJ})\to (E_{\tau}(F_I),E(u_{IJ}))$, $G$ is $\mathbb D^1$ invariant.
\item[(ii)']Let $\tau$ a topology on $\AnSp(\mathbb C)$. 
A complex $F=(F_I,u_{IJ})\in C(\AnSp(\mathbb C)^{2,(sm)pr}/(\tilde S_I))$ is said to be $\mathbb D^1$ local 
for the $\tau$ topology induced on $\AnSp(\mathbb C)^2/(\tilde S_I)$, 
if for an (hence every) $\tau$ local equivalence $k:F\to G$ with $k$ injective and 
$G=(G_I,u_{IJ})\in C(\AnSp(\mathbb C)^{2,(sm)pr}/(\tilde S_I))$ $\tau$ fibrant,
e.g. $k:(F_I,u_{IJ})\to (E_{\tau}(F_I),E(u_{IJ}))$, $G$ is $\mathbb D^1$ invariant.
\item[(iii)] A morphism $m=(m_I):(F_I,u_{IJ})\to (G_I,v_{IJ})$ with 
$(F_I,u_{IJ}),(G_I,v_{IJ})\in C(\AnSp(\mathbb C)^{2,(sm)}/(\tilde S_I))$ 
is said to an $(\mathbb D^1,usu)$ local equivalence 
if for all $H=(H_I,w_{IJ})\in C(\AnSp(\mathbb C)^{2,(sm)}/(\tilde S_I))$ which is $\mathbb D^1$ local for the usual topology
\begin{eqnarray*}
(\Hom(L(m_I),E_{et}(H_I))):\Hom(L(G_I,v_{IJ}),E_{et}(H_I,w_{IJ}))\to\Hom(L(F_I,u_{IJ}),E_{et}(H_I,w_{IJ})) 
\end{eqnarray*}
is a quasi-isomorphism.
Obviously, if a morphism $m=(m_I):(F_I,u_{IJ})\to (G_I,v_{IJ})$ with 
$(F_I,u_{IJ}),(G_I,u_{IJ})\in C(\AnSp(\mathbb C)^{2,(sm)}/(\tilde S_I))$ 
is an $(\mathbb D^1,usu)$ local equivalence, 
then all the $m_I:F_I\to G_I$ are $(\mathbb D^1,usu)$ local equivalence.
\item[(iii)'] A morphism $m=(m_I):(F_I,u_{IJ})\to (G_I,v_{IJ})$ with 
$(F_I,u_{IJ}),(G_I,v_{IJ})\in C(\AnSp(\mathbb C)^{2,(sm)pr}/(\tilde S_I))$ 
is said to an $(\mathbb D^1,usu)$ local equivalence 
if for all $(H_I,w_{IJ})\in C(\AnSp(\mathbb C)^{2,(sm)pr}/(\tilde S_I))$ which is $\mathbb D^1$ local for the usual topology
\begin{eqnarray*}
(\Hom(L(m_I),E_{et}(H_I))):\Hom(L(G_I,v_{IJ}),E_{et}(H_I,w_{IJ}))\to\Hom(L(F_I,u_{IJ}),E_{et}(H_I,w_{IJ})) 
\end{eqnarray*}
is a quasi-isomorphism.
Obviously, if a morphism $m=(m_I):(F_I,u_{IJ})\to (G_I,v_{IJ})$ with 
$(F_I,u_{IJ}),(G_I,u_{IJ})\in C(\AnSp(\mathbb C)^{2,(sm)pr}/(\tilde S_I))$ 
is an $(\mathbb D^1,usu)$ local equivalence, 
then all the $m_I:F_I\to G_I$ are $(\mathbb D^1,usu)$ local equivalence.
\item[(iv)] A morphism $m=(m_I):(F_I,u_{IJ})\to (G_I,v_{IJ})$ with 
$(F_I,u_{IJ}),(G_I,v_{IJ})\in C(\AnSp(\mathbb C)^{2,(sm)}/(\tilde S_I)^{op})$ 
is said to an $(\mathbb D^1,usu)$ local equivalence 
if for all $H=(H_I,w_{IJ})\in C(\AnSp(\mathbb C)^{2,(sm)}/(\tilde S_I))$ which is $\mathbb D^1$ local for the usual topology
\begin{eqnarray*}
(\Hom(L(m_I),E_{et}(H_I))):\Hom(L(G_I,v_{IJ}),E_{et}(H_I,w_{IJ}))\to\Hom(L(F_I,u_{IJ}),E_{et}(H_I,w_{IJ})) 
\end{eqnarray*}
is a quasi-isomorphism.
Obviously, if a morphism $m=(m_I):(F_I,u_{IJ})\to (G_I,v_{IJ})$ with 
$(F_I,u_{IJ}),(G_I,u_{IJ})\in C(\AnSp(\mathbb C)^{2,(sm)}/(\tilde S_I)^{op})$ 
is an $(\mathbb D^1,usu)$ local equivalence, 
then all the $m_I:F_I\to G_I$ are $(\mathbb D^1,usu)$ local equivalence.
\item[(iv)'] A morphism $m=(m_I):(F_I,u_{IJ})\to (G_I,v_{IJ})$ with 
$(F_I,u_{IJ}),(G_I,v_{IJ})\in C(\AnSp(\mathbb C)^{2,(sm)pr}/(\tilde S_I)^{op})$ 
is said to an $(\mathbb D^1,usu)$ local equivalence 
if for all $(H_I,w_{IJ})\in C(\AnSp(\mathbb C)^{2,(sm)pr}/(\tilde S_I))$ which is $\mathbb D^1$ local for the usual topology
\begin{eqnarray*}
(\Hom(L(m_I),E_{et}(H_I))):\Hom(L(G_I,v_{IJ}),E_{et}(H_I,w_{IJ}))\to\Hom(L(F_I,u_{IJ}),E_{et}(H_I,w_{IJ})) 
\end{eqnarray*}
is a quasi-isomorphism.
Obviously, if a morphism $m=(m_I):(F_I,u_{IJ})\to (G_I,v_{IJ})$ with 
$(F_I,u_{IJ}),(G_I,u_{IJ})\in C(\AnSp(\mathbb C)^{2,(sm)pr}/(\tilde S_I)^{op})$ 
is an $(\mathbb D^1,usu)$ local equivalence, 
then all the $m_I:F_I\to G_I$ are $(\mathbb D^1,usu)$ local equivalence.
\end{itemize}
\end{defi}

\begin{prop}\label{d1loceqprop12IJ}
Let $S\in\AnSp(\mathbb C)$. Let $S=\cup_{i=1}^l S_i$ an open affine cover and denote by $S_I=\cap_{i\in I} S_i$.
Let $i_i:S_i\hookrightarrow\tilde S_i$ closed embeddings, with $\tilde S_i\in\AnSp(\mathbb C)$.
\begin{itemize}
\item[(i)]A morphism $m:F\to G$ with $F,G\in C(\AnSp(\mathbb C)^{2,(sm)}/(\tilde S_I)^{op})$ 
is an $(\mathbb D^1,usu)$ local equivalence if and only if there exists 
\begin{eqnarray*}
\left\{((X_{1,\alpha,I},Z_{1,\alpha,I})/\tilde S_I,u^1_{IJ}),\alpha\in\Lambda_1\right\},\ldots,
\left\{((X_{r,\alpha,I},Z_{r,\alpha,I})/\tilde S_I,u^r_{IJ}),\alpha\in\Lambda_r\right\}
\subset\AnSp(\mathbb C)^{2,(sm)}/(\tilde S_I)^{op}
\end{eqnarray*}
with 
\begin{equation*}
u^l_{IJ}:(X_{l,\alpha,J},Z_{l,\alpha,J})/\tilde S_J\to 
(X_{l,\alpha,I}\times\tilde S_{J\backslash I},Z_{l,\alpha,I}\times\tilde S_{J\backslash I})/\tilde S_J
\end{equation*}
such that we have in $\Ho_{et}(C(\AnSp(\mathbb C)^{2,(sm)}/(\tilde S_I)^{op}))$
\begin{eqnarray*}
\Cone(m)\xrightarrow{\sim}\Cone( \\ \oplus_{\alpha\in\Lambda_1}
\Cone((\mathbb Z((X_{1,\alpha,I}\times\mathbb D^1,Z_{1,\alpha,I}\times\mathbb D^1)/\tilde S_I),\mathbb Z(u_{IJ}^1\times I))
\to(\mathbb Z((X_{1,\alpha,I},Z_{1,\alpha,I})/\tilde S_I),\mathbb Z(u_{IJ}^1))) \\
\to\cdots\to \\ \oplus_{\alpha\in\Lambda_r}
\Cone((\mathbb Z((X_{r,\alpha,I}\times\mathbb D^1,Z_{r,\alpha,I}\times\mathbb D^1)/\tilde S_I),\mathbb Z(u_{IJ}^r\times I))
\to(\mathbb Z((X_{r,\alpha,I},Z_{r,\alpha,I})/\tilde S_I),\mathbb Z(u^r_{IJ}))))
\end{eqnarray*}
\item[(i)']A morphism $m:F\to G$ with $F,G\in C(\AnSp(\mathbb C)^{2,(sm)pr}/(\tilde S_I)^{op})$ 
is an $(\mathbb D^1,usu)$ local equivalence if and only if there exists 
\begin{eqnarray*}
\left\{((Y_{1,\alpha,I}\times\tilde S_I,Z_{1,\alpha,I})/\tilde S_I,u^1_{IJ}),\alpha\in\Lambda_1\right\},\ldots,
\left\{((Y_{r,\alpha,I}\times\tilde S_I,Z_{r,\alpha,I})/\tilde S_I,u^r_{IJ}),\alpha\in\Lambda_r\right\} \\
\subset\AnSp(\mathbb C)^{2,(sm)pr}/(\tilde S_I)
\end{eqnarray*}
with 
\begin{equation*}
u^l_{IJ}:(Y_{l,\alpha,J}\times\tilde S_J,Z_{l,\alpha,J})/\tilde S_J\to 
(Y_{l,\alpha,I}\times\tilde S_J,Z_{l,\alpha,I}\times\tilde S_{J\backslash I})/\tilde S_J
\end{equation*}
such that we have in $\Ho_{et}(C(\AnSp(\mathbb C)^{2,(sm)}/(\tilde S_I)^{op}))$
\begin{eqnarray*}
\Cone(m)\xrightarrow{\sim}\Cone(\oplus_{\alpha\in\Lambda_1} \\
\Cone((\mathbb Z((Y_{1,\alpha,I}\times\mathbb D^1\times\tilde S_I,Z_{1,\alpha,I}\times\mathbb D^1)/\tilde S_I),
\mathbb Z(u_{IJ}^1\times I))
\to(\mathbb Z((Y_{1,\alpha,I}\times S,Z_{1,\alpha,I})/\tilde S_I),\mathbb Z(u_{IJ}))) \\
\to\cdots\to\oplus_{\alpha\in\Lambda_r} \\
\Cone((\mathbb Z((Y_{r,\alpha,I}\times\mathbb D^1\times\tilde S_I,Z_{r,\alpha,I}\times\mathbb D^1)/\tilde S_I),
\mathbb Z(u_{IJ}^r\times I))
\to(\mathbb Z((Y_{r,\alpha,I}\times\tilde S_I,Z_{r,\alpha})/\tilde S_I),\mathbb Z(u_{IJ}^r)))
\end{eqnarray*}
\item[(ii)] A similar statement holds for $(\mathbb D^1,usu)$ local equivalence
$m:F\to G$ with $F,G\in C(\AnSp(\mathbb C)^{2,(sm)}/(\tilde S_I))$.
\item[(ii)'] A similar statement holds for $(\mathbb D^1,usu)$ local equivalence
$m:F\to G$ with $F,G\in C(\AnSp(\mathbb C)^{2,(sm)pr}/(\tilde S_I))$.
\end{itemize}
\end{prop}

\begin{proof}
Standard. See Ayoub's thesis for example.
\end{proof}

\begin{defi}
\begin{itemize}
\item[(i)]A filtered complex $(G,F)\in C_{fil}(\AnSp(\mathbb C)^{2,(sm)}/S)$ 
is said to be $r$-filtered $\mathbb D^1$ homotopic if 
$\ad(p_a^*,p_{a*})(G,F):(G,F)\to p_{a*}p_a^*(G,F)$ is an $r$-filtered homotopy equivalence.
\item[(i)']A filtered complex $((G_I,F),u_{IJ})\in C_{fil}(\AnSp(\mathbb C)^{2,(sm)}/(\tilde S_I))$ 
is said to be $r$-filtered $\mathbb D^1$ homotopic if 
$\ad(p_a^*,p_{a*})((G_I,F),u_{IJ}):((G_I,F),u_{IJ})\to p_{a*}p_a^*((G_I,F),u_{IJ})$ is an $r$-filtered homotopy equivalence.
\item[(ii)]A filtered complex $(G,F)\in C_{fil}(\AnSp(\mathbb C)^{2,(sm)pr}/S)$ 
is said to be $r$-filtered $\mathbb D^1$ homotopic if 
$\ad(p_a^*,p_{a*})(G,F):(G,F)\to p_{a*}p_a^*(G,F)$ is an $r$-filtered homotopy equivalence.
\item[(ii)']A filtered complex $((G_I,F),u_{IJ})\in C_{fil}(\AnSp(\mathbb C)^{2,(sm)pr}/(\tilde S_I))$ 
is said to be $r$-filtered $\mathbb D^1$ homotopic if 
$\ad(p_a^*,p_{a*})((G_I,F),u_{IJ}):((G_I,F),u_{IJ})\to p_{a*}p_a^*((G_I,F),u_{IJ})$ is an $r$-filtered homotopy equivalence.
\end{itemize}
\end{defi}

We have the following canonical functor :

\begin{defi}\label{eta12defan}
\begin{itemize}
\item[(i)]For $S\in\AnSp(\mathbb C)$, we have the functor 
\begin{eqnarray*}
(-)^{\Gamma}:C(\AnSp(\mathbb C)^{sm}/S)\to C(\AnSp(\mathbb C)^{2,sm}/S), \\ 
F\longmapsto F^{\Gamma}:(((U,Z)/S)=((U,Z),h)\mapsto F^{\Gamma}((U,Z)/S):=(\Gamma^{\vee}_Zh^*LF)(U/U), \\ 
(g:((U',Z'),h')\to((U,Z),h))\mapsto \\
(F^{\Gamma}(g):(\Gamma^{\vee}_Zh^*LF)(U/U)\xrightarrow{i_{(\Gamma^{\vee}_Zh^*LF)(U/U)}}(g^*(\Gamma^{\vee}_Zh^*LF))(U'/U') \\ 
\xrightarrow{T(g,\gamma^{\vee})(h^*LF)(U'/U')}(\Gamma^{\vee}_{Z\times_U U'}g^*h^*LF)(U'/U') \\ 
\xrightarrow{T(Z'/Z\times_U U',\gamma^{\vee})(g^*h^*LF)(U'/U')}(\Gamma^{\vee}_{Z'}g^*h^*LF)(U'/U')))
\end{eqnarray*}
where $i_{(\Gamma^{\vee}_Zh^*LF)(U/U)}$ is the canonical arrow of the inductive limit.
Similarly, we have, for $S\in\AnSp(\mathbb C)$, the functor
\begin{eqnarray*}
(-)^{\Gamma}:C(\AnSp(\mathbb C)/S)\to C(\AnSp(\mathbb C)^2/S), \\ 
F\longmapsto F^{\Gamma}:(((X,Z)/S)=((X,Z),h)\mapsto F^{\Gamma}((X,Z)/S):=(\Gamma^{\vee}_Zh^*F)(X/X), \\ 
(g:((X',Z'),h')\to((X,Z),h))\mapsto 
(F^{\Gamma}(g):(\Gamma^{\vee}_Zh^*LF)(X/X)\to (\Gamma^{\vee}_{Z'}h^{'*}LF)(X'/X')))
\end{eqnarray*}
Note that for $S\in\AnSp(\mathbb C)$, $I(S/S):\mathbb Z((S,S)/S)\to\mathbb Z(S/S)^{\Gamma}$ given by
\begin{eqnarray*}
I(S/S)((U,Z),h):\mathbb Z((S,S)/S)(((U,Z),h))\xrightarrow{\gamma^{\vee}_Z(\mathbb Z(U/U))(U/U)} 
\mathbb Z(S/S)^{\Gamma}((U,Z),h):=(\Gamma_Z^{\vee}\mathbb Z(U/U))(U/U)
\end{eqnarray*} 
is an isomorphism.
\item[(ii)]Let $f:T\to S$ a morphism with $T,S\in\AnSp(\mathbb C)$.
For $F\in C(\AnSp(\mathbb C)^{sm}/S)$, we have the canonical morphism in $C(\AnSp(\mathbb C)^{2,sm}/T)$
\begin{eqnarray*}
T(f,\Gamma)(F):=T^*(f,\Gamma)(F):f^*(F^{\Gamma})\to(f^*F)^{\Gamma},\\
T(f,\Gamma)(F)((U',Z')/T=((U',Z'),h')): \\
f^*(F^{\Gamma})((U',Z'),h'):=
\lim_{((U',Z'),h')\xrightarrow{l}((U_T,Z_T),h_T)\xrightarrow{f_U}((U,Z),h)}(\Gamma^{\vee}_Zh^*LF)(U/U) \\
\xrightarrow{F^{\Gamma}(f_U\circ l)}(\Gamma^{\vee}_{Z'}l^*f_U^*h^*LF)(U'/U')=(\Gamma^{\vee}_{Z'}h^{'*}f^*LF)(U'/U') \\
\xrightarrow{(\Gamma^{\vee}_{Z'}h^{'*}T(f,L)(F))(U'/U')}(\Gamma^{\vee}_{Z'}h^{'*}Lf^*F)(U'/U')=:(f^*F)^{\Gamma}((U',Z'),h') 
\end{eqnarray*}
where $f_U:U_T:U\times_S T\to U$ and $h_T:U_T:=U\times_S T\to T$ are the base change maps,
the equality following from the fact that $h\circ f_U\circ l=f\circ h_T\circ l=f\circ h'$.
For $F\in C(\AnSp(\mathbb C)/S)$, we have similarly the canonical morphism in $C(\AnSp(\mathbb C)^2/T)$
\begin{equation*}
T(f,\Gamma)(F):f^*(F^{\Gamma})\to (f^*F)^{\Gamma}.
\end{equation*}
\item[(iii)]Let $h:U\to S$ a smooth morphism with $U,S\in\AnSp(\mathbb C)$. 
We have, for $F\in C(\AnSp(\mathbb C)^{sm}/U)$, the canonical morphism in $C(\AnSp(\mathbb C)^{2,sm}/S)$
\begin{eqnarray*}
T_{\sharp}(h,\Gamma)(F):h_{\sharp}(F^{\Gamma})\to(h_{\sharp}LF)^{\Gamma}, \\
T_{\sharp}(h,\Gamma)(F)((U',Z'),h'):h_{\sharp}(F^{\Gamma})((U',Z'),h'):=
\lim_{((U',Z'),h')\xrightarrow{l}((U,U),h)}(\Gamma^{\vee}_{Z'}l^*LF)(U'/U') \\
\xrightarrow{(\Gamma^{\vee}_{Z'}l^*\ad(h_{\sharp},h^*)(LF))(U'/U')}
(\Gamma^{\vee}_{Z'}l^*h^*h_{\sharp}LF)(U'/U')=:(h_{\sharp}LF)^{\Gamma}((U',Z')/h')
\end{eqnarray*}
\item[(iv)]Let $i:Z_0\hookrightarrow S$ a closed embedding with $Z_0,S\in\AnSp(\mathbb C)$. 
We have the canonical morphism in $C(\AnSp(\mathbb C)^{2,sm}/S)$
\begin{eqnarray*}
T_*(i,\Gamma)(\mathbb Z(Z_0/Z_0)):i_*((\mathbb Z(Z_0/Z_0))^{\Gamma}\to(i_*\mathbb Z(Z/Z))^{\Gamma}, \\
T_*(i,\Gamma)(\mathbb Z(Z_0/Z_0))((U,Z),h):
i_*((\mathbb Z(Z_0/Z_0))^{\Gamma}((U,Z),h):=(\Gamma^{\vee}_{Z\times_SZ_0}\mathbb Z(Z_0/Z_0))(U\times_S Z_0) \\
\xrightarrow{T(i_*,\gamma^{\vee})(\mathbb Z(Z_0/Z_0))(U\times_S Z_0)}
(\Gamma^{\vee}_Zi_*\mathbb Z(Z_0/Z_0))(U\times_S Z_0)=:(i_*\mathbb Z(Z/Z))^{\Gamma}((U,Z),h)
\end{eqnarray*}
\end{itemize}
\end{defi}

\begin{defi}\label{GrGammaan}
Let $S\in\AnSp(\mathbb C)$. We have for $F\in C(\AnSp(\mathbb C)^{sm}/S)$ the canonical map in $C(\AnSp(\mathbb C)^{sm}/S)$
\begin{eqnarray*}
\Gr(F):\Gr^{12}_{S*}\mu_{S*}F^{\Gamma}\to F, \\
\Gr(F)(U/S):\Gamma_U^{\vee}p^*F(U\times S/U\times S)\xrightarrow{\ad(l^*,l_*)(p^*F)(U\times S/U\times S)}h^*F(U/U)=F(U/S)
\end{eqnarray*}
where $h:U\to S$ is a smooth morphism with $U\in\AnSp(\mathbb C)$ and $h:U\xrightarrow{l}U\times S\xrightarrow{p}S$
is the graph factorization with $l$ the graph embedding and $p$ the projection.
\end{defi}

\begin{prop}\label{eta12propan}
Let $S\in\AnSp(\mathbb C)$. 
\begin{itemize}
\item[(i)] Then,
\begin{itemize}
\item if $m:F\to G$ with $F,G\in C(\AnSp(\mathbb C)^{sm}/S)$ is a quasi-isomorphism,
$m^{\Gamma}:F^{\Gamma}\to G^{\Gamma}$ is a quasi-isomorphism in $C(\AnSp(\mathbb C)^{2,sm}/S)$,
\item if $m:F\to G$ with $F,G\in C(\AnSp(\mathbb C)^{sm}/S)$ is an usu local equivalence,
$m^{\Gamma}:F^{\Gamma}\to G^{\Gamma}$ is an usu local equivalence in $C(\AnSp(\mathbb C)^{2,sm}/S)$,
\item if $m:F\to G$ with $F,G\in C(\AnSp(\mathbb C)^{sm}/S)$ is an $(\mathbb D^1,usu)$ local equivalence,
$m^{\Gamma}:F^{\Gamma}\to G^{\Gamma}$ is an $(\mathbb D^1,usu)$ local equivalence in $C(\AnSp(\mathbb C)^{2,sm}/S)$.
\end{itemize}
\item[(ii)] Then,
\begin{itemize}
\item if $m:F\to G$ with $F,G\in C(\AnSp(\mathbb C)/S)$ is a quasi-isomorphism,
$m^{\Gamma}:F^{\Gamma}\to G^{\Gamma}$ is a quasi-isomorphism in $C(\AnSp(\mathbb C)^2/S)$,
\item if $m:F\to G$ with $F,G\in C(\AnSp(\mathbb C)^{sm}/S)$ is an usu local equivalence,
$m^{\Gamma}:F^{\Gamma}\to G^{\Gamma}$ is an usu local equivalence in $C(\AnSp(\mathbb C)^{2,sm}/S)$,
\item if $m:F\to G$ with $F,G\in C(\AnSp(\mathbb C)^{sm}/S)$ is an $(\mathbb D^1,usu)$ local equivalence,
$m^{\Gamma}:F^{\Gamma}\to G^{\Gamma}$ is an $(\mathbb D^1,usu)$ local equivalence in $C(\AnSp(\mathbb C)^2/S)$.
\end{itemize}
\end{itemize}
\end{prop}

\begin{proof}
Similar to the proof of proposition \ref{eta12prop}.
\end{proof}

\subsection{The analytical functor for presheaves on the big Zariski or etale site 
and on the big Zariski or etale site of pairs}

We have for $f:T\to S$ a morphism with $T,S\in\Var(\mathbb C)$ the following commutative diagram of sites
\begin{equation}\label{AnVar}
Dia(S):=\xymatrix{
\AnSp(\mathbb C)/T^{an}\ar[rr]^{\An_T}\ar[dd]^{P(f)}\ar[rd]^{\rho_T} & \, & 
\Var(\mathbb C)/T\ar[dd]^{P(f)}\ar[rd]^{\rho_T} & \, \\  
 \, & \AnSp(\mathbb C)^{sm}/T^{an}\ar[rr]^{\An_T}\ar[dd]^{P(f)} & \, & \Var(\mathbb C)^{sm}/T\ar[dd]^{P(f)} \\  
\AnSp(\mathbb C)/S^{an}\ar[rr]^{\An_S}\ar[rd]^{\rho_S} & \, & \Var(\mathbb C)/S\ar[rd]^{\rho_S} & \, \\  
\, &  \AnSp(\mathbb C)^{sm}/S^{an}\ar[rr]^{\An_S} & \, & \Var(\mathbb C)^{sm}/S}  
\end{equation}
and
\begin{equation}\label{AnVar12}
Dia^{12}(S):=\xymatrix{
\AnSp(\mathbb C)^2/T^{an}\ar[rr]^{\An_T}\ar[dd]^{P(f)}\ar[rd]^{\rho_T} & \, & 
\Var(\mathbb C)^2/T\ar[dd]^{P(f)}\ar[rd]^{\rho_T} & \, \\  
 \, & \AnSp(\mathbb C)^{2,sm}/T^{an}\ar[rr]^{\An_T}\ar[dd]^{P(f)} & \, & \Var(\mathbb C)^{2,sm}/T\ar[dd]^{P(f)} \\  
\AnSp(\mathbb C)^2/S^{an}\ar[rr]^{\An_S}\ar[rd]^{\rho_S} & \, & \Var(\mathbb C)^2/S\ar[rd]^{\rho_S} & \, \\  
\, &  \AnSp(\mathbb C)^{2,sm}/S^{an}\ar[rr]^{\An_S} & \, & \Var(\mathbb C)^{2,sm}/S}  
\end{equation}
For $S\in\Var(\mathbb C)$ we have the following commutative diagrams of sites
\begin{equation*}
\xymatrix{\AnSp(\mathbb C)^2/S\ar[rr]^{\mu_S}\ar[dd]_{\An_S}\ar[rd]^{\rho_S} & \, & 
\AnSp(\mathbb C)^{2,pr}/S\ar[dd]^{\An_S}\ar[rd]^{\rho_S} & \, \\
\, & \AnSp(\mathbb C)^{2,sm}/S\ar[rr]^{\mu_S}\ar[dd]_{\An_S} & \, & \AnSp(\mathbb C)^{2,smpr}/S\ar[dd]^{\An_S} \\
\Var(\mathbb C)^2/S\ar[rr]^{\mu_S}\ar[rd]^{\rho_S} & \, &  \Var(\mathbb C)^{2,smpr}/S\ar[rd]^{\rho_S} & \, \\
\, & \Var(\mathbb C)^2/S\ar[rr]^{\mu_S} & \, &  \Var(\mathbb C)^{2,smpr}/S & \,}
\end{equation*}
and
\begin{equation}\label{GrAn}
\xymatrix{\AnSp(\mathbb C)^{2,pr}/S\ar[rr]^{\Gr_S^{12}}\ar[dd]_{\An_S}\ar[rd]^{\rho_S} & \, & 
\AnSp(\mathbb C)/S\ar[dd]^{\An_S}\ar[rd]^{\rho_S} & \, \\
\, & \AnSp(\mathbb C)^{2,smpr}/S\ar[rr]^{\Gr_S^{12}}\ar[dd]_{\An_S} & \, & \AnSp(\mathbb C)^{sm}/S\ar[dd]^{\An_S} \\
\Var(\mathbb C)^{2,pr}/S\ar[rr]^{\Gr_S^{12}}\ar[rd]^{\rho_S} & \, & \Var(\mathbb C)/S\ar[rd]^{\rho_S} & \, \\
\, & \Var(\mathbb C)^{2,sm}/S\ar[rr]^{\Gr_S^{12}} & \, & \Var(\mathbb C)^{sm}/S}.
\end{equation}
For $f:T\to S$ a morphism in $\Var(\mathbb C)$ the diagramm $Dia(S)$ and $Dia(T)$ commutes with the pullback functors :
we have $e(S)\circ P(f)=P(f)\circ e(T)$.

For $S\in\Var(\mathbb C)$, the analytical functor is 
\begin{equation*}
(-)^{an}:C_{O_S}(S)\to C_{O_{S^{an}}}, G\mapsto G^{an}:=\an_S^{*mod}G:=\an_S^*G\otimes_{\an_S^*O_S}O_{S^{an}}
\end{equation*}

Let $S\in\Var(\mathbb C)$.   
\begin{itemize}
\item As $\an_S^*:\PSh(S)\to\PSh(S^{an})$ preserve monomorphisms 
(the colimits involved being filtered colimits), we define, for $(G,F)\in C_{(2)fil}(S)$, 
$\an_S^*(G,F):=(\an_S^*G,\an_S^*F)\in C_{(2)fil}(S^{an})$.
\item As $(-)^{an}:=\an_S^{*mod}:\PSh_{O_S}(S)\to\PSh(S^{an})$ preserve monomorphisms 
($\an_S^*$ preserve monomorphism and $(-)\otimes_{O_S}O_{S^{an}}$ preserve monomorphism 
since $O_{S^{an}}$ is a flat $O_S$ module), we define, for $(G,F)\in C_{(2)fil}(S)$, 
$(G,F)^{an}:=(G^{an},\an_S^*F\otimes_{O_S}O_{S^{an}})\in C_{(2)fil}(S^{an})$.
\end{itemize}

Let $f:T\to S$ a morphism with $T,S\in\Var(\mathbb C)$. Then,
\begin{itemize}
\item the commutative diagrams of sites $D(\An,f):=(\An_S,f,\An_T,f=f^{an})$ gives,
for $G\in C(\Var(\mathbb C)^{sm}/T)$, the canonical map in $C(\AnSp(\mathbb C)^{am}/T)$
\begin{eqnarray*}
T(\An,f)(G):\An_S^*f_*G\xrightarrow{\ad(\An_T^*,\An_{T*})(G)}\An_S^*f_*\An_{T*}\An_T^*G=\An_S^*\An_{S*}f_*\An_T^*G \\
\xrightarrow{\ad(\An_S^*\An_{S*})(f_*\An_T^*G)}f_*\An_T^*G.
\end{eqnarray*} 
\item the commutative diagrams of sites $D(an,f):=(\an_S,f\an_T,f)$ gives, for $G\in C(T)$, the canonical map in $C(T^{an})$
\begin{eqnarray*}
T(an,f)(G):\an_S^*f_*G\xrightarrow{\ad(\an_T^*,\an_{T*})(G)}\an_S^*f_*\an_{T*}\an_T^*G=\an_S^*\an_{S*}f_*\an_T^*G \\
\xrightarrow{\ad(\an_S^*,\an_{S*})(f_*\an_T^*G)}f_*\an_T^*G
\end{eqnarray*}
and for $G\in C_{O_T}(T)$, the canonical map in $C_{O_{T^{an}}}(T^{an})$
\begin{eqnarray*}
T^{mod}(an,f)(G):(f_*G)^{an}:=\an_S^{*mod}f_*G\xrightarrow{\ad(\an_T^{*mod},\an_{T*})(G)} \\
\an_S^{*mod}f_*\an_{T*}\an_T^{*mod}G=\an_S^{*mod}\an_{S*}f_*\an_T^{*mod}G 
\xrightarrow{\an_S^{*mod},\an_{S*})(f_*\an_T^*G)}f_*\an_T^{*mod}G=:f_*G^{an}
\end{eqnarray*}
\end{itemize}

\begin{defiprop}\label{gamma4sect2}
Consider a closed embedding $i:Z\hookrightarrow S$ with $S,Z\in\Var(\mathbb C)$.
Then, for $G^{\bullet}\in C(\Var(\mathbb C)^{sm}/S)$, there exist a map in $C(\AnSp(\mathbb C)^{sm}/S)$
\begin{equation*}
T(\An,\gamma)(G):\An_S^*\Gamma_{Z}G\to\Gamma_{Z}\An_S^*G
\end{equation*}
unique up to homotopy, such that $\gamma_{Z}(\An_S^*G)\circ T(\An,\gamma)(G)=\An_S^*\gamma_{Z}G$. 
\end{defiprop}

\begin{proof}
Denote by $j:S\backslash Z\hookrightarrow S$ the open complementary embedding.
The map is given by $(I,T(\An,j)(j^*G)):\Cone(\An_S^*G\to\An_S^*j_*j^*G)\to(\An_S^*G\to j_*j^*\An_S^*G)$.
\end{proof}

\begin{defi}\label{projBMmotvaranprop}
Let $f:X\to S$ a morphism with $X,S\in\Var(\mathbb C)$. Assume that there exist a factorization
$f:X\xrightarrow{i} Y\times S\xrightarrow{p} S$, with $Y\in\SmVar(\mathbb C)$, 
$i:X\hookrightarrow Y$ is a closed embedding and $p$ the projection.
We then have the canonical isomorphism in $C(\AnSp(\mathbb C)^{sm}/S^{an})$
\begin{eqnarray*}
T(f,g,Q):=T_{\sharp}(\An,p)(-)^{-1}\circ T_{\sharp}(\An,j)(-)^{-1}: \\
\An_S^*Q(X/S):=An_S^*p_{\sharp}\Gamma^{\vee}_X\mathbb Z_{Y\times S}[d_Y]
\xrightarrow{=}p_{\sharp}\Gamma^{\vee}_{X^{an}}\mathbb Z_{Y\times S}[d_Y]=:Q(X^{an}/S^{an})
\end{eqnarray*}
with $j:Y\times S\backslash X\hookrightarrow Y\times S$ the closed embedding.
\end{defi}

\begin{defiprop}
Consider a closed embedding $i:Z\hookrightarrow S$  with $S\in\Var(\mathbb C)$. 
Then, for $G\in C_{O_S}(S)$, there is a canonical map in $C_{O_{S^{an}}}(S^{an})$
\begin{equation*}
T^{mod}(an,\gamma)(G):(\Gamma_{Z}G)^{an}\to\Gamma_{Z^{an}}G^{an}
\end{equation*}
unique up to homotopy, such that $\gamma_{Z^{an}}(G^{an})\circ T^{mod}(an,\gamma)(G)=g^*\gamma_{Z}G$.
\end{defiprop}

\begin{proof}
It is a particular case of definition-proposition \ref{gamma1sect2mod}(i).
\end{proof}

We recall the first GAGA theorem for coherent sheaf on the projective spaces :

\begin{thm}\label{GAGA1}
For $X\in\Var(\mathbb C)$ and $F\in C_{O_X}(X)$ denote by 
\begin{equation*}
a(F):\ad(\an_X^{*mod},\an(X)_*)(E(F)):E(F)\to\an_{X*}(E(F))^{an}=\an_{X*}E(F^{an}), 
\end{equation*}
the canonical morphism.
\begin{itemize}
\item[(i)] Let $X\in\PVar(\mathbb C)$ a proper complex algebraic variety. For $F\in\Coh_{O_X}(X)$ a coherent sheaf, the morphism 
\begin{equation*}
H^n\Gamma(X,a(F)):H^n(X,F)=H^n\Gamma(X,E(F))\to H^n(X,F^{an})=H^n\Gamma(X,E(F^{an}))
\end{equation*}
is an isomorphism for all $n\in\mathbb Z$.
\item[(ii)] Let $f:X\to S$ a proper morphism with $X,S\in\Var(\mathbb C)$. For $F\in\Coh_{O_X}(X)$ a coherent sheaf, the morphism 
\begin{equation*}
H^nf_*a(F):R^nf_*F=H^nf_*(E(F))\to R^nf_*F^{an}=H^nf_*E(F^{an})
\end{equation*}
is an isomorphism for all $n\in\mathbb Z$.
\end{itemize}
\end{thm}

\begin{proof}
See \cite{Serre}. (i) reduces to the case where $X$ is projective and (ii) to the case where $f$ is projective.
Hence, the theorem reduce to the case of a coherent sheaf $F\in\Coh_{O_{\mathbb P^N}}(\mathbb P^N)$ on $\mathbb P^N$.
\end{proof}

We have for $s:\mathcal I\to\mathcal J$ a functor with $\mathcal I,\mathcal J\in\Cat$ and
$f:T_{\bullet}\to S_{s(\bullet)}$ a morphism of diagram of algebraic varieties with 
$T_{\bullet}\in\Fun(\mathcal I,\Var(\mathbb C))$, $S_{\bullet}\in\Fun(\mathcal J,\Var(\mathbb C))$ 
the following commutative diagram of sites
\begin{equation}\label{AnVarIJ}
Dia(S):=\xymatrix{
\AnSp(\mathbb C)/T_{\bullet}^{an}\ar[rr]^{\An_{T_{\bullet}}}\ar[dd]^{P(f)}\ar[rd]^{\rho_{T_{\bullet}}} & \, & 
\Var(\mathbb C)/T_{\bullet}\ar[dd]^{P(f_{\bullet})}\ar[rd]^{\rho_{T_{\bullet}}} & \, \\  
 \, & \AnSp(\mathbb C)^{sm}/T_{\bullet}^{an}\ar[rr]^{\An_{T_{\bullet}}}\ar[dd]^{P(f_{\bullet})} & \, & 
\Var(\mathbb C)^{sm}/T_{\bullet}\ar[dd]^{P(f_{\bullet})} \\  
\AnSp(\mathbb C)/S_{\bullet}^{an}\ar[rr]^{\An_{S_{\bullet}}}\ar[rd]^{\rho_{S_{\bullet}}} & \, & 
\Var(\mathbb C)/S_{\bullet}\ar[rd]^{\rho_{S_{\bullet}}} & \, \\  
\, &  \AnSp(\mathbb C)^{sm}/S_{\bullet}^{an}\ar[rr]^{\An_{S_{\bullet}}} & \, & 
\Var(\mathbb C)^{sm}/S_{\bullet}}  
\end{equation}
and
\begin{equation}\label{AnVar12IJ}
Dia^{12}(S):=\xymatrix{
\AnSp(\mathbb C)^2/T_{\bullet}^{an}\ar[rr]^{\An_{T_{\bullet}}}\ar[dd]^{P(f_{\bullet})}\ar[rd]^{\rho_{T_{\bullet}}} & \, & 
\Var(\mathbb C)^2/T_{\bullet}\ar[dd]^{P(f_{\bullet})}\ar[rd]^{\rho_{T_{\bullet}}} & \, \\  
 \, & \AnSp(\mathbb C)^{2,sm}/T_{\bullet}^{an}\ar[rr]^{\An_{T_{\bullet}}}\ar[dd]^{P(f_{\bullet})} & \, & 
\Var(\mathbb C)^{2,sm}/T_{\bullet}\ar[dd]^{P(f_{\bullet})} \\  
\AnSp(\mathbb C)^2/S_{\bullet}^{an}\ar[rr]^{\An_{S_{\bullet}}}\ar[rd]^{\rho_{S_{\bullet}}} & \, & 
\Var(\mathbb C)^2/S_{\bullet}\ar[rd]^{\rho_{S_{\bullet}}} & \, \\  
\, &  \AnSp(\mathbb C)^{2,sm}/S_{\bullet}^{an}\ar[rr]^{\An_{S_{\bullet}}} & \, & 
\Var(\mathbb C)^{2,sm}/S_{\bullet}}.  
\end{equation}

\subsection{The De Rahm complexes of algebraic varieties and analytical spaces}

For $X\in\Var(\mathbb C)$, we denote by $\iota_X:\mathbb C_{X}\to\Omega^{\bullet}_{X}=:DR(X)$ the canonical inclusion map.
More generaly, for $f:X\to S$ a morphism with $X,S\in\Var(\mathbb C)$, 
we denote by $\iota_{X/S}:f^*O_S\to\Omega^{\bullet}_{X/S}=:DR(X/S)$ the canonical inclusion map.

For $X\in\AnSp(\mathbb C)$, we denote by $\iota_X:\mathbb C_{X}\to\Omega^{\bullet}_{X}=:DR(X)$ the canonical inclusion map.
More generaly, for $f:X\to S$ a morphism with $X,S\in\AnSp(\mathbb C)$, 
we denote by  $\iota_{X/S}:f^*O_S\to\Omega^{\bullet}_{X/S}=:DR(X/S)$ the canonical inclusion map.

Let $f:X\to S$ a morphism with $X,S\in\Var(\mathbb C)$. 
Then, the commutative diagram of site $(an,f):=(f,\An_S,f=f^{an},\an(X))$ 
gives the transformation map in $C_{O_{S^{an}}}(S^{an})$ (definition \ref{TDw})
\begin{eqnarray*}
T^O_{\omega}(an,f):(f_*E(\Omega^{\bullet}_{X/S},F_b))^{an}:=\an_S^{*mod}f_*E(\Omega^{\bullet}_{X/S},F_b)
\xrightarrow{T(an(X),E)(-)\circ T(an,f)(E(\Omega^{\bullet}_{X/S}))} \\
(f_*E(\an(X)^*(\Omega^{\bullet}_{X/S},F_b)))\otimes_{\an_S^*O_S}O_{S^{an}}
\xrightarrow{m\otimes E(\Omega_{(X^{an}/X)/(S^{an}/S)})}
f_*E(\Omega^{\bullet}_{X^{an}/S^{an}},F_b)
\end{eqnarray*}

We will give is this paper a relative version for all smooth morphisms of the following theorem of Grothendieck

\begin{thm}\label{GAGA1w}
Let $U\in\SmVar(\mathbb C)$. Denote by $a_U:U\to\left\{\pt\right\}$ the terminal map. Then the map
\begin{equation*}
T_{\omega}^O(a_U,an):\Gamma(U,E(\Omega^{\bullet}_U))\to\Gamma(U^{an},E(\Omega^{\bullet}_U))
\end{equation*}
is a quasi-isomorphism of complexes.
\end{thm}

\begin{proof}
Take a compactification $(X,D)$ of $U$, with $X\in\PSmVar(\mathbb C)$ and $D=X\backslash U$ a normal crossing divisor.
The proof then use proposition \ref{star}, the first GAGA theorem (theorem \ref{GAGA1} (i)) 
for the coherent sheaves $\Omega^p_U(nD)$ on $X$,
and the fact (which is specific of the De Rahm complex) that
$\Omega^{\bullet}_{U^{an}}(*D^{an})\to j_*E(\Omega_{U^{an}})$ is a quasi-isomorphism.
\end{proof}

We recall Poincare lemma for smooth morphisms of complex analytic spaces and in particular complex analytic manifold :
\begin{prop}\label{Poincarelem}
\begin{itemize}
\item[(i)] For $h:U\to S$ a smooth morphism with $U,S\in\AnSp(\mathbb C)$, the inclusion map 
$\iota_{X/S}:h^*O_S\to\Omega^{\bullet}_{U/S}=:DR(U/S)$ is a quasi-isomorphism.
\item[(ii)] For $X\in\AnSm(\mathbb C)$, the inclusion map $\iota_X:\mathbb C_{X}\to\Omega^{\bullet}_{X}$ is a quasi-isomorphism.
\end{itemize}
\end{prop}

\begin{proof}
Standard. (ii) is a particular case of (i) (the absolute case $S=\left\{\pt\right\}$).
\end{proof}

\begin{rem}\label{Poincaresing}
We do NOT have poincare lemma in general if $h:U\to S$ is not a smooth morphism.
Already in the absolute case, we can find $X\in\Var(\mathbb C)$ singular such that
the inclusion map $\iota_X:\mathbb C_{X^{an}}\to\Omega^{\bullet}_{X^{an}}$ is not a quasi-isomorphism.
Indeed, we can find exemple of $X\in\PVar(\mathbb C)$ projective singular where 
\begin{equation*}
H^p(c_X):H^p(X^{an},\mathbb C_{X^{an}})\xrightarrow{\sim}H^pC^{\bullet}_{sing}(X^{an})
\end{equation*}
$X^{an}$ being locally contractible since $X^{an}\in\CW$, have not the same dimension then the De Rham cohomology 
\begin{equation*}
H^p(T^O_{\omega}(an,a_X)):\mathbb H^p(X,E(\Omega^{\bullet}_X))\xrightarrow{\sim}\mathbb H^p(X^{an},E(\Omega^{\bullet}_{X^{an})})
\end{equation*}
$X$ being projective, that is are not isomorphic as vector spaces. Hence, in particular, the canonical map 
\begin{equation*}
H^p\iota_X:H^p(X^{an},\mathbb C_{X^{an}})\to\mathbb H^p(X^{an},E(\Omega^{\bullet}_{X^{an}}))
\end{equation*}
 is not an isomorphism.
\end{rem}

Consider a commutative diagram
\begin{equation*}
D_0=\xymatrix{f: X\ar[r]^{i} & Y\ar[r]^{p} & S \\
f':X'\ar[r]^{i'}\ar[u]^{g'} & Y'\ar[u]^{g''}\ar[r]^{p'} & T\ar[u]^{g} }
\end{equation*}
with $X,X',Y,Y',S,T\in\Var(\mathbb C)$ or $X,X',Y,Y',S,T\in\AnSp(\mathbb C)$, $i$, $i'$ being closed embeddings. 
Denote by $D$ the right square of $D_0$. The closed embedding $i':X'\hookrightarrow Y'$ factors through
$i':X'\xrightarrow{i'_1} X\times_Y Y'\xrightarrow{i'_0} Y'$ where $i'_1,i'_0$ are closed embeddings. 
Then, definition-proposition \ref{TDwgamma} say that
\begin{itemize}
\item there is a canonical map,
\begin{eqnarray*}
E(\Omega_{((Y'/Y))/(T/S)})\circ T(g'',E)(-)\circ T(g'',\gamma)(-):
g^{''*}\Gamma_{X}E(\Omega^{\bullet}_{Y/S},F_b)\to\Gamma_{X\times_Y Y'}E(\Omega^{\bullet}_{Y'/T},F_b)
\end{eqnarray*}
unique up to homotopy such that the following diagram in $C_{g^{''*}p^*O_Sfil}(Y')=C_{p^{'*}g^*O_Sfil}(Y')$ commutes
\begin{equation*}
\xymatrix{g^{''*}\Gamma_{X}E(\Omega^{\bullet}_{Y/S},F_b)
\ar[rrrr]^{E(\Omega_{((Y'/Y))/(T/S)})\circ T(g'',E)(-)\circ T(g'',\gamma)(-)}\ar[d]_{\gamma_X(-)} & \, & \, & \, &
\Gamma_{X\times_Y Y'}E(\Omega^{\bullet}_{Y'/T},F_b)\ar[d]^{\gamma_{X\times_Y Y'}(-)} \\
g^{''*}E(\Omega^{\bullet}_{Y/S},F_b)\ar[rrrr]^{E(\Omega_{((Y'/Y)/(T/S))}\circ T(g'',E)(-)} 
& \, & \, & \, &  E(\Omega^{\bullet}_{Y'/T},F_b)},
\end{equation*}
\item there is a canonical map,
\begin{equation*}
T^O_{\omega}(D)^{\gamma}:g^{*mod}L_Op_*\Gamma_{X}E(\Omega^{\bullet}_{Y/S},F_b)\to p'_*\Gamma_{X\times_Y Y'}E(\Omega^{\bullet}_{Y'/T},F_b)
\end{equation*}
unique up to homotopy such that the following diagram in $C_{O_{T}fil}(T)$ commutes
\begin{equation*}
\xymatrix{g^{*mod}L_Op_*\Gamma_{X}E(\Omega^{\bullet}_{Y/S})
\ar[rrr]^{T_{\omega}^O(D)^{\gamma}}\ar[d]_{\gamma_X(-)} & \, & \, & 
p'_*\Gamma_{X\times_Y Y'}E(\Omega^{\bullet}_{Y'/T})\ar[d]^{\gamma_{X\times_Y Y'}(-)} \\
g^{*mod}L_Op_*E(\Omega^{\bullet}_{Y/S})\ar[rrr]^{T^O_{\omega}(D)} & \, & \, &  
p'_*E(\Omega^{\bullet}_{Y'/T})}.
\end{equation*}
\item[(iii)] there is a map in $C_{f^{'*}O_Tfil}(Y')$
\begin{equation*}
T(X'/X\times_Y Y',\gamma)(E(\Omega^{\bullet}_{Y'/T},F_b)):\Gamma_{X'}E(\Omega^{\bullet}_{Y'/T},F_b)\to
\Gamma_{X\times_Y Y'}E(\Omega^{\bullet}_{Y'/T},F_b)
\end{equation*}
unique up to homotopy such that $\gamma_{X\times_Y Y'}(-)\circ T^(X'/X\times_Y Y',\gamma)(-)=\gamma_{X'}(-)$.
\end{itemize}

Let $h:Y\to S$ a morphism and $i:X\hookrightarrow Y$ a closed embedding with $S,Y,X\in\Var(\mathbb C)$. 
Then, definition-proposition \ref{TDwgamma} say that
\begin{itemize}
\item there is a canonical map
\begin{equation*}
E(\Omega_{(Y^{an}/Y)/(S^{an}/S)})\circ T(an,\gamma)(-):
\an(Y)^*\Gamma_XE(\Omega^{\bullet}_{Y/S},F_b)\to\Gamma_{X^{an}}E(\Omega^{\bullet}_{Y/S},F_b)
\end{equation*}
unique up to homotopy such that the following diagram in $C_{h^*O_Sfil}(Y^{an})$ commutes
\begin{equation*}
\xymatrix{
\an(Y)^*\Gamma_XE(\Omega^{\bullet}_{Y/S},F_b)
\ar[rrrr]^{E(\Omega_{(Y^{an}/Y)/(S^{an}/S)})\circ T(an,\gamma)(-)}\ar[d]_{\gamma_X(-)} & \, & \, & \, &  
\Gamma_{X^{an}}E(\Omega^{\bullet}_{Y/S},F_b)\ar[d]^{\gamma_{X^{an}}(-)} \\
\an(Y)^*E(\Omega^{\bullet}_{Y/S},F_b)\ar[rrrr]^{E(\Omega_{(Y^{an}/Y)/(S^{an}/S)})} & \, & \, & \, & E(\Omega^{\bullet}_{Y/S},F_b)}
\end{equation*}
\item there is a canonical map
\begin{equation*}
T_{\omega}^O(an,h)^{\gamma}:(h_*\Gamma_XE(\Omega^{\bullet}_{Y/S},F_b))^{an}\to h_*\Gamma_{X^{an}}E(\Omega^{\bullet}_{Y/S},F_b)
\end{equation*}
unique up to homotopy such that the following diagram in $C(Y)$commutes
\begin{equation*}
\xymatrix{
(h_*\Gamma_XE(\Omega^{\bullet}_{Y/S},F_b))^{an}\ar[rrrr]^{T_{\omega}^O(an,h)^{\gamma}}\ar[d]_{\gamma_X(-)} & \, & \, & \, & 
h_*\Gamma_{X^{an}}E(\Omega^{\bullet}_{Y/S},F_b)\ar[d]^{\gamma_{X^{an}}(-)} \\
(h_*E(\Omega^{\bullet}_{Y/S},F_b))^{an}\ar[rrrr]^{T_{\omega}^O(an,h)} & \, & \, & \, & h_*E(\Omega^{\bullet}_{Y/S},F_b)}
\end{equation*}
\end{itemize}

\subsection{The Corti-Hanamura resolution functors $\hat R^{CH}$ $R^{CH}$, $\hat R^{0CH}$, $R^{0CH}$
from complexes of representable presheaves on $\Var(\mathbb C)^{sm}/S$ with $S$ smooth, 
and the functorialities of these resolutions}

\begin{defi}\label{desVardef}
\begin{itemize}
\item[(i)]Let $X_0\in\Var(\mathbb C)$ and $Z\subset X_0$ a closed subset.
A desingularization of $(X_0,Z)$ is a pair of complex varieties $(X,D)\in\Var^2(\mathbb C))$,
together with a morphism of pair of varieties $\epsilon:(X,D)\to(X_0,\Delta)$ with $Z\subset\Delta$ such that 
\begin{itemize}
\item $X\in\SmVar(\mathbb C)$ and 
$D:=\epsilon^{-1}(\Delta)=\epsilon^{-1}(Z)\cup(\cup_i E_i)\subset X$ is a normal crossing divisor
\item $\epsilon:X\to X_0$ is a proper modification with discriminant $\Delta$, 
that is $\epsilon:X\to X_0$ is proper and $\epsilon:X\backslash D\xrightarrow{\sim}X\backslash\Delta$ is an isomorphism.
\end{itemize}
\item[(ii)]Let $X_0\in\Var(\mathbb C)$ and $Z\subset X_0$ a closed subset such that $X_0\backslash Z$ is smooth.
A strict desingularization of $(X_0,Z)$ is a pair of complex varieties $(X,D)\in\Var^2(\mathbb C))$,
together with a morphism of pair of varieties $\epsilon:(X,D)\to(X_0,Z)$ such that 
\begin{itemize}
\item $X\in\SmVar(\mathbb C)$ and $D:=\epsilon^{-1}(Z)\subset X$ is a normal crossing divisor
\item $\epsilon:X\to X_0$ is a proper modification with discriminant $Z$, 
that is $\epsilon:X\to X_0$ is proper and $\epsilon:X\backslash D\xrightarrow{\sim}X\backslash Z$ is an isomorphism.
\end{itemize}
\end{itemize}
\end{defi}

We have the following well known resolution of singularities of complex algebraic varieties
and their functorialities :

\begin{thm}\label{desVar}
\begin{itemize}
\item[(i)] Let $X_0\in\Var(\mathbb C)$ and $Z\subset X_0$ a closed subset.
There exists a desingularization of $(X_0,Z)$, that is a pair of complex varieties $(X,D)\in\Var^2(\mathbb C))$,
together with a morphism of pair of varieties $\epsilon:(X,D)\to(X_0,\Delta)$ with $Z\subset\Delta$ such that 
\begin{itemize}
\item $X\in\SmVar(\mathbb C)$ and 
$D:=\epsilon^{-1}(\Delta)=\epsilon^{-1}(Z)\cup(\cup_iE_i)\subset X$ is a normal crossing divisor
\item $\epsilon:X\to X_0$ is a proper modification with discriminant $\Delta$, 
that is $\epsilon:X\to X_0$ is proper and $\epsilon:X\backslash D\xrightarrow{\sim}X\backslash\Delta$ is an isomorphism.
\end{itemize}
\item[(ii)] Let $X_0\in\PVar(\mathbb C)$ and $Z\subset X_0$ a closed subset such that $X_0\backslash Z$ is smooth.
There exists a strict desingularization of $(X_0,Z)$, that is a pair of complex varieties $(X,D)\in\PVar^2(\mathbb C))$,
together with a morphism of pair of varieties $\epsilon:(X,D)\to(X_0,Z)$ such that 
\begin{itemize}
\item $X\in\PSmVar(\mathbb C)$ and $D:=\epsilon^{-1}(Z)\subset X$ is a normal crossing divisor
\item $\epsilon:X\to X_0$ is a proper modification with discriminant $Z$, 
that is $\epsilon:X\to X_0$ is proper and $\epsilon:X\backslash D\xrightarrow{\sim}X\backslash Z$ is an isomorphism.
\end{itemize}
\end{itemize}
\end{thm}

\begin{proof}
\noindent(i):Standard. See \cite{PS} for example.

\noindent(ii):Follows immediately from (i).
\end{proof}

We use this theorem to construct a resolution of a morphism by Corti-Hanamura morphisms,
we will need these resolution in the definition of the filtered De Rham realization functor :

\begin{defiprop}\label{RCHdef0}
\begin{itemize}
\item[(i)]Let $h:V\to S$ a morphism, with $V,S\in\Var(\mathbb C)$. 
Let $\bar{S}\in\PVar(\mathbb C)$ be a compactification of $S$.
\begin{itemize}
\item There exist a compactification $\bar X_0\in\PVar(\mathbb C)$ of $V$ such that $h:V\to S$
extend to a morphism $\bar f_0=\bar{h}_0:\bar X_0\to\bar{S}$. Denote by $\bar Z=\bar X_0\backslash V$.
We denote by $j:V\hookrightarrow\bar X_0$ the open embedding and by 
$i_0:\bar Z\hookrightarrow\bar X_0$ the complementary closed embedding.
We then consider $X_0:=\bar f_0^{-1}(S)\subset\bar X_0$ the open subset, $f_0:=\bar f_{0|X_0}:X_0\to S$, $Z=\bar Z\cap X_0$,
and we denote again $j:V\hookrightarrow X_0$ the open embedding and by $i_0:Z\hookrightarrow X_0$ the complementary closed embedding.
\item In the case $V$ is smooth, we take, using theorem \ref{desVar}(ii), a strict desingularization 
$\bar\epsilon:(\bar X,\bar D)\to(\bar X_0,\bar Z)$ of the pair $(\bar X_0,\bar Z)$,
with $\bar X\in\PSmVar(\mathbb C)$ and $\bar D=\cup_{i=1}^s\bar D_i\subset\bar X$ a normal crossing divisor.  
We denote by $i_{\bullet}:\bar D_{\bullet}\hookrightarrow\bar X=\bar X_{c(\bullet)}$ the morphism of simplicial varieties
given by the closed embeddings $i_I:\bar D_I=\cap_{i\in I}\bar D_i\hookrightarrow\bar X$.
Then the morphisms $\bar f:=\bar f_0\circ\bar\epsilon:\bar X\to\bar S$ and 
$\bar f_{D_{\bullet}}:=\bar f\circ i_{\bullet}:\bar D_{\bullet}\to\bar S$ are projective since
$\bar X$ and $\bar D_I$ are projective varieties.  
We then consider $(X,D):=\bar\epsilon^{-1}(X_0,Z)$, $\epsilon:=\bar\epsilon_{|X}:(X,D)\to(X_0,Z)$
We denote again by $i_{\bullet}:D_{\bullet}\hookrightarrow X=X_{c(\bullet)}$ the morphism of simplicial varieties
given by the closed embeddings $i_I:D_I=\cap_{i\in I}D_i\hookrightarrow X$.
Then the morphisms $f:=f_0\circ\epsilon:X\to S$ and 
$f_{D_{\bullet}}:=f\circ i_{\bullet}:D_{\bullet}\to S$ are projective since $\bar f:\bar X_0\to\bar S$ is projective.
\end{itemize}
\item[(ii)]Let $g:V'/S\to V/S$ a morphism, with $V'/S=(V',h'),V/S=(V,h)\in\Var(\mathbb C)/S$
\begin{itemize}
\item Take (see (i)) a compactification $\bar X_0\in\PVar(\mathbb C)$ of $V$ such that $h:V\to S$
extend to a morphism $\bar f_0=\bar{h}_0:\bar X_0\to\bar{S}$. Denote by $\bar Z=\bar X_0\backslash V$.
Then, there exist a compactification $\bar X'_0\in\PVar(\mathbb C)$ of $V'$ such that $h':V'\to S$
extend to a morphism $\bar f'_0=\bar{h}'_0:\bar X'_0\to\bar{S}$, $g:V'\to V$ extend to a morphism 
$\bar{g}_0:\bar X'_0\to\bar X_0$ and $\bar f_0\circ\bar{g}_0=\bar f'_0$ 
that is $\bar{g}_0$ is gives a morphism $\bar{g}_0:\bar X'_0/\bar{S}\to\bar X_0/\bar{S}$. 
Denote by $\bar Z'=\bar X'_0\backslash V'$. We then have the following commutative diagram
\begin{equation*}
\xymatrix{V\ar[r]^j & \bar X_0 & \, & \bar Z\ar[ll]_i \\
V'\ar[r]^{j'}\ar[u]^{\bar g} & \bar X'_0\ar[u]^{\bar{g}_0} & \bar Z'\ar[l]_{i'} & 
\bar{g_0}^{-1}(\bar Z)\ar[l]_{i''_{g,0}}\ar[u]^{\bar{g}'}:i'_{g,0}}
\end{equation*}
It gives the following commutative diagram
\begin{equation*}
\xymatrix{V\ar[r]^j & X_0:=\bar f_0^{-1}(S) & \, & Z\ar[ll]_i \\
V'\ar[r]^{j'}\ar[u]^g & X'_0:=\bar f_0^{'-1}(S)\ar[u]^{\bar{g}_0} & Z'\ar[l]_{i'} & 
\bar{g_0}^{-1}(Z)\ar[l]_{i''_{g,0}}\ar[u]^{\bar{g}'_0}:i'_{g,0}}
\end{equation*}
\item In the case $V$ and $V'$ are smooth, we take using theorem \ref{desVar} a strict desingularization 
$\bar\epsilon:(\bar X,\bar D)\to(\bar X_0,\bar Z)$ of $(\bar X_0,\bar Z)$.  
Then there exist a strict desingularization
$\bar\epsilon'_{\bullet}:(\bar X',\bar D')\to(\bar X'_0,\bar Z')$ of $(\bar X'_0,\bar Z')$
and a morphism $\bar g:\bar X'\to\bar X$ such that the following diagram commutes
\begin{equation*}
\xymatrix{\bar X'_0\ar[r]^{\bar{g}_0} & \bar X_0 \\
\bar X'\ar[u]^{\bar\epsilon'}\ar[r]^{\bar g} & \bar X\ar[u]^{\bar\epsilon}}.
\end{equation*}  
We then have the following commutative diagram in $\Fun(\Delta,\Var(\mathbb C))$
\begin{equation*}
\xymatrix{V=V_{c(\bullet)}\ar[r]^j & \bar X=\bar X_{c(\bullet)} & \, & \bar D_{s_g(\bullet)}\ar[ll]_{i_{\bullet}} \\
V'=V'_{c(\bullet)}\ar[r]^{j'}\ar[u]^g & \bar X'=\bar X'_{c(\bullet)}\ar[u]^{\bar{g}} & \bar D'_{\bullet}\ar[l]_{i'_{\bullet}} & 
\bar{g}^{-1}(\bar D_{s_g(\bullet)})\ar[l]_{i''_{g\bullet}}\ar[u]^{\bar{g}'_{\bullet}}:i'_{g\bullet}}
\end{equation*}
where $i_{\bullet}:\bar D_{\bullet}\hookrightarrow\bar X_{\bullet}$ the morphism of simplicial varieties
given by the closed embeddings $i_n:\bar D_n\hookrightarrow\bar X_n$,
and $i'_{\bullet}:\bar D'_{\bullet}\hookrightarrow\bar X'_{\bullet}$ the morphism of simplicial varieties
given by the closed embeddings $i'_n:\bar D'_n\hookrightarrow\bar X'_n$. 
It gives the commutative diagram in $\Fun(\Delta,\Var(\mathbb C))$
\begin{equation*}
\xymatrix{V=V_{c(\bullet)}\ar[r]^j & X:=\bar\epsilon^{-1}(X_0)=X_{c(\bullet)} & \, & D_{s_g(\bullet)}\ar[ll]_{i_{\bullet}} \\
V'=V'_{c(\bullet)}\ar[r]^{j'}\ar[u]^g & X':=\bar\epsilon^{',-1}(X'_0)=X'_{c(\bullet)}\ar[u]^{\bar{g}} & 
D'_{\bullet}\ar[l]_{i'_{\bullet}} & \bar{g}^{-1}(D_{s_g(\bullet)})\ar[l]_{i''_{g\bullet}}\ar[u]^{\bar{g}'_{\bullet}}:i'_{g\bullet}}
\end{equation*}
\end{itemize}
\end{itemize}
\end{defiprop}

\begin{proof}
\noindent(i): Let $\bar X_{00}\in\PVar(\mathbb C)$ be a compactification of $V$. 
Let $l_0:\bar X_0=\bar{\Gamma_h}\hookrightarrow\bar X_{00}\times\bar{S}$ be the closure of the graph of $h$ and
$\bar f_0:=p_{\bar{S}}\circ l_0:\bar X_0\hookrightarrow\bar X_{00}\times\bar{S}\to\bar{S}$, 
$\epsilon_{\bar X_0}:=p_{\bar X_{00}}\circ l_0:\bar X_0\hookrightarrow\bar X_{00}\times\bar{S}\to\bar X_{00}$
be the restriction to $\bar X_0$ of the projections.
Then, $\bar X\in\PVar(\mathbb C)$, $\epsilon_{\bar X_0}:\bar X_0\to\bar X_{00}$ is a proper modification which does not affect
the open subset $V\subset\bar X_0$, and $\bar f_0=\bar{h}_0:\bar X_0\to\bar{S}$ is a compactification of $h$.

\noindent(ii): There is two things to prove: 
\begin{itemize}
\item Let $\bar f_0:\bar X_0\to\bar S$ a compactification of $h:V\to S$ and 
$\bar f'_{00}:\bar X'_{00}\to\bar S$ a compactification of $h':V'\to S$ (see (i)). 
Let $l_0:\bar X'_0\hookrightarrow\bar{\Gamma_g}\subset\bar X'_{00}\times_{\bar S}\bar X_0$ be the closure of the graph of $g$,
$\bar f'_0:=(\bar f'_{00},\bar f_0)\circ l_0:\bar X'_0\hookrightarrow\bar X'_{00}\times_S\bar X_0\to\bar S$ and
$\bar{g}_0:=p_{\bar X_0}\circ l_0:\bar X'_0\hookrightarrow\bar X'_{00}\times_{\bar S}\bar X_0\to\bar X_0$, 
$\epsilon_{\bar X'_{00}}:=p_{\bar X'_0}\circ i:\bar X'_0\hookrightarrow\bar X'_{00}\times_{\bar S}\bar X_0\to\bar X'_{00}$
be the restriction to $X$ of the projections.
Then $\epsilon_{\bar X'_{00}}:\bar X'_0\to\bar X'_{00}$ is a proper modification which does not affect
the open subset $V'\subset\bar X'_0$, $\bar f'_0:\bar X'_0\to\bar S$ is an other compactification of $h':V'\to S$
and $\bar{g}_0:\bar X'_0\to\bar X_0$ is a compactification of $g$.
\item In the case $V$ and $V'$ are smooth, we take, using theorem \ref{desVar}, a strict desingularization 
$\bar\epsilon:(\bar X,\bar D)\to(\bar X_0,\bar Z)$ of the pair $(\bar X_0,\bar Z)$.
Take then, using theorem \ref{desVar}, a strict desingularization 
$\bar\epsilon'_1:(\bar X',\bar D')\to(\bar X\times_{\bar X_0}\bar X'_0,\bar X\times_{\bar X_0}\bar Z')$ 
of the pair $(\bar X\times_{\bar X_0}\bar X'_0,\bar X\times_{\bar X_0}\bar Z')$.
We consider then following commutative diagram whose square is cartesian :
\begin{equation*}
\xymatrix{ \, & \bar X'_0\ar[r]^{\bar g_0} &  X_0 \\
\, & \bar X\times_{\bar X_0}\bar X'_0\ar[u]^{\bar\epsilon'_0}\ar[r]^{\bar{g}'_0} & \bar X\ar[u]^{\epsilon} \\
\bar X'\ar[ruu]^{\epsilon'}\ar[rru]^{\bar{g}}\ar[ru]^{\bar\epsilon'_1} & \, & \,}
\end{equation*}  
and $\bar\epsilon':=\bar\epsilon'_0\circ\bar\epsilon'_1:(\bar X',\bar D')\to(\bar X'_0,\bar Z')$
is a strict desingularization of the pair $(\bar X\times_{\bar X_0}\bar X'_0, \bar X\times_{\bar X_0}\bar Z')$.
\end{itemize}
\end{proof}

Let $S\in\Var(\mathbb C)$. Recall we have the dual functor
\begin{eqnarray*}
\mathbb D_S:C(\Var(\mathbb C)/S)\to C(\Var(\mathbb C)/S), \; F\mapsto\mathbb D_S(F):=\mathcal Hom(F,E_{et}(\mathbb Z(S/S)))
\end{eqnarray*}
which induces the functor
\begin{eqnarray*}
L\mathbb D_S:C(\Var(\mathbb C)/S)\to C(\Var(\mathbb C)/S), \; 
F\mapsto L\mathbb D_S(F):=\mathbb D_S(LF):=\mathcal Hom(LF,E_{et}(\mathbb Z(S/S))).
\end{eqnarray*}

We will use the following resolutions of representable presheaves by Corti-Hanamura presheaves and
their the functorialities.

\begin{defi}\label{RCHdef} 
\begin{itemize}
\item[(i)]Let $h:U\to S$ a morphism, with $U,S\in\Var(\mathbb C)$ and $U$ smooth.
Take, see definition-proposition \ref{RCHdef0},
$\bar f_0=\bar{h}_0:\bar X_0\to\bar S$ a compactification of $h:U\to S$ and denote by $\bar Z=\bar X_0\backslash U$.
Take, using theorem \ref{desVar}(ii), a strict desingularization 
$\bar\epsilon:(\bar X,\bar D)\to(\bar X_0,\bar Z)$ of the pair $(\bar X_0,\bar Z)$, with $\bar X\in\PSmVar(\mathbb C)$ and 
$\bar D:=\epsilon^{-1}(\bar Z)=\cup_{i=1}^s\bar D_i\subset\bar X$ a normal crossing divisor.  
We denote by $i_{\bullet}:\bar D_{\bullet}\hookrightarrow\bar X=\bar X_{c(\bullet)}$ the morphism of simplicial varieties
given by the closed embeddings $i_I:\bar D_I=\cap_{i\in I}\bar D_i\hookrightarrow\bar X$
We denote by $j:U\hookrightarrow\bar X$ the open embedding and by $p_S:\bar X\times S\to S$ 
and $p_S:U\times S\to S$ the projections.
Considering the graph factorization $\bar f:\bar X\xrightarrow{\bar l}\bar X\times\bar S\xrightarrow{p_{\bar S}}\bar S$
of $\bar f:\bar X\to\bar S$, where $\bar l$ is the graph embedding and $p_{\bar S}$ the projection,
we get closed embeddings $l:=\bar l\times_{\bar S}S:X\hookrightarrow\bar X\times S$ and 
$l_{D_I}:=\bar D_I\times_{\bar X} l:D_I\hookrightarrow\bar D_I\times S$.
We then consider the following map in $C(\Var(\mathbb C)^2/S)$
\begin{eqnarray*}
r_{(\bar X,\bar D)/S}(\mathbb Z(U/S)):R_{(\bar X,\bar D)/S}(\mathbb Z(U/S)) \\
\xrightarrow{:=}p_{S*}E_{et}(\Cone(\mathbb Z(i_{\bullet}\times I):
(\mathbb Z((\bar D_{\bullet}\times S,D_{\bullet})/\bar X\times S),u_{IJ})\to
\mathbb Z((\bar X\times S,X)/\bar X\times S))) \\
\xrightarrow{p_{S*}E_{et}(0,k\circ\ad((j\times I)^*,(j\times I)_*)(\mathbb Z((\bar X\times S,X)/\bar X\times S)))}
p_{S*}E_{et}(\mathbb Z((U\times S,U)/U\times S)))=:\mathbb D_S^{12}(\mathbb Z(U/S)).
\end{eqnarray*}
Note that $\mathbb Z((\bar D_I\times S,D_I)/\bar X\times S)$ and $\mathbb Z((\bar X\times S,X)/\bar X\times S))$
are obviously $\mathbb A^1$ invariant.
Note that $r_{(X,D)/S}$ is NOT an equivalence $(\mathbb A^1,et)$ local by proposition \ref{rho12} since
$\rho_{\bar X\times S*}\mathbb Z((\bar D_{\bullet}\times S,D_{\bullet})/\bar X\times S)=0$, whereas 
$\rho_{\bar X\times S*}\ad((j\times I)^*,(j\times I)_*)(\mathbb Z((\bar X\times S,X)/\bar X\times S)))$ 
is not an equivalence $(\mathbb A^1,et)$ local.
\item[(ii)]Let $g:U'/S\to U/S$ a morphism, with $U'/S=(U',h'),U/S=(U,h)\in\Var(\mathbb C)/S$, with $U$ and $U'$ smooth.
Take, see definition-proposition \ref{RCHdef0}(ii),a compactification $\bar f_0=\bar h:\bar X_0\to\bar S$ of $h:U\to S$ 
and a compactification $\bar f'_0=\bar{h}':\bar X'_0\to\bar S$ of $h':U'\to S$ such that 
$g:U'/S\to U/S$ extend to a morphism $\bar g_0:\bar X'_0/\bar S\to\bar X_0/\bar S$. 
Denote $\bar Z=\bar X_0\backslash U$ and $\bar Z'=\bar X'_0\backslash U'$.
Take, see definition-proposition \ref{RCHdef0}(ii), a strict desingularization 
$\bar\epsilon:(\bar X,\bar D)\to(\bar X_0,\bar Z)$ of $(\bar X_0,\bar Z)$,
a strict desingularization $\bar\epsilon'_{\bullet}:(\bar X',\bar D')\to(\bar X'_0,\bar Z')$ of $(\bar X'_0,\bar Z')$
and a morphism $\bar g:\bar X'\to\bar X$ such that the following diagram commutes
\begin{equation*}
\xymatrix{\bar X'_0\ar[r]^{\bar{g}_0} & \bar X_0 \\
\bar X'\ar[u]^{\bar\epsilon'}\ar[r]^{\bar g} & \bar X\ar[u]^{\bar\epsilon}}.
\end{equation*} 
We then have, see definition-proposition \ref{RCHdef0}(ii),
the following commutative diagram in $\Fun(\Delta,\Var(\mathbb C))$
\begin{equation}\label{RCHdia}
\xymatrix{U=U_{c(\bullet)}\ar[r]^j & \bar X=\bar X_{c(\bullet)} & \, & \bar D_{s_g(\bullet)}\ar[ll]_{i_{\bullet}} \\
U'=U'_{c(\bullet)}\ar[r]^{j'}\ar[u]^g & \bar X'=\bar X'_{c(\bullet)}\ar[u]^{\bar{g}} & \bar D'_{\bullet}\ar[l]_{i'_{\bullet}} & 
\bar{g}^{-1}(\bar D_{s_g(\bullet)})\ar[l]_{i''_{g\bullet}}\ar[u]^{\bar{g}'_{\bullet}}:i'_{g\bullet}}
\end{equation}
Denote by $p_S:\bar X\times S\to S$ and $p'_S:\bar X'\times S\to S$ the projections
We then consider the following map in $C(\Var(\mathbb C)^2/S)$
\begin{eqnarray*}
R_S^{CH}(g):R_{(\bar X,\bar D)/S}(\mathbb Z(U/S))\xrightarrow{:=} \\
p_{S*}E_{et}(\Cone(\mathbb Z(i_{\bullet}\times I):
(\mathbb Z((\bar D_{s_g(\bullet)}\times S,D_{s_g(\bullet)})/\bar X\times S),u_{IJ})\to
\mathbb Z((\bar X\times S,X)/\bar X\times S))) \\
\xrightarrow{T((\bar g\times I),E)(-)\circ p_{S*}\ad((\bar{g}\times I)^*,(\bar{g}\times I)_*)(-)} \\
p'_{S*}E_{et}(\Cone(\mathbb Z(i'_{g\bullet}\times I): \\
(\mathbb Z((\bar{g}^{-1}(\bar D_{s_g(\bullet)})\times S,\bar g^{-1}(D_{s_g(\bullet)})/\bar X'\times S),u_{IJ})
\to\mathbb Z((\bar X'\times S,X')/\bar X'\times S)))) \\
\xrightarrow{p'_{S*}E_{et}(\mathbb Z(i''_{g\bullet}\times I),I)} \\
p'_{S*}E_{et}(\Cone(\mathbb Z(i'_{\bullet}\times I):
((\mathbb Z((\bar D'_{\bullet}\times S,D'_{\bullet)})/\bar X'\times S),u_{IJ})
\to\mathbb Z((\bar X'\times S,X')/\bar X'\times S))) \\
\xrightarrow{=:}R_{(\bar X',\bar D')/S}(\mathbb Z(U'/S))
\end{eqnarray*}
Then by the diagram (\ref{RCHdia}) and adjonction, the following diagram in $C(\Var(\mathbb C)^2/S)$ obviously commutes
\begin{equation*}
\xymatrix{R_{(\bar X,\bar D)/S}(\mathbb Z(U/S))\ar[rr]^{r_{(\bar X,\bar D)/S}(\mathbb Z(U/S))}\ar[d]_{R_S^{CH}(g)} & \, & 
p_{S*}E_{et}(\mathbb Z((U\times S,U)/U\times S))=:\mathbb D_S^{12}(\mathbb Z(U/S))
\ar[d]^{D^{12}_S(g):=T(g\times I,E)(-)\circ\ad((g\times I)^*,(g\times I)_*)(E_{et}(\mathbb Z((U\times S,U)/U\times S)))} \\
R_{(\bar X',\bar D')/S}(\mathbb Z(U'/S))\ar[rr]^{r_{(\bar X',\bar D')/S}(\mathbb Z(U'/S))} & \, & 
p'_{S*}E_{et}(\mathbb Z((U'\times S,U')/U'\times S))=:\mathbb D^{12}_S(\mathbb Z(U'/S))}.
\end{equation*}
\item[(iii)]For $g_1:U''/S\to U'/S$, $g_2:U'/S\to U/S$ two morphisms
with $U''/S=(U',h''),U'/S=(U',h'),U/S=(U,h)\in\Var(\mathbb C)/S$, with $U$, $U'$ and $U''$ smooth.
We get from (i) and (ii) 
a compactification $\bar f=\bar{h}:\bar X\to\bar S$ of $h:U\to S$,
a compactification $\bar f'=\bar{h}':\bar X'\to\bar S$ of $h':U'\to S$,
and a compactification $\bar f''=\bar{h}'':\bar X''\to\bar S$ of $h'':U''\to S$,
with $\bar X,\bar X',\bar X''\in\PSmVar(\mathbb C)$, 
$\bar D:=\bar X\backslash U\subset\bar X$ $\bar D':=\bar X'\backslash U'\subset\bar X'$, 
and $\bar D'':=\bar X''\backslash U''\subset\bar X''$ normal crossing divisors, 
such that $g_1:U''/S\to U'/S$ extend to $\bar g_1:\bar X''/\bar S\to\bar X'/\bar S$,
$g_2:U'/S\to U/S$ extend to $\bar g_2:\bar X'/\bar S\to\bar X/\bar S$, and
\begin{eqnarray*}
R_S^{CH}(g_2\circ g_1)=R_S^{CH}(g_1)\circ R_S^{CH}(g_2):R_{(\bar X,\bar D)/S}\to R_{(\bar X'',\bar D'')/S} 
\end{eqnarray*}
\item[(iv)] For  
\begin{eqnarray*}
Q^*:=(\cdots\to\oplus_{\alpha\in\Lambda^n}\mathbb Z(U^n_{\alpha}/S)
\xrightarrow{(\mathbb Z(g^n_{\alpha,\beta}))}\oplus_{\beta\in\Lambda^{n-1}}\mathbb Z(U^{n-1}_{\beta}/S)\to\cdots)
\in C(\Var(\mathbb C)/S)
\end{eqnarray*}
a complex of (maybe infinite) direct sum of representable presheaves with $U^*_{\alpha}$ smooth,
we get from (i), (ii) and (iii) the map in $C(\Var(\mathbb C)^2/S)$
\begin{eqnarray*}
r_S^{CH}(Q^*):R^{CH}(Q^*):=
(\cdots\to\oplus_{\beta\in\Lambda^{n-1}}\varinjlim_{(\bar X^{n-1}_{\beta},\bar D^{n-1}_{\beta})/S}
R_{(\bar X^{n-1}_{\beta},\bar D^{n-1}_{\beta})/S}(\mathbb Z(U^{n-1}_{\beta}/S)) \\
\xrightarrow{(R_S^{CH}(g^n_{\alpha,\beta}))}\oplus_{\alpha\in\Lambda^n}\varinjlim_{(\bar X^n_{\alpha},\bar D^n_{\alpha})/S}
R_{(\bar X^n_{\alpha},\bar D^n_{\alpha})/S}(\mathbb Z(U^n_{\alpha}/S))\to\cdots)
\to\mathbb D^{12}_S(Q^*),
\end{eqnarray*}
where for $(U^n_{\alpha},h^n_{\alpha})\in\Var(\mathbb C)/S$, the inductive limit run over all the compactifications 
$\bar f_{\alpha}:\bar X_{\alpha}\to\bar S$ of $h_{\alpha}:U_{\alpha}\to S$ with $\bar X_{\alpha}\in\PSmVar(\mathbb C)$
and $\bar D_{\alpha}:=\bar X_{\alpha}\backslash U_{\alpha}$ a normal crossing divisor.
For $m=(m^*):Q_1^*\to Q_2^*$ a morphism with 
\begin{eqnarray*}
Q_1^*:=(\cdots\to\oplus_{\alpha\in\Lambda^n}\mathbb Z(U^n_{1,\alpha}/S)
\xrightarrow{(\mathbb Z(g^n_{\alpha,\beta}))}\oplus_{\beta\in\Lambda^{n-1}}\mathbb Z(U^{n-1}_{1,\beta}/S)\to\cdots), \\
Q_2^*:=(\cdots\to\oplus_{\alpha\in\Lambda^n}\mathbb Z(U^n_{2,\alpha}/S)
\xrightarrow{(\mathbb Z(g^n_{\alpha,\beta}))}\oplus_{\beta\in\Lambda^{n-1}}\mathbb Z(U^{n-1}_{2,\beta}/S)\to\cdots)
\in C(\Var(\mathbb C)/S)
\end{eqnarray*}
complexes of (maybe infinite) direct sum of representable presheaves with $U^*_{1,\alpha}$ and $U^*_{2,\alpha}$ smooth,
we get again from (i), (ii) and (iii) a commutative diagram in $C(\Var(\mathbb C)^2/S)$
\begin{equation*}
\xymatrix{R^{CH}(Q_2^*)\ar[rr]^{r_S^{CH}(Q_2^*)}\ar[d]_{R_S^{CH}(m):=(R_S^{CH}(m^*))} & \, & 
\mathbb D^{12}_S(Q_2^*)\ar[d]^{\mathbb D^{12}_S(m):=(\mathbb D^{12}_S(m^*))} \\
R^{CH}(Q_1^*)\ar[rr]^{r_S^{CH}(Q_1^*)} & \, & \mathbb D^{12}_S(Q_1^*)}.
\end{equation*}
\end{itemize}
\end{defi}

\begin{itemize}
\item Let $S\in\Var(\mathbb C)$
For $(h,m,m')=(h^*,m^*,m^{'*}):Q_1^*[1]\to Q_2^*$ an homotopy with $Q_1^*,Q_2^*\in C(\Var(\mathbb C)/S)$
complexes of (maybe infinite) direct sum of representable presheaves with $U^*_{1,\alpha}$ and $U^*_{2,\alpha}$ smooth,
\begin{equation*}
(R_S^{CH}(h),R_S^{CH}(m),R_S^{CH}(m'))=(R_S^{CH}(h^*),R_S^{CH}(m^*),R_S^{CH}(m^{'*})):
R^{CH}(Q_2^*)[1]\to R^{CH}(Q_1^*)
\end{equation*}
is an homotopy in $C(\Var(\mathbb C)^2/S)$ using definition \ref{RCHdef} (iii). 
In particular if $m:Q_1^*\to Q_2^*$ with $Q_1^*,Q_2^*\in C(\Var(\mathbb C)/S)$
complexes of (maybe infinite) direct sum of representable presheaves with $U^*_{1,\alpha}$ and $U^*_{2,\alpha}$ smooth
is an homotopy equivalence, then $R_S^{CH}(m):R^{CH}(Q_2^*)\to R^{CH}(Q_1^*)$ is an homotopy equivalence.
\item Let $S\in\SmVar(\mathbb C)$. Let $F\in\PSh(\Var(\mathbb C)^{sm}/S)$. Consider 
\begin{eqnarray*}
q:LF:=(\cdots\to\oplus_{(U_{\alpha},h_{\alpha})\in\Var(\mathbb C)^{sm}/S}\mathbb Z(U_{\alpha}/S)
\xrightarrow{(\mathbb Z(g^n_{\alpha,\beta}))}
\oplus_{(U_{\alpha},h_{\alpha})\in\Var(\mathbb C)^{sm}/S}\mathbb Z(U_{\alpha}/S)\to\cdots)\to F
\end{eqnarray*}
the canonical projective resolution given in subsection 2.3.3.
Note that the $U_{\alpha}$ are smooth since $S$ is smooth and $h_{\alpha}$ are smooth morphism.
Definition \ref{RCHdef}(iv) gives in this particular case the map in $C(\Var(\mathbb C)^2/S)$
\begin{eqnarray*}
r_S^{CH}(\rho_S^*LF):R^{CH}(\rho_S^*LF):=
(\cdots\to\oplus_{(U_{\alpha},h_{\alpha})\in\Var(\mathbb C)^{sm}/S}\varinjlim_{(\bar X_{\alpha},\bar D_{\alpha})/S}
R_{(\bar X_{\alpha},\bar D_{\alpha})/S}(\mathbb Z(U_{\alpha}/S)) \\
\xrightarrow{(R_S^{CH}(g^n_{\alpha,\beta}))}
\oplus_{(U_{\alpha},h_{\alpha})\in\Var(\mathbb C)^{sm}/S}\varinjlim_{(\bar X_{\alpha},\bar D_{\alpha})/S}
R_{(\bar X_{\alpha},\bar D_{\alpha})/S}(\mathbb Z(U_{\alpha}/S))\to\cdots)
\to\mathbb D^{12}_S(\rho_S^*LF),
\end{eqnarray*}
where for $(U_{\alpha},h_{\alpha})\in\Var(\mathbb C)^{sm}/S$, the inductive limit run over all the compactifications 
$\bar f_{\alpha}:\bar X_{\alpha}\to\bar S$ of $h_{\alpha}:U_{\alpha}\to S$ with $\bar X_{\alpha}\in\PSmVar(\mathbb C)$
and $\bar D_{\alpha}:=\bar X_{\alpha}\backslash U_{\alpha}$ a normal crossing divisor.
Definition \ref{RCHdef}(iv) gives then by functoriality in particular, for $F=F^{\bullet}\in C(\Var(\mathbb C)^{sm}/S)$, 
the map in $C(\Var(\mathbb C)^2/S)$
\begin{eqnarray*}
r_S^{CH}(\rho_S^*LF)=(r_S^{CH}(\rho_S^*LF^*)):R^{CH}(\rho_S^*LF)\to\mathbb D^{12}_S(\rho_S^*LF).
\end{eqnarray*}

\item Let $g:T\to S$ a morphism with $T,S\in\SmVar(\mathbb C)$. Let $h:U\to S$ a smooth morphism with $U\in\Var(\mathbb C)$.
Consider the cartesian square
\begin{equation*}
\xymatrix{U_T\ar[r]^{h'}\ar[d]^{g'} & T\ar[d]^g \\
U\ar[r]^h & S}
\end{equation*}
Note that $U$ is smooth since $S$ and $h$ are smooth, and $U_T$ is smooth since $T$ and $h'$ are smooth. 
Take, see definition-proposition \ref{RCHdef0}(ii),a compactification $\bar f_0=\bar h:\bar X_0\to\bar S$ of $h:U\to S$ 
and a compactification $\bar f'_0=\bar{g\circ h'}:\bar X'_0\to\bar S$ of $g\circ h':U'\to S$ such that 
$g':U_T/S\to U/S$ extend to a morphism $\bar{g}'_0:\bar X'_0/\bar S\to \bar X_0/\bar S$. 
Denote $\bar Z=\bar X_0\backslash U$ and $\bar Z'=\bar X'_0\backslash U_T$.
Take, see definition-proposition \ref{RCHdef0}(ii), a strict desingularization 
$\bar\epsilon:(\bar X,\bar D)\to(\bar X_0,\bar Z)$ of $(\bar X_0,\bar Z)$,
a desingularization $\bar\epsilon'_{\bullet}:(\bar X',\bar D')\to(\bar X'_0,\bar Z')$ of $(\bar X'_0,\bar Z')$
and a morphism $\bar g':\bar X'\to\bar X$ such that the following diagram commutes
\begin{equation*}
\xymatrix{\bar X'_0\ar[r]^{\bar{g}'_0} & \bar X_0 \\
\bar X'\ar[u]^{\bar\epsilon'}\ar[r]^{\bar{g}'} & \bar X\ar[u]^{\bar\epsilon}}.
\end{equation*} 
We then have, see definition-proposition \ref{RCHdef0}(ii),
the following commutative diagram in $\Fun(\Delta,\Var(\mathbb C))$
\begin{equation*}
\xymatrix{U=U_{c(\bullet)}\ar[r]^j & \bar X=\bar X_{c(\bullet)} & \, & \bar D_{s_{g'}(\bullet)}\ar[ll]_{i_{\bullet}} \\
U_T=U_{T,c(\bullet)}\ar[r]^{j'}\ar[u]^{g'} & \bar X'=X'_{c(\bullet)}\ar[u]^{\bar{g}'} & \bar D'_{\bullet}\ar[l]_{i'_{\bullet}} & 
\bar{g}^{'-1}(\bar D_{s_{g'}(\bullet)})\ar[l]_{i''_{g'\bullet}}\ar[u]^{(\bar{g}')'_{\bullet}}:i'_{g\bullet}}
\end{equation*}
We then consider the following map in $C(\Var(\mathbb C)^2/T)$, see definition \ref{RCHdef}(ii)
\begin{eqnarray*}
T(g,R^{CH})(\mathbb Z(U/S)):g^*R_{(\bar X,\bar D)/S}(\mathbb Z(U/S)) \\ 
\xrightarrow{g^*R^{CH}_S(g')}g^*R_{(\bar X',\bar D')/S}(\mathbb Z(U_T/S))=g^*g_*R_{(\bar X',\bar D')/T}(\mathbb Z(U_T/T)) \\
\xrightarrow{\ad(g^*,g_*)(R_{(\bar X',\bar D')/T}(\mathbb Z(U_T/T)))}R_{(\bar X',\bar D')/T}(\mathbb Z(U_T/T))
\end{eqnarray*}
For  
\begin{eqnarray*}
Q^*:=(\cdots\to\oplus_{\alpha\in\Lambda^n}\mathbb Z(U^n_{\alpha}/S)
\xrightarrow{(\mathbb Z(g^n_{\alpha,\beta}))}\oplus_{\beta\in\Lambda^{n-1}}\mathbb Z(U^{n-1}_{\beta}/S)\to\cdots)
\in C(\Var(\mathbb C)/S)
\end{eqnarray*}
a complex of (maybe infinite) direct sum of representable presheaves with $h^n_{\alpha}:U^n_{\alpha}\to S$ smooth,
we get the map in $C(\Var(\mathbb C)^2/T)$
\begin{eqnarray*}
T(g,R^{CH})(Q^*):g^*R^{CH}(Q^*)=
(\cdots\to\oplus_{\alpha\in\Lambda^n}\varinjlim_{(\bar X^n_{\alpha},\bar D^n_{\alpha})/S}
g^*R_{(\bar X^n_{\alpha},\bar D^n_{\alpha})/S}(\mathbb Z(U^n_{\alpha}/S))\to\cdots) \\
\xrightarrow{(T(g,R^{CH})(\mathbb Z(U^n_{\alpha}/S)))} 
(\cdots\to\oplus_{\alpha\in\Lambda^n}\varinjlim_{(\bar X^{n'}_{\alpha},\bar D^{n'}_{\alpha})/T}
R_{(\bar X^{n'}_{\alpha},\bar D^{n'}_{\alpha})/T}(\mathbb Z(U^n_{\alpha,T}/S))\to\cdots)=:R^{CH}(g^*Q^*).
\end{eqnarray*}
Let $F\in\PSh(\Var(\mathbb C)^{sm}/S)$. Consider 
\begin{eqnarray*}
q:LF:=(\cdots\to\oplus_{(U_{\alpha},h_{\alpha})\in\Var(\mathbb C)^{sm}/S}\mathbb Z(U_{\alpha}/S)\to\cdots)\to F
\end{eqnarray*}
the canonical projective resolution given in subsection 2.3.3.
We then get in particular the map in $C(\Var(\mathbb C)^2/T)$
\begin{eqnarray*}
T(g,R^{CH})(\rho_S^*LF):g^*R^{CH}(\rho_S^*LF)= \\
(\cdots\to\oplus_{(U_{\alpha},h_{\alpha})\in\Var(\mathbb C)^{sm}/S}\varinjlim_{(\bar X_{\alpha},\bar D_{\alpha})/S}
g^*R_{(\bar X_{\alpha},\bar D_{\alpha})/S}(\mathbb Z(U_{\alpha}/S))\to\cdots) 
\xrightarrow{(T(g,R^{CH})(\mathbb Z(U_{\alpha}/S)))} \\ 
(\cdots\to\oplus_{(U_{\alpha},h_{\alpha})\in\Var(\mathbb C)^{sm}/S}\varinjlim_{(\bar X'_{\alpha},\bar D'_{\alpha})/T}
R_{(\bar X'_{\alpha},\bar D'_{\alpha})/T}(\mathbb Z(U_{\alpha,T}/S))\to\cdots)=:R^{CH}(\rho_T^*g^*LF). 
\end{eqnarray*}
By functoriality, we get in particular for $F=F^{\bullet}\in C(\Var(\mathbb C)^{sm}/S)$, the map in $C(\Var(\mathbb C)^2/T)$
\begin{eqnarray*}
T(g,R^{CH})(\rho_S^*LF):g^*R^{CH}(\rho_S^*LF)\to R^{CH}(\rho_T^*g^*LF).
\end{eqnarray*}

\item Let $S_1,S_2\in\SmVar(\mathbb C)$ and $p:S_1\times S_2\to S_1$ the projection. 
Let $h:U\to S_1$ a smooth morphism with $U\in\Var(\mathbb C)$. Consider the cartesian square
\begin{equation*}
\xymatrix{U\times S_2\ar[r]^{h\times I}\ar[d]^{p'} & S_1\times S_2\ar[d]^p \\
U\ar[r]^h & S_1}
\end{equation*}
Take, see definition-proposition \ref{RCHdef0}(i),a compactification $\bar f_0=\bar h:\bar X_0\to\bar S_1$ of $h:U\to S_1$.
Then $\bar f_0\times I:\bar X_0\times S_2\to\bar S_1\times S_2$ is a compactification of $h\times I:U\times S_2\to S_1\times S_2$ 
and $p':U\times S_2\to U$ extend to $\bar{p}'_0:=p_{X_0}:\bar X_0\times S_2\to\bar X_0$. Denote $Z=X_0\backslash U$. 
Take see theorem \ref{desVar}(i), a strict desingularization 
$\bar\epsilon:(\bar X,\bar D)\to(\bar X_0,\bar Z)$ of the pair $(\bar X_0,\bar Z)$.  
We then have the following commutative diagram in $\Fun(\Delta,\Var(\mathbb C))$ whose squares are cartesian
\begin{equation}\label{RCHdiap}
\xymatrix{U=U_{c(\bullet)}\ar[r]^j & \bar X & \bar D_{\bullet}\ar[l]_{i_{\bullet}} \\
U\times S_2=(U\times S_2)_{c(\bullet)}\ar[r]^{j\times I}\ar[u]^g & \bar X\times S_2\ar[u]^{\bar{p}':=p_{\bar X}} & 
\bar D_{\bullet}\times S_2\ar[l]_{i'_{\bullet}}\ar[u]^{\bar{p'}'_{\bullet}}}
\end{equation}
Then the map in $C(\Var(\mathbb C)^2/S_1\times S_2)$
\begin{eqnarray*}
T(p,R^{CH})(\mathbb Z(U/S_1)):p^*R_{(\bar X,\bar D)/S_1}(\mathbb Z(U/S_1))\xrightarrow{\sim} 
R_{(\bar X\times S_2,\bar D_{\bullet}\times S_2)/S_1\times S_2}(\mathbb Z(U\times S_2/S_1\times S_2))
\end{eqnarray*}
is an isomorphism.
Hence, for $Q^*\in C(\Var(\mathbb C)/S_1)$ a complex of (maybe infinite) direct sum of representable presheaves of smooth morphism,
the map in $C(\Var(\mathbb C)^2/S_1\times S_2)$
\begin{eqnarray*}
T(p,R^{CH})(Q^*):p^*R^{CH}(Q^*)\xrightarrow{\sim}R^{CH}(p^*Q^*)
\end{eqnarray*}
is an isomorphism.
In particular, for $F\in C(\Var(\mathbb C)^{sm}/S_1)$ the map in $C(\Var(\mathbb C)^2/S_1\times S_2)$
\begin{eqnarray*}
T(p,R^{CH})(\rho_{S_1}^*LF):p^*R^{CH}(\rho_{S_1}^*LF)\xrightarrow{\sim}R^{CH}(\rho_{S_1\times S_2}^*p^*LF)
\end{eqnarray*}
is an isomorphism.

\item Let $h_1:U_1\to S$, $h_2:U_2\to S$ two morphisms with $U_1,U_2,S\in\Var(\mathbb C)$, $U_1,U_2$ smooth.
Denote by $p_1:U_1\times_SU_2\to U_1$ and $p_2:U_1\times_S U_2\to U_2$ the projections.
Take, see definition-proposition \ref{RCHdef0}(i)),
a compactification $\bar f_{10}=\bar{h}_1:\bar X_{10}\to\bar S$ of $h_1:U_1\to S$
and a compactification $\bar f_{20}=\bar{h}_2:\bar X_{20}\to\bar S$ of $h_2:U_2\to S$. Then, 
\begin{itemize}
\item $\bar f_{10}\times\bar f_{20}:\bar X_{10}\times_{\bar S}\bar X_{20}\to S$ 
is a compactification of $h_1\times h_2:U_1\times_SU_2\to S$.
\item $\bar p_{10}:=p_{X_{10}}:\bar X_{10}\times_{\bar S} \bar X_{20}\to\bar X_{10}$ 
is a compactification of $p_1:U_1\times_SU_2\to U_1$.
\item $\bar p_{20}:=p_{X_{20}}:\bar X_{10}\times_{\bar S}\bar X_{20}\to\bar X_{20}$ 
is a compactification of $p_2:U_1\times_SU_2\to U_2$.
\end{itemize}
Denote $\bar Z_1=\bar X_{10}\backslash U_1$ and $\bar Z_2=\bar X_{20}\backslash U_2$.
Take, see theorem \ref{desVar}(i), a strict desingularization 
$\bar\epsilon_1:(\bar X_1,\bar D)\to(\bar X_{10},Z_1)$ of the pair $(\bar X_{10},\bar Z_1)$
and a  strictdesingularization 
$\bar\epsilon_2:(\bar X_2,\bar E)\to(\bar X_{20},Z_2)$ of the pair $(\bar X_{20},\bar Z_2)$. 
Take then a strict desingularization
\begin{equation*}
\bar\epsilon_{12}:((\bar X_1\times_{\bar S}\bar X_2)^N,\bar F)\to
(\bar X_1\times_{\bar S}\bar X_2,(D\times_{\bar S}\bar X_2)\cup(\bar X_1\times_{\bar S}\bar E)) 
\end{equation*}
of the pair $(\bar X_1\times_{\bar S}\bar X_2,(\bar D\times_{\bar S}\bar X_2)\cup(\bar X_1\times_{\bar S}\bar E))$. 
We have then the following commutative diagram 
\begin{equation*}
\xymatrix{\, & \bar X_1\ar[r]^{\bar f_1} & \bar S \\
\, & \bar X_1\times_{\bar S}\bar X_2\ar[r]^{\bar p_1}\ar[u]^{\bar p_2} &\bar X_2\ar[u]^{\bar f_2} \\
(\bar X_1\times_{\bar S}\bar X_2)^N\ar[ru]^{\bar\epsilon_{12}}\ar[rru]^{(\bar p_1)^N}\ar[ruu]^{(\bar p_2)^N} & \, & \,}
\end{equation*}
and
\begin{itemize}
\item $\bar f_1\times\bar f_2:\bar X_1\times_{\bar S}\bar X_2\to\bar S$ 
is a compactification of $h_1\times h_2:U_1\times_SU_2\to S$.
\item $(\bar p_1)^N:=\bar p_1\circ\epsilon_{12}:(\bar X_1\times_{\bar S}\bar X_2)^N\to\bar X_1$ is a compactification of 
$p_1:U_1\times_SU_2\to U_1$.
\item $(\bar p_2)^N:=\bar p_2\circ\epsilon_{12}:(\bar X_1\times_{\bar S}\bar X_2)^N\to\bar X_2$ is a compactification of 
$p_2:U_1\times_SU_2\to U_2$.
\end{itemize}
We have then the morphism in $C(\Var(\mathbb C)^2/S)$
\begin{eqnarray*}
T(\otimes,R_S^{CH})(\mathbb Z(U_1/S),\mathbb Z(U_2/S)):=R_S^{CH}(p_1)\otimes R_S^{CH}(p_2): \\
R_{(\bar X_1,\bar D)/S}(\mathbb Z(U_1/S))\otimes R_{(X_2,E))/S}(\mathbb Z(U_2/S)) 
\xrightarrow{\sim} R_{(\bar X_1\times_{\bar S}\bar X_2)^N,\bar F)/S}(\mathbb Z(U_1\times_S U_2/S))
\end{eqnarray*}
For   
\begin{eqnarray*}
Q_1^*:=(\cdots\to\oplus_{\alpha\in\Lambda^n}\mathbb Z(U^n_{1,\alpha}/S)
\xrightarrow{(\mathbb Z(g^n_{\alpha,\beta}))}\oplus_{\beta\in\Lambda^{n-1}}\mathbb Z(U^{n-1}_{1,\beta}/S)\to\cdots), \\
Q_2^*:=(\cdots\to\oplus_{\alpha\in\Lambda^n}\mathbb Z(U^n_{2,\alpha}/S)
\xrightarrow{(\mathbb Z(g^n_{\alpha,\beta}))}\oplus_{\beta\in\Lambda^{n-1}}\mathbb Z(U^{n-1}_{2,\beta}/S)\to\cdots)
\in C(\Var(\mathbb C)/S)
\end{eqnarray*}
complexes of (maybe infinite) direct sum of representable presheaves with $U^*_{\alpha}$ smooth, 
we get the morphism in $C(\Var(\mathbb C)^2/S)$
\begin{eqnarray*}
T(\otimes,R_S^{CH})(Q^*_1,Q^*_2): R^{CH}(Q^*_1)\otimes R^{CH}(Q^*_2) 
\xrightarrow{(T(\otimes,R_S^{CH})(\mathbb Z(U_{1,\alpha}^m),\mathbb Z(U_{2,\beta}^n))} R^{CH}(Q^*_1\otimes Q^*_2)).
\end{eqnarray*}
For $F_1,F_2\in C(\Var(\mathbb C)^{sm}/S)$, we get in particular the morphism in $C(\Var(\mathbb C)^2/S)$
\begin{eqnarray*}
T(\otimes,R_S^{CH})(\rho_S^*LF_1,\rho_S^*LF_2):R^{CH}(\rho_S^*LF_1)\otimes R^{CH}(\rho_S^*LF_2) 
\to R^{CH}(\rho_S^*(LF_1\otimes LF_2)).
\end{eqnarray*}

\end{itemize}

\begin{defi}\label{sharpstar}
Let $h:U\to S$ a morphism, with $U,S\in\Var(\mathbb C)$, $U$ irreducible.
Take, see definition-proposition \ref{RCHdef0},
$\bar f_0=\bar{h}_0:\bar X_0\to\bar S$ a compactification of $h:U\to S$ and denote by $\bar Z=\bar X_0\backslash U$.
Take, using theorem \ref{desVar}, a desingularization 
$\bar\epsilon:(\bar X,\bar D)\to(\bar X_0,\Delta)$ of the pair $(\bar X_0,\Delta)$, $\bar Z\subset\Delta$,
with $\bar X\in\PSmVar(\mathbb C)$ and 
$\bar D:=\bar\epsilon^{-1}(\Delta)=\cup_{i=1}^s\bar D_i\subset\bar X$ a normal crossing divisor.
Denote $d_X:=\dim(\bar X)=\dim(U)$.
\begin{itemize}
\item[(i)] The cycle $(\Delta_{\bar D_{\bullet}}\times S)\subset\bar D_{\bullet}\times\bar D_{\bullet}\times S$
induces by the diagonal $\Delta_{\bar D_{\bullet}}\subset\bar D_{\bullet}\times\bar D_{\bullet}$
gives the morphism in $C(\Var(\mathbb C)^2/S)$
\begin{eqnarray*}
[\Delta_{\bar D_{\bullet}}]\in\Hom(\mathbb Z^{tr}((\bar D_{\bullet}\times S,D_{\bullet})/S),
p_{S*}E_{et}(\mathbb Z((\bar D_{\bullet}\times S,D_{\bullet})/\bar X\times S)(d_X)[2d_X]))
\xrightarrow{\sim} \\
\Hom(\mathbb Z((\bar D_{\bullet}\times S\times\bar X,D_{\bullet})/\bar X\times S), \\
\mathbb Z^{tr}((\bar D_{\bullet}\times S\times\mathbb P^{d_X},D_{\bullet}\times\mathbb P^{d_X})/\bar X\times S)/
\mathbb Z^{tr}((-)\times\mathbb P^{d_X-1},(-)\times\mathbb P^{d_X-1})) \\
\subset H^0(\mathcal Z_{d_{D_{\bullet}}+d_S}(\square^*\times\bar D_{\bullet}\times\bar D_{\bullet}\times S),
\mbox{s.t.}\alpha_*(\times D_{\bullet})=D_{\bullet})
\end{eqnarray*}
\item[(ii)] The cycle $(\Delta_{\bar X}\times S)\subset\bar X\times\bar X\times S$
induces by the diagonal $\Delta_{\bar X}\subset\bar X\times\bar X$
gives the morphism in $C(\Var(\mathbb C)^2/S)$
\begin{eqnarray*}
[\Delta_{\bar X}]\in\Hom(\mathbb Z^{tr}((\bar X\times S,X)/S),
p_{S*}E_{et}(\mathbb Z((\bar X\times S,X)/\bar X\times S)(d_X)[2d_X]))
\xrightarrow{\sim} \\
\Hom(\mathbb Z((\bar X\times S\times\bar X,X)/\bar X\times S), \\
\mathbb Z^{tr}((\bar X\times S\times\mathbb P^{d_X},X\times\mathbb P^{d_X})/\bar X\times S)/
\mathbb Z^{tr}((-)\times\mathbb P^{d_X-1},(-)\times\mathbb P^{d_X-1})) \\
\subset H^0(\mathcal Z_{d_X+d_S}(\square^*\times\bar X\times\bar X\times S), \mbox{s.t.}\alpha_*(\times X)=X)
\end{eqnarray*}
\end{itemize}
Let $h:U\to S$ a morphism, with $U,S\in\Var(\mathbb C)$, $U$ smooth connected (hence irreducible by smoothness).
Take, see definition-proposition \ref{RCHdef0},
$\bar f_0=\bar{h}_0:\bar X_0\to\bar S$ a compactification of $h:U\to S$ and denote by $\bar Z=\bar X_0\backslash U$.
Take, using theorem \ref{desVar}(ii), a strict desingularization 
$\bar\epsilon:(\bar X,\bar D)\to(\bar X_0,\bar Z)$ of the pair $(\bar X_0,\bar Z)$
with $\bar X\in\PSmVar(\mathbb C)$ and 
$\bar D:=\bar\epsilon^{-1}(\bar Z)=\cup_{i=1}^s\bar D_i\subset\bar X$ a normal crossing divisor.
Denote $d_X:=\dim(\bar X)=\dim(U)$.
\begin{itemize}
\item[(iii)]We get from (i) and (ii) the morphism in $C(\Var(\mathbb C)^2/S)$
\begin{eqnarray*}
T(p_{S\sharp},p_{S*})(\mathbb Z((\bar D_{\bullet}\times S,D_{\bullet})/\bar X\times S),
\mathbb Z((\bar X\times S,X)/\bar X\times S)):=([\Delta_{\bar D_{\bullet}}],[\Delta_{\bar X}]): \\
\Cone(\mathbb Z(i_{\bullet}\times I):
(\mathbb Z^{tr}((\bar D_{\bullet}\times S,D_{\bullet})/S),u_{IJ})\to\mathbb Z^{tr}((\bar X\times S,X)/S))\to \\
p_{S*}E_{et}(\Cone(\mathbb Z(i_{\bullet}\times I):
(\mathbb Z((\bar D_{\bullet}\times S,D_{\bullet})/\bar X\times S),u_{IJ})\to \\
\mathbb Z((\bar X\times S,X)/\bar X\times S)))(d_X)[2d_X]=:R_{(\bar X,\bar D)/S}(\mathbb Z(U/S))(d_X)[2d_X]
\end{eqnarray*}
\item[(iii)']which gives the map in $C(\Var(\mathbb C)^{2,smpr}/S)$
\begin{eqnarray*}
T^{\mu,q}(p_{S\sharp},p_{S*})(\mathbb Z((\bar D_{\bullet}\times S,D_{\bullet})/\bar X\times S),
\mathbb Z((\bar X\times S,X)/\bar X\times S)): \\
\Cone(\mathbb Z(i_{\bullet}\times I):
(\mathbb Z^{tr}((\bar D_{\bullet}\times S,D_{\bullet})/S),u_{IJ})\to\mathbb Z^{tr}((\bar X\times S,X)/S))= \\
L\rho_{S*}\mu_{S*}\Cone(\mathbb Z(i_{\bullet}\times I):
(\mathbb Z^{tr}((\bar D_{\bullet}\times S,D_{\bullet})/S),u_{IJ})\to\mathbb Z^{tr}((\bar X\times S,X)/S)) \\
\xrightarrow{L\rho_{S*}\mu_{S*}T(p_{S\sharp},p_{S*})(\mathbb Z((\bar D_{\bullet}\times S,D_{\bullet})/\bar X\times S),
\mathbb Z((\bar X\times S,X)/\bar X\times S))}
L\rho_{S*}\mu_{S*}R_{(\bar X,\bar D)/S}(\mathbb Z(U/S))(d_X)[2d_X]
\end{eqnarray*}
\end{itemize}
\end{defi}

\begin{prop}\label{sharpstarprop}
Let $h:U\to S$ a morphism, with $U,S\in\Var(\mathbb C)$, $U$ irreducible.
Take, see definition-proposition \ref{RCHdef0},
$\bar f_0=\bar{h}_0:\bar X_0\to\bar S$ a compactification of $h:U\to S$ and denote by $\bar Z=\bar X_0\backslash U$.
Take, using theorem \ref{desVar}(ii), a desingularization 
$\bar\epsilon:(\bar X,\bar D)\to(\bar X_0,\Delta)$ of the pair $(\bar X_0,\Delta)$, $\bar Z\subset\Delta$
with $\bar X\in\PSmVar(\mathbb C)$ and 
$\bar D:=\bar\epsilon^{-1}(\Delta)=\cup_{i=1}^s\bar D_i\subset\bar X$ a normal crossing divisor.
Denote $d_X:=\dim(\bar X)=\dim(U)$.
\begin{itemize}
\item[(i)] The morphism
\begin{eqnarray*}
[\Delta_{\bar D_{\bullet}}]:\mathbb Z^{tr}((\bar D_{\bullet}\times S,D_{\bullet})/S),
\to p_{S*}E_{et}(\mathbb Z((\bar D_{\bullet}\times S,D_{\bullet})/\bar X\times S)(d_X)[2d_X])
\end{eqnarray*}
given in definition \ref{sharpstar}(i) is an equivalence $(\mathbb A^1,et)$ local.
\item[(ii)] The morphism
\begin{eqnarray*}
[\Delta_{\bar X}]:\mathbb Z^{tr}((\bar X\times S,X)/S),
\to p_{S*}E_{et}(\mathbb Z((\bar X\times S,X)/\bar X\times S)(d_X)[2d_X])
\end{eqnarray*}
given in definition \ref{sharpstar}(ii) is an equivalence $(\mathbb A^1,et)$ local.
\end{itemize}
Let $h:U\to S$ a morphism, with $U,S\in\Var(\mathbb C)$, $U$ smooth connected (hence irreducible by smoothness).
Take, see definition-proposition \ref{RCHdef0},
$\bar f_0=\bar{h}_0:\bar X_0\to\bar S$ a compactification of $h:U\to S$ and denote by $\bar Z=\bar X_0\backslash U$.
Take, using theorem \ref{desVar}(ii), a strict desingularization 
$\bar\epsilon:(\bar X,\bar D)\to(\bar X_0,\bar Z)$ of the pair $(\bar X_0,\bar Z)$,
with $\bar X\in\PSmVar(\mathbb C)$ and 
$\bar D:=\bar\epsilon^{-1}(Z)=\cup_{i=1}^s\bar D_i\subset\bar X$ a normal crossing divisor.
\begin{itemize}
\item[(iii)] The morphism
\begin{eqnarray*}
T(p_{S\sharp},p_{S*})(\mathbb Z((\bar D_{\bullet}\times S,D_{\bullet})/\bar X\times S),
\mathbb Z((\bar X\times S,X)/\bar X\times S)):=([\Delta_{\bar D_{\bullet}}],[\Delta_{\bar X}]): \\
\Cone(\mathbb Z(i_{\bullet}\times I):
(\mathbb Z^{tr}((\bar D_{\bullet}\times S,D_{\bullet})/S),u_{IJ})\to\mathbb Z^{tr}((\bar X\times S,X)/S))\to \\
p_{S*}E_{et}(\Cone(\mathbb Z(i_{\bullet}\times I):
(\mathbb Z((\bar D_{\bullet}\times S,D_{\bullet})/\bar X\times S),u_{IJ})\to \\
\mathbb Z((\bar X\times S,X)/\bar X\times S)))(d_X)[2d_X]=:R_{(\bar X,\bar D)/S}(\mathbb Z(U/S))(d_X)[2d_X]
\end{eqnarray*}
given in definition \ref{sharpstar}(iii)' is an equivalence $(\mathbb A^1,et)$ local.
\item[(iii)']The morphism
\begin{eqnarray*}
T^{\mu,q}(p_{S\sharp},p_{S*})(\mathbb Z((\bar D_{\bullet}\times S,D_{\bullet})/\bar X\times S),
\mathbb Z((\bar X\times S,X)/\bar X\times S)): \\
\Cone(\mathbb Z(i_{\bullet}\times I):
(\mathbb Z^{tr}((\bar D_{\bullet}\times S,D_{\bullet})/S),u_{IJ})\to\mathbb Z^{tr}((\bar X\times S,X)/S)) \\
\to L\rho_{S*}\mu_{S*}R_{(\bar X,\bar D)/S}(\mathbb Z(U/S))(d_X)[2d_X]
\end{eqnarray*}
given in definition \ref{sharpstar}(iii)' is an equivalence $(\mathbb A^1,et)$ local.
\end{itemize}
\end{prop}

\begin{proof}
\noindent(i):
By Yoneda lemma, it is equivalent to show that for every morphism $g:T\to S$
with $T\in\Var(\mathbb C)$ and every closed subset $E\subset T$, the composition morphism
\begin{eqnarray*}
[\Delta_{\bar D_{\bullet}}]:
\Hom^{\bullet}(\mathbb Z((T,E)/S),C_*\mathbb Z^{tr}((\bar D_{\bullet}\times S,D_{\bullet})/S))
\xrightarrow{\Hom^{\bullet}(\mathbb Z((T,E)/S),C_*\Delta_{\bar D_{\bullet}})} \\
\Hom^{\bullet}(\mathbb Z((T,E)/S),p_{S*}E_{et}(\mathbb Z((\bar D_{\bullet}\times S,D_{\bullet})/\bar X\times S)(d_X)[2d_X]))
\end{eqnarray*}
is a quasi-isomorphism of abelian groups.
But this map is the composite
\begin{eqnarray*}
\Hom^{\bullet}(\mathbb Z((T,E)/S),\mathbb Z^{tr}((\bar D_{\bullet}\times S,D_{\bullet})/S))
\xrightarrow{[\Delta_{\bar D_{\bullet}}]} \\ 
\Hom^{\bullet}(\mathbb Z((T,E)/S),p_{S*}E_{et}(\mathbb Z((\bar D_{\bullet}\times S,D_{\bullet})/\bar X\times S)(d_X)[2d_X]))
\xrightarrow{\sim} \\
\Hom^{\bullet}(\mathbb Z((T\times\bar X,E)/S\times\bar X), \\
C_*\mathbb Z^{tr}((\bar D_{\bullet}\times S\times\mathbb P^{d_X},D_{\bullet}\times\mathbb P^{d_X})/\bar X\times S)/
\mathbb Z^{tr}((-)\times\mathbb P^{d_X-1},(-)\times\mathbb P^{d_X-1}))
\end{eqnarray*}
which is clearly a quasi-isomorphism.

\noindent(ii): Similar to (i).

\noindent(iii):Follows from (i) and (ii).

\noindent(iii)':Follows from (iii) and the fact that 
$\mu_{S*}$ preserve $(\mathbb A^1,et)$ local equivalence (see proposition \ref{mu12}) and the fact that
$\rho_{S*}$ preserve $(\mathbb A^1,et)$ local equivalence (see proposition \ref{rho12}).
\end{proof}

\begin{defi}\label{tus}
\begin{itemize}
\item[(i)]Let $h:U\to S$ a morphism, with $U,S\in\Var(\mathbb C)$, $U$ smooth.
Take, see definition-proposition \ref{RCHdef0},
$\bar f_0=\bar{h}_0:\bar X_0\to\bar S$ a compactification of $h:U\to S$ and denote by $\bar Z=\bar X_0\backslash U$.
Take, using theorem \ref{desVar}(ii), a strict desingularization 
$\bar\epsilon:(\bar X,\bar D)\to(\bar X_0,\bar Z)$ of the pair $(\bar X_0,\bar Z)$,
with $\bar X\in\PSmVar(\mathbb C)$ and 
$\bar D:=\bar\epsilon^{-1}(\bar Z)=\cup_{i=1}^s\bar D_i\subset\bar X$ a normal crossing divisor.
We will consider the following canonical map in $C(\Var(\mathbb C)^{sm}/S)$
\begin{eqnarray*}
T_{(\bar X,\bar D)/S}(U/S):\Gr_{S*}^{12}L\rho_{S*}\mu_{S*}R_{(\bar X,\bar D)/S}(\mathbb Z(U/S)) 
\xrightarrow{q}
\Gr_{S*}^{12}\rho_{S*}\mu_{S*}R_{(\bar X,\bar D)/S}(\mathbb Z(U/S)) \\
\xrightarrow{r_{(\bar X,\bar D)/S}(\mathbb Z(U/S))} 
\Gr_{S*}^{12}\rho_{S*}\mu_{S*}p_{S*}E_{et}(\mathbb Z((U\times S,U)/U\times S)) 
\xrightarrow{l(U/S)}h_*E_{et}(\mathbb Z(U/U))=:\mathbb D_S^0(\mathbb Z(U/S))
\end{eqnarray*}
where, for $h':V\to S$ a smooth morphism with $V\in\Var(\mathbb C)$,
\begin{eqnarray*}
l^{00}(U/S)(V/S):\mathbb Z((U\times S,U)/U\times S)(V\times U\times S,V\times_SU/U\times S)
\to\mathbb Z(U/U)(V\times_S U), \alpha\mapsto\alpha_{|V\times_SU}  
\end{eqnarray*}
which gives
\begin{eqnarray*}
l^0(U/S)(V/S):E^0_{et}(\mathbb Z((U\times S,U)/U\times S))(V\times U\times S,V\times_SU/U\times S)
\to E^0_{et}(\mathbb Z(U/U))(V\times_S U), 
\end{eqnarray*}
and by induction 
\begin{equation*}
\tau^{\leq i}l(U/S):\Gr_{S*}^{12}\rho_{S*}\mu_{S*}p_{S*}E^{\leq i}_{et}(\mathbb Z((U\times S,U)/U\times S)) 
\to h_*E^{\leq i}_{et}(\mathbb Z(U/U))
\end{equation*}
where $\tau^{\leq i}$ is the cohomological truncation.
\item[(ii)]Let $g:U'/S\to U/S$ a morphism, with $U'/S=(U',h'),U/S=(U,h)\in\Var(\mathbb C)/S$, $U$,$U'$ smooth.
Take, see definition-proposition \ref{RCHdef0}(ii),a compactification $\bar f_0=\bar h:\bar X_0\to\bar S$ of $h:U\to S$ 
and a compactification $\bar f'_0=\bar{h}':\bar X'_0\to S$ of $h':U'\to S$ such that 
$g:U'/S\to U/S$ extend to a morphism $\bar g_0:\bar X'_0/\bar S\to\bar X_0/\bar S$. 
Denote $\bar Z=\bar X_0\backslash U$ and $\bar Z'=\bar X'_0\backslash U'$.
Take, see definition-proposition \ref{RCHdef0}(ii), a strict desingularization 
$\bar\epsilon:(\bar X,\bar D)\to(\bar X_0,\bar Z)$ of $(\bar X_0,\bar Z)$,
a strict desingularization $\bar\epsilon'_{\bullet}:(\bar X',\bar D')\to(\bar X'_0,\bar Z')$ of $(\bar X'_0,\bar Z')$
and a morphism $\bar g:\bar X'\to\bar X$ such that the following diagram commutes
\begin{equation*}
\xymatrix{\bar X'_0\ar[r]^{\bar{g}_0} & \bar X_0 \\
\bar X'\ar[u]^{\bar\epsilon'}\ar[r]^{\bar g} & \bar X\ar[u]^{\bar\epsilon}}.
\end{equation*}
Then by the diagram given in definition \ref{RCHdef}(ii), 
the following diagram in $C(\Var(\mathbb C)^{sm}/S)$ obviously commutes
\begin{equation*}
\xymatrix{\Gr_{S*}^{12}L\rho_{S*}\mu_{S*}R_{(\bar X,\bar D)/S}(\mathbb Z(U/S))
\ar[rr]^{T_{(\bar X,\bar D)/S}(U/S)}\ar[d]_{R_S^{CH}(g)} & \, & 
h_*E_{et}(\mathbb Z(U/U)):=\mathbb D_S^0(\mathbb Z(U/S))
\ar[d]^{T(g,E)(-)\circ\ad(g^*,g_*)(E_{et}(\mathbb Z(U/U))):=\mathbb D_S^0(g)} \\
\Gr_{S*}^{12}L\rho_{S*}\mu_{S*}R_{(\bar X',\bar D')/S}(\mathbb Z(U'/S))
\ar[rr]^{T_{(\bar X',\bar D')/S}(U'/S)} & \, & 
h'_*E_{et}(\mathbb Z(U'/U')):=\mathbb D_S^0(\mathbb Z(U'/S))}
\end{equation*}
where $l(U/S)$ are $l(U'/S)$ are the maps given in (i).
\item[(iii)]Let $S\in\SmVar(\mathbb C)$.
Let $F\in C(\Var(\mathbb C)^{sm}/S)$. We get from (i) and (ii) morphisms in $C(\Var(\mathbb C)^{sm}/S)$
\begin{eqnarray*}
T_S^{CH}(LF):\Gr_{S*}^{12}L\rho_{S*}\mu_{S*}R_{(\bar X^*,\bar D^*)/S}(\rho_S^*LF) \\
\xrightarrow{r_S^{CH}(LF)}\Gr_{S*}^{12}L\rho_{S*}\mu_{S*}\mathbb D^{12}_S(\rho_S^*LF)
\xrightarrow{l(L(F)}\mathbb D_S^0(L(F))
\end{eqnarray*}
\end{itemize}
\end{defi}

We will also have the following lemma

\begin{lem}\label{grRCH}
\begin{itemize}
\item[(i)]Let $h:U\to S$ a morphism, with $U,S\in\Var(\mathbb C)$, $U$ smooth.
Take, see definition-proposition \ref{RCHdef0},
$\bar f_0=\bar{h}_0:\bar X_0\to\bar S$ a compactification of $h:U\to S$ and denote by $\bar Z=\bar X_0\backslash U$.
Take, using theorem \ref{desVar}(ii), a  strict desingularization 
$\bar\epsilon:(\bar X,\bar D)\to(\bar X_0,\bar Z)$ of the pair $(\bar X_0,\bar Z)$,
with $\bar X\in\PSmVar(\mathbb C)$ and 
$\bar D:=\bar\epsilon^{-1}(\bar Z)=\cup_{i=1}^s\bar D_i\subset\bar X$ a normal crossing divisor.
Then the map in $C(\Var(\mathbb C)^{2,smpr}/S)$
\begin{eqnarray*}
\ad(\Gr_S^{12*},\Gr_{S*}^{12})(L\rho_{S*}\mu_{S*}R_{(\bar X,\bar D)/S}(\mathbb Z(U/S)))\circ q: \\
\Gr_S^{12*}L\Gr_{S*}^{12}L\rho_{S*}\mu_{S*}R_{(\bar X,\bar D)/S}(\mathbb Z(U/S))\to 
L\rho_{S*}\mu_{S*}R_{(\bar X,\bar D)/S}(\mathbb Z(U/S))
\end{eqnarray*}
is an equivalence $(\mathbb A^1,et)$ local.
\item[(ii)]Let $S\in\SmVar(\mathbb C)$.
Let $F\in C(\Var(\mathbb C)^{sm}/S)$. Then the map in $C(\Var(\mathbb C)^{2,smpr}/S)$
\begin{eqnarray*}
\ad(\Gr_S^{12*},\Gr_{S*}^{12})(L\rho_{S*}\mu_{S*}R_{(\bar X^*,\bar D^*)/S}(\rho_S^*LF))\circ q: \\
\Gr_S^{12*}L\Gr_{S*}^{12}L\rho_{S*}\mu_{S*}R_{(\bar X^*,\bar D^*)/S}(\rho_S^*LF)\to 
L\rho_{S*}\mu_{S*}R_{(\bar X^*,\bar D^*)/S}(\rho_S^*LF)
\end{eqnarray*}
is an equivalence $(\mathbb A^1,et)$ local.
\end{itemize}
\end{lem}

\begin{proof}
\noindent(i): Follows from proposition \ref{sharpstarprop}.

\noindent(ii): Follows from (i).
\end{proof}

\begin{defi}\label{RCHhatdef} 
\begin{itemize}
\item[(i)]Let $h:U\to S$ a morphism, with $U,S\in\Var(\mathbb C)$ and $U$ smooth.
Take, see definition-proposition \ref{RCHdef0},
$\bar f_0=\bar{h}_0:\bar X_0\to\bar S$ a compactification of $h:U\to S$ and denote by $\bar Z=\bar X_0\backslash U$.
Take, using theorem \ref{desVar}(ii), a strict desingularization 
$\bar\epsilon:(\bar X,\bar D)\to(\bar X_0,\bar Z)$ of the pair $(\bar X_0,\bar Z)$, with $\bar X\in\PSmVar(\mathbb C)$ and 
$\bar D:=\epsilon^{-1}(\bar Z)=\cup_{i=1}^s\bar D_i\subset\bar X$ a normal crossing divisor.  
We denote by $i_{\bullet}:\bar D_{\bullet}\hookrightarrow\bar X=\bar X_{c(\bullet)}$ the morphism of simplicial varieties
given by the closed embeddings $i_I:\bar D_I=\cap_{i\in I}\bar D_i\hookrightarrow\bar X$
We denote by $j:U\hookrightarrow\bar X$ the open embedding and by $p_S:\bar X\times S\to S$ 
and $p_S:U\times S\to S$ the projections.
Considering the graph factorization $\bar f:\bar X\xrightarrow{\bar l}\bar X\times\bar S\xrightarrow{p_{\bar S}}\bar S$
of $\bar f:\bar X\to\bar S$, where $\bar l$ is the graph embedding and $p_{\bar S}$ the projection,
we get closed embeddings $l:=\bar l\times_{\bar S}S:X\hookrightarrow\bar X\times S$ and 
$l_{D_I}:=\bar D_I\times_{\bar X} l:D_I\hookrightarrow\bar D_I\times S$.
We then consider the map in $C(\Var(\mathbb C)^{2,smpr}/S)$
\begin{eqnarray*}
T(\hat R^{CH},R^{CH})(\mathbb Z(U/S)):\hat R_{(\bar X,\bar D)/S}(\mathbb Z(U/S)) \\
\xrightarrow{:=}\Cone(\mathbb Z(i_{\bullet}\times I):
(\mathbb Z^{tr}((\bar D_{\bullet}\times S,D_{\bullet})/S),u_{IJ})\to\mathbb Z^{tr}((\bar X\times S,X)/S))(-d_X)[-2d_X] \\
\xrightarrow{T^{\mu,q}(p_{S\sharp},p_{S*})(\mathbb Z((\bar D_{\bullet}\times S,D_{\bullet})/\bar X\times S),
\mathbb Z((\bar X\times S,X)/\bar X\times S))(-d_X)[-2d_X]} \\
L\rho_{S*}\mu_{S*}R_{(\bar X,\bar D)/S}(\mathbb Z(U/S)).
\end{eqnarray*}
given in definition \ref{sharpstar}(iii).
\item[(ii)]Let $g:U'/S\to U/S$ a morphism, with $U'/S=(U',h'),U/S=(U,h)\in\Var(\mathbb C)/S$, with $U$ and $U'$ smooth.
Take, see definition-proposition \ref{RCHdef0}(ii),a compactification $\bar f_0=\bar h:\bar X_0\to\bar S$ of $h:U\to S$ 
and a compactification $\bar f'_0=\bar{h}':\bar X'_0\to\bar S$ of $h':U'\to S$ such that 
$g:U'/S\to U/S$ extend to a morphism $\bar g_0:\bar X'_0/\bar S\to\bar X_0/\bar S$. 
Denote $\bar Z=\bar X_0\backslash U$ and $\bar Z'=\bar X'_0\backslash U'$.
Take, see definition-proposition \ref{RCHdef0}(ii), a strict desingularization 
$\bar\epsilon:(\bar X,\bar D)\to(\bar X_0,\bar Z)$ of $(\bar X_0,\bar Z)$,
a strict desingularization $\bar\epsilon'_{\bullet}:(\bar X',\bar D')\to(\bar X'_0,\bar Z')$ of $(\bar X'_0,\bar Z')$
and a morphism $\bar g:\bar X'\to\bar X$ such that the following diagram commutes
\begin{equation*}
\xymatrix{\bar X'_0\ar[r]^{\bar{g}_0} & \bar X_0 \\
\bar X'\ar[u]^{\bar\epsilon'}\ar[r]^{\bar g} & \bar X\ar[u]^{\bar\epsilon}}.
\end{equation*} 
We then have, see definition-proposition \ref{RCHdef0}(ii),
the diagram (\ref{RCHdia}) in $\Fun(\Delta,\Var(\mathbb C))$
\begin{equation*}
\xymatrix{U=U_{c(\bullet)}\ar[r]^j & \bar X=\bar X_{c(\bullet)} & \, & \bar D_{s_g(\bullet)}\ar[ll]_{i_{\bullet}} \\
U'=U'_{c(\bullet)}\ar[r]^{j'}\ar[u]^g & \bar X'=\bar X'_{c(\bullet)}\ar[u]^{\bar{g}} & \bar D'_{\bullet}\ar[l]_{i'_{\bullet}} & 
\bar{g}^{-1}(\bar D_{s_g(\bullet)})\ar[l]_{i''_{g\bullet}}\ar[u]^{\bar{g}'_{\bullet}}:i'_{g\bullet}}
\end{equation*}
Consider 
\begin{eqnarray*}
[\Gamma_{\bar g}]^t\in\Hom(\mathbb Z^{tr}((\bar X\times S,X)/S)(-d_X)[-2d_X],
\mathbb Z^{tr}((\bar X'\times S,X')/S)(-d_{X'})[-2d_{X'}]) \\ 
\xrightarrow{\sim}\Hom(\mathbb Z^{tr}((\bar X\times\mathbb A^{d_{X'}}\times S,X\times\mathbb A^{d_{X'}})/S), \\
\mathbb Z_{tr}((\bar X'\times\mathbb P^{d_X}\times S,X'\times\mathbb P^{d_X})/S)/
\mathbb Z_{tr}((-)\times\mathbb P^{d_X-1},(-)\times\mathbb P^{d_X-1}))
\end{eqnarray*}
the morphism given by the transpose of the graph $\Gamma_g\subset X'\times_S X$ of $\bar g:\bar X'\to\bar X$. Then,  
since $i_{\bullet}\circ\bar g'_{\bullet}=\bar g\circ i''_{g\bullet}=\bar g\circ i'\circ\circ i'_{g\bullet}$, 
we have the factorization 
\begin{eqnarray*}
[\Gamma_g]^t\circ\mathbb Z(i_{\bullet}\times I): 
(\mathbb Z^{tr}((\bar D_{s_g(\bullet)}\times S,D_{s_g(\bullet)})/S),u_{IJ})(-d_X)[-2d_X] \\
\xrightarrow{[\Gamma_{\bar g'_{\bullet}}]^t}
(\mathbb Z^{tr}((\bar{g}^{-1}(\bar D_{s_g(\bullet)})\times S,\bar{g}^{-1}(D_{s_g(\bullet)}))/S),u_{IJ})(-d_{X'})[-2d_{X'}]) \\
\xrightarrow{\mathbb Z(i'_{g\bullet}\times I)} 
\mathbb Z^{tr}((\bar X'\times S,X')/S)(-d_{X'})[-2d_{X'}].
\end{eqnarray*}
with
\begin{eqnarray*}
[\Gamma_{\bar g'_{\bullet}}]^t\in
\Hom((\mathbb Z^{tr}((\bar D_{s_g(\bullet)}\times\mathbb A^{d_{X'}}\times S,
D_{s_g(\bullet)}\times\mathbb A^{d_{X'}})/S),u_{IJ}), \\
(\mathbb Z_{tr}((\bar{g}^{-1}(\bar D_{s_g(\bullet)})\times\mathbb P^{d_X}\times S,
\bar{g}^{-1}(D_{s_g(\bullet)})\times\mathbb P^{d_X})/S),u_{IJ})/
\mathbb Z^{tr}((-)\times\mathbb P^{d_X-1},(-)\times\mathbb P^{d_X-1})).
\end{eqnarray*}
We then consider the following map in $C(\Var(\mathbb C)^{2,pr}/S)$
\begin{eqnarray*}
\hat R_S^{CH}(g):\hat R_{(\bar X,\bar D)/S}(\mathbb Z(U/S))\xrightarrow{:=} \\
\Cone(\mathbb Z(i_{\bullet}\times I):
(\mathbb Z^{tr}((\bar D_{s_g(\bullet)}\times S,D_{s_g(\bullet)})/S),u_{IJ})\to
\mathbb Z^{tr}((\bar X\times S,X)/S))(-d_X)[-2d_X] \\
\xrightarrow{([\Gamma_{\bar g'_{\bullet}}]^t,[\Gamma_{\bar g}]^t)} \\
\Cone(\mathbb Z(i'_{g\bullet}\times I): \\
(\mathbb Z^{tr}((\bar{g}^{-1}(\bar D_{s_g(\bullet)})\times S,\bar g^{-1}(D_{s_g(\bullet)})/S),u_{IJ})
\to\mathbb Z^{tr}((\bar X'\times S,X')/S)))(-d_{X'})[-2d_{X'}] \\
\xrightarrow{(\mathbb Z(i''_{g\bullet}\times I),I)(-d_{X'})[-2d_{X'}]} \\
\Cone(\mathbb Z(i'_{\bullet}\times I):
((\mathbb Z^{tr}((\bar D'_{\bullet}\times S,D'_{\bullet)})/S),u_{IJ})\to
\mathbb Z^{tr}((\bar X'\times S,X')/S)))(-d_{X'})[-2d_{X'}] \\
\xrightarrow{=:}\hat R_{(\bar X',\bar D')/S}(\mathbb Z(U'/S))
\end{eqnarray*}
Then the following diagram in $C(\Var(\mathbb C)^{2,smpr}/S)$ commutes by definition
\begin{equation*}
\xymatrix{\hat R_{(\bar X,\bar D)/S}(\mathbb Z(U/S))
\ar[rr]^{T(\hat R^{CH},R^{CH})(\mathbb Z(U/S))}\ar[d]_{\hat R_S^{CH}(g)} & \, & 
L\rho_{S*}\mu_{S*}R_{(\bar X,\bar D)/S}(\mathbb Z(U/S))\ar[d]^{L\rho_{S*}\mu_{S*}R_S^{CH}(g)} \\
\hat R_{(\bar X',\bar D')/S}(\mathbb Z(U'/S))\ar[rr]^{T(\hat R^{CH},R^{CH})(\mathbb Z(U'/S))} & \, & 
L\rho_{S*}\mu_{S*}R_{(\bar X',\bar D')/S}(\mathbb Z(U'/S))}.
\end{equation*}
\item[(iii)] For $g_1:U''/S\to U'/S$, $g_2:U'/S\to U/S$ two morphisms
with $U''/S=(U',h''),U'/S=(U',h'),U/S=(U,h)\in\Var(\mathbb C)/S$, with $U$, $U'$ and $U''$ smooth.
We get from (i) and (ii) 
a compactification $\bar f=\bar{h}:\bar X\to\bar S$ of $h:U\to S$,
a compactification $\bar f'=\bar{h}':\bar X'\to\bar S$ of $h':U'\to S$,
and a compactification $\bar f''=\bar{h}'':\bar X''\to\bar S$ of $h'':U''\to S$,
with $\bar X,\bar X',\bar X''\in\PSmVar(\mathbb C)$, 
$\bar D:=\bar X\backslash U\subset\bar X$ $\bar D':=\bar X'\backslash U'\subset\bar X'$, 
and $\bar D'':=\bar X''\backslash U''\subset\bar X''$ normal crossing divisors, 
such that $g_1:U''/S\to U'/S$ extend to $\bar g_1:\bar X''/\bar S\to\bar X'/\bar S$,
$g_2:U'/S\to U/S$ extend to $\bar g_2:\bar X'/\bar S\to\bar X/\bar S$, and
\begin{eqnarray*}
\hat R_S^{CH}(g_2\circ g_1)=\hat R_S^{CH}(g_1)\circ \hat R_S^{CH}(g_2):
\hat R_{(\bar X,\bar D)/S}\to \hat R_{(\bar X'',\bar D'')/S} 
\end{eqnarray*}
\item[(iv)] For  
\begin{eqnarray*}
Q^*:=(\cdots\to\oplus_{\alpha\in\Lambda^n}\mathbb Z(U^n_{\alpha}/S)
\xrightarrow{(\mathbb Z(g^n_{\alpha,\beta}))}\oplus_{\beta\in\Lambda^{n-1}}\mathbb Z(U^{n-1}_{\beta}/S)\to\cdots)
\in C(\Var(\mathbb C)/S)
\end{eqnarray*}
a complex of (maybe infinite) direct sum of representable presheaves with $U^*_{\alpha}$ smooth,
we get from (i),(ii) and (iii) the map in $C(\Var(\mathbb C)^{2,smpr}/S)$
\begin{eqnarray*}
T(\hat R^{CH},R^{CH})(Q^*):\hat R^{CH}(Q^*):=
(\cdots\to\oplus_{\beta\in\Lambda^{n-1}}\varinjlim_{(\bar X^{n-1}_{\beta},\bar D^{n-1}_{\beta})/S}
\hat R_{(\bar X^{n-1}_{\beta},\bar D^{n-1}_{\beta})/S}(\mathbb Z(U^{n-1}_{\beta}/S)) \\
\xrightarrow{(\hat R_S^{CH}(g^n_{\alpha,\beta}))}\oplus_{\alpha\in\Lambda^n}\varinjlim_{(\bar X^n_{\alpha},\bar D^n_{\alpha})/S}
\hat R_{(\bar X^n_{\alpha},\bar D^n_{\alpha})/S}(\mathbb Z(U^n_{\alpha}/S))\to\cdots) 
\to L\rho_{S*}\mu_{S*}R^{CH}(Q^*),
\end{eqnarray*}
where for $(U^n_{\alpha},h^n_{\alpha})\in\Var(\mathbb C)/S$, the inductive limit run over all the compactifications 
$\bar f_{\alpha}:\bar X_{\alpha}\to\bar S$ of $h_{\alpha}:U_{\alpha}\to S$ with $\bar X_{\alpha}\in\PSmVar(\mathbb C)$
and $\bar D_{\alpha}:=\bar X_{\alpha}\backslash U_{\alpha}$ a normal crossing divisor.
For $m=(m^*):Q_1^*\to Q_2^*$ a morphism with 
\begin{eqnarray*}
Q_1^*:=(\cdots\to\oplus_{\alpha\in\Lambda^n}\mathbb Z(U^n_{1,\alpha}/S)
\xrightarrow{(\mathbb Z(g^n_{\alpha,\beta}))}\oplus_{\beta\in\Lambda^{n-1}}\mathbb Z(U^{n-1}_{1,\beta}/S)\to\cdots), \\
Q_2^*:=(\cdots\to\oplus_{\alpha\in\Lambda^n}\mathbb Z(U^n_{2,\alpha}/S)
\xrightarrow{(\mathbb Z(g^n_{\alpha,\beta}))}\oplus_{\beta\in\Lambda^{n-1}}\mathbb Z(U^{n-1}_{2,\beta}/S)\to\cdots)
\in C(\Var(\mathbb C)/S)
\end{eqnarray*}
complexes of (maybe infinite) direct sum of representable presheaves with $U^*_{1,\alpha}$ and $U^*_{2,\alpha}$ smooth,
we get again from (i),(ii) and (iii) a commutative diagram in $C(\Var(\mathbb C)^{2,smpr}/S)$
\begin{equation*}
\xymatrix{\hat R^{CH}(Q_2^*)\ar[rr]^{T(\hat R_S^{CH},R_S^{CH})(Q_2^*)}\ar[d]_{\hat R_S^{CH}(m):=(\hat R_S^{CH}(m^*))} & \, & 
L\rho_{S*}\mu_{S*}R^{CH}(Q_2^*)\ar[d]^{L\rho_{S*}\mu_{S*}R_S^{CH}(m):=L\rho_{S*}\mu_{S*}(R_S^{CH}(m^*))} \\
\hat R^{CH}(Q_1^*)\ar[rr]^{T(\hat R_S^{CH},R_S^{CH})(Q_1^*)} & \, & 
L\rho_{S*}\mu_{S*}R^{CH}(Q_1^*)}.
\end{equation*}
\end{itemize}
\end{defi}

\begin{itemize}
\item Let $S\in\Var(\mathbb C)$
For $(h,m,m')=(h^*,m^*,m^{'*}):Q_1^*[1]\to Q_2^*$ an homotopy with $Q_1^*,Q_2^*\in C(\Var(\mathbb C)/S)$
complexes of (maybe infinite) direct sum of representable presheaves with $U^*_{1,\alpha}$ and $U^*_{2,\alpha}$ smooth,
\begin{equation*}
(\hat R_S^{CH}(h),\hat R_S^{CH}(m),\hat R_S^{CH}(m'))=(\hat R_S^{CH}(h^*),\hat R_S^{CH}(m^*),\hat R_S^{CH}(m^{'*})):
R^{CH}(Q_2^*)[1]\to R^{CH}(Q_1^*)
\end{equation*}
is an homotopy in $C(\Var(\mathbb C)^{2,smpr}/S)$ using definition \ref{RCHhatdef} (iii). 
In particular if $m:Q_1^*\to Q_2^*$ with $Q_1^*,Q_2^*\in C(\Var(\mathbb C)/S)$
complexes of (maybe infinite) direct sum of representable presheaves with $U^*_{1,\alpha}$ and $U^*_{2,\alpha}$ smooth
is an homotopy equivalence, then $\hat R_S^{CH}(m):\hat R^{CH}(Q_2^*)\to\hat R^{CH}(Q_1^*)$ is an homotopy equivalence.
\item Let $S\in\SmVar(\mathbb C)$. Let $F\in\PSh(\Var(\mathbb C)^{sm}/S)$. Consider 
\begin{eqnarray*}
q:LF:=(\cdots\to\oplus_{(U_{\alpha},h_{\alpha})\in\Var(\mathbb C)^{sm}/S}\mathbb Z(U_{\alpha}/S)
\xrightarrow{(\mathbb Z(g^n_{\alpha,\beta}))}
\oplus_{(U_{\alpha},h_{\alpha})\in\Var(\mathbb C)^{sm}/S}\mathbb Z(U_{\alpha}/S)\to\cdots)\to F
\end{eqnarray*}
the canonical projective resolution given in subsection 2.3.3.
Note that the $U_{\alpha}$ are smooth since $S$ is smooth and $h_{\alpha}$ are smooth morphism.
Definition \ref{RCHhatdef}(iv) gives in this particular case the map in $C(\Var(\mathbb C)^2/S)$
\begin{eqnarray*}
T(\hat R_S^{CH},R_S^{CH})(\rho_S^*LF):\hat R^{CH}(\rho_S^*LF):=
(\cdots\to\oplus_{(U_{\alpha},h_{\alpha})\in\Var(\mathbb C)^{sm}/S}\varinjlim_{(\bar X_{\alpha},\bar D_{\alpha})/S}
\hat R_{(\bar X_{\alpha},\bar D_{\alpha})/S}(\mathbb Z(U_{\alpha}/S)) \\
\xrightarrow{(\hat R_S^{CH}(g^n_{\alpha,\beta}))}
\oplus_{(U_{\alpha},h_{\alpha})\in\Var(\mathbb C)^{sm}/S}\varinjlim_{(\bar X_{\alpha},\bar D_{\alpha})/S}
\hat R_{(\bar X_{\alpha},\bar D_{\alpha})/S}(\mathbb Z(U_{\alpha}/S))\to\cdots)
\to L\rho_{S*}\mu_{S*}R^{CH}(\rho_S^*LF),
\end{eqnarray*}
where for $(U_{\alpha},h_{\alpha})\in\Var(\mathbb C)^{sm}/S$, the inductive limit run over all the compactifications 
$\bar f_{\alpha}:\bar X_{\alpha}\to\bar S$ of $h_{\alpha}:U_{\alpha}\to S$ with $\bar X_{\alpha}\in\PSmVar(\mathbb C)$
and $\bar D_{\alpha}:=\bar X_{\alpha}\backslash U_{\alpha}$ a normal crossing divisor.
Definition \ref{RCHhatdef}(iv) gives then by functoriality in particular, for $F=F^{\bullet}\in C(\Var(\mathbb C)^{sm}/S)$, 
the map in $C(\Var(\mathbb C)^{2,smpr}/S)$
\begin{eqnarray*}
T(\hat R_S^{CH},R_S^{CH})(\rho_S^*LF):\hat R^{CH}(\rho_S^*LF)\to L\rho_{S*}\mu_{S*}R^{CH}(\rho_S^*LF).
\end{eqnarray*}

\item Let $g:T\to S$ a morphism with $T,S\in\SmVar(\mathbb C)$. Let $h:U\to S$ a smooth morphism with $U\in\Var(\mathbb C)$.
Consider the cartesian square
\begin{equation*}
\xymatrix{U_T\ar[r]^{h'}\ar[d]^{g'} & T\ar[d]^g \\
U\ar[r]^h & S}
\end{equation*}
Note that $U$ is smooth since $S$ and $h$ are smooth, and $U_T$ is smooth since $T$ and $h'$ are smooth. 
Take, see definition-proposition \ref{RCHdef0}(ii),a compactification $\bar f_0=\bar h:\bar X_0\to\bar S$ of $h:U\to S$ 
and a compactification $\bar f'_0=\bar{g\circ h'}:\bar X'_0\to\bar S$ of $g\circ h':U'\to S$ such that 
$g':U_T/S\to U/S$ extend to a morphism $\bar{g}'_0:\bar X'_0/\bar S\to \bar X_0/\bar S$. 
Denote $\bar Z=\bar X_0\backslash U$ and $\bar Z'=\bar X'_0\backslash U_T$.
Take, see definition-proposition \ref{RCHdef0}(ii), a strict desingularization 
$\bar\epsilon:(\bar X,\bar D)\to(\bar X_0,\bar Z)$ of $(\bar X_0,\bar Z)$,
a desingularization $\bar\epsilon'_{\bullet}:(\bar X',\bar D')\to(\bar X'_0,\bar Z')$ of $(\bar X'_0,\bar Z')$
and a morphism $\bar g':\bar X'\to\bar X$ such that the following diagram commutes
\begin{equation*}
\xymatrix{\bar X'_0\ar[r]^{\bar{g}'_0} & \bar X_0 \\
\bar X'\ar[u]^{\bar\epsilon'}\ar[r]^{\bar{g}'} & \bar X\ar[u]^{\bar\epsilon}}.
\end{equation*} 
We then have, see definition-proposition \ref{RCHdef0}(ii),
the following commutative diagram in $\Fun(\Delta,\Var(\mathbb C))$
\begin{equation*}
\xymatrix{U=U_{c(\bullet)}\ar[r]^j & \bar X=\bar X_{c(\bullet)} & \, & \bar D_{s_{g'}(\bullet)}\ar[ll]_{i_{\bullet}} \\
U_T=U_{T,c(\bullet)}\ar[r]^{j'}\ar[u]^{g'} & \bar X'=X'_{c(\bullet)}\ar[u]^{\bar{g}'} & \bar D'_{\bullet}\ar[l]_{i'_{\bullet}} & 
\bar{g}^{'-1}(\bar D_{s_{g'}(\bullet)})\ar[l]_{i''_{g'\bullet}}\ar[u]^{(\bar{g}')'_{\bullet}}:i'_{g\bullet}}
\end{equation*}
We then consider the following map in $C(\Var(\mathbb C)^{2,pr}/T)$, 
\begin{eqnarray*}
T(g,\hat R^{CH})(\mathbb Z(U/S)):g^*\hat R_{(\bar X,\bar D)/S}(\mathbb Z(U/S)) \\
\xrightarrow{:=}g^*\Cone(\mathbb Z(i_{\bullet}\times I):
(\mathbb Z^{tr}((\bar D_{\bullet}\times S,D_{\bullet})/S),u_{IJ})\to\mathbb Z^{tr}((\bar X\times S,X)/S))(-d_X)[-2d_X] \\ 
\xrightarrow{T(g,L)(-)\circ T(g,c)(-)} \\
\Cone(\mathbb Z(i'_{g\bullet}\times I): 
(\mathbb Z^{tr}((\bar D_{\bullet}\times T,\bar g^{-1}(D_{s_g(\bullet)})/T),u_{IJ})
\to\mathbb Z^{tr}((\bar X\times T,X')/T)))(-d_X)[-2d_X] \\
\xrightarrow{([\Gamma_{\bar g'_{\bullet}}]^t,[\Gamma_{\bar g}]^t)} \\
\Cone(\mathbb Z(i'_{g\bullet}\times I): \\
(\mathbb Z^{tr}((\bar{g}^{-1}(\bar D_{s_g(\bullet)})\times T,\bar g^{-1}(D_{s_g(\bullet)})/T),u_{IJ})
\to\mathbb Z^{tr}((\bar X'\times T,X')/T)))(-d_{X'})[-2d_{X'}] \\
\xrightarrow{(\mathbb Z(i''_{g\bullet}\times I),I)(-d_{X'})[-2d_{X'}]} \\
\Cone(\mathbb Z(i'_{\bullet}\times I):
((\mathbb Z^{tr}((\bar D'_{\bullet}\times T,D'_{\bullet)})/T),u_{IJ})\to
\mathbb Z^{tr}((\bar X'\times S,X')/T)))(-d_{X'})[-2d_{X'}] \\
\xrightarrow{=:}\hat R_{(\bar X',\bar D')/T}(\mathbb Z(U_T/T))
\end{eqnarray*}
For  
\begin{eqnarray*}
Q^*:=(\cdots\to\oplus_{\alpha\in\Lambda^n}\mathbb Z(U^n_{\alpha}/S)
\xrightarrow{(\mathbb Z(g^n_{\alpha,\beta}))}\oplus_{\beta\in\Lambda^{n-1}}\mathbb Z(U^{n-1}_{\beta}/S)\to\cdots)
\in C(\Var(\mathbb C)/S)
\end{eqnarray*}
a complex of (maybe infinite) direct sum of representable presheaves with $h^n_{\alpha}:U^n_{\alpha}\to S$ smooth,
we get the map in $C(\Var(\mathbb C)^{2,smpr}/T)$
\begin{eqnarray*}
T(g,\hat R^{CH})(Q^*):g^*\hat R^{CH}(Q^*)=
(\cdots\to\oplus_{\alpha\in\Lambda^n}\varinjlim_{(\bar X^n_{\alpha},\bar D^n_{\alpha})/S}
g^*\hat R_{(\bar X^n_{\alpha},\bar D^n_{\alpha})/S}(\mathbb Z(U^n_{\alpha}/S))\to\cdots) \\
\xrightarrow{(T(g,\hat R^{CH})(\mathbb Z(U^n_{\alpha}/S)))} 
(\cdots\to\oplus_{\alpha\in\Lambda^n}\varinjlim_{(\bar X^{n'}_{\alpha},\bar D^{n'}_{\alpha})/T}
\hat R_{(\bar X^{n'}_{\alpha},\bar D^{n'}_{\alpha})/T}(\mathbb Z(U^n_{\alpha,T}/S))\to\cdots)=:\hat R^{CH}(g^*Q^*)
\end{eqnarray*}
together with the commutative diagram in $C(\Var(\mathbb C)^{2,smpr}/T)$
\begin{equation*}
\xymatrix{g^*\hat R^{CH}(Q^*)\ar[rrr]^{T(g,\hat R^{CH})(Q^*)}\ar[d]_{g^*T(\hat R_S^{CH},R_S^{CH})(Q^*)} 
& \, & \, & \hat R^{CH}(g^*Q^*)\ar[d]^{T(\hat R_T^{CH},R_T^{CH})(g^*Q)} \\
g^*L\rho_{S*}\mu_{S*}R^{CH}(Q^*)\ar[rrr]^{T(g,R^{CH})(Q^*)\circ T(g,\mu)(-)\circ T(g,\rho)(-)\circ T(g,L)(-)} 
& \, & \, & L\rho_{T*}\mu_{T*}R^{CH}(g^*Q^*)}.
\end{equation*}
Let $F\in\PSh(\Var(\mathbb C)^{sm}/S)$. Consider 
\begin{eqnarray*}
q:LF:=(\cdots\to\oplus_{(U_{\alpha},h_{\alpha})\in\Var(\mathbb C)^{sm}/S}\mathbb Z(U_{\alpha}/S)\to\cdots)\to F
\end{eqnarray*}
the canonical projective resolution given in subsection 2.3.3.
We then get in particular the map in $C(\Var(\mathbb C)^{2,smpr}/T)$
\begin{eqnarray*}
T(g,\hat R^{CH})(\rho_S^*LF):g^*\hat R^{CH}(\rho_S^*LF)= \\
(\cdots\to\oplus_{(U_{\alpha},h_{\alpha})\in\Var(\mathbb C)^{sm}/S}\varinjlim_{(\bar X_{\alpha},\bar D_{\alpha})/S}
g^*\hat R_{(\bar X_{\alpha},\bar D_{\alpha})/S}(\mathbb Z(U_{\alpha}/S))\to\cdots) 
\xrightarrow{(T(g,\hat R^{CH})(\mathbb Z(U_{\alpha}/S)))} \\ 
(\cdots\to\oplus_{(U_{\alpha},h_{\alpha})\in\Var(\mathbb C)^{sm}/S}\varinjlim_{(\bar X'_{\alpha},\bar D'_{\alpha})/T}
\hat R_{(\bar X'_{\alpha},\bar D'_{\alpha})/T}(\mathbb Z(U_{\alpha,T}/S))\to\cdots)=:\hat R^{CH}(\rho_T^*g^*LF), 
\end{eqnarray*}
and by functoriality, we get in particular for $F=F^{\bullet}\in C(\Var(\mathbb C)^{sm}/S)$, 
the map in $C(\Var(\mathbb C)^{2,smpr}/T)$
\begin{eqnarray*}
T(g,\hat R^{CH})(\rho_S^*LF):g^*\hat R^{CH}(\rho_S^*LF)\to\hat R^{CH}(\rho_T^*g^*LF)
\end{eqnarray*}
together with the commutative diagram in $C(\Var(\mathbb C)^{2,smpr}/T)$
\begin{equation*}
\xymatrix{g^*\hat R^{CH}(\rho_S^*LF)\ar[rrr]^{T(g,\hat R^{CH})(\rho_S^*LF)}\ar[d]_{g^*T(\hat R_S^{CH},R_S^{CH})(\rho_S^*LF)} 
& \, & \, & \hat R^{CH}(\rho_T^*g^*LF)\ar[d]^{T(\hat R_T^{CH},R_T^{CH})(\rho_T^*g^*LF)} \\
g^*L\rho_{S*}\mu_{S*}R^{CH}(\rho_S^*LF)
\ar[rrr]^{L\rho_{T*}\mu_{T*}T(g,R^{CH})(\rho_S^*LF)\circ T(g,\mu)(-)\circ T(g,\rho)(-)\circ T(g,L)(-)} 
& \, & \, & L\rho_{T*}\mu_{T*}R^{CH}(\rho_T^*g^*LF)}.
\end{equation*}

\item Let $S_1,S_2\in\SmVar(\mathbb C)$ and $p:S_1\times S_2\to S_1$ the projection. 
Let $h:U\to S_1$ a smooth morphism with $U\in\Var(\mathbb C)$. Consider the cartesian square
\begin{equation*}
\xymatrix{U\times S_2\ar[r]^{h\times I}\ar[d]^{p'} & S_1\times S_2\ar[d]^p \\
U\ar[r]^h & S_1}
\end{equation*}
Take, see definition-proposition \ref{RCHdef0}(i),a compactification $\bar f_0=\bar h:\bar X_0\to\bar S_1$ of $h:U\to S_1$.
Then $\bar f_0\times I:\bar X_0\times S_2\to\bar S_1\times S_2$ is a compactification of $h\times I:U\times S_2\to S_1\times S_2$ 
and $p':U\times S_2\to U$ extend to $\bar{p}'_0:=p_{X_0}:\bar X_0\times S_2\to\bar X_0$. Denote $Z=X_0\backslash U$. 
Take see theorem \ref{desVar}(i), a strict desingularization 
$\bar\epsilon:(\bar X,\bar D)\to(\bar X_0,\bar Z)$ of the pair $(\bar X_0,\bar Z)$.  
We then have the commutative diagram (\ref{RCHdiap}) in $\Fun(\Delta,\Var(\mathbb C))$ whose squares are cartesian
\begin{equation*}
\xymatrix{U=U_{c(\bullet)}\ar[r]^j & \bar X & \bar D_{\bullet}\ar[l]_{i_{\bullet}} \\
U\times S_2=(U\times S_2)_{c(\bullet)}\ar[r]^{j\times I}\ar[u]^g & \bar X\times S_2\ar[u]^{\bar{p}':=p_{\bar X}} & 
\bar D_{\bullet}\times S_2\ar[l]_{i'_{\bullet}}\ar[u]^{\bar{p'}'_{\bullet}}}
\end{equation*}
Then the map in $C(\Var(\mathbb C)^{2,smpr}/S_1\times S_2)$
\begin{eqnarray*}
T(p,\hat R^{CH})(\mathbb Z(U/S_1)):p^*\hat R_{(\bar X,\bar D)/S_1}(\mathbb Z(U/S_1))\xrightarrow{\sim} 
\hat R_{(\bar X\times S_2,\bar D_{\bullet}\times S_2)/S_1\times S_2}(\mathbb Z(U\times S_2/S_1\times S_2))
\end{eqnarray*}
is an isomorphism.
Hence, for $Q^*\in C(\Var(\mathbb C)/S_1)$ a complex of (maybe infinite) direct sum of representable presheaves of smooth morphism,
the map in $C(\Var(\mathbb C)^{2,smpr}/S_1\times S_2)$
\begin{eqnarray*}
T(p,\hat R^{CH})(Q^*):p^*\hat R^{CH}(Q^*)\xrightarrow{\sim}\hat R^{CH}(p^*Q^*)
\end{eqnarray*}
is an isomorphism.
In particular, for $F\in C(\Var(\mathbb C)^{sm}/S_1)$ the map in $C(\Var(\mathbb C)^{2,smpr}/S_1\times S_2)$
\begin{eqnarray*}
T(p,\hat R^{CH})(\rho_{S_1}^*LF):p^*\hat R^{CH}(\rho_{S_1}^*LF)\xrightarrow{\sim}\hat R^{CH}(\rho_{S_1\times S_2}^*p^*LF)
\end{eqnarray*}
is an isomorphism.

\item Let $h_1:U_1\to S$, $h_2:U_2\to S$ two morphisms with $U_1,U_2,S\in\Var(\mathbb C)$, $U_1,U_2$ smooth.
Denote by $p_1:U_1\times_SU_2\to U_1$ and $p_2:U_1\times_S U_2\to U_2$ the projections.
Take, see definition-proposition \ref{RCHdef0}(i)),
a compactification $\bar f_{10}=\bar{h}_1:\bar X_{10}\to\bar S$ of $h_1:U_1\to S$
and a compactification $\bar f_{20}=\bar{h}_2:\bar X_{20}\to\bar S$ of $h_2:U_2\to S$. Then, 
\begin{itemize}
\item $\bar f_{10}\times\bar f_{20}:\bar X_{10}\times_{\bar S}\bar X_{20}\to S$ 
is a compactification of $h_1\times h_2:U_1\times_SU_2\to S$.
\item $\bar p_{10}:=p_{X_{10}}:\bar X_{10}\times_{\bar S} \bar X_{20}\to\bar X_{10}$ 
is a compactification of $p_1:U_1\times_SU_2\to U_1$.
\item $\bar p_{20}:=p_{X_{20}}:\bar X_{10}\times_{\bar S}\bar X_{20}\to\bar X_{20}$ 
is a compactification of $p_2:U_1\times_SU_2\to U_2$.
\end{itemize}
Denote $\bar Z_1=\bar X_{10}\backslash U_1$ and $\bar Z_2=\bar X_{20}\backslash U_2$.
Take, see theorem \ref{desVar}(i), a strict desingularization 
$\bar\epsilon_1:(\bar X_1,\bar D)\to(\bar X_{10},Z_1)$ of the pair $(\bar X_{10},\bar Z_1)$
and a  strictdesingularization 
$\bar\epsilon_2:(\bar X_2,\bar E)\to(\bar X_{20},Z_2)$ of the pair $(\bar X_{20},\bar Z_2)$. 
Take then a strict desingularization
\begin{equation*}
\bar\epsilon_{12}:((\bar X_1\times_{\bar S}\bar X_2)^N,\bar F)\to
(\bar X_1\times_{\bar S}\bar X_2,(D\times_{\bar S}\bar X_2)\cup(\bar X_1\times_{\bar S}\bar E)) 
\end{equation*}
of the pair $(\bar X_1\times_{\bar S}\bar X_2,(\bar D\times_{\bar S}\bar X_2)\cup(\bar X_1\times_{\bar S}\bar E))$. 
We have then the following commutative diagram 
\begin{equation*}
\xymatrix{\, & \bar X_1\ar[r]^{\bar f_1} & \bar S \\
\, & \bar X_1\times_{\bar S}\bar X_2\ar[r]^{\bar p_1}\ar[u]^{\bar p_2} &\bar X_2\ar[u]^{\bar f_2} \\
(\bar X_1\times_{\bar S}\bar X_2)^N\ar[ru]^{\bar\epsilon_{12}}\ar[rru]^{(\bar p_1)^N}\ar[ruu]^{(\bar p_2)^N} & \, & \,}
\end{equation*}
and
\begin{itemize}
\item $\bar f_1\times\bar f_2:\bar X_1\times_{\bar S}\bar X_2\to\bar S$ 
is a compactification of $h_1\times h_2:U_1\times_SU_2\to S$.
\item $(\bar p_1)^N:=\bar p_1\circ\epsilon_{12}:(\bar X_1\times_{\bar S}\bar X_2)^N\to\bar X_1$ is a compactification of 
$p_1:U_1\times_SU_2\to U_1$.
\item $(\bar p_2)^N:=\bar p_2\circ\epsilon_{12}:(\bar X_1\times_{\bar S}\bar X_2)^N\to\bar X_2$ is a compactification of 
$p_2:U_1\times_SU_2\to U_2$.
\end{itemize}
We have then the morphism in $C(\Var(\mathbb C)^{2,smpr}/S)$
\begin{eqnarray*}
T(\otimes,\hat R_S^{CH})(\mathbb Z(U_1/S),\mathbb Z(U_2/S)):=\hat R_S^{CH}(p_1)\otimes\hat R_S^{CH}(p_2): \\
\hat R_{(\bar X_1,\bar D)/S}(\mathbb Z(U_1/S))\otimes\hat R_{(X_2,E))/S}(\mathbb Z(U_2/S)) 
\xrightarrow{\sim}\hat R_{(\bar X_1\times_{\bar S}\bar X_2)^N,\bar F)/S}(\mathbb Z(U_1\times_S U_2/S))
\end{eqnarray*}
For   
\begin{eqnarray*}
Q_1^*:=(\cdots\to\oplus_{\alpha\in\Lambda^n}\mathbb Z(U^n_{1,\alpha}/S)
\xrightarrow{(\mathbb Z(g^n_{\alpha,\beta}))}\oplus_{\beta\in\Lambda^{n-1}}\mathbb Z(U^{n-1}_{1,\beta}/S)\to\cdots), \\
Q_2^*:=(\cdots\to\oplus_{\alpha\in\Lambda^n}\mathbb Z(U^n_{2,\alpha}/S)
\xrightarrow{(\mathbb Z(g^n_{\alpha,\beta}))}\oplus_{\beta\in\Lambda^{n-1}}\mathbb Z(U^{n-1}_{2,\beta}/S)\to\cdots)
\in C(\Var(\mathbb C)/S)
\end{eqnarray*}
complexes of (maybe infinite) direct sum of representable presheaves with $U^*_{\alpha}$ smooth, 
we get the morphism in $C(\Var(\mathbb C)^{2,smpr}/S)$
\begin{eqnarray*}
T(\otimes,\hat R_S^{CH})(Q^*_1,Q^*_2):\hat R^{CH}(Q^*_1)\otimes R^{CH}(Q^*_2) 
\xrightarrow{(T(\otimes,\hat R_S^{CH})(\mathbb Z(U_{1,\alpha}^m),\mathbb Z(U_{2,\beta}^n))}
\hat R^{CH}(Q^*_1\otimes Q^*_2))
\end{eqnarray*},
together with the commutative diagram in $C(\Var(\mathbb C)^{2,smpr}/S)$
\begin{equation*}
\xymatrix{\hat R^{CH}(Q^*_1)\otimes R^{CH}(Q^*_2)
\ar[rrr]^{T(\otimes,\hat R_S^{CH})(Q_1^*,Q_2)}
\ar[d]_{T(\hat R_S^{CH},R_S^{CH})(Q_1^*)\otimes T(\hat R_S^{CH},R_S^{CH})(Q_2^*)} 
& \, & \, & \hat R^{CH}(Q^*_1\times Q_2^*)\ar[d]^{T(\hat R_S^{CH},R_S^{CH})(Q^*_1\otimes Q^*_2)} \\
L\rho_{S*}\mu_{S*}(R^{CH}(Q_1^*)\otimes R^{CH}(Q_2^*))
\ar[rrr]^{L\rho_{S*}\mu_{S*}T(\otimes,R_S^{CH})(Q_1^*,Q_2)} 
& \, & \, & L\rho_{S*}\mu_{S*}R^{CH}(Q_1^*\otimes Q_2)}.
\end{equation*}
For $F_1,F_2\in C(\Var(\mathbb C)^{sm}/S)$, we get in particular the morphism in $C(\Var(\mathbb C)^2/S)$
\begin{eqnarray*}
T(\otimes,R_S^{CH})(\rho_S^*LF_1,\rho_S^*LF_2):R^{CH}(\rho_S^*LF_1)\otimes R^{CH}(\rho_S^*LF_2) 
\to R^{CH}(\rho_S^*(LF_1\otimes LF_2))
\end{eqnarray*}
together with the commutative diagram in $C(\Var(\mathbb C)^{2,smpr}/S)$
\begin{equation*}
\xymatrix{\hat R^{CH}(\rho_S^*LF_1)\otimes R^{CH}(\rho_S^*LF_2)
\ar[rrr]^{T(\otimes,\hat R_S^{CH})(\rho_S^*LF_1,\rho_S^*LF_2)}
\ar[d]_{T(\hat R_S^{CH},R_S^{CH})(\rho_S^*LF_1)\otimes T(\hat R_S^{CH},R_S^{CH})(\rho_S^*LF_2)} 
& \, & \, & \hat R^{CH}(\rho_S^*LF_1\times\rho_S^*LF_2)
\ar[d]^{T(\hat R_S^{CH},R_S^{CH})(\rho_S^*LF_1\otimes\rho_S^*LF_2)} \\
L\rho_{S*}\mu_{S*}(R^{CH}(\rho_S^*LF_1)\otimes R^{CH}(\rho_S^*LF_2))
\ar[rrr]^{L\rho_{S*}\mu_{S*}T(\otimes,R_S^{CH})(\rho_S^*LF_1,\rho_S^*LF_2)} 
& \, & \, & L\rho_{S*}\mu_{S*}R^{CH}(\rho_S^*LF_1\times\rho_S^*LF_2)}.
\end{equation*}

\end{itemize}

For $S\in\Var(\mathbb C)$, we will use rather the functors $R^{0CH}_S$ and $\hat R^{0CH}_S$
since we are working in the image of the graph functor $\Gr_S^{12}:\Var(\mathbb C)/S\to\Var(\mathbb C)^2/S$.
We have the full subcategory $\SmVar(\mathbb C)/S\subset\Var(\mathbb C)/S$ 
whose objects are morphisms $f:X\to S$ with $X\in\SmVar(\mathbb C)$. Then $\Gr_S^{12}(\SmVar(\mathbb C)/S)\subset\Var(\mathbb C)^{2,smpr}/S$.
If $S\in\SmVar(\mathbb C)$, we have the factorization of morphism of site 
\begin{equation*}
\Gr_S^{12}:\Var(\mathbb C)^{2,smpr}/S\xrightarrow{\Gr_S^{12}}\SmVar(\mathbb C)/S\xrightarrow{\rho_S}\Var(\mathbb C)^{sm}/S.
\end{equation*}

\begin{defi}\label{R0CHdef} 
\begin{itemize}
\item[(i)]Let $h:U\to S$ a morphism, with $U,S\in\Var(\mathbb C)$ and $U$ smooth.
Take, see definition-proposition \ref{RCHdef0},
$\bar f_0=\bar{h}_0:\bar X_0\to\bar S$ a compactification of $h:U\to S$ and denote by $\bar Z=\bar X_0\backslash U$.
Take, using theorem \ref{desVar}(ii), a strict desingularization 
$\bar\epsilon:(\bar X,\bar D)\to(\bar X_0,\bar Z)$ of the pair $(\bar X_0,\bar Z)$, with $\bar X\in\PSmVar(\mathbb C)$ and 
$\bar D:=\epsilon^{-1}(\bar Z)=\cup_{i=1}^s\bar D_i\subset\bar X$ a normal crossing divisor.  
We denote by $i_{\bullet}:\bar D_{\bullet}\hookrightarrow\bar X=\bar X_{c(\bullet)}$ the morphism of simplicial varieties
given by the closed embeddings $i_I:\bar D_I=\cap_{i\in I}\bar D_i\hookrightarrow\bar X$
We denote by $j:U\hookrightarrow\bar X$ the open embedding. We then consider the following map in $C(\Var(\mathbb C)/S)$
\begin{eqnarray*}
r^0_{(\bar X,\bar D)/S}(\mathbb Z(U/S)):R^0_{(\bar X,\bar D)/S}(\mathbb Z(U/S)) \\
\xrightarrow{:=}
\bar f_*E_{et}(\Cone(\mathbb Z(i_{\bullet}):(\mathbb Z((\bar D_{\bullet})/\bar X),u_{IJ})\to\mathbb Z((\bar X/\bar X)))) \\
\xrightarrow{\bar f_*E_{et}(0,k\circ\ad(j^*,j_*)(\mathbb Z(\bar X/\bar X)))}
h_*E_{et}(\mathbb Z(U/U))=:\mathbb D^0_S(\mathbb Z(U/S)).
\end{eqnarray*}
Note that $\mathbb Z(\bar D_I/\bar X)$ and $\mathbb Z(\bar X/\bar X)$ are obviously $\mathbb A^1$ invariant.
Note that $r_{(X,D)/S}$ is NOT an equivalence $(\mathbb A^1,et)$ local by proposition \ref{rho1} since
$\rho_{\bar X*}\mathbb Z(\bar D_{\bullet}/\bar X)=0$, 
and $\rho_{\bar X*}\ad(j^*,j_*)(\mathbb Z(\bar X/\bar X))$ is not an equivalence $(\mathbb A^1,et)$ local.
\item[(ii)]Let $g:U'/S\to U/S$ a morphism, with $U'/S=(U',h'),U/S=(U,h)\in\Var(\mathbb C)/S$, with $U$ and $U'$ smooth.
Take, see definition-proposition \ref{RCHdef0}(ii),a compactification $\bar f_0=\bar h:\bar X_0\to\bar S$ of $h:U\to S$ 
and a compactification $\bar f'_0=\bar{h}':\bar X'_0\to\bar S$ of $h':U'\to S$ such that 
$g:U'/S\to U/S$ extend to a morphism $\bar g_0:\bar X'_0/\bar S\to\bar X_0/\bar S$. 
Denote $\bar Z=\bar X_0\backslash U$ and $\bar Z'=\bar X'_0\backslash U'$.
Take, see definition-proposition \ref{RCHdef0}(ii), a strict desingularization 
$\bar\epsilon:(\bar X,\bar D)\to(\bar X_0,\bar Z)$ of $(\bar X_0,\bar Z)$,
a strict desingularization $\bar\epsilon'_{\bullet}:(\bar X',\bar D')\to(\bar X'_0,\bar Z')$ of $(\bar X'_0,\bar Z')$
and a morphism $\bar g:\bar X'\to\bar X$ such that the following diagram commutes
\begin{equation*}
\xymatrix{\bar X'_0\ar[r]^{\bar{g}_0} & \bar X_0 \\
\bar X'\ar[u]^{\bar\epsilon'}\ar[r]^{\bar g} & \bar X\ar[u]^{\bar\epsilon}}.
\end{equation*} 
We then have, see definition-proposition \ref{RCHdef0}(ii),
the commutative diagram (\ref{RCHdia}) in $\Fun(\Delta,\Var(\mathbb C))$
\begin{equation*}
\xymatrix{U=U_{c(\bullet)}\ar[r]^j & \bar X=\bar X_{c(\bullet)} & \, & \bar D_{s_g(\bullet)}\ar[ll]_{i_{\bullet}} \\
U'=U'_{c(\bullet)}\ar[r]^{j'}\ar[u]^g & \bar X'=\bar X'_{c(\bullet)}\ar[u]^{\bar{g}} & \bar D'_{\bullet}\ar[l]_{i'_{\bullet}} & 
\bar{g}^{-1}(\bar D_{s_g(\bullet)})\ar[l]_{i''_{g\bullet}}\ar[u]^{\bar{g}'_{\bullet}}:i'_{g\bullet}}
\end{equation*}
We then consider the following map in $C(\Var(\mathbb C)/S)$
\begin{eqnarray*}
R_S^{0CH}(g):R^0_{(\bar X,\bar D)/S}(\mathbb Z(U/S))\xrightarrow{:=} \\
\bar f_*E_{et}(\Cone(\mathbb Z(i_{\bullet}):(\mathbb Z((\bar D_{s_g(\bullet)})/\bar X),u_{IJ})\to\mathbb Z(\bar X/\bar X))) \\
\xrightarrow{T(\bar g,E)(-)\circ p_{S*}\ad(\bar{g}^*,\bar{g}_*)(-)} \\
\bar f'_*E_{et}(\Cone(\mathbb Z(i'_{g\bullet}): 
(\mathbb Z((\bar{g}^{-1}(\bar D_{s_g(\bullet)})/\bar X'),u_{IJ})\to\mathbb Z((\bar X'/\bar X'))) \\
\xrightarrow{\bar f'_*E_{et}(\mathbb Z(i''_{g\bullet}),I)} \\
\bar f'_*E_{et}(\Cone(\mathbb Z(i'_{\bullet}):(\mathbb Z(\bar D'_{\bullet}/\bar X'),u_{IJ})\to\mathbb Z(\bar X'/\bar X'))) \\
\xrightarrow{=:}R^0_{(\bar X',\bar D')/S}(\mathbb Z(U'/S))
\end{eqnarray*}
Then by the diagram (\ref{RCHdia}) and adjonction, the following diagram in $C(\Var(\mathbb C)/S)$ obviously commutes
\begin{equation*}
\xymatrix{R^0_{(\bar X,\bar D)/S}(\mathbb Z(U/S))\ar[rr]^{r_{(\bar X,\bar D)/S}(\mathbb Z(U/S))}\ar[d]_{R_S^{0CH}(g)} & \, & 
h_*E_{et}(\mathbb Z(U/U)=:\mathbb D^0_S(\mathbb Z(U/S))\ar[d]^{D_S(g):=T(g,E)(-)\circ\ad(g^*,g_*)(E_{et}(\mathbb Z(U/U)))} \\
R^0_{(\bar X',\bar D')/S}(\mathbb Z(U'/S))\ar[rr]^{r^0_{(\bar X',\bar D')/S}(\mathbb Z(U'/S))} & \, & 
h'_*E_{et}(\mathbb Z(U'/U'))=:\mathbb D^0_S(\mathbb Z(U'/S))}.
\end{equation*}
\item[(iii)]For $g_1:U''/S\to U'/S$, $g_2:U'/S\to U/S$ two morphisms
with $U''/S=(U',h''),U'/S=(U',h'),U/S=(U,h)\in\Var(\mathbb C)/S$, with $U$, $U'$ and $U''$ smooth.
We get from (i) and (ii) 
a compactification $\bar f=\bar{h}:\bar X\to\bar S$ of $h:U\to S$,
a compactification $\bar f'=\bar{h}':\bar X'\to\bar S$ of $h':U'\to S$,
and a compactification $\bar f''=\bar{h}'':\bar X''\to\bar S$ of $h'':U''\to S$,
with $\bar X,\bar X',\bar X''\in\PSmVar(\mathbb C)$, 
$\bar D:=\bar X\backslash U\subset\bar X$ $\bar D':=\bar X'\backslash U'\subset\bar X'$, 
and $\bar D'':=\bar X''\backslash U''\subset\bar X''$ normal crossing divisors, 
such that $g_1:U''/S\to U'/S$ extend to $\bar g_1:\bar X''/\bar S\to\bar X'/\bar S$,
$g_2:U'/S\to U/S$ extend to $\bar g_2:\bar X'/\bar S\to\bar X/\bar S$, and
\begin{eqnarray*}
R_S^{0CH}(g_2\circ g_1)=R_S^{0CH}(g_1)\circ R_S^{0CH}(g_2):R^0_{(\bar X,\bar D)/S}\to R^0_{(\bar X'',\bar D'')/S} 
\end{eqnarray*}
\item[(iv)] For  
\begin{eqnarray*}
Q^*:=(\cdots\to\oplus_{\alpha\in\Lambda^n}\mathbb Z(U^n_{\alpha}/S)
\xrightarrow{(\mathbb Z(g^n_{\alpha,\beta}))}\oplus_{\beta\in\Lambda^{n-1}}\mathbb Z(U^{n-1}_{\beta}/S)\to\cdots)
\in C(\Var(\mathbb C)/S)
\end{eqnarray*}
a complex of (maybe infinite) direct sum of representable presheaves with $U^*_{\alpha}$ smooth,
we get from (i), (ii) and (iii) the map in $C(\Var(\mathbb C)/S)$
\begin{eqnarray*}
r_S^{0CH}(Q^*):R^{0CH}(Q^*):=
(\cdots\to\oplus_{\beta\in\Lambda^{n-1}}\varinjlim_{(\bar X^{n-1}_{\beta},\bar D^{n-1}_{\beta})/S}
R^0_{(\bar X^{n-1}_{\beta},\bar D^{n-1}_{\beta})/S}(\mathbb Z(U^{n-1}_{\beta}/S)) \\
\xrightarrow{(R_S^{CH}(g^n_{\alpha,\beta}))}\oplus_{\alpha\in\Lambda^n}\varinjlim_{(\bar X^n_{\alpha},\bar D^n_{\alpha})/S}
R^0_{(\bar X^n_{\alpha},\bar D^n_{\alpha})/S}(\mathbb Z(U^n_{\alpha}/S))\to\cdots)
\to\mathbb D_S(Q^*),
\end{eqnarray*}
where for $(U^n_{\alpha},h^n_{\alpha})\in\Var(\mathbb C)/S$, the inductive limit run over all the compactifications 
$\bar f_{\alpha}:\bar X_{\alpha}\to\bar S$ of $h_{\alpha}:U_{\alpha}\to S$ with $\bar X_{\alpha}\in\PSmVar(\mathbb C)$
and $\bar D_{\alpha}:=\bar X_{\alpha}\backslash U_{\alpha}$ a normal crossing divisor.
For $m=(m^*):Q_1^*\to Q_2^*$ a morphism with 
\begin{eqnarray*}
Q_1^*:=(\cdots\to\oplus_{\alpha\in\Lambda^n}\mathbb Z(U^n_{1,\alpha}/S)
\xrightarrow{(\mathbb Z(g^n_{\alpha,\beta}))}\oplus_{\beta\in\Lambda^{n-1}}\mathbb Z(U^{n-1}_{1,\beta}/S)\to\cdots), \\
Q_2^*:=(\cdots\to\oplus_{\alpha\in\Lambda^n}\mathbb Z(U^n_{2,\alpha}/S)
\xrightarrow{(\mathbb Z(g^n_{\alpha,\beta}))}\oplus_{\beta\in\Lambda^{n-1}}\mathbb Z(U^{n-1}_{2,\beta}/S)\to\cdots)
\in C(\Var(\mathbb C)/S)
\end{eqnarray*}
complexes of (maybe infinite) direct sum of representable presheaves with $U^*_{1,\alpha}$ and $U^*_{2,\alpha}$ smooth,
we get again from (i), (ii) and (iii) a commutative diagram in $C(\Var(\mathbb C)/S)$
\begin{equation*}
\xymatrix{R^{0CH}(Q_2^*)\ar[rr]^{r_S^{0CH}(Q_2^*)}\ar[d]_{R_S^{0CH}(m):=(R_S^{0CH}(m^*))} & \, & 
\mathbb D^0_S(Q_2^*)\ar[d]^{\mathbb D_S(m):=(\mathbb D^0_S(m^*))} \\
R^{0CH}(Q_1^*)\ar[rr]^{r_S^{0CH}(Q_1^*)} & \, & \mathbb D^0_S(Q_1^*)}.
\end{equation*}
\item[(v)] Let  
\begin{eqnarray*}
Q^*:=(\cdots\to\oplus_{\alpha\in\Lambda^n}\mathbb Z(U^n_{\alpha}/S)
\xrightarrow{(\mathbb Z(g^n_{\alpha,\beta}))}\oplus_{\beta\in\Lambda^{n-1}}\mathbb Z(U^{n-1}_{\beta}/S)\to\cdots)
\in C(\Var(\mathbb C)/S)
\end{eqnarray*}
a complex of (maybe infinite) direct sum of representable presheaves with $U^*_{\alpha}$ smooth,
we have by definition 
\begin{equation*}
\Gr_S^{12*}R^{0CH}(Q^*)=R^{CH}(Q^*)\in C(\Var(\mathbb C)^2/S).
\end{equation*}
\end{itemize}
\end{defi}

\begin{itemize}
\item Let $S\in\Var(\mathbb C)$
For $(h,m,m')=(h^*,m^*,m^{'*}):Q_1^*[1]\to Q_2^*$ an homotopy with $Q_1^*,Q_2^*\in C(\Var(\mathbb C)/S)$
complexes of (maybe infinite) direct sum of representable presheaves with $U^*_{1,\alpha}$ and $U^*_{2,\alpha}$ smooth,
\begin{equation*}
(R_S^{0CH}(h),R_S^{0CH}(m),R_S^{0CH}(m'))=(R_S^{0CH}(h^*),R_S^{0CH}(m^*),R_S^{0CH}(m^{'*})):
R^{0CH}(Q_2^*)[1]\to R^{0CH}(Q_1^*)
\end{equation*}
is an homotopy in $C(\Var(\mathbb C)/S)$ using definition \ref{R0CHdef} (iii). 
In particular if $m:Q_1^*\to Q_2^*$ with $Q_1^*,Q_2^*\in C(\Var(\mathbb C)/S)$
complexes of (maybe infinite) direct sum of representable presheaves with $U^*_{1,\alpha}$ and $U^*_{2,\alpha}$ smooth
is an homotopy equivalence, then $R_S^{0CH}(m):R^{0CH}(Q_2^*)\to R^{0CH}(Q_1^*)$ is an homotopy equivalence.
\item Let $S\in\SmVar(\mathbb C)$. Let $F\in\PSh(\Var(\mathbb C)^{sm}/S)$. Consider 
\begin{eqnarray*}
q:LF:=(\cdots\to\oplus_{(U_{\alpha},h_{\alpha})\in\Var(\mathbb C)^{sm}/S}\mathbb Z(U_{\alpha}/S)
\xrightarrow{(\mathbb Z(g^n_{\alpha,\beta}))}
\oplus_{(U_{\alpha},h_{\alpha})\in\Var(\mathbb C)^{sm}/S}\mathbb Z(U_{\alpha}/S)\to\cdots)\to F
\end{eqnarray*}
the canonical projective resolution given in subsection 2.3.3.
Note that the $U_{\alpha}$ are smooth since $S$ is smooth and $h_{\alpha}$ are smooth morphism.
Definition \ref{R0CHdef}(iv) gives in this particular case the map in $C(\Var(\mathbb C)/S)$
\begin{eqnarray*}
r_S^{0CH}(\rho_S^*LF):R^{0CH}(\rho_S^*LF):=
(\cdots\to\oplus_{(U_{\alpha},h_{\alpha})\in\Var(\mathbb C)^{sm}/S}\varinjlim_{(\bar X_{\alpha},\bar D_{\alpha})/S}
R^0_{(\bar X_{\alpha},\bar D_{\alpha})/S}(\mathbb Z(U_{\alpha}/S)) \\
\xrightarrow{(R_S^{0CH}(g^n_{\alpha,\beta}))}
\oplus_{(U_{\alpha},h_{\alpha})\in\Var(\mathbb C)^{sm}/S}\varinjlim_{(\bar X_{\alpha},\bar D_{\alpha})/S}
R^0_{(\bar X_{\alpha},\bar D_{\alpha})/S}(\mathbb Z(U_{\alpha}/S))\to\cdots)
\to\mathbb D^0_S(\rho_S^*LF),
\end{eqnarray*}
where for $(U_{\alpha},h_{\alpha})\in\Var(\mathbb C)^{sm}/S$, the inductive limit run over all the compactifications 
$\bar f_{\alpha}:\bar X_{\alpha}\to\bar S$ of $h_{\alpha}:U_{\alpha}\to S$ with $\bar X_{\alpha}\in\PSmVar(\mathbb C)$
and $\bar D_{\alpha}:=\bar X_{\alpha}\backslash U_{\alpha}$ a normal crossing divisor.
Definition \ref{R0CHdef}(iv) gives then by functoriality in particular, for $F=F^{\bullet}\in C(\Var(\mathbb C)^{sm}/S)$, 
the map in $C(\Var(\mathbb C)/S)$
\begin{eqnarray*}
r_S^{0CH}(\rho_S^*LF)=(r_S^{0CH}(\rho_S^*LF^*)):R^{0CH}(\rho_S^*LF)\to\mathbb D^0_S(\rho_S^*LF).
\end{eqnarray*}

\item Let $g:T\to S$ a morphism with $T,S\in\SmVar(\mathbb C)$. Let $h:U\to S$ a smooth morphism with $U\in\Var(\mathbb C)$.
Consider the cartesian square
\begin{equation*}
\xymatrix{U_T\ar[r]^{h'}\ar[d]^{g'} & T\ar[d]^g \\
U\ar[r]^h & S}
\end{equation*}
Note that $U$ is smooth since $S$ and $h$ are smooth, and $U_T$ is smooth since $T$ and $h'$ are smooth. 
Take, see definition-proposition \ref{RCHdef0}(ii),a compactification $\bar f_0=\bar h:\bar X_0\to\bar S$ of $h:U\to S$. 
Take, see definition-proposition \ref{RCHdef0}(ii), a strict desingularization 
$\bar\epsilon:(\bar X,\bar D)\to(\bar X_0,\bar Z)$ of $(\bar X_0,\bar Z)$.
Then $\bar f'_0=\bar{g\circ h'}:\bar X_T\to\bar T$ is a compactification of $g\circ h':U_T\to S$ such that 
$g':U_T/S\to U/S$ extend to a morphism $\bar{g}'_0:\bar X_T/\bar S\to \bar X/\bar S$. 
Denote $\bar Z=\bar X_0\backslash U$ and $\bar Z'=\bar X_T\backslash U_T$.
Take, see definition-proposition \ref{RCHdef0}(ii), a strict desingularization 
$\epsilon'_{\bullet}:(\bar X',\bar D')\to(\bar X_T,\bar Z')$ of $(\bar X_T,\bar Z')$.
Denote $\bar g'=\bar g'_0\circ\epsilon'_{\bullet}:\bar X'\to\bar X$. 
We then have, see definition-proposition \ref{RCHdef0}(ii),
the following commutative diagram in $\Fun(\Delta,\Var(\mathbb C))$
\begin{equation*}
\xymatrix{U=U_{c(\bullet)}\ar[r]^j & \bar X=\bar X_{c(\bullet)} & \, & \bar D_{s_{g'}(\bullet)}\ar[ll]_{i_{\bullet}} \\
U_T=U_{T,c(\bullet)}\ar[r]^{j'}\ar[u]^{g'} & \bar X_T=X_{T,c(\bullet)}\ar[u]^{\bar{g}'_0} & \, & 
\bar{g}^{'-1}(\bar D_{s_{g'}(\bullet)})\ar[ll]_{i_{g\bullet}}\ar[u]^{(\bar{g}')'_{\bullet}} \\
U_T=U_{T,c(\bullet)}\ar[r]^{j'}\ar[u]^{g'} & \bar X'=X'_{c(\bullet)}\ar[u]^{\bar{g}'} & \bar D'_{\bullet}\ar[l]_{i'_{\bullet}} & 
\bar{g}^{'-1}(\bar D_{s_{g'}(\bullet)})\ar[l]_{i''_{g'\bullet}}\ar[u]^{\epsilon'_{\bullet}}:i'_{g\bullet}}
\end{equation*}
We then consider the following map in $C(\Var(\mathbb C)/T)$, see definition \ref{R0CHdef}(ii)
\begin{eqnarray*}
T(g,R^{0CH})(\mathbb Z(U/S)):g^*R^0_{(\bar X,\bar D)/S}(\mathbb Z(U/S)) \\ 
\xrightarrow{g^*R^{0CH}_S(g')}g^*R^0_{(\bar X',\bar D')/S}(\mathbb Z(U_T/S))=g^*g_*R^0_{(\bar X',\bar D')/T}(\mathbb Z(U_T/T)) \\
\xrightarrow{\ad(g^*,g_*)(R^0_{(\bar X',\bar D')/T}(\mathbb Z(U_T/T)))}R^0_{(\bar X',\bar D')/T}(\mathbb Z(U_T/T))
\end{eqnarray*}
For  
\begin{eqnarray*}
Q^*:=(\cdots\to\oplus_{\alpha\in\Lambda^n}\mathbb Z(U^n_{\alpha}/S)
\xrightarrow{(\mathbb Z(g^n_{\alpha,\beta}))}\oplus_{\beta\in\Lambda^{n-1}}\mathbb Z(U^{n-1}_{\beta}/S)\to\cdots)
\in C(\Var(\mathbb C)/S)
\end{eqnarray*}
a complex of (maybe infinite) direct sum of representable presheaves with $h^n_{\alpha}:U^n_{\alpha}\to S$ smooth,
we get the map in $C(\Var(\mathbb C)/T)$
\begin{eqnarray*}
T(g,R^{0CH})(Q^*):g^*R^{0CH}(Q^*)=
(\cdots\to\oplus_{\alpha\in\Lambda^n}\varinjlim_{(\bar X^n_{\alpha},\bar D^n_{\alpha})/S}
g^*R^0_{(\bar X^n_{\alpha},\bar D^n_{\alpha})/S}(\mathbb Z(U^n_{\alpha}/S))\to\cdots) \\
\xrightarrow{(T(g,R^{0CH})(\mathbb Z(U^n_{\alpha}/S)))} 
(\cdots\to\oplus_{\alpha\in\Lambda^n}\varinjlim_{(\bar X^{n'}_{\alpha},\bar D^{n'}_{\alpha})/T}
R^0_{(\bar X^{n'}_{\alpha},\bar D^{n'}_{\alpha})/T}(\mathbb Z(U^n_{\alpha,T}/S))\to\cdots)=:R^{CH}(g^*Q^*).
\end{eqnarray*}
Let $F\in\PSh(\Var(\mathbb C)^{sm}/S)$. Consider 
\begin{eqnarray*}
q:LF:=(\cdots\to\oplus_{(U_{\alpha},h_{\alpha})\in\Var(\mathbb C)^{sm}/S}\mathbb Z(U_{\alpha}/S)\to\cdots)\to F
\end{eqnarray*}
the canonical projective resolution given in subsection 2.3.3.
We then get in particular the map in $C(\Var(\mathbb C)/T)$
\begin{eqnarray*}
T(g,R^{0CH})(\rho_S^*LF):g^*R^{0CH}(\rho_S^*LF)= \\
(\cdots\to\oplus_{(U_{\alpha},h_{\alpha})\in\Var(\mathbb C)^{sm}/S}\varinjlim_{(\bar X_{\alpha},\bar D_{\alpha})/S}
g^*R^0_{(\bar X_{\alpha},\bar D_{\alpha})/S}(\mathbb Z(U_{\alpha}/S))\to\cdots) 
\xrightarrow{(T(g,R^{0CH})(\mathbb Z(U_{\alpha}/S)))} \\ 
(\cdots\to\oplus_{(U_{\alpha},h_{\alpha})\in\Var(\mathbb C)^{sm}/S}\varinjlim_{(\bar X'_{\alpha},\bar D'_{\alpha})/T}
R^0_{(\bar X'_{\alpha},\bar D'_{\alpha})/T}(\mathbb Z(U_{\alpha,T}/S))\to\cdots)=:R^{CH}(\rho_T^*g^*LF). 
\end{eqnarray*}
By functoriality, we get in particular for $F=F^{\bullet}\in C(\Var(\mathbb C)^{sm}/S)$, the map in $C(\Var(\mathbb C)/T)$
\begin{eqnarray*}
T(g,R^{0CH})(\rho_S^*LF):g^*R^{0CH}(\rho_S^*LF)\to R^{0CH}(\rho_T^*g^*LF).
\end{eqnarray*}

\item Let $S_1,S_2\in\SmVar(\mathbb C)$ and $p:S_1\times S_2\to S_1$ the projection. 
Let $h:U\to S_1$ a smooth morphism with $U\in\Var(\mathbb C)$. Consider the cartesian square
\begin{equation*}
\xymatrix{U\times S_2\ar[r]^{h\times I}\ar[d]^{p'} & S_1\times S_2\ar[d]^p \\
U\ar[r]^h & S_1}
\end{equation*}
Take, see definition-proposition \ref{RCHdef0}(i),a compactification $\bar f_0=\bar h:\bar X_0\to\bar S_1$ of $h:U\to S_1$.
Then $\bar f_0\times I:\bar X_0\times S_2\to\bar S_1\times S_2$ is a compactification of $h\times I:U\times S_2\to S_1\times S_2$ 
and $p':U\times S_2\to U$ extend to $\bar{p}'_0:=p_{X_0}:\bar X_0\times S_2\to\bar X_0$. Denote $Z=X_0\backslash U$. 
Take see theorem \ref{desVar}(i), a strict desingularization 
$\bar\epsilon:(\bar X,\bar D)\to(\bar X_0,\bar Z)$ of the pair $(\bar X_0,\bar Z)$.  
We then have the commutative diagram (\ref{RCHdiap}) in $\Fun(\Delta,\Var(\mathbb C))$ whose squares are cartesian
\begin{equation*}
\xymatrix{U=U_{c(\bullet)}\ar[r]^j & \bar X & \bar D_{\bullet}\ar[l]_{i_{\bullet}} \\
U\times S_2=(U\times S_2)_{c(\bullet)}\ar[r]^{j\times I}\ar[u]^g & \bar X\times S_2\ar[u]^{\bar{p}':=p_{\bar X}} & 
\bar D_{\bullet}\times S_2\ar[l]_{i'_{\bullet}}\ar[u]^{\bar{p'}'_{\bullet}}}
\end{equation*}
Then the map in $C(\Var(\mathbb C)/S_1\times S_2)$
\begin{eqnarray*}
T(p,R^{0CH})(\mathbb Z(U/S_1)):p^*R^0_{(\bar X,\bar D)/S_1}(\mathbb Z(U/S_1))\xrightarrow{\sim} 
R^0_{(\bar X\times S_2,\bar D_{\bullet}\times S_2)/S_1\times S_2}(\mathbb Z(U\times S_2/S_1\times S_2))
\end{eqnarray*}
is an isomorphism.
Hence, for $Q^*\in C(\Var(\mathbb C)/S_1)$ a complex of (maybe infinite) direct sum of representable presheaves of smooth morphism,
the map in $C(\Var(\mathbb C)/S_1\times S_2)$
\begin{eqnarray*}
T(p,R^{0CH})(Q^*):p^*R^{0CH}(Q^*)\xrightarrow{\sim}R^{0CH}(p^*Q^*)
\end{eqnarray*}
is an isomorphism.
In particular, for $F\in C(\Var(\mathbb C)^{sm}/S_1)$ the map in $C(\Var(\mathbb C)/S_1\times S_2)$
\begin{eqnarray*}
T(p,R^{0CH})(\rho_{S_1}^*LF):p^*R^{0CH}(\rho_{S_1}^*LF)\xrightarrow{\sim}R^{0CH}(\rho_{S_1\times S_2}^*p^*LF)
\end{eqnarray*}
is an isomorphism.

\item Let $h_1:U_1\to S$, $h_2:U_2\to S$ two morphisms with $U_1,U_2,S\in\Var(\mathbb C)$, $U_1,U_2$ smooth.
Denote by $p_1:U_1\times_SU_2\to U_1$ and $p_2:U_1\times_S U_2\to U_2$ the projections.
Take, see definition-proposition \ref{RCHdef0}(i)),
a compactification $\bar f_{10}=\bar{h}_1:\bar X_{10}\to\bar S$ of $h_1:U_1\to S$
and a compactification $\bar f_{20}=\bar{h}_2:\bar X_{20}\to\bar S$ of $h_2:U_2\to S$. Then, 
\begin{itemize}
\item $\bar f_{10}\times\bar f_{20}:\bar X_{10}\times_{\bar S}\bar X_{20}\to S$ 
is a compactification of $h_1\times h_2:U_1\times_SU_2\to S$.
\item $\bar p_{10}:=p_{X_{10}}:\bar X_{10}\times_{\bar S} \bar X_{20}\to\bar X_{10}$ 
is a compactification of $p_1:U_1\times_SU_2\to U_1$.
\item $\bar p_{20}:=p_{X_{20}}:\bar X_{10}\times_{\bar S}\bar X_{20}\to\bar X_{20}$ 
is a compactification of $p_2:U_1\times_SU_2\to U_2$.
\end{itemize}
Denote $\bar Z_1=\bar X_{10}\backslash U_1$ and $\bar Z_2=\bar X_{20}\backslash U_2$.
Take, see theorem \ref{desVar}(i), a strict desingularization 
$\bar\epsilon_1:(\bar X_1,\bar D)\to(\bar X_{10},Z_1)$ of the pair $(\bar X_{10},\bar Z_1)$
and a  strictdesingularization 
$\bar\epsilon_2:(\bar X_2,\bar E)\to(\bar X_{20},Z_2)$ of the pair $(\bar X_{20},\bar Z_2)$. 
Take then a strict desingularization
\begin{equation*}
\bar\epsilon_{12}:((\bar X_1\times_{\bar S}\bar X_2)^N,\bar F)\to
(\bar X_1\times_{\bar S}\bar X_2,(D\times_{\bar S}\bar X_2)\cup(\bar X_1\times_{\bar S}\bar E)) 
\end{equation*}
of the pair $(\bar X_1\times_{\bar S}\bar X_2,(\bar D\times_{\bar S}\bar X_2)\cup(\bar X_1\times_{\bar S}\bar E))$. 
We have then the following commutative diagram 
\begin{equation*}
\xymatrix{\, & \bar X_1\ar[r]^{\bar f_1} & \bar S \\
\, & \bar X_1\times_{\bar S}\bar X_2\ar[r]^{\bar p_1}\ar[u]^{\bar p_2} &\bar X_2\ar[u]^{\bar f_2} \\
(\bar X_1\times_{\bar S}\bar X_2)^N\ar[ru]^{\bar\epsilon_{12}}\ar[rru]^{(\bar p_1)^N}\ar[ruu]^{(\bar p_2)^N} & \, & \,}
\end{equation*}
and
\begin{itemize}
\item $\bar f_1\times\bar f_2:\bar X_1\times_{\bar S}\bar X_2\to\bar S$ 
is a compactification of $h_1\times h_2:U_1\times_SU_2\to S$.
\item $(\bar p_1)^N:=\bar p_1\circ\epsilon_{12}:(\bar X_1\times_{\bar S}\bar X_2)^N\to\bar X_1$ is a compactification of 
$p_1:U_1\times_SU_2\to U_1$.
\item $(\bar p_2)^N:=\bar p_2\circ\epsilon_{12}:(\bar X_1\times_{\bar S}\bar X_2)^N\to\bar X_2$ is a compactification of 
$p_2:U_1\times_SU_2\to U_2$.
\end{itemize}
We have then the morphism in $C(\Var(\mathbb C)/S)$
\begin{eqnarray*}
T(\otimes,R_S^{0CH})(\mathbb Z(U_1/S),\mathbb Z(U_2/S)):=R_S^{0CH}(p_1)\otimes R_S^{CH}(p_2): \\
R^0_{(\bar X_1,\bar D)/S}(\mathbb Z(U_1/S))\otimes R^0_{(X_2,E))/S}(\mathbb Z(U_2/S)) 
\xrightarrow{\sim} R^0_{(\bar X_1\times_{\bar S}\bar X_2)^N,\bar F)/S}(\mathbb Z(U_1\times_S U_2/S))
\end{eqnarray*}
For   
\begin{eqnarray*}
Q_1^*:=(\cdots\to\oplus_{\alpha\in\Lambda^n}\mathbb Z(U^n_{1,\alpha}/S)
\xrightarrow{(\mathbb Z(g^n_{\alpha,\beta}))}\oplus_{\beta\in\Lambda^{n-1}}\mathbb Z(U^{n-1}_{1,\beta}/S)\to\cdots), \\
Q_2^*:=(\cdots\to\oplus_{\alpha\in\Lambda^n}\mathbb Z(U^n_{2,\alpha}/S)
\xrightarrow{(\mathbb Z(g^n_{\alpha,\beta}))}\oplus_{\beta\in\Lambda^{n-1}}\mathbb Z(U^{n-1}_{2,\beta}/S)\to\cdots)
\in C(\Var(\mathbb C)/S)
\end{eqnarray*}
complexes of (maybe infinite) direct sum of representable presheaves with $U^*_{\alpha}$ smooth, 
we get the morphism in $C(\Var(\mathbb C)/S)$
\begin{eqnarray*}
T(\otimes,R_S^{0CH})(Q^*_1,Q^*_2): R^{0CH}(Q^*_1)\otimes R^{0CH}(Q^*_2) 
\xrightarrow{(T(\otimes,R_S^{0CH})(\mathbb Z(U_{1,\alpha}^m),\mathbb Z(U_{2,\beta}^n))} R^{0CH}(Q^*_1\otimes Q^*_2)).
\end{eqnarray*}
For $F_1,F_2\in C(\Var(\mathbb C)^{sm}/S)$, we get in particular the morphism in $C(\Var(\mathbb C)/S)$
\begin{eqnarray*}
T(\otimes,R_S^{0CH})(\rho_S^*LF_1,\rho_S^*LF_2):R^{0CH}(\rho_S^*LF_1)\otimes R^{0CH}(\rho_S^*LF_2) 
\to R^{0CH}(\rho_S^*(LF_1\otimes LF_2)).
\end{eqnarray*}

\end{itemize}

\begin{defi}\label{sharpstar0}
Let $h:U\to S$ a morphism, with $U,S\in\Var(\mathbb C)$, $U$ irreducible.
Take, see definition-proposition \ref{RCHdef0},
$\bar f_0=\bar{h}_0:\bar X_0\to\bar S$ a compactification of $h:U\to S$ and denote by $\bar Z=\bar X_0\backslash U$.
Take, using theorem \ref{desVar}, a desingularization 
$\bar\epsilon:(\bar X,\bar D)\to(\bar X_0,\Delta)$ of the pair $(\bar X_0,\Delta)$, $\bar Z\subset\Delta$,
with $\bar X\in\PSmVar(\mathbb C)$ and 
$\bar D:=\bar\epsilon^{-1}(\Delta)=\cup_{i=1}^s\bar D_i\subset\bar X$ a normal crossing divisor.
Denote $d_X:=\dim(\bar X)=\dim(U)$.
\begin{itemize}
\item[(i)] The diagonal $\Delta_{\bar D_{\bullet}}\subset\bar D_{\bullet}\times\bar D_{\bullet}$
induces the morphism in $C(\Var(\mathbb C)/S)$
\begin{eqnarray*}
[\Delta_{\bar D_{\bullet}}]\in\Hom(\mathbb Z^{tr}(\bar D_{\bullet}/S),
\bar f_*E_{et}(\mathbb Z(\bar D_{\bullet}/\bar X)(d_X)[2d_X]))
\xrightarrow{\sim} \\
\Hom(\mathbb Z(\bar D_{\bullet}\times_S\bar X/\bar X),
\mathbb Z^{tr}(\bar D_{\bullet}\times\mathbb P^{d_X}/\bar X)/
\mathbb Z^{tr}(\bar D_{\bullet}\times\mathbb P^{d_X-1}/\bar X)) \\
\subset H^0(\mathcal Z_{d_{D_{\bullet}}}(\square^*\times\bar D_{\bullet}\times_S\bar D_{\bullet}))
\end{eqnarray*}
\item[(ii)] The cycle $\Delta_{\bar X}\subset\bar X\times_S\bar X$ induces by the morphism in $C(\Var(\mathbb C)/S)$
\begin{eqnarray*}
[\Delta_{\bar X}]\in\Hom(\mathbb Z^{tr}(\bar X/S),
\bar f_*E_{et}(\mathbb Z(\bar X/\bar X)(d_X)[2d_X]))
\xrightarrow{\sim} \\
\Hom(\mathbb Z(\bar X\times_S\bar X/\bar X),
\mathbb Z^{tr}(\bar X\times\mathbb P^{d_X}/\bar X)/\mathbb Z^{tr}(\bar X\times\mathbb P^{d_X-1}/\bar X)) \\
\subset H^0(\mathcal Z_{d_X}(\square^*\times\bar X\times_S\bar X))
\end{eqnarray*}
\end{itemize}
Let $h:U\to S$ a morphism, with $U,S\in\Var(\mathbb C)$, $U$ smooth connected (hence irreducible by smoothness).
Take, see definition-proposition \ref{RCHdef0},
$\bar f_0=\bar{h}_0:\bar X_0\to\bar S$ a compactification of $h:U\to S$ and denote by $\bar Z=\bar X_0\backslash U$.
Take, using theorem \ref{desVar}(ii), a strict desingularization 
$\bar\epsilon:(\bar X,\bar D)\to(\bar X_0,\bar Z)$ of the pair $(\bar X_0,\bar Z)$
with $\bar X\in\PSmVar(\mathbb C)$ and 
$\bar D:=\bar\epsilon^{-1}(\bar Z)=\cup_{i=1}^s\bar D_i\subset\bar X$ a normal crossing divisor.
Denote $d_X:=\dim(\bar X)=\dim(U)$.
We get from (i) and (ii) the morphism in $C(\Var(\mathbb C)/S)$
\begin{eqnarray*}
T(\bar f_{\sharp},\bar f_*)(\mathbb Z(D_{\bullet}/\bar X),\mathbb Z(\bar X/\bar X))
:=([\Delta_{\bar D_{\bullet}}],[\Delta_{\bar X}]): \\
\Cone(\mathbb Z(i_{\bullet}):(\mathbb Z^{tr}(\bar D_{\bullet}/S),u_{IJ})\to\mathbb Z^{tr}(\bar X/S))\to \\
\bar f_*E_{et}(\Cone(\mathbb Z(i_{\bullet}):
(\mathbb Z(\bar D_{\bullet}/\bar X),u_{IJ})\to\mathbb Z(\bar X/\bar X)))(d_X)[2d_X] \\
=:R^0_{(\bar X,\bar D)/S}(\mathbb Z(U/S))(d_X)[2d_X].
\end{eqnarray*}
\end{defi}

\begin{defi}\label{R0CHhatdef} 
\begin{itemize}
\item[(i)]Let $h:U\to S$ a morphism, with $U,S\in\Var(\mathbb C)$ and $U$ smooth.
Take, see definition-proposition \ref{RCHdef0},
$\bar f_0=\bar{h}_0:\bar X_0\to\bar S$ a compactification of $h:U\to S$ and denote by $\bar Z=\bar X_0\backslash U$.
Take, using theorem \ref{desVar}(ii), a strict desingularization 
$\bar\epsilon:(\bar X,\bar D)\to(\bar X_0,\bar Z)$ of the pair $(\bar X_0,\bar Z)$, with $\bar X\in\PSmVar(\mathbb C)$ and 
$\bar D:=\epsilon^{-1}(\bar Z)=\cup_{i=1}^s\bar D_i\subset\bar X$ a normal crossing divisor.  
We denote by $i_{\bullet}:\bar D_{\bullet}\hookrightarrow\bar X=\bar X_{c(\bullet)}$ the morphism of simplicial varieties
given by the closed embeddings $i_I:\bar D_I=\cap_{i\in I}\bar D_i\hookrightarrow\bar X$
We denote by $j:U\hookrightarrow\bar X$ the open embedding. We then consider the map in $C(\Var(\mathbb C)/S)$
\begin{eqnarray*}
T(\hat R^{0CH},R^{0CH})(\mathbb Z(U/S)):\hat R^0_{(\bar X,\bar D)/S}(\mathbb Z(U/S)) \\
\xrightarrow{:=}\Cone(\mathbb Z(i_{\bullet}):
(\mathbb Z^{tr}(D_{\bullet}/S),u_{IJ})\to\mathbb Z^{tr}(X/S))(-d_X)[-2d_X] \\
\xrightarrow{T(\bar f_{\sharp},\bar f_*)(\mathbb Z(\bar D_{\bullet}/\bar X),\mathbb Z(\bar X/\bar X))(-d_X)[-2d_X]} \\
R^0_{(\bar X,\bar D)/S}(\mathbb Z(U/S)).
\end{eqnarray*}
given in definition \ref{sharpstar}(iii).
\item[(ii)]Let $g:U'/S\to U/S$ a morphism, with $U'/S=(U',h'),U/S=(U,h)\in\Var(\mathbb C)/S$, with $U$ and $U'$ smooth.
Take, see definition-proposition \ref{RCHdef0}(ii),a compactification $\bar f_0=\bar h:\bar X_0\to\bar S$ of $h:U\to S$ 
and a compactification $\bar f'_0=\bar{h}':\bar X'_0\to\bar S$ of $h':U'\to S$ such that 
$g:U'/S\to U/S$ extend to a morphism $\bar g_0:\bar X'_0/\bar S\to\bar X_0/\bar S$. 
Denote $\bar Z=\bar X_0\backslash U$ and $\bar Z'=\bar X'_0\backslash U'$.
Take, see definition-proposition \ref{RCHdef0}(ii), a strict desingularization 
$\bar\epsilon:(\bar X,\bar D)\to(\bar X_0,\bar Z)$ of $(\bar X_0,\bar Z)$,
a strict desingularization $\bar\epsilon'_{\bullet}:(\bar X',\bar D')\to(\bar X'_0,\bar Z')$ of $(\bar X'_0,\bar Z')$
and a morphism $\bar g:\bar X'\to\bar X$ such that the following diagram commutes
\begin{equation*}
\xymatrix{\bar X'_0\ar[r]^{\bar{g}_0} & \bar X_0 \\
\bar X'\ar[u]^{\bar\epsilon'}\ar[r]^{\bar g} & \bar X\ar[u]^{\bar\epsilon}}.
\end{equation*} 
We then have, see definition-proposition \ref{RCHdef0}(ii),
the diagram (\ref{RCHdia}) in $\Fun(\Delta,\Var(\mathbb C))$
\begin{equation*}
\xymatrix{U=U_{c(\bullet)}\ar[r]^j & \bar X=\bar X_{c(\bullet)} & \, & \bar D_{s_g(\bullet)}\ar[ll]_{i_{\bullet}} \\
U'=U'_{c(\bullet)}\ar[r]^{j'}\ar[u]^g & \bar X'=\bar X'_{c(\bullet)}\ar[u]^{\bar{g}} & \bar D'_{\bullet}\ar[l]_{i'_{\bullet}} & 
\bar{g}^{-1}(\bar D_{s_g(\bullet)})\ar[l]_{i''_{g\bullet}}\ar[u]^{\bar{g}'_{\bullet}}:i'_{g\bullet}}
\end{equation*}
Consider 
\begin{eqnarray*}
[\Gamma_{\bar g}]^t\in\Hom(\mathbb Z^{tr}(\bar X/S)(-d_X)[-2d_X],
\mathbb Z^{tr}(\bar X'/S)(-d_{X'})[-2d_{X'}]) \\ 
\xrightarrow{\sim}\Hom(\mathbb Z^{tr}(\bar X\times\mathbb P^{d_X}/S)/\mathbb Z_{tr}(\bar X\times\mathbb P^{d_X-1}/S), \\
\mathbb Z_{tr}(\bar X'\times\mathbb P^{d_{X'}}/S)/\mathbb Z_{tr}(\bar X'\times\mathbb P^{d_{X'}-1}/S)
\end{eqnarray*}
the morphism given by the transpose of the graph $\Gamma_g\subset X'\times_S X$ of $\bar g:\bar X'\to\bar X$. Then,  
since $i_{\bullet}\circ\bar g'_{\bullet}=\bar g\circ i''_{g\bullet}=\bar g\circ i'\circ\circ i'_{g\bullet}$, 
we have the factorization 
\begin{eqnarray*}
[\Gamma_g]^t\circ\mathbb Z(i_{\bullet}): 
(\mathbb Z^{tr}(\bar D_{s_g(\bullet)}/S),u_{IJ})(-d_X)[-2d_X] \\
\xrightarrow{[\Gamma_{\bar g'_{\bullet}}]^t}
(\mathbb Z^{tr}(\bar{g}^{-1}(\bar D_{s_g(\bullet)})/S),u_{IJ})(-d_{X'})[-2d_{X'}] \\
\xrightarrow{\mathbb Z(i'_{g\bullet})} 
\mathbb Z^{tr}(\bar X'/S)(-d_{X'})[-2d_{X'}].
\end{eqnarray*}
with
\begin{eqnarray*}
[\Gamma_{\bar g'_{\bullet}}]^t\in
\Hom((\mathbb Z^{tr}(\bar D_{s_g(\bullet)}\times\mathbb P^{d_X}/S),u_{IJ})/
(\mathbb Z^{tr}(\bar D_{s_g(\bullet)}\times\mathbb P^{d_{X-1}}/S),u_{IJ}), \\
(\mathbb Z_{tr}(\bar{g}^{-1}(\bar D_{s_g(\bullet)})\times\mathbb P^{d_{X'}}/S),u_{IJ})/
(\mathbb Z_{tr}(\bar{g}^{-1}(\bar D_{s_g(\bullet)})\times\mathbb P^{d_{X'-1}}/S),u_{IJ})).
\end{eqnarray*}
We then consider the following map in $C(\Var(\mathbb C)/S)$
\begin{eqnarray*}
\hat R_S^{0CH}(g):\hat R^0_{(\bar X,\bar D)/S}(\mathbb Z(U/S))\xrightarrow{:=} \\
\Cone(\mathbb Z(i_{\bullet}):
(\mathbb Z^{tr}(\bar D_{s_g(\bullet)}/S),u_{IJ})\to\mathbb Z^{tr}(\bar X/S)(-d_X)[-2d_X] \\
\xrightarrow{([\Gamma_{\bar g'_{\bullet}}]^t,[\Gamma_{\bar g}]^t)} \\
\Cone(\mathbb Z(i'_{g\bullet}): 
(\mathbb Z^{tr}(\bar{g}^{-1}(\bar D_{s_g(\bullet)})/S),u_{IJ})\to\mathbb Z^{tr}(\bar X'/S))(-d_{X'})[-2d_{X'}] \\
\xrightarrow{(\mathbb Z(i''_{g\bullet}),I)(-d_{X'})[-2d_{X'}]} \\
\Cone(\mathbb Z(i'_{\bullet}):
((\mathbb Z^{tr}(\bar D'_{\bullet}/S),u_{IJ})\to\mathbb Z^{tr}(\bar X'/S))(-d_{X'})[-2d_{X'}] \\
\xrightarrow{=:}\hat R^0_{(\bar X',\bar D')/S}(\mathbb Z(U'/S))
\end{eqnarray*}
Then the following diagram in $C(\Var(\mathbb C)/S)$ commutes by definition
\begin{equation*}
\xymatrix{\hat R^0_{(\bar X,\bar D)/S}(\mathbb Z(U/S))
\ar[rr]^{T(\hat R^{0CH},R^{0CH})(\mathbb Z(U/S))}\ar[d]_{\hat R_S^{0CH}(g)} & \, & 
R^0_{(\bar X,\bar D)/S}(\mathbb Z(U/S))\ar[d]^{R_S^{0CH}(g)} \\
\hat R^0_{(\bar X',\bar D')/S}(\mathbb Z(U'/S))\ar[rr]^{T(\hat R^{0CH},R^{0CH})(\mathbb Z(U'/S))} & \, & 
R^0_{(\bar X',\bar D')/S}(\mathbb Z(U'/S))}.
\end{equation*}
\item[(iii)] For $g_1:U''/S\to U'/S$, $g_2:U'/S\to U/S$ two morphisms
with $U''/S=(U',h''),U'/S=(U',h'),U/S=(U,h)\in\Var(\mathbb C)/S$, with $U$, $U'$ and $U''$ smooth.
We get from (i) and (ii) 
a compactification $\bar f=\bar{h}:\bar X\to\bar S$ of $h:U\to S$,
a compactification $\bar f'=\bar{h}':\bar X'\to\bar S$ of $h':U'\to S$,
and a compactification $\bar f''=\bar{h}'':\bar X''\to\bar S$ of $h'':U''\to S$,
with $\bar X,\bar X',\bar X''\in\PSmVar(\mathbb C)$, 
$\bar D:=\bar X\backslash U\subset\bar X$ $\bar D':=\bar X'\backslash U'\subset\bar X'$, 
and $\bar D'':=\bar X''\backslash U''\subset\bar X''$ normal crossing divisors, 
such that $g_1:U''/S\to U'/S$ extend to $\bar g_1:\bar X''/\bar S\to\bar X'/\bar S$,
$g_2:U'/S\to U/S$ extend to $\bar g_2:\bar X'/\bar S\to\bar X/\bar S$, and
\begin{eqnarray*}
\hat R_S^{0CH}(g_2\circ g_1)=\hat R_S^{0CH}(g_1)\circ \hat R_S^{0CH}(g_2):
\hat R^0_{(\bar X,\bar D)/S}\to \hat R^0_{(\bar X'',\bar D'')/S} 
\end{eqnarray*}
\item[(iv)] For  
\begin{eqnarray*}
Q^*:=(\cdots\to\oplus_{\alpha\in\Lambda^n}\mathbb Z(U^n_{\alpha}/S)
\xrightarrow{(\mathbb Z(g^n_{\alpha,\beta}))}\oplus_{\beta\in\Lambda^{n-1}}\mathbb Z(U^{n-1}_{\beta}/S)\to\cdots)
\in C(\Var(\mathbb C)/S)
\end{eqnarray*}
a complex of (maybe infinite) direct sum of representable presheaves with $U^*_{\alpha}$ smooth,
we get from (i),(ii) and (iii) the map in $C(\Var(\mathbb C)/S)$
\begin{eqnarray*}
T(\hat R^{0CH},R^{0CH})(Q^*):\hat R^{0CH}(Q^*):=
(\cdots\to\oplus_{\beta\in\Lambda^{n-1}}\varinjlim_{(\bar X^{n-1}_{\beta},\bar D^{n-1}_{\beta})/S}
\hat R^0_{(\bar X^{n-1}_{\beta},\bar D^{n-1}_{\beta})/S}(\mathbb Z(U^{n-1}_{\beta}/S)) \\
\xrightarrow{(\hat R_S^{0CH}(g^n_{\alpha,\beta}))}\oplus_{\alpha\in\Lambda^n}\varinjlim_{(\bar X^n_{\alpha},\bar D^n_{\alpha})/S}
\hat R^0_{(\bar X^n_{\alpha},\bar D^n_{\alpha})/S}(\mathbb Z(U^n_{\alpha}/S))\to\cdots) 
\to R^{0CH}(Q^*),
\end{eqnarray*}
where for $(U^n_{\alpha},h^n_{\alpha})\in\Var(\mathbb C)/S$, the inductive limit run over all the compactifications 
$\bar f_{\alpha}:\bar X_{\alpha}\to\bar S$ of $h_{\alpha}:U_{\alpha}\to S$ with $\bar X_{\alpha}\in\PSmVar(\mathbb C)$
and $\bar D_{\alpha}:=\bar X_{\alpha}\backslash U_{\alpha}$ a normal crossing divisor.
For $m=(m^*):Q_1^*\to Q_2^*$ a morphism with 
\begin{eqnarray*}
Q_1^*:=(\cdots\to\oplus_{\alpha\in\Lambda^n}\mathbb Z(U^n_{1,\alpha}/S)
\xrightarrow{(\mathbb Z(g^n_{\alpha,\beta}))}\oplus_{\beta\in\Lambda^{n-1}}\mathbb Z(U^{n-1}_{1,\beta}/S)\to\cdots), \\
Q_2^*:=(\cdots\to\oplus_{\alpha\in\Lambda^n}\mathbb Z(U^n_{2,\alpha}/S)
\xrightarrow{(\mathbb Z(g^n_{\alpha,\beta}))}\oplus_{\beta\in\Lambda^{n-1}}\mathbb Z(U^{n-1}_{2,\beta}/S)\to\cdots)
\in C(\Var(\mathbb C)/S)
\end{eqnarray*}
complexes of (maybe infinite) direct sum of representable presheaves with $U^*_{1,\alpha}$ and $U^*_{2,\alpha}$ smooth,
we get again from (i),(ii) and (iii) a commutative diagram in $C(\Var(\mathbb C)/S)$
\begin{equation*}
\xymatrix{\hat R^{0CH}(Q_2^*)\ar[rr]^{T(\hat R_S^{0CH},R_S^{0CH})(Q_2^*)}\ar[d]_{\hat R_S^{0CH}(m):=(\hat R_S^{0CH}(m^*))} & \, & 
R^{0CH}(Q_2^*)\ar[d]^{(R_S^{0CH}(m^*))} \\
\hat R^{0CH}(Q_1^*)\ar[rr]^{T(\hat R_S^{0CH},R_S^{0CH})(Q_1^*)} & \, & R^{0CH}(Q_1^*)}.
\end{equation*}
\item[(v)] Let  
\begin{eqnarray*}
Q^*:=(\cdots\to\oplus_{\alpha\in\Lambda^n}\mathbb Z(U^n_{\alpha}/S)
\xrightarrow{(\mathbb Z(g^n_{\alpha,\beta}))}\oplus_{\beta\in\Lambda^{n-1}}\mathbb Z(U^{n-1}_{\beta}/S)\to\cdots)
\in C(\Var(\mathbb C)/S)
\end{eqnarray*}
a complex of (maybe infinite) direct sum of representable presheaves with $U^*_{\alpha}$ smooth,
we have by definition 
\begin{equation*}
\Gr_S^{12*}\hat R^{0CH}(Q^*)=\hat R^{CH}(Q^*)\in C(\Var(\mathbb C)^{2,smpr}/S).
\end{equation*}
\end{itemize}
\end{defi}

\begin{itemize}
\item Let $S\in\Var(\mathbb C)$
For $(h,m,m')=(h^*,m^*,m^{'*}):Q_1^*[1]\to Q_2^*$ an homotopy with $Q_1^*,Q_2^*\in C(\Var(\mathbb C)/S)$
complexes of (maybe infinite) direct sum of representable presheaves with $U^*_{1,\alpha}$ and $U^*_{2,\alpha}$ smooth,
\begin{equation*}
(\hat R_S^{0CH}(h),\hat R_S^{0CH}(m),\hat R_S^{0CH}(m'))=(\hat R_S^{0CH}(h^*),\hat R_S^{0CH}(m^*),\hat R_S^{0CH}(m^{'*})):
R^{0CH}(Q_2^*)[1]\to R^{0CH}(Q_1^*)
\end{equation*}
is an homotopy in $C(\Var(\mathbb C)/S)$ using definition \ref{R0CHhatdef} (iii). 
In particular if $m:Q_1^*\to Q_2^*$ with $Q_1^*,Q_2^*\in C(\Var(\mathbb C)/S)$
complexes of (maybe infinite) direct sum of representable presheaves with $U^*_{1,\alpha}$ and $U^*_{2,\alpha}$ smooth
is an homotopy equivalence, then $\hat R_S^{0CH}(m):\hat R^{0CH}(Q_2^*)\to\hat R^{0CH}(Q_1^*)$ is an homotopy equivalence.
\item Let $S\in\SmVar(\mathbb C)$. Let $F\in\PSh(\Var(\mathbb C)^{sm}/S)$. Consider 
\begin{eqnarray*}
q:LF:=(\cdots\to\oplus_{(U_{\alpha},h_{\alpha})\in\Var(\mathbb C)^{sm}/S}\mathbb Z(U_{\alpha}/S)
\xrightarrow{(\mathbb Z(g^n_{\alpha,\beta}))}
\oplus_{(U_{\alpha},h_{\alpha})\in\Var(\mathbb C)^{sm}/S}\mathbb Z(U_{\alpha}/S)\to\cdots)\to F
\end{eqnarray*}
the canonical projective resolution given in subsection 2.3.3.
Note that the $U_{\alpha}$ are smooth since $S$ is smooth and $h_{\alpha}$ are smooth morphism.
Definition \ref{R0CHhatdef}(iv) gives in this particular case the map in $C(\Var(\mathbb C)/S)$
\begin{eqnarray*}
T(\hat R_S^{0CH},R_S^{0CH})(\rho_S^*LF):\hat R^{0CH}(\rho_S^*LF):=
(\cdots\to\oplus_{(U_{\alpha},h_{\alpha})\in\Var(\mathbb C)^{sm}/S}\varinjlim_{(\bar X_{\alpha},\bar D_{\alpha})/S}
\hat R^0_{(\bar X_{\alpha},\bar D_{\alpha})/S}(\mathbb Z(U_{\alpha}/S)) \\
\xrightarrow{(\hat R_S^{0CH}(g^n_{\alpha,\beta}))}
\oplus_{(U_{\alpha},h_{\alpha})\in\Var(\mathbb C)^{sm}/S}\varinjlim_{(\bar X_{\alpha},\bar D_{\alpha})/S}
\hat R^0_{(\bar X_{\alpha},\bar D_{\alpha})/S}(\mathbb Z(U_{\alpha}/S))\to\cdots)
\to R^{0CH}(\rho_S^*LF),
\end{eqnarray*}
where for $(U_{\alpha},h_{\alpha})\in\Var(\mathbb C)^{sm}/S$, the inductive limit run over all the compactifications 
$\bar f_{\alpha}:\bar X_{\alpha}\to\bar S$ of $h_{\alpha}:U_{\alpha}\to S$ with $\bar X_{\alpha}\in\PSmVar(\mathbb C)$
and $\bar D_{\alpha}:=\bar X_{\alpha}\backslash U_{\alpha}$ a normal crossing divisor.
Definition \ref{R0CHhatdef}(iv) gives then by functoriality in particular, for $F=F^{\bullet}\in C(\Var(\mathbb C)^{sm}/S)$, 
the map in $C(\Var(\mathbb C)/S)$
\begin{eqnarray*}
T(\hat R_S^{0CH},R_S^{0CH})(\rho_S^*LF):\hat R^{0CH}(\rho_S^*LF)\to R^{0CH}(\rho_S^*LF).
\end{eqnarray*}

\item Let $g:T\to S$ a morphism with $T,S\in\SmVar(\mathbb C)$. Let $h:U\to S$ a smooth morphism with $U\in\Var(\mathbb C)$.
Consider the cartesian square
\begin{equation*}
\xymatrix{U_T\ar[r]^{h'}\ar[d]^{g'} & T\ar[d]^g \\
U\ar[r]^h & S}
\end{equation*}
Note that $U$ is smooth since $S$ and $h$ are smooth, and $U_T$ is smooth since $T$ and $h'$ are smooth. 
Take, see definition-proposition \ref{RCHdef0}(ii),a compactification $\bar f_0=\bar h:\bar X_0\to\bar S$ of $h:U\to S$. 
Take, see definition-proposition \ref{RCHdef0}(ii), a strict desingularization 
$\bar\epsilon:(\bar X,\bar D)\to(\bar X_0,\bar Z)$ of $(\bar X_0,\bar Z)$.
Then $\bar f'_0=\bar{g\circ h'}:\bar X_T\to\bar T$ is a compactification of $g\circ h':U_T\to S$ such that 
$g':U_T/S\to U/S$ extend to a morphism $\bar{g}'_0:\bar X_T/\bar S\to \bar X/\bar S$. 
Denote $\bar Z=\bar X_0\backslash U$ and $\bar Z'=\bar X_T\backslash U_T$.
Take, see definition-proposition \ref{RCHdef0}(ii), a strict desingularization 
$\epsilon'_{\bullet}:(\bar X',\bar D')\to(\bar X_T,\bar Z')$ of $(\bar X_T,\bar Z')$.
Denote $\bar g'=\bar g'_0\circ\epsilon'_{\bullet}:\bar X'\to\bar X$. 
We then have, see definition-proposition \ref{RCHdef0}(ii),
the following commutative diagram in $\Fun(\Delta,\Var(\mathbb C))$
\begin{equation*}
\xymatrix{U=U_{c(\bullet)}\ar[r]^j & \bar X=\bar X_{c(\bullet)} & \, & \bar D_{s_{g'}(\bullet)}\ar[ll]_{i_{\bullet}} \\
U_T=U_{T,c(\bullet)}\ar[r]^{j'}\ar[u]^{g'} & \bar X_T=X_{T,c(\bullet)}\ar[u]^{\bar{g}'_0} & \, & 
\bar{g}^{'-1}(\bar D_{s_{g'}(\bullet)})\ar[ll]_{i_{g\bullet}}\ar[u]^{(\bar{g}')'_{\bullet}} \\
U_T=U_{T,c(\bullet)}\ar[r]^{j'}\ar[u]^{g'} & \bar X'=X'_{c(\bullet)}\ar[u]^{\bar{g}'} & \bar D'_{\bullet}\ar[l]_{i'_{\bullet}} & 
\bar{g}^{'-1}(\bar D_{s_{g'}(\bullet)})\ar[l]_{i''_{g'\bullet}}\ar[u]^{\epsilon'_{\bullet}}:i'_{g\bullet}}
\end{equation*}
We then consider the following map in $C(\Var(\mathbb C)/T)$, 
\begin{eqnarray*}
T(g,\hat R^{0CH})(\mathbb Z(U/S)):g^*\hat R^0_{(\bar X,\bar D)/S}(\mathbb Z(U/S)) \\
\xrightarrow{:=}g^*\Cone(\mathbb Z(i_{\bullet}):
(\mathbb Z^{tr}(\bar D_{\bullet}/S),u_{IJ})\to\mathbb Z^{tr}(\bar X/S))(-d_X)[-2d_X] \\ 
\xrightarrow{=} \\
\Cone(\mathbb Z(i_{g\bullet}):(\mathbb Z^{tr}(\bar g^{-1}(D_{s_g(\bullet)})/T),u_{IJ})
\to\mathbb Z^{tr}(\bar X_T/T))(-d_X)[-2d_X] \\
\xrightarrow{(\mathbb Z(i''_{g\bullet}),[\Gamma_{\epsilon'}]^t)} \\
\Cone(\mathbb Z(i'_{\bullet}):
((\mathbb Z^{tr}(\bar D'_{\bullet}/T),u_{IJ})\to\mathbb Z^{tr}((\bar X'/T)))(-d_{X'})[-2d_{X'}] \\
\xrightarrow{=:}\hat R^0_{(\bar X',\bar D')/T}(\mathbb Z(U_T/T))
\end{eqnarray*}
For  
\begin{eqnarray*}
Q^*:=(\cdots\to\oplus_{\alpha\in\Lambda^n}\mathbb Z(U^n_{\alpha}/S)
\xrightarrow{(\mathbb Z(g^n_{\alpha,\beta}))}\oplus_{\beta\in\Lambda^{n-1}}\mathbb Z(U^{n-1}_{\beta}/S)\to\cdots)
\in C(\Var(\mathbb C)/S)
\end{eqnarray*}
a complex of (maybe infinite) direct sum of representable presheaves with $h^n_{\alpha}:U^n_{\alpha}\to S$ smooth,
we get the map in $C(\Var(\mathbb C)/T)$
\begin{eqnarray*}
T(g,\hat R^{0CH})(Q^*):g^*\hat R^{0CH}(Q^*)=
(\cdots\to\oplus_{\alpha\in\Lambda^n}\varinjlim_{(\bar X^n_{\alpha},\bar D^n_{\alpha})/S}
g^*\hat R^0_{(\bar X^n_{\alpha},\bar D^n_{\alpha})/S}(\mathbb Z(U^n_{\alpha}/S))\to\cdots) \\
\xrightarrow{(T(g,\hat R^{0CH})(\mathbb Z(U^n_{\alpha}/S)))} 
(\cdots\to\oplus_{\alpha\in\Lambda^n}\varinjlim_{(\bar X^{n'}_{\alpha},\bar D^{n'}_{\alpha})/T}
\hat R_{(\bar X^{n'}_{\alpha},\bar D^{n'}_{\alpha})/T}(\mathbb Z(U^n_{\alpha,T}/S))\to\cdots)=:\hat R^{CH}(g^*Q^*)
\end{eqnarray*}
together with the commutative diagram in $C(\Var(\mathbb C)/T)$
\begin{equation*}
\xymatrix{g^*\hat R^{0CH}(Q^*)\ar[rrr]^{T(g,\hat R^{0CH})(Q^*)}\ar[d]_{g^*T(\hat R_S^{0CH},R_S^{0CH})(Q^*)} 
& \, & \, & \hat R^{0CH}(g^*Q^*)\ar[d]^{T(\hat R_T^{0CH},R_T^{0CH})(g^*Q)} \\
g^*R^{0CH}(Q^*)\ar[rrr]^{T(g,R^{0CH})(Q^*)} & \, & \, & R^{0CH}(g^*Q^*)}.
\end{equation*}
Let $F\in\PSh(\Var(\mathbb C)^{sm}/S)$. Consider 
\begin{eqnarray*}
q:LF:=(\cdots\to\oplus_{(U_{\alpha},h_{\alpha})\in\Var(\mathbb C)^{sm}/S}\mathbb Z(U_{\alpha}/S)\to\cdots)\to F
\end{eqnarray*}
the canonical projective resolution given in subsection 2.3.3.
We then get in particular the map in $C(\Var(\mathbb C)/T)$
\begin{eqnarray*}
T(g,\hat R^{0CH})(\rho_S^*LF):g^*\hat R^{0CH}(\rho_S^*LF)= \\
(\cdots\to\oplus_{(U_{\alpha},h_{\alpha})\in\Var(\mathbb C)^{sm}/S}\varinjlim_{(\bar X_{\alpha},\bar D_{\alpha})/S}
g^*\hat R_{(\bar X_{\alpha},\bar D_{\alpha})/S}(\mathbb Z(U_{\alpha}/S))\to\cdots) 
\xrightarrow{(T(g,\hat R^{0CH})(\mathbb Z(U_{\alpha}/S)))} \\ 
(\cdots\to\oplus_{(U_{\alpha},h_{\alpha})\in\Var(\mathbb C)^{sm}/S}\varinjlim_{(\bar X'_{\alpha},\bar D'_{\alpha})/T}
\hat R_{(\bar X'_{\alpha},\bar D'_{\alpha})/T}(\mathbb Z(U_{\alpha,T}/S))\to\cdots)=:\hat R^{0CH}(\rho_T^*g^*LF), 
\end{eqnarray*}
and by functoriality, we get in particular for $F=F^{\bullet}\in C(\Var(\mathbb C)^{sm}/S)$, 
the map in $C(\Var(\mathbb C)/T)$
\begin{eqnarray*}
T(g,\hat R^{0CH})(\rho_S^*LF):g^*\hat R^{0CH}(\rho_S^*LF)\to\hat R^{0CH}(\rho_T^*g^*LF)
\end{eqnarray*}
together with the commutative diagram in $C(\Var(\mathbb C)/T)$
\begin{equation*}
\xymatrix{g^*\hat R^{0CH}(\rho_S^*LF)\ar[rrr]^{T(g,\hat R^{0CH})(\rho_S^*LF)}\ar[d]_{g^*T(\hat R_S^{0CH},R_S^{0CH})(\rho_S^*LF)} 
& \, & \, & \hat R^{0CH}(\rho_T^*g^*LF)\ar[d]^{T(\hat R_T^{0CH},R_T^{0CH})(\rho_T^*g^*LF)} \\
g^*L\rho_{S*}\mu_{S*}R^{CH}(\rho_S^*LF)
\ar[rrr]^{T(g,R^{0CH})(\rho_S^*LF)} & \, & \, & R^{CH}(\rho_T^*g^*LF)}.
\end{equation*}

\item Let $S_1,S_2\in\SmVar(\mathbb C)$ and $p:S_1\times S_2\to S_1$ the projection. 
Let $h:U\to S_1$ a smooth morphism with $U\in\Var(\mathbb C)$. Consider the cartesian square
\begin{equation*}
\xymatrix{U\times S_2\ar[r]^{h\times I}\ar[d]^{p'} & S_1\times S_2\ar[d]^p \\
U\ar[r]^h & S_1}
\end{equation*}
Take, see definition-proposition \ref{RCHdef0}(i),a compactification $\bar f_0=\bar h:\bar X_0\to\bar S_1$ of $h:U\to S_1$.
Then $\bar f_0\times I:\bar X_0\times S_2\to\bar S_1\times S_2$ is a compactification of $h\times I:U\times S_2\to S_1\times S_2$ 
and $p':U\times S_2\to U$ extend to $\bar{p}'_0:=p_{X_0}:\bar X_0\times S_2\to\bar X_0$. Denote $Z=X_0\backslash U$. 
Take see theorem \ref{desVar}(i), a strict desingularization 
$\bar\epsilon:(\bar X,\bar D)\to(\bar X_0,\bar Z)$ of the pair $(\bar X_0,\bar Z)$.  
We then have the commutative diagram (\ref{RCHdiap}) in $\Fun(\Delta,\Var(\mathbb C))$ whose squares are cartesian
\begin{equation*}
\xymatrix{U=U_{c(\bullet)}\ar[r]^j & \bar X & \bar D_{\bullet}\ar[l]_{i_{\bullet}} \\
U\times S_2=(U\times S_2)_{c(\bullet)}\ar[r]^{j\times I}\ar[u]^g & \bar X\times S_2\ar[u]^{\bar{p}':=p_{\bar X}} & 
\bar D_{\bullet}\times S_2\ar[l]_{i'_{\bullet}}\ar[u]^{\bar{p'}'_{\bullet}}}
\end{equation*}
Then the map in $C(\Var(\mathbb C)/S_1\times S_2)$
\begin{eqnarray*}
T(p,\hat R^{0CH})(\mathbb Z(U/S_1)):p^*\hat R^0_{(\bar X,\bar D)/S_1}(\mathbb Z(U/S_1))\xrightarrow{\sim} 
\hat R^0_{(\bar X\times S_2,\bar D_{\bullet}\times S_2)/S_1\times S_2}(\mathbb Z(U\times S_2/S_1\times S_2))
\end{eqnarray*}
is an isomorphism.
Hence, for $Q^*\in C(\Var(\mathbb C)/S_1)$ a complex of (maybe infinite) direct sum of representable presheaves of smooth morphism,
the map in $C(\Var(\mathbb C)/S_1\times S_2)$
\begin{eqnarray*}
T(p,\hat R^{0CH})(Q^*):p^*\hat R^{0CH}(Q^*)\xrightarrow{\sim}\hat R^{0CH}(p^*Q^*)
\end{eqnarray*}
is an isomorphism.
In particular, for $F\in C(\Var(\mathbb C)^{sm}/S_1)$ the map in $C(\Var(\mathbb C)/S_1\times S_2)$
\begin{eqnarray*}
T(p,\hat R^{0CH})(\rho_{S_1}^*LF):p^*\hat R^{0CH}(\rho_{S_1}^*LF)\xrightarrow{\sim}\hat R^{0CH}(\rho_{S_1\times S_2}^*p^*LF)
\end{eqnarray*}
is an isomorphism.

\item Let $h_1:U_1\to S$, $h_2:U_2\to S$ two morphisms with $U_1,U_2,S\in\Var(\mathbb C)$, $U_1,U_2$ smooth.
Denote by $p_1:U_1\times_SU_2\to U_1$ and $p_2:U_1\times_S U_2\to U_2$ the projections.
Take, see definition-proposition \ref{RCHdef0}(i)),
a compactification $\bar f_{10}=\bar{h}_1:\bar X_{10}\to\bar S$ of $h_1:U_1\to S$
and a compactification $\bar f_{20}=\bar{h}_2:\bar X_{20}\to\bar S$ of $h_2:U_2\to S$. Then, 
\begin{itemize}
\item $\bar f_{10}\times\bar f_{20}:\bar X_{10}\times_{\bar S}\bar X_{20}\to S$ 
is a compactification of $h_1\times h_2:U_1\times_SU_2\to S$.
\item $\bar p_{10}:=p_{X_{10}}:\bar X_{10}\times_{\bar S} \bar X_{20}\to\bar X_{10}$ 
is a compactification of $p_1:U_1\times_SU_2\to U_1$.
\item $\bar p_{20}:=p_{X_{20}}:\bar X_{10}\times_{\bar S}\bar X_{20}\to\bar X_{20}$ 
is a compactification of $p_2:U_1\times_SU_2\to U_2$.
\end{itemize}
Denote $\bar Z_1=\bar X_{10}\backslash U_1$ and $\bar Z_2=\bar X_{20}\backslash U_2$.
Take, see theorem \ref{desVar}(i), a strict desingularization 
$\bar\epsilon_1:(\bar X_1,\bar D)\to(\bar X_{10},Z_1)$ of the pair $(\bar X_{10},\bar Z_1)$
and a  strictdesingularization 
$\bar\epsilon_2:(\bar X_2,\bar E)\to(\bar X_{20},Z_2)$ of the pair $(\bar X_{20},\bar Z_2)$. 
Take then a strict desingularization
\begin{equation*}
\bar\epsilon_{12}:((\bar X_1\times_{\bar S}\bar X_2)^N,\bar F)\to
(\bar X_1\times_{\bar S}\bar X_2,(D\times_{\bar S}\bar X_2)\cup(\bar X_1\times_{\bar S}\bar E)) 
\end{equation*}
of the pair $(\bar X_1\times_{\bar S}\bar X_2,(\bar D\times_{\bar S}\bar X_2)\cup(\bar X_1\times_{\bar S}\bar E))$. 
We have then the following commutative diagram 
\begin{equation*}
\xymatrix{\, & \bar X_1\ar[r]^{\bar f_1} & \bar S \\
\, & \bar X_1\times_{\bar S}\bar X_2\ar[r]^{\bar p_1}\ar[u]^{\bar p_2} &\bar X_2\ar[u]^{\bar f_2} \\
(\bar X_1\times_{\bar S}\bar X_2)^N\ar[ru]^{\bar\epsilon_{12}}\ar[rru]^{(\bar p_1)^N}\ar[ruu]^{(\bar p_2)^N} & \, & \,}
\end{equation*}
and
\begin{itemize}
\item $\bar f_1\times\bar f_2:\bar X_1\times_{\bar S}\bar X_2\to\bar S$ 
is a compactification of $h_1\times h_2:U_1\times_SU_2\to S$.
\item $(\bar p_1)^N:=\bar p_1\circ\epsilon_{12}:(\bar X_1\times_{\bar S}\bar X_2)^N\to\bar X_1$ is a compactification of 
$p_1:U_1\times_SU_2\to U_1$.
\item $(\bar p_2)^N:=\bar p_2\circ\epsilon_{12}:(\bar X_1\times_{\bar S}\bar X_2)^N\to\bar X_2$ is a compactification of 
$p_2:U_1\times_SU_2\to U_2$.
\end{itemize}
We have then the morphism in $C(\Var(\mathbb C)/S)$
\begin{eqnarray*}
T(\otimes,\hat R_S^{0CH})(\mathbb Z(U_1/S),\mathbb Z(U_2/S)):=\hat R_S^{0CH}(p_1)\otimes\hat R_S^{0CH}(p_2): \\
\hat R^0_{(\bar X_1,\bar D)/S}(\mathbb Z(U_1/S))\otimes\hat R^0_{(X_2,E))/S}(\mathbb Z(U_2/S)) 
\xrightarrow{\sim}\hat R^0_{(\bar X_1\times_{\bar S}\bar X_2)^N,\bar F)/S}(\mathbb Z(U_1\times_S U_2/S))
\end{eqnarray*}
For   
\begin{eqnarray*}
Q_1^*:=(\cdots\to\oplus_{\alpha\in\Lambda^n}\mathbb Z(U^n_{1,\alpha}/S)
\xrightarrow{(\mathbb Z(g^n_{\alpha,\beta}))}\oplus_{\beta\in\Lambda^{n-1}}\mathbb Z(U^{n-1}_{1,\beta}/S)\to\cdots), \\
Q_2^*:=(\cdots\to\oplus_{\alpha\in\Lambda^n}\mathbb Z(U^n_{2,\alpha}/S)
\xrightarrow{(\mathbb Z(g^n_{\alpha,\beta}))}\oplus_{\beta\in\Lambda^{n-1}}\mathbb Z(U^{n-1}_{2,\beta}/S)\to\cdots)
\in C(\Var(\mathbb C)/S)
\end{eqnarray*}
complexes of (maybe infinite) direct sum of representable presheaves with $U^*_{\alpha}$ smooth, 
we get the morphism in $C(\Var(\mathbb C)/S)$
\begin{eqnarray*}
T(\otimes,\hat R_S^{0CH})(Q^*_1,Q^*_2):\hat R^{0CH}(Q^*_1)\otimes R^{0CH}(Q^*_2) 
\xrightarrow{(T(\otimes,\hat R_S^{CH})(\mathbb Z(U_{1,\alpha}^m),\mathbb Z(U_{2,\beta}^n))}
\hat R^{0CH}(Q^*_1\otimes Q^*_2))
\end{eqnarray*},
together with the commutative diagram in $C(\Var(\mathbb C)/S)$
\begin{equation*}
\xymatrix{\hat R^{0CH}(Q^*_1)\otimes R^{0CH}(Q^*_2)
\ar[rrr]^{T(\otimes,\hat R_S^{0CH})(Q_1^*,Q_2)}
\ar[d]_{T(\hat R_S^{0CH},R_S^{CH})(Q_1^*)\otimes T(\hat R_S^{0CH},R_S^{0CH})(Q_2^*)} 
& \, & \, & \hat R^{0CH}(Q^*_1\times Q_2^*)\ar[d]^{T(\hat R_S^{0CH},R_S^{0CH})(Q^*_1\otimes Q^*_2)} \\
R^{0CH}(Q_1^*)\otimes R^{0CH}(Q_2^*)
\ar[rrr]^{T(\otimes,R_S^{0CH})(Q_1^*,Q_2)} & \, & \, & R^{0CH}(Q_1^*\otimes Q_2)}.
\end{equation*}
For $F_1,F_2\in C(\Var(\mathbb C)^{sm}/S)$, we get in particular the morphism in $C(\Var(\mathbb C)/S)$
\begin{eqnarray*}
T(\otimes,R_S^{0CH})(\rho_S^*LF_1,\rho_S^*LF_2):R^{0CH}(\rho_S^*LF_1)\otimes R^{0CH}(\rho_S^*LF_2) 
\to R^{0CH}(\rho_S^*(LF_1\otimes LF_2))
\end{eqnarray*}
together with the commutative diagram in $C(\Var(\mathbb C)/S)$
\begin{equation*}
\xymatrix{\hat R^{0CH}(\rho_S^*LF_1)\otimes R^{0CH}(\rho_S^*LF_2)
\ar[rrr]^{T(\otimes,\hat R_S^{0CH})(\rho_S^*LF_1,\rho_S^*LF_2)}
\ar[d]_{T(\hat R_S^{0CH},R_S^{0CH})(\rho_S^*LF_1)\otimes T(\hat R_S^{0CH},R_S^{0CH})(\rho_S^*LF_2)} 
& \, & \, & \hat R^{0CH}(\rho_S^*LF_1\times\rho_S^*LF_2)
\ar[d]^{T(\hat R_S^{0CH},R_S^{0CH})(\rho_S^*LF_1\otimes\rho_S^*LF_2)} \\
R^{0CH}(\rho_S^*LF_1)\otimes R^{0CH}(\rho_S^*LF_2)
\ar[rrr]^{T(\otimes,R_S^{0CH})(\rho_S^*LF_1,\rho_S^*LF_2)} 
& \, & \, & R^{0CH}(\rho_S^*LF_1\times\rho_S^*LF_2)}.
\end{equation*}

\end{itemize}

\subsection{The derived categories of filtered complexes of presheaves on a site or 
of filtered complexes of presheaves of modules on a ringed topos}

\begin{defi}
Let $\mathcal S\in\Cat$ a site endowed with topology $\tau$.
\begin{itemize}
\item[(i)]We denote by $D(\mathcal S):=\Ho_{\Top}C(\mathcal S)$
the localization of the category of complexes of presheaves on $S$ with respect to top local equivalence and by 
$D(\tau):C(\mathcal S)\to D(\mathcal S)$ the localization functor.
\item[(ii)]We denote for $r=1,\ldots\infty$, resp. $r=(1,\ldots\infty)^2$, 
\begin{equation*}
D_{fil,r}(\mathcal S):=K_{fil,r}(\mathcal S)([E_1]^{-1}) \; , \; 
D_{2fil,r,r'}(\mathcal S):=K_{2fil,r,r'}(\mathcal S)([E_1]^{-1}),
\end{equation*}
the localizations of the category of filtered complexes of presheaves on $\mathcal S$
whose filtration is biregular modulo $r$-filtered homotopies 
with respect to the classes of filtered $\tau$ local equivalence $[E_1]$. 
Note that the classes of filtered $\tau$ local equivalence constitute a right multiplicative system.
By definition, if $m:(G_1,F)\to(G_2,F)$ with $(G_1,F),(G_2,F)\in C_{fil}(\mathcal S)$ is an $r$-filtered $\tau$ local equivalence
then $m:=D(\tau)(m):(G_1,F)\to(G_2,F)$ is an isomorphism in $D_{fil,r}(\mathcal S)$.
By definition, we have sequences of functors
\begin{eqnarray*}
C_{fil}(\mathcal S)\to K_{fil}(\mathcal S)\to D_{fil}(\mathcal S)\to D_{fil,2}(\mathcal S)\to\cdots\to D_{fil,\infty}(\mathcal S).
\end{eqnarray*}
and commutative diagrams of functors
\begin{equation*}
\xymatrix{K_{fil}(\mathcal S)\ar[r]\ar[d] & D_{fil}(\mathcal S)\ar[d] \\
K_{fil,2}(\mathcal S)\ar[r] & D_{fil,2}(\mathcal S)} \; , \;
\xymatrix{K_{fil,r}(\mathcal S)\ar[r]\ar[d] & D_{fil,r}(\mathcal S)\ar[d] \\
K_{fil,r+1}(\mathcal S)\ar[r] & D_{fil,r}(\mathcal S)}.
\end{equation*}
Then, for $r=1$, $K_{fil}(\mathcal S)$ and $D_{fil}(\mathcal S)$ are in the canonical way triangulated categories.
However, for $r>1$, the categories $K_{fil,r}(\mathcal S)$ and $D_{fil,r}(\mathcal S)$ together with the canonical triangles
does NOT satisfy the 2 of 3 axiom of triangulated categories.
\end{itemize}
\end{defi}

\begin{defi}
Let $(\mathcal S,O_S)\in\RCat$ where $\mathcal S\in\Cat$ is a site endowed with topology $\tau$.
\begin{itemize}
\item[(i)]We denote by $D_{O_S}(\mathcal S):=\Ho_{\Top}C_{O_S}(\mathcal S)$
the localization of the category of complexes of presheaves on $S$ with respect to top local equivalence and by 
$D(\tau):C_{O_S}(\mathcal S)\to D_{O_S}(\mathcal S)$ the localization functor.
\item[(ii)]We denote for $r=1,\ldots\infty$, resp. $r=(1,\ldots\infty)^2$, 
\begin{equation*}
D_{O_Sfil,r}(\mathcal S):=K_{O_Sfil,r}(\mathcal S)([E_1]^{-1}) \; , \; 
D_{O_S2fil,r,r'}(\mathcal S):=K_{O2fil,r,r'}(\mathcal S)([E_1]^{-1}),
\end{equation*}
the localizations of the category of filtered complexes of presheaves on $\mathcal S$
whose filtration is biregular modulo $r$-filtered homotopies 
with respect to the classes of filtered $\tau$ local equivalence $[E_1]$. 
Note that the classes of filtered $\tau$ local equivalence constitute a right multiplicative system.
By definition, if $m:(G_1,F)\to(G_2,F)$ with $(G_1,F),(G_2,F)\in C_{O_Sfil}(\mathcal S)$ 
is an $r$-filtered $\tau$ local equivalence
then $m:=D(\tau)(m):(G_1,F)\to(G_2,F)$ is an isomorphism in $D_{O_Sfil,r}(\mathcal S)$.
By definition, we have sequences of functors
\begin{eqnarray*}
C_{O_Sfil}(\mathcal S)\to K_{O_Sfil}(\mathcal S)\to 
D_{O_Sfil}(\mathcal S)\to D_{O_Sfil,2}(\mathcal S)\to\cdots\to D_{O_Sfil,\infty}(\mathcal S).
\end{eqnarray*}
and commutative diagrams of functors
\begin{equation*}
\xymatrix{K_{O_Sfil}(\mathcal S)\ar[r]\ar[d] & D_{O_Sfil}(\mathcal S)\ar[d] \\
K_{O_Sfil,2}(\mathcal S)\ar[r] & D_{fil,2}(\mathcal S)} \; , \;
\xymatrix{K_{fil,r}(\mathcal S)\ar[r]\ar[d] & D_{O_Sfil,r}(\mathcal S)\ar[d] \\
K_{O_Sfil,r+1}(\mathcal S)\ar[r] & D_{O_Sfil,r}(\mathcal S)}.
\end{equation*}
Then, for $r=1$, $K_{O_Sfil}(\mathcal S)$ and $D_{O_Sfil}(\mathcal S)$ are in the canonical way triangulated categories.
However, for $r>1$, the categories $K_{O_Sfil,r}(\mathcal S)$ and $D_{O_Sfil,r}(\mathcal S)$ together with the canonical triangles
does NOT satisfy the 2 of 3 axiom of triangulated categories.
\end{itemize}
\end{defi}

Let $f:\mathcal T\to\mathcal S$ a morphism of presite, 
with $\mathcal S,\mathcal T\in\Cat$ endowed with a topology $\tau$.
If $f$ is a morphism of site, the adjonctions
\begin{equation*}
(f^*,f_*)=(f^{-1},f_*):C(\mathcal S)\leftrightarrows C(\mathcal T),  \;
(f^*,f_*)=(f^{-1},f_*):C_{(2)fil}(\mathcal S)\leftrightarrows C_{(2)fil}(\mathcal T). 
\end{equation*}
are Quillen adjonctions.
They induces respectively in the derived categories, 
for $r=(1,\ldots,\infty)$, resp. $r=(1,\ldots,\infty)$ (note that $f^*$ derive trivially)
\begin{equation*}
(f^*,Rf_*):D(\mathcal S)\leftrightarrows D(\mathcal T) \; , \; 
(f^*,Rf_*):D_{fil,r}(\mathcal S)\leftrightarrows D_{fil,r}(\mathcal T).
\end{equation*}
For $F^{\bullet}\in C(\mathcal S)$, we have the adjonction maps
\begin{equation*}
\ad(f^*,f_*)(F^{\bullet}):F^{\bullet}\to f_*f^*F^{\bullet} \; , \; \ad(f^*,f_*)(F^{\bullet}):f^*f_*F^{\bullet}\to F^{\bullet},
\end{equation*}
induces in the derived categories, for $(M,F)\in D_{fil}(\mathcal S)$ and $(N,F)\in D_{fil}(\mathcal T)$, the adjonction maps
\begin{equation*}
\ad(f^*,Rf_*)(M):(M,F)\to Rf_*f^*(M,F) \; , \; \ad(f^*,Rf_*)(N,F):f^*Rf_*(N,F)\to (N,F).
\end{equation*}

For a commutative diagram of sites : 
\begin{equation*}
D=\xymatrix{ 
\mathcal Y\ar[r]^{g_2}\ar[d]^{f_2} & \mathcal X\ar[d]^{f_1} \\
\mathcal T\ar[r]^{g_1} & \mathcal S},
\end{equation*}
with $\mathcal Y,\mathcal T,\mathcal S,\mathcal X\in\Cat$ with topology $\tau_Y,\tau_T,\tau_S,\tau X$,
the maps, for $F\in C(\mathcal X)$, 
\begin{equation*}
T(D)(F): g_1^*f_{1*} F\to f_{2*}g_2^*F
\end{equation*}
induce in the derived category the maps in $D_{fil,r}(\mathcal T)$,
given by, for $(G,F)\in D_{fil,r}(\mathcal X)$ with $(G,F)=D(\tau_X,r)((G,F))$,
\begin{equation*}
\xymatrix{
g_1^*Rf_{1*}(M,F)\ar[rr]^{T(D)(M,F)} & \, & Rf_{2*}g_2^*(M,F) \\
g_1^*f_{1*}(E(G,F))\ar[rr]^{k\circ T(D)(E(G,F))}\ar[u]^{=} & \, & f_{2*}E(g_2^*(E(G,F)))\ar[u]^{=}}.
\end{equation*}

Let $\mathcal S\in\Cat$ a site with topology $\tau$. 
The tensor product of complexes of abelian groups and the internal hom of presheaves on $\mathcal S$
\begin{equation*}
((\cdot\otimes\cdot),\mathcal Hom^{\bullet}(\cdot,\cdot)):C(\mathcal S)^2\to C(\mathcal S),
\end{equation*}
is a Quillen adjonction which induces in the derived category
\begin{equation*}
((\cdot\otimes^L\cdot),R\mathcal Hom^{\bullet}(\cdot,\cdot)):D_{fil,r}(\mathcal S)^2\to D_{fil,r}(\mathcal S), \;
R\mathcal Hom^{\bullet}((M,W),(N,W))=\mathcal Hom^{\bullet}((Q,W),E(G,F)),
\end{equation*}
where, $Q$ is projectively cofibrant such that $M=D(\tau)(Q^{\bullet})$ and $G$ such that $N=D(\tau)(G)$.

Let $i:Z\hookrightarrow S$ a closed embedding, with $S,Z\in\Top$. 
Denote by $j:S\backslash Z\hookrightarrow S$ the open embedding of the complementary subset.
The adjonction
\begin{equation*}
(i_*,i^!):=(i_*,i^{\bot}):C(Z)\to C(S), \; \mbox{with \, in \, this \, case} \; i^!F:=\ker(F\to j_*j^*F) 
\end{equation*}
is a Quillen adjonction. Since $i^!$ preserve monomorphisms, we have also Quillen adjonctions
\begin{equation*}
(i_*,i^!):C_{(2)fil}(Z)\to C_{(2)fil}(S), \; \mbox{with} \; i^!(G,F)=(i^!G,F).
\end{equation*}
which induces in the derived category ($i_*$ derive trivially)
\begin{equation*}
(i_*,Ri^!):D_{(2)fil}(Z)\to D_{(2)fil}(S), \; \mbox{with} \; Ri^!(G,F)=i^!E(G,F).
\end{equation*}

The 2-functor $S\in\Top\mapsto D(S)$ obviously satisfy the localization property, that is
for $i:Z\hookrightarrow S$ a closed embedding with $Z,S\in\Top$, 
denote by $j:S\backslash Z\hookrightarrow S$ the open complementary subset, 
we have for $K\in D(S)$ a distinguish triangle in $D(S)$
\begin{equation*}
j_{\sharp}j^*K\xrightarrow{\ad(j_{\sharp},j^*)(K)}K\xrightarrow{\ad(i^*,i_*)(K)}i_*i^*K\to j_{\sharp}j^*K[1]
\end{equation*}
equivalently, 
\begin{itemize}
\item the functor 
\begin{equation*}
(i^*,j^*): D(S)\xrightarrow{\sim} D(Z)\times D(S\backslash Z)
\end{equation*}
is conservative,
\item  and for $K\in C(Z)$, the adjonction map $\ad(i^*,i_*)(K):i^*i_*K\to K$ is an equivalence top local, 
hence for $K\in D(S)$, the induced map in the derived category 
\begin{equation*}
\ad(i^*,i_*)(K):i^*i_*K\xrightarrow{\sim} K 
\end{equation*}
is an isomorphism.
\end{itemize}

\section{Triangulated category of motives}

\subsection{Definition and the six functor formalism}

The category of motives is obtained by inverting the $(\mathbb A^1_S,et)$ equivalence.
Hence the $\mathbb A^1_S$ local complexes of presheaves plays a key role.
\begin{defi}
The derived category of motives of complex algebraic varieties over $S$ is
the category
\begin{equation*}
\DA(S):=\Ho_{\mathbb A_S^1,et}(C(\Var(\mathbb C)^{sm}/S)),
\end{equation*}
which is the localization of the category of complexes of presheaves on $\Var(\mathbb C)^{sm}/S$
with respect to $(\mathbb A_S^1,et)$ local equivalence and we denote by
\begin{equation*}
D(\mathbb A_S^1,et):=D(\mathbb A_S^1)\circ D(et):C(\Var(\mathbb C)^{sm}/S)\to\DA(S)
\end{equation*}
the localization functor. We have 
$\DA^-(S):=D(\mathbb A_S^1,et)(\PSh(\Var(\mathbb C)^{sm}/S,C^-(\mathbb Z)))\subset\DA(S)$
the full subcategory consisting of bounded above complexes.
\end{defi}

\begin{defi}
The stable derived category of motives of complex algebraic varieties over $S$ is the category
\begin{equation*}
\DA_{st}(S):=\Ho_{\mathbb A_S^1,et}(C_{\Sigma}(\Var(\mathbb C)^{sm}/S)),
\end{equation*}
which is the localization of the category of 
$\mathbb G_{mS}$-spectra ($\Sigma F^{\bullet}=F^{\bullet}\otimes\mathbb G_{mS}$) 
of complexes of presheaves on $\Var(\mathbb C)^{sm}/S$
with respect to $(\mathbb A_S^1,et)$ local equivalence. The functor 
\begin{equation*}
\Sigma^{\infty}:C(\Var(\mathbb C)^{sm}/S)\hookrightarrow
C_{\Sigma}(\Var(\mathbb C)^{sm}/S)
\end{equation*}
induces the functor $\Sigma^{\infty}:\DA(S)\to\DA_{st}(S)$.
\end{defi}

We have all the six functor formalism by \cite{C.D}. 
We give a list of the operation we will use : 

\begin{itemize}

\item For $f:T\to S$ a morphism with $S,T\in\Var(\mathbb C)$, the adjonction
\begin{equation*}
(f^*,f_*):C(\Var(\mathbb C)^{sm}/S)\leftrightarrows C(\Var(\mathbb C)^{sm}/T)
\end{equation*}
is a Quillen adjonction which induces in the derived categories ($f^*$ derives trivially),
$(f^*,Rf_*):\DA(S)\leftrightarrows\DA(T)$.

\item For $h:V\to S$ a smooth morphism with $V,S\in\Var(\mathbb C)$, the adjonction
\begin{equation*}
(h_{\sharp},h^*):C(\Var(\mathbb C)^{sm}/V)\leftrightarrows C(\Var(\mathbb C)^{sm}/S)
\end{equation*}
is a Quillen adjonction which induces in the derived categories ($h^*$ derive trivially) 
$(Lh_{\sharp},h^*)=:\DA(V)\leftrightarrows\DA(S)$.

\item For $i:Z\hookrightarrow S$ a closed embedding, with $Z,S\in\Var(\mathbb C)$,
\begin{equation*}
(i_*,i^!):=(i_*,i^{\bot}):C(\Var(\mathbb C)^{sm}/Z)\leftrightarrows C(\Var(\mathbb C)^{sm}/S)
\end{equation*}
is a Quillen adjonction, which induces in the derived categories ($i_*$ derive trivially) 
$(i_*,Ri^!):\DA(Z)\leftrightarrows\DA(S)$. The fact that $i_*$ derive trivially 
(i.e. send $(\mathbb A^1,et)$ local equivalence to $(\mathbb A^1,et)$ local equivalence is proved in \cite{AyoubT}.

\item For $S\in\Var(\mathbb C)$, the adjonction given by
the tensor product of complexes of abelian groups and the internal hom of presheaves
\begin{equation*}
((\cdot\otimes\cdot),\mathcal Hom^{\bullet}(\cdot,\cdot)):
C(\Var(\mathbb C)^{sm}/S)^2\to C(\Var(\mathbb C)^{sm}/S),
\end{equation*}
is a Quillen adjonction, which induces in the derived category
\begin{equation*},
((\cdot\otimes^L\cdot),R\mathcal Hom^{\bullet}(\cdot,\cdot)):\DA(S)^2\to\DA(S),
\end{equation*} 
\begin{itemize}
\item Let $M,N\in\DA(S)$, $Q^{\bullet}$ projectively cofibrant such that $M=D(\mathbb A^1,et)(Q^{\bullet})$, 
and $G^{\bullet}$ be $\mathbb A^1$ local for the etale topology such that $N=D(\mathbb A^1,et)(G^{\bullet})$. Then, 
\begin{equation}
R\mathcal Hom^{\bullet}(M,N)=\mathcal Hom^{\bullet}(Q^{\bullet},E(G^{\bullet}))\in\DA(S).
\end{equation}
This is well defined since if $s:Q_1\to Q_2$ is a etale local equivalence, 
\begin{equation*}
\mathcal Hom(s,E(G)):\mathcal Hom(Q_1,E(G))\to\mathcal Hom(Q_2,E(G))
\end{equation*}
is a etale local equivalence for $1\leq i\leq l$. 
\end{itemize}

\item For a commutative diagram in $\Var(\mathbb C)$ : 
\begin{equation*}
D=\xymatrix{ 
Y\ar[r]^{g_2}\ar[d]^{f_2} & X\ar[d]^{f_1} \\
T\ar[r]^{g_1} & S},
\end{equation*}
and $F\in C(\Var(\mathbb C)^{sm}/X)$, the transformation map $T(D)(F): g_1^*f_{1*} F\to f_{2*}g_2^*F$ induces 
in the derived category, for $M\in\DA(X)$, $M=D(\mathbb A^1,et)(F)$ with $F$ $\mathbb A^1$ local for the etale topology,
\begin{equation*}
\xymatrix{
g_1^*Rf_{1*} M\ar[rr]^{T(D)(M)} & \, & Rf_{2*}g_2^*M \\
g_1^*f_{1*}E(F)\ar[rr]^{k\circ c\circ T(D)(E(F))}\ar[u]^{=} & \, & f_{2*}E(C_*(g_2^*E(F)))\ar[u]_{=}}
\end{equation*}
If $D$ is cartesian with $f_1=f$, $g_1=g$ $f_2=f':X_T\to T$, $g':X_T\to X$, we denote
\begin{itemize}
\item $T(D)(F)=:T(f,g)(F): g^*f_*F\to f'_*g^{'*}F$,   
\item $T(D)(M)=:T(f,g)(M): g^*Rf_*M\to Rf'_*g^{*'}M$.
\end{itemize}

\end{itemize}

We get from the first point 2 functors :
\begin{itemize}
\item The 2-functor $C(\Var(\mathbb C)^{sm}/\cdot):\Var(\mathbb C)\to Ab\Cat$, given by
\begin{equation*}
S\mapsto C(\Var(\mathbb C)^{sm}/S) \;,\;
 (f:T\to S)\mapsto(f^*:C(\Var(\mathbb C)^{sm}/S)\to C(\Var(\mathbb C)^{sm}/T)). 
\end{equation*}
\item The 2-functor $\DA(\cdot):\Var(\mathbb C)\to\TriCat$, given by
\begin{equation*}
S\mapsto\DA(S) \;,\; (f:T\to S)\mapsto(f^*:\DA(S)\to\DA(T)). 
\end{equation*}
\end{itemize}
The main theorem is the following :
\begin{thm}\cite{AyoubT}\cite{C.D}\label{2functDM}
The 2-functor $\DA(\cdot):\Var(\mathbb C)\to\TriCat$, given by
\begin{equation*}
S\mapsto\DA(S) \;,\; (f:T\to S)\mapsto(f^*:\DA(S)\to\DA(T)) 
\end{equation*}
is a 2-homotopic functor (\cite{AyoubT})
\end{thm}

From theorem \ref{2functDM}, we get in particular
\begin{itemize}

\item For $f:T\to S$ a morphism with $T,S\in\Var(\mathbb C)$, there by theorem \ref{2functDM} is also a pair of adjoint functor
\begin{equation*}
(f_!,f^!):\DA(S)\leftrightarrows\DA(T) 
\end{equation*}
\begin{itemize}
\item with $f_!=Rf_*$ if $f$ is proper,
\item with $f^!=f^*[d]$ if $f$ is smooth of relative dimension $d$. 
\end{itemize}
For $h:U\to S$ a smooth morphism with $U,S\in\Var(\mathbb C)$ irreducible, have, for $M\in\DA(U)$, an isomorphism  
\begin{equation}\label{If}
Lh_{\sharp}M\to h_!M[d],
\end{equation}
in $\DA(S)$.

\item The 2-functor $S\in\Var(\mathbb C)\mapsto\DA(S)$ satisfy the localization property, that is
for $i:Z\hookrightarrow X$ a closed embedding with $Z,X\in\Var(\mathbb C)$, 
denote by $j:S\backslash Z\hookrightarrow S$ the open complementary subset, 
we have for $M\in\DA(S)$ a distinguish triangle in $\DA(S)$
\begin{equation*}
j_{\sharp}j^*M\xrightarrow{\ad(j_{\sharp},j^*)(M)}M\xrightarrow{\ad(i^*,i_*)(M)}i_*i^*M\to j_{\sharp}j^*M[1]
\end{equation*}
equivalently, 
\begin{itemize}
\item the functor 
\begin{equation*}
(i^*,j^*):\DA(S)\xrightarrow{\sim}\DA(Z)\times\DA(S\backslash Z)
\end{equation*}
is conservative,
\item  and for $F\in C(\Var(\mathbb C)^{sm}/Z)$,
the adjonction map $\ad(i^*,i_*)(F):i^*i_*F\to F$ is an equivalence Zariski local, 
hence for $M\in\DA(S)$, the induced map in the derived category 
\begin{equation*}
\ad(i^*,i_*)(M):i^*i_*M\xrightarrow{\sim} M 
\end{equation*}
is an isomorphism.
\end{itemize}

\item For $f:X\to S$ a proper map, $g:T\to S$ a morphism, with $T,X,S\in\Var(\mathbb C)$, and $M\in\DA(X)$,
\begin{equation*}
T(f,g)(M):g^*Rf_*M\to Rf'_*\tilde g^{'*}M 
\end{equation*}
is an isomorphism in $\DA(T)$ if $f$ is proper.

\end{itemize}

\begin{defi}
The derived category of extended motives of complex algebraic varieties over $S$ is
the category
\begin{equation*}
\underline{\DA}(S):=\Ho_{\mathbb A_S^1,et}(C(\Var(\mathbb C)/S)),
\end{equation*}
which is the localization of the category of complexes of presheaves on $\Var(\mathbb C)/S$
with respect to $(\mathbb A_S^1,et)$ local equivalence and we denote by
\begin{equation*}
D(\mathbb A_S^1,et):=D(\mathbb A_S^1)\circ D(et):C(\Var(\mathbb C)/S)\to\underline{\DA}(S)
\end{equation*}
the localization functor. We have 
$\underline{\DA}^-(S):=D(\mathbb A_S^1,et)(\PSh(\Var(\mathbb C)/S,C^-(\mathbb Z)))\subset\underline{\DA}(S)$
the full subcategory consisting of bounded above complexes.
\end{defi}

\begin{rem}\label{isharprem}
Let $i:Z\hookrightarrow S$ a closed embedding, with $Z,S\in\Var(\mathbb C)$.
\begin{itemize}
\item[(i)] By theorem \ref{2functDM}, for $X/S=(X,h)\in\Var(\mathbb C)^{sm}/S)$,
\begin{eqnarray*}
(0,\ad(i^*,i_*)(\mathbb Z(X/S))):\Gamma_Z^{\vee}\mathbb Z(X/S)\to i_*\mathbb Z(X_Z/Z)
\end{eqnarray*}
is an equivalence $(\mathbb A^1,et)$ local.
\item[(ii)] For $X/S=(X,f)\in\Var(\mathbb C)/S)$,
\begin{eqnarray*}
(0,\ad(i^*,i_*)(\mathbb Z(X/S))):\Gamma_Z^{\vee}\mathbb Z(X/S)\to i_*\mathbb Z(X_Z/Z)
\end{eqnarray*}
is NOT an equivalence $(\mathbb A^1,et)$ local in general, since for example if $f(X)=Z\subset S$,
$\rho_{S*}\mathbb Z(X/S)=0$ but 
$D(\mathbb A^1,et)(\rho_{S*}i_*\mathbb Z(X_Z/Z)=i_*\rho_{S*}\mathbb Z(X_Z/Z))\neq 0\in\underline{\DA}(S)$,
hence it is NOT an equivalence $(\mathbb A^1,et)$ local in this case by proposition \ref{rho1}.
In particular $\underline{\DA}(S)$ dos NOT satisfy the localization property.
\item[(ii)'] For $X/Z=(X,f)\in\Var(\mathbb C)/Z)$, the inclusion
\begin{eqnarray*}
T(i_{\sharp},i_*):i_{\sharp}\mathbb Z(X/Z)\hookrightarrow i_*\mathbb Z(X/Z)
\end{eqnarray*}
is NOT an equivalence $(\mathbb A^1,et)$ local by proposition \ref{rho1}
since $\rho_{S*}i_{\sharp}\mathbb Z(X/Z)=0$ but 
$D(\mathbb A^1,et)(\rho_{S*}i_*\mathbb Z(X/Z)=i_*\rho_{S*}\mathbb Z(X/Z))\neq 0\in\underline{\DA}(S)$.
\item[(iii)]Let $f:X\to S$ a smooth proper morphism with $X,S\in\Var(\mathbb C)$ of relative dimension $d=d_X-d_S$
and $X$ smooth.
Then, we have then by proposition \ref{sharpstarprop}(i) the equivalence $(\mathbb A^1,et)$ local in $C(\Var(\mathbb C)^{sm}/S)$
\begin{equation*}
T_0(f_{\sharp},f_*)(\mathbb Z(X/X)):=[\Delta_X]:f_{\sharp}\mathbb Z(X/X)=\mathbb Z(X/S)\to f_*E_{et}(\mathbb Z(X/X))(d)[2d]
\end{equation*}
given by the class of the diagonal $[\Delta_X]\in\Hom(f_{\sharp}\mathbb Z(X/X),f_*E_{et}(\mathbb Z(X/X))(d)[2d])$.
\item[(iii)']Let $f:X\to S$ a proper surjective morphism with $X,S\in\Var(\mathbb C)$ 
with equidimensional fiber of relative dimension $d=d_X-d_S$. Assume $X$ smooth.
Then, we have then by proposition \ref{sharpstarprop}(i) the equivalence $(\mathbb A^1,et)$ local in $C(\Var(\mathbb C)/S)$
\begin{equation*}
T_0(f_{\sharp},f_*)(\mathbb Z(X/X)):=[\Delta_X]:f_{\sharp}\mathbb Z(X/X)=\mathbb Z(X/S)\to f_*E_{et}(\mathbb Z(X/X))(d)[2d]
\end{equation*}
given by the class of the diagonal $[\Delta_X]\in\Hom(f_{\sharp}\mathbb Z(X/X),f_*E_{et}(\mathbb Z(X/X))(d)[2d])$.
\end{itemize}
\end{rem}

\subsection{Constructible motives and resolution of a motive by Corti-Hanamura motives}

We now give the definition of the motives of morphisms $f:X\to S$ which are constructible motives and
the definition of the category of Corti-Hanamura motives.

\begin{defi}
Let $S\in\Var(\mathbb C)$, 
\begin{itemize}
\item the homological motive functor is 
$M(/S):\Var(\mathbb C)/S\to\DA(S) \;, \; (f:X\to S)\mapsto M(X/S):=f_!f^!M(S/S)$,
\item the cohomological motive functor is 
$M^{\vee}(/S):\Var(\mathbb C)/S\to\DA(S) \;, \; (f:X\to S)\mapsto M(X/S)^{\vee}:=Rf_*M(X/X)=f_*E(\mathbb Z_X)$,
\item the Borel-Moore motive functor is 
$M^{BM}(/S):\Var(\mathbb C)/S\to\DA(S)\;, \; (f:X\to S)\mapsto M^{BM}(X/S):=f_!M(X/X)$,
\item the (homological) motive with compact support functor is
$M_c(/S):\Var(\mathbb C)/S\to\DA(S)\; , \; (f:X\to S)\mapsto M_c(X/S):=Rf_*f^!M(S/S)$.
\end{itemize}
Let $f:X\to S$ a morphism with $X,S\in\Var(\mathbb C)$. Assume that there exist a factorization
$f:X\xrightarrow{i} Y\times S\xrightarrow{p} S$, with $Y\in\SmVar(\mathbb C)$, 
$i:X\hookrightarrow Y$ is a closed embedding and $p$ the projection. Then, 
\begin{equation*}
Q(X/S):=p_{\sharp}\Gamma^{\vee}_X\mathbb Z_{Y\times S}\in C(\Var(\mathbb C)^{sm}/S)
\end{equation*}
(see definition \ref{projBMmotdef})is projective, admits transfert, and satisfy $D(\mathbb A^1_S,et)(Q(X/S))=M(X/S)$.
\end{defi}

\begin{defi}
\begin{itemize}
\item[(i)] Let $S\in\Var(\mathbb C)$. 
We define the full subcategory $\DA_c(S)\subset\DA(S)$ whose objects are constructible motives to
be the thick triangulated category generated by the motives of the form $M(X/S)$, 
with $f:X\to S$ a morphism, $X\in\Var(\mathbb C)$. 
\item[(ii)]Let $X,S\in\Var(\mathbb C)$. If $f:X\to S$ is proper (but not necessary smooth) and $X$ is smooth, 
$M(X/S)$ is said to be a Corti-Hanamura motive and we have by above in this case 
$M(X/S)=M^{BM}(X/S)[c]=M(X/S)^{\vee}[c]$, with $c=\codim(X,X\times S)$ where $f:X\hookrightarrow X\times S\to S$. 
We denote by 
\begin{equation*}
\mathcal{CH}(S)=\left\{M(X/S)\right\}_{\left\{X/S=(X,f),f \mbox{pr},X \mbox{sm}\right\}}^{pa}\subset DM(S) 
\end{equation*}
the full subcategory which is the pseudo-abelian completion of the full subcategory 
whose objects are Corti-Hanamura motives.
\item[(iii)] We denote by
\begin{equation*}
\mathcal {CH}^0(S)\subset\mathcal{CH}(S) 
\end{equation*}
the full subcategory which is the pseudo-abelian completion of the full subcategory 
whose objects are Corti-Hanamura motives $M(X/S)$ such that the morphism $f:X\to S$ is projective.
\end{itemize}
\end{defi}

For bounded above motives, we have

\begin{thm}\label{weightst}
Let $S\in\Var(\mathbb C)$.
\begin{itemize}
\item[(i)] There exists a unique weight structure $\omega$ on $\DA^-(S)$ such that $\DA^-(S)^{\omega=0}=\mathcal{CH}(S)$
\item[(ii)] There exist a well defined functor 
\begin{equation*}
W(S):\DA^-(S)\to K^-(\mathcal{CH}(S)) \; , \; W(S)(M)=[M^{(\bullet)}]
\end{equation*}
where $M^{(\bullet)}\in C^-(\mathcal{CH}(S))$ is a bounded above weight complex, such that if $m\in\mathbb Z$ is the highest weight,
we have a generalized distinguish triangle for all $i\leq m$
\begin{equation}\label{weightDT}
T_i:M^{(i)}[i]\to M^{(i+1)}[(i+1)]\to\cdots\to M^{(m)}[m]\to M^{w\geq i}
\end{equation}
Moreover the maps $w(M)^{(\geq i)}:M^{\geq i}\to M$ induce an isomorphism
 $w(M):holim_{i} M^{\geq i}\xrightarrow{\sim} M$
in $\DA^-(S)$.
\item[(iii)]Denote by $Chow(S)$ the category of Chow motives, which is the pseudo-abelian completion of the category 
\begin{itemize}
\item whose set of objects consist of the $X/S=(X,f)\in\Var(\mathbb C)/S$ such that $f$ is proper and $X$ is smooth, 
\item whose set of morphisms between $X_1/S$ and $X_2/S$ is $\CH^{d_1}(X_1\times_S X_2)$, 
and the composition law is given in \cite{C.H}.
\end{itemize} 
We have then a canonical functor $CH_S:Chow(S)\hookrightarrow\DA(S)$, with $CH_S(X/S):=M(X/S):=Rf_*\mathbb Z(X/X)$,
which is a full embedding whose image is the category $\mathcal{CH}(S)$.
\end{itemize}
\end{thm}

\begin{proof}
\noindent(i): The category $\DA(S)$ is clearly weakly generated by $\mathcal{CH}(S)$.
 Moreover $\mathcal{CH}(S)\subset\DA(S)$ is negative. Hence, the result follows from
\cite{Bondarko} theorem 4.3.2 III.

\noindent(ii): Follows from (i) by standard fact of weight structure on triangulated categories. 
See \cite{Bondarko} theorem 3.2.2 and theorem 4.3.2 V for example.

\noindent(iii): See \cite{Fangzhou}.
\end{proof}

\begin{thm}\label{weightst2}
Let $S\in\Var(\mathbb C)$.
\begin{itemize}
\item[(i)] There exists a unique weight structure $\omega$ on $\DA^-(S)$ such that $\DA^-(S)^{\omega=0}=\mathcal{CH}^0(S)$
\item[(ii)] There exist a well defined functor 
\begin{equation*}
W(S):\DA^-(S)\to K^-(\mathcal{CH}^0(S)) \; , \; W(S)(M)=[M^{(\bullet)}]
\end{equation*}
where $M^{(\bullet)}\in C^-(\mathcal{CH}^0(S))$ is a bounded above weight complex, such that if $m\in\mathbb Z$ is the highest weight,
we have a generalized distinguish triangle for all $i\leq m$
\begin{equation}\label{weightDT2}
T_i:M^{(i)}[i]\to M^{(i+1)}[(i+1)]\to\cdots\to M^{(m)}[m]\to M^{w\geq i}
\end{equation}
Moreover the maps $w(M)^{(\geq i)}:M^{\geq i}\to M$ induce an isomorphism 
$w(M):holim_{i} M^{\geq i}\xrightarrow{\sim} M$ in $\DA^-(S)$.
\end{itemize}
\end{thm}

\begin{proof}
Similar to the proof of theorem \ref{weightst}.
\end{proof}

\begin{cor}\label{weightst2Cor}
Let $S\in\Var(\mathbb C)$. Let $M\in\DA(S)$. 
Then there exist $(F,W)\in C_{fil}(\Var(\mathbb C)^{sm}/S)$ such that $D(\mathbb A^1,et)(F)=M$
and $D(\mathbb A^1,et)(\Gr^W_pF)\in\mathcal{CH}^0(S)$.
\end{cor}

\begin{proof}
By theorem \ref{weightst2}, there exist, by induction, for $i\in\mathbb Z$, a distinguish triangle in $\DA(S)$
\begin{equation}
T_i:M^{(i)}[i]\xrightarrow{m_i} M^{(i+1)}\xrightarrow{m_{i+1}}\cdots\xrightarrow{m_{m-1}}M^{(m)}[m]\to M^{w\geq i}
\end{equation}
with $M^{(j)}[j]\in\mathcal{CH}^0(S)$ and $w(M):holim_{i} M^{\geq i}\xrightarrow{\sim} M$ in $\DA^-(S)$.
For $i\in\mathbb Z$, take $(F_j)_{j\geq i},F_{w\geq i}\in C(\Var(\mathbb C)^{sm}/S)$ 
such that $D(\mathbb A^1,et)(F_j)=M^{(j)}[j]$, $D(\mathbb A^1,et)(F_{w\geq i})=M^{w\geq i}$ 
and such that we have in $C(\Var(\mathbb C)^{sm}/S)$, 
\begin{equation}
F_{w\geq i}=\Cone(F_i\xrightarrow{m_i} F_{i+1}\xrightarrow{m_{i+1}}\cdots\xrightarrow{m_{m-1}}F_m)
\end{equation}
where $m_j:F_j\to F_{j+1}$ are morphisms in $C(\Var(\mathbb C)^{sm}/S)$ such that $D(\mathbb A^1,et)(m_j)=m_j$.
Now set $F=\holim_iF_{w\geq i}\in C(\Var(\mathbb C)^{sm}/S)$ and $W_iF:=F_{w\geq i}\hookrightarrow F$,
so that $(F,W)\in C_{fil}(\Var(\mathbb C)^{sm}/S)$ satisfy $D(\mathbb A^1,et)(\Gr^W_pF)=M^{(p)}[p]\in\mathcal{CH}^0(S)$.
\end{proof}

\subsection{The restriction of relative motives to their Zariski sites}

Let $S\in\Var(\mathbb C)$. The adjonction
\begin{equation*}
(e(S)^*,e(S)_*):C(\Var(\mathbb C)^{sm}/S)\leftrightarrows C(S)
\end{equation*}
is a Quillen adjonction and induces in the derived category 
\begin{itemize}
\item $(e(S)^*,e(S)_*):\Ho_{zar}(\Var(\mathbb C)^{sm}/S)\leftrightarrows D(T):=\Ho_{zar}C(S)$,
since $e(S)_*$ sends Zariski local equivalence on the big site $\Var(\mathbb C)^{sm}/S$ to Zariski local equivalence
in the small Zariski site of $S$, 
\item $(e(S)^*,Re(S)_*):\DA(S)\leftrightarrows D(T):=\Ho_{zar}C(S)$. 
\end{itemize}

We will use in the defintion of the De Rahm realization functor on $\DA(S)$ the following proposition
concerning the restriction of the derived internal hom functor to the Zariski site :

\begin{prop}\label{eShom}
Let $M,N\in\DA(S)$ and $m:M\to N$ be a morphism. 
Let $F^{\bullet},G^{\bullet}\in\PSh(\Var(\mathcal C)^{sm}/S,C(\mathbb Z))$
such that $M=D(\mathbb A^1_S,et)(F^{\bullet})$ and $N=D(\mathbb A^1_S,et)(G^{\bullet})$.
If we take $G^{\bullet}$ $(\mathbb A^1_S,et)$ fibrant and admitting transfert, and $F^{\bullet}$ cofibrant 
for the projective model structure, we have
\begin{equation*}
Re(S)_*R\mathcal Hom^{\bullet}(M,N)=e(S)_*\mathcal Hom^{\bullet}(F^{\bullet},G^{\bullet})
\end{equation*} 
in D(S).
\end{prop}

\begin{proof}
Since $F^{\bullet}$ is projectively cofibrant and $G^{\bullet}$ is (projectively) $(\mathbb A^1_S,et)$ fibrant,
we have $R\mathcal Hom^{\bullet}(M,N)=\mathcal Hom^{\bullet}(F^{\bullet},G^{\bullet})$.
Then, $\mathcal Hom^{\bullet}(F^{\bullet},G^{\bullet})$ is $\mathbb A^1_S$ local and
admits transfert. On the other hand, we have   
\begin{equation*}
\mathcal L_{\mathbb A^1_S}D_{et}(\Cor\Var(\mathbb C)^{sm}/S)=\mathcal L_{\mathbb A^1_S}D_{zar}(\Cor\Var(\mathbb C)^{sm}/S)
\subset D(\Var(\mathbb C)^{sm}/S)
\end{equation*}
by theorem \ref{DDADM} (ii). This gives the equality of the proposition.
\end{proof}

We will also have :

\begin{prop}\label{keylem2Dmot}
For $f:T\to S$ a morphism and $i:Z\hookrightarrow S$ a closed embedding, 
with $Z,S,T\in\Var(\mathbb C)$, we have
\begin{itemize}
\item[(i)] $Re(S)_*Rf_*=Rf_*Re(T)_*$ and $e(S)^*Rf_*=Rf_*e(T)^*$ 
\item[(ii)] $Re(S)_*R\Gamma_Z=R\Gamma_ZRe(S)_*$.
\end{itemize}
\end{prop}

\begin{proof}
\noindent(i):Follows from proposition \ref{keylem2} (i) and the fact that $f_*$ preserve $(\mathbb A^1,et)$ fibrant complex of presheaves.

\noindent(ii):Follows from proposition \ref{keylem2} (ii) and the fact that $\Gamma_Z$ preserve $(\mathbb A^1,et)$ fibrant complex of presheaves.
\end{proof}

\subsection{Motives of complex analytic spaces}

The category of motives is obtained by inverting the $(\mathbb D^1_S,usu)$ local equivalence.
Hence the $\mathbb D^1_S$ local complexes of presheaves plays a key role.
\begin{defi}
The derived category of motives of complex algebraic varieties over $S$ is
the category
\begin{equation*}
\AnDA(S):=\Ho_{\mathbb D_S^1,usu}(C(\AnSp(\mathbb C)^{sm}/S)),
\end{equation*}
which is the localization of the category of complexes of presheaves on $\AnSp(\mathbb C)^{sm}/S$
with respect to $(\mathbb D_S^1,usu)$ local equivalence and we denote by
\begin{equation*}
D(\mathbb D_S^1,usu):=D(\mathbb A_S^1)\circ D(et):C(\AnSp(\mathbb C)^{sm}/S)\to\AnDA(S)
\end{equation*}
the localization functor. We have 
$\DA^-(S):=D(\mathbb D_S^1,usu)(\PSh(\AnSp(\mathbb C)^{sm}/S,C^-(\mathbb Z)))\subset\DA(S)$
the full subcategory consisting of bounded above complexes.
\end{defi}

\begin{thm}\label{AnDMthm}
Let $S\in\AnSp(\mathbb C)$.
The adjonction $(e(S)^*,e(S)_*):C(\AnSp(\mathbb C)^{sm}/S)\to C(S)$ induces an equivalence of categories
\begin{equation*}
(e(S)^*,Re(S)_*):\AnDM(S)\xrightarrow{\sim} D(S).
\end{equation*}
In particular, for $F\in C(\AnSp(\mathbb C)^{sm}/S)$, the adjonction map
$\ad(e(S)^*,e(S)_*)(F):e(S)^*e(S)_*F\to F$ is an equivalence $(\mathbb D^1,usu)$ local.
\end{thm}

\begin{proof}
See \cite{AyoubB}.
\end{proof}

We deduce from this theorem the following :

\begin{prop}\label{eShomAn}
Let $S\in\AnSp(\mathbb C)$. Let $F,G\in C(\AnSp(\mathbb C)^{sm}/S)$.
If $G$ is $\mathbb D^1$ local, then the canonical map
\begin{equation*}
T(e,hom)(F,G):e(S)_*\mathcal Hom(F,G)\to\mathcal Hom(e(S)_*F,e(S)_*G)
\end{equation*}
is an equivalence usu local.
\end{prop}

\begin{proof}
The map $T(e,hom)(F,G)$ is the composite
\begin{eqnarray*}
T(e,hom)(F,G):
e(S)_*\mathcal Hom(F,G)\xrightarrow{\mathcal Hom(\ad(e(S)^*,e(S)_*)(F),G)}e(S)_*\mathcal Hom(e(S)^*e(S)_*F,G) \\
\xrightarrow{I(e(S)^*,e(S)_*)(F,G)}\mathcal Hom(e(S)_*F,e(S)_*G)
\end{eqnarray*}
where the last map is the adjonction isomorphism. 
The first map is an isomorphism by theorem \ref{AnDMthm} since $G$ is $\mathbb D^1$ local.
\end{proof}

\section{The category of filtered D modules on commutative ringed topos, on commutative ringed spaces, 
complex algebraic varieties complex analytic spaces and the functorialities}

\subsection{The The category of filtered D modules on commutative ringed topos, on commutative ringed spaces, and the functorialities}

\subsubsection{Definitions et functorialities}

Let $(\mathcal S,O_S)\in\RCat$ with $O_S$ commutative.
Recall that $\Omega_{O_S}:=\mathbb D_S^O(\mathcal I_S/\mathcal I_S^2)\in\PSh_{O_S}(S)$ 
is the universal derivation $O_S$-module together with its derivation map $d:O_S\to\Omega_{O_S}$,
where $\mathcal I_S=\ker(s_S:O_S\otimes O_S\to O_S)\in\PSh_{O_S\times O_S}(\mathcal S)$ the diagonal ideal. 

In the particular case of a ringed space $(S,O_S)\in\RTop$, $s_S:O_S\otimes O_S=\Delta_S^*(p_1^*O_S\otimes p_2^*O_S)\to O_S$
is the structural morphism of diagonal embedding $\Delta_S:(S,O_S)\hookrightarrow(S\times S,p_1^*O_S\otimes p_2^*O_S)$,
$p_1:S\times S\to S$ and $p_2:S\times S\to S$ being the projections. 
More generally, for $k\in\mathbb N$, $k\geq 1$ we have the sheaf of $k$-jets $J^k(O_S):=\Delta_S^*\mathcal I_S/\mathcal I_S^{k+1}$
with in particular $J^1(O_S)=T_S$. 
We have, for $s\in S$, $J^k(O_S)_s=m_s/m_s^k$ where $m_s\subset O_{S,s}$ is the maximal ideal if $O_{S,s}$ is a local ring.

\begin{defi}
\begin{itemize}
\item[(i)]Let $(\mathcal S,O_S)\in\RCat$ with $O_S$ a commutative sheaf of ring and 
$\mathcal S$ is endowed with a topology $\tau$. We denote by 
\begin{equation*}
D(O_S)=<O_S,\Der_{O_S}(O_S,O_S)>\subset a_{\tau}\mathcal Hom(O_S,O_S)
\end{equation*}
the subsheaf of ring generated by $O_S$ and the subsheaf of derivations $\Der_{O_S}(O_S,O_S)=T_S:=\mathbb D^O_S\Omega_{O_S}$,
$a_{\tau}:\PSh(\mathcal S)\to\Shv(\mathcal S)$ being the sheaftification functor.
\item[(ii)] Let $f:\mathcal X\to \mathcal S$ be a morphism of site, 
with $\mathcal X,\mathcal S\in\Cat$ endowed with topology $\tau$, resp. $\tau'$, 
and $O_S\in\PSh(\mathcal S,\cRing)$ a commutative sheaf of ring.
We will note in this case by abuse $f^*O_S:=a_{\tau}f^*O_S$ and $f^*D(O_S):=a_{\tau}f^*D(O_S)$,
$a_{\tau}:\PSh(\mathcal X)\to\Shv(\mathcal X)$ being the sheaftification functor.
\end{itemize}
\end{defi}

Let $f:\mathcal X\to \mathcal S$ a morphism of site, 
with $\mathcal X,\mathcal S\in\Cat$ endowed with topology $\tau$, resp. $\tau'$, 
and $O_S\in\PSh(\mathcal S,\cRing)$ a commutative sheaf of ring. 
Consider the ringed space $(\mathcal X,f^*O_S):=(\mathcal X,a_{\tau}f^*O_S)\in\RCat$,
$a_{\tau}:\PSh(\mathcal X)\to\Shv(\mathcal X)$ being the sheaftification functor. Then, the map in $\PSh(\mathcal X)$
\begin{equation*}
T(f,hom)(O_S,O_S):f^*\mathcal Hom(O_S,O_S)\to\mathcal Hom(f^*O_S,f^*O_S)
\end{equation*} 
induces a canonical morphism of sheaf of rings 
\begin{equation*}
T(f,hom)(O_S,O_S):a_{\tau}f^*D(O_S)=:f^*D(O_S)\to D(a_{\tau}f^*O_S)=:D(f^*O_S). 
\end{equation*}

In the special case of ringed spaces, we have then :

\begin{prop}\label{DmodRTopiso}
Let $f:X\to S$ is a continous map, with $X,S\in\Top$ and $O_S\in\PSh(S,\cRing)$ a commutative sheaf of ring.
Consider the ringed space $(X,f^*O_S):=(X,a_{\tau}f^*O_S)\in\RTop$,
 $a_{\tau}:\PSh(X)\to\Shv(X)$ being the sheaftification functor.
Then, the map in $\PSh(X)$
\begin{equation*}
T(f,hom)(O_S,O_S):f^*\mathcal Hom(O_S,O_S)\to\mathcal Hom(f^*O_S,f^*O_S)
\end{equation*} 
induces a canonical isomorphism of sheaf of rings 
\begin{equation*}
T(f,hom)(O_S,O_S):f^*D(O_S):=a_{\tau}f^*D(O_S)\xrightarrow{\sim} D(a_{\tau}f^*O_S)=:D(f^*O_S). 
\end{equation*}
\end{prop}

\begin{proof}
For all $x\in X$,
\begin{equation*}
T(f,hom)(O_S,O_S)_x:(f^*D(O_S))_x\xrightarrow{\sim}D(O_{S,f(x)})\xrightarrow{\sim}(D(f^*O_S))_x. 
\end{equation*}
Hence, since $a_{\tau}f^*D(O_S)$ and $D(a_{\tau}f^*O_S)$ are sheaves,
\begin{equation*}
T(f,hom)(O_S,O_S):f^*D(O_S):=a_{\tau}f^*D(O_S)\xrightarrow{\sim} D(a_{\tau}f^*O_S)=:D(f^*O_S). 
\end{equation*}
is an isomorphism
\end{proof}

We will consider presheaves of $D(O_S)$ modules on a ringed topos $(\mathcal S,O_S)$ :

\begin{defi}\label{Dmodefinition}
Let $(\mathcal S,O_S)\in\RCat$ with $O_S$ a commutative sheaf of ring. 
\begin{itemize}
\item[(i)]We will consider $\PSh_{D(O_S)}(\mathcal S)$ the category of presheaves of (left) $D(O_S)$
modules on $S$ and $C_{D(O_S)}(\mathcal S):=C(\PSh_{D(O_S)}(\mathcal S))$ its category of complexes.
We will consider $\PSh_{D(O_S)^{op}}(\mathcal S)$ the category of presheaves of right $D(O_S)$
modules on $S$ and $C_{D(O_S)^{op}}(\mathcal S):=C(\PSh_{D(O_S)^{op}}(\mathcal S))$ its category of complexes.
We denote again by abuse 
\begin{equation*}
\PSh_{D(O_S)fil}(\mathcal S)=(\PSh_{D(O_S)}(\mathcal S),F):=(\PSh_{(D(O_S),F^{ord})}(\mathcal S),F) 
\end{equation*}
the category of filtered $(D(O_S),F^{ord})$-module, with, for $-p\leq 0$,
$F^{ord,-p}D(O_S)=\left\{P\in D(O_S),ord(P)\leq p\right\}$ and $F^{ord,p}D(O_S)=0$ for $p>0$, 
\begin{itemize}
\item whose objects are $(M,F)\in(\PSh_{O_S}(\mathcal S),F)$ such that $(M,F)$ is compatible with $(D(O_S),F^{ord})$
that is $F^{ord,-p}D(O_S)\cdot F^qM\subset F^{q-p}M$ (Griffitz transversality), this is to say that
the structural map $md:M\otimes_{O_S}D(O_S)\to M$ induces a filtered map of presheaves (i.e a map in $(\PSh_{O_S}(\mathcal S),F)$)
$md:(M,F)\otimes_{O_S}(D(O_S),F^{ord})\to(M,F)$,
\item whose morphism $\phi:(M_1,F)\to(M_2,F)$ are as usual the morphisms of presheaves $\phi:M_1\to M_2$
which are morphism of filtered presheaves (i.e. $\phi(F^pM_1)\subset F^pM_2$) and which are $D(O_S)$ linear (in particular $O_S$ linear).
\end{itemize}
Note that this a NOT the category of filtered $D(O_S)$ modules in the usual sense, that is the $(M,F)\in(\PSh_{O_S}(\mathcal S),F)$ 
together with a map $md:(M,F)\otimes_{O_S}D(O_S)\to(M,F)$ in $(\PSh_{O_S}(\mathcal S),F)$, since $F^{ord}$ is NOT the trivial filtration. 
More precisely the $O_S$ submodules $F^pM\subset M$ are NOT $D(O_S)$ submodules but satisfy Griffitz transversality. 
We denote by 
\begin{equation*}
\PSh_{D(O_S)0fil}(\mathcal S)\subset\PSh_{D(O_S)fil}(\mathcal S), \; 
\PSh_{D(O_S)(1,0)fil}(\mathcal S)\subset\PSh_{D(O_S)2fil}(\mathcal S):=(\PSh_{(D(O_S),F^{ord})}(\mathcal S),F,W)
\end{equation*}
the full subcategory consisting of filtered $D(O_S)$ module in the usual sense, resp.
the full subcategory such that $W$ is a filtration by $D(O_S)$ submodules. 
\item[(ii)] We denote again by 
\begin{equation*}
C_{D(O_S)fil}(\mathcal S)\subset C(\PSh_{D(O_S)}(\mathcal S),F) \; , \; 
C_{D(O_S)2fil}(\mathcal S)\subset C(\PSh_{D(O_S)}(\mathcal S),F,W)
\end{equation*}
the full subcategory of complexes such that the filtration(s) is (are) regular.
We will consider also
\begin{equation*}
C_{D(O_S)0fil}(\mathcal S)\subset C_{D(O_S)fil}(\mathcal S), C_{D(O_S)(1,0)fil}(\mathcal S)\subset C_{D(O_S)2fil}(\mathcal S)
\end{equation*}
the full subcategory consisting of complexes of filtered $D(O_S)$ modules in the usual sense (i.e. by $D(O_S)$ submodule),
respectively the full subcategory consisting of complexes of bifilterd $D(O_S)$ modules such that $W^pM\subset M$ are $D(O_S)$ submodules
i.e. the filtration $W$ is a filtration in the usual sense, but NOT $F$ wich satify only Griffitz transversality.
\end{itemize}
\end{defi}

\begin{prop}\label{RSpaceDmod}
Let $(\mathcal S,O_S)\in\RCat$ with a $O_S$ commutative sheaf of ring. 
\begin{itemize}
\item[(i)] Let $M\in\PSh_{O_S}(\mathcal S)$. 
Then, there is a one to one correspondence between 
\begin{itemize}
\item the $D(O_S)$ module structure on $M$ compatible with the $O_S$ module structure,
that is the maps $md:D(O_S)\otimes_{O_S} M\to M$ in $\PSh_{O_S}(\mathcal S)$ and  
\item the integrable connexions on $M$, that is the maps 
$\nabla:M\to\Omega_{O_S}\otimes_{O_S}M$ satisfying $\nabla\circ\nabla=0$
with $\nabla:\Omega_{O_S}\otimes_{O_S}M\to\Omega^2_{O_S}\otimes_{O_S}M$ 
given by $\nabla(\omega\otimes m)=(d\omega)\otimes m+\omega\wedge\nabla(m)$
\end{itemize}
\item[(ii)] Let $(M,F)\in\PSh_{O_Sfil}(\mathcal S)$. 
Then, there is a one to one correspondence between 
\begin{itemize}
\item the $D(O_S)$ module structure on $(M,F)$ compatible with the $O_S$ module structure,
that is the maps $md:(D(O_S),F^{ord})\otimes_{O_S} (M,F)\to (M,F)$ in $\PSh_{O_Sfil}(\mathcal S)$ and  
\item the integrable connexions on $M$, 
that is the maps $\nabla:(M,F)\to\Omega_{O_S}\otimes_{O_S}(M,F)$ satisfying $\nabla\circ\nabla=0$
and Griffith transversality (i.e. $\nabla(F^pM)\subset\Omega_{O_S}\otimes_{O_S}F^{p-1}M$).
\end{itemize}
\end{itemize}
\end{prop}

\begin{proof}
Standard.
\end{proof}

The following proposition tells that the $O$-tensor product of $D$ modules has a canonical structure of $D$ module.

\begin{defiprop}\label{RCatDOXOS}
\begin{itemize}
\item[(i)] Let $f:(\mathcal X,O_X)\to (\mathcal S,O_S)$ a morphism with $(\mathcal X,O_X),(\mathcal S,O_S)\in\RCat$
with commutative structural sheaf of ring. 
For $N\in\PSh_{O_X,D(f^*O_S)}(\mathcal X)$ and $M\in\PSh_{O_X,D(f^*O_S)}(\mathcal X)$, $N\otimes_{O_X}M$ 
has the canonical $D(f^*O_S)$ module structure given by, for $X^o\in\mathcal X$, 
\begin{equation*}
\gamma\in\Gamma(X^o,D(f^*O_S)), m\in\Gamma(X^o,M),n\in\Gamma(X^o,N), \; 
\gamma.(n\otimes m)=(\gamma.n)\otimes m+n\otimes(\gamma.m).
\end{equation*}
This gives the functor
\begin{eqnarray*}
\PSh_{O_X,D(f^*O_S)fil}(\mathcal X)\times\PSh_{O_X,D(f^*O_S)fil}(\mathcal X)\to\PSh_{O_X,D(f^*O_S)fil}(\mathcal X),  
((M,F),(N,F))\mapsto \\ (M,F)\otimes_{O_X}(N,F), 
F^p(M,F)\otimes_{O_X}(N,F):=\Im(\oplus_{q\in\mathbb Z} F^qM\otimes_{O_X}F^{p-q}N\to M\otimes_{O_X}N)
\end{eqnarray*}

\item[(ii)] Let $f:(\mathcal X,O_X)\to (\mathcal S,O_S)$ a morphism with $(\mathcal X,O_X),(\mathcal S,O_S)\in\RCat$
with commutative structural sheaf of ring. 
For $N\in C_{D(O_X),D(f^*O_S)}(\mathcal X)$ and $M\in C_{D(O_X)^{op}}(\mathcal X)$, $N\otimes_{D(O_X)}M$ 
has the canonical $f^*D(O_S)$ module structure given by, for $X^o\in\mathcal X$, 
\begin{equation*}
\gamma\in\Gamma(X^o,D(f^*O_S)), m\in\Gamma(X^o,M),n\in\Gamma(X^o,N), \; \gamma.(n\otimes m)=(\gamma.n)\otimes m.
\end{equation*}
This gives the functor
\begin{eqnarray*}
C_{D(O_X),D(f^*O_S)fil}(\mathcal X)\times C_{D(O_X)fil}(X)\to C_{D(f^*O_S)fil}(X),
((M,F),(N,F))\mapsto \\ (M,F)\otimes_{D(O_X)}(N,F),
F^p(M,F)\otimes_{D(O_X)}(N,F):=\Im(\oplus_{q\in\mathbb Z} F^qM\otimes_{D(O_X)}F^{p-q}N\to M\otimes_{D(O_X)}N)
\end{eqnarray*}
Note that, by definition, we have for $(M,F)\in(\PSh_{D(O_S)fil}(\mathcal S))$, the canonical isomorphism 
\begin{equation*}
(M,F)\otimes_{D(O_S)}(D(O_S),F^{ord})\xrightarrow{\sim}(M,F), \; m\otimes P\mapsto Pm, \;  m\mapsto (m\otimes 1)
\end{equation*}
\end{itemize}
\end{defiprop}

\begin{proof}
Immediate from definition.
\end{proof}

We now look at the functorialites for morphisms of ringed spaces, using proposition \ref{DmodRTopiso}. 
First note that for $f:(X,O_X)\to(S,O_S)$ a morphism,
with $(X,O_X),(S,O_S)\in\RTop$ with structural presheaves commutative sheaves of rings,
there is NO canonical morphism between $D(f^*O_S)=f^*D(O_S)$ (see proposition \ref{DmodRTopiso}) and $D(O_X)$. 

We have the pullback functor for (left) D-modules :

\begin{defiprop}\label{RSpacef}
\begin{itemize}
\item[(i)]Let $f:(X,O_X)\to(S,O_S)$ a morphism with $(X,O_X),(S,O_S)\in\RTop$ 
with structural presheaves commutative sheaves of rings. Recall that $f^*D(O_S)=D(f^*O_S)$ in this case. 
Then for $(M,F)\in\PSh_{D(O_S)fil}(S)$,
\begin{equation*}
f^{*mod}(M,F):=(O_X,F_b)\otimes_{f^*O_S}f^*(M,F)\in\PSh_{O_Xfil}(X) 
\end{equation*}
has a canonical structure of filtered $D(O_X)$ module given by 
\begin{eqnarray*}
\mbox{for} \; \gamma\in\Gamma(X^o,T_{O_X}), \, n\otimes m\in\Gamma(X^o,O_X\otimes_{f^*O_S}f^*M), \; 
\gamma.(n\otimes m):=(\gamma.n)\otimes m+n\otimes df(\gamma)(m)
\end{eqnarray*}
with $df:=\mathbb D^O_S\Omega_{O_X/f^*O_S}:T_{O_X}\to T_{f^*O_S}=f^*T_{O_S}$ and 
$f^*(M,F)\in\PSh_{f^*D(O_S)fil}(X)=\PSh_{D(f^*O_S)fil}(X)$.
\item[(ii)] More generally, let $f:(\mathcal X,O_X)\to(\mathcal S,O_S)$ a morphism 
with $(\mathcal X,O_X),(\mathcal S,O_S)\in\RCat$ with structural presheaves commutative sheaves of rings. 
Assume that the canonical morphism $T(f,\hom)(O_S,O_S):f^*D(O_S)\to D(f^*O_S)$ is an isomorphism of sheaves.
Then for $(M,F)\in\PSh_{D(O_S)fil}(\mathcal S)$, 
\begin{equation*}
f^{*mod}(M,F):=(O_X,F_b)\otimes_{f^*O_S}f^*(M,F)\in\PSh_{O_Xfil}(\mathcal X) 
\end{equation*}
has a canonical structure of filtered $D(O_X)$ module given by 
\begin{eqnarray*}
\mbox{for} \; \gamma\in\Gamma(X^o,T_{O_X}), \, n\otimes m\in\Gamma(X^o,O_X\otimes_{f^*O_S}f^*M), \; 
\gamma.(n\otimes m):=(\gamma.n)\otimes m+n\otimes df(\gamma)(m)
\end{eqnarray*}
with $df:=\mathbb D^O_S\Omega_{O_X/f^*O_S}:T_{O_X}\to T_{f^*O_S}=f^*T_{O_S}$ and 
$f^*(M,F)\in\PSh_{f^*D(O_S)fil}(\mathcal X)=\PSh_{D(f^*O_S)fil}(\mathcal X)$.
\end{itemize}
\end{defiprop}

\begin{proof}
Standard.
\end{proof}

\begin{rem}
\begin{itemize}
\item Let $f:(X,O_X)\to(S,O_S)$ a morphism with $(X,O_X),(S,O_S)\in\RTop$ with structural presheaves commutative sheaves of rings. 
Recall that $f^*D(O_S)=D(f^*O_S)$.Then by definition $f^{*mod}(O_S,F_b)=(O_X,F_b)$.
\item More generally, let $f:(\mathcal X,O_X)\to(\mathcal S,O_S)$ a morphism 
with $(\mathcal X,O_X),(\mathcal S,O_S)\in\RCat$ with structural presheaves commutative sheaves of rings. 
Assume that the canonical morphism $T(f,\hom)(O_S,O_S):f^*D(O_S)\to D(f^*O_S)$ is an isomorphism of sheaves.
Then by definition $f^{*mod}(O_S,F_b)=(O_X,F_b)$.
\end{itemize}
\end{rem}

For the definition of a push-forward functor for a right D module we use the transfert module

Let $f:(X,O_X)\to (S,O_S)$ be a morphism with $(X,O_X),(S,O_S)\in\RTop$ with structural presheaves commutative sheaves of rings.
Then, the transfer module is
\begin{equation*}
(D(O_X\to f^*O_S),F^{ord}):=f^{*mod}(D(O_S),F^{ord}):=f^*(D(O_S),F^{ord})\otimes_{f^*O_S}(O_X,F_b)
\end{equation*}
which is a left $D(O_X)$ module and a left and right $f^*D(O_S)=D(f^*O_S)$ module.

\begin{lem}\label{DcompRTop}
Let $f_1:(X,O_X)\to (Y,O_Y)$, $f_2:(Y,O_Y)\to (S,O_S)$ be two morphism with $(X,O_X),(Y,O_Y)(S,O_S)\in\RTop$.
We have in $C_{D(O_X),(f_2\circ f_1)^*D(O_S)fil}(X)$ 
\begin{equation*}
(D_{O_X\to (f_2\circ f_1)^*O_S},F^{ord})=f_1^*(D_{O_Y\to f_2^*O_S},F^{ord})\otimes_{f_1^*D(O_Y)}(D_{O_X\to f_1^*O_Y},F^{ord})
\end{equation*}
\end{lem}

\begin{proof}
Follows immediately from definition.
\end{proof}

For right D module, we have the direct image functor :

\begin{defi}\label{RSpacefd}
Let $f:(X,O_X)\to(S,O_S)$ a morphism with $(X,O_X),(S,O_S)\in\RTop$ with structural presheaves commutative sheaves of rings. 
Then for $(M,F)\in C_{D(O_X)^{op}fil}(X)$, we define
\begin{equation*}
f^{00}_{*mod}(M,F)=f_*((D_{O_X\to f^*O_S},F^{ord})\otimes_{D(O_X)}(M,F))\in C_{D(O_S)fil}(S)
\end{equation*}
\end{defi}

For a closed embedding of topological spaces, there is the $V$-filtration on the structural sheaf, 
it will play an important role in this article

\begin{defi}\label{Vfil}
\begin{itemize}
\item[(i)]Let $(S,O_S)\in\RTop$ a locally ringed space. Let $Z=V(\mathcal I_Z)\subset S$ a Zariski closed subset. 
We set, for $S^o\subset S$ an open subset, $p\in\mathbb Z$,
\begin{itemize}
\item  $V_{Z,p}O_S(S^o):=O_S(S^o)$ if $p>0$,
\item  $V_{Z,-q}O_S(S^o):=\mathcal I^q_Z(S^o)\subset O_S(S^o)$ for $p=-q\leq 0$.
\end{itemize}
We immediately check that, by definition, this filtration satisfy Griffitz transversality, that is
$(O_S,V_Z)\in\PSh_{D(O_S)fil}(S)$.
For a morphism $g:((T,O_T),Z')\to((S,O_S),Z)$ with $((T,O_T),Z),((S,O_S),Z)\in\RTop^2$ locally ringed spaces,
where $Z$ and $Z'$ are Zariski closed subsets, the structural morphism $a_g:g^*O_S\to O_T$ is a filtered morphism :
\begin{eqnarray*}
a_g:g^*(O_S,V_Z)\to (O_T,V_{Z'}), h\mapsto a_g(h)
\end{eqnarray*}
\item[(ii)]Let $(S,O_S)\in\RTop$. Let $i:Z\hookrightarrow S$ a closed embedding.
The $V_Z$-filtration on $O_S$ (see (i)) gives the filtration, given by for $p\in\mathbb Z$,
\begin{equation*}
V_{Z,p}\Hom(O_S,O_S):=\left\{P\in\Hom(O_S,O_S),\, \mbox{s.t.} P\mathcal I_Z^k\subset\mathcal I_Z^{k-p}\right\}
\end{equation*}
on $\Hom(O_S,O_S)$, which induces the filtration $V_{Z,p}D(O_S):=D(O_S)\cap V_{Z,p}\Hom(O_S,O_S)$ on $D(O_S)\subset\Hom(O_S,O_S)$. 
We get $(D(O_S),V_Z)\in\PSh_{fil}(S,\Ring)$ and we call it the $V_Z$-filtration on $D(O_S)$.
\item[(iii)]Let $(S,O_S)\in\RTop$ a locally ringed space. 
Let $i:Z=V(\mathcal I_Z)\hookrightarrow S$ a Zariski closed embedding and $O_Z:=i^*O_S/\mathcal I_Z$.
We say that $M\in\PSh_{D(O_S)}(S)$ is specializable on $Z$ if it admits an (increasing) filtration 
(called a $V_Z$-filtration) $(M,V)\in\PSh_{O_Sfil}(S)$ compatible with $(D_S,V_Z)$, 
that is $V_{Z,p}D_S\cdot V_{q}M\subset V_{p+q}M$, this is to say that the structural map $md:M\otimes_{O_S}D(O_S)\to M$ 
induces filtered map of presheaves $md:(M,V)\otimes_{O_S}(D(O_S),V_Z)\to(M,V)$.
For $(M,F)\in\PSh_{D(O_S)fil}(S)$ such that $M$ is specializable on $Z$, we thus get a filtered morphism
$md:(M,F,V)\otimes_{O_S}(D(O_S),F^{ord},V_Z)\to(M,F,V)$.
\item[(iii)'] Consider an injective morphism $m:M_1\hookrightarrow M_2$ with $M_1,M_2\in\PSh_{D(O_S)}(S)$.
If $M_2$ admits a $V_Z$ filtration $V_2$, then the filtration $V_{21}$ induced on $M_1$ 
(recall $V_{21,p}M_1:=V_{2,p}M_2\cap M_1$) is a $V_Z$ filtration.
Consider a surjective morphism $n:M_1\to M_2$ with $M_1,M_2\in\PSh_{D(O_S)}(S)$.
If $M_1$ admits a $V_Z$ filtration $V_1$, then the filtration $V_{12}$ induced on $M_2$ 
(recall $V_{12,p}M_2:=n(V_{1,p}M_1)$) is a $V_Z$ filtration.
\item[(iv)]Let $(S,O_S)\in\RTop$ a locally ringed space. 
Let $i:Z=V(\mathcal I_Z)\hookrightarrow S$ a Zariski closed embedding and $O_Z:=i^*O_S/\mathcal I_Z$.
For $(M,F)\in\PSh_{D(O_S)fil}(S)$ such that $M$ admits a $V_Z$ filtration $V$ so that
$(M,F,V)\in\PSh_{O_S2fil}(S)$, we will consider the quotient map in $\PSh_{O_Sfil}(S)$
\begin{equation*}
q_{V0}:(M,F)\to (M,F)/V_{-1}(M,F)=:Q_{V,0}(M,F). 
\end{equation*}
The quotient $i^*Q_{V,0}(M,F)$ has an action of $T_{O_Z}$ since for $S^o\subset S$ an open subset
and $\partial_z\in\Gamma(Z\cap S^o,T_{O_Z})\subset\Gamma(S^o,T_{O_S})$,
we have $\partial_z\in\Gamma(S,V_{Z,0}D(O_S))$ since for $f=\sum_{i=1}^r t_ih_i\in\Gamma(S^o,\mathcal I_Z)$,
where $(t_i)=\mathcal I_Z(S^o)$ are generators of the ideal $I_Z(S^o)\subset O_S(S^o)$ and $h_i\in\Gamma(S^o,O_S)$, we have
\begin{equation*}
\partial_z(\sum_{i=1}^r t_ih_i)=\sum_{i=1}^r(\partial_z(t_i)h_i+t_i\partial_z(h_i))=\sum_{i=1}^rt_i(\partial_z(h_i))
\in\Gamma(S,I_Z)
\end{equation*}
as $\partial_z(t_i)=0$ (only the vector field of $T_{O_S}$ which are transversal to $T_{O_Z}\subset T_{O_S}$ increase the grading),
Then, obviously, by definition, the map in $\PSh_{i^*O_Sfil}(Z)$
\begin{equation*}
i^*q_{V0}:i^*(M,F)/V_{-1}(M,F)=:i^*Q_{V,0}(M,F)
\end{equation*}
commutes with the action of $T_{O_Z}\subset i^*T_{O_S}$ and we call it the specialization map.
\end{itemize}
\end{defi}

\begin{defiprop}\label{Vfildefiprop}
Let $(S,O_S)\in\RTop$ a locally ringed space. Consider a commutative diagram
\begin{equation*}
\xymatrix{
Z_1=V(\mathcal I_1)\ar[r]^{i_1} & (S,O_S) \\
Z=V(\mathcal I)\ar[r]^{i'_1}\ar[u]^{i'_2} & Z_2=V(\mathcal I_2)\ar[u]^{i_2}}
\end{equation*}
where the maps are Zariski closed embeddings and which is cartesian 
(i.e. $\mathcal I=(\mathcal I_1,\mathcal I_2)$, in particular $Z=Z_1\cap Z_2$).
\begin{itemize}
\item[(i)] Let $(M,F)\in\PSh_{D(O_S)fil}(S)$ such that $M$ admits a $V_{Z_1}$-filtration $V_1$ and 
a $V_{Z_2}$-filtration $V_2$ (see definition \ref{Vfil}). Let $p,q\in\mathbb Z$. Then, we consider 
\begin{itemize}
\item the filtration $V_{21}$ on $V_{1,p}(M,F):=(V_{1,p}M,F)$ induced by $V_2$, 
\item the filtration $V_{12}$ on $Q_{V_2,p}(M,F):=(M/V_{2,p-1}M,F)$ induced by $V_1$. 
\end{itemize}
The quotient map in $\PSh_{i'_{1*}O_{Z_1}fil}(S)$
\begin{equation*}
q_{V_2,p}:V_{1,q}(M,F)\to V_{12,q}Q_{V_2,p}(M,F)
\end{equation*}
factors trough
\begin{equation*}
q_{V_2,p}:V_{1,q}(M,F)\xrightarrow{q_{V_{21},p}}
Q_{V_{21},p}V_{1,q}(M,F)\xrightarrow{Q^{p,q}_{V_1,V_2}(M,F)}V_{12,q}Q_{V_2,p}(M,F),
\end{equation*}
and the map $Q^{p,q}_{V_1,V_2}(M,F)$ in $\PSh_{i'_{1*}O_{Z_1}fil}(S)$ commute with the action of $T_{O_Z}$.
\item[(ii)]If $m:(M,F)\to(M',F)$ is a morphism with 
$(M,F),(M',F)\in\PSh_{D(O_S)fil}(S)$ admitting $V_{Z_1}$-filtration $V_1$ and $V'_1$ respectively
such that $m(V_{1,p}M)\subset V'_{1,p}M'$
and $V_{Z_2}$-filtration $V_2$ and $V'_2$ respectively such that $m(V_{1,p}M)\subset V'_{1,p}M'$. 
Then for all $p,q\in\mathbb Z$ the following diagram commutes
\begin{equation*}
\xymatrix{Q_{V_{21},p}V_{1,q}(M,F)\ar[rr]^{Q^{p,q}_{V_1,V_2}(M,F)}\ar[d]^{m} & \, &
V_{12,q}Q_{V_2,p}(M,F)\ar[d]^{m} \\
Q_{V'_{21},p}V'_{1,q}(M',F)\ar[rr]^{Q^{p,q}_{V'_1,V'_2}(M',F)} & \, & 
V_{12,q}Q_{V_2,p}(M',F)}
\end{equation*} 
\end{itemize}
Let $(S,O_S)\in\RTop$ a locally ringed space. Consider a commutative diagram
\begin{equation*}
\xymatrix{
Z'_1=V(\mathcal I'_1)\ar[r]^{l_1} & Z_1=V(\mathcal I_1)\ar[r]^{i_1} & (S,O_S) \\
Z'=V(\mathcal I')\ar[r]^{l'_1}\ar[u]^{i''_2} & Z=V(\mathcal I)\ar[r]^{i'_1}\ar[u]^{i'_2} & Z_2=V(\mathcal I_2)\ar[u]^{i_2}}
\end{equation*}
where the maps are Zariski closed embeddings and whose squares are cartesian
(i.e. $\mathcal I=(\mathcal I_1,\mathcal I_2)$ and $\mathcal I'=(\mathcal I'_1,\mathcal I)$,
in particular $Z=Z_1\cap Z_2$ and $Z'=Z'_1\cap Z$).
Let $(M,F)\in\PSh_{D(O_S)fil}(S)$ such that $M$ admits a $V_{Z_1}$-filtration $V_1$, a $V_{Z'_1}$-filtration $V'_1$, 
and a $V_{Z_2}$-filtration $V_2$(see definition \ref{Vfil}).
Then for all $p,q\in\mathbb Z$, denoting again $V'_{12}$ the filtration induced by $V'_1$ on $Q_{V_{21},p}V_{1,q}M$ 
and $V'_{21}$ the filtration induced by $V'_2$ on $V_{12,q}Q_{V_2,p}$
\begin{itemize}
\item $Q_{V'_{21},p}V'_{12,q}Q_{V_{21},p}V_{1,q}(M,F)=Q_{V'_{21},p}V'_{1,q}(M,F)$,
\item $V'_{12,q}Q_{V'_{21},p}V_{12,q}Q_{V_2,p}(M,F)=V'_{12,q}Q_{V_2,p}(M,F)$
\end{itemize}
and
\begin{equation*}
Q^{p,q}_{V'_{12},V_{21}}(Q^{p,q}_{V_1,V_2}(M,F))=Q^{p,q}_{V'_1,V_2}(M_2,F).
\end{equation*}
\end{defiprop}

\begin{proof}
Obvious.
\end{proof}

We will also consider the following categories
\begin{defi}\label{ODmod}
Let $(\mathcal X,O_X)\in\RCat$.
We denote by $C_{O_Xfil,D(O_X)}(\mathcal X)$ the category 
\begin{itemize}
\item whose objects $(M,F)\in C_{O_Xfil,D(O_X)}(\mathcal X)$ are filtered complexes of presheaves of $O_X$ modules
$(M,F)\in C_{O_Xfil}(\mathcal X)$ whose cohomology presheaves $H^n(M,F)\in\PSh_{O_Xfil}(\mathcal X)$ 
are emdowed with a structure of filtered $D(O_X)$ modules for all $n\in\mathbb Z$. 
\item whose set of morphisms 
$\Hom_{C_{O_Xfil,D(O_X)}(\mathcal X)}((M,F),(N,F))\subset\Hom_{C_{O_Xfil}(\mathcal X)}((M,F),(N,F))$ 
between $(M,F),(N,F)\in C_{O_Xfil,D(O_X)}(\mathcal X)$ are the morphisms of filtered complexes of $O_X$ modules 
$m:(M,F)\to(N,F)$ such that $H^nm:H^n(M,F)\to H^n(N,F)$ is $D(O_X)$ linear,
i.e. is a morphism of (filtered) $D(O_X)$ modules, for all $n\in\mathbb Z$.
\end{itemize}
\end{defi}

\subsubsection{The De Rham complex of a (left) filtered D-module and the Spencer complex of a right filtered D-module}

Using proposition \ref{RSpaceDmod}, we define the filtered De Rham complex of a complex of filtered (left) D-modules :

\begin{defi}\label{DOSDR}
\begin{itemize}
\item[(i)]Let $(\mathcal S,O_S)\in\RCat$ with $O_S$ commutative. 
Let $(M,F)\in C_{D(O_S)fil}(\mathcal S)$. By proposition \ref{RSpaceDmod}, we have the complex
\begin{equation*}
DR(O_S)(M,F):=(\Omega_{O_S}^{\bullet},F_b)\otimes_{O_S}(M,F)\in C_{fil}(\mathcal S)
\end{equation*}
whose differentials are $d(\omega\otimes m)=(d\omega)\otimes m+\omega\wedge(\nabla m)$.
\item[(ii)] More generally, let $f:(\mathcal X,O_X)\to (\mathcal S,O_S)$ with $(\mathcal X,O_X),(\mathcal S,O_S)\in\RCat$. 
The quotient map $q:\Omega_{O_X}\to\Omega_{O_X/f^*O_S}$ induce, for $G\in\PSh_{O_X}(\mathcal X)$ the quotient map 
\begin{equation*}
q^p(G):=\wedge^pq\otimes I\Omega^q_{O_X}\otimes_{O_X}G\to\Omega^q_{O_X/f^*O_S}\otimes_{O_X}G.
\end{equation*}
Let $(M,F)\in C_{D(O_X)fil}(\mathcal X)$. By proposition \ref{RSpaceDmod}, we have the relative De Rham complex
\begin{equation*}
DR(O_X/f^*O_S)(M,F):=(\Omega_{X/S}^{\bullet},F_b)\otimes_{O_X}(M,F)\in C_{f^*O_Sfil}(\mathcal X)
\end{equation*}
whose differentials are $d(q^p(M)(\omega\otimes m)):=q^{p+1}(M)((d\omega)\otimes m)+q^{p+1}(M)(\omega\otimes(\nabla m))$.
\item[(iii)] Let $(X,O_X)/\mathcal F\in FolRTop$, that is $(X,O_X)\in\RTop$ endowed with a foliation
with quotient map $q:\Omega_{O_X}\to\Omega_{O_X/\mathcal F}$.
Let $(M,F)\in C_{D(O_X)fil}(X)$. By proposition \ref{RSpaceDmod}, we have the foliated De Rham complex
\begin{equation*}
DR(O_X/\mathcal F)(M,F):=(\Omega_{O_X/\mathcal F}^{\bullet},F_b)\otimes_{O_X}(M,F)\in C_{fil}(X)
\end{equation*}
whose differentials are $d(q(M)(\omega\otimes m)):=q(M)((d\omega)\otimes m)+q(M)(\omega\otimes(\nabla m))$.
\end{itemize}
By definition,
\begin{itemize}
\item with the notation of (ii) if $\phi:(M_1,F)\to (M_2,F)$ is a morphism in $C_{D(O_X)fil}(\mathcal X)$, 
\begin{equation*}
DR(O_X/f^*O_S)(\phi):=(I\otimes\phi):(\Omega_{X/S}^{\bullet},F)\otimes_{O_X}(M_1,F)\to (\Omega_{X/S}^{\bullet},F)\otimes_{O_X}(M_2,F) 
\end{equation*}
is a morphism in $C_{f^*O_Sfil}(\mathcal X)$,
\item with the notation of (ii) $DR(O_X)(O_X)=DR(O_X)$ and more generally in the relative case $DR(O_X/f^*O_S)(O_X)=DR(O_X/f^*O_S)$, 
and with the notation of (iii) $DR(O_X/\mathcal F)(O_X)=DR(O_X/\mathcal F)$.
\end{itemize}
\end{defi}

Dually, we have the filtered Spencer complex of a complex of filtered right D-module :
\begin{defi}\label{DOSSP}
\begin{itemize}
\item[(i)]Let $(\mathcal S,O_S)\in\RCat$ with $O_S$ commutative. 
Let $(M,F)\in C_{D(O_S)^{op}fil}(\mathcal S)$. By proposition \ref{RSpaceDmod}, we have the complex
\begin{equation*}
SP(O_S)(M,F):=(T_{O_S}^{\bullet},F_b)\otimes_{O_S}(M,F)\in C_{fil}(\mathcal S)
\end{equation*}
whose differentials are, for $X\in\mathcal S$, and 
$\partial_1\wedge\cdots\wedge\partial_r\otimes m\in\Gamma(X,T^{r-1}_{O_S}\otimes_{O_S}M)$, 
\begin{eqnarray*}
d(\partial_1\wedge\cdots\wedge\partial_r\otimes m):(\omega\in\Gamma(X,\Omega^{r-1}_{O_S})\mapsto 
\sum_i\omega(\partial_1\wedge\cdots\wedge\partial_{\hat i}\cdots\partial_r)m-\sum_{i<j}\omega([\partial_i,\partial_j])m).
\end{eqnarray*}
\item[(ii)] More generally, let $f:(\mathcal X,O_X)\to (\mathcal S,O_S)$ with $(\mathcal X,O_X),(\mathcal S,O_S)\in\RCat$. 
The quotient map $q:\Omega_{O_X}\to\Omega_{O_X/f^*O_S}$ induce, for $G\in\PSh_{O_X}(\mathcal X)$ the injective map 
\begin{equation*}
q^{\vee,p}(G):=\wedge^pq^{\vee}\otimes I:T^q_{O_X/f^*O_S}\otimes_{O_X}G\to T^q_{O_X}\otimes_{O_X}G.
\end{equation*}
Let $(M,F)\in C_{D(O_X)^{op}fil}(\mathcal X)$. By proposition \ref{RSpaceDmod}, we have the relative Spencer complex
\begin{equation*}
SP(O_X/f^*O_S)(M,F):=(T_{X/S}^{\bullet},F_b)\otimes_{O_X}(M,F)\in C_{f^*O_Sfil}(\mathcal X)
\end{equation*}
whose differentials are the one of $SP(O_X)(M,F)$ given in (i) by the embedding 
$q^{\vee}:SP(O_X/f^*O_S)(M,F)\hookrightarrow SP(O_X)(M,F)$.
\item[(iii)] Let $(X,O_X)/\mathcal F\in FolRTop$, that is $(X,O_X)\in\RTop$ endowed with a foliation
with quotient map $q:\Omega_{O_X}\to\Omega_{O_X/\mathcal F}$.
Let $(M,F)\in C_{D(O_X)^{op}fil}(X)$. By proposition \ref{RSpaceDmod}, we have the foliated Spencer complex
\begin{equation*}
SP(O_X/\mathcal F)(M,F):=(T_{O_X/\mathcal F}^{\bullet},F_b)\otimes_{O_X}(M,F)\in C_{fil}(X)
\end{equation*}
whose differentials are of $SP(O_X)(M,F)$ given in (i) by the embedding 
$q^{\vee}:SP(O_X/\mathcal F)(M,F)\hookrightarrow SP(O_X)(M,F)$.
\end{itemize}
By definition, with the notation of (ii) if $\phi:(M_1,F)\to (M_2,F)$ is a morphism in $C_{D(O_X)^{op}fil}(\mathcal X)$, 
\begin{equation*}
SP(O_X/f^*O_S)(\phi):=(I\otimes\phi):(T_{X/S}^{\bullet},F_b)\otimes_{O_X}(M_1,F)\to 
(T_{X/S}^{\bullet},F_b)\otimes_{O_X}(M_2,F) 
\end{equation*}
is a morphism in $C_{f^*O_Sfil}(\mathcal X)$.
\end{defi}

\begin{prop}\label{otimesDmodsection2}
\begin{itemize}
\item[(i)] Let $f:(X,O_X)\to (S,O_S)$ a morphism with $(S,O_S),(X,O_X)\in\RTop$.
Assume that the canonical map $T(f,\hom)(O_X,O_X):f^*D(O_X)\to D(f^*O_X)$ is an isomorphism of sheaves.
For $(M,F)\in C_{D(O_X)^{op},f^*D(O_S)fil}(X)$ and $(M',F),(N,F)\in C_{D(O_X)fil}(X)$, 
we have canonical isomorphisms in $C_{f^*D(O_S)fil}(X)$ :
\begin{eqnarray*}
(M',F)\otimes_{O_X}(N,F)\otimes_{D(O_X)}(M,F)&=&(M',F)\otimes_{D(O_X)}((M,F)\otimes_{O_X}(N,F)) \\
&=&((M',F)\otimes_{O_X}(M,F))\otimes_{D(O_X)}(N,F)
\end{eqnarray*}
\item[(ii)] Let $f:(X,O_X)\to (S,O_S)$ a morphism with $(S,O_S),(X,O_X)\in\RTop$.
For $(M,F)\in C_{D(O_X)fil}(X)$, we have a canonical isomorphisms of filtered $f^*O_S$ modules, 
i.e. isomorphisms in $C_{f^*O_Sfil}(X)$,
\begin{eqnarray*}
(\Omega^{\bullet}_{O_X/f^*O_S},F_b)\otimes_{O_X}(M,F)=((\Omega^{\bullet}_{O_X/f^*O_S},F_b)\otimes_{O_X}D(O_X))\otimes_{D(O_X)}(M,F)
\end{eqnarray*}
\end{itemize}
\end{prop}

\begin{proof}
These are standard fact of algebra.
\end{proof}

\begin{defiprop}\label{TDwM}
Consider a commutative diagram in $\RCat$
\begin{equation*}
D=\xymatrix{(\mathcal X,O_X)\ar[r]^{f} & (\mathcal S,O_S) \\
(\mathcal X',O_{X'})\ar[u]^{g'}\ar[r]^{f'} & (\mathcal T,O_T)\ar[u]^{g}}.
\end{equation*}
with commutative structural sheaf of rings.
Assume that the canonical map $T(g',\hom)(O_X,O_X):g^{'*}D(O_X)\to D(g^{'*}O_X)$ is an isomorphism of sheaves.
\begin{itemize}
\item[(i)] For $(M,F)\in\PSh_{D(O_X)fil}(\mathcal X)$, the graded map in $(\PSh_{g^{'*}O_X}(\mathbb N\times\mathcal X'),F)$
\begin{eqnarray*}
\Omega_{(O_{X'}/g^{'*}O_X)/(O_T/g^*O_S)}(M,F):=m'\circ(\Omega_{(O_{X'}/g^{'*}O_X)/(O_T/g^*O_S)}\otimes I): \\
g^{'*}((\Omega^{\bullet}_{O_X/f^*O_S},F_b)\otimes_{O_X}(M,F))\to
(\Omega^{\bullet}_{O_{X'}/f^{'*}O_T},F_b)\otimes_{O_{X'}}g^{'*mod}(M,F) 
\end{eqnarray*}
given in degree $p\in\mathbb N$ by, for $X^{'o}\in\mathcal X'$ and $X^o\in\mathcal X$ such that $g^{'*}(X^o)\leftarrow X^{'o}$,   
\begin{eqnarray*}
\Omega^p_{(O_{X'}/g^{'*}O_X)/(O_T/g^*O_S)}(M)(X'^o):=m'\circ(\Omega^p_{(O_{X'}/g^{'*}O_X)/(O_T/g^*O_S)}\otimes I)(X'^o): \\
\omega\otimes m\in\Gamma(X^o,\Omega^p_{O_X}\otimes_{O_X}M)\mapsto\Omega_{O_{X'}/g^{'*}O_X}(\omega)\otimes(m\otimes 1) 
\end{eqnarray*}
is a map of complexes, that is a map in $C_{(f\circ g')^*O_Sfil}(\mathcal X')$.
\item[(ii)] For $(M,F)\in C_{D(O_X)fil}(\mathcal X)$, we get from (i) by functoriality,
the map in $C_{(f\circ g')^*O_Sfil}(\mathcal X')$
\begin{eqnarray*}
\Omega_{(O_{X'}/g^{'*}O_X)/(O_T/g^*O_S)}(M,F):=m'\circ(\Omega_{(O_{X'}/g^{'*}O_X)/(O_T/g^*O_S)}\otimes I): \\
g^{'*}((\Omega^{\bullet}_{O_X/f^*O_S},F_b)\otimes_{O_X}(M,F))\to
(\Omega^{\bullet}_{O_{X'}/f^{'*}O_T},F_b)\otimes_{O_{X'}}g^{'*mod}(M,F) 
\end{eqnarray*}
\item[(iii)]For $(M,F)\in C_{D(O_X)fil}(\mathcal X)$, we get from (ii) 
the canonical transformation map in $C_{O_Tfil}(\mathcal T)$
\begin{eqnarray*}
T^O_{\omega}(D)(M,F):g^{*mod}L_O(f_*E((\Omega^{\bullet}_{O_X/f^*O_S},F_b)\otimes_{O_X}(M,F)))
\xrightarrow{q} \\
(g^*f_*E((\Omega^{\bullet}_{O_X/f^*O_S},F_b)\otimes_{O_X}(M,F)))\otimes_{g^*O_S}O_T  
\xrightarrow{T(g',E)(-)\circ T(D)(E(\Omega^{\bullet}_{O_X/f^*O_S}\otimes_{O_X}(M,F)))} \\
(f'_*E(g'^*((\Omega^{\bullet}_{O_X/f^*O_S},F_b)\otimes_{O_X}(M,F))))\otimes_{g^*O_S}O_T 
\xrightarrow{E(\Omega_{(O_{X'}/g^{'*}O_X)/(O_T/g^*O_S})(M,F))} \\
(f'_*E((\Omega^{\bullet}_{O_{X'}/f^{'*}O_T},F_b)\otimes_{O_X}(M,F)))\otimes_{g^*O_S}O_T 
\xrightarrow{m}
f'_*E((\Omega^{\bullet}_{O_{X'}/f^{'*}O_T},F_b)\otimes_{O_{X'}}g^{'*mod}(M,F))
\end{eqnarray*}
with $m(n\otimes s)=s.n$.  
\end{itemize} 
\end{defiprop}

\begin{proof}
\noindent(i): We check that the map in $(\PSh_{g^{'*}O_X}(\mathbb N\times\mathcal X'),F)$
\begin{eqnarray*}
\Omega_{(O_{X'}/g^{'*}O_X)/(O_T/g^*O_S)}(M,F):=m'\circ(\Omega_{(O_{X'}/g^{'*}O_X)/(O_T/g^*O_S)}\otimes I): \\
g^{'*}((\Omega^{\bullet}_{O_X/f^*O_S},F_b)\otimes_{O_X}(M,F))\to
(\Omega^{\bullet}_{O_{X'}/f^{'*}O_T},F_b)\otimes_{O_{X'}}g^{'*mod}(M,F) 
\end{eqnarray*}
is a map in $C_{(f\circ g')^*O_Sfil}(\mathcal X')$. 
But we have, for $X^{'o}\in\mathcal X'$ the following equality in 
$\Gamma(X^{'o},\Omega^{p+1}_{O_{X'}}\otimes_{O_{X'}}g^{'*mod}M)$
\begin{eqnarray*}
d(\Omega^p_{(O_{X'}/g^{'*}O_X)(O_T/g^*O_S)}(M)(\omega\otimes m)):&=&d(\Omega^p_{O_{X'}/g^{'*}O_X}(\omega)\otimes (m\otimes 1)) \\
:&=&d(\Omega^p_{O_{X'}/g^{'*}O_X}(\omega))\otimes(m\otimes 1)+\Omega^p_{O_{X'}/g^{'*}O_X}(\omega)\otimes\nabla(m\otimes 1) \\
&=&\Omega^{p+1}_{O_{X'}/g^{'*}O_X}(d\omega)\otimes(m\otimes 1)+\Omega^p_{O_{X'}/g^{'*}O_X}(\omega)\otimes\nabla(m)\otimes 1 \\
&=&\Omega^{p+1}_{O_{X'}/g^{'*}O_X}(d\omega)\otimes (m\otimes 1)+\Omega^{p+1}_{O_{X'}/g^{'*}O_X}(\omega\otimes\nabla(m))\otimes 1 \\ 
&=&:\Omega^{p+1}_{(O_{X'}/g^{'*}O_X)(O_T/g^*O_S)}(M)(d(\omega)\otimes m+\omega\otimes\nabla(m)) \\
&=&:\Omega^{p+1}_{(O_{X'}/g^{'*}O_X)(O_T/g^*O_S)}(M)(d(\omega\otimes m))
\end{eqnarray*}
since for $\partial'\in T_{O_{X'}}(X^{'o})$,
\begin{equation*}
\nabla_{\partial'}(m\otimes 1)=\nabla_{dg'(\partial')}(m)\otimes 1+m\otimes\nabla_{\partial'} 1
=\nabla_{dg'(\partial')}(m)\otimes 1 : 
\end{equation*}
see in definition-proposition \ref{RSpacef} the definition of the $D(O_{X'})$ module structure 
on the $O_{X'}$ module $g^{'*mod}M:=g^{'*}M\otimes_{g^{'*}O_X}O_{X'}$.

\noindent(ii) and (iii): There is nothing to prove.
\end{proof}

\begin{rem}\label{TDwMrem}
Consider a commutative diagram in $\RCat$
\begin{equation*}
D=\xymatrix{(\mathcal X,O_X)\ar[r]^{f} & (\mathcal S,O_S) \\
(\mathcal X',O_{X'})\ar[u]^{g'}\ar[r]^{f'} & (\mathcal T,O_T)\ar[u]^{g}}.
\end{equation*}
Assume that the canonical map $T(g',\hom)(O_X,O_X):g^{'*}D(O_X)\to D(g^{'*}O_X)$ is an isomorphism of sheaves.
Under the canonical isomophism 
$(-)\otimes 1:(\Omega^{\bullet}_{O_X/f^*O_S},F_b)\xrightarrow{\sim}(\Omega^{\bullet}_{O_X/f^*O_S},F_b)\otimes_{O_X}(O_X,F_b)$,
we have (see definition-proposition \ref{TDwM} and definition \ref{TDw})
\begin{itemize}
\item $\Omega_{(O_{X'}/g^{'*}O_X)/(O_T/g^*O_S)}(O_X)=\Omega_{(O_{X'}/g^{'*}O_X)/(O_T/g^*O_S)}:
g^{'*}\Omega^{\bullet}_{O_X/f^*O_S}\to\Omega^{\bullet}_{O_{X'}/f^{'*}O_T}$ 
\item $T^O_{\omega}(D)(O_X)=T^O_{\omega}(D):
g^{*mod}L_O(f_*E(\Omega^{\bullet}_{O_X/f^*O_S},F_b))\to f'_*E(\Omega^{\bullet}_{O_{X'}/f^{'*}O_T},F_b)$.
\end{itemize}
\end{rem}

\begin{defi}\label{TDotimeswM}
Consider a commutative diagram in $\RCat$
\begin{equation*}
D=\xymatrix{(\mathcal X,O_X)\ar[r]^{f} & (\mathcal S,O_S) \\
(\mathcal X',O_{X'})\ar[u]^{g'}\ar[r]^{f'} & (\mathcal T,O_T)\ar[u]^{g}}.
\end{equation*}
with commutative structural sheaf of rings.
Assume that the canonical map $T(g',\hom)(O_X,O_X):g^{'*}D(O_X)\to D(g^{'*}O_X)$ is an isomorphism of sheaves.
For $(N,F)\in C_{D(O_{X'}),g^{'*}D(O_X)fil}(\mathcal X')$, 
we have by definition-proposition \ref{TDwM} the map in $C_{f^*O_Sfil}(\mathcal X)$ 
\begin{eqnarray*}
T^O_{\omega}(g',\otimes)(N,F):\Omega_{O_X/f^*O_S}^{\bullet}\otimes_{O_X}g'_*(N,F)\xrightarrow{\ad(g^{'*mod},g'_*)(-)} \\
g'_*(g^{'*}(\Omega_{O_X/f^*O_S}^{\bullet}\otimes_{O_X}g'_*N)\otimes_{g^{'*}O_X}O_{X'} 
\xrightarrow{m\circ\Omega_{(O_{X'}/g^{'*}O_X)/(O_T/g^*O_S)}(g'_*(N,F))} \\ 
g'_*(\Omega_{O_X/f^*O_S}^{\bullet}\otimes_{O_X}g^{'*mod}g'_*(N,F))\xrightarrow{\ad(g^{'*mod},g'_*)(N,F)}
g'_*(\Omega_{O_{X'}/f^*O_T}^{\bullet}\otimes_{O_{X'}}(N,F))
\end{eqnarray*}
with $m(n\otimes s)=s.n$ and $g'_*N\in C_{D(O_X)}(\mathcal X)$, 
the structure of $D(O_X)$ module being given by the canonical morphism
$\ad(g^{'*},g'_*)(D(O_X)):D(O_X)\to g'_*g^{'*}D(O_X)$ applied to $g'_*N\in C_{g'_*g^{'*}D(O_X)}(\mathcal X)$.
\end{defi}

We finish this subsection by a proposition for ringed spaces similar to proposition \ref{projformula}

\begin{prop}\label{projformulaw}
Let $f:(X,O_X)\to (S,O_S)$ a morphism with $(X,O_X),(S,O_S)\in\RTop$ with commutative sheaves of rings. 
Assume that $\Omega_{O_X/f^*O_S}\in\PSh_{O_X}(X)$ is a locally free $O_X$ module of finite rank.
\begin{itemize}
\item[(i)] If $\phi:(M,F)\to (N,F)$ is an $r$-filtered top local equivalence with $(M,F),(N,F)\in C_{D(O_X)fil}(X)$, then
\begin{equation*}
DR(O_X/f^*O_S)(\phi):(\Omega_{O_X/f^*O_S}^{\bullet},F_b)\otimes_{O_X}(M,F)\to
(\Omega_{O_X/f^*O_S}^{\bullet},F_b)\otimes_{O_T}(N,F) 
\end{equation*}
is an $r$-filtered top local equivalence. 
\item[(ii)] Consider a commutative diagram in $\RTop$
\begin{equation*}
D=\xymatrix{(X,O_X)\ar[r]^{f} & (S,O_S) \\
(X',O_{X'})\ar[u]^{g'}\ar[r]^{f'} & (T,O_T)\ar[u]^{g}}.
\end{equation*}
with commutative structural sheaf of rings. For $(N,F)\in C_{D(O_{X_T})fil}(X')$, the map in $C_{f^*O_Sfil}(X)$
\begin{eqnarray*}
k\circ T^O_{\omega}(g',\otimes)(E(N,F)):(\Omega_{O_X/f^*O_S}^{\bullet},F_b)\otimes_{O_X}g'_*E(N,F)\to  
g'_*E((\Omega_{O_{X'}/f^*O_T}^{\bullet},F_b)\otimes_{O_{X'}}E(N,F)) 
\end{eqnarray*}
is a filtered top local equivalence (see definition \ref{TDotimeswM}).
\end{itemize}
\end{prop}

\begin{proof}
\noindent(i):Follows from proposition \ref{projformula} (i) since $\Omega_{O_X/f^*O_S}^{\bullet}\in C^b(X)$ is then a bounded complex 
with $\Omega_{O_X/f^*O_S}^n\in\PSh_{O_X}(X)$ a locally free $O_X$ module of finite rank.

\noindent(ii):Follows from proposition \ref{projformula} (ii) since $\Omega_{O_X/f^*O_S}^{\bullet}\in C^b(X)$ is then a bounded complex 
with $\Omega_{O_X/f^*O_S}^n\in\PSh_{O_X}(X)$ a locally free $O_X$ module of finite rank.
\end{proof}

\subsubsection{The support section functor for D module on ringed spaces}

Let $(S,O_S)\in\RTop$ with $O_S$ commutative. Let $Z\subset S$ a closed subset.
Denote by $j:S\backslash Z\hookrightarrow S$ the open complementary embedding, 
\begin{itemize}
\item For $G\in C_{D(O_S)}(S)$, $\Gamma_ZG:=\Cone(\ad(j^*,j_*)(G):F\to j_*j^*G)[-1]$ has a (unique) structure of $D(O_S)$ module
such that $\gamma_Z(G):\Gamma_ZG\to G$ is a map in $C_{D(O_S)}(S)$. This gives the functor
\begin{equation*}
\Gamma_Z:C_{D(O_S)fil}(S)\to C_{D(O_S)fil}(S), \; (G,F)\mapsto\Gamma_Z(G,F) 
\end{equation*}
together with the canonical map $\gamma_Z(G,F):\Gamma_Z(G,F)\to (G,F)$.
Let $Z_2\subset Z$ a closed subset, then for $G\in C_{D(O_S)}(S)$, 
$T(Z_2/Z,\gamma)(G):\Gamma_{Z_2}G\to\Gamma_Z G$ is a map in $C_{D(O_S)}(S)$.
\item For $G\in C_{O_S}(S)$, $\Gamma^{\vee}_ZG:=\Cone(\ad(j_!,j^*)(G):j_!j^*G\to G)$ has a unique structure of $D(O_S)$ module,
such that $\gamma^{\vee}_Z(G):G\to\Gamma_Z^{\vee}G$ is a map in $C_{D(O_S)}(S)$. This gives the functor
\begin{equation*}
\Gamma^{\vee}_Z:C_{D(O_S)fil}(S)\to C_{D(O_S)fil}(S), \; (G,F)\mapsto\Gamma^{\vee}_Z(G,F), 
\end{equation*}
together with the canonical map $\gamma^{\vee}_Z(G,F):(G,F)\to\Gamma^{\vee}_Z(G,F)$.
Let $Z_2\subset Z$ a closed subset, then for $G\in C_{D(O_S)}(S)$, 
$T(Z_2/Z,\gamma^{\vee})(G):\Gamma_Z^{\vee}G\to\Gamma^{\vee}_{Z_2} G$ is a map in $C_{D(O_S)}(S)$.
\item For $G\in C_{D(O_S)}(S)$, 
\begin{eqnarray*}
\Gamma_Z^{\vee,h} G:&=&\mathbb D_S^OL_O\Gamma_ZE(\mathbb D^O_SG) \\
:&=&\Cone(\mathbb D^O_SL_O\ad(j_*,j^*)(E(\mathbb D^O_SG)):\mathbb D^O_SL_Oj_*j^*E(\mathbb D^O_SG)\to\mathbb D^O_SL_OE(\mathbb D^O_SG))
\end{eqnarray*}
has also canonical $D(O_S)$-module structure, and $\gamma_Z^{\vee,h}(G):G\to\Gamma_Z^{\vee,h}G$ is a map in $C_{D(O_S)}$.
This gives the functor
\begin{equation*}
\Gamma^{\vee,h}_Z:C_{D(O_S)fil}(S)\to C_{D(O_S)fil}(S), \; (G,F)\mapsto\Gamma^{\vee,h}_Z(G,F), 
\end{equation*}
together with the canonical map $\gamma^{\vee,h}_Z(G,F):(G,F)\to\Gamma^{\vee,h}_Z(G,F)$.
\item  Consider $\mathcal I^o_Z\subset O_S$ the ideal of vanishing function on $Z$
and $\mathcal I_Z\subset D_S$ the right ideal of $D_S$ generated by $\mathcal I^o_Z$.
We have then $\mathcal I^D_Z\subset\mathcal I_Z$, 
where $\mathcal I^D_Z\subset D_S$ is the left and right ideal consisting of sections which vanish on $Z$.
For $G\in\PSh_{D(O_S)}(S)$, we consider, $S^o\subset S$ being an open subset, 
\begin{equation*}
\mathcal I_ZG(S^o)=<\left\{f.m, m\in G(S^o),f\in\mathcal I_Z(S^o)\right\}>\subset G(S^o) 
\end{equation*}
the $D(O_S)$-submodule generated by the functions which vanish on $Z$ ($\mathcal I_Z$ is a right $D(O_S)$ ideal), 
This gives the functor,
\begin{eqnarray*}
\Gamma^{\vee,O}_Z:=\Gamma^{\vee,O,I_Z}_Z:C_{D(O_S)fil}(S)\to C_{D(O_S)fil}(S), \\ 
(G,F)\mapsto\Gamma_Z^{\vee,O}(G,F):=\Cone(b_Z(G,F):\mathcal I_Z(G,F)\to (G,F)), \, b_Z(-):=b_{I_Z}(-)
\end{eqnarray*}
together with the canonical map $\gamma_Z^{\vee,O}(G,F):(G,F)\to\Gamma_Z^{\vee,O}(G,F)$.
which factors through
\begin{equation*}
\gamma_Z^{\vee,O}(G):G\xrightarrow{\gamma_Z^{\vee}(G)}\Gamma^{\vee}_ZG\xrightarrow{b_{S/Z}(G)}\Gamma^{\vee,O}_ZG.
\end{equation*}
with $b_{S/Z}(-)=b_{S/Z}^I$and we have an homotopy equivalence $c_Z(G):=c_{I_Z}(G):\Gamma_Z^{\vee,O}G\to G/\mathcal I_ZG$.

\end{itemize}

\begin{lem}\label{gammaotimessection2}
Let $(Y,O_Y)\in\RTop$ and $i:X\hookrightarrow Y$ a closed embedding.
\begin{itemize}
\item[(i)] For $(M,F)\in C_{D(O_Y)fil}(Y)$ and 
$(N,F)\in\PSh_{D(O_Y)^{op}fil}(Y)$ such that $a_{\tau}N$ is a locally free $D(O_Y))$ module of finite rank, the canonical map
\begin{eqnarray*}
T(\gamma,\otimes)(E(M,F),(N,F)):=(I,T(j,\otimes)(E(M,F),(N,F))): \\
(\Gamma_XE(M,F))\otimes_{D(O_Y)}(N,F)\to\Gamma_XE((M,F)\otimes_{D(O_Y)}(N,F))
\end{eqnarray*}
is an equivalence top local.
\item[(ii)] For $(M,F)\in C_{D(O_Y)^{(op)}fil}(Y)$ and 
$(N,F)\in \PSh_{D(O_Y)^{(op)}fil}(Y)$ such that $a_{\tau}N$ is a locally free $O_Y$ module of finite rank, the canonical map 
\begin{eqnarray*}
T(\gamma,\otimes)(E(M,F),(N,F)):=(I,T(j,\otimes)(E(M,F),(N,F))): \\
(\Gamma_XE(M,F))\otimes_{O_Y}(N,F)\to\Gamma_XE((M,F)\otimes_{O_Y}(N,F))
\end{eqnarray*}
is a filtered top local equivalence.
\end{itemize}
\end{lem}

\begin{proof}
Follows from proposition \ref{projformula}. Also note that $T(j,\otimes)(-,-)=T^{mod}(j,\otimes)(-,-)$.
\end{proof}

We now look at the pullback map and the transformation map of De Rahm complexes together with the support section functor.
The follwoing is a generalization of definition-proposition \ref{TDwgamma} :

\begin{defiprop}\label{TDwMgamma}
Consider a commutative diagram in $\RTop$
\begin{equation*}
D_0=\xymatrix{f: (X,O_X)\ar[r]^{i} & (Y,O_Y)\ar[r]^{p} & (S,O_S) \\
f':(X',O_{X'})\ar[r]^{i'}\ar[u]^{g'} & (Y',O_{Y'})\ar[u]^{g''}\ar[r]^{p'} & (T,O_T)\ar[u]^{g} }
\end{equation*}
with $i$, $i'$ being closed embeddings. Denote by $D$ the right square of $D$. 
We have a factorization $i':X'\xrightarrow{i'_1}X\times Y Y'\xrightarrow{i'_0}Y'$, where $i'_0,i'_1$ are closed embedding. 
\begin{itemize}
\item[(i)] For $(M,F)\in C_{\mathcal D(O_Y)fil}(Y)$, the canonical map,
\begin{eqnarray*}
E(\Omega_{(O_{Y'}/g^{''*}O_Y)/(O_T/g^*O_S)}(M,F))\circ T(g'',E)(-)\circ T(g'',\gamma)(-): \\
g^{''*}\Gamma_{X}E((\Omega^{\bullet}_{O_Y/p^*O_S},F_b)\otimes_{O_Y} (M,F))\to 
\Gamma_{X\times Y Y'}E((\Omega^{\bullet}_{O_{Y'}/p^{'*}O_T},F_b)\otimes_{O_{Y'}}g^{''*mod}(M,F))
\end{eqnarray*}
unique up to homotopy such that the following diagram in $C_{g^{''*}p^*O_Sfil}(Y')=C_{p^{'*}g^*O_Sfil}(Y')$ commutes
\begin{equation*}
\xymatrix{g^{''*}\Gamma_{X}E((\Omega^{\bullet}_{O_Y/p^*O_S},F_b)\otimes_{O_Y}(M,F))
\ar[rr]^{E(\Omega_{(-)/(-)}(M,F))\circ T(g'',E)(-)\circ T(g'',\gamma)(-)}\ar[d]_{\gamma_X(-)} & \, & 
\Gamma_{X\times_Y Y'}E((\Omega^{\bullet}_{O_{Y'}/p^{'*}O_T},F_b)\otimes_{O_{Y'}}(g^{''*mod}(M,F)))
\ar[d]^{\gamma_{X\times Y Y'}(-)} \\
g^{''*}E((\Omega^{\bullet}_{O_Y/p^*O_S},F_b)\otimes_{O_Y}(M,F))
\ar[rr]^{E(\Omega_{(-)/(-)}(M,F)\circ T(g'',E)(-)} & \, &   
E((\Omega^{\bullet}_{O_{Y'}/p^{'*}O_T},F_b)\otimes_{O_{Y'}}g^{''*mod}(M,F))}.
\end{equation*}
\item[(ii)] For $M\in C_{\mathcal D}(Y)$, there is a canonical map
\begin{eqnarray*}
T^O_{\omega}(D)(M,F)^{\gamma}:g^{*mod}L_Op_*\Gamma_{X}E((\Omega^{\bullet}_{O_Y/p^*O_S},F_b)\otimes_{O_Y}(M,F))\to \\
p'_*\Gamma_{X\times_Y Y'}E((\Omega^{\bullet}_{O_{Y'}/p^{'*}O_T},F_b)\otimes_{O_{Y'}}g^{''*mod}(M,F))
\end{eqnarray*}
unique up to homotopy such that the following diagram in $C_{O_{T}fil}(T)$ commutes
\begin{equation*}
\xymatrix{g^{*mod}L_Op_*\Gamma_{X}E((\Omega^{\bullet}_{O_Y/p^*O_S},F_b)\otimes_{O_{Y}}(M,F))
\ar[r]^{T_{\omega}^O(D)(M,F)^{\gamma}}\ar[d]_{\gamma_X(-)} &  
p'_*\Gamma_{X\times_Y Y'}E((\Omega^{\bullet}_{O_{Y'}/p^{'*}O_T},F_b)\otimes_{O_{Y'}}(g^{''*mod}(M,F)))
\ar[d]^{\gamma_{X\times_Y Y'}(-)} \\
g^{*mod}L_Op_*E((\Omega^{\bullet}_{O_Y/p^*O_S},F_b)\otimes_{O_{Y}}(M,F))
\ar[r]^{T^O_{\omega}(D)(M,F)} &  
p'_*E((\Omega^{\bullet}_{O_{Y'}/p^{'*}O_T},F_b)\otimes_{O_{Y'}}g^{''*mod}(M,F))}.
\end{equation*}
\item[(iii)] For $N\in C_{\mathcal D}(Y\times T)$, the canonical map in $C_{h^{'*}O_Tfil}(Y')$
\begin{equation*}
T(X'/X\times_Y Y',\gamma)(-):\Gamma_{X'}E((\Omega^{\bullet}_{Y'/T},F_b)\otimes_{O_{Y'}}(N,F))\to
\Gamma_{X\times_Y Y'}E((\Omega^{\bullet}_{O_{Y'}/O_T},F_b)\otimes_{O_{Y'}}(N,F))
\end{equation*}
is unique up to homotopy such that $\gamma_{X\times_Y Y'}(-)\circ T(X'/X\times_ Y Y',\gamma)(-)=\gamma_{X'}(-)$.
\item[(iv)] For $M=O_Y$, we have $T_{\omega}^O(D)(O_{Y\times S})^{\gamma}=T_{\omega}^O(D)^{\gamma}$ 
and $T_{\omega}^O(X\times_Y Y'/Y')(O_{Y'})^{\gamma}=T_{\omega}^O(X\times_Y Y'/Y')^{\gamma}$ 
(see definition-proposition \ref{TDwgamma}).
\end{itemize}
\end{defiprop}

\begin{proof}
Immediate from definition. We take for the map of point (ii) the composite
\begin{eqnarray*}
T^O_{\omega}(D)(M,F)^{\gamma}: 
g^{*mod}L_Op_*\Gamma_{X}E((\Omega^{\bullet}_{O_Y/p^*O_S},F_b)\otimes_{O_Y}(M,F))
\xrightarrow{q} \\
g^*p_*\Gamma_{X}E((\Omega^{\bullet}_{O_Y/p^*O_S},F_b)\otimes_{O_Y}(M,F))\otimes_{g^*O_S}O_T  
\xrightarrow{T(g'',E)(-)\circ T(g'',\gamma)(-)\circ T(D)(E(\Omega^{\bullet}_{O_Y/p^*O_S},F_b))} \\   
(p'_*\Gamma_{X\times_Y Y'}E(g^{''}((\Omega^{\bullet}_{O_Y/p^*O_S},F_b)\otimes_{O_Y}(M,F))))\otimes_{g^*O_S}O_T  
\xrightarrow{E(\Omega_{(O_{Y'}/g^{''*}O_Y)/(O_T/g^*O_S}(M,F)))} \\
p'_*\Gamma_{X\times_Y Y'}E((\Omega^{\bullet}_{O_{Y'}/p^{'*}O_T},F_b)\otimes_{O_{Y'}}g^{''*mod}(M,F))\otimes_{g^*O_S}O_T 
\xrightarrow{m} \\
p'_*\Gamma_{X\times_Y Y'}E((\Omega^{\bullet}_{O_{Y'}/p^{'*}O_T},F_b)\otimes_{O_{Y'}}g^{''*mod}(M,F)),
\end{eqnarray*}
with $m(n\otimes s)=s.n$.
\end{proof}

Let $p:(Y,O_Y)\to (S,O_S)$ a morphism with $(Y,O_Y),(S,O_S)\in\RTop$. Let $i:X\hookrightarrow Y$ a closed embedding.
Denote by $j:Y\backslash X\hookrightarrow Y$ the complementary open embedding. 
Consider, for $(M,F)\in C_{D(O_Y)fil}(Y)$, the map in $C_{p^*O_Sfil}(Y)$ (see definition \ref{TDotimeswM}): 
\begin{eqnarray*}
k\circ T^O_w(j,\otimes)(E(M,F)):(\Omega^{\bullet}_{O_Y/p^*O_S},F_b)\otimes_{O_Y}j_*j^*E(M,F)) \\
\xrightarrow{DR(O_Y/p^*O_S)(\ad(j^*,j_*)(-))} \\
j_*j^*((\Omega^{\bullet}_{O_Y/p^*O_S},F_b)\otimes_{O_Y}j_*j^*E(M,F))=
j_*j^*(\Omega^{\bullet}_{O_Y/p^*O_S},F_b)\otimes_{O_Y}j^*j_*j^*E(M,F) \\
\xrightarrow{k\circ DR(O_Y/p^*O_S)(\ad(j^*,j_*)(j^*E(M)))} \\
j_*E(j^*(\Omega^{\bullet}_{O_Y/p^*O_S},F_b)\otimes_{O_Y}j^*E(M,F))=
j_*E(j^*((\Omega^{\bullet}_{O_Y/p^*O_S},F_b)\otimes_{O_Y}E(M,F)))
\end{eqnarray*}

\begin{defi}\label{mw0def}
Let $p:(Y,O_Y)\to (S,O_S)$ a morphism with $(Y,O_Y),(S,O_S)\in\RTop$. Let $i:X\hookrightarrow Y$ a closed embedding.  
Denote by $j:Y\backslash X\hookrightarrow Y$ the complementary open embedding. 
We consider, for $(M,F)\in C_{D(O_Y)fil}(Y)$ the canonical map in $C_{p^*O_Sfil}(Y)$ 
\begin{eqnarray*}
T^O_w(\gamma,\otimes)(M,F):=(I,k\circ T^O_w(j,\otimes)(E(M,F))): \\
(\Omega^{\bullet}_{O_Y/p^*O_S},F_b)\otimes_{O_Y}\Gamma_XE(M,F)\to
\Gamma_XE((\Omega^{\bullet}_{O_Y/p^*O_S},F_b)\otimes_{O_Y}E(M,F)).
\end{eqnarray*}
\end{defi}

\begin{prop}\label{mw0propsection2}
Let $p:(Y,O_Y)\to (S,O_S)$ a morphism with $(Y,O_Y),(S,O_S)\in\RTop$. Let $i:X\hookrightarrow Y$ a closed embedding.
Then, if $\Omega_{O_Y/p^*O_S}$ is a locally free $O_Y$ module, for $(M,F)\in C_{D(O_Y)fil}(Y)$
\begin{itemize}
\item[(i)] the map 
\begin{equation*}
T^O_w(\gamma,\otimes)(M,F):(\Omega^{\bullet}_{O_Y/p^*O_S},F_b)\otimes_{O_Y}\Gamma_XE(M,F)\to
\Gamma_XE((\Omega^{\bullet}_{O_Y/p^*O_S},F_b)\otimes_{O_Y}E(M,F))
\end{equation*}
is a $1$-filtered top local equivalence,
\item[(ii)] the map in $D_{p^*O_Sfil}(Y)$
\begin{eqnarray*}
T^O_w(\gamma,\otimes):=DR(O_Y/p^*O_S)(k)^{-1}\circ T^O_w(\gamma,\otimes)(M,F): \\
(\Omega^{\bullet}_{O_Y/p^*O_S},F_b)\otimes_{O_Y}\Gamma_XE(M,F)\to
\Gamma_XE((\Omega^{\bullet}_{O_Y/p^*O_S},F_b)\otimes_{O_Y}(M,F)) 
\end{eqnarray*}
is an isomorphism.
\end{itemize}
\end{prop}

\begin{proof}
By proposition \ref{projformulaw}, 
\begin{itemize}
\item $\Gr_F^p(k\circ T^O_w(j,\otimes)(E(M,F))):
\Omega^{\bullet}_{O_Y/p^*O_S}\otimes_{O_Y}j_*j^*F^{p-\bullet}E(M)\to 
j_*E(j^*(\Omega^{\bullet}_{O_Y/p^*O_S}\otimes_{O_Y}F^{p-\bullet}E(M)))$
is a top local equivalence and 
\item $DR(O_Y/p^*O_S)(k):\Omega^{\bullet}_{O_Y/p^*O_S}\otimes_{O_Y}(M,F)\to\Omega^{\bullet}_{O_Y/p^*O_S}\otimes_{O_Y}E(M,F)$
is a filtered top local equivalence.
\end{itemize}
\end{proof}

\subsection{The D-modules on smooth complex algebraic varieties and on complex analytic maninfold 
and their functorialities in the filtered case}

For convenience, we will work with and state the results for presheaves of D-modules. 
In this section, it is possible to assume that all the presheaves are sheaves and take the sheaftification functor after the
pullback functor $f^*$ for a morphism $f:X\to S$, $X,S\in\Var(\mathbb C)$ or $X,S\in\AnSp(\mathbb C)$, 
and after the internal hom functors and tensor products of presheaves of modules on $S\in\Var(\mathbb C)$ or $S\in\AnSp(\mathbb C)$.

For $S=(S,O_S)\in\SmVar(\mathbb C)$, resp. $S=(S,O_S)\in\AnSm(\mathbb C)$, we denote by
\begin{itemize}
\item $D_S:=D(O_S)\subset\mathcal Hom_{\mathbb C_S}(O_S,O_S)$ the subsheaf consisting of differential operators.
By a $D_S$ module, we mean a left $D_S$ module.
\item we denote by
\begin{itemize}
\item $\PSh_{\mathcal D}(S)$ the abelian category of Zariski (resp. usu) presheaves on $S$ with a structure of left $D_S$ module,
and by $\PSh_{\mathcal D,h}(S)\subset\PSh_{\mathcal D,c}(S)\subset\PSh_{\mathcal D}(S)$ 
the full subcategories whose objects are coherent, resp. holonomic, sheaves of left $D_S$ modules,
and by $\PSh_{\mathcal D,rh}(S)\subset\PSh_{\mathcal D,h}(S)$ 
the full subcategory of regular holonomic sheaves of left $D_S$ modules,
\item $\PSh_{\mathcal D^{op}}(S)$ the abelian category of Zariski (resp. usu) presheaves on $S$ with a structure of right $D_S$ module,
and by $\PSh_{\mathcal D^{op},h}(S)\subset\PSh_{\mathcal D^{op},c}(S)\subset\PSh_{\mathcal D^{op}}(S)$ 
the full subcategories whose objects are coherent, resp. holonomic, sheaves of right $D_S$ modules,
and by $\PSh_{\mathcal D^{op},rh}(S)\subset\PSh_{\mathcal D^{op},h}(S)$ 
the full subcategory of regular holonomic sheaves of right $D_S$ modules,
\end{itemize}
\item we denote by
\begin{itemize}
\item $C_{\mathcal D}(S)=C(\PSh_{\mathcal D}(S))$ the category of complexes of Zariski presheaves on $S$ with a structure of $D_S$ module,
\begin{equation*}
C_{\mathcal D,rh}(S)\subset C_{\mathcal D,h}(S)\subset C_{\mathcal D,c}(S)\subset C_{\mathcal D}(S) 
\end{equation*}
the full subcategories consisting of complexes of presheaves $M$ such that 
$a_{\tau}H^n(M)$ are coherent, resp. holonomic, resp. regular holonomic, sheaves of $D_S$ modules, 
$a_{\tau}$ being the sheaftification functor for the Zariski, resp. usual, topology,
\item $C_{\mathcal D^{op}}(S)=C(\PSh_{\mathcal D^{op}}(S))$ the category of complexes of Zariski presheaves on $S$
with a structure of right $D_S$ module,
\begin{equation*}
C_{\mathcal D^{op},rh}(S)\subset C_{\mathcal D^{op},h}(S)\subset C_{\mathcal D^{op},c}(S)\subset C_{\mathcal D^{op}}(S)  
\end{equation*}
the full subcategories consisting of complexes of presheaves $M$ such that 
$a_{\tau}H^n(M)$ are coherent, resp. holonomic, resp. regular holonomic, sheaves of right $D_S$ modules,
\end{itemize}
\item in the filtered case we have
\begin{itemize}
\item $C_{\mathcal D(2)fil}(S)\subset C(\PSh_{\mathcal D}(S),F,W):=C(\PSh_{D(O_S)}(S),F,W)$ 
the category of (bi)filtered complexes of algebraic (resp. analytic) $D_S$ modules such that the filtration is biregular 
(see definition \ref{Dmodefinition},
\begin{equation*}
C_{\mathcal D(2)fil,rh}(S)\subset C_{\mathcal D(2)fil,h}(S)\subset C_{\mathcal D(2)fil,c}(S)\subset C_{\mathcal D(2)fil}(S), 
\end{equation*}
the full subcategories consisting of filtered complexes of presheaves $(M,F)$ such that 
$a_{\tau}H^n(M,F)$ are filtered coherent, resp. filtered holonomic, resp filtered regular holonomic, sheaves of $D_S$ modules, 
that is $a_{\tau}H^n(M)$ are coherent, resp. holonomic, resp. regular holonomic, sheaves of $D_S$ modules and 
$F$ induces a good filtration on $a_{\tau}H^n(M)$ 
(in particular $F^pa_{\tau}H^n(M)\subset a_{\tau}H^n(M)$ are coherent sub $O_S$ modules),
\item $C_{\mathcal D0fil}(S)\subset C_{\mathcal Dfil}(S)$ the full subcategory such that
the filtration is a filtration by $D_S$ submodule (which is stronger then Griffitz transversality), 
$C_{\mathcal D(1,0)fil}(S)\subset C_{\mathcal D2fil}(S)$
the full subcategory such that $W$ is a filtration by $D_S$ submodules (see definition \ref{Dmodefinition}),
\begin{equation*}
C_{\mathcal D(1,0)fil,h}(S)=C_{\mathcal D2fil,h}(S)\cap C_{\mathcal D(1,0)fil}(S)\subset C_{\mathcal D2fil,h}(S), 
\end{equation*}
the full subcategory consisting of filtered complexes of presheaves $(M,F,W)$ 
such that $a_{\tau}H^n(M,F)$ are filtered holonomic sheaves of $D_S$ modules and such that $W^pM\subset M$ are $D_S$ submodules
(recall that the $O_S$ submodules $F^pM\subset M$ are NOT $D_S$ submodules but satisfy by definition 
$md:F^rD_S\otimes F^pM\subset F^{p+r}M$),
\begin{equation*}
C_{\mathcal D(1,0)fil,rh}(S)=C_{\mathcal D2fil,rh}(S)\cap C_{\mathcal D(1,0)fil}(S)\subset C_{\mathcal D2fil,rh}(S), 
\end{equation*}
the full subcategory consisting of filtered complexes of presheaves $(M,F,W)$ 
such that $a_{\tau}H^n(M,F)$ are filtered regular holonomic sheaves of $D_S$ modules 
and such that $W^pM\subset M$ are $D_S$ submodules
\item $C_{\mathcal D^{op}(2)fil}(S)\subset C(\PSh_{\mathcal D^{op}}(S),F,W):=C(\PSh_{D(O_S)^{op}}(S),F,W)$ 
the category of (bi)filtered complexes of algebraic (resp. analytic) right $D_S$ modules such that the filtration is biregular,
as in the left case we consider the subcategories
\begin{equation*}
C_{\mathcal D^{op}(2)fil,rh}(S)\subset C_{\mathcal D^{op}(2)fil,h}(S)\subset C_{\mathcal D^{op}(2)fil,c}(S)
\subset C_{\mathcal D^{op}(2)fil}(S), 
\end{equation*}
the full subcategories consisting of filtered complexes of presheaves $(M,F)$ 
such that $a_{\tau}H^n(M,F)$ are filtered coherent, resp. filtered holonomic, resp. filtered regular holonomic,
sheaves of right $D_S$ modules.
\end{itemize}
\end{itemize}

For $S=(S,O_S)\in\AnSm(\mathbb C)$, we have the natural extension 
$D_S\subset D_S^{\infty}\subset\mathcal Hom_{\mathbb C_S}(O_S,O_S)$ where $D_S^{\infty}\subset\mathcal Hom_{\mathbb C_S}(O_S,O_S)$ 
is the subsheaf of differential operators of possibly infinite order (see \cite{Kashiwara})
for the definition of the action of a differential operator of infinite order on $O_S$)
Similarly, we have 
\begin{itemize}
\item $C_{\mathcal D^{\infty}(2)fil}(S)\subset C(\PSh_{\mathcal D^{\infty}}(S),F,W):=C(\PSh_{D_S^{\infty}}(S),F,W)$ 
the category of (bi)filtered complexes of $D^{\infty}_S$ modules such that the filtration is biregular,
\begin{equation*}
C_{\mathcal D^{\infty}(2)fil,h}(S)\subset C_{\mathcal D^{\infty}(2)fil,c}(S)\subset C_{\mathcal D^{\infty}(2)fil}(S), 
\end{equation*}
the full subcategories consisting of filtered complexes of presheaves $(M,F)$ 
such that $a_{\tau}H^n(M,F)$ are filtered coherent (resp. holonomic) sheaves of $D^{\infty}_S$ modules
that is $a_{\tau}H^n(M)$ are coherent (resp. holonomic) sheaves of $D^{\infty}_S$ modules and 
$F$ induces a good filtration on $a_{\tau}H^n(M)$.
\item $C_{\mathcal D^{\infty}0fil}(S)\subset C_{\mathcal D^{\infty}fil}(S)$ the full subcategory such that
the filtration is a filtration by $D^{\infty}_S$ submodule, 
$C_{\mathcal D^{\infty}(1,0)fil}(S)\subset C_{\mathcal D^{\infty}2fil}(S)$
the full subcategory such that $W$ is a filtarion by $D^{\infty}_S$ submodules,
\begin{equation*}
C_{\mathcal D^{\infty}(1,0)fil,h}(S)=C_{\mathcal D^{\infty}2fil,h}(S)\cap C_{\mathcal D^{\infty}(1,0)fil}(S)
\subset C_{\mathcal D^{\infty}2fil,h}(S), 
\end{equation*}
the full subcategory consisting of filtered complexes of presheaves $(M,F,W)$ 
such that $a_{\tau}H^n(M,F)$ are filtered holonomic sheaves of $D^{\infty}_S$ modules and such that $W^pM\subset M$ are $D_S$ submodules
\item $C_{\mathcal D^{\infty,op}(2)fil}(S)\subset C(\PSh_{\mathcal D^{\infty,op}}(S),F,W):=C(\PSh_{D_S^{\infty,op}}(S),F,W)$ 
the category of (bi)filtered complexes of right $D^{\infty}_S$ modules such that the filtration is biregular,
\begin{equation*}
C_{\mathcal D^{\infty,op}(2)fil,h}(S)\subset C_{\mathcal D^{\infty,op}(2)fil,c}(S)\subset C_{\mathcal D^{\infty,op}(2)fil}(S), 
\end{equation*}
the full subcategories consisting of filtered complexes of presheaves $(M,F)$ 
such that $a_{\tau}H^n(M,F)$ are filtered coherent (resp. holonomic) sheaves of $D_S$ modules.
\end{itemize}

For $f:X\to S$ a morphism with $X,S\in\SmVar(\mathbb C)$ or with $(X,S)\in\AnSm(\mathbb C)$,
\begin{itemize}
\item we denote by
\begin{itemize}
\item $\PSh_{f^*\mathcal D}(X)$ the abelian category of Zariski (resp. usu) presheaves on $S$ with a structure of left $f^*D_S$ module,
and $C_{f^*\mathcal D}(X)=C(\PSh_{f^*\mathcal D}(X))$,
\item $\PSh_{\mathcal D,f^*\mathcal D}(X)$ the abelian category of Zariski (resp. usu) presheaves on $S$ with a structure 
of left $f^*D_S$ module and left $D_X$ module, 
and $C_{\mathcal D,f^*\mathcal D}(X)=C(\PSh_{\mathcal D,f^*\mathcal D}(X))$, 
\item $\PSh_{\mathcal D^{op},f^*\mathcal D}(X)$ the abelian category of Zariski (resp. usu) presheaves on $S$ with a structure  
of left $f^*D_S$ module and right $D_X$ module
and $C_{\mathcal D^{op},f^*\mathcal D}(X)=C(\PSh_{\mathcal D^{op},f^*\mathcal D}(X))$,  
\end{itemize}
\item we denote by
\begin{itemize}
\item $C_{f^*\mathcal Dfil}(X)\subset C(\PSh_{f^*\mathcal D}(X),F):=C(\PSh_{f^*D(O_S)}(X),F)$ 
the category of filtered complexes of algebraic (resp. analytic) $f^*D_S$ modules such that the filtration is biregular,
\item $C_{\mathcal D,f^*\mathcal Dfil}(X)\subset C(\PSh_{\mathcal D,f^*\mathcal D}(X),F)$ 
the category of filtered complexes of algebraic (resp. analytic) $(f^*D_S,D_X)$ modules such that the filtration is biregular,
\item $C_{\mathcal D^{op},f^*\mathcal Dfil}(X)\subset C(\PSh_{\mathcal D^{op},f^*\mathcal D}(X),F)$ 
the category of filtered complexes of algebraic (resp. analytic) $(f^*D_S,D^{op}_X)$ modules such that the filtration is biregular.
\end{itemize}
\end{itemize}

For $f:X\to S$ a morphism with $X,S\in\AnSm(\mathbb C)$, we denote by
\begin{itemize}
\item $C_{f^*\mathcal D^{\infty}fil}(X)\subset C(\PSh_{f^*\mathcal D^{\infty}}(X),F):=C(\PSh_{f^*D_S^{\infty}}(X),F)$ 
the category of filtered complexes of $f^*D^{\infty}_S$ modules such that the filtration is biregular,
\item $C_{\mathcal D^{\infty},f^*\mathcal D^{\infty}fil}(X)\subset C(\PSh_{\mathcal D^{\infty},f^*\mathcal D^{\infty}}(X),F)$ 
the category of filtered complexes of $(f^*D^{\infty}_S,D^{\infty}_X)$ modules such that the filtration is biregular,
\item $C_{\mathcal D^{\infty,op},f^*\mathcal D^{\infty}fil}(X)\subset C(\PSh_{\mathcal D^{\infty,op},f^*\mathcal D^{\infty}}(X),F)$ 
the category of filtered complexes of $(f^*D^{\infty}_S,D^{\infty,op}_X)$ modules such that the filtration is biregular.
\end{itemize}

For $S\in\AnSm(\mathbb C)$, we denote by
\begin{equation*}
J_S:C_{\mathcal D(2)fil}(S)\to C_{\mathcal D^{\infty}(2)fil}(S), \; 
(M,F)\mapsto J_S(M,F):=(M,F)\otimes_{D_S}(D_S^{\infty},F^{ord})
\end{equation*}
the natural functor. For $(M,F)\in C_{\mathcal D^{\infty}fil}(S)$, we will consider the map 
\begin{equation*}
\mathcal J_S(M,F):J_S(M,F):=(M,F)\otimes_{D_S}(D_S^{\infty},F^{ord})\to (M,F), m\otimes P\mapsto Pm\; 
\end{equation*}
Of course $J_S(C_{\mathcal D(1,0)fil}(S))\subset C_{\mathcal D^{\infty}(1,0)fil}(S)$.
More generally, for $f:X\to S$ a morphism with $X,S\in\AnSm(\mathbb C)$, we denote by
\begin{equation*}
J_{X/S}:C_{f^*\mathcal D(2)fil}(X)\to C_{f^*\mathcal D^{\infty}(2)fil}(X), \; 
(M,F)\mapsto J_{X/S}(M,F):=(M,F)\otimes_{f^*(D_S,F)}f^*(D_S^{\infty},F)
\end{equation*}
the natural functor, together with, for $(M,F)\in C_{f^*\mathcal D^{\infty}fil}(X)$, the map $\mathcal J_S(M,F):J_S(M,F)\to (M,F)$.

\begin{defi}\label{DmodsuppZ}
Let $S\in\SmVar(\mathbb C)$, resp. $S\in\AnSm(\mathbb C)$.
Let $Z\subset S$ a closed subset and denote by $j:S\backslash Z\hookrightarrow S$ the open embedding.
\begin{itemize}
\item[(i)]We denote by 
\begin{itemize}
\item $\PSh_{\mathcal D,Z}(S)\subset \PSh_{\mathcal D}(S)$, 
the full subcategory consisting of presheaves $M\in\PSh_{\mathcal D}(S)$, such that $j^*M=0$,
\item $C_{\mathcal D,Z}(S)\subset C_{\mathcal D}(S)$, 
the full subcategory consisting of complexes presheaves $M\in C_{\mathcal D}(S)$ 
such that $a_{\tau}j^*H^nM=0$ for all $n\in\mathbb Z$,
\item $C_{\mathcal D,Z,h}(S):=C_{\mathcal D,Z}(S)\cap C_{\mathcal D,h}(S)\subset C_{\mathcal D}(S)$
the full subcategory consising of $M\in C_{\mathcal D}(S)$ such that $a_{\tau}H^n(M)$ are holonomic 
and $a_{\tau}j^*H^nM=0$ for all $n\in\mathbb Z$,
\item $C_{\mathcal D,Z,c}(S):=C_{\mathcal D,Z}(S)\cap C_{\mathcal D,c}(S)\subset C_{\mathcal D}(S)$
the full subcategory consising of $M\in C_{\mathcal D}(S)$ such that $a_{\tau}H^n(M)$ are coherent and 
$a_{\tau}j^*H^nM=0$ for all $n\in\mathbb Z$.
\end{itemize}
\item[(ii)] We denote by
\begin{itemize}
\item $C_{\mathcal D(2)fil,Z}(S)\subset C_{\mathcal D(2)fil}(S)$, 
the full subcategory consisting of $(M,F)\in C_{\mathcal D(2)fil}(S)$ 
such that there exists $r\in\mathbb N$ and an $r$-filtered homotopy equivalence $\phi:(M,F)\to(N,F)$ 
with $(N,F)\in C_{\mathcal D(2)fil}(S)$ such that $a_{\tau}j^*H^n\Gr_F^p(M,F)=0$ for all $n,p\in\mathbb Z$,
note that by definition this $r$ does NOT depend on $n$ and $p$, 
\item $C_{\mathcal D(2)fil,Z,rh}(S):=C_{\mathcal D(2)fil,Z}(S)\cap C_{\mathcal D(2)fil,rh}(S)\subset C_{\mathcal D(2)fil}(S)$
the full subcategory consising of $(M,F)$ such that $a_{\tau}H^n(M,F)$ are filtered regular holonomic for all $n\in\mathbb Z$ and
such that there exists $r\in\mathbb N$ and an $r$-filtered homotopy equivalence $\phi:(M,F)\to(N,F)$ 
with $(N,F)\in C_{\mathcal D(2)fil}(S)$ such that $a_{\tau}j^*H^n\Gr_F^p(M,F)=0$ for all $n,p\in\mathbb Z$, 
\item $C_{\mathcal D(2)fil,Z,h}(S):=C_{\mathcal D(2)fil,Z}(S)\cap C_{\mathcal D(2)fil,h}(S)\subset C_{\mathcal D(2)fil}(S)$
the full subcategory consising of $(M,F)$ such that $a_{\tau}H^n(M,F)$ are filtered holonomic for all $n\in\mathbb Z$ and
such that there exists $r\in\mathbb N$ and an $r$-filtered homotopy equivalence $\phi:(M,F)\to(N,F)$ 
with $(N,F)\in C_{\mathcal D(2)fil}(S)$ such that $a_{\tau}j^*H^n\Gr_F^p(M,F)=0$ for all $n,p\in\mathbb Z$, 
\item $C_{\mathcal D(2)fil,Z,c}(S):=C_{\mathcal D(2)fil,Z}(S)\cap C_{\mathcal D(2)fil,c}(S)\subset C_{\mathcal D(2)fil}(S)$
the full subcategory consising of $(M,F)$ such that $a_{\tau}H^n(M,F)$ are filtered coherent for all $n\in\mathbb Z$ and
such that there exists $r\in\mathbb N$ and an $r$-filtered homotopy equivalence $\phi:(M,F)\to(N,F)$ 
with $(N,F)\in C_{\mathcal D(2)fil}(S)$ such that $a_{\tau}j^*H^n\Gr_F^p(M,F)=0$ for all $n,p\in\mathbb Z$. 
\end{itemize}
\item[(iii)] We have then the full subcategories 
\begin{itemize}
\item $C_{\mathcal D(1,0)fil,Z}(S)=C_{\mathcal D(1,0)fil}(S)\cap C_{\mathcal D2fil,Z}(S)\subset C_{\mathcal D2fil}(S)$,
\item $C_{\mathcal D(1,0)fil,Z,rh}(S)=C_{\mathcal D(1,0)fil}(S)\cap C_{\mathcal D2fil,Z,rh}(S)\subset C_{\mathcal D2fil}(S)$. 
\item $C_{\mathcal D(1,0)fil,Z,h}(S)=C_{\mathcal D(1,0)fil}(S)\cap C_{\mathcal D2fil,Z,h}(S)\subset C_{\mathcal D2fil}(S)$. 
\end{itemize}
\end{itemize}
\end{defi}

Similarly :

\begin{defi}\label{DmodsuppZinfty}
Let $S\in\AnSm(\mathbb C)$.
Let $Z\subset S$ a closed subset and denote by $j:S\backslash Z\hookrightarrow S$ the open embedding. 
\begin{itemize}
\item[(i)]We denote by
\begin{itemize}
\item $C_{\mathcal D^{\infty}(2)fil,Z}(S)\subset C_{\mathcal D^{\infty}(2)fil}(S)$. 
the full subcategory consisting of $(M,F)\in C_{\mathcal D^{\infty}}(S)$ such that $j^*M$ is acyclic 
\item  $C_{\mathcal D^{\infty}(2)fil,Z,h}(S):=C_{\mathcal D^{\infty}(2)fil,Z}(S)\cap C_{\mathcal D^{\infty}(2)fil,h}(S)
\subset C_{\mathcal D^{\infty}(2)fil}(S)$
the full subcategory consising of $(M,F)$ such that $a_{\tau}H^n(M)$ are holonomic and such that there exist $r\in\mathbb Z$ 
and an $r$-filtered homotopy equivalence $\phi:(M,F)\to(N,F)$ 
with $(N,F)\in C_{\mathcal D^{\infty}(2)fil}(S)$ such that $a_{\tau}j^*H^n\Gr_F^p(M,F)=0$. 
\end{itemize}
\item[(ii)] We have then the full subcategories
\begin{itemize}
\item $C_{\mathcal D^{\infty}(1,0)fil,Z}(S)=C_{\mathcal D^{\infty}(1,0)fil}(S)\cap C_{\mathcal D^{\infty}2fil,Z}(S)
\subset C_{\mathcal D^{\infty}2fil}(S)$,
\item $C_{\mathcal D^{\infty}(1,0)fil,Z,h}(S):=C_{\mathcal D^{\infty}(1,0)fil}(S)\cap C_{\mathcal D^{\infty}2fil,Z,h}(S)
\subset C_{\mathcal D^{\infty}(2)fil}(S)$.
\end{itemize}
\end{itemize}
\end{defi}

Let $S\in\SmVar(\mathbb C)$ or $S\in\AnSm(\mathbb C)$.
We recall (see section 2) that a morphism $m:(M,F)\to (N,F)$ with $(M,F),(N,F)\in C_{\mathcal Dfil}(S)$
is said to be an $r$-filtered zariski, resp. usu, local equivalence if there exist morphisms
$m_i:(C_i,F)\to(C_{i+1},F)$ with $(C_i,F),(C_{i+1},F)\in C_{\mathcal Dfil}(S)$, $0\leq i\leq s$,
$(C_0,F)=(M,F)$, $(C_s,F)=(N,F)$, such that
\begin{equation*}
m=m_s\circ\cdots\circ m_i\circ\cdots\circ m_0:(M,F)\to (N,F)
\end{equation*}
and $m_i:(C_i,F)\to (C_{i+1},F)$ is either a filtered zariski, resp. usu, local equivalence or an $r$-filtered homotopy.

\begin{defi}
\begin{itemize}
\item[(i)]Let $f:X\to S$ a morphism with $X,S\in\SmVar(\mathbb C)$, or with $X,S\in\AnSm(\mathbb C)$.
We have, for $r=1,\ldots\infty$, the homotopy categories 
\begin{equation*}
K_{\mathcal D,f^*\mathcal D(2)fil,r}(S):=\Ho_{r,0}C_{\mathcal D,f^*\mathcal D(2)fil}(S) \; , \; 
K_{\mathcal D^{op},f^*\mathcal D(2)fil,r}(S):=\Ho_{r,0}C_{\mathcal D^{op},f^*\mathcal D(2)fil}(S),
\end{equation*}
whose objects are those of $C_{\mathcal D,f^*\mathcal D(2)fil}(S)$, 
resp. those of $C_{\mathcal D^{op},f^*\mathcal D(2)fil}(S)$, and whose morphisms
are $r$-filtered for the first filtration (filtered for the second) homotopy classes of morphisms (see section 2.1). 
We have then
\begin{equation*}
D_{\mathcal D,f^*\mathcal D(2)fil,r}(S):=K_{\mathcal D,f^*\mathcal D(2)fil,r}(S)([E_1]^{-1}) \; , \; 
D_{\mathcal D^{op},f^*\mathcal D(2)fil,r}(S):=K_{\mathcal D^{op},f^*\mathcal D(2)fil,r}(S)([E_1]^{-1}),
\end{equation*}
the localizations with respect to the classes of filtered Zariski, resp. usu, local equivalences (see section 2). 
Note that the classes of filtered $\tau$ local equivalence constitute a right multiplicative system.
If $m:(M,F)\to (N,F)$ with $(M,F),(N,F)\in C_{\mathcal D,f^*\mathcal D(2)fil}(S)$
is an $r$-filtered zariski, resp. usu, local equivalence, 
then $m=D(top):(M,F)\to (N,F)$ is an isomorphism in $D_{\mathcal D,f^*\mathcal D(2)fil,r}(S)$.
\item[(ii)]Let $S\in\SmVar(\mathbb C)$ or $S\in\AnSm(\mathbb C)$. We denote by
\begin{equation*}
D_{\mathcal D(2)fil,\infty,rh}(S)\subset D_{\mathcal D(2)fil,\infty,h}(S)\subset D_{\mathcal D(2)fil,\infty}(S), 
\end{equation*}
the full subcategories  consisting of the image of 
$C_{\mathcal D(2)fil,h}(S)$, resp. $C_{\mathcal D(2)fil,rh}(S)$, by the localization functor
\begin{equation*}
D(top):C_{\mathcal D(2)fil}(S)\to D_{\mathcal D(2)fil,\infty}(S)
\end{equation*}
that is consisting of $(M,F)\in C_{\mathcal Dfil}(S)$ such that 
$a_{\tau}H^n(M,F)$ are filtered holonomic, resp. filtered regular holonomic for all $n\in\mathbb Z$,
\item[(iii)]Let $S\in\SmVar(\mathbb C)$ or $S\in\AnSm(\mathbb C)$. We denote by
\begin{equation*}
D_{\mathcal D(1,0)fil,\infty,rh}(S)\subset D_{\mathcal D(1,0)fil,\infty,h}(S)\subset D_{\mathcal D(2)fil,\infty}(S), 
\end{equation*}
the full subcategories  consisting of the image of 
$C_{\mathcal D(1,0)fil,h}(S)$, resp. $C_{\mathcal D(1,0)fil,rh}(S)$, by the localization functor
\begin{equation*}
D(top):C_{\mathcal D(2)fil}(S)\to D_{\mathcal D(2)fil,\infty}(S)
\end{equation*}
that is consisting of $(M,F,W)\in C_{\mathcal D2fil}(S)$ such that 
$a_{\tau}H^n(M,F)$ are filtered holonomic, resp. filtered regular holonomic,  
and $W^pM^n\subset M^n$ are $D_S$ submodules for all $n\in\mathbb Z$.
\end{itemize}
\end{defi}

Let $S\in\SmVar(\mathbb C)$ or $S\in\AnSm(\mathbb C)$. By definition (see section 2), we have sequences of functors
\begin{eqnarray*}
C_{\mathcal D(2)fil}(\mathcal S)\to K_{\mathcal D(2)fil}(S)\to 
D_{\mathcal D(2)fil}(\mathcal S)\to D_{\mathcal D(2)fil,2}(S)\to\cdots\to D_{\mathcal D(2)fil,\infty}(S).
\end{eqnarray*}
and commutative diagrams of functors
\begin{equation*}
\xymatrix{K_{\mathcal D(2)fil}(S)\ar[r]\ar[d] & D_{fil}(S)\ar[d] \\
K_{\mathcal D(2)fil,2}(\mathcal S)\ar[r] & D_{\mathcal D(2)fil,2}(S)} \; , \;
\xymatrix{K_{\mathcal D(2)fil,r}(S)\ar[r]\ar[d] & D_{\mathcal D(2)fil,r}(S)\ar[d] \\
K_{\mathcal D(2)fil,r+1}(S)\ar[r] & D_{\mathcal D(2)fil,r+1}(S)}.
\end{equation*}
Then, for $r=1$, $K_{\mathcal D(2)fil}(S)$ and $D_{\mathcal D(2)fil}(S)$ are in the canonical way triangulated categories.
However, for $r>1$, the categories $K_{\mathcal D(2)fil,r}(S)$ and $D_{\mathcal D(2)fil,r}(S)$
together with the canonical triangles does NOT satisfy the 2 of 3 axiom of triangulated categories.

Similarly,
\begin{defi}
\begin{itemize}
\item[(i)] Let $f:X\to S$ a morphism with $X,S\in\AnSm(\mathbb C)$.
We have, for $r=1,\ldots\infty$, the categories 
\begin{equation*}
K_{\mathcal D^{\infty},f^*\mathcal D^{\infty}(2)fil,r}(S):=\Ho_{r,0}C_{\mathcal D^{\infty},f^*\mathcal D^{\infty}(2)fil}(S) \; , \;
K_{\mathcal D^{\infty,op},f^*\mathcal D^{\infty}(2)fil,r}(S):=\Ho_{r,0}C_{\mathcal D^{\infty,op},f^*\mathcal D^{\infty}(2)fil}(S),
\end{equation*}
whose objects are those of $C_{\mathcal D^{\infty},f^*\mathcal D^{\infty}(2)fil}(S)$, 
resp. those of $C_{\mathcal D^{\infty,op},f^*\mathcal D^{\infty}(2)fil}(S)$, and whose morphisms
are $r$-filtered for the first filtration (filtered for the second) homotopy classes of morphisms (see section 2.1). 
We have the categories 
\begin{equation*}
D_{\mathcal D^{\infty},f^*\mathcal D^{\infty}(2)fil,r}(S):=
K_{\mathcal D^{\infty},f^*\mathcal D^{\infty}(2)fil,r}(S)([E_1]^{-1}) \; , \;
D_{\mathcal D^{\infty,op},f^*\mathcal D^{\infty}(2)fil,r}(S):=
K_{\mathcal D^{\infty,op},f^*\mathcal D^{\infty}(2)fil,r}(S)([E_1]^{-1}),
\end{equation*}
the localizations with respect to the classes of filtered usu local equivalence (see section 2.1). 
Note that the classes of filtered usu local equivalence constitute a right multiplicative system.
\item[(ii)]Let $S\in\AnSm(\mathbb C)$. We denote by
\begin{equation*}
D_{\mathcal D^{\infty}(2)fil,\infty,h}(S)\subset D_{\mathcal D^{\infty}(2)fil,\infty}(S), \;  
D_{\mathcal D^{\infty}(1,0)fil,\infty,h}(S)\subset D_{\mathcal D^{\infty}2fil,\infty}(S)
\end{equation*}
the full subcategories 
consisting of the image of $C_{\mathcal D^{\infty}(2)fil,h}(S)$, resp. $C_{\mathcal D^{\infty}(1,0)fil,h}(S)$, by the localization functor
\begin{equation*}
D(top):C_{\mathcal D^{\infty}(2)fil}(S)\to D_{\mathcal D^{\infty}(2)fil,\infty}(S).
\end{equation*}
\end{itemize}
\end{defi}

We begin this subsection by recalling the following well known facts
\begin{prop}
Let $S\in\SmVar(\mathbb C)$ or $S\in\AnSm(\mathbb C)$.
\begin{itemize}
\item[(i)] The sheaf of differential operators $D_S$ is a locally free sheaf of $O_S$ module. Hence,
a coherent $D_S$ module $M\in\Coh_{\mathcal D}(S)$ is a quasi-coherent sheaf of $O_S$ modules.  
\item[(ii)] A coherent sheaf $M\in\Coh_{O_S}(S)$ of $O_S$ module admits a $D_S$ module structure if and only if
it is locally free (of finite rank by coherency) and admits an integrable connexion. In particular
if $i:Z\hookrightarrow S$ is a closed embedding for the Zariski topology, then $i_*O_Z$ does NOT admit
a $D_S$ module structure since it is a coherent but not locally free $O_S$ module.
\end{itemize}
\end{prop}

\begin{proof}
Standard.
\end{proof}

In order to prove a version of the first GAGA theorem for coherent D modules, we will need to following.
We start by a definition (cf. \cite{LvDmod} definition 1.4.2) :
\begin{defi}\label{Daffine}
An $X\in\SmVar(\mathbb C)$ is said to be D-affine if the following two condition hold:
\begin{itemize}
\item[(i)] The global section functor $\Gamma(X,\cdot):\mathcal QCoh_{\mathcal D}(X)\to Mod(\Gamma(X,D_X))$ is exact.
\item[(ii)] If $\Gamma(X,M)=0$ for $M\in\mathcal QCoh_{\mathcal D}(X)$, then $M=0$.
\end{itemize} 
\end{defi}

\begin{prop}\label{Daffineprop}
If $X\in\SmVar(\mathbb C)$ is D-affine, then :
\begin{itemize}
\item[(i)] Any $M\in\mathcal QCoh_{\mathcal D}(X)$ is generated by its global sections.
\item[(ii)] The functor $\Gamma(X,\cdot):\mathcal QCoh_{\mathcal D}(X)\to Mod(\Gamma(X,D_X))$ is an equivalence
of category whose inverse is $L\in Mod(\Gamma(X,D_X))\mapsto D_X\otimes_{\Gamma(X,D_X)}L\in\mathcal QCoh_{\mathcal D}(X)$.
\item[(iii)] We have $\Gamma(X,\cdot)(\mathcal Coh_{\mathcal D}(X))=Mod(\Gamma(X,D_X))_f$, that is the global sections of a
coherent $D_X$ module is a finite module over the differential operators on $X$.
\end{itemize} 
\end{prop}

\begin{proof}
See \cite{LvDmod}.
\end{proof}

The following proposition is from Kashiwara.
\begin{prop}\label{Jprop}
Let $S\in\AnSm(\mathbb C)$.
\begin{itemize}
\item[(i)] For $K\in C_c(S)$ a complex of presheaves with constructible cohomology sheaves, 
we have $\mathcal Hom(L(K),E(O_S))\in C_{\mathcal D^{\infty},h}(S)$.
\item[(ii)] The functor $J_S:C_{\mathcal D(2)fil}(S)\to C_{\mathcal D^{\infty}(2)fil}(S)$ 
satisfy $J_S(C_{\mathcal D(2)fil,rh}(S))\subset C_{\mathcal D^{\infty}(2)fil,h}(S)$, derive trivially, and induce an equivalence of category
\begin{equation*}
J_S:D_{\mathcal D(2)fil,\infty,rh}(S)\to D_{\mathcal D^{\infty}(2)fil,\infty,h}(S).
\end{equation*}
whose inverse satify, for $(M,F)\in\mathcal Hol_{\mathcal D^{\infty}(2)fil}(S)$ a (filtered) holonomic $D^{\infty}_S$ module, 
that $J_S^{-1}(M,F)=(M_{reg},F)\subset (M,F)$ is the $D_S$ sub-module of $M$ which is the regular part.
\item[(iii)] We have $J_S(C_{\mathcal D(1,0)fil,rh}(S))\subset C_{\mathcal D^{\infty}(1,0)fil,h}(S)$ and 
$J_S(D_{\mathcal D(1,0)fil,\infty,rh}(S))=D_{\mathcal D^{\infty}(1,0)fil,\infty,h}(S)$.
\end{itemize}
\end{prop}

\begin{proof}
Follows from \cite{Kashiwara}.
\end{proof}

Let $S\in\SmVar(\mathbb C)$ or $S\in\AnSm(\mathbb C)$, and 
let $i:Z\hookrightarrow S$ a closed embedding and denote by $j:S\backslash Z\hookrightarrow S$ the open complementary. 
For $M\in\PSh_{\mathcal D}(S)$, we denote $\mathcal I_ZM\subset M$ the (left) $D_S$ submodule
given by, for $S^o\subset S$ an open subset, $\mathcal I_ZM(S^o)\subset M(S^o)$ is the (left) $D_S(S^o)$ submodule  
\begin{equation*}
\mathcal I_ZM(S^o)=<\left\{fm,f\in\mathcal I_Z(S^o),m\in M(S^o)\right\}>\subset M(S^o)
\end{equation*}
generated by the elements of the form $fm$. 
We denote by $b_Z(M):\mathcal I_ZM\to M$ the inclusion map and $c_Z(M):M\to M/\mathcal I_ZM$ the quotient map
of (left) $D_S$ modules.
For $M\in\PSh_{\mathcal D}(S)$, we denote $M\mathcal I_Z\subset M$ the right $D_S$ submodule
given by, for $S^o\subset S$ an open subset, $\mathcal I_ZM(S^o)\subset M(S^o)$ is the right $D_S(S^o)$ submodule  
\begin{equation*}
\mathcal I_ZM(S^o)=<\left\{mf,f\in\mathcal I_Z(S^o),m\in M(S^o)\right\}>\subset M(S^o)
\end{equation*}
generated by the elements of the form $mf$. 
We denote by $b_Z(M):\mathcal I_ZM\to M$ the inclusion map and $c_Z(M):M\to M/\mathcal I_ZM$ the quotient map
of right $D_S$ modules.

\subsubsection{Functorialities}

Let $f:X\to S$ be a morphism with $X,S\in\SmVar(\mathbb C)$,
or let $f:X\to S$ be a morphism with $X,S\in\AnSm(\mathbb C)$.
Then, we recall from section 4.1, the transfers modules 
\begin{itemize}
\item $(D_{X\to S},F^{ord}):=f^{*mod}(D_S,F^{ord}):=f^*(D_S,F^{ord})\otimes_{f^*O_S}(O_X,F_b)$ 
which is a left $D_X$ module and a left and right $f^*D_S$ module
\item $(D_{X\leftarrow S},F^{ord}):=(K_X,F_b)\otimes_{O_X}(D_{X\to S},F^{ord})\otimes_{f^*O_S}f^*(K_S,F_b)$.
which is a right $D_X$ module and a left and right $f^*D_S$ module.
\end{itemize}

Let $f:X\to S$ be a morphism with $X,S\in\AnSm(\mathbb C)$.
Then, the transfers modules of infite order are
\begin{itemize}
\item $(D^{\infty}_{X\to S},F^{ord}):=f^{*mod}(D^{\infty}_S,F^{ord}):=f^*(D^{\infty}_S,F^{ord})\otimes_{f^*O_S}(O_X,F_b)$ 
which is a left $D^{\infty}_X$ module and a left and right $f^*D^{\infty}_S$ module
\item $(D^{\infty}_{X\leftarrow S},F^{ord}):=(K_X,F_b){\otimes}_{O_X}(D^{\infty}_{X\to S},F^{ord})\otimes_{f^*O_S}f^*(K_S,F_b)$.
which is a right $D^{\infty}_X$ module and a left and right $f^*D^{\infty}_S$ module.
\end{itemize}

We have the following :

\begin{lem}\label{Dcomp}
Let $f_1:X\to Y$, $f_2:Y\to S$ be two morphism with $X,S,Y\in\SmVar(\mathbb C)$,
or let $f_1:X\to Y$, $f_2:Y\to S$ be two morphism with $X,S,Y\in\AnSm(\mathbb C)$.
\begin{itemize}
\item[(i)] We have  
$(D_{X\to S},F^{ord})=f_1^*(D_{Y\to S},F^{ord})\otimes_{f_1^*D_Y}(D_{X\to Y},F^{ord})$
in $C_{\mathcal D,(f_2\circ f_1)^*\mathcal Dfil}(X)$ and  
\begin{equation*}
(D_{X\to S},F^{ord})=f_1^*(D_{Y\to S},F^{ord})\otimes_{f_1^*D_Y}(D_{X\to Y},F^{ord})=
f_1^*(D_{Y\to S},F^{ord})\otimes^L_{f_1^*D_Y}(D_{X\to Y},F^{ord}). 
\end{equation*}
in $D_{\mathcal D,(f_2\circ f_1)^*\mathcal Dfil,r}(X)$.
\item[(ii)] We have  
$(D_{X\leftarrow S},F^{ord})=f_1^*(D_{Y\to S},F^{ord})\otimes_{f_1^*D_Y}(D_{X\leftarrow Y},F^{ord})$
in $C_{\mathcal D^{op},(f_2\circ f_1)^*\mathcal Dfil}(X)$ and 
\begin{equation*}
(D_{X\leftarrow S},F^{ord})=f_1^*(D_{Y\to S},F^{ord})\otimes_{f_1^*D_Y}(D_{X\leftarrow Y},F^{ord})=
f_1^*(D_{Y\leftarrow S},F^{ord})\otimes^L_{f_1^*D_Y}(D_{X\leftarrow Y},F^{ord}), 
\end{equation*}
in $D_{\mathcal D^{op},(f_2\circ f_1)^*\mathcal Dfil,r}(X)$. 
\end{itemize}
\end{lem}

\begin{proof}
Follows immediately from definition. The first assertions of (i) and (ii) are particular cases of lemma \ref{DcompRTop}.
See \cite{LvDmod} for example.
\end{proof}

In the analytical case we also have

\begin{lem}\label{Dcompinfty}
Let $f_1:X\to Y$, $f_2:Y\to S$ be two morphism with $X,S,Y\in\AnSm(\mathbb C)$.
\begin{itemize}
\item[(i)] We have  
$(D^{\infty}_{X\to S},F^{ord})=f_1^*(D^{\infty}_{Y\to S},F^{ord})\otimes_{f_1^*D^{\infty}_Y}(D^{\infty}_{X\to Y},F^{ord})$
in $C_{\mathcal D^{\infty},(f_2\circ f_1)^*\mathcal D^{\infty}fil}(X)$ and  
\begin{equation*}
(D^{\infty}_{X\to S},F^{ord})=f_1^*(D^{\infty}_{Y\to S},F^{ord})\otimes_{f_1^*D^{\infty}_Y}(D_{X\to Y},F^{ord})=
f_1^*(D^{\infty}_{Y\to S},F^{ord})\otimes^L_{f_1^*D^{\infty}_Y}(D^{\infty}_{X\to Y},F^{ord}). 
\end{equation*}
in $D_{\mathcal D^{\infty},(f_2\circ f_1)^*\mathcal D^{\infty}fil,r}(X)$.
\item[(ii)] We have  
$(D^{\infty}_{X\leftarrow S},F^{ord})=f_1^*(D^{\infty}_{Y\to S},F^{ord})\otimes_{f_1^*D^{\infty}_Y}(D^{\infty}_{X\leftarrow Y},F^{ord})$
in $C_{\mathcal D^{\infty,op},(f_2\circ f_1)^*\mathcal D^{\infty}fil}(X)$ and 
\begin{equation*}
(D^{\infty}_{X\leftarrow S},F^{ord})=f_1^*(D^{\infty}_{Y\to S},F^{ord})\otimes_{f_1^*D^{\infty}_Y}(D^{\infty}_{X\leftarrow Y},F^{ord})=
f_1^*(D^{\infty}_{Y\leftarrow S},F^{ord})\otimes^L_{f_1^*D^{\infty}_Y}(D^{\infty}_{X\leftarrow Y},F^{ord}), 
\end{equation*}
in $D_{\mathcal D^{\infty,op},(f_2\circ f_1)^*\mathcal D^{\infty}fil,r}(X)$. 
\end{itemize}
\end{lem}

\begin{proof}
Similar to the proof of lemma \ref{Dcomp}
\end{proof}

For closed embeddings, we have :

\begin{prop}
\begin{itemize}
\item[(i)] Let $i:Z\hookrightarrow S$ be a closed embedding with $Z,S\in\SmVar(\mathbb C)$.
Then, $D_{Z\to S}=i^*D_S/D_S\mathcal I_Z$ and it is a locally free (left) $D_Z$ module.
Similarly, $D_{Z\leftarrow S}=i^*D_S/\mathcal I_ZD_S$ and it is a locally free right $D_Z$ module.
\item[(ii)] Let $i:Z\to S$ be a closed embedding with $Z,S\in\AnSm(\mathbb C)$.
Then, $D_{Z\to S}=i^*D_S/D_S\mathcal I_Z$ and it is a locally free (left) $D_Z$ module.
Similarly, $D_{Z\leftarrow S}=i^*D_S/\mathcal I_ZD_S$ and it is a locally free right $D_Z$ module.
\item[(iii)] Let $i:Z\to S$ be a closed embedding with $Z,S\in\AnSm(\mathbb C)$.
Then, $D^{\infty}_{Z\to S}=i^*D^{\infty}_S/D^{\infty}_S\mathcal I_Z$ and it is a locally free (left) $D^{\infty}_Z$ module.
Similarly, $D^{\infty}_{Z\leftarrow S}=i^*D^{\infty}_S/\mathcal I_ZD^{\infty}_S$ and it is a locally free right $D^{\infty}_Z$ module.
\end{itemize}
\end{prop}

\begin{proof}
\noindent(i): See \cite{LvDmod}.

\noindent(ii):See \cite{PhamDmod}.

\noindent(iii):Similar to (ii).
\end{proof}

We now enumerate some functorialities we will use, 
all of them are particular case of the functoriality given in subsection 2.3 for any ringed spaces :

\begin{itemize}

\item Let $f:X\to S$ be a morphism with $X,S\in\Var(\mathbb C)$,
or let $f:X\to S$ be a morphism with $X,S\in\AnSp(\mathbb C)$.
Then, the inverse image functor
\begin{equation*}
f^{*mod}:\PSh_{O_S}(S)\to\PSh_{O_X}(X), \; \; M\mapsto f^{*mod}M:=O_X\otimes_{f^*O_S}f^*M 
\end{equation*}
is a Quillen adjonction which induces in the derived category the functor
\begin{equation*}
Lf^{*mod}:D_{O_S}(S)\to D_{O_X}(X), \; \; M\mapsto Lf^{*mod}M:=O_X\otimes^L_{f^*O_S}f^*M=O_X\otimes_{f^*O_S}f^*L_OM,
\end{equation*}
The adjonction $(f^{*mod},f_*):\PSh_{O_S}(S)\leftrightarrows\PSh_{O_X}(X)$ is a Quillen adjonction, the adjonction map are the maps
\begin{itemize}
\item for $M\in C_{O_S}(S)$,
$\ad(f^{*mod},f_*)(M):M\xrightarrow{\ad(f^*,f_*)(M)}f_*f^*M\xrightarrow{f_*m}f_*(f^*M\otimes_{f^*O_S}O_X)=f_*f^{*mod}M$
where $m(m)=m\otimes 1$,
\item for $M\in C_{O_X}(X)$,
$\ad(f^{*mod},f_*)(M):f^{*mod}f_*M=f^*f_*M\otimes_{f^*O_S}O_X\xrightarrow{\ad(f^*,f_*)(M)\otimes_{f^*O_S}O_X}M\otimes_{f^*O_S}O_X
\xrightarrow{n}M$, where $n(m\otimes h)=h.m$ is the multiplication map.
\end{itemize}

\item Let $S\in\SmVar(\mathbb C)$ or $S\in\AnSm(\mathbb C)$. 
\begin{itemize}
\item For $M\in C_{\mathcal D}(S)$, we have the canonical projective resolution $q:L_D(M)\to M$ of complexes
of $D_S$ modules.
\item For $M\in C_{\mathcal D}(S)$, there exist a unique strucure of $D_S$ module on the flasque presheaves $E^i(M)$ such that
$E(M)\in C_{\mathcal D}(S)$ (i.e. is a complex of $D_S$ modules) 
and that the map $k:M\to E(M)$ is a morphism of complexes of $D_S$ modules.
\end{itemize}
Let $S\in\AnSm(\mathbb C)$.
\begin{itemize}
\item For $M\in C_{\mathcal D^{\infty}}(S)$, we have the canonical projective resolution $q:L_{D^{\infty}}(M)\to M$ of complexes
of $D^{\infty}_S$ modules.
\item For $M\in C_{\mathcal D^{\infty}}(S)$, 
there exist a unique strucure of $D^{\infty}_S$ module on the flasque presheaves $E^i(M)$ such that
$E(M)\in C_{\mathcal D^{\infty}}(S)$ (i.e. is a complex of $D^{\infty}_S$ modules) and that the map $k:M\to E(M)$ 
is a morphism of complexes of $D^{\infty}_S$ modules.
\end{itemize}

\item Let $S\in\SmVar(\mathbb C)$ or let $S\in\AnSm(\mathbb C)$. 
For $M\in C_{\mathcal D^{(op)}}(S)$, $N\in C(S)$, we will consider the induced D module structure
(right $D_S$ module in the case one is a left $D_S$ module and the other one is a right one)
on the presheaf $M\otimes N:=M\otimes_{\mathbb Z_S}N$ (see section 2). We get the bifunctor
\begin{eqnarray*}
C(S)\times C_{\mathcal D}(S)\to C_{\mathcal D}(S), (M,N)\mapsto M\otimes N
\end{eqnarray*}
For $S\in\AnSm(\mathbb C)$, we also have the bifunctor
$C(S)\times C_{\mathcal D^{\infty}}(S)\to C_{\mathcal D^{\infty}}(S), \; (M,N)\mapsto M\otimes N$.

\item Let $S\in\SmVar(\mathbb C)$ or let $S\in\AnSm(\mathbb C)$. 
For $M,N\in C_{\mathcal D^{(op)}}(S)$, $M\otimes_{O_S} N$ (see section 2),
has a canonical structure of $D_S$ modules 
(right $D_S$ module in the case one is a left $D_S$ module and the other one is a right one)
given by (in the left case) for $S^o\subset S$ an open subset, 
\begin{equation*}
m\otimes n\in\Gamma(S^o,M\otimes_{O_S}N),\gamma\in\Gamma(S^o,D_S), \;
\gamma.(m\otimes n):=(\gamma.m)\otimes n-m\otimes\gamma.n
\end{equation*}
This gives the bifunctor
\begin{eqnarray*}
C_{\mathcal D^{(op)}}(S)^2\to C_{\mathcal D^{(op)}}(S), (M,N)\mapsto M\otimes_{O_S}N
\end{eqnarray*}
More generally, let $f:X\to S$ a morphism with $X,S\in\Var(\mathbb C)$ or with $X,S\in\AnSp(\mathbb C)$. Assume $S$ smooth. 
For $M,N\in C_{f^*\mathcal D^{(op)}}(X)$, $M\otimes_{f^*O_S} N$ (see section 2),
has a canonical structure of $f^*D_S$ modules 
(right $f^*D_S$ module in the case one is a left $f^*D_S$ module and the other one is a right one)
given by (in the left case) for $X^o\subset X$ an open subset, 
\begin{equation*}
m\otimes n\in\Gamma(X^o,M\otimes_{f^*O_S} N),\gamma\in\Gamma(X^o,f^*D_S), \;
\gamma.(m\otimes n):=(\gamma.m)\otimes n-m\otimes\gamma.n
\end{equation*}
This gives the bifunctor
\begin{eqnarray*}
C_{f^*\mathcal D^{(op)}}(X)^2\to C_{f^*\mathcal D^{(op)}}(X), (M,N)\mapsto M\otimes_{f^*O_S}N
\end{eqnarray*}
For $f:X\to S$ a morphism with $X,S\in\AnSp(\mathbb C)$ and $S$ smooth, we also have the bifunctor
$C_{f^*\mathcal D^{\infty,(op)}}(X)^2\to C_{f^*\mathcal D^{\infty,(op)}}(X), (M,N)\mapsto M\otimes_{f^*O_S}N$.

\item Let $S\in\SmVar(\mathbb C)$ or let $S\in\AnSm(\mathbb C)$. 
For $M\in C_{\mathcal D^{op}}(S)$ and $N\in C_{\mathcal D}(S)$, we have $M\otimes_{D_S} N\in C(S)$ (see section 2). 
This gives the bifunctor
\begin{eqnarray*}
C_{\mathcal D^{op}}(S)\times C_{\mathcal D}(S)\to C(S), (M,N)\mapsto M\otimes_{D_S}N
\end{eqnarray*}
For $S\in\AnSm(\mathbb C)$, we also have the bifunctor
$C_{\mathcal D^{\infty,op}}(S)\times C_{\mathcal D^{\infty}}(S)\to C(S), (M,N)\mapsto M\otimes_{D^{\infty}_S}N$.

\item Let $S\in\SmVar(\mathbb C)$ or let $S\in\AnSm(\mathbb C)$. 
The internal hom bifunctor 
\begin{equation*}
\mathcal Hom(\cdot,\cdot):=\mathcal Hom_{\mathbb Z_S}(\cdot,\cdot):C(S)^2\to C(S)
\end{equation*}
induces a bifunctor
\begin{equation*}
\mathcal Hom(\cdot,\cdot):=\mathcal Hom_{\mathbb Z_S}(\cdot,\cdot):C(S)\times C_{\mathcal D}(S)\to C_{\mathcal D}(S)
\end{equation*}
such that, for $F\in C(S)$ and $G\in C_{\mathcal D}(S)$, the $D_S$ structure on $\mathcal Hom^{\bullet}(F,G)$ is given by
\begin{equation*}
\gamma\in\Gamma(S^o,D_S)\longmapsto 
(\phi\in\Hom^{p}(F^{\bullet}_{|S^o},G_{|S^o})\mapsto
(\gamma\cdot\phi:\alpha\in F^{\bullet}(S^o)\mapsto\gamma\cdot\phi^p(S^o)(\alpha))
\end{equation*}
where $\phi^p(S^o)(\alpha)\in\Gamma(S^o,G)$. For $S\in\AnSm(\mathbb C)$, it also induce the bifunctor
\begin{equation*}
\mathcal Hom(\cdot,\cdot):=\mathcal Hom_{\mathbb Z_S}(\cdot,\cdot):
C(S)\times C_{\mathcal D^{\infty}}(S)\to C_{\mathcal D^{\infty}}(S)
\end{equation*}

\item Let $S\in\SmVar(\mathbb C)$ or let $S\in\AnSm(\mathbb C)$. 
For $M,N\in C_{\mathcal D}(S)$, $\mathcal Hom_{O_S}(M,N)$,
has a canonical structure of $D_S$ modules given by for $S^o\subset S$ an open subset and
$\phi\in\Gamma(S^o,\mathcal Hom(M,O_S))$, $\gamma\in\Gamma(S^o,D_S)$,
$(\gamma.\phi)(m):=\gamma.(\phi(m))-\phi(\gamma.m)$
This gives the bifunctor
\begin{eqnarray*}
\Hom^{\bullet}_{O_S}(-,-):C_{\mathcal D}(S)^2\to C_{\mathcal D}(S)^{op}, (M,N)\mapsto\mathcal Hom^{\bullet}_{O_S}(M,N)
\end{eqnarray*}
In particular, for $M\in C_{\mathcal D}(S)$, we get the dual 
\begin{equation*}
\mathbb D^{O}_S(M):=\mathcal Hom^{\bullet}_{O_S}(M,O_S)\in C_{\mathcal D}(S)
\end{equation*}
with respect to $O_S$, together with the canonical map $d(M):M\to\mathbb D^{O,2}_S(M)$.
Let $f:X\to S$ a morphism with $X,S\in\SmVar(\mathbb C)$ or with $X,S\in\AnSm(\mathbb C)$.
We have, for $M\in C_{\mathcal D}(S)$, the canonical transformation map
\begin{eqnarray*}
T(f,D^o)(M):f^{*mod}\mathbb D_S^OM=(f^*\mathcal Hom_{O_S}(M,O_S))\otimes_{f^*O_S}O_X \\
\xrightarrow{T^{mod}(f,hom)(M,O_S)}\mathcal Hom_{O_X}(f^*M\otimes_{f^*O_S}O_X,O_X)=:\mathbb D_X^O(f^{*mod}M).
\end{eqnarray*}

\item Let $S\in\SmVar(\mathbb C)$ or let $S\in\AnSm(\mathbb C)$. 
We have the bifunctors
\begin{itemize}
\item $\Hom^{\bullet}_{D_S}(-,-):C_{\mathcal D}(S)^2\to C(S)$, $(M,N)\mapsto\mathcal Hom^{\bullet}_{D_S}(M,N)$,
and if $N$ is a bimodule (i.e. has a right $D_S$ module structure whose opposite coincide with the left one),
$\mathcal Hom_{D_S}(M,N)\in C_{\mathcal D^{op}}(S)$ given by for $S^o\subset S$ an open subset and
$\phi\in\Gamma(S^o,\mathcal Hom(M,N))$, $\gamma\in\Gamma(S^o,D_S)$, $(\phi.\gamma)(m):=(\phi(m)).\gamma$
\item $\Hom_{D_S}(-,-):C_{\mathcal D^{op}}(S)^2\to C(S)$, $(M,N)\mapsto\mathcal Hom_{D_S}(M,N)$
and if $N$ is a bimodule, $\mathcal Hom_{D_S}(M,N)\in C_{\mathcal D}(S)$
\end{itemize}
For $M\in C_{\mathcal D}(S)$, we get in particular the dual with respect $\mathbb D_S$,
\begin{equation*}
\mathbb D_SM:=\mathcal Hom_{D_S}(M,D_S)\in C_{\mathcal D}(S) \; ; \;
\mathbb D^K_SM:=\mathcal Hom_{D_S}(M,D_S)\otimes_{O_S}\mathbb D^O_Sw(K_S)[d_S]\in C_{\mathcal D}(S)
\end{equation*}
and we have canonical map $d:M\to\mathbb D_S^2M$.
This functor induces in the derived category, for $M\in D_{\mathcal D}(S)$,
\begin{equation*}
L\mathbb D_SM:=R\mathcal Hom_{D_S}(L_DM,D_S)\otimes_{O_S}\mathbb D^O_Sw(K_S)[d_S]=\mathbb D^K_SL_DM\in D_{\mathcal D}(S).
\end{equation*}
where $\mathbb D^O_Sw(S):\mathbb D_S^Ow(K_S)\to\mathbb D^O_SK_S=K_S^{-1}$ 
is the dual of the Koczul resolution of the canonical bundle (proposition \ref{resw}), 
and the canonical map $d:M\to L\mathbb D_S^2M$.
For $S\in\AnSm(\mathbb C)$, we also have the bifunctors
\begin{itemize}
\item $\Hom^{\bullet}_{D^{\infty}_S}(-,-):C_{\mathcal D^{\infty}}(S)^2\to C(S)$, 
$(M,N)\mapsto\mathcal Hom^{\bullet}_{D^{\infty}_S}(M,N)$,
and if $N$ is a bimodule, $\mathcal Hom_{D^{\infty}_S}(M,N)\in C_{\mathcal D^{\infty}}(S)$,
\item $\Hom_{D_S}(-,-):C_{\mathcal D^{\infty,op}}(S)^2\to C(S)$, $(M,N)\mapsto\mathcal Hom_{D^{\infty}_S}(M,N)$
and if $N$ is a bimodule, $\mathcal Hom_{D^{\infty}_S}(M,N)\in C_{\mathcal D^{\infty,op}}(S)$
\end{itemize}
For $M\in C_{\mathcal D}(S)$, we get in particular the dual with respect $\mathbb D^{\infty}_S$,
\begin{equation*}
\mathbb D^{\infty}_SM:=\mathcal Hom_{D^{\infty}_S}(M,D^{\infty}_S)\in C_{\mathcal D^{\infty}}(S) \; , \;
\mathbb D^{\infty,K}_SM:=\mathcal Hom_{D^{\infty}_S}(M,D^{\infty}_S)\otimes_{O_S}\mathbb D^O_Sw(K_S)[d_S]
\in C_{\mathcal D^{\infty}}(S) 
\end{equation*}
and we have canonical maps $d:M\to\mathbb D^{\infty,2}_SM$, $d:M\to\mathbb D^{\infty,K,2}_SM$.
This functor induces in the derived category, for $M\in D_{\mathcal D^{\infty}}(S)$,
\begin{equation*}
L\mathbb D^{\infty}_SM:=R\mathcal Hom_{D^{\infty}_S}(M,D^{\infty}_S)\otimes_{O_S}\mathbb D_S^Ow(K_S)[d_S]
=\mathbb D^{\infty,K}_SL_{D^{\infty}}M\in D_{\mathcal D^{\infty}}(S).
\end{equation*}
and the canonical map $d:M\to L\mathbb D^{\infty,2}_SM$.

\item Let $f:X\to S$ a morphism with $X,S\in\SmVar(\mathbb C)$ or with $X,S\in\AnSm(\mathbb C)$. 
For $N\in C_{\mathcal D,f^*\mathcal D}(X)$ and $M\in C_{\mathcal D}(X)$, $N\otimes_{D_X}M$ 
has the canonical $f^*D_S$ module structure given by, for $X^o\subset X$ an open subset, 
\begin{equation*}
\gamma\in\Gamma(X^o,f^*D_S), m\in\Gamma(X^o,M),n\in\Gamma(X^o,N), \; \gamma.(n\otimes m)=(\gamma.n)\otimes m.
\end{equation*}
This gives the functor
\begin{equation*}
C_{\mathcal D,f^*\mathcal D}(X)\times C_{\mathcal D}(X)\to C_{f^*\mathcal D}(X), \; (M,N)\mapsto M\otimes_{D_X}N
\end{equation*}

\item Let $f:X\to S$ a morphism with $X,S\in\AnSm(\mathbb C)$. 
For $N\in C_{\mathcal D^{\infty},f^*\mathcal D^{\infty}}(X)$ and $M\in C_{\mathcal D^{\infty}}(X)$, $N\otimes_{D^{\infty}_X}M$ 
has the canonical $f^*D^{\infty}_S$ module structure given by, for $X^o\subset X$ an open subset, 
\begin{equation*}
\gamma\in\Gamma(X^o,f^*D_S), m\in\Gamma(X^o,M),n\in\Gamma(X^o,N), \; \gamma.(n\otimes m)=(\gamma.n)\otimes m.
\end{equation*}
This gives the functor
\begin{equation*}
C_{\mathcal D^{\infty},f^*\mathcal D^{\infty}}(X)\times C_{\mathcal D^{\infty}}(X)\to C_{f^*\mathcal D^{\infty}}(X), \; 
(M,N)\mapsto M\otimes_{D^{\infty}_X}N
\end{equation*}

\item Let $f:X\to S$ be a morphism with $X,S\in\SmVar(\mathbb C)$,
or let $f:X\to S$ be a morphism with $X,S\in\AnSm(\mathbb C)$.
Then, for $M\in C_{\mathcal D}(S)$, $O_X\otimes_{f^*O_S}f^*M$ has a canonical $D_X$ module structure given by
given by, for $X^o\subset X$ an open subset,
\begin{equation*}
m\otimes n\in\Gamma(X^o,O_X\otimes_{f^*O_S}f^*M), \gamma\in\Gamma(X^o,D_X), \; 
\gamma.(m\otimes n):=(\gamma.m)\otimes n-m\otimes df(\gamma).n.
\end{equation*}
This gives the inverse image functor
\begin{equation*}
f^{*mod}:\PSh_{\mathcal D}(S)\to\PSh_{\mathcal D}(X), \; \; 
M\mapsto f^{*mod}M:=O_X\otimes_{f^*O_S}f^*M=D_{X\to S}\otimes_{f^*D_S}f^*M 
\end{equation*}
which induces in the derived category the functor
\begin{equation*}
Lf^{*mod}:D_{\mathcal D}(S)\to D_{\mathcal D}(X), \; \;
M\mapsto Lf^{*mod}M:=O_X\otimes^L_{f^*O_S}f^*M=O_X\otimes_{f^*O_S}f^*L_DM,
\end{equation*}
We will also consider the shifted inverse image functor
\begin{equation*}
Lf^{*mod[-]}:=Lf^{*mod}[d_S-d_X]:D_{\mathcal D}(S)\to D_{\mathcal D}(X).
\end{equation*}

\item Let $f:X\to S$ be a morphism with $X,S\in\AnSm(\mathbb C)$.
Then, for $M\in C_{\mathcal D^{\infty}}(S)$, 
$O_X\otimes_{f^*O_S}f^*M$ has a canonical $D^{\infty}_X$ module structure induced by the finite order case.
This gives the inverse image functor
\begin{equation*}
f^{*mod}:\PSh_{\mathcal D^{\infty}}(S)\to\PSh_{\mathcal D^{\infty}}(X), \; \; 
M\mapsto f^{*mod}M:=O_X\otimes_{f^*O_S}f^*M=D_{X\to S}\otimes_{f^*D^{\infty}_S}f^*M 
\end{equation*}
which induces in the derived category the functor
\begin{equation*}
Lf^{*mod}:D_{\mathcal D^{\infty}}(S)\to D_{\mathcal D^{\infty}}(X), \; \;
M\mapsto Lf^{*mod}M:=O_X\otimes^L_{f^*O_S}f^*M=O_X\otimes_{f^*O_S}f^*L_{D^{\infty}}M,
\end{equation*}
We will also consider the shifted inverse image functor
\begin{equation*}
Lf^{*mod[-]}:=Lf^{*mod}[d_S-d_X]:D_{\mathcal D^{\infty}}(S)\to D_{\mathcal D^{\infty}}(X).
\end{equation*}

\item Let $f:X\to S$ be a morphism with $X,S\in\SmVar(\mathbb C)$,
or let $f:X\to S$ be a morphism with $X,S\in\AnSm(\mathbb C)$.
For $M\in C_{\mathcal D}(X)$, $D_{X\leftarrow S}\otimes_{D_X}M$ has the canonical $f^*D_S$ module structure given above. 
Then, the direct image functor 
\begin{equation*}
f^0_{*mod}:\PSh_{\mathcal D}(X)\to\PSh_{\mathcal D}(S), \; \; 
M\mapsto f_{*mod} M:=f_*(D_{X\leftarrow S}\otimes_{D_X}M)
\end{equation*}
induces in the derived category the functor
\begin{equation*}
\int_{f}=Rf_{*mod}:D_{\mathcal D}(X)\to D_{\mathcal D}(S), \; \; 
M\mapsto \int_{f}M=Rf_*(D_{X\leftarrow S}\otimes^L_{D_X}M).
\end{equation*}

\item Let $f:X\to S$ be a morphism with $X,S\in\AnSm(\mathbb C)$.
For $M\in C_{\mathcal D^{\infty}}(X)$, 
$D_{X\leftarrow S}\otimes_{D_X}M$ has the canonical $f^*D_S$ module structure given above. 
Then, the direct image functor 
\begin{equation*}
f^{00}_{*mod}:\PSh_{\mathcal D^{\infty}}(X)\to\PSh_{\mathcal D^{\infty}}(S), \; \; 
M\mapsto f_{*mod} M:=f_*(D^{\infty}_{X\leftarrow S}\otimes_{D^{\infty}_X}M)
\end{equation*}
induces in the derived category the functor
\begin{equation*}
\int_{f}=Rf_{*mod}:D_{\mathcal D^{\infty}}(X)\to D_{\mathcal D^{\infty}}(S), \; \;
M\mapsto \int_{f}M=Rf_*(D^{\infty}_{X\leftarrow S}\otimes^L_{D^{\infty}_X}M).
\end{equation*}

\item Let $f:X\to S$ be a morphism with $X,S\in\AnSm(\mathbb C)$.
The direct image functor with compact support 
\begin{equation*}
f^{00}_{!mod}:\PSh_{\mathcal D}(X)\to\PSh_{\mathcal D}(S), \; \; M\mapsto f_{!mod} M:=f_!(D_{S\leftarrow X}\otimes_{D_X}M)
\end{equation*}
induces in the derived category the functor
\begin{equation*}
\int_{f!}=Rf_{!mod}:D_{\mathcal D}(X)\to D_{\mathcal D}(S), \; \; M\mapsto \int_{f}M=Rf_!(D_{X\leftarrow S}\otimes^L_{D_X}M).
\end{equation*}

\item Let $f:X\to S$ be a morphism with $X,S\in\AnSm(\mathbb C)$.
The direct image functor with compact support 
\begin{equation*}
f^{00}_{!mod}:\PSh_{\mathcal D^{\infty}}(X)\to\PSh_{\mathcal D^{\infty}}(S), \; \; 
M\mapsto f_{!mod} M:=f_!(D^{\infty}_{S\leftarrow X}\otimes_{D^{\infty}_X}M)
\end{equation*}
induces in the derived category the functor
\begin{equation*}
\int_{f!}=Rf_{!mod}:D_{\mathcal D^{\infty}}(X)\to D_{\mathcal D^{\infty}}(S), \; \; 
M\mapsto \int_{f}M=Rf_!(D^{\infty}_{X\leftarrow S}\otimes^L_{D^{\infty}_X}M).
\end{equation*}

\item Let $S\in\SmVar(\mathbb C)$. The analytical functor of a $D_S$ modules has a canonical structure of $D_{S^{an}}$ module:
\begin{equation*}  
(-)^{an}:C_{\mathcal D}(S)\to C_{\mathcal D}(S^{an}), \; M\mapsto M^{an}:=\an_S^{*mod}M:=M\otimes_{\an_S^*O_S}O_{S^{an}}
\end{equation*}
which induces in the derived category 
\begin{equation*}
(-)^{an}:D_{\mathcal D}(S)\to D_{\mathcal D}(S^{an}), \; M\mapsto M^{an}:=\an_S^{*mod}M)
\end{equation*}
since $\an_S^{*mod}$ derive trivially.

\end{itemize}

The functorialities given above induce :

\begin{itemize}

\item Let $f:X\to S$ be a morphism with $X,S\in\Var(\mathbb C)$,
or let $f:X\to S$ be a morphism with $X,S\in\AnSp(\mathbb C)$.
The adjonction map induces
\begin{itemize}
\item for $(M,F)\in C_{O_Sfil}(S)$, the map in $D_{O_Sfil}(S)$
\begin{eqnarray*}
\ad(Lf^{*mod},Rf_*)(M,F):(M,F)\xrightarrow{k\circ\ad(f^*,f_*)(M,F)}f_*E(f^*(M,F))=Rf_*f^*(M,F) \\
\xrightarrow{f_*m}Rf_*(f^*(M,F)\otimes^L_{f^*O_S}O_X)=Rf_*f^*(M,F),
\end{eqnarray*}
where $m(m)=m\otimes 1$,
\item for $(M,F)\in C_{O_Xfil}(X)$, the map in $D_{O_Xfil}(X)$
\begin{eqnarray*}
\ad(Lf^{*mod},Rf_*)(M,F):Lf^{*mod}Rf_*(M,F)=f^*f_*E(M,F)\otimes^L_{f^*O_S}O_X \\
\xrightarrow{\ad(f^*,f_*)(E(M,F))\otimes^L_{f^*O_S}O_X}(M,F)\otimes^L_{f^*O_S}O_X\xrightarrow{n}(M,F), 
\end{eqnarray*}
where $n(m\otimes h)=h.m$ is the multiplication map.
\end{itemize}

\item For a commutative diagram in $\Var(\mathbb C)$ or in $\AnSp(\mathbb C)$ : 
\begin{equation*}
D=\xymatrix{ 
Y\ar[r]^{g_2}\ar[d]^{f_2} & X\ar[d]^{f_1} \\
T\ar[r]^{g_1} & S},
\end{equation*}
we have, for $(M,F)\in C_{O_Xfil}(X)$,
the canonical map in $D_{O_Tfil}(T)$ 
\begin{eqnarray*}
T^{mod}(D)(M,F): Lg_1^{*mod}f_{1*}(M,F)\xrightarrow{\ad(Lf_2^{*mod},Rf_{2*})(Lg_1^{*mod}f_{1*}E(M,F))} \\
Rf_{2*}Lf_2^{*mod}Lg_1^{*mod}Rf_{1*}(M,F)=Rf_{2*}Lg_2^{*mod}Lf_1^{*mod}Rf_{1*}(M,F) \\ 
\xrightarrow{\ad(Lf_1^{*mod},Rf_1)(M,F)}Rf_{2*}Lg_2^{*mod}(M,F) 
\end{eqnarray*}
the canonical transformation map given by the adjonction maps.

\item Let $S\in\SmVar(\mathbb C)$ or $S\in\AnSm(\mathbb C)$.
For $(M,F)\in C_{fil}(S)$ and $(N,F)\in C_{fil}(S)$, recall that (see section 2)
\begin{equation*}
F^p((M,F)\otimes(N,F)):=\Im(\oplus_qF^qM\otimes F^{p-q}N\to M\otimes N)
\end{equation*}
This gives the functor
\begin{equation*}
(\cdot,\cdot):C_{fil}(S)\times C_{\mathcal Dfil}(S)\to C_{\mathcal Dfil}(S) \; , 
\;((M,F),(N,F))\mapsto (M,F)\otimes(N,F).
\end{equation*}
It induces in the derived categories by taking r-projective resolutions the bifunctors, for $r=1,\ldots,\infty$,
\begin{equation*}
(\cdot,\cdot):D_{\mathcal Dfil,r}(S)\times D_{fil,r}(S)\to D_{\mathcal Dfil,r}(S) \; , \;
((M,F),(N,F))\mapsto (M,F)\otimes^L(N,F)=L_D(M,F)\otimes (N,F).
\end{equation*}
For $S\in\AnSm(\mathbb C)$, it gives the bifunctor
\begin{equation*}
(\cdot,\cdot):C_{fil}(S)\times C_{\mathcal D^{\infty}fil}(S)\to C_{\mathcal D^{\infty}fil}(S) \; , 
\; ((M,F),(N,F))\mapsto (M,F)\otimes(N,F),
\end{equation*}
and its derived functor.

\item Let $S\in\SmVar(\mathbb C)$ or $S\in\AnSm(\mathbb C)$ and $O'_S\in\PSh(S,\cRing)$ a sheaf of commutative ring.
For $(M,F)\in C_{O'_Sfil}(S)$ and $(N,F)\in C_{O'_Sfil}(S)$, recall that (see section 2)
\begin{equation*}
F^p((M,F)\otimes_{O'_S}(N,F)):=\Im(\oplus_qF^qM\otimes_{O'_S} F^{p-q}N\to M\otimes_{O'_S} N)
\end{equation*}
It induces in the derived categories by taking r-projective resolutions the bifunctors, for $r=1,\ldots,\infty$,
\begin{equation*}
(\cdot,\cdot):D_{\mathcal Dfil,r}(S)^2\to D_{\mathcal Dfil,r}(S) \; , \;((M,F),(N,F))\mapsto (M,F)\otimes_{O_S}^L(N,F).
\end{equation*}
More generally, let $f:X\to S$ a morphism with $X,S\in\Var(\mathbb C)$ or with $X,S\in\AnSp(\mathbb C)$. Assume $S$ smooth.
We have the bifunctors
\begin{equation*}
(\cdot,\cdot):D_{f^*\mathcal Dfil,r}(X)^2\to D_{f^*\mathcal Dfil,r}(X) \; , \;
((M,F),(N,F))\mapsto (M,F)\otimes_{f^*O_S}^L(N,F)=(M,F)\otimes_{f^*O_S}L_{f^*D}(N,F).
\end{equation*}

\item Let $S\in\SmVar(\mathbb C)$ or let $S\in\AnSm(\mathbb C)$. 
The hom functor induces the bifunctor
\begin{eqnarray*}
\Hom(-,-):C_{\mathcal Dfil}(S)\times C_{fil}(S)\to C_{\mathcal D(1,0)fil}(S), ((M,W),(N,F))\mapsto\mathcal Hom((M,W),(N,F)).
\end{eqnarray*}
For $S\in\AnSm(\mathbb C)$, the hom functor also induces the bifunctor
\begin{eqnarray*}
\Hom(-,-):C_{\mathcal D^{\infty}fil}(S)\times C_{fil}(S)\to C_{\mathcal D^{\infty}(1,0)fil}(S), 
((M,W),(N,F))\mapsto\mathcal Hom((M,W),(N,F)).
\end{eqnarray*}
Note that the filtration given by $W$ satisfy that the $W^p$ are $D_S$ submodule which is stronger than Griffitz transversality.

\item Let $S\in\SmVar(\mathbb C)$ or let $S\in\AnSm(\mathbb C)$. 
The hom functor induces the bifunctor
\begin{eqnarray*}
\Hom_{O_S}(-,-):C_{\mathcal Dfil}(S)^2\to C_{\mathcal D2fil}(S), ((M,W),(N,F))\mapsto\mathcal Hom_{O_S}((M,W),(N,F)).
\end{eqnarray*}

\item Let $S\in\SmVar(\mathbb C)$ or let $S\in\AnSm(\mathbb C)$. 
The hom functor induces the bifunctors
\begin{itemize}
\item $\Hom_{D_S}(-,-):C_{\mathcal Dfil}(S)^2\to C_{2fil}(S)$, $((M,W),(N,F))\mapsto\mathcal Hom_{D_S}((M,W),(N,F))$,
\item $\Hom_{D_S}(-,-):C_{\mathcal D^{op}fil}(S)^2\to C_{2fil}(S)$, $((M,W),(N,F))\mapsto\mathcal Hom_{D_S}((M,W),(N,F))$.
\end{itemize}
We get the filtered dual 
\begin{eqnarray*}
\mathbb D^K_S(\cdot):C_{\mathcal D(2)fil}(S)\to C_{\mathcal D(2)fil}(S)^{op}, \; 
(M,F)\mapsto\mathbb D^K_S(M,F):=\mathcal Hom_{D_S}((M,F),D_S)\otimes_{O_S}\mathbb D_S^Ow(K_S)[d_S]
\end{eqnarray*}
together with the canonical map $d(M,F):(M,F)\to\mathbb D^{2,K}_S(M,F)$.
Of course $\mathbb D^K_S(\cdot)(C_{\mathcal D(1,0)fil}(S))\subset C_{\mathcal D(1,0)fil}(S)$.
It induces in the derived categories $D_{\mathcal Dfil,r}(S)$, for $r=1,\ldots,\infty$, the functors
\begin{equation*}
L\mathbb D_S(\cdot):D_{\mathcal D(2)fil,r}(S)\to D_{\mathcal D(2)fil,r}(S)^{op}, 
(M,F)\mapsto L\mathbb D_S(M,F):=\mathbb D^K_S L_D(M,F).
\end{equation*}
together with the canonical map $d(M,F):L_D(M,F)\to \mathbb D^{2,K}_SL_D(M,F)$.

\item Let $S\in\AnSm(\mathbb C)$. 
The hom functor also induces the bifunctors
\begin{itemize}
\item $\Hom_{D^{\infty}_S}(-,-):C_{\mathcal D^{\infty}fil}(S)^2\to C_{2fil}(S)$, 
$((M,W),(N,F))\mapsto\mathcal Hom_{D^{\infty}_S}((M,W),(N,F))$,
\item $\Hom_{D_S}(-,-):C_{\mathcal D^{\infty,op}fil}(S)^2\to C_{2fil}(S)$, 
$((M,W),(N,F))\mapsto\mathcal Hom_{D^{\infty}_S}((M,W),(N,F))$.
\end{itemize}
We get the filtered dual 
\begin{eqnarray*}
\mathbb D^{\infty,K}_S(\cdot):C_{\mathcal D^{\infty}(2)fil}(S)\to C_{\mathcal D^{\infty}(2)fil}(S)^{op}, \; 
(M,F)\mapsto\mathbb D^{\infty,K}_S(M,F):=\mathcal Hom_{D^{\infty}_S}((M,F),D^{\infty}_S)\otimes_{O_S}\mathbb D_S^Ow(K_S)[d_S]
\end{eqnarray*}
together with the canonical map $d(M,F):(M,F)\to\mathbb D^{\infty,2}_S(M,F)$. 
Of course $\mathbb D^{\infty,K}_S(\cdot)(C_{\mathcal D^{\infty}(1,0)fil}(S))\subset C_{\mathcal D^{\infty}(1,0)fil}(S)$.
It induces in the derived categories $D_{\mathcal Dfil,r}(S)$, for $r=1,\ldots,\infty$, the functors
\begin{equation*}
L\mathbb D^{\infty}_S(\cdot):D_{\mathcal D^{\infty}(2)fil,r}(S)\to D_{\mathcal D^{\infty}(2)fil,r}(S)^{op}, 
(M,F)\mapsto L\mathbb D^{\infty}_S(M,F)=\mathbb D^{\infty,K}_S L_{D^{\infty}}(M,F).
\end{equation*}
together with the canonical map $d(M,F):(M,F)\to L\mathbb D^{\infty,2}_S(M,F)$.

\item Let $f:X\to S$ be a morphism with $X,S\in\SmVar(\mathbb C)$,
or let $f:X\to S$ be a morphism with $X,S\in\AnSm(\mathbb C)$.
Then, the inverse image functor
\begin{eqnarray*}
f^{*mod}:C_{\mathcal D(2)fil}(S)\to C_{\mathcal D(2)fil}(X), \\
(M,F)\mapsto f^{*mod}(M,F):=(O_X,F_b)\otimes_{f^*O_S}f^*(M,F)=(D_{X\to S},F^{ord})\otimes_{f^*D_S}f^*(M,F), 
\end{eqnarray*} 
induces in the derived categories the functors, for $r=1,\ldots,\infty$ (resp. $r\in(1,\ldots\infty)^2$),
\begin{eqnarray*}
Lf^{*mod}:D_{\mathcal D(2)fil,r}(S)\to D_{\mathcal D(2)fil,r}(X), \\  
(M,F)\mapsto Lf^{*mod}M:=(O_X,F_b)\otimes^L_{f^*O_S}f^*(M,F)=(O_X,F_b)\otimes_{f^*O_S}f^*L_D(M,F).
\end{eqnarray*}
Of course $f^{*mod}(C_{\mathcal D(1,0)fil}(S))\subset C_{\mathcal D(1,0)fil}(X)$. Note that
\begin{itemize}
\item If the $M$ is a complex of locally free $O_S$ modules, 
then $Lf^{*mod}(M,F)=f^{*mod}(M,F)$ in $D_{\mathcal D(2)fil,\infty}(S)$.
\item If the $\Gr^p_FM$ are complexes of locally free $O_S$ modules, 
then $Lf^{*mod}(M,F)=f^{*mod}(M,F)$ in $D_{\mathcal D(2)fil}(S)$.
\end{itemize}
We will consider also the shifted inverse image functors
\begin{equation*}
Lf^{*mod[-]}:=Lf^{*mod}[d_S-d_X]:D_{\mathcal D(2)fil,r}(S)\to D_{\mathcal D(2)fil,r}(X). 
\end{equation*}

\item Let $f:X\to S$ be a morphism with $X,S\in\AnSm(\mathbb C)$.
Then, the inverse image functor
\begin{eqnarray*}
f^{*mod}:C_{\mathcal D^{\infty}(2)fil}(S)\to C_{\mathcal D^{\infty}(2)fil}(X),  
(M,F)\mapsto f^{*mod}(M,F):=(O_X,F_b)\otimes_{f^*O_S}f^*(M,F), 
\end{eqnarray*} 
induces in the derived categories the functors, for $r=1,\ldots,\infty$ (resp. $r\in(1,\ldots\infty)^2$),
\begin{eqnarray*}
Lf^{*mod}:D_{\mathcal D^{\infty}(2)fil,r}(S)\to D_{\mathcal D^{\infty}(2)fil,r}(X), \\  
(M,F)\mapsto Lf^{*mod}M:=(O_X,F_b)\otimes^L_{f^*O_S}f^*(M,F)=(O_X,F_b)\otimes_{f^*O_S}f^*L_{D^{\infty}}(M,F).
\end{eqnarray*}
Of course $f^{*mod}(C_{\mathcal D^{\infty}(1,0)fil}(S))\subset C_{\mathcal D^{\infty}(1,0)fil}(X)$. Note that
We will consider also the shifted inverse image functors
\begin{equation*}
Lf^{*mod[-]}:=Lf^{*mod}[d_S-d_X]:D_{\mathcal D^{\infty}(2)fil,r}(S)\to D_{\mathcal D^{\infty}(2)fil,r}(X). 
\end{equation*}

\item Let $f:X\to S$ be a morphism with $X,S\in\SmVar(\mathbb C)$,
or let $f:X\to S$ be a morphism with $X,S\in\AnSm(\mathbb C)$.
Then,the direct image functor 
\begin{equation*}
f^{00}_{*mod}:(\PSh_{\mathcal D}(X),F)\to(\PSh_{\mathcal D}(S),F), \; \;
(M,F)\mapsto f_{*mod} (M,F):=f_*((D_{S\leftarrow X},F^{ord})\otimes_{D_X}(M,F))
\end{equation*}
induces in the derived categories by taking r-injective resolutions the functors, for $r=1,\ldots,\infty$,
\begin{eqnarray*}
\int_{f}=Rf_{*mod}:D_{\mathcal D(2)fil,r}(X)\to D_{\mathcal D(2)fil,r}(S), 
(M,F)\mapsto \int_{f}(M,F)=Rf_*((D_{S\leftarrow X},F^{ord})\otimes^L_{D_X}(M,F)).
\end{eqnarray*}
Let $f_1:X\to Y$ and $f_2:Y\to S$ two morphism with $X,Y,S\in\SmVar(\mathbb C)$ or with $X,Y,S\in\AnSm(\mathbb C)$.
We have, for $(M,F)\in C_{\mathcal Dfil}(X)$, the canonical transformation map in $D_{\mathcal D(2)fil,r}(S)$
\begin{eqnarray*}
T(\int_{f_2}\circ\int_{f_1},\int_{f_2\circ f_1})(M,F): \\
\int_{f_2}\int_{f_1}(M,F):=
Rf_{2*}((D_{Y\leftarrow S},F^{ord})\otimes^L_{D_Y}Rf_{1*}((D_{X\leftarrow Y},F^{ord})\otimes^L_{D_X}(M,F))) \\
\xrightarrow{T(f_1,\otimes)(-,-)}
Rf_{2*}Rf_{1*}(f_1^*(D_{Y\leftarrow S},F^{ord})\otimes^L_{D_Y}((D_{X\leftarrow Y},F^{ord})\otimes^L_{D_X}(M,F))) \\
\xrightarrow{\sim} 
Rf_{2*}Rf_{1*}((f_1^*(D_{Y\leftarrow S},F^{ord})\otimes^L_{D_Y}(D_{X\leftarrow Y},F^{ord}))\otimes^L_{D_X}(M,F)) \\
\xrightarrow{\sim}
Rf_{2*}Rf_{1*}((D_{X\leftarrow S},F^{ord})\otimes^L_{D_X}(M,F)):=\int_{f_2\circ f_1}(M,F)
\end{eqnarray*}

\item Let $f:X\to S$ be a morphism with $X,S\in\AnSm(\mathbb C)$. Then,the direct image functor 
\begin{equation*}
f^{00}_{*mod}:(\PSh_{\mathcal D^{\infty}}(X),F)\to(\PSh_{\mathcal D^{\infty}}(S),F), \; \;
(M,F)\mapsto f_{*mod} (M,F):=f_*((D^{\infty}_{S\leftarrow X},F^{ord})\otimes_{D^{\infty}_X}(M,F))
\end{equation*}
induces in the derived categories by taking r-injective resolutions the functors, for $r=1,\ldots,\infty$,
\begin{eqnarray*}
\int_{f}=Rf_{*mod}:D_{\mathcal D^{\infty}(2)fil,r}(X)\to D_{\mathcal D^{\infty}(2)fil,r}(S), 
(M,F)\mapsto \int_{f}(M,F)=Rf_*((D^{\infty}_{S\leftarrow X},F^{ord})\otimes^L_{D^{\infty}_X}(M,F)).
\end{eqnarray*}
We have, similarly, for $(M,F)\in C_{\mathcal D^{\infty}fil}(X)$, 
the canonical transformation map in $D_{\mathcal D^{\infty}(2)fil,r}(S)$
\begin{equation*}
T(\int_{f_2}\circ\int_{f_1},\int_{f_2\circ f_1})(M,F):\int_{f_2}\int_{f_1}(M,F)\to\int_{f_2\circ f_1}(M,F)
\end{equation*}

\item Let $f:X\to S$ be a morphism with $X,S\in\AnSm(\mathbb C)$.
Then,the direct image functor with compact support 
\begin{equation*}
f^{00}_{!mod}:(\PSh_{\mathcal D}(X),F)\to(\PSh_{\mathcal D}(S),F), \; \;
(M,F)\mapsto f^{00}_{!mod} (M,F):=f_!((D_{S\leftarrow X},F^{ord})\otimes_{D_X}(M,F))
\end{equation*}
induces in the derived categories by taking r-injective resolutions the functors, for $r=1,\ldots,\infty$,
\begin{equation*}
\int_{f!}=Rf_{!mod}:D_{\mathcal Dfil,r}(X)\to D_{\mathcal Dfil,r}(S), \; \;
(M,F)\mapsto \int_{f}(M,F)=Rf_!((D_{S\leftarrow X},F^{ord})\otimes^L_{D_X}(M,F)).
\end{equation*}
We have, similarly, for $(M,F)\in C_{\mathcal Dfil}(X)$, the canonical transformation map in $D_{\mathcal D(2)fil,r}(S)$
\begin{equation*}
T(\int_{f_2!}\circ\int_{f_1!},\int_{(f_2\circ f_1})!)(M,F):\int_{f_2!}\int_{f_1!}(M,F)\to\int_{(f_2\circ f_1)!}(M,F)
\end{equation*}

\item Let $f:X\to S$ be a morphism with $X,S\in\AnSm(\mathbb C)$.
Then,the direct image functor with compact support 
\begin{equation*}
f^{00}_{!mod}:(\PSh_{\mathcal D^{\infty}}(X),F)\to(\PSh_{\mathcal D^{\infty}}(S),F), \; \;
(M,F)\mapsto f^{00}_{!mod} (M,F):=f_!((D^{\infty}_{S\leftarrow X},F^{ord})\otimes_{D^{\infty}_X}(M,F))
\end{equation*}
induces in the derived categories by taking r-injective resolutions the functors, for $r=1,\ldots,\infty$,
\begin{eqnarray*}
\int_{f!}=Rf_{!mod}:D_{\mathcal D^{\infty}fil,r}(X)\to D_{\mathcal D^{\infty}fil,r}(S), 
(M,F)\mapsto \int_{f}(M,F)=Rf_!((D^{\infty}_{S\leftarrow X},F^{ord})\otimes^L_{D^{\infty}_X}(M,F)).
\end{eqnarray*}
We have, similarly, for $(M,F)\in C_{\mathcal D^{\infty}fil}(X)$, 
the canonical transformation map in $D_{\mathcal D^{\infty}(2)fil,r}(S)$
\begin{equation*}
T(\int_{f_2!}\circ\int_{f_1!},\int_{(f_2\circ f_1})!)(M,F):\int_{f_2!}\int_{f_1!}(M,F)\to\int_{(f_2\circ f_1)!}(M,F)
\end{equation*}

\item Let $S\in\SmVar(\mathbb C)$. 
The analytical functor for filtered $D_S$-modules is 
\begin{equation*}  
(\cdot)^{an}:C_{\mathcal D(2)fil}(S)\to C_{\mathcal D(2)fil}(S^{an}), \;  
(M,F)\mapsto (M,F)^{an}:=\an_S^*(M,F)\otimes_{\an_S^*O_S}(O_{S^{an}},F_b).
\end{equation*}
It induces in the derived categories the functors, for $r=1,\ldots,\infty$, 
\begin{equation*}
(\cdot)^{an}:D_{\mathcal D(2)fil,r}(S)\to D_{\mathcal D(2)fil,r}(S^{an}), \; 
(M,F)\mapsto (M,F)^{an}:=\an_S^*(M,F)\otimes^{L}_{\an_S^*O_S}(O_{S^{an}},F_b).
\end{equation*}

\item  Let $f:X\to S$ be a morphism with $X,S\in\SmVar(\mathbb C)$,
or let $f:X\to S$ be a morphism with $X,S\in\AnSm(\mathbb C)$.
Then the functor
\begin{eqnarray*}
f^{\hat{*}mod}:C_{\mathcal D2fil}(S)\to C_{\mathcal D2fil}(X), 
(M,F)\mapsto f^{\hat{*}mod}(M,F):=\mathbb D^K_XL_Df^{*mod}L_D\mathbb D^K_S(M,F) 
\end{eqnarray*} 
induces in the derived categories the exceptional inverse image functors, 
for $r=1,\ldots,\infty$ (resp. $r\in(1,\ldots\infty)^2$),
\begin{eqnarray*}
Lf^{\hat{*}mod}:D_{\mathcal D(2)fil,r}(S)\to D_{\mathcal D(2)fil,r}(X), \\
(M,F)\mapsto Lf^{\hat{*}mod}(M,F):=L\mathbb D_XLf^{*mod}L\mathbb D_S(M,F):=f^{\hat{*}mod}L_D(M,F).
\end{eqnarray*}
Of course $f^{\hat*mod}(C_{\mathcal D(1,0)fil}(S))\subset C_{\mathcal D(1,0)fil}(X)$.
We will also consider the shifted exceptional inverse image functors
\begin{equation*}
Lf^{\hat{*}mod[-]}:=Lf^{\hat{*}mod}[d_S-d_X]:D_{\mathcal D(2)fil,r}(S)\to D_{\mathcal D(2)fil,r}(X).
\end{equation*}

\item  Let $f:X\to S$ be a morphism with $X,S\in\AnSm(\mathbb C)$.
Then the functor
\begin{eqnarray*}
f^{\hat{*}mod}:C_{\mathcal D^{\infty}(2)fil}(S)\to C_{\mathcal D^{\infty}(2)fil}(X), 
(M,F)\mapsto f^{\hat{*}mod}(M,F):=\mathbb D^{K,\infty}_XL_Df^{*mod}L_D\mathbb D^{K,\infty}_S(M,F) 
\end{eqnarray*} 
induces in the derived categories the exceptional inverse image functors, 
for $r=1,\ldots,\infty$ (resp. $r\in(1,\ldots\infty)^2$),
\begin{eqnarray*}
Lf^{\hat{*}mod}:D_{\mathcal D^{\infty}(2)fil,r}(S)\to D_{\mathcal D^{\infty}(2)fil,r}(X), \\ 
(M,F)\mapsto Lf^{\hat{*}mod}(M,F):=L\mathbb D^{\infty}_XLf^{*mod}L\mathbb D^{\infty}_S(M,F):=f^{\hat{*}mod}(M,F).
\end{eqnarray*}
Of course $f^{\hat*mod}(C_{\mathcal D^{\infty}(1,0)fil}(S))\subset C_{\mathcal D^{\infty}(1,0)fil}(X)$.
We will also consider the shifted exceptional inverse image functors
\begin{equation*}
Lf^{\hat{*}mod[-]}:=Lf^{\hat{*}mod}[d_S-d_X]:D_{\mathcal D^{\infty}(2)fil,r}(S)\to D_{\mathcal D^{\infty}(2)fil,r}(X).
\end{equation*}

\item Let $S_1,S_2\in\SmVar(\mathbb C)$ or $S_1,S_2\in\AnSm(\mathbb C)$. Consider $p:S_1\times S_2\to S_1$ the projection. 
Since $p$ is a projection, we have a canonical embedding $p^*D_{S_1}\hookrightarrow D_{S_1\times S_2}$.
For $(M,F)\in C_{\mathcal D(2)fil}(S_1\times S_2)$, $(M,F)$ has a canonical $p^*D_{S_1}$ module structure. 
Moreover, with this structure, for $(M_1,F)\in C_{\mathcal D(2)fil}(S_1)$ 
\begin{equation*}
\ad(p^{*mod},p)(M_1,F):(M_1,F)\to p_*p^{*mod}(M_1,F)
\end{equation*}
is a map of complexes of $D_{S_1}$ modules, and for $(M_{12},F)\in C_{\mathcal D(2)fil}(S_1\times S_2))$
\begin{equation*}
\ad(p^{*mod},p)(M_{12},F):p^{*mod}p_*(M_{12},F)\to (M_{12},F)
\end{equation*}
is a map of complexes of $D_{S_1\times S_2}$ modules.
Indeed, for the first adjonction map, we note that $p^{*mod}(M_1,F)$ has a structure of $p^*D_{S_1}$ module,
hence $p_*p^{*mod}(M_1,F)$ has a structure of $p_*p^*D_{S_1}$ module, hence a structure of $D_{S_1}$ module using the adjonction 
map $\ad(p^*,p_*)(D_{S_1}):D_{S_1}\to p_*p^*D_{S_1}$.
For the second adjonction map, we note that $(M_{12},F)$ has a structure of $p^*D_{S_1}$ module, hence
$p_*(M_{12},F)$ has a structure of $p_*p^*D_{S_1}$, hence a structure of $D_{S_1}$ module using the adjonction 
map $\ad(p^*,p_*)(D_{S_1}):D_{S_1}\to p_*p^*D_{S_1}$.

\item Let $S_1,S_2\in\AnSm(\mathbb C)$. Consider $p:S_1\times S_2\to S_1$ the projection. 
Since $p$ is a projection, we have a canonical embedding $p^*D^{\infty}_{S_1}\hookrightarrow D^{\infty}_{S_1\times S_2}$.
For $(M,F)\in C_{\mathcal D^{\infty}(2)fil}(S_1\times S_2)$, $(M,F)$ has a canonical $p^*D^{\infty}_{S_1}$ module structure. 
Moreover, with this structure, for $(M_1,F)\in C_{\mathcal D^{\infty}(2)fil}(S_1)$ 
\begin{equation*}
\ad(p^{*mod},p)(M_1,F):(M_1,F)\to p_*p^{*mod}(M_1,F)
\end{equation*}
is a map of complexes of $D^{\infty}_{S_1}$ modules, and for $(M_{12},F)\in C_{\mathcal D^{\infty}(2)fil}(S_1\times S_2)$
\begin{equation*}
\ad(p^{*mod},p)(M_{12},F):p^{*mod}p_*(M_{12},F)\to (M_{12},F)
\end{equation*}
is a map of complexes of $D^{\infty}_{S_1\times S_2}$ modules, similarly to the finite order case.

\end{itemize}

We following proposition concern the commutativity of the inverse images functors 
and the commutativity of the direct image functors.

\begin{prop}\label{compDmod}
\begin{itemize}
\item[(i)] Let $f_1:X\to Y$ and $f_2:Y\to S$ two morphism with $X,Y,S\in\SmVar(\mathbb C)$. 
\begin{itemize}
\item Let $(M,F)\in C_{\mathcal D(2)fil,r}(S)$. Then $(f_2\circ f_1)^{*mod}(M,F)=f_1^{*mod}f_2^{*mod}(M,F)$.
\item Let $(M,F)\in D_{\mathcal D(2)fil,r}(S)$. Then $L(f_2\circ f_1)^{*mod}(M,F)=Lf_1^{*mod}(Lf_2^{*mod}(M,F))$.
\end{itemize}
\item[(ii)] Let $f_1:X\to Y$ and $f_2:Y\to S$ two morphism with $X,Y,S\in\SmVar(\mathbb C)$. 
Let $M\in D_{\mathcal D}(X)$. Then, 
\begin{equation*}
T(\int_{f_2}\circ\int_{f_1},\int_{f_2\circ f_1})(M):\int_{f_2}\int_{f_1}(M)\xrightarrow{\sim}\int_{f_2\circ f_1}(M)
\end{equation*}
is an isomorphism in $D_{\mathcal D}(S)$ (i.e. if we forget filtration).
\item[(iii)] Let $i_0:Z_2\hookrightarrow Z_1$ and $i_1:Z_1\hookrightarrow S$ two closed embedding, 
with $Z_2,Z_1,S\in\SmVar(\mathbb C)$. Let $(M,F)\in C_{\mathcal D(2)fil}(Z_2)$. 
Then, $(i_1\circ i_0)_{*mod}(M,F)=i_{1*mod}(i_{0*mod}(M,F))$ in  $C_{\mathcal D(2)fil}(S)$.
\end{itemize}
\end{prop}

\begin{proof}
\noindent(i): Obvious : we have
\begin{itemize}
\item $f_1^{*mod}f_2^{*mod}(M,F)=f_1^*(f_2^*(M,F)\otimes_{f_2^*O_S}O_Y)\otimes_{f_1^*O_Y}O_X
=f_1^*f_2^*(M,F)\otimes_{f_1^*f_2^*O_S}f_1^*O_Y\otimes_{f_1^*O_Y}O_X=(f_2\circ f_1)^{*mod}(M,F)$
\item $Lf_1^{*mod}Lf_2^{*mod}(M,F)=f_1^*(f_2^*(M,F)\otimes^L_{f_2^*O_S}O_Y)\otimes^L_{f_1^*O_Y}O_X
=f_1^*f_2^*(M,F)\otimes^L_{f_1^*f_2^*O_S}f_1^*O_Y\otimes^L_{f_1^*O_Y}O_X=L(f_2\circ f_1)^{*mod}(M,F)$
\end{itemize}.

\noindent(ii): See \cite{LvDmod} : we have by lemma \ref{Dcomp}
\begin{eqnarray*}
\int_{f_2\circ f_1}M:=Rf_{2*}Rf_{1*}(D_{X\leftarrow S}\otimes^L_{D_X}M) \\
\xrightarrow{=}
Rf_{2*}Rf_{1*}((f_1^*D_{Y\leftarrow S}\otimes^L_{f_1^*D_Y}D_{X\leftarrow Y})\otimes^L_{D_X}M) 
\xrightarrow{Rf_{2*}T(f_1,\otimes)(D_{Y\leftarrow S},D_{X\leftarrow Y}\otimes^L_{D_X}M)^{-1}} \\
Rf_{2*}(D_{Y\leftarrow S}\otimes^L_{D_Y}Rf_{1*}(D_{X\leftarrow Y}\otimes^L_{D_X}M))
=:\int_{f_2}\int{f_1}M
\end{eqnarray*}
where, $D_{Y\leftarrow S}$ being a quasi-coherent $D_Y$ module,
we used the fact that for $N\in C_{f_1^*\mathcal D}(X)$ and $N'\in C_{\mathcal D}(Y)$ 
\begin{equation*}
T(f_1,\otimes)(N',N):N'\otimes^L_{D_Y}Rf_{1*}N\to Rf_{1*}(f_1^*N'\otimes^L_{f_1^*D_Y}N)
\end{equation*}
is an isomorphism if $N'$ is quasi-coherent, 
which follows from the fact that $f_{1*}$ commutes with arbitrary (possibly infinite) direct sums (see \cite{LvDmod}).

\noindent(iii): Denote $i_2=i_1\circ i_0:Z_2\hookrightarrow S$. We have
\begin{eqnarray*}
i_{2*mod}(M,F)=i_{2*}((M,F)\otimes_{D_{Z_2}}(D_{Z_2\leftarrow S},F^{ord}))\xrightarrow{=} \\ 
i_{1*}i_{0*}((M,F)\otimes_{D_{Z_2}}(D_{Z_2\leftarrow Z_1},F^{ord})\otimes_{i_0^*D_{Z_1}}i_0^*(D_{Z_1\leftarrow S},F^{ord}))
\xrightarrow{i_{1*}T(i_0,\otimes)(-)^{-1}} i_{1*mod}i_{0*mod}((M,F))
\end{eqnarray*}
using proposition \ref{Dcomp} and proposition \ref{Tiotimes}.
\end{proof}

\begin{rem}\label{compDmodrem}
Let $f_1:X\to Y$ and $f_2:Y\to S$ two morphism with $X,Y,S\in\SmVar(\mathbb C)$.
Then, for $(M,F)\in D_{\mathcal D(2)fil,\infty}(X)$, 
$\int_{f_2}\int_{f_1}(M,F)$ is NOT isomorphic to $\int_{f_2\circ f_1}(M,F)$ in general, the filtrations
on the isomorphic cohomology sheaves may be different.
\end{rem}

\begin{prop}\label{compDmodh}
Let $f:X\to S$ a morphism with $X,S\in\SmVar(\mathbb C)$. Then,
\begin{itemize}
\item[(i)] For $(M,F)\in C_{\mathcal D(2)fil,h}(S)$, we have $Lf^{*mod}(M,F)\in D_{\mathcal D(2)fil,\infty,h}(X)$.
\item[(ii)] For $M\in C_{\mathcal D,h}(X)$, we have $\int_fM\in D_{\mathcal D,h}(S)$.
\item[(iii)] If $f$ is proper, for $(M,F)\in C_{\mathcal D(2)fil,h}(X)$, 
we have $\int_f(M,F)\in D_{\mathcal D(2)fil,\infty,h}(S)$.
\end{itemize}
\end{prop}

\begin{proof}
See \cite{LvDmod} for the non filtered case. 
The filtered case follows immediately from the non filtered case 
and the fact the pullback of a good filtration is a good filtration 
(since the pullback of a coherent $O_S$ module is coherent) 
and the direct image of a good filtration by a proper morphism is a good filtration
(since the pushforward of a coherent $O_X$ module by a proper morphism is coherent).
\end{proof}
 
The following easy proposition says that the analytical functor commutes we the pullback of $D$ modules and the tensor product. 
Again it is well known in the non filtered case. Note that for $S\in\SmVar(\mathbb C)$, $D_S^{an}=D_{S^{an}}$.

\begin{prop}\label{anmor} 
\begin{itemize}
\item[(i)] Let $f:T\to S$ a morphism with $T,S\in\SmVar(\mathbb C)$. 
\begin{itemize}
\item Let $(M,F)\in C_{\mathcal D(2)fil,r}(S)$. Then $(f^{*mod}(M,F))^{an}=f^{*mod}(M,F)^{an}$.
\item Let $(M,F)\in D_{\mathcal Dfil,r}(S)$, for $r=1,\ldots\infty$. Then, $(Lf^{*mod}(M,F))^{an}=Lf^{*mod}(M,F)^{an}$.
\end{itemize}
\item[(ii)] Let $S\in\SmVar(\mathbb C)$
\begin{itemize}
\item Let $(M,F),(N,F)\in C_{\mathcal Dfil}(S)$. 
Then, $((M,F)\otimes_{O_S}(N,F))^{an}=(M,F)^{an}\otimes_{O_{S^{an}}}(N,F)^{an}$. 
\item Let $(M,F),(N,F)\in D_{\mathcal Dfil,r}(S)$, for $r=1,\ldots\infty$. 
Then, $((M,F)\otimes^L_{O_S}(N,F))^{an}=(M,F)^{an}\otimes^L_{O_{S^{an}}}(N,F)^{an}$. 
\end{itemize}
\end{itemize}
\end{prop}

\begin{proof}
\noindent(i): For $(M,F)\in C_{\mathcal D,fil}(S)$, we have, since $f^*\an_S^*=\an_X^*f^{an*}$, 
\begin{eqnarray*}
(f^{*mod}(M,F))^{an}&=&\an(X)^*(f^*(M,F)\otimes_{f^*O_S}O_X)\otimes_{\an(X)^*O_X}O_{X^{an}} \\
&=&f^{an*}\an_S^*(M,F)\otimes_{f^{an*}O_{S^{an}}}\otimes O_{X^{an}}=:f^{an*mod}(M^{an},F)
\end{eqnarray*}
For $(M,F)\in D_{\mathcal D,fil,r}(S)$,
we take $(M,F)\in C_{\mathcal D,fil}(S)$ an r-projective $f^*O_S$ module such that $D_{top,r}(M,F)=(M,F)$ so that  
\begin{equation*}
(Lf^{*mod}(M,F))^{an}=(f^{*mod}(M,F))^{an}=f^{an*mod}(M^{an},F)=Lf^{an*mod}(M^{an},F)
\end{equation*}

\noindent(ii): For $(M,F),(N,F)\in C_{\mathcal D,fil}(S)$, we have
\begin{eqnarray*}
((M,F)\otimes_{O_S} (N,F))^{an}:&=&\an_S^*((M,F)\otimes_{O_S}(N,F))\otimes_{\an_S^*O_S}O_{S^{an}} \\
&=&\an_S^*(M,F)\otimes_{\an_S^*O_S}\an_S^*(N,F)\otimes_{\an_S^*O_S}O_{S^{an}} \\
&=&\an_S^*(M,F)\otimes_{\an_S^*O_S}\otimes O_{S^{an}}\otimes_{O_{S^{an}}}\an_S^*(N,F)\otimes_{\an_S^*O_S}O_{S^{an}} \\
&=&:(M^{an},F)\otimes_{O_{S^{an}}}(N^{an},F)
\end{eqnarray*}
It implies the isomorphism in the derived category by taking an r-projective resolution of $(M,F)$ (e.g $(L_D(M),F)=L_D(M,F)$). 
\end{proof}

\begin{prop}\label{compAnDmod}
\begin{itemize}
\item[(i)] Let $f_1:X\to Y$ and $f_2:Y\to S$ two morphism with $X,Y,S\in\AnSm(\mathbb C)$. 
\begin{itemize}
\item Let $(M,F)\in C_{\mathcal D(2)fil}(S)$ or let $(M,F)\in C_{\mathcal D^{\infty}(2)fil}(S)$. 
Then $(f_2\circ f_1)^{*mod}(M,F)=f_1^{*mod}f_2^{*mod}(M,F)$.
\item Let $(M,F)\in D_{\mathcal D(2)fil,r}(S)$ or let $(M,F)\in D_{\mathcal D^{\infty}(2)fil,r}(S)$.
Then $L(f_2\circ f_1)^{*mod}(M,F)=Lf_1^{*mod}(Lf_2^{*mod}(M,F))$.
\end{itemize}
\item[(ii)] Let $f_1:X\to Y$ and $f_2:Y\to S$ two morphisms with $X,Y,S\in\AnSm(\mathbb C)$. 
Let $M\in D_{\mathcal D}(X)$. 
If $f_1$ is proper, we have $\int_{f_2\circ f_1}M=\int_{f_2}(\int_{f_1}M)$.
\item[(ii)'] Let $f_1:X\to Y$ and $f_2:Y\to S$ two morphisms with $X,Y,S\in\AnSm(\mathbb C)$. 
Let $M\in D_{\mathcal D^{\infty}}(X)$. 
If $f_1$ is proper, we have $\int_{f_2\circ f_1}M=\int_{f_2}(\int_{f_1}M)$.
\item[(iii)] Let $f_1:X\to Y$ and $f_2:Y\to S$ two morphisms with $X,Y,S\in\AnSm(\mathbb C)$. 
Let $(M,F)\in D_{\mathcal D}(X)$. 
We have $\int_{(f_2\circ f_1)!}M=\int_{f_2!}(\int_{f_1!}M)$.
\item[(iii)'] Let $f_1:X\to Y$ and $f_2:Y\to S$ two morphisms with $X,Y,S\in\AnSm(\mathbb C)$. 
Let $M\in D_{\mathcal D^{\infty}}(X)$. 
We have $\int_{(f_2\circ f_1)!}M=\int_{f_2!}(\int_{f_1!}M)$.
\item[(iv)] Let $i_0:Z_2\hookrightarrow Z_1$ and $i_1:Z_1\hookrightarrow S$ two closed embedding,
 with $Z_2,Z_1,S\in\AnSm(\mathbb C)$. Let $(M,F)\in C_{\mathcal D(2)fil}(Z_2)$. 
Then, $(i_1\circ i_0)_{*mod}(M,F)=i_{1*mod}(i_{0*mod}(M,F))$ in  $C_{\mathcal D(2)fil}(S)$.
\item[(iv)'] Let $i_0:Z_2\hookrightarrow Z_1$ and $i_1:Z_1\hookrightarrow S$ two closed embedding,
with $Z_2,Z_1,S\in\AnSm(\mathbb C)$. Let $(M,F)\in C_{\mathcal D^{\infty}(2)fil}(Z_2)$. 
Then, $(i_1\circ i_0)_{*mod}(M,F)=i_{1*mod}(i_{0*mod}(M,F))$ in  $C_{\mathcal D^{\infty}(2)fil}(S)$.
\end{itemize}
\end{prop}

\begin{proof}
\noindent(i): Similar to the proof of proposition \ref{compDmod}(i).

\noindent(ii): Similar to the proof of proposition \ref{compDmod}(ii) : we use lemma \ref{Dcomp} and the fact that
for $N\in C_{f_1^*\mathcal D}(X)$ and $N'\in C_{\mathcal D}(Y)$, the canonical morphism 
\begin{equation*}
T(f_1,\otimes)(N',N):N'\otimes^L_{D_Y}Rf_{1*}N\to Rf_{1*}(f_1^*N'\otimes^L_{f_1^*D_Y}N)
\end{equation*}
is an isomorphism if $f_1$ is proper (in this case $f_{1!}=f_{1*}$).

\noindent(ii)': Similar to the proof of proposition \ref{compDmod}(ii) : we use lemma \ref{Dcompinfty} and the fact that
for $N\in C_{f_1^*\mathcal D^{\infty}}(X)$ and $N'\in C_{\mathcal D^{\infty}}(Y)$, the canonical morphism 
\begin{equation*}
T(f_1,\otimes)(N',N):N'\otimes^L_{D^{\infty}_Y}Rf_{1*}N\to Rf_{1*}(f_1^*N'\otimes^L_{f_1^*D^{\infty}_Y}N)
\end{equation*}
is an isomorphism if $f_1$ is proper (in this case $f_{1!}=f_{1*}$).

\noindent(iii): Similar to the proof of proposition \ref{compDmod}(ii) : we use lemma \ref{Dcomp} and the fact that
for $N\in C_{f_1^*\mathcal D}(X)$ and $N'\in C_{\mathcal D}(Y)$, the canonical morphism 
\begin{equation*}
T(f_1!,\otimes)(N',N):N'\otimes^L_{D_Y}Rf_{1!}N\to Rf_{1!}(f_1^*N'\otimes^L_{f_1^*D_Y}N)
\end{equation*}
is an isomorphism.

\noindent(iii)': Similar to the proof of proposition \ref{compDmod}(ii) : we use lemma \ref{Dcompinfty} and the fact that
for $N\in C_{f_1^*\mathcal D^{\infty}}(X)$ and $N'\in C_{\mathcal D^{\infty}}(Y)$, the canonical morphism 
\begin{equation*}
T(f_1!,\otimes)(N',N):N'\otimes^L_{D^{\infty}_Y}Rf_{1!}N\to Rf_{1!}(f_1^*N'\otimes^L_{f_1^*D^{\infty}_Y}N)
\end{equation*}
is an isomorphism

\noindent(iv): Similar to the proof of proposition \ref{compDmod}(iii) :we have
\begin{eqnarray*}
i_{2*mod}(M,F)=i_{2*}((M,F)\otimes_{D_{Z_2}}(D_{Z_2\leftarrow S}),F^{ord})\xrightarrow{=} \\
i_{1*}i_{0*}((M,F)\otimes_{D_{Z_2}}(D_{Z_2\leftarrow Z_1},F^{ord})\otimes_{i_0^*D_{Z_1}}i_0^*(D_{Z_1\leftarrow S}),F^{ord}) 
\xrightarrow{i_{1*}T(i_0,\otimes)(-)^{-1}} i_{1*mod}i_{0*mod}((M,F))
\end{eqnarray*}
using lemma \ref{Dcomp} and proposition \ref{Tiotimes}.

\noindent(iv)':Similar to (iv): we have
\begin{eqnarray*}
i_{2*mod}(M,F)=i_{2*}((M,F)\otimes_{D^{\infty}_{Z_2}}(D^{\infty}_{Z_2\leftarrow S},F^{ord}))
\xrightarrow{=} \\
i_{1*}i_{0*}((M,F)\otimes_{D^{\infty}_{Z_2}}(D^{\infty}_{Z_2\leftarrow Z_1},F^{ord})
\otimes_{i_0^*D^{\infty}_{Z_1}}i_0^*(D^{\infty}_{Z_1\leftarrow S},F^{ord})) 
\xrightarrow{i_{1*}T(i_0,\otimes)(-)^{-1}} i_{1*mod}i_{0*mod}((M,F))
\end{eqnarray*}
using lemma \ref{Dcompinfty} and proposition \ref{Tiotimes}.

\end{proof}

\begin{prop}\label{compAnDmodh}
\begin{itemize}
\item[(i)] Let $f:X\to S$ a morphism with $X,S\in\AnSm(\mathbb C)$.
For $(M,F)\in C_{\mathcal D(2)fil,h}(S)$, we have $Lf^{*mod}(M,F)\in D_{\mathcal D(2)fil,\infty,h}(X)$.
For $(M,F)\in C_{\mathcal D^{\infty}(2)fil,h}(S)$, we have $Lf^{*mod}(M,F)\in D_{\mathcal D^{\infty}(2)fil,\infty,h}(X)$.
\item[(ii)] Let $f:X\to S$ a proper morphism with $X,S\in\AnSm(\mathbb C)$. Then,
for $(M,F)\in C_{\mathcal D(2)fil,h}(X)$, we have $\int_f(M,F)\in D_{\mathcal D(2)fil,\infty,h}(S)$.
\item[(iii)] Let $f:X\to S$ a morphism with $X,S\in\AnSm(\mathbb C)$. Then,
for $(M,F)\in C_{\mathcal D^{\infty}(2)fil,h}(X)$, we have $\int_f(M,F)\in D_{\mathcal D^{\infty}(2)fil,\infty,h}(S)$.
\end{itemize}
\end{prop}

\begin{proof}
\noindent(i)and (ii):Follows imediately from the non filtered case since we look at the complex in the derived category with respect to
$\infty$-usu local equivalence. It says that the pullback and the proper pushforward of an holonomic D module is still holonomic.
See \cite{LvDmod} for the non filtered case. 

\noindent(iii):In the case the morphism is proper, it follows from the finite order case (ii). 
In the case of an open embedding, it follows from proposition \ref{Jprop}(i) : we have for $j:S^o\hookrightarrow S$ an open embedding,
\begin{equation*}
j_*E(O_{S^o})=j_*\mathcal Hom(\mathbb Z_{S^o},E(O_{S^o}))=\mathcal Hom(j_!\mathbb Z_{S^o},E(O_S))\in C_{\mathcal D^{\infty},h}(S).
\end{equation*}
and on the other hand
\begin{eqnarray*}
T(j,\otimes)(-,-)=T^{mod}(j,\otimes)(-,-):\int_j(M,F)=j_*E(M,F)=j_*E(j^*O_S\otimes_{O_{S^o}}(M,F)) \\
\xrightarrow{\sim}j_*E(O_{S^o})\otimes_{O_S}(M,F)
\end{eqnarray*}
is an isomorphism by proposition \ref{projformula}.
\end{proof}

For $X,Y\in\SmVar(\mathbb C)$ or $X,Y\in\AnSm(\mathbb C)$, we denote by
\begin{itemize}
\item $C_{O_X}(X)\times C_{O_Y}(Y)\to C_{O_{X\times Y}}(X\times Y),
(M,N)\mapsto M\cdot N:=O_{X\times Y}\otimes_{p_X^*O_X\otimes p_Y^*O_Y}p_X^*M\otimes p_Y^*N$,
\item $C_{\mathcal D}(X)\times C_{\mathcal D}(Y)\to C_{\mathcal D}(X\times Y),
(M,N)\mapsto M\cdot N:=O_{X\times Y}\otimes_{p_X^*O_X\otimes p_Y^*O_Y}p_X^*M\otimes p_Y^*N$
\end{itemize}
the natural functors which induces functors in the filtered cases and the derived categories,
$p_X:X\times Y\to X$ and $p_Y:X\times Y\to Y$ the projections.

We have then the following easy proposition :

\begin{prop}\label{otimesed}
For $X\in\SmVar(\mathbb C)$ or $X\in\AnSm(\mathbb C)$, 
we have for $(M,F),(N,F)\in C_{O_X,fil}(X)$ or $(M,F),(N,F)\in C_{\mathcal D,fil}(X)$,
\begin{equation*}
(M,F)\otimes_{O_X}(N,F)=\Delta_X^{*mod}(M,F)\cdot (N,F)
\end{equation*}
\end{prop}

\begin{proof}
Standard.
\end{proof}

\begin{defi}\label{Tfinfty}
Let $f:X\to S$ a morphism with $X,S\in\AnSm(\mathbb C)$. 
We have the canonical map in $C_{f^*\mathcal D,\mathcal D^{\infty}}(X)$ modules :
\begin{equation*}
T(f,\infty):(D_{X\to S},F^{ord})\otimes_{D_X}(D_X^{\infty},F^{ord})\to (D^{\infty}_{X\to S},F^{ord}), \; 
(h_X\otimes P_S)\otimes P_X\mapsto (P_X.h_X\otimes P_S+h_X\otimes df(P_X)P_S
\end{equation*}
where $h_X\in\Gamma(X^o,O_X)$, $P_S\in\Gamma(X^o,f^*D_S)$ and $P_X\in\Gamma(X^o,D_X^{\infty})$.
This gives, for $(M,F)\in C_{\mathcal D(2)fil}(S)$,  the following transformation map in $C_{\mathcal D^{\infty}(2)fil}(X)$
\begin{eqnarray*}
T(f,\infty)(M,F):J_X(f^{*mod}(M,F)):=f^*(M,F)\otimes_{f^*D_S}(D_{X\to S},F^{ord}\otimes_{D_X}(D_X^{\infty},F^{ord})
\xrightarrow{I\otimes T(f,\infty)} \\ f^*(M,F)\otimes_{f^*D_S}(D^{\infty}_{X\to S},F^{ord})=
f^*(M,F)\otimes_{f^*D_S}f^*D^{\infty}_S\otimes_{f^*D_S^{\infty}}(D^{\infty}_{X\to S},F^{ord})=:f^{*mod}J_S(M,F)
\end{eqnarray*}
where we recall that $J_S(M,F)=(M,F)\otimes_{D_S}(D_S^{\infty},F^{ord})$.
\end{defi}

We now look at some properies of the dual functor for D modules :
For complex of $D$ module with coherent cohomology we have the following: 

\begin{prop}\label{dmodinv}
\begin{itemize}
\item[(i)] Let $S\in\SmVar(\mathbb C)$. For $M\in C_{\mathcal D,c}(S)$, 
the canonical map $d(M):M\to\mathbb D^2_SL_DM$ is an equivalence Zariski local.
\item[(ii)] Let $S\in\AnSm(\mathbb C)$. For $M\in C_{\mathcal D}(S)$, 
the canonical map $d(M):M\to\mathbb D^2_SL_D(M)$ is an equivalence usu local.
\item[iii)] Let $S\in\AnSm(\mathbb C)$. For $(M,F)\in C_{\mathcal D^{\infty}}(S)$,
the canonical map $d(M):M\to\mathbb D^2_SL_{D^{\infty}}(M)$ is an equivalence usu local.
\end{itemize}
\end{prop}

\begin{proof}
Standard :follows from the definition of coherent sheaves. See \cite{LvDmod} for exemple.
\end{proof}

Let $S_1,S_2\in\SmVar(\mathbb C)$ or $S_1,S_2\in\AnSm(\mathbb C)$ and $p:S_{12}:=S_1\times S_2\to S_1$ the projection.
In this case we have a canonical embedding $D_{S_1}\hookrightarrow p_*D_{S_{12}}$.
This gives, for $(M,F)\in C_{\mathcal Dfil}(S_1\times S_2)$, the following transformation map in $C_{\mathcal Dfil}(S_1)$ 
\begin{eqnarray*}
T_*(p,D)(M,F):p_*\mathbb D_{S_{12}}^K(M,F):=
p_*\mathcal Hom_{D_{S_{12}}}((M,F),D_{S_{12}})\otimes_{O_{S_{12}}}\mathbb D_{S_{12}}^Ow(K_{S_{12}})[d_{S_{12}}] \\
\xrightarrow{T_*(p,hom)(-,-)}
\mathcal Hom_{p_*D_{S_{12}}}(p_*(M,F),p_*D_{S_{12}})\otimes_{p_*O_{S_{12}}}\mathbb D_{p_*S_{12}}^Ow(p_*K_{S_{12}})[d_{S_{12}}] \\
\xrightarrow{\sim}
\mathcal Hom_{D_{S_1}}(p_*(M,F),D_{S_1})\otimes_{O_{S_1}}\mathbb D_{S_1}^Ow(K_{S_1})[d_{S_1}]=:\mathbb D_{S_1}^Kp_*(M,F)
\end{eqnarray*}
We have the canonical map 
\begin{eqnarray*}
p(D):p^{*mod}D_{S_1}=p^*D_{S_1}\otimes_{p^*O_{S_1}}O_{S_{12}}\to D_{S_{12}}, \; \gamma\otimes f\mapsto f.\gamma
\end{eqnarray*}
induced by the embedding $p^*D_{S_1}\hookrightarrow D_{S_{12}}$.
This gives, for $(M,F)\in C_{\mathcal Dfil}(S_1)$, the following transformation map in $C_{\mathcal Dfil}(S_1\times S_2)$ 
\begin{eqnarray*}
T(p,D)(M,F):p^{*mod}\mathbb D_{S_1}^K(M,F):=
p^*\mathcal Hom_{D_{S_1}}((M,F),D_{S_1})\otimes_{p^*O_{S_1}}p^{*mod}\mathbb D_{S_1}^Ow(K_{S_1})[d_{S_1}] \\
\xrightarrow{T(p,hom)(-,-)\otimes I}
\mathcal Hom_{p^*D_{S_1}}(p^*(M,F),p^*D_{S_1})\otimes_{p^*O_{S_1}}p^{*mod}\mathbb D_{S_1}^OwK_{S_1})[d_{S_1}] \\
\xrightarrow{(\phi\mapsto\phi\otimes I_{O_{S_{12}}})\otimes I}
\mathcal Hom_{D_{S_{12}}}(p^{*mod}(M,F),p^{*mod}D_{S_1})\otimes_{p^*O_{S_1}}p^{*mod}\mathbb D_{S_1}^Ow(K_{S_1})[d_{S_1}] \\
\xrightarrow{I\otimes K^{-1}(S_1/S_{12})} 
\mathcal Hom_{D_{S_{12}}}(p^{*mod}(M,F),p^{*mod}D_{S_1})\otimes_{p^*O_{S_1}}\mathbb D_{S_{12}}^Ow(K_{S_{12}})[d_{S_{12}}] \\
\xrightarrow{q(p^*O_{S_1}/O_{S_{12}})}
\mathcal Hom_{D_{S_{12}}}(p^{*mod}(M,F),p^{*mod}D_{S_1})\otimes_{O_{S_{12}}}\mathbb D_{S_{12}}^Ow(K_{S_{12}})[d_{S_{12}}] \\
\xrightarrow{\mathcal Hom(p^{*mod}(M,F),p(D))}
\mathcal Hom_{D_{S_{12}}}(p^{*mod}(M,F),D_{S_{12}})\otimes_{O_{S_{12}}}\mathbb D_{S_{12}}^Ow(K_{S_{12}})[d_{S_{12}}]=:
\mathbb D^K_{S_{12}}(p^{*mod}(M,F))
\end{eqnarray*} 
whre $K^{-1}(S_1/S_{12})$ is given by the wedge product with a generator of 
$\wedge^{d_{S_2}}T_{S_{12}/S_1}\xrightarrow{\sim}K^{-1}_{S_2}$.

In the case $S_1,S_2\in\AnSm(\mathbb C)$, 
we also have the embedding $p^*D^{\infty}_{S_1}\hookrightarrow D^{\infty}_{S_{12}}$.
This gives in the same way, for $(M,F)\in C_{\mathcal Dfil}(S_1\times S_2)$, 
the following transformation map in $C_{\mathcal Dfil}(S_1)$ 
\begin{eqnarray*}
T_*(p,D^{\infty})(M,F):p_*\mathbb D^{\infty,K}_{S_{12}}(M,F)\to\mathbb D^{\infty,K}_{S_1}p_*(M,F).
\end{eqnarray*}
The map   
\begin{eqnarray*}
p(D^{\infty}):p^{*mod}D_{S_1}=p^*D^{\infty}_{S_1}\otimes_{p^*O_{S_1}}O_{S_{12}}\to D^{\infty}_{S_{12}}, \; 
\gamma\otimes f\mapsto f.\gamma
\end{eqnarray*}
induced by the embedding $p^*D^{\infty}_{S_1}\hookrightarrow D^{\infty}_{S_{12}}$, 
gives in the same way, for $(M,F)\in C_{\mathcal D^{\infty}fil}(S_1)$, 
the transformation map in $C_{\mathcal D^{\infty}fil}(S_1\times S_2)$ 
\begin{eqnarray*}
T(p,D^{\infty})(M,F):p^{*mod}\mathbb D_{S_1}^{\infty,K}(M,F):=
p^{*mod}(\mathbb D^{\infty}_{S_1}(M,F)\otimes_{O_{S_1}}\mathbb D^O_{S_1}w(K_{S_1})[d_{S_1}])\to \\
\mathcal Hom_{D^{\infty}_{S_{12}}}(p^{*mod}(M,F),
D^{\infty}_{S_{12}})\otimes_{O_{S_{12}}}\mathbb D_{S_{12}}^Ow(K_{S_{12}}[d_{S_{12}}]=:
\mathbb D^{\infty,K}_{S_{12}}(p^{*mod}(M,F)),
\end{eqnarray*} 
given in the same way then $T(p,D)(M,F)$.

\begin{prop}\label{dmodg}
\begin{itemize}
\item[(i)] Let $g:T\to S$ a morphism with $T,S\in\SmVar(\mathbb C)$.
We have, for $M\in D_{\mathcal D}(S)$ canonical maps
\begin{itemize}
\item $T'(g,D)(M):L\mathbb D_SLg^{*mod}M\to Lg^{*mod}L\mathbb D_SM$
\item $T'(g,D)(M):Lg^{*mod}L\mathbb D_SM\to L\mathbb D_SLg^{*mod}M$
\end{itemize}
Moreover, in the case where $g$ is non caracteristic with respect to $M$ (e.g if $g$ is smooth), these maps are isomorphism.
\item[(ii)] Let $S_1,S_2\in\SmVar(\mathbb C)$, $p:S_1\times S_2\to S_1$ the projection. For $M\in D_{\mathcal D}(S_1)$,
we have $T(p,D)(L_D(M))=T'(p,D)(M)$ in $D_{\mathcal Dfil}(S_1\times S_2)$ (c.f.(i)).
\end{itemize}
\end{prop}

\begin{proof}
\noindent(i):See \cite{LvDmod} for the first map. The second one follows from the first by proposition \ref{dmodinv}(i) and (iii).

\noindent(ii):See the proof of (i) in \cite{LvDmod}
\end{proof}

We have the followings :

\begin{prop}\label{compDmodhat}
Let $f_1:X\to Y$ and $f_2:Y\to S$ two morphism with $X,Y,S\in\SmVar(\mathbb C)$. 
Let $M\in C_{\mathcal D,h}(S)$. 
Then, we have $L(f_2\circ f_1)^{\hat{*}mod}M=Lf_1^{\hat{*}mod}(Lf_2^{\hat{*}mod}M)$ in $D_{\mathcal D,h}(X)$.
\end{prop}

\begin{proof}
Follows from proposition \ref{compDmod} (i), proposition \ref{compDmodh} and proposition \ref{dmodinv}.
\end{proof}

\begin{prop}\label{compAnDmodhat}
Let $f_1:X\to Y$ and $f_2:Y\to S$ two morphism with $X,Y,S\in\AnSm(\mathbb C)$. 
\begin{itemize}
\item[(i)]Let $(M)\in C_{\mathcal D,h}(S)$. 
Then, we have $L(f_2\circ f_1)^{\hat{*}mod}(M)=Lf_1^{\hat{*}mod}(Lf_2^{\hat{*}mod}(M))$ in $D_{\mathcal D,h}(X)$.
\item[(ii)]Let $M\in C_{\mathcal D^{\infty},h}(S)$. 
Then, we have $L(f_2\circ f_1)^{\hat{*}mod}M=Lf_1^{\hat{*}mod}(Lf_2^{\hat{*}mod}M)$ in $D_{\mathcal D^{\infty},h}(X)$.
\end{itemize}
\end{prop}

\begin{proof}
Follows from proposition \ref{compAnDmod} (i), proposition \ref{compAnDmodh} and proposition \ref{dmodinv}.
\end{proof}

In the analytic case, we have the following transformation map which we will use in subection 5.3:

\begin{defi}\label{TDinfty}
Let $S\in\AnSm(\mathbb C)$.
We have for $(M,F)\in C_{\mathcal Dfil}(S)$ the canonical transformation map in $C_{\mathcal D^{\infty}fil}(S)$ :
\begin{eqnarray*}
T(D,\infty)(M,F): \\
J_S(\mathbb D^K_S(M,F)):=\mathcal Hom_{D_S}((M,F),D_S)\otimes_{D_S}(D^{\infty}_S,F^{ord})\otimes_{O_S}\mathbb D^O_Sw(K_S)[d_S] 
\xrightarrow{ev_{D_S}(\hom,\otimes)(-,-,-)\otimes I} \\
\mathcal Hom_{D_S}(L_D(M,F),D^{\infty}_S)\otimes_{O_S}\mathbb D^O_Sw(K_S)[d_S] 
\xrightarrow{I(D^{\infty}_S/D_S)((M,F),D^{\infty}_S)\otimes I} \\
\mathcal Hom_{D^{\infty}_S}((M,F)\otimes_{D_S}(D^{\infty}_S,F^{ord}),D^{\infty}_S)\otimes_{O_S}\mathbb D^O_Sw(K_S)[d_S]
=:\mathbb D^{\infty,K}_SJ_S(M,F).
\end{eqnarray*}
which is an isomorphism.
\end{defi}

\subsubsection{The (relative) De Rahm of a (filtered) complex of a D-module and the filtered De Rham direct image}

Recall that for $f:X\to S$ a morphism with $X,S\in\Var(\mathbb C)$ or with $X,S\in\AnSp(\mathbb C)$,
\begin{equation*}
DR(X/S):=\Omega^{\bullet}_{X/S}\in C_{f^*O_S}(X) 
\end{equation*}
denotes (see section 2) the relative De Rham complex of the morphism of ringed spaces $f:(X,O_X)\to(S,O_S)$,
with $\Omega^{p}_{X/S}:=\wedge^p\Omega_{X/S}\in\PSh_{O_X}(X)$ and $\Omega_{X/S}:=\coker(f^*\Omega_S\to\Omega_X)\in\PSh_{O_X}(X)$. 
Recall that $\Omega^{\bullet}_{X/S}\in C_{f^*O_S}(S)$ is a complex of $f^*O_S$ modules, 
but is NOT a complex of $O_X$ module since the differential is a derivation hence NOT $O_X$ linear.
Recall that (see section 4.1), for $(M,F)\in C_{D(O_X)fil}(X)$, we have the relative (filtered) De rham complex of $(M,F)$
\begin{equation*}
DR(X/S)(M,F):=(\Omega^{\bullet}_{X/S},F_b)\otimes_{O_X}(M,F)\in C_{f^*O_Sfil}(X),
\end{equation*}
and that if $\phi:(M_1,F)\to(M_2,F)$ a morphism with $(M_1,F),(M_2,F)\in C_{D(O_X)fil}(X)$,
\begin{equation*}
(I\otimes\phi):DR(X/S)(M_1,F):=(\Omega^{\bullet}_{X/S},F)\otimes_{O_X}(M_1,F)\to DR(X/S)(M_2,F):=
(\Omega^{\bullet}_{X/S},F)\otimes_{O_X}(M_2,F)
\end{equation*}
is by definition a morphism of complexes, that is a morphism in $C_{f^*O_Sfil}(X)$.
For $(N,F)\in C_{D(O_X)^{op}fil}(X)$, we have the relative (filtered) Spencer complex of $(N,F)$
\begin{equation*}
SP(X/S)(N,F):=(T^{\bullet}_{X/S},F_b)\otimes_{O_X}(N,F)\in C_{f^*O_Sfil}(X),
\end{equation*}
and that if $\phi:(N_1,F)\to(N_2,F)$ a morphism with $(N_1,F),(N_2,F)\in C_{D(O_X)^{op}fil}(X)$,
\begin{equation*}
(I\otimes\phi):SP(X/S)(N_1,F):=(T^{\bullet}_{X/S},F)\otimes_{O_X}(N_1,F)\to SP(X/S)(N_2,F):=
(T^{\bullet}_{X/S},F)\otimes_{O_X}(N_2,F)
\end{equation*}
is by definition a morphism of complexes, that is a morphism in $C_{f^*O_Sfil}(X)$.
\begin{prop}
Let $f:X\to S$ a smooth morphism with $X,S\in\Var(\mathbb C)$ or with $X,S\in\AnSp(\mathbb C)$, denote $d=d_X-d_S$.
The inner product gives, for $(M,F)\in C_{D(O_X)fil}(X)$, an isomorphism in $C_{f^*O_Sfil}(X)$ and termwise $O_X$ linear
\begin{equation*}
T(DR,SP)(M,F):T^{\bullet}_{X/S}\otimes_{O_X}(M,F)\otimes_{O_X}K_{X/S}\xrightarrow{\sim}
\Omega^{d-\bullet}_{X/S}\otimes_{O_X}(M,F), \partial\otimes m\otimes\kappa\mapsto\iota(\partial)\kappa\otimes m
\end{equation*}
\end{prop}

\begin{proof}
Standard.
\end{proof}

For a commutative diagram in $\Var(\mathbb C)$ or in $\AnSp(\mathbb C)$ :
\begin{equation*}
D=\xymatrix{X\ar[r]^{f} & S \\
X'\ar[u]^{g'}\ar[r]^{f'} & T\ar[u]^{g}}
\end{equation*}
we have (see section 2) the relative differential map of $g'$ given by the pullback of differential forms:
\begin{eqnarray*}
\Omega_{(X'/X)/(T/S)}:g^{'*}\Omega_{X/S}\to\Omega_{X'/T}, \; \mbox{given by for} \; 
X^{'o}\subset X', \; X^o\supset g'(X^{'o}) (\mbox{i.e.} \, g^{'-1}(X^o)\supset X^{'o}), \\ 
\omega\in\Gamma(X^o,\Omega^p_{X/S})\mapsto \Omega_{(X'/X)/(T/S)}(X^{'o})(\omega):=[g^{'*}\omega]\in\Gamma(X^{'o},\Omega^p_{X'/T}).
\end{eqnarray*} 
Moreover, by definition-proposition \ref{TDwM} (section 4.1), for $(M,F)\in C_{D(O_X)fil}(X)$ the map
\begin{equation*}
\Omega_{(X'/X)/(T/S)}(M,F):g'^*(\Omega^{\bullet}_{X/S}\otimes_{O_X}(M,F))\to\Omega^{\bullet}_{X'/T}\otimes_{O_{X'}}g^{'*mod}(M,F)
\end{equation*}
given in degree $(p,i)$ by, for $X^{'o}\subset X'$ an open subset and $X^o\subset X$ an open subset such that
$g^{'-1}(X^o)\supset X^{'o}$ (i.e. $X^o\supset g'(X^{'o})$), $\omega\in\Gamma(X^o,\Omega^p_{X/S})$ and $m\in\Gamma(X^o,M^i)$,
\begin{equation*}
\Omega_{(X'/X)/(T/S)}(M,F)(\omega\otimes m)=g^{'*}\omega\otimes(m\otimes 1)
\end{equation*}
is a map of complexes, that is a map in $C_{f^*O_Sfil}(X')$.
This give, for $(M,F)\in C_{D(O_X)fil}(X)$, the following transformation map in $C_{O_Tfil}(T)$
\begin{eqnarray*}
T^O_{\omega}(D)(M,F):g^{*mod}L_O(f_*E(\Omega^{\bullet}_{X/S}\otimes_{O_X}(M,F)))
\xrightarrow{T(g,L_O)(-)} \\
(g^*f_*E(\Omega^{\bullet}_{X/S}\otimes_{O_X}(M,F)))\otimes_{g^*O_S}O_T  
\xrightarrow{T(g',E)(-)\circ T(D)(E(\Omega^{\bullet}_{X/S}\otimes_{O_X}M))} \\
(f'_*E(g'^*(\Omega^{\bullet}_{X/S}\otimes_{O_X}(M,F))))\otimes_{g^*O_S}O_T 
\xrightarrow{m\circ E(\Omega_{(Y/X)/(T/S)}(M))}
f'_*E(\Omega^{\bullet}_{X'/T}\otimes_{O_{X'}}g^{'*mod}(M,F)),
\end{eqnarray*}
with $m'(m)=m\otimes 1$. Under the canonical isomorphism 
$\Omega^{\bullet}_{X/S}\xrightarrow{\sim}\Omega^{\bullet}_{X/S}\otimes_{O_X}O_X$
given by $\omega\mapsto\omega\otimes 1$, we have (see remark \ref{TDwMrem}) 
\begin{equation*}
T^O_{\omega}(D)(O_X)=T^O_{\omega}(D):g^{*mod}L_O(f_*E(\Omega^{\bullet}_{X/S}))\to f'_*E(\Omega^{\bullet}_{X'/T}).
\end{equation*}

Let $f:X\to S$ a morphism with $X,S\in\Var(\mathbb C)$. 
Again by definition-proposition \ref{TDwM} (section 4.1), for $(M,F)\in C_{D(O_X)fil}(X)$ the map
\begin{equation*}
\Omega_{(X^{an}/X)/(S^{an}/S)}(M,F):an_X^*(\Omega^{\bullet}_{X/S}\otimes_{O_X}(M,F))\to
\Omega^{\bullet}_{X^{an}/S^{an}}\otimes_{O_{X^{an}}}M^{an}
\end{equation*}
given in degree $(p,i)$ by, for $X^o\subset X$ and $X^o\supset X^{oo}$ an open subsets of $X$ for the usual, resp. Zariski topology,
$\omega\in\Gamma(X^o,\Omega^p_{X/S})$ and $m\in\Gamma(X^o,M^i)$,
\begin{equation*}
\Omega_{(X^{an}/X)/(S^{an}/S)}(M,F)(\omega\otimes m=\omega\otimes(m\otimes 1)
\end{equation*}
is a map of complexes, that is a map in $C_{f^*O_{S^{an}}fil}(X^{an})$.
This gives, for $(M,F)\in C_{D(O_X)fil}(X)$, we have the  following transformation map in $C_{O_{S^{an}}fil}(S^{an})$
\begin{eqnarray*}
T^O_{\omega}(an,f)(M,F):(f_*E(\Omega^{\bullet}_{X/S}\otimes_{O_X}(M,F)))^{an}
:=\an_S^*(f_*E(\Omega^{\bullet}_{X/S}\otimes_{O_X}(M,F)))\otimes_{\an_S^*O_S}O_{S^{an}} \\
\xrightarrow{T(an(X),E)(-)\circ T(an,f)(E(\Omega^{\bullet}_{X/S}\otimes_{O_X}M))} 
(f_*E(\an_X^*(\Omega^{\bullet}_{X/S}\otimes_{O_X}(M,F))))\otimes_{\an_S^*O_S}O_{S^{an}} \\
\xrightarrow{m\circ E(\Omega_{(X^{an}/X)/(S^{an}/S)}(M,F))}
f_*E(\Omega^{\bullet}_{X/S}\otimes_{O_X^{an}}(M^{an},F))
\end{eqnarray*}
with $m(n\otimes s)=s.n$.Under the canonical isomorphism $\Omega^{\bullet}_{X/S}\xrightarrow{\sim}\Omega^{\bullet}_{X/S}\otimes_{O_X}O_X$
given by $\omega\mapsto\omega\otimes 1$, we have (see remark \ref{TDwMrem}) 
\begin{equation*}
T^O_{\omega}(an,f)(O_X)=T^O_{\omega}(an,f):(f_*E(\Omega^{\bullet}_{X/S}))^{an}\to f_*E(\Omega^{\bullet}_{X^{an}/S^{an}})
\end{equation*}

Let $f:X\to S$ a morphism with $X,S\in\Var(\mathbb C)$ or with $X,S\in\AnSp(\mathbb C)$.
In the case where $X$ is smooth, for $(M,F)=(M^{\bullet},F)\in C_{\mathcal Dfil}(X)$, 
the differential of the relative De Rham complex of $(M,F)$
\begin{equation*}
DR(X/S)(M,F):=(\Omega^{\bullet}_{X/S},F)\otimes_{O_X}(M,F)=
\Tot((\Omega^{\bullet}_{X/S},F)\otimes_{O_X}(M^{\bullet},F))\in C_{f^*O_Sfil}(X)
\end{equation*}
are given by
\begin{itemize}
\item $d_{p,p+1}:\Omega^p_{X/S}\otimes_{O_X} M^i\to\Omega^{p+1}_{X/S}\otimes_{O_X} M^i$, 
with for $X^o\subset X$ an open affine subset with $(x_1,\ldots, x_n)$ local coordinate (since $X$ is smooth, $T_X$ is locally free),
$m\in\Gamma(X^o,M^i)$ and $\omega\in\Gamma(X^o,\Omega^p_{X/S})$, 
\begin{equation*}
d_{p,p+1}(\omega\otimes m):=(d\omega)\otimes m+\sum_{i=1}^{n}(dx_i\wedge\omega)\otimes(\partial_i)m
\end{equation*}
\item $d_{i,i+1}:\Omega^p_{X/S}\otimes_{O_X} M^i\to\Omega^p_{X/S}\otimes_{O_X} M^{i+1}$,
with for $X^o\subset X$ an open subset, $m\in\Gamma(X^o,M^i)$ and $\omega\in\Gamma(X^o,\Omega^p_{X/S})$, 
$d_{i,i+1}(\omega\otimes m):=(\omega\otimes dm)$.
\end{itemize}
For $D_X$ only, the differential of its De Rahm complex $(\Omega^{\bullet}_{X/S},F)\otimes_{O_X}D_X$ are right linear, so that 
\begin{equation*}
(\Omega^{\bullet}_{X/S},F_b)\otimes_{O_X}(D_X,F^{ord})\in C_{\mathcal D^{op},f^*O_Sfil}(X)
\end{equation*}

In the particular case of a projection $p:Y\times S\to S$ with $Y,S\in\SmVar(\mathbb C)$ or with $Y,S\in\AnSm(\mathbb C)$ we have :

\begin{prop}\label{DRhUS}
Let $Y,S\in\SmVar(\mathbb C)$ or $Y,S\in\AnSm(\mathbb C)$. Let $p:Y\times S\to S$ the projection.
For $(M,F)\in C_{\mathcal Dfil}(Y\times S)$, 
\begin{equation*}
DR(Y\times S/S)(M,F):=(\Omega^{\bullet}_{Y\times S/S},F_b)\otimes_{O_{Y\times S}}(M,F)\in C_{p^*O_Sfil}(Y\times S)
\end{equation*}
is a naturally a complex of filtered $p^*D_S$ modules, that is 
\begin{equation*}
DR(Y\times S/S)(M,F):=(\Omega^{\bullet}_{Y\times S/S},F_b)\otimes_{O_{Y\times S}}(M,F)\in C_{p^*\mathcal Dfil}(Y\times S),
\end{equation*}
where the $p^*D_S$ module structure on $\Omega^p_{Y\times S/S}\otimes_{O_{Y\times S}}M^n$ is given by
for $(Y\times S)^o\subset Y\times S$ an open subset,
\begin{equation*}
(\gamma\in\Gamma((Y\times S)^o,T_{Y\times S}),
\hat\omega\otimes m\in\Gamma((Y\times S)^o,\Omega^p_{Y\times S/S}\otimes_{O_{Y\times S}}M^n))\mapsto 
\gamma.(\hat\omega\otimes m):=(\hat\omega\otimes(\gamma.m).
\end{equation*}
Moreover, if $\phi:(M_1,F)\to(M_2,F)$ a morphism with $(M_1,F),(M_2,F)\in C_{\mathcal Dfil}(Y\times S)$,
\begin{equation*}
DR(Y\times S/S)(\phi):=(I\otimes\phi):(\Omega^{\bullet}_{Y\times S/S},F_b)\otimes_{O_{Y\times S}}(M_1,F)
\to(\Omega^{\bullet}_{Y\times S/S},F_b)\otimes_{O_{Y\times S}}(M_2,F)
\end{equation*}
is a morphism in $C_{p^*\mathcal Dfil}(Y\times S)$.
\end{prop}

\begin{proof}
Standard.
\end{proof}

In the analytic case, we also have
\begin{prop}\label{DRhUSinfty}
Let $Y,S\in\AnSm(\mathbb C)$. Let $p:Y\times S\to S$ the projection.
For $(M,F)\in C_{\mathcal D^{\infty}fil}(Y\times S)$, 
\begin{equation*}
DR(Y\times S/S)(M,F):=(\Omega^{\bullet}_{Y\times S/S},F_b)\otimes_{O_{Y\times S}}(M,F)\in C_{p^*O_Sfil}(Y\times S)
\end{equation*}
is naturally a complex of filtered $p^*D^{\infty}_S$ modules, that is 
\begin{equation*}
DR(Y\times S/S)(M,F):=(\Omega^{\bullet}_{Y\times S/S},F_b)\otimes_{O_{Y\times S}}(M,F)\in C_{p^*\mathcal D^{\infty}fil}(Y\times S),
\end{equation*}
where the $p^*D^{\infty}_S$ module structure on $\Omega^p_{Y\times S/S}\otimes_{O_{Y\times S}}M^n$ is given by
for $(Y\times S)^o\subset Y\times S$ an open subset,
\begin{equation*}
(\gamma\in\Gamma((Y\times S)^o,T_{Y\times S}),
\hat\omega\otimes m\in\Gamma((Y\times S)^o,\Omega^p_{Y\times S/S}\otimes_{O_{Y\times S}}M^n))\mapsto 
\gamma.(\hat\omega\otimes m):=(\hat\omega\otimes(\gamma.m).
\end{equation*}
Moreover, if $\phi:(M_1,F)\to(M_2,F)$ a morphism with $(M_1,F),(M_2,F)\in C_{\mathcal D^{\infty}fil}(Y\times S)$,
\begin{equation*}
DR(Y\times S/S)(\phi):=(I\otimes\phi):(\Omega^{\bullet}_{Y\times S/S},F_b)\otimes_{O_{Y\times S}}(M_1,F)
\to(\Omega^{\bullet}_{Y\times S/S},F_b)\otimes_{O_{Y\times S}}(M_2,F)
\end{equation*}
is a morphism in $C_{p^*\mathcal D^{\infty}fil}(Y\times S)$.
\end{prop}

\begin{proof}
Standard : follows from the finite order case (proposition \ref{DRhUS}).
\end{proof}

We state on the one hand the commutativity of the tensor product with respect to $D_S$ and with respect to $O_S$, for
$S\in\SmVar(\mathbb C)$ or $S\in\AnSm(\mathbb C)$ in the filtered case, 
and on the other hand the commutativity between the tensor product with respect to $D_S$ by $D_S$ and the De Rahm complex :

\begin{prop}\label{otimesDmod}
\begin{itemize}
\item[(i)] Let $f:X\to S$ a morphism with $X,S\in\SmVar(\mathbb C)$ or with $X,S\in\AnSm(\mathbb C)$.
For $(M',F)\in C_{\mathcal D^{op}fil,f^*\mathcal D}(X)$ and $(M,F),(N,F)\in C_{\mathcal Dfil}(X)$.
we have canonical isomorphisms of filtered $f^*D_S$ modules, i.e. isomorphisms in $C_{f^*\mathcal D}(X)$,
\begin{eqnarray*}
(M',F)\otimes_{O_X}(N,F)\otimes_{D_X}(M,F)&=&(M',F)\otimes_{D_X}((M,F)\otimes_{O_X}(N,F)) \\
&=&((M',F)\otimes_{O_X}(M,F))\otimes_{D_X}(N,F)
\end{eqnarray*}
\item[(ii)] Let $f:X\to S$ a morphism with $X,S\in\Var(\mathbb C)$ or with $X,S\in\AnSp(\mathbb C)$.
For $(M,F)\in C_{D(O_X)fil}(X)$, we have a canonical isomorphisms of filtered $f^*O_S$ modules, i.e. isomorphisms in $C_{f^*O_Sfil}(X)$,
\begin{eqnarray*}
(\Omega^{\bullet}_{X/S},F_b)\otimes_{O_X}(M,F)=((\Omega^{\bullet}_{X/S},F_b)\otimes_{O_X}(D(O_X),F_b))\otimes_{D(O_X)}(M,F)
\end{eqnarray*}
\item[(iii)] Let $p:Y\times S\to S$ a projection with $Y,S\in\SmVar(\mathbb C)$ or with $X,S\in\AnSm(\mathbb C)$.
For $(M,F)\in C_{\mathcal Dfil}(Y\times S)$, the isomorphisms of filtered $p^*O_S$ modules of (ii)
\begin{eqnarray*}
(\Omega^{\bullet}_{Y\times S/S},F_b)\otimes_{O_{Y\times S}}(M,F)=
((\Omega^{\bullet}_{Y\times S/S},F_b)\otimes_{O_{Y\times S}}(D_{Y\times S},F_b)\otimes_{D_{Y\times S}}(M,F)
\end{eqnarray*}
are isomorphisms of filtered $p^*D_S$ modules, that is isomorphism in $C_{p^*\mathcal Dfil}(Y\times S)$.
\end{itemize}
\end{prop}

\begin{proof}
\noindent(i) and (ii) are particular case of proposition \ref{otimesDmodsection2}.

\noindent(iii): follows immediately by definition of the $p^*D_S$ module structure.  
\end{proof}

We now look at the functorialities of the relative De Rham complex of a smooth morphisms of smooth complex algebraic varieties :

\begin{prop}\label{TDhwM}
Consider a commutative diagram in $\SmVar(\mathbb C)$ or in $\AnSm(\mathbb C)$ :
\begin{equation*}
D=\xymatrix{Y\times S\ar[r]^{p} & S \\
Y'\times T\ar[u]^{g''=(g''_0\times g)}\ar[r]^{p'} & T\ar[u]^{g}}
\end{equation*}
with $p$ and $p'$ the projections. 
For $(M,F)\in C_{\mathcal Dfil}(Y\times S)$ the map in $C_{g^{''*}p^*O_Sfil}(Y'\times T)$
\begin{equation*}
\Omega_{(Y'\times T/Y\times S)/(T/S)}(M,F):g''^*((\Omega^{\bullet}_{Y\times S/S},F_b)\otimes_{O_{Y\times S}}(M,F))
\to(\Omega^{\bullet}_{Y'\times T/T},F_b)\otimes_{O_{Y'\times T}}g^{''*mod}(M,F)
\end{equation*}
given in definition-proposition \ref{TDwM} is a map in $C_{g^{''*}p^*\mathcal Dfil}(Y'\times T)$. 
Hence, for $(M,F)\in C_{\mathcal Dfil}(Y\times S)$, the map in $C_{O_Tfil}(T)$ (with $L_D$ instead of $L_O$)
\begin{eqnarray*}
T^O_{\omega}(D)(M):g^{*mod}L_D(p_*E((\Omega^{\bullet}_{Y\times S/S},F_b)\otimes_{O_{Y\times S}}(M,F)))\to
p'_*E((\Omega^{\bullet}_{Y'\times T/T},F_b)\otimes_{O_{Y'\times T}}g^{''*mod}(M,F)),
\end{eqnarray*}
is a map in $C_{\mathcal Dfil}(T)$.
\end{prop}

\begin{proof}
Follows imediately by definition. 
\end{proof}

In the analytic case, we also have

\begin{prop}\label{TDhwMinfty}
Consider a commutative diagram in $\AnSm(\mathbb C)$ :
\begin{equation*}
D=\xymatrix{Y\times S\ar[r]^{p} & S \\
Y'\times T\ar[u]^{g''=(g''_0\times g)}\ar[r]^{p'} & T\ar[u]^{g}}
\end{equation*}
with $p$ and $p'$ the projections. 
For $(M,F)\in C_{\mathcal D^{\infty}fil}(Y\times S)$ the map in $C_{g^{''*}p^*O_Sfil}(Y'\times T)$
\begin{equation*}
\Omega_{(Y'\times T/Y\times S)/(T/S)}(M,F):g''^*((\Omega^{\bullet}_{Y\times S/S},F_b)\otimes_{O_{Y\times S}}(M,F))
\to(\Omega^{\bullet}_{Y'\times T/T},F_b)\otimes_{O_{Y'\times T}}g^{''*mod}(M,F)
\end{equation*}
is a map in $C_{g^{''*}p^*\mathcal D^{\infty}fil}(Y'\times T)$. 
Hence, for $(M,F)\in C_{\mathcal D^{\infty}fil}(Y\times S)$, the map in $C_{O_Tfil}(T)$ (with $L_D$ instead of $L_O$)
\begin{eqnarray*}
T^O_{\omega}(D)(M):g^{*mod}L_D(p_*E((\Omega^{\bullet}_{Y\times S/S},F_b)\otimes_{O_{Y\times S}}(M,F)))\to
p'_*E((\Omega^{\bullet}_{Y'\times T/T},F_b)\otimes_{O_{Y'\times T}}g^{''*mod}(M,F)),
\end{eqnarray*}
is a map in $C_{\mathcal D^{\infty}fil}(T)$.
\end{prop}

\begin{proof}
Follows immediately by definition.
\end{proof}

Similarly, we have

\begin{prop}\label{TDhanwM}
Let $p:Y\times S\to S$ a projection with $Y,S\in\SmVar(\mathbb C)$. 
For $(M,F)\in C_{\mathcal Dfil}(Y\times S)$ the map in $C_{p^*O_{S^{an}}}(Y^{an}\times S^{an})$
\begin{equation*}
\Omega_{(Y^{an}\times S^{an}/Y\times S)/(S^{an}/S)}(M,F):
\an(Y\times S)^*((\Omega^{\bullet}_{Y\times S/S},F_b)\otimes_{O_{Y\times S}}(M,F))\to
(\Omega^{\bullet}_{Y^{an}\times S^{an}/S^{an}},F_b)\otimes_{O_{Y\times S}^{an}}(M^{an},F)
\end{equation*}
is a map in $C_{p^\mathcal Dfil}(Y^{an}\times S^{an})$.
For $(M,F)\in C_{\mathcal Dfil}(Y\times S)$, the map in $C_{O_{S^{an}}fil}(S^{an})$
\begin{eqnarray*}
T^O_{\omega}(an,h)(M,F):(p_*E((\Omega^{\bullet}_{Y\times S/S},F_b)\otimes_{O_{Y\times S}}(M,F)))^{an}
\to p_*E((\Omega^{\bullet}_{Y\times S/S},F_b)\otimes_{O_{Y\times S}^{an}}(M,F)^{an})
\end{eqnarray*}
is a map in $C_{\mathcal Dfil}(S^{an})$.
\end{prop}

\begin{proof}
Similar to the proof of proposition \ref{TDhwM}.
\end{proof}

\begin{prop}\label{projformulahw}
Let $p:Y\times S\to S$ a projection with $Y,S\in\SmVar(\mathbb C)$ or with $Y,S\in\AnSm(\mathbb C)$.
\begin{itemize}
\item[(i)] If $\phi:(M,F)\to (N,F)$ is an $r$-filtered Zariski, resp. usu, local equivalence with 
$(M_1,F),(M_2,F)\in C_{\mathcal Dfil}(Y\times S)$, then
\begin{equation*}
DR(Y\times S/S)(\phi):(\Omega_{Y\times S/S}^{\bullet},F_b)\otimes_{O_{Y\times S}}(M,F)\to
\Omega_{Y\times S/S}^{\bullet}\otimes_{O_{Y\times S}}(N,F) 
\end{equation*}
is an $r$-filtered equivalence Zariski, resp. usu, local in $C_{p^*\mathcal Dfil}(Y\times S)$. 
\item[(ii)] Consider a commutative diagram in $\SmVar(\mathbb C)$ or in $\AnSm(\mathbb C)$
\begin{equation*}
D=\xymatrix{Y\times S\ar[r]^{p} & S \\
V\ar[u]^{l}\ar[r]^{k} & S\ar[u]^{=}}.
\end{equation*}
with $p$ the projection. 
For $(N,F)\in C_{\mathcal D,l^*\mathcal Dfil}(V)$, the map in $C_{p^*O_S}(Y\times S)$ (see definition \ref{TDotimeswM})
\begin{eqnarray*}
k\circ T^O_{\omega}(l,\otimes)(E(N,F)):(\Omega_{Y\times S/S}^{\bullet},F_b)\otimes_{O_{Y\times S}}l_*E(N,F)\to  
l_*((\Omega_{V/S}^{\bullet},F_b)\otimes_{O_V}E(N,F)) \\
\to l_*E((\Omega_{V/S}^{\bullet},F_b)\otimes_{O_V}E(N,F)) 
\end{eqnarray*}
is a filtered equivalence Zariski, resp. usu, local in $C_{p^*\mathcal Dfil}(Y\times S)$.
\end{itemize}
\end{prop}

\begin{proof}
\noindent(i):Follows from proposition \ref{DRhUS}that it is a morphism of $p^*D_S$ module.
The fact that it is an equivalence Zariski, resp usu, local is a particular case of proposition \ref{projformulaw}(i).

\noindent(ii):Follows from proposition \ref{DRhUS} and the first part of proposition \ref{TDhwM} that it is a morphism of $h^*D_S$ module.
The fact that it is an equivalence Zariski, resp usu, local is a particular case of proposition \ref{projformulaw}(ii).
\end{proof}

In the analytical case, we also have
\begin{prop}\label{projformulahwinfty}
Let $p:Y\times S\to S$ a projection with $Y,S\in\SmVar(\mathbb C)$ or with $Y,S\in\AnSm(\mathbb C)$.
\begin{itemize}
\item[(i)] If $\phi:(M,F)\to (N,F)$ is an $r$-filtered usu local equivalence with 
$(M_1,F),(M_2,F)\in C_{\mathcal D^{\infty}fil}(Y\times S)$, then
\begin{equation*}
DR(Y\times S/S)(\phi):(\Omega_{Y\times S/S}^{\bullet},F_b)\otimes_{O_{Y\times S}}(M,F)\to
\Omega_{Y\times S/S}^{\bullet}\otimes_{O_{Y\times S}}(N,F) 
\end{equation*}
is an $r$-filtered equivalence usu local in $C_{p^*\mathcal D^{\infty}fil}(Y\times S)$. 
\item[(ii)] Consider a commutative diagram in $\AnSm(\mathbb C)$
\begin{equation*}
D=\xymatrix{Y\times S\ar[r]^{p} & S \\
V\ar[u]^{l}\ar[r]^{k} & S\ar[u]^{=}}.
\end{equation*}
with $p$ the projection. For $(N,F)\in C_{\mathcal D^{\infty},l^*\mathcal D^{\infty}fil}(V)$, 
the map in $C_{p^*O_S}(Y\times S)$ (see definition \ref{TDotimeswM})
\begin{eqnarray*}
k\circ T^O_{\omega}(l,\otimes)(E(N,F)):(\Omega_{Y\times S/S}^{\bullet},F_b)\otimes_{O_{Y\times S}}l_*E(N,F)\to  
l_*((\Omega_{V/S}^{\bullet},F_b)\otimes_{O_V}E(N,F)) \\
\to l_*E((\Omega_{V/S}^{\bullet},F_b)\otimes_{O_V}E(N,F)) 
\end{eqnarray*}
is a filtered equivalence usu local in $C_{p^*\mathcal D^{\infty}fil}(Y\times S)$.
\end{itemize}
\end{prop}

\begin{proof}
Follows from the finite order case : proposition \ref{projformulahw}.
\end{proof}

Dually of the De Rahm complex of a $D_S$ module $M$, we have the Spencer complex of $M$.
In the particular case of $D_S$, we have the following:

\begin{prop}\label{resw}
Let $S\in\SmVar(\mathbb C)$ or $S\in\AnSm(\mathbb C)$.
\begin{itemize}
\item We have the filtered resolutions of $K_S$ by the following complex of locally free right $D_S$ modules:
$\omega(S):\omega(K_S):=(\Omega^{\bullet}_{S},F_b)[d_S]\otimes_{O_S}(D_S,F_b)\to (K_S,F_b)$ and
$\omega(S):\omega(K_S,F^{ord}):=(\Omega^{\bullet}_{S},F_b)[d_S]\otimes_{O_S}(D_S,F^{ord})\to (K_S,F^{ord})$
\item Dually, we have the filtered resolution of $O_S$ by the following complex of locally free (left) $D_S$ modules:
$\omega^{\vee}(S):\omega(O_S):=(\wedge^{\bullet}T_S,F_b)[d_S]\otimes_{O_S}(D_S,F_b)\to (O_S,F_b)$ and
$\omega^{\vee}(S):\omega(O_S,F^{ord}):=(\wedge^{\bullet}T_S,F_b)[d_S]\otimes_{O_S}(D_S,F^{ord})\to (O_S,F^{ord})$.
\end{itemize}
Let $S_1,S_2\in\SmVar(\mathbb C)$ or $S_1,S_2\in\AnSm(\mathbb C)$. Consider the projection $p=p_1:S_1\times S_2\to S_1$.
\begin{itemize}
\item We have the filtered resolution of $D_{S_1\times S_2\to S_1}$ by the following complexes of (left)
$(p^*D_{S_1}$ and right $D_{S_1\times S_2})$ modules :
\begin{equation*}
\omega(S_1\times S_2/S_1):
(\Omega^{\bullet}_{S_1\times S_2/S_1}[d_{S_2}],F_b)\otimes_{O_{S_1\times S_2}}(D_{S_1\times S_2},F^{ord})
\to(D_{S_1\times S_2\leftarrow S_1},F^{ord}).
\end{equation*}
\item Dually, we have the filtered resolution of $D_{S_1\times S_2\to S_1}$ by the following complexes of (left)
$(p^*D_{S_1},D_{S_1\times S_2})$ modules :
\begin{equation*}
\omega^{\vee}(S_1\times S_2/S_1):
(\wedge^{\bullet}T_{S_1\times S_2/S_1}[d_{S_2}],F_b)\otimes_{O_{S_1\times S_2}}(D_{S_1\times S_2},F^{ord})
\to(D_{S_1\times S_2\to S_1},F^{ord}),
\end{equation*}
\end{itemize}
\end{prop}

\begin{proof} 
See \cite{LvDmod}.
\end{proof}

In the analytical case, we also have

\begin{prop}\label{reswan}
Let $S\in\AnSm(\mathbb C)$.
\begin{itemize}
\item We have the filtered resolutions of $K_S$ by the following complex of locally free right $D_S$ modules:
$\omega(S):\omega(K_S):=(\Omega^{\bullet}_{S},F_b)[d_S]\otimes_{O_S}(D^{\infty}_S,F^{ord})\to (K_S,F_b)$.
\item Dually, we have the filtered resolution of $O_S$ by the following complex of locally free (left) $D_S$ modules:
$\omega^{\vee}(S):\omega(O_S):=(\wedge^{\bullet}T_S,F_b)[d_S]\otimes_{O_S}(D^{\infty}_S,F^{ord})\to (O_S,F_b)$.
\end{itemize}
Let $S_1,S_2\in\AnSm(\mathbb C)$. Consider the projection $p=p_1:S_1\times S_2\to S_1$.
\begin{itemize}
\item We have the filtered resolution of $D^{\infty}_{S_1\times S_2\to S_1}$ by the following complexes of (left)
$(p^*D^{\infty}_{S_1}$ and right $D^{\infty}_{S_1\times S_2})$ modules :
\begin{equation*}
\omega(S_1\times S_2/S_1):
(\Omega^{\bullet}_{S_1\times S_2/S_1}[d_{S_2}],F_b)\otimes_{O_{S_1\times S_2}}(D^{\infty}_{S_1\times S_2},F^{ord})
\to(D^{\infty}_{S_1\times S_2\leftarrow S_1},F^{ord}).
\end{equation*}
\item Dually, we have the filtered resolution of $D^{\infty}_{S_1\times S_2\to S_1}$ by the following complexes of (left)
$(p^*D^{\infty}_{S_1},D^{\infty}_{S_1\times S_2})$ modules :
\begin{equation*}
\omega^{\vee}(S_1\times S_2/S_1):
(\wedge^{\bullet}T_{S_1\times S_2/S_1}[d_{S_2}],F_b)\otimes_{O_{S_1\times S_2}}(D^{\infty}_{S_1\times S_2},F^{ord})
\to(D^{\infty}_{S_1\times S_2\to S_1},F^{ord}),
\end{equation*}
\end{itemize}
\end{prop}

\begin{proof} 
Similar to the finite order case : the first map on the right is a surjection and the kernel are obtained by tensoring
$D_S^{\infty}$ with the kernel of the kozcul resolution of $K_S$ 
(note that $D_S^{\infty}$ is a locally free hence flat $O_S$ module).
\end{proof}

Motivated by these resolutions, we make the following definition

\begin{defi}\label{directmod}
\begin{itemize}
\item[(i)] Let $i:Z\hookrightarrow S$ be a closed embedding, with $Z,S\in\SmVar(\mathbb C)$ or with $Z,S\in\AnSm(\mathbb C)$.
Then, for $(M,F)\in C_{\mathcal Dfil}(Z)$, we set 
\begin{equation*}
i_{*mod}(M,F):=i^0_{*mod}(M,F):=i_*((M,F)\otimes_{D_Z}(D_{Z\leftarrow S},F^{ord}))\in C_{\mathcal Dfil}(S)
\end{equation*}
\item[(ii)] Let $S_1,S_2\in\SmVar(\mathbb C)$ or $S_1,S_2\in\AnSm(\mathbb C)$ and $p:S_1\times S_2\to S_1$ be the projection.
Then, for $(M,F)\in C_{\mathcal Dfil}(S_1\times S_2)$, we set
\begin{itemize}
\item $p^0_{*mod}(M,F):=p_*(DR(S_1\times S_2/S_1)(M,F)):=
p_*((\Omega^{\bullet}_{S_1\times S_2/S_1},F_b)\otimes_{O_{S_1\times S_2}}(M,F))[d_{S_2}]\in C_{\mathcal Dfil}(S_1)$,
\item $p_{*mod}(M,F):=p_*E(DR(S_1\times S_2/S_1)(M,F)):=
p_*E((\Omega^{\bullet}_{S_1\times S_2/S_1},F_b)\otimes_{O_{S_1\times S_2}}(M,F))[d_{S_2}]\in C_{\mathcal Dfil}(S_1)$.
\end{itemize}
\item[(iii)] Let $f:X\to S$ be a morphism, with $X,S\in\SmVar(\mathbb C)$ or $X,S\in\AnSm(\mathbb C)$.
Consider the factorization $f:X\xrightarrow{i} X\times S\xrightarrow{p_S}S$, 
where $i$ is the graph embedding and $p_S:X\times S\to S$ is the projection.
Then, for $(M,F)\in C_{\mathcal Dfil}(X)$ we set
\begin{itemize}
\item $f^{FDR}_{*mod}(M,F):=p_{S*mod}i_{*mod}(M,F)\in C_{\mathcal Dfil}(S)$,
\item $\int_f^{FDR}(M,F):=f^{FDR}_{*mod}(M,F):=p_{S*mod}i_{*mod}(M,F)\in D_{\mathcal Dfil,\infty}(S)$.
\end{itemize}
By proposition \ref{Pham} below, we have $\int_f^{FDR}M=\int_fM\in D_{\mathcal D}(X)$.
\item[(iii)] Let $f:X\to S$ be a morphism, with $X,S\in\SmVar(\mathbb C)$ or $X,S\in\AnSm(\mathbb C)$.
Consider the factorization $f:X\xrightarrow{i} X\times S\xrightarrow{p_S}S$, 
where $i$ is the graph embedding and $p_S:X\times S\to S$ is the projection.
Then, for $(M,F)\in C_{\mathcal Dfil}(X)$ we set
\begin{itemize}
\item $f^{FDR}_{!mod}(M,F):=\mathbb D_S^KL_Df^{FDR}_{*mod}\mathbb D_S^KL_D(M,F):=
\mathbb D_S^KL_Dp_{S*mod}i_{*mod}\mathbb D_{X\times S}^KL_D(M,F)\in C_{\mathcal Dfil}(S)$,
\item $\int_{f!}^{FDR}(M,F):=f^{FDR}_{!mod}(M,F):=
\mathbb D_S^KL_Dp_{S*mod}i_{*mod}\mathbb D_{X\times S}^KL_D(M,F)\in D_{\mathcal Dfil,\infty}(S)$.
\end{itemize}
\end{itemize}
\end{defi}

In the analytical case we also consider :

\begin{defi}\label{directmod'}
\begin{itemize}
\item[(i)] Let $i:Z\hookrightarrow S$ be a closed embedding with $Z,S\in\AnSm(\mathbb C)$.
Then, for $(M,F)\in C_{\mathcal D^{\infty}fil}(Z)$, we set 
\begin{equation*}
i_{*mod}(M,F):=i^0_{*mod}(M,F):=i_*((M,F)\otimes_{D^{\infty}_Z}(D^{\infty}_{Z\leftarrow S},F^{ord}))
\in C_{\mathcal D^{\infty}fil}(S)
\end{equation*}
\item[(ii)] Let $S_1,S_2\in\AnSm(\mathbb C)$ and $p:S_1\times S_2\to S_1$ be the projection.
For $(M,F)\in C_{\mathcal D^{\infty}fil}(S_1\times S_2)$, we consider
\begin{itemize}
\item $p^0_{*mod}(M,F):=p_*(DR(S_1\times S_2/S_1)(M,F)):=
p_*((\Omega^{\bullet}_{S_1\times S_2/S_1},F_b)\otimes_{O_{S_1\times S_2}}(M,F))[d_{S_2}]\in C_{\mathcal D^{\infty}fil}(S_1)$,
\item $p_{*mod}(M,F):=p_*E(DR(S_1\times S_2/S_1)(M,F)):=
p_*E((\Omega^{\bullet}_{S_1\times S_2/S_1},F_b)\otimes_{O_{S_1\times S_2}}(M,F))[d_{S_2}]\in C_{\mathcal D^{\infty}fil}(S_1)$.
\end{itemize}
\item[(iii)] Let $S_1,S_2\in\AnSm(\mathbb C)$ and $p:S_1\times S_2\to S_1$ be the projection.
For $(M,F)\in C_{\mathcal Dfil}(S_1\times S_2)$ or $(M,F)\in C_{\mathcal D^{\infty}fil}(S_1\times S_2)$, we set
\begin{itemize}
\item $p^0_{!mod}(M,F):=p_!(DR(S_1\times S_2/S_1)(M,F)):=
p_!((\Omega^{\bullet}_{S_1\times S_2/S_1},F_b)\otimes_{O_{S_1\times S_2}}(M,F))[d_{S_2}]\in C_{\mathcal Dfil}(S_1)$,
\item $p_{!mod}(M,F):=p_!E(DR(S_1\times S_2/S_1)(M,F)):=
p_!E((\Omega^{\bullet}_{S_1\times S_2/S_1},F_b)\otimes_{O_{S_1\times S_2}}(M,F))[d_{S_2}]\in C_{\mathcal Dfil}(S_1)$.
\end{itemize}
\item[(iv)] Let $f:X\to S$ be a morphism, with $X,S\in\AnSm(\mathbb C)$.
Consider the factorization $f:X\xrightarrow{i} X\times S\xrightarrow{p_S}S$, 
where $i$ is the graph embedding and $p_S:X\times S\to S$ is the projection. 
Then, for $(M,F)\in C_{\mathcal D^{\infty}fil}(X)$ we set
\begin{itemize}
\item  $f^{FDR}_{*mod}(M,F):=p_{S*mod}i_{*mod}(M,F)\in C_{\mathcal D^{\infty}fil}(S)$,
\item $\int_f^{FDR}(M,F):=f^{FDR}_{*mod}(M,F):=p_{S*mod}i_{*mod}(M,F)\in D_{\mathcal D^{\infty}fil,\infty}(S)$,
\item $f^{FDR}_{!mod}(M,F):=p_{S!mod}i_{*mod}(M,F)\in C_{\mathcal D^{\infty}fil}(S)$,
\item $\int_{f!}^{FDR}(M,F):=f^{FDR}_{!*mod}(M,F):=p_{S!mod}i_{*mod}(M,F)\in D_{\mathcal D^{\infty}fil,\infty}(S)$.
\end{itemize}
By proposition \ref{Pham'} below, we have $\int_{f!}^{FDR}M=\int_{f!}M\in D_{\mathcal D^{\infty}}(X)$
and $\int_f^{FDR}M=\int_fM\in D_{\mathcal D^{\infty}}(X)$.
\item[(v)] Let $f:X\to S$ be a morphism, with $X,S\in\AnSm(\mathbb C)$.
Consider the factorization $f:X\xrightarrow{i} X\times S\xrightarrow{p_S}S$, 
where $i$ is the graph embedding and $p_S:X\times S\to S$ is the projection. 
Then, for $(M,F)\in C_{\mathcal Dfil}(X)$ we set
\begin{itemize}
\item $f^{FDR}_{!mod}(M,F):=p_{S!mod}i_{*mod}(M,F)\in C_{\mathcal Dfil}(S)$,
\item $\int_{f!}^{FDR}(M,F):=f^{FDR}_{!*mod}(M,F):=p_{S!mod}i_{*mod}(M,F)\in D_{\mathcal Dfil,\infty}(S)$.
\end{itemize}
By proposition \ref{Pham'} below, we have $\int_{f!}^{FDR}M=\int_{f!}M\in D_{\mathcal D}(X)$.
\end{itemize}
\end{defi}

\begin{prop}\label{Pham}
\begin{itemize}
\item[(i)] Let $i:Z\hookrightarrow S$ a closed embedding with $S,Z\in\SmVar(\mathbb C)$ or with $S,Z\in\AnSm(\mathbb C)$.
Then for $(M,F)\in C_{\mathcal Dfil}(Z)$, we have
\begin{equation*}
\int_i(M,F):=Ri_*((M,F)\otimes^L_{D_Z}(D_{Z\leftarrow S},F^{ord})=i_*((M,F)\otimes_{D_Z}(D_{Z\leftarrow S},F^{ord}))=i_{*mod}(M,F).
\end{equation*}
\item[(ii)] Let $S_1,S_2\in\SmVar(\mathbb C)$ or $S_1,S_2\in\AnSm(\mathbb C)$ and $p:S_{12}:=S_1\times S_2\to S_1$ be the projection.
Then, for $(M,F)\in C_{\mathcal Dfil}(S_1\times S_2)$ we have
\begin{eqnarray*}
\int_p(M,F):&=&Rp_*((M,F)\otimes^L_{D_{S_1\times S_2}}(D_{S_1\times S_2\leftarrow S_1},F^{ord})) \\
&=&p_*E((\Omega^{\bullet}_{S_1\times S_2/S_1},F_b)\otimes_{O_{S_1\times S_2}}
(D_{S_1\times S_2},F^{ord})\otimes_{D_{S_1\times S_2}}(M,F))[d_{S_2}] \\
&=&p_*E((\Omega^{\bullet}_{S_1\times S_2/S_1},F_b)\otimes_{O_{S_1\times S_2}}(M,F))[d_{S_2}]
=:p_{*mod}(M,F).
\end{eqnarray*}
where the second equality follows from Griffitz transversality (the canonical isomorphism map respect by definition the filtration).
\item[(iii)] Let $f:X\to S$ be a morphism with $X,S\in\SmVar(\mathbb C)$ or with $X,S\in\AnSm(\mathbb C)$.
Then for $M\in C_{\mathcal D}(X)$, we have $\int^{FDR}_fM=\int_fM$.
\end{itemize}
\end{prop}

\begin{proof}

\noindent(i):Follows from the fact that $D_{Z\leftarrow S}$ is a locally free $D_Z$ module and that $i_*$ is an exact functor.

\noindent(ii): Since $\Omega^{\bullet}_{S_{12}/S_1}[d_{S_2}],F_b)\otimes_{O_{S_{12}}}D_{S_{12}}$
is a complex of locally free $D_{S_1\times S_2}$ modules, we have in $D_{fil}(S_1\times S_2)$, using proposition \ref{resw},
\begin{eqnarray*}
(D_{S_1\times S_2\leftarrow S_1},F^{ord})\otimes^L_{D_{S_1\times S_2}}(M,F)=
(\Omega^{\bullet}_{S_{12}/S_1}[d_{S_2}],F_b)\otimes_{O_{S_{12}}}(D_{S_{12}},F^{ord})\otimes_{D_{S_{12}}}(M,F).
\end{eqnarray*} 

\noindent(iii): Follows from (i) and (ii) by proposition \ref{compDmod} (ii) in the algebraic case
and by proposition \ref{compAnDmod}(ii) in the analytic case since a closed embedding is proper.

\end{proof}

In the analytical case, we also have :

\begin{prop}\label{Pham'}
\begin{itemize}
\item[(i)] Let $i:Z\hookrightarrow S$ a closed embedding  with $S,Z\in\AnSm(\mathbb C)$.
Then for $(M,F)\in C_{\mathcal D^{\infty}fil}(Z)$, we have $\int_i(M,F)=i_{*mod}(M,F)$.
\item[(ii)] Let $S_1,S_2\in\AnSm(\mathbb C)$ and $p:S_{12}:=S_1\times S_2\to S_1$ be the projection.
Then, for $(M,F)\in C_{\mathcal D^{\infty}fil}(S_1\times S_2)$ we have
\begin{eqnarray*}
\int_{p}(M,F):&=&Rp_*((M,F)\otimes^L_{D^{\infty}_{S_1\times S_2}}(D^{\infty}_{S_1\times S_2\leftarrow S_1},F^{ord})) \\
&=&p_*E((\Omega^{\bullet}_{S_1\times S_2/S_1},F_b)\otimes_{O_{S_1\times S_2}}
(D_{S_1\times S_2},F^{ord})\otimes_{D_{S_1\times S_2}}(M,F)[d_{S_2}]) \\
&=&p_*E((\Omega^{\bullet}_{S_1\times S_2/S_1},F_b)\otimes_{O_{S_1\times S_2}}(M,F)[d_{S_2}])
=:p_{*mod}(M,F).
\end{eqnarray*}
\item[(ii)'] Let $S_1,S_2\in\AnSm(\mathbb C)$ and $p:S_{12}:=S_1\times S_2\to S_1$ be the projection.
Then, for $(M,F)\in C_{\mathcal Dfil}(S_1\times S_2)$ or $(M,F)\in C_{\mathcal Dfil}(S_1\times S_2)$, we have
\begin{eqnarray*}
\int_{p!}(M,F):&=&Rp_!((M,F)\otimes^L_{D_{S_1\times S_2}}(D_{S_1\times S_2\leftarrow S_1},F^{ord}) \\
&=&p_!E((\Omega^{\bullet}_{S_1\times S_2/S_1},F_b)\otimes_{O_{S_1\times S_2}}
(D_{S_1\times S_2},F^{ord})\otimes_{D_{S_1\times S_2}}(M,F)[d_{S_2}]) \\
&=&p_!E((\Omega^{\bullet}_{S_1\times S_2/S_1},F_b)\otimes_{O_{S_1\times S_2}}(M,F)[d_{S_2}])
=:p_{!mod}(M,F).
\end{eqnarray*}
\item[(iii)] Let $f:X\to S$ be a morphism with $X,S\in\AnSm(\mathbb C)$. 
For $M\in C_{\mathcal D^{\infty}}(X)$, we have $\int^{FDR}_fM=\int_fM$ and $\int^{FDR}_{f!}M=\int_{f!}M$.
For $M\in C_{\mathcal D}(X)$, we have $\int^{FDR}_{f!}M=\int_{f!}M$.
\end{itemize}
\end{prop}

\begin{proof}

\noindent(i):Follows from the fact that $D^{\infty}_{Z\leftarrow S}$ is a locally free $D^{\infty}_Z$ module 
and that $i_*$ is an exact functor.

\noindent(ii): Similar to the proof of proposition \ref{Pham}(ii):follows from proposition \ref{reswan}.

\noindent(ii)': Similar to the proof of proposition \ref{Pham}(ii):follows from proposition \ref{reswan}.

\noindent(iii):The first assertion follows from (i), (ii) and (ii)' by proposition \ref{compAnDmod}.
The second one follows from proposition \ref{Pham}(i) and  (ii)' and by proposition \ref{compAnDmod}.

\end{proof}

\begin{prop}\label{compDmodDRd}
Let $f_1:X\to Y$ and $f_2:Y\to S$ two morphism with $X,Y,S\in\SmVar(\mathbb C)$. 
\begin{itemize}
\item[(i)]Let $(M,F)\in C_{\mathcal D(2)fil}(X)$. 
Then $\int^{FDR}_{f_2\circ f_1}(M,F)=\int_{f_1}^{FDR}\int_{f_2}^{FDR}(M,F)\in D_{\mathcal D(2)fil,\infty}(S)$.
\item[(ii)]Let $(M,F)\in C_{\mathcal D(2)fil,h}(X)$. 
Then $\int^{FDR}_{(f_2\circ f_1)!}(M,F)=\int_{f_1!}^{FDR}\int_{f_2!}^{FDR}(M,F)\in D_{\mathcal D(2)fil,\infty}(S)$.
\end{itemize}
\end{prop}

\begin{proof}
See \cite{Laumont}.
\end{proof}

\begin{prop}\label{compAnDmodDRd}
Let $f_1:X\to Y$ and $f_2:Y\to S$ two morphism with $X,Y,S\in\AnSm(\mathbb C)$. 
\begin{itemize}
\item[(i)] Let $(M,F)\in C_{\mathcal D^{\infty}(2)fil}(X)$. 
Then $\int_{(f_2\circ f_1)!}^{FDR}(M,F)=\int_{f_1!}^{FDR}\int_{f_2!}^{FDR}(M,F)$.
\item[(ii)] Let $(M,F)\in C_{\mathcal D^{\infty}(2)fil,h}(X)$. 
Then $\int_{f_2\circ f_1}^{FDR}(M,F)=\int_{f_1}^{FDR}\int_{f_2}^{FDR}(M,F)$.
\end{itemize}
\end{prop}

\begin{proof}
Similar to proposition \ref{compDmodDRd}.
\end{proof}

\begin{defi}\label{TDDR}
\begin{itemize}
\item[(i)]Let $f:X\to S$ be a morphism, with $X,S\in\SmVar(\mathbb C)$ or $X,S\in\AnSm(\mathbb C)$.
Consider the graph factorization $f:X\xrightarrow{l}X\times S\xrightarrow{p}S$, 
with $l$ the graph embedding and $p$ the projection.
We have the transformation map given by, for $(M,F)\in C_{\mathcal Dfil}(X)$,
\begin{eqnarray*}
T(\int_f^{FDR},\int_f)(M,F):\int_f^{FDR}(M,F):=\int_p\int_l(M,F)
\xrightarrow{T(\int_p\circ\int_l,\int_{p\circ l})(M,F)}\int_f(M,F)
\end{eqnarray*}
\item[(ii)]Let $j:S^o\hookrightarrow S$ an open embedding with $S\in\Var(\mathbb C)$.
Consider the graph factorization $j:S^o\xrightarrow{l}S^o\times S\xrightarrow{p} S$, 
with $l$ the graph embedding and $p$ the projection.
We have, for $(M,F)\in C_{\mathcal Dfil}(S^o)$, the canonical map in $C_{\mathcal Dfil}(S)$,
\begin{eqnarray*}
T(j^{FDR}_{*mod},j_*)(M,F):j^{FDR}_{*mod}(M,F):=
p_*E((\Omega^{\bullet}_{S^o\times S/S},F_b)\otimes_{O_{S^o\times S}}l_{*mod}(M,F))
\xrightarrow{k\circ\omega(S^o\times S/S)} \\
p_*E((D_{S^o\times S\leftarrow S},F^{ord})\otimes_{D_{S^o\times S}}l_*(D_{S^o\leftarrow S^o\times S}\otimes_{D_{S^o}}E(M,F))
\xrightarrow{T(l,\otimes)(-,-)} j_*E(M,F)
\end{eqnarray*}
We have, for $(M,F)\in C_{\mathcal Dfil}(S)$, the canonical map in $C_{\mathcal Dfil}(S)$,
\begin{eqnarray*}
\ad(j^*,j_{*mod}^{FDR})(M,F):(M,F)\xrightarrow{\ad(p^{*mod},p_*)(M,F)}
p_*E((\Omega^{\bullet}_{S^o\times S/S},F_b)\otimes_{O_{S^o\times S}}p^{*mod}(M,F))
\end{eqnarray*}
\end{itemize}
\end{defi}

\subsubsection{The support section functors for D modules and the graph inverse image}
 
Let $S\in\SmVar(\mathbb C)$ or $S\in\AnSm(\mathbb C)$. Let $i:Z\hookrightarrow S$ a closed embedding
and denote $j:S\backslash Z\hookrightarrow S$ the complementary open embedding. 
More generally, let $h:Y\to S$ a morphism with $Y,S\in\Var(\mathbb C)$ or $Y,S\in\AnSp(\mathbb C)$, $S$ smooth, and 
let $i:X\hookrightarrow Y$ a closed embedding and denote by $j:Y\backslash X\hookrightarrow Y$ the open complementary. 
We then get from section 2 the following functors :
\begin{itemize}
\item We get the functor
\begin{eqnarray*}
\Gamma_Z:C_{\mathcal D(2)fil}(S)\to C_{\mathcal D(2)fil}(S), \\ 
(M,F)\mapsto\Gamma_Z(M,F):=\Cone(\ad(j^*,j_*)((M,F)):(M,F)\to j_*j^*(M,F))[-1],
\end{eqnarray*}
together we the canonical map $\gamma_Z(M,F):\Gamma_Z(M,F)\to (M,F)$, and more generally the functor
\begin{eqnarray*}
\Gamma_X:C_{h^*\mathcal D(2)fil}(Y)\to C_{h^*\mathcal D(2)fil}(Y), \\ 
(M,F)\mapsto\Gamma_X(M,F):=\Cone(\ad(j^*,j_*)((M,F)):(M,F)\to j_*j^*(M,F))[-1],
\end{eqnarray*}
together we the canonical map $\gamma_X(M,F):\Gamma_X(M,F)\to (M,F)$.
\item We get the functor
\begin{eqnarray*}
\Gamma_Z^{\vee}:C_{\mathcal D(2)fil}(S)\to C_{\mathcal D(2)fil}(S), \\ 
(M,F)\mapsto\Gamma_Z^{\vee}(M,F):=\Cone(\ad(j_!,j^*)((M,F)):j_!j^*(M,F)\to(M,F)),
\end{eqnarray*}
together we the canonical map $\gamma_Z^{\vee}(M,F):(M,F)\to\Gamma_Z^{\vee}(M,F)$, 
and more generally the functor
\begin{eqnarray*}
\Gamma_X^{\vee}:C_{h^*\mathcal D(2)fil}(Y)\to C_{h^*\mathcal D(2)fil}(Y), \\ 
(M,F)\mapsto\Gamma_X^{\vee}(M,F):=\Cone(\ad(j_!,j^*)((M,F)):j_!j^*(M,F)\to(M,F)),
\end{eqnarray*}
together we the canonical map $\gamma_X^{\vee}(M,F):(M,F)\to\Gamma_X^{\vee}(M,F)$.
\item We get the functor
\begin{eqnarray*}
\Gamma_Z^{\vee,h}:C_{\mathcal D(2)fil}(S)\to C_{\mathcal D(2)fil}(S), \;  
(M,F)\mapsto\Gamma_Z^{\vee,h}(M,F):=\mathbb D_S^KL_D\Gamma_ZE(\mathbb D^K_S(M,F)),
\end{eqnarray*}
together with the factorization 
\begin{eqnarray*}
\gamma_Z^{\vee,h}(L_D(M,F)):L_D(M,F)\xrightarrow{\gamma_Z^{\vee}(L_D(M,F))}\Gamma_Z^{\vee}L_D(M,F)
\xrightarrow{k\circ\mathbb D^KI(j_!,j^*)(-)\circ d(-)}\Gamma_Z^{\vee,h}L_D(M,F),
\end{eqnarray*}
and more generally the functor
\begin{eqnarray*}
\Gamma_X^{\vee,h}:C_{h^*\mathcal D(2)fil}(Y)\to C_{h^*\mathcal D(2)fil}(Y), \;  
(M,F)\mapsto\Gamma_X^{\vee,h}(M,F):=\mathbb D_Y^{h^*D,K}L_{h^*D}\Gamma_XE(\mathbb D^{h^*D,K}_Y(M,F)),
\end{eqnarray*}
together with the factorization 
\begin{eqnarray*}
\gamma_X^{\vee,h}(L_{h^*D}(M,F)):L_{h^*D}(M,F)\xrightarrow{\gamma_X^{\vee}(L_{h^*D}(M,F))}\Gamma_X^{\vee}L_{h^*D}(M,F) \\
\xrightarrow{k\circ\mathbb D^{h^*D,K}I(j_!,j^*)(-)\circ d(-)}\Gamma_X^{\vee,h}L_{h^*D}(M,F).
\end{eqnarray*}
\item We get the functor
\begin{eqnarray*}
\Gamma_Z^{\vee,O}:C_{\mathcal D(2)fil}(S)\to C_{\mathcal D(2)fil}(S), \\ 
(M,F)\mapsto\Gamma_Z^{\vee,O}(M,F):=\Cone(b_Z((M,F)):\mathcal I_Z(M,F)\to(M,F)),
\end{eqnarray*}
together we the factorization 
\begin{eqnarray*}
\gamma_Z^{\vee,O}(M,F):(M,F)\xrightarrow{\gamma_Z^{\vee}(M,F)}\Gamma_Z^{\vee}(M,F)
\xrightarrow{b_{S/Z}(M,F)}\Gamma_Z^{\vee,O}(M,F).
\end{eqnarray*}
Since $M\mapsto M/\mathcal I_ZM$ is a right exact functor, $M\mapsto\Gamma^{\vee,O}_ZM$ send Zariski, resp. usu, local
equivalence between projective complexes of presheaves to Zariski, resp. usu local equivalence, and thus induces in the derived category
\begin{eqnarray*}
L\Gamma^{\vee,O}_Z:D_{\mathcal Dfil,\infty}(S)\to D_{\mathcal Dfil,\infty}(S), \\ 
(M,F)\mapsto\Gamma^{\vee,O}_ZL_D(M,F):=\Cone(b_Z(L_D(M,F)):\mathcal I_ZL_D(M,F)\to L_D(M,F)).
\end{eqnarray*}
\item We get the functor
\begin{eqnarray*}
\Gamma_Z^O:C_{\mathcal D(2)fil}(S)\to C_{\mathcal D(2)fil}(S), \\ 
(M,F)\mapsto\Gamma_Z^O(M,F):=\Cone(b'_Z((M,F)):(M,F)\to(M,F)\otimes_{O_S}\mathbb D_S^K(\mathcal I_ZD_S)),
\end{eqnarray*}
together we the factorization 
\begin{eqnarray*}
\gamma_Z^O(M,F):(M,F)\Gamma_Z^O\xrightarrow{b'_{S/Z}(M,F)}\Gamma_Z(M,F)\xrightarrow{\gamma_Z(M,F)}(M,F).
\end{eqnarray*}
\item We have, for $(M,F)\in C_{\mathcal Dfil}(S)$, a canonical isomorphism 
\begin{equation*}
I(D,\gamma^O)(M,F):\mathbb D^K_S\Gamma_Z^{\vee,O}(M,F)\xrightarrow{\sim}\Gamma_Z^O\mathbb D^K_S(M,F)
\end{equation*}
which gives the transformation map in $C_{\mathcal Dfil}(S)$
\begin{eqnarray*}
T(D,\gamma^O)(M,F):\Gamma_Z^{\vee,O}\mathbb D^K_S(M,F)\to\mathbb D^K_S\Gamma_Z^O(M,F)
\end{eqnarray*}
\end{itemize}
Let $S\in\AnSm(\mathbb C)$. Let $i:Z\hookrightarrow S$ a closed embedding
and denote $j:S\backslash Z\hookrightarrow S$ the complementary open embedding. 
More generally, let $h:Y\to S$ a morphism with $Y,S\in\AnSp(\mathbb C)$, $S$ smooth, and 
let $i:X\hookrightarrow Y$ a closed embedding and denote by $j:Y\backslash X\hookrightarrow Y$ the open complementary. 
\begin{itemize}
\item We get the functor
\begin{eqnarray*}
\Gamma_Z:C_{\mathcal D^{\infty}(2)fil}(S)\to C_{\mathcal D^{\infty}(2)fil}(S), \\ 
(M,F)\mapsto\Gamma_Z(M,F):=\Cone(\ad(j^*,j_*)((M,F)):(M,F)\to j_*j^*(M,F))[-1],
\end{eqnarray*}
together we the canonical map $\gamma_Z(M,F):\Gamma_Z(M,F)\to (M,F)$, and more generally the functor
\begin{eqnarray*}
\Gamma_X:C_{h^*\mathcal D^{\infty}(2)fil}(Y)\to C_{h^*\mathcal D^{\infty}(2)fil}(Y), \\ 
(M,F)\mapsto\Gamma_X(M,F):=\Cone(\ad(j^*,j_*)((M,F)):(M,F)\to j_*j^*(M,F))[-1],
\end{eqnarray*}
together we the canonical map $\gamma_X(M,F):\Gamma_X(M,F)\to (M,F)$.
\item We get the functor
\begin{eqnarray*}
\Gamma_Z^{\vee}:C_{\mathcal D^{\infty}(2)fil}(S)\to C_{\mathcal D^{\infty}(2)fil}(S), \\ 
(M,F)\mapsto\Gamma_Z^{\vee}(M,F):=\Cone(\ad(j_!,j^*)((M,F)):j_!j^*(M,F)\to(M,F)),
\end{eqnarray*}
together we the canonical map $\gamma_Z^{\vee}(M,F):(M,F)\to\Gamma_Z^{\vee}(M,F)$, 
and more generally the functor
\begin{eqnarray*}
\Gamma_X^{\vee}:C_{h^*\mathcal D^{\infty}(2)fil}(Y)\to C_{h^*\mathcal D^{\infty}(2)fil}(Y), \\ 
(M,F)\mapsto\Gamma_X^{\vee}(M,F):=\Cone(\ad(j_!,j^*)((M,F)):j_!j^*(M,F)\to(M,F)),
\end{eqnarray*}
together we the canonical map $\gamma_X^{\vee}(M,F):(M,F)\to\Gamma_X^{\vee}(M,F)$. 
\item We get the functor
\begin{eqnarray*}
\Gamma_Z^{\vee,h}:C_{\mathcal D^{\infty}(2)fil}(S)\to C_{\mathcal D^{\infty}(2)fil}(S), \;  
(M,F)\mapsto\Gamma_Z^{\vee,h}(M,F):=\mathbb D_S^{\infty,K}L_{D^{\infty}}\Gamma_ZE(\mathbb D^{\infty,K}_S(M,F)),
\end{eqnarray*}
together with the factorization 
\begin{eqnarray*}
\gamma_Z^{\vee,h}(L_{D^{\infty}}(M,F)):L_{D^{\infty}}(M,F)\xrightarrow{\gamma_Z^{\vee}(L_{D^{\infty}}(M,F))}
\Gamma_Z^{\vee}L_{D^{\infty}}(M,F)
\xrightarrow{k\circ\mathbb D^{\infty}I(j_!,j^*)(-)\circ d(-)}\Gamma_Z^{\vee,h}L_{D^{\infty}}(M,F),
\end{eqnarray*}
and more generally the functor
\begin{eqnarray*}
\Gamma_X^{\vee,h}:C_{h^*\mathcal D^{\infty}(2)fil}(Y)\to C_{h^*\mathcal D^{\infty}(2)fil}(Y), \;  
(M,F)\mapsto\Gamma_X^{\vee,h}(M,F):=\mathbb D_Y^{h^*\infty,K}L_{h^*D^{\infty}}\Gamma_XE(\mathbb D^{h^*\infty,K}_Y(M,F)),
\end{eqnarray*}
together with the factorization 
\begin{eqnarray*}
\gamma_X^{\vee,h}(L_{h^*D^{\infty}}(M,F)):L_{h^*D^{\infty}}(M,F)\xrightarrow{\gamma_X^{\vee}(L_{h^*D^{\infty}}(M,F))}
\Gamma_X^{\vee}L_{h^*D^{\infty}}(M,F) \\ 
\xrightarrow{k\circ\mathbb D^{h^*D^{\infty,K}}I(j_!,j^*)(-)\circ d(-)}\Gamma_X^{\vee,h}L_{h^*D^{\infty}}(M,F).
\end{eqnarray*}
\item We get the functor
\begin{eqnarray*}
\Gamma_Z^{\vee,O}:C_{\mathcal D^{\infty}(2)fil}(S)\to C_{\mathcal D^{\infty}(2)fil}(S), \\ 
(M,F)\mapsto\Gamma_Z^{\vee,O}(M,F):=\Cone(b_Z((M,F)):\mathcal I_Z(M,F)\to(M,F)),
\end{eqnarray*}
together we the factorization 
\begin{eqnarray*}
\gamma_Z^{\vee,O}(M,F):(M,F)\xrightarrow{\gamma_Z^{\vee}(M,F)}\Gamma_Z^{\vee}(M,F)
\xrightarrow{b_{S/Z}(M,F)}\Gamma_Z^{\vee,O}(M,F).
\end{eqnarray*}
\item We get the functor
\begin{eqnarray*}
\Gamma_Z^O:C_{\mathcal D^{\infty}(2)fil}(S)\to C_{\mathcal D^{\infty}(2)fil}(S), \\ 
(M,F)\mapsto\Gamma_Z^O(M,F):=\Cone(b'_Z((M,F)):(M,F)\to(M,F)\otimes_{O_S}\mathbb D_S^K(\mathcal I_ZD_S)),
\end{eqnarray*}
together we the factorization 
\begin{eqnarray*}
\gamma_Z^O(M,F):(M,F)\Gamma_Z^O\xrightarrow{b'_{S/Z}(M,F)}\Gamma_Z(M,F)\xrightarrow{\gamma_Z(M,F)}(M,F).
\end{eqnarray*}
\item We have, for $(M,F)\in C_{\mathcal D^{\infty}fil}(S)$, a canonical isomorphism 
\begin{equation*}
I(D,\gamma^O)(M,F):\mathbb D^{K,\infty}_S\Gamma_Z^{\vee,O}(M,F)\xrightarrow{\sim}\Gamma_Z^O\mathbb D^{K,\infty}_S(M,F)
\end{equation*}
which gives the transformation map in $C_{\mathcal D^{\infty}fil}(S)$
\begin{eqnarray*}
T(D,\gamma^O)(M,F):\Gamma_Z^{\vee,O}\mathbb D^{\infty,K}_S(M,F)\to\mathbb D^{\infty,K}_S\Gamma_Z^O(M,F)
\end{eqnarray*}
\end{itemize}

In the analytic case, we have
\begin{defi}\label{Tgammainfty}
Let $S\in\AnSm(\mathbb C)$. For $(M,F)\in C_{\mathcal Dfil}(S)$, we have the map in $C_{\mathcal D^{\infty}fil}(S)$
\begin{eqnarray*}
T(\infty,\gamma)(M,F):=(I,T(j,\otimes)(-,-)): \\
J_S(\Gamma_Z(M,F)):=\Gamma_Z(M,F)\otimes_{D_S}(D_S^{\infty},F^{ord})\to
\Gamma_Z((M,F)\otimes_{D_S}(D_S^{\infty},F^{ord}))=:\Gamma_ZJ_S(M,F)
\end{eqnarray*}
\end{defi}

Let $i:Z\hookrightarrow S$ a closed embedding, with $Z,S\in\SmVar(\mathbb C)$ or $Z,S\in\AnSm(\mathbb C)$. 
We have the functor
\begin{eqnarray*}
i^{\sharp}:C_{\mathcal Dfil}(S)\to C_{\mathcal Dfil}(Z), 
(M,F)\mapsto i^{\sharp}(M,F):=\mathcal Hom_{i^*D_S}((D_{S\leftarrow Z},F^{ord}),i^*(M,F))
\end{eqnarray*}
where the (left) $D_Z$ module structure on $i^{\sharp}M$ comes from the right module structure on $D_{S\leftarrow Z}$, resp. $O_Z$.
We denote by
\begin{itemize}
\item for $(M,F)\in C_{\mathcal Dfil}(S)$, the canonical map in $C_{\mathcal Dfil}(S)$
\begin{eqnarray*}
\ad(i_{*mod},i^{\sharp})(M,F): i_{*mod}i^{\sharp}(M,F):=
i_*(\mathcal Hom_{i^*D_S}((D_{S\leftarrow Z},F^{ord}),i^*(M,F))\otimes_{D_Z}(D_{S\leftarrow Z},F^{ord})) \\
\to (M,F), \phi\otimes P\mapsto\phi(P)
\end{eqnarray*}
\item for $(N,F)\in C_{\mathcal Dfil}(Z)$, the canonical map in $C_{\mathcal Dfil}(Z)$
\begin{eqnarray*}
\ad(i_{*mod},i^{\sharp})(N,F):(N,F)\to i^{\sharp}i_{*mod}(N,F):=
\mathcal Hom_{i^*D_S}(D_{S\leftarrow Z},i^*i_*((N,F)\otimes_{D_Z}(D_{S\leftarrow Z},F^{ord}))) \\
n\mapsto(P\mapsto n\otimes P)
\end{eqnarray*}
\end{itemize}
The functor $i^{\sharp}$ induces in the derived category the functor :
\begin{eqnarray*}
Ri^{\sharp}:D_{\mathcal D(2)fil,r}(S)\to D_{\mathcal D(2)fil,r}(Z), \\  
(M,F)\mapsto Ri^{\sharp}(M,F):=R\mathcal Hom_{i^*D_S}((D_{Z\leftarrow S},F^{ord}),i^*(M,F))
=\mathcal Hom_{i^*D_S}((D_{Z\leftarrow S},F^{ord}),E(i^*(M,F))).
\end{eqnarray*}

\begin{prop}\label{isharp}
Let $i:Z\hookrightarrow S$ a closed embedding, with $Z,S\in\SmVar(\mathbb C)$ or $Z,S\in\AnSm(\mathbb C)$. 
The functor  $i_{*mod}:C_{\mathcal D}(Z)\to C_{\mathcal D}(S)$ admit a right adjoint
which is the functor $i^{\sharp}:C_{\mathcal D}(S)\to C_{\mathcal D}(Z)$ and
\begin{equation*}
\ad(i_{*mod},i^{\sharp})(N):N\to i^{\sharp}i_{*mod}N \mbox{\; and \;} \ad(i_{*mod},i^{\sharp})(M):i_{*mod}i^{\sharp}M\to M
\end{equation*}
are the adjonction maps. 
\end{prop}

\begin{proof}
See \cite{LvDmod} for the algebraic case. The analytic case is completely analogue.
\end{proof}

One of the main results in D modules is Kashiwara equivalence :
\begin{thm}\label{Keq}
\begin{itemize}
\item[(i)]Let $i:Z\hookrightarrow S$ a closed embedding with $Z,S\in\SmVar(\mathbb C)$. 
\begin{itemize}
\item The functor 
$i_{*mod}:\mathcal QCoh_{\mathcal D}(Z)\to\mathcal QCoh_{\mathcal D}(S)$ is an equivalence of category whose inverse is 
$i^{\sharp}:=a_{\tau}i^{\sharp}:\mathcal QCoh_{\mathcal D}(S)\to\mathcal QCoh_{\mathcal D}(Z)$.
That is, for $M\in QCoh_{\mathcal D}(S)$ and $N\in QCoh_{\mathcal D}(Z)$, the adjonction maps
\begin{equation*}
\ad(i_{*mod},i^{\sharp})(M):i_{*mod}i^{\sharp}M\xrightarrow{\sim} M \; , \; 
\ad(i_{*mod},i^{\sharp})(N):i^{\sharp}i_{*mod}N\xrightarrow{\sim} N
\end{equation*}
are isomorphisms.
\item The functor 
$\int_i=i_{*mod}:D_{\mathcal D}(Z)\to D_{\mathcal D}(S)$ is an equivalence of category whose inverse is 
$Ri^{\sharp}:D_{\mathcal D}(S)\to D_{\mathcal D}(Z)$.
That is, for $M\in D_{\mathcal D}(S)$ and $N\in D_{\mathcal D}(Z)$, the adjonction maps
\begin{equation*}
\ad(\int_i,Ri^{\sharp})(M):\int_iRi^{\sharp}M\xrightarrow{\sim} M \; , \; 
\ad(\int_i,Ri^{\sharp})(N):Ri^{\sharp}\int_i N\xrightarrow{\sim} N
\end{equation*}
are isomorphisms.
\end{itemize}
\item[(ii)]Let $i:Z\hookrightarrow S$ a closed embedding with $Z,S\in\AnSm(\mathbb C)$. 
\begin{itemize}
\item The functor 
$i_{*mod}:\mathcal Coh_{\mathcal D}(Z)\to\mathcal Coh_{\mathcal D}(S)$ is an equivalence of category whose inverse is 
$i^{\sharp}:=a_{\tau}i^{\sharp}:\mathcal Coh_{\mathcal D}(S)\to\mathcal Coh_{\mathcal D}(Z)$.
That is, for $M\in Coh_{\mathcal D}(S)$ and $N\in Coh_{\mathcal D}(Z)$, the adjonction maps
\begin{equation*}
\ad(i_{*mod},i^{\sharp})(M):i_{*mod}i^{\sharp}M\xrightarrow{\sim} M \; , \; 
\ad(i_{*mod},i^{\sharp})(N):i^{\sharp}i_{*mod}N\xrightarrow{\sim} N
\end{equation*}
are isomorphisms.
\item The functor 
$\int_i=i_{*mod}:D_{\mathcal D,c}(Z)\to D_{\mathcal D,c}(S)$ is an equivalence of category whose inverse is 
$Ri^{\sharp}:D_{\mathcal D,c}(S)\to D_{\mathcal D,c}(Z)$.
That is, for $M\in D_{\mathcal D,c}(S)$ and $N\in D_{\mathcal D,c}(Z)$, the adjonction maps
\begin{equation*}
\ad(\int_i,Ri^{\sharp})(M):\int_iRi^{\sharp}M\xrightarrow{\sim} M \; , \; 
\ad(\int_i,Ri^{\sharp})(N):Ri^{\sharp}\int_i N\xrightarrow{\sim} N
\end{equation*}
are isomorphisms.
\end{itemize}
\end{itemize}
\end{thm}

\begin{proof}
\noindent(i):Standard. Note that the second point follows from the first.

\noindent(ii):Standard. Note that the second point follows from the first.
\end{proof}

We have a canonical embedding of rings $D_Z\hookrightarrow D_{Z\to S}:=i^*D_S\otimes_{i^*O_S}O_Z$.
We denote by $C_{i^*\mathcal D,Z}(Z)$ the category whose objects are complexes of presheaves $M$ of $i^*D_S$ modules on $Z$
such that the cohomology presheaves $H^nM$ have an induced structure of $D_Z$ modules. We denote by
\begin{equation*}
q_K:K_{O_S}(i_*O_Z)\to i_*O_Z
\end{equation*} 
the Kozcul complex which is a resolution of the $O_S$ module $i_*O_Z$ of lenght $c=\codim(Z,S)$ by locally free sheaves of finite rank. 
The fact that it is a locally free resolution of finite rank comes from the fact that
$Z$ is a locally complete intersection in $S$ since both $Z$ and $S$ are smooth. We denote again 
\begin{equation*}
q_K=i^*q_K:K_{i^*O_S}(O_Z):=i^*K_{O_S}(i_*O_Z)\to i^*i_*O_Z=O_Z
\end{equation*}
We denote by $K^{\vee}_{i^*O_S}(O_Z):=\mathcal Hom_{i^*O_S}(K_{i^*O_S}(O_Z),i^*O_S)$ its dual, so that we have a canonical map 
\begin{equation*}
q_K^{\vee}:K^{\vee}_{i^*O_S}(O_Z)\to O_Z[-c].
\end{equation*}
Let $M\in C_{\mathcal D}(S)$. The $i^*D_S$ module structure on 
$\mathcal Hom_{i^*O_S}(K_{i^*O_S}(O_Z),i^*M)$ and $K_{i^*O_S}(O_Z)\otimes_{i^*O_S}i^*M$ 
induce a canonical $D_Z$ module structure on the cohomology groups 
$H^n\mathcal Hom_{i^*O_S}(K_{i^*O_S}(O_Z),i^*M)$ and $H^n(K_{i^*O_S}(O_Z)\otimes_{i^*O_S}i^*M)$ for all $n\in\mathbb Z$.

The projection formula for ringed spaces (proposition \ref{projformula}) implies the following lemma :

\begin{lem}\label{imodj}
Let $i:Z\hookrightarrow S$ a closed embedding with $Z,S\in\Var(\mathbb C)$ or with $Z,S\in\AnSp(\mathbb C)$. 
Denote by $j:U:=S\backslash Z\hookrightarrow Z$ the open complementary embedding.
Then, if $i$ is a locally complete intersection embedding (e.g. if $Z,S$ are smooth),
we have for $M\in C_{O_U}(U)$, $Li^{*mod}Rj_*M=0$.
\end{lem}

\begin{proof}
We have
\begin{eqnarray*}
i_*Li^{*mod}Rj_*M:=i_*(i^*L_O(j_*E(M))\otimes_{i^*O_S}O_Z)\xrightarrow{T(i,\otimes)(L_O(j_*E(M)),O_Z)^{-1}}
L_O(j_*E(M))\otimes_{O_S}i_*O_Z \\ \xrightarrow{q\circ (i_*q_K)^{-1}}(j_*E(M))\otimes_{O_S}i_*K_{i^*O_S}(O_Z)
\xrightarrow{T(j,\otimes)(E(M),K_{O_S}(i_*O_Z))}j_*(E(M)\otimes_{O_U}j^*K_{O_S}(i_*O_Z)),
\end{eqnarray*}
$T(i,\otimes)(L_O(j_*E(M)),O_Z)$ being an equivalence Zariski, resp. usu, local by proposition \ref{Tiotimes} and 
follows from the fact that $j^*K_{O_S}(i_*O_Z)$ is acyclic. But 
\begin{equation*}
T(j,\otimes)(E(M),K_{O_S}(i_*O_Z)):(j_*E(M))\otimes_{O_S}K_{O_S}(i_*O_Z)\to j_*(E(M)\otimes_{O_U}j^*i_*K_{i^*O_S}(O_Z))
\end{equation*}
is an equivalence Zariski, resp. usu, local by proposition \ref{projformula}
since $K_{O_S}(i_*O_Z)$ is a finite complex of locally free $O_S$ modules of finite rank.
\end{proof}

We deduce from theorem \ref{Keq}(i) and lemma \ref{imodj} the localization for $D$-modules 
for a closed embedding of smooth algebraic varieties:

\begin{thm}\label{KZS}
Let $i:Z\hookrightarrow S$ a closed embedding with $Z,S\in\SmVar(\mathbb C)$. Denote by $c=\codim(Z,S)$.
Then, for $M\in C_{\mathcal D}(S)$, we have by Kashiwara equivalence the following map in $C_{\mathcal D}(S)$ :
\begin{eqnarray*}
\mathcal K_{Z/S}(M):\Gamma_ZE(M)
\xrightarrow{\ad(i_{*mod},i^{\sharp})(-)^{-1}}i_{*mod}i^{\sharp}\Gamma_ZE(M) \\
\xrightarrow{\gamma_Z(-)}i_{*mod}i^{\sharp}(E(M))
\xrightarrow{\mathcal Hom(q_K,E(i^*M))\circ\mathcal Hom(O_Z,T(i,E)(M)}
i_{*mod}K^{\vee}_{i^*O_S}(O_Z)\otimes_{i^*O_S} M 
\end{eqnarray*}
which is an equivalence Zariski local. It gives the isomorphism in $D_{\mathcal D}(S)$
\begin{equation*}
\mathcal K_{Z/S}(M):R\Gamma_ZM\to i_{*mod}K^{\vee}_{i^*O_S}(O_Z)=i_{*mod}Li^{*mod}M[c]
\end{equation*}
\end{thm}

\begin{proof}
Follows from theorem \ref{Keq}(i) and lemma \ref{imodj} : see \cite{LvDmod} for example.
\end{proof}

\begin{defi}\label{inverseGamma}
Let $f:X\to S$ be a morphism, with $X,S\in\SmVar(\mathbb C)$ or $X,S\in\AnSm(\mathbb C)$.
Consider the factorization $f:X\xrightarrow{i} X\times S\xrightarrow{p_S}S$, 
where $i$ is the graph embedding and $p_S:X\times S\to S$ is the projection.
\begin{itemize}
\item[(i)]Then, for $(M,F)\in C_{\mathcal D(2)fil}(S)$ we set
\begin{equation*}
f^{*mod[-],\Gamma}(M,F):=\Gamma_XE(p_S^{*mod[-]}(M,F))\in C_{\mathcal D(2)fil,\infty}(X\times S), 
\end{equation*}
It induces in the derived category
\begin{equation*}
Rf^{*mod[-],\Gamma}(M,F):=f^{*mod[-],\Gamma}(M,F):=\Gamma_XE(p_S^{*mod[-]}(M,F))\in D_{\mathcal D(2)fil,\infty}(X\times S), 
\end{equation*}
By definition-proposition \ref{KZS}, we have in the algebraic case 
$Li^{*mod}f^{*mod,\Gamma}M=Lf^{*mod}M\in D_{\mathcal D}(X)$.
\item[(ii)]Then, for $(M,F)\in C_{\mathcal D(2)fil}(S)$ we set
\begin{equation*}
Lf^{\hat*mod[-],\Gamma}(M,F):=\Gamma^{\vee,h}_XL_Dp_S^{*mod[-]}(M,F)
:=\mathbb D_S^KL_D\Gamma_XE(\mathbb D_S^KL_Dp_S^{*mod[-]}(M,F))
\in D_{\mathcal D(2)fil,\infty}(X\times S). 
\end{equation*}
\end{itemize}
\end{defi}

In the analytical case we also have

\begin{defi}\label{inverseGammainfty}
Let $f:X\to S$ be a morphism, with $X,S\in\AnSm(\mathbb C)$.
Consider the factorization $f:X\xrightarrow{i} X\times S\xrightarrow{p_S}S$, 
where $i$ is the graph embedding and $p_S:X\times S\to S$ is the projection.
\begin{itemize}
\item[(i)]Then, for $(M,F)\in C_{\mathcal D^{\infty}(2)fil}(S)$ we set
\begin{equation*}
f^{*mod[-],\Gamma}(M,F):=\Gamma_XE(p_S^{*mod[-]}(M,F))\in C_{\mathcal D^{\infty}(2)fil,\infty}(X\times S), 
\end{equation*}
It induces in the derived category
\begin{equation*}
Rf^{*mod[-],\Gamma}(M,F):=f^{*mod[-],\Gamma}(M,F):=\Gamma_XE(p_S^{*mod[-]}(M,F))\in D_{\mathcal D^{\infty}(2)fil,\infty}(X\times S), 
\end{equation*}
\item[(ii)]Then, for $(M,F)\in C_{\mathcal D^{\infty}(2)fil}(S)$ we set
\begin{equation*}
Lf^{\hat*mod[-],\Gamma}(M,F):=\Gamma^{\vee,h}_XL_Dp_S^{*mod[-]}(M,F)
:=\mathbb D_S^KL_D\Gamma_XE(\mathbb D_S^KL_Dp_S^{*mod[-]}(M,F))
\in D_{\mathcal D^{\infty}(2)fil,\infty}(X\times S). 
\end{equation*}
\end{itemize}
\end{defi}

\subsubsection{The 2 functors and transformations maps for D modules on the smooth complex algebraic varieties 
and the complex analytic manifolds}

By the definitions and 
the propositions \ref{compDmod}, \ref{compDmodh}, \ref{compDmodDRd}, for the algebraic case,  
and the propositions \ref{compAnDmod}, \ref{compAnDmodh}, \ref{compAnDmodDRd}, for the analytic case,
\begin{itemize}
\item we have the 2 functors on $\SmVar(\mathbb C)$ :
\begin{itemize}
\item $C_{\mathcal D(2)fil}(\cdot):\SmVar(\mathbb C)\to C_{\mathcal D(2)fil}(\cdot), \; 
S\mapsto C_{\mathcal D(2)fil}(S), \, (f:T\to S)\mapsto f^{*mod}$, $(f:T\to S)\mapsto f^{*mod[-]}$, 
\item $D_{\mathcal D(2)fil,r}(\cdot):\SmVar(\mathbb C)\to D_{\mathcal D(2)fil,r}(\cdot), \;
S\mapsto D_{\mathcal D(2)fil,r}(S), \, (f:T\to S)\mapsto Lf^{*mod}$, $(f:T\to S)\mapsto Lf^{*mod[-]}$,
\item $D_{\mathcal D(2)fil,\infty}(\cdot):\SmVar(\mathbb C)\to D_{\mathcal D(2)fil,\infty}(\cdot), \; 
S\mapsto D_{\mathcal D(2)fil,\infty}(S), \, (f:T\to S)\mapsto \int^{FDR}_f$, 
\end{itemize}
\item we have the 2 functors on $\AnSm(\mathbb C)$ :
\begin{itemize}
\item $C_{\mathcal D(2)fil}(\cdot):\AnSm(\mathbb C)\to C_{\mathcal D(2)fil}(\cdot), \; 
S\mapsto C_{\mathcal D(2)fil}(S), \, (f:T\to S)\mapsto f^{*mod}$, $(f:T\to S)\mapsto f^{*mod[-]}$, 
\item $D_{\mathcal D(2)fil,r}(\cdot):\AnSm(\mathbb C)\to D_{\mathcal D(2)fil,r}(\cdot), \;
S\mapsto D_{\mathcal D(2)fil,r}(S), \, (f:T\to S)\mapsto Lf^{*mod}$, $(f:T\to S)\mapsto Lf^{*mod[-]}$,
\item $D_{\mathcal D(2)fil,\infty}(\cdot):\AnSm(\mathbb C)\to D_{\mathcal D(2)fil,\infty}(\cdot), \; 
S\mapsto D_{\mathcal D(2)fil,\infty}(S), \, (f:T\to S)\mapsto \int^{FDR}_{f!}$, 
\end{itemize}
\item we have also the 2 functors on $\AnSm(\mathbb C)$ :
\begin{itemize}
\item $C_{\mathcal D^{\infty}(2)fil}(\cdot):\AnSm(\mathbb C)\to C_{\mathcal D^{\infty}(2)fil}(\cdot), \; 
S\mapsto C_{\mathcal D^{\infty}(2)fil}(S), \, (f:T\to S)\mapsto f^{*mod}$, $(f:T\to S)\mapsto f^{*mod[-]}$, 
\item $D_{\mathcal D^{\infty}(2)fil,r}(\cdot):\AnSm(\mathbb C)\to D_{\mathcal D^{\infty}(2)fil,r}(\cdot), \;
S\mapsto D_{\mathcal D^{\infty}(2)fil,r}(S), \, (f:T\to S)\mapsto Lf^{*mod}$, $(f:T\to S)\mapsto Lf^{*mod[-]}$,
\item $D_{\mathcal D^{\infty}(2)fil,r}(\cdot):\AnSm(\mathbb C)\to D_{\mathcal D^{\infty}(2)fil,r}(\cdot), \; 
S\mapsto D_{\mathcal D^{\infty}(2)fil,r}(S), \, (f:T\to S)\mapsto \int^{FDR}_{f!}$,
\end{itemize}
\end{itemize}
inducing the following commutative diagrams of functors :
\begin{eqnarray*}
\xymatrix{\SmVar(\mathbb C)\ar[rr]^{f\mapsto f^{*mod}}\ar[d]^{\An} & \, & 
C_{\mathcal D(2)fil}(\cdot)\ar[d]^{an}\ar[rd] \\
\AnSm(\mathbb C)\ar[rr]^{f\mapsto f^{*mod}} & \, & 
C_{\mathcal D(2)fil}(\cdot)\ar[r]^{J_{\cdot}} & C_{\mathcal D^{\infty}(2)fil}(\cdot)}, \;
\xymatrix{\SmVar(\mathbb C)\ar[rr]^{f\mapsto Lf^{*mod}}\ar[d]^{\An} & \, & 
D_{\mathcal D(2)fil,r}(\cdot)\ar[d]^{an}\ar[rd] \\
\AnSm(\mathbb C)\ar[rr]^{f\mapsto Lf^{*mod}} & \, & D_{\mathcal D(2)fil,r}(\cdot)\ar[r]^{J_{\cdot}}  & 
D_{\mathcal D^{\infty}(2)fil,r}(\cdot)}, 
\end{eqnarray*}
where, for $S\in\AnSm(\mathbb C)$, 
\begin{itemize}
\item $D_{\mathcal D(2)fil,\infty,rh}(S)\subset D_{\mathcal D(2)fil,\infty,h}(S)$ is the full subcategory
consisting of filtered complexes of $D_S$ module whose cohomology sheaves are regular holonomic, 
\item $J:C_{\mathcal D(2)fil}(S)\to C_{\mathcal D^{\infty}(2)fil}(S)$ 
is the functor $(M,F)\mapsto J(M,F):=(M,F)\otimes_{D_S}D_S^{\infty}$, which derive trivially.
\end{itemize}

We first look at the pullback map and the transformation map of De Rahm complexes 
(see definition \ref{TDwM} and definition-proposition \ref{TDwMgamma}) together with the support section functor :

\begin{prop}\label{TDhwMgamma}
Consider a commutative diagram and a factorization
\begin{equation*}
D_0=\xymatrix{X\ar[r]^{f} & S \\
X'\ar[u]^{g'}\ar[r]^{f'} & T\ar[u]^{g} }
D_0=\xymatrix{f: X\ar[r]^{i} & Y\times S\ar[r]^{p} & S \\
f':X'\ar[r]^{i'}\ar[u]^{g'} & Y\times T\ar[u]^{g''=I\times g}\ar[r]^{p'} & T\ar[u]^{g} }
\end{equation*}
with $X,X',Y,S,T\in\Var(\mathbb C)$ or $X,X',Y,S,T\in\AnSp(\mathbb C)$, $i$, $i'$ being closed embeddings,
and $p$, $p'$ the projections. Denote by $D$ the right square of $D_0$. 
We have a factorization $i':X'\xrightarrow{i'_1}X_T=X\times_{Y\times S}Y\times T\xrightarrow{i'_0}Y\times T$, 
where $i'_0,i'_1$ are closed embedding. Assume $S,T,Y,Y'$ are smooth.
\begin{itemize}
\item[(i)] For $(M,F)\in C_{\mathcal Dfil}(Y\times S)$, the canonical map in $C_{p^{'*}O_Tfil}(Y\times T)$ 
(c.f. definition-proposition \ref{TDwMgamma}),
\begin{eqnarray*}
E(\Omega_{((Y'\times T)/(X\times S))/(T/S)}(M,F))\circ T(g'',E)(-)\circ T(g'',\gamma)(-): \\
g^{''*}\Gamma_{X}E((\Omega^{\bullet}_{Y\times S},F_b)\otimes_{O_{Y\times S}}(M,F))\to 
\Gamma_{X_T}E((\Omega^{\bullet}_{Y\times T/T},F_b)\otimes_{O_{Y\times T}}g^{''*mod}(M,F))
\end{eqnarray*}
is a map in $C_{p^{'*}\mathcal Dfil}(Y\times T)$.
\item[(ii)] For $(M,F)\in C_{\mathcal Dfil}(Y\times S)$, the canonical map in $C_{O_{T}fil}(T)$ 
(c.f. definition-proposition \ref{TDwMgamma} with $L_D$ instead of $L_O$)
\begin{equation*}
T^O_{\omega}(D)(M,F)^{\gamma}:g^{*mod}L_Dp_*\Gamma_{X}E(\Omega^{\bullet}_{Y\times S}\otimes_{O_{Y\times S}}(M,F))\to 
p'_*\Gamma_{X_T}E(\Omega^{\bullet}_{Y\times T/T}\otimes_{O_{Y\times T}}g^{''*mod}(M,F))
\end{equation*}
is a map in $C_{\mathcal Dfil}(T)$.
\item[(iii)] For $(N,F)\in C_{\mathcal Dfil}(Y\times T)$, the canonical map in $C_{p^{'*}O_Tfil}(Y\times T)$
\begin{equation*}
T(X'/X_T,\gamma)(-):\Gamma_{X'}E((\Omega^{\bullet}_{Y\times T/T},F_b)\otimes_{O_{Y\times T}}(N,F))\to
\Gamma_{X_T}E((\Omega^{\bullet}_{Y\times T/T},F_b)\otimes_{O_{Y\times T}}(N,F))
\end{equation*}
is a map in $C_{p^{'*}\mathcal D0fil}(Y\times T)$ .
\item[(iv)] For $M=O_Y$, we have $T_{\omega}^O(D)(O_{Y\times S})^{\gamma}=T_{\omega}^O(D)^{\gamma}$ as complexes of $D_T$ modules 
and $T_{\omega}^O(X_T/Y\times T)(O_{Y\times T})^{\gamma}=T_{\omega}^O(X_T/Y\times T)^{\gamma}$.
as complexes of $p^{'*}D_T$ modules.
\end{itemize}
\end{prop}

\begin{proof}
Follows by definition from proposition \ref{TDhwM}.
\end{proof}

In the analytical case, we also have

\begin{prop}\label{TDhwMgammainfty}
Consider a commutative diagram and a factorization
\begin{equation*}
D_0=\xymatrix{X\ar[r]^{f} & S \\
X'\ar[u]^{g'}\ar[r]^{f'} & T\ar[u]^{g} }
D_0=\xymatrix{f: X\ar[r]^{i} & Y\times S\ar[r]^{p} & S \\
f':X'\ar[r]^{i'}\ar[u]^{g'} & Y\times T\ar[u]^{g''=I\times g}\ar[r]^{p'} & T\ar[u]^{g} }
\end{equation*}
with $X,X',Y,S,T\in\AnSp(\mathbb C)$, $i$, $i'$ being closed embeddings,
and $p$, $p'$ the projections. Denote by $D$ the right square of $D_0$. 
We have a factorization $i':X'\xrightarrow{i'_1}X_T=X\times_{Y\times S}Y\times T\xrightarrow{i'_0}Y\times T$, 
where $i'_0,i'_1$ are closed embedding. Assume $S,T,Y,Y'$ are smooth.
\begin{itemize}
\item[(i)] For $(M,F)\in C_{\mathcal D^{\infty}fil}(Y\times S)$, the canonical map in $C_{p^{'*}O_Tfil}(Y\times T)$ 
(c.f. definition-proposition \ref{TDwMgamma}),
\begin{eqnarray*}
E(\Omega_{((Y'\times T)/(X\times S))/(T/S)}(M,F))\circ T(g'',E)(-)\circ T(g'',\gamma)(-): \\
g^{''*}\Gamma_{X}E((\Omega^{\bullet}_{Y\times S},F_b)\otimes_{O_{Y\times S}}(M,F))\to 
\Gamma_{X_T}E((\Omega^{\bullet}_{Y\times T/T},F_b)\otimes_{O_{Y\times T}}g^{''*mod}(M,F))
\end{eqnarray*}
is a map in $C_{p^{'*}\mathcal D^{\infty}fil}(Y\times T)$.
\item[(ii)] For $(M,F)\in C_{\mathcal D^{\infty}fil}(Y\times S)$, the canonical map in $C_{O_{T}fil}(T)$ 
(c.f. definition-proposition \ref{TDwMgamma} with $L_{D^{\infty}}$ instead of $L_O$)
\begin{equation*}
T^O_{\omega}(D)(M,F)^{\gamma}:
g^{*mod}L_{D^{\infty}}p_*\Gamma_{X}E((\Omega^{\bullet}_{Y\times S},F_b)\otimes_{O_{Y\times S}}(M,F))\to 
p'_*\Gamma_{X_T}E((\Omega^{\bullet}_{Y\times T/T},F_b)\otimes_{O_{Y\times T}}g^{''*mod}(M,F))
\end{equation*}
is a map in $C_{\mathcal D^{\infty}fil}(T)$.
\item[(iii)] For $(N,F)\in C_{\mathcal D^{\infty}fil}(Y\times T)$, the canonical map in $C_{p^{'*}O_Tfil}(Y\times T)$
\begin{equation*}
T(X'/X_T,\gamma)(-):\Gamma_{X'}E(\Omega^{\bullet}_{Y\times T/T}\otimes_{O_{Y\times T}}(N,F))\to
\Gamma_{X_T}E(\Omega^{\bullet}_{Y\times T/T}\otimes_{O_{Y\times T}}(N,F))
\end{equation*}
is a map in $C_{p^{'*}\mathcal D^{\infty}fil}(Y\times T)$ .
\item[(iv)] For $M=O_Y$, we have $T_{\omega}^O(D)(O_{Y\times S})^{\gamma}=T_{\omega}^O(D)^{\gamma}$ 
as complexes of $D^{\infty}_T$ modules 
and $T_{\omega}^O(X_T/Y\times T)(O_{Y\times T})^{\gamma}=T_{\omega}^O(X_T/Y\times T)^{\gamma}$.
as complexes of $p^{'*}D^{\infty}_T$ modules.
\end{itemize}
\end{prop}

\begin{proof}
Follows from proposition \ref{TDhwMgamma}.
\end{proof}

Similarly, we have :

\begin{prop}\label{TDhanwMgamma}
Let $p:Y\times S\to S$ a projection and $i:X\hookrightarrow Y\times S$ a closed embedding with $S,Y\in\SmVar(\mathbb C)$. 
\begin{itemize}
\item[(i)] For $(M,F)\in C_{\mathcal Dfil}(Y\times S)$ the canonical map in $C_{p^*O_Sfil}(Y^{an}\times S^{an})$
 (see definition-proposition \ref{TDwMgamma})
\begin{eqnarray*}
E(\Omega_{(Y^{an}\times S^{an}/Y\times S)/(S^{an}/S)}(M,F))\circ T(an,\gamma)(-): \\
(\Gamma_XE((\Omega^{\bullet}_{Y\times S/S},F_b)\otimes_{O_Y}(M,F)))^{an}\to
\Gamma_{X^{an}}E((\Omega^{\bullet}_{Y\times S/S},F_b)\otimes_{O_{Y^{an}\times S^{an}}}(M,F)^{an})
\end{eqnarray*}
is a map in $C_{h^*\mathcal Dfil}(Y^{an}\times S)$.
\item[(ii)] For $(M,F)\in C_{\mathcal Dfil}(Y\times S)$ the canonical map in $C_{O_Sfil}(S^{an})$ 
(see definition-proposition \ref{TDwMgamma})
\begin{equation*}
T_{\omega}^O(an,p)(M,F)^{\gamma}:(p_*\Gamma_XE((\Omega^{\bullet}_{Y\times S/S},F_b)\otimes_{O_Y}(M,F)))^{an}\to
p_*\Gamma_{X^{an}}E((\Omega^{\bullet}_{Y\times S/S},F_b)\otimes_{O_{Y^{an}}}(M,F)^{an})
\end{equation*}
is a map in $C_{\mathcal Dfil}(S^{an})$.
\item[(iii)] For $M=O_Y$, we have $T_{\omega}^O(an,h)(O_Y)^{\gamma}=T_{\omega}^O(an,h)^{\gamma}$ as complexes of $D_S$ modules
\end{itemize}
\end{prop}

\begin{proof}
Follows by definition from proposition \ref{TDhanwM}
\end{proof}

Let $p:Y\times S\to S$ a projection with $Y,S\in\SmVar(\mathbb C)$ or with $Y,S\in\AnSm(\mathbb C)$. 
Let $j:V\hookrightarrow Y\times S$ an open embedding.
Consider (see proposition \ref{projformulahw}), for $(M,F)\in C_{\mathcal Dfil}(Y\times S)$,
the canonical transformation map in $C_{p^*\mathcal Dfil}(Y\times S)$ 
\begin{eqnarray*}
k\circ T^O_w(j,\otimes)(E(M,F)):(\Omega^{\bullet}_{Y\times S/S},F_b)\otimes_{O_Y}j_*j^*E(M,F) \\
\xrightarrow{DR(Y\times S/S)(\ad(j^*,j_*)(-))} \\
j_*j^*((\Omega^{\bullet}_{Y\times S/S},F_b)\otimes_{O_Y}j_*j^*E(M,F))=
j_*j^*(\Omega^{\bullet}_{Y\times S/S},F_b)\otimes_{O_Y}j^*j_*j^*E(M,F) \\
\xrightarrow{k\circ DR(Y\times S/S)(\ad(j^*,j_*)(E(M,F)))} \\
j_*E(j^*(\Omega^{\bullet}_{Y\times S/S},F_b)\otimes_{O_{Y\times S}}j^*E(M,F))
=j_*E(j^*((\Omega^{\bullet}_{Y\times S/S},F_b)\otimes_{O_{Y\times S}}E(M,F)))
\end{eqnarray*}
We have then :

\begin{prop}\label{mw0prop}
Let $p:Y\times S\to S$ a projection with $Y,S\in\SmVar(\mathbb C)$ or with $Y,S\in\AnSm(\mathbb C)$. 
Let $i:X\hookrightarrow Y\times S$ a closed embedding. Then, for $(M,F)\in C_{\mathcal Dfil}(Y\times S)$
\begin{itemize}
\item[(i)] the canonical map in $C_{p^*\mathcal Dfil}(Y)$ (definition \ref{mw0def}) 
\begin{eqnarray*}
T^O_w(\gamma,\otimes)(M,F):=(I,k\circ T^O_w(j,\otimes)(E(M,F))):\\
(\Omega^{\bullet}_{Y\times S/S},F_b)\otimes_{O_{Y\times S}}\Gamma_XE(M,F)\to
\Gamma_XE((\Omega^{\bullet}_{Y\times S/S},F_b)\otimes_{O_{Y\times S}}E(M,F)),
\end{eqnarray*}
is a ($1$-)filtered Zariski, resp usu, local equivalence.
\item[(ii)] the map of point (i) gives the following canonical isomorphism in $D_{p^*\mathcal Dfil}(Y)$ 
\begin{eqnarray*}
T^O_w(\gamma,\otimes)(M,F):(\Omega^{\bullet}_{Y\times S/S},F_b)\otimes_{O_{Y\times S}}\Gamma_XE(M,F)
\xrightarrow{T^O_w(\gamma,\otimes)(M,F)} \\
\Gamma_XE((\Omega^{\bullet}_{Y\times S/S},F_b)\otimes_{O_{Y\times S}}E(M,F)) 
\xrightarrow{DR(Y\times S/S)(k)^{-1}} 
\Gamma_XE((\Omega^{\bullet}_{Y\times S/S},F_b)\otimes_{O_{Y\times S}}(M,F)).
\end{eqnarray*}
\end{itemize}
\end{prop}

\begin{proof}
\noindent(i): By proposition \ref{projformulahw} 
\begin{itemize}
\item $T^O_w(j,\otimes)(M,F):(\Omega^{\bullet}_{Y\times S/S},F_b)\otimes_{O_{Y\times S}}j_*j^*E(M,F)\to 
j_*E(j^*((\Omega^{\bullet}_{Y\times S/S},F_b)\otimes_{O_{Y\times S}}E(M,F)))$
is a filtered Zariski, resp usu, local equivalence in $C_{p^*\mathcal Dfil}(Y\times S)$ and
\item $DR(Y\times S/S)(k):(\Omega^{\bullet}_{Y\times S/S},F_b)\otimes_{O_{Y\times S}}(M,F)\to
(\Omega^{\bullet}_{Y\times S/S},F_b)\otimes_{O_{Y\times S}}E(M,F) $
is a filtered Zariski, resp usu, local equivalence in $C_{p^*\mathcal Dfil}(Y\times S)$.
\end{itemize}

\noindent(ii): Follows from (i).
\end{proof}

In the analytic case, we also have
\begin{prop}\label{mw0propinfty}
Let $p:Y\times S\to S$ a projection with $Y,S\in\AnSm(\mathbb C)$. 
Let $i:X\hookrightarrow Y$ a closed embedding. Then, for $(M,F)\in C_{\mathcal D^{\infty}fil}(Y\times S)$ 
\begin{itemize}
\item[(i)] the canonical map in $C_{p^*\mathcal Dfil}(Y)$ 
\begin{eqnarray*}
T^O_w(\gamma,\otimes)(M,F):=(I,T^O_w(j,\otimes)(E(M,F))):\\
(\Omega^{\bullet}_{Y\times S/S},F_b)\otimes_{O_{Y\times S}}\Gamma_XE(M,F)\to
\Gamma_XE((\Omega^{\bullet}_{Y\times S/S},F_b)\otimes_{O_{Y\times S}}E(M,F)).
\end{eqnarray*}
is a map in $C_{p^*\mathcal D^{\infty}fil}(Y\times S)$. Proposition \ref{mw0prop} says that it is a filtered equivalence usu local,
\item[(ii)] the map of point (i) gives the following canonical isomorphism in $D_{p^*\mathcal D^{\infty},fil}(Y\times S)$ 
\begin{eqnarray*}
T^O_w(\gamma,\otimes)(M,F):(\Omega^{\bullet}_{Y\times S/S},F_b)\otimes_{O_{Y\times S}}\Gamma_XE(M,F)
\xrightarrow{T^O_w(\gamma,\otimes)(M,F)} \\
\Gamma_XE((\Omega^{\bullet}_{Y\times S/S},F_b)\otimes_{O_{Y\times S}}E(M,F)) 
\xrightarrow{DR(Y\times S/S)(k)^{-1}} 
\Gamma_XE((\Omega^{\bullet}_{Y\times S/S},F_b)\otimes_{O_{Y\times S}}(M,F)).
\end{eqnarray*}
\end{itemize}
\end{prop}

\begin{proof}
\noindent(i): By proposition \ref{projformulahwinfty} 
\begin{itemize}
\item $T^O_w(j,\otimes)(M):\Omega^{\bullet}_{Y\times S/S}\otimes_{O_{Y\times S}}j_*j^*E(M)\to 
j_*E(j^*(\Omega^{\bullet}_{Y\times S/S}\otimes_{O_{Y\times S}}E(M)))$
is an equivalence usu local in $C_{p^*\mathcal D^{\infty}}(Y\times S)$ and
\item $DR(Y\times S/S)(k):\Omega^{\bullet}_{Y\times S/S}\otimes_{O_{Y\times S}} M\to
\Omega^{\bullet}_{Y\times S/S}\otimes_{O_{Y\times S}}E(M) $
is an equivalence usu local in $C_{p^*\mathcal D^{\infty}}(Y\times S)$.
\end{itemize}

\noindent(ii): Follows from (i).
\end{proof}

In the projection case, we consider the following canonical maps :
Let $S_1,S_2\in\SmVar(\mathbb C)$ or let $S_1,S_2\in\AnSm(\mathbb C)$.
Denote by $p=p_1:S_{12}=S_1\times S_2\to S_1$ and $p_2:S_{12}=S_1\times S_2\to S_1$ the projection. We consider
\begin{itemize}
\item $p(M_1,F):(M_1,F)\to p_{*mod}p^{*mod[-]}(M_1,F)$ in $C_{\mathcal D(2)fil}(S_1)$, 
for $(M_1,F)\in C_{\mathcal Dfil(2)}(S_1)$, which is the composite 
\begin{eqnarray*}
p(M_1,F):(M_1,F)\xrightarrow{\ad(p^*,p_*)(M_1)}p_*p^*(M_1,F)
\xrightarrow{m_1}p_*((\Omega^{\bullet}_{S_{12}/S_1},F_b)\otimes_{p^*O_{S_1}}p^*(M_1,F)) \\
\xrightarrow{=}p_*((\Omega^{\bullet}_{S_{12}/S_1},F_b)\otimes_{O_{S_{12}}}p^{*mod}(M_1,F)=p_{*mod}p^{*mod[-]}(M_1,F)
\end{eqnarray*}
where $m_1:p^*M_1\to p^*M_1\otimes_{p^*O_{S_1}}\Omega^{\bullet}_{S_{12}/S_1}$ is given by $m_1(m)=m\otimes 1$,
\item $p(M_{12},F):p^{*mod[-]}p_{*mod}(M_{12},F)\to (M_{12},F)$ in $C_{\mathcal Dfil}(S_1\times S_2)$, 
for $(M_{12},F)\in C_{\mathcal D}(S_1\times S_2)$, which is the composite
\begin{eqnarray*}
p(M_{12},F):p^{*mod[-]}p_{*mod}(M_{12},F)=
p^*p_*((M_{12,F)}\otimes_{O_{S_{12}}}(\Omega^{\bullet}_{S_{12}/S_1},F_b))\otimes_{p^*O_{S_1}}O_{S_{12}} \\
\xrightarrow{\ad(p^*,p_*)(-)\otimes_{p^*O_{S_1}}I}  
(M_{12},F)\otimes_{O_{S_{12}}}\Omega^{\bullet}_{S_{12}/S_1}\otimes_{p^*O_{S_1}}O_{S_{12}}
=(M_{12},F)\otimes_{p^*O_{S_1}}\Omega^{\bullet}_{S_{12}/S_1}\xrightarrow{m_{12}}(M_{12},F)
\end{eqnarray*}
where $m_{12}:M_{12}\otimes_{p^*O_{S_1}}\Omega^{\bullet}_{S_{12}/S_1}\to M_{12}$ is the multiplication map:
\begin{itemize}
\item $m_{12}(M_{12}\otimes_{p^*O_{S_1}}\Omega^p_{S_{12}/S_1})=0$ for $p\neq 0$ and 
\item $m_{12}:M_{12}\otimes_{p^*O_{S_1}}\Omega^0_{S_{12}/S_1}=M_{12}\otimes_{p^*O_{S_1}}O_{S_{12}}\to M_{12}$ is 
given by $m_{12}(m\otimes f)=fm$ 
\end{itemize}
\end{itemize}
We have then $p(p^{*mod[-]}(M_1,F))\circ p^{*mod[-]}p(M_1,F)=I_{p^{*mod[-]}(M_1,F)}$.  
It gives the following maps
\begin{itemize}
\item $p_!(M_{12}):(M_{12},F)\to p^{*mod[-]}\int_{p!}(M_{12,F)}$ in $D_{\mathcal D(2)fil}(S_1\times S_2)$, 
for $(M_{12},F)\in C_{\mathcal Dfil,h}(S_1\times S_2)$, given by
\begin{eqnarray*}
p_!(M_{12}): (M_{12},F)\xrightarrow{d(M_{12},F)}L\mathbb D_S^{2}(M_{12},F)
\xrightarrow{L\mathbb D_S(p(-)\circ q)} \mathbb D_S^KL_Dp^{*mod[-]}p_{*mod}E(\mathbb D^K_SL_D(M_{12},F)) \\
\xrightarrow{T(p,D)(-)^{-1}}p^{*mod}\mathbb D_S^KL_Dp_{*mod}E(\mathbb D_S^KL_D(M_{12},F))=p^{*mod[-]}\int_{p_!}(M_{12},F)
\end{eqnarray*}
\item $p_!(M_1,F):\int_{p!}p^{*mod[-]}(M_1,F)\to (M_1,F)$ in $D_{\mathcal Dfil}(S_1)$, 
for $M_1\in C_{\mathcal Dfil,h}(S_1)$, given by
\begin{eqnarray*}
p_!(M_1,F):\int_{p!}p^{*mod[-]}(M_1,F)=\mathbb D_S^KL_Dp_{*mod}E(\mathbb D_S^Kp^{*mod[-]}L_D(M_1,F))
\xrightarrow{(\mathbb D_S^Kk)\circ T(p,D)(-)^{-1}} \\
\mathbb D_S^Kp_{*mod}p^{*mod[-]}\mathbb D_S^KL_D(M_1,F)\xrightarrow{\mathbb D_S^Kp(\mathbb D_S^KL_D(M_1,F))}
\mathbb D_S^{K,2}L_D(M_1,F)\xrightarrow{d(M_1,F)^{-1}}(M_1,F)
\end{eqnarray*}
\end{itemize}
so that $p^{*mod[-]}(p_!(M_1,F))\circ p_!(p^{*mod[-]}(M_1,F))=I_{p^{*mod[-]}(M_1,F)}$.

\begin{defi}\label{TDmodlem}
\begin{itemize}
\item[(i)]Consider a commutative diagram in $\SmVar(\mathbb C)$ which is cartesian, together with its factorization
\begin{equation*}
D=(f,q)=\xymatrix{
X\times T\ar[r]^{f'}\ar[d]_{g'} & T\times S\ar[d]^{q} \\ X\ar[r]^{f} & S}
D=(f,q)=\xymatrix{
f'':X\times T\ar[r]^{i''}\ar[d]^{q'} & X\times T\times S\ar[r]^{p''}\ar[d]^{q''=I\times q} & T\times S\ar[d]^{q} \\
f:X\ar[r]^{i} & X\times S\ar[r]^{p} &  S},
\end{equation*}
where the squares are cartesian, $f=p\circ i$ being the graph factorization and $q$, $q'$ the projections.
We have, for $(M,F)\in C_{\mathcal D(2)fil,\infty}(X)$, 
the following transformation map in $C_{\mathcal D(2)fil,\infty}(T\times S)$ : 
\begin{eqnarray*}
T^{\mathcal Dmod}(f,q)(M,F):
q^{*mod}p_*E((\Omega_{X\times S/S},F_b)\otimes_{O_{X\times S}}i_{*mod}(M,F)) \\
\xrightarrow{T^O_{\omega}(q,p)(i_{*mod}(M,F))} 
p''_*E((\Omega_{X\times T\times S/T\times S},F_b)\otimes_{O_{X\times T\times S}}q^{''*mod}i_{*mod}(M,F)) \\
\xrightarrow{p_*E(T^{\mathcal Dmod}(i,q'')(M,F)\otimes I)}
p''_*E((\Omega_{X\times T\times S/T\times S},F_b)\otimes_{O_{X\times T\times S}}i''_{*mod}q^{'*mod}(M,F))
\end{eqnarray*}
where 
\begin{eqnarray*}
T^{\mathcal Dmod}(i,q'')(M,F):q^{''*mod}i_{*mod}(M,F):=q^{''*mod}i_*((M,F)\otimes_{D_X}i^{*mod}(D_{X\times S},F^{ord})) \\
\xrightarrow{T^{mod}(q'',i)(-)} 
i_*q^{'*mod}((M,F)\otimes_{D_X}i^{*mod}(D_{X\times S},F^{ord})) \\
\xrightarrow{=}
i_*(q^{'*mod}(M,F)\otimes_{q^{'*}D_X}q^{'*}i^{*mod}(D_{X\times S},F^{ord})) 
\xrightarrow{=} \\
i_*(q^{'*mod}(M,F)\otimes_{D_{X\times T}}i^{'*mod}(D_{X\times S\times T},F^{ord}))=:i'_{*mod}q^{'*mod}(M,F)
\end{eqnarray*}
\item[(ii)]Consider a commutative diagram in $\SmVar(\mathbb C)$ which is cartesian, together with its factorization
\begin{equation*}
D=(f,g)=\xymatrix{
X_T\ar[r]^{f'}\ar[d]_{g'} & T\ar[d]^{g} \\ X\ar[r]^{f} & S}
D=(f,g)=\xymatrix{ 
f':X_T\ar[r]^{i'}\ar[d]^{l'} & X\times T\ar[r]^{p'}\ar[d]^{l''=I\times l} & T\ar[d]^{l} \\
f'':X\times T\ar[r]^{i''}\ar[d]^{q'} & X\times T\times S\ar[r]^{p''}\ar[d]^{q''=I\times q} & T\times S\ar[d]^{q} \\
f:X\ar[r]^{i} & X\times S\ar[r]^{p} &  S},
\end{equation*}
where the squares are cartesian, $f=p\circ i$, $g=q\circ l$ being the graph factorizations.
We have, for $(M,F)\in D_{\mathcal D(2)fil,\infty}(X)$, 
the following transformation map in $D_{\mathcal D(2)fil,\infty}(T\times S)$ : 
\begin{eqnarray*}
T^{\mathcal Dmod}(f,g)((M,F)): \\
Rg^{*mod,\Gamma}(M,F)\int^{FDR}_{f}(M,F):=
\Gamma_TE(q^{*mod}p_*E((\Omega_{X\times S/S},F_b)\otimes_{O_{X\times S}}i_{*mod}(M,F))) \\ 
\xrightarrow{\Gamma_TE(T^{\mathcal Dmod}(f,q)(M,F))}
\Gamma_TE(p''_*E((\Omega_{X\times T\times S/T\times S},F_b)\otimes_{O_{X\times T\times S}}i''_{*mod}q^{'*mod}(M,F))) \\
\xrightarrow{=}
p''_*\Gamma_{X\times T}E((\Omega_{X\times T\times S/T\times S},F_b)\otimes_{O_{X\times T\times S}}i''_{*mod}q^{'*mod}(M,F)) \\
\xrightarrow{T^O_w(\gamma,\otimes)(-)}
p''_*E((\Omega_{X\times T\times S/T\times S},F_b)\otimes_{O_{X\times T\times S}}\Gamma_{X\times T}E(i''_{*mod}q^{'*mod}(M,F))) \\
\xrightarrow{=}
p''_*E((\Omega_{X\times T\times S/T\times S},F_b)\otimes_{O_{X\times T\times S}}(i''_{*mod}q^{'*mod}\Gamma_{X_T}E(M,F)))
=:\int^{FDR}_{f'}Rg^{'*mod,\Gamma}(M,F) 
\end{eqnarray*}
\item[(ii)'] We have, for $M\in D_{\mathcal D}(X)$, the following transformation map in $D_{\mathcal D}(T)$ : 
\begin{eqnarray*}
T^{\mathcal Dmod}(f,g)(M): \\
g^{*mod}(M,F)\int_{f}(M)=l^{*mod}q^{*mod}\int_fM
\xrightarrow{q'(M)}l^{*mod}q^{*mod}\int_fq'_{*mod}q^{'*mod}M \\
\xrightarrow{=}l^{*mod}q^{*mod}q_{*mod}\int_{f''}q^{'*mod}M
\xrightarrow{q(-)}l^{*mod}\int_{f''}q^{'*mod}M \\
\xrightarrow{l^{*mod}\ad(l^{'\sharp},l'_{*mod})(-)^{-1}}l^{*mod}\int_{f''}l'_{*mod}l^{'\sharp}q^{'*mod}M
\xrightarrow{=}l^{\sharp}l_{*mod}\int_{f'}l^{'*mod}q^{*mod}M \\
\xrightarrow{\ad(l^{\sharp},l_{*mod})(-)}
\int_{f'}l^{'*mod}q^{*mod}M=:\int_{f'}g^{'*mod}(M) 
\end{eqnarray*}
where $l^{*mod}\ad(l^{'\sharp},l'_{*mod})(-)$ is an isomorphism by lemma \ref{imodj}.
\end{itemize}
\end{defi}

In the analytic case, we have :

\begin{defi}\label{TDmodleman}
Consider a commutative diagram in $\AnSm(\mathbb C)$ which is cartesian together with a factorization
\begin{equation*}
D=(f,g)=\xymatrix{
X_T\ar[r]^{f'}\ar[d]_{g'} & T\ar[d]^{g} \\
X\ar[r]^{f} & S}
D=(f,g)=\xymatrix{ 
f':X_T\ar[r]^{i'}\ar[d]^{g'} & Y\times T\ar[r]^{p'}\ar[d]^{g''=I\times g} & T\ar[d]^{g} \\
f:X\ar[r]^{i} & Y\times S\ar[r]^{p} &  S},
\end{equation*}
where $Y\in\AnSm(\mathbb C)$, $i$, $i'$ are closed embeddings and $p$, $p'$ the projections.
\begin{itemize}
\item[(i)]We have, for $(M,F)\in D_{\mathcal D(2)fil,\infty,h}(X)$, 
the following transformation map in $D_{\mathcal D(2)fil,\infty}(T\times S)$ 
\begin{eqnarray*}
T^{\mathcal Dmod}(f,g)((M,F)):Rg^{*mod,\Gamma}\int^{FDR}_{f}(M,F)\to\int^{FDR}_{f'}Rg^{'*mod,\Gamma}(M,F) 
\end{eqnarray*}
define in the same way as in definition \ref{TDmodlem}
\item[(ii)] For $(M,F)\in D_{\mathcal D^{\infty}(2)fil,\infty}(X)$,
 the following transformation map in $D_{\mathcal D^{\infty}(2)fil,\infty}(T\times S)$ 
\begin{equation*}
T^{\mathcal Dmod}(f,g)((M,F)):Rg^{*mod,\Gamma}\int^{FDR}_{f}(M,F)\to\int^{FDR}_{f'}Rg^{'*mod,\Gamma}(M,F)
\end{equation*}
is defined in the same way as in (ii) : see definition \ref{TDmodlem}. 
\end{itemize}
\end{defi}

In the algebraic case, we have the following proposition:

\begin{prop}\label{PDmod1}
Consider a cartesian square in $\SmVar(\mathbb C)$ 
\begin{equation*}
D=\xymatrix{ 
X_T\ar[r]^{g'}\ar[d]_{f'} & X\ar[d]^{f} \\
T\ar[r]^{g} & S}
\end{equation*}
\begin{itemize}
\item[(i)]For $(M,F)\in D_{\mathcal D(2)fil,\infty,c}(X)$, 
\begin{equation*}
T^{\mathcal Dmod}(f,g)((M,F)):Rg^{*mod,\Gamma}\int^{FDR}_{f}(M,F)\xrightarrow{\sim}\int^{FDR}_{f'}Rg^{'*mod,\Gamma}(M,F)
\end{equation*}
is an isomorphism in $D_{\mathcal D(2)fil,\infty}(T)$.
\item[(ii)]For $M\in D_{\mathcal D,c}(X)$, 
\begin{equation*}
T^{\mathcal Dmod}(f,g)(M):g^{*mod}\int_{f}M\xrightarrow{\sim}\int_{f'}g^{'*mod}M
\end{equation*}
is an isomorphism in $D_{\mathcal D}(T)$.
\end{itemize}
\end{prop}

\begin{proof}
Follows from the projection case and the closed embedding case.
\end{proof}

In the analytic case, we have similarly:

\begin{prop}\label{PDmod1an}
Consider a cartesian square in $\AnSm(\mathbb C)$ 
\begin{equation*}
D=\xymatrix{ 
X_T\ar[r]^{g'}\ar[d]_{f'} & X\ar[d]^{f} \\
T\ar[r]^{g} & S}
\end{equation*}
\begin{itemize}
\item[(i)] Assume that $f$, hence $f'$ is proper. For $(M,F)\in D_{\mathcal D(2)fil,\infty,h}(X)$, 
\begin{equation*}
T^{\mathcal Dmod}(f,g)((M,F)):Rg^{*mod,\Gamma}\int^{FDR}_{f}(M,F)\xrightarrow{\sim}\int^{FDR}_{f'}Rg^{'*mod,\Gamma}(M,F)
\end{equation*}
is an isomorphism in $D_{\mathcal D(2)fil,\infty}(T)$.
\item[(ii)] For $(M,F)\in D_{\mathcal D^{\infty}(2)fil,\infty,h}(X)$, 
\begin{equation*}
T^{\mathcal Dmod}(f,g)((M,F)):Rg^{*mod,\Gamma}\int^{FDR}_{f}(M,F)\xrightarrow{\sim}\int^{FDR}_{f'}Rg^{'*mod,\Gamma}(M,F)
\end{equation*}
is an isomorphism in $D_{\mathcal D^{\infty}(2)fil,\infty}(T)$.
\end{itemize}
\end{prop}

\begin{proof}
\noindent(i):Similar to the proof of proposition \ref{PDmod1}.

\noindent(ii):Similar to the proof of proposition \ref{PDmod1}.
\end{proof}

\begin{defi}\label{projDmod}
Let $f:X\to S$ a morphism with $X,S\in\SmVar(\mathbb C)$.
\begin{itemize}
\item[(i)] We have, for $(M,F)\in C_{\mathcal Dfil}(S)$ and $(N,F)\in C_{\mathcal Dfil}(X)$, 
we have the map in $C_{\mathcal Dfil}(S)$
\begin{eqnarray*}
T^{\mathcal Dmod,0}(\otimes,f)((M,F),(N,F)): \\ 
(M,F)\otimes_{O_S}f^0_{*mod}(N,F):=(M,F)\otimes_{O_S}f_*((N,F)\otimes_{D_X}(D_{X\leftarrow S},F^{ord})) \\
\xrightarrow{T(\otimes,f)(-,-)}  
f_*(f^*(M,F)\otimes_{f^*O_S}(N,F)\otimes_{D_X}(D_{X\leftarrow S},F^{ord}))
\xrightarrow{=} \\
f_*(f^{*mod}(M,F)\otimes_{O_X}(N,F)\otimes_{D_X}(D_{X\leftarrow S},F^{ord}))=:f^0_{*mod}(f^{*mod}(M,F)\otimes_{O_X}(N,F))
\end{eqnarray*}
\item[(ii)]Consider the cartesian square
\begin{equation*}
D=\xymatrix{ 
X\ar[r]^{i}\ar[d]^{f} & X\times S\ar[d]^{f\times I_S} \\
S\ar[r]^{\Delta_S} & S\times S},
\end{equation*}
where $i_f=(f\times I_S)\circ\Delta_X:X\hookrightarrow X\times S$ is the graph embedding. Then,
for $(M,F)\in C_{\mathcal D(2)fil}(S)$ and $(N,F)\in C_{\mathcal Dfil}(X)$, we have the map in $D_{\mathcal D(2)fil,r}(S)$
\begin{eqnarray*} 
T^{\mathcal Dmod}(\otimes,f)((M,F),(N,F)):
\int^{FDR}_f((N,F)\otimes_{O_X}f_{FDR}^{*mod}(M,F))=\int^{FDR}_{f}i_{f,FDR}^{*mod}(p_X^*N\otimes p_S^*M) \\
\xrightarrow{T^{\mathcal Dmod}(\Delta_S,f\times I_S)(-)}
\Delta_{S,FDR}^{*mod}\int^{FDR}_{(f\times I_S)}(p_X^*N\otimes p_S^*M)=(\int_f(N,F))\otimes^L_{O_S}(M,F). 
\end{eqnarray*}
\end{itemize}
Clearly if $i:Z\hookrightarrow S$ is a closed embedding with $Z,S\in\SmVar(\mathbb C)$ or with $Z,S\in\AnSm(\mathbb C)$,
then $T^{\mathcal D,0}(\otimes,i)(M,N)=T^{\mathcal D}(\otimes,i)(M,N)$ in $D_{\mathcal D(2)fil,\infty}(S)$.  
\end{defi}

We have then the following :

\begin{prop}\label{PDmod2}
\begin{itemize}
\item[(i)] Let $i:Z\hookrightarrow S$ is a closed embedding with $Z,S\in\SmVar(\mathbb C)$,
then for $(M,F)\in C_{\mathcal Dfil}(S)$ and $(N,F)\in C_{\mathcal Dfil}(Z)$
\begin{equation*}
T^{\mathcal D,0}(\otimes,i)((M,F),(N,F)):
(M,F)\otimes_{O_S}i_{*mod}(N,F)\xrightarrow{\sim} i_{*mod}(i^{*mod}(M,F)\otimes_{O_Z}(N,F))  
\end{equation*}
is an isomorphism in $C_{\mathcal Dfil}(S)$.
\item[(ii)] Let $f:X\to S$ a morphism with $X,S\in\SmVar(\mathbb C)$.
Then, for $(M,F)\in C_{\mathcal D(2)fil}(X)$ and $(N,F)\in C_{\mathcal D(2)fil}(S)$,
\begin{equation*}
T^{\mathcal Dmod}(\otimes,f)((M,F),(N,F)):
\int^{FDR}_{f}((M,F)\otimes^L_{O_X}f_{FDR}^{*mod}(N,F))\xrightarrow{\sim}(\int^{FDR}_{f}(M,F))\otimes^L_{O_Y}(N,F) 
\end{equation*}
is an isomorphism in $D_{\mathcal D(2)fil,\infty}(S)$.
\end{itemize}
\end{prop}

\begin{proof}
\noindent(i): Follows from proposition \ref{Tiotimes}.

\noindent(ii):Follows from proposition \ref{PDmod1}(i).
\end{proof}

Let $f:X\to S$ a morphism with $X,S\in\SmVar(\mathbb C)$.
Consider the graph embedding $f:X\xrightarrow{i}X\times S\xrightarrow{p}S$, with  $X,Y,S\in\SmVar(\mathbb C)$.  
We have, for $(M,F)\in C_{\mathcal Dfil}(X)$, the canonical isomorphism in $C_{\mathcal D(2)fil}(S^{an})$
\begin{eqnarray*}
\an_X^{*mod}i^{*mod}L_D(p^{*mod}(M,F)\otimes_{O_{X\times S}}(O_{X\times S},\mathcal V_X)) 
\xrightarrow{=} \\
i^{*mod}L_Dp^{*mod}((M,F)^{an}\otimes_{O_{X^{an}\times S^{an}}}(O_{X^{an}\times S^{an}},\mathcal V_{X^{an}}))
\end{eqnarray*}

We then define and study the transformation map between the direct image functor and the analytical functor for D-modules :

\begin{defi}
Let $f:X\to S$ a morphism with $X,S\in\SmVar(\mathbb C)$.
\begin{itemize}
\item[(i)] We have for $(M,F)\in C_{\mathcal D(2)fil}(X)$ the canonical map in $C_{\mathcal D(2)fil}(S^{an})$
\begin{eqnarray*}
T^{\mathcal Dmod0}(an,f)(M,F):\an_S^{*mod}(f_*E((D_{X\leftarrow S},F^{ord})\otimes_{D_X} L_D(M,F)))
\xrightarrow{T^{mod}(an,f)(-)} \\
f_*(E((D_{X\leftarrow S},F^{ord})\otimes_{D_X} L_D(M,F)))^{an}\xrightarrow{=}
f_*E(D_{X^{an}\leftarrow S^{an}}\otimes_{D_{X^{an}}} L_D(M^{an},F))
\end{eqnarray*}
\item[(ii)] Consider the graph embedding $f:X\xrightarrow{i}X\times S\xrightarrow{p}S$, with  $X,Y,S\in\SmVar(\mathbb C)$.  
We have, for $(M,F)\in C_{\mathcal Dfil}(X)$, the canonical map in $C_{\mathcal D(2)fil}(S^{an})$
\begin{eqnarray*}
T^{\mathcal Dmod}(an,f)(M,F):
\an_S^{*mod}(p_*E((\Omega^{\bullet}_{Y\times S/S},F_b)\otimes_{O_{Y\times S}}i_{*mod}(M,F))) \\
\xrightarrow{T^O_{\omega}(an,p)(i_{*mod}(M,F))} 
p_*E((\Omega^{\bullet}_{Y\times S/S},F_b)\otimes_{O_{Y^{an}\times S^{an}}}(i_{*mod}(M,F))^{an}) \\
\xrightarrow{p_{*mod}T^{\mathcal Dmod0}(an,i)((M,F))}
p_*E((\Omega^{\bullet}_{Y\times S/S},F_b)\otimes_{O_{Y^{an}\times S^{an}}}i_{*mod}((M,F)^{an})).
\end{eqnarray*}
\end{itemize}
\end{defi}

In order to prove that this map gives an isomorphism in the derived category in the non filtered case 
if $f$ is proper and $M$ coherent, we will need the following (c.f.\cite{LvDmod}):
\begin{thm}\label{Daffinethm}
A product $X\times S$ of a smooth projective variety $X$ and a smooth affine variety $S$ is D-affine.
\end{thm}

\begin{proof}
See \cite{LvDmod} theorem 1.6.5.
\end{proof}

A main result is that we have the following version of the first GAGA theorem for coherent D-modules :

\begin{thm}\label{GAGADmod}
Let $f:X\to S$ a morphism with $X,S\in\SmVar(\mathbb C)$. Let $M\in D_{\mathcal D(2)fil,c}(X)$, for $r=1,\ldots\infty$.
If $f$ is proper, 
\begin{equation*}
T^{\mathcal Dmod}(an,f)(M,F):(\int_fM)^{an}\xrightarrow{\sim}\int_{f^{an}}(M^{an})
\end{equation*}
is an isomorphism.
\end{thm}

\begin{proof}
We may assume that $f$ is projective, so that we have a factorization 
$f:X\xrightarrow{i}\mathbb P^N\times S\xrightarrow{p} S$ where $i$ is a closed embedding and $p$ the projection.
The question being local on $S$, we may assume that $S$ is affine.
Since $\mathbb P^N\times S$ is D-affine by theorem \ref{Daffinethm}, we have by proposition \ref{Daffineprop}(iii) a 
complex $F\in C_{\mathcal D}(\mathbb P^N\times S)$ such that $i_{*mod}M=F\simeq F\in D_{\mathcal D,r}(\mathbb P^N\times S)$
and each $F^n$ is a direct summand of a free $D_{\mathbb P^N\times S}$ module of finite rank.
The theorem now follows from the fact that $\int_pD_{\mathbb P^N\times S}\simeq D_S[-N]$ and the fact that $(D_S)^{an}=D_{S^{an}}$.
\end{proof}

We also have

\begin{defi}\label{TDH0}
\begin{itemize}
\item[(i)] Let $f:X\to S$ a morphism with $X,S\in\SmVar(\mathbb C)$ or with $X,S\in\AnSm(\mathbb C)$.
We have, for $M,N\in C_{\mathcal Dfil}(X)$, the canonical transformation map in $D_{\mathcal Dfil,\infty}(S)$
\begin{eqnarray*}
T^{0,\mathcal D}(f,hom)((M,F),(N,F)):Rf_*R\mathcal Hom_{f^*D_S}((M,F),(N,F))\to \\ 
Rf_*R\mathcal Hom_{D_X}((M,F)\otimes_{D_X}L_D(D_{X\leftarrow S},F^{ord}),(N,F)\otimes_{D_X}L_D(D_{X\leftarrow S},F^{ord})) 
\xrightarrow{T^0(f,hom)(E(-),E(-))} \\
R\mathcal Hom_{D_X}(Rf_*((M,F)\otimes_{D_X}L_D(D_{X\leftarrow S},F^{ord})),Rf_*((N,F)\otimes_{D_X}L_D(D_{X\leftarrow S},F^{ord})))= \\ 
R\mathcal Hom_{D_X}(\int_f(M,F),\int_f(N,F))
\end{eqnarray*}
\item[(ii)] Let $f:X\to S$ a morphism with $X,S\in\AnSm(\mathbb C)$.
We have, for $(M,F),(N,F)\in C_{\mathcal Dfil}(X)$, the canonical transformation map in $D_{\mathcal Dfil,\infty}(S)$
\begin{eqnarray*}
T^{0,\mathcal D}(f_!,hom)((M,F),(N,F)):Rf_*\mathcal Hom_{f^*D_S}((M,F),(N,F))\to \\
Rf_*\mathcal Hom_{D_X}((M,F)\otimes_{D_X}L_D(D_{X\leftarrow S},F^{ord}),(N,F)\otimes_{D_X}L_D(D_{X\leftarrow S},F^{ord})) 
\xrightarrow{T^0(f_!,hom)(E(-),E(-))} \\
R\mathcal Hom_{D_X}(Rf_!((M,F)\otimes_{D_X}L_D(D_{X\leftarrow S},F^{ord})),Rf_!((N,F)\otimes_{D_X}L_D(D_{X\leftarrow S},F^{ord})))= \\
R\mathcal Hom_{D_X}(\int_{f!}(M,F),\int_{f!}(N,F))
\end{eqnarray*}
\end{itemize}
\end{defi}

\begin{defi}\label{TDH}
Let $f:X\to S$ a morphism with $X,S\in\SmVar(\mathbb C)$ or with $X,S\in\AnSm(\mathbb C)$.
We have, for $(M,F),(N,F)\in C_{\mathcal Dfil}(S)$, the canonical transformation map in $C_{\mathcal Dfil}(X)$
\begin{eqnarray*}
T^{\mathcal D}(f,hom)((M,F),(N,F)):f^*\mathcal Hom_{D_S}((M,F),(N,F)) \\
\xrightarrow{T(f,hom)((M,F),(N,F))}\mathcal Hom_{f^*D_S}(f^*(M,F),(f^*(N,F)) \\ 
\to\mathcal Hom_{D_X}(f^*(M,F)\otimes_{f^*D_S}L_{f^*D}(D_{X\to S},F^{ord}),f^*(N,F)\otimes_{f^*D_S}L_{f^*D}(D_{X\to S},F^{ord})) \\
\xrightarrow{=}\mathcal Hom_{D_X}(f^{*mod}(M,F),f^{*mod}(N,F))
\end{eqnarray*}
which is the one given by Kashiwara (see \cite{KS}).
\end{defi}

In the algebraic case, we have, in the non filtered case, the six functor formalism for holonomic D-modules :

\begin{thm}\label{Dmodad}
Let $f:X\to S$ a morphism with $X,S\in\SmVar(\mathbb C)$.
\begin{itemize}
\item[(i)] We have, for $M\in D_{\mathcal D,h}(X)$ and $N\in D_{\mathcal D,h}(S)$ 
a canonical isomorphism in $D_{\mathcal D}(S)$
\begin{eqnarray*}
I^{\mathcal Dmod}(Lf^{\hat{*}mod[-]},\int_f)(M,N):  
Rf_*R\mathcal Hom_{D_X}(Lf^{\hat{*}mod[-]}N,M)\xrightarrow{\sim} R\mathcal Hom_{D_S}(N,\int_fM).
\end{eqnarray*}
\item[(ii)] We have, for $M\in D_{\mathcal D,h}(X)$ and $N\in D_{\mathcal D,h}(S)$ 
a canonical isomorphism in $D_{\mathcal D}(X)$
\begin{eqnarray*}
I^{\mathcal Dmod}(\int_{f!},Lf^{*mod[-]})(M,N):
R\mathcal Hom_{D_X}(\int_{f!}M,N)\xrightarrow{\sim} Rf_*R\mathcal Hom_{D_S}(M,Lf^{*mod[-]}N).
\end{eqnarray*}
\end{itemize}
\end{thm}

\begin{proof}
Follows from the projection case and the closed embedding case.
\end{proof}

\begin{cor}
Let $f:X\to S$ a morphism with $X,S\in\SmVar(\mathbb C)$. Then,
\begin{itemize}
\item $(Lf^{\hat{*}mod[-]},\int_f):D_{\mathcal D,h}(S)\to D_{\mathcal D,h}(X)$ is a pair of adjoint functors.
\item $(\int_{f!},Lf^{*mod[-]}):D_{\mathcal D,h}(S)\to D_{\mathcal D,h}(X)$ is a pair of adjoint functors.
\end{itemize}
\end{cor}

\begin{proof}
Follows immediately from theorem \ref{Dmodad} by taking global sections.
\end{proof}

Consider a commutative diagram in $\SmVar(\mathbb C)$, 
\begin{equation*}
D=\xymatrix{
X'\ar[r]^{f'}\ar[d]_{g'} & T\ar[d]^{g} \\
X\ar[r]^{f} & S}.
\end{equation*}
We have, for $M\in C_{\mathcal D,h}(X)$, the following transformation maps
\begin{eqnarray*}
T_1^{\mathcal Dmod}(D)(M):Lg^{\hat*mod[-]}\int_fM\xrightarrow{\ad(Lf^{'\hat*mod[-]},\int_{f'})(-)} 
\int_{f'}Lf^{'\hat*mod[-]}Lg^{\hat*mod[-]}\int_fM\xrightarrow{=} \\ 
\int_{f'}Lg_{FDR}^{'\hat*mod[-]}Lf^{\hat*mod[-]}\int_fM\xrightarrow{\ad(Lf^{\hat*mod[-]},\int_{f})(M)}
\int_{f'}Lg_{FDR}^{'\hat*mod[-]}M
\end{eqnarray*}
and
\begin{eqnarray*}
T_2^{\mathcal Dmod}(D)(M,F):\int_{f'!}Lg^{'*mod[-]}M\xrightarrow{\ad(\int_{f!},Lf^{*mod[-]})(-)}
\int_{f'!}Lg^{'*mod[-]}Lf^{*mod[-]}\int_{f!}M\xrightarrow{=} \\
\int_{f'!}Lf^{'*mod[-]}Lg^{*mod[-]}M\int_{f!}\xrightarrow{\ad(\int_{f'!},Lf^{'*mod[-]})(-)}
Lg^{*mod[-]}\int^{FDR}_{f!}M
\end{eqnarray*}

\begin{prop}
Consider a cartesian square in $\SmVar(\mathbb C)$ 
\begin{equation*}
D=\xymatrix{ 
X_T\ar[r]^{g'}\ar[d]_{f'} & X\ar[d]^{f} \\
T\ar[r]^{g} & S}
\end{equation*}
Assume that $f$ (and hence $f'$) is proper. Then, for $(M,F)\in D_{\mathcal D(2)fil,\infty,h}(X)$,
\begin{itemize}
\item $T_1^{\mathcal Dmod}(f,g)(M):Lg^{\hat*mod[-]}\int_fM\xrightarrow{\sim}\int_{f'}Lg^{'\hat*mod[-]}M$ and
\item $T_2^{\mathcal Dmod}(f,g)(M):\int_{f'!}Lg^{'*mod[-]}M\xrightarrow{\sim}Lg^{*mod[-]}\int_{f!}M$
\end{itemize}
are isomorphisms in $D_{\mathcal D}(T)$.
\end{prop}

\begin{proof}
Follows from proposition \ref{PDmod1} and the fact that
the map $T_1^{\mathcal Dmod}(f,g)(M)$ is given by the composite 
\begin{eqnarray*}
T_1^{\mathcal Dmod}(f,g)(M)[d_T-d_S]:Lg^{\hat*mod}\int_f(M,F)=L\mathbb D_TLg^{*mod}L\mathbb D_S\int_fM
\xrightarrow{T(f_*,f_!)(-)} \\ L\mathbb D_TLg^{*mod}\int_fL\mathbb D_XM
\xrightarrow{(L\mathbb D_TT^{\mathcal Dmod}(f,g)(\mathbb D_XM))^{-1}}
\mathbb D_T\int_{f'}Lg^{'*mod}\mathbb D_XM \\ \xrightarrow{T(f'_!,f'_*)(-)}
\int_{f'}L\mathbb D_{X_T}Lg^{'\hat*mod}\mathbb D_X(M,F)=\int_{f'}Lg^{'\hat*mod}M
\end{eqnarray*}
and the map $T_2^{\mathcal Dmod}(f,g)(M,F)$ is given by the composite 
\begin{eqnarray*}
T_2^{\mathcal Dmod}(f,g)(M)[d_T-d_S]:\int_{f'!}Lg^{'*mod}M=L\mathbb D_T\int_{f'}\mathbb D_{X_T}Lg^{'*mod}M
\xrightarrow{d(-)\circ T(f_*,f_!)(-)} \\ \int_{f'}Lg^{'*mod}M 
\xrightarrow{T^{\mathcal Dmod}(f,g)(M)^{-1}}
L\mathbb D_TLg^{*mod}\int_fL\mathbb D_XM=Lg^{*mod}(M,F)\int_{f!}M
\end{eqnarray*}
\end{proof}

\subsection{The D modules on singular algebraic varieties and singular complex analytic spaces}

In this subsection by defining the category of complexes of filtered D-modules in the singular case and
there functorialities. 

\subsubsection{Definition}

In all this subsection, we fix the notations:
\begin{itemize}
\item For $S\in\Var(\mathbb C)$, we denote by $S=\cup_iS_i$ an open cover
such that there exits closed embeddings $i_iS_i\hookrightarrow\tilde S_i$ with $\tilde S_i\in\SmVar(\mathbb C)$.
We have then closed embeddings $i_I:S_I:=\cap_{i\in I} S_i\hookrightarrow\tilde S_I:=\Pi_{i\in I}\tilde S_I$.
Then for $I\subset J$, we denote by $j_{IJ}:S_J\hookrightarrow S_I$ the open embedding 
and $p_{IJ}:\tilde S_J\to\tilde S_I$ the projection, so that $p_{IJ}\circ i_J=i_I\circ j_{IJ}$.
This gives the diagram of algebraic varieties $(\tilde S_I)\in\Fun(\mathcal P(\mathbb N),\Var(\mathbb C))$ which
gives the diagram of sites $(\tilde S_I):=\Ouv(\tilde S_I)\in\Fun(\mathcal P(\mathbb N),\Cat)$.
It also gives the diagram of sites $(\tilde S_I)^{op}:=\Ouv(\tilde S_I)^{op}\in\Fun(\mathcal P(\mathbb N),\Cat)$.
For $I\subset J$, we denote by $m:\tilde S_I\backslash(S_I\backslash S_J)\hookrightarrow\tilde S_I$ the open embedding.
\item For $S\in\AnSp(\mathbb C)$ we denote by $S=\cup_iS_i$ an open cover
such that there exist closed embeddings $i_i:S_i\hookrightarrow\tilde S_i$ with $\tilde S_i\in\AnSm(\mathbb C)$.
We have then closed embeddings $i_I:S_I=\cap_{i\in I} S_i\hookrightarrow\tilde S_I=\Pi_{i\in I}\tilde S_I$.
Then for $I\subset J$, we denote by $j_{IJ}:S_J\hookrightarrow S_I$ the open embedding 
and $p_{IJ}:\tilde S_J\to\tilde S_I$ the projection, so that $p_{IJ}\circ i_J=i_I\circ j_{IJ}$.
This gives the diagram of analytic spaces $(\tilde S_I)\in\Fun(\mathcal P(\mathbb N),\AnSp(\mathbb C))$ which
gives the diagram of sites $(\tilde S_I):=\Ouv(\tilde S_I)\in\Fun(\mathcal P(\mathbb N),\Cat)$.
It also gives the diagram of sites $(\tilde S_I)^{op}:=\Ouv(\tilde S_I)^{op}\in\Fun(\mathcal P(\mathbb N),\Cat)$.
For $I\subset J$, we denote by $m:\tilde S_I\backslash(S_I\backslash S_J)\hookrightarrow\tilde S_I$ the open embedding.
\end{itemize}

The first definition is from \cite{Saito} remark 2.1.20, where we give a shifted version to have compatibility
with perverse sheaves. 

\begin{defi}
Let $S\in\Var(\mathbb C)$ and let $S=\cup_iS_i$ an open cover
such that there exist closed embeddings $i_iS_i\hookrightarrow\tilde S_i$ with $\tilde S_i\in\SmVar(\mathbb C)$ ;
or let $S\in\AnSp(\mathbb C)$ and let $S=\cup_iS_i$ an open cover such that there exit closed embeddings
$i_i:S_i\hookrightarrow\tilde S_i$ with $\tilde S_i\in\AnSm(\mathbb C)$. 
Then, $\PSh_{\mathcal D(2)fil}(S/(\tilde S_I))\subset\PSh_{\mathcal D(2)fil}((\tilde S_I))$ is the full subcategory
\begin{itemize}
\item whose objects are $(M,F)=((M_I,F)_{I\subset\left[1,\cdots l\right]},s_{IJ})$, with
\begin{itemize}
\item $(M_I,F)\in\PSh_{\mathcal D(2)fil}(\tilde S_I)$ such that $\mathcal I_{S_I}M_I=0$, 
in particular $(M_I,F)\in\PSh_{\mathcal D(2)fil,S_I}(\tilde S_I)$
\item $s_{IJ}:m^*(M_I,F)\xrightarrow{\sim} m^*p_{IJ*}(M_J,F)[d_{\tilde S_I}-d_{\tilde S_J}]$ 
for $I\subset J$, are isomorphisms, $p_{IJ}:\tilde S_J\to\tilde S_I$ being the projection,
satisfying for $I\subset J\subset K$, $p_{IJ*}s_{JK}\circ s_{IJ}=s_{IK}$ ; 
\end{itemize}
\item the morphisms $m:(M,F)\to(N,F)$ between  
$(M,F)=((M_I,F)_{I\subset\left[1,\cdots l\right]},s_{IJ})$ and $(N,F)=((N_I,F)_{I\subset\left[1,\cdots l\right]},r_{IJ})$
are by definition a family of morphisms of complexes,  
\begin{equation*}
m=(m_I:(M_I,F)\to (N_I,F))_{I\subset\left[1,\cdots l\right]}
\end{equation*}
such that $r_{IJ}\circ m_J=p_{IJ*}m_J\circ s_{IJ}$ in $C_{\mathcal D,S_J}(\tilde S_J)$.
\end{itemize}
We denote by 
\begin{equation*}
\PSh_{\mathcal D(2)fil,rh}(S/(\tilde S_I))\subset
\PSh_{\mathcal D(2)fil,h}(S/(\tilde S_I))\subset\PSh_{\mathcal D(2)fil,c}(S/(\tilde S_I))
\subset\PSh_{\mathcal D(2)fil}(S/(\tilde S_I))
\end{equation*}
the full subcategory consisting of $((M_I,F),s_{IJ})$ 
such that $(M_I,F)$ is filtered coherent, resp. filtered holonomic, resp. filtered regular holonomic, i.e. 
$M_I$ are coherent, resp. holonomic,resp. filtered regular holonomic, 
sheaves of $D_{\tilde S_I}$ modules and $F$ is a good filtration.
We have the full subcategories
\begin{eqnarray*} 
\PSh_{\mathcal D(1,0)fil,rh}(S/(\tilde S_I))\subset\PSh_{\mathcal D2fil,rh}(S/(\tilde S_I)), \;
\PSh_{\mathcal D(1,0)fil,h}(S/(\tilde S_I))\subset\PSh_{\mathcal D2fil,h}(S/(\tilde S_I)), \\
\PSh_{\mathcal D(1,0)fil,h}(S/(\tilde S_I))\subset\PSh_{\mathcal D2fil,h}(S/(\tilde S_I)), 
\end{eqnarray*}
consisting of $((M_I,F,W),s_{IJ})$ such that $W^pM_I$ are $D_{\tilde S_I}$ submodules.
\end{defi}

We recall from section 2 the following
\begin{itemize}
\item A morphism $m=(m_I):((M_I),s_{IJ})\to((N_I),r_{IJ})$ in $C(\PSh_{\mathcal D}(S/(\tilde S_I)))$ 
is a Zariski, resp. usu, local equivalence if and only if all the $m_I$ are Zariski, resp. usu, local equivalences.
\item A morphism $m=(m_I):((M_I,F),s_{IJ}\to((N_I,F),r_{IJ}))$ in $C(\PSh_{\mathcal D(2)fil}(S/(\tilde S_I)))$
is a filtered Zariski, resp. usu, local equivalence 
if and only if all the $m_I$ are filtered Zariski, resp. usu, local equivalence.
\item By definition, a morphism $m=(m_I):((M_I,F),s_{IJ})\to((N_I,F),r_{IJ})$ in $C(\PSh_{\mathcal D(2)fil}(S/(\tilde S_I)))$
is an $r$-filtered Zariski, resp. usu, local equivalence 
if there exist $m_i:(C_{iI},F),s_{iIJ})\to(C_{(i+1)I},F),s_{(i+1)IJ})$, $0\leq i\leq s$, 
with $(C_{iI},F),s_{iIJ})\in C(\PSh_{\mathcal D(2)fil}(S/(\tilde S_I)))$, 
$(C_{0I},F),s_{iIJ})=(M_I,F),s_{IJ})$, $(C_{sI},F),s_{sIJ})=(N_I,F),r_{IJ})$ such that
\begin{equation*}
m=m_s\circ\cdots\circ m_i\circ\cdots\circ m_0:((M_I,F),s_{IJ}\to((N_I,F),r_{IJ}))
\end{equation*}
with $m_i:(C_{iI},F),s_{iIJ})\to(C_{(i+1)I},F),s_{(i+1)IJ})$ either filtered Zariski, resp. usu, local equivalence
or $r$-filtered homotopy equivalence.
\end{itemize}

Let $S\in\Var(\mathbb C)$ or $S\in\AnSp(\mathbb C)$. 
\begin{itemize}
\item If $S\in\Var(\mathbb C)$, let $S=\cup_{i=1}^lS_i$ an open cover such that there exist closed embeddings 
$i_i:S_i\hookrightarrow\tilde S_i$ with $\tilde S_i\in\SmVar(\mathbb C)$,
and let $S=\cup_{i'=1}^{l'}S_{i'}$ an other open cover such that there exist closed embeddings 
$i_{i'}:S_{i'}\hookrightarrow\tilde S_{i'}$ with $\tilde S_{i'}\in\SmVar(\mathbb C)$.
\item If $S\in\AnSp(\mathbb C)$, let $S=\cup_{i=1}^lS_i$ an open cover such that there exist closed embeddings 
$i_i:S_i\hookrightarrow\tilde S_i$ with $\tilde S_i\in\AnSm(\mathbb C)$,
and let $S=\cup_{i'=1}^{l'}S_{i'}$ an other open cover such that there exist closed embeddings 
$i_{i'}:S_{i'}\hookrightarrow\tilde S_{i'}$ with $\tilde S_{i'}\in\AnSm(\mathbb C)$.
\end{itemize}
Denote $L=[1,\ldots,l]$, $L'=[1,\ldots,l']$ and $L'':=[1,\ldots,l]\sqcup[1,\ldots,l']$.
We have then the refined open cover $S=\cup_{k\in L} S_k$ and we denote for $I\sqcup I'\subset L''$, 
$S_{I\sqcup I'}:=\cap_{k\in I\sqcup I'} S_k$ and $\tilde S_{I\sqcup I'}:=\Pi_{k\in I\sqcup I'} \tilde S_k$,
so that we have a closed embedding $i_{I\sqcup I'}:S_{I\sqcup I'}\hookrightarrow\tilde S_{I\sqcup I'}$.
For $I\sqcup I'\subset J\sqcup J'$, denote by $p_{I\sqcup I',J\sqcup J'}:\tilde S_{J\sqcup J'}\to\tilde S_{I\sqcup I'}$ the projection.
We then have a natural transfer map
\begin{eqnarray*}
T^{L/L'}_{S}:\PSh_{\mathcal Dfil}(S/(S_I))\to\PSh_{\mathcal Dfil}(S/(S_{I'})), \\ 
((M_I,F),s_{IJ})\mapsto (\ho\lim_{I\in L}p_{I'(I\sqcup I')*}
(p_{I(I\sqcup I')}^{*mod}(M_I,F))/\mathcal I_{S_{I\sqcup I'}},s_{I'J'}),
\end{eqnarray*}
with, in the homotopy limit, the natural transition morphisms 
\begin{eqnarray*}
p_{I'(I\sqcup I')*}\ad(p_{IJ}^{*mod},p_{IJ*})(p_{I(I\sqcup I')}^{*mod[-]}(M_I,F)): \\
p_{I'(J\sqcup I')*}(p_{J(J\sqcup I')}^{*mod[-]}(M_J,F))/\mathcal I_{S_{J\sqcup I'}}\to
p_{I'(I\sqcup I')*}(p_{I(I\sqcup I')}^{*mod[-]}(M_I,F))/\mathcal I_{S_{I\sqcup I'}}
\end{eqnarray*}
for $J\subset I$, and 
\begin{eqnarray*}
s_{I'J'}:\holim_{I\in L}m^*p_{I'(I\sqcup I')*}(p_{I(I\sqcup I')}^{*mod[-]}(M_I,F)/\mathcal I_{S_{I\sqcup I'}})\to \\
\holim_{I\in L}p_{I'J'*}(p_{I'J'}^{*mod[-]}m^*p_{I'(I\sqcup I')*}
p_{I(I\sqcup I')}^{*mod[-]}((M_I,F)/\mathcal I_{S_{I\sqcup I'}}))/\mathcal I_{S_{J'}} \\
\to\holim_{I\in L} p_{I'J'*}p_{J'(I\sqcup J')*}(p_{I(I\sqcup J')}^{*mod[-]}(M_I,F)\mathcal I_{S_{I\sqcup I'}})
\end{eqnarray*}

\begin{defiprop}\label{defprops} 
Let $S\in\Var(\mathbb C)$ and let $S=\cup_iS_i$ an open cover
such that there exist closed embeddings $i_iS_i\hookrightarrow\tilde S_i$ with $\tilde S_i\in\SmVar(\mathbb C)$ ;
or let $S\in\AnSp(\mathbb C)$ and let $S=\cup_iS_i$ an open cover such that there exist closed embeddings
$i_iS_i\hookrightarrow\tilde S_i$ with $\tilde S_i\in\AnSm(\mathbb C)$.
Then $\PSh_{\mathcal D(2)fil}(S/(\tilde S_I))$ does not depend on the open covering of $S$ and the closed embeddings and we set
\begin{equation*}
\PSh_{\mathcal D(2)fil}(S):=\PSh_{\mathcal D(2)fil}(S/(\tilde S_I))
\end{equation*}
We denote by $C^0_{\mathcal D(2)fil}(S):=C(\PSh_{\mathcal D(2)fil}(S/(\tilde S_I)))$ and by 
$D^0_{\mathcal D(2)fil,r}(S):=K^0_{\mathcal D(2)fil,r}(S)([E_1]^{-1})$
the localization of the $r$-filtered homotopy category with respect to 
the classes of filtered Zariski, resp. usu, local equivalences.
\end{defiprop}

\begin{proof}
It is obvious that 
$T_S^{L/L'}:\PSh_{\mathcal Dfil}(S/(S_I))\to\PSh_{\mathcal D}(S/(S_{I'}))$ 
is an equivalence of category with inverse
$T_S^{L'/L}:\PSh_{\mathcal Dfil}(S/(S_{I'}))\to\PSh_{\mathcal D}(S/(S_I))$.
\end{proof}

We now give the definition of our category :

\begin{defi}\label{Dmodsingdef}
Let $S\in\Var(\mathbb C)$ and let $S=\cup_iS_i$ an open cover
such that there exist closed embeddings $i_i:S_i\hookrightarrow\tilde S_i$ with $\tilde S_i\in\SmVar(\mathbb C)$ ;
or let $S\in\AnSp(\mathbb C)$ and let $S=\cup_iS_i$ an open cover 
such that there exist closed embeddings $i_i:S_i\hookrightarrow\tilde S_i$ with $\tilde S_i\in\AnSm(\mathbb C)$. 
Then, $C_{\mathcal D(2)fil}(S/(\tilde S_I))\subset C_{\mathcal D(2)fil}((\tilde S_I))$ is the full subcategory
\begin{itemize}
\item whose objects are $(M,F)=((M_I,F)_{I\subset\left[1,\cdots l\right]},u_{IJ})$, with
\begin{itemize}
\item $(M_I,F)\in C_{\mathcal D(2)fil,S_I}(\tilde S_I)$ (see definition \ref{DmodsuppZ}),
\item $u_{IJ}:m^*(M_I,F)\to m^*p_{IJ*}(M_J,F)[d_{\tilde S_I}-d_{\tilde S_J}]$ 
for $J\subset I$, are morphisms, $p_{IJ}:\tilde S_J\to\tilde S_I$ being the projection, 
satisfying for $I\subset J\subset K$, $p_{IJ}*u_{JK}\circ u_{IJ}=u_{IK}$ in $C_{\mathcal Dfil}(\tilde S_I)$ ;
\end{itemize}
\item the morphisms $m:((M_I,F),u_{IJ})\to((N_I,F),v_{IJ})$ between  
$(M,F)=((M_I,F)_{I\subset\left[1,\cdots l\right]},u_{IJ})$ and $(N,F)=((N_I,F)_{I\subset\left[1,\cdots l\right]},v_{IJ})$
being a family of morphisms of complexes,  
\begin{equation*}
m=(m_I:(M_I,F)\to (N_I,F))_{I\subset\left[1,\cdots l\right]}
\end{equation*}
such that $v_{IJ}\circ m_I=p_{IJ*}m_J\circ u_{IJ}$ in $C_{\mathcal Dfil}(\tilde S_I)$.
\end{itemize}
We denote by $C^{\sim}_{\mathcal D(2)fil}(S/(\tilde S_I))\subset C_{\mathcal D(2)fil}(S/(\tilde S_I))$ the full subcategory 
consisting of objects $((M_I,F),u_{IJ})$ such that the $u_{IJ}$ are $\infty$-filtered Zariski, resp. usu, local equivalences.
\end{defi}

Let $S\in\Var(\mathbb C)$ and let $S=\cup_iS_i$ an open cover
such that there exist closed embeddings $i_i:S_i\hookrightarrow\tilde S_i$ with $\tilde S_i\in\SmVar(\mathbb C)$ ;
or let $S\in\AnSp(\mathbb C)$ and let $S=\cup_iS_i$ an open cover 
such that there exist closed embeddings $i_i:S_i\hookrightarrow\tilde S_i$ with $\tilde S_i\in\AnSm(\mathbb C)$. Then,
We denote by
\begin{equation*}
C_{\mathcal D(2)fil,rh}(S/(\tilde S_I))\subset 
C_{\mathcal D(2)fil,h}(S/(\tilde S_I))\subset C_{\mathcal D(2)fil,c}(S/(\tilde S_I))
\subset C^{\sim}_{\mathcal D(2)fil}(S/(\tilde S_I))
\end{equation*}
the full subcategories consisting of those $((M_I,F),u_{IJ})\in C^{\sim}_{\mathcal D(2)fil}(S/(\tilde S_I))$ 
such that $(M_I,F)\in C_{\mathcal D(2)fil,S_I,c}(\tilde S_I)$, 
that is such that $a_{\tau}H^n(M_I,F)$ are filtered coherent for all $n\in\mathbb Z$ and all $I\subset\left[1,\cdots l\right]$
(i.e. $a_{\tau}H^n(M_I)$ are coherent sheaves of $D_{\tilde S_I}$ modules and $F$ induces a good filtration on $a_{\tau}H^n(M_I)$),
resp. such that $(M_I,F)\in C_{\mathcal D(2)fil,S_I,h}(\tilde S_I)$, 
that is such that $a_{\tau}H^n(M_I,F)$ are filtered holonomic for all $n\in\mathbb Z$ and all $I\subset\left[1,\cdots l\right]$
(i.e. $a_{\tau}H^n(M_I)$ are holonomic sheaves of $D_{\tilde S_I}$ modules and $F$ induces a good filtration on $a_{\tau}H^n(M_I)$),
resp. such that $(M_I,F)\in C_{\mathcal D(2)fil,S_I,rh}(\tilde S_I)$, 
that is such that $a_{\tau}H^n(M_I,F)$ are filtered regular holonomic for all $n\in\mathbb Z$ and all $I\subset\left[1,\cdots l\right]$
(i.e. $a_{\tau}H^n(M_I)$ are regular holonomic sheaves of $D_{\tilde S_I}$ modules 
and $F$ induces a good filtration on $a_{\tau}H^n(M_I)$).

We denote by
\begin{eqnarray*}
C_{\mathcal D(1,0)fil,h}(S/(\tilde S_I))\subset C_{\mathcal D2fil,h}(S/(\tilde S_I)), \;
C_{\mathcal D(1,0)fil,rh}(S/(\tilde S_I))\subset C_{\mathcal D2fil,rh}(S/(\tilde S_I)), \\
C^{\sim}_{\mathcal D(1,0)fil}(S/(\tilde S_I))\subset C^{\sim}_{\mathcal D2fil}(S/(\tilde S_I)), 
\end{eqnarray*}
the full subcategories consisting of those $((M_I,F,W),u_{IJ})\in C^{\sim}_{\mathcal D2fil}(S/(\tilde S_I))$
such that $W^pM_I$ are $D_{\tilde S_I}$ submodules (resp. and $a_{\tau}H^n(M_I,F)$ are filtered holonomic).

We recall from section 2 the following
\begin{itemize}
\item A morphism $m=(m_I):((M_I),u_{IJ})\to((N_I),v_{IJ})$ in $C_{\mathcal D}(S/(\tilde S_I))$ 
is a Zariski, resp. usu, local equivalence if and only if all the $m_I$ are Zariski, resp. usu, local equivalences.
\item A morphism $m=(m_I):((M_I,F),u_{IJ}\to((N_I,F),v_{IJ}))$ in $C_{\mathcal D(2)fil}(S/(\tilde S_I))$
is a filtered Zariski, resp. usu, local equivalence 
if and only if all the $m_I$ are filtered Zariski, resp. usu, local equivalence.
\item Let $r\in\mathbb N$.
By definition, a morphism $m=(m_I):((M_I,F),u_{IJ})\to((N_I,F),v_{IJ})$ in $C_{\mathcal D(2)fil}(S/(\tilde S_I))$
is an $r$-filtered Zariski, resp. usu, local equivalence 
if there exist $m_i:(C_{iI},F),u_{iIJ})\to(C_{(i+1)I},F),u_{(i+1)IJ})$, $0\leq i\leq s$, 
with $(C_{iI},F),u_{iIJ})\in C_{\mathcal D(2)fil}(S/(\tilde S_I))$, 
$(C_{0I},F),u_{iIJ})=(M_I,F),u_{IJ})$, $(C_{sI},F),u_{sIJ})=(N_I,F),v_{IJ})$ such that
\begin{equation*}
m=m_s\circ\cdots\circ m_i\circ\cdots\circ m_0:((M_I,F),u_{IJ}\to((N_I,F),v_{IJ}))
\end{equation*}
with $m_i:(C_{iI},F),u_{iIJ})\to(C_{(i+1)I},F),u_{(i+1)IJ})$ either filtered Zariski, resp. usu, local equivalence
or $r$-filtered homotopy equivalence (i.e. $r$-filtered for the first filtration and filtered for the second filtration).
\end{itemize}

In the analytic case, we also define in the same way :

\begin{defi}\label{DmodsingVarinftydef}
Let $S\in\AnSp(\mathbb C)$ and let $S=\cup_iS_i$ an open cover 
such that there exist closed embeddings $i_i:S_i\hookrightarrow\tilde S_i$ with $\tilde S_i\in\AnSm(\mathbb C)$. Then,
$C_{\mathcal D^{\infty}(2)fil}(S/(\tilde S_I))\subset C_{\mathcal D^{\infty}(2)fil}((\tilde S_I))$ is the full subcategory
\begin{itemize}
\item whose objects are $(M,F)=((M_I,F)_{I\subset\left[1,\cdots l\right]},u_{IJ})$, with
\begin{itemize}
\item $(M_I,F)\in C_{\mathcal D^{\infty}fil,S_I}(\tilde S_I)$ (see definition \ref{DmodsuppZinfty}),
\item $u_{IJ}:m^*(M_I,F)\to m^*p_{IJ*}(M_I,F)[d_{\tilde S_I}-d_{\tilde S_J}]$, 
for $J\subset I$, are morphisms, $p_{IJ}:\tilde S_J\to\tilde S_I$ being the projection, 
satisfying for $I\subset J\subset K$, $p_{IJ*}u_{JK}\circ u_{IJ}=u_{IK}$ in $C_{\mathcal D^{\infty}fil}(\tilde S_I)$ ;
\end{itemize}
\item the morphisms $m:((M_I,F),u_{IJ})\to((N_I,F),v_{IJ})$ between  
$(M,F)=((M_I,F)_{I\subset\left[1,\cdots l\right]},u_{IJ})$ and $(N,F)=((N_I,F)_{I\subset\left[1,\cdots l\right]},v_{IJ})$
being a family of morphisms of complexes,  
\begin{equation*}
m=(m_I:(M_I,F)\to (N_I,F))_{I\subset\left[1,\cdots l\right]}
\end{equation*}
such that $v_{IJ}\circ m_I=p_{IJ*}m_J\circ u_{IJ}$ in $C_{\mathcal D^{\infty}fil}(\tilde S_I)$.
\end{itemize}
We denote by $C^{\sim}_{\mathcal D^{\infty}(2)fil}(S/(\tilde S_I))\subset C_{\mathcal D^{\infty}(2)fil}(S/(\tilde S_I))$ 
the full subcategory consisting of objects $((M_I,F),u_{IJ})$ such that the $u_{IJ}$ are $\infty$-filtered usu local equivalence.
\end{defi}

Let $S\in\AnSp(\mathbb C)$ and let $S=\cup_iS_i$ an open cover 
such that there exist closed embeddings $i_i:S_i\hookrightarrow\tilde S_i$ with $\tilde S_i\in\AnSm(\mathbb C)$. 
We denote by
\begin{equation*}
C_{\mathcal D^{\infty}(2)fil,h}(S/(\tilde S_I))\subset C_{\mathcal D^{\infty}(2)fil,c}(S/(\tilde S_I))
\subset C^{\sim}_{\mathcal D^{\infty}(2)fil}(S/(\tilde S_I))
\end{equation*}
the full subcategories consisting of 
$((M_I,F),u_{IJ})\in C^{\sim}_{\mathcal D^{\infty}(2)fil}(S/(\tilde S_I))$ 
such that $(M_I,F)\in C_{\mathcal D^{\infty}(2)fil,S_I,c}(\tilde S_I)$, 
that is such that $a_{\tau}H^n(M_I,F)$ are filtered coherent for all $n\in\mathbb Z$ and all $I\subset\left[1,\cdots l\right]$,
resp. such that $(M_I,F)\in C_{\mathcal D^{\infty}(2)fil,S_I,h}(\tilde S_I)$, 
that is such that $a_{\tau}H^n(M_I,F)$ are filtered holonomic for all $n\in\mathbb Z$ and all $I\subset\left[1,\cdots l\right]$.
We denote by
\begin{equation*}
C_{\mathcal D^{\infty}(1,0)fil}(S/(\tilde S_I))\subset C^{\sim}_{\mathcal D^{\infty}2fil}(S/(\tilde S_I)), \;
C_{\mathcal D^{\infty}(1,0)fil,h}(S/(\tilde S_I))\subset C^{\sim}_{\mathcal D^{\infty}2fil,h}(S/(\tilde S_I))
\end{equation*}
the full subcategories consisting of those $((M_I,F,W),u_{IJ})\in C^{\sim}_{\mathcal D2fil}(S/(\tilde S_I))$
such that $W^pM_I$ are $D_{\tilde S_I}$ submodules (resp. and $a_{\tau}H^n(M_I,F)$ filtered holonomic).

We recall from section 2 the following
\begin{itemize}
\item A morphism $m=(m_I):((M_I),u_{IJ})\to((N_I),v_{IJ})$ in $C_{\mathcal D^{\infty}}(S/(\tilde S_I))$ 
is a usu local equivalence if and only if all the $m_I$ are usu local equivalences.
\item A morphism $m=(m_I):((M_I,F),u_{IJ}\to((N_I,F),v_{IJ}))$ in $C_{\mathcal D^{\infty}(2)fil}(S/(\tilde S_I))$
is a filtered usu local equivalence if and only if all the $m_I$ are filtered usu local equivalence.
\item Let $r\in\mathbb N$.
By definition, a morphism $m=(m_I):((M_I,F),u_{IJ})\to((N_I,F),v_{IJ})$ in $C_{\mathcal D^{\infty}(2)fil}(S/(\tilde S_I))$
is an $r$-filtered usu local equivalence 
if there exist $m_i:(C_{iI},F),u_{iIJ})\to(C_{(i+1)I},F),u_{(i+1)IJ})$, $0\leq i\leq s$, 
with $(C_{iI},F),u_{iIJ})\in C_{\mathcal D(2)fil}(S/(\tilde S_I))$, 
$(C_{0I},F),u_{iIJ})=(M_I,F),u_{IJ})$, $(C_{sI},F),u_{sIJ})=(N_I,F),v_{IJ})$ such that
\begin{equation*}
m=m_s\circ\cdots\circ m_i\circ\cdots\circ m_0:((M_I,F),u_{IJ}\to((N_I,F),v_{IJ}))
\end{equation*}
with $m_i:(C_{iI},F),u_{iIJ})\to(C_{(i+1)I},F),u_{(i+1)IJ})$ either filtered usu local equivalence
or $r$-filtered homotopy equivalence (i.e. $r$-filtered for the first filtration and filtered for the second filtration).
\end{itemize}

\begin{defi}\label{defpropK} 
Let $S\in\Var(\mathbb C)$ and let $S=\cup_iS_i$ an open cover
such that there exist closed embeddings $i_i:S_i\hookrightarrow\tilde S_i$ with $\tilde S_i\in\SmVar(\mathbb C)$ ;
or let $S\in\AnSp(\mathbb C)$ and let $S=\cup_iS_i$ an open cover such that there exist closed embeddings
$i_i:S_i\hookrightarrow\tilde S_i$ with $\tilde S_i\in\AnSm(\mathbb C)$. 
\begin{itemize}
\item[(i)]We have, for $r=1,\ldots,\infty$, (see section 2.1) the $\infty$-filtered homotopy category
\begin{equation*}
K_{\mathcal D(2)fil,r}(S/(\tilde S_I)):=\Ho_{r}(C_{\mathcal D(2)fil}(S/(\tilde S_I)))
\end{equation*}
whose objects are those of $C_{\mathcal D(2)fil}(S/(\tilde S_I))$ and
whose morphism are $r$-filtered homotopy classes of morphisms 
($r$-filtered for the first filtration and filtered for the second), and
\begin{equation*}
D_{\mathcal D(2)fil,r}(S/(\tilde S_I)):=K_{\mathcal D(2)fil,r}(S/(\tilde S_I)))([E_1]^{-1})
\end{equation*}
the localization with respect to the classes of filtered Zariski, resp. usu, local equivalences, together with
the localization functor
\begin{equation*}
D(top):C_{\mathcal D(2)fil}(S/(\tilde S_I))\to K_{\mathcal D(2)fil,r}(S/(\tilde S_I))
\to D_{\mathcal D(2)fil,r}(S/(\tilde S_I)).
\end{equation*}
Note that the classes of filtered $\tau$ local equivalence constitute a right multiplicative system.
By definition if $m:((M_I,F),u_{IJ})\to((N_I,F),v_{IJ})$ in $C_{\mathcal D(2)fil}(S/(\tilde S_I))$ is
an $r$-filtered Zariski, resp. usu, local equivalence, then 
$m=D(\top)(m):((M_I,F),u_{IJ})\to((N_I,F),v_{IJ})$ is an isomorphism in $D_{\mathcal D(2)fil,r}(S/(\tilde S_I))$
\item[(ii)]We have 
\begin{equation*}
D_{\mathcal D(1,0)fil,\infty,h}(S/(\tilde S_I))\subset D_{\mathcal D2fil,\infty,h}(S/(\tilde S_I))
\subset D_{\mathcal D2fil,\infty}(S/(\tilde S_I)) 
\end{equation*}
the full subcategories which are the image of $C_{\mathcal D2fil,h}(S/(\tilde S_I))$, 
resp. of $C_{\mathcal D(1,0)fil,h}(S/(\tilde S_I))$, by the localization functor
$D(top):C_{\mathcal D(2)fil}(S/(\tilde S_I))\to D_{\mathcal D(2)fil,\infty}(S/(\tilde S_I))$.
\end{itemize}
\end{defi}

In the analytic case, we also have
\begin{defi}\label{defpropKinfty} 
Let $S\in\AnSp(\mathbb C)$ and let $S=\cup_iS_i$ an open cover such that there exist closed embeddings
$i_i:S_i\hookrightarrow\tilde S_i$ with $\tilde S_i\in\AnSm(\mathbb C)$.
\begin{itemize}
\item[(i)]We have, for $r=1,\ldots,\infty$, (see section 2.1) the $\infty$-filtered homotopy category
\begin{equation*}
K_{\mathcal D^{\infty}(2)fil,r}(S/(\tilde S_I)):=\Ho_{r}(C_{\mathcal D^{\infty}(2)fil}(S/(\tilde S_I)))
\end{equation*}
whose objects are those of $C_{\mathcal D^{\infty}(2)fil}(S/(\tilde S_I))$ and
whose morphism are $r$-filtered homotopy classes of morphisms 
($r$-filtered for the first filtration and filtered for the second), and
\begin{equation*}
D_{\mathcal D^{\infty}(2)fil,r}(S/(\tilde S_I)):=K_{\mathcal D^{\infty}(2)fil,r}(S/(\tilde S_I)))([E_1]^{-1})
\end{equation*}
the localization with respect to the classes of filtered usu local equivalences, together with
the localization functor
\begin{equation*}
D(usu):C_{\mathcal D^{\infty}(2)fil}(S/(\tilde S_I))\to K_{\mathcal D^{\infty}(2)fil,r}(S/(\tilde S_I))
\to D_{\mathcal D^{\infty}(2)fil,r}(S/(\tilde S_I)).
\end{equation*}
Note that the classes of filtered usu local equivalence constitute a right multiplicative system.
By definition if $m:((M_I,F),u_{IJ})\to((N_I,F),v_{IJ})$ in $C_{\mathcal D^{\infty}(2)fil}(S/(\tilde S_I))$ is
an $r$-filtered Zariski, resp. usu, local equivalence, then 
$m=D(\top)(m):((M_I,F),u_{IJ})\to((N_I,F),v_{IJ})$ is an isomorphism in $D_{\mathcal D^{\infty}(2)fil,r}(S/(\tilde S_I))$
\item[(ii)]We have then 
\begin{equation*}
D_{\mathcal D^{\infty}(1,0)fil,\infty,h}(S/(\tilde S_I))\subset D_{\mathcal D^{\infty}2fil,\infty,h}(S/(\tilde S_I))
\subset D_{\mathcal D^{\infty}2fil,\infty}(S/(\tilde S_I)) 
\end{equation*}
the full subcategories wich are the image of $C_{\mathcal D^{\infty}2fil,h}(S/(\tilde S_I))$, 
resp. $C_{\mathcal D^{\infty}(1,0)fil,h}(S/(\tilde S_I))$, by the localization functor 
$D(usu):C_{\mathcal D^{\infty}(2)fil}(S/(\tilde S_I))\to D_{\mathcal D^{\infty}(2)fil,\infty}(S/(\tilde S_I))$.
\end{itemize}
\end{defi}

\begin{defi}\label{iota0S}
Let $S\in\Var(\mathbb C)$ and let $S=\cup_iS_i$ an open cover
such that there exist closed embeddings $i_iS_i\hookrightarrow\tilde S_i$ with $\tilde S_i\in\SmVar(\mathbb C)$.
Or let $S\in\AnSp(\mathbb C)$ and let $S=\cup_iS_i$ an open cover such that there exist closed embeddings
$i_iS_i\hookrightarrow\tilde S_i$ with $\tilde S_i\in\AnSm(\mathbb C)$. 
\begin{itemize}
\item[(i)]We denote by
\begin{equation*}
C_{\mathcal D(2)fil}(S/(\tilde S_I))^0\subset C_{\mathcal D(2)fil}(S/(\tilde S_I))
\end{equation*}
the full subcategory consisting of $((M_I,F),u_{IJ})\in C_{\mathcal D(2)fil}(S/(\tilde S_I))$ such that 
\begin{equation*}
H^n((M_I,F),u_{IJ})=(H^n(M_I,F),H^nu_{IJ})\in\PSh^0_{\mathcal D(2)fil}(S/(\tilde S_I))
\end{equation*}
that is such that the $H^nu_{IJ}$ are isomorphism.
We denote by 
$D_{\mathcal D(2)fil}(S/(\tilde S_I))^0:=D(top)(C_{\mathcal D(2)fil}(S/(\tilde S_I))^0)$
its image by the localization functor.
\item[(ii)]We have the full embedding functor
\begin{eqnarray*}
\iota^0_{S/(\tilde S_I)}:C^0_{\mathcal D(2)fil}(S):=C^0_{\mathcal D(2)fil}(S/(\tilde S_I))
\hookrightarrow C_{\mathcal D(2)fil}(S/(\tilde S_I)), \\
((M_I,F),s_{IJ})\mapsto ((M_I,F),s_{IJ})
\end{eqnarray*}
By definition, 
$\iota^0_{S/(\tilde S_I)}(C^0_{\mathcal D(1,0)fil}(S/(\tilde S_I)))\subset C^{\sim}_{\mathcal D(1,0)fil}(S/(\tilde S_I))$.
This full embedding induces in the derived category the functor
\begin{eqnarray*}
\iota^0_{S/(\tilde S_I)}:D^0_{\mathcal D(2)fil,\infty}(S):=D^0_{\mathcal D(2)fil,\infty}(S/(\tilde S_I))\to 
D_{\mathcal D(2)fil,\infty}(S/(\tilde S_I)), \\
((M_I,F),s_{IJ})\mapsto ((M_I,F),s_{IJ}).
\end{eqnarray*}
\end{itemize}
\end{defi}

\begin{prop}\label{iota0Sprop}
Let $S\in\Var(\mathbb C)$ and let $S=\cup_iS_i$ an open cover
such that there exist closed embeddings $i_iS_i\hookrightarrow\tilde S_i$ with $\tilde S_i\in\SmVar(\mathbb C)$.
Or let $S\in\AnSp(\mathbb C)$ and let $S=\cup_iS_i$ an open cover such that there exist closed embeddings
$i_iS_i\hookrightarrow\tilde S_i$ with $\tilde S_i\in\AnSm(\mathbb C)$. Then, 
\begin{eqnarray*}
\iota^0_{S/(\tilde S_I)}:D^0_{\mathcal D(2)fil,\infty}(S)\to D_{\mathcal D(2)fil,\infty}(S/(\tilde S_I))
\end{eqnarray*}
is a full embedding whose image is $D_{\mathcal D(2)fil,\infty}(S/(\tilde S_I))^0$, that is 
consists of $((M_I,F),s_{IJ})\in C_{\mathcal D(2)fil}(S/(\tilde S_I))$ such that 
\begin{equation*}
H^n((M_I,F),s_{IJ}):=(H^n(M_I,F),H^n(s_{IJ}))\in\PSh^0_{\mathcal D}(S/(\tilde S_I)).
\end{equation*}
and 
\begin{eqnarray*}
\iota^0_S:=\iota^0_{S/(\tilde S_I)}:D^0_{\mathcal D(2)fil,\infty}(S)\xrightarrow{\sim}
D_{\mathcal D(2)fil,\infty}(S/(\tilde S_I))^0
\end{eqnarray*}
the induced equivalence of categories.
\end{prop}

\begin{proof}
Follows immediately from the fact that for $((M_I,F),s_{IJ}),((N_I,F),r_{IJ})\in C^0_{\mathcal D(2)fil}(S)$
\begin{eqnarray*}
\Hom_{D^0_{\mathcal D(2)fil,\infty}(S)}(((M_I,F),s_{IJ}),((N_I,F),r_{IJ}))&=&
\Hom_{K_{\mathcal D(2)fil,\infty}(S/(\tilde S_I))}((L(M_I,F),s^q_{IJ}),(E(N_I,F),r^k_{IJ})) \\
&=&\Hom_{D_{\mathcal D(2)fil,\infty}(S/(\tilde S_I))}(((M_I,F),s_{IJ}),((N_I,F),r_{IJ}))
\end{eqnarray*}

\end{proof}

We finish this subsection by the statement a result of kashiwara in the singular case.

\begin{defi}
Let $S\in\AnSp(\mathbb C)$ and $S=\cup_{i=1}^l S_i$ an open cover such that there exist closed embeddings
$i_i:S_i\hookrightarrow\tilde S_i$ with $\tilde S_i\in\AnSm(\mathbb C)$.
We will consider the functor
\begin{eqnarray*}
J_S:C_{\mathcal D(2)fil}(S/(\tilde S_I))\to C_{\mathcal D^{\infty}(2)fil}(S/(\tilde S_I)), \\
((M_I,F),u_{IJ})\mapsto J_S((M_I,F),u_{IJ}):=(J_{\tilde S_I}(M_I,F),J(u_{IJ})):=((M_I\otimes_{D_S}D_S^{\infty},F),J(u_{IJ}))
\end{eqnarray*}
with, denoting for short $d_{IJ}:=d_{\tilde S_J}-d_{\tilde S_I}$, 
\begin{eqnarray*}
J(u_{IJ}):J(M_I,F)\xrightarrow{J(u_{IJ})}J(p_{IJ*}(M_J,F)[d_{IJ}])\xrightarrow{T_*(p_{IJ},J)(-)} p_{IJ*}J(M_J,F)[d_{IJ}].
\end{eqnarray*}
Of course $J_S(C_{\mathcal D(1,0)fil}(S/(\tilde S_I)))\subset C_{\mathcal D^{\infty}(1,0)fil}(S/(\tilde S_I))$.
\end{defi}

\begin{prop}\label{Jpropsing}
Let $S\in\AnSp(\mathbb C)$ and $S=\cup_{i=1}^l S_i$ an open cover such that there exist closed embeddings
$i_i:S_i\hookrightarrow\tilde S_i$ with $\tilde S_i\in\AnSm(\mathbb C)$. Then the functor
\begin{equation*}
J_S:C_{\mathcal D(2)fil}(S/(\tilde S_I))\to C_{\mathcal D^{\infty}(2)fil}(S/(\tilde S_I)),
\end{equation*}
satisfy $J_S:C^{\sim}_{\mathcal D(2)fil}(S/(\tilde S_I))\subset C^{\sim}_{\mathcal D^{\infty}(2)fil}(S/(\tilde S_I))$ and
induces an equivalence of category
\begin{equation*}
J_S:D_{\mathcal D(2)fil,\infty,rh}(S/(\tilde S_I))\to D_{\mathcal D^{\infty}(2)fil,\infty,h}(S/(\tilde S_I)).
\end{equation*}
and $J_S(D_{\mathcal D(1,0)fil,\infty,rh}(S/(\tilde S_I))\subset D_{\mathcal D^{\infty,h}(1,0)fil,\infty}(S/(\tilde S_I))$.
\end{prop}

\begin{proof}
Follows immediately from the smooth case (proposition \ref{Jprop}).
\end{proof}

\subsubsection{Duality in the singular case}

The definition of Saito's category comes with a dual functor :

\begin{defi}\label{dualsing}
Let $S\in\Var(\mathbb C)$ and let $S=\cup S_i$ an open cover such that there exist closed embeddings
$i_i:S_i\hookrightarrow\tilde S_i$ with $\tilde S_i\in\SmVar(\mathbb C)$ ;
or let $S\in\AnSp(\mathbb C)$ and $S=\cup S_i$ an open cover such that there exist
closed embedding $i_i:S_i\hookrightarrow\tilde S_i$ with $\tilde S_i\in\AnSm(\mathbb C)$.
We have the dual functor : 
\begin{equation*}
\mathbb D_S^K:C^0_{\mathcal Dfil}(S/(\tilde S_I))\to C^0_{\mathcal Dfil}(S/(\tilde S_I)), \;
((M_I,F),s_{IJ})\mapsto(\mathbb D^K_{\tilde S_I}(M_I,F),s^d_{IJ}),
\end{equation*}
with, denoting for short $d_{IJ}:=d_{\tilde S_J}-d_{\tilde S_I}$,
\begin{eqnarray*}
u^q_{IJ}:\mathbb D^K_{\tilde S_I}(M_I,F)\xrightarrow{\mathbb D^K(s_{IJ}^{-1})}
\mathbb D^K_{\tilde S_I}p_{IJ*}(M_J,F)[d_{IJ}]\xrightarrow{T_*(p_{IJ},D)(-)}
p_{IJ*}\mathbb D^K_{\tilde S_J}(M_J,F)[d_{IJ}]
\end{eqnarray*}
It induces in the derived category the functor
\begin{equation*}
L\mathbb D_S^K:D^0_{\mathcal Dfil}(S/(\tilde S_I))\to D^0_{\mathcal Dfil}(S/(\tilde S_I)), \;
((M_I,F),s_{IJ})\mapsto\mathbb D^K_SQ((M_I,F),s_{IJ}),
\end{equation*}
with $q:Q((M_I,F),s_{IJ})\to((M_I,F),s_{IJ})$ a projective resolution.
\end{defi}

In the analytic case we also define

\begin{defi}\label{dualsinginfty}
Let $S\in\AnSp(\mathbb C)$ and $S=\cup S_i$ an open cover such that there exist
closed embedding $i_i:S_i\hookrightarrow\tilde S_i$ with $\tilde S_i\in\AnSm(\mathbb C)$.
We have the dual functor : 
\begin{equation*}
\mathbb D_S^{K,\infty}:C_{\mathcal D^{\infty}fil}(S/(\tilde S_I))\to C_{\mathcal D^{\infty}fil}(S/(\tilde S_I)), \;
((M_I,F),u_{IJ})\mapsto(\mathbb D^{K,\infty}_{\tilde S_I}(M_I,F),u^d_{IJ}),
\end{equation*}
with $u_{IJ}^d$ defined similarly as in definition \ref{dualsing}.
It induces in the derived category the functor
\begin{equation*}
L\mathbb D_S^{K,\infty}:D_{\mathcal D^{\infty}fil}(S/(\tilde S_I))\to D_{\mathcal D^{\infty}fil}(S/(\tilde S_I)), \;
((M_I,F),u_{IJ})\mapsto(\mathbb D^{K,\infty}_SQ((M_I,F),u^{q,d}_{IJ}),
\end{equation*}
with $q:Q((M_I,F),s_{IJ})\to((M_I,F),s_{IJ})$ a projective resolution.
\end{defi}

\subsubsection{Inverse image in the singular case}

We give in this subsection the inverse image functors between our categories.

Let $n:S^o\hookrightarrow S$ be an open embedding with $S\in\Var(\mathbb C)$ and
let $S=\cup_i S_i$ an open cover such that there exist closed embeddings
$i_i:S_i\hookrightarrow\tilde S_i$ with $\tilde S_i\in\SmVar(\mathbb C)$ ; or let 
$n:S^o\hookrightarrow S$ be an open embedding with $S\in\AnSp(\mathbb C)$
and let $S=\cup_i S_i$ an open cover such that there exist closed embeddings
$i_i:S_i\hookrightarrow\tilde S_i$ with $\tilde S_i\in\AnSm(\mathbb C)$. 
Denote $S^o_I:=n^{-1}(S_I)=S_I\cap S^o$ and $n_I:=n_{|S^o_I}:S^o_I\hookrightarrow S^o$ the open embeddings. 
Consider open embeddings $\tilde n_I:\tilde S^o_I\hookrightarrow\tilde S_I$ such that $\tilde S^o_I\cap S_I=S^o_I$,
that is which are lift of $n_I$. We have the functor
\begin{eqnarray*}
n^*:C_{\mathcal Dfil}(S/(\tilde S_I))\to C_{\mathcal Dfil}(S^o/(\tilde S^o_I)), \\ 
(M,F)=((M_I,F),u_{IJ})\mapsto n^*(M,F):=(\tilde n_I)^*(M,F):=(\tilde n_I^*(M_I,F),n^*u_{IJ})
\end{eqnarray*}
which derive trivially.

Let $f:X\to S$ be a morphism, with $X,S\in\Var(\mathbb C)$, 
such that there exist a factorization $f;X\xrightarrow{l}Y\times S\xrightarrow{p_S}S$
with $Y\in\SmVar(\mathbb C)$, $l$ a closed embedding and $p_S$ the projection,
and consider $S=\cup_{i=1}^lS_i$ an open cover such that there exist closed embeddings 
$i_i:S_i\hookrightarrow\tilde S_i$, with $\tilde S_i\in\SmVar(\mathbb C)$ ;
or let $f:X\to S$ be a morphism, with $X,S\in\AnSp(\mathbb C)$, 
such that there exist a factorization $f;X\xrightarrow{l}Y\times S\xrightarrow{p_S}S$
with $Y\in\AnSm(\mathbb C)$, $l$ a closed embedding and $p_S$ the projection
and consider $S=\cup_{i=1}^lS_i$ an open cover such that there exist closed embeddings 
$i_i:S_i\hookrightarrow\tilde S_i$, with $\tilde S_i\in\AnSm(\mathbb C)$.
Then, $X=\cup^l_{i=1} X_i$ with $X_i:=f^{-1}(S_i)$.
Denote by $p_{IJ}:\tilde S_J\to\tilde S_I$ and $p'_{IJ}:Y\times\tilde S_J\to Y\times\tilde S_I$ the projections and by
\begin{equation*}
E_{IJ}=\xymatrix{\tilde S_J\backslash S_J\ar[r]^{m_J}\ar[d]^{p_{IJ}} & \tilde S_J\ar[d]^{p_{IJ}} \\
\tilde S_I\backslash(S_I\backslash S_J)\ar[r]^{m=m_{IJ}} & \tilde S_I}, \;
E'_{IJ}=\xymatrix{Y\times\tilde S_J\backslash X_J\ar[r]^{m'_J}\ar[d]^{p'_{IJ}} & Y\times\tilde S_J\ar[d]^{p'_{IJ}} \\
Y\times\tilde S_I\backslash(X_I\backslash X_J)\ar[r]^{m'=m'_{IJ}} & Y\times\tilde S_I}, \;
E_{fIJ}\xymatrix{\tilde X_J\ar[r]^{\tilde f_J}\ar[d]^{p'_{IJ}} & \tilde S_J\ar[d]^{p_{IJ}} \\
Y\times\tilde S_I\ar[r]^{\tilde f_I} & \tilde S_I}
\end{equation*}
the commutative diagrams. The (graph) inverse image functors is :
\begin{eqnarray*}
f^{*mod[-],\Gamma}:C_{\mathcal Dfil}(S/(\tilde S_I))\to C_{\mathcal Dfil}(X/(Y\times\tilde S_I)), \\  
(M,F)=((M_I,F),u_{IJ})\mapsto  
f^{*mod[-],\Gamma}(M,F):=(\Gamma_{X_I}E(p_{\tilde S_I}^{*mod[-]}(M_I,F)),\tilde f_J^{*mod[-]}u_{IJ})
\end{eqnarray*}
with, denoting for short $d_{IJ}:=d_{\tilde S_I}-d_{\tilde S_J}$,
\begin{eqnarray*}
\tilde f_J^{*mod[-]}u_{IJ}:\Gamma_{X_I}E(p_{\tilde S_I}^{*mod[-]}(M_I,F))
\xrightarrow{\Gamma_{X_I}E(p_{\tilde S_I}^{*mod[-]})(u_{IJ})}\Gamma_{X_I}E(p_{\tilde S_I}^{*mod[-]}p_{IJ*}(M_J,F)[d_{IJ}]) \\
\xrightarrow{\Gamma_{X_I}E(T(p_{IJ}^{*mod},p_{\tilde S_I})(-)^{-1})[d_Y+d_{IJ}]}
\Gamma_{X_I}E(p'_{IJ*}p_{\tilde S_J}^{*mod}(M_J,F)[d_Y+d_{IJ}]) \\
\xrightarrow{=}p'_{IJ*}\Gamma_{X_J}E(p_{\tilde S_J}^{*mod[-]}(M_J,F))[d_{IJ}].
\end{eqnarray*}
It induces in the derived categories the functor
\begin{eqnarray*}
Rf^{*mod[-],\Gamma}:D_{\mathcal D(2)fil,\infty}(S/(\tilde S_I))\to D_{\mathcal D(2)fil,\infty}(X/(Y\times\tilde S_I)), \\   
(M,F)=((M_I,F),u_{IJ})\mapsto  
f^{*mod[-],\Gamma}(M,F):=(\Gamma_{X_I}E(p_{\tilde S_I}^{*mod[-]}(M_I,F)),\tilde f_J^{*mod[-]}u_{IJ}).
\end{eqnarray*}
It gives by duality the functor
\begin{eqnarray*}
Lf^{\hat*mod[-],\Gamma}:D_{\mathcal D(2)fil,\infty}(S/(\tilde S_I))^0\to 
D_{\mathcal D(2)fil,\infty}(X/(Y\times\tilde S_I))^0, \\  
(M,F)=((M_I,F),u_{IJ})\mapsto 
Lf^{\hat*mod[-],\Gamma}(M,F):=L\mathbb D^K_SRf^{*mod[-],\Gamma}L\mathbb D^K_S\iota_S^{0,-1}(M,F).
\end{eqnarray*}
where 
$\iota^0_S:D^0_{\mathcal D(2)fil,\infty}(S/(\tilde S_I))\xrightarrow{\sim}D_{\mathcal D(2)fil,\infty}(S/(\tilde S_I))^0$
is the isomorphism of definition \ref{iota0S}.

Let $f:X\to S$ be a morphism, with $X,S\in\AnSp(\mathbb C)$, 
such that there exist a factorization $f;X\xrightarrow{l}Y\times S\xrightarrow{p_S}S$
with $Y\in\AnSm(\mathbb C)$, $l$ a closed embedding and $p_S$ the projection
and consider $S=\cup_{i=1}^lS_i$ an open cover such that there exist closed embeddings 
$i_i:S_i\hookrightarrow\tilde S_i$, with $\tilde S_i\in\AnSm(\mathbb C)$.
Then, $X=\cup^l_{i=1} X_i$ with $X_i:=f^{-1}(S_i)$.
We have also the functors,
\begin{eqnarray*}
f^{*mod[-],\Gamma}:C_{\mathcal D^{\infty}fil}(S/(\tilde S_I))\to C_{\mathcal D^{\infty}fil}(X/(Y\times\tilde S_I)), \\ 
(M,F)=((M_I,F),u_{IJ})\mapsto  
f^{*mod[-],\Gamma}(M,F):=(\Gamma_{X_I}(p_{\tilde S_I}^{*mod[-]}(M_I,F),\tilde f_J^{*mod[-]}u_{IJ}))
\end{eqnarray*}
with, denoting for short $d_{IJ}:=d_{\tilde S_J}-d_{\tilde S_I}$,
\begin{eqnarray*}
\tilde f_J^{*mod[-]}u_{IJ}:\Gamma_{X_I}E(p_{\tilde S_I}^{*mod[-]}(M_I,F))
\xrightarrow{\Gamma_{X_I}E(p_{\tilde S_I}^{*mod[-]})(u_{IJ})}\Gamma_{X_I}E(p_{\tilde S_I}^{*mod[-]}p_{IJ*}(M_J,F)[d_{IJ}]) \\
\xrightarrow{\Gamma_{X_I}E(T(p_{IJ}^{*mod},p_{\tilde S_I})(-)^{-1})[d_Y+d_{IJ}]}
\Gamma_{X_I}E(p'_{IJ*}p_{\tilde S_J}^{*mod}(M_J,F)[d_Y+d_{IJ}]) \\
\xrightarrow{=}p'_{IJ*}\Gamma_{X_J}E(p_{\tilde S_J}^{*mod[-]}(M_J,F))[d_{IJ}].
\end{eqnarray*}
It induces in the derived categories, the functor
\begin{eqnarray*}
Rf^{*mod[-],\Gamma}:D_{\mathcal D^{\infty}(2)fil,\infty}(S/(\tilde S_I))\to 
D_{\mathcal D^{\infty}(2)fil,\infty}(X/(Y\times\tilde S_I)), \\  
(M,F)=((M_I,F),u_{IJ})\mapsto  
Rf^{*mod[-],\Gamma}(M,F):=(\Gamma_{X_I}E(p_{\tilde S_I}^{*mod[-]}(M_I,F)),\tilde f_J^{*mod[-]}u_{IJ}).
\end{eqnarray*}
It gives by duality the functor
\begin{eqnarray*}
Lf^{\hat*mod[-],\Gamma}:D_{\mathcal D^{\infty}(2)fil,\infty}(S/(\tilde S_I))^0\to 
D_{\mathcal D^{\infty}(2)fil,\infty}(X/(Y\times\tilde S_I))^0, \\  
(M,F)=((M_I,F),u_{IJ})\mapsto Lf^{\hat*mod[-],\Gamma}(M,F):=
L\mathbb D^{K,\infty}Rf^{*mod[-],\Gamma}L\mathbb D^{K,\infty}\iota_S^{0,-1}(M,F).
\end{eqnarray*}
where 
$\iota^0_S:D^0_{\mathcal D^{\infty}(2)fil,\infty}(S/(\tilde S_I))\xrightarrow{\sim}
D_{\mathcal D^{\infty}(2)fil,\infty}(S/(\tilde S_I))^0$
is the isomorphism of definition \ref{iota0S}.

The following proposition are then easy :

\begin{prop}\label{compDImod}
Let $f_1:X\to Y$ and $f_2:Y\to S$ two morphism with $X,Y,S\in\Var(\mathbb C)$. 
Assume there exist factorizations $f_1:X\xrightarrow{l_1} Y'\times Y\xrightarrow{p_Y} Y$ and 
$f_2:Y\xrightarrow{l_2} Y''\times S\xrightarrow{p_S} S$ with $Y',Y''\in\SmVar(\mathbb C)$,
$l_1,l_2$ closed embeddings and $p_S,p_Y$ the projections. We have then the factorization 
\begin{equation*}
f_2\circ f_1:X\xrightarrow{(l_2\circ I_{Y'})\circ l_1}Y'\times Y''\times S\xrightarrow{p_S} S.
\end{equation*}
We have, for $(M,F)\in C^{\sim}_{\mathcal D(2)fil}(S/(\tilde S_I))$,
$R(f_2\circ f_1)^{*mod[-],\Gamma}(M,F)=Rf_{2}^{*mod[-],\Gamma}\circ Rf_{1}^{*mod[-],\Gamma}(M,F)$.
\end{prop}

\begin{proof}
Follows from the the fact that 
for $(M,F)=((M_I,F),u_{IJ})\in C^{\sim}_{\mathcal D(2)fil}(S/(\tilde S_I))$,
\begin{eqnarray*}
(\Gamma_{X_I}E(\tilde f_{1I}^{*mod[-]}\Gamma_{Y_I}E(\tilde f_{2I}^{*mod[-]}(M_I,F))),
\tilde f_{1J}^{*mod[-]}(\tilde f_{2J}^{*mod[-]}u_{IJ}))
\xrightarrow{=} \\
(\Gamma_{X_I}E((\tilde f_{1I}\circ\tilde f_{2I})^{*mod[-]}(M_I,F)),(\tilde f_{1J}\circ\tilde f_{2J})^{*mod[-]}u^q_{IJ})
\end{eqnarray*}
by proposition \ref{compDmod}(i) and the fact that $X_I\subset\tilde f_{1I}^{-1}(Y_I)$.
\end{proof}

\begin{prop}\label{compDImodan}
Let $f_1:X\to Y$ and $f_2:Y\to S$ two morphism with $X,Y,S\in\Var(\mathbb C)$. 
Assume there exist factorizations $f_1:X\xrightarrow{l_1} Y'\times Y\xrightarrow{p_Y} Y$ and 
$f_2:Y\xrightarrow{l_2} Y''\times S\xrightarrow{p_S} S$ with $Y',Y''\in\SmVar(\mathbb C)$,
$l_1,l_2$ closed embeddings and $p_S,p_Y$ the projections. We have then the factorization 
\begin{equation*}
f_2\circ f_1:X\xrightarrow{(l_2\circ I_{Y'})\circ l_1}Y'\times Y''\times S\xrightarrow{p_S} S.
\end{equation*}
We have, for $(M,F)\in C^{\sim}_{\mathcal D(2)fil}(S/(\tilde S_I))$ or
$(M,F)\in C^{\sim}_{\mathcal D^{\infty}(2)fil}(S/(\tilde S_I))$,
$R(f_2\circ f_1)^{*mod[-],\Gamma}(M,F)=Rf_{2}^{*mod[-],\Gamma}\circ Rf_{1}^{*mod[-],\Gamma}(M,F)$.
\end{prop}

\begin{proof}
Similar to the proof of proposition \ref{compDImod}.
\end{proof}

\subsubsection{Direct image functor in the singular case}

We define the direct image functors between our category.

Let $f:X\to S$ be a morphism with $X,S\in\Var(\mathbb C)$,
and assume there exist a factorization $f:X\xrightarrow{l}Y\times S\xrightarrow{p_S}S$ 
with $Y\in\SmVar(\mathbb C)$, $l$ a closed embedding and $p_S$ a the projection ; 
or let $f:X\to S$ be a morphism with $X,S\in\AnSp(\mathbb C)$, and
assume there exist a factorization $f:X\xrightarrow{l}Y\times S\xrightarrow{p_S}S$ 
with $Y\in\AnSm(\mathbb C)$, $l$ a closed embedding and $p_S$ a the projection. 
Let $S=\cup_{i=1}^l S_i$ an open cover such that there exist closed embeddings
$i_i:S_i\hookrightarrow\tilde S_i$ with $\tilde S_i\in\SmVar(\mathbb C)$ ; resp. 
let $S=\cup_{i=1}^l S_i$ an open cover such that there exist closed embeddings
$i_i:S_i\hookrightarrow\tilde S_i$ with $\tilde S_i\in\AnSm(\mathbb C)$.
Then $X=\cup_{i=1}^lX_i$ with $X_i:=f^{-1}(S_i)$. 
Denote, for $I\subset\left[1,\cdots l\right]$, $S_I=\cap_{i\in I} S_i$ and $X_I=\cap_{i\in I}X_i$.
For $I\subset\left[1,\cdots l\right]$, denote by $\tilde S_I=\Pi_{i\in I}\tilde S_i$,
We then have, for $I\subset\left[1,\cdots l\right]$, closed embeddings $i_I:S_I\hookrightarrow\tilde S_I$
and the following commutative diagrams which are cartesian (we take $Y=\mathbb P^{N,o}$ in the algebraic case) 
\begin{equation*}
\xymatrix{
f_I=f_{|X_I}:X_I\ar[r]^{l_I}\ar[rd] & Y\times S_I\ar[r]^{p_{S_I}}\ar[d]^{i'_I} & S_I\ar[d]^{i_I} \\
\, & Y\times\tilde S_I\ar[r]^{p_{\tilde S_I}} & \tilde S_I} \;, \;
\xymatrix{Y\times\tilde S_J\ar[r]^{p_{\tilde S_J}}\ar[d]_{p'_{IJ}} & \tilde S_J\ar[d]^{p_{IJ}} \\
Y\times\tilde S_I\ar[r]^{p_{\tilde S_I}} & \tilde S_I}
\end{equation*}
with $l_I:l_{|X_I}$, $i'_I=I\times i_I$, $p_{S_I}$ and $p_{\tilde S_I}$ are the projections and $p'_{IJ}=I\times p_{IJ}$. 
Then $\tilde f_I:=p_{\tilde S_I}:Y\times\tilde S_I\to\tilde S_I$ is a lift of $f_I=f_{|X_I}$. 
We define the direct image functor on our category by
\begin{eqnarray*}
f^{FDR}_{*mod}:C_{\mathcal D(2)fil}(X/(Y\times\tilde S_I))\to C_{\mathcal D(2)fil}(S/(\tilde S_I)), \\ 
((M_I,F),u_{IJ})\mapsto (\tilde f^{FDR}_{I*mod}(M_I,F),f^k(u_{IJ})):=
(p_{\tilde S_I*}E((\Omega^{\bullet}_{Y\times\tilde S_I/\tilde S_I},F_b)\otimes_{O_{Y\times\tilde S_I}}(M_I,F)[d_Y]),f^k(u_{IJ}))
\end{eqnarray*}
with, denoting for short $d_{IJ}:=d_{\tilde S_J}-d_{\tilde S_I}$,  
\begin{eqnarray*}
f^k(u_{IJ})[d_Y]:
p_{\tilde S_I*}E((\Omega^{\bullet}_{Y\times\tilde S_I/\tilde S_I},F_b)\otimes_{O_{Y\times\tilde S_I}}(M_I,F)) \\
\xrightarrow{p_{\tilde S_I*}E(DR(Y\times\tilde S_I/\tilde S_I)(u_{IJ}))}
p_{\tilde S_I*}E((\Omega^{\bullet}_{Y\times\tilde S_I/\tilde S_I},F_b)\otimes_{O_{Y\times\tilde S_I}}p'_{IJ*}(M_J,F)[d_{IJ}]) \\
\xrightarrow{T^O_w(p_{IJ},\otimes)(M_I,F)}
p_{\tilde S_I*}E(p_{IJ*}(\Omega^{\bullet}_{Y\times\tilde S_J/\tilde S_J},F_b)\otimes_{O_{Y\times\tilde S_J}}(M_J,F)[d_{IJ}]) \\ 
\xrightarrow{=}
p_{\tilde S_J*}E((\Omega^{\bullet}_{Y\times\tilde S_J/\tilde S_J},F_b)\otimes_{O_{Y\times\tilde S_J}}(M_J,F))[d_{IJ}].
\end{eqnarray*}
It induces in the derived categories the functor
\begin{eqnarray*}
\int^{FDR}_f:D_{\mathcal D(2)fil,\infty}(X)\to D_{\mathcal D(2)fil,\infty}(S), \;  
((M_I,F),u_{IJ})\mapsto (\tilde f^{FDR}_{I*mod}(M_I,F),f^k(u_{IJ}))
\end{eqnarray*}

Let $f:X\to S$ be a morphism with $X,S\in\AnSp(\mathbb C)$, and
assume there exist a factorization $f:X\xrightarrow{l}Y\times S\xrightarrow{p_S}S$ 
with $Y\in\AnSm(\mathbb C)$, $l$ a closed embedding and $p_S$ a the projection. 
Let $S=\cup_{i=1}^l S_i$ an open cover such that there exist closed embeddings
$i_i:S_i\hookrightarrow\tilde S_i$ with $\tilde S_i=\in\AnSm(\mathbb C)$.
Then $X=\cup_{i=1}^lX_i$ with $X_i:=f^{-1}(S_i)$. 
We also have the functors
\begin{eqnarray*}
f^{FDR}_{*mod}:C_{\mathcal D^{\infty}(2)fil}(X/(Y\times\tilde S_I))\to C_{\mathcal D^{\infty}(2)fil}(S/(\tilde S_I)), \\ 
((M_I,F),u_{IJ})\mapsto (\tilde f^{FDR}_{I*mod}(M_I,F),f^k(u_{IJ})):=
(p_{\tilde S_I*}((\Omega^{\bullet}_{Y\times\tilde S_I/\tilde S_I},F_b)\otimes_{O_{Y\times\tilde S_I}}(M_I,F)[d_Y]),f^k(u_{IJ}))
\end{eqnarray*}
where  $f^k(u_{IJ})[d_Y]$ is given as above,
\begin{eqnarray*}
\int^{FDR}_f:D_{\mathcal D^{\infty}(2)fil,\infty}(X)\to D_{\mathcal D^{\infty}(2)fil,\infty}(S), \\ 
((M_I,F),u_{IJ})\mapsto (\tilde f^{FDR}_{I*mod}(M_I,F),f^k(u_{IJ})):=
(p_{\tilde S_I*}E((\Omega^{\bullet}_{Y\times\tilde S_I/\tilde S_I},F_b)\otimes_{O_{Y\times\tilde S_I}}(M_I,F)[d_Y]),f^k(u_{IJ}))
\end{eqnarray*}
where $f^k(u_{IJ})[d_Y]$ is given as above.

In the algebraic case, we have the followings:

\begin{prop}\label{compDmodiihsing}
Let $f_1:X\to Y$ and $f_2:Y\to S$ two morphism with $X,Y,S\in\QPVar(\mathbb C)$ quasi-projective. 
Then there exist factorizations $f_1:X\xrightarrow{l_1} Y'\times Y\xrightarrow{p_Y} Y$ and 
$f_2:Y\xrightarrow{l_2} Y''\times S\xrightarrow{p_S} S$ 
with $Y'=\mathbb P^{N,o}\subset\mathbb P^N$,$Y''=\mathbb P^{N',o}\subset\mathbb P^{N'}$ open subsets,
$l_1,l_2$ closed embeddings and $p_S,p_Y$ the projections.
We have then the factorization 
$f_2\circ f_1:X\xrightarrow{(l_2\circ I_{Y'})\circ l_1}Y'\times Y''\times S\xrightarrow{p_S} S$.
Let $i:S\hookrightarrow\tilde S$ a closed embedding with $\tilde S=\mathbb P^{n,o}\subset\mathbb P^n$ an open subset.
\begin{itemize}
\item[(i)]Let $(M,F)\in C_{\mathcal D(2)fil}(X/(Y'\times Y''\times\tilde S))$. 
Then, we have $\int^{FDR}_{f_2\circ f_1}(M,F)=\int^{FDR}_{f_2}(\int^{FDR}_{f_1}(M,F))$ 
in $D_{\mathcal D(2)fil,\infty}(S/(\tilde S_I))$.
\item[(ii)]Let $(M,F)\in C_{\mathcal D(2)fil,h}(X/(Y'\times Y''\times\tilde S))$.
Then, we have $\int^{FDR}_{(f_2\circ f_1)!}(M,F)=\int^{FDR}_{f_2!}(\int^{FDR}_{f_1!}(M,F))$ 
in $D_{\mathcal D(2)fil,\infty,h}(S/(\tilde S_I))$.
\end{itemize}
\end{prop}

\begin{proof}
\noindent(i):By the smooth case : proposition \ref{compDmodDRd}, we have en isomorphism
\begin{eqnarray*}
\int^{FDR}_{f_2}\int^{FDR}_{f_1}(M,F):=\int^{FDR}_{p_{\tilde S}}\int^{FDR}_{p_{Y'\times\tilde S}}(M,F)
\xrightarrow{\sim}\int^{FDR}_{p_{\tilde S}}(M,F):=\int^{FDR}_{(f_2\circ f_1)}(M,F).
\end{eqnarray*}

\noindent(ii):Follows from (i).
\end{proof}

In the analytic case, we have the followings:

\begin{prop}\label{compAnDmodiihsing}
Let $f_1:X\to Y$ and $f_2:Y\to S$ two morphism with $X,Y,S\in\AnSp(\mathbb C)$ quasi-projective. 
Then there exist factorizations $f_1:X\xrightarrow{l_1} Y'\times Y\xrightarrow{p_Y} Y$ and 
$f_2:Y\xrightarrow{l_2} Y''\times S\xrightarrow{p_S} S$ 
with $Y'=\mathbb P^{N,o}\subset\mathbb P^N$,$Y''=\mathbb P^{N',o}\subset\mathbb P^{N'}$ open subsets,
$l_1,l_2$ closed embeddings and $p_S,p_Y$ the projections.
We have then the factorization 
$f_2\circ f_1:X\xrightarrow{(l_2\circ I_{Y'})\circ l_1}Y'\times Y''\times S\xrightarrow{p_S} S$.
Let $i:S\hookrightarrow\tilde S$ a closed embedding with $\tilde S=\mathbb P^{n,o}\subset\mathbb P^n$ an open subset.
\begin{itemize}
\item[(i)]Let $(M,F)\in C_{\mathcal D^{\infty}(2)fil,h}(X/(Y'\times Y''\times\tilde S))$. 
Then, we have $\int^{FDR}_{f_2\circ f_1}(M,F)=\int^{FDR}_{f_2}(\int^{FDR}_{f_1}(M,F))$ 
in $D_{\mathcal D^{\infty}(2)fil,\infty}(S/(\tilde S_I))$.
\item[(ii)]Let $(M,F)\in C_{\mathcal D^{\infty}(2)fil,h}(X/(Y'\times Y''\times\tilde S))$.
Then, we have $\int^{FDR}_{(f_2\circ f_1)!}(M,F)=\int^{FDR}_{f_2!}(\int^{FDR}_{f_1!}(M,F))$ 
in $D_{\mathcal D^{\infty}(2)fil,\infty}(S/(\tilde S_I))$.
\end{itemize}
\end{prop}

\begin{proof}
\noindent(i): By the smooth case : proposition \ref{compAnDmodDRd}, we have en isomorphism
\begin{eqnarray*}
\int^{FDR}_{f_2}\int^{FDR}_{f_1}(M,F):=\int^{FDR}_{p_{\tilde S}}\int^{FDR}_{p_{Y'\times\tilde S}}(M,F)
\xrightarrow{\sim}\int^{FDR}_{p_{\tilde S}}(M,F):=\int^{FDR}_{(f_2\circ f_1)}(M,F).
\end{eqnarray*}

\noindent(ii):Follows from (i).
\end{proof}

\subsubsection{Tensor product in the singular case}

Let $S\in\Var(\mathbb C)$ and let $S=\cup S_i$ an open cover such that there exist closed embeddings
$i_i:S_i\hookrightarrow\tilde S_i$ with $\tilde S_i\in\SmVar(\mathbb C)$ ;
or let $S\in\AnSp(\mathbb C)$ and $S=\cup S_i$ an open cover such that there exist
closed embeddings $i_i:S_i\hookrightarrow\tilde S_i$ with $\tilde S_i\in\AnSm(\mathbb C)$.
We have the tensor product functors
\begin{eqnarray*}
(-)\otimes^{[-]}_{O_S}(-):C^2_{\mathcal Dfil}(S/(\tilde S_I))\to C_{\mathcal Dfil}(S/(\tilde S_I)), \\
(((M_I,F),u_{IJ}),((N_I,F),v_{IJ}))\mapsto((M_I,F)\otimes_{O_{\tilde S_I}}(N_I,F)[d_{\tilde S_I}],u_{IJ}\otimes v_{IJ}),
\end{eqnarray*}
with, denoting for short $d_{IJ}:=d_{\tilde S_J}-d_{\tilde S_I}$ and $d_I:=d_{\tilde S_I}$,
\begin{eqnarray*}
u_{IJ}\otimes v_{IJ}:(M_I,F)\otimes_{O_{\tilde S_I}}(N_I,F)[d_I]
\xrightarrow{T(p_{IJ}^{*mod},p_{IJ})(-)[d_I]}p_{IJ*}p_{IJ}^{*mod}((M_I,F)\otimes_{O_{\tilde S_I}}(N_I,F))[d_I] \\
\xrightarrow{=}p_{IJ*}(p_{IJ}^{*mod}(M_I,F)\otimes_{O_{\tilde S_J}}p_{IJ}^{*mod}(N_I,F))[d_I] \\
\xrightarrow{I(p_{IJ}^{*mod},p_{IJ})(-,-)(u_{IJ})\otimes I(p_{IJ}^{*mod},p_{IJ})(-,-)(v_{IJ})[d_I]}
p_{IJ*}((M_J,F)\otimes_{O_{\tilde S_J}}(N_J,F))[d_J+d_{IJ}].
\end{eqnarray*}

Let $S\in\AnSp(\mathbb C)$ and $S=\cup S_i$ an open cover such that there exist
closed embeddings $i_i:S_i\hookrightarrow\tilde S_i$ with $\tilde S_i\in\AnSm(\mathbb C)$.
We have the tensor product functors
\begin{eqnarray*}
(-)\otimes^{[-]}_{O_S}(-):C^2_{\mathcal D^{\infty}fil}(S/(\tilde S_I))\to C_{\mathcal D^{\infty}fil}(S/(\tilde S_I)), \\
(((M_I,F),u_{IJ}),((N_I,F),v_{IJ}))\mapsto((M_I,F)\otimes_{O_{\tilde S_I}}(N_I,F),u_{IJ}\otimes v_{IJ}),
\end{eqnarray*}
with $u_{IJ}\otimes v_{IJ}$ as above.

\begin{prop}\label{otimesedsing}
Let $S\in\Var(\mathbb C)$. Denote $\Delta_S:S\hookrightarrow S\times S$ the diagonal embedding.
Let $S=\cup S_i$ an open cover such that there exist
closed embeddings $i_i:S_i\hookrightarrow\tilde S_i$ closed embedding with $\tilde S_i\in\SmVar(\mathbb C)$ ;
or let $S\in\AnSp(\mathbb C)$ and $S=\cup S_i$ an open cover such that there exist
closed embedding $i_i:S_i\hookrightarrow\tilde S_i$ with $\tilde S_i\in\AnSm(\mathbb C)$.
We have, for $((M_I,F),u_{IJ}),((N_I,F),v_{IJ})\in C_{\mathcal Dfil}(S/(\tilde S_I))$,
\begin{equation*}
((M_I,F),u_{IJ})\otimes^{[-]}_{O_{S_I}}((N_I,F),v_{IJ})=\Delta_S^{*mod}(((M_I,F),u_{IJ}).((N_I,F),v_{IJ}))
\end{equation*}
\end{prop}

\begin{proof}
Follows from proposition \ref{otimesed}. 
\end{proof}

\subsubsection{The 2 functors of D modules on the category of complex algebraic varieties 
and on the category of complex analytic spaces, and the transformation maps}

\begin{defi}\label{TDmodlem0sing}
Consider a commutative diagram in $\Var(\mathbb C)$ which is cartesian :
\begin{equation*}
D=\xymatrix{X_T\ar[r]^{f'}\ar[d]^{g'} & T\ar[d]^{g} \\
X\ar[r]^{f} & S}.
\end{equation*}
Assume there exist factorizations $f:X\xrightarrow{l_1}Y_1\times S\xrightarrow{p_S}S$, 
$g:T\xrightarrow{l_2}Y_2\times S\xrightarrow{p_S}S$, with $Y_1,Y_2\in\SmVar(\mathbb C)$, 
$l_1,l_2$ closed embeddings and $p_S$, $p_S$ the projections.
Then, the above commutative diagram factors through
\begin{equation*}
D=\xymatrix{f':X_T\ar[r]^{l'_1}\ar[d]^{l'_2} & Y_1\times T\ar[d]^{l''_2=I\times l_2}\ar[r]^{p_T} & T\ar[d]^{l_2} \\
f'':X\times Y_2\ar[r]^{l_1''=I\times l}\ar[d]^{p_X} & Y_1\times Y_2\times S\ar[d]^{p_{Y_1\times S}}\ar[r]^{p_{Y_2\times S}} & 
Y_2\times S\ar[d]^{p_S} \\
f:X\ar[r]^{l} & Y_1\times S\ar[r]^{p_S} & S}.
\end{equation*}
whose squares are cartesian.
Let $S=\cup_i S_i$ be an open cover such that there exist closed embeddings $i_i:S_i\hookrightarrow\tilde S_i$ 
with $\tilde S_i\in\SmVar(\mathbb C)$. 
Then $X=\cup_iX_i$ and $T=\cup_iT_i$ with $X_i:=f^{-1}(S_i)$ and $T_i:=f^{-1}(S_i)$.
Moreover, $f_i=f_{|X_i}:X_i\to S_i$ lift to $\tilde f_i:=p_{\tilde S_i}:Y_1\times\tilde S_i\to\tilde S_i$
and $g_i=g_{|T_i}:T_i\to S_i$ lift to $\tilde g_i:=p_{\tilde S_i}:Y_2\times\tilde S_i\to\tilde S_i$.
We then have the following commutative diagram whose squares are cartesian
\begin{equation*}
\xymatrix{f':X_{IT}\ar[r]^{l'}\ar[d]^{l'_{2I}} & Y_1\times T_I\ar[d]^{l''_{2I}}\ar[r]^{p_T} & T_I\ar[d]^{l_{2I}} \\
X_I\times Y_2\ar[r]^{l_1''=I\times l_1}\ar[d]^{p_X} & 
Y_1\times Y_2\times\tilde S_I\ar[d]^{p_{Y_1\times\tilde S_I}}\ar[r]^{p_{Y_2\times\tilde S_I}} & 
Y_2\times\tilde S_I\ar[d]^{\tilde g_I} \\
X_I\ar[r]^{i_I\circ l_I} & Y_1\times\tilde S_I\ar[r]^{\tilde f_I} & \tilde S_I}
\end{equation*}
We then define, for $(M,F)=((M_I,F),u_{IJ})\in C_{\mathcal D(2)fil}(X/(Y_1\times\tilde S_I))$,
the following canonical transformation map in $D_{\mathcal D(2)fil,\infty}(T/(Y_2\times\tilde S_I))$,
using proposition \ref{mw0prop},
\begin{eqnarray*}
T^{\mathcal Dmod}(f,g)(M,F): \\
Rg^{*mod,\Gamma}\int^{FDR}_f(M,F):=(\Gamma_{T_I}E(\tilde g_I^{*mod}p_{\tilde S_I*}
E((\Omega^{\bullet}_{Y_1\times\tilde S_I/\tilde S_I},F_b)\otimes_{O_{Y_1\times\tilde S_I}}(M_I,F))),
\tilde g_J^{*mod}f^k(u_{IJ})) \\
\xrightarrow{(T^O_{\omega}(p_{\tilde S_I},\tilde g_I)(M_I,F))} \\
(\Gamma_{T_I}E(p_{Y_2\times\tilde S_I*}E((\Omega^{\bullet}_{Y_1\times Y_2\times\tilde S_I/Y_2\times\tilde S_I},F_b)
\otimes_{O_{Y_1\times Y_2\times\tilde S_I}}p_{Y_1\times\tilde S_I}^{*mod}(M_I,F))),
f^{'k}(p_{Y_1\times\tilde S_J}^{*mod}(u_{IJ}))) \\
\xrightarrow{(T^O_{\omega}(\gamma,\otimes)(p_{Y_1\times\tilde S_I}^{*mod}(M_I,F)))^{-1}} \\
(p_{Y_2\times\tilde S_I*}E((\Omega^{\bullet}_{Y_1\times Y_2\times\tilde S_I/Y_2\times\tilde S_I},F_b)
\otimes_{O_{Y_1\times Y_2\times\tilde S_I}}\Gamma_{Y_1\times T_I}E(p_{Y_1\times\tilde S_I}^{*mod}((M_I,F)))),
f^{'k}(\tilde g_J^{''*mod}(u^q_{IJ}))) \\
=:\int^{FDR}_{f'}Rg^{'*mod,\Gamma}(M,F).
\end{eqnarray*}
\end{defi}

In the analytic case, we have
\begin{defi}\label{TDmodlem0singan}
Consider a commutative diagram in $\AnSp(\mathbb C)$ which is cartesian :
\begin{equation*}
D=(f,g)=\xymatrix{X_T\ar[r]^{f'}\ar[d]^{g'} & T\ar[d]^{g} \\
X\ar[r]^{f} & S}.
\end{equation*}
Assume there exist factorizations $f:X\xrightarrow{l_1}Y_1\times S\xrightarrow{p_S}S$, 
$g:T\xrightarrow{l_2}Y_2\times S\xrightarrow{p_S}S$, with $Y_1,Y_2\in\AnSm(\mathbb C)$, 
$l_1,l_2$ closed embeddings and $p_S$, $p_S$ the projections.
\begin{itemize}
\item[(i)]We have, for $(M,F)\in D_{\mathcal D(2)fil,\infty,h}(X/(Y_1\times\tilde S_I))$, 
the following transformation map in $D_{\mathcal D(2)fil,\infty}(T/(Y_2\times\tilde S_I))$ 
\begin{eqnarray*}
T^{\mathcal Dmod}(f,g)((M,F)):Rg^{*mod,\Gamma}\int^{FDR}_{f}(M,F)\to\int^{FDR}_{f'}Rg^{'*mod,\Gamma}(M,F) 
\end{eqnarray*}
define in the same way as in definition \ref{TDmodlem0sing}
\item[(ii)] For $(M,F)\in D_{\mathcal D^{\infty}(2)fil,\infty}(X/(Y_1\times\tilde S_I))$,
 the following transformation map in $D_{\mathcal D^{\infty}(2)fil,\infty}(T/(Y_2\times\tilde S_I))$ 
\begin{equation*}
T^{\mathcal Dmod}(f,g)((M,F)): Rg^{*mod,\Gamma}\int^{FDR}_{f}(M,F)\to\int^{FDR}_{f'}Rg^{'*mod,\Gamma}(M,F)
\end{equation*}
is defined in the same way as in (ii) : see definition \ref{TDmodlem0sing}. 
\end{itemize}
\end{defi}

In the algebraic case, we have the following :

\begin{prop}\label{PDmod1sing}
Consider a commutative diagram in $\Var(\mathbb C)$
\begin{equation*}
D=(f,g)=\xymatrix{X_T\ar[r]^{f'}\ar[d]^{g'} & T\ar[d]^{g} \\
X\ar[r]^{f} & S}.
\end{equation*}
which is cartesian. Assume there exist factorizations $f:X\xrightarrow{l_1}Y_1\times S\xrightarrow{p_S}S$, 
$g:T\xrightarrow{l_2}Y_2\times S\xrightarrow{p_S}S$, with $Y_1,Y_2\in\SmVar(\mathbb C)$, 
$l_1,l_2$ closed embeddings and $p_S$, $p_S$ the projections.
For $(M,F)=((M_I,F),u_{IJ})\in C_{\mathcal D(2)fil,c}(X/(Y\times\tilde S_I))$,
\begin{equation*}
T^{\mathcal Dmod}(f,g):Rg^{*mod,\Gamma}\int^{FDR}_f(M,F)\to\int^{FDR}_{f'}Rg^{'*mod,\Gamma}(M,F)
\end{equation*}
is an isomorphism in $D_{\mathcal D(2)fil,\infty}(T/(Y_2\times\tilde S_I))$.
\end{prop}

\begin{proof}
Similar to the proof of proposition \ref{PDmod1}: the maps 
\begin{eqnarray*}
T^O_{\omega}(p_{\tilde S_I},\tilde g_I)(M_I,F):\tilde g_I^{*mod}p_{\tilde S_I*}
E((\Omega^{\bullet}_{Y\times\tilde S_I/\tilde S_I},F_b)\otimes_{O_{Y\times\tilde S_I}}(M_I,F))\to \\
p_{\tilde T_I*}E((\Omega_{Y\times\tilde T_I/\tilde T_I},F_b)\otimes_{O_{Y\times\tilde T_I}}\tilde g_I^{''*mod}(M_I,F))
\end{eqnarray*}
are $\infty$-filtered Zariski local equivalences since $\tilde g_I:Y_2\times\tilde S_I\to\tilde S_I$ are projections.
\end{proof}

\begin{prop}\label{PDmod1singan}
Consider a commutative diagram in $\AnSp(\mathbb C)$
\begin{equation*}
D=(f,g)=\xymatrix{X_T\ar[r]^{f'}\ar[d]^{g'} & T\ar[d]^{g} \\
X\ar[r]^{f} & S}.
\end{equation*}
which is cartesian. Assume that $f$ (hence $f'$) is proper and that 
there exist factorizations $f:X\xrightarrow{l_1}Y_1\times S\xrightarrow{p_S}S$, 
$g:T\xrightarrow{l_2}Y_2\times S\xrightarrow{p_S}S$, with $Y_1,Y_2\in\AnSm(\mathbb C)$, 
$l_1,l_2$ closed embeddings and $p_S$, $p_S$ the projections.
\begin{itemize}
\item[(i)] For $(M,F)=((M_I,F),u_{IJ})\in C_{\mathcal D(2)fil,h}(X/(Y_1\times\tilde S_I))$
\begin{equation*}
T^{\mathcal Dmod}(f,g):Rg^{*mod,\Gamma}\int^{FDR}_f(M,F)\to \int^{FDR}_{f'}Rg^{'*mod,\Gamma}(M,F)
\end{equation*}
is an isomorphism in $D_{\mathcal D(2)fil,\infty}(T/Y_2\times\tilde S_I)$.
\item[(ii)] For $(M,F)=((M_I,F),u_{IJ})\in C_{\mathcal D^{\infty}(2)fil,h}(X/(Y_1\times\tilde S_I))$
\begin{equation*}
T^{\mathcal Dmod}(f,g):Rg^{*mod,\Gamma}\int^{FDR}_f(M,F)\to \int^{FDR}_{f'}Rg^{'*mod,\Gamma}(M,F)
\end{equation*}
is an isomorphism in $D_{\mathcal D^{\infty}(2)fil,\infty}(T/(Y_2\times\tilde S_I))$.
\end{itemize}
\end{prop}

\begin{proof}
\noindent(i):Similar to the proof of proposition \ref{PDmod1sing}.

\noindent(ii):Similar to the proof of proposition \ref{PDmod1sing}.
\end{proof}

\begin{defi}\label{Tfansing} 
Let $f:X\to S$ be a morphism, with $X,S\in\Var(\mathbb C)$, 
such that there exist a factorization $f;X\xrightarrow{l}Y\times S\xrightarrow{p_S}S$
with $Y\in\SmVar(\mathbb C)$, $l$ a closed embedding and $p_S$ the projection,
and consider $S=\cup_{i=1}^lS_i$ an open cover such that there exist closed embeddings 
$i_i:S_i\hookrightarrow\tilde S_i$, with $\tilde S_i\in\SmVar(\mathbb C)$ ;
Then, $X=\cup^l_{i=1} X_i$ with $X_i:=f^{-1}(S_i)$.
We have, for $(M,F)=((M_I,F),u_{IJ})\in C_{\mathcal D(2)fil}(S/(\tilde S_I))$, 
the canonical transformation map in $D_{\mathcal D(2)fil}(T^{an}/(\tilde T^{an}_I))$
\begin{eqnarray*}
T^{mod}(an,\gamma_T)(M,F): \\
f^{*mod[-],\Gamma}(M,F))^{an}:=((\Gamma_{T_I}E(p_{\tilde S_I}^{*mod[-]}(M_I,F)))^{an},(f^{*mod[-]}u_{IJ})^{an}) \\
\xrightarrow{(T^{mod}(an,\gamma_{T_I})(-))}
(\Gamma_{T^{an}_I}E((p_{\tilde S_I}^{*mod[-]}(M_I,F))^{an}),f^{*mod[-]}u^{an}_{IJ}) \\
\xrightarrow{=}
(\Gamma_{T^{an}_I}E(p_{\tilde S_I}^{*mod[-]}(M_I^{an},F)),f^{*mod[-]}u^{an}_{IJ})=:f^{*mod[-],\Gamma}((M,F)^{an})
\end{eqnarray*}
where the equality is obvious (see proposition \ref{anmor}). 
\end{defi}

\begin{defi}
Let $f:X\to S$ a morphism with $X,S\in\Var(\mathbb C)$.
Assume there exist a factorization $f:X\xrightarrow{l}Y\times S\xrightarrow{p_S}S$ 
with $Y\in\SmVar(\mathbb C)$, $l$ a closed embedding and $p_S$ the projection. 
Let $S=\cup_{i=1}^l S_i$ be an open cover such that there exist closed embeddings
$i_i:S_i\hookrightarrow\tilde S_i$ closed embeddings with $\tilde S_i\in\SmVar(\mathbb C)$.
We have, for $(M,F)=((M_I,F),u_{IJ})\in C_{\mathcal Dfil}(X/Y\times\tilde S_I)$,
the following transformation map in $D_{\mathcal Dfil}(X^{an}/(Y\times\tilde S_I)^{an})$
\begin{eqnarray*}
T^{\mathcal Dmod}(an,f)(M,F):(\int^{FDR}_f(M,F))^{an}=
(p_{\tilde S_I*}E((\Omega^{\bullet}_{Y\times\tilde S_I/\tilde S_I},F_b)\otimes_{O_{Y\times\tilde S_I}}L_D(M_I,F)))^{an},
(f^k(u^q_{IJ}))^an) \\
\xrightarrow{(T^O_{\omega}(p_{\tilde S_I},an)(M_I,F))}
(p_{\tilde T_I*}E((\Omega^{\bullet}_{Y\times\tilde T_I/\tilde T_I},F_b)\otimes_{O_{(Y\times\tilde T_I)^{an}}}L_D(M_I,F)^{an}),
f^{'k}((u^q_{IJ})^{an}))=:\int^{FDR}_{f^{an}}(M,F)^{an}
\end{eqnarray*}
\end{defi}

\begin{thm}\label{GAGADmodsing}
Let $f:X\to S$ a morphism with $X,S\in\Var(\mathbb C)$.
Assume there exist a factorization $f:X\xrightarrow{l}Y\times S\xrightarrow{p_S}S$ 
with $Y\in\SmVar(\mathbb C)$, $l$ a closed embedding and $p_S$ the projection. 
Let $S=\cup_{i=1}^l S_i$ be an open cover such that there exist closed embeddings
$i_i:S_i\hookrightarrow\tilde S_i$ closed embeddings with $\tilde S_i\in\SmVar(\mathbb C)$.
Let $M\in D_{\mathcal Dfil,c}(X/Y\times\tilde S_I)$. If $f$ is proper, 
\begin{equation*}
T^{\mathcal D}(an,f)(M):(\int_{f}M)^{an}\xrightarrow{\sim}\int_{f^{an}}(M)^{an}
\end{equation*}
is an isomorphism.
\end{thm}

\begin{proof}
By theorem \ref{GAGADmod}, $T^O_{\omega}(p_{\tilde S_I},an)(M_I)$ are usu local equivalences.
\end{proof}

In the analytic case, we have the following canonical transformation maps

\begin{defi}\label{Tfinftysing}
Let $f:X\to S$ be a morphism, with $X,S\in\AnSp(\mathbb C)$, 
such that there exist a factorization $f;X\xrightarrow{l}Y\times S\xrightarrow{p_S}S$
with $Y\in\AnSm(\mathbb C)$, $l$ a closed embedding and $p_S$ the projection,
and consider $S=\cup_{i=1}^lS_i$ an open cover such that there exist closed embeddings 
$i_i:S_i\hookrightarrow\tilde S_i$, with $\tilde S_i\in\AnSm(\mathbb C)$ ;
Then, $X=\cup^l_{i=1} X_i$ with $X_i:=f^{-1}(S_i)$.
We have, for $(M,F)=((M_I,F),u_{IJ})\in C_{\mathcal D(2)fil}(S/(\tilde S_I))$, 
the canonical transformation map in $D_{\mathcal D^{\infty}fil}(T/(\tilde T_I))$ 
obtained by the canonical maps given in definition \ref{Tfinfty} and definition \ref{Tgammainfty} :
\begin{eqnarray*}
T(f,\infty)(M,F):
J_T(f^{*mod[-],\Gamma}(M,F)):=(J_{\tilde T_I}(\Gamma_{T_I}E(p_{\tilde S_I}^{*mod[-]}(M_I,F))),J(f^{*mod[-]}u_{IJ})) \\
\xrightarrow{(T(\infty,\gamma_{T_I})(-))}
(\Gamma_{T_I}E(J_{\tilde T_I}(p_{\tilde S_I}^{*mod[-]}(M_I,F))),N_{IJ}) \\
\xrightarrow{(T(p_{\tilde S_I},\infty)(-))}
(\Gamma_{T_I}E(p_{\tilde S_I}^{*mod[-]}J_{\tilde S_I}(M_I,F)),f^{*mod[-]}J(u_{IJ}))=:f^{*mod[-],\Gamma}(J_S(M,F))
\end{eqnarray*}
\end{defi}

\subsection{The category of complexes of quasi-coherent sheaves whose cohomology sheaves has a structure of D-modules}

\subsubsection{Definition on a smooth complex algebraic variety or smooth complex analytic space
and the functorialities}

Let $X\in\SmVar(\mathbb C)$ or let $X\in\AnSm(\mathbb C)$.
Recall that (see definition \ref{ODmod} section 4.1) $C_{O_Xfil,\mathcal D}(X)$ is the category 
\begin{itemize}
\item whose objects $(M,F)\in C_{O_Xfil,\mathcal D}(X)$ are filtered complexes of presheaves of $O_X$ modules
$(M,F)\in C_{O_Xfil}(X)$ whose cohomology presheaves $H^n(M,F)\in\PSh_{O_Xfil}(X)$ are emdowed with a structure
of filtered $D_X$ modules for all $n\in\mathbb Z$. 
\item whose set of morphisms $\Hom_{C_{O_Xfil,\mathcal D}(X)}((M,F),(N,F))\subset\Hom_{C_{O_Xfil}(X)}((M,F),(N,F))$ 
between $(M,F),(N,F)\in C_{O_Xfil,\mathcal D}(X)$
are the morphisms of filtered complexes of $O_X$ modules 
$m:(M,F)\to(N,F)$ such that $H^nm:H^n(M,F)\to H^n(N,F)$ is $D_X$ linear,
i.e. is a morphism of (filtered) $D_X$ modules, for all $n\in\mathbb Z$.
\end{itemize}
More generally, let $h:X\to S$ a morphism with $X,S\in\SmVar(\mathbb C)$ or with $X,S\in\AnSm(\mathbb C)$.
Then, $C_{h^*O_Sfil,h^*\mathcal D}(X)$ the category 
\begin{itemize}
\item whose objects $(M,F)\in C_{h^*O_Sfil,h^*\mathcal D}(X)$ are filtered complexes of presheaves of $h^*O_S$ modules
$(M,F)\in C_{h^*O_Sfil}(X)$ whose cohomology presheaves $H^n(M,F)\in\PSh_{h^*O_Sfil}(X)$ are emdowed with a structure
of filtered $h^*D_S$ modules for all $n\in\mathbb Z$. 
\item whose set of morphisms 
$\Hom_{C_{h^*O_Sfil,h^*\mathcal D}(X)}((M,F),(N,F))\subset\Hom_{C_{h^*O_Sfil}(X)}((M,F),(N,F))$ 
between $(M,F),(N,F)\in C_{h^*O_Sfil,h^*\mathcal D}(X)$
are the morphisms of filtered complexes of $h^*O_S$ modules $m:(M,F)\to(N,F)$ 
such that $H^nm:H^n(M,F)\to H^n(N,F)$ is $h^*D_S$ linear,
i.e. is a morphism of (filtered) $h^*D_S$ modules, for all $n\in\mathbb Z$.
\end{itemize}

\begin{defi}
Let $S\in\SmVar(\mathbb C)$ or $S\in\AnSm(\mathbb C)$. Let $Z\subset S$ a closed subset.
Denote by $j:S\backslash Z\hookrightarrow S$ the open complementary embedding.
\begin{itemize}
\item[(i)] We denote by $C_{O_S,\mathcal D,Z}(S)\subset C_{O_S,\mathcal D}(S)$
the full subcategory consisting of $M\in C_{O_S,\mathcal D}(S)$ such that
such that $j^*H^nM=0$ for all $n\in\mathbb Z$.
\item[(ii)] We denote by $C_{O_Sfil,\mathcal D,Z}(S)\subset C_{O_Sfil,\mathcal D}(S)$
the full subcategory consisting of $(M,F)\in C_{O_Sfil,\mathcal D}(S)$ such that there exist $r\in\mathbb N$
and an $r$-filtered homotopy equivalence $m:(M,F)\to(M',F)$ with $(M',F)\in C_{O_Sfil,\mathcal D}(S)$
such that $j^*H^n\Gr_F^p(M',F)=0$ for all $n,p\in\mathbb Z$.
\end{itemize}
\end{defi}

\begin{defi}
Let $S\in\SmVar(\mathbb C)$ or let $S\in\AnSm(\mathbb C)$.
We have then (see section 2), for $r=1,\cdots,\infty$, the homotopy category 
$K_{O_Sfil,\mathcal D,r}(S)=\Ho_r(C_{O_Sfil,\mathcal D}(S))$ whose objects are those of $C_{O_Sfil,\mathcal D}(S)$
and whose morphisms are $r$-filtered homotopy classes of morphism, and
its localization $D_{O_Sfil,\mathcal D,r}(S)=K_{O_Sfil,\mathcal D,r}(S)([E_1]^{-1})$ with respect to
filtered zariski, resp. usu local equivalence.
Note that the classes of filtered $\tau$ local equivalence constitute a right multiplicative system.
\end{defi}

We look at functoriality

\begin{itemize}

\item Let $S\in\SmVar(\mathbb C)$ or $S\in\AnSm(\mathbb C)$. Let $(M,F)\in C_{O_Sfil,\mathcal D}(S)$.
Then, the canonical morphism $q:L_O(M,F)\to(M,F)$ in $C_{O_Sfil}(S)$ being a quasi-isomorphism of $O_S$ modules, 
we get in a unique way $L_O(M,F)\in C_{O_Sfil,\mathcal D}(S)$ such that $q:L_O(M,F)\to(M,F)$ is a morphism
in $C_{O_Sfil,\mathcal D}(S)$

\item Let $f:X\to S$ be a morphism with $X,S\in\SmVar(\mathbb C)$,
or let $f:X\to S$ be a morphism with $X,S\in\AnSm(\mathbb C)$.
Let $(M,F)\in C_{O_Sfil,\mathcal D}(S)$. Then, $f^{*mod}H^n(M,F):=(O_X,F_b)\otimes_{f^*O_S}f^*H^n(M,F)$
is canonical a filtered $D_X$ module (see section 4.1 or 4.2).
Consider the canonical surjective map $q(f):H^nf^{*mod}(M,F)\to f^{*mod}H^n(M,F)$.
Then, $q(f)$ is an isomorphism if $f$ is smooth.
Let $h:U\to S$ be a smooth morphism with $U,S\in\SmVar(\mathbb C)$,
or let $h:U\to S$ be a smooth morphism with $U,S\in\AnSm(\mathbb C)$. We get the functor
\begin{eqnarray*}
h^{*mod}:C_{O_Sfil,\mathcal D}(S)\to C_{O_Ufil,\mathcal D}(U),
(M,F)\mapsto h^{*mod}(M,F), 
\end{eqnarray*} 

\item Let $S\in\SmVar(\mathbb C)$ or $S\in\AnSm(\mathbb C)$, and 
let $i:Z\hookrightarrow S$ a closed embedding and denote by $j:S\backslash Z\hookrightarrow S$ the open complementary. 
For $M\in C_{O_S,\mathcal D}(S)$, the cohomology presheaves of 
\begin{equation*}
\Gamma_Z M:=\Cone(\ad(j^*,j_*)(M):M\to j_*j^*M)[-1]
\end{equation*}
has a canonical $D_S$-module structure 
(as $j^*H^nM$ is a $j^*D_S$ module, $H^nj_*j^*M=j_*j^*H^nM$ has an induced structure of $D_S$ module),
and $\gamma_Z(M):\Gamma_Z M\to M$ is a map in $C_{O_S,\mathcal D}(S)$.
For $Z_2\subset Z$ a closed subset and $M\in C_{O_S,\mathcal D}(S)$, 
$T(Z_2/Z,\gamma)(M):\Gamma_{Z_2}M\to\Gamma_Z M$ is a map in $C_{O_S,\mathcal D}(S)$.
We get the functor
\begin{eqnarray*}
\Gamma_Z:C_{O_Sfil,\mathcal D}(S)\to C_{O_Sfil,\mathcal D}(S), \\ 
(M,F)\mapsto\Gamma_Z(M,F):=\Cone(\ad(j^*,j_*)((M,F)):(M,F)\to j_*j^*(M,F))[-1],
\end{eqnarray*}
together we the canonical map $\gamma_Z(M,F):\Gamma_Z(M,F)\to (M,F)$

More generally, let $h:Y\to S$ a morphism with $Y,S\in\Var(\mathbb C)$ or $Y,S\in\AnSp(\mathbb C)$, $S$ smooth, and 
let $i:X\hookrightarrow Y$ a closed embedding and denote by $j:Y\backslash X\hookrightarrow Y$ the open complementary. 
For $M\in C_{h^*O_S,h^*\mathcal D}(Y)$, 
\begin{equation*}
\Gamma_X M:=\Cone(\ad(j^*,j_*)(M):M\to j_*j^*M)[-1]
\end{equation*}
has a canonical $h^*D_S$-module structure,
(as $j^*H^nM$ is a $j^*h^*D_S$ module, $H^nj_*j^*M=j_*j^*H^nM$ has an induced structure of $j^*h^*D_S$ module),
and $\gamma_X(M):\Gamma_X M\to M$ is a map in $C_{h^*O_S,h^*\mathcal D}(Y)$.
For $X_2\subset X$ a closed subset and $M\in C_{h^*O_S,h^*\mathcal D}(Y)$, 
$T(Z_2/Z,\gamma)(M):\Gamma_{X_2}M\to\Gamma_X M$ is a map in $C_{h^*O_S,h^*\mathcal D}(Y)$. We get the functor
\begin{eqnarray*}
\Gamma_X:C_{h^*O_Sfil,h^*\mathcal D}(Y)\to C_{h^*O_Sfil,h^*\mathcal D}(Y), \\ 
(M,F)\mapsto\Gamma_X(M,F):=\Cone(\ad(j^*,j_*)((M,F)):(M,F)\to j_*j^*(M,F))[-1],
\end{eqnarray*}
together we the canonical map $\gamma_X(M,F):\Gamma_X(M,F)\to (M,F)$

\item Let $f:X\to S$ be a morphism with $X,S\in\SmVar(\mathbb C)$,
or let $f:X\to S$ be a morphism with $X,S\in\AnSm(\mathbb C)$. 
Consider the factorization $f:X\xrightarrow{l}X\times S\xrightarrow{p}S$,
where $l$ is the graph embedding and $p$ the projection. We get from the two preceding points the functor
\begin{eqnarray*}
f^{*mod,\Gamma}:C_{O_Sfil,\mathcal D}(S)\to C_{O_Xfil,\mathcal D}(X\times S),
(M,F)\mapsto f^{*mod,\Gamma}(M,F):=\Gamma_Xp^{*mod}(M,F), 
\end{eqnarray*} 
and
\begin{eqnarray*}
f^{*mod[-],\Gamma}:C_{O_Sfil,\mathcal D}(S)\to C_{O_Xfil,\mathcal D}(X\times S),
(M,F)\mapsto f^{*mod[-],\Gamma}(M,F):=\Gamma_XE(p^{*mod}(M,F))[-d_X], 
\end{eqnarray*} 
which induces in the derived categories the functor
\begin{eqnarray*}
Rf^{*mod[-],\Gamma}:D_{O_Sfil,\mathcal D}(S)\to D_{O_Xfil,\mathcal D}(X\times S),
(M,F)\mapsto Rf^{*mod[-],\Gamma}(M,F):=\Gamma_XE(p^{*mod[-]}(M,F)). 
\end{eqnarray*} 
For $(M,F)\in C_{O_Sfil,\mathcal D}(S)$ or $(M,F)\in C_{O_Sfil}(S)$, the canonical map in $C_{O_Xfil}(X\times S)$
\begin{eqnarray*}
\ad(i^{*mod},i_*)(-):L_O\Gamma_XE(p^{*mod}(M,F))\to i_*i^{*mod}L_O\Gamma_XE(p^{*mod}(M,F))
\end{eqnarray*}
gives in the derived category, the canonical map in $D_{O_Xfil,\infty}(X\times S)$
\begin{eqnarray*}
I(f^{*mod,\Gamma})(M,F):Rf^{*mod,\Gamma}(M,F)=L_O\Gamma_XE(p^{*mod}(M,F))\xrightarrow{\ad(i^{*mod},i_*)(-)} \\
i_*i^{*mod}L_O\Gamma_XE(p^{*mod}(M,F))\xrightarrow{\sim}i_*i^{*mod}L_O(p^{*mod}(M,F))=Lf^{*mod}(M,F)
\end{eqnarray*}
where the isomorphism is given by lemma \ref{imodj}.

\item Let $S\in\SmVar(\mathbb C)$. We have the analytical functor :
\begin{equation*}  
(-)^{an}:C_{O_Sfil,\mathcal D}(S)\to C_{O_Sfil,\mathcal D}(S^{an}), \; 
(M,F)\mapsto (M,F)^{an}:=\an_S^{*mod}(M,F):=(M,F)\otimes_{\an_S^*O_S}O_{S^{an}}
\end{equation*}
which induces in the derived category 
\begin{equation*}
(-)^{an}:D_{O_Sfil,\mathcal D,r}(S)\to D_{O_Sfil,\mathcal D,r}(S^{an}), \; ((M,F)\mapsto (M,F)^{an}:=\an_S^{*mod}(M,F))
\end{equation*}
since $\an_S^{*mod}$ is an exact functor.

\end{itemize}

We have, for $f:T\to S$ with $T,S\in\SmVar(\mathbb C)$ or with $T,S\in\AnSm(\mathbb C)$, the commutative diagrams of functors
\begin{equation*}
\xymatrix{C_{\mathcal Dfil}(S)\ar[r]^{o_S}\ar[d]^{f^{*mod[-],\Gamma}} & 
C_{Ofil,\mathcal D}(S)\ar[d]^{f^{*mod[-],\Gamma}} \\
C_{\mathcal Dfil}(T)\ar[r]^{o_T} & C_{Ofil,\mathcal D}(T)} \; , \; 
\xymatrix{D_{\mathcal Dfil,r}(S)\ar[r]^{o_S}\ar[d]^{Rf^{*mod[-],\Gamma}} & 
D_{Ofil,\mathcal D,r}(S)\ar[d]^{Rf^{*mod[-],\Gamma}} \\
D_{\mathcal Dfil,r}(T)\ar[r]^{o_T} & D_{Ofil,\mathcal D,r}(T)}
\end{equation*}
where $o_S$ and $o_T$ are the forgetfull functors.

\subsubsection{Definition on a singular complex algebraic variety or singular complex analytic space
and the functorialities}

\begin{defi}\label{DmodsingVardef}
Let $S\in\Var(\mathbb C)$ and let $S=\cup_iS_i$ an open cover
such that there exist closed embeddings $i_i:S_i\hookrightarrow\tilde S_i$ with $\tilde S_i\in\SmVar(\mathbb C)$ ;
or let $S\in\AnSp(\mathbb C)$ and let $S=\cup_iS_i$ an open cover 
such that there exist closed embeddings $i_i:S_i\hookrightarrow\tilde S_i$ with $\tilde S_i\in\AnSm(\mathbb C)$. Then,
$C_{Ofil,\mathcal D}(S/(\tilde S_I))$ is the category
\begin{itemize}
\item whose objects are $(M,F)=((M_I,F)_{I\subset\left[1,\cdots l\right]},u_{IJ})$, with
\begin{itemize}
\item $(M_I,F)\in C_{O_{\tilde S_I}fil\mathcal D,S_I}(\tilde S_I)$,
\item $u_{IJ}:m^*(M_I,F)\to m^*p_{IJ*}(M_J,F)[d_{\tilde S_J}-d_{\tilde S_I}]$ 
for $J\subset I$, are morphisms, $p_{IJ}:\tilde S_J\to\tilde S_I$ being the projection, 
satisfying for $I\subset J\subset K$, $p_{IJ*}u_{JK}\circ u_{IJ}=u_{IK}$ in $C_{O_{\tilde S_I}fil,\mathcal D}(\tilde S_I)$ ;
\end{itemize}
\item whose morphisms $m:((M_I,F),u_{IJ})\to((N_I,F),v_{IJ})$ between  
$(M,F)=((M_I,F)_{I\subset\left[1,\cdots l\right]},u_{IJ})$ and $(N,F)=((N_I,F)_{I\subset\left[1,\cdots l\right]},v_{IJ})$
are a family of morphisms of complexes,  
\begin{equation*}
m=(m_I:(M_I,F)\to (N_I,F))_{I\subset\left[1,\cdots l\right]}
\end{equation*}
such that $v_{IJ}\circ m_I=p_{IJ*}m_J\circ u_{IJ}$ in $C_{O_{\tilde S_I}fil,\mathcal D}(\tilde S_I)$.
\end{itemize}
We denote by $C^{\sim}_{Ofil,\mathcal D}(S/(\tilde S_I))\subset C_{Ofil,\mathcal D}(S/(\tilde S_I))$ the full subcategory 
consisting of objects $((M_I,F),u_{IJ})$ such that the $u_{IJ}$ are $\infty$-filtered Zariski, resp. usu, local equivalences.
\end{defi}

\begin{defi}
Let $S\in\Var(\mathbb C)$ and let $S=\cup_iS_i$ an open cover
such that there exist closed embeddings $i_i:S_i\hookrightarrow\tilde S_i$ with $\tilde S_i\in\SmVar(\mathbb C)$ ;
or let $S\in\AnSp(\mathbb C)$ and let $S=\cup_iS_i$ an open cover 
such that there exist closed embeddings $i_i:S_i\hookrightarrow\tilde S_i$ with $\tilde S_i\in\AnSm(\mathbb C)$.
We have then (see section 2), for $r=1,\cdots,\infty$, the homotopy category 
\begin{equation*}
K_{Ofil,\mathcal D,r}(S/(\tilde S_I)):=\Ho_r(C_{Ofil,\mathcal D}(S/(\tilde S_I))) 
\end{equation*}
whose objects are those of $C_{Ofil,\mathcal D}(S/(\tilde S_I))$
and whose morphisms are $r$-filtered homotopy classes of morphism, and its localization 
\begin{equation*}
D_{fil,\mathcal D,r}(S/(\tilde S_I)):=K_{Ofil,\mathcal D,r}(S/(\tilde S_I))([E_1]^{-1})
\end{equation*}
with respect to the classes of filtered zariski, resp. usu local equivalence.
Note that the classes of filtered $\tau$ local equivalence constitute a right multiplicative system.
\end{defi}

Let $f:X\to S$ be a morphism, with $X,S\in\Var(\mathbb C)$, 
such that there exist a factorization $f;X\xrightarrow{l}Y\times S\xrightarrow{p_S}S$
with $Y\in\SmVar(\mathbb C)$, $l$ a closed embedding and $p_S$ the projection,
and consider $S=\cup_{i=1}^lS_i$ an open cover such that there exist closed embeddings 
$i_i:S_i\hookrightarrow\tilde S_i$, with $\tilde S_i\in\SmVar(\mathbb C)$ ;
or let $f:X\to S$ be a morphism, with $X,S\in\AnSp(\mathbb C)$, 
such that there exist a factorization $f;X\xrightarrow{l}Y\times S\xrightarrow{p_S}S$
with $Y\in\AnSm(\mathbb C)$, $l$ a closed embedding and $p_S$ the projection
and consider $S=\cup_{i=1}^lS_i$ an open cover such that there exist closed embeddings 
$i_i:S_i\hookrightarrow\tilde S_i$, with $\tilde S_i\in\AnSm(\mathbb C)$.
Then, $X=\cup^l_{i=1} X_i$ with $X_i:=f^{-1}(S_i)$.
Denote by $p_{IJ}:\tilde S_J\to\tilde S_I$ and $p'_{IJ}:Y\times\tilde S_J\to Y\times\tilde S_I$ the projections and by
\begin{equation*}
E_{IJ}=\xymatrix{\tilde S_J\backslash S_J\ar[r]^{m_J}\ar[d]^{p_{IJ}} & \tilde S_J\ar[d]^{p_{IJ}} \\
\tilde S_I\backslash(S_I\backslash S_J)\ar[r]^{m=m_{IJ}} & \tilde S_I}, \;
E'_{IJ}=\xymatrix{Y\times\tilde S_J\backslash X_J\ar[r]^{m'_J}\ar[d]^{p'_{IJ}} & Y\times\tilde S_J\ar[d]^{p'_{IJ}} \\
Y\times\tilde S_I\backslash(X_I\backslash X_J)\ar[r]^{m'=m'_{IJ}} & Y\times\tilde S_I}, \;
E_{fIJ}\xymatrix{\tilde X_J\ar[r]^{\tilde f_J}\ar[d]^{p'_{IJ}} & \tilde S_J\ar[d]^{p_{IJ}} \\
Y\times\tilde S_I\ar[r]^{\tilde f_I} & \tilde S_I}
\end{equation*}
the commutative diagrams. We then have the filtered De Rham the inverse image functor :
\begin{eqnarray*}
f^{*mod[-],\Gamma}:C_{Ofil,\mathcal D}(S/(\tilde S_I))\to C_{Ofil,\mathcal D}(X/(Y\times\tilde S_I)), \; \; 
(M,F)=((M_I,F),u_{IJ})\mapsto \\
f^{*mod[-],\Gamma}(M,F):=(\Gamma_{X_I}E(p_{\tilde S_I}^{*mod[-]}(M_I,F))),\tilde f_J^{*mod[-]}u_{IJ})
\end{eqnarray*}
with, denoting for short $d_{IJ}:=d_{\tilde S_J}-d_{\tilde S_I}$
\begin{eqnarray*}
\tilde f_J^{*mod[-]}u_{IJ}:\Gamma_{X_I}E(p_{\tilde S_I}^{*mod[-]}(M_I,F))
\xrightarrow{\Gamma_{X_I}E(p_{\tilde S_I}^{*mod[-]})(u_{IJ})}\Gamma_{X_I}E(p_{\tilde S_I}^{*mod[-]}p_{IJ*}(M_J,F)[d_{IJ}]) \\
\xrightarrow{\Gamma_{X_I}E(T(p_{IJ}^{*mod},p_{\tilde S_I})(-)^{-1})[d_Y+d_{IJ}]}
\Gamma_{X_I}E(p'_{IJ*}p_{\tilde S_J}^{*mod}(M_J,F)[d_Y+d_{IJ}]) \\
\xrightarrow{=}p'_{IJ*}\Gamma_{X_J}E(p_{\tilde S_J}^{*mod[-]})(M_J,F)[d_{IJ}].
\end{eqnarray*}
It induces in the derived categories, the functor
\begin{eqnarray*}
Rf^{*mod[-],\Gamma}:D_{Ofil,\mathcal D,\infty}(S/(\tilde S_I)\to D_{Ofil,\mathcal D,\infty}(X/(Y\times\tilde S_I)), \\  
(M,F)=((M_I,F),u_{IJ})\mapsto \\
Rf^{*mod[-],\Gamma}:=f^{*mod[-],\Gamma}(M,F):=(\Gamma_{X_I}E(p_{\tilde S_I}^{*mod[-]}(M_I,F)),\tilde f_J^{*mod[-]}u_{IJ}).
\end{eqnarray*}

By definition, for $f:T\to S$ with $T,S\in\QPVar(\mathbb C)$ or with $T,S\in\AnSp(\mathbb C)^{QP}$, 
after considering a factorization $f:T\xrightarrow{l}Y\times S\xrightarrow{p_S}S$ with $Y\in\SmVar(\mathbb C)$,
$l$ a closed embedding and $p_S$ the projection, the commutative diagrams of functors
\begin{equation*}
\xymatrix{C_{\mathcal Dfil}(S/(\tilde S_I))\ar[r]^{o_S}\ar[d]^{f^{*mod[-],\Gamma}} & 
C_{Ofil,\mathcal D}(S/(\tilde S_I))\ar[d]^{f^{*mod[-],\Gamma}} \\
C_{\mathcal Dfil}(T/(Y\times\tilde S_I))\ar[r]^{o_T} & C_{Ofil,\mathcal D}(T/(Y\times\tilde S_I))} \; , \; 
\xymatrix{D_{\mathcal Dfil,r}(S/(\tilde S_I))\ar[r]^{o_S}\ar[d]^{Rf^{*mod[-],\Gamma}} & 
D_{Ofil,\mathcal D,\infty}(S/(\tilde S_I))\ar[d]^{Rf^{*mod[-],\Gamma}} \\
D_{\mathcal Dfil,\infty}(T/(Y\times\tilde S_I))\ar[r]^{o_T} & D_{Ofil,\mathcal D,r}(T/(Y\times\tilde S_I))}
\end{equation*}
where $o_S$ and $o_T$ are the forgetful functors.

Let $f:X\to S$ be a morphism, with $X,S\in\Var(\mathbb C)$, 
such that there exist a factorization $f;X\xrightarrow{l}Y\times S\xrightarrow{p_S}S$
with $Y\in\SmVar(\mathbb C)$, $l$ a closed embedding and $p_S$ the projection,
and consider $S=\cup_{i=1}^lS_i$ an open cover such that there exist closed embeddings 
$i_i:S_i\hookrightarrow\tilde S_i$, with $\tilde S_i\in\SmVar(\mathbb C)$ ;
Then, $X=\cup^l_{i=1} X_i$ with $X_i:=f^{-1}(S_i)$.
We have, for $(M,F)=((M_I,F),u_{IJ})\in C_{Ofil,\mathcal D}(S/(\tilde S_I))^{\vee}$, 
the canonical transformation map in $D_{Ofil,\mathcal D}(T^{an}/(\tilde T^{an}_I))^{\vee}$
\begin{eqnarray*}
T^{mod}(an,\gamma_T)(M,F): \\
f^{*mod[-],\Gamma}(M,F))^{an}:=((\Gamma_{T_I}E(p_{\tilde S_I}^{*mod[-]}(M_I,F)))^{an},(f^{*mod}[-]u_{IJ})^{an}) \\
\xrightarrow{(T^{mod}(an,\gamma_{T_I})(-))}
(\Gamma_{T^{an}_I}E((p_{\tilde S_I}^{*mod[-]}(M_I,F))^{an}),f^{*mod[-]}u^{an}_{IJ}) \\
\xrightarrow{=}
(\Gamma_{T^{an}_I}E(p_{\tilde S_I}^{*mod[-]}(M_I^{an},F)),f^{*mod[-]}u^{an}_{IJ})=:f^{*mod[-],\Gamma}((M,F)^{an})
\end{eqnarray*}
where the equality is obvious.

\section{The category of mixed Hodge modules on complex algebraic varieties and complex analytic spaces and
the functorialities}

For $S\in\Top$ a topological space endowed with a stratification $S=\sqcup_{k=1}^dS_k$ by locally closed subsets $S_k$
together with the perversity $p(S_k)$,
we denote by $\mathcal P(S,W)\subset D_{fil}(S)$ the category of filtered perverse sheaves of abelian groups.
For a locally compact (hence Hausdorf) topological space, we denote by $D_c(S)\subset D(S)$
the full subcategory of complexes of presheaves whose cohomology sheaves are constructible.

\subsection{The De Rahm functor for D modules on a complex analytic space}

Let $S\in\AnSm(\mathbb C)$. Recall we have the dual functor
\begin{eqnarray*}
\mathbb D_S:C(S)\to C(S), \; K\mapsto\mathbb D_S(K):=\mathcal Hom(K,E(\mathbb Z_S))
\end{eqnarray*}
which induces the functor
\begin{eqnarray*}
L\mathbb D_S:D(S)\to D(S), \; K\mapsto L\mathbb D_S(K):=\mathbb D_S(LK):=\mathcal Hom(LK,E_{et}(\mathbb Z_S)).
\end{eqnarray*}

Let $S\in\AnSm(\mathbb C)$.
\begin{itemize}
\item The functor 
\begin{equation*}
M\in\PSh_{\mathcal D}(S)\mapsto DR(S)(M):=\Omega^{\bullet}_S\otimes_{O_S}M\in C_{\mathbb C_S}(S) 
\end{equation*}
which sends a $D_S$ module to its De Rham complex (see section 4) induces, after
shifting by $d_S$ in order to send holonomic module (degree zero) to perverse sheaves, 
in the derived category the functor
\begin{eqnarray*}
DR(S)^{[-]}:D_{\mathcal D}(S)\to D_{\mathbb C_S}(S), M\mapsto \\ 
DR(S)^{[-]}(M):=DR(S)(M)[d_S]:=\Omega^{\bullet}_S\otimes_{O_S}M[d_S]
\simeq K_S\otimes_{D_S}^LM\simeq\mathcal Hom_{D_S}(\mathbb D_SL_DM,E(O_S))[d_S]
\end{eqnarray*} 
and, by functoriality, the functor
\begin{eqnarray*}
DR(S)^{[-]}:D_{\mathcal D0fil,\infty}(S)\to D_{\mathbb C_Sfil,\infty}(S), \\ 
(M,W)\mapsto DR(S)^{[-]}(M,W):=(\Omega^{\bullet}_S,F_b)\otimes_{O_S}(M,W)[d_S]=K_S\otimes_{D_S}^L(M,W)
\end{eqnarray*} 
\item On the other hand, we have the functor
\begin{equation*}
C_{\mathbb C_S}(S)\to C_{\mathcal D^{\infty}}(S), \; 
K\mapsto\mathcal Hom_{\mathbb C_S}(L_{\mathbb C}\mathbb D_S(LK),E(O_S))[-d_S]
\end{equation*}
together with, for $K\in C_{\mathbb C_S}(S)$, the canonical map 
\begin{eqnarray*}
s(K):K\to DR(S)^{[-]}(J_S^{-1}\mathcal Hom_{\mathbb C_S}(L_{\mathbb C}(\mathbb D_SLK),E(O_S))[-d_S]) \\
\xrightarrow{=}
\mathcal Hom_{D_S}(\mathbb D_S^KL_DJ_S^{-1}\mathcal Hom_{\mathbb C_S}(L_{\mathbb C}(\mathbb D_SLK),E(O_S)),E(O_S)), \\ 
c\in\Gamma(S^o,L(K))\longmapsto s(K)(c)=(\phi\in\Gamma(S^{oo},L_D\mathcal Hom(L_{\mathbb C}(K),E(O_S)))\mapsto\phi(c))
\end{eqnarray*}
where $S^{oo}\subset S^o\subset S$ are open subsets.
\end{itemize}

The main result is Riemann-Hilbert equivalence :
\begin{thm}\label{DRK}
Let $S\in\AnSm(\mathbb C)$.
\begin{itemize}
\item[(i)] The functor $J_S:D_{\mathcal D,rh}(S)\to D_{\mathcal D^{\infty},h}(S)$ is an equivalence of category.
Moreover, for $K\in C(S)$, we have $\mathcal Hom(L(K),E(O_S))\in C_{\mathcal D^{\infty},h}(S)$. 
\item[(ii)] The restriction of the De Rahm functor to the full subcategory $D_{\mathcal D,rh}(S)\subset D_{\mathcal D}(S)$
is an equivalence of category 
\begin{equation*}
DR(S)^{[-]}:D_{\mathcal D,rh}(S)\xrightarrow{\sim} D_{\mathbb C_S,c}(S)
\end{equation*}
whose inverse is the functor 
\begin{equation*}
K\in C_{\mathbb C_S,c}(S)\mapsto J^{-1}\mathcal Hom_{\mathbb C_S}(\mathbb D_SL(K),E(O_S))[-d_S],
\end{equation*}
the map 
$s(K):K\xrightarrow{\sim}DR(S)^{[-]}(J^{-1}\mathcal Hom_{\mathbb C_S}(L_{\mathbb C}\mathbb D_SL(K),E(O_S)))$ 
being an isomorphism.
\item[(iii)] The De Rahm functor $DR(S)^{[-]}$ sends regular holonomic modules to perverse sheaves.
\end{itemize}
\end{thm}

\begin{proof}
See \cite{Kashiwara}.
\end{proof}
 
Let $S\in\AnSp(\mathbb C)$ and $S=\cup_{i=1}^lS_i$ an open cover 
such that there exists closed embeddings $i_i:S_i\hookrightarrow\tilde S_i$.
\begin{itemize}
\item The De Rham functor is in this case 
\begin{eqnarray*}
DR(S)^{[-]}:D_{\mathcal D0fil,\infty}(S)\to D_{\mathbb C_Sfil,\infty}(S), M=((M_I,W),u_{IJ})\mapsto \\
DR(S)^{[-]}(M,W):=(DR(\tilde S_I)^{[-]}(M_I,W),DR^{[-]}(u_{IJ})):=
(\Omega^{\bullet}_{\tilde S_I}\otimes_{O_{\tilde S_I}}(M_I,W),DR^{[-]}(u_{IJ}))
\end{eqnarray*} 
with, denoting for short $d_I=d_{\tilde S_I}$ 
\begin{eqnarray*}
DR^{[-]}(u_{IJ}):\Omega^{\bullet}_{\tilde S_I}\otimes_{O_{\tilde S_I}}(M_I,W)[d_I]
\xrightarrow{\ad(p_{IJ},p_{IJ*})(-))}
p_{IJ*}p_{IJ}^*\Omega^{\bullet}_{\tilde S_I}\otimes_{O_{\tilde S_I}}(M_I,W)[d_I] \\
\xrightarrow{p_{IJ*}\Omega_{\tilde S_J/\tilde S_I}[d_I]}
p_{IJ*}\Omega^{\bullet}_{\tilde S_J}\otimes_{O_{\tilde S_J}}p_{IJ}^{*mod}(M_I,W)[d_I] \\
\xrightarrow{p_{IJ*}I(p_{IJ}^{*mod},p_{IJ})(-,-)(u_{IJ})[d_I]}
p_{IJ*}\Omega^{\bullet}_{\tilde S_J}\otimes_{O_{\tilde S_J}}(M_J,W)[d_J+d_{IJ}]
\end{eqnarray*}
\item Considering the diagrams
\begin{equation*}
D_{IJ}=\xymatrix{\tilde S_J\ar[r]^{p_{IJ}} &  \tilde S_I \\
S_J\ar[u]^{i_J}\ar[r]^{j_{IJ}} &  S_I\ar[u]^{i_I}}
\end{equation*}
we get the functor
\begin{eqnarray*}
C_{\mathbb C_Sfil}(S)\xrightarrow{T(S/(S_I))}C_{\mathbb C_Sfil}(S/(S_I))\to C_{\mathcal D0fil}(S/(S_I)),\\ 
(K,W)\mapsto 
(\mathcal Hom_{\mathbb C_{\tilde S_I}}(L_{\mathbb C}\mathbb D_{\tilde S_I}(Li_{I*}j_I^*(K,W)),
E(O_{\tilde S_I}))[-d_{\tilde S_I}],u_{IJ}(K,W))
\end{eqnarray*}
where 
\begin{eqnarray*}
u_{IJ}(K,W):\mathcal Hom_{\mathbb C_{\tilde S_I}}(L_{\mathbb C}\mathbb D_{\tilde S_I}L(i_{I*}j_I^*(K,W)),
E(O_{\tilde S_I}))[-d_{\tilde S_I}] \\
\xrightarrow{\ad(p_{IJ}^{*mod},p_{IJ})(-)}
p_{IJ*}p_{IJ}^{*mod[-]}\mathcal Hom_{\mathbb C_{\tilde S_I}}(L_{\mathbb C}\mathbb D_{\tilde S_I}L(i_{I*}j_I^*(K,W)),
E(O_{\tilde S_I}))[-d_{\tilde S_J}] \\
\xrightarrow{\mathcal Hom(-,Eo(p_{IJ}))\circ T(p_{IJ},\hom)(-,-)}
p_{IJ*}\mathcal Hom_{\mathbb C_{\tilde S_I}}(p_{IJ}^*L_{\mathbb C}\mathbb D_{\tilde S_I}L(i_{I*}j_I^*(K,W)),
E(O_{\tilde S_J}))[-d_{\tilde S_J}] \\
\xrightarrow{\mathcal Hom(T(p_{IJ},\mathbb D)(-)^{-1},-)) }
p_{IJ*}\mathcal Hom_{\mathbb C_{\tilde S_I}}(L_{\mathbb C}\mathbb D_{\tilde S_J}p_{IJ}^*L(i_{I*}j_I^*(K,W)),
E(O_{\tilde S_J}))[-d_{\tilde S_J}] \\
\xrightarrow{\mathcal Hom(L_{\mathbb C}\mathbb D_{\tilde S_J}T^q(D_{IJ})(j_I^*(K,W)),E(O_{\tilde S_J}))}
p_{IJ*}\mathcal Hom_{\mathbb C_{\tilde S_J}}(L_{\mathbb C}\mathbb D_{\tilde S_J}L(i_{J*}j_J^*(K,W)),
E(O_{\tilde S_J}))[-d_{\tilde S_J}]. 
\end{eqnarray*}
Moreover, for $(K,W)\in C_{fil}(S)$, we have 
\begin{equation*}
(\mathcal Hom_{\mathbb C_{\tilde S_I}}(L_{\mathbb C}\mathbb D_{\tilde S_I}L(i_{I*}j_I^*(K,W)),
E(O_{\tilde S_I}))[-d_{\tilde S_I}],u_{IJ}(K,W))\in C_{\mathcal D^{\infty}0fil,h}(S)^0
\end{equation*}
and a canonical map in $D_{fil}(S)=D_{fil}(S/(\tilde S_I))$
\begin{eqnarray*}
s(K):T(S/(S_I))(K,W):=(L(i_{I*}j_I^*(K,W)),I)\to \\
DR(S)^{[-]}(J_S^{-1}\mathcal Hom_{\mathbb C_{\tilde S_I}}(L_{\mathbb C}\mathbb D_{\tilde S_I}L(i_{I*}j_I^*(K,W)),
E(O_{\tilde S_I}))[-d_{\tilde S_I}],u_{IJ}(K,W))
\end{eqnarray*}
\end{itemize}

\begin{cor}\label{DRKsing}
Let $S\in\AnSp(\mathbb C)$. 
Let $S=\cup_{i=1}^lS_i$ an open cover such that there exists closed embeddings $i_i:S_i\hookrightarrow\tilde S_i$.
The restriction of the De Rahm functor to the full subcategory $D^0_{\mathcal D,rh}(S)\subset D^0_{\mathcal D}(S)$
is an equivalence of category 
\begin{equation*}
DR(S)^{[-]}:D^0_{\mathcal D,rh}(S)\xrightarrow{\sim} D_{\mathbb C_S,c}(S)
\end{equation*}
whose inverse is the functor 
\begin{equation*}
K\mapsto J_S^{-1}(\mathcal Hom_{\mathbb C_{\tilde S_I}}(L_{\mathbb C}\mathbb D_{\tilde S_I}L(i_{I*}j_I^*K),
E(O_{\tilde S_I}))[-d_{\tilde S_I}],u_{IJ}(K))
\end{equation*}
the map 
\begin{eqnarray*}
s(K):T(S/(S_I))(K,W):=(L(i_{I*}j_I^*(K,W)),I)\to \\
DR(S)^{[-]}(J_S^{-1}\mathcal Hom_{\mathbb C_{\tilde S_I}}(L_{\mathbb C}\mathbb D_{\tilde S_I}L(i_{I*}j_I^*(K,W)),
E(O_{\tilde S_I}))[-d_{\tilde S_I}],u_{IJ}(K,W))
\end{eqnarray*} 
being an isomorphism.
\end{cor}

\begin{proof}
Follows from theorem \ref{DRK}(ii), see \cite{Saito}.
\end{proof}

\begin{prop}\label{DRdual}
\begin{itemize}
\item[(i)] Let $S\in\AnSm(\mathbb C)$.Then, for $M\in C_{\mathcal D,c}(S)$, there is a canonical isomorphism
\begin{eqnarray*}
T(D,DR)(M):\mathbb D^{\mathbb C}_S DR(S)^{[-]}(M)\xrightarrow{\sim}DR(S)^{[-]}(\mathbb D_S^KL_DM)
\end{eqnarray*}
\item[(ii)] Let $S\in\AnSp(\mathbb C)$. 
Let $S=\cup_{i=1}^lS_i$ an open cover such that there exists closed embeddings $i_i:S_i\hookrightarrow\tilde S_i$.
Then, for $M=(M_I,u_{IJ})\in C^0_{\mathcal D,c}(S/(\tilde S_I))$, there is a canonical isomorphism
\begin{eqnarray*}
T(D,DR)(M):\mathbb D^{\mathbb C}_S DR(S)^{[-]}(M)\xrightarrow{\sim}DR(S)^{[-]}(L\mathbb D_S^KM)
\end{eqnarray*}
\end{itemize}
\end{prop}

\begin{proof}
\noindent(i):See \cite{LvDmod}.

\noindent(ii):Follows from (i), see \cite{Saito}.
\end{proof}

We have the following transformation maps :
\begin{itemize}
\item Let $g:T\to S$ a morphism with $T,S\in\AnSm(\mathbb C)$. 
We have, for $(M,W)\in C_{\mathcal D0fil}(S)$, the canonical transformation map in $D_{fil}(T)$ :
\begin{eqnarray*}
T(g,DR)(M,W):g^*DR(S)^{[-]}(M,W):=g^*(\Omega_S^{\bullet}\otimes_{O_S}L_D(M,W))[d_S] \\
\xrightarrow{\Omega_{T/S}(L_D(M,W))[d_S]}\Omega_T^{\bullet}\otimes_{O_T}g^{*mod}L_D(M,W)[d_S] \\
\xrightarrow{=}\Omega_T^{\bullet}\otimes_{O_T}g^{*mod[-]}L_D(M,W)[d_S]=:DR(T)^{(-)}(Lg^{*mod[-]}(M,W))
\end{eqnarray*}
Note that this transformation map is NOT an isomorphism in general. It is an isomorphism if $g$ is a smooth morphism.
If $g$ is a closed embedding, it is an isomorphism for $M$ non caracteristic with respect to $g$.

\item Let $j:S^o\hookrightarrow S$ an open embedding with $S\in\AnSm(\mathbb C)$.
We have, for $(M,W)\in C_{\mathcal D0fil}(S^o)$, the canonical transformation map in $D_{fil}(S)$ :
\begin{eqnarray*}
T_*(j,DR)(M,W):DR(S)^{[-]}(j_*(M,W)):=\Omega_S^{\bullet}\otimes_{O_S}j_*(M,W)[d_S] \\
\xrightarrow{T^O_w(j,\otimes)(L_D(M,W))[d_S]}
j_*(\Omega_{S^o}^{\bullet}\otimes_{O_{S^o}}L_D(M,W))[d_S]=:j_*DR(S)^{[-]}(M,W) 
\end{eqnarray*}
which is an isomorphism (see proposition \ref{mw0prop}).

\item  Let $g:T\to S$ a morphism with $T,S\in\AnSp(\mathbb C)$.
Assume there exist a factorization $g:T\xrightarrow{l}Y\times S\xrightarrow{p_S}S$ with $Y\in\AnSm(\mathbb C)$,
$l$ a closed embedding and $p_S$ the projection.
Let $S=\cup_iS_i$ an open covers such that there exist closed embeddings 
$i_i:S_i\hookrightarrow\tilde S_i$ with $\tilde S_i\in\AnSm(\mathbb C)$.
We have, for $M=(M_I,u_{IJ})\in C_{\mathcal D}(S/(\tilde S_I))$, 
the canonical transformation map in $D_{fil}(T/(Y\times\tilde S_I))$
\begin{eqnarray*}
T^!(g,DR)(M):T(T/(Y\times\tilde S_I))(g^!DR(S)^{[-]}(M,W)) \\
\xrightarrow{:=}(\Gamma_{T_I}E(\tilde g_I^*(\Omega_{\tilde S_I}^{\bullet}\otimes_{O_{\tilde S_I}}L_D(M_I,W))),
\tilde g_I^*DR(u_{IJ})) \\
\xrightarrow{(\Omega_{(Y\times\tilde S_I/\tilde S_I)}(L_D(M_I,W)))} 
(\Gamma_{T_I}E(\Omega_{Y\times\tilde S_I}^{\bullet}\otimes_{O_{Y\times\tilde S_I}}\tilde g_I^{*mod}(M_I,W)),
DR(\tilde g_I^{*mod}u_{IJ})) \\
\xrightarrow{(T^O_w(\gamma,\otimes)(\tilde g_I^{*mod}L_D(M,W)))}
(\Omega_{Y\times\tilde S_I}^{\bullet}\otimes_{O_{Y\times\tilde S_I}}\Gamma_{T_I}E(\tilde g_I^{*mod}(M_I,W)),
DR(\tilde g_I^{*mod}u_{IJ})) \\
\xrightarrow{=:}DR(T)^{[-]}(Rg^{*mod[-],\Gamma}(M,W))
\end{eqnarray*}
which is an isomorphism.

\item Let $f:T\to S$ a morphism with $T,S\in\Var(\mathbb C)$.
Assume there exist a factorization $f:T\xrightarrow{l}Y\times S\xrightarrow{p_S}S$ with $Y\in\SmVar(\mathbb C)$,
$l$ a closed embedding and $p_S$ the projection.
Let $S=\cup_iS_i$ an open covers such that there exist closed embeddings 
$i_i:S_i\hookrightarrow\tilde S_i$ with $\tilde S_i\in\SmVar(\mathbb C)$.
We have, for $M=(M_I,u_{IJ})\in C_{\mathcal D}(S/(\tilde S_I))$, the canonical map in $D_{fil}(T/(Y\times\tilde S_I))$ 
\begin{eqnarray*}
DR(T)^{[-]}(T^{mod}(an,\gamma_T)(M)): \\
DR(T)^{[-]}((Rf^{*mod[-],\Gamma}M)^{an}):=
DR(T)^{[-]}(((\Gamma_{T_I}E(\tilde f_I^{*mod[-]}(M_I)))^{an},(f^{*mod[-]}u^q_{IJ})^{an})) \\
\to DR(T)^{[-]}(Rf^{*mod[-],\Gamma}(M^{an})):=
DR(T)^{[-]}((\Gamma_{T^{an}_I}E(\tilde f_I^{*mod[-]}(M_I^{an})),f^{*mod[-]}u^{q,an}_{IJ})).
\end{eqnarray*}
\end{itemize}

\begin{prop}\label{DRTangamma}
Let $f:T\to S$ a morphism with $T,S\in\Var(\mathbb C)$.
Assume there exist a factorization $f:T\xrightarrow{l}Y\times S\xrightarrow{p_S}S$ with $Y\in\SmVar(\mathbb C)$,
$l$ a closed embedding and $p_S$ the projection.
Let $S=\cup_iS_i$ an open covers such that there exist closed embeddings 
$i_i:S_i\hookrightarrow\tilde S_i$ with $\tilde S_i\in\SmVar(\mathbb C)$.
Then, for $M=(M_I,u_{IJ})\in C_{\mathcal D,rh}(S/(\tilde S_I))$, the map in $D_{fil}(T/(Y\times\tilde S_I))$ 
\begin{eqnarray*}
DR(T)^{[-]}(T^{mod}(an,\gamma_T)(M)): \\
DR(T)^{[-]}((Rf^{*mod[-],\Gamma}M)^{an}):=
DR(T)^{[-]}(((\Gamma_{T_I}E(\tilde f_I^{*mod[-]}(M_I)))^{an},(f^{*mod[-]}u^q_{IJ})^{an})) \\
\to DR(T)^{[-]}(Rf^{*mod[-],\Gamma}(M^{an})):=
DR(T)^{[-]}((\Gamma_{T^{an}_I}E(\tilde f_I^{*mod[-]}(M_I^{an})),f^{*mod[-]}u^{q,an}_{IJ}))
\end{eqnarray*}
given above is an isomorphism.
\end{prop}

\begin{proof}
See \cite{LvDmod}.
\end{proof}

In the algebraic case, we have, by proposition \ref{DRTangamma},
for complexes of $D$-modules whose cohomology sheaves are regular holonomic the following canonical isomorphisms:

\begin{defi}\label{DRf}
\begin{itemize}
\item[(i)] Let $f:T\to S$ a morphism with $T,S\in\Var(\mathbb C)$.
Assume there exist a factorization $f:T\xrightarrow{l}Y\times S\xrightarrow{p_S}S$ with $Y\in\SmVar(\mathbb C)$,
$l$ a closed embedding and $p_S$ the projection.
We have, for $M=(M_I,u_{IJ})\in C_{\mathcal D,rh}(S/(\tilde S_I))^0$, the canonical map
\begin{eqnarray*}
T^!(f,DR)(M):f^!DR(S)^{[-]}(M^{an})
\xrightarrow{T^!(f,DR)(M^{an})}DR(T)^{[-]}(Rf^{*mod[-],\Gamma}(M^{an})) \\
\xrightarrow{DR(T)^{[-]}(T^{mod}(an,\gamma_T)(M))}
DR(T)^{[-]}((Rf^{*mod[-],\Gamma}M)^{an})=:DR(T)^{[-]}((Rf^{*mod[-],\Gamma}M)^{an}).
\end{eqnarray*}
which is an isomorphism by proposition \ref{DRTangamma}.

\item[(ii)] Let $f:T\to S$ a morphism with $T,S\in\Var(\mathbb C)$.
Assume there exist a factorization $f:T\xrightarrow{l}Y\times S\xrightarrow{p_S}S$ with $Y\in\SmVar(\mathbb C)$,
$l$ a closed embedding and $p_S$ the projection.
We have, for $M=(M_I,u_{IJ})\in C_{\mathcal D,rh}(S/(\tilde S_I))^0$, the canonical transformation map
\begin{eqnarray*}
T(f,DR)(M):DR(T)^{[-]}((Lf^{\hat*mod[-],\Gamma}M)^{an}):=DR(T)^{[-]}((L\mathbb D^K_TRf^{*mod[-],\Gamma}L\mathbb D^K_SM)^{an}) \\
\xrightarrow{T(D,DR)(-)}L\mathbb D^K_TDR(T)^{[-]}((Rf^{*mod[-],\Gamma}L\mathbb D^K_SM)^{an}))
\xrightarrow{L\mathbb D^K_TT^!(f,DR)(-)}L\mathbb D^K_Tf^!DR(S)^{[-]}(L\mathbb D^K_SM^{an}) \\
\xrightarrow{L\mathbb D^K_Tf^!T(D,DR)(-)^{-1}}L\mathbb D^K_Tf^!L\mathbb D^K_SDR(S)^{[-]}(M^{an})=f^*DR(S)^{[-]}(M^{an})
\end{eqnarray*}
which is an isomorphism by (i) and proposition \ref{DRdual}.

\item[(iii)] Let $f:T\to S$ a morphism with $T,S\in\QPVar(\mathbb C)$.
Consider a factorization $f:T\xrightarrow{l}Y\times S\xrightarrow{p_S}S$ with $Y=\mathbb P^{N,o}\subset\mathbb P^N$
an open subset, $l$ a closed embedding, and $p_S$ the projection.
We have, for $M\in C_{\mathcal D,rh}(T/Y\times\tilde S)^0$, the canonical transformation map
\begin{eqnarray*}
T_*(f,DR)(M):DR(S)^{[-]}((\int_fM)^{an})\xrightarrow{\ad(f^*,Rf_*)(-)}
Rf_*f^*DR(S)^{[-]}((\int_fM)^{an})\xrightarrow{Rf_*T(f,DR)((\int_fM))} \\
Rf_*DR(T)^{[-]}((Lf^{\hat*mod[-],\Gamma}\int_fM)^{an})\xrightarrow{Rf_*DR(T)^{[-]}((\ad(Lf^{\hat*mod[-],\Gamma},\int_f)(M))^{an})}
Rf_*DR(T)^{[-]}(M^{an})
\end{eqnarray*}
which is an isomorphism by GAGA in the proper case and by the open embedding case (c.f. proposition \ref{DRTangamma}).

\item[(iv)] Let $f:T\to S$ a morphism with $T,S\in\QPVar(\mathbb C)$.
Consider a factorization $f:T\xrightarrow{l}Y\times S\xrightarrow{p_S}S$ with $Y=\mathbb P^{N,o}\subset\mathbb P^N$
an open subset, $l$ a closed embedding, and $p_S$ the projection.
We have, for $M\in C_{\mathcal D,rh}(T)$, the canonical transformation map
\begin{eqnarray*}
T_!(f,DR)(M):Rf_!DR(T)^{[-]}(M^{an}) \\
\xrightarrow{Rf_!DR(T)^{[-]}(\ad(\int_{f!}Rf^{*mod[-],\Gamma})(M)^{an})}
Rf_!DR(T)^{[-]}((Rf^{*mod[-],\Gamma}\int_{f!}(M))^{an}) \\
\xrightarrow{T^!(f,DR)(\int_{f!}M)}Rf_!f^!DR(S)^{[-]}(\int_{f!}(M))  
\xrightarrow{\ad(Rf_!,f^!)(-)}DR(S)^{[-]}(\int_{f!}M)
\end{eqnarray*}
which is an isomorphism by (iii) and proposition \ref{DRdual}.
\end{itemize}
\end{defi}

In the filtered case, we will consider the weight monodromy filtration for open embeddings :

\begin{defi}\label{jw}
Let $k\subset\mathbb C$ a subfield.
\begin{itemize}
\item[(i)] Let $S\in\Var(\mathbb C)$ and $j:S^o\hookrightarrow S$ an open embedding such that 
$D:=S\backslash S^o=V(s)\subset S$ is a Cartier divisor.
\begin{itemize}
\item For $(K,W)\in P_{fil}(S^{o,an})$, we consider as in \cite{Saito}
\begin{equation*}
j_{*w}(K,W):=(Rj_*K,W)\in P_{fil}(S^{an}), \; W_kRj_*K:=<Rj_*W_kK,W(N)_kK>\subset Rj_*K 
\end{equation*}
so that $j^*j_{*w}(K,W)=(K,W)$,
where $W_kRj_*K\subset Rj_*K$ is given by $W$ and the weight monodromy filtration $W(N)$ of 
the universal cover $\pi:\tilde S^{o,an}\to S^{o,an}$.
Note that a stratitification of $W_kRj_*K$ is given by the closure of a stratification of $W_kK$ 
and $D:=S\backslash S^o$. 
\item For $(K,W)\in P_{fil}(S^{o,an})$, we consider
\begin{equation*}
j_{!w}(K,W):=\mathbb D_S^vj_{*w}\mathbb D_S^v(K,W)\in P_{fil}(S^{an}) 
\end{equation*}
so that $j^*j_{!w}(K,W)=(K,W)$.
\end{itemize}
For $(K',W)\in P_{fil}(S^{an})$, there is, by construction,
\begin{itemize}
\item a canonical map $\ad(j^*,j_{*w})(K',W)=\ad(j^*,j_*)(K'):(K',W)\to j_{*w}j^*(K',W)$ 
in $P_{fil}(S^{an})$,
\item a canonical map $\ad(j_{!w},j^*)(K',W)=\ad(j_!,j^*)(K'):j_{!w}j^*(K',W)\to(K',W)$ 
in $P_{fil}(S^{an})$.
\end{itemize}
\item[(ii)] Let $S\in\Var(\mathbb C)$. Let $j:S^o:=S\backslash Z\hookrightarrow S$ an open embedding with
$Z=V(\mathcal I)\subset S$ an arbitrary closed subset, $\mathcal I\subset O_S$ being an ideal subsheaf. 
Taking generators $\mathcal I=(s_1,\ldots,s_r)$, we get $Z=V(s_1,\ldots,s_r)=\cap^r_{i=1}Z_i\subset S$ with 
$Z_i=V(s_i)\subset S$, $s_i\in\Gamma(S,\mathcal L_i)$ and $L_i$ a line bundle. 
Note that $Z$ is an arbitrary closed subset, $d_Z\geq d_X-r$ needing not be a complete intersection. Denote by 
$j_I:S^{o,I}:=\cap_{i\in I}(S\backslash Z_i)=S\backslash(\cup_{i\in I}Z_i)\xrightarrow{j_I^o}S^o\xrightarrow{j} S$ 
the open complementary embeddings, where $I\subset\left\{1,\cdots,r\right\}$. Denote 
\begin{equation*}
\mathcal D(Z/S):=\left\{(Z_i)_{i\in[1,\ldots r]},Z_i\subset S,\cap Z_i=Z\right\},Z'_i\subset Z_i
\end{equation*}
the flag category. For $(K,W)\in C(P_{fil}(S^{o,an}))$, we define by (i) 
\begin{itemize}
\item the (bi)-filtered complex of $D_S$-modules
\begin{equation*}
j_{*w}(K,W):=\varinjlim_{\mathcal D(Z/S)}\Tot_{card I=\bullet}(j_{I*w}j_I^{o*}(K,W))\in C(P_{fil}(S^{an})) 
\end{equation*}
where the horizontal differential are given by, 
if $I\subset J$, $d_{IJ}:=\ad(j^*_{IJ},j_{IJ*w})(j_I^{o*}(K,W))$, 
$j_{IJ}:S^{oJ}\hookrightarrow S^{oI}$ being the open embedding, 
and $d_{IJ}=0$ if $I\notin J$,
\item the (bi)-filtered complex of $D_S$-modules
\begin{eqnarray*}
j_{!w}(K,W):=\varprojlim_{\mathcal D(Z/S)}\Tot_{card I=-\bullet}(j_{I!w}j_I^{o*}(K,W)) 
=\mathbb D_S^vj_{*w}\mathbb D_S^v(K,W)\in C(P_{fil}(S^{an})),
\end{eqnarray*}
where the horizontal differential are given by, 
if $I\subset J$, $d_{IJ}:=\ad(j_{IJ!w},j^*_{IJ})(j_I^{o*}(K,W))$, 
$j_{IJ}:S^{oJ}\hookrightarrow S^{oI}$ being the open embedding, and $d_{IJ}=0$ if $I\notin J$.
\end{itemize}
By definition, we have for $(K,W)\in C(P_{fil}(S^{o,an}))$, 
$j^*j_{*w}(K,W)=(K,W)$ and $j^*j_{!w}(K,W)=(K,W)$.
For $(K',W)\in C(P_{fil}(S^{an}))$, there is, by (i),
\begin{itemize}
\item a canonical map $\ad(j^*,j_{*w})(K',W):(K',W)\to j_{*w}j^*(K',W)$ in $C(P_{fil}(S^{an}))$,
\item a canonical map $\ad(j_{!w},j^*)(K',W):j_{!w}j^*(K',W)\to(K',W)$ in $C(P_{fil}(S^{an}))$.
\end{itemize}
\end{itemize}
\end{defi}

\begin{defi}\label{gammaw}
Let $S\in\Var(\mathbb C)$. Let $Z\subset S$ a closed subset.
Denote by $j:S\backslash Z\hookrightarrow S$ the complementary open embedding. 
\begin{itemize}
\item[(i)] We define using definition \ref{jw}, the filtered Hodge support section functor
\begin{eqnarray*}
\Gamma^w_Z:C(P_{fil}(S^{an}))\to C(P_{fil}(S^{an})), \\ 
(K,W)\mapsto\Gamma^w_Z(K,W):=\Cone(\ad(j^*,j_{*w})(K,W):(K,W)\to j_{*w}j^*(K,W))[-1],
\end{eqnarray*}
together we the canonical map $\gamma^w_Z(K,W):\Gamma^w_Z(K,W)\to (K,W)$.
\item[(i)'] Since $j_{*w}:C(P_{fil}(S^{o,an}))\to C(P_{fil}(S^{an}))$ is an exact functor, 
$\Gamma^w_Z$ induces the functor
\begin{eqnarray*}
\Gamma^w_Z:D_{fil,c}(S^{an})\to D_{fil,c}(S^{an}), \; (K,W)\mapsto\Gamma^w_Z(K,W)
\end{eqnarray*}
\item[(ii)] We define using definition \ref{jw}, the dual filtered Hodge support section functor
\begin{eqnarray*}
\Gamma^{\vee,w}_Z:C(P_{fil}(S^{an}))\to C(P_{fil}(S^{an})), \\ 
(K,W)\mapsto\Gamma^{\vee,w}_Z(K,W):=\Cone(\ad(j_{!w},j^*)(K,W):j_{!w},j^*(K,W)\to (K,W)),
\end{eqnarray*}
together we the canonical map $\gamma^{\vee,Hdg}_Z(K,W):(K,W)\to\Gamma_Z^{\vee,w}(K,W)$.
\item[(ii)'] Since $j_{!w}:C(P_{fil}(S^{o,an}))\to C(P_{fil}(S^{an}))$ is an exact functor, 
$\Gamma^{\vee,w}_Z$ induces the functor
\begin{eqnarray*}
\Gamma^{\vee,w}_Z:D_{fil,c}(S^{an})\to D_{fil,c}(S^{an}), \; (K,W)\mapsto\Gamma^{\vee,w}_Z(K,W)
\end{eqnarray*}
\end{itemize}
\end{defi}

Let $S\in\Var(\mathbb C)$ and $D=V(s)\subset S$ a Cartier divisor. 
Denote $i:D\hookrightarrow S$ the closed embedding and $j:S^o:=S\backslash D\hookrightarrow S$ the open embedding.
Let $\pi:\tilde S^{o,an}\to S^{o,an}$ the universal covering.
We then consider, for $(K,W)\in D_{fil,c}(S^{o,an})=\Ho(C(P_{fil}(S^{o,an}))$, 
\begin{itemize}
\item the filtered nearby cycle functor 
\begin{equation*}
\psi_D(K,W):=(\psi_DK,W)\in D_{fil,c}(D^{an}), \; W_k(\psi_D(K,W)):=<W_k\psi_DK,W(N)_k\psi_DK>\subset\psi_DK,
\end{equation*}
\item the vanishing cycle functor 
\begin{equation*}
\phi_D(K,W):=\Cone(\ad(j\circ\pi^*,j\circ\pi_*)(K):i^*(K,W)\to\psi_D(K,W))\in D_{fil,c}(D^{an}),
\end{equation*}
\item the canonical morphisms in $D_{fil,c}(D^{an})$
\begin{eqnarray*}
can(K,W):=c(\phi_D(K,W)):\psi_D(K,W)\to\phi_D(K,W), \\ 
var(K,W):=\mathbb D_S^vc(\phi_D\mathbb D_S^vD(K,W)):\phi_D(K,W)\to\psi_D(K,W).
\end{eqnarray*}
\end{itemize}

\begin{defi}\label{fw}
Let $k\subset\mathbb C$ a subfield.
\begin{itemize}
\item[(i)]Let $f:X\to S$ a morphism with $S,X\in\Var(\mathbb C)$. 
Consider the graph factorization $f:X\xrightarrow{l}X\times S\xrightarrow{p}S$ of $f$ 
where $l$ the the graph closed embedding and $p$ is the projection.
We have, using definition \ref{gammaw},
\begin{itemize}
\item the inverse image functor
\begin{eqnarray*}
f^{*w}:D_{fil,c}(S^{an})\to D_{fil,c}(X^{an}), \;
(K,W)\mapsto f^{*w}(K,W):=l^*\Gamma_X^{\vee,w}p^*(K,W)
\end{eqnarray*}
\item the exceptional inverse image functor
\begin{eqnarray*}
f^{*w}:D_{fil,c}(S^{an})\to D_{fil,c}(X^{an}), \;
(K,W)\mapsto f^{!w}(K,W):=l^*\Gamma_X^wp^*(K,W).
\end{eqnarray*}
\end{itemize}
\item[(ii)]Let $f:X\to S$ a morphism with $S,X\in\Var(\mathbb C)$.
Consider a compactification $f:X\hookrightarrow{j}\bar X\xrightarrow{\bar f}S$ of $f$
with $\bar X\in\Var(\mathbb C)$, $j$ an open embedding and $\bar f$ a proper morphism.
We have, using definition \ref{jw},
\begin{itemize}
\item the direct image functor
\begin{eqnarray*}
Rf_{*w}:D_{fil,c}(X^{an})\to D_{fil,c}(S^{an}), \;
(K,W)\mapsto Rf_{*w}(K,W):=R\bar f_*j_{*w}(K,W)
\end{eqnarray*}
\item the proper direct image functor
\begin{eqnarray*}
Rf_{!w}:D_{fil,c}(X^{an})\to D_{fil,c}(S^{an}), \;
(K,W)\mapsto Rf_{!w}(K,W):=R\bar f_*j_{!w}(K,W).
\end{eqnarray*}
\end{itemize}
\end{itemize}
\end{defi}

\subsection{The filtered Hodge direct image, the filtered Hodge inverse image, 
and the hodge support section functors for mixed hodge modules}

\begin{itemize}
\item Let $S\in\SmVar(\mathbb C)$.
The category $C_{\mathcal D(1,0)fil,rh}(S)\times_I D_{fil}(S^{an})$ is the category 
\begin{itemize}
\item whose set of objects is the set of triples $\left\{((M,F,W),(K,W),\alpha)\right\}$ with 
\begin{eqnarray*}
(M,F,W)\in C_{\mathcal D(1,0)fil,rh}(S), \, (K,W)\in D_{fil}(S^{an}), \; 
\alpha:(K,W)\otimes\mathbb C_{S^{an}}\to DR(S)^{[-]}((M,W)^{an})
\end{eqnarray*}
where $DR(S)^{[-]}:=DR(S)^{[-]}(S^{an}):C_{\mathcal D(1,0)fil,rh}(S^{an})\to C_{fil}(S^{an})$ 
is the De Rahm functor and $\alpha$ is an morphism in $D_{fil}(S^{an})$,
\item and whose set of morphisms are 
\begin{equation*}
\phi=(\phi_D,\phi_C,[\theta]):((M_1,F,W),(K_1,W),\alpha_1)\to((M_2,F,W),(K_2,W),\alpha_2)
\end{equation*}
where $\phi_D:(M_1,F,W)\to(M_2,F,W)$ and $\phi_C:(K_1,W)\to (K_2,W)$ are morphisms, and
\begin{eqnarray*}
\theta=(\theta^{\bullet},I(DR(S)(\phi^{an}_D))\circ I(\alpha_1),I(\alpha_2)\circ I(\phi_C\otimes I)):
I(K_1,W)\otimes\mathbb C_{S^{an}}[1]\to I(DR(S)(M^{an}_2,W)) 
\end{eqnarray*}
is an homotopy, i.e. for all $i\in\mathbb Z$, 
\begin{eqnarray*}
\theta^i\circ\partial^i-\partial^{i+1}\circ\theta^i=
(I(DR(S)(\phi_D))\circ I(\alpha_1))^i-(I(\alpha_2)\circ I(\phi_C))^i,
\end{eqnarray*}
$I:C_{fil}(S^{an})\to K_{fil}(S^{an})$ being the injective resolution functor : for
$(K,W)\in C_{fil}(S^{an})$, we take an injective resolution $k:(K,W)\to I(K,W)$ 
with $I(K,W)\in C_{fil}(S^{an})$ which is unique modulo homotopy, 
and the class $[\theta]$ of $\theta$ does NOT depend of the injective resolution ;
in particular, we have 
\begin{equation*}
DR(S)^{[-]}(\phi^{an}_D)\circ\alpha_1=\alpha_2\circ(\phi_C\otimes I) 
\end{equation*}
in $D_{fil}(S^{an})$ ; and for
\begin{itemize}
\item $\phi=(\phi_D,\phi_C,[\theta]):((M_1,F,W),(K_1,W),\alpha_1)\to((M_2,F,W),(K_2,W),\alpha_2)$
\item $\phi'=(\phi'_D,\phi'_C,[\theta']):((M_2,F,W),(K_2,W),\alpha_2)\to((M_3,F,W),(K_3,W),\alpha_3)$
\end{itemize}
the composition law is given by 
\begin{eqnarray*}
\phi'\circ\phi:=(\phi'_D\circ\phi_D,\phi'_C\circ\phi_C,
I(DR(S)(\phi^{'an}_D))\circ[\theta]+[\theta']\circ I(\phi_C\otimes I)[1]): \\
((M_1,F,W),(K_1,W),\alpha_1)\to((M_3,F,W),(K_3,W),\alpha_3),
\end{eqnarray*}
in particular for $((M,F,W),(K,W),\alpha)\in C_{\mathcal D(1,0)fil,rh}(S)\times_I D_{fil}(S^{an})$,
\begin{equation*}
I_{((M,F,W),(K,W),\alpha)}=(I_M,I_K,0).
\end{equation*}
\end{itemize}
We have then the full embedding
\begin{eqnarray*}
\PSh_{\mathcal D(1,0)fil,rh}(S)\times_I P_{fil}(S^{an})
\hookrightarrow C_{\mathcal D(1,0)fil,rh}(S)\times_I D_{fil}(S^{an})
\end{eqnarray*}
where $\PSh_{\mathcal D(1,0)fil,rh}(S)\times_I P_{fil}(S^{an})$ is the category
\begin{itemize}
\item whose set of objects is the set of triples $\left\{((M,F,W),(K,W),\alpha)\right\}$ with 
\begin{eqnarray*}
(M,F,W)\in\PSh_{\mathcal D(1,0)fil,rh}(S), \, (K,W)\in P_{fil}(S^{an}), \; 
\alpha:(K,W)\otimes\mathbb C_{S^{an}}\to DR(S)^{[-]}((M,W)^{an})
\end{eqnarray*}
where $DR(S)^{[-]}$ is the De Rahm functor and $\alpha$ is an isomorphism in $D_{fil}(S^{an})$,
\item and whose set of morphisms are 
\begin{equation*}
\phi=(\phi_D,\phi_C)=(\phi_D,\phi_C,0):((M_1,F,W),(K_1,W),\alpha_1)\to((M_2,F,W),(K_2,W),\alpha_2)
\end{equation*}
where $\phi_D:(M_1,F,W)\to(M_2,F,W)$ and $\phi_C:(K_1,W)\to (K_2,W)$ are morphisms (of filtered sheaves)
and $DR(S)^{[-]}(\phi^{an}_D)\circ\alpha_1=\alpha_2\circ(\phi_C\otimes I)$ in $P_{fil}(S^{an})$.
\end{itemize}
\item Let $S\in\Var(\mathbb C)$. Let $S=\cup_{i\in I}S_i$ an open cover such that there
exists closed embeddings $i_i:S_i\hookrightarrow\tilde S_i$ with $\tilde S_I\in\SmVar(\mathbb C)$.
The category $C_{\mathcal D(1,0)fil,rh}(S/(\tilde S_I))\times_I D_{fil}(S^{an})$ is the category 
\begin{itemize}
\item whose set of objects is the set of triples $\left\{(((M_I,F,W),u_{IJ}),(K,W),\alpha)\right\}$ with
\begin{eqnarray*} 
((M_I,F,W),u_{IJ})\in C_{\mathcal D(1,0)fil,rh}(S/(\tilde S_I)), \, (K,W)\in D_{fil}(S^{an}), \\ 
\alpha:T(S/(\tilde S_I))(K,W)\otimes\mathbb C_{S^{an}}\to DR(S)^{[-]}(((M_I,W),u_{IJ})^{an})
\end{eqnarray*}
where 
\begin{equation*}
DR(S)^{[-]}:=DR(S^{an})^{[-]}:C_{\mathcal D(1,0)fil,rh}(S^{an}/(\tilde S_I^{an}))
\to C_{fil}(S^{an}/(\tilde S_I^{an})) 
\end{equation*}
is the De Rahm functor and $\alpha$ is a morphism in $D_{fil}(S^{an}/(\tilde S_I^{an}))$,
\item and whose set of morphisms consists of 
\begin{equation*}
\phi=(\phi_D,\phi_C,[\theta]):(((M_{1I},F,W),u_{IJ}),(K_1,W),\alpha_1)\to(((M_{2I},F,W),u_{IJ}),(K_2,W),\alpha_2)
\end{equation*}
where $\phi_D:((M_1,F,W),u_{IJ})\to((M_2,F,W),u_{IJ})$ and $\phi_C:(K_1,W)\to (K_2,W)$ 
are morphisms, and
\begin{eqnarray*}
\theta=(\theta^{\bullet},I(DR(S)(\phi^{an}_D))\circ I(\alpha_1),I(\alpha_2)\circ I(\phi_C\otimes I)): \\
I(T(S/(\tilde S_I))(K_1,W))\otimes\mathbb C_{S^{an}}[1]\to I(DR(S)(((M_{2I},W),u_{IJ})^{an}))  
\end{eqnarray*}
is an homotopy,  
$I:C_{fil}(S^{an}/(\tilde S_I^{an}))\to K_{fil}(S^{an}/(\tilde S_I^{an}))$
being the injective resolution functor : for
$((K_I,W),t_{IJ})\in C_{fil}(S^{an}/(\tilde S_I^{an}))$, 
we take an injective resolution 
\begin{equation*}
k:((K_I,W),t_{IJ})\to I((K_I,W),t_{IJ}) 
\end{equation*}
with $I((K,W),t_{IJ})\in C_{fil}(S^{an}/(\tilde S_I^{an}))$ 
which is unique modulo homotopy, and the class $[\theta]$ of $\theta$ does NOT depend of the injective resolution ;
in particular we have 
\begin{equation*}
DR(S)^{[-]}(\phi^{an}_D)\circ\alpha_1=\alpha_2\circ(\phi_C\otimes I) 
\end{equation*}
in $D_{fil}(S^{an}/(\tilde S_I^{an}))$ ; and for
\begin{itemize}
\item $\phi=(\phi_D,\phi_C,[\theta]):(((M_{1I},F,W),u_{IJ}),(K_1,W),\alpha_1)\to(((M_{2I},F,W),u_{IJ}),(K_2,W),\alpha_2)$
\item $\phi'=(\phi'_D,\phi'_C,[\theta']):(((M_{2I},F,W),u_{IJ}),(K_2,W),\alpha_2)\to(((M_{3I},F,W),u_{IJ}),(K_3,W),\alpha_3)$
\end{itemize}
the composition law is given by 
\begin{eqnarray*}
\phi'\circ\phi:=(\phi'_D\circ\phi_D,\phi'_C\circ\phi_C,
I(DR(S)(\phi^{'an}_D))\circ[\theta]+[\theta']\circ I(\phi_C\otimes I)[1]): \\
(((M_{1I},F,W),u_{IJ}),(K_1,W),\alpha_1)\to(((M_{3I},F,W),u_{IJ}),(K_3,W),\alpha_3)
\end{eqnarray*}
in particular for 
$(((M_I,F,W),u_{IJ}),(K,W),\alpha)\in C_{\mathcal D(1,0)fil,rh}(S/(\tilde S_I))\times_I D_{fil}(S^{an})$,
\begin{equation*}
I_{(((M_I,F,W),u_{IJ}),(K,W),\alpha)}=((I_{M_I}),I_K,0).
\end{equation*}
\end{itemize}
We have then full embeddings
\begin{eqnarray*}
\PSh^0_{\mathcal D(1,0)fil,rh}(S/(\tilde S_I))\times_I P_{fil}(S^{an})
\hookrightarrow C^0_{\mathcal D(1,0)fil,rh}(S/(\tilde S_I))\times_I D_{fil}(S^{an}) \\
\xrightarrow{\iota^0_{S/\tilde S_I}} C_{\mathcal D(1,0)fil,rh}(S/(\tilde S_I))^0\times_I D_{fil}(S^{an})
\hookrightarrow C_{\mathcal D(1,0)fil,rh}(S/(\tilde S_I))\times_I D_{fil}(S^{an})
\end{eqnarray*}
where $\PSh^0_{\mathcal D(1,0)fil,rh}(S/(\tilde S_I))\times_I P_{fil}(S^{an})$ is the category 
\begin{itemize}
\item whose set of objects is the set of triples $\left\{(((M_I,F,W),u_{IJ}),(K,W),\alpha)\right\}$ with
\begin{eqnarray*} 
((M_I,F,W),u_{IJ})\in\PSh^0_{\mathcal D(1,0)fil,rh}(S/(\tilde S_I)), \, (K,W)\in P_{fil}(S^{an}), \\ 
\alpha:T(S/(\tilde S_I))(K,W)\otimes\mathbb C_{S^{an}}\to DR(S)^{[-]}(((M_I,W),u_{IJ})^{an})
\end{eqnarray*}
where $DR(S)^{[-]}$ is the De Rahm functor and $\alpha$ is an isomorphism in $D_{fil}(S^{an}/(\tilde S_I^{an}))$,
\item and whose set of morphisms are 
\begin{equation*}
\phi=(\phi_D,\phi_C)=(\phi_D,\phi_C,0):
(((M_{1I},F,W),u_{IJ}),(K_1,W),\alpha_1)\to(((M_{2I},F,W),u_{IJ}),(K_2,W),\alpha_2)
\end{equation*}
where $\phi_D:((M_1,F,W),u_{IJ})\to((M_2,F,W),u_{IJ})$ and $\phi_C:(K_1,W)\to (K_2,W)$ are morphisms (of filtered sheaves)
such that $\phi_D^{an}\circ\alpha_1=\alpha_2\circ\phi_C$ in $P_{fil}(S^{an})$.
\end{itemize}
\end{itemize}

For holonomic $D$-modules on a smooth variety $S\in\SmVar(\mathbb C)$, there exist for a closed embedding $Z\subset S$
with $Z$ smooth, a $V_Z$-filtration (see definition \ref{Vfil}) satisfying further hypothesis so that it is unique: 

\begin{defi}\label{VfilKM}
Let $S\in\SmVar(\mathbb C)$ or $S\in\AnSm(\mathbb C)$. 
\begin{itemize}
\item[(i)] Let $D=V(s)\subset S$ be a smooth (Cartier) divisor, 
where $s\in\Gamma(S,L)$ is a section of the line bundle $L=L_D$ associated to $D$. 
Let $M\in\PSh_{\mathcal D}(S)$. 
A $V_D$-filtration $V$ for $M$ (see definition \ref{Vfil}) is called a Kashiwara-Malgrange $V_D$-filtration for $M$ if
\begin{itemize}
\item $V_kM$ are coherent $V_{Dk}D_S$ modules for all $k\in\mathbb Z$, that is $V$ is a good filtration,
\item $sV_kM=V_{k-1}M$ for $k<<0$,
\item all eigenvalues of $s\partial_s:\Gr_{V,k}:=V_kM/V_{k-1}M\to\Gr_{V,k}M:=V_kM/V_{k-1}M$ have real part between $k-1$ and $k$.
\end{itemize}
Almost by definition, a Kashiwara-Malgrange $V_D$-filtration for $M$ if it exists is unique (see \cite{Sa2})
so that we denote it by $(M,V_D)\in\PSh_{O_Sfil}(S)$ and $(M,V_D)$ is strict.
In particular if $m:(M_1,F)\to (M_2,F)$ a morphism with $(M_1,F),(M_2,F)\in\PSh_{\mathcal D(2)fil}(S)$
such that $M_1$ and $M_2$ admit the Kashiwara-Malgrange filtration for $D\subset S$,
we have $m(V_{D,q}F^pM_1)\subset V_{D,q}F^pM_2$, that is we get $m:(M_1,F,V_D)\to(M_2,F,V_D)$ a filtered morphism,
and if $0\to M'\to M\to M''\to 0$ is an exact sequence, 
$0\to (M',V_D)\to (M,V_D)\to (M'',V_D)\to 0$ is an exact sequence (strictness).
\item[(ii)]More generally, let $Z=V(s_1,\ldots,s_r)=D_1\cap\cdots\cap D_r\subset S$ be a smooth Zariski closed subset,
where $s_i\in\Gamma(S,L_i)$ is a section of the line bundle $L=L_{D_i}$ associated to $D_i$. 
Let $M\in\PSh_{\mathcal D}(S)$. 
A $V_Z$-filtration $V$ for $M$ (see definition \ref{Vfil}) is called a Kashiwara-Malgrange $V_Z$-filtration for $M$ if
\begin{itemize}
\item $V_kM$ are coherent $O_S$ modules for all $k\in\mathbb Z$,
\item $\sum_{i=1}^rs_iV_kM=V_{k-1}M$ for $k<<0$,
\item all eigenvalues of $\sum_{i=1}^rs_i\partial_{s_i}:\Gr_{V,k}M:=V_kM/V_{k-1}M\to\Gr^V_kM:=V_kM/V_{k-1}M$ 
have real part between $k-1$ and $k$.
\end{itemize}
Almost by definition, a Kashiwara-Malgrange $V_Z$-filtration for $M$ if it exists is unique (see \cite{Sa2})
so that we denote it by $(M,V_Z)\in\PSh_{O_Sfil}(S)$ and $(M,V_Z)$ is strict.
In particular if $m:(M_1,F)\to (M_2,F)$ a morphism with $(M_1,F),(M_2,F)\in\PSh_{\mathcal D(2)fil}(S)$
such that $M_1$ and $M_2$ admit the Kashiwara-Malgrange filtration for $D\subset S$,
we have $m(V_{Z,q}F^pM_1)\subset V_{Z,q}F^pM_2$, that is we get $m:(M_1,F,V_Z)\to(M_2,F,V_Z)$ a filtered morphism,
and if $0\to M'\to M\to M''\to 0$ is an exact sequence, 
$0\to (M',V_Z)\to (M,V_Z)\to (M'',V_Z)\to 0$ is an exact sequence (strictness).
\end{itemize}
\end{defi}

\begin{prop}\label{VfilVarAn}
\begin{itemize}
\item[(i)]Let $S\in\AnSm(\mathbb C)$. 
\begin{itemize}
\item Let $D=V(s)\subset S$ a smooth (Cartier) divisor, 
where $s\in\Gamma(S,L)$ is a section of the line bundle $L=L_D$ associated to $D$. 
If $M\in\PSh_{\mathcal D,rh}(S)$, the Kashiwara-Malgrange $V_D$-filtration for $M$ (see definition \ref{VfilKM}) exist
so that we denote it by $(M,V_D)\in\PSh_{O_Sfil}(S)$.
\item More generally, let $Z=V(s_1,\ldots,s_r)=D_1\cap\cdots\cap D_r\subset S$ be a smooth Zariski closed subset,
where $s_i\in\Gamma(S,L_i)$ is a section of the line bundle $L=L_{D_i}$ associated to $D_i$.
If $M\in\PSh_{\mathcal D(2)rh}(S)$, the Kashiwara-Malgrange $V_Z$-filtration for $M$ (see definition \ref{VfilKM}) exist
so that we denote it by $(M,V_Z)\in\PSh_{O_Sfil}(S)$.
\end{itemize}
\item[(ii)]Let $S\in\SmVar(\mathbb C)$. 
\begin{itemize}
\item Let $D=V(s)\subset S$ a smooth (Cartier) divisor, 
where $s\in\Gamma(S,L)$ is a section of the line bundle $L=L_D$ associated to $D$. 
If $M\in\PSh_{\mathcal D,rh}(S)$, the Kashiwara-Malgrange $V_D$-filtration for $M$ (see definition \ref{VfilKM}) exist
so that we denote it by $(M,V_D)\in\PSh_{O_Sfil}(S)$.
\item More generally, let $Z=V(s_1,\ldots,s_r)=D_1\cap\cdots\cap D_r\subset S$ be a smooth Zariski closed subset,
where $s_i\in\Gamma(S,L_i)$ is a section of the line bundle $L=L_{D_i}$ associated to $D_i$.
If $M\in\PSh_{\mathcal D(2)rh}(S)$, the Kashiwara-Malgrange $V_Z$-filtration for $M$ (see definition \ref{VfilKM}) exist
so that we denote it by $(M,V_Z)\in\PSh_{O_Sfil}(S)$.
\end{itemize}
\end{itemize}
\end{prop}

\begin{proof}
\noindent (i):Follows from the work of Kashiwara. 
Note that the second point is a particular case of the first by induction.
\noindent(ii): Take a compactification $\bar S\in\PSmVar(\mathbb C)$ of $S$ and denote by $\bar D\subset\bar S$
the closure of $D$. Using the closed embedding $i:\bar S\hookrightarrow L_{\bar D}$ given by the zero section,
we may assume that $\bar D$ is smooth. Denote by $j:\bar S\backslash\bar D\hookrightarrow\bar S$ the open complementary.
Then, $j_*M\in\PSh_{\mathcal D,rh}(\bar S)$ is regular holonomic. 
The result then follows by (i) and GAGA for $j_*M$ and we get $(j_*M,V_D)\in\PSh_{O_{\bar S}fil}(\bar S)$
and $(M,V_D)=(j^*j_*M,j^*V_D)\in\PSh{O_Sfil}(S)$.
We can also prove the algebraic case directly using the theory of meromorphic connexions 
since a simple holonomic $D_S$-module with support $Z\subset S$ is an integrable connexion on $Z^o=Z\cap S^o$,
$S^o\subset S$ being an open subset.
\end{proof}

We have from Kashiwara or Malgrange the following which relates the graded piece of the Kashiwara-Magrange $V$-filtration $V_D$
of a $D_S$ module $M\in\PSh_{\mathcal D,rh}(S)$ along a smooth divisor $D$ 
with the nearby and vanishing cycle functors of $DR(S)(M)$ with respect to $D$ :

\begin{thm}\label{VfilKMDthm}
Let $S\in\AnSm(\mathbb C)$. Let $D=V(s)\subset S$ be a smooth (Cartier) divisor, 
where $s\in\Gamma(S,L)$ is a section of the line bundle $L=L_D$ associated to $D$. 
Denote by $j:S^o:=S\backslash D\hookrightarrow S$ the open complementary embedding and by 
$k:\tilde S^o\xrightarrow{k}S^o\xrightarrow{j} S$ with $k$ the universal covering of $S^o$
For $M\in\PSh_{\mathcal D,rh}(S)$ a regular holonomic $D_S$ module, 
consider $(M,V_D)\in\PSh_{O_Sfil}(S)$ it together with its $V_D$ filtration. Then,
\begin{itemize}
\item there is canonical isomorphism 
\begin{equation*}
T(\psi_D,DR)(M):DR(S)(\oplus_{-1\leq\alpha<0}\Gr_{V_D,\alpha}M)\xrightarrow{\sim}\psi_D(DR(S)(M)):=Rk_*k^*DR(S)(M)
\end{equation*}
\item there is canonical isomorphism 
\begin{eqnarray*}
T(\phi_D,DR)(M):DR(S)(\oplus_{-1<\alpha\leq 0}\Gr_{V_D,\alpha}M)\xrightarrow{\sim} \\
\phi_D(DR(S)(M)):=\Cone(DR(S)(M)\xrightarrow{\ad(k^*,Rk_*)(-)}\psi_DDR(S)(M))
\end{eqnarray*}
\item $DR(S)(\partial_s)=can(DR(S)(M))\circ T(\psi_D,DR)(M)$, 
with $can(DR(S)(M)):\psi_DDR(S)(M)\to\phi_D(DR(S)(M))$ the structural embedding of complexes of the cone, 
\item $DR(S)(s\partial_s)=T\circ T(\psi_D,DR)(M)$, with $T:\psi_D(DR(S)(M))\to\psi_D(DR(S)(M))$ the monodromy morphism.
\item $T(\phi_D,DR)(M):DR(S)(s)\simeq var(DR(S)(M))$ with 
$var(DR(S)(M)):=\mathbb D^vcan(\mathbb D^vDR(S)M):\phi_D(DR(S)(M))\to\psi_D(DR(S)(M))$.
\end{itemize}
\end{thm}

\begin{proof}
See \cite{Sa2}.
\end{proof}

The main tool is the nearby and vanishing cycle functors for Cartier divisors.
We need for the definition of Hodge modules on a smooth complex algebraic variety $S$
to extend the V-filtration associated to a smooth Cartier divisor $D\subset S$ of regular holonomic
$D_S$ module $M$ such that the monodromy morphism $T:\psi_D(DR(S)(M))\to\psi_D(DR(S))$ is quasi-unipotent
by a rational V-filtration (i.e. indexed by rational numbers).

\begin{defi}\label{phipsiMHM}
Let $S\in\SmVar(\mathbb C)$ or $S\in\AnSm(\mathbb C)$. 
Let $D=V(s)\subset S$ a (Cartier) divisor, where $s\in\Gamma(S,L)$ 
is a section of the line bundle $L=L_D$ associated to $D$. We then have the zero section embedding
$i:S\hookrightarrow L$. We denote $L_0=i(S)$ and $j:L^o:=L\backslash L_0\hookrightarrow L$ the open complementary subset.
We denote $\PSh_{\mathcal D(2)fil,rh}(S)^{sp_D0}\subset\PSh_{\mathcal D(2)fil,rh}(S)$ 
the full subcategory consiting of objects such that 
the monodromy operator $T:\psi_D(DR(S)(M^{(an)}))\to\psi_D(DR(S)(M^{(an)}))$ is quasi-unipotent.
\begin{itemize}
\item[(i)] Let $(M,F)\in\PSh_{\mathcal D(2)fil,rh}(S)^{sp_{D0}}$
By proposition \ref{VfilVarAn}, we have the Kashiwara-Malgrange $V_S$-filtration on $i_{*mod}M$.
We refine it to all rational numbers as follows : 
for $\alpha=k-1+r/q\in\mathbb Q$, $k,q,r\in\mathbb Z$, $q\leq 0$, $0\leq r\leq q-1$, we set
\begin{eqnarray*}
V_{S,\alpha}M:=q_{V,k}^{-1}(\oplus_{k-1<\beta\leq\alpha}\Gr_{k,\beta}^{V_S}M\subset V_{S,k}M
\end{eqnarray*}
with $\Gr_{k,\beta}^{V_S}M:=\ker(\partial_ss-\beta I)\subset\Gr_k^{V_S}M$ and $q_{V,k}:V_{S,k}M\to\Gr_k^{V_S}M$
is the projection.
We set similarly
\begin{eqnarray*}
V_{S,<\alpha}M:=q_{V,k}^{-1}(\oplus_{k-1<\beta <\alpha}\Gr_{k,\beta}^{V_S}M\subset V_{S,k}M
\end{eqnarray*}
The Hodge filtration induced on $\Gr^V_{\alpha}M$ is 
\begin{eqnarray*}
F^p\Gr^{V_S}_{\alpha}M:=(F^pM\cap V_{S,\alpha}M)/(F^pM\cap V_{S,<\alpha}M)
\end{eqnarray*}
\item[(ii)] we have using (i) the nearby cycle functor
\begin{eqnarray*}
\psi_D:\PSh_{\mathcal Dfil,rh}(S)^{sp_D0}\to\PSh_{\mathcal Dfil,rh}(D/(S)), \; 
(M,F)\mapsto\psi_D(M,F):=\oplus_{-1\leq\alpha<0}\Gr_{V_S,\alpha}i_{*mod}(M,F), 
\end{eqnarray*}
and vanishing cycle functors
\begin{eqnarray*}
\phi_D:\PSh_{\mathcal D(2)fil,rh}(S)^{sp_D0}\to\PSh_{\mathcal D(2)fil,rh}(D/(S)), \; 
(M,F)\mapsto\phi_D(M,F):=:=\oplus_{-1<\alpha\leq 0}\Gr_{V_S,\alpha}i_{*mod}(M,F), \\
\phi_{D1}:\PSh_{\mathcal D(2)fil,rh}(S)^{sp_D0}\to\PSh_{\mathcal D(2)fil,rh}(D/(S)), \; 
(M,F)\mapsto\phi_{D1}(M,F):=\Gr_{V_S,0}i_{*mod}(M,F).
\end{eqnarray*}
\item[(iii)] This induces, by theorem \ref{VfilKMDthm}, the nearby cycle functor
\begin{eqnarray*}
\psi_D:\PSh_{\mathcal D(1,0)fil,rh}(S)^{sp_D0}\times_I P_{fil}(S^{(an)})\to
\PSh_{\mathcal D(1,0)fil,rh,D}(S)\times_I P_{fil,D}(S^{an}), \\
((M,F,W),(K,W),\alpha)\mapsto\psi_D((M,F,W),(K,W),\alpha):=(\psi_D(M,F,W),\psi_D(K,W),\psi_D(\alpha))
\end{eqnarray*}
and the vanishing cycle functor
\begin{eqnarray*}
\phi_{D1}:\PSh_{\mathcal D(1,0)fil,rh}(S)^{sp_D0}\times_I D_{fil}(S^{(an)})\to 
\PSh_{\mathcal D(1,0)fil,rh,D}(S)\times_I P_{fil,D}(S^{an}), \\ 
((M,F,W),(K,W),\alpha)\mapsto\phi_{D1}((M,F,W),(K,W),\alpha):=(\phi_{D1}(M,F,W),\phi_D(K,W),\phi_D(\alpha))
\end{eqnarray*}
\end{itemize}
\end{defi}

We have the category of mixed Hodge modules over a complex algebraic variety or a complex analytic space $S$ 
defined by, for $S$ smooth, by induction on dimension of $S$,
and for $S$ singular using embeddings into smooth complex algebraic varieties, resp. smooth complex analytic spaces: 

\begin{defi}\cite{Saito}
\begin{itemize}
\item[(i)]Let $S\in\SmVar(\mathbb C)$ or $S\in\AnSm(\mathbb C)$. 
Denote $\PSh_{\mathcal Dfil,rh}(S)^{sp,ssd}\subset\PSh_{\mathcal Dfil,rh}(S)$
the full subcategory consisting of objects $(M,F)$ 
\begin{itemize}
\item such that for all Cartier divisor $D=V(s)\subset S$, $s\in\Gamma(S,L)$, 
denoting $i:S\hookrightarrow L$ the closed embedding
the monodromy morphism $T:\psi_D(DR(S)(M^{(an)}))\to\psi_D(DR(S)(M^{(an)}))$ is quasi-unipotent,
$sF^pV_{S,\alpha}i_{*mod}M=F^pV_{S,\alpha-1}i_{*mod}M$ for $\alpha<0$,
$\partial_sF^p\Gr^{V_S}_{\alpha}i_{*mod}M=\Gr^{V_S}_{\alpha+1}i_{*mod}M$ for $\alpha>-1$,
the filtration induced by $F$ on $\Gr^{V_S}_{\alpha}i_{*mod}M$ is good,
\item which admits a decomposition with $D_S$ module with strict support on closed irreducible subvarieties.
\end{itemize}
The category of Hodge modules over $S$ of weight $w$ is the full subcategory 
\begin{equation*}
\iota_S:HM(S,w)=\oplus_{d\in\mathbb N}HM_{\leq d}(S,w)\hookrightarrow\PSh_{\mathcal Dfil,rh}(S)^{sp,ssd}\times_I P(S^{(an)}), 
\hookrightarrow\PSh_{\mathcal Dfil,rh}(S)\times_I P(S^{(an)})
\end{equation*}
given inductively by, $d$ being the dimension of the support of the $D_S$ modules,
\begin{itemize}
\item for $i_0:s_0\hookrightarrow S$ a closed point,
$i_{0*}\iota_{\pt}:HM_{s_0}(S,w)=HS\hookrightarrow\PSh_{\mathcal Dfil,rh,s_0}(S)\times_I P_{s_0}(S^{(an)})$ 
consist of Hodge structures of weight $w$, this gives $HM_0(S,w)$
\item for $Z\subset S$ an irreducible closed subvariety of dimension $d$,
$((M,F),K,\alpha)\in\PSh_{\mathcal Dfil,rh}(S)\times_I P(S^{(an)})$ belongs to $HM_Z(S,w)$ if and only if
$M$ has strict support $Z$ (i.e. $\supp(M)=Z$ and for all non trivial subobject $N$ or quotient of $M$ $\supp(N)=Z$),
and for all proper maps $f:S^o\to\mathbb A^1$ such that $f_{|Z\cap S^o}\neq 0$, $j:S^o\hookrightarrow S$ being an open subset, 
\begin{equation*}
\Gr^{W(N)}_k\psi_{f^{-1}(0)}(j^*(M,F),j^*K,j^*\alpha)
\in HM_{\leq d-1}(S^o,w-1+k)\hookrightarrow\PSh_{\mathcal Dfil,rh,f^{-1}(0)}(S^o)\times_I P_{f^{-1}(0)}({S^o}^{(an)})
\end{equation*}
for all $k\in\mathbb Z$, see definition \ref{phipsiMHM}, $W(N)$ being the weight filtration associated
to the monodromy morphism $T:\psi_{f^{-1}(0)}(DR(S)(M^{(an)}))\to\psi_{f^{-1}(0)}(DR(S)(M^{(an)}))$,
we then set $HM_{\leq d}(S,w):=\oplus_{Z\subset S,\dim(Z)=d}HM_Z(S,w)$.
\end{itemize}
\item[(ii)]Let $S\in\SmVar(\mathbb C)$ or $S\in\AnSm(\mathbb C)$. 
The category of mixed Hodge modules over $S$ is the full subcategory 
\begin{equation*}
\iota_S:MHM(S)\hookrightarrow MHW(S)\hookrightarrow\PSh_{\mathcal D(1,0)fil,rh}(S)\times_I P_{fil}(S^{(an)}), 
\end{equation*}
where the full subcategory $MHW(S)$ consists of objects 
$((M,F,W),(K,W),\alpha)\in\PSh_{\mathcal D(1,0)fil,rh}(S)\times_I P_{fil}(S^{(an)})$ satisfy 
\begin{equation*}
(\Gr^W_i((M,F,W),\Gr^W_i(K,W),\Gr^W_i\alpha)\in HM(S).
\end{equation*}
and the objects of $MHM(S)$ satisfy in addition an admissibility condition 
(in particular the three filtration $F,W,V_Z$ are compatible).
As usual, for $Z\subset S$ a closed subset and $j:S\backslash Z\hookrightarrow S$ the open complementary subset,
we denote $MHM_Z(S)\subset MHM(S)$ the full subcategory consisting of $((M,F,W),(K,W),\alpha)\in MHM(S)$ such that 
\begin{equation*}
j^*((M,F,W),(K,W),\alpha):=(j^*(M,F,W),j^*(K,W),j^*\alpha)=0.
\end{equation*}
\item[(iii)]Let $S\in\Var(\mathbb C)$ or $S\in\AnSp(\mathbb C)$ non smooth. Take an open
cover $S=\cup_iS_i$ so that there are closed embedding $S_I\hookrightarrow\tilde S_I$, with $S_I\in\SmVar(\mathbb C)$,
resp $S_I\in\AnSm(\mathbb C)$. The category of mixed Hodge modules over $S$ is the full subcategory
\begin{equation*}
\iota_S:MHM(S)\hookrightarrow MHW(S)\hookrightarrow\PSh_{\mathcal D(1,0)fil,rh}(S/(\tilde S_I))\times_I P_{fil}(S^{(an)}) 
\end{equation*}
consisting of objects 
\begin{equation*}
(((M_I,F,W),u_{IJ}),(K,W),\alpha)\in\PSh^0_{\mathcal D(1,0)fil,rh}(S/(\tilde S_I))\times_I P_{fil}(S^{(an)})
\end{equation*}
such that $((M_I,F,W),T(S/(\tilde S_I))(K,W),\alpha)\in(MHM_{S_I}(\tilde S_I))$ (see (ii)).
The category $MHM(S)$ does NOT depend on the open cover an the closed embedding by proposition \ref{TSIHdg}.
\item[(iv)]Let $S\in\Var(\mathbb C)$. We get from (iii) $D(MHM(S)):=(\Ho_{zar},I)(C(MHM(S)))$.
By induction, using the result for mixed hodge structure 
and the strictness of the Kashiwara-Malgrange $V$-filtration for morphism of $D$-module, 
the morphism of $MHM(S)$ are strict for $F$ and $W$ (see \cite{Saito}).
\item[(iv)']Let $S\in\AnSp(\mathbb C)$. We get from (iii) $D(MHM(S)):=(\Ho_{usu},I)(C(MHM(S)))$.
By induction, using the result for mixed hodge structure 
and the strictness of the Kashiwara-Malgrange $V$-filtration for morphism of $D$-module, 
the morphism of $MHM(S)$ are strict for $F$ and $W$ (see \cite{Saito}).
\end{itemize}
\end{defi}

\begin{itemize}
\item Let $S\in\SmVar(\mathbb C)$. We consider the canonical functor
\begin{eqnarray*}
\pi_S:C(MHW(S))\xrightarrow{\iota_S} C_{\mathcal D(1,0)fil}(S)\times_I D_{fil}(S^{an})
\xrightarrow{p_S} C_{\mathcal D(1,0)fil}(S), \\ 
((M,F,W),(K,W),\alpha)\mapsto (M,F,W)
\end{eqnarray*}
where $p_S$ is the projection functor. 
Then $\pi_S(MHW(S))\subset\PSh_{\mathcal D(1,0)fil}(S)$ is the
subcategory consisting of $(M,F,W)\in\PSh_{\mathcal D(1,0)fil}(S)$ such that
$((M,F,W),(K,W),\alpha)\in MHW(S)$ is a W filtered Hodge module for some $(K,W)\in C_{fil}(S)$. 
It induces in the derived category the functor 
\begin{eqnarray*}
\pi_S:D(MHW(S))\xrightarrow{\iota_S}D_{\mathcal D(1,0)fil}(S)\times_I D_{fil}(S^{an})
\xrightarrow{p_S} D_{\mathcal D(1,0)fil}(S), \\ 
((M,F,W),(K,W),\alpha)\mapsto (M,F,W)
\end{eqnarray*}
after localization with respect to filtered Zariski and usu local equivalence.

\item Let $S\in\Var(\mathbb C)$ non smooth. 
Take an open cover $S=\cup_iS_i$ such that there are closed embedding $S_I\hookrightarrow\tilde S_I$, 
with $S_I\in\SmVar(\mathbb C)$. We consider the canonical functor
\begin{eqnarray*}
\pi_S:C(MHW(S))\hookrightarrow C_{\mathcal D(1,0)fil}(S/(\tilde S_I))\times_I D_{fil}(S^{an})
\xrightarrow{p_S} C_{\mathcal D(1,0)fil}(S/(\tilde S_I)), \\ 
(((M_I,F,W),u_{IJ}),(K,W),\alpha)\mapsto ((M_I,F,W),u_{IJ})
\end{eqnarray*}
where $p_S$ is the projection functor. 
Then $\pi_S(MHW(S))\subset\PSh_{\mathcal D(1,0)fil}(S/(\tilde S_I))$ is the
subcategory consisting of $((M,F,W),u_{IJ})\in\PSh_{\mathcal D(1,0)fil}(S/(\tilde S_I))$ such that
$(((M_I,F,W),u_{IJ}),(K,W),\alpha)\in MHW(S)$ is a W filtered Hodge module for some $(K,W)\in C_{fil}(S)$. 
It induces in the derived category the functor 
\begin{eqnarray*}
\pi_S:D(MHW(S))\xrightarrow{\iota_S} D_{\mathcal D(1,0)fil}(S/(\tilde S_I))\times_I D_{fil}(S^{an})
\xrightarrow{p_S} D_{\mathcal D(1,0)fil}(S/(\tilde S_I)), \\ 
(((M_I,F,W),u_{IJ},(K,W),\alpha)\mapsto ((M,F,W),u_{IJ})
\end{eqnarray*}
after localization with respect to filtered Zariski and usu local equivalence.
\end{itemize}

Let $S\in\Var(\mathbb C)$ or $S\in\AnSp(\mathbb C)$. 
\begin{itemize}
\item If $S\in\Var(\mathbb C)$, let $S=\cup_{i=1}^lS_i$ an open cover such that there exist closed embeddings 
$i_i:S_i\hookrightarrow\tilde S_i$ with $\tilde S_i\in\SmVar(\mathbb C)$,
and let $S=\cup_{i'=1}^{l'}S_{i'}$ an other open cover such that there exist closed embeddings 
$i_{i'}:S_{i'}\hookrightarrow\tilde S_{i'}$ with $\tilde S_{i'}\in\SmVar(\mathbb C)$.
\item If $S\in\AnSp(\mathbb C)$, let $S=\cup_{i=1}^lS_i$ an open cover such that there exist closed embeddings 
$i_i:S_i\hookrightarrow\tilde S_i$ with $\tilde S_i\in\AnSm(\mathbb C)$,
and let $S=\cup_{i'=1}^{l'}S_{i'}$ an other open cover such that there exist closed embeddings 
$i_{i'}:S_{i'}\hookrightarrow\tilde S_{i'}$ with $\tilde S_{i'}\in\AnSm(\mathbb C)$.
\end{itemize}
Denote $L=[1,\ldots,l]$, $L'=[1,\ldots,l']$ and $L'':=[1,\ldots,l]\sqcup[1,\ldots,l']$.
We have then the refined open cover $S=\cup_{k\in L} S_k$ and we denote for $I\sqcup I'\subset L''$, 
$S_{I\sqcup I'}:=\cap_{k\in I\sqcup I'} S_k$ and $\tilde S_{I\sqcup I'}:=\Pi_{k\in I\sqcup I'} \tilde S_k$,
so that we have a closed embedding $i_{I\sqcup I'}:S_{I\sqcup I'}\hookrightarrow\tilde S_{I\sqcup I'}$.
Consider $\pi^L_S(MHM(S))\subset\PSh_{\mathcal Dfil}(S/(S_I))$ and $\pi^{L'}_S(MHM(S))\subset\PSh_{\mathcal Dfil}(S/(S_{I'}))$.
For $I\sqcup I'\subset J\sqcup J'$, denote by 
$p_{I\sqcup I',J\sqcup J'}:\tilde S_{J\sqcup J'}\to\tilde S_{I\sqcup I'}$ the projection.
We then have a natural transfer map
\begin{eqnarray*}
T^{L/L'}_{S}:\pi^L_S(MHM(S))\to\pi^{L'}_S(MHM(S)), \\ 
((M_I,F,W),s_{IJ})\mapsto (\ho\lim_{I\in L}p_{I'(I\sqcup I')*}
\Gr_{V_{I\sqcup I'}}p_{I(I\sqcup I')}^{*mod}(M_I,F)),s_{I'J'}),
\end{eqnarray*}
with, in the homotopy limit, the natural transition morphisms 
\begin{eqnarray*}
p_{I'(I\sqcup I')*}\ad(p_{IJ}^{*mod},p_{IJ*})(p_{I(I\sqcup I')}^{*mod[-]}(M_I,F)): \\
p_{I'(J\sqcup I')*}(\Gr_{V_{J\sqcup I'}}p_{J(J\sqcup I')}^{*mod[-]}(M_J,F))\to
p_{I'(I\sqcup I')*}(\Gr_{V_{I\sqcup I'}}p_{I(I\sqcup I')}^{*mod[-]}(M_I,F))
\end{eqnarray*}
for $J\subset I$, and 
\begin{eqnarray*}
s_{I'J'}:\holim_{I\in L}m^*p_{I'(I\sqcup I')*}(\Gr_{V_{I\sqcup I'}}p_{I(I\sqcup I')}^{*mod[-]}(M_I,F)\to \\
\holim_{I\in L}p_{I'J'*}\Gr_{V_{J'}}(p_{I'J'}^{*mod[-]}m^*p_{I'(I\sqcup I')*}
\Gr_{V_{I\sqcup I'}}p_{I(I\sqcup I')}^{*mod[-]}((M_I,F))) \\
\to\holim_{I\in L} p_{I'J'*}p_{J'(I\sqcup J')*}\Gr_{V_{I\sqcup I'}}p_{I(I\sqcup J')}^{*mod[-]}(M_I,F)
\end{eqnarray*}

\begin{prop}\label{TSIHdg}
\begin{itemize}
\item[(i)]Let $S\in\Var(\mathbb C)$ and let $S=\cup_iS_i$ an open cover
such that there exist closed embeddings $i_iS_i\hookrightarrow\tilde S_i$ with $\tilde S_i\in\SmVar(\mathbb C)$.
Then $\pi_S(MHM(S)\subset\PSh_{\mathcal D(2)fil}(S/(\tilde S_I))$ 
does not depend on the open covering of $S$ and the closed embeddings.
More precisely, let $S=\cup_{i'=1}^{l'}S_{i'}$ an other open cover such that there exist closed embeddings 
$i_{i'}:S_{i'}\hookrightarrow\tilde S_{i'}$ with $\tilde S_{i'}\in\SmVar(\mathbb C)$. Then,
\begin{equation*}
T^{L/L'}_{S}:\pi^L_S(MHM(S))\to\pi^{L'}_S(MHM(S)), 
\end{equation*} 
is an equivalence of category with inverse is $T_S^{L'/L}:\pi^{L'}_S(MHM(S))\to\pi^L_S(MHM(S))$.
\item[(ii)]Let $S\in\AnSp(\mathbb C)$ and let $S=\cup_iS_i$ an open cover
such that there exist closed embeddings $i_iS_i\hookrightarrow\tilde S_i$ with $\tilde S_i\in\AnSm(\mathbb C)$.
Then $\pi_S(MHM(S))\subset\PSh_{\mathcal D(2)fil}(S/(\tilde S_I))$ 
does not depend on the open covering of $S$ and the closed embeddings.
More precisely, let $S=\cup_{i'=1}^{l'}S_{i'}$ an other open cover such that there exist closed embeddings 
$i_{i'}:S_{i'}\hookrightarrow\tilde S_{i'}$ with $\tilde S_{i'}\in\AnSm(\mathbb C)$. Then,
\begin{equation*}
T^{L/L'}_{S}:\pi^L_S(MHM(S))\to\pi^{L'}_S(MHM(S)), 
\end{equation*} 
is an equivalence of category with inverse is $T_S^{L'/L}:\pi^{L'}_S(MHM(S))\to\pi^L_S(MHM(S))$.
\end{itemize}
\end{prop}

\begin{proof}
Follows from the definition of the Hodge filtration which use the $V$-filtration : see \cite{Saito}.
\end{proof}

The main results of Saito, which implies in the algebraic case the six functor formalism on $DMHM(-)$ are the followings

\begin{defi}\label{VMHSdef}
Let $S\in\SmVar(\mathbb C)$ or $S\in\AnSm(\mathbb C)$.
We denote by $VMHS(S)\subset\PSh_{\mathcal D(1,0)fil,rh}(S)\times_I P_{fil}(S^{an})$ the full subcategory
consisting of variation of mixed Hodge structure, whose objects consist of
\begin{equation*}
(((L_S,W)\otimes O_S,F),(L_S,W),\alpha)\subset\PSh_{\mathcal D(1,0)fil,rh}(S)\times_I P_{fil}(S^{an})
\end{equation*}
with 
\begin{itemize}
\item $L_S\in\PSh(S^{an})$ a local system,
\item the $D_S$ module structure on $(L_S,W)\otimes O_S$ is given by the flat connection associated to the local system $L_S$,
\item $F^p(W^qL_S\otimes O_S)\subset(W^qL_S\otimes O_S)$ are locally free $O_S$ subbundle satisfying Griffitz transversality
for the $D_S$ module structure (i.e. for the flat connection).
\item $\alpha:(L_S,W)\to DR(S)^{[-]}((L_S,W)\otimes O_S)$ is the isomorphism given by theorem \ref{DRK}.
\end{itemize}
\end{defi}

\begin{thm}\label{VMHSthm}
Let $S\in\SmVar(\mathbb C)$ or $S\in\AnSm(\mathbb C)$.
\begin{itemize}
\item[(i)] A variation of mixed Hodge structure $(((L_S,W)\otimes O_S,F),(L_S,W),\alpha)\in VMHS(S)$ (see definition \ref{VMHSdef})
is a mixed module. That is $VMHS(S)\subset MHM(S)$.
\item[(ii)] For $((M,F,W),(K,W),\alpha)\in MHM(S)$ a mixed Hodge module with support $\supp M=Z$,
there exist an open subset $j:S^o\hookrightarrow S$,
such that $j^*((M,F,W),(K,W),\alpha):=(j^*(M,F,W),j^*(K,W),j^*\alpha)\in VMHS(Z\cap S^o)$.
That is a mixed Hodge module is generically a variation of mixed Hodge structure.
\end{itemize}
\end{thm}

\begin{proof}
See \cite{Saito}.
\end{proof}

\begin{thm}\label{Sa1}
\begin{itemize}
\item[(i)]Let $f:X\to S$ a projective morphism with $X,S\in\AnSp(\mathbb C)$, where projective means that there exist a factorization
$f:X\xrightarrow{l}\mathbb P^N\times S\xrightarrow{p_S}S$ with $l$ a closed embedding and $p_S$ the projection.
Let $S=\cup_{i=1}^s S_i$ an open cover such that there exits closed embeddings 
$i_I:S_i\hookrightarrow\tilde S_i$ with $\tilde S_i\in\AnSm(\mathbb C)$. 
For $I\subset[1,\ldots,s]$, recall that we denote $S_I:=\cap_{i\in I} S_i$ and  $X_I:=f^{-1}(S_I)$.
We have then the following commutative diagram
\begin{equation*}
\xymatrix{X_I\ar[r]^{i_I\circ l_I} & \mathbb P^N\times\tilde S_I\ar[r]^{p_{\tilde S_I}} & \tilde S_I \\
X_J\ar[u]^{j'_{IJ}}\ar[r]^{i_J\circ l_J} & \mathbb P^N\times\tilde S_J\ar[r]^{p_{\tilde S_J}}\ar[u]^{p'_{IJ}} & 
\tilde S_J\ar[u]^{p_{IJ}}}
\end{equation*}
whose right square is cartesian (see section 5).Then, for 
\begin{equation*}
((M,F,W),(K,W),\alpha)=(((M_I,F,W),u_{IJ}),(K,W),\alpha)\in MHM(X), 
\end{equation*}
where $((M_I,F,W),u_{IJ})\in C_{\mathcal D2fil}(X_I/(\mathbb P^N\times\tilde S_I))$, $(K,W)\in C_{fil}(X)$, 
we have for all $n\in\mathbb Z$,
\begin{equation*}
(H^n\int^{FDR}_f((M_I,F,W),u_{IJ}),R^nf_*(K,W),H^nf_*(\alpha))\in MHM(S)
\end{equation*}
and for all $p\in\mathbb Z$, the differentials of $\Gr_F^p\int^{FDR}_f((M_I,F,W),u_{IJ})$ 
are strict for the the Hodge filtration $F$.
\item[(ii)]Let $f:X\to S$ a projective morphism with $X,S\in\Var(\mathbb C)$, where projective means that there exist a factorization
$f:X\xrightarrow{l}\mathbb P^N\times S\xrightarrow{p_S}S$ with $l$ a closed embedding and $p_S$ the projection.
Let $S=\cup_{i=1}^s S_i$ an open cover such that there exits closed embeddings 
$i_I:S_i\hookrightarrow\tilde S_i$ with $\tilde S_i\in\SmVar(\mathbb C)$. 
For $I\subset[1,\ldots,s]$, recall that we denote $S_I:=\cap_{i\in I} S_i$ and  $X_I:=f^{-1}(S_I)$.
We have then the following commutative diagram
\begin{equation*}
\xymatrix{X_I\ar[r]^{i_I\circ l_I} & \mathbb P^N\times\tilde S_I\ar[r]^{p_{\tilde S_I}} & \tilde S_I \\
X_J\ar[u]^{j'_{IJ}}\ar[r]^{i_J\circ l_J} & \mathbb P^N\times\tilde S_J\ar[r]^{p_{\tilde S_J}}\ar[u]^{p'_{IJ}} & 
\tilde S_J\ar[u]^{p_{IJ}}}
\end{equation*}
whose right square is cartesian (see section 5). Then, for 
\begin{equation*}
((M,F,W),(K,W),\alpha)=(((M_I,F,W),u_{IJ}),(K,W),\alpha)\in D(MHM(X)), 
\end{equation*}
where $((M_I,F,W),u_{IJ})\in C_{\mathcal D2fil}(X_I/(\mathbb P^N\times\tilde S_I))$, $(K,W)\in C_{fil}(X^{an})$, we have
\begin{eqnarray*}
H^n(\int^{FDR}_f((M_I,F,W),u_{IJ}),Rf_*(K,W),f_*(\alpha))\in MHM(S)
\end{eqnarray*}
for all $n\in\mathbb Z$, and for all $p\in\mathbb Z$, 
the differentials of $\Gr_F^p\int^{FDR}_f((M_I,F,W),u_{IJ})$ are strict for the the Hodge filtration $F$.
\end{itemize}
\end{thm}

\begin{proof}
\noindent(i):See \cite{Saito}.

\noindent(ii): By (i) $(H^n\int_f((M,F,W)^{an}),R^nf_*(K,W),H^nf_*(\alpha))\in MHM(S^{an})$ for all $n\in\mathbb Z$.
On the other hand, 
$T^{\mathcal D}(an,f)(M,F,W):(\int_f(M,F,W))^{an}\xrightarrow{\sim}\int_f((M,F,W)^{an})$
is an isomorphism since $f$ is proper by theorem GAGA for mixed hodge modules : see \cite{Saito}.
\end{proof}

\begin{thm}\label{Sa2}
\begin{itemize}
\item[(i)]Let $S\in\AnSp(\mathbb C)$. Let $Y\in\AnSm(\mathbb C)$ and $p_S:Y\times S\to S$ the projection. 
Let $S=\cup_{i=1}^s S_i$ an open cover such that there exits closed embeddings 
$i_I:S_i\hookrightarrow\tilde S_i$ with $\tilde S_i\in\AnSm(\mathbb C)$. 
For $I\subset[1,\ldots,s]$, recall that we denote $S_I:=\cap_{i\in I} S_i$.
We have then the following commutative diagram
\begin{equation*}
\xymatrix{Y\times\tilde S_I\ar[r]^{p_{\tilde S_I}} & \tilde S_I \\
Y\times\tilde S_J\ar[r]^{p_{\tilde S_J}}\ar[u]^{p'_{IJ}} & \tilde S_J\ar[u]^{p_{IJ}}}
\end{equation*}
which is cartesian (see section 5).
Then, for 
\begin{equation*}
((M,F,W),(K,W),\alpha)=(((M_I,F,W),u_{IJ}),(K,W),\alpha)\in MHM(S), 
\end{equation*}
where $((M_I,F,W),u_{IJ})\in C_{\mathcal D2fil}(S_I/(\tilde S_I))$, $(K,W)\in C_{fil}(S)$,
\begin{itemize}
\item $(p_S^{*mod[-]}(M,F,W),p_S^*(K,W),p_S^*(\alpha)):=
((p_{\tilde S_I}^{*mod[-]}(M_I,F,W),p_{\tilde S_J}^{*mod[-]}u_{IJ}),p_S^*(K,W),p_S^*(\alpha))\in MHM(S)$
\item $(p_S^{\hat*mod[-]}(M,F,W),p_S^*(K,W),p_S^*(\alpha)):=
((p_{\tilde S_I}^{\hat*mod[-]}(M_I,F,W),p_{\tilde S_J}^{*mod[-]}u_{IJ}),p_S^*(K,W),p_S^*(\alpha))\in MHM(S)$
\end{itemize}
\item[(ii)]Let $S\in\Var(\mathbb C)$. Let $Y\in\SmVar(\mathbb C)$ and $p_S:Y\times S\to S$ the projection. 
Let $S=\cup_{i=1}^s S_i$ an open cover such that there exits closed embeddings 
$i_I:S_i\hookrightarrow\tilde S_i$ with $\tilde S_i\in\SmVar(\mathbb C)$. 
For $I\subset[1,\ldots,s]$, recall that we denote $S_I:=\cap_{i\in I} S_i$.
We have then the following commutative diagram
\begin{equation*}
\xymatrix{Y\times\tilde S_I\ar[r]^{p_{\tilde S_I}} & \tilde S_I \\
Y\times\tilde S_J\ar[r]^{p_{\tilde S_J}}\ar[u]^{p'_{IJ}} & \tilde S_J\ar[u]^{p_{IJ}}}
\end{equation*}
which is cartesian (see section 5).Then, for 
\begin{equation*}
((M,F,W),(K,W),\alpha)=(((M_I,F,W),u_{IJ}),(K,W),\alpha)\in D(MHM(S)) 
\end{equation*}
where $((M_I,F,W),u_{IJ})\in C_{\mathcal D2fil}(S_I/(\tilde S_I))$, $(K,W)\in C_{fil}(S^{an})$, we have
\begin{itemize}
\item $(p_S^{*mod[-]},p_S^!)((M,F,W),(K,W),\alpha):=
((p_{\tilde S_I}^{*mod[-]}(M_I,F,W),p_{\tilde S_J}^{*mod[-]}u_{IJ}),p_S^*(K,W),p_S^*(\alpha))\in D(MHM(S))$
\item $(p_S^{\hat*mod[-]},p_S^*)((M,F,W),(K,W),\alpha):=
((p_{\tilde S_I}^{\hat*mod[-]}(M_I,F,W),p_{\tilde S_J}^{*mod[-]}u_{IJ}),p_S^*(K,W),p_S^*(\alpha))\in D(MHM(S))$.
\end{itemize}
\end{itemize}
\end{thm}

\begin{proof}
\noindent(i):See \cite{Saito}.

\noindent(ii):Follows immediately from (i) since 
$(p_{\tilde S_I}^{*mod[-]}(M_I,F,W))^{an}=p_{\tilde S_I}^{*mod[-]}((M_I,F,W)^{an})$.
\end{proof}

We have, by the results of Saito, the following key definition.

\begin{defi}\label{DHdgj}
\begin{itemize}
\item[(i)] Let $S\in\SmVar(\mathbb C)$ or $S\in\AnSm(\mathbb C)$.
Let $D=V(s)\subset S$ a divisor with $s\in\Gamma(S,L)$ and $L$ a line bundle ($S$ being smooth, $D$ is Cartier).
Denote by $j:S^o:=S\backslash D\hookrightarrow S$ the open complementary embedding.
Let $(M,F,W)\in\pi_{S^o}(MHW(S^o))$. Consider the $V_S$-filtration on $i_{*mod}M$ (see proposition \ref{VfilVarAn}).
If $(M,F,W)$ is extendable (which is always the case in the algebraic case), then,
\begin{itemize}
\item we have 
\begin{eqnarray*}
j_{*Hdg}(M,F,W):=(j_*M,F,W)\in\pi_S(MHM(S)), \\
\mbox{with} \; F^pj_*M:=\sum_{k\in\mathbb N}\partial_s^kF^{p+k}V_{S,0}j_*M\subset j_*M, \; 
W_kj_*M:=<j_*W_kM,W(N)_k>\subset j_*M 
\end{eqnarray*}
which is unique such that $j^*j_{*Hdg}(M,F,W)=(M,F,W)$ and $DR(S)(j_{*Hdg}(M,F,W))=j_*DR(S^o)(M,W)$,
\item there exist 
\begin{eqnarray*}
j_{!Hdg}(M,F,W):=\mathbb D_S^{Hdg}j_{*Hdg}\mathbb D_S^{Hdg}(M,F,W)\in\pi_S(MHM(S)) 
\end{eqnarray*}
unique such that $j^*j_{!Hdg}(M,F,W)=(M,F,W)$ and $DR(S)(j_{!Hdg}(M,F,W))=j_!DR(S^o)(M,W)$.
\end{itemize}
Moreover for $(M',F,W)\in\pi_S(MHM(S))$, by proposition \ref{Saprop}
\begin{itemize}
\item there is a canonical map $\ad(j^*,j_{*Hdg})(M',F,W):(M',F,W)\to j_{*Hdg}j^*(M',F,W)$ in $\pi_S(MHM(S))$,
\item there is a canonical map $\ad(j_{!Hdg},j^*)(M',F,W):j_{!Hdg}j^*(M',F,W)\to(M',F,W)$ in $\pi_S(MHM(S))$.
\end{itemize}
\item[(ii)] Let $S\in\SmVar(\mathbb C)$.
Let $Z=V(\mathcal I)\subset S$ an arbitrary closed subset, $\mathcal I\subset O_S$ being an ideal subsheaf. 
Taking generators $\mathcal I=(s_1,\ldots,s_r)$, we get $Z=V(s_1,\ldots,s_r)=\cap^r_{i=1}Z_i\subset S$ with 
$Z_i=V(s_i)\subset S$, $s_i\in\Gamma(S,\mathcal L_i)$ and $L_i$ a line bundle. 
Note that $Z$ is an arbitrary closed subset, $d_Z\geq d_X-r$ needing not be a complete intersection. 
Denote by $j:S^o:=S\backslash Z\hookrightarrow S$,
$j_I:S^{o,I}:=\cap_{i\in I}(S\backslash Z_i)=S\backslash(\cup_{i\in I}Z_i)\xrightarrow{j_I^o}S^o\xrightarrow{j} S$ 
the open complementary embeddings, where $I\subset\left\{1,\cdots,r\right\}$.
For $(M,F,W)\in \pi_{S^o}(C(MHM(S^o)))$, we define by (i) 
\begin{itemize}
\item the (bi)-filtered complex of $D_S$-modules
\begin{equation*}
j_{*Hdg}(M,F,W):=\varinjlim_{\left\{(Z_i)_{i\in[1,\ldots r]},Z_i\subset S,\cap Z_i=Z\right\},Z'_i\subset Z_i}
\Tot_{card I=\bullet}(j_{I*}^{Hdg}j_I^{o*}(M,F,W))\in\pi_S(C(MHM(S))), 
\end{equation*}
where the horizontal differential are given by, 
if $I\subset J$, $d_{IJ}:=\ad(j^*_{IJ},j^{Hdg}_{IJ*})(j_I^{o*}(M,F,W))$, 
$j_{IJ}:S^{oJ}\hookrightarrow S^{oI}$ being the open embedding, 
and $d_{IJ}=0$ if $I\notin J$,
\item the (bi)-filtered complex of $D_S$-modules
\begin{eqnarray*}
j_{!Hdg}(M,F,W):&=&\varprojlim_{\left\{(Z_i)_{i\in[1,\ldots r]},Z_i\subset S,\cap Z_i=Z\right\},Z'_i\subset Z_i}
\Tot_{card I=-\bullet}(j_{I!}^{Hdg}j_I^{o*}(M,F,W)) \\
&=&\mathbb D_S^{Hdg}j_{*Hdg}\mathbb D_S^{Hdg}(M,F,W)\in\pi_S(C(MHM(S))),
\end{eqnarray*}
where the horizontal differential are given by, 
if $I\subset J$, $d_{IJ}:=\ad(j^{Hdg}_{IJ!},j^*_{IJ})(j_I^{o*}(M,F,W))$, 
$j_{IJ}:S^{oJ}\hookrightarrow S^{oI}$ being the open embedding, and $d_{IJ}=0$ if $I\notin J$.
\end{itemize}
By definition, we have for $(M,F,W)\in\pi_{S^o}(C(MHM(S^o)))$,
$j^*j_{*Hdg}(M,F,W)=(M,F,W)$ and $j^*j_{!Hdg}(M,F,W)=(M,F,W)$.
For $(M',F,W)\in\pi_S(C(MHM(S)))$, there is, by construction,
\begin{itemize}
\item a canonical map $\ad(j^*,j_{*Hdg})(M',F,W):(M',F,W)\to j_{*Hdg}j^*(M',F,W)$,
\item a canonical map $\ad(j_{!Hdg},j^*)(M',F,W):j_{!Hdg}j^*(M',F,W)\to(M',F,W)$.
\end{itemize}
For $(M,F,W)\in\pi_{S^o}(C(MHM(S^o)))$, 
\begin{itemize}
\item we have the canonical map in $C_{\mathcal D(1,0)fil}(S)$
\begin{equation*}
T(j_{*Hdg},j_*)(M,F,W):=k\circ\ad(j^*,j_*)(j_{*Hdg}(M,F,W)):j_{*Hdg}(M,F,W)\to j_*E(M,F,W),
\end{equation*}
\item we have the canonical map in $C_{\mathcal D(1,0)fil}(S)$
\begin{eqnarray*}
T(j_!,j_{!Hdg})(M,F,W):=\mathbb D^K_SL_D(k\circ\ad(j^*,j_*)(-)): \\
j_!(M,F,W):=\mathbb D_S^KL_Dj_*E(\mathbb D_S^K(M,F,W))\to\mathbb D_S^KL_Dj_{*Hdg}\mathbb D_S^K(M,F,W)=j_{!Hdg}(M,F,W)
\end{eqnarray*}
\end{itemize}
the canonical maps.
\end{itemize}
\end{defi}

\begin{rem}\label{remHdgkey}
Let $j:S^o\hookrightarrow S$ an open embedding, with $S\in\SmVar(\mathbb C)$.
Then, for $((M,F,W),(K,W),\alpha)\in MHM(S^o)$, 
\begin{itemize}
\item the map  $T(j_!,j_{!Hdg})(M,W):j_{!w}(M,W)\to j_{!Hdg}(M,W)$ in $C_{\mathcal D0fil}(S)$
is a filtered quasi-isomorphism (apply the functor $DR^{[-]}(S^o)$ and use theorem \ref{DRK} and theorem \ref{DRf}).
\item the map $T(j_{*Hdg},j_*)(M,W):j_{*Hdg}(M,W)\to j_{*w}(M,W)$ in $C_{\mathcal D0fil}(S)$
is a filtered quasi-isomorphism (apply the functor $DR^{[-]}(S^o)$ and use theorem \ref{DRK} and theorem \ref{DRf}).
\end{itemize}
Hence, for $((M,F,W),(K,W),\alpha)\in MHM(S^o)$,  
\begin{itemize}
\item we get, for all $p,n\in\mathbb N$, monomorphisms
\begin{equation*}
F^pH^nT(j_!,j_{!Hdg})(M,F,W):F^pH^nj_{!w}(M,F,W)\hookrightarrow F^pH^nj_{!Hdg}(M,F,W)
\end{equation*}
in $\PSh_{O_S}(S)$, but $F^pH^nj_{!w}(M,F,W)\neq F^pH^nj_{!Hdg}(M,F,W)$ (it leads to different F-filtrations),
since $F^pH^nj_!(M,F)\subset H^nj_!M$ are sub $D_S$ module 
while the F-filtration on $H^nj_{!Hdg}(M,F)$ is given by Kashiwara-Malgrange $V$-filtrations,
hence satisfy a non trivial Griffith transversality property, 
thus $H^nj_!(M,F)$ and $H^nj_{!Hdg}(M,F)$ are isomorphic as $D_S$-modules but NOT isomorphic as filtered $D_S$-modules.
\item we get, for all $p,n\in\mathbb N$, monomorphisms
\begin{equation*}
T(j_{*Hdg},j_*)(M,F,W):F^pH^nj_{*Hdg}(M,F,W)\hookrightarrow F^pH^nj_*E(M,F,W)
\end{equation*}
in $\PSh_{O_S}(S)$, but $F^pH^nj_{*Hdg}(M,F,W)\neq F^pH^nj_{*w}(M,F,W)$ (it leads to different F-filtrations),
since $F^pH^nj_*E(M,F)\subset H^nj_*E(M)$ are sub $D_S$ module 
while the F-filtration on $H^nj_{*Hdg}(M,F)$ is given by Kashiwara-Malgrange $V$-filtrations,
hence satisfy a non trivial Griffith transversality property, 
thus $H^nj_*E(M,F)$ and $H^nj_{*Hdg}(M,F)$ are isomorphic as $D_S$-modules but NOT isomorphic as filtered $D_S$-modules.
\end{itemize}
\end{rem}

We make the following key definition

\begin{defi}\label{gammaHdg}
Let $S\in\SmVar(\mathbb C)$. Let $Z\subset S$ a closed subset.
Denote by $j:S\backslash Z\hookrightarrow S$ the complementary open embedding. 
\begin{itemize}
\item[(i)] We define using definition \ref{DHdgj}, the filtered Hodge support section functor
\begin{eqnarray*}
\Gamma^{Hdg}_Z:\pi_S(C(MHM(S))\to\pi_S(C(MHM(S)), \\ 
(M,F,W)\mapsto\Gamma^{Hdg}_Z(M,F,W):=\Cone(\ad(j^*,j_{*Hdg})(M,F):(M,F)\to j_{*Hdg}j^*(M,F))[-1],
\end{eqnarray*}
together we the canonical map $\gamma^{Hdg}_Z(M,F,W):\Gamma^{Hdg}_Z(M,F,W)\to (M,F,W)$.
We then have the canonical map in $C_{\mathcal D(2)fil}(S)$
\begin{eqnarray*}
T(\Gamma^{Hdg}_Z,\Gamma_Z)(M,F,W):=(I,T(j_{*Hdg},j_*)(M,F,W)):\Gamma_Z^{Hdg}(M,F,W)\to\Gamma_ZE(M,F,W)
\end{eqnarray*}
unique up to homotopy such that $\gamma_Z^{Hdg}(M,F,W)=\gamma_Z(E(M,F,W))\circ T(\Gamma^{Hdg}_Z,\Gamma_Z)(M,F,W)$.
\item[(i)'] Since $j_{*Hdg}:\pi_{S^o}(C(MHM(S^o))\to\pi_S(C(MHM(S))$ is an exact functor, 
$\Gamma^{Hdg}_Z$ induces the functor
\begin{eqnarray*}
\Gamma^{Hdg}_Z:\pi_S(D(MHM(S))\to\pi_S(D(MHM(S)), \; (M,F,W)\mapsto\Gamma^{Hdg}_Z(M,F,W)
\end{eqnarray*}
\item[(ii)] We define using definition \ref{DHdgj}, the dual filtered Hodge support section functor
\begin{eqnarray*}
\Gamma^{\vee,Hdg}_Z:\pi_S(C(MHM(S))\to\pi_S(C(MHM(S)), \\ 
(M,F,W)\mapsto\Gamma^{\vee,Hdg}_Z(M,F,W):=\Cone(\ad(j_{!Hdg},j^*)(M,F,W):j_{!Hdg},j^*(M,F,W)\to (M,F,W)),
\end{eqnarray*}
together we the canonical map $\gamma^{\vee,Hdg}_Z(M,F,W):(M,F,W)\to\Gamma_Z^{\vee,Hdg}(M,F)$.
We then have the canonical map in $C_{\mathcal D(2)fil}(S)$
\begin{eqnarray*}
T(\Gamma^{\vee,h}_Z,\Gamma^{\vee,Hdg}_Z)(M,F,W):=(I,T(j_!,j_{!Hdg})(M,F,W)):
\Gamma^{\vee,h}_Z(M,F,W)\to\Gamma^{\vee,Hdg}_Z(M,F,W)
\end{eqnarray*}
unique up to homotopy such that 
$\gamma_Z^{\vee,Hdg}(M,F)=T(\Gamma^{\vee,h}_Z,\Gamma^{\vee,Hdg}_Z)(M,F,W)\circ\gamma_Z^{\vee,h}(M,F,W)$.
\item[(ii)'] Since $j_{!Hdg}:\pi_{S^o}(C(MHM(S^o))\to \pi_S(C(MHM(S))$ is an exact functor, 
$\Gamma^{Hdg,\vee}_Z$ induces the functor
\begin{eqnarray*}
\Gamma^{\vee,Hdg}_Z:\pi_S(D(MHM(S))\to\pi_S(D(MHM(S)), \; (M,F,W)\mapsto\Gamma^{\vee,Hdg}_Z(M,F,W)
\end{eqnarray*}
\end{itemize}
\end{defi}

We now give the definition of the filtered Hodge inverse image functor :

\begin{defi}\label{inverseHdg}
\begin{itemize}
\item[(i)] Let $i:Z\hookrightarrow S$ be a closed embedding, with $Z,S\in\SmVar(\mathbb C)$.
Then, for $(M,F,W)\in\pi_S(C(MHM(S))$, we set 
\begin{equation*}
i_{Hdg}^{*mod}(M,F,W):=i^*\mathcal S_Z^{-1}\Gamma_Z^{Hdg}(M,F,W)\in\pi_Z(D(MHM(Z))
\end{equation*}
and
\begin{equation*}
i_{Hdg}^{\hat*mod}(M,F,W):=i^*\mathcal S_Z^{-1}\Gamma_Z^{\vee,Hdg}(M,F,W)\in\pi_Z(D(MHM(Z))
\end{equation*}
using the fact that $\mathcal S_Z:\pi_Z(D(MHM(Z))\to\pi_S(D(MHM_Z(S))$ is an equivalence of category 
since $\mathcal S_Z:D(MHM_Z(S))\to D(MHM_Z(S))$ is an equivalence of category by \cite{Saito}.
\item[(ii)] Let $f:X\to S$ be a morphism, with $X,S\in\SmVar(\mathbb C)$.
Consider the factorization $f:X\xrightarrow{i} X\times S\xrightarrow{p_S}S$, 
where $i$ is the graph embedding and $p_S:X\times S\to S$ is the projection.
\begin{itemize}
\item For $(M,F,W)\in\pi_S(C(MHM(S))$ we set
\begin{equation*}
f_{Hdg}^{*mod}(M,F,W):=i_{Hdg}^{*mod}p_S^{*mod[-]}(M,F,W)(d_X)[2d_X]\in\pi_X(D(MHM(X)), 
\end{equation*}
\item For $(M,F,W)\in\pi_S(C(MHM(S))$ we set
\begin{equation*}
f_{Hdg}^{\hat*mod}(M,F,W):=i_{Hdg}^{\hat*mod}p_S^{*mod[-]}(M,F,W)\pi_X(D(MHM(X)), 
\end{equation*}
\end{itemize}
If $j:S^o\hookrightarrow S$ is a closed embedding, we have (see \cite{Saito}), for $(M,F,W)\in\pi_S(C(MHM(S)))$, 
\begin{equation*}
j_{Hdg}^{*mod}(M,F,W)=j_{Hdg}^{\hat*mod}(M,F,W)=j^*(M,F,W)\in \pi_{S^o}(D(MHM(S^o)))
\end{equation*}
\item[(iii)] Let $f:X\to S$ be a morphism, with $X,S\in\SmVar(\mathbb C)$ or $X,S\in\AnSm(\mathbb C)$.
Consider the factorization $f:X\xrightarrow{i} X\times S\xrightarrow{p_S}S$, 
where $i$ is the graph embedding and $p_S:X\times S\to S$ is the projection.
\begin{itemize}
\item For $(M,F,W)\in\pi_S(C(MHM(S))$ we set
\begin{equation*}
f_{Hdg}^{*mod}(M,F,W):=\Gamma_X^{Hdg}p_S^{*mod[-]}(M,F,W)(d_X)[2d_X]\in\pi_{X\times S}(C(MHM(X\times S))), 
\end{equation*}
We have for $(M,F,W)\in\pi_S(C(MHM(S))$, the canonical map in $C_{\mathcal D(1,0)fil}(X\times S)$
\begin{eqnarray*}
T(f_{Hdg}^{*mod},f^{*mod,\Gamma})(M,F,W):f_{Hdg}^{*mod,\Gamma}(M,F,W):=\Gamma_X^{Hdg}p_S^{*mod[-]}(M,F,W) \\
\xrightarrow{T(\Gamma_X^{Hdg},\Gamma_X)(-)}\Gamma_XE(p_S^{*mod[-]}(M,F,W))=:f^{*mod[-],\Gamma}(M,F,W)
\end{eqnarray*}
\item For $(M,F,W)\in\pi_S(C(MHM(S)))$ we set
\begin{equation*}
f_{Hdg}^{\hat*mod}(M,F,W):=\Gamma_X^{\vee,Hdg}p_S^{*mod[-]}(M,F,W)\in\pi_{X\times S}(C(MHM(X\times S))), 
\end{equation*}
We have for $(M,F,W)\in\pi_S(C(MHM(S))$, the canonical map in $C_{\mathcal D(1,0)fil}(X\times S)$
\begin{eqnarray*}
T(f^{\hat*mod,\Gamma},f_{Hdg}^{\hat*mod})(M,F,W):f^{\hat*mod,\Gamma}(M,F,W):=\Gamma_X^{\vee,h}p_S^{*mod[-]}(M,F,W) \\
\xrightarrow{T(\Gamma_X^{\vee,h},\Gamma_X^{\vee,Hdg})(-)}\Gamma_X^{\vee,Hdg}p_S^{*mod[-]}(M,F,W)=:f_{Hdg}^{\hat*mod}(M,F,W)
\end{eqnarray*}
\end{itemize}
\end{itemize}
\end{defi}

\begin{defiprop}\label{TgammaHdg}
\begin{itemize}
\item[(i)] Let $g:S'\to S$ a morphism with $S',S\in\SmVar(\mathbb C)$ and $i:Z\hookrightarrow S$ a closed subset. 
Then, for $(M,F,W)\in\pi_S(C(MHM(S)))$, there is a canonical map in $\pi_S(C(MHM_{S'}(S'\times S)))$
\begin{equation*}
T^{Hdg}(g,\gamma)(M,F,W):g_{Hdg}^{*mod,\Gamma}\Gamma^{Hdg}_{Z}(M,F,W)\to\Gamma^{Hdg}_{Z\times_S S'}g_{Hdg}^{*mod,\Gamma}(M,F,W)
\end{equation*}
unique up to homotopy such that 
\begin{equation*}
\gamma^{Hdg}_{Z\times_S S'}(g_{Hdg}^{*mod,\Gamma}(M,F,W))\circ T^{Hdg}(g,\gamma)(M,F,W)=
g^{*mod,\Gamma}_{Hdg}\gamma^{Hdg}_{Z}(M,F,W).
\end{equation*}
\item[(i)'] Let $g:S'\to S$ a morphism with $S',S\in\SmVar(\mathbb C)$ and $i:Z\hookrightarrow S$ a closed subset. 
Then, for $(M,F,W)\in\pi_S(C(MHM(S)))$, 
there is a canonical isomorphism in $\pi_S(C(MHM_{S'}(S'\times S)))$
\begin{equation*}
T^{Hdg}(g,\gamma^{\vee})(M,F,W):\Gamma^{Hdg}_{Z\times_S S'}g_{Hdg}^{\hat*mod,\Gamma}(M,F,W)\xrightarrow{\sim} 
g_{Hdg}^{\hat*mod,\Gamma}\Gamma^{Hdg}_{Z}(M,F,W)
\end{equation*}
unique up to homotopy such that 
\begin{equation*}
\gamma^{\vee,Hdg}_{Z\times_S S'}(g_{Hdg}^{\hat*mod,\Gamma}(M,F,W))
\circ g^{\hat*mod,\Gamma}_{Hdg}\gamma^{\vee,Hdg}_{Z}(M,F,W)=T^{Hdg}(g,\gamma)(M,F,W).
\end{equation*}
\item[(ii)] Let $S\in\SmVar(\mathbb C)$ and $i_1:Z_1\hookrightarrow S$, $i_2:Z_2\hookrightarrow Z_1$ be closed embeddings.
Then, for $(M,F,W)\in\pi_S(C(MHM(S)))$, 
\begin{itemize}
\item there is a canonical map 
$T(Z_2/Z_1,\gamma^{Hdg})(M,F,W):\Gamma^{Hdg}_{Z_2}(M,F,W)\to\Gamma^{Hdg}_{Z_1}(M,F,W)$ 
in $\pi_S(C(MHM(S)))$ unique up to homotopy such that 
\begin{equation*}
\gamma^{Hdg}_{Z_1}(G,F)\circ T(Z_2/Z_1,\gamma^{Hdg})(G,F)=\gamma^{Hdg}_{Z_2}(G,F) 
\end{equation*}
together with a distinguish triangle in $K(\pi_S(MHM(S)))$
\begin{eqnarray*}
\Gamma^{Hdg}_{Z_2}(M,F,W)\xrightarrow{T(Z_2/Z_1,\gamma^{Hdg})(M,F,W)}\Gamma^{Hdg}_{Z_1}(M,F,W) \\
\xrightarrow{\ad(j_2^*,j^{Hdg}_{2*})(\Gamma^{Hdg}_{Z_1}(G,F))}
\Gamma^{Hdg}_{Z_1/\backslash Z_2}(G,F)\to\Gamma^{Hdg}_{Z_2}(G,F)[1]
\end{eqnarray*} 
\item there is a canonical map 
$T(Z_2/Z_1,\gamma^{\vee,Hdg})(M,F,W):\Gamma_{Z_1}^{\vee,Hdg}(M,F,W)\to\Gamma_{Z_2}^{\vee,Hdg}(M,F,W)$ 
in $\pi_S(C(MHM(S)))$ unique up to homotopy such that 
\begin{equation*}
\gamma^{\vee,Hdg}_{Z_2}(M,F,W)=T(Z_2/Z_1,\gamma^{\vee,Hdg})(M,F,W)\circ\gamma^{\vee,Hdg}_{Z_1}(M,F,W). 
\end{equation*}
together with a distinguish triangle in $K(\pi_S((MHM(S))))$
\begin{eqnarray*}
\Gamma_{Z_1\backslash Z_2}^{\vee,Hdg}(M,F,W)\xrightarrow{\ad(j^{Hdg}_{2!},j_2^*)(M,F,W)}
\Gamma_{Z_1}^{\vee,Hdg}(M,F,W) \\ 
\xrightarrow{T(Z_2/Z_1,\gamma^{\vee,Hdg})(M,F,W))}
\Gamma^{\vee,Hdg}_{Z_2}(M,F,W)\to\Gamma_{Z_2\backslash Z_1}^{\vee,Hdg}(M,F,W)[1]
\end{eqnarray*} 
\end{itemize}
\end{itemize}
\end{defiprop}

\begin{proof}
Follows from the projection case and the closed embedding case using the adjonction maps.
\end{proof}

The definitions \ref{gammaHdg} and \ref{inverseHdg} immediately extends to the non smooth case :

\begin{defi}\label{gammaHdgsing}
Let $S\in\Var(\mathbb C)$. Let $Z\subset S$ a closed subset.
Let $S=\cup_iS_i$ an open cover such that there exist closed embeddings
$i_i:S_i\hookrightarrow\tilde S_i$ with $\tilde S_i\in\SmVar(\mathbb C)$. Denote $Z_I:=Z\cap S_I$. 
Denote by $j:S\backslash Z\hookrightarrow S$ and $\tilde j_I:\tilde S_I\backslash Z_I\hookrightarrow\tilde S_I$ 
the complementary open embeddings. 
\begin{itemize}
\item[(i)] We define using definition \ref{DHdgj}, the filtered Hodge support section functor
\begin{eqnarray*}
\Gamma^{Hdg}_Z:\pi(C(MHM(S)))\to\pi(C(MHM(S))), \; 
((M_I,F,W),u_{IJ})\mapsto\Gamma^{Hdg}_Z((M_I,F,W),u_{IJ}):= \\
\Cone(\ad(j^*,j_{*Hdg})((M_I,F,W),u_{IJ}):((M_I,F,W),u_{IJ})\to(\tilde j^{Hdg}_{I*}\tilde j_I^*(M_I,F,W)),\tilde j_J(u_{IJ}))[-1],
\end{eqnarray*}
together with the canonical map $\gamma^{Hdg}_Z((M_I,F,W),u_{IJ}):\Gamma^{Hdg}_Z((M_I,F,W),u_{IJ})\to((M_I,F,W),u_{IJ})$.
We then have the canonical map in $C_{\mathcal D(2)fil}(S/(\tilde S_I)$
\begin{eqnarray*}
T(\Gamma^{Hdg}_Z,\Gamma_Z)((M_I,F,W),u_{IJ}):=(I,T(j_{*Hdg},j_*)(M,F,W)): \\
\Gamma_Z^{Hdg}((M_I,F,W),u_{IJ})\to(\Gamma_ZE(M_I,F,W),\Gamma(u_{IJ}))
\end{eqnarray*}
unique up to homotopy such that 
\begin{equation*}
\gamma_Z^{Hdg}((M_I,F,W),u_{IJ})=(\gamma_{Z_I}(E(M_I,F,W)))\circ T(\Gamma^{Hdg}_Z,\Gamma_Z)((M_I,F,W),u_{IJ}).
\end{equation*}
\item[(i)'] Since 
$\tilde j_{I*}^{Hdg}:\pi_{\tilde S_I}(C(MHM(\tilde S_I\backslash S_I)))\to\pi_{\tilde S_I}(C(MHM(\tilde S_I)))$
are exact functors, $\Gamma^{Hdg}_Z$ induces the functor
\begin{eqnarray*}
\Gamma^{Hdg}_Z:\pi_S(D(MHM(S))\to\pi_S(D(MHM(S)),  ((M_I,F,W),u_{IJ})\mapsto\Gamma^{Hdg}_Z((M_I,F,W),u_{IJ})
\end{eqnarray*}
\item[(ii)] We define using definition \ref{DHdgj}, the dual filtered Hodge support section functor
\begin{eqnarray*}
\Gamma^{\vee,Hdg}_Z:\pi(C(MHM(S)))\to\pi(C(MHM(S))), \; 
((M_I,F,W),u_{IJ})\mapsto\Gamma^{\vee,Hdg}_Z((M_I,F,W),u_{IJ}):= \\
\Cone(\ad(j_{!Hdg},j^*)((M_I,F,W),u_{IJ}):
(\tilde j_{I!}^{Hdg}\tilde j_I^*(M_I,F,W),(\tilde j_J(u^d_{IJ}))^d)\to ((M_I,F,W),u_{IJ})),
\end{eqnarray*}
together we the canonical map 
$\gamma^{\vee,Hdg}_Z((M_I,F,W),u_{IJ}):((M_I,F,W),u_{IJ})\to\Gamma_Z^{\vee,Hdg}((M_I,F,W),u_{IJ})$.
We then have the canonical map in $C_{\mathcal D(2)fil}(S/(\tilde S_I))$
\begin{eqnarray*}
T(\Gamma^{\vee,h}_Z,\Gamma^{\vee,Hdg}_Z)((M_I,F,W),u_{IJ}):=(I,T(j_!,j_{!Hdg})((M_I,F,W),u_{IJ})): \\
(\Gamma^{\vee,h}_Z(M_I,F,W),\Gamma_Z^{\vee,h}(u_{IJ}))\to\Gamma^{\vee,Hdg}_Z((M_I,F,W),u_{IJ})
\end{eqnarray*}
unique up to homotopy such that 
\begin{equation*}
\gamma_Z^{\vee,Hdg}((M_I,F,W),u_{IJ})=T(\Gamma^{\vee,h}_Z,\Gamma^{\vee,Hdg}_Z)((M_I,F,W),u_{IJ})
\circ(\gamma_{Z_I}^{\vee,h}(M_I,F,W)).
\end{equation*}
\item[(ii)'] Since 
$\tilde j^{Hdg}_{I!}:\pi_{\tilde S_I}(C(MHM(\tilde S_I\backslash S_I)))\to\pi_{\tilde S_I}(C(MHM(\tilde S_I)))$ 
are exact functors, $\Gamma^{Hdg,\vee}_Z$ induces the functor
\begin{eqnarray*}
\Gamma^{\vee,Hdg}_Z:\pi_S(D(MHM(S))\to\pi_S(D(MHM(S)), \; ((M_I,F,W),u_{IJ})\mapsto\Gamma^{\vee,Hdg}_Z((M_I,F,W),u_{IJ})
\end{eqnarray*}
\end{itemize}
\end{defi}

\begin{defi}\label{inverseHdgsing}
Let $f:X\to S$ a morphism with $X,S\in Var(\mathbb C)$.
Assume there exist a factorization $f:X\xrightarrow{l}Y\times S\xrightarrow{p_S}S$ 
with $Y\in\SmVar(\mathbb C)$, $l$ a closed embedding and $p_S$ the projection.
Let $S=\cup_{i\in I}$ an open cover such that there exist closed embeddings
$i:S_i\hookrightarrow\tilde S_i$ with $\tilde S_i\in\SmVar(\mathbb C)$. 
Denote $X_I:=f^{-1}(S_I)$. We have then $X=\cup_{i\in I}X_i$ and the commutative diagrams
\begin{equation*}
\xymatrix{f:X_I\ar[r]^{l_I}\ar[rd] & Y\times S_I\ar[r]^{p_{S_I}}\ar[d]^{i_I':=(I\times i_I)} & S_I\ar[d]^{i_I} \\ 
\, & Y\times\tilde S_I\ar[r]^{p_{\tilde S_I}=:\tilde f_I} & \tilde S_I} 
\end{equation*}
\begin{itemize}
\item[(i)] For $((M_I,F,W),u_{IJ})\in\pi_S(C(MHM(S))$ we set (see definition \ref{gammaHdgsing} for $l$)
\begin{equation*}
f_{Hdg}^{*mod}((M_I,F,W),u_{IJ}):=\Gamma_X^{Hdg}(p_{\tilde S_I}^{*mod[-]}(M_I,F,W),u_{IJ})(d_Y)[2d_Y]\in\pi_X(C(MHM(X))), 
\end{equation*}
We have for $((M_I,F,W),u_{IJ})\in\pi_S(C(MHM(S))$, the canonical map in $C_{\mathcal D(1,0)fil}(X/(Y\times\tilde S_I))$
\begin{eqnarray*}
T(f_{Hdg}^{*mod},f^{*mod,\Gamma})((M_I,F,W),u_{IJ}):
f_{Hdg}^{*mod}((M_I,F,W),u_{IJ}):=\Gamma_X^{Hdg}(p_{\tilde S_I}^{*mod[-]}(M_I,F,W),p_{\tilde S_I}^{*mod[-]}u_{IJ}) \\
\xrightarrow{T(\Gamma_X^{Hdg},\Gamma_X)(-)}
(\Gamma_XE(p_{\tilde S_I}^{*mod[-]}(M_I,F,W)),\tilde f_I^{*mod[-]}u_{IJ}))=:f^{*mod[-],\Gamma}(M,F,W)
\end{eqnarray*}
\item[(ii)] For $((M_I,F,W),u_{IJ})\in\pi_S(C(MHM(S)))$ we set (see definition \ref{gammaHdgsing} for $l$)
\begin{equation*}
f_{Hdg}^{\hat*mod}(M,F,W):=\Gamma_X^{\vee,Hdg}(p_{\tilde S_I}^{*mod[-]}(M_I,F,W),p_{\tilde S_I}^{*mod[-]}u_{IJ})
\in\pi_X(C(MHM(X)), 
\end{equation*}
We have for $(M,F,W)\in\pi_S(C(MHM(S))$, the canonical map in $C_{\mathcal D(1,0)fil}(X/(Y\times\tilde S_I))$
\begin{eqnarray*}
T(f^{\hat*mod,\Gamma},f_{Hdg}^{\hat*mod[-]})((M_I,F,W),u_{IJ}):
f^{\hat*mod[-],\Gamma}(M,F,W):=\mathbb D_S^Kf^{*mod[-],\Gamma}\mathbb D_X^K((M_I,F,W),u_{IJ}) \\
\xrightarrow{\mathbb D_S^KT(\Gamma_X^{Hdg},\Gamma_X)(-)}
\Gamma_X^{\vee,Hdg}(p_S^{*mod[-]}(M_I,F,W),p_{\tilde S_I}^{*mod[-]}u_{IJ})=:f_{Hdg}^{\hat*mod[-]}(M,F,W)
\end{eqnarray*}
\end{itemize}
\end{defi}

\begin{prop}\label{compDmodDRHdg}
Let $f_1:X\to Y$ and $f_2:Y\to S$ two morphism with $X,Y,S\in\QPVar(\mathbb C)$ or with $X,Y,S\in\SmVar(\mathbb C)$. 
\begin{itemize}
\item[(i)]Let $(M,F,W)\in\pi_S(C(MHM(S)))$. Then, 
\begin{equation*}
(f_2\circ f_1)_{Hdg}^{*mod}(M,F)=f_{1Hdg}^{*mod}f_{2Hdg}^{*mod}(M,F)\in\pi_X(D(MHM(X))).
\end{equation*}
\item[(ii)]Let $(M,F,W)\in\pi_S(C(MHM(S)))$. Then,
\begin{equation*}
(f_2\circ f_1)_{Hdg}^{\hat*mod}(M,F)=f_{1Hdg}^{\hat*mod}f_{2Hdg}^{\hat*mod}(M,F)\in\pi_X(D(MHM(X)))
\end{equation*}
\end{itemize}
\end{prop}

\begin{proof}
\noindent(i):Follows from the unicity of the functor $j_{*Hdg}$.

\noindent(ii):Follows from the unicity of the functor $j_{!Hdg}$.
\end{proof}

\begin{itemize}
\item Let $f:X\to S$ a morphism with $S,X\in\SmVar(\mathbb C)$. 
Let $f:X\xrightarrow{j}\bar X\xrightarrow{\bar f}S$ a compactification of $f$ with $\bar X\in\SmVar(\mathbb C)$
and $j$ the open embedding. We set, for $(M,F,W)\in\pi_X(D(MHM(X))$, using definition \ref{DHdgj}
\begin{equation*}
\int_f^{Hdg}(M,F,W):=\int_{\bar f}^{FDR}j_{*Hdg}(M,F,W)\in D_{\mathcal D(1,0)fil}(S)
\end{equation*}
and
\begin{equation*}
\int_{f!}^{Hdg}(M,F,W):=\int_{\bar f}^{FDR}j_{!Hdg}(M,F,W)\in D_{\mathcal D(1,0)fil}(S)
\end{equation*}
\item Let $f:X\to S$ a morphism with $S,X\in\QPVar(\mathbb C)$.  
Consider a factorization $f:X\xrightarrow{l}Y\times S\xrightarrow{p}S$
with $Y\in\SmVar(\mathbb C)$, $l$ a closed embedding and $p_S$ the projection.
Let $\bar Y\in\PSmVar(\mathbb C)$ a smooth compactification of $Y$ with $j:Y\hookrightarrow\bar Y$ the open embedding.
Then $\bar f:\bar X\xrightarrow{\bar l}\bar Y\times_S\xrightarrow{\bar p}S$ is a compactification of $f$,
with $\bar X\subset\bar Y\times S$ the closure of $X$ and $\bar l$ the closed embedding.
Let $S=\cup_i S_i$ an open affine cover and 
$i_i:S_i\hookrightarrow\tilde S_i$ closed embedding with $\tilde S_i\in\SmVar(\mathbb C)$
For $((M_I,F,W),u_{IJ})\in \pi_X(C(MHM(X)))\subset C_{\mathcal D(1,0)fil,rh}(X/(Y\times\tilde S_I))$, 
we set using definition \ref{DHdgj}
\begin{eqnarray*}
\int_f^{Hdg}((M_I,F,W),u_{IJ}):=\int_{\bar f}^{FDR}((I\times j)_{*Hdg}((M_I,F,W),u_{IJ}))
\in D_{\mathcal D(1,0)fil}(S/(\tilde S_I))
\end{eqnarray*}
and
\begin{eqnarray*}
\int_{f!}^{Hdg}((M_I,F,W),u_{IJ}):=\int_{\bar f}^{FDR}((I\times j)_{!Hdg}((M_I,F,W),u_{IJ}))
\in D_{\mathcal D(1,0)fil}(S/(\tilde S_I))
\end{eqnarray*}
\end{itemize}

From the D-module case on algebraic varieties and the constructible sheaves case on CW complexes, we get :

\begin{defi}\label{falpha}
\begin{itemize}
\item[(i)] Let $f:X\to S$ a morphism with $S,X\in\SmVar(\mathbb C)$. 
Let $f:X\xrightarrow{j}\bar X\xrightarrow{\bar f}S$ a compactification of $f$ with $\bar X\in\SmVar(\mathbb C)$
and $j$ the open embedding. Let 
\begin{equation*}
\alpha:(K,W)\otimes\mathbb C_{X^{an}}\to DR(X)((M,W)^{an}) 
\end{equation*}
a morphism in $D_{fil}(X^{an})$, with $(M,F,W)\in \pi_X(C(MHM(X)))$ and $(K,W)\in D_{fil}(X^{an})$. 
We then consider the maps in $D_{fil}(S^{an})$
\begin{eqnarray*}
f_*\alpha:Rf_{*w}(K,W)\otimes\mathbb C_{S^{an}}:=
R\bar f_*Rj_{*w}(K,W)\otimes\mathbb C_{S^{an}} \\
\xrightarrow{R\bar f_*j_*\alpha}R\bar f_*Rj_{*w}DR(X)((M,W)^{an}) 
\xrightarrow{T^w(j,\otimes)(-)^{-1}}R\bar f_*DR(\bar X)(j_{*Hdg}(M,W)^{an}) \\
\xrightarrow{T(\bar f,DR)(-)^{-1}}
DR(S)((\int_{\bar f}(j_{*Hdg}(M,W))^{an})=DR(S)((\int^{Hdg}_f(M,W))^{an})
\end{eqnarray*}
and
\begin{eqnarray*}
f_!\alpha:Rf_{!w}(K,W)\otimes\mathbb C_{S^{an}}:=
R\bar f_*Rj_{!w}(K,W)\otimes\mathbb C_{S^{an}} \\
\xrightarrow{R\bar f_*j_!\alpha}R\bar f_*Rj_{!w}DR(X)((M,W)^{an}) 
\xrightarrow{\mathbb DT^w(j,\otimes)(-)}R\bar f_*DR(\bar X)(j_{!Hdg}(M,W)^{an}) \\
\xrightarrow{T(\bar f,DR)(-)^{-1}}
DR(S)((\int_{\bar f}(j_{!Hdg}(M,W))^{an})=DR(S)((\int^{Hdg}_{f!}(M,W))^{an}),
\end{eqnarray*}
see definition \ref{fw} and definition \ref{DHdgj} .
\item[(ii)] Let $f:X\to S$ a morphism with $S,X\in\QPVar(\mathbb C)$.  
Consider a factorization $f:X\xrightarrow{l}Y\times S\xrightarrow{p}S$
with $Y\in\SmVar(\mathbb C)$, $l$ a closed embedding and $p_S$ the projection.
Let $\bar Y\in\PSmVar(\mathbb C)$ a smooth compactification of $Y$ with $j:Y\hookrightarrow\bar Y$ the open embedding.
Then $\bar f:\bar X\xrightarrow{\bar l}\bar Y\times_S\xrightarrow{\bar p}S$ is a compactification of $f$,
with $\bar X\subset\bar Y\times S$ the closure of $X$ and $\bar l$ the closed embedding.
Let $S=\cup_i S_i$ an open affine cover and 
$i_i:S_i\hookrightarrow\tilde S_i$ closed embedding with $\tilde S_i\in\SmVar(\mathbb C)$. Let 
\begin{equation*}
\alpha:T(X/(Y\times\tilde S_I))((K,W)\otimes\mathbb C_{X^{an}})\to DR(X)((M_I,W)^{an},u^{an}_{IJ}) 
\end{equation*}
a morphism in in $D_{fil}(X^{an}/({Y\times\tilde S_I}^{an}))$, with
$((M_I,W),u_{IJ})\in \pi_X(C(MHM(X)))\subset C_{\mathcal D0fil,rh}(X/(Y\times\tilde S_I))$ 
and $(K,W)\in D_{fil}(X^{an})$. 
We then consider the maps in $D_{fil}(S^{an}/(\tilde S_I^{an}))$
\begin{eqnarray*}
f_*\alpha=f_*(\alpha):T(S/\tilde S_I)(Rf_{*w}(K,W))\otimes\mathbb C_{S^{an}} \\
\xrightarrow{:=}T(S/\tilde S_I)(R\bar p_*(I\times j)_{*w}(K,W))\otimes\mathbb C_{S^{an}}
\xrightarrow{=}R\bar p_*(I\times j)_{*w}T(X/(Y\times\tilde S_I))((K,W)\otimes\mathbb C_{X^{an}}) \\
\xrightarrow{Rp_*\alpha}R\bar p_*(I\times j)_{*w}DR(X)(((M_I,W),u_{IJ})^{an}) \\ 
\xrightarrow{(T^w(I\times j,\otimes)(-))}R\bar p_*DR(X)((I\times j)_{*Hdg}((M_I,W),u_{IJ})^{an}) \\
\xrightarrow{T(\bar f,DR)(-)}
DR(S)((\int_{\bar f}(I\times j)_{*Hdg}((M_I,W),u_{IJ}))^{an})=DR(S)((\int^{Hdg}_f((M_I,W),u_{IJ}))^{an})
\end{eqnarray*}
and
\begin{eqnarray*}
f_!\alpha=f_!(\alpha):T(S/\tilde S_I)(Rf_{!w}(K,W))\otimes\mathbb C_{S^{an}} \\
\xrightarrow{:=}T(S/\tilde S_I)(R\bar p_*(I\times j)_{*w}(K,W))\otimes\mathbb C_{S^{an}}
\xrightarrow{=}R\bar p_*(I\times j)_{!w}T(X/(Y\times\tilde S_I))((K,W)\otimes\mathbb C_{X^{an}}) \\
\xrightarrow{R\bar p_*\mathbb D^vR(I\times j)_*\mathbb D^v\alpha}
R\bar p_*(I\times j)_{!w}DR(X)(((M_I,W),u_{IJ})^{an}) \\
\xrightarrow{T(D,DR)(-)\circ(\mathbb DT^w(I\times j,\otimes)(-))\circ T(D,DR)(-)}
R\bar p_*DR(X)((I\times j)_{!Hdg}((M_I,W),u_{IJ})^{an}) \\
\xrightarrow{T(\bar f,DR)(-)}
DR(S)((\int_{\bar f}(I\times j)_{!Hdg}((M_I,W),u_{IJ}))^{an})=DR(S)((\int^{Hdg}_{f!}((M_I,W),u_{IJ}))^{an}), 
\end{eqnarray*}
see definition \ref{fw} and definition \ref{DHdgj} .
\item[(iii)] Let $l:S^o\hookrightarrow S$ an open embedding with $S\in\Var(\mathbb C)$ and denote $Z=S\backslash S^o$.
Let $S=\cup_i S_i$ an open affine cover and 
$i_i:S_i\hookrightarrow\tilde S_i$ closed embedding with $\tilde S_i\in\SmVar(\mathbb C)$. 
Let $l_I:\tilde S^o_I\hookrightarrow\tilde S_I$ open embeddings such that $\tilde S^o_I\cap S=S^o\cap S_I$. Let 
\begin{equation*}
\alpha:T(S/(\tilde S_I))((K,W)\otimes\mathbb C_{S^{an}})\to DR(S)(((M_I,W),u_{IJ})^{an}) 
\end{equation*}
a morphism in $D_{fil}(S^{an}/(\tilde S_I^{an}))$, with
$((M_I,W),u_{IJ})\in \pi_S(CDRM(S))\subset C_{\mathcal D0fil,rh}(S/(\tilde S_I))$ and $(K,W)\in D_{fil}(S^{an})$. 
We then consider the maps in $D_{fil}(S^{an}/(\tilde S_I^{an}))$
\begin{eqnarray*}
\Gamma_Z(\alpha):T(S/(\tilde S_I))(\Gamma^w_Z(K,W)\otimes\mathbb C_{S^{an}})
\xrightarrow{=}\Gamma^w_ZT(S/(\tilde S_I))((K,W)\otimes\mathbb C_{S^{an}}) \\
\xrightarrow{R\Gamma_Z\alpha}\Gamma^w_ZDR(S)(((M_I,W),u_{IJ})^{an}) \\
\xrightarrow{(T(\gamma_Z,DR)((M_I,W),u_{IJ})^{-1}:=((I,T^w(l_I,\otimes)(M_I,W))\circ (T^w(an,\otimes)(M_I,W)))^{-1}} 
DR(S)((\Gamma^{Hdg}_Z((M_I,W),u_{IJ}))^{an})
\end{eqnarray*}
and
\begin{eqnarray*}
\Gamma^{\vee}_Z(\alpha):T(S/(\tilde S_I))(\Gamma^{\vee,w}_Z(K,W)\otimes\mathbb C_{S^{an}})
\xrightarrow{=}\Gamma^{\vee}_ZT(S/(\tilde S_I))((K,W)\otimes\mathbb C_{S^{an}}) \\
\xrightarrow{\Gamma^{\vee}_Z\alpha}\Gamma^{\vee,w}_ZDR(S)(((M_I,W),u_{IJ})^{an}) \\
\xrightarrow{T(\gamma^{\vee}_Z,DR)((M_I,W),u_{IJ}):=
(\mathbb D(I,T^w(l_I,\otimes)(\mathbb D(M_I,W)))\circ (\mathbb DT^w(an,\otimes)(\mathbb D(M_I,W))))} 
DR(S)((\Gamma^{\vee,w}_Z((M_I,W),u_{IJ}))^{an}),
\end{eqnarray*}
see definition \ref{gammaw} and definition \ref{gammaHdgsing}.
\item[(iv)] Let $f:X\to S$ a morphism with $S,X\in\QPVar(\mathbb C)$. 
Consider a factorization $f:X\hookrightarrow Y\times S\xrightarrow{p}S$ with $Y\in\SmVar(\mathbb C)$. 
Let $S=\cup_i S_i$ an open affine cover and 
$i_i:S_i\hookrightarrow\tilde S_i$ closed embedding with $\tilde S_i\in\SmVar(\mathbb C)$. Let 
\begin{equation*}
\alpha:T(S/(\tilde S_I))((K,W)\otimes\mathbb C_{S^{an}})\to DR(S)(((M_I,W),u_{IJ})^{an}) 
\end{equation*}
a morphism in $D_{fil}(S^{an}/(\tilde S_I^{an}))$, with
$((M_I,W),u_{IJ})\in C(DRM(S))\subset C_{\mathcal D0fil,rh}(S/(\tilde S_I))$ and $(K,W)\in D_{fil}(S^{an})$. 
We then consider, see (iii), the maps in $D_{fil}(X^{an}/(Y\times\tilde S_{I,\mathbb C})^{an})$
\begin{eqnarray*}
f^!\alpha=f^!(\alpha):T(X/(Y\times\tilde S_I))(f^{!w}(K,W)\otimes\mathbb C_{X^{an}})
\xrightarrow{:=}T(X/(Y\times\tilde S_I))(\Gamma^w_Xp^*(K,W))\otimes\mathbb C_{X^{an}}) \\
\xrightarrow{=}
(\Gamma^w_Xp_{\tilde S_I}^*T(S/\tilde S_I)((K,W)_I\otimes\mathbb C_{S^{an}}),\Gamma^w_Xp^*T(D_{IJ})(-)) \\
\xrightarrow{R\Gamma_Xp^*\alpha}
(\Gamma^w_Xp_{\tilde S_I}^*DR(S)((M_I,W),u_{IJ})_I,\Gamma^w_Xp^*DR(u_{IJ})) \\
\xrightarrow{T^!(f,DR)(-):=T(\gamma_X,DR)(-)^{-1}\circ T^!(p,DR)(-)}DR(X)(f^{*mod}_{Hdg}((M_I,W),u_{IJ})^{an}) 
\end{eqnarray*}
and
\begin{eqnarray*}
f^*\alpha=f^*(\alpha):T(X/(Y\times\tilde S_I))(f^{*w}(K,W)\otimes\mathbb C_{X^{an}})
\xrightarrow{:=}T(X/(Y\times\tilde S_I))(\Gamma^{\vee,w}_Xp^*(K,W)\otimes\mathbb C_{X^{an}}) \\
\xrightarrow{=}
(\Gamma^{\vee,w}_Xp_{\tilde S_I}^*T(S/(\tilde S_I))((K,W)\otimes\mathbb C_{S^{an}})_I,
\Gamma_X^{\vee,w}p^*T(D_{IJ})(-)) \\
\xrightarrow{\Gamma^{\vee}_Xp^*\alpha}
(\Gamma^{\vee,w}_Xp_{\tilde S_I}^*DR(S)((M_I,W),u_{IJ})_I,\Gamma_X^{\vee,w}p^*DR(u_{IJ})) \\
\xrightarrow{T^*(f,DR)(-):=T(\gamma_X^{\vee},DR)(-)\circ T^*(p,DR)(-)}DR(X)(f_{Hdg}^{\hat*mod}((M_I,W),u_{IJ})^{an}) 
\end{eqnarray*}
\item[(v)] Let $S\in\Var(\mathbb C)$.Let $S=\cup_i S_i$ an open affine cover and 
$i_i:S_i\hookrightarrow\tilde S_i$ closed embedding with $\tilde S_i\in\SmVar(\mathbb C)$. Let 
\begin{eqnarray*}
\alpha:T(S/(\tilde S_I))((K,W)\otimes\mathbb C_{S^{an}})\to DR(S)(((M_I,W),u_{IJ})^{an}), \\
\alpha':T(S/(\tilde S_I))((K',W)\otimes\mathbb C_{S^{an}})\to DR(S)(((M'_I,W),v_{IJ})^{an}) 
\end{eqnarray*}
two morphism in $D_{fil}(S^{an}/(\tilde S_I^{an}))$, with
$((M_I,W),u_{IJ}),((M'_I,W),u_{IJ})\in \pi_S(C(MHM(S)))\subset C_{\mathcal D0fil,rh}(S/(\tilde S_I))$ and 
$(K,W),(K',W)\in D_{fil}(S^{an})$. We then consider the map
\begin{eqnarray*}
\alpha\otimes\alpha':T(S/(\tilde S_I))((K,W)\otimes(K',W)\otimes\mathbb C_{S^{an}}) \\
\xrightarrow{=}T(S/(\tilde S_I))((K,W)\otimes\mathbb C_{S^{an}})
\otimes T(S/(\tilde S_I))((K',W)\otimes\mathbb C_{S^{an}}) \\
\xrightarrow{\alpha\otimes\alpha'}DR(S)(((M_I,W),u_{IJ})^{an})\otimes DR(S)(((M'_I,W),v_{IJ})^{an}) \\
\xrightarrow{T(\otimes,DR)(-,-)}DR(S)((((M_I,W),u_{IJ})\otimes^L_{O_S}((M'_I,W),v_{IJ}))^{an}) \\
\xrightarrow{=}DR(S)((((M_I,W),u_{IJ})\otimes^{Hdg}_{O_S}((M'_I,W),v_{IJ}))^{an})
\end{eqnarray*}
in $D_{fil}(S^{an}/(\tilde S_I^{an}))$.
\end{itemize}
\end{defi}

\begin{defi}\label{DHdgjalpha}
Let $k\subset\mathbb C$ a subfield. Let $S\in\SmVar(\mathbb C)$. Let $j:S^o\hookrightarrow S$ an open embedding.
Let $Z:=S\backslash S^o=V(\mathcal I)\subset S$ an the closed complementary subset, 
$\mathcal I\subset O_S$ being an ideal subsheaf. 
Taking generators $\mathcal I=(s_1,\ldots,s_r)$, we get $Z=V(s_1,\ldots,s_r)=\cap^r_{i=1}Z_i\subset S$ with 
$Z_i=V(s_i)\subset S$, $s_i\in\Gamma(S,\mathcal L_i)$ and $L_i$ a line bundle. 
Note that $Z$ is an arbitrary closed subset, $d_Z\geq d_X-r$ needing not be a complete intersection. 
Denote by $j_I:S^{o,I}:=\cap_{i\in I}(S\backslash Z_i)=S\backslash(\cup_{i\in I}Z_i)\xrightarrow{j_I^o}S^o\xrightarrow{j} S$ 
the open embeddings.
Let $(M,F,W)\in MHM(S^o))$. We then define, using definition \ref{DHdgj} and definition \ref{jw}
\begin{itemize}
\item the canonical extension 
\begin{eqnarray*}
j_{*Hdg}((M,F,W),(K,W),\alpha):=(j_{*Hdg}(M,F,W),j_{*w}(K,W),j_*\alpha) \\
:=\Tot((j_{I*Hdg}j_I^*(M,F,W),j_{I*w}j_I^*(K,W),j_{I*}\alpha))\in MHM(S), 
\end{eqnarray*}
so that $j^*(j_{*Hdg}((M,F,W),(K,W),\alpha))=((M,F,W),(K,W),\alpha)$,
\item the canonical extension 
\begin{eqnarray*}
j_{!Hdg}((M,F,W),(K,W),\alpha):=(j_{!Hdg}(M,F,W),j_{!w}(K,W),j_!\alpha) \\
:=\Tot((j_{I!Hdg}j_I^*(M,F,W),j_{I!w}j_I^*(K,W),j_{I!}\alpha))\in MHM(S),  
\end{eqnarray*}
so that $j^*(j_{!Hdg}((M,F,W),(K,W),\alpha))=((M,F,W),(K,W),\alpha)$.
\end{itemize}
Moreover for $((M',F,W),(K',W),\alpha')\in MHM(S)$,
\begin{itemize}
\item there is a canonical map in $MHM(S)$
\begin{equation*}
\ad(j^*,j_{*Hdg})((M',F,W),(K',W),\alpha'):((M',F,W),(K',W),\alpha')\to j_{*Hdg}j^*((M',F,W),(K',W),\alpha'), 
\end{equation*}
\item there is a canonical map in $MHM(S)$
\begin{equation*}
\ad(j_{!Hdg},j^*)((M',F,W),(K',W),\alpha'):j_{!Hdg}j^*((M',F,W),(K',W),\alpha')\to((M',F,W),(K',W),\alpha').
\end{equation*}
\end{itemize}
\end{defi}

\begin{defi}\label{gammaHdgalpha}
Let $S\in\SmVar(\mathbb C)$. Let $Z\subset S$ a closed subset.
Denote by $j:S\backslash Z\hookrightarrow S$ the complementary open embedding. 
\begin{itemize}
\item[(i)] We define using definition \ref{gammaHdg}, definition \ref{gammaw} and definition \ref{falpha}(iii), 
the filtered Hodge support section functor
\begin{eqnarray*}
\Gamma^{Hdg}_Z:C(MHM(S))\to C(MHM(S)), \; \; ((M,F,W),(K,W),\alpha)\mapsto \\
\Gamma^{Hdg}_Z((M,F,W),(K,W),\alpha):=(\Gamma_Z^{Hdg}(M,F,W),\Gamma_Z^w(K,W),\Gamma_Z(\alpha)) \\
=\Cone(\ad(j^*,j_{*Hdg})(-):j_{*Hdg},j^*((M,F,W),(K,W),\alpha)\to((M,F,W),(K,W),\alpha)[-1]
\end{eqnarray*}
see definition \ref{DHdgjalpha} for the last equality, together we the canonical map 
\begin{eqnarray*}
\gamma^{Hdg}_Z((M,F,W),(K,W),\alpha):\Gamma^{Hdg}_Z((M,F,W),(K,W),\alpha)\to((M,F,W),(K,W),\alpha).
\end{eqnarray*}
\item[(i)'] Since $j_{*Hdg}:C(MHM(S^o))\to C(MHM(S))$ is an exact functor, 
$\Gamma^{Hdg}_Z$ induces the functor
\begin{eqnarray*}
\Gamma^{Hdg}_Z:D(MHM(S))\to D(MHM(S)), \\ 
((M,F,W),(K,W),\alpha)\mapsto\Gamma^{Hdg}_Z((M,F,W),(K,W),\alpha)
\end{eqnarray*}
\item[(ii)] We define using definition \ref{gammaHdg}, definition \ref{gammaw} and definition \ref{falpha}(iii) 
the dual filtered Hodge support section functor
\begin{eqnarray*}
\Gamma^{\vee,Hdg}_Z:C(MHM(S))\to C(MHM(S)), \; \; ((M,F,W),(K,W),\alpha)\mapsto \\
\Gamma^{\vee,Hdg}_Z((M,F,W),(K,W),\alpha):=(\Gamma_Z^{\vee,Hdg}(M,F,W),\Gamma_Z^{\vee,w}(K,W),\Gamma_Z^{\vee}(\alpha)) \\
=\Cone(\ad(j_{!Hdg},j^*)(-):j_{!Hdg},j^*((M,F,W),(K,W),\alpha) \to ((M,F,W),(K,W),\alpha))
\end{eqnarray*}
see definition \ref{DHdgjalpha} for the last equality, together we the canonical map 
\begin{eqnarray*}
\gamma^{\vee,Hdg}_Z((M,F,W),(K,W),\alpha):((M,F,W),(K,W),\alpha)\to\Gamma_Z^{\vee,Hdg}((M,F,W),(K,W),\alpha).
\end{eqnarray*}
\item[(ii)'] Since $j_{!Hdg}:C(MHM(S^o))\to C(MHM(S))$ is an exact functor, 
$\Gamma^{Hdg,\vee}_Z$ induces the functor
\begin{eqnarray*}
\Gamma^{\vee,Hdg}_Z:D(MHM(S))\to D(MHM(S)), \\ 
((M,F,W),(K,W),\alpha)\mapsto\Gamma^{\vee,Hdg}_Z((M,F,W),(K,W),\alpha)
\end{eqnarray*}
\end{itemize}
\end{defi}

In the singular case it gives :

\begin{defi}\label{gammaHdgsingalpha}
Let $S\in\Var(\mathbb C)$. Let $Z\subset S$ a closed subset.
Let $S=\cup_{i=1}^sS_i$ an open cover such that there exist closed embeddings
$i_i:S_i\hookrightarrow\tilde S_i$ with $\tilde S_i\in\SmVar(\mathbb C)$. Denote $Z_I:=Z\cap S_I$. 
Denote by $n:S\backslash Z\hookrightarrow S$ and $\tilde n_I:\tilde S_I\backslash Z_I\hookrightarrow\tilde S_I$ 
the complementary open embeddings. 
\begin{itemize}
\item[(i)] We define using definition \ref{gammaHdgsing}, definition \ref{gammaw} and definition \ref{falpha}(iii)
the filtered Hodge support section functor
\begin{eqnarray*}
\Gamma^{Hdg}_Z:C(MHM(S))\to C(MHM(S)), \\ 
(((M_I,F,W),u_{IJ}),(K,W),\alpha)\mapsto\Gamma^{Hdg}_Z(((M_I,F,W),u_{IJ}),(K,W),\alpha):= \\
:=(\Gamma_Z^{Hdg}((M_I,F,W),u_{IJ}),\Gamma_Z^w(K,W),\Gamma_Z(\alpha))
\end{eqnarray*}
together with the canonical map
\begin{eqnarray*}
\gamma^{Hdg}_Z(((M_I,F,W),u_{IJ}),(K,W),\alpha): \\
\Gamma^{Hdg}_Z(((M_I,F,W),u_{IJ}),(K,W),\alpha)\to(((M_I,F,W),u_{IJ}),(K,W),\alpha).
\end{eqnarray*}
\item[(i)'] By exactness of $\Gamma_Z^{Hdg}$ and $\Gamma_Z^w$ it induces the functor
\begin{eqnarray*}
\Gamma^{Hdg}_Z: D(MHM(S))\to D(MHM(S)), \\  
(((M_I,F,W),u_{IJ}),(K,W),\alpha)\mapsto\Gamma^{Hdg}_Z(((M_I,F,W),u_{IJ}),(K,W),\alpha)
\end{eqnarray*}
\item[(ii)] We define using definition \ref{gammaHdgsing}, definition \ref{gammaw} and definition \ref{falpha}(iii)
the dual filtered Hodge support section functor
\begin{eqnarray*}
\Gamma^{\vee,Hdg}_Z:C(MHM(S))\to C(MHM(S)), \; \; 
(((M_I,F,W),u_{IJ}),(K,W),\alpha)\mapsto \\
\Gamma^{\vee,Hdg}_Z(((M_I,F,W),u_{IJ}),(K,W),\alpha):=
(\Gamma_Z^{\vee,Hdg}((M_I,F,W),u_{IJ}),\Gamma_Z^{\vee,w}(K,W),\Gamma^{\vee}_Z(\alpha)),
\end{eqnarray*}
together we the canonical map 
\begin{eqnarray*}
\gamma^{\vee,Hdg}_Z(((M_I,F,W),u_{IJ}),(K,W),\alpha): \\
(((M_I,F,W),u_{IJ}),(K,W),\alpha)\to\Gamma_Z^{\vee,Hdg}(((M_I,F,W),u_{IJ}),(K,W),\alpha).
\end{eqnarray*}
\item[(ii)'] By exactness of $\Gamma_Z^{\vee,Hdg}$ and $\Gamma_Z^{\vee,w}$, it induces the functor
\begin{eqnarray*}
\Gamma^{\vee,Hdg}_Z:D(MHM(S))\to D(MHM(S)), \\
(((M_I,F,W),u_{IJ}),(K,W),\alpha)\mapsto\Gamma^{\vee,Hdg}_Z(((M_I,F,W),u_{IJ}),(K,W),\alpha) \\
:=(\Gamma_Z^{\vee,Hdg}((M_I,F,W),u_{IJ}),\Gamma_Z^{\vee,w}(K,W),\Gamma^{\vee}_Z(\alpha))
\end{eqnarray*}
\end{itemize}
\end{defi}

For $X\in\SmVar(\mathbb C)$, we have, by definition 
\begin{equation*}
\mathbb Z^{Hdg}_X:=a_X^{*Hdg}\mathbb Z^{Hdg}_{\pt}:=((O_X,F_b)[d_X],\mathbb Z_X,\alpha(X))\in D(MHM(X)),
\end{equation*}
with $\alpha(X):\mathbb C_X\hookrightarrow(0\to O_X\to\Omega_X\to\cdots K_X)$. If $X\in\SmVar(\mathbb C)$, 
\begin{equation*}
\mathbb Z^{Hdg}_X:=a_X^{*Hdg}\mathbb Z^{Hdg}_{\pt}:=((O_X,F_b)[d_X],\mathbb Z_{X^{an}},\alpha(X^{an}))\in D(MHM(X)).
\end{equation*}
Let $X\in\Var(\mathbb C)$ non smooth. 
Take an open cover $X=\cup_{i=1}^lX_i$ such that there exists closed
embeddings $i_i:X_i\hookrightarrow\tilde X_i$ with $\tilde X_i\in\SmVar(\mathbb C)$. 
Then, by definition  
\begin{eqnarray*}
\mathbb Z^{Hdg}_X:=a_X^*\mathbb Z^{Hdg}_{\pt}:=
((\Gamma^{\vee,Hdg}_{X_I}(O_{\tilde X_I},F_b)[d_{\tilde X_I}],o_{\tilde S_J/\tilde S_I}),(\mathbb Z_{X^{an}},W),
\alpha(X/\tilde X_I))\in D(MHM(X)).
\end{eqnarray*}
with
\begin{equation*}
\alpha(X/\tilde X_I):(\Gamma^{\vee,w}_{X_I}\alpha(\tilde X_I)):
T(X/(\tilde X_I))(\mathbb Z_{X^{an}}):=(i_{I*}\mathbb Z_{X_I^{an}},I)\to 
DR(X)^{[-]}((\Gamma^{\vee,Hdg}_{X_I}(O_{\tilde X_I})[d_{\tilde X_I}]),o_{\tilde S_J/\tilde S_I})
\end{equation*}

We have from \cite{Saito} the following proposition which shows us how to construct inductively mixed Hodge modules,
as we do for perverse sheaves :

\begin{prop}\label{Saprop}
\begin{itemize}
\item[(i)]Let $S\in\AnSm(\mathbb C)$. Let $D=V(s)\subset S$ a (Cartier) divisor, where $s\in\Gamma(S,L)$ 
is a section of the line bundle $L=L_D$ associated to $D$. We then have the zero section embedding
$i:S\hookrightarrow L$. We denote $L_0=i(S)$ and $j:L^o:=L\backslash L_0\hookrightarrow L$ the open complementary subset.
We denote by $MHW(S\backslash D)^{ex}\times_J MHW(D)$ the category whose set of objects consists of
\begin{equation*}
\left\{(\mathcal M,\mathcal N,a,b),\mathcal M\in MHW(S\backslash D)^{ex},\mathcal N\in MHW(D),
a:\psi_{D1}\mathcal M\to N,b:N\to\psi_{D1}M \right\}
\end{equation*}
where $MHW(S\backslash D)^{ex}\subset MHW(S\backslash D)$ is the full subcategory of extendable objects.
The functor (see definition \ref{phipsiMHM})
\begin{eqnarray*}
(j^*,\phi_{D1},c,v):MHW(S)\to MHW(S\backslash D)^{ex}\times_J MHW(D), \\
((M,F,W),(K,W),\alpha)\mapsto((j^*(M,F,W),j^*(K,W),j^*\alpha),\phi_{D1}((M,F,W),(K,W),\alpha), can,var)
\end{eqnarray*}
is an equivalence of category.
\item[(ii)]Let $S\in\SmVar(\mathbb C)$. Let $D=V(s)\subset S$ a (Cartier) divisor, where $s\in\Gamma(S,L)$ 
is a section of the line bundle $L=L_D$ associated to $D$. We then have the zero section embedding
$i:S\hookrightarrow L$. We denote $L_0=i(S)$ and $j:L^o:=L\backslash L_0\hookrightarrow L$ the open complementary subset.
We denote by $MHW(S\backslash D)\times_J MHW(D)$ the category whose set of objects consists of
\begin{equation*}
\left\{(\mathcal M,\mathcal N,a,b),\mathcal M\in MHW(S\backslash D),\mathcal N\in MHW(D),
a:\psi_{D1}\mathcal M\to N,b:N\to\psi_{D1}M \right\}
\end{equation*}
The functor (see definition \ref{phipsiMHM})
\begin{eqnarray*}
(j^*,\phi_{D1},c,v):MHW(S)\to MHW(S\backslash D)\times_J MHW(D), \\
((M,F,W),(K,W),\alpha)\mapsto((j^*(M,F,W),j^*(K,W),j^*\alpha),\phi_{D1}((M,F,W),(K,W),\alpha), can,var)
\end{eqnarray*}
is an equivalence of category.
\end{itemize}
\end{prop}

\begin{proof}
See \cite{Saito}.
\end{proof}

\begin{thm}
Let $S\in\Var(\mathbb C)$. The category of  mixed Hodge modules is the full subcategory
\begin{equation*}
\iota_S:MHM(S)\hookrightarrow MHW(S)\hookrightarrow\PSh_{\mathcal D(1,0)fil,rh}(S)\times_I P_{fil}(S^{an})
\end{equation*}
consisting of objects 
\begin{equation*}
((M,F,W),(K,W),\alpha)=(((M_I,F,W),u_{IJ}),(K,W),\alpha)\in\PSh_{\mathcal D(1,0)fil,rh}(S)\times_I P_{fil}(S^{an})
\end{equation*}
such that 
$((M,F,W)^{an},(K,W),\alpha):=(((M_I^{an},F,W),u_{IJ}),(K,W),\alpha)\in MHM(S^{an})$.
\end{thm}

\begin{proof}
Follows from GAGA and the extendableness in the algebraic case (proposition \ref{Saprop}).
\end{proof}

Let $S\in\Var(\mathbb C)$. Let $S=\cup_{i\in I}S_i$ an open cover such that there
exists closed embeddings $i_i:S_i\hookrightarrow\tilde S_i$ with $\tilde S_I\in\SmVar(\mathbb C)$.
We have the category $D_{\mathcal D(1,0)fil,rh}(S/(\tilde S_I))\times_I D_{fil}(S^{an})$  
\begin{itemize}
\item whose set of objects is the set of triples $\left\{(((M_I,F,W),u_{IJ}),(K,W),\alpha)\right\}$ with
\begin{eqnarray*} 
((M_I,F,W),u_{IJ})\in D_{\mathcal D(1,0)fil,rh}(S/(\tilde S_I)), \, (K,W)\in D_{fil}(S^{an}), \\ 
\alpha:T(S/(\tilde S_I))(K,W)\otimes\mathbb C_{S^{an}}\to DR(S)^{[-]}(((M_I,W),u_{IJ})^{an})
\end{eqnarray*}
where $\alpha$ is an morphism in $D_{fil}(S^{an}/(\tilde S_{I}^{an}))$,
\item and whose set of morphisms consists of 
\begin{equation*}
\phi=(\phi_D,\phi_C,[\theta]):(((M_{1I},F,W),u_{IJ}),(K_1,W),\alpha_1)\to(((M_{2I},F,W),u_{IJ}),(K_2,W),\alpha_2)
\end{equation*}
where $\phi_D:((M_1,F,W),u_{IJ})\to((M_2,F,W),u_{IJ})$ and $\phi_C:(K_1,W)\to (K_2,W)$ 
are morphisms and
\begin{eqnarray*}
\theta=(\theta^{\bullet},I(DR(S)(\phi^{an}_D))\circ I(\alpha_1),I(\alpha_2)\circ I(\phi_C\otimes I)): \\
I(T(S/(\tilde S_I))(K_1,W))\otimes\mathbb C_{S^{an}}[1]\to I(DR(S)(((M_{2I},W),u_{IJ})^{an}))  
\end{eqnarray*}
is an homotopy,  
$I:D_{fil}(S^{an}/(\tilde S^{an}_{I}))\to K_{fil}(S^{an}/(\tilde S^{an}_{I}))$
being the injective resolution functor, and for
\begin{itemize}
\item $\phi=(\phi_D,\phi_C,[\theta]):(((M_{1I},F,W),u_{IJ}),(K_1,W),\alpha_1)\to(((M_{2I},F,W),u_{IJ}),(K_2,W),\alpha_2)$
\item $\phi'=(\phi'_D,\phi'_C,[\theta']):(((M_{2I},F,W),u_{IJ}),(K_2,W),\alpha_2)\to(((M_{3I},F,W),u_{IJ}),(K_3,W),\alpha_3)$
\end{itemize}
the composition law is given by 
\begin{eqnarray*}
\phi'\circ\phi:=(\phi'_D\circ\phi_D,\phi'_C\circ\phi_C,
I(DR(S)(\phi^{'an}_D))\circ[\theta]+[\theta']\circ I(\phi_C\otimes I)[1]): \\
(((M_{1I},F,W),u_{IJ}),(K_1,W),\alpha_1)\to(((M_{3I},F,W),u_{IJ}),(K_3,W),\alpha_3),
\end{eqnarray*}
in particular for 
$(((M_I,F,W),u_{IJ}),(K,W),\alpha)\in D_{\mathcal D(1,0)fil,rh}(S/(\tilde S_I))\times_I D_{fil}(S^{an})$,
\begin{equation*}
I_{(((M_I,F,W),u_{IJ}),(K,W),\alpha)}=((I_{M_I}),I_K,0),
\end{equation*}
\end{itemize}
and also the category 
$D_{\mathcal D(1,0)fil,rh,\infty}(S/(\tilde S_I))\times_I D_{fil}(S^{an})$ defined in the same way,
together with the localization functor
\begin{eqnarray*}
(D(zar),I):C_{\mathcal D(1,0)fil,rh}(S/(\tilde S_I))\times_I D_{fil}(S^{an})
\to D_{\mathcal D(1,0)fil,rh}(S/(\tilde S_I))\times_I D_{fil}(S^{an}) \\
\to D_{\mathcal D(1,0)fil,rh,\infty}(S/(\tilde S_I))\times_I D_{fil}(S^{an}).
\end{eqnarray*}
Moreover,
\begin{itemize}
\item For $(((M_{I},F,W),u_{IJ}),(K,W),\alpha)D_{\mathcal D(1,0)fil,rh}(S/(\tilde S_I))\times_I D_{fil}(S^{an})$, we set
\begin{equation*}
(((M_{I},F,W),u_{IJ}),(K,W),\alpha)[1]:=(((M_{I},F,W),u_{IJ})[1],(K,W)[1],\alpha[1]).
\end{equation*}
\item For 
\begin{equation*}
\phi=(\phi_D,\phi_C,[\theta]):(((M_{1I},F,W),u_{IJ}),(K_1,W),\alpha_1)\to(((M_{2I},F,W),u_{IJ}),(K_2,W),\alpha_2)
\end{equation*}
a morphism in $D_{\mathcal D(1,0)fil,rh}(S/(\tilde S_I))\times_I D_{fil}(S^{an})$, we set (see \cite{CG} definition 3.12)
\begin{eqnarray*}
\Cone(\phi):=(\Cone(\phi_D),\Cone(\phi_C),((\alpha_1,\theta),(\alpha_2,0)))
\in D_{\mathcal D(1,0)fil,rh}(S/(\tilde S_I))\times_I D_{fil}(S^{an}).
\end{eqnarray*}
together with the canonical maps
\begin{itemize}
\item $c_1(-)=(c_1(\phi_D),c_1(\phi_C),0):(((M_{2I},F,W),u_{IJ}),(K_2,W),\alpha_2)\to\Cone(\phi)$
\item $c_2(-)=(c_2(\phi_D),c_2(\phi_C),0):\Cone(\phi)\to (((M_{1I},F,W),u_{IJ}),(K_1,W),\alpha_1)[1]$.
\end{itemize}
\end{itemize}

We now state and prove the following key theorem :

\begin{thm}\label{Be}
\begin{itemize}
\item[(i)]Let $S\in\Var(\mathbb C)$. Let $S=\cup_{i\in I}S_i$ an open cover such that there exists
closed embedding $i_i:S_i\hookrightarrow\tilde S_i$ with $\tilde S_i\in\SmVar(\mathbb C)$. Then the full embedding
\begin{eqnarray*}
\iota_S:MHM(S)\hookrightarrow\PSh^0_{\mathcal D(1,0)fil,rh}(S/(\tilde S_I))\times_I P_{fil}(S^{an}) 
\hookrightarrow C_{\mathcal D(1,0)fil,rh}(S/(\tilde S_I))\times_I D_{fil}(S^{an}) 
\end{eqnarray*}
induces a full embedding
\begin{equation*}
\iota_S:D(MHM(S))\hookrightarrow D_{\mathcal D(1,0)fil,rh}(S/(\tilde S_I))\times_I D_{fil}(S^{an}) 
\end{equation*}
whose image consists of 
$(((M_I,F,W),u_{IJ}),(K,W),\alpha)\in D_{\mathcal D(1,0)fil,rh}(S/(\tilde S_I))\times_I D_{fil}(S^{an})$ such that 
\begin{equation*}
((H^n(M_I,F,W),H^n(u_{IJ})),H^n(K,W),H^n\alpha)\in MHM(S) 
\end{equation*}
for all $n\in\mathbb Z$ and such that for all $p\in\mathbb Z$,
the differentials of $\Gr_W^p(M_I,F)$ are strict for the filtrations $F$.
\item[(ii)]Let $S\in\Var(\mathbb C)$. Let $S=\cup_{i\in I}S_i$ an open cover such that there exists
closed embedding $i_i:S_i\hookrightarrow\tilde S_i$ with $\tilde S_i\in\SmVar(\mathbb C)$. Then the full embedding
\begin{eqnarray*}
\iota_S:MHM(S)\hookrightarrow\PSh^0_{\mathcal D(1,0)fil,rh}(S/(\tilde S_I))\times_I P_{fil}(S^{an})
\hookrightarrow C_{\mathcal D(1,0)fil,rh}(S/(\tilde S_I))\times_I D_{fil}(S^{an}) 
\end{eqnarray*}
induces a full embedding
\begin{equation*}
\iota_S:D(MHM(S))\hookrightarrow D_{\mathcal D(1,0)fil,\infty,rh}(S/(\tilde S_I))\times_I D_{fil}(S^{an}) 
\end{equation*}
whose image consists of 
$(((M_I,F,W),u_{IJ}),(K,W),\alpha)\in D_{\mathcal D(1,0)fil,\infty,rh}(S/(\tilde S_I))\times_I D_{fil}(S^{an})$ such that 
\begin{equation*}
((H^n(M_I,F,W),H^n(u_{IJ})),H^n(K,W),H^n\alpha)\in MHM(S) 
\end{equation*}
for all $n\in\mathbb Z$ and such that there exist $r\in\mathbb Z$ and an $r$-filtered homotopy equivalence
$((M_I,F,W),u_{IJ})\to ((M'_I,F,W),u_{IJ})$ such that for all $p\in\mathbb Z$
the differentials of $\Gr_W^p(M'_I,F)$ are strict for the filtrations $F$.
\end{itemize}
\end{thm}

\begin{proof}
\noindent(i): We first show that $\iota_S$ is fully faithfull, that is for all
$\mathcal M=(((M_I,F,W),u_{IJ}),(K,W),\alpha),\mathcal M'=(((M'_I,F,W),u_{IJ}),(K',W),\alpha')\in MHM(S)$ 
and all $n\in\mathbb Z$,
\begin{eqnarray*}
\iota_S:\Ext_{D(MHM(S))}^n(\mathcal M,\mathcal M'):=\Hom_{D(MHM(S))}(\mathcal M,\mathcal M'[n]) \\
\to\Ext_{\mathcal D(S)}^n(\mathcal M,\mathcal M')
:=\Hom_{\mathcal D(S):=D_{\mathcal D(1,0)fil,rh}(S/(\tilde S_I))\times_I D_{fil}(S^{an})}(\mathcal M,\mathcal M'[n])
\end{eqnarray*}
For this it is enough to assume $S$ smooth. We then proceed by induction on $max(\dim\supp(M),\dim\supp(M'))$. 
\begin{itemize}
\item For $\supp(M)=\supp(M')=\left\{s\right\}$, 
it is the theorem for mixed hodge complexes or absolute Hodge complexes, see \cite{CG}. 
If $\supp(M)=\left\{s\right\}$ and $\supp(M')=\left\{s'\right\}$, then by the localization exact sequence
\begin{equation*}
\Ext_{D(MHM(S))}^n(\mathcal M,\mathcal M')=0=\Ext_{\mathcal D(S)}^n(\mathcal M,\mathcal M')
\end{equation*}
\item Denote $\supp(M)=Z\subset S$ and $\supp(M')=Z'\subset S$.
There exist an open subset $S^o\subset S$ such that $Z^o:=Z\cap S^o$ and $Z^{'o}:=Z'\cap S^o$ are smooth,
and $\mathcal M_{|Z^o}:=((i^*\Gr_{V_{Z^o},0}M_{|S^o},F,W),i^*j^*(K,W),\alpha^*(i))\in MHM(Z^o)$ and 
$\mathcal M'_{|Z^{'o}}:=((i^{'*}\Gr_{V_{Z^{'o}},0}M'_{|S^o},F,W),i^{'*}j^*K,\alpha^*(i'))\in MHM(Z^{'o})$ 
are variation of mixed Hodge structure, where $j:S^o\hookrightarrow S$ is the open embedding, and
$i:Z^o\hookrightarrow S^o$, $i:Z^o\hookrightarrow S^o$ the closed embeddings.
Considering the connected components of $Z^o$ and $Z^{'o}$, we way assume that $Z^o$ and $Z^{'o}$ are connected.
Shrinking $S^o$ if necessary, we may assume that either $Z^o=Z^{'o}$ or $Z^o\cap Z^{'o}=\emptyset$,
We denote $D=S\backslash S^o$. Shrinking $S^o$ if necessary, 
we may assume that $D$ is a divisor and denote by $l:S\hookrightarrow L_D$ the zero section embedding.
\begin{itemize}
\item If $Z^o=Z^{'o}$, denote $i:Z^o\hookrightarrow S^o$ the closed embedding.
We have then the following commutative diagram
\begin{equation*}
\xymatrix{\Ext_{D(MHM(S^o))}^n(\mathcal M_{|S^o},\mathcal M'_{|S^o})
\ar[rr]^{\iota_{S^o}}\ar[d]_{(i^*\Gr_{V_{Z^o},0},i^*,\alpha^*(i))} & \, & 
\Ext_{\mathcal D(S^o)}^n(\mathcal M_{|S^o},\mathcal M'_{|S^o})\ar[d]^{(i^*\Gr_{V_{Z^o},0},i^*,\alpha^*(i))} \\
\Ext_{D(MHM(Z^o))}^n(\mathcal M_{|Z^o},\mathcal M'_{|Z^o})\ar[rr]^{\iota_{Z^o}}\ar[u]_{(i_{*mod},i_*,\alpha_*(i))} & \, &
\Ext_{\mathcal D(Z^o)}^n(\mathcal M_{|Z^o},\mathcal M'_{|Z^o})\ar[u]^{(i_{*mod},i_*,\alpha_*(i))}}
\end{equation*}
Now we prove that $\iota_{Z^o}$ is an isomorphism similarly to the proof the the generic case of \ref{Be}.
On the other hand the left and right colummn are isomorphisms.
Hence $\iota_{S^o}$ is an isomorphism by the diagram.
\item If $Z^o\cap Z^{'o}=\emptyset$, we consider the following commutative diagram
\begin{equation*}
\xymatrix{\Ext_{D(MHM(S^o))}^n(\mathcal M_{|S^o},\mathcal M'_{|S^o})
\ar[rr]^{\iota_{S^o}}\ar[d]_{(i^*\Gr_{V_{Z^o},0},i^*,\alpha^*(i))} & \, & 
\Ext_{\mathcal D(S^o)}^n(\mathcal M_{|S^o},\mathcal M'_{|S^o})\ar[d]^{(i^*\Gr_{V_{Z^o},0},i^*,\alpha^*(i))} \\
\Ext_{D(MHM(Z^o))}^n(\mathcal M_{|Z^o},0)=0\ar[rr]^{\iota_{Z^o}}\ar[u]_{(i_{*mod},i_*,\alpha_*(i))} & \, &
\Ext_{\mathcal D(Z^o)}^n(\mathcal M_{|Z^o},0)=0\ar[u]^{(i_{*mod},i_*,\alpha_*(i))}}
\end{equation*}
where the left and right column are isomorphism by strictness of the $V_{Z^o}$ filtration
(use a bi-filtered injective resolution with respect to $F$ and $V_{Z^o}$ for the right column).
\end{itemize}
\item We consider now the following commutative diagram in $C(\mathbb Z)$ where we denote for short $H:=D(MHM(S))$
\begin{equation*}
\xymatrix{0\ar[r] & \Hom_{H}^{\bullet}(\Gamma^{\vee,Hdg}_D\mathcal M,\Gamma^{Hdg}_D\mathcal M')
\ar[r]^{\Hom(-,\gamma^{Hdg}_D(\mathcal M'))}\ar[d]^{\iota_S} &
\Hom_{H}^{\bullet}(\Gamma^{\vee,Hdg}_D\mathcal M,\mathcal M')
\ar[r]^{\Hom(-,\ad(j^*,j_{*Hdg})(\mathcal M'))}\ar[d]^{\iota_S} &
\Hom_{H}^{\bullet}(\Gamma^{\vee,Hdg}_D\mathcal M,j_{*Hdg}j^*\mathcal M')\ar[r]\ar[d]^{\iota_S} & 0 \\
0\ar[r] & \Hom_{\mathcal D(S)}^{\bullet}(\Gamma^{\vee,Hdg}_D\mathcal M,\Gamma^{Hdg}_D\mathcal M')
\ar[r]^{\Hom(-,\gamma^{Hdg}_D(\mathcal M'))} &
\Hom_{\mathcal D(S)}^{\bullet}(\Gamma^{\vee,Hdg}_D\mathcal M,\mathcal M')\ar[r]^{\Hom(-,\ad(j^*,j_{*Hdg})(\mathcal M'))} &
\Hom_{\mathcal D(S)}^{\bullet}(\Gamma^{\vee,Hdg}_D\mathcal M,j_{*Hdg}j^*\mathcal M')\ar[r] & 0}
\end{equation*}
whose lines are exact sequence. We have on the one hand,
\begin{equation*}
\Hom_{D(MHM(S))}^{\bullet}(\Gamma^{\vee,Hdg}_D\mathcal M,j_{*Hdg}j^*\mathcal M')=0=
\Hom_{\mathcal D(S)}^{\bullet}(\Gamma^{\vee,Hdg}_D\mathcal M,j_{*Hdg}j^*\mathcal M')
\end{equation*}
On the other hand by induction hypothesis
\begin{equation*}
\iota_S:\Hom_{D(MHM(S))}^{\bullet}(\Gamma^{\vee,Hdg}_D\mathcal M,\Gamma^{Hdg}_D\mathcal M')\to
\Hom_{\mathcal D(S)}^{\bullet}(\Gamma^{\vee,Hdg}_D\mathcal M,\Gamma^{Hdg}_D\mathcal M')
\end{equation*}
is a quasi-isomorphism. Hence, by the diagram
\begin{equation*}
\iota_S:\Hom_{D(MHM(S))}^{\bullet}(\Gamma^{\vee,Hdg}_D\mathcal M,\mathcal M')\to
\Hom_{\mathcal D(S)}^{\bullet}(\Gamma^{\vee,Hdg}_D\mathcal M,\mathcal M')
\end{equation*}
is a quasi-isomorphism.
\item We consider now the following commutative diagram in $C(\mathbb Z)$ where we denote for short $H:=D(MHM(S))$
\begin{equation*}
\xymatrix{0\ar[r] & \Hom_{H}^{\bullet}(\Gamma^{\vee,Hdg}_D\mathcal M,\mathcal M')
\ar[r]^{\Hom(\gamma^{\vee,Hdg}_D(\mathcal M),-)}\ar[d]^{\iota_S} &
\Hom_{H}^{\bullet}(\mathcal M,\mathcal M')
\ar[r]^{\Hom(\ad(j_{!Hdg},j^*)(\mathcal M'),-)}\ar[d]^{\iota_S} &
\Hom_{H}^{\bullet}(j_{!Hdg}j^*\mathcal M,\mathcal M')\ar[r]\ar[d]^{\iota_S} & 0 \\
0\ar[r] & \Hom_{\mathcal D(S)}^{\bullet}(\Gamma^{\vee,Hdg}_D\mathcal M,\mathcal M')
\ar[r]^{\Hom(\gamma^{\vee,Hdg}_D(\mathcal M),-)} &
\Hom_{\mathcal D(S)}^{\bullet}(\mathcal M,\mathcal M')\ar[r]^{\Hom(\ad(j_{!Hdg},j^*)(\mathcal M),-)} &
\Hom_{\mathcal D(S)}^{\bullet}(j_{!Hdg}j^*\mathcal M,\mathcal M')\ar[r] & 0}
\end{equation*}
whose lines are exact sequence. On the one hand, the commutative diagram
\begin{equation*}
\xymatrix{\Hom_{D(MHM(S))}^{\bullet}(j_{!Hdg}j^*\mathcal M,\mathcal M')\ar[r]^{j^*}\ar[d]^{\iota_{S}} &
\Hom_{D(MHM(S^o))}^{\bullet}(j^*\mathcal M,j^*\mathcal M')\ar[d]^{\iota_{S^o}} \\
\Hom_{\mathcal D(S)}^{\bullet}(j_{!Hdg}j^*\mathcal M,\mathcal M')\ar[r]^{j^*} &
\Hom_{\mathcal D(S^o)}^{\bullet}(j^*\mathcal M,j^*\mathcal M')}
\end{equation*}
together with the fact that the horizontal arrows $j^*$ are quasi-isomorphism 
by the functoriality given the uniqueness of the $V_S$ filtration for the embedding $l:S\hookrightarrow L_D$, 
(use a bi-filtered injective resolution with respect to $F$ and $V_S$ for the lower arrow)
and the fact that $\iota_{S^o}$ is a quasi-isomorphism by the first two point, show that
\begin{equation*}
\iota_S:\Hom_{D(MHM(S))}^{\bullet}(j_{!Hdg}j^*\mathcal M,\mathcal M')\to
\Hom_{\mathcal D(S)}^{\bullet}(j_{!Hdg}j^*\mathcal M,\mathcal M')
\end{equation*}
is a quasi-isomorphism. On the other hand, by the third point
\begin{equation*}
\iota_S:\Hom_{D(MHM(S))}^{\bullet}(\Gamma^{\vee,Hdg}_D\mathcal M,\mathcal M')\to
\Hom_{\mathcal D(S)}^{\bullet}(\Gamma^{\vee,Hdg}_D\mathcal M,\mathcal M')
\end{equation*}
is a quasi-isomorphism. Hence, by the diagram
\begin{equation*}
\iota_S:\Hom_{D(MHM(S))}^{\bullet}(\Gamma^{\vee,Hdg}_D\mathcal M,\mathcal M')\to
\Hom_{\mathcal D(S)}^{\bullet}(\Gamma^{\vee,Hdg}_D\mathcal M,\mathcal M')
\end{equation*}
is a quasi-isomorphism.
\end{itemize}
This shows the fully faithfulness. We now prove the essential surjectivity : let
\begin{equation*}
(((M_I,F,W),u_{IJ}),(K,W),\alpha)\in C_{\mathcal D(1,0)fil,rh}(S/(\tilde S_I))\times_I D_{fil}(S^{an}) 
\end{equation*}
such that the cohomology are mixed hodge modules and such that the differential are strict.
We proceed by induction on $card\left\{n\in\mathbb Z\right\}, \, \mbox{s.t.} H^n(M_I,F,W)\neq 0$ by taking for 
the cohomological troncation 
\begin{equation*}
\tau^{\leq n}(((M_I,F,W),u_{IJ}),(K,W),\alpha):=
((\tau^{\leq n}(M_I,F,W),\tau^{\leq n}u_{IJ}),\tau^{\leq n}(K,W),\tau^{\leq n}\alpha)
\end{equation*}
and using the fact that the differential are strict for the filtration $F$ and the fully faithfullness.

\noindent(ii): Follows from (i) and the fact that the image of the embedding given by (i) 
consists of classes of complexes whose differential are strict for $F$.
\end{proof}

\begin{defi}\label{DHdgalpha}
Let $f:X\to S$ a morphism with $X,S\in\Var(\mathbb C)$. 
Assume there exist a factorization $f:X\xrightarrow{l}Y\times S\xrightarrow{p_S}S$
with $Y\in\SmVar(\mathbb C)$, $l$ a closed embedding and $p_S$ the projection.
Let $\bar Y\in\PSmVar(\mathbb C)$ a smooth compactification of $Y$ with $n:Y\hookrightarrow\bar Y$ the open embedding.
Then $\bar f:\bar X\xrightarrow{\bar l}\bar Y\times_S\xrightarrow{\bar p_S}S$ is a compactification of $f$,
with $\bar X\subset\bar Y\times S$ the closure of $X$ and $\bar l$ the closed embedding,
and we denote by $n':X\hookrightarrow\bar X$ the open embedding so that $f=\bar f\circ n'$.
\begin{itemize}
\item[(i)]For $(((M_I,F,W),u_{IJ}),(K,W),\alpha)\in C(MHM(X))$, 
we define, using definition \ref{DHdgjalpha} and theorem \ref{Be}
\begin{eqnarray*}
Rf_{*Hdg}(((M_I,F,W),u_{IJ}),(K,W),\alpha):=
\iota_S^{-1}(\int_{\bar f}^{FDR}(n\times I)_{*Hdg}((M_I,F,W),u_{IJ}),Rf_{*w}(K,W),f_*(\alpha)) \\
\in D(MHM(S))
\end{eqnarray*}
where $f_*(\alpha)$ is given in definition \ref{falpha}, and since 
\begin{itemize} 
\item by definition 
\begin{equation*}
H^i(\int_{\bar f}^{FDR}\Gr_W^k(I\times n)_{Hdg*}((M_I,F,W),u_{IJ}),
R\bar f_*\Gr_W^kn'_{*w}(K,W),\bar f_*\Gr_W^kn'_*\alpha)\in HM(S) 
\end{equation*}
for all $i,k\in\mathbb Z$, hence by the spectral sequence for the filtered complexes 
$\int_{\bar f}^{FDR}(I\times n)_{Hdg*}((M_I,F,W),u_{IJ})$ and $R\bar f_*((I\times n)_{*w}(K,W))$
\begin{eqnarray*}
\Gr_W^kH^i(\int_f^{Hdg}((M_I,F,W),u_{IJ}),Rf_{*w}(K,W),f_*\alpha)):= \\
(\Gr_W^kH^i\int_{\bar f}^{FDR}(I\times n)_{Hdg*}((M_I,F,W),u_{IJ}),
\Gr_W^kH^iR\bar f_*n'_{*w}(K,W),\Gr_W^kH^if_*\alpha)\in HM(S) 
\end{eqnarray*}
this gives by definition 
$H^i(\int_f^{Hdg}((M_I,F,W),u_{IJ}),Rf_{*w}(K,W),f_*(\alpha))\in MHM(S)$ for all $i\in\mathbb Z$. 
\item $\int_{f}^{Hdg}((M_I,F,W),u_{IJ})$ is the class of a complex such that the differential are strict for $F$
by theorem \ref{Sa1}.
\end{itemize}
\item[(ii)]For $(((M_I,F,W),u_{IJ}),(K,W),\alpha)\in C(MHM(X))$, 
we define, using definition \ref{DHdgjalpha} and theorem \ref{Be}, 
\begin{eqnarray*}
Rf_{!Hdg}(((M_I,F,W),u_{IJ}),(K,W),\alpha):=
\iota_S^{-1}(\int_{\bar f}^{FDR}(n\times I)_{!Hdg}((M_I,F,W),u_{IJ}),Rf_{!w}(K,W),f_!(\alpha)) \\
\in D(MHM(S))
\end{eqnarray*}
where $f_!(\alpha)$ is given in definition \ref{falpha}, and since 
\begin{itemize} 
\item by definition 
\begin{equation*}
H^i(\int_{\bar f}^{FDR}\Gr_W^k(n\times I)_{!Hdg}((M_I,F,W),u_{IJ}),
R\bar f_*\Gr_W^kn'_{!w}(K,W),\Gr_W^kf_!\alpha)\in HM(S) 
\end{equation*}
for all $i,k\in\mathbb Z$, 
hence by the spectral sequence for the filtered complexes 
$\int_{\bar f}^{FDR}(n\times I)_{!Hdg}((M_I,F,W),u_{IJ})$ and $R\bar f_*(n\times I)_{!w}(K,W)$
\begin{eqnarray*}
\Gr_W^kH^i(\int_{f!}^{Hdg}((M_I,F,W),u_{IJ}),Rf_{!w}(K,W),f_!\alpha):= \\
(\Gr_W^kH^i\int_{\bar f}^{FDR}(n\times I)_{!Hdg}((M_I,F,W),u_{IJ}),
\Gr_W^kH^iR\bar f_*n'_{!w}(K,W),\Gr_W^kH^if_!\alpha)\in HM_{gm,k,\mathbb C}(S) 
\end{eqnarray*}
this gives by definition 
$H^i(\int_{f!}^{Hdg}((M_I,F,W),u_{IJ}),Rf_{!w}(K,W),f_!(\alpha))\in MHM(S)$ for all $i\in\mathbb Z$. 
\item $\int_{f!}^{Hdg}((M_I,F,W),u_{IJ})$ is the class of a complex such that the differential are strict for $F$
by theorem \ref{Sa1}.
\end{itemize}
\end{itemize}
\end{defi}

We have the six functors formalism for mixed Hodge modules on quasi-projective varieties :

\begin{defi}\label{2functsHdg} 
Let $f:X\to S$ a morphism with $X,S\in\QPVar(\mathbb C)$.
Then, since $X$ is quasi-projective, there exist a factorization $f:X\xrightarrow{l}\mathbb P^{N,o}\times S\xrightarrow{p_S}S$ 
with $n_0:\mathbb P^{N,o}\hookrightarrow\mathbb P^N$ an open subset, $l$ a closed embedding and $p_S$ the projection.
Since $S$ is quasi-projective, there exist a closed embedding $i:S\hookrightarrow\tilde S$ with $\tilde S\in\SmVar(\mathbb C)$.
We have then the commutative diagram
\begin{equation*}
\xymatrix{f:X\ar[r]^{l}\ar[rd] & \mathbb P^{N,o}\times S\ar[r]^{p_S}\ar[d]^{i':=(I\times i)} & S\ar[d]^{i} \\ 
\, & \mathbb P^{N,o}\times\tilde S\ar[d]^{n:=(n_0\times I)}\ar[r]^{p_{\tilde S}} & \tilde S\ar[d]^{=} \\ 
\, & \mathbb P^N\times\tilde S\ar[r]^{\bar p_{\tilde S}} & \tilde S} 
\end{equation*}
\begin{itemize}
\item[(i)] For $((M,F,W),(K,W),\alpha)\in D(MHM(X))$,
where $(M,F,W)\in C_{\mathcal D(1,0)fil}(X/\mathbb P^{N,o}\times\tilde S)$ and $(K,W)\in C_{fil}(X^{an})$,
we define, using definition \ref{DHdgalpha},
\begin{eqnarray*}
f_{*Hdg}((M,F,W),(K,W),\alpha):&=&(Rf_*^{Hdg}(M,F,W),Rf_{*w}(K,W),f_*(\alpha)) \\
:&=&\iota_S^{-1}(\int^{FDR}_{p_{\tilde S}}n^{Hdg}_*(M,F,W),Rf_{*w}(K,W),f_*(\alpha))\in D(MHM(S))
\end{eqnarray*} 
with $f_*(\alpha)$ given in definition \ref{falpha}. 
\item[(ii)] For $((M,F,W),(K,W),\alpha)\in D(MHM(X))$,
where $(M,F,W)\in C_{\mathcal D(1,0)fil}(X/\mathbb P^{N,o}\times\tilde S)$ and $(K,W)\in C_{fil}(X^{an})$,
we define, using definition \ref{DHdgalpha}
\begin{eqnarray*}
f_{!Hdg}((M,F,W),(K,W),\alpha):&=&(Rf^{Hdg}_!(M,F,W),Rf_{!w}(K,W),f_!(\alpha)) \\
:&=&\iota_S^{-1}(\int^{FDR}_{p_{\tilde S}}n^{Hdg}_!(M,F,W)),Rf_{!w}(K,W),f_!(\alpha))\in D(MHM(S))
\end{eqnarray*}
with $f_!(\alpha)$ given in definition \ref{falpha}. 
\item[(iii)] For $((M,F,W),(K,W),\alpha)\in D(MHM(S))$,
where $(M,F,W)\in C_{\mathcal D(1,0)fil}(S/(\tilde S))$, $(K,W)\in C_{fil}(S^{an})$,
we define, using definition \ref{inverseHdgsing} 
(see theorem \ref{Sa2}(ii) for $p_S$ and definition \ref{gammaHdg} for $i\circ l$),
\begin{eqnarray*}
f^{*Hdg}((M,F,W),(K,W),\alpha):&=&(f^{\hat*mod}_{Hdg}(M,F,W),f^{*w}(K,W),f^*\alpha) \\
:&=&(\Gamma_X^{\vee,Hdg}p_{\tilde S}^{\hat*mod[-]}(M,F,W),\Gamma_X^{\vee,w}p_{S}^*(K,W),f^*(\alpha))\in D(MHM(X))
\end{eqnarray*}
with $f^*\alpha$ given in definition \ref{falpha}. 
For $j:S^o\hookrightarrow S$ an open embedding and $((M,F,W),(K,W),\alpha)\in D(MHM(S))$, we have (see \cite{Saito})
\begin{equation*}
j^{*Hdg}((M,F,W),(K,W),\alpha)=(j^*(M,F,W),j^*(K,W),j^*\alpha)\in D(MHM(S^o)).
\end{equation*}
\item[(iv)] For $((M,F,W),(K,W),\alpha)\in D(MHM(S))$,
where $(M,F,W)\in C_{\mathcal D(1,0)fil}(S/(\tilde S))$, $(K,W)\in C_{fil}(S^{an})$,
we define, using definition \ref{inverseHdgsing}
(see theorem \ref{Sa2}(ii) for $p_S$ and definition \ref{gammaHdg} for $i\circ l$),
\begin{eqnarray*}
f^{!Hdg}((M,F,W),(K,W),\alpha):&=&(f^{*mod}_{Hdg}(M,F,W),f^{!w}(K,W),f^!\alpha) \\
:&=&(\Gamma_X^{Hdg}p_{\tilde S_I}^{*mod[-]}(M,F,W)(d_X)[2d_X], \Gamma_X^wp_S^*(K,W),f^!(\alpha))\in D(MHM(X))
\end{eqnarray*}
with $f^!\alpha$ given in definition \ref{falpha}. 
For $j:S^o\hookrightarrow S$ an open embedding and $((M,F,W),(K,W),\alpha)\in D(MHM(S))$, we have (see \cite{Saito})
\begin{equation*}
j^{!Hdg}((M,F,W),(K,W),\alpha)=(j^*(M,F,W),j^*(K,W),j^*\alpha)\in D(MHM(S^o)).
\end{equation*}
\end{itemize}
Using the unicity of proposition \ref{Saprop}, we see that these definitions does NOT depends on the choice of the factorization
$f:X\xrightarrow{l}Y\times S\xrightarrow{p_S}S$ of $f$.
Moreover, using the unicity of proposition \ref{Saprop} and proposition \ref{compDmod}, 
we see that they are 2 functors on the category of quasi-projective complex algebraic varieties $(\Var(\mathbb C))^{QP}$.
\item[(v)]Let $S\in\Var(\mathbb C)$. Take an open cover $S=\cup_{i=1}^lS_i$ such that there exists closed
embeddings $i_i:S_i\hookrightarrow\tilde S_i$ with $\tilde S_i\in\SmVar(\mathbb C)$.
We have the functor
\begin{eqnarray*}
((-)\otimes^{Hdg}_{O_S}(-),(-)\otimes(-)):C(MHM(S)^2\to C(MHM(S)), \\
((((M_I,F,W),u_{IJ}),(K,W),\alpha),(((N_I,F,W),v_{IJ}),(K',W),\alpha'))\mapsto \\
(((M_I,F,W),u_{IJ})\otimes^{Hdg}_{O_S}((N_I,F,W),v_{IJ}),(K,W)\otimes(K',W),\alpha\otimes\alpha')
\end{eqnarray*}
where 
\begin{eqnarray*}
(((M,F,W),u_{IJ}),((N,F,W),v_{IJ}))\mapsto ((M,F,W),u_{IJ})\otimes_{O_S}^{Hdg}((N,F,W),v_{IJ}):= \\
\Delta_{S,Hdg}^{*mod}(p_{1I}^{*mod}(M_I,F,W)\otimes_{O_{\tilde S_I\times\tilde S_I}}p_{2I}^{*mod}(N_I,F,W),
p_{1I}^{*mod}u_{IJ}\otimes p_{2I}^{*mod}v_{IJ}):= \\
\Delta_S^*\Gr_{V_{\Delta_S,0}}\Gamma_{\Delta_S}^{\vee,Hdg}
(p_{1I}^{*mod}(M_I,F,W)\otimes_{O_{\tilde S_I\times\tilde S_I}}p_{2I}^{*mod}(N_I,F,W),
p_{1I}^{*mod}u_{IJ}\otimes p_{2I}^{*mod}v_{IJ})
\end{eqnarray*}
and the map $\alpha\otimes\alpha'$ is given in definition \ref{falpha}.
\end{defi}

We have then the following

\begin{thm}\label{adHdg}
Let $f:X\to S$ a morphism with $X,S\in\Var(\mathbb C)$, $X$ quasi-projective. Then,
\begin{itemize}
\item[(i)] $(f^{*Hdg},f_{*Hdg}):D(MHM(S))\to D(MHM(X))$ is a pair of adjoint functors,
\item[(ii)]$(f^{*Hdg},f_{*Hdg}):D(MHM(S))\to D(MHM(X))$ is a pair of adjoint functors.
\end{itemize}
\end{thm}

\begin{proof}
For the projection case see section 4. For the open embedding see definition \ref{DHdgj}. 
\end{proof}
 
Definition \ref{2functsHdg} gives the following 2 functors :
\begin{itemize}
\item We have the following 2 functor on the category of complex algebraic varieties
\begin{eqnarray*}
D(MHW(\cdot)):\Var(\mathbb C)\to\TriCat, \; S\mapsto D(MHW(S)), \\
(f:T\to S)\longmapsto (f^{*Hdg}:((M,F,W),(K,W),\alpha)\mapsto \\
f^{*Hdg}((M,F,W),(K,W),\alpha):=(f_{Hdg}^{\hat*mod}(M,F,W),f^{*w}(K,W),f^*(\alpha))).
\end{eqnarray*}
\item We have the following 2 functor on the category of complex quasi-projective algebraic varieties
\begin{eqnarray*}
D(MHW(\cdot)):\QPVar(\mathbb C)\to\TriCat, \; S\mapsto D(MHW(S)), \\
(f:T\to S)\longmapsto (f_{*Hdg}:((M,F,W),(K,W),\alpha)\mapsto \\
f_{*Hdg}((M,F,W),(K,W),\alpha):=(Rf^{Hdg}_*(M,F,W),Rf_{*w}(K,W),f_*(\alpha))).
\end{eqnarray*}
\item We have the following 2 functor on the category of complex quasi-projective algebraic varieties
\begin{eqnarray*}
D(MHW(\cdot)):\QPVar(\mathbb C)\to\TriCat, \; S\mapsto D(MHW(S)), \\
(f:T\to S)\longmapsto (f_{!Hdg}:((M,F,W),(K,W),\alpha)\mapsto \\
f_{!Hdg}((M,F,W),(K,W),\alpha):=(Rf_!^{Hdg}(M,F,W),Rf_{!w}(K,W),f_!(\alpha))).
\end{eqnarray*}
\item We have the following 2 functor on the category of complex algebraic varieties
\begin{eqnarray*}
D(MHW(\cdot)):\Var(\mathbb C)\to\TriCat, \; S\mapsto D(MHW(S)), \\
(f:T\to S)\longmapsto (f^{!Hdg}:((M,F,W),(K,W),\alpha)\mapsto \\
f^{!Hdg}((M,F,W),(K,W),\alpha):=(f_{Hdg}^{*mod}(M,F,W),f^{!w}(K,W),f^!(\alpha))).
\end{eqnarray*}
\end{itemize}

For a commutative diagram in $\Var(\mathbb C)$
\begin{equation*}
D=\xymatrix{X\ar[r]^{f} & S \\
X'\ar[r]^{f'}\ar[u]^{g'} & T\ar[u]^{g}}
\end{equation*}
with $S,T,X',X$ quasi-projective, 
we have, for $((M,F,W),(K,W),\alpha)\in D(MHM(X))$ using theorem \ref{adHdg}, the following transformations maps
\begin{eqnarray*}
T_1^{Hdg}(D)((M,F,W),(K,W),\alpha): \\
g^{*Hdg}f_{*Hdg}((M,F,W),(K,W),\alpha)\xrightarrow{\ad(f^{'*Hdg},f_{*Hdg})(-)} 
f'_{*Hdg}f^{'*Hdg}g^{*Hdg}f_{*Hdg}((M,F,W),(K,W),\alpha) \\
\xrightarrow{=} f'_{*Hdg}g^{'*Hdg}f^{*Hdg}f_{*Hdg}((M,F,W),(K,W),\alpha) 
\xrightarrow{\ad(f^{*Hdg},f_{*Hdg})(-)}f'_{*Hdg}g^{'*Hdg}((M,F,W),(K,W),\alpha)
\end{eqnarray*}
and
\begin{eqnarray*}
T_2^{Hdg}(D)((M,F,W),(K,W),\alpha): \\
f'_{!Hdg}g^{'!Hdg}((M,F,W),(K,W),\alpha)\xrightarrow{\ad(f^{!Hdg},f_{!Hdg})(-)} 
f'_{!Hdg}g^{'!Hdg}f^{!Hdg}f_{!Hdg}((M,F,W),(K,W),\alpha) \\
\xrightarrow{=} f'_{!Hdg}f^{'!Hdg}g^{!Hdg}f_{!Hdg}((M,F,W),(K,W),\alpha) 
\xrightarrow{\ad(f^{'*Hdg},f'_{*Hdg})(-)}g^{!Hdg}f_{!Hdg}((M,F,W),(K,W),\alpha)
\end{eqnarray*}

One consequence of the unicity of proposition \ref{Saprop} is the following :

\begin{prop}\label{PDmod1Hdg}
For a commutative diagram in $\Var(\mathbb C)$
\begin{equation*}
D=\xymatrix{X\ar[r]^{f} & S \\
X_T\ar[r]^{f'}\ar[u]^{g'} & T\ar[u]^{g}}
\end{equation*}
which is cartesian, with $S,T,X',X$ quasi-projective and $f$ (hence $f'$ proper), and $((M,F,W),(K,W),\alpha)\in D(MHM(X))$
\begin{eqnarray*}
T_1^{Hdg}(f,g):((M,F,W),(K,W),\alpha): \\
g^{*Hdg}f_{*Hdg}((M,F,W),(K,W),\alpha)\xrightarrow{\sim}f'_{*Hdg}g^{'*Hdg}((M,F,W),(K,W),\alpha)
\end{eqnarray*}
is an isomorphism.
\end{prop}

\begin{proof}
See \cite{Saito}.
\end{proof}

We have the following proposition 

\begin{prop}\label{jD}
Let $Y\in\PSmVar(\mathbb C)$ and $i:Z\hookrightarrow S$ a closed embedding with $Z$ smooth. 
Denote by $j:U:=S\backslash Z\hookrightarrow Y$ the complementary open subset. 
\begin{itemize}
\item[(i)]We have
\begin{eqnarray*}
a_{UHdg!}\mathbb Z_U^{Hdg}:=a_{UHdg!}((O_U,F_b),\mathbb Z_{U^{an}},\alpha(U))  
\xrightarrow{=}(\int^{FDR}_{a_Y}j_{!Hdg}(O_U,F_b),(Ra_{U!}\mathbb Z_{U^{an}},W),a_{U*}\alpha(U)) \\ 
\xrightarrow{\sim}(\int^{FDR}_{a_Y}\Cone(\mathbb D^K_S\ad(i_{*mod},i^{\sharp})(-):(O_Y,F_b)\to i_{*mod}(O_Z,F_b)),
(Ra_{U!}\mathbb Z_{U^{an}},W),a_{U!}\alpha(U)) \\
\xrightarrow{\sim}(\Cone(E(\Omega_{D/Y})(D):\Gamma(Y,E(\Omega^{\bullet}_Y,F_b))\to\Gamma(Z,E(\Omega^{\bullet}_Z,F_b)),W),
(Ra_{U!}\mathbb Z_{U^{an}},W),a_{U!}\alpha(U))
\end{eqnarray*}
\item[(ii)] We have
\begin{eqnarray*}
a_{UHdg*}\mathbb Z_U^{Hdg}:=a_{UHdg*}((O_U,F),\mathbb Z_{U^{an}},\alpha(U))  
\xrightarrow{=}(\int^{FDR}_{a_Y}j_{*Hdg}(O_U,F,W),(Ra_{U*}\mathbb Z_{U^{an}},W),a_{U!}\alpha(U)) \\ 
\xrightarrow{\sim}(\int^{FDR}_{a_Y}\Cone(\ad(i_{*mod},i^{\sharp})(-):i_{*mod}(O_Z,F_b)[c]\to (O_Y,F_b)),
(Ra_{U*}\mathbb Z_{U^{an}},W),a_{U*}\alpha(U)) \\
\xrightarrow{\sim}(\Cone(i_{Z*}:\Gamma(Z,E(\Omega^{\bullet}_Z,F_b))(-c)[-2c]\to\Gamma(Y,E(\Omega^{\bullet}_Y,F_b),
(Ra_{U*}\mathbb Z_{U^{an}},W),a_{U*}\alpha(U))
\end{eqnarray*}
\end{itemize}
\end{prop}

\begin{proof}
See \cite{Saito}.
\end{proof}

In the case where $D=\cup D_i\subset Y$ is a normal crossing divisor, proposition \ref{jD} gives
\begin{eqnarray*}
a_{HdgU*}\mathbb Z_U^{Hdg}\xrightarrow{\sim}(\Gamma(Y,E(\Omega^{\bullet}_Y(\log D),F,W)),
(Ra_{U*}\mathbb Z_{U^{an}},W),a_{U*}\alpha(U))
\end{eqnarray*}
and
\begin{eqnarray*}
a_{HdgU!}\mathbb Z_U^{Hdg}&:=&(\Gamma(Y,E(\Omega^{\bullet}_Y(\nul D),F,W)),
(Ra_{U!}\mathbb Z_{U^{an}},W),a_{U!}\alpha(U))
\end{eqnarray*}

We recall the definition of the Deligne complex of a complex manifold and 
the Deligne cohomology class of an algebraic cycle of a complex algebraic variety.
\begin{defi}\label{Delkdef}
\begin{itemize}
\item[(i)] Let $X\in\AnSm(\mathbb C)$. We have for $d\in\mathbb Z$ the Deligne complex
\begin{equation*}
\mathbb Z_{\mathcal D,X}(d):=(\mathbb Z_X(d)\hookrightarrow\tau^{\leq d}DR(X))
=(\mathbb Z(d)\hookrightarrow (O_X\to\cdots\to\Omega_X^{d-1})\in C(X)
\end{equation*}
Let $D\subset X$ a normal crossing divisor. We have for $d\in\mathbb Z$ the Deligne complexes
\begin{equation*}
\mathbb Z_{\mathcal D,(X,D)}(d):=(\mathbb Z_X(d)\hookrightarrow(O_X\to\cdots\to\Omega_X^{d-1}(\log D)))\in C(X)
\end{equation*}
and
\begin{equation*}
\mathbb Z_{\mathcal D,(X,D)}(d)^{\vee}:=(\mathbb Z_X(d)\hookrightarrow(O_X\to\cdots\to\Omega_X^{d-1}(\null D)))\in C(X).
\end{equation*}
Moreover we have (see \cite{DP}) canonical products 
\begin{itemize}
\item $(-)\cdot(-):\mathbb Z_{\mathcal D,(X,D)}(d)\otimes\mathbb Z_{\mathcal D,(X,D)}(d')\to\mathbb Z_{\mathcal D,(X,D)}(d+d')$
\item $(-)\cdot(-):\mathbb Z_{\mathcal D,(X,D)}(d)^{\vee}\otimes\mathbb Z_{\mathcal D,(X,D)}(d')^{\vee}
\to\mathbb Z_{\mathcal D,(X,D)}(d+d')^{\vee}$
\end{itemize}
\item[(ii)] Let $X\in\AnSm(\mathbb C)$. We have for $d\in\mathbb Z$ the Deligne complex
\begin{eqnarray*}
C^{\bullet}_{\mathcal D}(X,\mathbb Z(d)):=
\Cone(\mathbb Z\Hom_{Diff(\mathbb R)}(\Delta^{\bullet},X)\oplus\Gamma(X,F^d\mathcal D^{\bullet}_X)
\hookrightarrow\Gamma(X,\mathcal D^{\bullet}_X))\in C(\mathbb C)
\end{eqnarray*}
Let $D\subset X$ a normal crossing divisor. Denote $U:=X\backslash D$. 
We have for $d\in\mathbb Z$ the Deligne complexes
\begin{eqnarray*}
C^{\bullet}_{\mathcal D}((X,D),\mathbb Z(d)):=
\Cone(\mathbb Z\Hom_{Diff(\mathbb R)}(\Delta^{\bullet},U)\oplus\Gamma(X,F^d\mathcal D^{\bullet}_X(\log D))
\hookrightarrow\Gamma(X,\mathcal D^{\bullet}_X(\log D)))\in C(\mathbb C)
\end{eqnarray*}
and
\begin{eqnarray*}
C^{\bullet}_{\mathcal D}(X,D,\mathbb Z(d)):=
\Cone(\mathbb Z\Hom_{Diff(\mathbb R)}(\Delta^{\bullet},(X,D))\oplus\Gamma(X,F^d\mathcal D^{\bullet}_X(\null D))
\hookrightarrow\Gamma(X,\mathcal D^{\bullet}_X(\null D)))\in C(\mathbb C).
\end{eqnarray*}
\item[(iii)] Let $X\in\PSmVar(\mathbb C)$. We have, for $k\in\mathbb Z$ and $d\in\mathbb Z$, the Deligne cohomology
\begin{eqnarray*}
H_{\mathcal D}^k(X^{an},\mathbb Z(d)):=\mathbb H^k(X^{an},\mathbb Z_{X,\mathcal D}(d))
=H^kC^{\bullet}_{\mathcal D}(X^{an},D,\mathbb Z(d))^{\vee}
\end{eqnarray*}
Let $U\in\SmVar(\mathbb C)$. Let $X\in\PSmVar(\mathbb C)$ a compactification of $U$ with $D:=X\backslash U$ a normal crossing divisor.
We have, for $k\in\mathbb Z$ and $d\in\mathbb Z$, the Deligne cohomology
\begin{eqnarray*}
H_{\mathcal D}^k(U^{an},\mathbb Z(d)):=
\mathbb H^k(X,\mathbb Z_{(X^{an},D^{an}),\mathcal D}(d))
=H^kC^{\bullet}_{\mathcal D}((X^{an},D^{an}),\mathbb Z(d))^{\vee}
\end{eqnarray*}
and
\begin{eqnarray*}
H_{\mathcal D}^k(X,D,\mathbb Z(d)):=
\mathbb H^k(X^{an},\mathbb Z_{(X^{an},D^{an}),\mathcal D}(d)^{\vee})
=H^kC^{\bullet}_{\mathcal D}(X^{an},D^{an},\mathbb Z(d))^{\vee}.
\end{eqnarray*}
\item[(iv)] Let $U\in\SmVar(\mathbb C)$. 
Let $X\in\PSmVar(\mathbb C)$ a compactification of $U$ with $D:=X\backslash U$ a normal crossing divisor.
We define the Deligne cohomology of a (higher) cycle $Z\in\mathcal Z^d(U,n)^{\partial=0}$ by
\begin{eqnarray*}
[Z]_{\mathcal D}:=\Im(H^{2d-n}(\gamma_{\supp(Z)})([Z])), \\ 
H^k(\gamma_{\supp(Z)}):\mathbb H^{2d-n}_{\supp(Z)}(X^{an},\mathbb Z_{X^{an},D^{an}}(d))
\to\mathbb H^{2d-n}(X^{an},\mathbb Z_{X^{an},D^{an}}(d)) 
\end{eqnarray*}
with $\supp(Z):=p_X(\supp(Z))\subset X$, where $\supp(Z)\subset X\times\square^n$ is the support of $Z$.
\item[(v)]Let $U\in\SmVar(\mathbb C)$. 
Let $X\in\PSmVar(\mathbb C)$ a compactification of $U$ with $D:=X\backslash U$ a normal crossing divisor.
We have for $d\in\mathbb Z$ the morphism of complexes
\begin{eqnarray*}
\mathcal R^d_U:\mathcal Z^d(U,\bullet)\to C^{\bullet}_{\mathcal D}(X^{an},D^{an},\mathbb Z(d)), \; 
Z\mapsto\mathcal R^d_U(Z):=(T_{\bar Z},\Omega_{\bar Z},R_{\bar Z})
\end{eqnarray*}
which gives for $Z\in\mathcal Z^d(U,n)^{\partial=0}$, 
\begin{equation*}
[\mathcal R^d_U(Z)]=[Z]_{\mathcal D}\in H^{2d-n}(U^{an},\mathbb Z(d))
\end{equation*}
\end{itemize}
\end{defi}

\begin{thm}\label{Delk}
\begin{itemize}
\item[(i)] Let $U\in\SmVar(\mathbb C)$. Denote by $a_U:U\to\pt$ the terminal map. 
Let $X\in\PSmVar(\mathbb C)$ a compactification of $U$ with $D:=X\backslash U$ a normal crossing divisor.
The embedding (see theorem \ref{Be})
\begin{equation*}
\iota:D(MHM(\left\{\pt\right\}))\to D_{fil}(\mathbb C)\times_ID(\mathbb Z) 
\end{equation*}
induces for $k\in\mathbb Z$ and $d\in\mathbb Z$, canonical isomorphisms
\begin{eqnarray*}
\iota(a_{U!Hdg}\mathbb Z^{Hdg}_U):H^k(a_{U!Hdg}\mathbb Z^{Hdg}_U)
\xrightarrow{\sim}H^k_{\mathcal D}(X^{an},D^{an},\mathbb Z(d)), 
\; \mbox{and} \\
\iota(a_{U*Hdg}\mathbb Z^{Hdg}_U):H^k(a_{U*Hdg}\mathbb Z^{Hdg}_U)\xrightarrow{\sim}
H^k_{\mathcal D}(U^{an},\mathbb Z(d)).
\end{eqnarray*}
\item[(ii)] Let $h:U\to S$ and $h':U'\to S$ two morphism with $S,U,U'\in\SmVar(\mathbb C)$.
Let $X\in\PSmVar(\mathbb C)$ a compactification of $U$ with $D:=X\backslash U$ a normal crossing divisor
such that $h:U\to S$ extend to $f:X\to\bar S$.
Let $X'\in\PSmVar(\mathbb C)$ a compactification of $U'$ with $D':=X'\backslash U'$ a normal crossing divisor
such that $h':U'\to S$ extend to $f':X'\to\bar S$.
The embedding $\iota:D(MHM(\pt))\to D_{fil}(k)\times_ID(\mathbb Z)$ (see theorem \ref{Be}) 
induces for $k\in\mathbb Z$ and $d\in\mathbb Z$ a canonical isomorphism
\begin{eqnarray*}
\iota(a_{U'\times_SU!Hdg}\mathbb Z^{Hdg}_{U'\times_SU}):
\Hom_{D(MHM(S))}(h_{U'!Hdg}\mathbb Z^{Hdg}_{U'},h_{U!Hdg}\mathbb Z^{Hdg}_U(d)[k]) \\
\xrightarrow{RI(-,-)}\Hom_{D(MHM_{gm,k}(\pt))}(\mathbb Z^{Hdg}_{\pt},a_{U'\times_SU!Hdg}\mathbb Z^{Hdg}_{U'\times_SU}(d)[k])
=H^k(a_{U'\times_SU!Hdg}\mathbb Z^{Hdg}_{U'\times_SU}(d)) \\
\xrightarrow{\sim}
H^k_{\mathcal D}((X'\times_SX)^{an},((X'\times_SU)\cup(U'\times_SX))^{an},\mathbb Z(d)).
\end{eqnarray*}
\item[(iii)]Let $U\in\SmVar(\mathbb C)$.  
Let $X\in\PSmVar(\mathbb C)$ a compactification of $U$ with $D:=X\backslash U$ a normal crossing divisor.
For $[Z]\in\CH^d(U,n)$ and $[Z']\in\CH^{d'}(U,n')$, we have
\begin{equation*}
([Z]\cdot[Z'])_{\mathcal D}=[Z]_{\mathcal D}\cdot[Z']_{\mathcal D}\in H^{2d+2d'-n-n'}(U^{an},\mathbb Z(d+d'))
\end{equation*}
where the product on the left is the intersection of higher Chow cycle which is well defined modulo boundary 
(they intersect properly modulo boundary) while the right product of Deligne cohomology classes is induced by 
the product of Deligne complexes 
$(-)\cdot(-):\mathbb Z_{\mathcal D,(X,D)}(d)\otimes\mathbb Z_{\mathcal D,(X,D)}(d')\to\mathbb Z_{\mathcal D,(X,D)}(d+d')$.
\item[(iv)]Let $h:U\to S$,$h':U'\to S$, $h'':U''\to S$ three morphism with $S,U,U',U''\in\SmVar(\mathbb C)$.
Let $X\in\PSmVar(\mathbb C)$ a compactification of $U$ with $D:=X\backslash U$ a normal crossing divisor
such that $h:U\to S$ extend to $f:X\to\bar S$.
Let $X'\in\PSmVar(\mathbb C)$ a compactification of $U'$ with $D':=X'\backslash U'$ a normal crossing divisor
such that $h':U'\to S$ extend to $f':X'\to\bar S$.
Let $X'\in\PSmVar(\mathbb C)$ a compactification of $U'$ with $D':=X'\backslash U'$ a normal crossing divisor
such that $h':U'\to S$ extend to $f':X'\to\bar S$.
For $[Z]\in\CH^d(U\times_SU',n)$ and $[Z']\in\CH^{d'}(U'\times_SU'',n')$, we have
\begin{equation*}
([Z]\circ[Z'])_{\mathcal D}=[Z]_{\mathcal D}\circ[Z']_{\mathcal D}
\in H^{d''-n''}((U\times_SU'')^{an},\mathbb Z(d''-n''))
\end{equation*}
where the composition on the left is the composition of higher correspondence modulo boundary
while the composition on the right is given by (ii).
\end{itemize}
\end{thm}

\begin{proof}
\noindent(i):Standard.

\noindent(ii):Follows on the one hand from (i) and 
on the other hand the six functor formalism on the 2-functor 
$D(MHM(-)):\SmVar(\mathbb C)\to\TriCat$  gives the isomorphism $RI(-,-)$.

\noindent(iii):Standard.

\noindent(iv):Follows from (iii).

\end{proof}

\section{The algebraic and analytic filtered De Rham realizations for Voevodsky relative motives} 

\subsection{The algebraic filtered De Rham realization functor}

\subsubsection{The algebraic Gauss-Manin filtered De Rham realization functor and its transformation map with pullbacks}

Consider, for $S\in\Var(\mathbb C)$, the following composition of morphism in $\RCat$ (see section 2)
\begin{eqnarray*}
\tilde e(S):(\Var(\mathbb C)/S,O_{\Var(\mathbb C)/S})\xrightarrow{\rho_S}(\Var(\mathbb C)^{sm}/S,O_{\Var(\mathbb C)^{sm}/S})
\xrightarrow{e(S)}(S,O_S)
\end{eqnarray*}
with, for $X/S=(X,h)\in\Var(\mathbb C)/S$,
\begin{itemize}
\item $O_{\Var(\mathbb C)/S}(X/S):=O_X(X)$, 
\item $(\tilde e(S)^*O_S(X/S)\to O_{\Var(\mathbb C)/S}(X/S)):=(h^*O_S\to O_X)$.
\end{itemize}
and $O_{\Var(\mathbb C)^{sm}/S}:=\rho_{S*}O_{\Var(\mathbb C)/S}$, that is, 
for $U/S=(U,h)\in\Var(\mathbb C)^{sm}/S$, $O_{\Var(\mathbb C)^{sm}/S}(U/S):=O_{\Var(\mathbb C)/S}(U/S):=O_U(U)$

\begin{defi}\label{wtildewGM}
\begin{itemize}
\item[(i)]For $S\in\Var(\mathbb C)$, we consider the complexes of presheaves 
\begin{equation*} 
\Omega^{\bullet}_{/S}:=
\coker(\Omega_{O_{\Var(\mathbb C)/S}/\tilde e(S)^*O_S}:
\Omega^{\bullet}_{\tilde e(S)^*O_S}\to\Omega^{\bullet}_{O_{\Var(\mathbb C)/S}})
\in C_{O_S}(\Var(\mathbb C)/S) 
\end{equation*}
which is by definition given by 
\begin{itemize}
\item for $X/S$ a morphism $\Omega^{\bullet}_{/S}(X/S)=\Omega^{\bullet}_{X/S}(X)$
\item for $g:X'/S\to X/S$ a morphism, 
\begin{eqnarray*}
\Omega^{\bullet}_{/S}(g):=\Omega_{(X'/X)/(S/S)}(X'):
\Omega^{\bullet}_{X/S}(X)\to g^*\Omega_{X/S}(X')\to\Omega^{\bullet}_{X'/S}(X') \\
\omega\mapsto\Omega_{(X'/X)/(S/S)}(X')(\omega):=g^*(\omega):(\alpha\in\wedge^{k}T_{X'}(X')\mapsto\omega(dg(\alpha)))
\end{eqnarray*}
\end{itemize}
\item[(ii)] For $S\in\Var(\mathbb C)$, we consider the complexes of presheaves 
\begin{equation*} 
\Omega^{\bullet}_{/S}:=\rho_{S*}\tilde\Omega^{\bullet}_{/S}=
\coker(\Omega_{O_{\Var(\mathbb C)^{sm}/S}/e(S)^*O_S}:\Omega^{\bullet}_{e(S)^*O_S}\to\Omega^{\bullet}_{O_{\Var(\mathbb C)^{sm}/S}})
\in C_{O_S}(\Var(\mathbb C)^{sm}/S) 
\end{equation*}
which is by definition given by 
\begin{itemize}
\item for $U/S$ a smooth morphism $\Omega^{\bullet}_{/S}(U/S)=\Omega^{\bullet}_{U/S}(U)$
\item for $g:U'/S\to U/S$ a morphism, 
\begin{eqnarray*}
\Omega^{\bullet}_{/S}(g):=\Omega_{(U'/U)/(S/S)}(U'):
\Omega^{\bullet}_{U/S}(U)\to g^*\Omega_{U/S}(U')\to\Omega^{\bullet}_{U'/S}(U') \\
\omega\mapsto\Omega_{(U'/U)/(S/S)}(U')(\omega):=g^*(\omega):(\alpha\in\wedge^{k}T_{U'}(U')\mapsto\omega(dg(\alpha)))
\end{eqnarray*}
\end{itemize}
\end{itemize}
\end{defi}

\begin{rem}
For $S\in\Var(\mathbb C)$, $\Omega^{\bullet}_{/S}\in C(\Var(\mathbb C)/S)$ 
is by definition a natural extension of $\Omega^{\bullet}_{/S}\in C(\Var(\mathbb C)^{sm}/S)$. 
However $\Omega^{\bullet}_{/S}\in C(\Var(\mathbb C)/S)$ does NOT satisfy cdh descent.
\end{rem}

For a smooth morphism $h:U\to S$ with $S,U\in\SmVar(\mathbb C)$, 
the cohomology presheaves  $H^n\Omega^{\bullet}_{U/S}$ of the relative De Rham complex
\begin{equation*}
DR(U/S):=\Omega^{\bullet}_{U/S}:=\coker(h^*\Omega_S\to\Omega_U)\in C_{h^*O_S}(U)
\end{equation*}
for all $n\in\mathbb Z$, have a canonical structure of a complex of $h^*D_S$ modules given by the Gauss Manin connexion : 
for $S^o\subset S$ an open subset, $U^o=h^{-1}(S^o)$, 
$\gamma\in\Gamma(S^o,T_S)$ a vector field and $\hat\omega\in\Omega^{p}_{U/S}(U^o)^c$ a closed form, the action is given by
\begin{equation*}
\gamma\cdot[\hat\omega]=[\widehat{\iota(\tilde\gamma)\partial\omega}], 
\end{equation*}
$\omega\in\Omega^p_U(U^o)$ being a representative of $\hat\omega$ and 
$\tilde\gamma\in\Gamma(U^o,T_U)$ a relevement of $\gamma$ ($h$ is a smooth morphism), so that 
\begin{equation*}
DR(U/S):=\Omega^{\bullet}_{U/S}:=\coker(h^*\Omega_S\to\Omega_U)\in C_{h^*O_S,h^*\mathcal D}(U)
\end{equation*}
with this $h^*D_S$ structure. Hence we get $h_*\Omega^{\bullet}_{U/S}\in C_{O_S,\mathcal D}(S)$ considering this structure.
Since $h$ is a smooth morphism, $\Omega^p_{U/S}$ are locally free $O_U$ modules.

The point (ii) of the definition \ref{wtildew} above gives 
the object in $\DA(S)$ which will, for $S$ smooth, represent the algebraic Gauss-Manin De Rham realisation.
It is the class of an explicit complex of presheaves on $\Var(\mathbb C)^{sm}/S$. 

\begin{prop}\label{aetfibGM}
Let $S\in\Var(\mathbb C)$.
\begin{itemize}
\item[(i)] For $U/S=(U,h)\in\Var(\mathbb C)^{sm}/S$, we have $e(U)_*h^*\Omega^{\bullet}_{/S}=\Omega^{\bullet}_{U/S}$.
\item[(ii)] The complex of presheaves $\Omega^{\bullet}_{/S}\in C_{O_S}(\Var(\mathbb C)^{sm}/S)$ 
is $\mathbb A^1$ homotopic, in particular $\mathbb A^1$ invariant.
Note that however, for $p>0$, the complexes of presheaves $\Omega^{\bullet\geq p}$ are NOT $\mathbb A^1$ local.
On the other hand, $(\Omega^{\bullet}_{/S},F_b)$ admits transferts 
(recall that means $\Tr(S)_*\Tr(S)^*\Omega^p_{/S}=\Omega^p_{/S}$).
\item[(iii)] If $S$ is smooth, we get $(\Omega^{\bullet}_{/S},F_b)\in C_{O_Sfil,D_S}(\Var(\mathbb C)^{sm}/S)$
with the structure given by the Gauss Manin connexion. Note that however the $D_S$ structure on the cohomology groups
given by Gauss Main connexion does NOT comes from a structure of $D_S$ module structure on the filtered complex of $O_S$ module.
The $D_S$ structure on the cohomology groups satisfy a non trivial Griffitz transversality (in the non projection cases),
whereas the filtration on the complex is the trivial one.
\end{itemize}
\end{prop}

\begin{proof} 

\noindent(i): Let $h':V\to U$ a smooth morphism with $V\in\Var(\mathbb C)$.
We have then
\begin{equation*}
h^*\Omega^p_{/S}(V\xrightarrow{h'}U)=\Omega^p_{/S}(V\xrightarrow{h'}U\xrightarrow{h}S).
\end{equation*}
Hence, if $h':V\hookrightarrow U$ is in particular an open embedding, 
$h^*\Omega^p_{/S}(V\xrightarrow{h'}U)=\Omega^p_{U/S}(V)$.
This proves the equality.

\noindent(ii): We prove that $E_{et}(\Omega^{\bullet}_{/S},F_b)\in C_{O_Sfil}(\Var(\mathbb C)^{sm}/S)$ 
is $2$-filtered $\mathbb A^1_S$ invariant. We follow \cite{LW}.
Consider the map in $C(\Var(\mathbb C)^{sm}/S)$
\begin{equation*}
\phi:=\ad(p_a^*,p_{a*})(-):\Omega^{\bullet}_{/S}\to p_{a*}p_a^*\Omega^{\bullet}_{/S}
\end{equation*} 
which is given, for $U/S\in\Var(\mathbb C)^{sm}/S$ by 
\begin{equation*}
\ad(p_a^*,p_{a*})(-)(U/S)=\Omega_{(U\times\mathbb A^1/U)/(S/S)}(U\times\mathbb A^1): \\
\Omega^{\bullet}_{U/S}(U)\to\Omega^{\bullet}_{U\times\mathbb A^1/S}(U\times\mathbb A^1), \;
\omega\mapsto p^*\omega
\end{equation*} 
where $p:U\times\mathbb A^1\to U$ is the projection.
On the other hand consider the map in $C(\Var(\mathbb C)^{sm}/S)$
\begin{equation*}
\psi:=I_0^*:p_{a*}p_a^*\Omega^{\bullet}_{/S}\to\Omega^{\bullet}_{/S}
\end{equation*} 
given, for $U/S\in\Var(\mathbb C)^{sm}/S$ by 
\begin{equation*}
I_0^*(U/S):\Omega^{\bullet}_{U\times\mathbb A^1/S}(U\times\mathbb A^1)\to\Omega^{\bullet}_{U/S}(U), \;
\omega\mapsto i_0^*\omega
\end{equation*} 
where $i_0:U\hookrightarrow:U\times\mathbb A^1$ is closed embedding given by $i_0(x):=(x,0)$.
Then,
\begin{itemize}
\item we have $\phi\circ\psi=I$
\item considering the map in $\PSh(\mathbb N\times\Var(\mathbb C)^{sm}/S)$
\begin{eqnarray*}
H:p_{a*}p_a^*\Omega^{\bullet}_{/S}[1]\to p_{a*}p_a^*\Omega^{\bullet}_{/S}
\end{eqnarray*}
given for $U/S\in\Var(\mathbb C)^{sm}/S$ by 
\begin{eqnarray*}
H(U/S)\Omega^p_{U\times\mathbb A^1/S}(U\times\mathbb A^1)\to\Omega^{p-1}_{U\times\mathbb A^1/S}(U\times\mathbb A^1), \\
H(U/S)(p^*\omega\wedge q^*(f(s)ds))=(\int_0^tf(s)ds)p^*\omega, \, H(U/S)(p^*\omega\wedge q^*f)=0,
\end{eqnarray*}
note that $g(t)=\int_0^tf(s)ds$ is algebraic since $f\in O_{\mathbb A^1}(\mathbb A^1)$ is a polynomial, 
we have $\psi\circ\phi-I=\partial H+H\partial$.
\end{itemize}
This shows that
\begin{equation*}
\ad(p_a^*,p_{a*})(-):\Omega^{\bullet}_{/S}\to p_{a*}p_a^*\Omega^{\bullet}_{/S}
\end{equation*} 
is an homotopy equivalence whose inverse is $I^*_0$. Hence,
\begin{equation*}
\ad(p_a^*,p_{a*})(-):(\Omega^{\bullet}_{/S},F_b)\to p_{a*}p_a^*(\Omega^{\bullet}_{/S},F_b)
\end{equation*}
is a $2$-filtered homotopy equivalence whose inverse is
\begin{equation*}
I_0^*:p_{a*}p_a^*(\Omega^{\bullet}_{/S},F_b)\to(\Omega^{\bullet}_{/S},F_b).
\end{equation*} 

\noindent(iii):For $h:U\to S$ a smooth morphism with $U,S\in\SmVar(\mathbb C)$,
recall that the $h^*D_S(U)=D_S(h(U))$ structure on $H^p\Omega^{\bullet}_{/S}(U/S):=H^p\Omega^{\bullet}_{U/S}(U)$ is given by,
for $\hat\omega\in\Omega^p_{U/S}(U)^c$, $\gamma\cdot[\hat\omega]=[\widehat{\iota(\tilde\gamma)\partial\omega}]$, 
$\omega\in\Omega^p_U(U^o)$ being a representative of $\hat\omega$ and 
$\tilde\gamma\in\Gamma(U^o,T_U)$ a relevement of $\gamma$ ($h$ is a smooth morphism).
Now, if $g:V/S\to U/S$ is a morphism, where $h':V\to S$ is a smooth morphism with $V\in\SmVar(\mathbb C)$, we have 
\begin{equation*}
g^*(\gamma\cdot\hat\omega)=\widehat{g^*(\iota(\tilde\gamma)\partial\omega)}=
\widehat{\iota(\tilde\gamma)\partial g^*\omega}=\gamma\cdot(g^*\hat\omega)
\end{equation*}
that is $H^p\Omega^{\bullet}_{/S}(g):H^p\Omega^{\bullet}(U/S)\to H^p\Omega^{\bullet}(V/S)$ is a map of $D_S(h(U))$ modules.
\end{proof}

We have the following canonical transformation map given by the pullback of (relative) differential forms:

Let $g:T\to S$ a morphism with $T,S\in\Var(\mathbb C)$.
Consider the following commutative diagram in $\RCat$ :
\begin{equation*}
D(g,e):\xymatrix{
(\Var(\mathbb C)^{sm}/T,O_{\Var(\mathbb C)^{sm}/T})\ar[rr]^{P(g)}\ar[d]^{e(T)} & \, & 
(\Var(\mathbb C)^{sm}/S,O_{\Var(\mathbb C)^{sm}/S})\ar[d]^{e(S)} \\
(T,O_T)\ar[rr]^{P(g)} & \, & (S,O_S)}
\end{equation*}
It gives (see section 2) the canonical morphism in $C_{g^*O_Sfil}(\Var(\mathbb C)^{sm}/T)$ 
\begin{eqnarray*}
\Omega_{/(T/S)}:=\Omega_{(O_{\Var(\mathbb C)^{sm}/T}/g^*O_{\Var(\mathbb C)^{sm}/S})/(O_T/g^*O_S}): \\
g^*(\Omega^{\bullet}_{/S},F_b)=\Omega^{\bullet}_{g^*O_{\Var(\mathbb C)^{sm}/S}/g^*e(S)^*O_S}\to
(\Omega^{\bullet}_{/T},F_b)=\Omega^{\bullet}_{O_{\Var(\mathbb C)^{sm}/T}/e(T)^*O_T}
\end{eqnarray*}
which is by definition given by the pullback on differential forms : for $(V/T)=(V,h)\in\Var(\mathbb C)^{sm}/T$,
\begin{eqnarray*}
\Omega_{/(T/S)}(V/T): 
g^*(\Omega^{\bullet}_{/S})(V/T):=\lim_{(h':U\to S \mbox{sm},g':V\to U,h,g)}\Omega^{\bullet}_{U/S}(U)
\xrightarrow{\Omega_{(V/U)/(T/S)}(V/T)}\Omega^{\bullet}_{V/T}(V)=:\Omega^{\bullet}_{/T}(V/T) \\
\hat\omega\mapsto\Omega_{(V/U)/(T/S)}(V/T)(\omega):=\hat{g^{'*}\omega}.
\end{eqnarray*}
If $S$ and $T$ are smooth, $\Omega_{/(T/S)}:g^*(\Omega^{\bullet}_{/S},F_b)\to(\Omega^{\bullet}_{/T},F_b)$
is a map in $C_{g^*O_Sfil,g^*D_S}(\Var(\mathbb C)^{sm}/T)$
It induces the canonical morphisms in $C_{g^*O_Sfil,g^*D_S}(\Var(\mathbb C)^{sm}/T)$:
\begin{eqnarray*}
E\Omega_{/(T/S)}:g^*E_{et}(\Omega^{\bullet}_{/S},F_b)\xrightarrow{T(g,E_{et})(\Omega^{\bullet}_{/S},F_b)}
E_{et}(g^*(\Omega^{\bullet}_{/S},F_b))\xrightarrow{E_{et}(\Omega_{/(T/S)})}E_{et}(\Omega^{\bullet}_{/T},F_b). 
\end{eqnarray*}
and
\begin{eqnarray*}
E\Omega_{/(T/S)}:g^*E_{zar}(\Omega^{\bullet}_{/S},F_b)\xrightarrow{T(g,E_{zar})(\Omega^{\bullet}_{/S},F_b)}
E_{zar}(g^*(\Omega^{\bullet}_{/S},F_b))\xrightarrow{E_{zar}(\Omega_{/(T/S)})}E_{zar}(\Omega^{\bullet}_{/T},F_b). 
\end{eqnarray*}

\begin{defi}\label{TgDRGM}
\begin{itemize}
\item[(i)]Let $g:T\to S$ a morphism with $T,S\in\Var(\mathbb C)$.
We have, for $F\in C(\Var(\mathbb C)^{sm}/S)$, the canonical transformation in $C_{O_Tfil}(T)$ :
\begin{eqnarray*}
T^O(g,\Omega_{/\cdot})(F): 
g^{*mod}L_Oe(S)_*\mathcal Hom^{\bullet}(F,E_{et}(\Omega^{\bullet}_{/S},F_b)) \\
\xrightarrow{:=}
(g^*L_Oe(S)_*\mathcal Hom^{\bullet}(F,E_{et}(\Omega^{\bullet}_{/S},F_b)))\otimes_{g^*O_S}O_T \\
\xrightarrow{T(e,g)(-)\circ T(g,L_O)(-)}  
L_O(e(T)_*g^*\mathcal Hom^{\bullet}(F,E_{et}(\Omega^{\bullet}_{/S},F))\otimes_{g^*O_S}O_T) \\ 
\xrightarrow{T(g,hom)(F,E_{et}(\Omega^{\bullet}_{/S}))\otimes I} 
L_O(e(T)_*\mathcal Hom^{\bullet}(g^*F,g^*E_{et}(\Omega^{\bullet}_{/S},F_b))\otimes_{g^*O_S}O_T) \\
\xrightarrow{ev(hom,\otimes)(-,-,-)} 
L_Oe(T)_*\mathcal Hom^{\bullet}(g^*F,g^*E_{et}(\Omega^{\bullet}_{/S},F_b)\otimes_{g^*e(S)^*O_S}e(T)^*O_T) \\
\xrightarrow{\mathcal Hom^{\bullet}(g^*F,E\Omega_{/(T/S)}\otimes I)} 
L_Oe(T)_*\mathcal Hom^{\bullet}(g^*F,E_{et}(\Omega^{\bullet}_{/T},F_b)\otimes_{g^*e(S)^*O_S}e(T)^*O_T) \\
\xrightarrow{m}
L_Oe(T)_*\mathcal Hom^{\bullet}(g^*F,E_{et}(\Omega^{\bullet}_{/T},F_b)
\end{eqnarray*}
where $m(\alpha\otimes h):=h.\alpha$ is the multiplication map.
\item[(ii)] Let $g:T\to S$ a morphism with $T,S\in\Var(\mathbb C)$, $S$ smooth.
Assume there is a factorization $g:T\xrightarrow{l}Y\times S\xrightarrow{p_S}S$
with $Y\in\SmVar(\mathbb C)$, $l$ a closed embedding and $p_S$ the projection.
We have, for $F\in C(\Var(\mathbb C)^{sm}/S)$, the canonical transformation in $C_{O_Tfil}(Y\times S)$ :
\begin{eqnarray*}
T(g,\Omega_{/\cdot})(F):
g^{*mod,\Gamma}e(S)_*\mathcal Hom^{\bullet}(F,E_{et}(\Omega^{\bullet}_{/S},F_b)) \\
\xrightarrow{:=} 
\Gamma_TE_{zar}(p_S^{*mod}e(S)_*\mathcal Hom^{\bullet}(F,E_{et}(\Omega^{\bullet}_{/S},F_b))) \\
\xrightarrow{T^O(p_S,\Omega_{/\cdot})(F)} 
\Gamma_TE_{zar}(e(T\times S)_*\mathcal Hom^{\bullet}(p_S^*F,E_{et}(\Omega^{\bullet}_{/Y\times S},F_b))) \\
\xrightarrow{=}
e(T\times S)_*\Gamma_T(\mathcal Hom^{\bullet}(p_S^*F,E_{et}(\Omega^{\bullet}_{/Y\times S},F_b))) \\
\xrightarrow{I(\gamma,\hom)(-,-)}
e(T\times S)_*\mathcal Hom^{\bullet}(\Gamma^{\vee}_Tp_S^*F,E_{et}(\Omega^{\bullet}_{/Y\times S},F_b)).
\end{eqnarray*}
For $Q\in Proj\PSh(\Var(\mathbb C)^{sm}/S)$,
\begin{eqnarray*}
T(g,\Omega_{/\cdot})(Q):
g^{*mod,\Gamma}e(S)_*\mathcal Hom^{\bullet}(Q,E_{et}(\Omega^{\bullet}_{/S},F_b))\to
e(T\times S)_*\mathcal Hom^{\bullet}(\Gamma^{\vee}_Tp_S^*Q,E_{et}(\Omega^{\bullet}_{/Y\times S},F_b))
\end{eqnarray*}
is a map in $C_{O_Tfil,\mathcal D}(Y\times S)$.
\end{itemize}
\end{defi}

The following easy lemma describe these transformation map on representable presheaves :

\begin{lem}\label{homQomegaGM}
Let $g:T\to S$ a morphism with $T,S\in\Var(\mathbb C)$ and $h:U\to S$ is a smooth morphism with $U\in\Var(\mathbb C)$.
Consider a commutative diagram whose square are cartesian :
\begin{equation*}
\xymatrix{g:T\ar[r]^l & S\times Y\ar[r]^{p_S} & S \\
g':U_T\ar[r]^{l'}\ar[u]^{h'} & U\times Y\ar[r]^{p_U}\ar[u]^{h'':=h\times I} & U\ar[u]^h}
\end{equation*} 
with $l$, $l'$ the graph embeddings and $p_S$, $p_U$ the projections. 
Then $g^*\mathbb Z(U/S)=\mathbb Z(U_T/T)$ and
\begin{itemize}
\item[(i)]we have the following commutative diagram in $C_{O_Tfil}(T)$ 
(see definition \ref{TDw} and definition \ref{TgDRGM}(i)) :
\begin{equation*}
\xymatrix{g^{*mod}L_Oe(S)_*\mathcal Hom^{\bullet}(\mathbb Z(U/S),E_{et}(\Omega^{\bullet}_{/S},F_b))
\ar[rrr]^{T(g,\Omega_{/\cdot})(\mathbb Z(U/S))} & \, & \, &
e(T)_*\mathcal Hom^{\bullet}(\mathbb Z(U_T/T),E_{et}(\Omega^{\bullet}_{/T},F_b)) \\
g^{*mod}L_Oe(S)_*\mathcal Hom^{\bullet}(\mathbb Z(U/S),E_{zar}(\Omega^{\bullet}_{/S},F_b))\ar[u]_k\ar[d]_{=}
\ar[rrr]^{T(g,\Omega_{/\cdot})(\mathbb Z(U/S))} & \, & \, &
e(T)_*\mathcal Hom^{\bullet}(\mathbb Z(U_T/T),E_{zar}(\Omega^{\bullet}_{/T},F_b))\ar[u]^k\ar[d]^{=} \\
g^{*mod}L_Oh_*E_{zar}(\Omega^{\bullet}_{U/S},F_b)\ar[rrr]^{T_{\omega}^{mod}(g,h)} & \, & \, & 
h'_*E_{zar}(\Omega^{\bullet}_{U_T/T},F_b)}
\end{equation*}
\item[(ii)]if $Y,S\in\SmVar(\mathbb C)$, we have the following commutative diagram in $C_{O_Tfil,\mathcal D}(Y\times S)$ 
(see definition \ref{TDw} and definition \ref{TgDRGM}(ii)) :
\begin{equation*}
\xymatrix{g^{*mod,\Gamma}e(S)_*\mathcal Hom^{\bullet}(\mathbb Z(U/S),E_{et}(\Omega^{\bullet}_{/S},F_b))
\ar[rrr]^{T(g,\Omega_{/\cdot})(\mathbb Z(U/S))} & \, & \, &
e(Y\times S)_*\mathcal Hom^{\bullet}(\Gamma_{U_T}^{\vee}\mathbb Z(U\times Y/S\times Y),
E_{et}(\Omega^{\bullet}_{/Y\times S},F_b)) \\
g^{*mod,\Gamma}L_Oe(S)_*\mathcal Hom^{\bullet}(\mathbb Z(U/S),E_{zar}(\Omega^{\bullet}_{/S},F_b))\ar[u]_k\ar[d]_{=}
\ar[rrr]^{T(g,\Omega_{/\cdot})(\mathbb Z(U/S))} & \, & \, &
e(Y\times S)_*\mathcal Hom^{\bullet}(\Gamma^{\vee}_{U_T}\mathbb Z(U\times Y/S\times Y),
E_{zar}(\Omega^{\bullet}_{/Y\times S},F_b))\ar[u]^k\ar[d]^{=} \\
g^{*mod,\Gamma}h_*E_{zar}(\Omega^{\bullet}_{U/S},F_b)\ar[rrr]^{T_{\omega}^O(p_S,h)(-)}& \, & \, & 
h''_*\Gamma_{U_T}E_{zar}(\Omega^{\bullet}_{U\times Y/S\times Y},F_b)}
\end{equation*}
where $j:T\backslash T\times S\hookrightarrow T\times S$ is the open complementary embedding,
\end{itemize}
with 
\begin{equation*}
k:E_{zar}(h^*\Omega^{\bullet}_{/S},F_b)\to E_{et}(E_{zar}(h^*\Omega^{\bullet}_{/S},F_b))=E_{et}(h^*\Omega^{\bullet}_{/S},F_b).
\end{equation*}
which is a filtered Zariski local equivalence.
\end{lem}

\begin{proof}
The commutative diagram follows from Yoneda lemma and proposition \ref{aetfibGM}(i).
On the other hand, $k:E_{zar}(\Omega^{\bullet}_{/S},F_b)\to E_{et}(\Omega^{\bullet}_{/S},F_b)$
is a ($1$-)filtered Zariski local equivalence by theorem \ref{DDADM} and proposition \ref{aetfibGM}(ii)
\end{proof}

In the projection case, we have the following :

\begin{prop}\label{TpDRQGM}
Let $p:S_{12}\to S_1$ is a smooth morphism  with $S_1,S_{12}\in\AnSp(\mathbb C)$. 
Then if $Q\in C(\Var(\mathbb C)^{sm}/S_1)$ is projective,
\begin{equation*}
T(p,\Omega_{/\cdot})(Q):
p^{*mod}e(S_1)_*\mathcal Hom^{\bullet}(Q,E_{et}(\Omega^{\bullet}_{/S_1},F_b))\to
e(S_{12})_*\mathcal Hom^{\bullet}(p^*Q,E_{et}(\Omega^{\bullet}_{/S_{12}},F_b))
\end{equation*}
is an isomorphism.
\end{prop}

\begin{proof}
Follows from lemma \ref{homQomegaGM} and base change by smooth morphisms of quasi-coherent sheaves.
\end{proof}

Let $S\in\Var(\mathbb C)$ and $h:U\to S$ a morphism with $U\in\Var(\mathbb C)$.
We then have the canonical map given by the wedge product
\begin{equation*}
w_{U/S}:\Omega^{\bullet}_{U/S}\otimes_{O_S}\Omega^{\bullet}_{U/S}\to\Omega^{\bullet}_{U/S}; 
\alpha\otimes\beta\mapsto\alpha\wedge\beta.
\end{equation*}
Let $S\in\Var(\mathbb C)$ and $h_1:U_1\to S$, $h_2:U_2\to S$ two morphisms with $U_1,U_2\in\Var(\mathbb C)$.
Denote $h_{12}:U_{12}:=U_1\times_S U_2\to S$ and $p_{112}:U_1\times_S U_2\to U_1$, $p_{212}:U_1\times_S U_2\to U_2$ the projections.
We then have the canonical map given by the wedge product
\begin{equation*}
w_{(U_1,U_2)/S}:p_{112}^*\Omega^{\bullet}_{U_1/S}\otimes_{O_S}p_{212}^*\Omega^{\bullet}_{U_2/S}
\to \Omega^{\bullet}_{U_{12}/S}; \alpha\otimes\beta\mapsto p_{112}^*\alpha\wedge p_{212}^*\beta
\end{equation*}
which gives the map
\begin{eqnarray*}
Ew_{(U_1,U_2)/S}:h_{1*}E_{zar}(\Omega^{\bullet}_{U_1/S})\otimes_{O_S}h_{2*}E_{zar}(\Omega^{\bullet}_{U_2/S}) \\
\xrightarrow{\ad(p_{112}^*,p_{112*})(-)\otimes\ad(p_{212}^*,p_{212*})(-)} 
(h_{1*}p_{112*}p_{112}^*E_{zar}(\Omega^{\bullet}_{U_1/S}))\otimes_{O_S}(h_{2*}p_{212*}p_{212}^*E_{zar}(\Omega^{\bullet}_{U_2/S})) \\
\xrightarrow{=}
h_{12*}(p_{112}^*E_{zar}(\Omega^{\bullet}_{U_1/S})\otimes_{h_{12}^*O_S}p_{212}^*E_{zar}(\Omega^{\bullet}_{U_2/S}) \\
\xrightarrow{T(\otimes,E)(-)\circ(T(p_{112},E)(-)\otimes T(p_{212},E)(-))}
h_{12*}E_{zar}(p_{112}^*\Omega^{\bullet}_{U_1/S}\otimes_{O_S}p_{212}^*\Omega^{\bullet}_{U_2/S})
\end{eqnarray*}

Let $S\in\Var(\mathbb C)$. We have the canonical map in $C_{O_Sfil}(\Var(\mathbb C)^{sm}/S)$
\begin{eqnarray*}
w_S:(\Omega^{\bullet}_{/S},F_b)\otimes_{O_S}(\Omega^{\bullet}_{/S},F_b)\to(\Omega^{\bullet}_{/S},F_b)
\end{eqnarray*}
given by for $h:U\to S\in\Var(\mathbb C)^{sm}/S$
\begin{eqnarray*}
w_S(U/S):(\Omega^{\bullet}_{U/S},F_b)\otimes_{h^*O_S}(\Omega^{\bullet}_{U/S},F_b)(U)\xrightarrow{w_{U/S}(U)}
(\Omega^{\bullet}_{U/S},F_b)(U)
\end{eqnarray*}
It gives the map
\begin{eqnarray*}
Ew_S:E_{et}(\Omega^{\bullet}_{/S},F_b)\otimes_{O_S}E_{et}(\Omega^{\bullet}_{/S},F_b)\xrightarrow{=}
E_{et}((\Omega^{\bullet}_{/S},F_b)\otimes_{O_S}(\Omega^{\bullet}_{/S},F_b))\xrightarrow{E_{et}(w_S)}
E_{et}(\Omega^{\bullet}_{/S},F_b)
\end{eqnarray*}
If $S\in\SmVar(\mathbb C)$, 
\begin{eqnarray*}
w_S:(\Omega^{\bullet}_{/S},F_b)\otimes_{O_S}(\Omega^{\bullet}_{/S},F_b)\to(\Omega^{\bullet}_{/S},F_b)
\end{eqnarray*}
is a map in $C_{O_Sfil,D_S}(\Var(\mathbb C)^{sm}/S)$.

\begin{defi}\label{TotimesDRGM}
Let $S\in\Var(\mathbb C)$.
We have, for $F,G\in C(\Var(\mathbb C)^{sm}/S)$, the canonical transformation in $C_{O_Sfil}(S)$ :
\begin{eqnarray*}
T(\otimes,\Omega)(F,G):
e(S)_*\mathcal Hom(F,E_{et}(\Omega^{\bullet}_{/S},F_b))\otimes_{O_S}e(S)_*\mathcal Hom(G,E_{et}(\Omega^{\bullet}_{/S},F_b)) \\ 
\xrightarrow{=}
e(S)_*(\mathcal Hom(F,E_{et}(\Omega^{\bullet}_{/S},F_b))\otimes_{O_S}\mathcal Hom(G,E_{et}(\Omega^{\bullet}_{/S},F_b))) \\
\xrightarrow{e(S)_*T(\mathcal Hom,\otimes)(-)}
e(S)_*\mathcal Hom(F\otimes G,E_{et}(\Omega^{\bullet}_{/S},F_b)\otimes_{O_S}E_{et}(\Omega^{\bullet}_{/S},F_b)) \\
\xrightarrow{\mathcal Hom(F\otimes G,Ew_S)}
e(S)_*\mathcal Hom(F\otimes G,E_{et}(\Omega^{\bullet}_{/S},F_b))
\end{eqnarray*}
If $S\in\SmVar(\mathbb C)$, $T(\otimes,\Omega)(F,G)$ is a map in $C_{O_Sfil,\mathcal D}(S)$. 
\end{defi}

\begin{lem}\label{homQomega2GM}
Let $S\in\Var(\mathbb C)$ and $h_1:U_1\to S$, $h_2:U_2\to S$ two smooth morphisms with $U_1,U_2\in\Var(\mathbb C)$.
Denote $h_{12}:U_{12}:=U_1\times_S U_2\to S$ and 
$p_{112}:U_1\times_S U_2\to U_1$, $p_{212}:U_1\times_S U_2\to U_2$ the projections.
We then have the following commutative diagram
\begin{equation*}
\xymatrix{
e(S)_*\mathcal Hom(F,E_{et}(\Omega^{\bullet}_{/S},F_b))\otimes_{O_S}e(S)_*\mathcal Hom(G,E_{et}(\Omega^{\bullet}_{/S},F_b))
\ar[rr]^{T(\otimes,\Omega)(F,G)} & \, & 
e(S)_*\mathcal Hom(F\otimes G,E_{et}(\Omega^{\bullet}_{/S},F)) \\
h_{1*}E_{zar}(\Omega^{\bullet}_{U_1/S},F_b)\otimes_{O_S}h_{2*}E_{zar}(\Omega^{\bullet}_{U_2/S},F_b)
\ar[rr]^{Ew_{(U_1,U_2)/S}}\ar[u]^{k} & \, & h_{12*}E_{zar}(\Omega^{\bullet}_{U_{12}/S},F_b)\ar[u]^{k}}
\end{equation*}
with 
\begin{equation*}
k:E_{zar}(\Omega^{\bullet}_{/S},F_b)\to E_{et}(E_{zar}(\Omega^{\bullet}_{/S},F_b))=E_{et}(\Omega^{\bullet}_{/S},F_b).
\end{equation*}
which is a filtered Zariski local equivalence.
\end{lem}

\begin{proof}
Follows from Yoneda lemma.
\end{proof}

Let $S\in\Var(\mathbb C)$ and $S=\cup_{i=1}^l S_i$ an open cover such that there exist 
closed embeddings $i_i:S_i\hookrightarrow\tilde S_i$ with $\tilde S_i\in\SmVar(\mathbb C)$. 
For $I\subset\left[1,\cdots l\right]$, denote by $S_I:=\cap_{i\in I} S_i$ and $j_I:S_I\hookrightarrow S$ the open embedding.
We then have closed embeddings $i_I:S_I\hookrightarrow\tilde S_I:=\Pi_{i\in I}\tilde S_i$.
Consider, for $I\subset J$, the following commutative diagram
\begin{equation*}
D_{IJ}=\xymatrix{ S_I\ar[r]^{i_I} & \tilde S_I \\
S_J\ar[u]^{j_{IJ}}\ar[r]^{i_J} & \tilde S_J\ar[u]^{p_{IJ}}}  
\end{equation*}
and $j_{IJ}:S_J\hookrightarrow S_I$ is the open embedding so that $j_I\circ j_{IJ}=j_J$.
Considering the factorization of the diagram $D_{IJ}$ by the fiber product :
\begin{equation*}
D_{IJ}=\xymatrix{
\tilde S_J=\tilde S_I\times\tilde S_{J\backslash I}\ar[rr]^{p_{IJ}} & \, & \tilde S_I \\
\, & S_I\times\tilde S_{J\backslash I}\ar[lu]^{i_I\times I}\ar[rd]^{p_{IJ}^0} & \, \\
S_J\ar[uu]^{i_J}\ar[ru]^{l_J}\ar[rr]^{j_{IJ}} & \, & S_I\ar[uu]^{i_I}}
\end{equation*}
the square of this factorization being cartesian, 
we have for $F\in C(\Var(\mathbb C)^{sm}/S)$ the canonical map in $C(\Var(\mathbb C)^{sm}/\tilde S_J)$
\begin{eqnarray*}
S(D_{IJ})(F):Li_{J*}j_J^*F\xrightarrow{q}i_{J*}j_J^*F=(i_I\times I)*l_{J*}j_J^*F 
\xrightarrow{(i_I\times I)_*\ad(p_{IJ\sharp}^o,p_{IJ}^{o*})(-)} \\
(i_I\times I)_*p_{IJ}^{o*}p_{IJ\sharp}^0l_{J*}j_J^*F 
\xrightarrow{T(p_{IJ},i_I)(-)^{-1}} p_{IJ}^*i_{I*}p_{IJ\sharp}^0l_{J*}j_I^*F=p_{IJ}^*i_{I*}j_I^*F
\end{eqnarray*}
which factors through
\begin{eqnarray*}
S(D_{IJ})(F):Li_{J*}j_I^*F\xrightarrow{S^q(D_{IJ})(F)}p_{IJ}^*Li_{I*}j_I^*F\xrightarrow{q}p_{IJ}^*i_{I*}j_I^*F
\end{eqnarray*}

\begin{defi}\label{DRalgdefFunctGM}
\begin{itemize}
\item[(i)]Let $S\in\SmVar(\mathbb C)$. We have the functor
\begin{eqnarray*}
C(\Var(\mathbb C)^{sm}/S)^{op}\to C_{Ofil,\mathcal D}(S), \; \;
F\mapsto e(S)_*\mathcal Hom^{\bullet}(L(i_{I*}j_I^*F),E_{et}(\Omega^{\bullet}_{/S},F_b))[-d_S].
\end{eqnarray*}
\item[(ii)]Let $S\in\Var(\mathbb C)$ and $S=\cup_{i=1}^l S_i$ an open cover such that there exist closed embeddings
$i_i:S_i\hookrightarrow\tilde S_i$  with $\tilde S_i\in\SmVar(\mathbb C)$. 
For $I\subset\left[1,\cdots l\right]$, denote by $S_I:=\cap_{i\in I} S_i$ and $j_I:S_I\hookrightarrow S$ the open embedding.
We then have closed embeddings $i_I:S_I\hookrightarrow\tilde S_I:=\Pi_{i\in I}\tilde S_i$.
We have the functor
\begin{eqnarray*}
C(\Var(\mathbb C)^{sm}/S)^{op}\to C_{Ofil,\mathcal D}(S/(\tilde S_I)), \; \;
F\mapsto(e(\tilde S_I)_*\mathcal Hom^{\bullet}(L(i_{I*}j_I^*F),
E_{et}(\Omega^{\bullet}_{/\tilde S_I},F_b))[-d_{\tilde S_I}],u^q_{IJ}(F))
\end{eqnarray*}
where
\begin{eqnarray*}
u^q_{IJ}(F)[d_{\tilde S_J}]:
e(\tilde S_I)_*\mathcal Hom^{\bullet}(L(i_{I*}j_I^*F),E_{et}(\Omega^{\bullet}_{/\tilde S_I},F_b)) \\
\xrightarrow{\ad(p_{IJ}^{*mod},p_{IJ*})(-)}
p_{IJ*}p_{IJ}^{*mod}e(\tilde S_I)_*\mathcal Hom^{\bullet}(L(i_{I*}j_I^*F),E_{et}(\Omega^{\bullet}_{/\tilde S_I},F_b)) \\
\xrightarrow{p_{IJ*}T(p_{IJ},\Omega_{\cdot})(L(i_{I*}j_I^*F))}
p_{IJ*}e(\tilde S_J)_*\mathcal Hom^{\bullet}(p_{IJ}^*L(i_{I*}j_I^*F),E_{et}(\Omega^{\bullet}_{/\tilde S_J},F_b)) \\
\xrightarrow{p_{IJ*}e(\tilde S_J)_*\mathcal Hom(S^q(D_{IJ})(F),E_{et}(\Omega_{/\tilde S_J}^{\bullet,\Gamma},F_b))}  
p_{IJ*}e(\tilde S_J)_*\mathcal Hom^{\bullet}(L(i_{J*}j_J^*F),E_{et}(\Omega^{\bullet}_{/\tilde S_J},F_b)). 
\end{eqnarray*}
For $I\subset J\subset K$, we have obviously $p_{IJ*}u_{JK}(F)\circ u_{IJ}(F)=u_{IK}(F)$.
\end{itemize}
\end{defi}

We then have the following key proposition

\begin{prop}\label{projwachGM}
\begin{itemize}
\item[(i)]Let $S\in\Var(\mathbb C)$. 
Let $m:Q_1\to Q_2$ be an equivalence $(\mathbb A^1,et)$ local in $C(\Var(\mathbb C)^{sm}/S)$
with $Q_1,Q_2$ complexes of projective presheaves. Then,
\begin{equation*}
e(S)_*\mathcal Hom(m,E_{et}(\Omega^{\bullet}_{/S},F_b)):
e(S)_*\mathcal Hom^{\bullet}(Q_2,E_{et}(\Omega^{\bullet}_{/S},F_b))\to 
e(S)_*\mathcal Hom^{\bullet}(Q_1,E_{et}(\Omega^{\bullet}_{/S},F_b))
\end{equation*}
is an $2$-filtered quasi-isomorphism. 
It is thus an isomorphism in $D_{O_Sfil,\mathcal D,\infty}(S)$ if $S$ is smooth.
\item[(ii)]Let $S\in\Var(\mathbb C)$. Let $S=\cup_{i=1}^l S_i$ an open cover such that there exist closed embeddings
$i_i:S_i\hookrightarrow\tilde S_i$  with $\tilde S_i\in\SmVar(\mathbb C)$.
Let $m=(m_I):(Q_{1I},s^1_{IJ})\to (Q_{2I},s^2_{IJ})$ be an equivalence $(\mathbb A^1,et)$ local 
in $C(\Var(\mathbb C)^{sm}/(\tilde S_I)^{op})$ with $Q_{1I},Q_{2I}$ complexes of projective presheaves. Then,
\begin{eqnarray*}
(e(\tilde S_I)_*\mathcal Hom(m_I,E_{et}(\Omega^{\bullet}_{/\tilde S_I},F_b))): \\
(e(\tilde S_I)_*\mathcal Hom^{\bullet}(Q_{2I},E_{et}(\Omega^{\bullet}_{/\tilde S_I},F_b)),u_{IJ}(Q_{2I},s^2_{IJ}))\to 
(e(\tilde S_I)_*\mathcal Hom^{\bullet}(Q_{1I},E_{et}(\Omega^{\bullet}_{/\tilde S_I},F_b)),u_{IJ}(Q_{1I},s^1_{IJ}))
\end{eqnarray*}
is an $2$-filtered quasi-isomorphism. It is thus an isomorphism in $D_{O_Sfil,\mathcal D,\infty}((\tilde S_I))$.
\end{itemize}
\end{prop}

\begin{proof}
\noindent(i): By definition of an $(\mathbb A^1,et)$ local equivalence (see proposition \ref{ca1Var}), there exist
\begin{equation*}
\left\{X_{1,\alpha}/S,\alpha\in\Lambda_1\right\},\ldots,\left\{X_{s,\alpha}/S,\alpha\in\Lambda_s\right\}
\subset\Var(\mathbb C)^{(sm)}/S 
\end{equation*}
such that we have in $\Ho_{et}(C(\Var(\mathbb C)^{(sm)}/S))$
\begin{eqnarray*}
\Cone(m)\xrightarrow{\sim}
\Cone(\oplus_{\alpha\in\Lambda_1}\Cone(\mathbb Z(X_{1,\alpha}\times\mathbb A^1/S)\to\mathbb Z(X_{1,\alpha}/S)) \\
\to\cdots\to\oplus_{\alpha\in\Lambda_s}\Cone(\mathbb Z(X_{s,\alpha}\times\mathbb A^1/S)\to\mathbb Z(X_{s,\alpha}/S)))
\end{eqnarray*}
This gives in $D_{fil}(\mathbb Z):=\Ho_{fil}(\mathbb Z)$,
\begin{eqnarray*}
\Cone(\Hom(m,E_{et}(\Omega^{\bullet}_{/S},F_b)))\xrightarrow{\sim}
\Cone(\oplus_{\alpha\in\Lambda_1}
\Cone(E_{et}(\Omega^{\bullet}_{/S},F_b)(X_{1,\alpha}/S)\to E_{et}(\Omega^{\bullet}_{/S},F_b)(X_{1,\alpha}\times\mathbb A^1/S)) \\
\to\cdots\to\oplus_{\alpha\in\Lambda_s}
\Cone(E_{et}(\Omega^{\bullet}_{/S},F_b)(X_{s,\alpha}/S)\to E_{et}(\Omega^{\bullet}_{/S},F_b)(X_{s,\alpha}\times\mathbb A^1/S)))
\end{eqnarray*}
Since $\Omega^{\bullet}_{/S}\in C(\Var(\mathbb C)^{sm}/S)$ is $\mathbb A^1$ homotopic,
for all $1\leq i\leq s$ and all $\alpha\in\Lambda_i$,
\begin{equation*}
\Cone(E_{et}(\Omega^{\bullet}_{/S})(X_{i,\alpha}/S)\to E_{et}(\Omega^{\bullet}_{/S})(X_{i,\alpha}\times\mathbb A^1/S))\to 0
\end{equation*}
are homotopy equivalence. 
Hence $\Cone(\Hom(m,E_{et}(G,F))\to 0$ is a $2$-filtered quasi-isomorphism.

\noindent(ii):Similar to (i) : see proposition \ref{ca1VarIJ}
\end{proof}

\begin{defi}\label{DRalgdefsingGM}
\begin{itemize}
\item[(i)] We define, using definition \ref{DRalgdefFunctGM}, by proposition \ref{projwachGM}, 
the filtered algebraic Gauss-Manin realization functor defined as
\begin{eqnarray*}
\mathcal F_S^{GM}:\DA_c(S)^{op}\to D_{O_Sfil,\mathcal D,\infty}(S), \; \;
M\mapsto\mathcal F_S^{GM}(M):=e(S)_*\mathcal Hom^{\bullet}(L(F),E_{et}(\Omega^{\bullet}_{/S},F_b))[-d_S]
\end{eqnarray*}
where $F\in C(\Var(\mathbb C)^{sm}/S)$ is such that $M=D(\mathbb A^1,et)(F)$,
\item[(ii)]Let $S\in\Var(\mathbb C)$ and $S=\cup_{i=1}^l S_i$ an open cover such that there exist closed embeddings
$i_i:S_i\hookrightarrow\tilde S_i$  with $\tilde S_i\in\SmVar(\mathbb C)$. 
For $I\subset\left[1,\cdots l\right]$, denote by $S_I=\cap_{i\in I} S_i$ and $j_I:S_I\hookrightarrow S$ the open embedding.
We then have closed embeddings $i_I:S_I\hookrightarrow\tilde S_I:=\Pi_{i\in I}\tilde S_i$.
We define, using definition \ref{DRalgdefFunctGM} and proposition \ref{projwachGM}
the filtered algebraic Gauss-Manin realization functor defined as
\begin{eqnarray*}
\mathcal F_S^{GM}:\DA_c(S)^{op}\to D_{Ofil,\mathcal D,\infty}(S/(\tilde S_I)), \; M\mapsto \\
\mathcal F_S^{GM}(M):=((e(\tilde S_I)_*\mathcal Hom^{\bullet}(L(i_{I*}j_I^*F),
E_{et}(\Omega^{\bullet}_{/\tilde S_I}),F_b))[-d_{\tilde S_I}],u^q_{IJ}(F))
\end{eqnarray*}
where $F\in C(\Var(\mathbb C)^{sm}/S)$ is such that $M=D(\mathbb A^1,et)(F)$.
\end{itemize}
\end{defi}

\begin{prop}\label{FDRwelldefGM}
For $S\in\Var(\mathbb C)$, the functor $\mathcal F_S^{GM}$ is well defined.
\end{prop}

\begin{proof}
Let $S\in\Var(\mathbb C)$ and $S=\cup_{i=1}^l S_i$ an open cover such that there exist closed embeddings
$i_i:S_i\hookrightarrow\tilde S_i$ with $\tilde S_i\in\SmVar(\mathbb C)$.
Denote, for $I\subset\left[1,\cdots, l\right]$, $S_I=\cap_{i\in I} S_i$ and $j_I:S_I\hookrightarrow S$ the open embedding.
We then have closed embeddings $i_I:S_I\hookrightarrow\tilde S_I:=\Pi_{i\in I}\tilde S_i$.
Let $M\in\DA(S)$. Let $F,F'\in C(\Var(\mathbb C)^{sm}/S)$ such that $M=D(\mathbb A_1,et)(F)=D(\mathbb A_1,et)(F')$.
Then there exist by definition a sequence of morphisms in $C(\Var(\mathbb C)^{sm}/S)$ :
\begin{equation*}
F=F_1\xrightarrow{s_1} F_2\xleftarrow{s_2} F_3\xrightarrow{s_3}F_4\to\cdots\xrightarrow{s_l} F'=F_s  
\end{equation*}
where, for $1\leq k\leq s$, and $s_k$ are $(\mathbb A^1,et)$ local equivalence.
But if $s:F_1\to F_2$ is an equivalence $(\mathbb A^1,et)$ local, 
\begin{equation*}
L(i_{I*}j_I^*s):L(i_{I*}j_I^*F_1)\to L(i_{I*}j_I^*F_2) 
\end{equation*}
is an equivalence $(\mathbb A^1,et)$ local, hence
\begin{eqnarray*}  
\mathcal Hom(L(i_{I*}j_I^*s),E_{et}(\Omega_{/\tilde S_I}^{\bullet},F_b)): 
(e(\tilde S_I)_*\mathcal Hom(L(i_{I*}j_I^*F_2),E_{et}(\Omega_{/\tilde S_I}^{\bullet},F_b)),u_{IJ}(F_2)) \\ 
\to(e(\tilde S_I)_*\mathcal Hom(L(i_{I*}j_I^*F_1),E_{et}(\Omega_{/\tilde S_I}^{\bullet},F_b)),u_{IJ}(F_1))
\end{eqnarray*}
is an $\infty$-filtered quasi-isomorphism by proposition \ref{projwachGM}.
\end{proof}

Let $f:X\to S$ a  morphism with $S,X\in\Var(\mathbb C)$. Assume that there is a factorization 
\begin{equation*}
f:X\xrightarrow{l}Y\times S\xrightarrow{p_S} S
\end{equation*}
of $f$, with $Y\in\SmVar(\mathbb C)$, $l$ a closed embedding and $p_S$ the projection.
Let $S=\cup_{i=1}^l S_i$ an open cover such that there exist closed embeddings
$i_i:S_i\hookrightarrow\tilde S_i$ with $\tilde S_i\in\SmVar(\mathbb C)$. 
We have $X=\cup_{i=1}^lX_i$ with $X_i:=f^{-1}(S_i)$. 
Denote, for $I\subset\left[1,\cdots l\right]$, $S_I=\cap_{i\in I} S_i$ and $X_I=\cap_{i\in I}X_i$.
For $I\subset\left[1,\cdots l\right]$, denote by $\tilde S_I=\Pi_{i\in I}\tilde S_i$,
We then have, for $I\subset\left[1,\cdots l\right]$, closed embeddings $i_I:S_I\hookrightarrow\tilde S_I$
and the following commutative diagrams which are cartesian 
\begin{equation*}
\xymatrix{
f_I=f_{|X_I}:X_I\ar[r]^{l_I}\ar[rd] & Y\times S_I\ar[r]^{p_{S_I}}\ar[d]^{i'_I} & S_I\ar[d]^{i_I} \\
\, & Y\times\tilde S_I\ar[r]^{p_{\tilde S_I}} & \tilde S_I} \;, \;
\xymatrix{Y\times\tilde S_J\ar[r]^{p_{\tilde S_J}}\ar[d]_{p'_{IJ}} & \tilde S_J\ar[d]^{p_{IJ}} \\
Y\times\tilde S_I\ar[r]^{p_{\tilde S_I}} & \tilde S_I}
\end{equation*}
with $l_I:l_{|X_I}$, $i'_I=I\times i_I$, $p_{S_I}$ and $p_{\tilde S_I}$ are the projections and $p'_{IJ}=I\times p_{IJ}$, 
and we recall that we denote by 
$j_I:\tilde S_I\backslash S_I\hookrightarrow\tilde S_I$ and $j'_I:Y\times\tilde S_I\backslash X_I\hookrightarrow Y\times S_I$
the open complementary embeddings.
We then have the commutative diagrams
\begin{equation*}
D_{IJ}=\xymatrix{S_J\ar[r]^{i_J}\ar[d]^{j_{IJ}} & \tilde S_J\ar[d]^{p_{IJ}} \\
S_I\ar[r]^{i_I} & \tilde S_I} \; , \;
D'_{IJ}=\xymatrix{X_J\ar[rr]^{i'_J\circ l_J}\ar[d]^{j'_{IJ}} & \, & Y\times\tilde S_J\ar[d]^{p'_{IJ}} \\
X_I\ar[rr]^{i'_I\circ l_I} & \, & Y\times\tilde S_I}.
\end{equation*}
and the factorization of $D'_{IJ}$ by the fiber product:
\begin{equation}
D'_{IJ}=
\xymatrix{ X_J\ar[r]^{i'_I\circ l_I}\ar[d]^{j'_{IJ}} & Y\times\tilde S_J\ar[d]^{p'_{IJ}} \\ 
X_I\ar[r]^{i'_I\circ l_I} & Y\times\tilde S_I}, \; \;
D'_{IJ}=
\xymatrix{ X_J\ar[rrr]^{i'_I\circ l_I}\ar[dd]^{j'_{IJ}}\ar[rrd]^{\hat l_J} & \, & \, & 
Y\times\tilde S_J\ar[dd]^{p'_{IJ}} \\ 
\,  & \, & X_I\times_{Y\times\tilde S_I}Y\times\tilde S_J=X_I\times\tilde S_{J\backslash I}
\ar[lld]^{\hat p_{IJ}}\ar[ru]^{\hat {il}_I}  & \, \\ 
X_I\ar[rrr]^{i'_I\circ l_I} & \, & \, & Y\times\tilde S_I}
\end{equation}
where $j'_{IJ}:X_J\hookrightarrow X_I$ is the open embedding.
Consider
\begin{equation*}
F(X/S):=p_{S,\sharp}\Gamma^{\vee}_X\mathbb Z(Y\times S/Y\times S)\in C(\Var(\mathbb C)^{sm}/S)
\end{equation*}
so that $D(\mathbb A^1,et)(F(X/S))=M(X/S)$. Then, by definition,
\begin{equation*}
\mathcal F_S^{GM}(M^{BM}(X/S)):=
(e(\tilde S_I)_*\mathcal Hom(L(i_{I*}j_I^*F(X/S)),E_{et}(\Omega_{/\tilde S_I},F_b))[-d_{\tilde S_I}],
u^q_{IJ}(F(X/S)))
\end{equation*}
On the other hand, let  
\begin{equation*}
Q(X_I/\tilde S_I):=p_{\tilde S_I,\sharp}\Gamma^{\vee}_{X_I}\mathbb Z(Y\times\tilde S_I/Y\times\tilde S_I)
\in C(\Var(\mathbb C)^{sm}/\tilde S_I), 
\end{equation*}
see definition \ref{projBMmotdef}.
We have then for $I\subset[1,l]$ the following map in $C(\Var(\mathbb C)^{sm}/\tilde S_J)$ :
\begin{eqnarray}\label{HYS}
N_I(X/S):Q(X_I/\tilde S_I)=p_{\tilde S_I\sharp}\Gamma^{\vee}_{X_I}\mathbb Z(Y\times\tilde S_I/Y\times\tilde S_I) 
\xrightarrow{p_{\tilde S_I\sharp}\ad(i^{'*}_I,i'_{I*})(-)} \\
p_{\tilde S_I\sharp}i'_{I*}i_I^{'*}\Gamma^{\vee}_{X_I}\mathbb Z(Y\times\tilde S_I/Y\times\tilde S_I)[d_Y] 
\xrightarrow{p_{\tilde S\sharp}(T(i'_I,\gamma^{\vee})(-))^{-1}}
p_{\tilde S_I\sharp}i'_{I*}\Gamma^{\vee}_{X_I}\mathbb Z(Y\times S_I/Y\times S_I)[d_Y] \\
\xrightarrow{\hat{T}_{\sharp}(p_{S_I},i_I)(-)}
i_{I*}p_{S_I\sharp}\Gamma^{\vee}_{X_I}\mathbb Z(Y\times S_I/Y\times S_I)[d_Y]=i_{I*}j_I^*F(X/S)
\end{eqnarray}
We have then for $I\subset J$ the following commutative diagram in $C(\Var(\mathbb C)^{sm}/\tilde S_J)$ :
\begin{equation}\label{HYSIJ}
\xymatrix{
p_{IJ}^*p_{\tilde S_I\sharp}\Gamma^{\vee}_{X_I}\mathbb Z(Y\times\tilde S_I/Y\times\tilde S_I)
\ar[rrr]^{p_{IJ}^*N_I(X/S)}
 & \, & \, & p_{IJ}^*(i_{I*}j_I^*F(X/S)) \\
p_{\tilde S_J\sharp}\Gamma^{\vee}_{X_J}\mathbb Z(Y\times\tilde S_J/Y\times\tilde S_J)
\ar[rrr]^{N_J(X/S)}\ar[u]^{H_{IJ}}
 & \, & \, & i_{J*}j_J^*F(X/S)\ar[u]_{S(D_{IJ})(F(X/S))}}
\end{equation}
with 
\begin{eqnarray*}
H_{IJ}:p_{\tilde S_J\sharp}\Gamma^{\vee}_{X_J}\mathbb Z(Y\times\tilde S_J/Y\times\tilde S_J) \\ 
\xrightarrow{=}
p_{\tilde S_J\sharp}p_{X\sharp}\Gamma^{\vee}_{X_J}p_{IJ}^{'*}\mathbb Z(Y\times\tilde S_I/Y\times\tilde S_I)
\xrightarrow{\Cone(\ad(p'_{IJ\sharp},p_{IJ}^{'*})(-),I)}
p_{\tilde S_J\sharp}\Gamma^{\vee}_{X_I\times\tilde S_{J\backslash I}}p_{IJ}^{'*}\mathbb Z(Y\times\tilde S_I/Y\times\tilde S_I) \\
\xrightarrow{T(p_{IJ},\gamma^{\vee})(-)}
p_{\tilde S_J\sharp}p_{IJ}^*\Gamma^{\vee}_{X_I}\mathbb Z(Y\times\tilde S_I/Y\times\tilde S_I) 
\xrightarrow{T_{\sharp}(p_{IJ},p_{\tilde S_I})(-)}
p_{IJ}^*p_{\tilde S_I\sharp}\Gamma^{\vee}_{X_I}\mathbb Z(Y\times\tilde S_I/Y\times\tilde S_I).
\end{eqnarray*}
This say that the maps $N_I(X/S)$ induces a map in $C(\Var(\mathbb C)^{sm}/(S/\tilde S_I))$
\begin{equation*}
(N_I(X/S)):(Q(X_I/\tilde S_I),H_{IJ})\to (i_{I*}j_I^*F(X/S),S(D_{IJ})(F(X/S))).
\end{equation*} 
We denote by $v_{IJ}^q(F(X/S))$ the composite
\begin{eqnarray*}
v_{IJ}^q(F(X/S))[d_{\tilde S_J}]:
e(\tilde S_I)_*\mathcal Hom(Q(X_I/\tilde S_I),E_{et}(\Omega^{\bullet}_{/\tilde S_I},F_b)) \\
\xrightarrow{\ad(p_{IJ}^{*mod},p_{IJ})(-)}
p_{IJ*}p_{IJ}^{*mod}e(\tilde S_I)_*\mathcal Hom(Q(X_I/\tilde S_I),E_{et}(\Omega^{\bullet}_{/\tilde S_I},F_b)) \\
\xrightarrow{p_{IJ*}T(p_{IJ},\Omega_{\cdot})(Q(X_I/\tilde S_I))}
p_{IJ*}e(\tilde S_J)_*\mathcal Hom(p_{IJ}^*Q(X_I/\tilde S_I),E_{et}(\Omega^{\bullet}_{/\tilde S_J},F_b)) \\
\xrightarrow{\mathcal Hom(H_{IJ},E_{et}(\Omega^{\bullet}_{/\tilde S_J},F_b))}
p_{IJ*}e(\tilde S_J)_*\mathcal Hom(Q(X_J/\tilde S_J),E_{et}(\Omega^{\bullet}_{/\tilde S_J},F_b)).
\end{eqnarray*} 
On the other hand, we have the following map in $C_{Ofil,\mathcal D,S_J}(\tilde S_J)$
\begin{eqnarray*}
w_{IJ}(X/S)[d_{\tilde S_J}]:
p_{\tilde S_I*}\Gamma_{X_I}E_{zar}(\Omega^{\bullet}_{Y\times\tilde S_I/\tilde S_I},F_b) 
\xrightarrow{\ad(p_{IJ}^{*mod},p_{IJ})(-)}
p_{IJ*}p_{IJ}^{*mod}p_{\tilde S_I*}\Gamma_{X_I}E_{zar}(\Omega^{\bullet}_{Y\times\tilde S_I/\tilde S_I},F_b) \\
\xrightarrow{T^O_w(p_{IJ},p_{\tilde S_I})^{\gamma}}
p_{IJ*}p_{\tilde S_J*}\Gamma_{X_I\times\tilde S_{J\backslash I}}E_{zar}(\Omega^{\bullet}_{Y\times\tilde S_I/\tilde S_I},F_b)
\xrightarrow{\Cone(I,\ad(p_{IJ}^{'*},p_{IJ*})(-))}
p_{\tilde S_J*}\Gamma_{X_J}E_{zar}(\Omega^{\bullet}_{Y\times\tilde S_J/\tilde S_J},F_b).
\end{eqnarray*}

\begin{lem}\label{keyalgsinglemGM}
\begin{itemize}
\item[(i)] The map in $C(\Var(\mathbb C)^{sm}/(S/\tilde S_I))$
\begin{equation*}
(N_I(X/S)):(Q(X_I/\tilde S_I),H_{IJ})\to(L(i_{I*}j_I^*F(X/S)),S^q(D_{IJ})(F(X/S))).
\end{equation*}
is an equivalence $(\mathbb A^1,et)$ local.
\item[(ii)] The maps $(N_I(X/S))$ induces an $\infty$-filtered quasi-isomorphism 
in $C_{Ofil,\mathcal D}(S/(\tilde S_I))$
\begin{eqnarray*}
(\mathcal Hom(N_I(X/S),E_{et}(\Omega^{\bullet}_{/\tilde S_I},F_b))): \\
(e(\tilde S_I)_*\mathcal Hom(L(i_{I*}j_I^*F(X/S)),E_{et}(\Omega^{\bullet}_{/\tilde S_I},F_b))[-d_{\tilde S_I}]
,u^q_{IJ}(F(X/S)))\to \\
(e(\tilde S_I)_*\mathcal Hom(Q(X_I/\tilde S_I),E_{et}(\Omega^{\bullet}_{/\tilde S_I},F_b))[-d_{\tilde S_I}],
v^q_{IJ}(F(X/S)))
\end{eqnarray*}
\item[(iii)] The maps $(I(\gamma,\hom)(-,-))$ and 
$(k:E_{zar}(p_{\tilde S_I}^*\Omega^{\bullet}_{/\tilde S_I},F_b)\to E_{et}(p_{\tilde S_I}^*\Omega^{\bullet}_{/\tilde S_I},F_b))$ 
induce an ($1$-)filtered Zariski local equivalence in $C_{Ofil,\mathcal D}(S/(\tilde S_I))$
\begin{eqnarray*}
(k\circ I(\gamma,\hom)(-,-)):
(p_{\tilde S_I*}\Gamma_{X_I}E_{zar}(\Omega^{\bullet}_{Y\times\tilde S_I/\tilde S_I},F_b)[-d_{\tilde S_I}],w_{IJ}(X/S)) \\
\to 
(e(\tilde S_I)_*\mathcal Hom(Q(X_I/\tilde S_I),E_{et}(\Omega^{\bullet}_{/\tilde S_I},F_b))[-d_{\tilde S_I}],v^q_{IJ}(F(X/S))) 
\end{eqnarray*}
\end{itemize}
\end{lem}

\begin{proof}
\noindent(i): Follows from theorem \ref{2functDM}.

\noindent(ii): These maps induce a morphism in $C_{Ofil,\mathcal D}(S/(\tilde S_I))$ by construction.
The fact that it is an $\infty$-filtered quasi-isomorphism follows from (i) and proposition \ref{projwachGM}.

\noindent(iii):These maps induce a morphism in $C_{Ofil,\mathcal D}(S/(\tilde S_I))$ by construction.
\end{proof}

\begin{prop}\label{keyalgsing1GM}
Let $f:X\to S$ a morphism with $S,X\in\Var(\mathbb C)$.
Let $S=\cup_{i=1}^l S_i$ an open cover such that there exist closed embeddings
$i_i:S_i\hookrightarrow\tilde S_i$ with $\tilde S_i\in\SmVar(\mathbb C)$. 
Then $X=\cup_{i=1}^lX_i$ with $X_i:=f^{-1}(S_i)$.
Denote, for $I\subset\left[1,\cdots l\right]$, $S_I=\cap_{i\in I} S_i$ and $X_I=\cap_{i\in I}X_i$.
Assume there exist a factorization 
\begin{equation*}
f:X\xrightarrow{l}Y\times S\xrightarrow{p_S} S
\end{equation*}
of $f$ with $Y\in\SmVar(\mathbb C)$, $l$ a closed embedding and $p_S$ the projection.
We then have, for $I\subset\left[1,\cdots l\right]$, the following commutative diagrams which are cartesian 
\begin{equation*}
\xymatrix{
f_I=f_{|X_I}:X_I\ar[r]^{l_I}\ar[rd] & Y\times S_I\ar[r]^{p_{S_I}}\ar[d]^{i'_I} & S_I\ar[d]^{i_I} \\
\, & Y\times\tilde S_I\ar[r]^{p_{\tilde S_I}} & \tilde S_I} \;, \;
\xymatrix{Y\times\tilde S_J\ar[r]^{p_{\tilde S_J}}\ar[d]_{p'_{IJ}} & \tilde S_J\ar[d]^{p_{IJ}} \\
Y\times\tilde S_I\ar[r]^{p_{\tilde S_I}} & \tilde S_I}
\end{equation*}
Let $F(X/S):=p_{S,\sharp}\Gamma_X^{\vee}\mathbb Z(Y\times S/Y\times S)$. 
The transformations maps $(N_I(X/S):Q(X_I/\tilde S_I)\to i_{I*}j_I^*F(X/S))$ and $(k\circ I(\gamma,\hom)(-,-))$, 
for $I\subset\left[1,\cdots,l\right]$, induce an isomorphism in $D_{Ofil,\mathcal D,\infty}(S/(\tilde S_I))$ 
\begin{eqnarray*} 
I^{GM}(X/S): \\
\mathcal F_S^{GM}(M(X/S)):=
(e(\tilde S_I)_*\mathcal Hom(L(i_{I*}j_I^*F(X/S)),E_{et}(\Omega^{\bullet}_{/\tilde S_I},F_b))[-d_{\tilde S_I}],
u^q_{IJ}(F(X/S))) \\
\xrightarrow{(e(\tilde S_I)_*\mathcal Hom(LN_I(X/S),E_{et}(\Omega^{\bullet}_{/\tilde S_I},F_b)))} 
(e(\tilde S_I)_*\mathcal Hom(Q(X_I/\tilde S_I),E_{et}(\Omega^{\bullet}_{/\tilde S_I},F_b))[-d_{\tilde S_I}],
v^q_{IJ}(F(X/S))) \\
\xrightarrow{(k\circ I(\gamma,\hom)(-,-))^{-1}}
(p_{\tilde S_I*}\Gamma_{X_I}E_{zar}(\Omega^{\bullet}_{Y\times\tilde S_I/\tilde S_I},F_b)[-d_{\tilde S_I}],
w_{IJ}(X/S)).
\end{eqnarray*}
\end{prop}

\begin{proof}
Follows from lemma \ref{keyalgsinglemGM}.
\end{proof}

We now define the functorialities of $\mathcal F_S^{GM}$ with respect to $S$ 
which makes $\mathcal F^{-}_{GM}$ a morphism of 2-functor.

\begin{defi}\label{TgDRdefGM}
Let $g:T\to S$ a morphism with $T,S\in\SmVar(\mathbb C)$.
Consider the factorization $g:T\xrightarrow{l}T\times S\xrightarrow{p_S}S$
where $l$ is the graph embedding and $p_S$ the projection.
Let $M\in\DA_c(S)$ and $F\in C(\Var(\mathbb C)^{sm}/S)$ such that $M=D(\mathbb A^1_S,et)(F)$. 
Then, $D(\mathbb A^1_T,et)(g^*F)=g^*M$. 
\begin{itemize}
\item[(i)]We have then the canonical transformation in $D_{Ofil,\mathcal D,\infty}(T\times S)$
(see definition \ref{TgDRGM}) :
\begin{eqnarray*}
T(g,\mathcal F^{GM})(M):Rg^{*mod[-],\Gamma}\mathcal F_S^{GM}(M):=
g^{*mod,\Gamma}e(S)_*\mathcal Hom^{\bullet}(LF,E_{et}(\Omega^{\bullet}_{/S},F_b)))[-d_T] \\
\xrightarrow{T(g,\Omega_{/\cdot})(LF)} \\ 
e(T\times S)_*\mathcal Hom^{\bullet}(\Gamma_T^{\vee}p_S^*LF,E_{et}(\Omega^{\bullet}_{/T\times S},F_b))[-d_T]
=:\mathcal F_{T\times S}^{GM}(l_*g^*(M,W)). 
\end{eqnarray*}
\item[(ii)]We have then the canonical transformation in $D_{Ofil,\infty}(T)$(see definition \ref{TgDRGM}) :
\begin{eqnarray*}
T^O(g,\mathcal F^{GM})(M,W):Lg^{*mod[-]}\mathcal F_S^{GM}(M):=
g^{*mod}e(S)_*\mathcal Hom^{\bullet}(LF,E_{et}(\Omega^{\bullet}_{/S},F_b)))[-d_T] \\
\xrightarrow{T^O(g,\Omega_{/\cdot})(LF)} \\ 
e(T)_*\mathcal Hom^{\bullet}(g^*LF,E_{et}(\Omega^{\bullet}_{/T},F_b))[-d_T]=:\mathcal F_T^{GM}(g^*M). 
\end{eqnarray*}
\end{itemize}
\end{defi}

We give now the definition in the non smooth case
Let $g:T\to S$ a morphism with $T,S\in\Var(\mathbb C)$.
Assume we have a factorization $g:T\xrightarrow{l}Y\times S\xrightarrow{p_S}S$
with $Y\in\SmVar(\mathbb C)$, $l$ a closed embedding and $p_S$ the projection.
Let $S=\cup_{i=1}^lS_i$ be an open cover such that 
there exists closed embeddings $i_i:S_i\hookrightarrow\tilde S_i$ with $\tilde S_i\in\SmVar(\mathbb C)$
Then, $T=\cup^l_{i=1} T_i$ with $T_i:=g^{-1}(S_i)$
and we have closed embeddings $i'_i:=i_i\circ l:T_i\hookrightarrow Y\times\tilde S_i$,
Moreover $\tilde g_I:=p_{\tilde S_I}:Y\times\tilde S_I\to\tilde S_I$ is a lift of $g_I:=g_{|T_I}:T_I\to S_I$.
We recall the commutative diagram :
\begin{equation*}
E_{IJg}=\xymatrix{(Y\times\tilde S_I)\backslash T_I\ar[d]^{p_{\tilde S_I}}\ar[r]^{m'_I} & 
Y\times\tilde S_J\ar[d]^{\tilde g_I} \\
\tilde S_I\backslash S_I\ar[r]^{m_I}& \tilde S_I}, \,
E_{IJ}=\xymatrix{\tilde S_J\backslash S_J\ar[d]^{p_{IJ}}\ar[r]^{m_J} & \tilde S_J\ar[d]^{p_{IJ}} \\
\tilde S_I\backslash(S_I\backslash S_J)\ar[r]^{m=m_{IJ}}& \tilde S_I}
E'_{IJ}=\xymatrix{(Y\times\tilde S_J)\backslash T_J\ar[d]^{p'_{IJ}}\ar[r]^{m'_J} & Y\times\tilde S_J\ar[d]^{p'_{IJ}} \\
(Y\times\tilde S_I)\backslash(T_I\backslash T_J)\ar[r]^{m'=m'_{IJ}}& Y\times\tilde S_I}
\end{equation*}
For $I\subset J$, denote by $p_{IJ}:\tilde S_J\to\tilde S_I$ and 
$p'_{IJ}:=I_Y\times p_{IJ}:Y\times\tilde S_J\to Y\times\tilde S_I$ the projections,
so that $\tilde g_I\circ p'_{IJ}=p_{IJ}\circ\tilde g_J$.
Consider, for $I\subset J\subset\left[1,\ldots,l\right]$,
resp. for each $I\subset\left[1,\ldots,l\right]$, the following commutative diagrams in $\Var(\mathbb C)$
\begin{equation*}
D_{IJ}=\xymatrix{ S_I\ar[r]^{i_I} & \tilde S_I \\
S_J\ar[u]^{j_{IJ}}\ar[r]^{i_J} & \tilde S_J\ar[u]^{p_{IJ}}} \; , \;  
D'_{IJ}=\xymatrix{ T_I\ar[r]^{i'_I} & Y\times\tilde S_I \\
T_J\ar[u]^{j'_{IJ}}\ar[r]^{i'_J} & Y\times\tilde S_J\ar[u]^{p'_{IJ}}}  
D_{gI}=\xymatrix{ S_I\ar[r]^{i_I} & \tilde S_I \\
T_I\ar[u]^{g_I}\ar[r]^{i'_I} & Y\times\tilde S_I\ar[u]^{\tilde g_I}} \; , \;  
\end{equation*}
and $j_{IJ}:S_J\hookrightarrow S_I$ is the open embedding so that $j_I\circ j_{IJ}=j_J$.
Let $F\in C(\Var(\mathbb C)^{sm}/S)$. 
Recall (see section 2) that since $j_I^{'*}i'_{I*}j_I^{'*}g^*F=0$,
the morphism $T(D_{gI})(j_I^*F):\tilde g_I^*i_{I*}j_I^*F\to i'_{I*}j_I^{'*}g^*F$ factors trough
\begin{equation*}
T(D_{gI})(j_I^*F):\tilde g_I^*i_{I*}j_I^*F\xrightarrow{\gamma_{X_I}^{\vee}(-)}
\Gamma_{X_I}^{\vee}\tilde g_I^*i_{I*}j_I^*F\xrightarrow{T^{\gamma}(D_{gI})(j_I^*F)}i'_{I*}j_I^{'*}g^*F
\end{equation*}
We then have, for each $I\subset\left[1,\ldots,l\right]$, the morphism 
\begin{eqnarray*}
T^{q,\gamma}(D_{gI})(j_I^*F):\Gamma^{\vee}_{T_I}\tilde g_I^*L(i_{I*}j_I^*F)\xrightarrow{T(\tilde g_I,L)(-)} \\
\Gamma^{\vee}_{T_I}L(\tilde g_I^*i_{I*}j_I^*F)=L(\Gamma^{\vee}_{T_I}\tilde g_I^*i_{I*}j_I^*F) 
\xrightarrow{L(T^{\gamma}(D_{gI})(j_I^*F))} L(i'_{I*}j_I^{'*}g^*F)
\end{eqnarray*}
and the following diagram in $C(\Var(\mathbb C)^{sm}/Y\times\tilde S_I)$ commutes
\begin{equation*}
\xymatrix{
\Gamma^{\vee}_{T_I}\tilde g_I^*L(i_{I*}j_I^*F)
\ar[rr]^{\Gamma^{\vee}_{T_I}\tilde g_I^*q_I(F)}\ar[d]_{T^{q,\gamma}(D_{gI})(j_I^*F)} & \, &
\Gamma^{\vee}_{T_I}\tilde g_I^*i_{I*}j_I^*F\ar[d]^{T^{\gamma}(D_{gI})(j_I^*F)} \\
L(i'_{I*}j_I^{'*}g^*F)\ar[rr]^{q_I(g^*F)} & \, & i'_{I*} g_I^*j_I^*F=i'_{I*}j_I^{'*}g^*F}
\end{equation*}  
We have the following commutative diagram in $C(\Var(\mathbb C)^{sm}/Y\times\tilde S_J)$
\begin{equation}\label{DgDI}
\xymatrix{
p^{'*}_{IJ}\tilde g_I^*i_{I*}j_I^*F=\tilde g_J^*p_{IJ}^*i_{I*}j_I^*F
\ar[rr]^{p^{'*}_{IJ}T(D_{gI})(j_I^*F)} & \, & 
p^{'*}_{IJ}i'_{I*}g_I^*j_I^*F=p^{'*}_{IJ}i'_{I*}j^{'*}_Ig^*F \\
\tilde g_J^*i_{J*}j_{IJ}^*j_I^*F=\tilde g_J^*i_{J*}j_J^*F
\ar[rr]^{T(D_{gJ})(j_J^*F)}\ar[u]^{\tilde g_J^*S(D_{IJ})(F)} & \, &
i'_{J*}g_J^*j_J^*F=i'_{J*}j^{'*}_{IJ}j^{'*}_Ig^*F=i'_{J*}j^{'*}_Jg^*F
\ar[u]_{S(D'_{IJ})(g^*F)}}
\end{equation}
This gives, after taking the functor $L$, the following commutative diagram in $C(\Var(\mathbb C)^{sm}/Y\times\tilde S_J)$
\begin{equation}\label{DgDIq}
\xymatrix{
\Gamma^{\vee}_{T_J}p^{'*}_{IJ}\Gamma^{\vee}_{T_I}\tilde g_I^*L(i_{I*}j_I^*F)=\Gamma^{\vee}_{T_J}\tilde g_J^*p_{IJ}^*L(i_{I*}j_I^*F)
\ar[rrr]^{p^{'*}_{IJ}T^{q,\gamma}(D_{gI})(j_I^*F)} & \, & \, & 
\Gamma^{\vee}_{T_J}p^{'*}_{IJ}L(i'_{I*}j_I^{'*}g^*F) \\
\Gamma^{\vee}_{T_J}\tilde g_J^*L(i_{J*}j_J^*g^*F)
\ar[rrr]^{T^{q,\gamma}(D_{gJ})(j_J^*F)}\ar[u]^{\tilde g_J^*S^q(D_{IJ})(F)} & \, & \, & 
L(i'_{J*}j_J^{'*}g^*F)\ar[u]_{S^q(D'_{IJ})(g^*F)}}
\end{equation}
The fact that the diagrams (\ref{DgDIq}) commutes says that the maps $T^{q,\gamma}(D_{gI})(j_I^*F)$
define a morphism in $C(\Var(\mathbb C)^{sm}/(T/(Y\times\tilde S_I)))$
\begin{equation*}
(T^{q,\gamma}(D_{gI})(j_I^*F)):(\Gamma^{\vee}_{T_I}\tilde g_I^*L(i_{I*}j_I^*F),\tilde g_J^*S^q(D_{IJ})(F))
\to (L(i'_{I*}j_I^{'*}g^*F),S^q(D'_{IJ})(g^*F))
\end{equation*}
We denote by $\tilde g_J^*u^q_{IJ}(F)_1$ the composite
\begin{eqnarray*}
\tilde g_J^*u^q_{IJ}(F)_1[d_Y+d_{\tilde S_I}]:
e(Y\times\tilde S_I)_*\Gamma_{T_I}\mathcal Hom(\tilde g_I^*L(i_{I*}j_I^*F),
E_{et}(\Omega^{\bullet}_{/Y\times\tilde S_I},F_b)) \\
\xrightarrow{\ad(p'_{IJ},p'_{IJ*})(-)} 
p'_{IJ*}p^{'*mod}_{IJ}e(Y\times\tilde S_I)_*\Gamma_{T_I}\mathcal Hom(\tilde g_I^*L(i_{I*}j_I^*F),
E_{et}(\Omega^{\bullet}_{/Y\times\tilde S_I},F_b)) \\
\xrightarrow{T^{mod}(p'_{IJ},\gamma)(-)}
p'_{IJ*}e(Y\times\tilde S_J)_*\Gamma_{T_I\times\tilde S_{J\backslash I}}
p_{IJ}^{'*mod}\mathcal Hom(\tilde g_J^*L(i_{J*}j_J^*F),E_{et}(\Omega^{\bullet}_{/Y\times\tilde S_J},F_b)) \\
\xrightarrow{\Cone(\ad(p'_{IJ\sharp},p_{IJ}^{'*})(-),I)}
p'_{IJ*}e(Y\times\tilde S_J)_*\Gamma_{T_J}p_{IJ}^{'*mod}\mathcal Hom(\tilde g_J^*L(i_{J*}j_J^*F),
E_{et}(\Omega^{\bullet}_{/Y\times\tilde S_J},F_b)) \\
\xrightarrow{T(p'_{IJ},\Omega_{/\cdot})(-)}
p'_{IJ*}e(Y\times\tilde S_J)_*\Gamma_{T_J}\mathcal Hom(\tilde g_J^*p_{IJ}^*L(i_{I*}j_I^*F),
E_{et}(\Omega^{\bullet}_{/Y\times\tilde S_J},F_b)) \\  
\xrightarrow{\mathcal Hom(\tilde g_J^*(S^q(D_{IJ})(F)),
E_{et}(\Omega^{\bullet}_{/Y\times\tilde S_J},F_b))}
p'_{IJ*}e(Y\times\tilde S_J)_*\Gamma_{T_J}\mathcal Hom(\tilde g_J^*L(i_{J*}j_J^*F),
E_{et}(\Omega^{\bullet}_{/Y\times\tilde S_J},F_b)) 
\end{eqnarray*}
We denote by $\tilde g_J^*u^q_{IJ}(F)_2$ the composite
\begin{eqnarray*}
\tilde g_J^*u^q_{IJ}(F)_2[d_Y+d_{\tilde S_I}]:
e(\tilde T_I)_*\mathcal Hom(\Gamma^{\vee}_{T_I}\tilde g_I^*L(i_{I*}j_I^*F),
E_{et}(\Omega^{\bullet}_{/Y\times\tilde S_I},F_b)) \\
\xrightarrow{\ad(p_{IJ}^{'*mod},p'_{IJ*})(-)}
p'_{IJ*}p^{'*mod}_{IJ}e(\tilde T_I)_*\mathcal Hom(\Gamma^{\vee}_{T_I}\tilde g_I^*L(i_{I*}j_I^*F),
E_{et}(\Omega^{\bullet}_{/Y\times\tilde S_I},F_b)) \\
\xrightarrow{T(p'_{IJ},\Omega_{/\cdot})(-)}
p'_{IJ*}e(\tilde T_I)_*\mathcal Hom(p_{IJ}^{'*}\Gamma^{\vee}_{T_I}\tilde g_I^*L(i_{I*}j_I^*F),
E_{et}(\Omega^{\bullet}_{/Y\times\tilde S_J},F_b)) \\
\xrightarrow{\mathcal Hom(T(p'_{IJ},\gamma^{\vee})(-),E_{et}(\Omega^{\bullet}_{/Y\times\tilde S_J},F_b))}   
p'_{IJ*}e(\tilde T_I)_*\mathcal Hom(\Gamma^{\vee}_{T_I\times\tilde S_{J\backslash I}}
p_{IJ}^{'*}\tilde g_I^*L(i_{I*}j_I^*F),E_{et}(\Omega^{\bullet}_{/Y\times\tilde S_J},F_b)) \\
\xrightarrow{\Cone(\ad(p'_{IJ\sharp},p_{IJ}^{'*})(-),I)}   
p'_{IJ*}e(\tilde T_I)_*\mathcal Hom(\Gamma^{\vee}_{T_J}p_{IJ}^{'*}\tilde g_I^*L(i_{I*}j_I^*F),
E_{et}(\Omega^{\bullet}_{/Y\times\tilde S_J},F_b)) \\
\xrightarrow{\mathcal Hom(\Gamma^{\vee}_{T_J}\tilde g_J^*(S^q(D_{IJ})(F)),
E_{et}(\Omega^{\bullet}_{/Y\times\tilde S_J},F_b))}
p'_{IJ*}e(Y\times\tilde S_J)_*\mathcal Hom(\Gamma^{\vee}_{T_J}\tilde g_J^*L(i_{J*}j_J^*F),
E_{et}(\Omega^{\bullet}_{/Y\times\tilde S_J},F_b)) 
\end{eqnarray*}

We then have then the following lemma :

\begin{lem}\label{TglemGM}
\begin{itemize}
\item[(i)]The morphism in $C(\Var(\mathbb C)^{sm}/(T/(Y\times\tilde S_I)))$
\begin{equation*}
(T^{q,\gamma}(D_{gI})(j_I^*F)):(\Gamma^{\vee}_{T_I}L\tilde g_I^*i_{I*}j_I^*F,\tilde g_J^*S^q(D_{IJ})(F))
\to (i'_{I*}j_I^{'*}g^*F,S^q(D'_{IJ})(F)(j_I^{'*}g^*F))
\end{equation*}
is an equivalence $(\mathbb A^1,et)$ local.
\item[(ii)] Denote for short $d_{YI}=-d_Y-d_{\tilde S_I}$. The maps 
$\mathcal Hom((T^{q,\gamma}(D_{gI})(j_I^*F)),E_{et}(\Omega^{\bullet}_{/Y\times\tilde S_I}),F_b))$
induce an $\infty$-filtered quasi-isomorphism in $C_{Ofil,\mathcal D}(T/(Y\times\tilde S_I))$
\begin{eqnarray*}
(\mathcal Hom(T^{q,\gamma}(D_{gI})(j_I^*F),E_{et}(\Omega^{\bullet}_{/Y\times\tilde S_I},F_b))): \\
(e(Y\times\tilde T_I)_*\mathcal Hom(L(i'_{I*}j_I^{'*}g^*F),
E_{et}(\Omega^{\bullet}_{/Y\times\tilde S_I},F_b))[d_{YI}],u^q_{IJ}(g^*F))\to \\ 
(e(Y\times\tilde T_I)_*\mathcal Hom((\Gamma^{\vee}_{T_I}L\tilde g_I^*i_{I*}j_I^*F),
E_{et}(\Omega^{\bullet}_{/Y\times\tilde S_I},F_b))[d_{YI}],\tilde g_J^*u^q_{IJ}(F)_2) 
\end{eqnarray*}
\item[(iii)] The maps $T(\tilde g_I,\Omega_{\cdot})(L(i_{I*}j_I^*F))$ (see definition \ref{TgDRGM})
induce a morphism in $C_{Ofil,\mathcal D}(T/(Y\times\tilde S_I))$
\begin{eqnarray*}
(T(\tilde g_I,\Omega_{/\cdot})(L(i_{I*}j_I^*F))): \\ 
(\Gamma_{T_I}E_{zar}(\tilde g_I^{*mod[-]}e(\tilde S_I)_*\mathcal Hom^{\bullet}(L(i_{I*}j_I^*F),
E_{et}(\Omega^{\bullet}_{/\tilde S_I},F_b)))[d_{YI}],\tilde g_J^{*mod}u^q_{IJ}(F))\to \\  
(\Gamma_{T_I}(e(Y\times\tilde S_I)_*\mathcal Hom(\tilde g_I^*L(i_{I*}j_I^*F),
E_{et}(\Omega^{\bullet}_{/Y\times\tilde S_I},F_b)))[d_{YI}],\tilde g_J^*u^q_{IJ}(F)_1). 
\end{eqnarray*}
\end{itemize}
\end{lem}

\begin{proof}

\noindent(i): Follows from theorem \ref{2functDM}.

\noindent(ii): These maps induce a morphism in $C_{Ofil,\mathcal D}(T/(Y\times\tilde S_I))$ by construction.
The fact that this morphism is an $\infty$-filtered equivalence Zariski local follows from (i) 
and proposition \ref{projwachGM}.

\noindent(iii): These maps induce a morphism in $C_{Ofil,\mathcal D}(T/(Y\times\tilde S_I))$ by construction.

\end{proof}

\begin{defi}\label{TgDRdefsingGM}
Let $g:T\to S$ a morphism with $T,S\in\Var(\mathbb C)$.
Assume we have a factorization $g:T\xrightarrow{l}Y\times S\xrightarrow{p_S}S$
with $Y\in\SmVar(\mathbb C)$, $l$ a closed embedding and $p_S$ the projection.
Let $S=\cup_{i=1}^lS_i$ be an open cover such that 
there exists closed embeddings $i_i:S_i\hookrightarrow\tilde S_i$ with $\tilde S_i\in\SmVar(\mathbb C)$
Then, $T=\cup^l_{i=1} T_i$ with $T_i:=g^{-1}(S_i)$
and we have closed embeddings $i'_i:=i_i\circ l:T_i\hookrightarrow Y\times\tilde S_i$,
Moreover $\tilde g_I:=p_{\tilde S_I}:Y\times\tilde S_I\to\tilde S_I$ is a lift of $g_I:=g_{|T_I}:T_I\to S_I$.
Denote for short $d_{YI}:=d_Y+d_{\tilde S_I}$.
Let $M\in\DA_c(S)$ and $F\in C(\Var(\mathbb C)^{sm}/S)$ such that $M=D(\mathbb A^1_S,et)(F)$.
Then, $D(\mathbb A^1_T,et)(g^*F)=g^*M$. 
We have, by lemma \ref{TglemGM}, the canonical transformation in $D_{Ofil,\mathcal D,\infty}(T/(Y\times\tilde S_I))$
\begin{eqnarray*}
T(g,\mathcal F^{GM})(M):Rg^{*mod[-],\Gamma}\mathcal F_S^{GM}(M):= \\
(\Gamma_{T_I}E_{zar}(\tilde g_I^{*mod}e(\tilde S_I)_*\mathcal Hom^{\bullet}(L(i_{I*}j_I^*F),
E_{et}(\Omega^{\bullet}_{/\tilde S_I},F_b)))[-d_Y-d_{\tilde S_I}],\tilde g_J^{*mod}u^q_{IJ}(F)) \\
\xrightarrow{(\Gamma_{T_I}E(T(\tilde g_I,\Omega_{/\cdot})(L(i_{I*}j_I^*(F,W)))))} \\
(\Gamma_{T_I}e(Y\times\tilde S_I)_*\mathcal Hom^{\bullet}(\tilde g_I^*L(i_{I*}j_I^*F),
E_{et}(\Omega^{\bullet}_{/Y\times\tilde S_I},F_b))[-d_Y-d_{\tilde S_I}],\tilde g_J^*u^q_{IJ}(F)_1) \\ 
\xrightarrow{(I(\gamma,\hom(-,-)))} \\
(e(Y\times\tilde S_I)_*\mathcal Hom^{\bullet}(\Gamma_{T_I}^{\vee}\tilde g_I^*L(i_{I*}j_I^*F),
E_{et}(\Omega^{\bullet}_{/Y\times\tilde S_I},F_b))[-d_Y-d_{\tilde S_I}],\tilde g_J^*u^q_{IJ}(F)_2) \\ 
\xrightarrow{(e(Y\times\tilde S_I)_*\mathcal Hom(T^{q,\gamma}(D_{gI})(j_I^*F),
E_{et}(\Omega^{\bullet}_{/Y\times\tilde S_I},F_b)))^{-1}} \\
(e(Y\times\tilde S_I)_*\mathcal Hom^{\bullet}(L(i'_{I*}j_I^{'*}g^*F),
E_{et}(\Omega^{\bullet}_{/Y\times\tilde S_I},F_b))[-d_Y-d_{\tilde S_I}],u^q_{IJ}(g^*F))=:\mathcal F_T^{GM}(g^*M).
\end{eqnarray*} 
\end{defi}

\begin{prop}\label{TgGMprop}
\begin{itemize}
\item[(i)]Let $g:T\to S$ a morphism with $T,S\in\Var(\mathbb C)$.
Assume we have a factorization $g:T\xrightarrow{l}Y_2\times S\xrightarrow{p_S}S$
with $Y_2\in\SmVar(\mathbb C)$, $l$ a closed embedding and $p_S$ the projection.
Let $S=\cup_{i=1}^lS_i$ be an open cover such that 
there exists closed embeddings $i_i:S_i\hookrightarrow\tilde S_i$ with $\tilde S_i\in\SmVar(\mathbb C)$
Then, $T=\cup^l_{i=1} T_i$ with $T_i:=g^{-1}(S_i)$
and we have closed embeddings $i'_i:=i_i\circ l:T_i\hookrightarrow Y_2\times\tilde S_i$,
Moreover $\tilde g_I:=p_{\tilde S_I}:Y\times\tilde S_I\to\tilde S_I$ is a lift of $g_I:=g_{|T_I}:T_I\to S_I$.
Let $f:X\to S$ a  morphism with $X\in\Var(\mathbb C)$. Assume that there is a factorization 
$f:X\xrightarrow{l}Y_1\times S\xrightarrow{p_S} S$, with $Y_1\in\SmVar(\mathbb C)$, 
$l$ a closed embedding and $p_S$ the projection. We have then the following commutative diagram
whose squares are cartesians
\begin{equation*}
\xymatrix{f':X_T\ar[r]\ar[d] & Y_1\times T\ar[d]\ar[r] & T\ar[d] \\
f''=f\times I:Y_2\times X\ar[r]\ar[d] & Y_1\times Y_2\times S\ar[r]\ar[d] & Y_2\times S\ar[d] \\
f:X\ar[r] & Y_1\times S\ar[r] & S}
\end{equation*} 
Consider $F(X/S):=p_{S,\sharp}\Gamma_X^{\vee}\mathbb Z(Y_1\times S/Y_1\times S)$ and
the isomorphism in $C(\Var(\mathbb C)^{sm}/S)$
\begin{eqnarray*}
T(f,g,F(X/S)):g^*F(X/S):=g^*p_{S,\sharp}\Gamma_X^{\vee}\mathbb Z(Y_1\times S/Y_1\times S)\xrightarrow{\sim} \\
p_{T,\sharp}\Gamma_{X_T}^{\vee}\mathbb Z(Y_1\times T/Y_1\times T)=:F(X_T/T).
\end{eqnarray*}
which gives in $\DA(S)$ the isomorphism $T(f,g,F(X/S)):g^*M(X/S)\xrightarrow{\sim}M(X_T/T)$. 
Then, the following diagram in $D_{Ofil,\mathcal D,\infty}(T/(Y_2\times\tilde S_I))$ commutes
\begin{equation*}
\begin{tikzcd}
Rg^{*mod,\Gamma}\mathcal F_S^{GM}(M(X/S))
\ar[rrr,"T(g{,}\mathcal F^{GM})(M(X/S))"]\ar[d,"I^{GM}(X/S)"] & \, & \, & 
\mathcal F_T^{GM}(M(X_T/T))\ar[d,"I^{GM}(X_T/T)"] \\
\shortstack{$g^{*mod[-],\Gamma}(p_{\tilde S_I*}\Gamma_{X_I}
E_{zar}(\Omega^{\bullet}_{Y_1\times\tilde S_I/\tilde S_I},F_b)[-d_{\tilde S_I}],$ \\ $w_{IJ}(X/S))$}
\ar[rrr,"(T(\tilde g_I\times I{,}\gamma)(-)\circ T^O_w(\tilde g_I{,}p_{\tilde S_I}))"] & \, & \, &
\shortstack{$(p_{Y_2\times\tilde S_I*}\Gamma_{X_{T_I}}
E_{zar}(\Omega^{\bullet}_{Y_2\times Y_1\times\tilde S_I/Y_2\times\tilde S_I},F_b)[-d_{Y_2}-d_{\tilde S_I}],$ \\ $w_{IJ}(X_T/T))$}.
\end{tikzcd}
\end{equation*} 
\item[(ii)] Let $g:T\to S$ a morphism with $T,S\in\SmVar(\mathbb C)$.
Let $f:X\to S$ a  morphism with $X\in\Var(\mathbb C)$. Assume that there is a factorization 
$f:X\xrightarrow{l}Y\times S\xrightarrow{p_S} S$, with $Y\in\SmVar(\mathbb C)$, $l$ a closed embedding and $p_S$ the projection.
Consider $F(X/S):=p_{S,\sharp}\Gamma_X^{\vee}\mathbb Z(Y\times S/Y\times S)$ and
the isomorphism in $C(\Var(\mathbb C)^{sm}/S)$
\begin{eqnarray*}
T(f,g,F(X/S)):g^*F(X/S):=g^*p_{S,\sharp}\Gamma_X^{\vee}\mathbb Z(Y\times S/Y\times S)\xrightarrow{\sim} \\
p_{T,\sharp}\Gamma_{X_T}^{\vee}\mathbb Z(Y\times T/Y\times T)=:F(X_T/T).
\end{eqnarray*}
which gives in $\DA(S)$ the isomorphism $T(f,g,F(X/S)):g^*M(X/S)\xrightarrow{\sim}M(X_T/T)$. 
Then, the following diagram in $D_{Ofil,\infty}(T)$ commutes
\begin{equation*}
\begin{tikzcd}
Lg^{*mod[-]}\mathcal F_S^{GM}(M(X/S))
\ar[rrr,"T^O(g{,}\mathcal F^{GM})(M(X/S))"]\ar[d,"I^{GM}(X/S)"] & \, & \, & 
\mathcal F_T^{GM}(M(X_T/T))\ar[d,"I^{GM}(X_T/T)"] \\
g^{*mod}L_O(p_{S*}\Gamma_{X}E_{zar}(\Omega^{\bullet}_{Y\times S/S},F_b)[-d_T]
\ar[d,"T_w(\otimes{,}\gamma)(O_{Y\times S})"]\ar[rrr,"(T(g\times I{,}\gamma)(-)\circ T^O_w(g{,}p_S))"] & \, & \, &
p_{Y\times T*}\Gamma_{X_T}E_{zar}(\Omega^{\bullet}_{Y\times T/T},F_b)[-d_T]
\ar[d,"T_w(\otimes{,}\gamma)(O_{Y\times T})"] \\
Lg^{*mod}\int^{FDR}_{p_S}\Gamma_XE(O_{Y\times S},F_b)[-d_Y-d_T]
\ar[rrr,"T^{\mathcal Dmod}(g{,}f)(\Gamma_XE(O_{Y\times S}{,}F_b))"] & \, & \, & 
\int^{FDR}_{p_T}\Gamma_{X_T}E(O_{Y\times T},F_b)[-d_Y-d_T].
\end{tikzcd}
\end{equation*}
\end{itemize}
\end{prop}

\begin{proof}
Follows immediately from definition.
\end{proof}

We have the following theorem:

\begin{thm}\label{mainthmGM}
\begin{itemize}
\item[(i)]Let $g:T\to S$ is a morphism with $T,S\in\Var(\mathbb C)$. 
Assume there exist a factorization $g:T\xrightarrow{l}Y\times S\xrightarrow{p_S}$
with $Y\in\SmVar(\mathbb C)$, $l$ a closed embedding and $p_S$ the projection. 
Let $S=\cup_{i=1}^lS_i$ be an open cover such that 
there exists closed embeddings $i_i:S_i\hookrightarrow\tilde S_i$ with $\tilde S_i\in\SmVar(\mathbb C)$.
Then, for $M\in\DA_c(S)$
\begin{eqnarray*}
T(g,\mathcal F^{GM})(M):Rg^{*mod[-],\Gamma}\mathcal F_S^{GM}(M)\xrightarrow{\sim}\mathcal F_T^{GM}(g^*M)
\end{eqnarray*}
is an isomorphism in $D_{O_Tfil,\mathcal D,\infty}(T/(Y\times\tilde S_I))$.
\item[(ii)]Let $g:T\to S$ is a morphism with $T,S\in\SmVar(\mathbb C)$. Then, for $M\in\DA_c(S)$
\begin{eqnarray*}
T^O(g,\mathcal F^{GM})(M):Lg^{*mod[-]}\mathcal F_S^{GM}(M)\xrightarrow{\sim}\mathcal F_T^{GM}(g^*M)
\end{eqnarray*}
is an isomorphism in $D_{O_T}(T)$.
\item[(iii)] A base change theorem for algebraic De Rham cohomology :
Let $g:T\to S$ is a morphism with $T,S\in\SmVar(\mathbb C)$. 
Let $h:U\to S$ a smooth morphism with $U\in\Var(\mathbb C)$. Then the map (see definition \ref{TDw})
\begin{equation*}
T^O_w(g,h):Lg^{*mod}Rh_*(\Omega^{\bullet}_{U/S},F_b)\xrightarrow{\sim} Rh'_*(\Omega^{\bullet}_{U_T/T},F_b)
\end{equation*}
is an isomorphism in $D_{O_T}(T)$.
\end{itemize}
\end{thm}

\begin{proof} 
\noindent(i):Follows from proposition \ref{TpDRQGM}.

\noindent(ii): Follows from proposition \ref{TgGMprop}(ii) 
and the base change for algebraic D modules (proposition \ref{PDmod1}).

\noindent(iii):Follows from (ii) and lemma \ref{homQomegaGM}.
\end{proof}

We finish this subsection by some remarks on the absolute case and on a particular case of the relative case:

\begin{prop}\label{FFXD}
\begin{itemize}
\item[(i)]Let $X\in\PSmVar(\mathbb C)$ and $D=\cup D_i\subset X$ a normal crossing divisor. 
Consider the open embedding $j:U:=X\backslash D\hookrightarrow X$. Then, 
\begin{itemize}
\item The map in $D_{fil,\infty}(\mathbb C)$ 
\begin{eqnarray*}
\Hom(L\mathbb D(\mathbb Z(U)),k)^{-1}\circ\Hom((0,\ad(j^*,j_*)(\mathbb Z(X/X))),E_{et}(\Omega^{\bullet},F_b)): \\
\mathcal F^{GM}(\mathbb D(\mathbb Z(U))):=\Hom(L\mathbb D(\mathbb Z(U)),E_{et}(\Omega^{\bullet},F_b)) \\
\xrightarrow{\sim}
\Hom(\Cone(\mathbb Z(D)\to\mathbb Z(X)),E_{zar}(\Omega^{\bullet},F_b))=\Gamma(X,E_{zar}(\Omega^{\bullet}_X(\nul D),F_b)).
\end{eqnarray*}
is an isomorphism, where we recall $\mathbb D(\mathbb Z(U):=a_{X*}j_*E_{et}(\mathbb Z(U/U))=a_{U*}E_{et}(\mathbb Z(U/U))$, 
\item $\mathcal F^{GM}(\mathbb Z(U))=\Gamma(U,E_{zar}(\Omega_U^{\bullet},F_b))\in D_{fil,\infty}(\mathbb C)$
is NOT isomorphic to $\Gamma(X,E_{zar}(\Omega^{\bullet}_X(\log D),F_b))$ in $D_{fil,\infty}(\mathbb C)$ in general.
For exemple $U$ is affine, then $H^n(U,\Omega_U^p)=0$ for all $p\in\mathbb N$, $p\neq 0$, 
so that the $E^{p,q}_{\infty}(\Gamma(U,E_{zar}(\Omega_U^{\bullet},F_b)))$ are NOT isomorphic to 
$E^{p,q}_{\infty}(\Gamma(X,E_{zar}(\Omega^{\bullet}_X(\log D),F_b)))$ in this case. 
In particular, the map,  
\begin{equation*}
j^*:=\ad(j^*,j_*)(-):H^n\Gamma(X,E_{zar}(\Omega^{\bullet}_X(\log D)))\xrightarrow{\sim}H^n\Gamma(U,E_{zar}(\Omega^{\bullet}_U)) 
\end{equation*}
which is an isomorphism in $D(\mathbb C)$ (i.e. if we forgot filtrations), gives embeddings  
\begin{eqnarray*}
j^*:=\ad(j^*,j_*)(-):F^pH^n(U,\mathbb C):=F^pH^n\Gamma(X,E_{zar}(\Omega^{\bullet}_X(\log D),F_b))
\hookrightarrow F^pH^n\Gamma(U,E_{zar}(\Omega^{\bullet}_U,F_b))
\end{eqnarray*}
which are NOT an isomorphism in general for $n,p\in\mathbb Z$. Note that, since $a_U:U\to\left\{\pt\right\}$ is not proper,
\begin{equation*}
[\Delta_U]:\mathbb Z(U)\to a_{U*}E_{et}(\mathbb Z(U/U))[2d_U]
\end{equation*}
is NOT an equivalence $(\mathbb A^1,et)$ local.
\item Let $Z\subset X$ a smooth subvariety and denote $U:=X\backslash Z$ the open complementary.
Denote $M_Z(X):=\Cone(M(U)\to M(X))\in\DA(\mathbb C)$. The map in $D_{fil,\infty}(\mathbb C)$
\begin{eqnarray*}
\Hom(G(X,Z),E_{et}(\Omega^{\bullet},F_b))^{-1}\circ\Hom(a_{X\sharp}\Gamma^{\vee}_Z\mathbb Z(X/X)),k)^{-1}: \\
\mathcal F^{GM}(M_Z(X)):=\Hom(a_{X\sharp}\Gamma^{\vee}_Z\mathbb Z(X/X)),E_{et}(\Omega^{\bullet},F_b))\xrightarrow{\sim} \\
\Gamma(X,\Gamma_ZE_{zar}(\Omega_X^{\bullet},F_b))=\Gamma_Z(X,E_{zar}(\Omega_X^{\bullet},F_b)) \\
\xrightarrow{\sim}\mathcal F^{GM}(M(Z)(c)[2c])=\Gamma(Z,E_{zar}(\Omega_Z^{\bullet},F_b))(-c)[-2c]
\end{eqnarray*}
is an isomorphism,
where $c=\codim(Z,X)$ and $G(X,Z):a_{X\sharp}\Gamma^{\vee}_Z\mathbb Z(X/X)\to\mathbb Z(Z)(c)[2c]$ is the Gynsin morphism.
\item Let $D\subset X$ a smooth divisor and denote $U:=X\backslash Z$ the open complementary
Note that the canonical distinguish triangle in $\DA(\mathbb C)$
\begin{equation*}
M(U)\xrightarrow{\ad(j_{\sharp},j^*)(\mathbb Z(X/X))}M(X)\xrightarrow{\gamma^{\vee}_Z(\mathbb Z(X/X))} M_D(X)\to M(U)[1]
\end{equation*}
give a canonical triangle in $D_{fil,\infty}(\mathbb C)$
\begin{eqnarray*}
\mathcal F^{GM}(M_D(X))\xrightarrow{\mathcal F^{GM}(\gamma^{\vee}_Z(\mathbb Z(X/X)))}
\mathcal F^{GM}(M(X))\xrightarrow{\mathcal F^{GM}(\ad(j_{\sharp},j^*)(\mathbb Z(X/X))} 
\mathcal F^{GM}(M(U))\to\mathcal F^{GM}(M_D(X))[1],
\end{eqnarray*}
which is NOT the image of a distinguish triangle in $\pi(D(MHM(\mathbb C)))$, 
as $\mathcal F^{GM}(M(U))\notin\pi(D(MHM(\mathbb C)))$ since the morphism
\begin{equation*}
\ad(j^*,j_*):H^n(X,E_{zar}(\Omega^{\bullet}_X(\log D),F_b)\to H^n(U,E_{zar}(\Omega^{\bullet}_U,F_b))
\end{equation*}
is not strict. Note that if $U:=S\backslash D$ is affine, then by the exact sequence in $C(\mathbb Z)$
\begin{equation*}
0\to \Gamma_Z(X,E_{zar}(\Omega^p_X))\to\Gamma(X,E_{zar}(\Omega^p_X))\to\Gamma(U,E_{zar}(\Omega^p_U))\to 0
\end{equation*}
we have $H^q\Gamma_Z(X,E_{zar}(\Omega^p_X))=H^q(\Gamma(X,E_{zar}(\Omega^p_X)))$.
In particular, the map,  
\begin{equation*}
j^*:=\ad(j^*,j_*)(-):\Gamma(X,E_{zar}(\Omega^{\bullet}_X(\log D),F_b))\to\Gamma(U,E_{zar}(\Omega^{\bullet}_U,F_b)) 
\end{equation*}
and hence the map
\begin{eqnarray*}
j^*:=\ad(j^*,j_*)(-):\Cone(\Gamma(X,\Omega^{\bullet}_X,F_b)\to\Gamma(X,E_{zar}(\Omega^{\bullet}_X(\log D),F_b))\to \\
\Cone(\Gamma(X,\Omega^{\bullet}_X,F_b)\to\Gamma(U,E_{zar}(\Omega^{\bullet}_U,F_b)))
=:\Gamma(X,\Gamma_ZE_{zar}(\Omega_X^{\bullet},F_b)) 
\end{eqnarray*}
are quasi-isomorphisms (i.e. if we forgot filtrations), but the first one is NOT an $\infty$-filtered quasi-isomorphism
whereas the second one is an $\infty$-filtered quasi-isomorphism (recall that for $r>1$ the $r$-filtered quasi-isomorphisms
does NOT satisfy the 2 of 3 property for morphism of canonical triangles : see section 2.1).
\end{itemize}
\item[(ii)] More generally, let $f:X\to S$ a smooth projective morphism with $S,X\in\SmVar(\mathbb C)$.
Let $D=\cup D_i\subset X$ a normal crossing divisor such that $f_{|D_I}:=f\circ i_I:D_I\to S$ are SMOOTH morphisms
(note that it is a very special case), with $i_I:D_I\hookrightarrow X$ the closed embeddings.
Consider the open embedding $j:U:=X\backslash D\hookrightarrow X$ and $h:=f\circ j:U\to S$. 
\begin{itemize}
\item The map in $D_{\mathcal Dfil,\infty}(S)$ 
\begin{eqnarray*}
\Hom(L\mathbb D(\mathbb Z(U)),k)^{-1}\circ\Hom(\ad(j^*,j_*)(\mathbb Z(X/X)),E_{et}(\Omega^{\bullet}_{/S},F_b)): \\
\mathcal F_S^{GM}(\mathbb D(\mathbb Z(U/S))):=\Hom(L\mathbb D(\mathbb Z(U/S)),E_{et}(\Omega^{\bullet}_{/S},F_b)) \\
\xrightarrow{\sim}
\Hom(\Cone(\mathbb Z(D)\to\mathbb Z(X)),E_{zar}(\Omega^{\bullet}_{/S},F_b))=f_*E_{zar}(\Omega^{\bullet}_{X/S}(\nul D),F_b).
\end{eqnarray*}
is an isomorphism, where we recall $\mathbb D(\mathbb Z(U):=f_*j_*E_{et}(\mathbb Z(U/U))=h_*E_{et}(\mathbb Z(U/U))$, 
\item $\mathcal F_S^{GM}(\mathbb Z(U/S))=h_*E_{zar}\Omega_{U/S}^{\bullet},F_b)\in D_{\mathcal Dfil,\infty}(S)$
is NOT isomorphic to $f_*E_{zar}(\Omega^{\bullet}_{X/S}(\log D),F_b)$ in $D_{\mathcal Dfil,\infty}(S)$ in general.
In particular, the map,  
\begin{equation*}
j^*:=\ad(j^*,j_*)(-):H^nf_*E_{zar}(\Omega^{\bullet}_{X/S}(\log D))\xrightarrow{\sim}H^nh_*E_{zar}(\Omega^{\bullet}_{U/S}) 
\end{equation*}
which is an isomorphism in $D_{\mathcal D}(S)$ (i.e. if we forgot filtrations), gives embeddings  
\begin{eqnarray*}
j^*:=\ad(j^*,j_*)(-):F^pH^nh_*\mathbb C_U:=F^pH^nf_*E_{zar}(\Omega^{\bullet}_{X/S}(\log D),F_b)
\hookrightarrow F^pH^nh_*E_{zar}(\Omega^{\bullet}_{U/S},F_b)
\end{eqnarray*}
which are NOT an isomorphism in general for $n,p\in\mathbb Z$. Note that, since $a_U:U\to\left\{\pt\right\}$ is not proper,
\begin{equation*}
[\Delta_U]:\mathbb Z(U/S)\to h_*E_{et}(\mathbb Z(U/U))[2d_U]
\end{equation*}
is NOT an equivalence $(\mathbb A^1,et)$ local.
\item Let $Z\subset X$ a subvariety and denote $U:=X\backslash Z$ the open complementary. 
Denote $M_Z(X/S):=\Cone(M(U/S)\to M(X/S))\in\DA(S)$. 
If $f_{|Z}:=f\circ i_Z:Z\to S$ is a SMOOTH morphism, the map in $D_{\mathcal Dfil,\infty}(S)$
\begin{eqnarray*}
\Hom(G(X,Z),E_{et}(\Omega^{\bullet},F_b))\circ\Hom(\Gamma^{\vee}_Z\mathbb Z(X/X)),k)^{-1}: \\
\mathcal F_S^{GM}(M_Z(X/S)):=\Hom(f_{\sharp}\Gamma^{\vee}_Z\mathbb Z(X/X)),E_{et}(\Omega^{\bullet}_{/S},F_b)) 
\xrightarrow{\sim}f_*\Gamma_ZE_{zar}(\Omega_{X/S}^{\bullet},F_b)) \\
\xrightarrow{\sim}\mathcal F_S^{GM}(M(Z/S)(c)[2c])=f_{Z*}E_{zar}(\Omega_{Z/S}^{\bullet},F_b)(-c)[-2c]
\end{eqnarray*}
is an isomorphism,
where $c=\codim(Z,X)$ and $G(X,Z):f_{\sharp}\Gamma^{\vee}_Z\mathbb Z(X/X)\to\mathbb Z(Z/S)(c)[2c]$ is the Gynsin morphism.
\item Let $D\subset X$ a smooth divisor and denote $U:=X\backslash Z$ the open complementary
Note that the canonical distinguish triangle in $\DA(S)$
\begin{equation*}
M(U/S)\xrightarrow{\ad(j_{\sharp},j^*)(\mathbb Z(X/X))}M(X/S)\xrightarrow{\gamma^{\vee}_Z(\mathbb Z(X/X))} M_D(X/S)\to M(U/S)[1]
\end{equation*}
give a canonical triangle in $D_{\mathcal Dfil,\infty}(S)$
\begin{eqnarray*}
\mathcal F_S^{GM}(M_D(X/S))\xrightarrow{\mathcal F^{GM}(\gamma^{\vee}_Z(\mathbb Z(X/X)))}
\mathcal F_S^{GM}(M(X/S))\xrightarrow{\mathcal F^{GM}(\ad(j_{\sharp},j^*)(\mathbb Z(X/X))} 
\mathcal F_S^{GM}(M(U/S)) \\
\to\mathcal F_S^{GM}(M_D(X/S))[1],
\end{eqnarray*}
which is NOT the image of a distinguish triangle in $D(MHM(S))$.
\end{itemize}
\end{itemize}
\end{prop}

\begin{proof}
\noindent(i):For simplicity, we may assume that $i:D\hookrightarrow X$ is a smooth divisor.
Then, by theorem \ref{2functDM}, the map
\begin{equation*}
(0,\ad(j_{\sharp},j_*)(\mathbb Z(X/X)):\mathbb Z(D)\to\mathbb Z(X))\to\mathbb D(\mathbb Z(U/U))
\end{equation*}
is an equivalence $(\mathbb A^1,et)$ local in $C(\SmVar(\mathbb C))$.
The result then follows from proposition \ref{aetfibGM}. 
By theorem \ref{2functDM}, we have an equivalence $(\mathbb A^1,et)$ local in $C(\SmVar(\mathbb C))$
\begin{equation*}
G(X,Z):a_{X\sharp}\Gamma^{\vee}\mathbb Z(X/X)\to\mathbb Z(Z)(c)[2c]
\end{equation*}
The result then follows from proposition \ref{aetfibGM}.

\noindent(ii):For simplicity, we may assume that $i:D\hookrightarrow X$ is a smooth divisor.
Then, by theorem \ref{2functDM}, the map
\begin{equation*}
(0,\ad(j_{\sharp},j_*)(\mathbb Z(X/X)):\mathbb Z(D/S)\to\mathbb Z(X/S))\to\mathbb D(\mathbb Z(U/U))
\end{equation*}
is an equivalence $(\mathbb A^1,et)$ local in $C(\Var(\mathbb C)^{sm}/S)$.
The result then follows from proposition \ref{aetfibGM}.
By theorem \ref{2functDM}, we have an equivalence $(\mathbb A^1,et)$ local in $C(\SmVar(\mathbb C))$
\begin{equation*}
G(X,Z):f_{\sharp}\Gamma^{\vee}\mathbb Z(X/X)\to\mathbb Z(Z/S)(c)[2c]
\end{equation*}
The result then follows from proposition \ref{aetfibGM}.
\end{proof}

\begin{defi}
Let $S\in\SmVar(\mathbb C)$. We have, for $M,N\in\DA(S)$ and $F,G\in C(\Var(\mathbb C)^{sm}/S))$ projective
such that $M=D(\mathbb A^1,et)(F)$ and $N=D(\mathbb A^1,et)(G)$, the following transformation map in $D_{Ofil,\mathcal D}(S)$
\begin{eqnarray*}
T(\mathcal F_S^{GM},\otimes)(M,N):
\mathcal F_S^{GM}(M)\otimes^{L[-]}_{O_S}\mathcal F_S^{GM}(N):= \\
(e(S)_*\mathcal Hom(F,E_{et}(\Omega^{\bullet}_{/S},F_b)))\otimes_{O_S}
(e(S)_*\mathcal Hom(G,E_{et}(\Omega^{\bullet}_{/S},F_b)))[-d_S] \\
\xrightarrow{T(\otimes,\Omega_{/S})(F,G)}
e(S)_*\mathcal Hom(F\otimes G,E_{et}(\Omega^{\bullet}_{/S},F_b))[-d_S] \\ 
\xrightarrow{=}
e(S)_*\mathcal Hom(F\otimes G,E_{et}(\Omega^{\bullet}_{/S},F_b))[-d_S]
=:\mathcal F_S^{GM}(M\otimes N)
\end{eqnarray*}
\end{defi}

We now give the definition in the non smooth case :

\begin{defi}
Let $S\in\Var(\mathbb C)$ and $S=\cup_{i=1}^l S_i$ an open affine covering and denote, 
for $I\subset\left[1,\cdots l\right]$, $S_I=\cap_{i\in I} S_i$ and $j_I:S_I\hookrightarrow S$ the open embedding.
Let $i_i:S_i\hookrightarrow\tilde S_i$ closed embeddings, with $\tilde S_i\in\SmVar(\mathbb C)$. 
We have, for $M,N\in\DA(S)$ and $F,G\in C(\Var(\mathbb C)^{sm}/S)$ such that 
$M=D(\mathbb A^1,et)(F)$ and $N=D(\mathbb A^1,et)(G)$, 
the following transformation map in $D_{Ofil,\mathcal D}(S/(\tilde S_I))$
\begin{eqnarray*}
T(\mathcal F_S^{GM},\otimes)(M,N):
\mathcal F_S^{GM}(M)\otimes^{L[-]}_{O_S}\mathcal F_S^{GM}(N):= \\
(e(\tilde S_I)_*\mathcal Hom(L(i_{I*}j_I^*F),E_{et}(\Omega^{\bullet}_{/\tilde S_I},F_b))[-d_{\tilde S_I}],u_{IJ}(F))
\otimes^{[-]}_{O_S} \\
(e(\tilde S_I)_*\mathcal Hom(L(i_{I*}j_I^*G),E_{et}(\Omega^{\bullet}_{/\tilde S_I},F_b))[-d_{\tilde S_I}],u_{IJ}(G)) \\
\xrightarrow{=}
((e(\tilde S_I)_*\mathcal Hom(L(i_{I*}j_I^*F),E_{et}(\Omega^{\bullet}_{/\tilde S_I},F_b))\otimes_{O_{\tilde S_I}} \\
e(\tilde S_I)_*\mathcal Hom(L(i_{I*}j_I^*G),E_{et}(\Omega^{\bullet}_{/\tilde S_I},F_b)))[-d_{\tilde S_I}],
u_{IJ}(F)\otimes u_{IJ}(G)) \\
\xrightarrow{(T(\otimes,\Omega_{/\tilde S_I})(L(i_{I*}j_I^*F),L(i_{I*}j_I^*G)))} \\
(e(\tilde S_I)_*\mathcal Hom(L(i_{I*}j_I^*F)\otimes L(i_{I*}j_I^*G),
E_{et}(\Omega^{\bullet}_{/\tilde S_I},F_{DR}))[-d_{\tilde S_I}],v_{IJ}(F\otimes G)) \\
\xrightarrow{=} 
(e(\tilde S_I)_*\mathcal Hom(L(i_{I*}j_I^*(F\otimes G),E_{et}(\Omega^{\bullet}_{/\tilde S_I},F_b)))[-d_{\tilde S_I}],
u_{IJ}(F\otimes G))=:\mathcal F_S^{GM}(M\otimes N)
\end{eqnarray*}
\end{defi}

\begin{prop}
Let $f_1:X_1\to S$, $f_2:X_2\to S$ two morphism with $X_1,X_2,S\in\Var(\mathbb C)$. 
Assume that there exist factorizations 
$f_1:X_1\xrightarrow{l_1}Y_1\times S\xrightarrow{p_S} S$, $f_2:X_2\xrightarrow{l_2}Y_2\times S\xrightarrow{p_S} S$
with $Y_1,Y_2\in\SmVar(\mathbb C)$, $l_1,l_2$ closed embeddings and $p_S$ the projections.
We have then the factorization
\begin{equation*}
f_1\times f_2:X_{12}:=X_1\times_S X_2\xrightarrow{l_1\times l_2}Y_1\times Y_2\times S\xrightarrow{p_S} S
\end{equation*}
Let $S=\cup_{i=1}^l S_i$ an open affine covering and denote, 
for $I\subset\left[1,\cdots l\right]$, $S_I=\cap_{i\in I} S_i$ and $j_I:S_I\hookrightarrow S$ the open embedding.
Let $i_i:S_i\hookrightarrow\tilde S_i$ closed embeddings, with $\tilde S_i\in\SmVar(\mathbb C)$. 
We have, for $M,N\in\DA(S)$ and $F,G\in C(\Var(\mathbb C)^{sm}/S)$ such that 
$M=D(\mathbb A^1,et)(F)$ and $N=D(\mathbb A^1,et)(G)$, 
the following commutative diagram in $D_{Ofil,\mathcal D}(S/(\tilde S_I))$
\begin{equation*}
\begin{tikzcd}
\mathcal F_S^{GM}(M(X_1/S))\otimes^L_{O_S}\mathcal F_S^{GM}(M(X_2/S))
\ar[d,"I^{GM}(X_1/S)\otimes I^{GM}(X_2/S)"]
\ar[rr,"T(\mathcal F_S^{GM}{,}\otimes)(M(X_1/S){,}M(X_2/S))"]  & \, &
\shortstack{$\mathcal F_S^{GM}(M(X_1/S)\otimes M(X_2/S))$ \\ $=\mathcal F_S^{GM}(M(X_1\times_S X_2/S))$}
\ar[d,"I^{GM}(X_{12}/S)"] \\
\shortstack{$(p_{\tilde S_I*}\Gamma_{X_{1I}}
E_{zar}(\Omega^{\bullet}_{Y_1\times\tilde S_I/\tilde S_I},F_b)[-d_{\tilde S_I}],w_{IJ}(X_1/S))\otimes_{O_S}$ \\
$(p_{\tilde S_I*}\Gamma_{X_{2I}}
E_{zar}(\Omega^{\bullet}_{Y_2\times\tilde S_I/\tilde S_I},F_b)[-d_{\tilde S_I}],w_{IJ}(X_2/S))$}
\ar[rr,"(Ew_{(Y_1\times\tilde S_I{,}Y_2\times\tilde S_I)/\tilde S_I})"] & \, &
\shortstack{$(p_{\tilde S_I*}\Gamma_{X_{12I}}
E_{zar}(\Omega^{\bullet}_{Y_1\times Y_2\times\tilde S_I/\tilde S_I},F_b)[-d_{\tilde S_I}],$ \\ $w_{IJ}(X_{12}/S))$}.
\end{tikzcd}
\end{equation*}
\end{prop}

\begin{proof}
Immediate from definition.
\end{proof}

\subsubsection{The algebraic filtered De Rham realization functor and the commutativity with the six operation}

We recall (see section 2), for $f:T\to S$ a morphism with $T,S\in\Var(\mathbb C)$, 
the commutative diagrams of sites (\ref{muf}) and (\ref{Grf})
\begin{equation*}
\xymatrix{
\Var(\mathbb C)^2/T\ar[rd]^{\rho_T}\ar[rr]^{\mu_T}\ar[dd]_{P(f)} & \, & 
\Var(\mathbb C)^{2,pr}/T\ar[dd]^{P(f)}\ar[rd]^{\rho_T} & \, \\
\, & \Var(\mathbb C)^{2,sm}/T\ar[rr]^{\mu_T}\ar[dd]_{P(f)} & \, & \Var(\mathbb C)^{2,smpr}/T\ar[dd]^{P(f)} \\
\Var(\mathbb C)^2/S\ar[rd]^{\rho_S}\ar[rr]^{\mu_S} & \, & \Var(\mathbb C)^{2,pr}/S\ar[rd]^{\rho_S} \\
\, & \Var(\mathbb C)^{2,sm}/S\ar[rr]^{\mu_S} & \, & \Var(\mathbb C)^{2,smpr}/S}
\end{equation*}
and
\begin{equation*}
\xymatrix{\Var(\mathbb C)^{2,pr}/T\ar[rr]^{\Gr_T^{12}}\ar[dd]_{P(f)}\ar[rd]^{\rho_T} & \, & 
\Var(\mathbb C)/T\ar[dd]^{P(f)}\ar[rd]^{\rho_T} & \, \\
\, & \Var(\mathbb C)^{2,smpr}/T\ar[rr]^{\Gr_T^{12}}\ar[dd]_{P(f)} & \, & \Var(\mathbb C)^{sm}/T\ar[dd]^{P(f)} \\
\Var(\mathbb C)^{2,pr}/S\ar[rr]^{\Gr_S^{12}}\ar[rd]^{\rho_S} & \, & \Var(\mathbb C)/S\ar[rd]^{\rho_S} & \, \\
\, & \Var(\mathbb C)^{2,sm}/S\ar[rr]^{\Gr_S^{12}} & \, & \Var(\mathbb C)^{sm}/S}.
\end{equation*}
Let $S\in\Var(\mathbb C)$. We have for $F\in C(\Var(\mathbb C)^{sm}/S)$ the canonical map in $C(\Var(\mathbb C)^{sm}/S)$
\begin{eqnarray*}
\Gr(F):\Gr^{12}_{S*}\mu_{S*}F^{\Gamma}\to F, \\
\Gr(F)(U/S):\Gamma_U^{\vee}p^*F(U\times S/U\times S)\xrightarrow{\ad(l^*,l_*)(p^*F)(U\times S/U\times S)}h^*F(U/U)=F(U/S)
\end{eqnarray*}
where $h:U\to S$ is a smooth morphism with $U\in\Var(\mathbb C)$ and $h:U\xrightarrow{l}U\times S\xrightarrow{p}S$
is the graph factorization with $l$ the graph embedding and $p$ the projection.

For $s:\mathcal I\to\mathcal J$ a functor, with $\mathcal I,\mathcal J\in\Cat$, and
$f_{\bullet}:T_{\bullet}\to S_{s(\bullet)}$ a morphism with 
$T_{\bullet}\in\Fun(\mathcal J,\Var(\mathbb C))$ and $S_{\bullet}\in\Fun(\mathcal I,\Var(\mathbb C))$, 
we have then the commutative diagrams of sites (\ref{mufIJ}) and (\ref{GrfIJ})
\begin{equation*}
\xymatrix{\Var(\mathbb C)^2/T_{\bullet}\ar[rr]^{\mu_{T_{\bullet}}}\ar[dd]_{P(f_{\bullet})}\ar[rd]^{\rho_{T_{\bullet}}} & \, & 
\Var(\mathbb C)^{2,pr}/T_{\bullet}\ar[dd]^{P(f_{\bullet})}\ar[rd]^{\rho_{T_{\bullet}}} & \, \\
\, & \Var(\mathbb C)^{2,sm}/T_{\bullet}\ar[rr]^{\mu_{T_{\bullet}}}\ar[dd]_{P(f_{\bullet})} & \, & 
\Var(\mathbb C)^{2,smpr}/T_{\bullet}\ar[dd]^{P(f_{\bullet})} \\
\Var(\mathbb C)^2/S_{\bullet}\ar[rr]^{\mu_{S_{\bullet}}}\ar[rd]^{\rho_{S_{\bullet}}} & \, & 
\Var(\mathbb C)^{2,pr}/S_{\bullet}\ar[rd]^{\rho_{S_{\bullet}}} & \, \\
\, & \Var(\mathbb C)^{2,sm}/S_{\bullet}\ar[rr]^{\mu_{S_{\bullet}}} & \, & \Var(\mathbb C)^{2,smpr}/S_{\bullet}}.
\end{equation*}
and
\begin{equation*}
\xymatrix{\Var(\mathbb C)^{2,pr}/T_{\bullet}
\ar[rr]^{\Gr_{T_{\bullet}}^{12}}\ar[dd]_{P(f_{\bullet})}\ar[rd]^{\rho_{T_{\bullet}}} & \, & 
\Var(\mathbb C)/T\ar[dd]^{P(f_{\bullet})}\ar[rd]^{\rho_{T_{\bullet}}} & \, \\
\, & \Var(\mathbb C)^{2,smpr}/T_{\bullet}\ar[rr]^{\Gr_{T}^{12}}\ar[dd]_{P(f_{\bullet})} & \, & 
\Var(\mathbb C)^{sm}/T_{\bullet}\ar[dd]^{P(f_{\bullet})} \\
\Var(\mathbb C)^{2,pr}/S_{\bullet}\ar[rr]^{\Gr_{S_{\bullet}}^{12}}\ar[rd]^{\rho_{S_{\bullet}}} & \, & 
\Var(\mathbb C)/S_{\bullet}\ar[rd]^{\rho_{S_{\bullet}}} & \, \\
\, & \Var(\mathbb C)^{2,sm}/S_{\bullet}\ar[rr]^{\Gr_{S_{\bullet}}^{12}} & \, & 
\Var(\mathbb C)^{sm}/S_{\bullet}}.
\end{equation*}

We will use the following map from the property of mixed Hodge module (see section 5)
together with the specialization map of a filtered D module for a closed embedding (see definition \ref{Vfil}) :

\begin{defiprop}\label{VfilOlem}
\begin{itemize}
\item[(i)]Let $l:Z\hookrightarrow S$ a closed embedding with $S,Z\in\SmVar(\mathbb C)$.
Consider an open embedding $j:S^o\hookrightarrow S$. We then have the cartesian square
\begin{equation*}
\xymatrix{S^o\ar[r]^j & S \\ 
Z^o:=Z\times_SS^o\ar[r]^{j'}\ar[u]^{l'} & Z\ar[u]^l}
\end{equation*}
where $j'$ is the open embedding given by base change. Using proposition \ref{VfilVarAn}, 
the morphisms $Q^{p,0}_{V_Z,V_D}(O_S,F_b)$ for $D\subset S$ a closed subset
of definition-proposition \ref{Vfildefiprop} induces a canonical morphism in $C_{l^*O_Sfil}(Z)$
\begin{equation*}
Q(Z,j_!)(O_S,F_b):l^*Q_{V_Z,0}j_{!Hdg}(O_{S^o},F_b)\to j'_{!Hdg}(O_{Z^o},F_b),
\end{equation*}
where $V_Z$ is the Kashiwara-Malgrange $V_Z$-filtration and $V_D$ is the Kashiwara-Malgrange $V_D$-filtration,
which commutes with the action of $T_Z$.
\item[(ii)]Let $l:Z\hookrightarrow S$ and $k:Z'\hookrightarrow Z$ be closed embeddings with $S,Z,Z'\in\SmVar(\mathbb C)$.
Consider an open embedding $j:S^o\hookrightarrow S$. We then have the commutative diagram whose squares are cartesian.
\begin{equation*}
\xymatrix{S^o\ar[r]^j & S \\ 
Z^o:=Z\times_SS^o\ar[r]^{j'}\ar[u]^{l'} & Z\ar[u]^l \\
Z^{'o}:=Z'\times_SS^o\ar[r]^{j''}\ar[u]^{k'} & Z'\ar[u]^{k}}
\end{equation*}
where $j'$ is the open embedding given by base change. Then, 
\begin{eqnarray*}
Q(Z',j_!)(O_S,F_b)=Q(Z',j'_!)(O_Z,F_b)\circ (k^*Q_{V_{Z'},0}Q(Z,j_!)(O_S,F_b)): \\
k^*Q_{V_{Z'},0}l^*Q_{V_Z,0}j_{!Hdg}(O_{S^o},F_b)
\xrightarrow{k^*Q_{V_{Z'},0}Q(Z,j_!)(O_S,F_b)}k^*Q_{V_{Z'},0}j'_{!Hdg}(O_{Z^o},F_b) \\
\xrightarrow{Q(Z',j'_!)(O_Z,F_b)}j^{''Hdg}_!(O_{Z^{'o}},F_b)
\end{eqnarray*}
in $C_{k^*l^*O_Sfil}(Z')$ which commutes with the action of $T_{Z'}$.
\item[(iii)] Consider a commutative diagram whose squares are cartesian
\begin{equation*}
\xymatrix{S^{oo}\ar[r]^{j_2} & S^o\ar[r]^{j_1} & S \\ 
Z^{oo}:=Z\times_SS^{oo}\ar[r]^{j'_2}\ar[u]^{l''} & Z^o:=Z\times_SS^o\ar[r]^{j'_1}\ar[u]^{l'} & Z\ar[u]^l}
\end{equation*}
where $j_1$, $j_2$, and hence $j'_1$,$j'_2$ are open embeddings.
We have then the following commutative diagram
\begin{equation*}
\xymatrix{l^*Q_{V_Z,0}j^{Hdg}_{1!}(O_{S^o},F_b)
\ar[rr]^{\ad(j_{2!Hdg},j_2^*)(O_{S^o},F_b)}\ar[d]^{Q(Z,j_!)(O_S,F_b)} & \, &
l^*Q_{V_Z,0}(j_1\circ j_2)_{!Hdg}(O_{S^{oo}},F_b)\ar[d]^{Q(Z,(j_1\circ j_2)_!)(O_S,F_b)} \\
j^{'Hdg}_{1!}(O_{Z^o},F_b)\ar[rr]^{\ad(j'_{2!Hdg},j_2^{'*})(O_{Z^o},F_b)} & \, &
(j'_1\circ j'_2)_{!Hdg}(O_{Z^{oo}},F_b)}
\end{equation*}
in $C_{l^*O_Sfil}(Z)$ which commutes with the action of $T_{Z}$.
\end{itemize}
\end{defiprop}

\begin{proof}
\noindent(i): By definition of $j_{!Hdg}:\pi_{S^o}(MHM(S^o))\to \pi_S(C(MHM(S)))$,
we have to construct the isomorphism for each complement of a (Cartier) divisor
$j=j_D:S^o=S\backslash D\hookrightarrow S$.
In this case, we have the closed embedding $i:S\hookrightarrow L$ given by the zero section of the line bundle
$L=L_D$ associated to $D$. 
We have then, using definition-proposition \ref{Vfildefiprop},
the canonical morphism in $PSh_{l^*O_Sfil}(Z)$ which commutes with the action of $T_Z$
\begin{eqnarray*}
Q(Z,j_!)(O_S,F_b):l^*Q_{V_Z,0}j_{!Hdg}(O_{S^o},F_b)\xrightarrow{T_!(l,j)(-)^{-1}}
j_{!Hdg}Q_{V_{Z^o},0}(O_{S^o},F_b)=j'_{!Hdg}(O_{Z^o},F_b).
\end{eqnarray*}
and $V_Z^pT_!(l,j)(-)^{-1}=Q^{p,0}_{V_Z,V_S}(i_{*mod}(O_S,F_b))$.
Now for $j:S^o=S\backslash R\hookrightarrow S$ an arbitrary open embedding, we set
\begin{eqnarray*}
Q(Z,j_!)(O_S,F_b):=\varprojlim_{(D_i),R\subset D_i\subset S} (Q(Z,j_{D_J!})(j_{D_I}^*(O_S,F_b)): 
l^*Q_{V_Z,0}j_{!Hdg}(O_{S^o},F_b)\xrightarrow{\sim}j'_{!Hdg}(O_{Z^o},F_b)
\end{eqnarray*}

\noindent(ii): Follows from definition-proposition \ref{Vfildefiprop}.

\noindent(iii): Follows from definition-proposition \ref{Vfildefiprop}.
\end{proof}

Using definition-proposition \ref{TgammaHdg} in the projection case, and
the specialization map given in definition \ref{Vfil} and the isomorphism of definition-proposition \ref{VfilOlem},
in the closed embedding case, we have the following canonical map :

\begin{defi}\label{TgHdgO}
Consider a commutative diagram in $\SmVar(\mathbb C)$ whose square are cartesian
\begin{equation*}
\xymatrix{Z_T\ar[r]^{i'}\ar[dd]^{g'}\ar[rd]^{l'} & T\ar[dd]^g\ar[rd]^l  & 
T\backslash Z_T\ar[l]^{j'}\ar[dd]^{g}\ar[rd]^{l} \\ 
 \, & T\times Z\ar[r]^{I\times i}\ar[ld]^{p_Z} & T\times S\ar[ld]^{p_S} & 
T\times S\backslash (T\times Z)\ar[l]^{I\times j}\ar[ld]^{p_S}\\
Z\ar[r]^i & S & S\backslash Z\ar[l]^j}
\end{equation*}
where $i$ and hence $I\times i$ and $i'$, are closed embeddings,
$j$, $I\times j$, $j'$ are the complementary open embeddings
and $g:T\xrightarrow{l}T\times S\xrightarrow{p_S}S$ 
is the graph factorization, where $l$ is the graph embedding and $p_S$ the projection. 
Then, the map in $C_{l^*O_{T\times S}fil}(T)$
\begin{eqnarray*}
sp_{V_T}(\Gamma_{T\times Z}^{\vee,Hdg}(O_{T\times S},F_b)):
l^*\Gamma_{T\times Z}^{\vee,Hdg}(O_{T\times S},F_b)\xrightarrow{q_{V_T,0}}
l^*Q_{V_T,0}(\Gamma_{T\times Z}^{\vee,Hdg}(O_{T\times S},F_b)) \\
\xrightarrow{Q(T,(I\times j)_!)(O_{T\times S},F_b):=T_!(l,(I\times j))(-)}
\Gamma_{Z_T}^{\vee,Hdg}(O_T,F_b)
\end{eqnarray*}
which commutes with the action of $T_T$, where the first map is given in definition \ref{Vfil} 
and the last map is studied definition-proposition \ref{VfilOlem}, factors through
\begin{eqnarray*}
sp_{V_T}(\Gamma_{T\times Z}^{\vee,Hdg}(O_{T\times S},F_b)):
l^*\Gamma_{T\times Z}^{\vee,Hdg}(O_{T\times S},F_b)\xrightarrow{n}
l^{*mod}\Gamma_{T\times Z}^{\vee,Hdg}(O_{T\times S},F_b) \\
\xrightarrow{\bar{sp}_{V_T}(\Gamma_{T\times Z}^{\vee,Hdg}(O_{T\times S},F_b))}
\Gamma_{Z_T}^{\vee,Hdg}(O_T,F_b),
\end{eqnarray*}
with for $U\subset T\times S$ an open subset, $m\in\Gamma(U,O_{T\times S})$ and $h\in\Gamma(U_T,O_T)$,
$n(m):=n\otimes 1$ and $\bar{sp}_{V_T}(-)(m\otimes h)=h\cdot sp_{V_T}(m)$ ;
see definition-proposition \ref{TgammaHdg}, proposition \ref{VfilVarAn} and theorem \ref{VMHSthm}. Then,
\begin{eqnarray*}
\bar{sp}_{V_T}(\Gamma_{T\times Z}^{\vee,Hdg}(O_{T\times S},F_b)):
l^{*mod}\Gamma_{T\times Z}^{\vee,Hdg}(O_{T\times S},F_b)\to\Gamma_{Z_T}^{\vee,Hdg}(O_T,F_b),
\end{eqnarray*}
is a map in $C_{\mathcal D(1,0)fil}(T)$, i.e. is $D_T$ linear.
We then consider the canonical map in $C_{\mathcal D(1,0)fil}(T)$ 
\begin{eqnarray*}
a(g,Z)(O_S,F_b):g^{*mod}\Gamma_Z^{\vee,Hdg}(O_S,F_b)=l^{*mod}p_S^{*mod}\Gamma_Z^{\vee,Hdg}(O_S,F_b) 
\xrightarrow{l^{*mod}T^{Hdg}(p,\gamma^{\vee})(O_S,F_b)^{-1}} \\
l^{*mod}\Gamma_{T\times Z}^{\vee,Hdg}(O_{T\times S},F_b) 
\xrightarrow{\bar{sp}_{V_T}(\Gamma_{T\times Z}^{\vee,Hdg}(O_{T\times S},F_b))}
\Gamma_{Z_T}^{\vee,Hdg}(O_T,F_b).
\end{eqnarray*}
\end{defi}

\begin{lem}\label{TgHdgOlem}
\begin{itemize}
\item[(i)]For $g:T\to S$ and $g:T'\to T$ two morphism with $S,T,T'\in\SmVar(\mathbb C)$,
considering the commutative diagram whose squares are cartesian
\begin{equation*}
\xymatrix{Z_{T'}\ar[r]^{i''}\ar[d]^{g'} & T'\ar[d]^{g'}  & T'\backslash Z_{T'}\ar[l]^{j''}\ar[d]^{g'} \\ 
Z_T\ar[r]^{i'}\ar[d]^g & T\ar[d]^{g} & T\backslash Z_T\ar[l]^{j'}\ar[d]^{g} \\
Z\ar[r]^i & S & S\backslash Z\ar[l]^j}
\end{equation*}
we have then
\begin{eqnarray*}
a(g\circ g',Z)(O_S,F_b)=a(g',Z_T)(O_T,F_b)\circ (g^{'*mod}a(g,Z)(O_S,F_b)): \\
(g\circ g')^{*mod}\Gamma_Z^{\vee,Hdg}(O_S,F_b)=g^{'*mod}g^{*mod}\Gamma_Z^{\vee,Hdg}(O_S,F_b)
\xrightarrow{g^{'*mod}a(g,Z)(O_S,F_b)}g^{'*mod}\Gamma_{Z_T}^{\vee,Hdg}(O_T,F_b) \\
\xrightarrow{a(g',Z_T)(O_T,F_b)}\Gamma_{Z_{T'}}^{\vee,Hdg}(O_{T'},F_b).
\end{eqnarray*}
\item[(ii)]For $g:T\to S$ a morphism with $S,T\in\SmVar(\mathbb C)$,
considering the commutative diagram whose squares are cartesian
\begin{equation*}
\xymatrix{Z'_T\ar[r]^{k'}\ar[d]^g & Z_T\ar[r]^{i'}\ar[d]^{g} & T\ar[d]^{g} \\
Z'\ar[r]^{k} & Z\ar[r]^i & S }
\end{equation*}
we have then the following commutative diagram
\begin{equation*}
\xymatrix{g^{*mod}\Gamma_Z^{\vee,Hdg}(O_S,F_b)
\ar[d]^{a(g,Z)(O_S,F_b)}\ar[rr]^{g^{*mod}T(Z'/Z,\gamma^{\vee,Hdg})(O_S,F_b)} & \, &
g^{*mod}\Gamma_{Z'}^{\vee,Hdg}(O_S,F_b)\ar[d]^{a(g,Z')(O_S,F_b)} \\
\Gamma_{Z_T}^{\vee,Hdg}(O_T,F_b)\ar[rr]^{T(Z'_T/Z_T,\gamma^{\vee,Hdg})(O_T,F_b)} & \, &
\Gamma_{Z'_T}^{\vee,Hdg}(O_T,F_b)}
\end{equation*}
\end{itemize}
\end{lem}

\begin{proof}
\noindent(i):Follows from definition-proposition \ref{VfilOlem} (ii)

\noindent(ii):Follows from definition-proposition \ref{VfilOlem} (iii)
\end{proof}

We can now define the main object :

\begin{defi}\label{wtildew}
\begin{itemize}
\item[(i)]For $S\in\SmVar(\mathbb C)$, we consider the filtered complexes of presheaves 
\begin{eqnarray*} 
(\Omega^{\bullet,\Gamma,pr}_{/S},F_{DR})\in C_{D_Sfil}(\Var(\mathbb C)^{2,smpr}/S) 
\end{eqnarray*}
given by, 
\begin{itemize}
\item for $(Y\times S,Z)/S=((Y\times S,Z),p)\in\Var(\mathbb C)^{2,smpr}/S$, 
\begin{eqnarray*}
(\Omega^{\bullet,\Gamma,pr}_{/S}((Y\times S,Z)/S),F_{DR}):=
((\Omega^{\bullet}_{Y\times S/S},F_b)\otimes_{O_{Y\times S}}\Gamma^{\vee,Hdg}_Z(O_{Y\times S},F_b))(Y\times S) 
\end{eqnarray*}
with the structure of $p^*D_S$ module given by proposition \ref{DRhUS},
\item for $g:(Y_1\times S,Z_1)/S=((Y_1\times S,Z_1),p_1)\to (Y\times S,Z)/S=((Y\times S,Z),p)$ 
a morphism in $\Var(\mathbb C)^{2,smpr}/S$, denoting for short $\hat{Z}:=Z\times_{Y\times S}(Y_1\times S)$,
\begin{eqnarray*}
\Omega^{\bullet,\Gamma,pr}_{/S}(g):
((\Omega^{\bullet}_{Y\times S/S},F_b)\otimes_{O_{Y\times S}}\Gamma^{\vee,Hdg}_Z(O_{Y\times S},F_b))(Y\times S) \\ 
\xrightarrow{i_{-}} 
g^*((\Omega^{\bullet}_{Y\times S/S},F_b)\otimes_{O_{Y\times S}}\Gamma^{\vee,Hdg}_Z(O_{Y\times S},F_b))(Y_1\times S) \\  
\xrightarrow{\Omega_{(Y_1\times S/Y\times S)/(S/S)}(\Gamma^{\vee,Hdg}_Z(O_{Y\times S},F_b))(Y_1\times S)} 
(\Omega^{\bullet}_{Y_1\times S/S},F_b)\otimes_{O_{Y_1\times S}}g^{*mod}\Gamma^{\vee,Hdg}_Z(O_{Y\times S},F_b))(Y_1\times S) \\ 
\xrightarrow{DR(Y_1\times S/S)(a(g,Z)(O_{Y\times S},F_b))(Y_1\times S)} 
(\Omega^{\bullet}_{Y_1\times S/S},F_b)\otimes_{O_{Y_1\times S}}\Gamma^{\vee,Hdg}_{\hat{Z}}(O_{Y_1\times S},F_b))(Y_1\times S) \\
\xrightarrow{DR(Y_1\times S/S)(T(Z_1/\hat{Z},\gamma^{\vee,Hdg})(O_{Y_1\times S},F_b))(Y_1\times S)} 
(\Omega^{\bullet}_{Y_1\times S/S},F_b)\otimes_{O_{Y_1\times S}}\Gamma^{\vee,Hdg}_{Z_1}(O_{Y_1\times S},F_b))(Y_1\times S), 
\end{eqnarray*}
where 
\begin{itemize}
\item $i_{-}$ is the arrow of the inductive limit,
\item we recall that
\begin{eqnarray*}
\Omega_{(Y_1\times S/Y\times S)/(S/S)}(\Gamma^{\vee,Hdg}_Z(O_{Y\times S},F_b)):
g^*((\Omega^{\bullet}_{Y\times S/S},F_b)\otimes_{O_{Y\times S}}\Gamma^{\vee,Hdg}_Z(O_{Y\times S},F_b)) \\  
\to(\Omega^{\bullet}_{Y_1\times S/S},F_b)\otimes_{O_{Y_1\times S}}g^{*mod}\Gamma^{\vee,Hdg}_Z(O_{Y\times S},F_b))
\end{eqnarray*}
is the map given in definition-proposition \ref{TDwM}, which is $p_1^*D_S$ linear by proposition \ref{TDhwM},
\item the map
\begin{equation*}
a(g,Z)(O_{Y\times S},F_b):g^{*mod}\Gamma^{\vee,Hdg}_{Z}(O_{Y\times S},F_b)\to\Gamma^{\vee,Hdg}_{\hat{Z}}(O_{Y_1\times S},F_b) 
\end{equation*}
is the map given in definition \ref{TgHdgO} 
\item the map
\begin{eqnarray*}
T(Z_1/\hat{Z},\gamma^{\vee,Hdg})(O_{Y_1\times S},F_b):
\Gamma^{\vee,Hdg}_{\hat{Z}}(O_{Y_1\times S},F_b)\to\Gamma^{\vee,Hdg}_{Z_1}(O_{Y_1\times S},F_b) 
\end{eqnarray*}
is given in definition-proposition \ref{TgammaHdg}.
\end{itemize}
\end{itemize}
For $g:((Y_1\times S,Z_1),p_1)\to ((Y\times S,Z),p)$ and $g':((Y'_1\times S,Z'_1),p_1)\to ((Y_1\times S,Z_1),p)$
two morphisms in $\Var(\mathbb C)^{2,smpr}/S$, we have
\begin{eqnarray*}
\Omega^{\bullet,\Gamma,pr}_{/S}(g\circ g')=
\Omega^{\bullet,\Gamma,pr}_{/S}(g')\circ\Omega^{\bullet,\Gamma,pr}_{/S}(g):
((\Omega^{\bullet}_{Y\times S/S},F_b)\otimes_{O_{Y\times S}}\Gamma^{\vee,Hdg}_Z(O_{Y\times S},F_b))(Y\times S) \\ 
\xrightarrow{\Omega^{\bullet,\Gamma,pr}_{/S}(g)} 
(\Omega^{\bullet}_{Y_1\times S/S},F_b)\otimes_{O_{Y_1\times S}}\Gamma^{\vee,Hdg}_{Z_1}(O_{Y_1\times S},F_b))(Y_1\times S) \\ 
\xrightarrow{\Omega^{\bullet,\Gamma,pr}_{/S}(g')} 
(\Omega^{\bullet}_{Y'_1\times S/S},F_b)\otimes_{O_{Y'_1\times S}}\Gamma^{\vee,Hdg}_{Z'_1}(O_{Y'_1\times S},F_b))(Y'_1\times S),
\end{eqnarray*}
since, denoting for short $\hat{Z}:=Z\times_{Y\times S}(Y_1\times S)$ and $\hat{Z}':=Z\times_{Y\times S}(Y'_1\times S)$
\begin{itemize}
\item we have by lemma \ref{TgHdgOlem}(i)
\begin{equation*}
a(g\circ g',\hat Z')(O_{Y\times S},F_b)=a(g',\hat Z)(O_{Y_1\times S},F_b)\circ g^{'*mod}a(g,Z)(O_{Y\times S},F_b),
\end{equation*}
\item we have by lemma \ref{TgHdgOlem}(ii)
\begin{eqnarray*}
T(Z'_1/\hat{Z}',\gamma^{\vee,Hdg})(O_{Y'_1\times S},F_b)\circ a(g',\hat Z)(O_{Y_1\times S},F_b) \\
=a(g',Z_1)(O_{Y_1\times S},F_b)\circ g^{'*mod}T(Z_1/\hat{Z},\gamma^{\vee,Hdg})(O_{Y_1\times S},F_b).
\end{eqnarray*}
\end{itemize}
\item[(ii)]For $S\in\SmVar(\mathbb C)$, we have the canonical map $C_{O_Sfil,D_S}(\Var(\mathbb C)^{sm}/S)$ 
\begin{eqnarray*} 
\Gr(\Omega_{/S}):\Gr_{S*}^{12}(\Omega^{\bullet,\Gamma,pr}_{/S},F_{DR})\to(\Omega^{\bullet}_{/S},F_b) 
\end{eqnarray*}
given by, for $U/S=(U,h)\in\Var(\mathbb C)^{sm}/S$
\begin{eqnarray*}
\Gr(\Omega_{/S})(U/S):\Gr_{S*}^{12}(\Omega^{\bullet,\Gamma,pr}_{/S},F_{DR})(U/S):=
((\Omega^{\bullet}_{U\times S/S},F_b)\otimes_{O_{U\times S}}\Gamma^{\vee,Hdg}_{U}(O_{U\times S},F_b))(U\times S) \\
\xrightarrow{\ad(i_U^*,i_{U*})(-)(U\times S)}
i^*((\Omega^{\bullet}_{U\times S/S},F_b)\otimes_{O_{U\times S}}\Gamma^{\vee,Hdg}_{U}(O_{U\times S},F_b))(U) \\
\xrightarrow{\Omega_{(U/U\times S)/(S/S)}(-)(U)}
((\Omega^{\bullet}_{U/S},F_b)\otimes_{O_{U}}i_U^{*mod}\Gamma^{\vee,Hdg}_{U}(O_{U\times S},F_b))(U) \\
\xrightarrow{DR(U/S)(a(i_U,U))(U)}
(\Omega^{\bullet}_{U/S},F_b)(U)=:(\Omega^{\bullet}_{/S},F_b)(U/S)
\end{eqnarray*}
where $h:U\xrightarrow{i_U}U\times S\xrightarrow{p_S}S$ is the graph factorization 
with $i_U$ the graph embedding and $p_S$ the projection, note that $a(i_U,U)$ is an isomorphism since for
$j_U:U\times S\backslash U\hookrightarrow U\times S$ the open complementary $i_U^{*mod}j^{Hdg}_{U!}(M,F,W)=0$.
\end{itemize}
\end{defi}

\begin{defi}\label{wtildewT}
For $S\in\SmVar(\mathbb C)$, we have the canonical map $C_{O_Sfil,D_S}(\Var(\mathbb C)^{2,smpr}/S)$ 
\begin{eqnarray*} 
T(\Omega^{\Gamma}_{/S}):\mu_{S*}(\Omega^{\bullet,\Gamma}_{/S},F_b)\to(\Omega^{\bullet,\Gamma,pr}_{/S},F_{DR}) 
\end{eqnarray*}
given by, for $(Y\times S,X)/S=((Y\times S,Z),p)\in\Var(\mathbb C)^{2,smpr}/S$
\begin{eqnarray*} 
T(\Omega^{\Gamma}_{/S})((Y\times S,Z)/S): \\
(\Omega^{\bullet,\Gamma}_{/S},F_b)((Y\times S,Z)/S):=
\mathbb D_{p^*O_S}L_{p^*O}\Gamma_{Z}E_{zar}(\mathbb D_{p^*O_S}L_{p^*O}(\Omega^{\bullet}_{Y\times S/S},F_b))(Y\times S) \\
\xrightarrow{DR(Y\times S/S)(\gamma^{\vee,Hdg}_Z(O_{Y\times S}))(Y\times S)} \\
\mathbb D_{p^*O_S}L_{p^*O}\Gamma_{Z}E_{zar}(\mathbb D_{p^*O_S}L_{p^*O}
((\Omega^{\bullet}_{Y\times S/S},F_b)\otimes_{O_{Y\times S}}\Gamma^{\vee,Hdg}_{Z}(O_{Y\times S},F_b))(Y\times S)
\xrightarrow{=} \\
((\Omega^{\bullet}_{Y\times S/S},F_b)\otimes_{O_{Y\times S}}\Gamma^{\vee,Hdg}_{Z}(O_{Y\times S},F_b))(Y\times S)=:
(\Omega^{\bullet,\Gamma,pr}_{/S},F_{DR})((Y\times S,Z)/S).
\end{eqnarray*}
By definition $\Gr(\Omega_{/S})\circ\Gr_{S*}T(\Omega^{\Gamma}_{/S})=\Gr^O(\Omega_{/S})$.
\end{defi}

\begin{rem} 
\begin{itemize}
\item[(i)] Let $S\in\Var(\mathbb C)$. We have by definition  
$o_{12*}(\Omega^{\bullet,\Gamma}_{/S},F_b)=(\Omega^{\bullet}_{/S},F_b)\in C_{O_Sfil}(\Var(\mathbb C)^{sm}/S)$.
Moreover, if $S\in\SmVar(\mathbb C)$,  
$o_{12*}(\Omega^{\bullet,\Gamma}_{/S},F_b)=(\Omega^{\bullet}_{/S},F_b)\in C_{O_Sfil,D_S}(\Var(\mathbb C)^{sm}/S)$.
\item[(ii)] Let $S\in\Var(\mathbb C)$.
Then, $(\Omega^{\bullet,\Gamma}_{/S},F_b)\in C_{O_Sfil}(\Var(\mathbb C)^2/S)$ is a natural extension of 
\begin{equation*}
(\Omega^{\bullet,\Gamma}_{/S},F_b):=\rho_{S*}(\Omega^{\bullet,\Gamma}_{/S},F_b)
\in C_{O_Sfil}(\Var(\mathbb C)^{2,sm}/S),
\end{equation*}
but does NOT satisfy cdh descent.
\end{itemize}
\end{rem}

We have the following canonical transformation map given by the pullback of (relative) differential forms:

\begin{itemize}
\item Let $g:T\to S$ a morphism with $T,S\in\Var(\mathbb C)$.
We have the canonical morphism in $C_{g^*O_Sfil,g^*D_S}(\Var(\mathbb C)^{2,sm}/T)$ 
\begin{eqnarray*}
\Omega^{\Gamma}_{/(T/S)}:g^*(\Omega^{\bullet,\Gamma}_{/S},F_b)\to (\Omega^{\bullet,\Gamma}_{/T},F_b)
\end{eqnarray*}
induced by the pullback of differential forms : for $((V,Z_1)/T)=((V,Z_1),h)\in\Var(\mathbb C)^{2,sm}/T$,
\begin{eqnarray*}
\Omega^{\Gamma}_{/(T/S)}((V,Z_1)/T): \\
g^*\Omega^{\bullet,\Gamma}_{/S}((V,Z_1)/T):= 
\lim_{(h:(U,Z)\to S \mbox{sm},g_1:(V,Z_1)\to (U_T,Z_T),h,g)}\Omega^{\bullet,\Gamma}_{/S}((U,Z)/S) \\
\xrightarrow{\Omega^{\bullet,\Gamma}_{/S}(g'\circ g_1)}\Omega^{\bullet,\Gamma}_{/S}((V,Z_1)/S)
\xrightarrow{\Gamma^{\vee,h}_{Z_1}q(Y_1\times T)}\Omega^{\bullet,\Gamma}_{/T}((V,Z_1)/T),
\end{eqnarray*}
where $g':U_T:=U\times_S T\to U$ is the base change map and 
$q:\Omega^{\bullet}_{Y_1\times T/S}\to\Omega^{\bullet}_{Y_1\times T/T}$ is the quotient map.
If $T,S\in\SmVar(\mathbb C)$,
\begin{eqnarray*}
\Omega^{\Gamma}_{/(T/S)}:g^*(\Omega^{\bullet,\Gamma}_{/S},F_b)\to (\Omega^{\bullet,\Gamma}_{/T},F_b)
\end{eqnarray*}
is a morphism in $C_{g^*O_Sfil,g^*D_S}(\Var(\mathbb C)^{2,sm}/T)$ 
It induces the canonical morphisms in $C_{g^*O_Sfil,g^*D_S}(\Var(\mathbb C)^{2,sm}/T)$ :
\begin{eqnarray*}
E\Omega^{\Gamma}_{/(T/S)}:g^*E_{et}(\Omega^{\bullet,\Gamma}_{/S},F_b)
\xrightarrow{T(g,E_{et})(\Omega^{\bullet,\Gamma}_{/S},F_b)}
E_{et}(g^*(\Omega^{\bullet,\Gamma}_{/S},F_b))
\xrightarrow{E_{et}(\Omega^{\Gamma}_{/(T/S)})}E_{et}(\Omega^{\bullet,\Gamma}_{/T},F_b)
\end{eqnarray*}
and
\begin{eqnarray*}
E\Omega^{\Gamma}_{/(T/S)}:g^*E_{zar}(\Omega^{\bullet,\Gamma}_{/S},F_b)
\xrightarrow{T(g,E_{zar})(\Omega^{\bullet,\Gamma}_{/S},F_b)}
E_{zar}(g^*(\Omega^{\bullet,\Gamma}_{/S},F_b))
\xrightarrow{E_{zar}(\Omega^{\Gamma}_{/(T/S)})}E_{zar}(\Omega^{\bullet,\Gamma}_{/T},F_b).
\end{eqnarray*}
\item Let $g:T\to S$ a morphism with $T,S\in\SmVar(\mathbb C)$.
We have the canonical morphism in $C_{g^*D_Sfil}(\Var(\mathbb C)^{2,smpr}/T)$ 
\begin{eqnarray*}
\Omega^{\Gamma,pr}_{/(T/S)}:g^*(\Omega^{\bullet,\Gamma,pr}_{/S},F_{DR})\to (\Omega^{\bullet,\Gamma,pr}_{/T},F_{DR})
\end{eqnarray*}
induced by the pullback of differential forms : 
for $((Y_1\times T,Z_1)/T)=((Y_1\times T,Z_1),p)\in\Var(\mathbb C)^{2,smpr}/T$,
\begin{eqnarray*}
\Omega^{\Gamma,pr}_{/(T/S)}((Y_1\times T,Z_1)/T): \\
g^*\Omega^{\bullet,\Gamma,pr}_{/S}((Y_1\times T,Z_1)/T):= 
\lim_{(h:(Y\times S,Z)\to S, \, g_1:(Y_1\times T,Z_1)\to (Y\times T,Z_T),h,g)}
\Omega^{\bullet,\Gamma,pr}_{/S}((Y\times T,Z)/S) \\
\xrightarrow{\Omega^{\bullet,\Gamma,pr}_{/S}(g'\circ g_1)}\Omega^{\bullet,\Gamma,pr}_{/S}((Y_1\times T,Z_1)/S)
\xrightarrow{q(-)(Y_1\times T)}\Omega^{\bullet,\Gamma,pr}_{/T}((Y_1\times T,Z_1)/T), 
\end{eqnarray*}
where $g'=(I_Y\times g):Y\times T\to Y\times S$ is the base change map and 
$q(M):\Omega_{Y_1\times T/S}\otimes_{O_{Y_1\times T}}(M,F)\to\Omega_{Y_1\times T/T}\otimes_{O_{Y_1\times T}}(M,F)$
is the quotient map. 
It induces the canonical morphisms in $C_{g^*D_Sfil}(\Var(\mathbb C)^{2,smpr}/T)$ :
\begin{eqnarray*}
E\Omega^{\Gamma,pr}_{/(T/S)}:g^*E_{et}(\Omega^{\bullet,\Gamma,pr}_{/S},F_{DR})
\xrightarrow{T(g,E)(-)} E_{et}(g^*(\Omega^{\bullet,\Gamma,pr}_{/S},F_{DR}))
\xrightarrow{E_{et}(\Omega^{\Gamma,pr}_{/(T/S)})}E_{et}(\Omega^{\bullet,\Gamma,pr}_{/T},F_{DR})
\end{eqnarray*}
and
\begin{eqnarray*}
E\Omega^{\Gamma,pr}_{/(T/S)}:g^*E_{zar}(\Omega^{\bullet,\Gamma,pr}_{/S},F_{DR})
\xrightarrow{T(g,E)(-)} E_{zar}(g^*(\Omega^{\bullet,\Gamma,pr}_{/S},F_{DR}))
\xrightarrow{E_{zar}(\Omega^{\Gamma,pr}_{/(T/S)})}E_{zar}(\Omega^{\bullet,\Gamma,pr}_{/T},F_{DR}).
\end{eqnarray*}
\end{itemize}

\begin{defi}\label{TgDR}
Let $g:T\to S$ a morphism with $T,S\in\SmVar(\mathbb C)$.
We have, for $F\in C(\Var(\mathbb C)^{2,smpr}/S)$, the canonical transformation in $C_{\mathcal Dfil}(T)$ :
\begin{eqnarray*}
T(g,\Omega^{\Gamma,pr}_{/\cdot})(F): 
g^{*mod}L_De(S)_*\Gr^{12}_{S*}\mathcal Hom^{\bullet}(F,E_{et}(\Omega^{\bullet,\Gamma,pr}_{/S},F_{DR})) \\
\xrightarrow{=} 
(g^*L_De(S)_*\mathcal Hom^{\bullet}(F,E_{et}(\Omega^{\bullet,\Gamma,pr}_{/S},F_{DR})))\otimes_{g^*O_S}O_T \\
\xrightarrow{T(g,\Gr^{12})(-)\circ T(e,g)(-)\circ q}  
e(T)_*\Gr^{12}_{T*}g^*\mathcal Hom^{\bullet}(F,E_{et}(\Omega^{\bullet,\Gamma,pr}_{/S},F_{DR}))\otimes_{g^*O_S}O_T \\ 
\xrightarrow{(T(g,hom)(-,-)\otimes I)}
e(T)_*\Gr^{12}_{T*}\mathcal Hom^{\bullet}(g^*F,g^*E_{et}(\Omega^{\bullet,\Gamma,pr}_{/S},F_{DR}))\otimes_{g^*O_S}O_T \\
\xrightarrow{ev(hom,\otimes)(-,-,-)} 
e(T)_*\Gr^{12}_{T*}\mathcal Hom^{\bullet}(g^*F,g^*E_{et}(\Omega^{\bullet,\Gamma,pr}_{/S},F_{DR}))
\otimes_{g^*e(S)^*O_S}e(T)^*O_T \\
\xrightarrow{\mathcal Hom^{\bullet}(g^*F,(E\Omega^{\Gamma,pr}_{/(T/S)}\otimes m))}
e(T)_*\Gr^{12}_{T*}\mathcal Hom^{\bullet}(g^*F,E_{et}(\Omega^{\bullet,\Gamma,pr}_{/T},F_{DR}))
\end{eqnarray*}
\end{defi}

Let $S\in\Var(\mathbb C)$. Recall that for and $h:U\to S$ a morphism with $U\in\Var(\mathbb C)$,
we have the canonical map given by the wedge product
\begin{equation*}
w_{U/S}:\Omega^{\bullet}_{U/S}\otimes_{O_S}\Omega^{\bullet}_{U/S}\to\Omega^{\bullet}_{U/S}; 
\alpha\otimes\beta\mapsto\alpha\wedge\beta.
\end{equation*}
Let $S\in\Var(\mathbb C)$ and $h_1:U_1\to S$, $h_2:U_2\to S$ two morphisms with $U_1,U_2\in\Var(\mathbb C)$.
Denote $h_{12}:U_{12}:=U_1\times_S U_2\to S$ and $p_{112}:U_1\times_S U_2\to U_1$, $p_{212}:U_1\times_S U_2\to U_2$ the projections.
Recall we have the canonical map given by the wedge product
\begin{equation*}
w_{(U_1,U_2)/S}:p_{112}^*\Omega^{\bullet}_{U_1/S}\otimes_{O_S}p_{212}^*\Omega^{\bullet}_{U_2/S}
\to \Omega^{\bullet}_{U_{12}/S}; \alpha\otimes\beta\mapsto p_{112}^*\alpha\wedge p_{212}^*\beta
\end{equation*}
which gives the map
\begin{eqnarray*}
Ew_{(U_1,U_2)/S}:h_{1*}E_{zar}(\Omega^{\bullet}_{U_1/S})\otimes_{O_S}h_{2*}E_{zar}(\Omega^{\bullet}_{U_2/S})
\to h_{12*}E_{zar}(p_{112}^*\Omega^{\bullet}_{U_1/S}\otimes_{O_S}p_{212}^*\Omega^{\bullet}_{U_2/S})
\end{eqnarray*}

Let $S\in\SmVar(\mathbb C)$. 
\begin{itemize}
\item The complex of presheaves $(\Omega^{\bullet,\Gamma}_{/S},F_b)\in C_{O_Sfil,D_S}(\Var(\mathbb C)^{2,sm}/S)$ 
have a monoidal structure given by the wedge product of differential forms: 
for $h:(U,Z)\to S\in\Var(\mathbb C)^2/S$, the map 
\begin{eqnarray*}
DR(-)(\gamma^{\vee,h}_Z(-))\circ w_{U/S}:(\Omega^{\bullet}_{U/S},F_b)\otimes_{p^*O_S}(\Omega^{\bullet}_{U/S},F_b)
\to\Gamma_Z^{\vee,h}L_{h^*O_S}(\Omega^{\bullet}_{U/S},F_b)
\end{eqnarray*}
factors trough
\begin{eqnarray*}
DR(-)(\gamma^{\vee,Hdg}_Z(-))\circ w_{U/S}:
(\Omega^{\bullet}_{U/S},F_b)\otimes_{p^*O_S}(\Omega^{\bullet}_{U/S},F_b) \\
\xrightarrow{DR(-)(\gamma_Z^{\vee,h}(-))\otimes DR(-)(\gamma_Z^{\vee,h}(-))} 
\Gamma_Z^{\vee,h}L_{h^*O_S}(\Omega^{\bullet}_{U/S},F_b)\otimes_{p^*O_S}
\Gamma_Z^{\vee,h}L_{h^*O_S}(\Omega^{\bullet}_{U/S},F_b) \\
\xrightarrow{(DR(-)(\gamma^{\vee,h}_Z(-))\circ w_{U/S})^{\gamma}}
\Gamma_Z^{\vee,h}L_{h^*O_S}(\Omega^{\bullet}_{U/S},F_b)
\end{eqnarray*}
unique up to homotopy, giving the map in $C_{O_Sfil,D_S}(\Var(\mathbb C)^{2,smpr}/S)$:
\begin{eqnarray*}
w_S:(\Omega^{\bullet,\Gamma}_{/S},F_b)\otimes_{O_S}(\Omega^{\bullet,\Gamma}_{/S},F_b)
\to(\Omega^{\bullet,\Gamma}_{/S},F_b)
\end{eqnarray*}
given by for $h:(U,Z)\to S\in\Var(\mathbb C)^{2,sm}/S$,
\begin{eqnarray*}
w_S((U,Z)/S):(\Gamma_Z^{\vee,h}L_{h^*O_S}(\Omega^{\bullet}_{U/S},F_b)\otimes_{p^*O_S}
\Gamma_Z^{\vee,h}L_{h^*O_S}(\Omega^{\bullet}_{U/S},F_b))(U) \\
\xrightarrow{(DR(-)(\gamma^{\vee,h}_Z(-))\circ w_{U/S})^{\gamma}(U)}
\Gamma_Z^{\vee,h}L_{h^*O_S}(\Omega^{\bullet}_{U/S},F_b)(U)
\end{eqnarray*}
which induces the map in $C_{O_Sfil,D_S}(\Var(\mathbb C)^{2,sm}/S)$
\begin{eqnarray*}
Ew_S:E_{et}(\Omega^{\bullet,\Gamma}_{/S},F_b)\otimes_{O_S} E_{et}(\Omega^{\bullet,\Gamma}_{/S},F_b)
\xrightarrow{=}
E_{et}((\Omega^{\bullet,\Gamma}_{/S},F_b)\otimes_{O_S}(\Omega^{\bullet,\Gamma}_{/S},F_b))
\xrightarrow{E_{et}(w_S)} E_{et}(\Omega^{\bullet,\Gamma}_{/S},F_b)
\end{eqnarray*}
given by the functoriality of the Godement resolution (see section 2).
\item The complex of presheaves 
$(\Omega^{\bullet,\Gamma,pr}_{/S},F_{DR})\in C_{D_Sfil}(\Var(\mathbb C)^{2,smpr}/S)$ 
have a monoidal structure given by the wedge product of differential forms: 
for $p:(Y\times S,Z)\to S\in\Var(\mathbb C)^{2,smpr}/S$, the map 
\begin{eqnarray*}
DR(-)(\gamma^{\vee,Hdg}_Z(-))\circ w_{Y\times S/S}:
(\Omega^{\bullet}_{Y\times S/S}\otimes_{O_{Y\times S}}(O_{Y\times S},F_b))\otimes_{p^*O_S}
(\Omega^{\bullet}_{Y\times S/S}\otimes_{O_{Y\times S}}(O_{Y\times S},F_b)) \\
\to\Omega^{\bullet}_{Y\times S/S}\otimes_{O_{Y\times S}}\Gamma_Z^{\vee,Hdg}(O_{Y\times S},F_b)
\end{eqnarray*}
factors trough
\begin{eqnarray*}
DR(-)(\gamma^{\vee,Hdg}_Z(-))\circ w_{Y\times S/S}: \\
(\Omega^{\bullet}_{Y\times S/S}\otimes_{O_{Y\times S}}(O_{Y\times S},F_b))\otimes_{p^*O_S}
(\Omega^{\bullet}_{Y\times S/S}\otimes_{O_{Y\times S}}(O_{Y\times S},F_b)) \\
\xrightarrow{DR(-)(\gamma_Z^{\vee,Hdg}(-))\otimes DR(-)(\gamma_Z^{\vee,Hdg}(-))} \\
(\Omega^{\bullet}_{Y\times S/S}\otimes_{O_{Y\times S}}\Gamma_Z^{\vee,Hdg})(O_{Y\times S},F_b)\otimes_{p^*O_S}
\Omega^{\bullet}_{Y\times S/S}\otimes_{O_{Y\times S}}\Gamma_Z^{\vee,Hdg}(O_{Y\times S},F_b) \\
\xrightarrow{(DR(-)(\gamma^{\vee,Hdg}_Z(-))\circ w_{Y\times S/S})^{\gamma}}
\Omega^{\bullet}_{Y\times S/S}\otimes_{O_{Y\times S}}\Gamma_Z^{\vee,Hdg}(O_{Y\times S},F_b)
\end{eqnarray*}
unique up to homotopy, giving the map in $C_{D_Sfil}(\Var(\mathbb C)^{2,smpr}/S)$:
\begin{eqnarray*}
w_S:(\Omega^{\bullet,\Gamma,pr}_{/S},F_{DR})\otimes_{O_S}(\Omega^{\bullet,\Gamma,pr}_{/S},F_{DR})
\to(\Omega^{\bullet,\Gamma,pr}_{/S},F_{DR})
\end{eqnarray*}
given by for $p:(Y\times S,Z)\to S\in\Var(\mathbb C)^{2,smpr}/S$,
\begin{eqnarray*}
w_S((Y\times S,Z)/S): \\
(((\Omega^{\bullet}_{Y\times S/S}\otimes_{O_{Y\times S}}\Gamma_Z^{\vee,Hdg})(O_{Y\times S},F_b))\otimes_{p^*O_S}
(\Omega^{\bullet}_{Y\times S/S}\otimes_{O_{Y\times S}}\Gamma_Z^{\vee,Hdg}(O_{Y\times S},F_b)))(Y\times S) \\
\xrightarrow{(DR(-)(\gamma^{\vee,Hdg}_Z(-))\circ w_{Y\times S/S})^{\gamma}(Y\times S)}
(\Omega^{\bullet}_{Y\times S/S}\otimes_{O_{Y\times S}}\Gamma_Z^{\vee,Hdg}(O_{Y\times S},F_b))(Y\times S)
\end{eqnarray*}
which induces the map in $C_{D_Sfil}(\Var(\mathbb C)^{2,smpr}/S)$
\begin{eqnarray*}
Ew_S:E_{et}(\Omega^{\bullet,\Gamma,pr}_{/S},F_{DR})\otimes_{O_S} E_{et}(\Omega^{\bullet,\Gamma,pr}_{/S},F_{DR})
\xrightarrow{=} \\
E_{et}((\Omega^{\bullet,\Gamma,pr}_{/S},F_{DR})\otimes_{O_S}(\Omega^{\bullet,\Gamma,pr}_{/S},F_{DR}))
\xrightarrow{E_{et}(w_S)} E_{et}(\Omega^{\bullet,\Gamma,pr}_{/S},F_{DR})
\end{eqnarray*}
by the functoriality of the Godement resolution (see section 2).
\end{itemize}

\begin{defi}\label{TotimesDR}
Let $S\in\SmVar(\mathbb C)$.
We have, for $F,G\in C(\Var(\mathbb C)^{2,smpr}/S)$, the canonical transformation in $C_{\mathcal Dfil}(S)$ :
\begin{eqnarray*}
T(\otimes,\Omega)(F,G): \\
e(S)_*\Gr^{12}_{S*}\mathcal Hom(F,E_{et}(\Omega^{\bullet,\Gamma,pr}_{/S},F_{DR}))\otimes_{O_S}
e(S)_*\Gr^{12}_{S*}\mathcal Hom(G,E_{et}(\Omega^{\bullet,\Gamma,pr}_{/S},F_{DR})) \\
\xrightarrow{=} 
e(S)_*\Gr^{12}_{S*}(\mathcal Hom(F,E_{et}(\Omega^{\bullet,\Gamma,pr}_{/S},F_{DR}))\otimes_{O_S}
\mathcal Hom(G,E_{et}(\Omega^{\bullet,\Gamma,pr}_{/S},F_{DR}))) \\ 
\xrightarrow{T(\mathcal Hom,\otimes)(-)} 
e(S)_*\Gr^{12}_{S*}\mathcal Hom(F\otimes G,E_{et}(\Omega^{\bullet,\Gamma,pr}_{/S},F_{DR})\otimes_{O_S}
E_{et}(\Omega^{\bullet,\Gamma,pr}_{/S},F_{DR}))) \\
\xrightarrow{=} 
e(S)_*\Gr^{12}_{S*}\mathcal Hom(F\otimes G,(E_{et}(\Omega^{\bullet,\Gamma,pr}_{/S},F_{DR})
\otimes_{O_S}E_{et}(\Omega^{\bullet,\Gamma,pr}_{/S},F_{DR}))) \\
\xrightarrow{\mathcal Hom(F\otimes G,Ew_S)} 
e(S)_*\Gr^{12}_{S*}\mathcal Hom(F\otimes G,E_{et}(\Omega^{\bullet,\Gamma,pr}_{/S},F_{DR})).
\end{eqnarray*}
\end{defi}

Let $S\in\Var(\mathbb C)$. Let $S=\cup_{i=1}^l S_i$ an open affine cover and denote by $S_I=\cap_{i\in I} S_i$.
Let $i_i:S_i\hookrightarrow\tilde S_i$ closed embeddings, with $\tilde S_i\in\Var(\mathbb C)$. 
For $I\subset\left[1,\cdots l\right]$, denote by $\tilde S_I=\Pi_{i\in I}\tilde S_i$.
We then have closed embeddings $i_I:S_I\hookrightarrow\tilde S_I$ and for $J\subset I$ the following commutative diagram
\begin{equation*}
D_{IJ}=\xymatrix{ S_I\ar[r]^{i_I} & \tilde S_I \\
S_J\ar[u]^{j_{IJ}}\ar[r]^{i_J} & \tilde S_J\ar[u]^{p_{IJ}}}  
\end{equation*}
where $p_{IJ}:\tilde S_J\to\tilde S_I$ is the projection
and $j_{IJ}:S_J\hookrightarrow S_I$ is the open embedding so that $j_I\circ j_{IJ}=j_J$.
This gives the diagram of algebraic varieties $(\tilde S_I)\in\Fun(\mathcal P(\mathbb N),\Var(\mathbb C))$ which
the diagram of sites $\Var(\mathbb C)^{2,smpr}/(\tilde S_I)\in\Fun(\mathcal P(\mathbb N),\Cat)$. 
This gives also the diagram of algebraic varieties $(\tilde S_I)^{op}\in\Fun(\mathcal P(\mathbb N)^{op},\Var(\mathbb C))$ which
the diagram of sites $\Var(\mathbb C)^{2,smpr}/(\tilde S_I)^{op}\in\Fun(\mathcal P(\mathbb N)^{op},\Cat)$. 
We then get
\begin{eqnarray*}
((\Omega^{\bullet,\Gamma,pr}_{/(\tilde S_I)},F_{DR})[-d_{\tilde S_I}],T_{IJ})
\in C_{D_{(\tilde S_I)}fil}(\Var(\mathbb C)^{2,smpr}/(\tilde S_I))
\end{eqnarray*}
with
\begin{eqnarray*}
T_{IJ}:(\Omega^{\bullet,\Gamma,pr}_{/\tilde S_I},F_{DR})[-d_{\tilde S_I}]
\xrightarrow{\ad(p_{IJ}^{*mod[-]},p_{IJ*}(-)}
p_{IJ*}p_{IJ}^*(\Omega^{\bullet,\Gamma,pr}_{/\tilde S_I},F_{DR})
\otimes_{p_{IJ}^*O_{\tilde S_I}}O_{\tilde S_J}[-d_{\tilde S_J}] \\
\xrightarrow{m\circ p_{IJ*}\Omega^{\Gamma,pr}_{/(\tilde S_J/\tilde S_I)}[-d_{\tilde S_J}]}
p_{IJ*}(\Omega^{\bullet,\Gamma,pr}_{/\tilde S_J},F_{DR})[-d_{\tilde S_J}].
\end{eqnarray*}
For $(G_I,K_{IJ})\in C(\Var(\mathbb C)^{2,smpr}/(\tilde S_I)^{op})$, we denote (see section 2)
\begin{eqnarray*}
e'((\tilde S_I))_*\mathcal Hom((G_I,K_{IJ}),(E_{zar}(\Omega^{\bullet,\Gamma,pr}_{/(\tilde S_I)},F_{DR})[-d_{\tilde S_I}],T_{IJ})):= \\
(e'(\tilde S_I)_*\mathcal Hom(G_I,E_{zar}(\Omega^{\bullet,\Gamma,pr}_{/\tilde S_I},F_{DR}))[-d_{\tilde S_I}],u_{IJ}((G_I,K_{IJ})))
\in C_{\mathcal Dfil}((\tilde S_I))
\end{eqnarray*}
with
\begin{eqnarray*}
u_{IJ}((G_I,K_{IJ})):e'(\tilde S_I)_*\mathcal Hom(G_I,E_{zar}(\Omega^{\bullet,\Gamma,pr}_{/\tilde S_I},F_{DR}))[-d_{\tilde S_I}] \\
\xrightarrow{\ad(p_{IJ}^{*mod[-]},p_{IJ*})(-)\circ T(p_{IJ},e)(-)}
p_{IJ*}e'(\tilde S_J)_*p_{IJ}^*\mathcal Hom(G_I,E_{zar}(\Omega^{\bullet,\Gamma,pr}_{/\tilde S_I},F_{DR}))
\otimes_{p_{IJ}^*O_{\tilde S_I}}O_{\tilde S_J}[-d_{\tilde S_J}] \\
\xrightarrow{T(p_{IJ},hom)(-,-)}
p_{IJ*}e'(\tilde S_J)_*\mathcal Hom(p_{IJ}^*G_I,p_{IJ}^*E_{zar}(\Omega^{\bullet,\Gamma,pr}_{/\tilde S_I},F_{DR})) 
\otimes_{p_{IJ}^*O_{\tilde S_I}}O_{\tilde S_J}[-d_{\tilde S_J}] \\
\xrightarrow{m\circ\mathcal Hom(p_{IJ}^*G_I,T_{IJ})}
p_{IJ*}e'(\tilde S_J)_*\mathcal Hom(p_{IJ}^*G_I,E_{zar}(\Omega^{\bullet,\Gamma,pr}_{/\tilde S_J},F_{DR}))[-d_{\tilde S_J}] \\
\xrightarrow{\mathcal Hom(K_{IJ},E_{zar}(\Omega^{\bullet,\Gamma,pr}_{/\tilde S_J},F_{DR}))}
p_{IJ*}e'(\tilde S_J)_*\mathcal Hom(G_J,E_{zar}(\Omega^{\bullet,\Gamma,pr}_{/\tilde S_J},F_{DR}))[-d_{\tilde S_J}].
\end{eqnarray*}
This gives in particular
\begin{eqnarray*}
(\Omega^{\bullet,\Gamma,pr}_{/(\tilde S_I)},F_{DR})[-d_{\tilde S_I}],T_{IJ})\in 
C_{D_{(\tilde S_I)}fil}(\Var(\mathbb C)^{2,(sm)pr}/(\tilde S_I)^{op}).
\end{eqnarray*}

We now define the filtered De Rahm realization functor.

\begin{defi}\label{DRalgdefFunct}
\begin{itemize}
\item[(i)]Let $S\in\SmVar(\mathbb C)$.
We have, using definition \ref{wtildew} and definition \ref{RCHhatdef}, the functor
\begin{eqnarray*}
\mathcal F_S^{FDR}:C(\Var(\mathbb C)^{sm}/S)\to C_{\mathcal Dfil}(S), \; F\mapsto \\ 
\mathcal F_S^{FDR}(F):=e(S)_*\Gr^{12}_{S*}\mathcal Hom^{\bullet}(\hat R^{CH}(\rho_S^*L(F)),
E_{zar}(\Omega^{\bullet,\Gamma,pr}_{/S},F_{DR}))[-d_S]
\end{eqnarray*}
\item[(ii)]Let $S\in\Var(\mathbb C)$ and $S=\cup_{i=1}^l S_i$ an open cover such that there exist closed embeddings
$i_i:S_i\hookrightarrow\tilde S_i$  with $\tilde S_i\in\SmVar(\mathbb C)$. 
For $I\subset\left[1,\cdots l\right]$, denote by $S_I:=\cap_{i\in I} S_i$ and $j_I:S_I\hookrightarrow S$ the open embedding.
We then have closed embeddings $i_I:S_I\hookrightarrow\tilde S_I:=\Pi_{i\in I}\tilde S_i$.
Consider, for $I\subset J$, the following commutative diagram
\begin{equation*}
D_{IJ}=\xymatrix{ S_I\ar[r]^{i_I} & \tilde S_I \\
S_J\ar[u]^{j_{IJ}}\ar[r]^{i_J} & \tilde S_J\ar[u]^{p_{IJ}}}  
\end{equation*}
and $j_{IJ}:S_J\hookrightarrow S_I$ is the open embedding so that $j_I\circ j_{IJ}=j_J$.
We have, using definition \ref{wtildew} and definition \ref{RCHhatdef}, the functor
\begin{eqnarray*}
\mathcal F_S^{FDR}:C(\Var(\mathbb C)^{sm}/S)\to C_{\mathcal Dfil}(S/(\tilde S_I)), \; F\mapsto \\
\mathcal F_S^{FDR}(F):=
e'((\tilde S_I))_*\mathcal Hom^{\bullet}(
(\hat R^{CH}(\rho_{\tilde S_I}^*L(i_{I*}j_I^*F)),\hat R^{CH}_{\tilde S_J}(T^q(D_{IJ})(j_I^*F))), \\
(E_{zar}(\Omega^{\bullet,\Gamma,pr}_{/(\tilde S_I)},F_{DR})[-d_{\tilde S_I}],T_{IJ})) \\
:=(e'(\tilde S_I)_*\mathcal Hom^{\bullet}(\hat R^{CH}(\rho_{\tilde S_I}^*L(i_{I*}j_I^*F)), 
E_{zar}(\Omega^{\bullet,\Gamma,pr}_{/\tilde S_I},F_{DR}))[-d_{\tilde S_I}],u^q_{IJ}(F))
\end{eqnarray*}
where we have denoted for short $e'(\tilde S_I)=e(\tilde S_I)\circ\Gr^{12}_{\tilde S_I}$, and
\begin{eqnarray*}
u^q_{IJ}(F)[d_{\tilde S_J}]:
e'(\tilde S_I)_*\mathcal Hom^{\bullet}(\hat R^{CH}(\rho_{\tilde S_I}^*L(i_{I*}j_I^*F)),
E_{zar}(\Omega^{\bullet,\Gamma,pr}_{/\tilde S_I},F_{DR})) \\
\xrightarrow{\ad(p_{IJ}^{*mod},p_{IJ})(-)} 
p_{IJ*}p_{IJ}^{*mod}e'(\tilde S_I)_*\mathcal Hom^{\bullet}(\hat R^{CH}(\rho_{\tilde S_I}^*L(i_{I*}j_I^*F)),
E_{zar}(\Omega^{\bullet,\Gamma,pr}_{/\tilde S_I},F_{DR})) \\
\xrightarrow{p_{IJ*}T(p_{IJ},\Omega^{\gamma,pr}_{\cdot})(-)}  
p_{IJ*}e'(\tilde S_J)_*\mathcal Hom^{\bullet}(p_{IJ}^*\hat R^{CH}(\rho_{\tilde S_I}^*L(i_{I*}j_I^*F)),
E_{zar}(\Omega^{\bullet,\Gamma,pr}_{/\tilde S_J},F_{DR})) \\
\xrightarrow{\mathcal Hom(T(p_{IJ},\hat R^{CH})(Li_{I*}j_I^*F)^{-1},E_{zar}(\Omega_{/\tilde S_J}^{\bullet,\Gamma,pr},F_{DR}))} \\
p_{IJ*}e'(\tilde S_J)_*\mathcal Hom^{\bullet}(\hat R^{CH}(\rho_{\tilde S_J}^*p_{IJ}^*L(i_{I*}j_I^*F)),
E_{zar}(\Omega^{\bullet,\Gamma,pr}_{/\tilde S_J},F_{DR})) \\
\xrightarrow{\mathcal Hom(\hat R^{CH}_{\tilde S_J}(T^q(D_{IJ})(j_I^*F)),
E_{zar}(\Omega_{/\tilde S_J}^{\bullet,\Gamma,pr},F_{DR}))} \\
p_{IJ*}e'(\tilde S_J)_*\mathcal Hom^{\bullet}(\hat R^{CH}(\rho_{\tilde S_J}^*L(i_{J*}j_J^*F)),
E_{zar}(\Omega^{\bullet,\Gamma,pr}_{/\tilde S_J},F_{DR})).
\end{eqnarray*}
For $I\subset J\subset K$, we have obviously $p_{IJ*}u_{JK}(F)\circ u_{IJ}(F)=u_{IK}(F)$.
\end{itemize}
\end{defi}

Recall, see section 2, that we have the projection morphisms of sites 
$p_a:\Var(\mathbb C)^{2,smpr}/(\tilde S_I)^{op}\to\Var(\mathbb C)^{2,smpr}/(\tilde S_I)^{op}$
given by the functor 
\begin{eqnarray*}
p_a:\Var(\mathbb C)^{2,smpr}/(\tilde S_I)^{op}\to\Var(\mathbb C)^{2,smpr}/(\tilde S_I)^{op}, \\
p_a((Y_I\times\tilde S_I,Z_I)/\tilde S_I,s_{IJ}):=
((Y_I\times\mathbb A^1\times\tilde S_I,Z_I\times\mathbb A^1)/\tilde S_I,s_{IJ}\times I), \\ 
p_a((g_I):((Y'_I\times\tilde S_I,Z'_I)/\tilde S_I,s'_{IJ})\to((Y_I\times\tilde S_I,Z_I)/\tilde S_I,s_{IJ}))= \\
(g_I\times I):((Y'_I\times\mathbb A^1\times\tilde S_I,Z'_I\times\mathbb A^1)/\tilde S_I,s'_{IJ}\times I)
\to((Y_I\times\mathbb A^1\times\tilde S_I,Z_I\times\mathbb A^1)/\tilde S_I),s_{IJ}\times I)).
\end{eqnarray*}

We have the following key proposition :

\begin{prop}\label{aetfib}
\begin{itemize}
\item[(i)]Let $S\in\Var(\mathbb C)$. Let $S=\cup_{i=1}^l S_i$ an open cover such that there exist closed embeddings
$i_i:S_i\hookrightarrow\tilde S_i$ with $\tilde S_i\in\SmVar(\mathbb C)$.
The complex of presheaves
$(\Omega^{\bullet,\Gamma,pr}_{/(\tilde S_I)},F_{DR})\in C_{D_{(\tilde S_I)}fil}(\Var(\mathbb C)^{2,smpr}/(\tilde S_I)^{op})$ 
is $2$-filtered $\mathbb A^1$ homotopic, that is
\begin{equation*}
\ad(p_a^*,p_{a*})(\Omega^{\bullet,\Gamma,pr}_{/(\tilde S_I)},F_{DR}):
(\Omega^{\bullet,\Gamma,pr}_{/S},F_{DR})\to p_{a*}p_a^*(\Omega^{\bullet,\Gamma,pr}_{/(\tilde S_I)},F_{DR})
\end{equation*}
is a $2$-filtered homotopy.
\item[(i2)]Let $S\in\SmVar(\mathbb C)$. The complex of presheaves
$(\Omega^{\bullet,\Gamma,pr}_{/S},F_{DR})\in C_{D_Sfil}(\Var(\mathbb C)^{2,smpr}/S)$ 
admits transferts, i.e. 
\begin{equation*}
\Tr(S)_*\Tr(S)^*(\Omega^{\bullet,\Gamma,pr}_{/S},F_{DR}=(\Omega^{\bullet,\Gamma,pr}_{/S},F_{DR}).
\end{equation*}
\item[(ii)]Let $S\in\Var(\mathbb C)$. Let $S=\cup_{i=1}^l S_i$ an open cover such that there exist closed embeddings
$i_i:S_i\hookrightarrow\tilde S_i$ with $\tilde S_i\in\SmVar(\mathbb C)$.
Let $m=(m_I):(Q_{1I},K^1_{IJ})\to(Q_{2I},K^2_{IJ})$ be an equivalence $(\mathbb A^1,et)$ local with 
$(Q_{1I},K_{IJ})\to(Q_{2I},K_{IJ})\in C(\Var(\mathbb C)^{smpr}/(\tilde S_I)^{op})$ complexes of representable presheaves. 
Then, the map in $C_{\mathcal Dfil}((\tilde S_I))$
\begin{eqnarray*} 
M:=(e(\tilde S_I)_*\mathcal Hom^{\bullet}(m_I,E_{zar}(\Omega^{\bullet,\Gamma,pr}_{/\tilde S_I},F_{DR})[-d_{\tilde S_I}])): \\  
e'((\tilde S_I))_*\mathcal Hom^{\bullet}((Q_{2I},K^1_{IJ}),
(E_{zar}(\Omega^{\bullet,\Gamma,pr}_{/(\tilde S_I)},F_{DR})[-d_{\tilde S_I}],T_{IJ})) \\
\to e'((\tilde S_I))_*\mathcal Hom^{\bullet}((Q_{1I},K^1_{IJ}),
(E_{zar}(\Omega^{\bullet,\Gamma,pr}_{/(\tilde S_I)},F_{DR})[-d_{\tilde S_I}],T_{IJ}))
\end{eqnarray*}
is a $2$-filtered quasi-isomorphism. It is thus an isomorphism in $D_{\mathcal Dfil,\infty}((\tilde S_I))$.
\end{itemize}
\end{prop}

\begin{proof}  
\noindent(i1):Let $(Y\times S,Z)/S=((Y\times S,Z),p)\in\Var(\mathbb C)^{2,smpr}/S$ so that
$p_a:(Y\times\mathbb A^1\times S,Z\times\mathbb A^1)\to(Y\times S,Z)$ is the projection, and
$i_0:(Y\times S,Z)\to(Y\times\mathbb A^1\times S,Z\times\mathbb A^1)$ the closed embedding. 
Then,
\begin{eqnarray*}
a(p_a,Z):p_a^{*mod}\Gamma_Z^{\vee,Hdg}(O_{Y\times\times S},F_b))\to
\Gamma_{Z\times\mathbb A^1}^{\vee,Hdg}(O_{Y\times\mathbb A^1\times S},F_b)).
\end{eqnarray*}
a quasi-isomorphism in $\pi_{Y\times\mathbb A^1\times S}(C(MHM(Y\times\mathbb A^1\times S)))$.
Since a morphism of mixed Hodge module is strict for the F-filtration,
\begin{eqnarray*}
a(p_a,Z):p_a^{*mod}\Gamma_Z^{\vee,Hdg}(O_{Y\times\times S},F_b)\to
\Gamma_{Z\times\mathbb A^1}^{\vee,Hdg}(O_{Y\times\mathbb A^1\times S},F_b).
\end{eqnarray*}
is a filtered quasi-isomorphism in $C_{\mathcal Dfil}(Y\times\mathbb A^1\times S)$. 
Hence, as
\begin{eqnarray*} 
I(p_a^*,p_{a*})(-,-)(\Omega_{(Y\times\mathbb A^1\times S/Y\times S)(S/S)}(\Gamma_Z^{\vee,Hdg}(O_{Y\times S},F_b))): \\ 
(\Omega^{\bullet}_{Y\times S/S},F_b)\otimes_{O_{Y\times S}}\Gamma_Z^{\vee,Hdg}(O_{Y\times S},F_b) \\
\to p_{a*}((\Omega^{\bullet}_{Y\times\mathbb A^1\times S/S},F_b)\otimes_{O_{Y\times\mathbb A^1\times S}}
p_a^{*mod}\Gamma_Z^{\vee,Hdg}(O_{Y\times\times S},F_b))
\end{eqnarray*} 
is a $2$-filtered homotopy equivalence whose inverse is 
\begin{eqnarray*} 
p_{a*}I(i_0^*,i_{0*})(-,-)
(\Omega_{(Y\times S/Y\times\mathbb A^1\times S)(S/S)}(p_a^{*mod}(\Gamma_Z^{\vee,Hdg}(O_{Y\times S},F_b)))): \\ 
p_{a*}((\Omega^{\bullet}_{Y\times\mathbb A^1\times S/S},F_b)\otimes_{O_{Y\times\mathbb A^1\times S}} 
p_a^{*mod}(\Gamma_Z^{\vee,Hdg}(O_{Y\times\times S},F_b))) \\
\to p_{a*}i_{0*}(\Omega^{\bullet}_{Y\times S/S}\otimes_{O_{Y\times S}}
(i_0^{*mod}p_a^{*mod}(\Gamma_Z^{\vee,Hdg}(O_{Y\times S},F_b)) \\
\xrightarrow{=}(\Omega^{\bullet}_{Y\times S/S},F_b)\otimes_{O_{Y\times S}}(\Gamma_Z^{\vee,Hdg}(O_{Y\times S},F_b))
\end{eqnarray*} 
(see the proof of proposition \ref{aetfibGM}), the map
\begin{eqnarray*} 
\ad(p_a^*,p_{a*})(\Omega^{\bullet,\Gamma,pr}_{/S},F_{DR}):
(\Omega^{\bullet,\Gamma,pr}_{/S},F_{DR})\to p_{a*}p_a^*(\Omega^{\bullet,\Gamma,pr}_{/S},F_{DR})
\end{eqnarray*}
is an homotopy equivalence whose inverse is
\begin{eqnarray*} 
\ad(i_0^*,i_{0*})(p_{a*}p_a^*(\Omega^{\bullet,\Gamma,pr}_{/S},F_{DR})):
p_{a*}p_a^*(\Omega^{\bullet,\Gamma,pr}_{/S},F_{DR})\to (\Omega^{\bullet,\Gamma,pr}_{/S},F_{DR}).
\end{eqnarray*}

\noindent (ii2):Let us shows that $\Omega^{\bullet,\Gamma,pr}_{/S}\in C_{\mathcal D_Sfil}(\Var(\mathbb C)^{2,smpr}/S)$
admits transferts. Let $\alpha\in\Cor(\Var(\mathbb C)^{2,smpr}/S)((Y_1\times S,Z_1)/S,(Y_2\times S,Z_2)/S)$ irreducible.
Denote by $i:\alpha\hookrightarrow Y_1\times Y_2\times S$ the closed embedding, 
and $p_1:Y_1\times Y_2\times S\to Y_1\times S$, $p_2:Y_1\times Y_2\times S\to Y_2\times S$ the projections.
The morphism $p_1\circ i:\alpha\to Y_1\times S$ is then finite surjective and 
$(Z_1\times Y_2)\cap\alpha\subset Y_1\times Z_2$ (i.e. $p_2(p_1^{-1}(Z_2)\cap\alpha)\subset Z_2$).
Then, the transfert map is given by
\begin{eqnarray*}
\Omega^{\bullet,\Gamma,pr}_{/S}(\alpha):
((\Omega^{\bullet}_{Y_2\times S/S},F_b)\otimes_{O_{Y_2\times S}}
\Gamma_{Z_2}^{\vee,Hdg}(O_{Y_2\times S},F_b))(Y_2\times S) \\
\xrightarrow{i_{-}}
p_2^*((\Omega^{\bullet}_{Y_2\times S/S},F_b)\otimes_{O_{Y_2\times S}}
\Gamma_{Z_2}^{\vee,Hdg}(O_{Y_2\times S},F_b))(Y_1\times Y_2\times S) \\
\xrightarrow{\Omega_{(Y_1\times Y_2\times S/Y_2\times S)/(S/S)}(-)(-)}
((\Omega^{\bullet}_{Y_1\times Y_2\times S/S},F_b)\otimes_{O_{Y_1\times Y_2\times S}}
\Gamma_{Y_1\times Z_2}^{\vee,Hdg}(O_{Y_1\times Y_2\times S},F_b))(Y_1\times Y_2\times S) \\
\xrightarrow{DR(-)(T((Z_1\times Y_2)\cap\alpha/Y_1\times Z_2,\gamma^{\vee,Hdg})(-)(-)} \\
((\Omega^{\bullet}_{Y_1\times Y_2\times S/S},F_b)\otimes_{O_{Y_1\times Y_2\times S}}
\Gamma_{(Z_1\times Y_2)\cap\alpha}^{\vee,Hdg}(O_{Y_1\times Y_2\times S},F_b))(Y_1\times Y_2\times S) \\
\xrightarrow{i_{-}}
i^*((\Omega^{\bullet}_{Y_1\times Y_2\times S/S},F_b)\otimes_{O_{Y_1\times Y_2\times S}}
\Gamma_{(Z_1\times Y_2)\cap\alpha}^{\vee,Hdg}(O_{Y_1\times Y_2\times S},F_b))(\alpha) \\
\xrightarrow{\Omega_{(\alpha/Y_1\times Y_2\times S)/(S/S)}(-)(-)}
((\Omega^{\bullet}_{\alpha/S},F_b)\otimes_{O_{\alpha}}
i^{*mod}\Gamma_{(Z_1\times Y_2)\cap\alpha}^{\vee,Hdg}(O_{Y_1\times Y_2\times S},F_b))(\alpha) \\
\xrightarrow{\Omega_{(\alpha/Y_1\times S)(S/S)}(-)(-)^{tr}}
((\Omega^{\bullet}_{Y_1\times S/S},F_b)\otimes_{O_{Y_1\times S}}
\Gamma_{Z_1}^{\vee,Hdg}(O_{Y_1\times S},F_b))(Y_1\times S).
\end{eqnarray*}

\noindent (ii): Two morphisms of complexes of representable presheaves in $C(\Var(\mathbb C)^{smpr}/(\tilde S_I)^{op})$
\begin{eqnarray*}
m=(m^*_{\alpha}),m'=(m^{'*}_{\alpha}):
\oplus_{\alpha\in\Lambda_*}(\mathbb Z((Y^{1,*}_{\alpha}\times\tilde S_I,Z^{1,*}_{\alpha})/\tilde S_I),s^{1,*}_{IJ})
\to\oplus_{\alpha\in\Lambda_*}(\mathbb Z((Y^{2,*}_{\alpha}\times\tilde S_I,Z^{2,*}_{\alpha})/\tilde S_I),s^{2,*}_{IJ})
\end{eqnarray*}
are said to induce a twisted homotopy if there exist morphisms in $C_{D(\tilde S_I)fil}(\mathbb N\times\mathbb Z)$
\begin{eqnarray*}
h=(h^{p,*}_{\alpha}):
(\Omega^{p}_{Y^{2,*}_{\alpha}\times\tilde S_I/\tilde S_I}\otimes_{O_{Y^{2,*}_{\alpha}\times\tilde S_I}}
\Gamma_{Z^{2,*}_{\alpha}}^{\vee,Hdg}(O_{Y^{2,*}_{\alpha}\times\tilde S_I},F_b)(Y^{2,*}_{\alpha}\times\tilde S_I),u_{IJ}(-)^p) \\
\to(\Omega^{p-1}_{Y^{1,*}_{\alpha}\times\tilde S_I/\tilde S_I}\otimes_{O_{Y^{1,*}_{\alpha}\times\tilde S_I}}
\Gamma_{Z^{1,*}_{\alpha}}^{\vee,Hdg}(O_{Y^{1,*}_{\alpha}\times\tilde S_I},F_b)(Y^{1,*}_{\alpha}\times\tilde S_I),u_{IJ}(-)^{p-1})
\end{eqnarray*}
where $p\in\mathbb N$, such that 
\begin{eqnarray*}
\Omega^{\Gamma,pr}(m)-\Omega^{\Gamma,pr}(m')=d_1\circ h+h\circ d_2.
\end{eqnarray*}
As for homotopy, twisted homotopy satisfy the 2 of 3 property for morphism of canonical triangles.
Moreover, by definition, if two morphisms of complexes of representable presheaves in 
$C(\Var(\mathbb C)^{smpr}/(\tilde S_I)^{op})$ $m,m':Q_{1I},K^1_{IJ})\to (Q_{2I},K^2_{IJ})$
induce a twisted homotopy, then 
\begin{eqnarray*} 
e((\tilde S_I))_*\mathcal Hom^{\bullet}((h,m,m'),
(E_{zar}(\Omega^{\bullet,\Gamma,pr}_{/\tilde S_I},F_{DR})[-d_{\tilde S_I}],T_{IJ})): \\  
e'((\tilde S_I))_*\mathcal Hom^{\bullet}((Q_{2I},K^1_{IJ}),
(E_{zar}(\Omega^{\bullet,\Gamma,pr}_{/(\tilde S_I)},F_{DR})[-d_{\tilde S_I}],T_{IJ}))[1] \\
\to e'((\tilde S_I))_*\mathcal Hom^{\bullet}((Q_{1I},K^1_{IJ}),
(E_{zar}(\Omega^{\bullet,\Gamma,pr}_{/(\tilde S_I)},F_{DR})[-d_{\tilde S_I}],T_{IJ}))
\end{eqnarray*}
is a $2$-filtered homotopy.
Now, by proposition \ref{ca1Var12IJ}, there exists 
\begin{eqnarray*}
\left\{((Y_{1I,\alpha}\times\tilde S_I,Z_{1I,\alpha})/\tilde S_I,s^{1,\alpha}_{IJ}),\alpha\in\Lambda_1\right\},\ldots,
\left\{((Y_{sI,\alpha}\times\tilde S_I,Z_{sI,\alpha})/\tilde S_I,s^{s,\alpha}_{IJ}),\alpha\in\Lambda_s\right\} \\
\subset\Var(\mathbb C)^{2,(sm)pr}/(\tilde S_I)^{op}
\end{eqnarray*}
such that we have in $\Ho_{et}(C(\Var(\mathbb C)^{2,(sm)pr}/(\tilde S_I)^{op})$
\begin{eqnarray*}
\Cone(m)\xrightarrow{\sim} \\
\Cone(\oplus_{\alpha\in\Lambda_1}
\Cone((\mathbb Z((Y_{1I,\alpha}\times\mathbb A^1\times\tilde S_I,
Z_{1I,\alpha}\times\mathbb A^1)/\tilde S_I),s^{1,\alpha}_{IJ}\times I)\to
(\mathbb Z((Y_{1I,\alpha}\times S,Z_{1I,\alpha})/\tilde S_I),s^{1,\alpha}_{IJ})) \\
\to\cdots\to\oplus_{\alpha\in\Lambda_s}
\Cone((\mathbb Z((Y_{sI,\alpha}\times\mathbb A^1\times\tilde S_I,
Z_{sI,\alpha}\times\mathbb A^1)/\tilde S_I),s^{s,\alpha}_{IJ}\times I)\to
(\mathbb Z((Y_{sI,\alpha}\times S,Z_{sI,\alpha})/\tilde S_I),s^{s,\alpha}_{IJ})))
\end{eqnarray*}
This gives in $D_{fil}((\tilde S_I)):=\Ho_{zar,fil}((\tilde S_I))$,
\begin{eqnarray*}
\Cone(M)\xrightarrow{\sim} \\
\Cone(\oplus_{\alpha\in\Lambda_1}
\Cone(e'(\tilde S_I)_*\mathcal Hom^{\bullet}((\mathbb Z((Y_{1I,\alpha}\times S,Z_{1,\alpha})/\tilde S_I),s^1_{IJ}),
(E_{zar}(\Omega^{\bullet,\Gamma,pr}_{/\tilde S_I},F_{DR})[-d_{\tilde S_I}],T_{IJ}))\to \\
e'(\tilde S_I)_*\mathcal Hom^{\bullet}((\mathbb Z((Y_{1I,\alpha}\times\mathbb A^1\times\tilde S_I,
Z_{1I,\alpha}\times\mathbb A^1)/\tilde S_I),s^1_{IJ}\times I),
(E_{zar}(\Omega^{\bullet,\Gamma,pr}_{/\tilde S_I},F_{DR})[-d_{\tilde S_I}],T_{IJ})) \\
\to\cdots\to\oplus_{\alpha\in\Lambda_s}
\Cone(e'(\tilde S_I)_*\mathcal Hom^{\bullet}((\mathbb Z((Y_{sI,\alpha}\times\tilde S_I,Z_{sI,\alpha})/\tilde S_I),s^s_{IJ}),
(E_{zar}(\Omega^{\bullet,\Gamma,pr}_{/\tilde S_I},F_{DR})[-d_{\tilde S_I}],T_{IJ}))\to \\
e'(\tilde S_I)_*\mathcal Hom^{\bullet}((\mathbb Z((Y_{sI,\alpha}\times\mathbb A^1\times\tilde S_I,
Z_{sI,\alpha}\times\mathbb A^1)/\tilde S_I),s^s_{IJ}\times I),
(E_{zar}(\Omega^{\bullet,\Gamma,pr}_{/\tilde S_I},F_{DR})[-d_{\tilde S_I}],T_{IJ})))
\end{eqnarray*}
Then by (i1),(i2) and theorem \ref{DDADM12fil}, for all $1\leq i\leq s$ and all $\alpha\in\Lambda_i$
\begin{eqnarray*}
\Cone(e'(\tilde S_I)_*\mathcal Hom^{\bullet}((\mathbb Z((Y_{iI,\alpha}\times\tilde S_I,Z_{iI,\alpha})/\tilde S_I),s^i_{IJ}),
(E_{zar}(\Omega^{\bullet,\Gamma,pr}_{/\tilde S_I},F_{DR})[-d_{\tilde S_I}],T_{IJ}))\to \\
e'(\tilde S_I)_*\mathcal Hom^{\bullet}((\mathbb Z((Y_{iI,\alpha}\times\mathbb A^1\times\tilde S_I,
Z_{iI,\alpha}\times\mathbb A^1)/\tilde S_I),s^1_{IJ}\times I),
(E_{zar}(\Omega^{\bullet,\Gamma,pr}_{/\tilde S_I},F_{DR})[-d_{\tilde S_I}],T_{IJ})))\to 0
\end{eqnarray*}
are twisted homotopy equivalence. 
Hence, we have in $C_{\mathcal Dfil}((\tilde S_I))$
\begin{eqnarray*} 
M=m_1\circ\cdots\circ m_l:  
e'((\tilde S_I))_*\mathcal Hom^{\bullet}((Q_{2I},K^1_{IJ}),
(E_{zar}(\Omega^{\bullet,\Gamma,pr}_{/(\tilde S_I)},F_{DR})[-d_{\tilde S_I}],T_{IJ})) \\
\to e'((\tilde S_I))_*\mathcal Hom^{\bullet}((Q_{1I},K^1_{IJ}),
(E_{zar}(\Omega^{\bullet,\Gamma,pr}_{/(\tilde S_I)},F_{DR})[-d_{\tilde S_I}],T_{IJ}))
\end{eqnarray*}
with $m_i$ either a filtered quasi-isomorphism or a twisted homotopy equivalence.
This proves (ii).
\end{proof}

We deduce the following :

\begin{prop}\label{projwach}
Let $S\in\Var(\mathbb C)$.Let $S=\cup_{i=1}^l S_i$ an open cover such that there exist closed embeddings
$i_i:S_i\hookrightarrow\tilde S_i$ with $\tilde S_i\in\SmVar(\mathbb C)$. 
\begin{itemize}
\item[(i)]Let $m=(m_I):(Q_{1I},K^1_{IJ})\to (Q_{2I},K^2_{IJ})$ be an etale local equivalence local 
with $(Q_{1I},K^1_{IJ}),(Q_{2I},K^2_{IJ})\in C(\Var(\mathbb C)^{sm}/(\tilde S_I))$
complexes of projective presheaves. Then,
\begin{eqnarray*} 
(e'(\tilde S_I)_*\mathcal Hom^{\bullet}(\hat R_S^{CH}(m_I),
E_{zar}(\Omega^{\bullet,\Gamma,pr}_{/\tilde S_I},F_{DR}))[-d_{\tilde S_I}]): \\ 
e'(\tilde S_I)_*\mathcal Hom^{\bullet}((\hat R^{CH}(\rho_S^*Q_{1I}),\hat R^{CH}(K^1_{IJ})),
(E_{zar}(\Omega^{\bullet,\Gamma,pr}_{/\tilde S_I},F_{DR})[-d_{\tilde S_I}],T_{IJ})) \\
\to e'(\tilde S_I)_*\mathcal Hom^{\bullet}((\hat R^{CH}(\rho_S^*Q_{2I}),\hat R^{CH}(K^2_{IJ})),
(E_{zar}(\Omega^{\bullet,\Gamma,pr}_{/\tilde S_I},F_{DR})[-d_{\tilde S_I}],T_{IJ}))
\end{eqnarray*}
is a filtered quasi-isomorphism. It is thus an isomorphism in $D_{\mathcal Dfil}((\tilde S_I))$.
\item[(ii)]Let $m=(m_I):(Q_{1I},K^1_{IJ})\to (Q_{2I},K^2_{IJ})$ be an equivalence $(\mathbb A^1,et)$ local equivalence local 
with $(Q_{1I},K^1_{IJ}),(Q_{2I},K^2_{IJ})\in C(\Var(\mathbb C)^{sm}/(\tilde S_I))$
complexes of projective presheaves. Then,
\begin{eqnarray*} 
(e'(\tilde S_I)_*\mathcal Hom^{\bullet}(\hat R_S^{CH}(m_I),
E_{zar}(\Omega^{\bullet,\Gamma,pr}_{/\tilde S_I},F_{DR}))[-d_{\tilde S_I}]): \\ 
e'(\tilde S_I)_*\mathcal Hom^{\bullet}((\hat R^{CH}(\rho_S^*Q_{1I}),\hat R^{CH}(K^1_{IJ})),
(E_{zar}(\Omega^{\bullet,\Gamma,pr}_{/\tilde S_I},F_{DR})[-d_{\tilde S_I}],T_{IJ})) \\
\to e'(\tilde S_I)_*\mathcal Hom^{\bullet}((\hat R^{CH}(\rho_S^*Q_{2I}),\hat R^{CH}(K^2_{IJ})),
(E_{zar}(\Omega^{\bullet,\Gamma,pr}_{/\tilde S_I},F_{DR})[-d_{\tilde S_I}],T_{IJ}))
\end{eqnarray*}
is a filtered quasi-isomorphism. It is thus an isomorphism in $D_{\mathcal Dfil}((\tilde S_I))$.
\end{itemize}
\end{prop}

\begin{proof}
Follows from proposition \ref{aetfib} and the fact that the differentials of the complexes are stricts for the
filtration by theorem \ref{Sa1}.
\end{proof}

\begin{defi}\label{DRalgdefsing}
\begin{itemize}
\item[(i)] Let $S\in\SmVar(\mathbb C)$.
We define using definition \ref{DRalgdefFunct}(i) and proposition \ref{projwach}(ii)
the filtered algebraic De Rahm realization functor defined as
\begin{eqnarray*}
\mathcal F_S^{FDR}:\DA_c(S)\to D_{\mathcal Dfil}(S), \; M\mapsto \\ 
\mathcal F_S^{FDR}(M):=e(S)_*\Gr^{12}_{S*}\mathcal Hom^{\bullet}(\hat R^{CH}(\rho_S^*L(F)),
E_{zar}(\Omega^{\bullet,\Gamma,pr}_{/S},F_{DR}))[-d_S] 
\end{eqnarray*}
where $F\in C(\Var(\mathbb C)^{sm}/S)$ is such that $M=D(\mathbb A^1,et)(F)$.
\item[(i)'] For the Corti-Hanamura weight structure $W$ on $\DA_c(S)^-$, 
we define using definition \ref{DRalgdefFunct}(i) and proposition \ref{projwach}(ii)
\begin{eqnarray*}
\mathcal F_S^{FDR}:\DA_c^-(S)\to D_{\mathcal D(1,0)fil,\infty}^-(S), \; M\mapsto \\
\mathcal F_S^{FDR}((M,W)):= e(S)_*\Gr^{12}_{S*}\mathcal Hom^{\bullet}(\hat R^{CH}(\rho_S^*L(F,W)),
E_{zar}(\Omega^{\bullet,\Gamma,pr}_{/S},F_{DR}))[-d_S] 
\end{eqnarray*}
where $(F,W)\in C_{fil}(\Var(\mathbb C)^{sm}/S)$ 
is such that $M=D(\mathbb A^1,et)((F,W))$ using corollary \ref{weightst2Cor}. 
Note that the filtration induced by $W$ is a filtration by sub $D_S$ module,
which is a stronger property then Griffitz transversality.
Of course, the filtration induced by $F$ satisfy only Griffitz transversality in general.
\item[(ii)]Let $S\in\Var(\mathbb C)$ and $S=\cup_{i=1}^l S_i$ an open cover such that there exist closed embeddings
$i_i:S_i\hookrightarrow\tilde S_i$  with $\tilde S_i\in\SmVar(\mathbb C)$. 
For $I\subset\left[1,\cdots l\right]$, denote by $S_I=\cap_{i\in I} S_i$ and $j_I:S_I\hookrightarrow S$ the open embedding.
We then have closed embeddings $i_I:S_I\hookrightarrow\tilde S_I:=\Pi_{i\in I}\tilde S_i$.
We define, using definition \ref{DRalgdefFunct}(ii) and proposition \ref{projwach}(ii),
the filtered algebraic De Rahm realization functor defined as
\begin{eqnarray*}
\mathcal F_S^{FDR}:\DA_c(S)\to D_{\mathcal Dfil}(S/(\tilde S_I)), \; M\mapsto \\
\mathcal F_S^{FDR}(M):= 
(e'(\tilde S_I)_*\mathcal Hom^{\bullet}(\hat R^{CH}(\rho_{\tilde S_I}^*L(i_{I*}j_I^*F)),
E_{zar}(\Omega^{\bullet,\Gamma,pr}_{/\tilde S_I},F_{DR}))[-d_{\tilde S_I}],u^q_{IJ}(F))
\end{eqnarray*}
where $F\in C(\Var(\mathbb C)^{sm}/S)$ is such that $M=D(\mathbb A^1,et)(F)$, 
see definition \ref{DRalgdefFunct}.
\item[(ii)'] For the Corti-Hanamura weight structure $W$ on $\DA_c^-(S)$, 
using definition \ref{DRalgdefFunct}(ii) and proposition \ref{projwach}(ii), 
\begin{eqnarray*}
\mathcal F_S^{FDR}:\DA_c(S)\to D_{\mathcal D(1,0)fil}(S/(\tilde S_I)), \; 
M\mapsto\mathcal F_S^{FDR}((M,W)):= \\ 
(e'(\tilde S_I)_*\mathcal Hom^{\bullet}(\hat R^{CH}(\rho_{\tilde S_I}^*L(i_{I*}j_I^*(F,W))),
E_{zar}(\Omega^{\bullet,\Gamma,pr}_{/\tilde S_I},F_{DR}))[-d_{\tilde S_I}],u^q_{IJ}(F,W))
\end{eqnarray*}
where $(F,W)\in C_{fil}(\Var(\mathbb C)^{sm}/S)$ 
is such that $(M,W)=D(\mathbb A^1,et)(F,W)$ using corollary \ref{weightst2Cor}.
Note that the filtration induced by $W$ is a filtration by sub $D_{\tilde S_I}$-modules,
which is a stronger property then Griffitz transversality.
Of course, the filtration induced by $F$ satisfy only Griffitz transversality in general.
\end{itemize}
\end{defi}

\begin{prop}\label{FDRwelldef}
For $S\in\Var(\mathbb C)$ and $S=\cup_{i=1}^l S_i$ an open cover such that there exist closed embeddings
$i_i:S_i\hookrightarrow\tilde S_i$ with $\tilde S_i\in\SmVar(\mathbb C)$, the functor $\mathcal F_S^{FDR}$ is well defined. 
\end{prop}

\begin{proof}
Let $S\in\Var(\mathbb C)$ and $S=\cup_{i=1}^l S_i$ an open cover such that there exist closed embeddings
$i_i:S_i\hookrightarrow\tilde S_i$ with $\tilde S_i\in\SmVar(\mathbb C)$.
Denote, for $I\subset\left[1,\cdots, l\right]$, $S_I=\cap_{i\in I} S_i$ and $j_I:S_I\hookrightarrow S$ the open embedding.
We then have closed embeddings $i_I:S_I\hookrightarrow\tilde S_I:=\Pi_{i\in I}\tilde S_i$.
Let $M\in\DA(S)$. Let $F,F'\in C(\Var(\mathbb C)^{sm}/S)$ such that $M=D(\mathbb A_1,et)(F)=D(\mathbb A_1,et)(F')$.
Then there exist by definition a sequence of morphisms in $C(\Var(\mathbb C)^{sm}/S)$ :
\begin{equation*}
F=F_1\xrightarrow{s_1} F_2\xleftarrow{s_2} F_3\xrightarrow{s_3}F_4\to\cdots\xrightarrow{s_l} F'=F_s  
\end{equation*}
where, for $1\leq k\leq s$, and $s_k$ are $(\mathbb A^1,et)$ local equivalence.
But if $s:F_1\to F_2$ is an equivalence $(\mathbb A^1,et)$ local, 
\begin{equation*}
L(i_{I*}j_I^*s):L(i_{I*}j_I^*F_1)\to L(i_{I*}j_I^*F_2) 
\end{equation*}
is an equivalence $(\mathbb A^1,et)$ local, hence
\begin{eqnarray*}  
\mathcal Hom(\hat R^{CH}_{\tilde S_I}(L(i_{I*}j_I^*s)),E_{zar}(\Omega_{/\tilde S_I}^{\bullet,\Gamma,pr},F_{DR})): \\ 
(e'(\tilde S_I)_*\mathcal Hom(\hat R^{CH}(\rho_{\tilde S_I}^*L(i_{I*}j_I^*F_1)),
E_{zar}(\Omega_{/\tilde S_I}^{\bullet,\Gamma,pr},F_{DR})),u^q_{IJ}(F_1)) \\ 
\to(e'(\tilde S_I)_*\mathcal Hom(\hat R^{CH}(\rho_{\tilde S_I}^*L(i_{I*}j_I^*F_2)),
E_{zar}(\Omega_{/\tilde S_I}^{\bullet,\Gamma,pr},F_{DR})),u^q_{IJ}(F_2))
\end{eqnarray*}
is a filtered quasi-isomorphism by proposition \ref{projwach}.
\end{proof}

Let $f:X\to S$ a morphism with $S,X\in\Var(\mathbb C)$. Assume there exists a factorization 
\begin{equation*}
f:X\xrightarrow{l}Y\times S\xrightarrow{p_S} S
\end{equation*}
of $f$, with $l$ a closed embedding, $Y\in\SmVar(\mathbb C)$ and $p_S$ the projection.
Let $\bar Y\in\PSmVar(\mathbb C)$ a smooth compactification of $Y$ with $\bar Y\backslash Y=D$ a normal crossing divisor,
denote $k:D\hookrightarrow \bar Y$ the closed embedding and $n:Y\hookrightarrow\bar Y$ the open embedding.
Denote $\bar X\subset\bar Y\times S$ the closure of $X\subset\bar Y\times S$.
We have then the following commutative diagram in $\Var(\mathbb C)$
\begin{equation*}
\xymatrix{X\ar[r]^l\ar[d] & Y\times S\ar[rd]^{p_S}\ar[d]^{(n\times I)} & \, \\
\bar X\ar[r]^l & \bar Y\times S\ar[r]^{\bar p_S} & S \\
Z:=\bar X\backslash X\ar[ru]_{l_Z}\ar[u]\ar[r] & D\times S\ar[ru]\ar[u]_{(k\times I)} & \, }.
\end{equation*}
Let $S=\cup_{i=1}^l S_i$ an open cover such that there exist closed embeddings
$i_i:S_i\hookrightarrow\tilde S_i$ with $\tilde S_i\in\SmVar(\mathbb C)$. 
We have $X=\cup_{i=1}^lX_i$ with $X_i:=f^{-1}(S_i)$. 
Denote, for $I\subset\left[1,\cdots l\right]$, $S_I=\cap_{i\in I} S_i$ and $X_I=\cap_{i\in I}X_i$.
For $I\subset\left[1,\cdots l\right]$, denote by $\tilde S_I=\Pi_{i\in I}\tilde S_i$.
We then have, for $I\subset\left[1,\cdots l\right]$, closed embeddings $i_I:S_I\hookrightarrow\tilde S_I$
and for $I\subset J$, the following commutative diagrams which are cartesian 
\begin{equation*}
\xymatrix{
f_I=f_{|X_I}:X_I\ar[r]^{l_I}\ar[rd] & Y\times S_I\ar[r]^{p_{S_I}}\ar[d]^{i'_I} & S_I\ar[d]^{i_I} \\
\, & Y\times\tilde S_I\ar[r]^{p_{\tilde S_I}} & \tilde S_I} \;, \;
\xymatrix{Y\times\tilde S_J\ar[r]^{p_{\tilde S_J}}\ar[d]_{p'_{IJ}} & \tilde S_J\ar[d]^{p_{IJ}} \\
Y\times\tilde S_I\ar[r]^{p_{\tilde S_I}} & \tilde S_I},
\end{equation*}
with $l_I:l_{|X_I}$, $i'_I=I\times i_I$, $p_{S_I}$ and $p_{\tilde S_I}$ are the projections and $p'_{IJ}=I\times p_{IJ}$, 
and we recall that we denote by 
$j_I:\tilde S_I\backslash S_I\hookrightarrow\tilde S_I$ and $j'_I:Y\times\tilde S_I\backslash X_I\hookrightarrow Y\times S_I$
the open complementary embeddings.
We then have the commutative diagrams
\begin{equation*}
D_{IJ}=\xymatrix{S_J\ar[r]^{i_J}\ar[d]^{j_{IJ}} & \tilde S_J\ar[d]^{p_{IJ}} \\
S_I\ar[r]^{i_I} & \tilde S_I} \; , \;
D'_{IJ}=\xymatrix{X_J\ar[rr]^{i'_J\circ l_J}\ar[d]^{j'_{IJ}} & \, & Y\times\tilde S_J\ar[d]^{p'_{IJ}} \\
X_I\ar[rr]^{i'_I\circ l_I} & \, & Y\times\tilde S_I}.
\end{equation*}
and the factorization of $D'_{IJ}$ by the fiber product:
\begin{equation}
D'_{IJ}=
\xymatrix{ X_J\ar[r]^{i'_I\circ l_I}\ar[d]^{j'_{IJ}} & Y\times\tilde S_J\ar[d]^{p'_{IJ}} \\ 
X_I\ar[r]^{i'_I\circ l_I} & Y\times\tilde S_I}, \; \;
D'_{IJ}=
\xymatrix{ X_J\ar[rrr]^{i'_I\circ l_I}\ar[dd]^{j'_{IJ}}\ar[rrd]^{\hat l_J} & \, & \, & 
Y\times\tilde S_J\ar[dd]^{p'_{IJ}} \\ 
\,  & \, & X_I\times_{Y\times\tilde S_I}Y\times\tilde S_J=X_I\times\tilde S_{J\backslash I}
\ar[lld]^{\hat p_{IJ}}\ar[ru]^{\hat {il}_I}  & \, \\ 
X_I\ar[rrr]^{i'_I\circ l_I} & \, & \, & Y\times\tilde S_I}
\end{equation}
where $j'_{IJ}:X_J\hookrightarrow X_I$ is the open embedding.
Consider
\begin{equation*}
F(X/S):=p_{S,\sharp}\Gamma^{\vee}_X\mathbb Z(Y\times S/Y\times S)\in C(\Var(\mathbb C)^{sm}/S)
\end{equation*}
so that $D(\mathbb A^1,et)(F(X/S))=M(X/S)$ since $Y$ is smooth. Then, by definition,
\begin{eqnarray*}
\mathcal F_S^{FDR}(M(X/S)):=
(e'(\tilde S_I)_*\mathcal Hom(\hat R^{CH}(\rho_{\tilde S_I}^*L(i_{I*}j_I^*F(X/S))), \\
E_{zar}(\Omega^{\bullet,\Gamma,pr}_{/\tilde S_I},F_{DR}))[-d_{\tilde S_I}],u^q_{IJ}(F(X/S)))
\end{eqnarray*}
On the other hand, let  
\begin{equation*}
Q(X_I/\tilde S_I):=p_{\tilde S_I,\sharp}\Gamma^{\vee}_{X_I}\mathbb Z(Y\times\tilde S_I/Y\times\tilde S_I).
\in C(\Var(\mathbb C)^{sm}/\tilde S_I), 
\end{equation*}
We have then for $I\subset[1,l]$ the map (\ref{HYS}) in $C(\Var(\mathbb C)^{sm}/\tilde S_J)$ :
\begin{eqnarray*}
N_I(X/S):Q(X_I/\tilde S_I)=p_{\tilde S_I\sharp}\Gamma^{\vee}_{X_I}\mathbb Z(Y\times\tilde S_I/Y\times\tilde S_I)
\xrightarrow{\ad(i^{'*}_I,i'_{I*})(-)} \\
p_{\tilde S_I\sharp}i'_{I*}i_I^{'*}\Gamma^{\vee}_{X_I}\mathbb Z(Y\times S_I/Y\times S_I) 
\xrightarrow{T(i'_I,\gamma^{\vee})(-)}
p_{\tilde S_I\sharp}i'_{I*}\Gamma^{\vee}_{X_I}\mathbb Z(Y\times S_I/Y\times S_I) \\
\xrightarrow{\hat{T}_{\sharp}(p_{S_I},i_I)(-)}
i_{I*}p_{S_I\sharp}\Gamma^{\vee}_{X_I}\mathbb Z(Y\times S_I/Y\times S_I)=i_{I*}j_I^*F(X/S).
\end{eqnarray*}
We then have the commutative diagram in $C(\Var(\mathbb C)^{sm}/\tilde S_J)$
\begin{equation}\label{HYSIJq}
\xymatrix{p_{IJ}^*Q(X_I/\tilde S_I)\ar[rr]^{p_{IJ}^*N_I(X/S)}\ar[d]_{H_{IJ}} & \, &
p_{IJ}^*Li_{I*}j_I^*F(X/S)\ar[d]^{T^q(D_{IJ})(j_I^*F(X/S))} \\
Q(X_J/\tilde S_J)\ar[rr]^{N_J(X/S)} & \, & Li_{J*}j_J^*F(X/S)}
\end{equation}
with 
\begin{eqnarray*}
H_{IJ}:p_{IJ}^*p_{\tilde S_I\sharp}\Gamma^{\vee}_{X_I}\mathbb Z(Y\times\tilde S_I/Y\times\tilde S_I) 
\xrightarrow{T_{\sharp}(p_{IJ},p_{\tilde S_I})(-)^{-1}}
p_{\tilde S_J\sharp}p_{IJ}^{'*}\Gamma^{\vee}_{X_I}\mathbb Z(Y\times\tilde S_I/Y\times\tilde S_I) \\
\xrightarrow{p_{\tilde S_J\sharp}T(p_{IJ},\gamma^{\vee})(-)}
p_{\tilde S_J\sharp}\Gamma^{\vee}_{X_I\times\tilde S_{J\backslash I}}\mathbb Z(Y\times\tilde S_J/Y\times\tilde S_J) 
\xrightarrow{p_{\tilde S_J\sharp}T(X_J/X_I\times\tilde S_{J\backslash I},\gamma^{\vee})(-)}
p_{\tilde S_J\sharp}\Gamma^{\vee}_{X_J}\mathbb Z(Y\times\tilde S_J/Y\times\tilde S_J). 
\end{eqnarray*}
The diagram \ref{HYSIJq} say that the maps $N_I(X/S)$ induces a map in $C(\Var(\mathbb C)^{sm}/(S/\tilde S_I))$
\begin{eqnarray*}
(N_I(X/S)):(Q(X_I/\tilde S_I),I(p_{IJ}^*,p_{IJ*})(-,-)(H_{IJ})) \\
\to(Li_{I*}j_I^*F(X/S),I(p_{IJ}^*,p_{IJ*})(-,-)(T^q(D_{IJ})(j_I^*F(X/S)))).
\end{eqnarray*} 
Denote $\bar X_I:=\bar X\cap(\bar Y\times S_I)\subset\bar Y\times\tilde S_I$ 
the closure of $X_I\subset\bar Y\times\tilde S_I$,
and $Z_I:=Z\cap(\bar Y\times S_I)=\bar X_I\backslash X_I$.
Consider for $I\subset\left[1,\cdots l\right]$ and $I\subset J$ the following commutative diagrams in $\Var(\mathbb C)$
\begin{equation*}
\xymatrix{X_I\ar[r]^{l_I}\ar[d] & Y\times \tilde S_I\ar[rd]^{p_{\tilde S_I}}\ar[d]^{(n\times I)} & \, \\
\bar X_I\ar[r]^{l_I} & \bar Y\times\tilde S_I\ar[r]^{\bar p_{\tilde S_I}} & \tilde S_I \\
Z_I=\bar X_I\backslash X_I\ar[ru]_{l_{Z_I}}\ar[u]\ar[r] & D\times\tilde S_I\ar[ru]\ar[u]_{(k\times I)} & \, } \, , \,
\xymatrix{\bar Y\times\tilde S_J\ar[r]^{p_{\tilde S_J}}\ar[d]_{p'_{IJ}} & \tilde S_J\ar[d]^{p_{IJ}} \\
\bar Y\times\tilde S_I\ar[r]^{p_{\tilde S_I}} & \tilde S_I}
\end{equation*}
Let $\epsilon_1:((\bar Y\times\tilde S_I)_1,E_1)\to(\bar Y\times\tilde S_I,Z_I)$ 
a strict desingularization of the pair $(\bar Y\times\tilde S_I,Z_I)$, 
$\epsilon_2:((\bar Y\times\tilde S_I)_2,E_2)\to(\bar Y\times\tilde S_I,\bar X_I)$
a strict desingularization of the pair $(\bar Y\times\tilde S_I,\bar X_I)$
and a morphism $\epsilon_{12}:(\bar Y\times\tilde S_I)_2\to(\bar Y\times\tilde S_I)_1$
such that the following diagram commutes (see definition-proposition \ref{RCHdef0}) :
\begin{equation*}
\xymatrix{(\bar Y\times\tilde S_I)_2\ar[r]^{\epsilon_{12}}\ar[d]^{\epsilon_2} & 
(\bar Y\times\tilde S_I)_1\ar[d]^{\epsilon_1} \\
\bar Y\times\tilde S_I\ar[r]^{=} & \bar Y\times\tilde S_I}
\end{equation*}
and we denote by $l'_{Z_I}:E_1\hookrightarrow(\bar Y\times\tilde S_I)_2$
$l'_I:E_2\hookrightarrow(\bar Y\times\tilde S_I)_2$ the closed embeddings.
We have then the canonical map in $C(\Var(\mathbb C)^{2,smpr}/(\tilde S_I)^{op})$
\begin{eqnarray*}
(I_{\delta}((\bar X_I,Z_I)/\tilde S_I)):
(\hat R^{CH}(\rho_{\tilde S_I}^*Q(X_I/\tilde S_I)),\hat R^{CH}(H_{IJ})) \\
:=(\Cone(\Cone(
(\mathbb Z^{tr}((E_{1\bullet}\times\tilde S_I,E_{1\bullet})/(\bar Y\times\tilde S_I)_1\times\tilde S_I,u_{IJ})\to \\
\mathbb Z^{tr}(((\bar Y\times\tilde S_I)_1\times\tilde S_I,(\bar Y\times\tilde S_I)_1)/(\bar Y\times\tilde S_I)_1)\times\tilde S_I)) 
\to\Cone((\mathbb Z^{tr}((E_{2\bullet}\times\tilde S_I,E_{2\bullet})/(\bar Y\times\tilde S_I)_2\times\tilde S_I),u_{IJ})\to \\
\mathbb Z^{tr}(((\bar Y\times\tilde S_I)_2\times\tilde S_I,(\bar Y\times\tilde S_I)_2)/(\bar Y\times\tilde S_I)_2\times\tilde S_I))),
\hat R^{CH}(H_{IJ})) \\
\xrightarrow{(([(\Gamma_{\epsilon_1\circ l'_{Z_I}})],[\Gamma_{\epsilon_1}]),
([\Gamma_{\epsilon_2\circ l'_I}],[\Gamma_{\epsilon_2}])))} \\
(\Cone(\mathbb Z^{tr}((\bar Y\times\tilde S_I,Z_I)/\tilde S_I)\to\mathbb Z^{tr}((\bar Y\times\tilde S_I,\bar X_I)/\tilde S_I))
(-d_Y-d_{\tilde S_I})[-2d_Y-2d_{\tilde S_I}],\mathbb Z^{tr}(I\times p_{IJ}))
\end{eqnarray*}
We denote by $v_{IJ}^q(F(X/S))$ the composite
\begin{eqnarray*}
v_{IJ}^q(F(X/S))[d_{\tilde S_J}]:
e'(\tilde S_I)_*\mathcal Hom(\hat R_{((\bar Y\times\tilde S_I)^*,E^*)/\tilde S_I}(\rho_{\tilde S_I}^*Q(X_I/\tilde S_I)),
E_{zar}(\Omega^{\bullet,\Gamma,pr}_{/\tilde S_I},F_{DR})) \\
\xrightarrow{p_{IJ*}T(p_{IJ},\Omega^{\Gamma,pr}_{\cdot})(-)\circ\ad(p_{IJ}^{*mod},p_{IJ})(-)} \\
e'(\tilde S_J)_*\mathcal Hom(
p_{IJ}^*\hat R_{((\bar Y\times\tilde S_I)^*,E^*)/\tilde S_I}(\rho_{\tilde S_I}^*Q(X_I/\tilde S_I),
E_{zar}(\Omega^{\bullet,\Gamma,pr}_{/\tilde S_I},F_{DR}))) \\
\xrightarrow{\mathcal Hom(T(p_{IJ},R^{CH})(Q(X_I/\tilde S_I))^{-1},
E_{zar}(\Omega^{\bullet,\Gamma,pr}_{/\tilde S_I},F_{DR}))} \\
p_{IJ*}e'(\tilde S_J)_*\mathcal Hom(\hat R_{((\bar Y\times\tilde S_I)^*\times\tilde S_{J\backslash I},
E^*\times\tilde S_{J\backslash I})/\tilde S_J}
(\rho_{\tilde S_J}^*p_{IJ}^*Q(X_I/\tilde S_I)),E_{zar}(\Omega^{\bullet,\Gamma,pr}_{/\tilde S_J},F_{DR})) \\
\xrightarrow{\mathcal Hom(\hat R_{\tilde S_J}^{CH}(H_{IJ}),
E_{zar}(\Omega^{\bullet,\Gamma,pr}_{/\tilde S_J},F_{DR}))} \\
p_{IJ*}e'(\tilde S_J)_*\mathcal Hom(\hat R_{((\bar Y\times\tilde S_J)^*,E^{'*})/\tilde S_J}
(\rho_{\tilde S_J}^*Q(X_J/\tilde S_J)),E_{zar}(\Omega^{\bullet,\Gamma,pr}_{/\tilde S_J},F_{DR})).
\end{eqnarray*} 
On the other hand, we have in $\pi_{\bar X}(C(MHM(\bar X)))\subset C_{\mathcal Dfil}(\bar X/(\bar Y\times\tilde S_I))$
\begin{eqnarray*}
(\Cone(T(Z_I/\bar X_I,\gamma^{\vee,Hdg})(-))):
(\Gamma^{\vee,Hdg}_{\bar X_I}(O_{\bar Y\times\tilde S_I},F_b),x_{IJ}(\bar X/S)) 
\to(\Gamma^{\vee,Hdg}_{Z_I}(O_{\bar Y\times\tilde S_I},F_b),x_{IJ}(Z/S)) \\
\xrightarrow{=}(n\times I)_{!Hdg}(\Gamma^{\vee,Hdg}_{X_I}(O_{Y\times\tilde S_I},F_b),x_{IJ}(X/S))
\end{eqnarray*}
with 
\begin{itemize}
\item for the closed embedding $\bar X\subset\bar Y\times S$
we consider the map in $\pi_{\bar Y\times\tilde S_J}(C(MHM(\bar Y\times\tilde S_J)))$
\begin{eqnarray*}
x_{IJ}(\bar X/S):\Gamma^{\vee,Hdg}_{\bar X_I}(O_{\bar Y\times\tilde S_I},F_b)
\xrightarrow{\ad(p_{IJ}^{'*mod},p'_{IJ*})(-)}
p'_{IJ*}p_{IJ}^{'*mod}\Gamma^{\vee,Hdg}_{\bar X_I}(O_{\bar Y\times\tilde S_I},F_b) \\
\xrightarrow{T(\bar X_J/p_{IJ}^{'-1}(\bar X_I),\gamma^{\vee})(-)\circ T(p'_{IJ},\gamma^{\vee})(-)} 
\Gamma^{\vee,Hdg}_{\bar X_J}(O_{\bar Y\times\tilde S_J},F_b),
\end{eqnarray*} 
\item for the closed embedding $Z\subset\bar Y\times S$
we consider the map in $\pi_{\bar Y\times\tilde S_J}(C(MHM(\bar Y\times\tilde S_J)))$
\begin{eqnarray*}
x_{IJ}(Z/S):\Gamma^{\vee,Hdg}_{Z_I}(O_{\bar Y\times\tilde S_I},F_b)
\xrightarrow{\ad(p_{IJ}^{'*mod},p'_{IJ*})(-)}
p'_{IJ*}p_{IJ}^{'*mod}\Gamma^{\vee,Hdg}_{Z_I}(O_{\bar Y\times\tilde S_I},F_b) \\
\xrightarrow{T(Z_J/p_{IJ}^{'-1}(Z_I),\gamma^{\vee})(-)\circ T(p'_{IJ},\gamma^{\vee})(-)} 
\Gamma^{\vee,Hdg}_{Z_J}(O_{\bar Y\times\tilde S_J},F_b),
\end{eqnarray*}
\item for the closed embedding $X\subset Y\times S$
we consider the map in $\pi_{Y\times\tilde S_J}(C(MHM(Y\times\tilde S_J)))$
\begin{eqnarray*}
x_{IJ}(X/S)(-d_Y)[-2d_Y]:\Gamma^{\vee,Hdg}_{X_I}(O_{Y\times\tilde S_I},F_b)
\xrightarrow{\ad(p_{IJ}^{'*mod},p'_{IJ*})(-)}
p'_{IJ*}p_{IJ}^{'*mod}\Gamma^{\vee,Hdg}_{X_I}(O_{Y\times\tilde S_I},F_b) \\ 
\xrightarrow{T(X_J/X_I\times\tilde S_{J\backslash I},\gamma^{\vee,Hdg})(-)\circ T(p'_{IJ},\gamma^{\vee})(-)} 
\Gamma^{\vee,Hdg}_{X_J}(O_{Y\times\tilde S_J},F_b).
\end{eqnarray*}
\end{itemize}
The maps $x_{IJ}(X/S)$ gives the following maps in $C_{\mathcal Dfil,S_J}(\tilde S_J)$
\begin{eqnarray*}
w_{IJ}(X/S)(-d_Y)[-2d_Y]:
p_{\tilde S_I*}E_{zar}((\Omega^{\bullet}_{\bar Y\times\tilde S_I/\tilde S_I},F_b)\otimes_{O_{\bar Y\times\tilde S_I}}
(n\times I)_!^{Hdg}\Gamma^{\vee,Hdg}_{X_I}(O_{Y\times\tilde S_I},F_b)) \\
\xrightarrow{\ad(p_{IJ}^{*mod},p_{IJ})(-)} 
p_{IJ*}p_{IJ}^{*mod}p_{\tilde S_I*}E_{zar}((\Omega^{\bullet}_{\bar Y\times\tilde S_I/\tilde S_I},F_b)
\otimes_{O_{Y\times\tilde S_I}}(n\times I)_!^{Hdg}\Gamma^{\vee,Hdg}_{X_I}(O_{Y\times\tilde S_I},F_b) \\
\xrightarrow{p_{IJ*}T^O_w(p_{IJ},p_{\tilde S_I})(-)}
p_{\tilde S_J*}E_{zar}((\Omega^{\bullet}_{\bar Y\times\tilde S_J/\tilde S_J},F_b)\otimes_{O_{\bar Y\times\tilde S_J}}
p_{IJ}^{'*mod}(n\times I)_{!Hdg}\Gamma^{\vee,Hdg}_{X_I}(O_{Y\times\tilde S_I},F_b)) \\
\xrightarrow{(p_{\tilde S_J*}E(DR(-)(x_{IJ}(X/S))}
p_{\tilde S_J*}E_{zar}((\Omega^{\bullet}_{\bar Y\times\tilde S_J/\tilde S_J},F_b)\otimes_{O_{\bar Y\times\tilde S_J}}
(n\times I)_{!Hdg}\Gamma^{\vee,Hdg}_{X_J}(O_{\bar Y\times\tilde S_J},F_b)).
\end{eqnarray*}
We have then the following lemma

\begin{lem}\label{keyalgsinglem}
\begin{itemize}
\item[(i)] The map in $C(\Var(\mathbb C)^{sm}/(S/\tilde S_I))$
\begin{equation*}
(N_I(X/S)):(Q(X_I/\tilde S_I),H_{IJ})\to(L(i_{I*}j_I^*F(X/S)),T^q(D_{IJ})(F(X/S))).
\end{equation*}
is an equivalence $(\mathbb A^1,et)$ local.
\item[(ii)] The maps $(N_I(X/S))$ induces a filtered quasi-isomorphism 
in $C_{\mathcal Dfil}(S/(\tilde S_I))$
\begin{eqnarray*}
(\mathcal Hom(\hat R^{CH}_{\tilde S_I}(N_I(X/S)),E_{zar}(\Omega^{\bullet,\Gamma,pr}_{/\tilde S_I},F_b))): \\
(e'(\tilde S_I)_*\mathcal Hom(\hat R^{CH}(\rho_{\tilde S_I}^*L(i_{I*}j_I^*F(X/S))),
E_{zar}(\Omega^{\bullet,\Gamma,pr}_{/\tilde S_I},F_{DR}))[-d_{\tilde S_I}],u^q_{IJ}(F(X/S))) \\
\to (e'(\tilde S_I)_*\mathcal Hom(\hat R^{CH}(\rho_{\tilde S_I}^*Q(X_I/\tilde S_I)),
E_{zar}(\Omega^{\bullet,\Gamma,pr}_{/\tilde S_I},F_{DR}))[-d_{\tilde S_I}],v^q_{IJ}(F(X/S)))
\end{eqnarray*}
\item[(iii)] The map $(I_{\delta}((\bar X_I,Z_I)/\tilde S_I))$  
induces a filtered Zariski local equivalence in $C_{\mathcal Dfil}(S/(\tilde S_I))$
\begin{eqnarray*}
(I((\bar X_I,Z_I)/\tilde S_I)):= \\
\mathcal Hom((I_{\delta}((\bar X_I,Z_I)/\tilde S_I)),-)\circ
(DR(Y\times\tilde S_I/\tilde S_I)(\ad((n\times I)_{!Hdg},(n\times I)^{!Hdg}(O_{Y\times\tilde S_I},F_b))),0): \\
(p_{\tilde S_I*}E_{zar}((\Omega^{\bullet}_{\bar Y\times\tilde S_I/\tilde S_I},F_b)\otimes_{O_{\bar Y\times\tilde S_I}}
(n\times I)_{!Hdg}(\Gamma^{\vee,Hdg}_{X_I}(O_{Y\times\tilde S_I},F_b))(d_Y)[2d_Y+d_{\tilde S_I}],w_{IJ}(X/S)) \\
\to (e'(\tilde S_I)_*\mathcal Hom(\hat R^{CH}(\rho_{\tilde S_I}^*Q(X_I/\tilde S_I)),
E_{zar}(\Omega^{\bullet,\Gamma,pr}_{/\tilde S_I},F_{DR}))[-d_{\tilde S_I}],v^q_{IJ}(F(X/S)))    
\end{eqnarray*}
\end{itemize}
\end{lem}

\begin{proof}
\noindent(i): See lemma \ref{keyalgsinglemGM}(i)

\noindent(ii): These maps induce a morphism in $C_{\mathcal D}(S/(\tilde S_I))$ by construction.
It is a filtered quasi-isomorphism by (i) and proposition \ref{aetfib}.

\noindent(iii): See section 5.2.
\end{proof}

\begin{prop}\label{keyalgsing1}
Let $f:X\to S$ a morphism with $S,X\in\Var(\mathbb C)$. Assume there exist a factorization 
\begin{equation*}
f:X\xrightarrow{l}Y\times S\xrightarrow{p_S} S
\end{equation*}
of $f$ with $Y\in\SmVar(\mathbb C)$, $l$ a closed embedding and $p_S$ the projection.
Let $\bar Y\in\PSmVar(\mathbb C)$ a compactification of $Y$ with $\bar Y\backslash Y=D$ a normal crossing divisor,
denote $k:D\hookrightarrow \bar Y$ the closed embedding and $n:Y\hookrightarrow\bar Y$ the open embedding.
Denote $\bar X\subset\bar Y\times S$ the closure of $X\subset\bar Y\times S$.
We have then the following commutative diagram in $\Var(\mathbb C)$
\begin{equation*}
\xymatrix{X\ar[r]^l\ar[d] & Y\times S\ar[rd]^{p_S}\ar[d]^{(n\times I)} & \, \\
\bar X\ar[r]^l & \bar Y\times S\ar[r]^{\bar p_S} & S \\
Z:=\bar X\backslash X\ar[ru]^{l_Z}\ar[u]\ar[r] & D\times S\ar[ru]\ar[u]_{(k\times I)} & \, }.
\end{equation*}
Let $S=\cup_{i=1}^l S_i$ an open cover such that there exist closed embeddings
$i_i:S_i\hookrightarrow\tilde S_i$ with $\tilde S_i\in\SmVar(\mathbb C)$. 
Then $X=\cup_{i=1}^lX_i$ with $X_i:=f^{-1}(S_i)$.
Denote, for $I\subset\left[1,\cdots l\right]$, $S_I=\cap_{i\in I} S_i$ and $X_I=\cap_{i\in I}X_i$.
Denote $\bar X_I:=\bar X\cap(\bar Y\times S_I)\subset\bar Y\times\tilde S_I$ 
the closure of $X_I\subset\bar Y\times\tilde S_I$, 
and $Z_I:=Z\cap(\bar Y\times S_I)=\bar X_I\backslash X_I\subset\bar Y\times\tilde S_I$.
We have then for $I\subset\left[1,\cdots l\right]$, the following commutative diagram in $\Var(\mathbb C)$
\begin{equation*}
\xymatrix{X_I\ar[r]^{l_I}\ar[d] & Y\times \tilde S_I\ar[rd]^{p_{\tilde S_I}}\ar[d]^{(n\times I)} & \, \\
\bar X_I\ar[r]^{l_I} & \bar Y\times\tilde S_I\ar[r]^{\bar p_{\tilde S_I}} & \tilde S_I \\
Z_I=\bar X_I\backslash X_I\ar[ru]^{l_{Z_I}}\ar[u]\ar[r] & D\times\tilde S_I\ar[ru]\ar[u]_{(k\times I)} & \, }.
\end{equation*}
Let $F(X/S):=p_{S,\sharp}\Gamma_X^{\vee}\mathbb Z(Y\times S/Y\times S)\in C(\Var(\mathbb C)^{sm}/S)$. 
We have then the following isomorphism in $D_{\mathcal Dfil}(S/(\tilde S_I))$ 
\begin{eqnarray*} 
I(X/S):\mathcal F_S^{FDR}(M(X/S))\xrightarrow{:=} \\
(e'(\tilde S_I)_*\mathcal Hom(\hat R^{CH}(\rho_{\tilde S_I}^*L(i_{I*}j_I^*F(X/S))), 
E_{zar}(\Omega^{\bullet,\Gamma,pr}_{/\tilde S_I},F_{DR}))[-d_{\tilde S_I}],u^q_{IJ}(F(X/S))) \\ 
\xrightarrow{(\mathcal Hom(\hat R^{CH}_{\tilde S_I}(N_I(X/S)), 
E_{zar}(\Omega^{\bullet,\Gamma,pr}_{/\tilde S_I},F_{DR})))} \\
(e'(\tilde S_I)_*\mathcal Hom(\hat R^{CH}(\rho_{\tilde S_I}^*Q(X_I/\tilde S_I)), 
E_{zar}(\Omega^{\bullet,\Gamma,pr}_{/\tilde S_I},F_{DR}))[-d_{\tilde S_I}],v^q_{IJ}(F(X/S))) \\
\xrightarrow{I((\bar X_I,Z_I)/\tilde S_I)} \\
(p_{\tilde S_I*}E_{zar}((\Omega^{\bullet}_{\bar Y\times\tilde S_I/\tilde S_I},F_b)\otimes_{O_{\bar Y\times\tilde S_I}}
(n\times I)_{!Hdg}\Gamma^{\vee,Hdg}_{X_I}(O_{Y\times\tilde S_I},F_b))(d_Y+d_{\tilde S_I})[2d_Y+d_{\tilde S_I}],w_{IJ}(X/S)) \\ 
\xrightarrow{=:}\iota_SRf^{Hdg}_!(\Gamma^{\vee,Hdg}_{X_I}(O_{Y\times\tilde S_I},F_b)(d_Y)[2d_Y],x_{IJ}(X/S)).
\xrightarrow{=:}\iota_SRf^{Hdg}_!f^{*mod}_{Hdg}\mathbb Z^{Hdg}_S.
\end{eqnarray*}
\end{prop}

\begin{proof}
Follows from lemma \ref{keyalgsinglem}, proposition \ref{aetfib} and proposition \ref{sharpstarprop}.
\end{proof}

\begin{cor}\label{FDRMHM}
Let $S\in\Var(\mathbb C)$ and $S=\cup_{i=1}^l S_i$ an open cover such that there exist closed embeddings
$i_i:S_i\hookrightarrow\tilde S_i$ with $\tilde S_i\in\SmVar(\mathbb C)$. 
Then, for $F\in C(\Var(\mathbb C)^{sm}/S)$ such that $M=D(\mathbb A^1,et)(F)\in\DA_c(S)$, 
\begin{eqnarray*}
H^i\mathcal F_S^{FDR}(M,W):= 
(a_{zar}H^ie'(\tilde S_I)_*\mathcal Hom^{\bullet}(L\rho_{\tilde S_I*}\mu_{\tilde S_I*}
R^{CH}(\rho_{\tilde S_I}^*L(i_{I*}j_I^*(F,W))), \\
E_{zar}(\Omega^{\bullet,\Gamma,pr}_{/\tilde S_I},F_{DR}))[-d_{\tilde S_I}],H^iu^q_{IJ}(F,W)) 
\in\pi_S(MHM(S))
\end{eqnarray*}
for all $i\in\mathbb Z$, and for all $p\in\mathbb Z$, 
\begin{eqnarray*}
\mathcal F_S^{FDR}(M,W)\in D_{\mathcal D(1,0)fil}(S/(\tilde S_I))
\end{eqnarray*}
is the class of a complex $\mathcal F_S^{FDR}(F,W)\in C_{\mathcal D(1,0)fil}(S/(\tilde S_I))$
such that for all $k\in\mathbb Z$, 
the differentials of $\Gr_k^W\mathcal F_S^{FDR}(F,W)$ are strict for the filtration $F$.
\end{cor}

\begin{proof}
Follows from theorem \ref{Sa1}:
Indeed, for $S$ smooth and $g:U'/S\to U/S$ a morphism 
with $U/S=(U,h),U'/S=(U',h')\in\Var(\mathbb C)^{sm}/S$, $U$,$U'$ connected, hence irreducible by smoothness, the complex
\begin{eqnarray*} 
\hat{\mathcal F_S^{FDR}}(g):=
e'(S)_*\mathcal Hom^{\bullet}(\hat R^{CH}(\mathbb Z(U'/S)),E_{zar}(\Omega^{\bullet,\Gamma,pr}_{/S},F_{DR}))[-d_S] \\
\xrightarrow{\mathcal Hom(\hat R^{CH}(g),-)} 
e'(S)_*\mathcal Hom^{\bullet}(\hat R^{CH}(\mathbb Z(U/S)),E_{zar}(\Omega^{\bullet,\Gamma,pr}_{/S},F_{DR})[-d_S])
\in C_{\mathcal Dfil}(S) 
\end{eqnarray*}
satisfy $H^i\mathcal F_S^{FDR}(g)\in\pi_S(MHM(S))$ and is the class of a complex such that the differentials are strict for $F$.
Let $U\subset\bar X$ a compactification of $U$ and $U'\subset\bar X'$ a compactification of $U'$, 
$S\subset\bar S$ a compactification of $S$ with $\bar X,\bar X'\in\PSmVar(\mathbb C)$,
$i:\bar D=\cup_i \bar D_i:=\bar X\backslash U\hookrightarrow\bar X$,
$i:\bar D'=\cup_i \bar D'_i:=\bar X'\backslash U'\hookrightarrow\bar X'$ normal crossing divisors, 
such that $g:U'/S\to U/S$ extend to $\bar g:\bar X'/\bar S\to\bar X/\bar S$ : see section 2. 
Denote $n:U\hookrightarrow\bar X$, $n':U'\hookrightarrow\bar X'$ and
$n'':U'\hookrightarrow\bar X'\backslash\bar g^{-1}(\bar D)=:U''$ the open embeddings. 
We get the following map in $C_{\mathcal Dfil}(S)$
\begin{eqnarray*} 
\Omega^{\Gamma,pr}(\Gamma_{\bar g}^t):=(p_{S*}E(\Omega^{\bullet}_{\bar X'\times S/S}\otimes_{O_{\bar X'\times S}}
(n'\times I)_{!Hdg}\Gamma_{U'}^{\vee,Hdg}(O_{U'\times S},F_b))(d')[2d'] \\
\xrightarrow{\Omega^{\Gamma,pr}([\Gamma_{\bar g}^t])\circ 
DR(\bar X'\times S/S)(\ad(n''_{!Hdg},n^{''}!_{Hdg})(O_{U''\times S},F_b))} \\
p_{S*}E(\Omega^{\bullet}_{\bar X\times S/S}\otimes_{O_{\bar X\times S}}
(n\times I)_{!Hdg}\Gamma_U^{\vee,Hdg}(O_{U\times S},F_b))(d)[2d])) \\
\xrightarrow{(DR(\bar X'\times S/S)(\ad((n'\times I)_{!Hdg},(n'\times I)^{!Hdg})(O_{\bar X'\times S},F_b)),0)
(DR(\bar X\times S/S)(\ad((n\times I)_{!Hdg},(n\times I)^{!Hdg})(O_{\bar X\times S},F_b)),0)} \\
e'(S)_*\mathcal Hom((\Cone(\mathbb Z(i_{\bullet}\times I):
\mathbb Z^{tr}((\bar D_{\bullet}\times S,D_{\bullet})/S)\to \mathbb Z^{tr}((\bar X\times S,X)/S))(-d)[-2d] \\
\xrightarrow{([\Gamma_{\bar g}\times\square^*]^t,[\Gamma_{\bar g'_{\bullet}}\times\square^*]^t)} \\
\Cone(\mathbb Z(i_{\bullet}\times I):\mathbb Z^{tr}((\bar D'_{\bullet}\times S,D'_{\bullet})/S)
\to\mathbb Z^{tr}((\bar X'\times S,X')/S))(-d')[-2d']),
E_{zar}(\Omega^{\bullet,\Gamma,pr}_{/S},F_{DR})):=\hat{\mathcal F_S^{FDR}}(g)
\end{eqnarray*}
which is a filtered quasi-isomorphism by proposition \ref{aetfib}.
Note that in the particular case where $\bar g=i:Z\hookrightarrow X$ is a closed embedding 
of codimension $c=\dim(X)-\dim(Z)$,
\begin{eqnarray*}
p_{S*}E(\Omega^{\bullet}_{Z\times S/S}\otimes_{O_{Z\times S}}\Gamma^{\vee,Hdg}_{Z}(O_{Z\times S},F_b))(-c)[-2c] \\
\xrightarrow{\Omega^{\Gamma,pr}([\Gamma_i^t])=(\alpha\mapsto\alpha\wedge\gamma}
p_{S*}E(\Omega^{\bullet}_{X\times S/S}\otimes_{O_{X\times S}}\Gamma^{\vee,Hdg}_X(O_{X\times S},F_b)))
\end{eqnarray*}
is the Gynsin morphism, where $\gamma\in H^0\Gamma(X,N^{\bullet}_{Z/X}[2c])$ is the De Rham fundamental class
of the conormal bundle, also note that denoting by $j:X\backslash Z\hookrightarrow X$ the open complementary
the composite in $C_{\mathcal D,fil}(S)$ 
\begin{eqnarray*}
p_{S*}E(\Omega^{\bullet}_{Z\times S/S}\otimes_{O_{Z\times S}}\Gamma^{\vee,Hdg}_{Z}(O_{Z\times S},F_b))(-c)[-2c] \\
\xrightarrow{\Omega^{\Gamma,pr}([\Gamma_i^t])=(\alpha\mapsto\alpha\wedge\gamma}
p_{S*}E(\Omega^{\bullet}_{X\times S/S}\otimes_{O_{X\times S}}\Gamma^{\vee,Hdg}_X(O_{X\times S},F_b))) \\
\xrightarrow{\Omega^{\Gamma,pr}(j)}
p_{S*}E(\Omega^{\bullet}_{(X\backslash Z)\times S/S}\otimes_{O_{(X\backslash Z)\times S}}
\Gamma^{\vee,Hdg}_{(X\backslash Z)}(O_{(X\backslash Z)\times S},F_b)))
\end{eqnarray*}
is NOT equal to zero, and the composite in $C(\Var(\mathbb C)^{2,smpr}/S)$
\begin{eqnarray*}
\mathbb Z^{tr}(((X\backslash Z)\times S,X\backslash Z)/S)\xrightarrow{\mathbb Z(j\times I)}
\mathbb Z^{tr}((X\times S,X)/S)\xrightarrow{[\Gamma_i^t]}\mathbb Z^{tr}((Z\times S,Z)/S)(c)[2c]
\end{eqnarray*}
is NOT equal to zero.
Now, for $g:U'/\tilde S_I\to U/\tilde S_I$ the commutative diagram
\begin{equation*}
\xymatrix{p_{IJ}^{*mod}\hat{\mathcal F^{FDR}}_{\tilde S_I}(g) 
\ar[rr]^{T(p_{IJ},\Omega^{\gamma,pr})(-)}\ar[d]_{p_{IJ}^{*mod}A_I} & \, &
\hat{\mathcal F^{FDR}}_{\tilde S_J}(g\times I)
\ar[d]^{A_J} \\
p_{IJ}^{*mod}\Omega^{\Gamma,pr}(\Gamma_{\bar g}^t)\ar[rr]^{T^O_w(p_{IJ},p_{\tilde S_I})(-)} & \, & 
\Omega^{\Gamma,pr}(\Gamma_{\bar g\times I}^t)}
\end{equation*}
with $g\times I:U'\times\tilde S_{J\backslash I}/\tilde S_J\to U\times\tilde S_{J\backslash I}/\tilde S_J$ and
\begin{eqnarray*}
A_I:=((DR(\bar X'\times\tilde S_I/\tilde S_I)(\ad((n'\times I)_{!Hdg},(n'\times I)^{!Hdg})(O_{\bar X'\times\tilde S_I},F_b)),0), \\
(DR(\bar X\times\tilde S_I/\tilde S_I)(\ad((n\times I)_{!Hdg},(n\times I)^{!Hdg})(O_{\bar X\times\tilde S_I},F_b)),0)).
\end{eqnarray*}
\end{proof}

\begin{prop}\label{FDRHdgwelldef}
For $S\in\Var(\mathbb C)$ not smooth, the functor (see corollary \ref{FDRMHM}) 
\begin{equation*}
\iota_S^{-1}\mathcal F_S^{FDR}:\DA_c^-(S)^{op}\to\pi_S(D(MHM(S))
\end{equation*}
does not depend on the choice of the open cover $S=\cup_iS_i$
and the closed embeddings $i_i:S_i\hookrightarrow\tilde S_i$ with $\tilde S_i\in\SmVar(\mathbb C)$.
\end{prop}

\begin{proof}
Let $S=\cup_{i=l+1}^{i=l'}S_i$ is an other open cover together with closed embeddings $i_i:S_i\hookrightarrow\tilde S_i$
with $\tilde S_i\in\SmVar(\mathbb C)$ for $l+1\leq i\leq l'$. Then, 
for $J'\subset I'\subset[l+1,\ldots,l']=L'$ and $J\subset I\subset L=[1,\ldots,l]$, 
\begin{eqnarray*}
T_S^{L/L'}(\iota_S^{-1}(e'(\tilde S_I)_*\mathcal Hom(\hat R^{CH}(\rho_{\tilde S_I}^*L(i_{I*}j_I^*F)),
E_{zar}(\Omega_{/\tilde S_I}^{\bullet,\Gamma,pr},F_{DR}))[-d_{\tilde S_I}],u_{IJ}(F))) \\
:=((\ho\lim_{I\in L}p_{I'(I\sqcup I')*}\Gamma^{\vee,Hdg}_{S_{I\sqcup I'}}p_{I(I\sqcup I')}^{*mod[-]}
(\iota_S^{-1}(e'(\tilde S_I)_*\mathcal Hom(\hat R^{CH}(\rho_{\tilde S_I}^*L(i_{I*}j_I^*F)), \\ 
E_{zar}(\Omega_{/\tilde S_I}^{\bullet,\Gamma,pr},F_{DR}))[-d_{\tilde S_I}],T_S^{L/L'}(u_{IJ}(F))))_I) 
\xrightarrow{(\ho\lim_{I\in L} u_{I(I\sqcup I')}(F))} \\
(\ho\lim_{I\in L}p_{I'(I\sqcup I')*}\Gamma^{\vee,Hdg}_{S_{I\sqcup I'}}p_{I(I\sqcup I')}^{*mod[-]}p_{I(I\sqcup I')*}
e'(\tilde S_{(I\sqcup I')})_*\mathcal Hom(\hat R^{CH}(\rho_{\tilde S_I}^*L(i_{(I\sqcup I')*}j_{(I\sqcup I')}^*F)), \\
E_{zar}(\Omega_{/\tilde S_{I\sqcup I'}}^{\bullet,\Gamma,pr},F_{DR}))[-d_{\tilde S_{I\sqcup I'}}], 
u_{(I\sqcup I')(I\sqcup J')}(F)) 
\xrightarrow{\ad(p_{I(I\sqcup I')}^{*mod},p_{I(I\sqcup I')*})(-)\circ\gamma_{S_{I\sqcup I'}}^{\vee,Hdg}(-)} \\
(\ho\lim_{I\in L}p_{I'(I\sqcup I')*}e'(\tilde S_{(I\sqcup I')})_*
\mathcal Hom(\hat R^{CH}(\rho_{\tilde S_I}^*L(i_{(I\sqcup I')*}j_{(I\sqcup I')}^*F)), \\ 
E_{zar}(\Omega_{/\tilde S_{I\sqcup I'}}^{\bullet,\Gamma,pr},F_{DR}))[-d_{\tilde S_{I\sqcup I'}}], 
u_{(I\sqcup I')(I\sqcup J')}(F)) \\
\xleftarrow{(\ho\lim_{I\in L}u_{I'(I\sqcup I')}(F))} 
(e'(\tilde S_{I'})_*\mathcal Hom(\hat R^{CH}(\rho_{\tilde S_{I'}}^*L(i_{I'*}j_{I'}^*F)),
E_{zar}(\Omega_{/\tilde S_{I'}}^{\bullet,\Gamma,pr},F_{DR}))[-d_{\tilde S_{I'}}],u_{I'J'}(F))
\end{eqnarray*}
is a filtered Zariski local equivalence, 
since all the morphisms are filtered Zariski local equivalences by proposition \ref{keyalgsing1} 
and proposition \ref{TSIHdg}.
\end{proof}

We have the canonical transformation map between the filtered De Rham realization functor and the Gauss-Manin realization functor :

\begin{defi}\label{GMFDRdef}
Let $S\in\Var(\mathbb C)$ and $S=\cup_{i=1}^l S_i$ an open cover such that there exist closed embeddings
$i_i:S_i\hookrightarrow\tilde S_i$ with $\tilde S_i\in\SmVar(\mathbb C)$. 
Let $M\in\DA_c(S)$ and $F\in C(\Var(\mathbb C)^{sm}/S)$ such that $M=D(\mathbb A^1,et)(F)$.
We have, using definition \ref{wtildew}(ii), the canonical map in $D_{O_Sfil,\mathcal D,\infty}(S/(\tilde S_I))$
\begin{eqnarray*}
T(\mathcal F^{GM}_S,\mathcal F^{FDR}_S)(M): \\
\mathcal F_S^{GM}(L\mathbb D_SM):=
(e(\tilde S_I)_*\mathcal Hom^{\bullet}(L(i_{I*}j_I^*\mathbb D_SLF),
E_{zar}(\Omega^{\bullet}_{/\tilde S_I},F_b))[-d_{\tilde S_I}],u^q_{IJ}(F)) \\
\xrightarrow{\sim} 
(e(\tilde S_I)_*\mathcal Hom^{\bullet}(L\mathbb D^0_{\tilde S_I}L(i_{I*}j_I^*F),
E_{zar}(\Omega^{\bullet}_{/\tilde S_I},F_b))[-d_{\tilde S_I}],u^{q,d}_{IJ}(F)) \\
\xrightarrow{(\mathcal Hom(-,\Gr(\Omega_{\tilde S_I})))^{-1}}
(e(\tilde S_I)_*\mathcal Hom^{\bullet}(L\mathbb D^0_{\tilde S_I}L(i_{I*}j_I^*F),
\Gr^{12}_{\tilde S_I*}E_{zar}(\Omega^{\bullet,\Gamma,pr}_{/\tilde S_I},F_{DR}))[-d_{\tilde S_I}],u^{q,d}_{IJ}(F)) \\
\xrightarrow{(\mathcal Hom^{\bullet}(r^{0CH}_{\tilde S_I}(L(i_{I*}j_I^*F))\circ T(\hat R^{0CH},R^{0CH})(L(i_{I*}j_I^*F))),-)} \\
(e(\tilde S_I)_*\mathcal Hom^{\bullet}(\hat R^{CH}(\rho_{\tilde S_I}^*L(i_{I*}j_I^*F)),
\Gr^{12}_{\tilde S_I*}E_{zar}(\Omega^{\bullet,\Gamma,pr}_{/\tilde S_I},F_{DR}))[-d_{\tilde S_I}],u^{q,d}_{IJ}(F)) \\
\xrightarrow{(I(\Gr^{12*}_{\tilde S_I},\Gr^{12}_{\tilde S_I*})(-,-))} 
(e(\tilde S_I)_*\mathcal Hom^{\bullet}(\Gr^{12*}_{\tilde S_I}\hat R^{0CH}(\rho_{\tilde S_I}^*L(i_{I*}j_I^*F)),
E_{zar}(\Omega^{\bullet,\Gamma,pr}_{/\tilde S_I},F_{DR}))[-d_{\tilde S_I}],u^{q,d}_{IJ}(F)) \\
\xrightarrow{=}
(e'(\tilde S_I)_*\mathcal Hom(\hat R^{CH}(\rho_{\tilde S_I}^*L(i_{I*}j_I^*F)),
E_{zar}(\Omega^{\bullet,\Gamma,pr}_{/\tilde S_I},F_{DR}))[-d_{\tilde S_I}],u^q_{IJ}(F))=:\mathcal F^{FDR}_S(M)
\end{eqnarray*}
\end{defi}

\begin{prop}\label{TGMFDRprop}
Let $S\in\Var(\mathbb C)$ and $S=\cup_{i=1}^l S_i$ an open cover such that there exist closed embeddings
$i_i:S_i\hookrightarrow\tilde S_i$ with $\tilde S_i\in\SmVar(\mathbb C)$. 
\begin{itemize}
\item[(i)] For $M\in\DA_c(S)$ the map in $D_{O_S,\mathcal D}(S/(\tilde S_I))=D_{O_S,\mathcal D}(S)$
\begin{equation*}
o_{fil}T(\mathcal F^{GM}_S,\mathcal F^{FDR}_S)(M):
o_{fil}\mathcal F_S^{GM}(L\mathbb D_SM)\xrightarrow{\sim} o_{fil}\mathcal F^{FDR}_S(M)
\end{equation*}
given in definition \ref{GMFDRdef} is an isomorphism if we forgot the Hodge filtration $F$.
\item[(ii)]For $M\in\DA_c(S)$ and all $n,p\in\mathbb Z$, the map in $\PSh_{O_S,\mathcal D}(S/(\tilde S_I))$
\begin{equation*}
F^pH^nT(\mathcal F^{GM}_S,\mathcal F^{FDR}_S)(M):F^pH^n\mathcal F_S^{GM}(L\mathbb D_SM)\hookrightarrow F^pH^n\mathcal F^{FDR}_S(M)
\end{equation*}
given in definition \ref{GMFDRdef} is a monomorphism.
Note that $F^pH^nT(\mathcal F^{GM}_S,\mathcal F^{FDR}_S)(M)$ is NOT an isomorphism in general :
take for example $M(S^o/S)^{\vee}=D(\mathbb A^1,et)(j_*E_{zar}(\mathbb Z(S^o/S)))$ 
for an open embedding $j:S^o\hookrightarrow S$, then 
\begin{equation*}
H^n\mathcal F_S^{GM}(L\mathbb D_SM(S^o/S)^{\vee})=\mathcal F_S^{GM}(\mathbb Z(S^o/S))=j_*E(O_{S^o},F_b)\notin\pi_S(MHM(S)) 
\end{equation*}
and hence is NOT isomorphic to $H^n\mathcal F_S^{FDR}(L\mathbb D_SM(S^o/S)^{\vee})\in\pi_S(MHM(S))$ 
as filtered $D_S$-modules (see remark \ref{remHdgkey}). 
It is an isomorphism in the very particular cases where $M=D(\mathbb A^1,et)(\mathbb Z(X/S))$
or $M=D(\mathbb A^1,et)(\mathbb Z(X^o/S))$ for $f:X\to S$ is a smooth proper morphism and $n:X^o\hookrightarrow X$ is an open subset
such that $X\backslash X^o=\cup D_i$ is a normal crossing divisor 
and such that $f_{|D_i}=f\circ i_i:D_i\to X$ are SMOOTH morphism with $i_i:D_i\hookrightarrow X$ the closed embedding and
considering $f_{|X^o}=f\circ n:X^o\to S$ (see proposition \ref{FFXD}).
\end{itemize}
\end{prop}

\begin{proof}
\noindent(i):Follows from the computation for a Borel-Moore motive.

\noindent(ii):Follows from (i).
\end{proof}

We now define the functorialities of $\mathcal F_S^{FDR}$ with respect to $S$ 
which makes $\mathcal F^{-}_{FDR}$ a morphism of 2 functor.

\begin{defi}\label{TGammaFDR}
Let $S\in\Var(\mathbb C)$. Let $Z\subset S$ a closed subset.
Let $S=\cup_{i=1}^l S_i$ an open cover such that there exist closed embeddings
$i_i:S_i\hookrightarrow\tilde S_i$ with $\tilde S_i\in\SmVar(\mathbb C)$.
Denote $Z_I:=Z\cap S_I$. We then have closed embeddings $Z_I\hookrightarrow S_I\hookrightarrow\tilde S_I$.
\begin{itemize}
\item[(i)]For $F\in C(\Var(\mathbb C)^{sm}/S)$, we will consider the following canonical map 
in $\pi_S(D(MHM(S)))\subset D_{\mathcal D(1,0)fil}(S/(\tilde S_I))$
\begin{eqnarray*}
T(\Gamma_Z^{\vee,Hdg},\Omega^{\Gamma,pr}_{/S})(F,W):\\
\Gamma_Z^{\vee,Hdg}\iota_S^{-1}(e'_*\mathcal Hom^{\bullet}(\hat R^{CH}(\rho_{\tilde S_I}^*L(i_{I*}j_I^*(F,W))),  
E_{zar}(\Omega^{\bullet,\Gamma,pr}_{/\tilde S_I},F_{DR}))[-d_{\tilde S_I}],u^q_{IJ}(F,W)) \\ 
\xrightarrow{\mathcal Hom^{\bullet}(\hat R_{\tilde S_I}^{CH}(\gamma^{\vee,Z_I}(L(i_{I*}j_I^*(F,W)))),
E_{zar}(\Omega^{\bullet,\Gamma,pr}_{/\tilde S_I},F_{DR}))} \\ 
\Gamma_Z^{\vee,Hdg}\iota_S^{-1}(e'_*\mathcal Hom^{\bullet}(
\hat R^{CH}(\rho_{\tilde S_I}^*\Gamma^{\vee}_{Z_I}L(i_{I*}j_I^*(F,W))),  
E_{zar}(\Omega^{\bullet,\Gamma,pr}_{/\tilde S_I},F_{DR}))[-d_{\tilde S_I}],u^{q,Z}_{IJ}(F,W)) \\
\xrightarrow{=}  
\iota_S^{-1}(e'_*\mathcal Hom^{\bullet}(\hat R^{CH}(\rho_{\tilde S_I}^*\Gamma^{\vee}_{Z_I}L(i_{I*}j_I^*(F,W))),  
E_{zar}(\Omega^{\bullet,\Gamma,pr}_{/\tilde S_I},F_{DR}))[-d_{\tilde S_I}],u^{q,Z}_{IJ}(F,W)).
\end{eqnarray*}
with
\begin{eqnarray*}
u^{q,Z}_{IJ}(F)[d_{\tilde S_I}]:
e'(\tilde S_I)_*\mathcal Hom^{\bullet}(\hat R^{CH}(\rho_{\tilde S_I}^*\Gamma^{\vee}_{Z_I}L(i_{I*}j_I^*F)),
E_{zar}(\Omega^{\bullet,\Gamma,pr}_{/\tilde S_I},F_{DR})) \\
\xrightarrow{p_{IJ*}T(p_{IJ},\Omega^{\gamma,pr}_{\cdot})(-)\circ \ad(p_{IJ}^{*mod},p_{IJ})(-)} \\ 
p_{IJ*}e'(\tilde S_J)_*\mathcal Hom^{\bullet}(p_{IJ}^*\hat R^{CH}(\rho_{\tilde S_I}^*\Gamma^{\vee}_{Z_I}L(i_{I*}j_I^*F)),
E_{zar}(\Omega^{\bullet,\Gamma,pr}_{/\tilde S_J},F_{DR})) \\
\xrightarrow{\mathcal Hom(\hat T(p_{IJ},R^{CH})(Li_{I*}j_I^*F)^{-1},
E_{zar}(\Omega_{/\tilde S_J}^{\bullet,\Gamma,pr},F_{DR}))} \\
p_{IJ*}e'(\tilde S_J)_*\mathcal Hom^{\bullet}(\hat R^{CH}(\rho_{\tilde S_J}^*p_{IJ}^*\Gamma^{\vee}_{Z_I}L(i_{I*}j_I^*F)),
E_{zar}(\Omega^{\bullet,\Gamma,pr}_{/\tilde S_J},F_{DR})) \\
\xrightarrow{\mathcal Hom(\hat R^{CH}_{\tilde S_J}(T^q(D_{IJ})(j_I^*F)
\circ T(Z_J/Z_I\times\tilde S_{J\backslash I},\gamma^{\vee})(-)\circ T(p_{IJ},\gamma)(-)),
E_{zar}(\Omega_{/\tilde S_J}^{\bullet,\Gamma,pr},F_{DR}))} \\
p_{IJ*}e'(\tilde S_J)_*\mathcal Hom^{\bullet}(\hat R^{CH}(\rho_{\tilde S_J}^*\Gamma^{\vee}_{Z_J}L(i_{J*}j_J^*F)),
E_{zar}(\Omega^{\bullet,\Gamma,pr}_{/\tilde S_J},F_{DR})).
\end{eqnarray*}
\item[(ii)]For $F\in C(\Var(\mathbb C)^{sm}/S)$, we have also the following canonical map 
in $\pi_S(D(MHM(S)))\subset D_{\mathcal D(1,0)fil}(S/(\tilde S_I))$
\begin{eqnarray*}
T(\Gamma_Z^{Hdg},\Omega^{\Gamma,pr}_{/S})(F,W):\\
\iota_S^{-1}(e'_*\mathcal Hom^{\bullet}(\hat R^{CH}(\rho_{\tilde S_I}^*L\Gamma_{Z_I}E(i_{I*}j_I^*\mathbb D_S(F,W))),  
E_{zar}(\Omega^{\bullet,\Gamma,pr}_{/\tilde S_I},F_{DR}))[-d_{\tilde S_I}],u^{q,Z,d}_{IJ}(F,W)) 
\xrightarrow{=} \\ 
\Gamma_Z^{Hdg}\iota_S^{-1}(e'_*\mathcal Hom^{\bullet}(
\hat R^{CH}(\rho_{\tilde S_I}^*L\Gamma_{Z_I}E(i_{I*}j_I^*\mathbb D_S(F,W))),  
E_{zar}(\Omega^{\bullet,\Gamma,pr}_{/\tilde S_I},F_{DR}))[-d_{\tilde S_I}],u^{q,Z,d}_{IJ}(F,W)) \\
\xrightarrow{\mathcal Hom^{\bullet}(\hat R_{\tilde S_I}^{CH}(\gamma^{Z_I}(-)),
E_{zar}(\Omega^{\bullet,\Gamma,pr}_{/\tilde S_I},F_{DR}))} \\
\Gamma_Z^{Hdg}\iota_S^{-1}(e'_*\mathcal Hom^{\bullet}(\hat R^{CH}(\rho_{\tilde S_I}^*L(i_{I*}j_I^*\mathbb D_S(F,W))), 
E_{zar}(\Omega^{\bullet,\Gamma,pr}_{/\tilde S_I},F_{DR}))[-d_{\tilde S_I}],u^q_{IJ}(F,W)) 
\end{eqnarray*}
with
\begin{eqnarray*}
u^{q,Z}_{IJ}(F)[d_{\tilde S_I}]:e'(\tilde S_I)_*\mathcal Hom^{\bullet}(L\rho_{\tilde S_I*}\mu_{\tilde S_I*}
R^{CH}(\rho_{\tilde S_I}^*L\Gamma_{Z_I}E(i_{I*}j_I^*\mathbb D_SLF)),
E_{zar}(\Omega^{\bullet,\Gamma,pr}_{/\tilde S_I},F_{DR})) \\
\xrightarrow{p_{IJ*}T(p_{IJ},\Omega^{\gamma,pr}_{\cdot})(-)\circ \ad(p_{IJ}^{*mod},p_{IJ})(-)} \\ 
p_{IJ*}e'(\tilde S_J)_*\mathcal Hom^{\bullet}(
p_{IJ}^*\hat R^{CH}(\rho_{\tilde S_I}^*L\Gamma_{Z_I}E(i_{I*}j_I^*\mathbb D_SF)),
E_{zar}(\Omega^{\bullet,\Gamma,pr}_{/\tilde S_J},F_{DR})) \\
\xrightarrow{\mathcal Hom(T(p_{IJ},R^{CH})(Li_{I*}j_I^*F)^{-1},E_{zar}(\Omega_{/\tilde S_J}^{\bullet,\Gamma,pr},F_{DR}))} \\
p_{IJ*}e'(\tilde S_J)_*\mathcal Hom^{\bullet}(\hat R^{CH}(\rho_{\tilde S_J}^*Lp_{IJ}^*\Gamma_{Z_I}E(i_{I*}j_I^*\mathbb D_SLF)),
E_{zar}(\Omega^{\bullet,\Gamma,pr}_{/\tilde S_J},F_{DR})) \\
\xrightarrow{\mathcal Hom(\hat R^{CH}_{\tilde S_J}(\mathbb D_{\tilde S_J}S^q(D_{IJ})(\mathbb D_SLF)),
E_{zar}(\Omega_{/\tilde S_J}^{\bullet,\Gamma,pr},F_{DR}))} \\
p_{IJ*}e'(\tilde S_J)_*\mathcal Hom^{\bullet}(
\hat R^{CH}(\rho_{\tilde S_J}^*\Gamma_{Z_J}E(i_{J*}j_J^*\mathbb D_SLF)),E_{zar}(\Omega^{\bullet,\Gamma,pr}_{/\tilde S_J},F_{DR})).
\end{eqnarray*}
\end{itemize}
\end{defi}

This transformation map will, with the projection case, gives the transformation between the pullback functor :

\begin{defi}\label{TgDRdef}
Let $g:T\to S$ a morphism with $T,S\in\SmVar(\mathbb C)$.
Consider the factorization $g:T\xrightarrow{l}T\times S\xrightarrow{p_S}S$
where $l$ is the graph embedding and $p_S$ the projection.
Let $M\in\DA_c(S)^-$ and $(F,W)\in C_{fil}(\Var(\mathbb C)^{sm}/S)$ such that $(M,W)=D(\mathbb A^1_S,et)(F,W)$. 
Then, $D(\mathbb A^1_T,et)(g^*F)=g^*M$ and there exist 
$(F',W)\in C_{fil}(\Var(\mathbb C)^{sm}/T\times S)$ and an equivalence $(\mathbb A^1,et)$ local 
$e:\Gamma^{\vee}p_S^*(F,W)\to(F',W)$ such that $D(\mathbb A^1_{T\times S},et)(F',W)=(\Gamma^{\vee}p_S^*M,W)$. 
We have then the canonical transformation in $\pi_T(D(MHM(T))$
using definition \ref{TgDR} and definition \ref{TGammaFDR}(i) :
\begin{eqnarray*}
T(g,\mathcal F^{FDR})(M):g^{\hat*mod,Hdg}\mathcal F_S^{FDR}(M):= \\ 
\Gamma^{\vee,Hdg}_T\iota_T^{-1}(p_S^{*mod[-]}(
e(S)_*\Gr^{12}_{S*}\mathcal Hom^{\bullet}(\hat R^{CH}(\rho_S^*L(F,W)),E_{zar}(\Omega^{\bullet,\Gamma,pr}_{/S},F_{DR}))[-d_S])) \\ 
\xrightarrow{T(p_S,\Omega^{\Gamma,pr}_{/\cdot})(-)} \\ 
\Gamma^{\vee,Hdg}_T(e'(T\times S)_*\mathcal Hom^{\bullet}(p_S^*\hat R^{CH}(\rho_S^*L(F)), 
E_{zar}(\Omega^{\bullet,\Gamma,pr}_{/T\times S},F_{DR}))[-d_S])) \\
\xrightarrow{\mathcal Hom(T(p_S,\hat R^{CH})(L(F,W))^{-1},E_{zar}(\Omega^{\bullet,\Gamma,pr}_{/S},F_{DR}))[-d_S]} \\
\Gamma^{\vee,Hdg}_T\iota_T^{-1}(e'(T\times S)_*\mathcal Hom^{\bullet}(\hat R^{CH}(\rho_{T\times S}^*p_S^*L(F,W)), 
E_{zar}(\Omega^{\bullet,\Gamma,pr}_{/T\times S},F_{DR}))[-d_S]) 
\xrightarrow{=} \\
\Gamma^{\vee,Hdg}_T\iota_T^{-1}(e'(T\times S)_*\mathcal Hom^{\bullet}(\hat R^{CH}(\rho_{T\times S}^*p_S^*L(F,W)), 
E_{zar}(\Omega^{\bullet,\Gamma,pr}_{/T\times S},F_{DR}))[-d_S]) \\
\xrightarrow{T(\Gamma_T^{\vee,Hdg},\Omega^{\Gamma,pr}_{/T\times S})(F,W)} \\
(e'(T\times S)_*\mathcal Hom^{\bullet}(\hat R^{CH}(\rho_{T\times S}^*\Gamma^{\vee}_Tp_S^*L(F,W)), 
E_{zar}(\Omega^{\bullet,\Gamma,pr}_{/T\times S},F_{DR}))[-d_S]) \\ 
\xrightarrow{\mathcal Hom(\hat R^{CH}_{T\times S}(e),-)} \\
(e'(T\times S)_*\mathcal Hom^{\bullet}(\hat R^{CH}(\rho_{T\times S}^*L(F',W)), 
E_{zar}(\Omega^{\bullet,\Gamma,pr}_{/T\times S},F_{DR}))[-d_S]) \\
=:\mathcal F_{T\times S}^{FDR}(l_*g^*M)=\mathcal F_T^{FDR}(g^*M)
\end{eqnarray*}
where the last equality follows from proposition \ref{FDRHdgwelldef}.
\end{defi}

We give now the definition in the non smooth case
Let $g:T\to S$ a morphism with $T,S\in\Var(\mathbb C)$.
Assume we have a factorization $g:T\xrightarrow{l}Y\times S\xrightarrow{p_S}S$
with $Y\in\SmVar(\mathbb C)$, $l$ a closed embedding and $p_S$ the projection.
Let $S=\cup_{i=1}^lS_i$ be an open cover such that 
there exists closed embeddings $i_i:S_i\hookrightarrow\tilde S_i$ with $\tilde S_i\in\SmVar(\mathbb C)$
Then, $T=\cup^l_{i=1} T_i$ with $T_i:=g^{-1}(S_i)$
and we have closed embeddings $i'_i:=i_i\circ l:T_i\hookrightarrow Y\times\tilde S_i$,
Moreover $\tilde g_I:=p_{\tilde S_I}:Y\times\tilde S_I\to\tilde S_I$ is a lift of $g_I:=g_{|T_I}:T_I\to S_I$.
We recall the commutative diagram :
\begin{equation*}
E_{IJg}=\xymatrix{(Y\times\tilde S_I)\backslash T_I\ar[d]^{p_{\tilde S_I}}\ar[r]^{m'_I} & 
Y\times\tilde S_J\ar[d]^{\tilde g_I} \\
\tilde S_I\backslash S_I\ar[r]^{m_I}& \tilde S_I}, \,
E_{IJ}=\xymatrix{\tilde S_J\backslash S_J\ar[d]^{p_{IJ}}\ar[r]^{m_J} & \tilde S_J\ar[d]^{p_{IJ}} \\
\tilde S_I\backslash(S_I\backslash S_J)\ar[r]^{m=m_{IJ}}& \tilde S_I}
E'_{IJ}=\xymatrix{(Y\times\tilde S_J)\backslash T_J\ar[d]^{p'_{IJ}}\ar[r]^{m'_J} & Y\times\tilde S_J\ar[d]^{p'_{IJ}} \\
(Y\times\tilde S_I)\backslash(T_I\backslash T_J)\ar[r]^{m'=m'_{IJ}}& Y\times\tilde S_I}
\end{equation*}
For $I\subset J$, denote by $p_{IJ}:\tilde S_J\to\tilde S_I$ and 
$p'_{IJ}:=I_Y\times p_{IJ}:Y\times\tilde S_J\to Y\times\tilde S_I$ the projections,
so that $\tilde g_I\circ p'_{IJ}=p_{IJ}\circ\tilde g_J$.
Consider, for $I\subset J\subset\left[1,\ldots,l\right]$,
resp. for each $I\subset\left[1,\ldots,l\right]$, the following commutative diagrams in $\Var(\mathbb C)$
\begin{equation*}
D_{IJ}=\xymatrix{ S_I\ar[r]^{i_I} & \tilde S_I \\
S_J\ar[u]^{j_{IJ}}\ar[r]^{i_J} & \tilde S_J\ar[u]^{p_{IJ}}} \; , \;  
D'_{IJ}=\xymatrix{ T_I\ar[r]^{i'_I} & Y\times\tilde S_I \\
T_J\ar[u]^{j'_{IJ}}\ar[r]^{i'_J} & Y\times\tilde S_J\ar[u]^{p'_{IJ}}}  
D_{gI}=\xymatrix{ S_I\ar[r]^{i_I} & \tilde S_I \\
T_I\ar[u]^{g_I}\ar[r]^{i'_I} & Y\times\tilde S_I\ar[u]^{\tilde g_I}} \; , \;  
\end{equation*}
and $j_{IJ}:S_J\hookrightarrow S_I$ is the open embedding so that $j_I\circ j_{IJ}=j_J$.
Let $F\in C(\Var(\mathbb C)^{sm}/S)$. 
Recall (see section 2) that since $j_I^{'*}i'_{I*}j_I^{'*}g^*F=0$,
the morphism $T(D_{gI})(j_I^*F):\tilde g_I^*i_{I*}j_I^*F\to i'_{I*}j_I^{'*}g^*F$ factors trough
\begin{equation*}
T(D_{gI})(j_I^*F):\tilde g_I^*i_{I*}j_I^*F\xrightarrow{\gamma_{X_I}^{\vee}(-)}
\Gamma_{X_I}^{\vee}\tilde g_I^*i_{I*}j_I^*F\xrightarrow{T^{\gamma}(D_{gI})(j_I^*F)}i'_{I*}j_I^{'*}g^*F.
\end{equation*}  
and that the fact that the diagrams (\ref{DgDIq}) commutes says that the maps $T^{q,\gamma}(D_{gI})(j_I^*F)$
define a morphism in $C(\Var(\mathbb C)^{sm}/(T/(Y\times\tilde S_I)))$
\begin{eqnarray*}
(T^{q,\gamma}(D_{gI})(j_I^*F)):
(\Gamma^{\vee}_{T_I}\tilde g_I^*L(i_{I*}j_I^*F),
T^q(D_{IJ})(j_I^*F)\circ T(T_I/T_I\times\tilde S_{J\backslash I},\gamma^{\vee})(-)\circ T(p'_{IJ},\gamma^{\vee})(-))) \\
\to (L(i'_{I*}j_I^{'*}g^*F),T^q(D'_{IJ})(j_I^{'*}g^*F))
\end{eqnarray*}
Denote for short $d_{YI}:=-d_Y-d_{\tilde S_I}$.
We denote by $\tilde g_J^*u^q_{IJ}(F)$ the composite
\begin{eqnarray*}
\tilde g_J^*u^q_{IJ}(F)[-d_{YJ}]: \\
e'(Y\times\tilde S_I)_*\mathcal Hom(\hat R^{CH}(\rho_{Y\times\tilde S_I}^*\tilde g_I^*L(i_{I*}j_I^*F)),
E_{zar}(\Omega^{\bullet,\Gamma,pr}_{/Y\times\tilde S_I},F_{DR})) \\
\xrightarrow{p'_{IJ*}T(p'_{IJ},\Omega^{\Gamma,pr})(-)\circ\ad(p_{IJ}^{'*mod},p'_{IJ*})(-)} \\
p'_{IJ*}e'(Y\times\tilde S_I)_*\mathcal Hom(p_{IJ}^{'*}\hat R^{CH}(\rho_{Y\times\tilde S_I}^*\tilde g_I^*L(i_{I*}j_I^*F)),
E_{zar}(\Omega^{\bullet,\Gamma,pr}_{/Y\times\tilde S_J},F_{DR})) \\
\xrightarrow{\mathcal Hom(T(p'_{IJ},\hat R^{CH})()^{-1},
E_{zar}(\Omega^{\bullet,\Gamma,pr}_{/Y\times\tilde S_J},F_{DR})} \\
p'_{IJ*}e'(Y\times\tilde S_I)_*\mathcal Hom(\hat R^{CH}(\rho_{Y\times\tilde S_J}^*p_{IJ}^{'*}\tilde g_I^*L(i_{I*}j_I^*F)),
E_{zar}(\Omega^{\bullet,\Gamma,pr}_{/Y\times\tilde S_J},F_{DR})) \\
\xrightarrow{\mathcal Hom(\hat R^{CH}_{Y\times\tilde S_J}(T^q(D_{IJ})(j_I^*F)),
E_{zar}(\Omega^{\bullet,\Gamma,pr}_{/Y\times\tilde S_J},F_{DR}))} \\
e'(Y\times\tilde S_J)_*\mathcal Hom(\hat R^{CH}(\rho_{Y\times\tilde S_J}^*\tilde g_J^*L(i_{J*}j_J^*F)),
E_{zar}(\Omega^{\bullet,\Gamma,pr}_{/Y\times\tilde S_I},F_{DR})).
\end{eqnarray*}
We denote by $\tilde g_J^{*,\gamma}u^q_{IJ}(F)$ the composite
\begin{eqnarray*}
\tilde g_J^{*,\gamma}u^q_{IJ}(F)[-d_{YJ}]: \\
e'(Y\times\tilde S_I)_*\mathcal Hom(\hat R^{CH}(\rho_{Y\times\tilde S_I}^*\Gamma^{\vee}_{T_I}\tilde g_I^*L(i_{I*}j_I^*F)),
E_{zar}(\Omega^{\bullet,\Gamma,pr}_{/Y\times\tilde S_I},F_{DR})) \\
\xrightarrow{p'_{IJ*}T(p'_{IJ},\Omega^{\Gamma,pr})(-)\circ\ad(p_{IJ}^{'*mod},p'_{IJ*})(-)} \\ 
p'_{IJ*}e'(Y\times\tilde S_I)_*\mathcal Hom(p_{IJ}^{'*}
\hat R^{CH}(\rho_{Y\times\tilde S_I}^*\Gamma^{\vee}_{T_I}\tilde g_I^*L(i_{I*}j_I^*F)),
E_{zar}(\Omega^{\bullet,\Gamma,pr}_{/Y\times\tilde S_J},F_{DR})) \\
\xrightarrow{\mathcal Hom(T(p'_{IJ},\hat R^{CH})()^{-1},
E_{zar}(\Omega^{\bullet,\Gamma,pr}_{/Y\times\tilde S_J},F_{DR})} \\
p'_{IJ*}e'(Y\times\tilde S_I)_*\mathcal Hom(
\hat R^{CH}(\rho_{Y\times\tilde S_J}^*p_{IJ}^{'*}\Gamma^{\vee}_{T_I}\tilde g_I^*L(i_{I*}j_I^*F)),
E_{zar}(\Omega^{\bullet,\Gamma,pr}_{/Y\times\tilde S_J},F_{DR})) \\
\xrightarrow{\mathcal Hom(\hat R^{CH}_{Y\times\tilde S_J}
(T^q(D_{IJ})(j_I^*F)\circ T(T_I/T_I\times\tilde S_{J\backslash I},\gamma^{\vee})(-)\circ T(p'_{IJ},\gamma^{\vee})(-)),
E_{zar}(\Omega^{\bullet,\Gamma,pr}_{/Y\times\tilde S_J},F_{DR}))} \\
e'(Y\times\tilde S_J)_*\mathcal Hom(\hat R^{CH}(\rho_{Y\times\tilde S_J}^*\Gamma^{\vee}_{T_J}\tilde g_J^*L(i_{J*}j_J^*F)),
E_{zar}(\Omega^{\bullet,\Gamma,pr}_{/Y\times\tilde S_I},F_{DR})).
\end{eqnarray*}

We then have then the following lemma :

\begin{lem}\label{Tglem}
\begin{itemize}
\item[(i)]The morphism in $C(\Var(\mathbb C)^{sm}/(T/(Y\times\tilde S_I)))$
\begin{eqnarray*}
(T^{q,\gamma}(D_{gI})(j_I^*F)):
(\Gamma^{\vee}_{T_I}L\tilde g_I^*i_{I*}j_I^*F,
T^q(D_{IJ})(j_I^*F)\circ T(T_I/T_I\times\tilde S_{J\backslash I},\gamma^{\vee})(-)\circ T(p'_{IJ},\gamma^{\vee})(-)) \\
\to (Li'_{I*}j_I^{'*}g^*F,T^q(D'_{IJ})(j_I^{'*}g^*F))
\end{eqnarray*}
is an equivalence $(\mathbb A^1,et)$ local.
\item[(ii)] The maps 
$\mathcal Hom((T^{q,\gamma}(D_{gI})(j_I^*F)),E_{zar}(\Omega^{\bullet,\Gamma,pr}_{/Y\times\tilde S_I}),F_{DR}))$
induce a filtered quasi-isomorphism in $C_{\mathcal Dfil}(T/(Y\times\tilde S_I))$
\begin{eqnarray*}
(\mathcal Hom(\hat R^{CH}_{Y\times\tilde S_I}(T^{q,\gamma}(D_{gI})(j_I^*F)),
E_{zar}(\Omega^{\bullet,\Gamma,pr}_{/Y\times\tilde S_I},F_{DR}))[d_{YI}]): \\
(e'(-)_*\mathcal Hom(\hat R^{CH}(\rho_{Y\times\tilde S_I}^*\Gamma^{\vee}_{T_I}\tilde g_I^*L(i_{I*}j_I^*F)),
E_{zar}(\Omega^{\bullet,\Gamma,pr}_{/Y\times\tilde S_I},F_{DR}))[d_{YI}],\tilde g_J^{*,\gamma}u^q_{IJ}(F)) \\
\to (e'(-)_*\mathcal Hom(\hat R^{CH}(\rho_{Y\times\tilde S_I}^*L(i'_{I*}j^{'*}_Ig^*F)),
E_{zar}(\Omega^{\bullet,\Gamma,pr}_{/Y\times\tilde S_I},F_{DR}))[d_{YI}],u^q_{IJ}(g^*F))
\end{eqnarray*}
\item[(iii)] The maps $T(\tilde g_I,\Omega^{\Gamma,pr}_{\cdot})(-)$ (see definition \ref{TgDR}),
induce a morphism in $C_{\mathcal Dfil}(T/(Y\times\tilde S_I))$
\begin{eqnarray*}
T(\tilde g_I,\Omega^{\Gamma,pr}_{/\cdot})(-))[d_{YI}]: \\ 
(\tilde g_I^{*mod}e'(\tilde S_I)_*\mathcal Hom^{\bullet}(\hat R^{CH}(\rho_{\tilde S_I}^*(L(i_{I*}j_I^*F)),
E_{zar}(\Omega^{\bullet,\Gamma,pr}_{/\tilde S_I},F_{DR}))[d_{YI}],\tilde g_J^{*mod}u^q_{IJ}(F)) \\
\to (e'(-)_*\mathcal Hom(\tilde g_I^*\hat R^{CH}(\rho_{\tilde S_I}^*L(i_{I*}j_I^*F)),
E_{zar}(\Omega^{\bullet,\Gamma,pr}_{/Y\times\tilde S_I},F_{DR}))[d_{YI}],\tilde g_J^*u^q_{IJ}(F)).
\end{eqnarray*}
\end{itemize}
\end{lem}

\begin{proof}
\noindent(i):Follows from theorem \ref{2functDM}.

\noindent(ii):These morphism induce a morphism in $C_{\mathcal Dfil}(T/(Y\times\tilde S_I))$ by construction.
The fact that this morphism is an filtered equivalence Zariski local follows from (i) and proposition \ref{projwach}.

\noindent(iii):These morphism induce a morphism in $C_{\mathcal Dfil}(T/(Y\times\tilde S_I))$ by construction.
\end{proof}

\begin{defi}\label{TgDRdefsing}
Let $g:T\to S$ a morphism with $T,S\in\Var(\mathbb C)$.
Assume we have a factorization $g:T\xrightarrow{l}Y\times S\xrightarrow{p_S}S$
with $Y\in\SmVar(\mathbb C)$, $l$ a closed embedding and $p_S$ the projection.
Let $S=\cup_{i=1}^lS_i$ be an open cover such that 
there exists closed embeddings $i_i:S_i\hookrightarrow\tilde S_i$ with $\tilde S_i\in\SmVar(\mathbb C)$
Then, $T=\cup^l_{i=1} T_i$ with $T_i:=g^{-1}(S_i)$
and we have closed embeddings $i'_i:=i_i\circ l:T_i\hookrightarrow Y\times\tilde S_i$,
Moreover $\tilde g_I:=p_{\tilde S_I}:Y\times\tilde S_I\to\tilde S_I$ is a lift of $g_I:=g_{|T_I}:T_I\to S_I$.
Let $M\in\DA_c(S)^-$ and $(F,W)\in C_{fil}(\Var(\mathbb C)^{sm}/S)$ such that $(M,W)=D(\mathbb A^1_S,et)(F,W)$.
Then, $D(\mathbb A^1_T,et)(g^*F)=g^*M$ and there exist 
$(F',W)\in C_{fil}(\Var(\mathbb C)^{sm}/S)$ and an equivalence $(\mathbb A^1,et)$ local $e:g^*(F,W)\to(F',W)$ 
such that $D(\mathbb A^1_T,et)(F',W)=(g^*M,W)$.Denote for short $d_{YI}:=-d_Y-d_{\tilde S_I}$.  
We have, using definition \ref{TgDR} and definition \ref{TGammaFDR}(i), by lemma \ref{Tglem}, 
the canonical map in $\pi_T(D(MHM(T)))\subset D_{\mathcal D(1,0)fil}(T/(Y\times\tilde S_I))$
\begin{eqnarray*}
T(g,\mathcal F^{FDR})(M):g^{\hat*mod}_{Hdg}\iota_S^{-1}\mathcal F_S^{FDR}(M):= \\
\Gamma^{\vee,Hdg}_T\iota_T^{-1}(\tilde g_I^{*mod}
(e'_*\mathcal Hom^{\bullet}(\hat R^{CH}(\rho_{\tilde S_I}^*(L(i_{I*}j_I^*(F,W))),   
E_{zar}(\Omega^{\bullet,\Gamma,pr}_{/\tilde S_I},F_{DR})))[d_{YI}],\tilde g_J^{*mod}u^q_{IJ}(F,W))) \\
\xrightarrow{(T(\tilde g_I,\Omega^{\Gamma,pr}_{/\cdot})(-))} \\
\Gamma^{\vee,Hdg}_T\iota_T^{-1}(e'_*\mathcal Hom^{\bullet}(\tilde g_I^*\hat R^{CH}(\rho_{\tilde S_I}^*L(i_{I*}j_I^*(F,W))), 
E_{zar}(\Omega^{\bullet,\Gamma,pr}_{/Y\times\tilde S_I},F_{DR}))[d_{YI}],\tilde g_J^*u^q_{IJ}(F,W)) \\
\xrightarrow{\mathcal Hom(T(\tilde g_I,\hat R^{CH})(-)^{-1},-)} \\ 
\Gamma^{\vee,Hdg}_T\iota_T^{-1}(e'_*\mathcal Hom^{\bullet}(\hat R^{CH}(\rho_{Y\times\tilde S_I}^*\tilde g_I^*L(i_{I*}j_I^*(F,W))), 
E_{zar}(\Omega^{\bullet,\Gamma,pr}_{/Y\times\tilde S_I},F_{DR}))[d_{YI}],\tilde g_J^*u^q_{IJ}(F,W)) \\  
\xrightarrow{T(\Gamma_T^{\vee,Hdg},\Omega^{\Gamma,pr}_{/S})(F,W)}\\
\iota_T^{-1}(e'_*\mathcal Hom^{\bullet}(\hat R^{CH}(\rho_{Y\times\tilde S_I}^*\Gamma^{\vee}_{T_I}\tilde g_I^*L(i_{I*}j_I^*(F,W))), 
E_{zar}(\Omega^{\bullet,\Gamma,pr}_{/Y\times\tilde S_I},F_{DR}))[d_{YI}],\tilde g_J^{*,\gamma}u^q_{IJ}(F,W)) \\ 
\xrightarrow{(\mathcal Hom(\hat R^{CH}_{Y\times\tilde S_I}(T^{q,\gamma}(D_{gI})(j_I^*(F,W))), 
E_{zar}(\Omega^{\bullet,\Gamma,pr}_{/Y\times\tilde S_I},F_{DR}))[d_{YI}])} \\
\iota_T^{-1}(e'_*\mathcal Hom^{\bullet}(\hat R^{CH}(\rho_{Y\times\tilde S_I}^*L(i'_{I*}j^{'*}_Ig^*(F,W))), 
E_{zar}(\Omega^{\bullet,\Gamma,pr}_{/Y\times\tilde S_I},F_{DR}))[d_{YI}],u^q_{IJ}(g^*(F,W))) \\ 
\xrightarrow{\mathcal Hom(\hat R^{CH}_{Y\times\tilde S_I}(Li'_{I*}j^{'*}_I(e)),-)} \\
\iota_T^{-1}(e'_*\mathcal Hom^{\bullet}(\hat R^{CH}(\rho_{Y\times\tilde S_I}^*L(i'_{I*}j^{'*}_I(F',W))), 
E_{zar}(\Omega^{\bullet,\Gamma,pr}_{/Y\times\tilde S_I},F_{DR}))[d_{YI}],u^q_{IJ}(F',W)) \\
\xrightarrow{=:}\mathcal F_T^{FDR}(g^*M)
\end{eqnarray*}
\end{defi}

\begin{defi}\label{SixTalg}
\begin{itemize}
\item Let $f:X\to S$ a morphism with $X,S\in\Var(\mathbb C)$. Assume there exist a factorization
$f:X\xrightarrow{l}Y\times S\xrightarrow{p_S}S$ with $Y\in\SmVar(\mathbb C)$, $l$ a closed embedding and $p_S$ the projection.
We have, for $M\in\DA_c(X)$, the following transformation map in $\pi_S(D(MHM(S)))$
\begin{eqnarray*}
T_*(f,\mathcal F^{FDR})(M):\mathcal F_S^{FDR}(Rf_*M)
\xrightarrow{\ad(f_{Hdg}^{\hat*mod},Rf^{Hdg}_*)(-)}Rf^{Hdg}_*f^{\hat*mod}_{Hdg}\mathcal F_S^{FDR}(Rf_*M) \\
\xrightarrow{T(f,\mathcal F^{FDR})(Rf_*M)}Rf^{Hdg}_*\mathcal F_X^{FDR}(f^*Rf_*M)
\xrightarrow{\mathcal F_X^{FDR}(\ad(f^*,Rf_*)(M))}Rf^{Hdg}_*\mathcal F_X^{FDR}(M)
\end{eqnarray*}
Clearly, for $p:Y\times S\to S$ a projection with $Y\in\PSmVar(\mathbb C)$, we have, for $M\in\DA_c(Y\times S)$,
$T_*(p,\mathcal F^{FDR})(M)=T_!(p,\mathcal F^{FDR})(M)[d_Y]$

\item Let $S\in\Var(\mathbb C)$. Let $Y\in\SmVar(\mathbb C)$ and $p:Y\times S\to S$ the projection.
We have then, for $M\in\DA(Y\times S)$ the following transformation map in $\pi_S(D(MHM(S)))$
\begin{eqnarray*}
T_!(p,\mathcal F^{FDR})(M):p^{Hdg}_!\mathcal F_{Y\times S}^{FDR}(M)
\xrightarrow{\mathcal F_{Y\times S}^{FDR}(\ad(Lp_{\sharp},p^*)(M))}
Rp^{Hdg}_!\mathcal F_{Y\times S}^{FDR}(p^*Lp_{\sharp}M) \\
\xrightarrow{T(p,\mathcal F^{FDR})(Lp_{\sharp}(M,W))}Rp^{Hdg}_!p^{\hat*mod[-]}\mathcal F_S^{FDR}(Lp_{\sharp}M)
\xrightarrow{T(p^{*mod},p^{\hat*mod})(-)}p^{Hdg}_!p^{*mod[-]} \\
\mathcal F_S^{FDR}(Lp_{\sharp}M)\xrightarrow{\ad(Rp^{Hdg}_!,p^{*mod[-]})(\mathcal F_S^{FDR}(Lp_{\sharp}M))}
\mathcal F_S^{FDR}(Lp_{\sharp}M)
\end{eqnarray*}

\item Let $f:X\to S$ a morphism with $X,S\in\Var(\mathbb C)$. Assume there exist a factorization
$f:X\xrightarrow{l}Y\times S\xrightarrow{p_S}S$ with $Y\in\SmVar(\mathbb C)$, $l$ a closed embedding and $p_S$ the projection.
We have then, using the second point, for $M\in\DA(X)$ the following transformation map in $\pi_S(D(MHM(S)))$
\begin{eqnarray*}
T_!(f,\mathcal F^{FDR})(M):
Rp^{Hdg}_!\mathcal F_X^{FDR}(M,W):=Rp^{Hdg}_!\mathcal F_{Y\times S}^{FDR}(l_*M) \\
\xrightarrow{T_!(p,\mathcal F^{FDR})(l_*M)}\mathcal F_S^{FDR}(Lp_{\sharp}l_*M)
\xrightarrow{=}\mathcal F_S^{FDR}(Rf_!M)
\end{eqnarray*}

\item Let $f:X\to S$ a morphism with $X,S\in\Var(\mathbb C)$. Assume there exist a factorization
$f:X\xrightarrow{l}Y\times S\xrightarrow{p_S}S$ with $Y\in\SmVar(\mathbb C)$, $l$ a closed embedding and $p_S$ the projection.
We have, using the third point, for $M\in\DA(S)$, the following transformation map in in $\pi_X(D(MHM(X)))$
\begin{eqnarray*}
T^!(f,\mathcal F^{FDR})(M):\mathcal F_X^{FDR}(f^!M)
\xrightarrow{\ad(Rf^{Hdg}_!,f^{*mod}_{Hdg})(\mathcal F_X^{FDR}(f^!M))}
f^{*mod}_{Hdg}Rf^{Hdg}_!\mathcal F_X^{FDR}(f^!M) \\
\xrightarrow{T_!(p_S,\mathcal F^{FDR})(\mathcal F^{FDR}(f^!M))}f^{*mod}_{Hdg}\mathcal F_S^{FDR}(Rf_!f^!M)
\xrightarrow{\mathcal F_S^{FDR}(\ad(Rf_!,f^!)(M))}f^{*mod}_{Hdg}\mathcal F_S^{FDR}(M)
\end{eqnarray*}

\item Let $S\in\Var(\mathbb C)$. Let $S=\cup_{i=1}^l S_i$ an open cover such that there exist closed embeddings 
$i_i:S_i\hookrightarrow\tilde S_i$ with $\tilde S_i\in\SmVar(\mathbb C)$. 
We have, using the preceding point, 
denoting $\Delta_S:S\hookrightarrow S$ the diagonal closed embedding 
and $p_1:S\times S\to S$, $p_2:S\times S\to S$ the projections, 
for $M,N\in\DA(S)$ and $(F,W),(G,W)\in C_{fil}(\Var(\mathbb C)^{sm}/S))$ 
such that $(M,W)=D(\mathbb A^1,et)(F,W)$ and $(N,W)=D(\mathbb A^1,et)(G,W)$, 
the following transformation map in $\pi_S(D(MHM(S)))$
\begin{eqnarray*}
T(\mathcal F_S^{FDR},\otimes)(M,N):\mathcal F_S^{FDR}(M)\otimes^{Hdg}_{O_S}\mathcal F_S^{FDR}(N) 
:=\Delta_S^{!Hdg}(p_1^{*mod}\mathcal F_S^{FDR}(M)\otimes_{O_{S\times S}}p_2^{*mod}\mathcal F_S^{FDR}(N) \\
\xrightarrow{T^!(p_1,\mathcal F_S^{FDR})(M)\otimes T^!(p_1,\mathcal F_S^{FDR})(M)}
\Delta_S^{!Hdg}(\mathcal F_{S\times S}^{FDR}(p_1^!M)\otimes_{O_{S\times S}}\mathcal F_{S\times S}^{FDR}(p_2^!N) \\
\xrightarrow{(T(\otimes,\Omega)(\hat R^{CH}(\rho_{\tilde S_I\times\tilde S_J}^*L(i_I\times i_J)_*(j_I\times j_J)^*p_1^*F[2d_S]),
\hat R^{CH}(\rho_{\tilde S_I\times\tilde S_J}^*L(i_I\times i_J)_*(j_I\times j_J)^*p_2^*F[2d_S])))} \\
\Delta_S^{!Hdg}(\mathcal F_{S\times S}^{FDR}(p_1^!M\otimes p_2^!N)
\xrightarrow{T^!(\Delta_S,\mathcal F^{FDR})(p_1^!M\otimes p_2^!N)}
\mathcal F_S^{FDR}(\Delta_S^!(p_1^!M\otimes p_2^!N))=\mathcal F_S^{FDR}(M\otimes N)
\end{eqnarray*}
where the last equality follows from the equality in $\DA(S)$
\begin{equation*}
\Delta_S^!(p_1^!M\otimes p_2^!N)=\Delta_S^!p_1^!M\otimes\Delta_S^!p_2^!N=M\otimes N
\end{equation*}
\end{itemize}
\end{defi}

\begin{prop}\label{Tgprop}
Let $g:T\to S$ a morphism with $T,S\in\Var(\mathbb C)$.
Assume we have a factorization $g:T\xrightarrow{l}Y_2\times S\xrightarrow{p_S}S$
with $Y_2\in\SmVar(\mathbb C)$, $l$ a closed embedding and $p_S$ the projection.
Let $S=\cup_{i=1}^lS_i$ be an open cover such that 
there exists closed embeddings $i_i:S_i\hookrightarrow\tilde S_i$ with $\tilde S_i\in\SmVar(\mathbb C)$
Then, $T=\cup^l_{i=1} T_i$ with $T_i:=g^{-1}(S_i)$
and we have closed embeddings $i'_i:=i_i\circ l:T_i\hookrightarrow Y_2\times\tilde S_i$,
Moreover $\tilde g_I:=p_{\tilde S_I}:Y\times\tilde S_I\to\tilde S_I$ is a lift of $g_I:=g_{|T_I}:T_I\to S_I$.
Let $f:X\to S$ a  morphism with $X\in\Var(\mathbb C)$ such that there exists a factorization 
$f:X\xrightarrow{l}Y_1\times S\xrightarrow{p_S} S$, with $Y_1\in\SmVar(\mathbb C)$, 
$l$ a closed embedding and $p_S$ the projection. We have then the following commutative diagram
whose squares are cartesians
\begin{equation*}
\xymatrix{f':X_T\ar[r]\ar[rd]\ar[dd]_{g'} & Y_1\times T\ar[rd]\ar[r] & T\ar[rd]\ar[dd]^{g} & \, \\
\, & Y_1\times X\ar[r]\ar[ld] & Y_1\times Y_2\times S\ar[r]\ar[ld] & Y_2\times S\ar[ld] \\
f:X\ar[r] & Y_1\times S\ar[r] & S & \,}
\end{equation*} 
Take a smooth compactification $\bar Y_1\in\PSmVar(\mathbb C)$ of $Y_1$, denote
$\bar X_I\subset\bar Y_1\times\tilde S_I$ the closure of $X_I$, and $Z_I:=\bar X_I\backslash X_I$.
Consider $F(X/S):=p_{S,\sharp}\Gamma_X^{\vee}\mathbb Z(Y_1\times S/Y_1\times S)\in C(\Var(\mathbb C)^{sm}/S)$ and
the isomorphism in $C(\Var(\mathbb C)^{sm}/T)$
\begin{eqnarray*}
T(f,g,F(X/S)):g^*F(X/S):=g^*p_{S,\sharp}\Gamma_X^{\vee}\mathbb Z(Y_1\times S/Y_1\times S)\xrightarrow{\sim} \\
p_{T,\sharp}\Gamma_{X_T}^{\vee}\mathbb Z(Y_1\times T/Y_1\times T)=:F(X_T/T).
\end{eqnarray*}
which gives in $\DA(T)$ the isomorphism $T(f,g,F(X/S)):g^*M(X/S)\xrightarrow{\sim}(X_T/T)$.
Then the following diagram in $\pi_T(D(MHM(T)))\subset D_{\mathcal D(1,0)fil}(T/(Y_2\times\tilde S_I))$, 
where the horizontal maps are given by proposition \ref{keyalgsing1}, commutes
\begin{equation*}
\begin{tikzcd}
g^{\hat*mod}_{Hdg}\iota_S^{-1}\mathcal F_S^{FDR}(M(X/S))
\ar[dd,"'T(g{,}\mathcal F^{FDR})(M(X/S))"']\ar[rr,"g^{\hat*mod}_{Hdg}I(X/S)"] & \, &
g^{\hat*mod}_{Hdg}Rf^{Hdg}_!(\Gamma^{\vee,Hdg}_{X_I}(O_{Y_1\times\tilde S_I},F_b)(d_{Y_1})[2d_{Y_1}],x_{IJ}(X/S))
\ar[d,"T(p_{\tilde S_I}{,}\gamma^{\vee{,}Hdg})(-)"] \\
\, & \, & Rf^{'Hdg}_!g^{'\hat*mod}_{Hdg}(\Gamma^{\vee,Hdg}_{X_I}(O_{Y_1\times\tilde S_I},F_b)(d_{Y_1})[2d_{Y_1}],x_{IJ}(X/S))
\ar[d,"T(p_{Y_1\times Y_2\times\tilde S_I{,}Hdg}^{\hat*mod}{,}p_{Y_1\times Y_2\times\tilde S_I{,}Hdg}^{*mod})(-)"] \\
\iota_T^{-1}\mathcal F_T^{FDR}(M(X_T/T))\ar[rr,"I(X_T/T)"] & \, & 
Rf^{'Hdg}_!(\Gamma^{\vee,Hdg}_{X_{T_I}}(O_{Y_2\times Y_1\times\tilde S_I},F_b)(d_{Y_{12}})[2d_{Y_{12}}],x_{IJ}(X_T/T)).
\end{tikzcd}
\end{equation*} 
with $d_{Y_{12}}=d_{Y_1}+d_{Y_2}$.
\end{prop}

\begin{proof}
Follows immediately from definition.
\end{proof}

\begin{prop}\label{mainthmprop2}
Let $S\in\Var(\mathbb C)$. Let $Y\in\SmVar(\mathbb C)$ and $p:Y\times S\to S$ the projection.
Let $S=\cup_{i=1}^l S_i$ an open cover such that there exist closed embeddings 
$i^o_i:S_i\hookrightarrow\tilde S_i$ with $\tilde S_i\in\SmVar(\mathbb C)$.
For $I\subset\left[1,\cdots l\right]$, we denote by $S_I=\cap_{i\in I} S_i$, $j^o_I:S_I\hookrightarrow S$ and
$j_I:Y\times S_I\hookrightarrow Y\times S$ the open embeddings. 
We then have closed embeddings $i_I:Y\times S_I\hookrightarrow Y\times\tilde S_I$.
and we denote by $p_{\tilde S_I}:Y\times\tilde S_I\to\tilde S_I$ the projections.
Let $f':X'\to Y\times S$ a morphism, with $X'\in\Var(\mathbb C)$ such that there exists a factorization
$f':X'\xrightarrow{l'}Y'\times Y\times S\xrightarrow{p'} Y\times S$ 
with $Y'\in\SmVar(\mathbb C)$, $l'$ a closed embedding and $p'$ the projection. 
Denoting $X'_I:=f^{'-1}(Y\times S_I)$, we have closed embeddings $i'_I:X'_I\hookrightarrow Y'\times Y\times\tilde S_I$
Consider 
\begin{eqnarray*}
F(X'/Y\times S):=p_{Y\times S,\sharp}\Gamma_{X'}^{\vee}\mathbb Z(Y'\times Y\times S/Y'\times Y\times S)
\in C(\Var(\mathbb C)^{sm}/Y\times S)
\end{eqnarray*}
and $F(X'/S):=p_{\sharp}F(X'/Y\times S)\in C(\Var(\mathbb C)^{sm}/S)$, 
so that $Lp_{\sharp}M(X'/Y\times S)[-2d_Y]=:M(X'/S)$. 
Then, the following diagram in $\pi_S(D(MHM(S)))\subset D_{\mathcal D(1,0)fil}(S/(Y\times\tilde S_I))$,
where the vertical maps are given by proposition \ref{keyalgsing1}, commutes
\begin{equation*}
\xymatrix{
Rp^{Hdg}!\mathcal F_{Y\times S}^{FDR}(M(X'/Y\times S))
\ar[rr]^{T_!(p,\mathcal F^{FDR})(M(X'/Y\times S))} & \, & 
\mathcal F_S^{FDR}(M(X'/S)) \\
Rp^{Hdg}!Rf^{'Hdg}_!f^{'*mod}_{Hdg}\mathbb Z_{Y\times S}^{Hdg} 
\ar[rr]^{=}\ar[u]^{T(p_{Hdg}^{\hat*mod}{,}p_{Hdg}^{*mod})(-)\circ Rp^{Hdg}!(I(X'/Y\times S))} & \, &
Rf^{Hdg}_!f^{*mod}_{Hdg}\mathbb Z_S^{Hdg}\ar[u]_{I(X'/S)}}. 
\end{equation*}
\end{prop}

\begin{proof}
Immediate from definition.
\end{proof}

\begin{prop}\label{TotimesDRprop}
Let $f_1:X_1\to S$, $f_2:X_2\to S$ two morphism with $X_1,X_2,S\in\Var(\mathbb C)$. 
Assume that there exist factorizations 
$f_1:X_1\xrightarrow{l_1}Y_1\times S\xrightarrow{p_S} S$, $f_2:X_2\xrightarrow{l_2}Y_2\times S\xrightarrow{p_S} S$
with $Y_1,Y_2\in\SmVar(\mathbb C)$, $l_1,l_2$ closed embeddings and $p_S$ the projections.
We have then the factorization
\begin{equation*}
f_{12}:=f_1\times f_2:X_{12}:=X_1\times_S X_2\xrightarrow{l_1\times l_2}Y_1\times Y_2\times S\xrightarrow{p_S} S
\end{equation*}
Let $S=\cup_{i=1}^l S_i$ an open affine covering and denote, 
for $I\subset\left[1,\cdots l\right]$, $S_I=\cap_{i\in I} S_i$ and $j_I:S_I\hookrightarrow S$ the open embedding.
Let $i_i:S_i\hookrightarrow\tilde S_i$ closed embeddings, with $\tilde S_i\in\SmVar(\mathbb C)$. 
We have then the following commutative diagram in $\pi_S(DMHM(S))\subset D_{\mathcal D(1,0)fil}(S/(\tilde S_I))$
where the vertical maps are given by proposition \ref{keyalgsing1}
\begin{equation*}
\begin{tikzcd}
\mathcal F_S^{FDR}(M(X_1/S))\otimes^{L[-]}_{O_S}\mathcal F_S^{FDR}(M(X_2/S))
\ar[rr,"I(X_1/S)\otimes I(X_2/S)"] \ar[d,"T(\mathcal F_S^{FDR}{,}\otimes)(M(X_1/S){,}M(X_2/S))"]  & \, &
\shortstack{$Rf_{1!}^{Hdg}(\Gamma^{\vee,Hdg}_{X_{1I}}(O_{Y_1\times\tilde S_I},F_b)(d_2)[2d_1],x_{IJ}(X_1/S))\otimes_{O_S}$ \\ 
$Rf_{2!}^{Hdg}(\Gamma^{\vee,Hdg}_{X_{2I}}(O_{Y_2\times\tilde S_I},F_b)(d_1)[2d_2],x_{IJ}(X_2/S))$}
\ar[d,"(Ew_{(Y_1\times\tilde S_I{,}Y_2\times\tilde S_I)/\tilde S_I})"] \\
\mathcal F_S^{FDR}(M(X_1/S)\otimes M(X_2/S)=M(X_1\times_S X_2/S))\ar[rr,"I(X_{12}/S)"] & \, &
Rf_{12!}^{Hdg}(\Gamma^{\vee,Hdg}_{X_{1I}\times_S X_{2I}}(O_{Y_1\times Y_2\times\tilde S_I},F_b)(d_{12})[2d_{12}],x_{IJ}(X_1/S)).
\end{tikzcd}
\end{equation*}
with $d_1=d_{Y_1}$, $d_2=d_{Y_2}$ and $d_{12}=d_{Y_1}+d_{Y_2}$.
\end{prop}

\begin{proof}
Immediate from definition.
\end{proof}

\begin{thm}\label{mainthm}
\begin{itemize}
\item[(i)]Let $g:T\to S$ a morphism, with $S,T\in\Var(\mathbb C)$. 
Assume we have a factorization $g:T\xrightarrow{l}Y\times S\xrightarrow{p_S}S$
with $Y\in\SmVar(\mathbb C)$, $l$ a closed embedding and $p_S$ the projection. Let $M\in\DA_c(S)$. 
Then map in $\pi_T(D(MHM(T)))$ 
\begin{eqnarray*}
T(g,\mathcal F^{FDR})(M):g^{\hat*mod}_{Hdg}\mathcal F_S^{FDR}(M)\xrightarrow{\sim}\mathcal F_T^{FDR}(g^*M)
\end{eqnarray*} 
given in definition \ref{TgDRdefsing} is an isomorphism.
\item[(ii)] Let $f:X\to S$ a morphism with $X,S\in\Var(\mathbb C)$. Assume there exist a factorization
$f:X\xrightarrow{l}Y\times S\xrightarrow{p_S}S$ with $Y\in\SmVar(\mathbb C)$, $l$ a closed embedding and $p_S$ the projection.
Then, for $M\in\DA_c(X)$, the map given in definition \ref{SixTalg}
\begin{equation*}
T_!(f,\mathcal F^{FDR})(M):Rf^{Hdg}_!\mathcal F_X^{FDR}(M)\xrightarrow{\sim}\mathcal F_S^{FDR}(Rf_!M)
\end{equation*}
is an isomorphism in $\pi_S(D(MHM(S))$.
\item[(iii)] Let $f:X\to S$ a morphism with $X,S\in\Var(\mathbb C)$, $S$ quasi-projective. Assume there exist a factorization
$f:X\xrightarrow{l}Y\times S\xrightarrow{p_S}S$ with $Y\in\SmVar(\mathbb C)$, $l$ a closed embedding and $p_S$ the projection.
We have, for $M\in\DA_c(X)$, the map given in definition \ref{SixTalg}
\begin{equation*}
T_*(f,\mathcal F^{FDR})(M):\mathcal F_S^{FDR}(Rf_*M)\xrightarrow{\sim}Rf^{Hdg}_*\mathcal F_X^{FDR}(M)
\end{equation*}
is an isomorphism in $\pi_S(D(MHM(S))$.
\item[(iv)] Let $f:X\to S$ a morphism with $X,S\in\Var(\mathbb C)$, $S$ quasi-projective. Assume there exist a factorization
$f:X\xrightarrow{l}Y\times S\xrightarrow{p_S}S$ with $Y\in\SmVar(\mathbb C)$, $l$ a closed embedding and $p_S$ the projection.
Then, for $M\in\DA_c(S)$, the map given in definition \ref{SixTalg}
\begin{equation*}
T^!(f,\mathcal F^{FDR})(M):\mathcal F_X^{FDR}(f^!M)\xrightarrow{\sim}f^{*mod}_{Hdg}\mathcal F_S^{FDR}(M)
\end{equation*}
is an isomorphism in $\pi_X(D(MHM(X))$. 
\item[(v)]Let $S\in\Var(\mathbb C)$ and $S=\cup_{i=1}^l S_i$ an open affine covering and denote, 
for $I\subset\left[1,\cdots l\right]$, $S_I=\cap_{i\in I} S_i$ and $j_I:S_I\hookrightarrow S$ the open embedding.
Let $i_i:S_i\hookrightarrow\tilde S_i$ closed embeddings, with $\tilde S_i\in\SmVar(\mathbb C)$. 
Then, for $M,N\in\DA_c(S)$, the map in $\pi_S(D(MHM(S)))$ 
\begin{eqnarray*}
T(\mathcal F_S^{FDR},\otimes)(M,N): 
\mathcal F_S^{FDR}(M)\otimes^{Hdg}_{O_S}\mathcal F_S^{FDR}(N)\xrightarrow{\sim}\mathcal F_S^{FDR}(M\otimes N)
\end{eqnarray*}
given in definition \ref{SixTalg} is an isomorphism.
\end{itemize}
\end{thm}

\begin{proof}
\noindent(i):Follows from proposition \ref{Tgprop} and proposition \ref{keyalgsing1}.

\noindent(ii): By proposition \ref{mainthmprop2} and proposition \ref{keyalgsing1},
for $S\in\Var(\mathbb C)$, $Y\in\SmVar(\mathbb C)$, $p:Y\times S\to S$ the projection and $M\in\DA_c(Y\times S)$, 
\begin{equation*}
T_!(p,\mathcal F^{FDR})(M):Rp^{Hdg}_!\mathcal F_{Y\times S}^{FDR}(M)\to\mathcal F_S^{FDR}(Rp_!M)
\end{equation*}
is an isomorphism.

\noindent(iii): 
Consider first an open embedding $n:S^o\hookrightarrow S$ with $S\in\Var(\mathbb C)$ so that there exist
a closed embedding $i:S\hookrightarrow\tilde S$ with $\tilde S\in\SmVar(\mathbb C)$. Then, since 
\begin{equation*}
n^*:C(\Var(\mathbb C)^{sm}/S)\to C(\Var(\mathbb C)^{sm}/S^o)
\end{equation*}
is surjective, $n^*:\DA(S)\to\DA(S^o)$ is surjective. 
Denote by $i:Z=S\backslash S^o\hookrightarrow S$ the complementary closed embedding.
By \cite{AyoubB}, $\DA_c(S)$ is generated by motives of the form 
\begin{eqnarray*}
\DA_c(S)=<M(X'/S)=f'_*E(\mathbb Z_{X'}), \; f':X'\to S \; \mbox{proper with} \; X'\in\SmVar(\mathbb C), \\ 
\mbox{s.t.} \; f^{'-1}(Z)=X' \; \mbox{or} \; f^{'-1}(Z)=\cup D_i=D\subset X'>. 
\end{eqnarray*}
If $f^{'-1}(Z)=X'$, $n^*M(X'/S)=0$. So let consider the case $f^{'-1}(Z)=\cup_{i=1}^l D_i=D\subset X'$ is a normal crossing divisor.
Denote $f'_D:f'_{|D}:D\to Z$, $D_I=\cap_{i\in I} D_i$ and $i'_I:D_I\hookrightarrow X'$,
$n':X^{'o}:=X'\backslash D\hookrightarrow X'$ the complementary open embedding and $f^{'o}:f'_{|X^{'o}}:X^{'o}\to S^o$.
Denote $L=[1,\ldots,l]$. We have then a generalized distinguish triangle in $\DA(X')$
\begin{eqnarray*}
a(n',i'):n'_*n^{'*}E_{et}(\mathbb Z_{X'})&\xrightarrow{\sim}&
\Cone(\gamma_D(-):\Gamma_DE_{et}(\mathbb Z_{X'})\to E_{et}(\mathbb Z_{X'})) \\
&\xrightarrow{\sim}&
\Cone(\Gamma_{D_L}E_{et}(\mathbb Z_{X'})\to\cdots\to\bigoplus_{i=1}^lE_{et}(\mathbb Z_{X'}) \\
&\xrightarrow{\sim}&
\Cone(i'_{L*}i_L^{'!}E_{et}(\mathbb Z_{X'})\to\cdots\to\bigoplus_{i=1}^li'_{i*}i_i^{'!}E_{et}(\mathbb Z_{X'})
\xrightarrow{\ad(i'_{i*},i_i^{'!})(E_{et}(\mathbb Z_{X'}))}E_{et}(\mathbb Z_{X'})) \\
&\xrightarrow{\sim}&
\Cone(i'_{L*}\mathbb Z_{D_{L}}[-l]\to\cdots\bigoplus_{i=1}^li'_{i*}\mathbb Z_{D_i}[-1]\to\mathbb Z_{X'})
\end{eqnarray*}
where the first isomorphism is the image of an homotopy equivalence by definition, 
the second one is the image of an homotopy equivalence by definition-proposition \ref{gamma2sect2}(ii), 
the third one follows by the localization property (see section 3, theorem \ref{2functDM})
and the last one follows from purity since the $D_I$ and $X'$ are smooth (see section 3, theorem \ref{2functDM}).
Similarly, we have a generalized distinguish triangle in $D(MHM(X'))$
\begin{eqnarray*}
a^{mod}(n',i'):n^{'Hdg}_*n^{'*}(O_{X'},F_b)&\xrightarrow{\sim}&
\Cone(\gamma^{Hdg}_D(O_{X'},F_b):\Gamma^{Hdg}_D(O_{X'},F_b)\to(O_{X'},F_b)) \\
&\xrightarrow{\sim}&
\Cone(\Gamma^{Hdg}_{D_L}(O_{X'},F_b)\to\cdots\to\bigoplus_{i=1}^l\Gamma^{Hdg}_{D_i}(O_{X'},F_b)
\xrightarrow{\oplus_i\gamma^{Hdg}_{D_i}(-)}(O_{X'},F_b)) \\
&\xrightarrow{\sim}&
\Cone(i'_{L*mod}(O_{D_L},F_b)[-l]\to\cdots\to\bigoplus_{i=1}^li'_{i*mod}(O_{D_i},F_b)[-1] \\
& \; & \xrightarrow{\ad(i'_{i*mod},i_i^{'\sharp})(O_{X'},F_b)}\cdots\to(O_{X'},F_b))
\end{eqnarray*}
where the first isomorphism is the image of an homotopy equivalence by definition, 
the second one is the image of an homotopy equivalence by definition-proposition \ref{TgammaHdg}, 
and the third one follows by the localization property of mixed Hodge modules (see section 5).
Consider $n^*M(X'/S)=M(X^{'o}/S^o)$. We have then the following commutative diagram in $\pi_S(MHM(S))$
\begin{equation*}
\xymatrix{\mathcal F_S^{FDR}(Rn_*n^*M(X'/S))=\mathcal F_S^{FDR}(Rn_*n^*Rf'_*\mathbb Z_{X'})
\ar[r]^{T_*(n,\mathcal F^{FDR})(n^*M(X'/S))}\ar[d]_{\mathcal F_S^{FDR}(Rn_*T(n,f')(\mathbb Z_{X'}))}  &
n_{*Hdg}\mathcal F_{S^o}^{FDR}(n^*M(X'/S))=n_{*Hdg}\mathcal F_{S^o}^{FDR}(n^*Rf'_*\mathbb Z_{X'}) 
\ar[d]^{n_{*Hdg}\mathcal F_{S^o}^{FDR}(T(n,f')(\mathbb Z_{X'}))} \\
\mathcal F_S^{FDR}(Rn_*Rf^{'o}_*n^{'*}\mathbb Z_{X'})=\mathcal F_S^{FDR}(Rf'_*Rn'_*n^{'*}\mathbb Z_{X'})
\ar[r]^{T_*(n,\mathcal F^{FDR})(M(X^{'o}/S))}\ar[d]_{T_*(f',\mathcal F^{FDR})(Rn'_*n^{'*}\mathbb Z_{X'})}  &
n_{*Hdg}\mathcal F_{S^o}^{FDR}(Rf^{'o}_*n^{'*}\mathbb Z_{X'})
\ar[d]^{n_{*Hdg}T_*(f^{'o},\mathcal F^{FDR})(n^{'*}\mathbb Z_{X'})} \\
Rf'_{*Hdg}\mathcal F_{X'}^{FDR}(Rn'_*n^{'*}\mathbb Z_{X'})
\ar[r]^{T_*(n,\mathcal F^{FDR})(n^{'*}\mathbb Z_{X'})}\ar[d]_{Rf'_{*Hdg}\mathcal F_{X'}^{FDR}(a(n',i'))}  &
n_{*Hdg}Rf^{'o}_{*Hdg}\mathcal F_{X^{'o}}^{FDR}(n^{'*}\mathbb Z_{X'})=Rf'_{*Hdg}n'_{*Hdg}n^{'*}\mathbb Z^{Hdg}_{X'})
\ar[d]^{Rf'_{*Hdg}a^{mod}(n',i')} \\
Rf^{'Hdg}_*\mathcal F_{X'}^{FDR}(\Cone(i'_{L*}\mathbb Z_{D_{L}}[-l]\to\cdots\to\mathbb Z_{X'}))
\ar[r]^{(T_*(i'_I,\mathcal F^{FDR})(-)\circ T(i'_I,\mathcal F^{FDR})(-))}  &
Rf^{'Hdg}_*\Cone(i'_{L*mod}(O_{D_L},F_b)[-l]\to\cdots\to(O_{X'},F_b))}
\end{equation*}
Since all the morphism involved are isomorphisms, $T_*(n,\mathcal F^{FDR})(n^*M(X'/S))$ is an isomorphism.
Hence, $T_*(n,\mathcal F^{FDR})(M)$ is an isomorphism for all $M\in\DA(S^o)$. 
Consider now the case of a general morphism $f:X\to S$, $X,S\in\Var(\mathbb C)$, $S$ quasi-projective,
which factors trough $f:X\xrightarrow{l}Y\times S\xrightarrow{p_S}S$ with some $Y\in\SmVar(\mathbb C)$.
By definition, for $M\in\DA_c(X)$ 
\begin{eqnarray*}
T_*(f,\mathcal F^{FDR})(M):\mathcal F_S^{FDR}(Rf_*M)=\mathcal F_S^{FDR}(Rp_{S*}l_*M) \\ 
\xrightarrow{T_*(p_S,\mathcal F^{FDR})(l_*M)} 
Rp^{Hdg}_{S*}\mathcal F_{Y\times S}^{FDR}(l_*M)=:Rf^{Hdg}_*\mathcal F_X^{FDR}(M).
\end{eqnarray*}
Hence, we have to show that for $S\in\Var(\mathbb C)$, $Y\in\SmVar(\mathbb C)$, $p:Y\times S\to S$ the projection,
and $M\in\DA_c(Y\times S)$, 
\begin{equation*}
T_*(p,\mathcal F^{FDR})(M):\mathcal F_S^{FDR}(Rp_*M)\to Rp^{Hdg}_*\mathcal F_{Y\times S}^{FDR}(M)
\end{equation*}
is an isomorphism.
Take a smooth compactification $\bar Y\in\PSmVar(\mathbb C)$ of $Y$. Denote
by $n_0:Y\hookrightarrow\bar Y$ and $n:=n_0\times I_S:Y\times S\hookrightarrow\bar Y\times S$ the open embeddings
and by $\bar p:\bar Y\times S\to S$ the projection. We have $p=\bar p\circ n:Y\times S\to S$, which gives the factorization
\begin{eqnarray*}
T_*(p,\mathcal F^{FDR})(M):\mathcal F_S^{FDR}(Rp_*M)=\mathcal F_S^{FDR}(R\bar p_*Rn_*M)
\xrightarrow{T_*(\bar p,\mathcal F^{FDR})(Rn_*M)} R\bar{p}^{Hdg}_*\mathcal F_{\bar Y\times S}^{FDR}(Rn_*M) \\
\xrightarrow{T_*(n,\mathcal F^{FDR})(M)} 
R\bar{p}^{Hdg}_*n^{Hdg}_*\mathcal F_{Y\times S}^{FDR}(M)=Rp^{Hdg}_*\mathcal F_{Y\times S}^{FDR}(M).
\end{eqnarray*}
By the open embedding case $T_*(n,\mathcal F^{FDR})(M)$ is an isomorphism.
On the other hand, since $\bar p$ is proper, 
$T_*(\bar p,\mathcal F^{FDR})(Rn_*M)=T_!(\bar p,\mathcal F^{FDR})(Rn_*M)$ is an isomorphism by (i).

\noindent(iv): Denote by $n:Y\times S\backslash X\hookrightarrow Y\times S$ the complementary open embedding.
We have, for $M\in\DA_c(S)$, the factorization
\begin{eqnarray*}
T^!(f,\mathcal F^{FDR})(M):\mathcal F_X^{FDR}(f^!M)=\mathcal F_{Y\times S}^{FDR}(l_*l^!p_S^!M) 
\xrightarrow{\mathcal F_{Y\times S}^{FDR}(a(n,l))} 
\mathcal F_{Y\times S}^{FDR}(\Cone(p_S^!M\to Rn_*n^*p_S^!M)[-1]) \\
\xrightarrow{=} 
\Cone(\mathcal F_{Y\times S}^{FDR}(p_S^!M)\to\mathcal F_{Y\times S}^{FDR}(Rn_*n^*p_S^!M))[-1] \\
\xrightarrow{(T(n,\mathcal F^{FDR})(p_S^!M)\circ T^!(p_S,\mathcal F^{FDR}(M)),T^!(p_S,\mathcal F^{FDR}(M)))} \\
\Cone(p_S^{*mod[-]}\mathcal F_S^{FDR}(M)\to n^{Hdg}_*n^*p_S^{*mod[-]}\mathcal F_S^{FDR}(M))[-1]
\xrightarrow{\sim}f^{*mod}_{Hdg}\mathcal F_S^{FDR}(M).
\end{eqnarray*}
By (ii),$T(n,\mathcal F^{FDR})(p_S^!M)$ is an isomorphism.
On the other hand, since $p_S$ is a smooth morphism, 
$T^!(p_S,\mathcal F^{FDR}(M))=T(p_S,\mathbb D\mathcal F^{FDR}(M))[d_Y]$ ; 
hence, $T^!(p_S,\mathcal F^{FDR}(M))$ is an isomorphism by (i).
\end{proof}

\begin{lem}\label{TgDRexlem}
Let $g:T\to S$ a morphism with $T,S\in\Var(\mathbb C)$.
Assume we have a factorization $g:T\xrightarrow{l}Y\times S\xrightarrow{p_S}S$
with $Y\in\SmVar(\mathbb C)$, $l$ a closed embedding and $p_S$ the projection.
Let $S=\cup_{i=1}^lS_i$ be an open cover such that 
there exists closed embeddings $i_i:S_i\hookrightarrow\tilde S_i$ with $\tilde S_i\in\SmVar(\mathbb C)$
Then, $T=\cup^l_{i=1} T_i$ with $T_i:=g^{-1}(S_i)$
and we have closed embeddings $i'_i:=i_i\circ l:T_i\hookrightarrow Y\times\tilde S_i$,
Moreover $\tilde g_I:=p_{\tilde S_I}:Y\times\tilde S_I\to\tilde S_I$ is a lift of $g_I:=g_{|T_I}:T_I\to S_I$.
Let $M\in\DA_c(S)^{-}$ and $(F,W)\in C_{fil}(\Var(\mathbb C)^{sm}/S)$ such that $(M,W)=D(\mathbb A^1_S,et)(F,W)$.
Then, $g^!M=L\mathbb D_Sg^*L\mathbb D_SM$, $D(\mathbb A^1_T,et)(g^*\mathbb D_SLF)=g^*L\mathbb D_SM$
and there exist $(F',W)\in C_{fil}(\Var(\mathbb C)^{sm}/S)$ and an equivalence $(\mathbb A^1,et)$ local 
$e:(F',W)\to(g^*\mathbb D_SL(F,W))$ such that $D(\mathbb A^1_T,et)(F',W)=g^*L\mathbb D_S(M,W)$
and, using definition \ref{TgDR} and definition \ref{TGammaFDR}(ii) and lemma \ref{Tglem},
the map in $\pi_T(D(MHM(T)))\subset D_{\mathcal D(1,0)fil}(T/(Y\times\tilde S_I))$
\begin{equation*}
T^!(g,\mathcal F^{FDR})(M):\mathcal F_T^{FDR}(g^!M)\to g^{*mod}_{Hdg}\mathcal F_S^{FDR}(M) 
\end{equation*}
given in definition \ref{SixTalg} is the inverse of the following map
\begin{eqnarray*}
T^{!,-1}(g,\mathcal F^{FDR})(M):g^{*mod}_{Hdg}\iota_S^{-1}\mathcal F_S^{FDR}(M) 
\xrightarrow{:=} \\ 
(\Gamma^{Hdg}_T\iota_T^{-1}(\tilde g_I^{*mod}
(e'_*\mathcal Hom^{\bullet}(\hat R^{CH}(\rho_{\tilde S_I}^*L(i_{I*}j_I^*\mathbb D_SL(F,W))),  
E_{zar}(\Omega^{\bullet,\Gamma,pr}_{/\tilde S_I},F_{DR})))[-d_{YI}],\tilde g_J^{*mod}u^q_{IJ}(\mathbb D_SL(F,W))) \\  
\xrightarrow{(T(\tilde g_I,\Omega^{\Gamma,pr}_{/\cdot})(-))} \\
\Gamma^{Hdg}_T\iota_T^{-1}(e'_*\mathcal Hom^{\bullet}(
\tilde g_I^*\hat R^{CH}(\rho_{\tilde S_I}^*L(i_{I*}j_I^*\mathbb D_SL(F,W))), 
E_{zar}(\Omega^{\bullet,\Gamma,pr}_{/Y\times\tilde S_I},F_{DR})),\tilde g_J^*u^q_{IJ}(\mathbb D_SL(F,W))) \\
\xrightarrow{\mathcal Hom(T(\tilde g_I,\hat R^{CH})(-)^{-1},-)} \\
\Gamma^{Hdg}_T\iota_T^{-1}(e'_*\mathcal Hom^{\bullet}(
\hat R^{CH}(\rho_{Y\times\tilde S_I}^*\tilde g_I^*L(i_{I*}j_I^*\mathbb D_SL(F,W))), 
E_{zar}(\Omega^{\bullet,\Gamma,pr}_{/Y\times\tilde S_I},F_{DR}))[d_{YI}],\tilde g_J^*u^q_{IJ}\mathbb D_SL(F,W)) \\
\xrightarrow{T(\Gamma_T^{Hdg},\Omega^{\Gamma,pr}_{/S})(F,W)^{-1}} \\
\iota_T^{-1}(e'_*\mathcal Hom^{\bullet}(
\hat R^{CH}(\rho_{Y\times\tilde S_I}^*L\Gamma_{T_I}E(\tilde g_I^*L(i_{I*}j_I^*\mathbb D_SL(F,W))), 
E_{zar}(\Omega^{\bullet,\Gamma,pr}_{/Y\times\tilde S_I},F_{DR}))[d_{YI}],\tilde g_J^{*,\gamma,d}u^q_{IJ}(F,W)) \\
\xrightarrow{(\mathcal Hom(
\hat R^{CH}_{Y\times\tilde S_I}(\mathbb D_{Y\times\tilde S_I}T^{q,\gamma}(D_{gI})(j_I^*\mathbb D_SL(F,W))), 
E_{zar}(\Omega^{\bullet,\Gamma,pr}_{/Y\times\tilde S_I},F_{DR}))[d_{YI}])^{-1}} \\
\iota_T^{-1}(e'_*\mathcal Hom^{\bullet}(
\hat R^{CH}(\rho_{Y\times\tilde S_I}^*L\mathbb D_{Y\times\tilde S_I}L(i'_{I*}j^{'*}_Ig^*\mathbb D_SL(F,W))), 
E_{zar}(\Omega^{\bullet,\Gamma,pr}_{/Y\times\tilde S_I},F_{DR}))[d_{YI}],u^{q,d}_{IJ}(g^*\mathbb D_SL(F,W))) \\
\xrightarrow{\mathcal Hom(\hat R^{CH}_{Y\times\tilde S_I}(\mathbb D_{Y\times\tilde S_I}Li'_{I*}j_I^{'*}(e)),-))} \\
\iota_T^{-1}(e'_*\mathcal Hom(
\hat R^{CH}(\rho_{Y\times\tilde S_I}^*L\mathbb D_{Y\times\tilde S_I}L(i'_{I*}j^{'*}_Ig^*\mathbb D_SL(F,W))),  
E_{zar}(\Omega^{\bullet,\Gamma,pr}_{/Y\times\tilde S_I},F_{DR}))[d_{YI}],u^{q,d}_{IJ}(L(F',W))) \\
\xrightarrow{=:}\mathcal F_T^{FDR}(g^!M)
\end{eqnarray*}
\end{lem}

We have the following proposition :

\begin{prop}\label{TgDRGMdia}
Let $g:T\to S$ a morphism with $T,S\in\Var(\mathbb C)$.
Assume we have a factorization $g:T\xrightarrow{l}Y\times S\xrightarrow{p_S}S$
with $Y\in\SmVar(\mathbb C)$, $l$ a closed embedding and $p_S$ the projection.
Let $S=\cup_{i=1}^lS_i$ be an open cover such that 
there exists closed embeddings $i_i:S_i\hookrightarrow\tilde S_i$ with $\tilde S_i\in\SmVar(\mathbb C)$
Then, $T=\cup^l_{i=1} T_i$ with $T_i:=g^{-1}(S_i)$
and we have closed embeddings $i'_i:=i_i\circ l:T_i\hookrightarrow Y\times\tilde S_i$,
Moreover $\tilde g_I:=p_{\tilde S_I}:Y\times\tilde S_I\to\tilde S_I$ is a lift of $g_I:=g_{|T_I}:T_I\to S_I$.
Let $M\in\DA_c(S)$ and $F\in C(\Var(\mathbb C)^{sm}/S)$ such that $M=D(\mathbb A^1_S,et)(F)$.
Then, $D(\mathbb A^1_T,et)(g^*F)=g^*M$. 
Then the following diagram in $D_{Ofil,\mathcal D,\infty}(T/(Y\times\tilde S_I))$ commutes 
\begin{equation*}
\xymatrix{Rg^{*mod[-],\Gamma}\mathcal F_S^{GM}(L\mathbb D_SM)
\ar[d]^{T(g,\mathcal F^{GM})(L\mathbb D_SM)}\ar[rr]^{Rg^{*mod[-],\Gamma}T(\mathcal F_S^{GM},\mathcal F_S^{FDR})(M)} & \, &
Rg^{*mod[-],\Gamma}\mathcal F_S^{FDR}(M)\ar[rrd] & \, &
g^{*mod}_{Hdg}\mathcal F_S^{FDR}(M)
\ar[ll]_{T(g^{*mod}_{Hdg},Rg^{*mod[-],\Gamma})(\mathcal F_S^{FDR}(M))}\ar[d]^{T^!(g,\mathcal F^{FDR})(M)^{-1}} \\ 
\mathcal F_T^{GM}(g^*L\mathbb D_SM=L\mathbb D_Tg^!M)
\ar[rrrr]^{T(\mathcal F_T^{GM},\mathcal F_T^{FDR})(g^*M)} & \, & \, & \, &
\mathcal F_T^{FDR}(g^!M)}
\end{equation*}
\end{prop}

\begin{proof}
Follows from lemma \ref{TgDRexlem}.
\end{proof}

We have the following easy proposition

\begin{prop}
Let $S\in\Var(\mathbb C)$ and $S=\cup_{i=1}^l S_i$ an open affine covering and denote, 
for $I\subset\left[1,\cdots l\right]$, $S_I=\cap_{i\in I} S_i$ and $j_I:S_I\hookrightarrow S$ the open embedding.
Let $i_i:S_i\hookrightarrow\tilde S_i$ closed embeddings, with $\tilde S_i\in\SmVar(\mathbb C)$. 
We have, for $M,N\in\DA(S)$ and $F,G\in C(\Var(\mathbb C)^{sm}/S)$ such that 
$M=D(\mathbb A^1,et)(F)$ and $N=D(\mathbb A^1,et)(G)$, 
the following commutative diagram in $D_{O_Sfil,\mathcal D,\infty}(S/(\tilde S_I))$
\begin{equation*}
\xymatrix{\mathcal F_S^{GM}(L\mathbb D_SM)\otimes^L_{O_S}\mathcal F_S^{GM}(L\mathbb D_SN)
\ar[d]^{T(\mathcal F_S^{GM},\otimes)(L\mathbb D_SM,L\mathbb D_SN)}
\ar[rrrr]^{T(\mathcal F_S^{GM},\mathcal F_S^{FDR})(M)\otimes T(\mathcal F_S^{GM},\mathcal F_S^{FDR})(N)} & \, & \, & \, &
\mathcal F_S^{FDR}(M)\otimes^{Hdg}_{O_S}\mathcal F_S^{FDR}(N)\ar[d]^{T(\mathcal F_S^{FDR},\otimes)(M,N)} \\
\mathcal F_S^{GM}(L\mathbb D_S(M\otimes N))\ar[rrrr]^{T(\mathcal F_S^{GM},\mathcal F_S^{FDR})(M\otimes N)} & \, & \, & \, & 
\mathcal F_S^{FDR}(M\otimes N)}
\end{equation*}
\end{prop}

\begin{proof}
Immediate from definition.
\end{proof}

\subsection{The analytic filtered De Rahm realization functor}

On $\AnSp(\mathbb C)$ the usual topology is equivalent to the etale topology since
a morphism $r:U'\to U$ is etale (which means non ramified and flat, see section 2)
if and only if for all $x\in U'$ there exist an open neighborhood $U'_x\subset U$
of $x$ such that $r$ induces an isomorphism $r_{|U'_x}:U'_x\xrightarrow{\sim}r(U'_x)$. 
We note $\tau=et$ the etale topology.

\subsubsection{The analytic Gauss-Manin filtered De Rham realization functor and its transformation map with pullbacks}

Consider, for $S\in\AnSp(\mathbb C)$, the following composition of morphism in $\RCat$ (see section 2)
\begin{eqnarray*}
\tilde e(S):(\AnSp(\mathbb C)/S,O_{\AnSp(\mathbb C)/S})\xrightarrow{\rho_S}(\AnSp(\mathbb C)^{sm}/S,O_{\AnSp(\mathbb C)^{sm}/S})
\xrightarrow{e(S)}(S,O_S)
\end{eqnarray*}
with, for $X/S=(X,h)\in\AnSp(\mathbb C)/S$,
\begin{itemize}
\item $O_{\AnSp(\mathbb C)/S}(X/S):=O_X(X)$, 
\item $(\tilde e(S)^*O_S(X/S)\to O_{\AnSp(\mathbb C)/S}(X/S)):=(h^*O_S\to O_X)$.
\end{itemize}
and $O_{\AnSp(\mathbb C)^{sm}/S}:=\rho_{S*}O_{\AnSp(\mathbb C)/S}$, that is, 
for $U/S=(U,h)\in\AnSp(\mathbb C)^{sm}/S$, $O_{\AnSp(\mathbb C)^{sm}/S}(U/S):=O_{\AnSp(\mathbb C)/S}(U/S):=O_U(U)$

\begin{defi}\label{wtildewGMan}
\begin{itemize}
\item[(i)]For $S\in\Var(\mathbb C)$, we consider the complexes of presheaves 
\begin{equation*} 
\Omega^{\bullet}_{/S}:=
\coker(\Omega_{O_{\AnSp(\mathbb C)/S}/\tilde e(S)^*O_S}:
\Omega^{\bullet}_{\tilde e(S)^*O_S}\to\Omega^{\bullet}_{O_{\AnSp(\mathbb C)/S}})
\in C_{O_S}(\AnSp(\mathbb C)/S) 
\end{equation*}
which is by definition given by 
\begin{itemize}
\item for $X/S$ a morphism $\Omega^{\bullet}_{/S}(X/S)=\Omega^{\bullet}_{X/S}(X)$
\item for $g:X'/S\to X/S$ a morphism, 
\begin{eqnarray*}
\Omega^{\bullet}_{/S}(g):=\Omega_{(X'/X)/(S/S)}(X'):
\Omega^{\bullet}_{X/S}(X)\to g^*\Omega_{X/S}(X')\to\Omega^{\bullet}_{X'/S}(X') \\
\omega\mapsto\Omega_{(X'/X)/(S/S)}(X')(\omega):=g^*(\omega):(\alpha\in\wedge^{k}T_{X'}(X')\mapsto\omega(dg(\alpha)))
\end{eqnarray*}
\end{itemize}
\item[(ii)] For $S\in\AnSp(\mathbb C)$, we consider the complexes of presheaves 
\begin{equation*} 
\Omega^{\bullet}_{/S}:=\rho_{S*}\tilde\Omega^{\bullet}_{/S}=
\coker(\Omega_{O_{\AnSp(\mathbb C)^{sm}/S}/e(S)^*O_S}:\Omega^{\bullet}_{e(S)^*O_S}\to\Omega^{\bullet}_{O_{\AnSp(\mathbb C)^{sm}/S}})
\in C_{O_S}(\AnSp(\mathbb C)^{sm}/S) 
\end{equation*}
which is by definition given by 
\begin{itemize}
\item for $U/S$ a smooth morphism $\Omega^{\bullet}_{/S}(U/S)=\Omega^{\bullet}_{U/S}(U)$
\item for $g:U'/S\to U/S$ a morphism, 
\begin{eqnarray*}
\Omega^{\bullet}_{/S}(g):=\Omega_{(U'/U)/(S/S)}(U'):
\Omega^{\bullet}_{U/S}(U)\to g^*\Omega_{U/S}(U')\to\Omega^{\bullet}_{U'/S}(U') \\
\omega\mapsto\Omega_{(U'/U)/(S/S)}(U')(\omega):=g^*(\omega):(\alpha\in\wedge^{k}T_{U'}(U')\mapsto\omega(dg(\alpha)))
\end{eqnarray*}
\end{itemize}
\end{itemize}
\end{defi}

\begin{rem}
For $S\in\AnSp(\mathbb C)$, $\Omega^{\bullet}_{/S}\in C(\AnSp(\mathbb C)/S)$ 
is by definition a natural extension of $\Omega^{\bullet}_{/S}\in C(\AnSp(\mathbb C)^{sm}/S)$. 
However $\Omega^{\bullet}_{/S}\in C(\AnSp(\mathbb C)/S)$ does NOT satisfy cdh descent.
\end{rem}

For a smooth morphism $h:U\to S$ with $S,U\in\AnSm(\mathbb C)$, 
the cohomology presheaves  $H^n\Omega^{\bullet}_{U/S}$ of the relative De Rham complex
\begin{equation*}
DR(U/S):=\Omega^{\bullet}_{U/S}:=\coker(h^*\Omega_S\to\Omega_U)\in C_{h^*O_S}(U)
\end{equation*}
for all $n\in\mathbb Z$, have a canonical structure of a complex of $h^*D^{\infty}_S$ modules given by the Gauss Manin connexion : 
for $S^o\subset S$ an open subset, $U^o=h^{-1}(S^o)$, 
$\gamma\in\Gamma(S^o,T_S)$ a vector field and $\hat\omega\in\Omega^{p}_{U/S}(U^o)^c$ a closed form, the action is given by
\begin{equation*}
\gamma\cdot[\hat\omega]=[\widehat{\iota(\tilde\gamma)\partial\omega}], 
\end{equation*}
$\omega\in\Omega^p_U(U^o)$ being a representative of $\hat\omega$ and 
$\tilde\gamma\in\Gamma(U^o,T_U)$ a relevement of $\gamma$ ($h$ is a smooth morphism), so that 
\begin{equation*}
DR(U/S):=\Omega^{\bullet}_{U/S}:=\coker(h^*\Omega_S\to\Omega_U)\in C_{h^*O_S,h^*\mathcal D^{\infty}}(U)
\end{equation*}
with this $h^*D^{\infty}_S$ structure. 
Hence we get $h_*\Omega^{\bullet}_{U/S}\in C_{O_S,\mathcal D^{\infty}}(S)$ considering this structure.
Since $h$ is a smooth morphism, $\Omega^p_{U/S}$ are locally free $O_U$ modules.

The point (ii) of the definition \ref{wtildewan} above gives 
the object in $\DA(S)$ which will, for $S$ smooth, represent the analytic Gauss-Manin De Rham realisation.
It is the class of an explicit complex of presheaves on $\AnSp(\mathbb C)^{sm}/S$. 

\begin{prop}\label{ausufibGM}
Let $S\in\Var(\mathbb C)$.
\begin{itemize}
\item[(i)] For $U/S=(U,h)\in\AnSp(\mathbb C)^{sm}/S$, we have $e(U)_*h^*\Omega^{\bullet}_{/S}=\Omega^{\bullet}_{U/S}$.
\item[(ii)] The complex of presheaves $\Omega^{\bullet}_{/S}\in C_{O_S}(\AnSp(\mathbb C)^{sm}/S)$ is $\mathbb D^1$ homotopic.
Note that however, for $p>0$, the complexes of presheaves $\Omega^{\bullet\geq p}$ are NOT $\mathbb D^1_S$ local.
On the other hand, $\Omega^{\bullet}_{/S}$ admits transferts (recall that means $\Tr(S)_*\Tr(S)^*\Omega^p_{/S}=\Omega^p_{/S}$).
\item[(iii)] If $S$ is smooth, we get $(\Omega^{\bullet}_{/S},F_b)\in C_{O_Sfil,D^{\infty}_S}(\Var(\mathbb C)^{sm}/S)$
with the structure given by the Gauss Manin connexion. Note that however the $D^{\infty}_S$ structure on the cohomology groups
given by Gauss Main connexion does NOT comes from a structure of $D^{\infty}_S$ module structure on the filtered complex of $O_S$ module.
The $D_S$ structure on the cohomology groups satisfy a non trivial Griffitz transversality (in the non projection cases),
whereas the filtration on the complex is the trivial one.
\end{itemize}
\end{prop}

\begin{proof}
Similar to the proof of proposition \ref{aetfibGM}.
\end{proof}

We have the following canonical transformation map given by the pullback of (relative) differential forms:

Let $g:T\to S$ a morphism with $T,S\in\AnSp(\mathbb C)$.
Consider the following commutative diagram in $\RCat$ :
\begin{equation*}
D(g,e):\xymatrix{
(\AnSp(\mathbb C)^{sm}/T,O_{\AnSp(\mathbb C)^{sm}/T})\ar[rr]^{P(g)}\ar[d]^{e(T)} & \, & 
(\AnSp(\mathbb C)^{sm}/S,O_{\AnSp(\mathbb C)^{sm}/S})\ar[d]^{e(S)} \\
(T,O_T)\ar[rr]^{P(g)} & \, & (S,O_S)}
\end{equation*}
It gives (see section 2) the canonical morphism in $C_{g^*O_Sfil}(\Var(\mathbb C)^{sm}/T)$ 
\begin{eqnarray*}
\Omega_{/(T/S)}:=\Omega_{(O_{\AnSp(\mathbb C)^{sm}/T}/g^*O_{\AnSp(\mathbb C)^{sm}/S})/(O_T/g^*O_S}): \\
g^*(\Omega^{\bullet}_{/S},F_b)=\Omega^{\bullet}_{g^*O_{\AnSp(\mathbb C)^{sm}/S}/g^*e(S)^*O_S}\to
(\Omega^{\bullet}_{/T},F_b)=\Omega^{\bullet}_{O_{\AnSp(\mathbb C)^{sm}/T}/e(T)^*O_T}
\end{eqnarray*}
which is by definition given by the pullback on differential forms : for $(V/T)=(V,h)\in\Var(\mathbb C)^{sm}/T$,
\begin{eqnarray*}
\Omega_{/(T/S)}(V/T): 
g^*(\Omega^{\bullet}_{/S})(V/T):=\lim_{(h':U\to S \mbox{sm},g':V\to U,h,g)}\Omega^{\bullet}_{U/S}(U)
\xrightarrow{\Omega_{(V/U)/(T/S)}(V/T)}\Omega^{\bullet}_{V/T}(V)=:\Omega^{\bullet}_{/T}(V/T) \\
\hat\omega\mapsto\Omega_{(V/U)/(T/S)}(V/T)(\omega):=\hat{g^{'*}\omega}.
\end{eqnarray*}
If $S$ and $T$ are smooth, $\Omega_{/(T/S)}:g^*(\Omega^{\bullet}_{/S},F_b)\to(\Omega^{\bullet}_{/T},F_b)$
is a map in $C_{g^*O_Sfil,g^*D^{\infty}_S}(\AnSp(\mathbb C)^{sm}/T)$
It induces the canonical morphism in $C_{g^*O_Sfil,g^*D^{\infty}_S}(\AnSp(\mathbb C)^{sm}/T)$:
\begin{eqnarray*}
E\Omega_{/(T/S)}:g^*E_{usu}(\Omega^{\bullet}_{/S},F_b)\xrightarrow{T(g,E_{usu})(\Omega^{\bullet}_{/S},F_b)}
E_{usu}(g^*(\Omega^{\bullet}_{/S},F_b))\xrightarrow{E_{usu}(\Omega_{/(T/S)})}E_{usu}(\Omega^{\bullet}_{/T},F_b). 
\end{eqnarray*}

\begin{defi}\label{TgDRGMan}
\begin{itemize}
\item[(i)]Let $g:T\to S$ a morphism with $T,S\in\AnSp(\mathbb C)$.
We have, for $F\in C(\AnSp(\mathbb C)^{sm}/S)$, the canonical transformation in $C_{O_Tfil}(T)$ :
\begin{eqnarray*}
T^O(g,\Omega_{/\cdot})(F): 
g^{*mod}L_Oe(S)_*\mathcal Hom^{\bullet}(F,E_{usu}(\Omega^{\bullet}_{/S},F_b)) \\
\xrightarrow{:=}
(g^*L_Oe(S)_*\mathcal Hom^{\bullet}(F,E_{et}(\Omega^{\bullet}_{/S},F_b)))\otimes_{g^*O_S}O_T \\
\xrightarrow{T(e,g)(-)\circ T(g,L_O)(-)}  
L_O(e(T)_*g^*\mathcal Hom^{\bullet}(F,E_{usu}(\Omega^{\bullet}_{/S},F))\otimes_{g^*O_S}O_T) \\ 
\xrightarrow{T(g,hom)(F,E_{et}(\Omega^{\bullet}_{/S}))\otimes I} 
L_O(e(T)_*\mathcal Hom^{\bullet}(g^*F,g^*E_{usu}(\Omega^{\bullet}_{/S},F_b))\otimes_{g^*O_S}O_T) \\
\xrightarrow{ev(hom,\otimes)(-,-,-)}
L_Oe(T)_*\mathcal Hom^{\bullet}(g^*F,g^*E_{usu}(\Omega^{\bullet}_{/S},F_b)\otimes_{g^*e(S)^*O_S}e(T)^*O_T) \\
\xrightarrow{\mathcal Hom^{\bullet}(g^*F,E\Omega_{/(T/S)}\otimes I)} 
L_Oe(T)_*\mathcal Hom^{\bullet}(g^*F,E_{usu}(\Omega^{\bullet}_{/T},F_b)\otimes_{g^*e(S)^*O_S}e(T)^*O_T) \\
\xrightarrow{m}
L_Oe(T)_*\mathcal Hom^{\bullet}(g^*F,E_{usu}(\Omega^{\bullet}_{/T},F_b)
\end{eqnarray*}
where $m(\alpha\otimes h):=h.\alpha$ is the multiplication map.
\item[(ii)] Let $g:T\to S$ a morphism with $T,S\in\AnSp(\mathbb C)$, $S$ smooth.
Assume there is a factorization $g:T\xrightarrow{l}Y\times S\xrightarrow{p_S}S$
with $l$ a closed embedding, $Y\in\AnSm(\mathbb C)$ and $p_S$ the projection.
We have, for $F\in C(\AnSp(\mathbb C)^{sm}/S)$, the canonical transformation in $C_{O_Tfil,\mathcal D^{\infty}}(Y\times S)$ :
\begin{eqnarray*}
T(g,\Omega_{/\cdot})(F):
g^{*mod,\Gamma}e(S)_*\mathcal Hom^{\bullet}(F,E_{usu}(\Omega^{\bullet}_{/S},F_b)) \\
\xrightarrow{:=} 
\Gamma_TE_{usu}(p_S^{*mod}e(S)_*\mathcal Hom^{\bullet}(F,E_{usu}(\Omega^{\bullet}_{/S},F_b))) \\
\xrightarrow{T^O(p_S,\Omega_{/\cdot})(F)} 
\Gamma_TE_{usu}(e(T\times S)_*\mathcal Hom^{\bullet}(p_S^*F,E_{usu}(\Omega^{\bullet}_{/T\times S},F_b))) \\
\xrightarrow{=}
e(T\times S)_*\Gamma_T(\mathcal Hom^{\bullet}(p_S^*F,E_{usu}(\Omega^{\bullet}_{/T\times S},F_b))) \\
\xrightarrow{I(\gamma,\hom)(-,-)}
e(T\times S)_*\mathcal Hom^{\bullet}(\Gamma^{\vee}_Tp_S^*F,E_{usu}(\Omega^{\bullet}_{/T\times S},F_b)).
\end{eqnarray*}
For $Q\in Proj\PSh(\AnSp(\mathbb C)^{sm}/S)$,
\begin{eqnarray*}
T(g,\Omega_{/\cdot})(Q):
g^{*mod,\Gamma}e(S)_*\mathcal Hom^{\bullet}(Q,E_{usu}(\Omega^{\bullet}_{/S},F_b))\to
e(T\times S)_*\mathcal Hom^{\bullet}(\Gamma^{\vee}_Tp_S^*Q,E_{usu}(\Omega^{\bullet}_{/Y\times S},F_b))
\end{eqnarray*}
is a map in $C_{O_Tfil,\mathcal D^{\infty}}(Y\times S)$.
\end{itemize}
\end{defi}

The following easy lemma describe these transformation map on representable presheaves :

\begin{lem}\label{homQomegaGMan}
Let $g:T\to S$ a morphism with $T,S\in\AnSp(\mathbb C)$ and $h:U\to S$ is a smooth morphism with $U\in\AnSp(\mathbb C)$.
Consider a commutative diagram whose square are cartesian :
\begin{equation*}
\xymatrix{g:T\ar[r]^l & S\times Y\ar[r]^{p_S} & S \\
g':U_T\ar[r]^{l'}\ar[u]^{h'} & U\times Y\ar[r]^{p_U}\ar[u]^{h'':=h\times I} & U\ar[u]^h}
\end{equation*} 
with $l$, $l'$ the graph embeddings and $p_S$, $p_U$ the projections. 
Then $g^*\mathbb Z(U/S)=\mathbb Z(U_T/T)$ and 
\begin{itemize}
\item[(i)] we have the following commutative diagram in $C_{O_Tfil}(T)$ 
(see definition \ref{TDw} and definition \ref{TgDRGMan}(i)) :
\begin{equation*}
\xymatrix{g^{*mod}L_Oe(S)_*\mathcal Hom^{\bullet}(\mathbb Z(U/S),E_{usu}(\Omega^{\bullet}_{/S},F_b))
\ar[rrr]^{T(g,\Omega_{/\cdot})(\mathbb Z(U/S))}\ar[d]_{=} & \, & \, &
e(T)_*\mathcal Hom^{\bullet}(\mathbb Z(U_T/T),E_{usu}(\Omega^{\bullet}_{/T},F_b))\ar[d]^{=} \\
g^{*mod}L_Oh_*E_{usu}(\Omega^{\bullet}_{U/S},F_b)\ar[rrr]^{T_{\omega}^{mod}(g,h)} & \, & \, & 
h'_*E_{usu}(\Omega^{\bullet}_{U_T/T},F_b)}
\end{equation*}
\item[(ii)] if $Y,S\in\AnSm(\mathbb C)$, we have the following commutative diagram in $C_{O_Tfil,\mathcal D^{\infty}}(T)$ 
(see definition \ref{TDw} and definition \ref{TgDRGMan}(ii)) :
\begin{equation*}
\xymatrix{g^{*mod,\Gamma}e(S)_*\mathcal Hom^{\bullet}(\mathbb Z(U/S),E_{usu}(\Omega^{\bullet}_{/S},F_b))
\ar[rrr]^{T(g,\Omega_{/\cdot})(\mathbb Z(U/S))}\ar[d]_{=} & \, & \, &
e(T)_*\mathcal Hom^{\bullet}(\mathbb Z(U_T/T),E_{usu}(\Omega^{\bullet}_{/T},F_b))\ar[d]^{=} \\
g^{*mod,\Gamma}h_*E_{usu}(\Omega^{\bullet}_{U/S},F_b)
\ar[rrr]^{T_{\omega}^O(\gamma,\otimes)(-)\circ T_{\omega}^O(p_S,h)(-)} & \, & \, & 
h'_*E_{usu}(\Omega^{\bullet}_{U_T/T},F_b)}
\end{equation*}
where $j:T\backslash T\times S\hookrightarrow T\times S$ is the open complementary embedding,
\end{itemize}
\end{lem}

\begin{proof}
Obvious.
\end{proof}

\begin{prop}\label{TpDRQGMan}
Let $p:S_{12}\to S_1$ is a smooth morphism  with $S_1,S_{12}\in\AnSp(\mathbb C)$. 
Then if $Q\in C(\AnSp(\mathbb C)^{sm}/S_1)$ is projective,
\begin{equation*}
T(p,\Omega_{/\cdot})(Q):
p^{*mod}e(S_1)_*\mathcal Hom^{\bullet}(Q,E_{usu}(\Omega^{\bullet}_{/S_1},F_b))\to
e(S_{12})_*\mathcal Hom^{\bullet}(p^*Q,E_{usu}(\Omega^{\bullet}_{/S_{12}},F_b))
\end{equation*}
is an isomorphism.
\end{prop}

\begin{proof}
Similar to the proof of proposition \ref{TpDRQGM}.
\end{proof}

Let $S\in\AnSp(\mathbb C)$ and $h:U\to S$ a morphism with $U\in\AnSp(\mathbb C)$.
We then have the canonical map given by the wedge product
\begin{equation*}
w_{U/S}:\Omega^{\bullet}_{U/S}\otimes_{O_S}\Omega^{\bullet}_{U/S}\to\Omega^{\bullet}_{U/S}; 
\alpha\otimes\beta\mapsto\alpha\wedge\beta.
\end{equation*}
Let $S\in\Var(\mathbb C)$ and $h_1:U_1\to S$, $h_2:U_2\to S$ two morphisms with $U_1,U_2\in\AnSp(\mathbb C)$.
Denote $h_{12}:U_{12}:=U_1\times_S U_2\to S$ and $p_{112}:U_1\times_S U_2\to U_1$, $p_{212}:U_1\times_S U_2\to U_2$ the projections.
We then have the canonical map given by the wedge product
\begin{equation*}
w_{(U_1,U_2)/S}:p_{112}^*\Omega^{\bullet}_{U_1/S}\otimes_{O_S}p_{212}^*\Omega^{\bullet}_{U_2/S}
\to \Omega^{\bullet}_{U_{12}/S}; \alpha\otimes\beta\mapsto p_{112}^*\alpha\wedge p_{212}^*\beta
\end{equation*}
which gives the map
\begin{eqnarray*}
Ew_{(U_1,U_2)/S}:h_{1*}E_{usu}(\Omega^{\bullet}_{U_1/S})\otimes_{O_S}h_{2*}E_{usu}(\Omega^{\bullet}_{U_2/S}) \\
\xrightarrow{\ad(p_{112}^*,p_{112*})(-)\otimes\ad(p_{212}^*,p_{212*})(-)} 
(h_{1*}p_{112*}p_{112}^*E_{usu}(\Omega^{\bullet}_{U_1/S}))\otimes_{O_S}(h_{2*}p_{212*}p_{212}^*E_{usu}(\Omega^{\bullet}_{U_2/S})) \\
\xrightarrow{=} 
h_{12*}(p_{112}^*E_{usu}(\Omega^{\bullet}_{U_1/S})\otimes_{h_{12}^*O_S}p_{212}^*E_{usu}(\Omega^{\bullet}_{U_2/S}) \\
\xrightarrow{T(\otimes,E)(-)\circ(T(p_{112},E)(-)\otimes T(p_{212},E)(-))}
h_{12*}E_{zar}(p_{112}^*\Omega^{\bullet}_{U_1/S}\otimes_{O_S}p_{212}^*\Omega^{\bullet}_{U_2/S})
\end{eqnarray*}

Let $S\in\AnSp(\mathbb C)$. We have the canonical map in $C_{O_Sfil}(\AnSp(\mathbb C)^{sm}/S)$
\begin{eqnarray*}
w_S:(\Omega^{\bullet}_{/S},F_b)\otimes_{O_S}(\Omega^{\bullet}_{/S},F_b)\to(\Omega^{\bullet}_{/S},F_b)
\end{eqnarray*}
given by for $h:U\to S\in\AnSp(\mathbb C)^{sm}/S$
\begin{eqnarray*}
w_S(U/S):(\Omega^{\bullet}_{U/S},F_b)\otimes_{h^*O_S}(\Omega^{\bullet}_{U/S},F_b)(U)\xrightarrow{w_{U/S}(U)}
(\Omega^{\bullet}_{U/S},F_b)(U)
\end{eqnarray*}
It gives the map
\begin{eqnarray*}
Ew_S:E_{usu}(\Omega^{\bullet}_{/S},F_b)\otimes_{O_S}E_{usu}(\Omega^{\bullet}_{/S},F_b)\xrightarrow{=}
E_{usu}((\Omega^{\bullet}_{/S},F_b)\otimes_{O_S}(\Omega^{\bullet}_{/S},F_b))\xrightarrow{E_{usu}(w_S)}
E_{usu}(\Omega^{\bullet}_{/S},F_b)
\end{eqnarray*}
If $S\in\AnSm(\mathbb C)$, 
\begin{eqnarray*}
w_S:(\Omega^{\bullet}_{/S},F_b)\otimes_{O_S}(\Omega^{\bullet}_{/S},F_b)\to(\Omega^{\bullet}_{/S},F_b)
\end{eqnarray*}
is a map in $C_{O_Sfil,D^{\infty}_S}(\Var(\mathbb C)^{sm}/S)$.

\begin{defi}\label{TotimesDRGMan}
Let $S\in\AnSp(\mathbb C)$.
We have, for $F,G\in C(\AnSp(\mathbb C)^{sm}/S)$, the canonical transformation in $C_{O_Sfil}(S)$ :
\begin{eqnarray*}
T(\otimes,\Omega)(F,G):
e(S)_*\mathcal Hom(F,E_{usu}(\Omega^{\bullet}_{/S},F_b))\otimes_{O_S}e(S)_*\mathcal Hom(G,E_{usu}(\Omega^{\bullet}_{/S},F_b)) \\ 
\xrightarrow{=}
e(S)_*(\mathcal Hom(F,E_{usu}(\Omega^{\bullet}_{/S},F_b))\otimes_{O_S}\mathcal Hom(G,E_{usu}(\Omega^{\bullet}_{/S},F_b))) \\
\xrightarrow{e(S)_*T(\mathcal Hom,\otimes)(-)}
e(S)_*\mathcal Hom(F\otimes G,E_{usu}(\Omega^{\bullet}_{/S},F_b)\otimes_{O_S}E_{usu}(\Omega^{\bullet}_{/S},F_b)) \\
\xrightarrow{\mathcal Hom(F\otimes G,Ew_S)}
e(S)_*\mathcal Hom(F\otimes G,E_{usu}(\Omega^{\bullet}_{/S},F_b))
\end{eqnarray*}
If $S\in\AnSm(\mathbb C)$, $T(\otimes,\Omega)(F,G)$ is a map in $C_{O_Sfil,\mathcal D^{\infty}}(S)$. 
\end{defi}

\begin{lem}\label{homQomega2GMan}
Let $S\in\AnSp(\mathbb C)$ and $h_1:U_1\to S$, $h_2:U_2\to S$ two smooth morphisms with $U_1,U_2\in\AnSp(\mathbb C)$.
Denote $h_{12}:U_{12}:=U_1\times_S U_2\to S$ and 
$p_{112}:U_1\times_S U_2\to U_1$, $p_{212}:U_1\times_S U_2\to U_2$ the projections.
We then have the following commutative diagram
\begin{equation*}
\xymatrix{
e(S)_*\mathcal Hom(F,E_{usu}(\Omega^{\bullet}_{/S},F_b))\otimes_{O_S}e(S)_*\mathcal Hom(G,E_{usu}(\Omega^{\bullet}_{/S},F_b))
\ar[rr]^{T(\otimes,\Omega)(F,G)}\ar[d]^{=} & \, & 
e(S)_*\mathcal Hom(F\otimes G,E_{usu}(\Omega^{\bullet}_{/S},F))\ar[d]^{=} \\
h_{1*}E_{usu}(\Omega^{\bullet}_{U_1/S},F_b)\otimes_{O_S}h_{2*}E_{usu}(\Omega^{\bullet}_{U_2/S},F_b)
\ar[rr]^{Ew_{(U_1,U_2)/S}} & \, & h_{12*}E_{zar}(\Omega^{\bullet}_{U_{12}/S},F_b)}
\end{equation*}
\end{lem}

\begin{proof}
Follows from Yoneda lemma.
\end{proof}

We now define the analytic Gauss Manin De Rahm realization functor.

Let $S\in\AnSp(\mathbb C)$ and $S=\cup_{i=1}^l S_i$ an open cover such that there exist 
closed embeddings $i_i:S_i\hookrightarrow\tilde S_i$ with $\tilde S_i\in\AnSm(\mathbb C)$ an affine space. 
For $I\subset\left[1,\cdots l\right]$, denote by $S_I:=\cap_{i\in I} S_i$ and $j_I:S_I\hookrightarrow S$ the open embedding.
We then have closed embeddings $i_I:S_I\hookrightarrow\tilde S_I:=\Pi_{i\in I}\tilde S_i$.
Consider, for $I\subset J$, the following commutative diagram
\begin{equation*}
D_{IJ}=\xymatrix{ S_I\ar[r]^{i_I} & \tilde S_I \\
S_J\ar[u]^{j_{IJ}}\ar[r]^{i_J} & \tilde S_J\ar[u]^{p_{IJ}}}  
\end{equation*}
and $j_{IJ}:S_J\hookrightarrow S_I$ is the open embedding so that $j_I\circ j_{IJ}=j_J$.
Considering the factorization of the diagram $D_{IJ}$ by the fiber product :
\begin{equation*}
D_{IJ}=\xymatrix{
\tilde S_J=\tilde S_I\times\tilde S_{J\backslash I}\ar[rr]^{p_{IJ}} & \, & \tilde S_I \\
\, & S_I\times\tilde S_{J\backslash I}\ar[lu]^{i_I\times I}\ar[rd]^{p_{IJ}^0} & \, \\
S_J\ar[uu]^{i_J}\ar[ru]^{l_J}\ar[rr]^{j_{IJ}} & \, & S_I\ar[uu]^{i_I}}
\end{equation*}
the square of this factorization being cartesian, 
we have for $F\in C(\AnSp(\mathbb C)^{sm}/S)$ the canonical map in $C(\AnSp(\mathbb C)^{sm}/\tilde S_J)$
\begin{eqnarray*}
S(D_{IJ})(F):Li_{J*}j_J^*F\xrightarrow{q}i_{J*}j_J^*F=(i_I\times I)*l_{J*}j_J^*F 
\xrightarrow{(i_I\times I)_*\ad(p_{IJ\sharp}^o,p_{IJ}^{o*})(-)} \\
(i_I\times I)_*p_{IJ}^{o*}p_{IJ\sharp}^0l_{J*}j_J^*F 
\xrightarrow{T(p_{IJ},i_I)(-)^{-1}} p_{IJ}^*i_{I*}p_{IJ\sharp}^0l_{J*}j_I^*F=p_{IJ}^*i_{I*}j_I^*F
\end{eqnarray*}
which factors through
\begin{eqnarray*}
S(D_{IJ})(F):Li_{J*}j_I^*F\xrightarrow{S^q(D_{IJ})(F)}p_{IJ}^*Li_{I*}j_I^*F\xrightarrow{q}p_{IJ}^*i_{I*}j_I^*F
\end{eqnarray*}

\begin{defi}\label{DRalgdefFunctGMan}
\begin{itemize}
\item[(i)]Let $S\in\AnSm(\mathbb C)$. We have the functor
\begin{equation*}
\mathcal Hom^{\bullet}(\cdot, E_{usu}(\Omega^{\bullet}_{/S},F_b)):
C(\AnSp(\mathbb C)^{sm}/S)\to C_{O_Sfil,D^{\infty}_S}(S), \;
F\mapsto e(S)_*\mathcal Hom^{\bullet}(LF,E_{usu}(\Omega^{\bullet}_{/S},F_b)).
\end{equation*}
\item[(ii)]Let $S\in\Var(\mathbb C)$ and $S=\cup_{i=1}^l S_i$ an open cover such that there exist closed embeddings
$i_i:S_i\hookrightarrow\tilde S_i$  with $\tilde S_i\in\SmVar(\mathbb C)$. 
For $I\subset\left[1,\cdots l\right]$, denote by $S_I:=\cap_{i\in I} S_i$ and $j_I:S_I\hookrightarrow S$ the open embedding.
We then have closed embeddings $i_I:S_I\hookrightarrow\tilde S_I:=\Pi_{i\in I}\tilde S_i$.
We have the functor
\begin{eqnarray*}
C(\Var(\mathbb C)^{sm}/S)^{op}\to C_{Ofil,\mathcal D^{\infty}}(S/(\tilde S_I)), \; \;
F\mapsto(e(\tilde S_I)_*\mathcal Hom^{\bullet}(\An_{\tilde S_I}^*L(i_{I*}j_I^*F),
E_{usu}(\Omega^{\bullet}_{/\tilde S_I},F_b))[-d_{\tilde S_I}],u^q_{IJ}(F))
\end{eqnarray*}
where
\begin{eqnarray*}
u^q_{IJ}(F)[d_{\tilde S_J}]:
e(\tilde S_I)_*\mathcal Hom^{\bullet}(\An_{\tilde S_I}^*L(i_{I*}j_I^*F),E_{usu}(\Omega^{\bullet}_{/\tilde S_I},F_b)) \\
\xrightarrow{\ad(p_{IJ}^{*mod},p_{IJ*})(-)}
p_{IJ*}p_{IJ}^{*mod}e(\tilde S_I)_*\mathcal Hom^{\bullet}(\An_{\tilde S_I}^*L(i_{I*}j_I^*F),
E_{usu}(\Omega^{\bullet}_{/\tilde S_I},F_b)) \\
\xrightarrow{p_{IJ*}T(p_{IJ},\Omega_{\cdot})(L(i_{I*}j_I^*F))}
p_{IJ*}e(\tilde S_J)_*\mathcal Hom^{\bullet}(\An_{\tilde S_J}^*p_{IJ}^*L(i_{I*}j_I^*F),
E_{usu}(\Omega^{\bullet}_{/\tilde S_J},F_b)) \\
\xrightarrow{p_{IJ*}e(\tilde S_J)_*\mathcal Hom(\An_{\tilde S_J}^*S^q(D_{IJ})(F),
E_{usu}(\Omega_{/\tilde S_J}^{\bullet,\Gamma},F_b))}  
p_{IJ*}e(\tilde S_J)_*\mathcal Hom^{\bullet}(\An_{\tilde S_J}^*L(i_{J*}j_J^*F),E_{usu}(\Omega^{\bullet}_{/\tilde S_J},F_b)). 
\end{eqnarray*}
For $I\subset J\subset K$, we have obviously $p_{IJ*}u_{JK}(F)\circ u_{IJ}(F)=u_{IK}(F)$.
\end{itemize}
\end{defi}

\begin{prop}\label{projwachGMan}
\begin{itemize}
\item[(i)]Let $S\in\AnSp(\mathbb C)$. 
Let $m:Q_1\to Q_2$ be an equivalence $(\mathbb D^1,et)$ local in $C(\AnSp(\mathbb C)^{sm}/S)$
with $Q_1,Q_2$ complexes of projective presheaves. Then,
\begin{equation*}
e(S)_*\mathcal Hom(m,E_{et}(\Omega^{\bullet}_{/S},F_b)):
e(S)_*\mathcal Hom^{\bullet}(Q_2,E_{et}(\Omega^{\bullet}_{/S},F_b))\to 
e(S)_*\mathcal Hom^{\bullet}(Q_1,E_{et}(\Omega^{\bullet}_{/S},F_b))
\end{equation*}
is an $2$-filtered quasi-isomorphism. 
It is thus an isomorphism in $D_{O_Sfil,\mathcal D^{\infty},\infty}(S)$ if $S$ is smooth.
\item[(ii)]Let $S\in\AnSp(\mathbb C)$. Let $S=\cup_{i=1}^l S_i$ an open cover such that there exist closed embeddings
$i_i:S_i\hookrightarrow\tilde S_i$  with $\tilde S_i\in\AnSm(\mathbb C)$.
Let $m=(m_I):(Q_{1I},s^1_{IJ})\to (Q_{2I},s^2_{IJ})$ be an equivalence $(\mathbb D^1,usu)$ local 
in $C(\AnSp(\mathbb C)^{sm}/(\tilde S_I)^{op})$ with $Q_{1I},Q_{2I}$ complexes of projective presheaves. Then,
\begin{eqnarray*}
(e(\tilde S_I)_*\mathcal Hom(m_I,E_{et}(\Omega^{\bullet}_{/\tilde S_I},F_b))): \\
(e(\tilde S_I)_*\mathcal Hom^{\bullet}(Q_{2I},E_{et}(\Omega^{\bullet}_{/\tilde S_I},F_b)),u_{IJ}(Q_{2I},s^2_{IJ}))\to 
(e(\tilde S_I)_*\mathcal Hom^{\bullet}(Q_{1I},E_{et}(\Omega^{\bullet}_{/\tilde S_I},F_b)),u_{IJ}(Q_{1I},s^1_{IJ}))
\end{eqnarray*}
is an $2$-filtered quasi-isomorphism. 
It is thus an isomorphism in $D_{O_Sfil,\mathcal D^{\infty},\infty}(S/(\tilde S_I))$.
\end{itemize}
\end{prop}

\begin{proof}
Similar to the proof of proposition \ref{projwachGM}.
\end{proof}

\begin{defi}\label{DRalgdefsingGMan}
\begin{itemize}
\item[(i)]We define, according to proposition \ref{projwachGMan}, 
the filtered analytic Gauss-Manin realization functor defined as
\begin{eqnarray*}
\mathcal F_{S,an}^{GM}:\DA_c(S)^{op}\to D_{O_Sfil,\mathcal D^{\infty},\infty}(S), \; \; M\mapsto \\
\mathcal F_{S,an}^{GM}(M):=e(S)_*\mathcal Hom^{\bullet}(\An_S^*L(F),E_{usu}(\Omega^{\bullet}_{/S},F_b))[-d_S] \\
=e(S)_*\mathcal Hom^{\bullet}(L(F),\An_{S*}E_{usu}(\Omega^{\bullet}_{/S},F_b))[-d_S]
\end{eqnarray*}
where $F\in C(\Var(\mathbb C)^{sm}/S)$ is such that $M=D(\mathbb A^1,et)(F)$,
\item[(ii)]Let $S\in\Var(\mathbb C)$ and $S=\cup_{i=1}^l S_i$ an open cover such that there exist closed embeddings
$i_i:S_i\hookrightarrow\tilde S_i$  with $\tilde S_i\in\SmVar(\mathbb C)$. 
For $I\subset\left[1,\cdots l\right]$, denote by $S_I=\cap_{i\in I} S_i$ and $j_I:S_I\hookrightarrow S$ the open embedding.
We then have closed embeddings $i_I:S_I\hookrightarrow\tilde S_I:=\Pi_{i\in I}\tilde S_i$.
We define the filtered analytic Gauss-Manin realization functor defined as
\begin{eqnarray*}
\mathcal F_{S,an}^{GM}:\DA_c(S)^{op}\to D_{Ofil,\mathcal D^{\infty},\infty}(S/(\tilde S_I)), \; M\mapsto \\
\mathcal F_{S,an}^{GM}(M):=((e(\tilde S_I)_*\mathcal Hom^{\bullet}(\An_{\tilde S_I}^*L(i_{I*}j_I^*F),
E_{usu}(\Omega^{\bullet}_{/\tilde S_I}),F_b))[-d_{\tilde S_I}],u^q_{IJ}(F)) \\
=((e(\tilde S_I)_*\mathcal Hom^{\bullet}(L(i_{I*}j_I^*F),
\An_{\tilde S_I*}E_{usu}(\Omega^{\bullet}_{/\tilde S_I}),F_b))[-d_{\tilde S_I}],u^q_{IJ}(F))
\end{eqnarray*}
where $F\in C(\Var(\mathbb C)^{sm}/S)$ is such that $M=D(\mathbb A^1,et)(F)$, 
see definition \ref{DRalgdefFunctGMan}.
\end{itemize}
\end{defi}

\begin{prop}\label{FDRwelldefGMan}
For $S\in\Var(\mathbb C)$, the functor $\mathcal F_S^{GM}$ is well defined.
\end{prop}

\begin{proof}
Similar to the proof of proposition \ref{FDRwelldefGM}.
\end{proof}

\begin{prop}\label{keyalgsing1GMan}
Let $f:X\to S$ a morphism with $S,X\in\Var(\mathbb C)$.
Let $S=\cup_{i=1}^l S_i$ an open cover such that there exist closed embeddings
$i_i:S_i\hookrightarrow\tilde S_i$ with $\tilde S_i\in\SmVar(\mathbb C)$. 
Then $X=\cup_{i=1}^lX_i$ with $X_i:=f^{-1}(S_i)$.
Denote, for $I\subset\left[1,\cdots l\right]$, $S_I=\cap_{i\in I} S_i$ and $X_I=\cap_{i\in I}X_i$.
Assume there exist a factorization 
\begin{equation*}
f:X\xrightarrow{l}Y\times S\xrightarrow{p_S} S
\end{equation*}
of $f$ with $Y\in\SmVar(\mathbb C)$, $l$ a closed embedding and $p_S$ the projection.
We then have, for $I\subset\left[1,\cdots l\right]$, the following commutative diagrams which are cartesian 
\begin{equation*}
\xymatrix{
f_I=f_{|X_I}:X_I\ar[r]^{l_I}\ar[rd] & Y\times S_I\ar[r]^{p_{S_I}}\ar[d]^{i'_I} & S_I\ar[d]^{i_I} \\
\, & Y\times\tilde S_I\ar[r]^{p_{\tilde S_I}} & \tilde S_I} \;, \;
\xymatrix{Y\times\tilde S_J\ar[r]^{p_{\tilde S_J}}\ar[d]_{p'_{IJ}} & \tilde S_J\ar[d]^{p_{IJ}} \\
Y\times\tilde S_I\ar[r]^{p_{\tilde S_I}} & \tilde S_I}
\end{equation*}
Let $F(X/S):=p_{S,\sharp}\Gamma_X^{\vee}\mathbb Z(Y\times S/Y\times S)$. 
The transformations maps $(N_I(X/S):Q(X_I/\tilde S_I)\to i_{I*}j_I^*F(X/S))$ and $(k\circ I(\gamma,\hom)(-,-))$, 
for $I\subset\left[1,\cdots,l\right]$, induce an isomorphism in $D_{Ofil,\mathcal D^{\infty},\infty}(S/(\tilde S_I))$ 
\begin{eqnarray*} 
I^{GM}(X/S): \\
\mathcal F_{S,an}^{GM}(M(X/S)):=
(e(\tilde S_I)_*\mathcal Hom(\An_{\tilde S_I}^*L(i_{I*}j_I^*F(X/S)),
E_{usu}(\Omega^{\bullet}_{/\tilde S_I},F_b))[-d_{\tilde S_I}],u^q_{IJ}(F(X/S))) \\
\xrightarrow{(e(\tilde S_I)_*\mathcal Hom(\An_{\tilde S_I}^*N_I(X/S),E_{usu}(\Omega^{\bullet}_{/\tilde S_I},F_b)))} 
(e(\tilde S_I)_*\mathcal Hom(\An_{\tilde S_I}^*Q(X_I/\tilde S_I),
E_{usu}(\Omega^{\bullet}_{/\tilde S_I},F_b))[-d_{\tilde S_I}],v^q_{IJ}(F(X/S))) \\
\xrightarrow{(e(\tilde S_I)_*\mathcal Hom(T(\An,\gamma^{\vee})(-)^{-1},E_{usu}(\Omega^{\bullet}_{/\tilde S_I},F_b)))} 
(e(\tilde S_I)_*\mathcal Hom(Q(X^{an}_I/\tilde S^{an}_I),
E_{usu}(\Omega^{\bullet}_{/\tilde S_I},F_b))[-d_{\tilde S_I}],v^q_{IJ}(F(X/S))) \\
\xrightarrow{(I(\gamma,\hom)(-,-))^{-1}}
(p_{\tilde S_I*}\Gamma_{X_I}E_{usu}(\Omega^{\bullet}_{Y\times\tilde S_I/\tilde S_I},F_b)[-d_{\tilde S_I}],
w_{IJ}(X/S)).
\end{eqnarray*}
\end{prop}

\begin{proof}
Similar to the proof of proposition  \ref{keyalgsing1GM}.
\end{proof}

We now define the functorialities of $\mathcal F_S^{GM}$ with respect to $S$ 
which makes $\mathcal F^{-}_{GM}$ a morphism of 2-functor.
\begin{defi}\label{TgDRdefGMan}
Let $g:T\to S$ a morphism with $T,S\in\SmVar(\mathbb C)$.
Consider the factorization $g:T\xrightarrow{l}T\times S\xrightarrow{p_S}S$
where $l$ is the graph embedding and $p_S$ the projection.
Let $M\in\DA_c(S)$ and $F\in C(\Var(\mathbb C)^{sm}/S)$ such that $M=D(\mathbb A^1_S,et)(F)$. 
Then, $D(\mathbb A^1_T,et)(g^*F)=g^*M$. 
\begin{itemize}
\item[(i)]We have then the canonical transformation in $D_{\mathcal D^{\infty}fil,\infty}(T\times S)$
(see definition \ref{TgDRGMan}) :
\begin{eqnarray*}
T(g,\mathcal F^{GM})(M):Rg^{*mod[-],\Gamma}\mathcal F_{S,an}^{GM}(M):=
g^{*mod,\Gamma}e(S)_*\mathcal Hom^{\bullet}(\An_S^*L(F),E_{usu}(\Omega^{\bullet}_{/S},F_b)))[-d_T] \\
\xrightarrow{T(g,\Omega_{/\cdot})(\An_S^*L(F))}  
e(T\times S)_*\mathcal Hom^{\bullet}(\Gamma_T^{\vee}p_S^*\An_S^*L(F),E_{usu}(\Omega^{\bullet}_{/Y\times S},F_b))[-d_T] \\
\xrightarrow{\mathcal Hom(T(\An,\gamma^{\vee})(p_S^*LF)^{-1},-)}
\mathcal F_{T\times S,an}^{GM}(l_*g^*M). 
\end{eqnarray*}
where the last isomorphism in the derived category comes from proposition \ref{FDRwelldefGMan}.
\item[(ii)]We have then the canonical transformation in $D_{Ofil,\infty}(T)$
(see definition \ref{TgDRGMan}) :
\begin{eqnarray*}
T^O(g,\mathcal F^{GM})(M):Lg^{*mod[-]}\mathcal F_{S,an}^{GM}(M):=
g^{*mod}L_Oe(S)_*\mathcal Hom^{\bullet}(\An_S^*L(F),E_{usu}(\Omega^{\bullet}_{/S},F_b)))[-d_T] \\
\xrightarrow{T(g,\Omega_{/\cdot})(\An_S^*L(F))}  
e(T\times S)_*\mathcal Hom^{\bullet}(g^*\An_S^*L(F),E_{usu}(\Omega^{\bullet}_{/Y\times S},F_b))[-d_T]
=:\mathcal F_{T,an}^{GM}(g^*M). 
\end{eqnarray*}
\end{itemize}
\end{defi}

We give now the definition in the non smooth case
Let $g:T\to S$ a morphism with $T,S\in\Var(\mathbb C)$.
Assume we have a factorization $g:T\xrightarrow{l}Y\times S\xrightarrow{p_S}S$
with $Y\in\SmVar(\mathbb C)$, $l$ a closed embedding and $p_S$ the projection.
Let $S=\cup_{i=1}^lS_i$ be an open cover such that 
there exists closed embeddings $i_i:S_i\hookrightarrow\tilde S_i$ with $\tilde S_i\in\SmVar(\mathbb C)$
Then, $T=\cup^l_{i=1} T_i$ with $T_i:=g^{-1}(S_i)$
and we have closed embeddings $i'_i:=i_i\circ l:T_i\hookrightarrow Y\times\tilde S_i$,
Moreover $\tilde g_I:=p_{\tilde S_I}:Y\times\tilde S_I\to\tilde S_I$ is a lift of $g_I:=g_{|T_I}:T_I\to S_I$.
We recall the commutative diagram :
\begin{equation*}
E_{IJg}=\xymatrix{(Y\times\tilde S_I)\backslash T_I\ar[d]^{p_{\tilde S_I}}\ar[r]^{m'_I} & 
Y\times\tilde S_J\ar[d]^{\tilde g_I} \\
\tilde S_I\backslash S_I\ar[r]^{m_I}& \tilde S_I}, \,
E_{IJ}=\xymatrix{\tilde S_J\backslash S_J\ar[d]^{p_{IJ}}\ar[r]^{m_J} & \tilde S_J\ar[d]^{p_{IJ}} \\
\tilde S_I\backslash(S_I\backslash S_J)\ar[r]^{m=m_{IJ}}& \tilde S_I}
E'_{IJ}=\xymatrix{(Y\times\tilde S_J)\backslash T_J\ar[d]^{p'_{IJ}}\ar[r]^{m'_J} & Y\times\tilde S_J\ar[d]^{p'_{IJ}} \\
(Y\times\tilde S_I)\backslash(T_I\backslash T_J)\ar[r]^{m'=m'_{IJ}}& Y\times\tilde S_I}
\end{equation*}
For $I\subset J$, denote by $p_{IJ}:\tilde S_J\to\tilde S_I$ and 
$p'_{IJ}:=I_Y\times p_{IJ}:Y\times\tilde S_J\to Y\times\tilde S_I$ the projections,
so that $\tilde g_I\circ p'_{IJ}=p_{IJ}\circ\tilde g_J$.
Consider, for $I\subset J\subset\left[1,\ldots,l\right]$,
resp. for each $I\subset\left[1,\ldots,l\right]$, the following commutative diagrams in $\Var(\mathbb C)$
\begin{equation*}
D_{IJ}=\xymatrix{ S_I\ar[r]^{i_I} & \tilde S_I \\
S_J\ar[u]^{j_{IJ}}\ar[r]^{i_J} & \tilde S_J\ar[u]^{p_{IJ}}} \; , \;  
D'_{IJ}=\xymatrix{ T_I\ar[r]^{i'_I} & Y\times\tilde S_I \\
T_J\ar[u]^{j'_{IJ}}\ar[r]^{i'_J} & Y\times\tilde S_J\ar[u]^{p'_{IJ}}}  
D_{gI}=\xymatrix{ S_I\ar[r]^{i_I} & \tilde S_I \\
T_I\ar[u]^{g_I}\ar[r]^{i'_I} & Y\times\tilde S_I\ar[u]^{\tilde g_I}} \; , \;  
\end{equation*}
and $j_{IJ}:S_J\hookrightarrow S_I$ is the open embedding so that $j_I\circ j_{IJ}=j_J$.
Let $F\in C(\Var(\mathbb C)^{sm}/S)$. 
The fact that the diagrams (\ref{DgDIq}) commutes says that the maps $T^{q,\gamma}(D_{gI})(j_I^*F)$
define a morphism in $C(\Var(\mathbb C)^{sm}/(T/(Y\times\tilde S_I)))$
\begin{equation*}
(T^{q,\gamma}(D_{gI})(j_I^*F)):(\Gamma^{\vee}_{T_I}\tilde g_I^*L(i_{I*}j_I^*F),\tilde g_J^*S^q(D_{IJ})(F))
\to (L(i'_{I*}j_I^{'*}g^*F),S^q(D'_{IJ})(g^*F))
\end{equation*}

We then have then the following lemma :

\begin{lem}\label{TglemGMan}
\begin{itemize}
\item[(i)]The morphism in $C(\Var(\mathbb C)^{sm}/(T/(Y\times\tilde S_I)))$
\begin{equation*}
(T^{q,\gamma}(D_{gI})(j_I^*F)):(\Gamma^{\vee}_{T_I}L\tilde g_I^*i_{I*}j_I^*F,\tilde g_J^*S^q(D_{IJ})(F))
\to (i'_{I*}j_I^{'*}g^*F,S^q(D'_{IJ})(g^*F))
\end{equation*}
is an equivalence $(\mathbb A^1,et)$ local.
\item[(ii)] Denote for short $d_{YI}:=-d_Y-d_{\tilde S_I}$. The maps 
$\mathcal Hom(\An_{Y\times\tilde S_I}^*(T^{q,\gamma}(D_{gI})(j_I^*F)),E_{usu}(\Omega^{\bullet}_{/Y\times\tilde S_I},F_b))$
induce an $\infty$-filtered quasi-isomorphism in $C_{Ofil,\mathcal D^{\infty}}(T/(Y\times\tilde S_I))$
\begin{eqnarray*}
(\mathcal Hom(\An_{Y\times\tilde S_I}^*T^{q,\gamma}(D_{gI})(j_I^*F),E_{usu}(\Omega^{\bullet}_{/Y\times\tilde S_I},F_b))): \\
(e(Y\times\tilde T_I)_*\mathcal Hom(\An_{Y\times\tilde S_I}^*L(i'_{I*}j_I^{'*}g^*F),
E_{usu}(\Omega^{\bullet}_{/Y\times\tilde S_I},F_b))[d_{YI}],u^q_{IJ}(g^*F))\to \\ 
(e(Y\times\tilde T_I)_*\mathcal Hom(\An_{Y\times\tilde S_I}^*\Gamma^{\vee}_{T_I}L(\tilde g_I^*i_{I*}j_I^*F),
E_{usu}(\Omega^{\bullet}_{/Y\times\tilde S_I},F_b))[d_{YI}],\tilde g_J^*u^q_{IJ}(F)_2) 
\end{eqnarray*}
\item[(iii)] The maps $T(\tilde g_I,\Omega_{\cdot})(L(i_{I*}j_I^*F))$ (see definition \ref{TgDRGMan})
induce a morphism in $C_{Ofil,\mathcal D^{\infty}}(T/(Y\times\tilde S_I))$
\begin{eqnarray*}
(T(\tilde g_I,\Omega_{/\cdot})(L(i_{I*}j_I^*F))): \\ 
(\Gamma_{T_I}E_{zar}(\tilde g_I^{*mod[-]}e(\tilde S_I)_*\mathcal Hom^{\bullet}(\An_{\tilde S_I}^*L(i_{I*}j_I^*F),
E_{usu}(\Omega^{\bullet}_{/\tilde S_I},F_b)))[d_{YI}],\tilde g_J^{*mod}u^q_{IJ}(F))\to \\  
(\Gamma_{T_I}(e(Y\times\tilde S_I)_*\mathcal Hom(\An_{Y\times\tilde S_I}^*\tilde g_I^*L(i_{I*}j_I^*F),
E_{usu}(\Omega^{\bullet}_{/Y\times\tilde S_I},F_b)))[d_{YI}],\tilde g_J^*u^q_{IJ}(F)_1). 
\end{eqnarray*}
\end{itemize}
\end{lem}

\begin{proof}

\noindent(i):Follows from theorem \ref{2functDM}

\noindent(ii): Similar to lemma \ref{TglemGM}(ii).

\noindent(iii):Similar to lemma \ref{TglemGM}(iii).

\end{proof}

\begin{defi}\label{TgDRdefsingGMan}
Let $g:T\to S$ a morphism with $T,S\in\Var(\mathbb C)$.
Assume we have a factorization $g:T\xrightarrow{l}Y\times S\xrightarrow{p_S}S$
with $Y\in\SmVar(\mathbb C)$, $l$ a closed embedding and $p_S$ the projection.
Let $S=\cup_{i=1}^lS_i$ be an open cover such that 
there exists closed embeddings $i_i:S_i\hookrightarrow\tilde S_i$ with $\tilde S_i\in\SmVar(\mathbb C)$
Then, $T=\cup^l_{i=1} T_i$ with $T_i:=g^{-1}(S_i)$
and we have closed embeddings $i'_i:=i_i\circ l:T_i\hookrightarrow Y\times\tilde S_i$,
Moreover $\tilde g_I:=p_{\tilde S_I}:Y\times\tilde S_I\to\tilde S_I$ is a lift of $g_I:=g_{|T_I}:T_I\to S_I$.
Denote for short $d_{YI}:=-d_Y-d_{\tilde S_I}$.
Let $M\in\DA_c(S)$ and $F\in C(\Var(\mathbb C)^{sm}/S)$ such that $M=D(\mathbb A^1_S,et)(F)$.
Then, $D(\mathbb A^1_T,et)(g^*F)=g^*M$. 
We have, by lemma \ref{TglemGM}, the canonical transformation in $D_{Ofil,\mathcal D^{\infty},\infty}(T/(Y\times\tilde S_I))$
\begin{eqnarray*}
T(g,\mathcal F^{GM})(M):Rg^{*mod[-],\Gamma}\mathcal F_{S,an}^{GM}(M):= \\
(\Gamma_{T_I}E_{zar}(\tilde g_I^{*mod}e(\tilde S_I)_*\mathcal Hom^{\bullet}(\An_{\tilde S_I}^*L(i_{I*}j_I^*F),
E_{usu}(\Omega^{\bullet}_{/\tilde S_I},F_b)))[d_{YI}],\tilde g_J^{*mod}u^q_{IJ}(F)) \\
\xrightarrow{(\Gamma_{T_I}E(T(\tilde g_I,\Omega_{/\cdot})(\An_{\tilde S_I}^*L(i_{I*}j_I^*F))))} \\
(\Gamma_{T_I}e(Y\times\tilde S_I)_*\mathcal Hom^{\bullet}(\An_{Y\times\tilde S_I}^*\tilde g_I^*L(i_{I*}j_I^*F),
E_{usu}(\Omega^{\bullet}_{/Y\times\tilde S_I},F_b))[d_{YI}],\tilde g_J^*u^q_{IJ}(F)_1) \\ 
\xrightarrow{(I(\gamma,\hom(-,-)))} \\
(e(Y\times\tilde S_I)_*\mathcal Hom^{\bullet}(\Gamma_{T_I}^{\vee}\An_{Y\times\tilde S_I}^*\tilde g_I^*L(i_{I*}j_I^*F),
E_{usu}(\Omega^{\bullet}_{/Y\times\tilde S_I},F_b))[d_{YI}],\tilde g_J^*u^q_{IJ}(F)_2) \\ 
\xrightarrow{(e(Y\times\tilde S_I)_*\mathcal Hom^{\bullet}(T(\An,\gamma^{\vee})(\tilde g_I^*L(i_{I*}j_I^*F))^{-1},
E_{usu}(\Omega^{\bullet}_{/Y\times\tilde S_I},F_b)))} \\
(e(Y\times\tilde S_I)_*\mathcal Hom^{\bullet}(\An_{Y\times\tilde S_I}^*\Gamma_{T_I}^{\vee}\tilde g_I^*L(i_{I*}j_I^*F),
E_{usu}(\Omega^{\bullet}_{/Y\times\tilde S_I},F_b))[d_{YI}],\tilde g_J^*u^q_{IJ}(F)_2) \\ 
\xrightarrow{(e(Y\times\tilde S_I)_*\mathcal Hom(\An_{Y\times\tilde S_I}^*T^{q,\gamma}(D_{gI})(j_I^*F),
E_{usu}(\Omega^{\bullet}_{/Y\times\tilde S_I},F_b)))^{-1}} \\
(e(Y\times\tilde S_I)_*\mathcal Hom^{\bullet}(\An_{Y\times\tilde S_I}^*L(i'_{I*}j_I^{'*}g^*F),
E_{usu}(\Omega^{\bullet}_{/Y\times\tilde S_I},F_b))[d_{YI}],u^q_{IJ}(g^*F))=:\mathcal F_{T,an}^{GM}(g^*M).
\end{eqnarray*} 
\end{defi}

\begin{prop}\label{TgGMpropan}
\begin{itemize}
\item[(i)]Let $g:T\to S$ a morphism with $T,S\in\Var(\mathbb C)$.
Assume we have a factorization $g:T\xrightarrow{l}Y_2\times S\xrightarrow{p_S}S$
with $Y_2\in\SmVar(\mathbb C)$, $l$ a closed embedding and $p_S$ the projection.
Let $S=\cup_{i=1}^lS_i$ be an open cover such that 
there exists closed embeddings $i_i:S_i\hookrightarrow\tilde S_i$ with $\tilde S_i\in\SmVar(\mathbb C)$
Then, $T=\cup^l_{i=1} T_i$ with $T_i:=g^{-1}(S_i)$
and we have closed embeddings $i'_i:=i_i\circ l:T_i\hookrightarrow Y_2\times\tilde S_i$,
Moreover $\tilde g_I:=p_{\tilde S_I}:Y\times\tilde S_I\to\tilde S_I$ is a lift of $g_I:=g_{|T_I}:T_I\to S_I$.
Let $f:X\to S$ a  morphism with $X\in\Var(\mathbb C)$. Assume that there is a factorization 
$f:X\xrightarrow{l}Y_1\times S\xrightarrow{p_S} S$, with $Y_1\in\SmVar(\mathbb C)$, 
$l$ a closed embedding and $p_S$ the projection. We have then the following commutative diagram
whose squares are cartesians
\begin{equation*}
\xymatrix{f':X_T\ar[r]\ar[d] & Y_1\times T\ar[d]\ar[r] & T\ar[d] \\
f''=f\times I:Y_2\times X\ar[r]\ar[d] & Y_1\times Y_2\times S\ar[r]\ar[d] & Y_2\times S\ar[d] \\
f:X\ar[r] & Y_1\times S\ar[r] & S}
\end{equation*} 
Consider $F(X/S):=p_{S,\sharp}\Gamma_X^{\vee}\mathbb Z(Y_1\times S/Y_1\times S)[d_{Y_1}]$ and
the isomorphism in $C(\Var(\mathbb C)^{sm}/S)$
\begin{eqnarray*}
T(f,g,F(X/S)):g^*F(X/S):=g^*p_{S,\sharp}\Gamma_X^{\vee}\mathbb Z(Y_1\times S/Y_1\times S)\xrightarrow{\sim} \\
p_{T,\sharp}\Gamma_{X_T}^{\vee}\mathbb Z(Y_1\times T/Y_1\times T)=:F(X_T/T).
\end{eqnarray*}
which gives in $\DA(S)$ the isomorphism $T(f,g,F(X/S)):g^*M(X/S)\xrightarrow{\sim}M(X_T/T)$. 
Then, the following diagram in $D_{Ofil,\mathcal D^{\infty},\infty}(T/(Y_2\times\tilde S_I))$ commutes
\begin{equation*}
\begin{tikzcd}
Rg^{*mod,\Gamma}\mathcal F_{S,an}^{GM}(M(X/S))
\ar[rrr,"T(g{,}\mathcal F^{GM})(M(X/S))"]\ar[d,"I^{GM}(X/S)"] & \, & \, & 
\mathcal F_{T,an}^{GM}(M(X_T/T))\ar[d,"I^{GM}(X_T/T)"] \\
\shortstack{$g^{*mod[-],\Gamma}(p_{\tilde S_I*}\Gamma_{X_I}
E_{usu}(\Omega^{\bullet}_{Y_1\times\tilde S_I/\tilde S_I},F_b)[-d_{\tilde S_I}],$ \\ $w_{IJ}(X/S))$}
\ar[rrr,"(T(\tilde g_I\times I{,}\gamma)(-)\circ T^O_w(\tilde g_I{,}p_{\tilde S_I}))"] & \, & \, &
\shortstack{$(p_{Y_2\times\tilde S_I*}\Gamma_{X_{T_I}}
E_{usu}(\Omega^{\bullet}_{Y_2\times Y_1\times\tilde S_I/Y_2\times\tilde S_I},F_b)[-d_{Y_2}-d_{\tilde S_I}],$ \\ $w_{IJ}(X_T/T))$}.
\end{tikzcd}
\end{equation*} 
\item[(ii)] Let $g:T\to S$ a morphism with $T,S\in\SmVar(\mathbb C)$.
Let $f:X\to S$ a  morphism with $X\in\Var(\mathbb C)$. Assume that there is a factorization 
$f:X\xrightarrow{l}Y\times S\xrightarrow{p_S} S$, with $Y\in\SmVar(\mathbb C)$, $l$ a closed embedding and $p_S$ the projection.
Consider $F(X/S):=p_{S,\sharp}\Gamma_X^{\vee}\mathbb Z(Y\times S/Y\times S)$ and
the isomorphism in $C(\Var(\mathbb C)^{sm}/S)$
\begin{eqnarray*}
T(f,g,F(X/S)):g^*F(X/S):=g^*p_{S,\sharp}\Gamma_X^{\vee}\mathbb Z(Y\times S/Y\times S)\xrightarrow{\sim} \\
p_{T,\sharp}\Gamma_{X_T}^{\vee}\mathbb Z(Y\times T/Y\times T)[d_{Y}]=:F(X_T/T).
\end{eqnarray*}
which gives in $\DA(S)$ the isomorphism $T(f,g,F(X/S)):g^*M(X/S)\xrightarrow{\sim}M(X_T/T)$. 
Then, the following diagram in $D_{Ofil,\infty}(T)$ commutes
\begin{equation*}
\begin{tikzcd}
Lg^{*mod[-]}\mathcal F_{S,an}^{GM}(M(X/S))
\ar[rrr,"T^O(g{,}\mathcal F^{GM})(M(X/S))"]\ar[d,"I^{GM}(X/S)"] & \, & \, & 
\mathcal F_{T,an}^{GM}(M(X_T/T))\ar[d,"I^{GM}(X_T/T)"] \\
g^{*mod}L_O(p_{S*}\Gamma_{X}E_{usu}(\Omega^{\bullet}_{Y\times S/S},F_b)[-d_T]
\ar[d,"T_w(\otimes{,}\gamma)(O_{Y\times S})"]\ar[rrr,"(T(g\times I{,}\gamma)(-)\circ T^O_w(g{,}p_S))"] & \, & \, &
p_{Y\times T*}\Gamma_{X_T}E_{usu}(\Omega^{\bullet}_{Y\times T/T},F_b)[-d_T]
\ar[d,"T_w(\otimes{,}\gamma)(O_{Y\times T})"] \\
Lg^{*mod}\int^{FDR}_{p_S}(\Gamma_XE_{usu}(O_{Y\times S},F_b)[-d_Y-d_T]\ar[rrr,"T^{\mathcal Dmod}(p_S{,}f)(-)"] & \, & \, & 
\int^{FDR}_{p_T}(\Gamma_{X_T}E_{usu}(O_{Y\times T},F_b))[-d_Y-d_T].
\end{tikzcd}
\end{equation*}
\end{itemize}
\end{prop}

\begin{proof}
Follows immediately from definition.
\end{proof}

We have the following theorem:

\begin{thm}\label{mainthmGMan}
\begin{itemize}
\item[(i)]Let $g:T\to S$ is a morphism with $T,S\in\Var(\mathbb C)$. 
Assume there exist a factorization $g:T\xrightarrow{l}Y\times S\xrightarrow{p_S}$
with $Y\in\SmVar(\mathbb C)$, $l$ a closed embedding and $p_S$ the projection. 
Let $S=\cup_{i=1}^lS_i$ be an open cover such that 
there exists closed embeddings $i_i:S_i\hookrightarrow\tilde S_i$ with $\tilde S_i\in\SmVar(\mathbb C)$.
Then, for $M\in\DA_c(S)$
\begin{eqnarray*}
T(g,\mathcal F_{an}^{GM})(M):Rg^{*mod[-],\Gamma}\mathcal F_{S,an}^{GM}(M)\to\mathcal F_{T,an}^{GM}(g^*M)
\end{eqnarray*}
is an isomorphism in $D_{O_Tfil,\mathcal D^{\infty},\infty}(T/(Y\times\tilde S_I))$.
\item[(ii)]Let $g:T\to S$ is a morphism with $T,S\in\SmVar(\mathbb C)$. Then, for $M\in\DA_c(S)$
\begin{eqnarray*}
T(g,\mathcal F_{an}^{GM})(M):Lg^{*mod[-]}\mathcal F_{S,an}^{GM}(M)\to\mathcal F_{T,an}^{GM}(g^*M)
\end{eqnarray*}
is an isomorphism in $D_{O_T}(T)$.
\end{itemize}
\end{thm}

\begin{proof}  
\noindent(i):Follows from proposition \ref{TpDRQGMan}.

\noindent(ii):
\noindent: \textbf{First proof} : Follows from proposition \ref{TgGMpropan}, proposition \ref{keysing1} 
and proposition \ref{PDmod1singan}.

\noindent: \textbf{Second proof} : In the analytic case only, we can give a direct proof of this proposition :
Indeed, let $g:T\to S$ is a morphism with $T,S\in\AnSp(\mathbb C)$ and let $h:U\to S$ a smooth morphism with $U\in\AnSp(\mathbb C$, then,
\begin{equation*}
T^O_{\omega}(g,h):g^{*mod}L_{D^{\infty}}h_*E(\Omega^{\bullet}_{U/S},F)\to h'_*E(\Omega^{\bullet}_{U_T/T},F)
\end{equation*}
is an equivalence usu local : consider the following commutative diagram
\begin{equation*}
\xymatrix{g^{*mod}L_O(h_*E(\mathbb Z_U)\otimes O_S)\ar[rrr]^{T(g,h)(E(\mathbb Z_{U}))}\ar[d]_{g^{*mod}L_OT(h,\otimes)(-,-)} 
& \, & \, & h'_*E(\mathbb Z_{U_T})\otimes O_T\ar[d]^{T(h',\otimes)(-,-)} \\
g^{*mod}L_Oh_*E(h^*O_S)\ar[rrr]^{T^{mod}(g,h)(h^*O_S)}\ar[d]_{g^{*mod}L_Oh_*E(\iota_{U/S})} 
& \, & \, & h'_*E(h^{'*}O_T)\ar[d]^{h'_*E(\iota_{U_T/T})} \\
g^{*mod}L_Oh_*E(\Omega^{\bullet}_{U/S})\ar[rrr]^{T^O_{\omega}(g,h)} 
& \, & \, & h'_*E(\Omega^{\bullet}_{U_T/T})},
\end{equation*}
then,
\begin{itemize}
\item the maps $T(h',\otimes)(-,-)$ and $T(h,\otimes)(-,-)$ are usu local equivalence by proposition \ref{projformula},
\item since $h:U\to S$ is a smooth morphism, the inclusion $\iota_{U/S}:h^*O_S\to\Omega^{\bullet}_{U/S}$ is a quasi-isomorphism,
\item since $h':U_T\to T$ is a smooth morphism, the inclusion $\iota_{U_T/T}:h^*O_T\to\Omega^{\bullet}_{U_T/T}$ is a quasi-isomorphism,
\item since $U$, $U_T$, $S$, $T$ are paracompact topological spaces (in particular Hausdorf),
$T(g,h)(E(\mathbb Z_{U}):g^*h_*E(\mathbb Z_U)\to h'_*E(\mathbb Z_{U_T})$ is a quasi-isomorphism.
\end{itemize}
This fact, together with lemma \ref{homQomegaGMan}, proves the proposition.
\end{proof}

We finish this subsection by some remarks on the absolute case and on a particular case of the relative case:

\begin{prop}\label{FFXDan}
\begin{itemize}
\item[(i)]Let $X\in\PSmVar(\mathbb C)$ and $D=\cup D_i\subset X$ a normal crossing divisor. 
Consider the open embedding $j:U:=X\backslash D\hookrightarrow X$. Then, 
\begin{itemize}
\item The map in $D_{fil,\infty}(\mathbb C)$ 
\begin{eqnarray*}
\Hom((0,\ad(j^*,j_*)(\mathbb Z(X/X)),E_{usu}(\Omega^{\bullet},F_b)): \\
\mathcal F_{an}^{GM}(\mathbb D(\mathbb Z(U))):=\Hom(L\mathbb D(\mathbb Z(U)),E_{usu}(\Omega^{\bullet},F_b)) \\
\xrightarrow{\sim}
\Hom(\Cone(\mathbb Z(D)\to\mathbb Z(X)),E_{zar}(\Omega^{\bullet},F_b))=\Gamma(X,E_{usu}(\Omega^{\bullet}_X(\nul D),F_b)).
\end{eqnarray*}
is an isomorphism, where we recall $\mathbb D(\mathbb Z(U):=a_{X*}j_*E_{et}(\mathbb Z(U/U))=a_{U*}E_{et}(\mathbb Z(U/U))$, 
\item $\mathcal F_{an}^{GM}(\mathbb Z(U))=\Gamma(U,E_{usu}\Omega_U^{\bullet},F_b)\in D_{fil,\infty}(\mathbb C)$
is NOT isomorphic to $\Gamma(X,E_{usu}(\Omega^{\bullet}_X(\log D),F_b)$ in $D_{fil,\infty}(\mathbb C)$ in general.
For exemple $U$ is affine, then $U^{an}$ is Stein so that $H^n(U,\Omega_U^p)=0$ for all $p\in\mathbb N$, $p\neq 0$, 
so that the $E^{p,q}_{\infty}(\Gamma(U,E_{usu}(\Omega_U^{\bullet},F_b)))$ are NOT isomorphic to 
$E^{p,q}_{\infty}(\Gamma(X,E_{usu}(\Omega^{\bullet}_X(\log D),F_b)))$ in this case. 
In particular, the map,  
\begin{equation*}
j^*:=\ad(j^*,j_*)(-):H^n\Gamma(X,E_{usu}(\Omega^{\bullet}_X(\log D)))\xrightarrow{\sim}H^n\Gamma(U,E_{zar}(\Omega^{\bullet}_U)) 
\end{equation*}
which is an isomorphism in $D(\mathbb C)$ (i.e. if we forgot filtrations), gives embeddings  
\begin{eqnarray*}
j^*:=\ad(j^*,j_*)(-):F^pH^n(U,\mathbb C):=F^pH^n\Gamma(X,E_{usu}(\Omega^{\bullet}_X(\log D),F_b))
\hookrightarrow F^pH^n\Gamma(U,E_{usu}(\Omega^{\bullet}_U,F_b))
\end{eqnarray*}
which are NOT an isomorphism in general for $n,p\in\mathbb Z$. Note that, since $a_U:U\to\left\{\pt\right\}$ is not proper,
\begin{equation*}
[\Delta_U]:\mathbb Z(U)\to a_{U*}E_{et}(\mathbb Z(U/U))[2d_U]
\end{equation*}
is NOT an equivalence $(\mathbb A^1,et)$ local.
\item Let $Z\subset X$ a smooth subvariety and denote $U:=X\backslash Z$ the open complementary.
Denote $M_Z(X):=\Cone(M(U)\to M(X))\in\DA(\mathbb C)$. The map in $D_{fil,\infty}(\mathbb C)$
\begin{eqnarray*}
\Hom(G(X,Z),E_{usu}(\Omega^{\bullet},F_b))^{-1}: 
\mathcal F_{an}^{GM}(M_Z(X)):=\Hom(a_{X\sharp}\Gamma^{\vee}_Z\mathbb Z(X/X)),E_{usu}(\Omega^{\bullet},F_b))\xrightarrow{\sim} \\
\Gamma(X,\Gamma_ZE_{usu}(\Omega_X^{\bullet},F_b))=\Gamma_Z(X,E_{usu}(\Omega_X^{\bullet},F_b)) \\
\xrightarrow{\sim}\mathcal F_{an}^{GM}(M(Z)(c)[2c])=\Gamma(Z,E_{usu}(\Omega_Z^{\bullet},F_b))(-c)[-2c]
\end{eqnarray*}
is an isomorphism,
where $c=\codim(Z,X)$ and $G(X,Z):a_{X\sharp}\Gamma^{\vee}_Z\mathbb Z(X/X)\to\mathbb Z(Z)(c)[2c]$ is the Gynsin morphism.
\item Let $D\subset X$ a smooth divisor and denote $U:=X\backslash Z$ the open complementary
Note that the canonical distinguish triangle in $\DA(\mathbb C)$
\begin{equation*}
M(U)\xrightarrow{\ad(j_{\sharp},j^*)(\mathbb Z(X/X))}M(X)\xrightarrow{\gamma^{\vee}_Z(\mathbb Z(X/X))} M_D(X)\to M(U)[1]
\end{equation*}
give a canonical triangle in $D_{fil,\infty}(\mathbb C)$
\begin{eqnarray*}
\mathcal F^{GM}(M_D(X))\xrightarrow{\mathcal F^{GM}(\gamma^{\vee}_Z(\mathbb Z(X/X)))}
\mathcal F^{GM}(M(X))\xrightarrow{\mathcal F^{GM}(\ad(j_{\sharp},j^*)(\mathbb Z(X/X))} 
\mathcal F^{GM}(M(U))\to\mathcal F^{GM}(M_D(X))[1],
\end{eqnarray*}
which is NOT the image of a distinguish triangle in $\pi(D(MHM(\mathbb C)))$, as
$\mathcal F^{GM}(M(U))\notin\pi(D(MHM(\mathbb C)))$ since the morphism
\begin{eqnarray*}
j^*:=\ad(j^*,j_*)(-):H^n\Gamma(X,E_{usu}(\Omega^{\bullet}_X(\log D),F_b))\to H^n\Gamma(U,E_{usu}(\Omega^{\bullet}_U,F_b)) 
\end{eqnarray*}
are not strict.
Note that if $U:=S\backslash D$ is affine, then by the exact sequence in $C(\mathbb Z)$
\begin{equation*}
0\to \Gamma_Z(X,E_{usu}(\Omega^p_X))\to\Gamma(X,E_{usu}(\Omega^p_X))\to \Gamma(U,E_{usu}(\Omega^p_U))\to 0
\end{equation*}
we have $H^q\Gamma_Z(X,E_{usu}(\Omega^p_X))=H^q(\Gamma(X,E_{usu}(\Omega^p_X)))$.
\end{itemize}
In particular, the maps
\begin{eqnarray*}
j^*:=\ad(j^*,j_*)(-):\Gamma(X,E_{usu}(\Omega^{\bullet}_X(\log D),F_b))\to\Gamma(U,E_{usu}(\Omega^{\bullet}_U,F_b)) 
\end{eqnarray*}
and
\begin{eqnarray*}
j^*:=\ad(j^*,j_*)(-):
\Cone(\Gamma(X,E_{usu}(\Omega^{\bullet}_X,F_b))\to\Gamma(X,E_{usu}(\Omega^{\bullet}_X(\log D),F_b)))\to \\
\Cone(\Gamma(X,E_{usu}(\Omega^{\bullet}_X,F_b))\to\Gamma(U,E_{usu}(\Omega^{\bullet}_U,F_b)))=: 
\Gamma(X,\Gamma_ZE_{usu}(\Omega_X^{\bullet},F_b))
\end{eqnarray*}
are quasi-isomorphism (i.e. if we forgot filtrations),
but the first one is NOT an $\infty$-filtered quasi-isomorphism
whereas the second one is an $\infty$-filtered quasi-isomorphism (recall that for $r>1$ the $r$-filtered quasi-isomorphism
does NOT satisfy the 2 of 3 property for morphism of canonical triangles : see section 2.1).
\item[(ii)] More generally, let $f:X\to S$ a smooth projective morphism with $S,X\in\SmVar(\mathbb C)$.
Let $D=\cup D_i\subset X$ a normal crossing divisor such that $f_{|D_I}:=f\circ i_I:D_I\to S$ are SMOOTH morphisms
(note that it is a very special case), with $i_I:D_I\hookrightarrow X$ the closed embeddings.
Consider the open embedding $j:U:=X\backslash D\hookrightarrow X$ and $h:=f\circ j:U\to S$. 
\begin{itemize}
\item The map in $D_{\mathcal Dfil,\infty}(S)$ 
\begin{eqnarray*}
\Hom((0,\ad(j^*,j_*)(\mathbb Z(X/X))),E_{et}(\Omega^{\bullet}_{/S},F_b)): \\
\mathcal F_{S,an}^{GM}(\mathbb D(\mathbb Z(U/S))):=\Hom(L\mathbb D(\mathbb Z(U/S)),E_{usu}(\Omega^{\bullet}_{/S},F_b)) \\
\xrightarrow{\sim}
\Hom(\Cone(\mathbb Z(D)\to\mathbb Z(X)),E_{zar}(\Omega^{\bullet}_{/S},F_b))=f_*E_{usu}(\Omega^{\bullet}_{X/S}(\nul D),F_b).
\end{eqnarray*}
is an isomorphism, where we recall $\mathbb D(\mathbb Z(U):=f_*j_*E_{et}(\mathbb Z(U/U))=h_*E_{usu}(\mathbb Z(U/U))$, 
\item $\mathcal F_{S,an}^{GM}(\mathbb Z(U/S))=h_*E_{usu}\Omega_{U/S}^{\bullet},F_b)\in D_{\mathcal Dfil,\infty}(S)$
is NOT isomorphic to $f_*E_{usu}(\Omega^{\bullet}_{X/S}(\log D),F_b)$ in $D_{\mathcal Dfil,\infty}(S)$ in general.
In particular, the map,  
\begin{equation*}
j^*:=\ad(j^*,j_*)(-):H^nf_*E_{usu}(\Omega^{\bullet}_{X/S}(\log D))\xrightarrow{\sim}H^nh_*E_{usu}(\Omega^{\bullet}_{U/S}) 
\end{equation*}
which is an isomorphism in $D_{\mathcal D}(S)$ (i.e. if we forgot filtrations), gives embeddings  
\begin{eqnarray*}
j^*:=\ad(j^*,j_*)(-):F^pH^nh_*\mathbb C_U:=F^pH^nf_*E_{usu}(\Omega^{\bullet}_{X/S}(\log D),F_b)
\hookrightarrow F^pH^nh_*E_{usu}(\Omega^{\bullet}_{U/S},F_b)
\end{eqnarray*}
which are NOT an isomorphism in general for $n,p\in\mathbb Z$. Note that, since $a_U:U\to\left\{\pt\right\}$ is not proper,
\begin{equation*}
[\Delta_U]:\mathbb Z(U/S)\to h_*E_{usu}(\mathbb Z(U/U))[2d_U]
\end{equation*}
is NOT an equivalence $(\mathbb A^1,et)$ local.
\item Let $Z\subset X$ a subvariety and denote $U:=X\backslash Z$ the open complementary. 
Denote $M_Z(X/S):=\Cone(M(U/S)\to M(X/S))\in\DA(S)$. 
If $f_{|Z}:=f\circ i_Z:Z\to S$ is a SMOOTH morphism, the map in $D_{\mathcal Dfil,\infty}(S)$
\begin{eqnarray*}
\Hom(G(X,Z),E_{usu}(\Omega^{\bullet},F_b)): \\
\mathcal F_{S,an}^{GM}(M_Z(X/S)):=\Hom(f_{\sharp}\Gamma^{\vee}_Z\mathbb Z(X/X)),E_{usu}(\Omega^{\bullet}_{/S},F_b)) 
\xrightarrow{\sim}f_*\Gamma_ZE_{usu}(\Omega_{X/S}^{\bullet},F_b)) \\
\xrightarrow{\sim}\mathcal F_S^{GM}(M(Z/S)(c)[2c])=f_{Z*}E_{usu}(\Omega_{Z/S}^{\bullet},F_b)(-c)[-2c]
\end{eqnarray*}
is an isomorphism,
where $c=\codim(Z,X)$ and $G(X,Z):f_{\sharp}\Gamma^{\vee}_Z\mathbb Z(X/X)\to\mathbb Z(Z/S)(c)[2c]$ is the Gynsin morphism.
\item Let $D\subset X$ a smooth divisor and denote $U:=X\backslash Z$ the open complementary
Note that the canonical distinguish triangle in $\DA(S)$
\begin{equation*}
M(U/S)\xrightarrow{\ad(j_{\sharp},j^*)(\mathbb Z(X/X))}M(X/S)\xrightarrow{\gamma^{\vee}_Z(\mathbb Z(X/X))} M_D(X/S)\to M(U/S)[1]
\end{equation*}
give a canonical triangle in $D_{\mathcal Dfil,\infty}(S)$
\begin{eqnarray*}
\mathcal F_{S,an}^{GM}(M_D(X/S))\xrightarrow{\mathcal F^{GM}(\gamma^{\vee}_Z(\mathbb Z(X/X)))}
\mathcal F_{S,an}^{GM}(M(X/S))\xrightarrow{\mathcal F^{GM}(\ad(j_{\sharp},j^*)(\mathbb Z(X/X))}  
\mathcal F_{S,an}^{GM}(M(U/S)) \\
\to\mathcal F_{S,an}^{GM}(M_D(X/S))[1],
\end{eqnarray*}
which is NOT the image of a distinguish triangle in $\pi_S(D(MHM(S)))$.
\end{itemize}
\end{itemize}
\end{prop}

\begin{proof}
Similar to the proof of theorem \ref{FFXD}.
\end{proof}

\begin{defi}
Let $S\in\Var(\mathbb C)$ and $S=\cup_{i=1}^l S_i$ an open affine covering and denote, 
for $I\subset\left[1,\cdots l\right]$, $S_I=\cap_{i\in I} S_i$ and $j_I:S_I\hookrightarrow S$ the open embedding.
Let $i_i:S_i\hookrightarrow\tilde S_i$ closed embeddings, with $\tilde S_i\in\SmVar(\mathbb C)$. 
We have, for $M,N\in\DA(S)$ and $F,G\in C(\Var(\mathbb C)^{sm}/S)$ such that 
$M=D(\mathbb A^1,et)(F)$ and $N=D(\mathbb A^1,et)(G)$, 
the following transformation map in $D_{Ofil,\mathcal D^{\infty}}(S/(\tilde S_I))$
\begin{eqnarray*}
T(\mathcal F_{S,an}^{GM},\otimes)(M,N): \\
\mathcal F_{S,an}^{GM}(M)\otimes^L_{O_S}\mathcal F_{S,an}^{GM}(N):= 
(e(\tilde S_I)_*\mathcal Hom(\An_{\tilde S_I}^*L(i_{I*}j_I^*F),
E_{usu}(\Omega^{\bullet}_{/\tilde S_I},F_b)),u_{IJ}(F))\otimes_{O_S} \\
(e(\tilde S_I)_*\mathcal Hom(\An_{\tilde S_I}^*L(i_{I*}j_I^*G),
E_{usu}(\Omega^{\bullet}_{/\tilde S_I},F_b)),u_{IJ}(G)) \\
\xrightarrow{=}
((e(\tilde S_I)_*\mathcal Hom(\An_{\tilde S_I}^*L(i_{I*}j_I^*F),
E_{usu}(\Omega^{\bullet}_{/\tilde S_I},F_b))\otimes_{O_{\tilde S_I}} \\
e(\tilde S_I)_*\mathcal Hom(\An_{\tilde S_I}^*L(i_{I*}j_I^*G),
E_{usu}(\Omega^{\bullet}_{/\tilde S_I},F_b))),u_{IJ}(F)\otimes u_{IJ}(G)) \\
\xrightarrow{(T(\otimes,\Omega_{/\tilde S_I})(L(i_{I*}j_I^*F),L(i_{I*}j_I^*G)))} \\
(e(\tilde S_I)_*\mathcal Hom(\An_{\tilde S_I}^*L(i_{I*}j_I^*F)\otimes\An_{\tilde S_I}^*L(i_{I*}j_I^*G),
E_{et}(\Omega^{\bullet}_{/\tilde S_I},F_b)),v_{IJ}(F\otimes G)) \\
\xrightarrow{=} 
(e(\tilde S_I)_*\mathcal Hom(\An_{\tilde S_I}^*L(i_{I*}j_I^*(F\otimes G),E_{usu}(\Omega^{\bullet}_{/\tilde S_I},F_b))),
u_{IJ}(F\otimes G))=:\mathcal F_{S,an}^{GM}(M\otimes N)
\end{eqnarray*}
\end{defi}

We have in the analytical case the following :
\begin{prop}
Let $S\in\Var(\mathbb C)$. Then, for $M,N\in\DA_c(S)$
\begin{eqnarray*}
T(\otimes,\mathcal F_{S,an}^{GM})(M,N):
\mathcal F_{S,an}^{GM}(M\otimes N)\xrightarrow{\sim}\mathcal F_{S,an}^{GM}(M)\otimes^L_{O_S}\mathcal F_{S,an}^{GM}(N)
\end{eqnarray*}
is an isomorphism.
\end{prop}

\begin{proof}
Asumme first that $S$ is smooth. 
Let $h_1:U_1\to S$ and $h_2:U_2\to S$ smooth morphisms with $U_1,U_2\in\Var(\mathbb C)$ and consider $h_{12}:U_1\times_S U_2\to S$.
We then have by lemma \ref{homQomega2GMan} the following commutative diagram
\begin{equation*}
\xymatrix{
e(S)_*\mathcal Hom(\mathbb Z(U_1/S)\otimes\mathbb Z(U_2/S),E(\Omega^{\bullet}_{/S},F_b))
\ar[r]^{T(\otimes,\mathcal F_S^{GM}(M(U_1/S),M(U_2/S)))}\ar[d]^{=} & 
e(S)_*\mathcal Hom(\mathbb Z(U_1/S),E(\Omega^{\bullet}_{/S},F_b))\otimes_{O_S}
e(S)_*\mathcal Hom(\mathbb Z(U_2/S),E(\Omega^{\bullet}_{/S},F_b))\ar[d]^{=} \\
h_{12*}E(\Omega_{U_1\times_SU_2/S},F)\ar[r]^{Ew_{(U_1,U_2)/S}}\ar[d]^{\iota\otimes\iota} & 
h_{1*}E(\omega_{U_1/S},F_b)\otimes_{O_S}h_{2*}E(\Omega_{U_2/S},F_b)\ar[d]^{\iota} \\
h_{12*}E(h_{12}^*O_S)\ar[r]^{Ew_{(U_1,U_2)/S}}\ar[d]^{T(h_{12},\otimes)(O_S,\mathbb Z_{U_{12}})} & 
h_{1*}E(h_1^*O_S)\otimes_{O_S}h_{2*}E(h_2^*O_S)
\ar[d]^{T(h_1,\otimes)(O_S,\mathbb Z_{U_1})\otimes T(h_2,\otimes)(O_S,\mathbb Z_{U_1})} \\
h_{12*}E(\mathbb Z_{U_{12}})\otimes O_S\ar[r]^{Ew_{(U_1,U_2)/S}} & 
(h_{1*}E(\mathbb Z_{U_1})\otimes_{O_S} h_{2*}E(\mathbb Z_{U_2}))}
\end{equation*}
Since $U_1,U_2\in\AnSp(\mathbb C)$ are locally contractible topological spaces, 
the lower row is an equivalence usu local by Kunneth formula for topological spaces (see section 2). 
This proves the proposition in the case $S$ is smooth.
Let $S\in\Var(\mathbb C)$ and $S=\cup_{i=1}^l S_i$ an open cover such that there exist closed embeddings 
$i_i:S_i\hookrightarrow\tilde S_i$ with $\tilde S_i\in\SmVar(\mathbb C)$.
By definition, for $F,G\in C(\Var(\mathbb C)^{sm}/S)$ such that $M=D(\mathbb A^1,et)(F)$ and $N=D(\mathbb A^1,et)(G)$, 
\begin{eqnarray*}
T(\otimes,\mathcal F_{S,an}^{GM}(M,N)): \\
e(\tilde S_I)_*\mathcal Hom(\An_{\tilde S_I}^*L(i_{I*}j_I^*(F\otimes G),
E_{usu}(\Omega^{\bullet}_{/\tilde S_I},F_b))),u_{IJ}(F\otimes G)) \\
\xrightarrow{(T(\otimes,\Omega_{/\tilde S_I})(\An_{\tilde S_I}^*L(i_{I*}j_I^*F),\An_{\tilde S_I}^*L(i_{I*}j_I^*G)))} \\
(e(\tilde S_I)_*\mathcal Hom(\An_{\tilde S_I}^*L(i_{I*}j_I^*F),
E_{usu}(\Omega^{\bullet}_{/\tilde S_I},F_b)),u_{IJ}(F))\otimes_{O_S} \\
(e(\tilde S_I)_*\mathcal Hom(\An_{\tilde S_I}^*L(i_{I*}j_I^*G),
E_{usu}(\Omega^{\bullet}_{/\tilde S_I},F_b)),u_{IJ}(G))
\end{eqnarray*}
Since $L(i_{I*}j_I^*F),L(i_{I*}j_I^*G)\in\DA_c(\tilde S_I)$, by the smooth case applied to $\tilde S_I$ for each $I$,
$T(\otimes,\mathcal F_{S,an}^{FDR}(M,N))$ is an equivalence usu local.
\end{proof}

\subsubsection{The analytic filtered De Rham realization functor}

Recall from section 2 that, for $S\in\Var(\mathbb C)$ we have the following commutative diagrams of sites
\begin{equation*}
\xymatrix{\AnSp(\mathbb C)^2/S\ar[rr]^{\mu_S}\ar[dd]_{\An_S}\ar[rd]^{\rho_S} & \, & 
\AnSp(\mathbb C)^{2,pr}/S\ar[dd]^{\An_S}\ar[rd]^{\rho_S} & \, \\
\, & \AnSp(\mathbb C)^{2,sm}/S\ar[rr]^{\mu_S}\ar[dd]_{\An_S} & \, & \AnSp(\mathbb C)^{2,smpr}/S\ar[dd]^{\An_S} \\
\Var(\mathbb C)^2/S\ar[rr]^{\mu_S}\ar[rd]^{\rho_S} & \, &  \Var(\mathbb C)^{2,smpr}/S\ar[rd]^{\rho_S} & \, \\
\, & \Var(\mathbb C)^2/S\ar[rr]^{\mu_S} & \, &  \Var(\mathbb C)^{2,smpr}/S & \,}
\end{equation*}
and
\begin{equation}
\xymatrix{\AnSp(\mathbb C)^{2,pr}/S\ar[rr]^{\Gr_S^{12}}\ar[dd]_{\An_S}\ar[rd]^{\rho_S} & \, & 
\AnSp(\mathbb C)/S\ar[dd]^{\An_S}\ar[rd]^{\rho_S} & \, \\
\, & \AnSp(\mathbb C)^{2,smpr}/S\ar[rr]^{\Gr_S^{12}}\ar[dd]_{\An_S} & \, & \AnSp(\mathbb C)^{sm}/S\ar[dd]^{\An_S} \\
\Var(\mathbb C)^{2,pr}/S\ar[rr]^{\Gr_S^{12}}\ar[rd]^{\rho_S} & \, & \Var(\mathbb C)/S\ar[rd]^{\rho_S} & \, \\
\, & \Var(\mathbb C)^{2,sm}/S\ar[rr]^{\Gr_S^{12}} & \, & \Var(\mathbb C)^{sm}/S},
\end{equation}
and that for $f:T\to S$ a morphism with $T,S\in\Var(\mathbb C)$ we have the following commutative diagrams of site,
\begin{equation*}
\xymatrix{
\AnSp(\mathbb C)^2/T^{an}\ar[rr]^{\An_T}\ar[dd]^{P(f)}\ar[rd]^{\rho_T} & \, & 
\Var(\mathbb C)^2/T\ar[dd]^{P(f)}\ar[rd]^{\rho_T} & \, \\  
 \, & \AnSp(\mathbb C)^{2,sm}/T^{an}\ar[rr]^{\An_T}\ar[dd]^{P(f)} & \, & \Var(\mathbb C)^{2,sm}/T\ar[dd]^{P(f)} \\  
\AnSp(\mathbb C)^2/S^{an}\ar[rr]^{\An_S}\ar[rd]^{\rho_S} & \, & \Var(\mathbb C)^2/S\ar[rd]^{\rho_S} & \, \\  
\, &  \AnSp(\mathbb C)^{2,sm}/S^{an}\ar[rr]^{\An_S} & \, & \Var(\mathbb C)^{2,sm}/S}.  
\end{equation*}

For $s:\mathcal I\to\mathcal J$ a functor, with $\mathcal I,\mathcal J\in\Cat$, and
$f_{\bullet}:T_{\bullet}\to S_{s(\bullet)}$ a morphism with 
$T_{\bullet}\in\Fun(\mathcal J,\Var(\mathbb C))$ and $S_{\bullet}\in\Fun(\mathcal I,\Var(\mathbb C))$, 
we have then the commutative diagram of sites \ref{mufIJan}
\begin{equation*}
\xymatrix{\AnSp(\mathbb C)^2/T_{\bullet}\ar[rr]^{\mu_{T_{\bullet}}}\ar[dd]_{P(f_{\bullet})}\ar[rd]^{\rho_{T_{\bullet}}} & \, & 
\AnSp(\mathbb C)^{2,pr}/T_{\bullet}\ar[dd]^{P(f_{\bullet})}\ar[rd]^{\rho_{T_{\bullet}}} & \, \\
\, & \AnSp(\mathbb C)^{2,sm}/T_{\bullet}\ar[rr]^{\mu_{T_{\bullet}}}\ar[dd]_{P(f_{\bullet})} & \, & 
\AnSp(\mathbb C)^{2,smpr}/T_{\bullet}\ar[dd]^{P(f_{\bullet})} \\
\AnSp(\mathbb C)^2/S_{\bullet}\ar[rr]^{\mu_{S_{\bullet}}}\ar[rd]^{\rho_{S_{\bullet}}} & \, & 
\AnSp(\mathbb C)^{2,pr}/S_{\bullet}\ar[rd]^{\rho_{S_{\bullet}}} & \, \\
\, & \AnSp(\mathbb C)^{2,sm}/S_{\bullet}\ar[rr]^{\mu_{S_{\bullet}}} & \, & \AnSp(\mathbb C)^{2,smpr}/S_{\bullet}}.
\end{equation*}
and the following commutative diagrams of site,
\begin{equation*}
\xymatrix{
\AnSp(\mathbb C)^2/T_{\bullet}^{an}\ar[rr]^{\An_{T_{\bullet}}}\ar[dd]^{P(f_{\bullet})}\ar[rd]^{\rho_{T_{\bullet}}} & \, & 
\Var(\mathbb C)^2/T_{\bullet}\ar[dd]^{P(f_{\bullet})}\ar[rd]^{\rho_{T_{\bullet}}} & \, \\  
 \, & \AnSp(\mathbb C)^{2,sm}/T_{\bullet}^{an}\ar[rr]^{\An_{T_{\bullet}}}\ar[dd]^{P(f_{\bullet})} & \, & 
\Var(\mathbb C)^{2,sm}/T_{\bullet}\ar[dd]^{P(f_{\bullet})} \\  
\AnSp(\mathbb C)^2/S_{\bullet}^{an}\ar[rr]^{\An_{S_{\bullet}}}\ar[rd]^{\rho_{S_{\bullet}}} & \, & 
\Var(\mathbb C)^2/S_{\bullet}\ar[rd]^{\rho_{S_{\bullet}}} & \, \\  
\, &  \AnSp(\mathbb C)^{2,sm}/S_{\bullet}^{an}\ar[rr]^{\An_{S_{\bullet}}} & \, & \Var(\mathbb C)^{2,sm}/S_{\bullet}}.  
\end{equation*}

\begin{defi}\label{wtildewan}
\begin{itemize}
\item[(i)] For $S\in\SmVar(\mathbb C)$, we consider, using definition \ref{wtildew}(i),
the filtered complexes of presheaves 
\begin{eqnarray*} 
(\Omega^{\bullet,\Gamma,pr}_{/S^{an}},F_{DR})\in C_{D^{\infty}_Sfil}(\Var(\mathbb C)^{2,smpr}/S) 
\end{eqnarray*}
given by, 
\begin{itemize}
\item for $(Y\times S,Z)/S=((Y\times S,Z),p)\in\Var(\mathbb C)^{2,smpr}/S$, 
\begin{eqnarray*}
(\Omega^{\bullet,\Gamma,pr}_{/S^{an}}((Y\times S,Z)/S),F_{DR}):=
((\Omega^{\bullet}_{(Y\times S)^{an}/S^{an}},F_b)\otimes_{O_{(Y\times S)^{an}}}
(\Gamma^{\vee,Hdg}_Z(O_{Y\times S},F_b))^{an})((Y\times S)^{an}) 
\end{eqnarray*}
with the structure of $p^*D_S$ module given by proposition \ref{DRhUS}.
\item for $g:(Y_1\times S,Z_1)/S=((Y_1\times S,Z_1),p_1)\to (Y\times S,Z)/S=((Y\times S,Z),p)$ 
a morphism in $\Var(\mathbb C)^{2,smpr}/S$, 
\begin{eqnarray*}
\Omega^{\bullet,\Gamma,pr}_{/S^{an}}(g):=(\Omega^{\bullet,\Gamma,pr}_{/S}(g))^{an}:
((\Omega^{\bullet}_{(Y\times S)^{an}/S^{an}},F_b)\otimes_{O_{(Y\times S)^{an}}}
(\Gamma^{\vee,Hdg}_Z(O_{Y\times S},F_b))^{an})((Y\times S)^{an})\to \\ 
((\Omega^{\bullet}_{(Y_1\times S)^{an}/S^{an}},F_b)\otimes_{O_{(Y_1\times S)^{an}}}
(\Gamma^{\vee,Hdg}_{Z_1}(O_{Y_1\times S},F_b))^{an})((Y_1\times S)^{an}). 
\end{eqnarray*}
\end{itemize}
For $S\in\SmVar(\mathbb C)$, we get the filtered complexes of presheaves 
\begin{eqnarray*} 
(\Omega^{\bullet,\Gamma,pr,an}_{/S^{an}},F_{DR}):=\An_S^{*mod}(\Omega^{\bullet,\Gamma,pr}_{/S^{an}},F_{DR})
:=\An_S^*(\Omega^{\bullet,\Gamma,pr}_{/S^{an}},F_{DR})\otimes_{O_S}O_{S^{an}}
\in C_{D^{\infty}_Sfil}(\AnSp(\mathbb C)^{2,smpr}/S). 
\end{eqnarray*}
\item[(ii)] For $S\in\SmVar(\mathbb C)$, we have the canonical map $C_{O_Sfil,D^{\infty}_S}(\Var(\mathbb C)^{sm}/S)$ 
\begin{eqnarray*} 
\Gr(\Omega_{/S^{an}}):\Gr_{S*}^{12}(\Omega^{\bullet,\Gamma,pr}_{/S^{an}},F_b)\to\An_{S*}(\Omega^{\bullet}_{/S},F_b) 
\end{eqnarray*}
given by
\begin{eqnarray*} 
\Gr(\Omega_{/S^{an}})(U/S):=(\Gr(\Omega_{/S})(U/S))^{an}\otimes m: \\
J_S((\Omega^{\bullet}_{(U\times S)^{an}/S^{an}},F_b)\otimes_{O_{(U\times S)^{an}}}
(\Gamma^{\vee,Hdg}_U(O_{U\times S},F_b))^{an})((U\times S)^{an})
\to(\Omega^{\bullet}_{U^{an}/S^{an}},F_b), 
\end{eqnarray*}
where $\Gr(\Omega_{/S^{an}})(U/S)(\omega\otimes m\otimes P):=P(\Gr(\Omega_{/S})(U/S)(\omega\otimes m))$
with $P\in\Gamma(S,D_S^{\infty})$,
see definition \ref{wtildew}(ii), which gives by adjonction
\begin{eqnarray*} 
\Gr(\Omega_{/S^{an}}):=I(\An_S^{*mod},\An_S)(\Gr(\Omega_{/S^{an}})): 
J_S(\Gr_{S*}^{12}(\Omega^{\bullet,\Gamma,pr,an}_{/S^{an}},F_b))\to(\Omega^{\bullet}_{/S},F_b) 
\end{eqnarray*}
in $C_{O_Sfil,D^{\infty}_S}(\AnSp(\mathbb C)^{sm}/S)$.
\end{itemize}
\end{defi}

\begin{defi}\label{wtildewTan}
For $S\in\SmVar(\mathbb C)$, we have the canonical map in $C_{O_Sfil,D^{\infty}_S}(\Var(\mathbb C)^{2,smpr}/S)$ 
\begin{eqnarray*} 
T(\Omega^{\Gamma}_{/S^{an}}):\An_{S*}\mu_{S*}(\Omega^{\bullet,\Gamma}_{/S^{an}},F_b)\to
(\Omega^{\bullet,\Gamma,pr}_{/S^{an}},F_{DR}) 
\end{eqnarray*}
given by, for $(Y\times S,X)/S=((Y\times S,Z),p)\in\Var(\mathbb C)^{2,smpr}/S$
\begin{eqnarray*} 
T(\Omega^{\Gamma}_{/S^{an}})((Y\times S,Z)/S):=(T(\Omega^{\Gamma}_{/S})((Y\times S,Z)/S))^{an}: \\
(\Omega^{\bullet,\Gamma}_{/S},F_b)(((Y\times S)^{an},Z^{an})/S):=
\mathbb D_{p^*O_S}L_{p^*O}\Gamma_{Z}E_{usu}(\mathbb D_{p^*O_S}L_{p^*O}
(\Omega^{\bullet}_{(Y\times S)^{an}/S^{an}},F_b))((Y\times S)^{an})\to \\
((\Omega^{\bullet}_{(Y\times S)^{an}/S^{an}},F_b)\otimes_{O_{(Y\times S)^{an}}}
(\Gamma^{\vee,Hdg}_{Z}(O_{Y\times S},F_b))^{an})((Y\times S)^{an})=:
(\Omega^{\bullet,\Gamma,pr}_{/S^{an}},F_{DR})((Y\times S,Z)/S),
\end{eqnarray*}
see definition \ref{wtildewT}.
By definition we have $\Gr^O(\Omega_{/S^{an}})=\Gr(\Omega_{/S^{an}})\circ T(\Omega^{\Gamma}_{/S^{an}})$.
\end{defi}

We have the following canonical transformation map given by the pullback of (relative) differential forms:

Let $g:T\to S$ a morphism with $T,S\in\AnSm(\mathbb C)$.
\begin{itemize}
\item We have the canonical morphism in $C_{g^*O_Sfil,g^*D^{\infty}_S}(\AnSp(\mathbb C)^{2,sm}/T)$ 
\begin{eqnarray*}
\Omega^{\Gamma}_{/(T/S)}:g^*(\Omega^{\bullet,\Gamma}_{/S},F_b)\to (\Omega^{\bullet,\Gamma}_{/T},F_b)
\end{eqnarray*}
induced by the pullback of differential forms : for $((V,Z_1)/T)=((V,Z_1),h)\in\AnSp(\mathbb C)^{2,sm}/T$,
\begin{eqnarray*}
\Omega^{\Gamma}_{/(T/S)}((V,Z_1)/T): \\
g^*\Omega^{\bullet,\Gamma}_{/S}((V,Z_1)/T):= 
\lim_{(h:(U,Z)\to S \mbox{sm},g_1:(V,Z_1)\to (U_T,Z_T),h,g)}\Omega^{\bullet,\Gamma}_{/S}((U,Z)/S) \\
\xrightarrow{\Omega^{\bullet,\Gamma}_{/S}(g'\circ g_1)}\Omega^{\bullet,\Gamma}_{/S}((V,Z_1)/S)
\xrightarrow{\Gamma^{\vee,h}_{Z_1}q(Y_1\times T)}\Omega^{\bullet,\Gamma}_{/T}((V,Z_1)/T),
\end{eqnarray*}
where $g':U_T:=U\times_S T\to U$ is the base change map and 
$q:\Omega^{\bullet}_{Y_1\times T/S}\to\Omega^{\bullet}_{Y_1\times T/T}$ is the quotient map.
It induces the canonical morphisms in $C_{g^*O_Sfil,g^*D^{\infty}_S}(\AnSp(\mathbb C)^{2,sm}/T)$ :
\begin{eqnarray*}
E\Omega^{\Gamma}_{/(T/S)}:g^*E_{et}(\Omega^{\bullet,\Gamma}_{/S},F_b)
\xrightarrow{T(g,E_{et})(\Omega^{\bullet,\Gamma}_{/S},F_b)}
E_{et}(g^*(\Omega^{\bullet,\Gamma}_{/S},F_b))
\xrightarrow{E_{et}(\Omega^{\Gamma}_{/(T/S)})}E_{et}(\Omega^{\bullet,\Gamma}_{/T},F_b)
\end{eqnarray*}
\item We have the canonical morphism in $C_{g^*D^{\infty}_Sfil}(\Var(\mathbb C)^{2,smpr}/T)$ 
\begin{eqnarray*}
\Omega^{\Gamma,pr}_{/(T/S)^{an}}:g^*(\Omega^{\bullet,\Gamma,pr}_{/S^{an}},F_{DR})\to 
(\Omega^{\bullet,\Gamma,pr}_{/T^{an}},F_{DR})
\end{eqnarray*}
induced by the pullback of differential forms : 
for $((Y_1\times T,Z_1)/T)=((Y_1\times T,Z_1),p)\in\Var(\mathbb C)^{2,smpr}/T$,
\begin{eqnarray*}
\Omega^{\Gamma,pr}_{/(T/S)^{an}}((Y_1\times T,Z_1)/T): \\
g^*\Omega^{\bullet,\Gamma,pr}_{/S^{an}}((Y_1\times T,Z_1)/T):= 
\lim_{(h:(Y\times S,Z)\to S, \, g_1:(Y_1\times T,Z_1)\to (Y\times T,Z_T),h,g)}
\Omega^{\bullet,\Gamma,pr}_{/S^{an}}((Y\times T,Z)/S) \\
\xrightarrow{\Omega^{\bullet,\Gamma,pr}_{/S^{an}}(g'\circ g_1)}\Omega^{\bullet,\Gamma,pr}_{/S^{an}}((Y_1\times T,Z_1)/S)
\xrightarrow{q(-)((Y_1\times T)^{an})}\Omega^{\bullet,\Gamma,pr}_{/T^{an}}((Y_1\times T,Z_1)/T), 
\end{eqnarray*}
where $g'=(I_Y\times g):Y\times T\to Y\times S$ is the base change map and 
\begin{equation*}
q(M):\Omega_{(Y_1\times T)^{an}/S^{an}}\otimes_{O_{(Y_1\times T)^{an}}}(M,F)\to
\Omega_{(Y_1\times T)^{an}/T^{an}}\otimes_{O_{(Y_1\times T)^{an}}}(M,F)
\end{equation*}
is the quotient map. 
It induces the canonical morphisms in $C_{g^*D_Sfil}(\Var(\mathbb C)^{2,smpr}/T)$ :
\begin{eqnarray*}
E\Omega^{\Gamma,pr}_{/(T/S)}:g^*E_{et}(\Omega^{\bullet,\Gamma,pr}_{/S^{an}},F_{DR})
\xrightarrow{T(g,E)(-)}E_{et}(g^*(\Omega^{\bullet,\Gamma,pr}_{/S^{an}},F_{DR}))
\xrightarrow{E_{et}(\Omega^{\Gamma,pr}_{/(T/S)^{an}})}E_{et}(\Omega^{\bullet,\Gamma,pr}_{/T^{an}},F_{DR})
\end{eqnarray*}
and
\begin{eqnarray*}
E\Omega^{\Gamma,pr}_{/(T/S)^{an}}:g^*E_{zar}(\Omega^{\bullet,\Gamma,pr}_{/S^{an}},F_{DR})
\xrightarrow{T(g,E)(-)} E_{zar}(g^*(\Omega^{\bullet,\Gamma,pr}_{/S^{an}},F_{DR}))
\xrightarrow{E_{zar}(\Omega^{\Gamma,pr}_{/(T/S)^{an}})}E_{zar}(\Omega^{\bullet,\Gamma,pr}_{/T^{an}},F_{DR}).
\end{eqnarray*}
\end{itemize}

\begin{defi}\label{TgDRan}
Let $g:T\to S$ a morphism with $T,S\in\SmVar(\mathbb C)$.
We have, for $F\in C(\Var(\mathbb C)^{2,smpr}/S)$, the canonical transformation in $C_{\mathcal D^{\infty}fil}(T)$ :
\begin{eqnarray*}
T(g,\Omega^{\Gamma,pr}_{/\cdot})(F): 
g^{*mod}L_De(S)_*\Gr^{12}_{S*}\mathcal Hom^{\bullet}(\An_S^*F,E_{et}(\Omega^{\bullet,\Gamma,pr,an}_{/S^{an}},F_{DR})) \\
\xrightarrow{:=} 
(g^*L_De(S)_*\mathcal Hom^{\bullet}(F,E_{et}(\Omega^{\bullet,\Gamma,pr,an}_{/S^{an}},F_{DR})))\otimes_{g^*O_S}O_T \\
\xrightarrow{T(g,\Gr^{12})(-)\circ T(e,g)(-)\circ q}  
e(T)_*\Gr^{12}_{T*}g^*\mathcal Hom^{\bullet}(\An_S^*F,E_{et}(\Omega^{\bullet,\Gamma,pr}_{/S^{an}},F_{DR}))\otimes_{g^*O_S}O_T \\ 
\xrightarrow{(T(g,hom)(-,-)\otimes I)}
e(T)_*\Gr^{12}_{T*}\mathcal Hom^{\bullet}(\An_T^*g^*F,g^*E_{et}(\Omega^{\bullet,\Gamma,pr,an}_{/S^{an}},F_{DR}))
\otimes_{g^*O_S}O_T \\
\xrightarrow{ev(hom,\otimes)(-,-,-)}  
e(T)_*\Gr^{12}_{T*}\mathcal Hom^{\bullet}(\An_T^*g^*F,g^*E_{et}(\Omega^{\bullet,\Gamma,pr,an}_{/S^{an}},F_{DR}))
\otimes_{g^*e(S)^*O_S}e(T)^*O_T \\
\xrightarrow{\mathcal Hom^{\bullet}(\An_T^*g^*F,(E\Omega^{\Gamma,pr}_{/(T/S)}\otimes m))}
e(T)_*\Gr^{12}_{T*}\mathcal Hom^{\bullet}(\An_T^*g^*F,E_{et}(\Omega^{\bullet,\Gamma,pr,an}_{/T^{an}},F_{DR}))
\end{eqnarray*}
\end{defi}

\begin{itemize}
\item Let $S\in\AnSm(\mathbb C)$. We have the map in $C_{O_Sfil,D_S}(\Var(\mathbb C)^{2,smpr}/S)$:
\begin{eqnarray*}
w_S:(\Omega^{\bullet,\Gamma}_{/S},F_b)\otimes_{O_S}(\Omega^{\bullet,\Gamma}_{/S},F_b)
\to(\Omega^{\bullet,\Gamma}_{/S},F_b): \\
\end{eqnarray*}
given by for $h:(U,Z)\to S\in\Var(\mathbb C)^{2,sm}/S$,
\begin{eqnarray*}
w_S((U,Z)/S):(\Gamma_Z^{\vee,h}L_{h^*O_S}(\Omega^{\bullet}_{U/S},F_b)\otimes_{p^*O_S}
\Gamma_Z^{\vee,h}L_{h^*O_S}(\Omega^{\bullet}_{U/S},F_b))(U) \\
\xrightarrow{(DR(-)(\gamma^{\vee,h}_Z(-))\circ w_{U/S})^{\gamma}(U)}
\Gamma_Z^{\vee,h}L_{h^*O_S}(\Omega^{\bullet}_{U/S},F_b)(U)
\end{eqnarray*}
which induces the map in $C_{O_Sfil,D_S}(\Var(\mathbb C)^{2,sm}/S)$
\begin{eqnarray*}
Ew_S:E_{et}(\Omega^{\bullet,\Gamma}_{/S},F_b)\otimes_{O_S} E_{et}(\Omega^{\bullet,\Gamma}_{/S},F_b)
\xrightarrow{=}
E_{et}((\Omega^{\bullet,\Gamma}_{/S},F_b)\otimes_{O_S}(\Omega^{\bullet,\Gamma}_{/S},F_b))
\xrightarrow{E_{et}(w_S)} E_{et}(\Omega^{\bullet,\Gamma}_{/S},F_b).
\end{eqnarray*}
\item Let $S\in\SmVar(\mathbb C)$. We have the map in $C_{D^{\infty}_Sfil}(\Var(\mathbb C)^{2,smpr}/S)$:
\begin{eqnarray*}
w_S:(\Omega^{\bullet,\Gamma,pr}_{/S^{an}},F_{DR})\otimes_{O_S}(\Omega^{\bullet,\Gamma,pr}_{/S^{an}},F_{DR})
\to(\Omega^{\bullet,\Gamma,pr}_{/S^{an}},F_{DR})
\end{eqnarray*}
given by for $p:(Y\times S,Z)\to S\in\Var(\mathbb C)^{2,smpr}/S$,
\begin{eqnarray*}
w_S((Y\times S,Z)/S): \\
(((\Omega^{\bullet}_{Y\times S/S}\otimes_{O_{Y\times S}}\Gamma_Z^{\vee,Hdg})(O_{Y\times S},F_b))\otimes_{p^*O_S}
(\Omega^{\bullet}_{Y\times S/S}\otimes_{O_{Y\times S}}\Gamma_Z^{\vee,Hdg}(O_{Y\times S},F_b)))(Y\times S) \\
\xrightarrow{(DR(-)(\gamma^{\vee,Hdg}_Z(-))\circ w_{Y\times S/S})^{\gamma}(Y\times S)}
(\Omega^{\bullet}_{Y\times S/S}\otimes_{O_{Y\times S}}\Gamma_Z^{\vee,Hdg}(O_{Y\times S},F_b))(Y\times S)
\end{eqnarray*}
which induces the map in $C_{D^{\infty}_Sfil}(\Var(\mathbb C)^{2,smpr}/S)$
\begin{eqnarray*}
Ew_S:E_{et}(\Omega^{\bullet,\Gamma,pr}_{/S^{an}},F_{DR})\otimes_{O_S} E_{et}(\Omega^{\bullet,\Gamma,pr}_{/S^{an}},F_{DR})
\xrightarrow{=} \\
E_{et}((\Omega^{\bullet,\Gamma,pr}_{/S^{an}},F_{DR})\otimes_{O_S}(\Omega^{\bullet,\Gamma,pr}_{/S^{an}},F_{DR}))
\xrightarrow{E_{et}(w_S)} E_{et}(\Omega^{\bullet,\Gamma,pr}_{/S^{an}},F_{DR})
\end{eqnarray*}
by the functoriality of the Godement resolution (see section 2).
\end{itemize}

\begin{defi}\label{TotimesDRan}
Let $S\in\SmVar(\mathbb C)$.
We have, for $F,G\in C(\Var(\mathbb C)^{2,smpr}/S)$, the canonical transformation in $C_{\mathcal D^{\infty}fil}(S^{an})$ :
\begin{eqnarray*}
T(\otimes,\Omega)(F,G): \\
e(S)_*\Gr^{12}_{S*}\mathcal Hom(\An_S^*F,E_{et}(\Omega^{\bullet,\Gamma,pr,an}_{/S},F_{DR}))\otimes_{O_S}
e(S)_*\Gr^{12}_{S*}\mathcal Hom(\An_S^*G,E_{et}(\Omega^{\bullet,\Gamma,pr,an}_{/S},F_{DR})) \\
\xrightarrow{=} 
e(S)_*\Gr^{12}_{S*}(\mathcal Hom(\An_S^*F,E_{et}(\Omega^{\bullet,\Gamma,pr,an}_{/S},F_{DR}))\otimes_{O_S}
\mathcal Hom(\An_S^*G,E_{et}(\Omega^{\bullet,\Gamma,pr,an}_{/S},F_{DR}))) \\ 
\xrightarrow{T(\mathcal Hom,\otimes)(-)} 
e(S)_*\Gr^{12}_{S*}\mathcal Hom(\An_S^*F\otimes \An_S^*G,E_{et}(\Omega^{\bullet,\Gamma,pr}_{/S},F_{DR})\otimes_{O_S}
E_{usu}(\Omega^{\bullet,\Gamma,pr}_{/S},F_{DR}))) \\
\xrightarrow{=} 
e(S)_*\Gr^{12}_{S*}\mathcal Hom(\An_S^*(F\otimes G),E_{et}(\Omega^{\bullet,\Gamma,pr}_{/S},F_{DR})
\otimes_{O_S}E_{et}(\Omega^{\bullet,\Gamma,pr}_{/S},F_{DR})) \\
\xrightarrow{\mathcal Hom(F\otimes G,\An_S^{*mod}Ew_S)} 
e(S)_*\Gr^{12}_{S*}\mathcal Hom(F\otimes G,E_{et}(\Omega^{\bullet,\Gamma,pr,an}_{/S},F_{DR})).
\end{eqnarray*}
\end{defi}

Let $S\in\Var(\mathbb C)$. Let $S=\cup_{i=1}^l S_i$ an open affine cover and denote by $S_I=\cap_{i\in I} S_i$.
Let $i_i:S_i\hookrightarrow\tilde S_i$ closed embeddings, with $\tilde S_i\in\Var(\mathbb C)$. 
For $I\subset\left[1,\cdots l\right]$, denote by $\tilde S_I=\Pi_{i\in I}\tilde S_i$.
We then have closed embeddings $i_I:S_I\hookrightarrow\tilde S_I$ and for $J\subset I$ the following commutative diagram
\begin{equation*}
D_{IJ}=\xymatrix{ S_I\ar[r]^{i_I} & \tilde S_I \\
S_J\ar[u]^{j_{IJ}}\ar[r]^{i_J} & \tilde S_J\ar[u]^{p_{IJ}}}  
\end{equation*}
where $p_{IJ}:\tilde S_J\to\tilde S_I$ is the projection
and $j_{IJ}:S_J\hookrightarrow S_I$ is the open embedding so that $j_I\circ j_{IJ}=j_J$.
This gives the diagram of algebraic varieties $(\tilde S_I)\in\Fun(\mathcal P(\mathbb N),\Var(\mathbb C))$ which
the diagram of sites $\Var(\mathbb C)^{2,smpr}/(\tilde S_I)\in\Fun(\mathcal P(\mathbb N),\Cat)$. 
This gives also the diagram of algebraic varieties $(\tilde S_I)^{op}\in\Fun(\mathcal P(\mathbb N)^{op},\Var(\mathbb C))$ which
the diagram of sites $\Var(\mathbb C)^{2,smpr}/(\tilde S_I)^{op}\in\Fun(\mathcal P(\mathbb N)^{op},\Cat)$. 
We then get
\begin{eqnarray*}
((\Omega^{\bullet,\Gamma,pr,an}_{/(\tilde S_I)},F_{DR})[-d_{\tilde S_I}],T_{IJ})
\in C_{D^{\infty}_{(\tilde S_I)}fil}(\AnSp(\mathbb C)^{2,smpr}/(\tilde S_I))
\end{eqnarray*}
with
\begin{eqnarray*}
T_{IJ}:(\Omega^{\bullet,\Gamma,pr,an}_{/\tilde S_I},F_{DR})[-d_{\tilde S_I}]
\xrightarrow{\ad(p_{IJ}^{*mod[-]},p_{IJ*}(-)}
p_{IJ*}p_{IJ}^*(\Omega^{\bullet,\Gamma,pr,an}_{/\tilde S_I},F_{DR})
\otimes_{p_{IJ}^*O_{\tilde S_I}}O_{\tilde S_J}[-d_{\tilde S_J}] \\
\xrightarrow{m\circ p_{IJ*}\Omega^{\Gamma,pr,an}_{/(\tilde S_J/\tilde S_I)}[-d_{\tilde S_J}]}
p_{IJ*}(\Omega^{\bullet,\Gamma,pr,an}_{/\tilde S_J},F_{DR})[-d_{\tilde S_J}].
\end{eqnarray*}
For $(G_I,K_{IJ})\in C(\AnSp(\mathbb C)^{2,smpr}/(\tilde S_I)^{op})$, we denote (see section 2)
\begin{eqnarray*}
e'((\tilde S_I))_*\mathcal Hom((G_I,K_{IJ}),
(E_{usu}(\Omega^{\bullet,\Gamma,pr,an}_{/(\tilde S_I)},F_{DR})[-d_{\tilde S_I}],T_{IJ})):= \\
(e'(\tilde S_I)_*\mathcal Hom(G_I,E_{et}(\Omega^{\bullet,\Gamma,pr,an}_{/\tilde S_I},F_{DR}))[-d_{\tilde S_I}],
u_{IJ}((G_I,K_{IJ})))\in C_{\mathcal D^{\infty}fil}((\tilde S_I))
\end{eqnarray*}
with
\begin{eqnarray*}
u_{IJ}((G_I,K_{IJ})):e'(\tilde S_I)_*\mathcal Hom(G_I,E_{usu}(\Omega^{\bullet,\Gamma,pr,an}_{/\tilde S_I},F_{DR}))[-d_{\tilde S_I}] \\
\xrightarrow{\ad(p_{IJ}^{*mod[-]},p_{IJ*})(-)\circ T(p_{IJ},e)(-)}
p_{IJ*}e'(\tilde S_J)_*p_{IJ}^*\mathcal Hom(G_I,E_{usu}(\Omega^{\bullet,\Gamma,pr,an}_{/\tilde S_I},F_{DR}))
\otimes_{p_{IJ}^*O_{\tilde S_I}}O_{\tilde S_J}[-d_{\tilde S_J}] \\
\xrightarrow{T(p_{IJ},hom)(-,-)}
p_{IJ*}e'(\tilde S_J)_*\mathcal Hom(p_{IJ}^*G_I,p_{IJ}^*E_{usu}(\Omega^{\bullet,\Gamma,pr,an}_{/\tilde S_I},F_{DR})) 
\otimes_{p_{IJ}^*O_{\tilde S_I}}O_{\tilde S_J}[-d_{\tilde S_J}] \\
\xrightarrow{m\circ\mathcal Hom(p_{IJ}^*G_I,T_{IJ})}
p_{IJ*}e'(\tilde S_J)_*\mathcal Hom(p_{IJ}^*G_I,E_{usu}(\Omega^{\bullet,\Gamma,pr,an}_{/\tilde S_J},F_{DR}))[-d_{\tilde S_J}] \\
\xrightarrow{\mathcal Hom(K_{IJ},E_{usu}(\Omega^{\bullet,\Gamma,pr,an}_{/\tilde S_J},F_{DR}))}
p_{IJ*}e'(\tilde S_J)_*\mathcal Hom(G_J,E_{usu}(\Omega^{\bullet,\Gamma,pr}_{/\tilde S_J},F_{DR}))[-d_{\tilde S_J}].
\end{eqnarray*}
This gives in particular
\begin{eqnarray*}
(\Omega^{\bullet,\Gamma,pr,an}_{/(\tilde S_I)},F_{DR})[-d_{\tilde S_I}],T_{IJ})\in 
C_{D^{\infty}_{(\tilde S_I)}fil}(\AnSp(\mathbb C)^{2,(sm)pr}/(\tilde S_I)^{op}).
\end{eqnarray*}

We now define the filtered analytic De Rahm realization functor.

\begin{defi}\label{DRalgdefFunctan}
\begin{itemize}
\item[(i)]Let $S\in\SmVar(\mathbb C)$.
We have, using definition \ref{wtildewan} and definition \ref{RCHhatdef}, the functor
\begin{eqnarray*}
\mathcal F_{S,an}^{FDR}:C(\Var(\mathbb C)^{sm}/S)\to C_{\mathcal D^{\infty}fil}(S^{an}), \; F\mapsto \\ 
\mathcal F_{S,an}^{FDR}(F):=e'(S)_*\mathcal Hom^{\bullet}(\An_S^*\hat R^{CH}(\rho_S^*L(F)),
E_{et}(\Omega^{\bullet,\Gamma,pr,an}_{/S},F_{DR}))[-d_S] \\
=e'(S)_*\mathcal Hom^{\bullet}(\hat R^{CH}(\rho_S^*L(F)),
\An_{S*}E_{et}(\Omega^{\bullet,\Gamma,pr,an}_{/S},F_{DR}))[-d_S] 
\end{eqnarray*}
denoting for short $e'(S)=e(S)\circ\Gr^{12}_S$.
\item[(ii)]Let $S\in\Var(\mathbb C)$ and $S=\cup_{i=1}^l S_i$ an open cover such that there exist closed embeddings
$i_i:S_i\hookrightarrow\tilde S_i$  with $\tilde S_i\in\SmVar(\mathbb C)$. 
For $I\subset\left[1,\cdots l\right]$, denote by $S_I:=\cap_{i\in I} S_i$ and $j_I:S_I\hookrightarrow S$ the open embedding.
We then have closed embeddings $i_I:S_I\hookrightarrow\tilde S_I:=\Pi_{i\in I}\tilde S_i$.
Consider, for $I\subset J$, the following commutative diagram
\begin{equation*}
D_{IJ}=\xymatrix{ S_I\ar[r]^{i_I} & \tilde S_I \\
S_J\ar[u]^{j_{IJ}}\ar[r]^{i_J} & \tilde S_J\ar[u]^{p_{IJ}}}  
\end{equation*}
and $j_{IJ}:S_J\hookrightarrow S_I$ is the open embedding so that $j_I\circ j_{IJ}=j_J$.
We have, using definition \ref{wtildewan} and definition \ref{RCHhatdef}, the functor
\begin{eqnarray*}
\mathcal F_{S,an}^{FDR}:C(\Var(\mathbb C)^{sm}/S)\to C_{\mathcal D^{\infty}fil}(S^{an}/(\tilde S^{an}_I)), \; 
F\mapsto \\
\mathcal F_{S,an}^{FDR}(F):=e'(\tilde S_I)_*\mathcal Hom^{\bullet}((\An_{\tilde S_I}^*
\hat R^{CH}(\rho_{\tilde S_I}^*L(i_{I*}j_I^*F)),\hat R^{CH}(T^q(D_{IJ})(j_I^*F))), \\ 
(E_{usu}(\Omega^{\bullet,\Gamma,pr,an}_{/\tilde S_I},F_{DR})[-d_{\tilde S_I}],T_{IJ})) \\
(e'(\tilde S_I)_*\mathcal Hom^{\bullet}(\An_{\tilde S_I}^*\hat R^{CH}(\rho_{\tilde S_I}^*L(i_{I*}j_I^*F)), 
E_{et}(\Omega^{\bullet,\Gamma,pr,an}_{/\tilde S_I},F_{DR}))[-d_{\tilde S_I}],u^q_{IJ}(F)) \\
=(e'(\tilde S_I)_*\mathcal Hom^{\bullet}(\hat R^{CH}(\rho_{\tilde S_I}^*L(i_{I*}j_I^*F)), 
\An_{\tilde S_I*}E_{usu}(\Omega^{\bullet,\Gamma,pr,an}_{/\tilde S_I},F_{DR}))[-d_{\tilde S_I}],u^q_{IJ}(F)) \\
\end{eqnarray*}
where we have denoted for short $e'(\tilde S_I)=e(\tilde S_I)\circ\Gr^{12}_{\tilde S_I}$, and
\begin{eqnarray*}
u^q_{IJ}(F)[d_{\tilde S_J}]:
e'(\tilde S_I)_*\mathcal Hom^{\bullet}(\An_{\tilde S_I}^*
\hat R^{CH}(\rho_{\tilde S_I}^*L(i_{I*}j_I^*F)),E_{et}(\Omega^{\bullet,\Gamma,pr,an}_{/\tilde S_I},F_{DR})) \\
\xrightarrow{\ad(p_{IJ}^{*mod},p_{IJ})(-)} 
p_{IJ*}p_{IJ}^{*mod}e'(\tilde S_I)_*\mathcal Hom^{\bullet}(\An_{\tilde S_I}^*\hat R^{CH}(\rho_{\tilde S_I}^*L(i_{I*}j_I^*F)),
E_{usu}(\Omega^{\bullet,\Gamma,pr,an}_{/\tilde S_I},F_{DR})) \\
\xrightarrow{p_{IJ*}T(p_{IJ},\Omega^{\gamma,pr}_{\cdot})(-)}  
p_{IJ*}e'(\tilde S_J)_*\mathcal Hom^{\bullet}(\An_{\tilde S_I}^*p_{IJ}^*\hat R^{CH}(\rho_{\tilde S_I}^*L(i_{I*}j_I^*F)),
E_{usu}(\Omega^{\bullet,\Gamma,pr,an}_{/\tilde S_J},F_{DR})) \\
\xrightarrow{\mathcal Hom(T(p_{IJ},\hat R^{CH})(Li_{I*}j_I^*F)^{-1},E_{et}(\Omega_{/\tilde S_J}^{\bullet,\Gamma,pr,an},F_{DR}))} \\
p_{IJ*}e'(\tilde S_J)_*\mathcal Hom^{\bullet}(\An_{\tilde S_I}^*\hat R^{CH}(\rho_{\tilde S_J}^*p_{IJ}^*L(i_{I*}j_I^*F)),
E_{usu}(\Omega^{\bullet,\Gamma,pr}_{/\tilde S_J},F_{DR})) \\
\xrightarrow{\mathcal Hom(\hat R^{CH}_{\tilde S_J}(T^q(D_{IJ})(j_I^*F)),
E_{usu}(\Omega_{/\tilde S_J}^{\bullet,\Gamma,pr,an},F_{DR}))} \\
p_{IJ*}e'(\tilde S_J)_*\mathcal Hom^{\bullet}(\An_{\tilde S_I}^*\hat R^{CH}(\rho_{\tilde S_J}^*L(i_{J*}j_J^*F)),
E_{usu}(\Omega^{\bullet,\Gamma,pr,an}_{/\tilde S_J},F_{DR})).
\end{eqnarray*}
For $I\subset J\subset K$, we have obviously $p_{IJ*}u_{JK}(F)\circ u_{IJ}(F)=u_{IK}(F)$.
\end{itemize}
\end{defi}

We have the following key proposition :

Recall, see section 2, that we have the projection morphisms of sites 
$p_a:\AnSp(\mathbb C)^{2,smpr}/(\tilde S_I)^{op}\to\AnSp(\mathbb C)^{2,smpr}/(\tilde S_I)^{op}$
given by the functor 
\begin{eqnarray*}
p_a:\AnSp(\mathbb C)^{2,smpr}/(\tilde S_I)^{op}\to\AnSp(\mathbb C)^{2,smpr}/(\tilde S_I)^{op}, \\
p_a((Y_I\times\tilde S_I,Z_I)/\tilde S_I,s_{IJ}):=
((Y_I\times\mathbb D^1\times\tilde S_I,Z_I\times\mathbb D^1)/\tilde S_I,s_{IJ}\times I), \\ 
p_a((g_I):((Y'_I\times\tilde S_I,Z'_I)/\tilde S_I,s'_{IJ})\to((Y_I\times\tilde S_I,Z_I)/\tilde S_I,s_{IJ}))= \\
(g_I\times I):((Y'_I\times\mathbb D^1\times\tilde S_I,Z'_I\times\mathbb D^1)/\tilde S_I,s'_{IJ}\times I)
\to((Y_I\times\mathbb D^1\times\tilde S_I,Z_I\times\mathbb D^1)/\tilde S_I),s_{IJ}\times I)).
\end{eqnarray*}

We have the following key proposition :

\begin{prop}\label{ausufib}
\begin{itemize}
\item[(i1)]Let $S\in\Var(\mathbb C)$. Let $S=\cup_{i=1}^l S_i$ an open cover such that there exist closed embeddings
$i_i:S_i\hookrightarrow\tilde S_i$ with $\tilde S_i\in\SmVar(\mathbb C)$.
The complex of presheaves
$(\Omega^{\bullet,\Gamma,pr,an}_{/(\tilde S_I)},F_{DR})\in 
C_{D^{\infty}_{(\tilde S_I)}fil}(\AnSp(\mathbb C)^{2,smpr}/(\tilde S_I)^{op})$ 
is $2$-filtered $\mathbb D^1$ homotopic, that is
\begin{equation*}
\ad(p_a^*,p_{a*})(\Omega^{\bullet,\Gamma,pr,an}_{/(\tilde S_I})^{an},F_{DR}):
(\Omega^{\bullet,\Gamma,pr,an}_{/S},F_{DR})\to p_{a*}p_a^*(\Omega^{\bullet,\Gamma,pr}_{/(\tilde S_I})^{an},F_{DR})
\end{equation*}
is a $2$-filtered homotopy.
\item[(i2)]Let $S\in\SmVar(\mathbb C)$. The complex of presheaves
$(\Omega^{\bullet,\Gamma,pr,an}_{/S},F_{DR})\in C_{D^{\infty}_Sfil}(\AnSp(\mathbb C)^{2,smpr}/S)$ 
admits transferts, i.e. 
\begin{equation*}
\Tr(S)_*\Tr(S)^*(\Omega^{\bullet,\Gamma,pr,an}_{/S},F_{DR}=(\Omega^{\bullet,\Gamma,pr,an}_{/S},F_{DR}).
\end{equation*}
\item[(ii)]Let $S\in\Var(\mathbb C)$. Let $S=\cup_{i=1}^l S_i$ an open cover such that there exist closed embeddings
$i_i:S_i\hookrightarrow\tilde S_i$ with $\tilde S_i\in\SmVar(\mathbb C)$.
Let $m=(m_I):(Q_{1I},K^1_{IJ})\to(Q_{2I},K^2_{IJ})$ be an equivalence $(\mathbb D^1,usu)$ local with 
$(Q_{1I},K_{IJ})\to(Q_{2I},K_{IJ})\in C(\AnSp(\mathbb C)^{smpr}/(\tilde S_I)^{op})$ complexes of representable presheaves. 
Then, the map in $C_{\mathcal D^{\infty}fil}((\tilde S_I))$
\begin{eqnarray*} 
M:=(e(\tilde S_I)_*\mathcal Hom^{\bullet}(m_I,E_{usu}(\Omega^{\bullet,\Gamma,pr,an}_{/\tilde S_I},F_{DR})[-d_{\tilde S_I}])): \\  
e'((\tilde S_I))_*\mathcal Hom^{\bullet}((Q_{2I},K^1_{IJ}),
(E_{usu}(\Omega^{\bullet,\Gamma,pr,an}_{/\tilde S_I},F_{DR})[-d_{\tilde S_I}],T_{IJ})) \\
\to e'((\tilde S_I))_*\mathcal Hom^{\bullet}((Q_{1I},K^1_{IJ}),
(E_{usu}(\Omega^{\bullet,\Gamma,pr,an}_{/\tilde S_I},F_{DR})[-d_{\tilde S_I}],T_{IJ}))
\end{eqnarray*}
is a $2$-filtered quasi-isomorphism. It is thus an isomorphism in $D_{\mathcal D^{\infty}fil,\infty}((\tilde S_I))$.
\end{itemize}
\end{prop}

\begin{proof} 
Similar to proposition \ref{aetfib}.
\end{proof}

We deduce the following:

\begin{prop}\label{projwachan}
Let $S\in\Var(\mathbb C)$.Let $S=\cup_{i=1}^l S_i$ an open cover such that there exist closed embeddings
$i_i:S_i\hookrightarrow\tilde S_i$ with $\tilde S_i\in\SmVar(\mathbb C)$. 
\begin{itemize}
\item[(i)]Let $m=(m_I):(Q_{1I},K^1_{IJ})\to (Q_{2I},K^2_{IJ})$ be an etale local equivalence local 
with $(Q_{1I},K^1_{IJ}),(Q_{2I},K^2_{IJ})\in C(\Var(\mathbb C)^{sm}/(\tilde S_I))$
complexes of projective presheaves. Then,
\begin{eqnarray*} 
(e'(\tilde S_I)_*\mathcal Hom^{\bullet}(\An_{\tilde S_I}^*\hat R_S^{CH}(m_I),
E_{usu}(\Omega^{\bullet,\Gamma,pr,an}_{/\tilde S_I},F_{DR}))[-d_{\tilde S_I}]): \\ 
e'(\tilde S_I)_*\mathcal Hom^{\bullet}((\An_{\tilde S_I}^*\hat R^{CH}(\rho_S^*Q_{1I}),R^{CH}(K^1_{IJ})),
(E_{usu}(\Omega^{\bullet,\Gamma,pr,an}_{/\tilde S_I},F_{DR})[-d_{\tilde S_I}],T_{IJ})) \\
\to e'(\tilde S_I)_*\mathcal Hom^{\bullet}((\An_{\tilde S_I}^*L\hat R^{CH}(\rho_S^*Q_{2I}),R^{CH}(K^2_{IJ})),
(E_{usu}(\Omega^{\bullet,\Gamma,pr,an}_{/\tilde S_I},F_{DR})[-d_{\tilde S_I}],T_{IJ}))
\end{eqnarray*}
is a filtered quasi-isomorphism. It is thus an isomorphism in $D_{\mathcal D^{\infty}fil}((\tilde S_I))$.
\item[(ii)]Let $m=(m_I):(Q_{1I},K^1_{IJ})\to (Q_{2I},K^2_{IJ})$ be an equivalence $(\mathbb A^1,et)$ local equivalence local 
with $(Q_{1I},K^1_{IJ}),(Q_{2I},K^2_{IJ})\in C(\Var(\mathbb C)^{sm}/(\tilde S_I))$
complexes of projective presheaves. Then,
\begin{eqnarray*} 
(e'(\tilde S_I)_*\mathcal Hom^{\bullet}(\An_{\tilde S_I}^*\hat R_S^{CH}(m_I),
E_{usu}(\Omega^{\bullet,\Gamma,pr}_{/\tilde S_I},F_{DR}))[-d_{\tilde S_I}]): \\ 
e'(\tilde S_I)_*\mathcal Hom^{\bullet}((\An_{\tilde S_I}^*\hat R^{CH}(\rho_S^*Q_{1I}),R^{CH}(K^1_{IJ})),
(E_{usu}(\Omega^{\bullet,\Gamma,pr}_{/\tilde S_I},F_{DR})[-d_{\tilde S_I}],T_{IJ})) \\
\to e'(\tilde S_I)_*\mathcal Hom^{\bullet}((\An_{\tilde S_I}^*\hat R^{CH}(\rho_S^*Q_{2I}),R^{CH}(K^2_{IJ})),
(E_{usu}(\Omega^{\bullet,\Gamma,pr}_{/\tilde S_I},F_{DR})[-d_{\tilde S_I}],T_{IJ}))
\end{eqnarray*}
is an filtered quasi-isomorphism. It is thus an isomorphism in $D_{\mathcal D^{\infty}fil}((\tilde S_I))$.
\end{itemize}
\end{prop}

\begin{proof}
Similar to the proof of proposition \ref{projwach}.
\end{proof}

\begin{defi}\label{DRalgdefsingan}
\begin{itemize}
\item[(i)] Let $S\in\SmVar(\mathbb C)$.
We define using definition \ref{DRalgdefFunctan}(i) and proposition \ref{projwachan}(ii)
the filtered algebraic De Rahm realization functor defined as
\begin{eqnarray*}
\mathcal F_{S,an}^{FDR}:\DA_c(S)\to D_{\mathcal D^{\infty}fil}(S^{an}), 
M\mapsto\mathcal F_S^{FDR}(M):= \\ 
e'(S)_*\mathcal Hom^{\bullet}(\An_S^*\hat R^{CH}(\rho_S^*L(F)), E_{usu}(\Omega^{\bullet,\Gamma,pr,an}_{/S},F_{DR}))[-d_S] 
\end{eqnarray*}
where $F\in C(\Var(\mathbb C)^{sm}/S)$ is such that $M=D(\mathbb A^1,et)(F)$.
\item[(i)'] For the Corti-Hanamura weight structure $W$ on $\DA_c(S)^-$, 
we define using definition \ref{DRalgdefFunctan}(i) and proposition \ref{projwachan}(ii)
\begin{eqnarray*}
\mathcal F_{S,an}^{FDR}:\DA_c^-(S)\to D_{\mathcal D^{\infty}(1,0)fil}^-(S^{an}), 
M\mapsto\mathcal F_S^{FDR}((M,W)):= \\ 
e'(S)_*\mathcal Hom^{\bullet}(\An_S^*\hat R^{CH}(\rho_S^*L(F,W)),E_{usu}(\Omega^{\bullet,\Gamma,pr,an}_{/S},F_{DR}))[-d_S] 
\end{eqnarray*}
where $(F,W)\in C_{fil}(\Var(\mathbb C)^{sm}/S)$ 
is such that $M=D(\mathbb A^1,et)((F,W))$ using corollary \ref{weightst2Cor}. 
Note that the filtration induced by $W$ is a filtration by sub $D_S$ module,
which is a stronger property then Griffitz transversality.
Of course, the filtration induced by $F$ satisfy only Griffitz transversality in general.
\item[(ii)]Let $S\in\Var(\mathbb C)$ and $S=\cup_{i=1}^l S_i$ an open cover such that there exist closed embeddings
$i_i:S_i\hookrightarrow\tilde S_i$  with $\tilde S_i\in\SmVar(\mathbb C)$. 
For $I\subset\left[1,\cdots l\right]$, denote by $S_I=\cap_{i\in I} S_i$ and $j_I:S_I\hookrightarrow S$ the open embedding.
We then have closed embeddings $i_I:S_I\hookrightarrow\tilde S_I:=\Pi_{i\in I}\tilde S_i$.
We define, using definition \ref{DRalgdefFunctan}(ii) and proposition \ref{projwachan}(ii),
the filtered algebraic De Rahm realization functor defined as
\begin{eqnarray*}
\mathcal F_{S,an}^{FDR}:\DA_c(S)\to D_{\mathcal D^{\infty}fil}(S^{an}/(\tilde S^{an}_I)), 
M\mapsto\mathcal F_S^{FDR}(M):= \\
(e'(\tilde S_I)_*\mathcal Hom^{\bullet}(\An_{\tilde S_I}^*\hat R^{CH}(\rho_{\tilde S_I}^*L(i_{I*}j_I^*F)),
E_{usu}(\Omega^{\bullet,\Gamma,pr,an}_{/\tilde S_I},F_{DR}))[-d_{\tilde S_I}],u^q_{IJ}(F))
\end{eqnarray*}
where $F\in C(\Var(\mathbb C)^{sm}/S)$ is such that $M=D(\mathbb A^1,et)(F)$, 
see definition \ref{DRalgdefFunct} .
\item[(ii)'] For the Corti-Hanamura weight structure $W$ on $\DA_c^-(S)$, 
using definition \ref{DRalgdefFunct}(ii) and proposition \ref{projwach}(ii), 
\begin{eqnarray*}
\mathcal F_{S,an}^{FDR}:\DA_c^-(S)\to D_{\mathcal D^{\infty}(1,0)fil}^-(S^{an}/(\tilde S^{an}_I)), 
M\mapsto\mathcal F_S^{FDR}((M,W)):= \\
(e'(\tilde S_I)_*\mathcal Hom^{\bullet}(\An_{\tilde S_I}^*\hat R^{CH}(\rho_{\tilde S_I}^*L(i_{I*}j_I^*(F,W))),
E_{usu}(\Omega^{\bullet,\Gamma,pr,an}_{/\tilde S_I},F_{DR}))[-d_{\tilde S_I}],u^q_{IJ}(F,W))
\end{eqnarray*}
where $(F,W)\in C_{fil}(\Var(\mathbb C)^{sm}/S)$ 
is such that $(M,W)=D(\mathbb A^1,et)(F,W)$ using corollary \ref{weightst2Cor}.
Note that the filtration induced by $W$ is a filtration by sub $D_{\tilde S_I}$-modules,
which is a stronger property then Griffitz transversality.
Of course, the filtration induced by $F$ satisfy only Griffitz transversality in general.
\end{itemize}
\end{defi}

\begin{prop}\label{FDRwelldefan}
For $S\in\Var(\mathbb C)$ and $S=\cup_{i=1}^l S_i$ an open cover such that there exist closed embeddings
$i_i:S_i\hookrightarrow\tilde S_i$ with $\tilde S_i\in\SmVar(\mathbb C)$, the functor $\mathcal F_{S,an}^{FDR}$ is well defined. 
\end{prop}

\begin{proof}
Similar to the proof of proposition \ref{FDRwelldef}.
\end{proof}

\begin{prop}\label{keysing1}
Let $f:X\to S$ a morphism with $S,X\in\Var(\mathbb C)$. Assume there exist a factorization 
\begin{equation*}
f:X\xrightarrow{l}Y\times S\xrightarrow{p_S} S
\end{equation*}
of $f$ with $Y\in\SmVar(\mathbb C)$, $l$ a closed embedding and $p_S$ the projection.
Let $\bar Y\in\PSmVar(\mathbb C)$ a compactification of $Y$ with $\bar Y\backslash Y=D$ a normal crossing divisor,
denote $k:D\hookrightarrow \bar Y$ the closed embedding and $n:Y\hookrightarrow\bar Y$ the open embedding.
Denote $\bar X\subset\bar Y\times S$ the closure of $X\subset\bar Y\times S$.
We have then the following commutative diagram in $\Var(\mathbb C)$
\begin{equation*}
\xymatrix{X\ar[r]^l\ar[d] & Y\times S\ar[rd]^{p_S}\ar[d]^{(n\times I)} & \, \\
\bar X\ar[r]^l & \bar Y\times S\ar[r]^{\bar p_S} & S \\
Z:=\bar X\backslash X\ar[ru]^{l_Z}\ar[u]\ar[r] & D\times S\ar[ru]\ar[u]^{(k\times I)} & \, }.
\end{equation*}
Let $S=\cup_{i=1}^l S_i$ an open cover such that there exist closed embeddings
$i_i:S_i\hookrightarrow\tilde S_i$ with $\tilde S_i\in\SmVar(\mathbb C)$. 
Then $X=\cup_{i=1}^lX_i$ with $X_i:=f^{-1}(S_i)$.
Denote, for $I\subset\left[1,\cdots l\right]$, $S_I=\cap_{i\in I} S_i$ and $X_I=\cap_{i\in I}X_i$.
Denote $\bar X_I:=\bar X\cap(\bar Y\times S_I)\subset\bar Y\times\tilde S_I$ 
the closure of $X_I\subset\bar Y\times\tilde S_I$, 
and $Z_I:=Z\cap(\bar Y\times S_I)=\bar X_I\backslash X_I\subset\bar Y\times\tilde S_I$.
We have then for $I\subset\left[1,\cdots l\right]$, the following commutative diagram in $\Var(\mathbb C)$
\begin{equation*}
\xymatrix{X_I\ar[r]^{l_I}\ar[d] & Y\times \tilde S_I\ar[rd]^{p_{\tilde S_I}}\ar[d]^{(n\times I)} & \, \\
\bar X_I\ar[r]^{l_I} & \bar Y\times\tilde S_I\ar[r]^{\bar p_{\tilde S_I}} & \tilde S_I \\
Z_I=\bar X_I\backslash X_I\ar[ru]^{l_{Z_I}}\ar[u]\ar[r] & D\times\tilde S_I\ar[ru]\ar[u]_{(k\times I)} & \, }.
\end{equation*}
Let $F(X/S):=p_{S,\sharp}\Gamma_X^{\vee}\mathbb Z(X\times S/X\times S)$. 
We have then the following isomorphism in $D_{\mathcal Dfil}(S/(\tilde S_I))$ 
\begin{eqnarray*} 
I(X/S):\mathcal F_{S,an}^{FDR}(M(X/S))\xrightarrow{:=} \\
(e'_*\mathcal Hom(\An_{\tilde S_I}^*\hat R^{CH}(\rho_{\tilde S_I}^*L(i_{I*}j_I^*F(X/S))), 
E_{usu}(\Omega^{\bullet,\Gamma,pr,an}_{/\tilde S_I},F_{DR}))[-d_{\tilde S_I}],u^q_{IJ}(F(X/S))) \\ 
\xrightarrow{(\mathcal Hom(\An_{\tilde S_I}^*\hat R^{CH}_{\tilde S_I}(N_I(X/S)),
E_{usu}(\Omega^{\bullet,\Gamma,pr,an}_{/\tilde S_I},F_{DR})))}  \\
(e'_*\mathcal Hom(\An_{\tilde S_I}^*\hat R^{CH}(\rho_{\tilde S_I}^*Q(X_I/\tilde S_I)), 
E_{usu}(\Omega^{\bullet,\Gamma,pr,an}_{/\tilde S_I},F_{DR}))[-d_{\tilde S_I}],v^q_{IJ}(F(X/S))) \\ 
\xrightarrow{I((\bar X_I,Z_I)/\tilde S_I)} \\
(\bar p_{\tilde S_I*}E_{usu}((\Omega^{\bullet}_{\bar Y\times\tilde S_I/\tilde S_I},F_b)\otimes_{O_{Y\times\tilde S_I}}
(n\times I)_{!Hdg}\Gamma^{\vee,Hdg}_{X_I}(O_{(Y\times\tilde S_I)^{an}},F_b))(d_Y+d_{\tilde S_I})[2d_Y+d_{\tilde S_I}],
w_{IJ}(X/S)) \\ 
\xrightarrow{=:}\iota_SRf^{Hdg}_!(\Gamma^{\vee,Hdg}_{X_I}(O_{(Y\times\tilde S_I)^{an}},F_b)(d_Y)[2d_Y],x_{IJ}(X/S))
\xrightarrow{=:}\iota_SRf^{Hdg}_!f^{*mod}_{Hdg}\mathbb Z_{S^{an}}^{Hdg}
\end{eqnarray*}
\end{prop}

\begin{proof}
Similar to the proof of proposition \ref{keyalgsing1}.
\end{proof}

\begin{cor}\label{FDRMHMan}
Let $S\in\Var(\mathbb C)$ and $S=\cup_{i=1}^l S_i$ an open cover such that there exist closed embeddings
$i_i:S_i\hookrightarrow\tilde S_i$ with $\tilde S_i\in\SmVar(\mathbb C)$. 
Then, for $F\in C(\Var(\mathbb C)^{sm}/S)$ such that $M=D(\mathbb A^1,et)(F)\in\DA_c(S)$, 
\begin{eqnarray*}
H^i\mathcal F_{S,an}^{FDR}(M,W):= 
(a_{usu}H^ie'(\tilde S_I)_*\mathcal Hom^{\bullet}(\An_{\tilde S_I}^*\hat R^{CH}(\rho_{\tilde S_I}^*L(i_{I*}j_I^*(F,W))), \\
E_{usu}(\Omega^{\bullet,\Gamma,pr,an}_{/\tilde S_I},F_{DR}))[-d_{\tilde S_I}],H^iu^q_{IJ}(F,W)) 
\in\pi_S(MHM(S^{an}))
\end{eqnarray*}
for all $i\in\mathbb Z$, and for all $p\in\mathbb Z$, 
\begin{eqnarray*}
\mathcal F_{S,an}^{FDR}(M,W)\in D_{\mathcal D^{\infty}(1,0)fil}(S/(\tilde S_I))
\end{eqnarray*}
is the class of a complex $\mathcal F_{S,an}^{FDR}(F,W)\in C_{\mathcal D^{\infty}(1,0)fil}(S/(\tilde S_I))$
such that for all $k\in\mathbb Z$,
the differentials of $\Gr_k^W\mathcal F_{S,an}^{FDR}(F,W)$ are strict for the filtration $F$.
\end{cor}

\begin{proof}
Similar to the proof of corollary \ref{FDRMHM}.
\end{proof}

\begin{prop}\label{FDRHdgwelldefan}
For $S\in\Var(\mathbb C)$ not smooth, the functor (see corollary \ref{FDRMHM}) 
\begin{equation*}
\iota_S^{-1}\mathcal F_{S,an}^{FDR}:\DA_c^-(S)^{op}\to\pi_S(D(MHM(S^{an}))
\end{equation*}
does not depend on the choice of the open cover $S=\cup_iS_i$
and the closed embeddings $i_i:S_i\hookrightarrow\tilde S_i$ with $\tilde S_i\in\SmVar(\mathbb C)$.
\end{prop}

\begin{proof}
Similar to the proof of proposition \ref{FDRHdgwelldef}.
\end{proof}

We have the canonical transformation map between the filtered analytic De Rham realization functor and 
the analytic Gauss-Manin realization functor :

\begin{defi}
Let $S\in\Var(\mathbb C)$ and $S=\cup_{i=1}^l S_i$ an open cover such that there exist closed embeddings
$i_i:S_i\hookrightarrow\tilde S_i$ with $\tilde S_i\in\SmVar(\mathbb C)$. 
Let $M\in\DA_c(S)$ and $F\in C(\Var(\mathbb C)^{sm}/S)$ such that $M=D(\mathbb A^1,et)(F)$.
We have, using definition \ref{wtildewan}(ii), 
the canonical map in $D_{O_Sfil,\mathcal D^{\infty},\infty}(S^{an}/(\tilde S^{an}_I))$
\begin{eqnarray*}
T(\mathcal F^{GM}_{S,an},\mathcal F^{FDR}_{S,an})(M): \\
\mathcal F_{S,an}^{GM}(L\mathbb D_SM):=
(e(\tilde S_I)_*\mathcal Hom^{\bullet}(\An_{\tilde S_I}^*L(i_{I*}j_I^*\mathbb D_SLF),
E_{usu}(\Omega^{\bullet}_{/\tilde S_I},F_b)),u^q_{IJ}(F)) \\
\xrightarrow{\sim}
(e(\tilde S_I)_*\mathcal Hom^{\bullet}(\An_{\tilde S_I}^*L\mathbb D^0_{\tilde S_I}(L(i_{I*}j_I^*F)),
E_{usu}(\Omega^{\bullet}_{/\tilde S_I},F_b)),u^{q,d}_{IJ}(F)) \\
\xrightarrow{\mathcal Hom(-,\Gr(\Omega_{\tilde S^{an}_I}))^{-1}}
J_S(e(\tilde S_I)_*\mathcal Hom^{\bullet}(\An_{\tilde S_I}^*L\mathbb D^0_{\tilde S_I}(L(i_{I*}j_I^*F)),
\Gr^{12}_{\tilde S_I*}E_{usu}(\Omega^{\bullet,\Gamma,pr,an}_{/\tilde S_I},F_{DR})),u^{q,d}_{IJ}(F)) \\
\xrightarrow{(\mathcal Hom^{\bullet}(\An_{\tilde S_I}^*L(r^{0CH}_{\tilde S_I}(L(i_{I*}j_I^*F))\circ 
T(\hat R^{0CH},R^{0CH})(L(i_{I*}j_I^*F))),-))} \\
J_S(e(\tilde S_I)_*\mathcal Hom^{\bullet}(\An_{\tilde S_I}^*\hat R^{0CH}(\rho_{\tilde S_I}^*L(i_{I*}j_I^*F)),
\Gr^{12}_{\tilde S_I*}E_{usu}(\Omega^{\bullet,\Gamma,pr,an}_{/\tilde S_I},F_{DR})),u^{q,d}_{IJ}(F)) \\
\xrightarrow{I(\Gr_{\tilde S_I}^{12*},\Gr^{12}_{\tilde S_I*})(-,-)} \\
J_S(e'(\tilde S_I)_*\mathcal Hom(\An_{\tilde S_I}^*\hat R^{CH}(\rho_{\tilde S_I}^*L(i_{I*}j_I^*F)),
E_{usu}(\Omega^{\bullet,\Gamma,pr,an}_{/\tilde S_I},F_{DR}))[-d_{\tilde S_I}],u^q_{IJ}(F)) 
=:J_S(\mathcal F^{FDR}_{S,an}(M))
\end{eqnarray*}
\end{defi}

We now define the functorialities of $\mathcal F_S^{FDR}$ with respect to $S$ 
which makes $\mathcal F^{-}_{FDR}$ a morphism of 2 functor.

\begin{defi}\label{TGammaFDRan}
Let $S\in\Var(\mathbb C)$. Let $Z\subset S$ a closed subset.
Let $S=\cup_{i=1}^l S_i$ an open cover such that there exist closed embeddings
$i_i:S_i\hookrightarrow\tilde S_i$ with $\tilde S_i\in\SmVar(\mathbb C)$.
Denote $Z_I:=Z\cap S_I$. We then have closed embeddings $Z_I\hookrightarrow S_I\hookrightarrow\tilde S_I$.
\begin{itemize}
\item[(i)]For $F\in C(\Var(\mathbb C)^{sm}/S)$, we will consider the following canonical map 
in $\pi_S(D(MHM(S^{an})))\subset D_{\mathcal D(1,0)fil}(S^{an}/(\tilde S^{an}_I))$
\begin{eqnarray*}
T(\Gamma_Z^{\vee,Hdg},\Omega^{\Gamma,pr,an}_{/S})(F,W):\\
\Gamma_Z^{\vee,Hdg}\iota_S^{-1}(
e'(\tilde S_I)_*\mathcal Hom^{\bullet}(\An_{\tilde S_I}^*\hat R^{CH}(\rho_{\tilde S_I}^*L(i_{I*}j_I^*(F,W))), 
E_{usu}(\Omega^{\bullet,\Gamma,pr,an}_{/\tilde S_I},F_{DR}))[-d_{\tilde S_I}], u^q_{IJ}(F,W)) \\
\xrightarrow{\mathcal Hom^{\bullet}(\An_{\tilde S_I}^*\hat R_{\tilde S_I}^{CH}(\gamma^{\vee,Z_I}(L(i_{I*}j_I^*(F,W)))),
E_{usu}(\Omega^{\bullet,\Gamma,pr,an}_{/\tilde S_I},F_{DR}))} \\
\Gamma_Z^{\vee,Hdg}\iota_S^{-1}(
e'(\tilde S_I)_*\mathcal Hom^{\bullet}(\An_{\tilde S_I}^*\hat R^{CH}(\rho_{\tilde S_I}^*\Gamma^{\vee}_{Z_I}L(i_{I*}j_I^*(F,W))),  
E_{usu}(\Omega^{\bullet,\Gamma,pr,an}_{/\tilde S_I},F_{DR}))[-d_{\tilde S_I}], u^{q,Z}_{IJ}(F,W)) \\
\xrightarrow{=} 
\iota_S^{-1}(e'(\tilde S_I)_*\mathcal Hom^{\bullet}(
\An_{\tilde S_I}^*\hat R^{CH}(\rho_{\tilde S_I}^*\Gamma^{\vee}_{Z_I}L(i_{I*}j_I^*(F,W))),  
E_{usu}(\Omega^{\bullet,\Gamma,pr}_{/\tilde S_I},F_{DR}))[-d_{\tilde S_I}], u^{q,Z}_{IJ}(F,W)),
\end{eqnarray*}
with $u^{q,Z}_{IJ}(F)$ given as in definition \ref{TGammaFDR}(i).
\item[(ii)]For $F\in C(\Var(\mathbb C)^{sm}/S)$, we have also the following canonical map 
in $\pi_S(D(MHM(S^{an})))\subset D_{\mathcal D(1,0)fil}(S^{an}/(\tilde S^{an}_I))$
\begin{eqnarray*}
T(\Gamma_Z^{Hdg},\Omega^{\Gamma,pr}_{/S})(F,W):\\
\iota_S^{-1}(e'(\tilde S_I)_*\mathcal Hom^{\bullet}(
\An_{\tilde S_I}^*\hat R^{CH}(\rho_{\tilde S_I}^*L\Gamma_{Z_I}E(i_{I*}j_I^*\mathbb D_S(F,W))), 
E_{usu}(\Omega^{\bullet,\Gamma,pr}_{/\tilde S_I},F_{DR}))[-d_{\tilde S_I}], u^{q,Z,d}_{IJ}(F,W)) \\
\xrightarrow{=}
\Gamma_Z^{Hdg}\iota_S^{-1}(e'(\tilde S_I)_*\mathcal Hom^{\bullet}(
\An_{\tilde S_I}^*\hat R^{CH}(\rho_{\tilde S_I}^*L\Gamma_{Z_I}E(i_{I*}j_I^*\mathbb D_S(F,W))),  
E_{usu}(\Omega^{\bullet,\Gamma,pr}_{/\tilde S_I},F_{DR}))[-d_{\tilde S_I}], u^{q,Z,d}_{IJ}(F,W)) \\
\xrightarrow{\mathcal Hom^{\bullet}(\An_{\tilde S_I}^*\hat R_{\tilde S_I}^{CH}(\gamma^{Z_I}(-)),
E_{usu}(\Omega^{\bullet,\Gamma,pr}_{/\tilde S_I},F_{DR}))} \\
\Gamma_Z^{Hdg}\iota_S^{-1}(e'(\tilde S_I)_*\mathcal Hom^{\bullet}(
\An_{\tilde S_I}^*\hat R^{CH}(\rho_{\tilde S_I}^*L(i_{I*}j_I^*\mathbb D_S(F,W))), 
E_{usu}(\Omega^{\bullet,\Gamma,pr}_{/\tilde S_I},F_{DR}))[-d_{\tilde S_I}], u^q_{IJ}(F,W)) 
\end{eqnarray*}
with $u^{q,Z}_{IJ}(F)$ given as in definition \ref{TGammaFDR}(ii).
\end{itemize}
\end{defi}

\begin{defi}\label{TgDRdefsingan}
Let $g:T\to S$ a morphism with $T,S\in\Var(\mathbb C)$.
Assume we have a factorization $g:T\xrightarrow{l}Y\times S\xrightarrow{p_S}S$
with $Y\in\SmVar(\mathbb C)$, $l$ a closed embedding and $p_S$ the projection.
Let $S=\cup_{i=1}^lS_i$ be an open cover such that 
there exists closed embeddings $i_i:S_i\hookrightarrow\tilde S_i$ with $\tilde S_i\in\SmVar(\mathbb C)$
Then, $T=\cup^l_{i=1} T_i$ with $T_i:=g^{-1}(S_i)$
and we have closed embeddings $i'_i:=i_i\circ l:T_i\hookrightarrow Y\times\tilde S_i$,
Moreover $\tilde g_I:=p_{\tilde S_I}:Y\times\tilde S_I\to\tilde S_I$ is a lift of $g_I:=g_{|T_I}:T_I\to S_I$.
Let $M\in\DA_c(S)$ and $(F,W)\in C_{fil}(\Var(\mathbb C)^{sm}/S)$ such that $(M,W)=D(\mathbb A^1_S,et)(F,W)$.
Then, $D(\mathbb A^1_T,et)(g^*F)=g^*M$ and there exist 
$(F',W)\in C_{fil}(\Var(\mathbb C)^{sm}/S)$ and an equivalence $(\mathbb A^1,et)$ local $e:g^*(F,W)\to (F',W)$ 
such that $D(\mathbb A^1_T,et)(F',W)=(g^*M,W)$.  
We have, using definition \ref{TgDRan} and definition \ref{TGammaFDRan}(i), 
the canonical map in $\pi_T(D(MHM(T^{an})))\subset D_{\mathcal D(1,0)fil}(T^{an}/(Y^{an}\times\tilde S^{an}_I))$
\begin{eqnarray*}
T(g,\mathcal F_{an}^{FDR})(M):g^{\hat*mod}_{Hdg}\iota_S^{-1}\mathcal F_{S,an}^{FDR}(M):= \\
(\Gamma^{\vee,Hdg}_T\iota_T^{-1}(\tilde g_I^{*mod}
(e'(\tilde S_I)_*\mathcal Hom^{\bullet}(\An_{\tilde S_I}^*\hat R^{CH}(\rho_{\tilde S_I}^*(L(i_{I*}j_I^*(F,W))),  
E_{usu}(\Omega^{\bullet,\Gamma,pr,an}_{/\tilde S_I},F_{DR})))[d_{YI}],\tilde g_J^{*mod}u^q_{IJ}(F,W)) \\
\xrightarrow{(T(\tilde g_I,\Omega^{\Gamma,pr,an}_{/\cdot})(-)} \\
\Gamma^{\vee,Hdg}_T\iota_T^{-1}(e'(-)_*
\mathcal Hom(\An_{Y\times\tilde S_I}^*\tilde g_I^*\hat R^{CH}(\rho_{\tilde S_I}^*L(i_{I*}j_I^*(F,W))), 
E_{usu}(\Omega^{\bullet,\Gamma,pr,an}_{/Y\times\tilde S_I},F_{DR}))[d_{YI}],\tilde g_J^*u^q_{IJ}(F,W)) \\
\xrightarrow{\mathcal Hom(T(\tilde g_I,\hat R^{CH})(-)^{-1},-)} \\ 
\Gamma^{\vee,Hdg}_T\iota_T^{-1}(e'_*\mathcal Hom(
\An_{Y\times\tilde S_I}^*\hat R^{CH}(\rho_{Y\times\tilde S_I}^*\tilde g_I^*L(i_{I*}j_I^*(F,W))), 
E_{usu}(\Omega^{\bullet,\Gamma,pr,an}_{/Y\times\tilde S_I},F_{DR}))[d_{YI}],\tilde g_J^*u^q_{IJ}(F,W)) \\
\xrightarrow{T(\Gamma_T^{\vee,Hdg},\Omega^{\Gamma,pr,an}_{/S})(F,W)}\\
\iota_T^{-1}(e'_*\mathcal Hom(
\An_{Y\times\tilde S_I}^*\hat R^{CH}(\rho_{Y\times\tilde S_I}^*\Gamma^{\vee}_{T_I}\tilde g_I^*L(i_{I*}j_I^*(F,W))), 
E_{usu}(\Omega^{\bullet,\Gamma,pr,an}_{/Y\times\tilde S_I},F_{DR}))[d_{YI}],\tilde g_J^{*,\gamma}u^q_{IJ}(F,W)) \\
\xrightarrow{(\mathcal Hom(\An_{Y\times\tilde S_I}^*\hat R^{CH}_{Y\times\tilde S_I}(T^{q,\gamma}(D_{gI})(j_I^*(F,W))),
E_{usu}(\Omega^{\bullet,\Gamma,pr,an}_{/Y\times\tilde S_I},F_{DR}))[d_{YI}])} \\
\iota_T^{-1}(e'_*\mathcal Hom(\An_{Y\times\tilde S_I}^*\hat R^{CH}(\rho_{Y\times\tilde S_I}^*L(i'_{I*}j^{'*}_Ig^*(F,W))), 
E_{et}(\Omega^{\bullet,\Gamma,pr,an}_{/Y\times\tilde S_I},F_{DR}))[d_{YI}],u^q_{IJ}(g^*(F,W))) \\
\xrightarrow{\mathcal Hom(\hat R^{CH}_{Y\times\tilde S_I}(Li'_{I*}j^{'*}_I(e)),)} \\
\iota_T^{-1}(e'_*\mathcal Hom(\An_{Y\times\tilde S_I}^*\hat R^{CH}(\rho_{Y\times\tilde S_I}^*L(i'_{I*}j^{'*}_I(F',W))), 
E_{usu}(\Omega^{\bullet,\Gamma,pr,an}_{/Y\times\tilde S_I},F_{DR}))[d_{YI}],u^q_{IJ}(F',W))
\xrightarrow{=:}\mathcal F_{T,an}^{FDR}(g^*M)
\end{eqnarray*}
\end{defi}

\begin{defi}\label{SixTalgan}
\begin{itemize}
\item Let $f:X\to S$ a morphism with $X,S\in\Var(\mathbb C)$. Assume there exist a factorization
$f:X\xrightarrow{l}Y\times S\xrightarrow{p_S}S$ with $Y\in\SmVar(\mathbb C)$, $l$ a closed embedding and $p_S$ the projection.
We have, for $M\in\DA_c(X)$, the following transformation map in $\pi_S(D(MHM(S^{an})))$
\begin{eqnarray*}
T_*(f,\mathcal F_{an}^{FDR})(M):\mathcal F_{S,an}^{FDR}(Rf_*M)
\xrightarrow{\ad(f_{Hdg}^{\hat*mod},Rf^{Hdg}_*)(-)}Rf^{Hdg}_*f^{\hat*mod}_{Hdg}\mathcal F_{S,an}^{FDR}(Rf_*M) \\
\xrightarrow{T(f,\mathcal F_{an}^{FDR})(Rf_*M)}Rf^{Hdg}_*\mathcal F_{X,an}^{FDR}(f^*Rf_*M)
\xrightarrow{\mathcal F_X^{FDR}(\ad(f^*,Rf_*)(M))}Rf^{Hdg}_*\mathcal F_{X,an}^{FDR}(M)
\end{eqnarray*}
Clearly, for $p:Y\times S\to S$ a projection with $Y\in\PSmVar(\mathbb C)$, we have, for $M\in\DA_c(Y\times S)$,
$T_*(p,\mathcal F^{FDR})(M)=T_!(p,\mathcal F^{FDR})(M)[2d_Y]$

\item Let $S\in\Var(\mathbb C)$. Let $Y\in\SmVar(\mathbb C)$ and $p:Y\times S\to S$ the projection.
We have then, for $M\in\DA(Y\times S)$ the following transformation map in $\pi_S(D(MHM(S^{an})))$
\begin{eqnarray*}
T_!(p,\mathcal F_{an}^{FDR})(M):p^{Hdg}_!\mathcal F_{Y\times S,an}^{FDR}(M)
\xrightarrow{\mathcal F_{Y\times S,an}^{FDR}(\ad(Lp_{\sharp},p^*)(M))}
Rp^{Hdg}_!\mathcal F_{Y\times S,an}^{FDR}(p^*Lp_{\sharp}(M)) \\
\xrightarrow{T(p,\mathcal F^{FDR})(Lp_{\sharp}(M,W))}Rp^{Hdg}_!p^{\hat*mod[-]}\mathcal F_{S,an}^{FDR}(Lp_{\sharp}M)
\xrightarrow{T(p^{*mod},p^{\hat*mod})(-)}p^{Hdg}_!p^{*mod[-]} \\
\mathcal F_{S,an}^{FDR}(Lp_{\sharp}M)\xrightarrow{\ad(Rp^{Hdg}_!,p^{*mod[-]})(\mathcal F_{S,an}^{FDR}(Lp_{\sharp}M))}
\mathcal F_{S,an}^{FDR}(Lp_{\sharp}M)
\end{eqnarray*}

\item Let $f:X\to S$ a morphism with $X,S\in\Var(\mathbb C)$. Assume there exist a factorization
$f:X\xrightarrow{l}Y\times S\xrightarrow{p_S}S$ with $Y\in\SmVar(\mathbb C)$, $l$ a closed embedding and $p_S$ the projection.
We have then, using the second point, for $M\in\DA(X)$ the following transformation map in $\pi_S(D(MHM(S^{an})))$
\begin{eqnarray*}
T_!(f,\mathcal F_{an}^{FDR})(M):
Rp^{Hdg}_!\mathcal F_X^{FDR}(M):=Rp^{Hdg}_!\mathcal F_{Y\times S,an}^{FDR}(l_*M) \\
\xrightarrow{T_!(p,\mathcal F_{an}^{FDR})(l_*M)}\mathcal F_{S,an}^{FDR}(Lp_{\sharp}l_*M)
\xrightarrow{=}\mathcal F_{S,an}^{FDR}(Rf_!M)
\end{eqnarray*}

\item Let $f:X\to S$ a morphism with $X,S\in\Var(\mathbb C)$. Assume there exist a factorization
$f:X\xrightarrow{l}Y\times S\xrightarrow{p_S}S$ with $Y\in\SmVar(\mathbb C)$, $l$ a closed embedding and $p_S$ the projection.
We have, using the third point, for $M\in\DA(S)$, the following transformation map in in $\pi_X(D(MHM(X^{an})))$
\begin{eqnarray*}
T^!(f,\mathcal F_{an}^{FDR})(M):\mathcal F_{X,an}^{FDR}(f^!(M,W))
\xrightarrow{\ad(Rf^{Hdg}_!,f^{*mod}_{Hdg})(\mathcal F_{X,an}^{FDR}(f^!M))}
f^{*mod}_{Hdg}Rf^{Hdg}_!\mathcal F_{X,an}^{FDR}(f^!M) \\
\xrightarrow{T_!(p_S,\mathcal F_{an}^{FDR})(\mathcal F_{an}^{FDR}(f^!M))}
f^{*mod}_{Hdg}\mathcal F_{S,an}^{FDR}(Rf_!f^!(M,W))
\xrightarrow{\mathcal F_{S,an}^{FDR}(\ad(Rf_!,f^!)(M))}f^{*mod}_{Hdg}\mathcal F_{S,an}^{FDR}(M)
\end{eqnarray*}

\item Let $S\in\Var(\mathbb C)$. Let $S=\cup_{i=1}^l S_i$ an open cover such that there exist closed embeddings 
$i_i:S_i\hookrightarrow\tilde S_i$ with $\tilde S_i\in\SmVar(\mathbb C)$. 
We have, using the preceding point, 
denoting $\Delta_S:S\hookrightarrow S$ the diagonal closed embedding 
and $p_1:S\times S\to S$, $p_2:S\times S\to S$ the projections, 
for $M,N\in\DA(S)$ and $(F,W),(G,W)\in C_{fil}(\Var(\mathbb C)^{sm}/S))$ 
such that $(M,W)=D(\mathbb A^1,et)(F,W)$ and $(N,W)=D(\mathbb A^1,et)(G,W)$, 
the following transformation map in $\pi_S(D(MHM(S)))$
\begin{eqnarray*}
T(\mathcal F_{S,an}^{FDR},\otimes)(M,N):\mathcal F_{S,an}^{FDR}(M)\otimes^{Hdg}_{O_S}\mathcal F_{S,an}^{FDR}(N) 
:=\Delta_S^{!Hdg}(p_1^{*mod}\mathcal F_{S,an}^{FDR}(M)\otimes_{O_{S\times S}}p_2^{*mod}\mathcal F_{S,an}^{FDR}(N) \\
\xrightarrow{T^!(p_1,\mathcal F_{S,an}^{FDR})(M)\otimes T^!(p_1,\mathcal F_{S,an}^{FDR})(M)}
\Delta_S^{!Hdg}(\mathcal F_{S\times S,an}^{FDR}(p_1^!M)\otimes_{O_{S\times S}}\mathcal F_{S\times S,an}^{FDR}(p_2^!N) \\
\xrightarrow{(T(\otimes,\Omega)(\hat R^{CH}(\rho_{\tilde S_I\times\tilde S_J}^*L(i_I\times i_J)_*(j_I\times j_J)^*p_1^*F[2d_S]),
\hat R^{CH}(\rho_{\tilde S_I\times\tilde S_J}^*L(i_I\times i_J)_*(j_I\times j_J)^*p_2^*F[2d_S])))} \\
\Delta_S^{!Hdg}(\mathcal F_{S\times S,an}^{FDR}(p_1^!M\otimes p_2^!N)
\xrightarrow{T^!(\Delta_S,\mathcal F_{an}^{FDR})(p_1^!M\otimes p_2^!N)}
\mathcal F_{S,an}^{FDR}(\Delta_S^!(p_1^!M\otimes p_2^!N))=\mathcal F_{S,an}^{FDR}(M\otimes N)
\end{eqnarray*}
\end{itemize}
\end{defi}

\begin{prop}\label{Tgpropan}
Let $g:T\to S$ a morphism with $T,S\in\Var(\mathbb C)$.
Assume we have a factorization $g:T\xrightarrow{l}Y_2\times S\xrightarrow{p_S}S$
with $Y_2\in\SmVar(\mathbb C)$, $l$ a closed embedding and $p_S$ the projection.
Let $S=\cup_{i=1}^lS_i$ be an open cover such that 
there exists closed embeddings $i_i:S_i\hookrightarrow\tilde S_i$ with $\tilde S_i\in\SmVar(\mathbb C)$
Then, $T=\cup^l_{i=1} T_i$ with $T_i:=g^{-1}(S_i)$
and we have closed embeddings $i'_i:=i_i\circ l:T_i\hookrightarrow Y_2\times\tilde S_i$,
Moreover $\tilde g_I:=p_{\tilde S_I}:Y\times\tilde S_I\to\tilde S_I$ is a lift of $g_I:=g_{|T_I}:T_I\to S_I$.
Let $f:X\to S$ a  morphism with $X\in\Var(\mathbb C)$ such that there exists a factorization 
$f:X\xrightarrow{l}Y_1\times S\xrightarrow{p_S} S$, with $Y_1\in\SmVar(\mathbb C)$, 
$l$ a closed embedding and $p_S$ the projection. We have then the following commutative diagram
whose squares are cartesians
\begin{equation*}
\xymatrix{f':X_T\ar[r]\ar[rd]\ar[dd]_{g'} & Y_1\times T\ar[rd]\ar[r] & T\ar[rd]\ar[dd]^{g} & \, \\
\, & Y_1\times X\ar[r]\ar[ld] & Y_1\times Y_2\times S\ar[r]\ar[ld] & Y_2\times S\ar[ld] \\
f:X\ar[r] & Y_1\times S\ar[r] & S & \,}
\end{equation*} 
Take a smooth compactification $\bar Y_1\in\PSmVar(\mathbb C)$ of $Y_1$, denote
$\bar X_I\subset\bar Y_1\times\tilde S_I$ the closure of $X_I$, and $Z_I:=\bar X_I\backslash X_I$.
Consider $F(X/S):=p_{S,\sharp}\Gamma_X^{\vee}\mathbb Z(Y_1\times S/Y_1\times S)\in C(\Var(\mathbb C)^{sm}/S)$ and
the isomorphism in $C(\Var(\mathbb C)^{sm}/T)$
\begin{eqnarray*}
T(f,g,F(X/S)):g^*F(X/S):=g^*p_{S,\sharp}\Gamma_X^{\vee}\mathbb Z(Y_1\times S/Y_1\times S)\xrightarrow{\sim} \\
p_{T,\sharp}\Gamma_{X_T}^{\vee}\mathbb Z(Y_1\times T/Y_1\times T)=:F(X_T/T).
\end{eqnarray*}
which gives in $\DA(T)$ the isomorphism $T(f,g,F(X/S)):g^*M(X/S)\xrightarrow{\sim}M(X_T/T)$.
Then the following diagram in $\pi_T(D(MHM(T)))\subset D_{\mathcal D(1,0)fil,\infty}(T/(Y_2\times\tilde S_I))$, 
where the horizontal maps are given by proposition \ref{keysing1}, commutes
\begin{equation*}
\begin{tikzcd}
g^{\hat*mod}_{Hdg}\iota_S^{-1}\mathcal F_{S,an}^{FDR}(M(X/S))
\ar[dd,"'T(g{,}\mathcal F^{FDR})(M(X/S))"']\ar[rr,"g^{\hat*mod}_{Hdg}I(X/S)"] & \, &
g^{\hat*mod}_{Hdg}Rf^{Hdg}_!(\Gamma^{\vee,Hdg}_{X_I}(O_{(Y_1\times\tilde S_I)^{an}},F_b)(d_{Y_1})[2d_{Y_1}],x_{IJ}(X/S))
\ar[d,"T(p_{\tilde S_I}{,}\gamma^{\vee{,}Hdg})(-)"] \\
\, & \, & 
Rf^{'Hdg}_!g^{'\hat*mod}_{Hdg}(\Gamma^{\vee,Hdg}_{X_I}(O_{(Y_1\times\tilde S_I)^{an}},F_b)(d_{Y_1})[2d_{Y_1}],x_{IJ}(X/S))
\ar[d,"T(p_{Y_1\times Y_2\times\tilde S_I{,}Hdg}^{\hat*mod}{,}p_{Y_1\times Y_2\times\tilde S_I{,}Hdg}^{*mod})(-)"] \\
\iota_T^{-1}\mathcal F_{T,an}^{FDR}(M(X_T/T))\ar[rr,"I(X_T/T)"] & \, & 
Rf^{'Hdg}_!(\Gamma^{\vee,Hdg}_{X_{T_I}}(O_{(Y_2\times Y_1\times\tilde S_I)^{an}},F_b)(d_{Y_{12}})[2d_{Y_{12}}],x_{IJ}(X_T/T)).
\end{tikzcd}
\end{equation*} 
with $d_{Y_{12}}=d_{Y_1}+d_{Y_2}$.
\end{prop}

\begin{proof}
Follows immediately from definition.
\end{proof}

\begin{prop}\label{mainthmprop2an}
Let $S\in\Var(\mathbb C)$. Let $Y\in\SmVar(\mathbb C)$ and $p:Y\times S\to S$ the projection.
Let $S=\cup_{i=1}^l S_i$ an open cover such that there exist closed embeddings 
$i^o_i:S_i\hookrightarrow\tilde S_i$ with $\tilde S_i\in\SmVar(\mathbb C)$.
For $I\subset\left[1,\cdots l\right]$, we denote by $S_I=\cap_{i\in I} S_i$, $j^o_I:S_I\hookrightarrow S$ and
$j_I:Y\times S_I\hookrightarrow Y\times S$ the open embeddings. 
We then have closed embeddings $i_I:Y\times S_I\hookrightarrow Y\times\tilde S_I$.
and we denote by $p_{\tilde S_I}:Y\times\tilde S_I\to\tilde S_I$ the projections.
Let $f':X'\to Y\times S$ a morphism, with $X'\in\Var(\mathbb C)$ such that there exists a factorization
$f':X'\xrightarrow{l'}Y'\times Y\times S\xrightarrow{p'} Y\times S$ 
with $Y'\in\SmVar(\mathbb C)$, $l'$ a closed embedding and $p'$ the projection. 
Denoting $X'_I:=f^{'-1}(Y\times S_I)$, we have closed embeddings $i'_I:X'_I\hookrightarrow Y'\times Y\times\tilde S_I$
Consider 
\begin{eqnarray*}
F(X'/Y\times S):=p_{Y\times S,\sharp}\Gamma_{X'}^{\vee}\mathbb Z(Y'\times Y\times S/Y'\times Y\times S)
\in C(\Var(\mathbb C)^{sm}/Y\times S)
\end{eqnarray*}
and $F(X'/S):=p_{\sharp}F(X'/Y\times S)\in C(\Var(\mathbb C)^{sm}/S)$, 
so that $Lp_{\sharp}M(X'/Y\times S)[-2d_Y]=:M(X'/S)$. 
Then, the following diagram in 
$\pi_S(D(MHM(S^{an})))\subset D_{\mathcal D(1,0)fil}(S^{an}/(Y^{an}\times\tilde S^{an}_I))$,
where the vertical maps are given by proposition \ref{keysing1}, commutes
\begin{equation*}
\xymatrix{
Rp^{Hdg}!\mathcal F_{Y\times S,an}^{FDR}(M(X'/Y\times S))
\ar[rr]^{T_!(p,\mathcal F_{an}^{FDR})(M(X'/Y\times S))} & \, & 
\mathcal F_{S,an}^{FDR}(M(X'/S)) \\
Rp^{Hdg}!Rf^{'Hdg}_!f_{Hdg}^{'*mod}\mathbb Z_{(Y\times S)^{an}}^{Hdg} 
\ar[rr]^{=}\ar[u]^{T(p_{Hdg}^{\hat*mod}{,}p_{Hdg}^{*mod})(-)\circ Rp^{Hdg}!(I(X'/Y\times S))} & \, &
Rf^{Hdg}_!f_{Hdg}^{*mod}\mathbb Z_{S^{an}}^{Hdg}
\ar[u]_{I(X'/S)}} 
\end{equation*}
\end{prop}

\begin{proof}
Immediate from definition.
\end{proof}

\begin{prop}\label{TotimesDRpropan}
Let $f_1:X_1\to S$, $f_2:X_2\to S$ two morphism with $X_1,X_2,S\in\Var(\mathbb C)$. 
Assume that there exist factorizations 
$f_1:X_1\xrightarrow{l_1}Y_1\times S\xrightarrow{p_S} S$, $f_2:X_2\xrightarrow{l_2}Y_2\times S\xrightarrow{p_S} S$
with $Y_1,Y_2\in\SmVar(\mathbb C)$, $l_1,l_2$ closed embeddings and $p_S$ the projections.
We have then the factorization
\begin{equation*}
f_{12}:=f_1\times f_2:X_{12}:=X_1\times_S X_2\xrightarrow{l_1\times l_2}Y_1\times Y_2\times S\xrightarrow{p_S} S
\end{equation*}
Let $S=\cup_{i=1}^l S_i$ an open affine covering and denote, 
for $I\subset\left[1,\cdots l\right]$, $S_I=\cap_{i\in I} S_i$ and $j_I:S_I\hookrightarrow S$ the open embedding.
Let $i_i:S_i\hookrightarrow\tilde S_i$ closed embeddings, with $\tilde S_i\in\SmVar(\mathbb C)$. 
We have then the following commutative diagram in 
$\pi_S(DMHM(S^{an}))\subset D_{\mathcal D(1,0)fil}(S^{an}/(\tilde S^{an}_I))$
where the vertical maps are given by proposition \ref{keysing1}
\begin{equation*}
\begin{tikzcd}
\mathcal F_{S,an}^{FDR}(M(X_1/S))\otimes^{Hdg}_{O_S}\mathcal F_{S,an}^{FDR}(M(X_2/S))
\ar[rr,"I(X_1/S)\otimes I(X_2/S)"] \ar[d,"T(\mathcal F_{S,an}^{FDR}{,}\otimes)(M(X_1/S){,}M(X_2/S))"]  & \, &
\shortstack{$Rf_{1!}^{Hdg}(\Gamma^{\vee,Hdg}_{X_{1I}}(O_{(Y_1\times\tilde S_I)^{an}},F_b)(d_1)[2d_1],x_{IJ}(X_1/S))\otimes_{O_S}$ \\ 
$Rf_{2!}^{Hdg}(\Gamma^{\vee,Hdg}_{X_{2I}}(O_{(Y_2\times\tilde S_I)^{an}},F_b)(d_2)[2d_2],x_{IJ}(X_2/S))$}
\ar[d,"(Ew_{(Y_1\times\tilde S_I{,}Y_2\times\tilde S_I)/\tilde S_I})"] \\
\mathcal F_{S,an}^{FDR}(M(X_1/S)\otimes M(X_2/S)=M(X_1\times_S X_2/S))\ar[rr,"I(X_{12}/S)"] & \, &
Rf_{12!}^{Hdg}(\Gamma^{\vee,Hdg}_{X_{1I}\times_S X_{2I}}(O_{(Y_1\times Y_2\times\tilde S_I)^{an}},F_b)(d_{12})[2d_{12}],
x_{IJ}(X_1/S)).
\end{tikzcd}
\end{equation*}
\end{prop}

\begin{proof}
Immediate from definition.
\end{proof}

\begin{thm}\label{mainthmCoran}
\begin{itemize}
\item[(i)]Let $g:T\to S$ a morphism, with $S,T\in\Var(\mathbb C)$. 
Assume we have a factorization $g:T\xrightarrow{l}Y\times S\xrightarrow{p_S}S$
with $Y\in\SmVar(\mathbb C)$, $l$ a closed embedding and $p_S$ the projection. Let $M\in\DA_c(S)$. 
Then map in $\pi_T(D(MHM(T^{an})))$ 
\begin{eqnarray*}
T(g,\mathcal F_{an}^{FDR})(M):g^{\hat*mod}_{Hdg}\mathcal F_{S,an}^{FDR}(M)
\xrightarrow{\sim}\mathcal F_{T,an}^{FDR}(g^*M)
\end{eqnarray*} 
given in definition \ref{TgDRdefsingan} is an isomorphism.
\item[(ii)] Let $f:X\to S$ a morphism with $X,S\in\Var(\mathbb C)$. Assume there exist a factorization
$f:X\xrightarrow{l}Y\times S\xrightarrow{p_S}S$ with $Y\in\SmVar(\mathbb C)$, $l$ a closed embedding and $p_S$ the projection.
Then, for $M\in\DA_c(X)$,
\begin{equation*}
T_!(f,\mathcal F_{an}^{FDR})(M):Rf^{Hdg}_!\mathcal F_{X,an}^{FDR}(M)\xrightarrow{\sim}\mathcal F_{S,an}^{FDR}(Rf_!M)
\end{equation*}
is an isomorphism in $\pi_S(D(MHM(S^{an}))$
\item[(iii)] Let $f:X\to S$ a morphism with $X,S\in\Var(\mathbb C)$, $S$ quasi-projective. Assume there exist a factorization
$f:X\xrightarrow{l}Y\times S\xrightarrow{p_S}S$ with $Y\in\SmVar(\mathbb C)$, $l$ a closed embedding and $p_S$ the projection.
We have, for $M\in\DA_c(X)$,
\begin{equation*}
T_*(f,\mathcal F_{an}^{FDR})(M):\mathcal F_{S,an}^{FDR}(Rf_*M)\xrightarrow{\sim}Rf^{Hdg}_*\mathcal F_{X,an}^{FDR}(M)
\end{equation*}
is an isomorphism in $\pi_S(D(MHM(S^{an}))$.
\item[(iv)] Let $f:X\to S$ a morphism with $X,S\in\Var(\mathbb C)$, $S$ quasi-projective. Assume there exist a factorization
$f:X\xrightarrow{l}Y\times S\xrightarrow{p_S}S$ with $Y\in\SmVar(\mathbb C)$, $l$ a closed embedding and $p_S$ the projection.
Then, for $M\in\DA_c(S)$
\begin{equation*}
T^!(f,\mathcal F_{an}^{FDR})(M):\mathcal F_{X,an}^{FDR}(f^!M)\xrightarrow{\sim}f^{*mod}_{Hdg}\mathcal F_{S,an}^{FDR}(M)
\end{equation*}
is an isomorphism in $\pi_X(D(MHM(X^{an}))$. 
\item[(v)]Let $S\in\Var(\mathbb C)$ and $S=\cup_{i=1}^l S_i$ an open affine covering and denote, 
for $I\subset\left[1,\cdots l\right]$, $S_I=\cap_{i\in I} S_i$ and $j_I:S_I\hookrightarrow S$ the open embedding.
Let $i_i:S_i\hookrightarrow\tilde S_i$ closed embeddings, with $\tilde S_i\in\SmVar(\mathbb C)$. 
Then, for $M,N\in\DA_c(S)$, the map in $\pi_S(D(MHM(S^{an})))$ 
\begin{eqnarray*}
T(\mathcal F_{S,an}^{FDR},\otimes)(M,N): 
\mathcal F_{S,an}^{FDR}(M)\otimes^{Hdg}_{O_S}\mathcal F_{S,an}^{FDR}(N)\xrightarrow{\sim}\mathcal F_{S,an}^{FDR}(M\otimes N)
\end{eqnarray*}
given in definition is an isomorphism $\pi_S(D(MHM(S^{an}))$.
\end{itemize}
\end{thm}

\begin{proof}
\noindent(i):Follows from proposition \ref{Tgpropan} and proposition \ref{keysing1}.

\noindent(ii),(iii),(iv):Similar to the proof of theorem \ref{mainthm}.

\noindent(v):Follows from proposition \ref{TotimesDRpropan}.
\end{proof}

\begin{prop}\label{TgDRGMdiaan}
Let $g:T\to S$ a morphism with $T,S\in\Var(\mathbb C)$.
Assume we have a factorization $g:T\xrightarrow{l}Y\times S\xrightarrow{p_S}S$
with $Y\in\SmVar(\mathbb C)$, $l$ a closed embedding and $p_S$ the projection.
Let $S=\cup_{i=1}^lS_i$ be an open cover such that 
there exists closed embeddings $i_i:S_i\hookrightarrow\tilde S_i$ with $\tilde S_i\in\SmVar(\mathbb C)$
Then, $T=\cup^l_{i=1} T_i$ with $T_i:=g^{-1}(S_i)$
and we have closed embeddings $i'_i:=i_i\circ l:T_i\hookrightarrow Y\times\tilde S_i$,
Moreover $\tilde g_I:=p_{\tilde S_I}:Y\times\tilde S_I\to\tilde S_I$ is a lift of $g_I:=g_{|T_I}:T_I\to S_I$.
Let $M\in\DA_c(S)$ and $F\in C(\Var(\mathbb C)^{sm}/S)$ such that $M=D(\mathbb A^1_S,et)(F)$.
Then, $D(\mathbb A^1_T,et)(g^*F)=g^*M$. 
Then the following diagram in $D_{Ofil,\mathcal D^{\infty},\infty}(T^{an}/(Y^{an}\times\tilde S^{an}_I))$ commutes 
\begin{equation*}
\xymatrix{Rg^{*mod[-],\Gamma}\mathcal F_{S,an}^{GM}(L\mathbb D_SM)
\ar[d]^{T(g,\mathcal F_{an}^{GM})(L\mathbb D_SM)}
\ar[rr]^{Rg^{*mod[-],\Gamma}T(\mathcal F_{S,an}^{GM},\mathcal F_{S,an}^{FDR})(M)} & \, &
Rg^{*mod[-],\Gamma}\mathcal F_{S,an}^{FDR}(M)\ar[rrd] & \, &
g^{*mod}_{Hdg}\mathcal F_{S,an}^{FDR}(M)
\ar[ll]_{T(g^{*mod}_{Hdg},Rg^{*mod[-],\Gamma})(\mathcal F_{S,an}^{FDR}(M))}
\ar[d]^{T^!(g,\mathcal F_{an}^{FDR})(M)^{-1}} \\ 
\mathcal F_{T,an}^{GM}(g^*L\mathbb D_SM=L\mathbb D_Tg^!M)
\ar[rrrr]^{T(\mathcal F_{T,an}^{GM},\mathcal F_{T,an}^{FDR})(g^*M)} & \, & \, & \, &
\mathcal F_{T,an}^{FDR}(g^!M)}
\end{equation*}
\end{prop}

\begin{proof}
Similar to the proof of proposition \ref{TgDRGMdia}.
\end{proof}

We have the following easy proposition

\begin{prop}
Let $S\in\Var(\mathbb C)$ and $S=\cup_{i=1}^l S_i$ an open affine covering and denote, 
for $I\subset\left[1,\cdots l\right]$, $S_I=\cap_{i\in I} S_i$ and $j_I:S_I\hookrightarrow S$ the open embedding.
Let $i_i:S_i\hookrightarrow\tilde S_i$ closed embeddings, with $\tilde S_i\in\SmVar(\mathbb C)$. 
We have, for $M,N\in\DA(S)$ and $F,G\in C(\Var(\mathbb C)^{sm}/S)$ such that 
$M=D(\mathbb A^1,et)(F)$ and $N=D(\mathbb A^1,et)(G)$, 
the following commutative diagram in $D_{O_Sfil,\mathcal D^{\infty},\infty}(S^{an}/(\tilde S^{an}_I))$
\begin{equation*}
\xymatrix{\mathcal F_{S,an}^{GM}(L\mathbb D_SM)\otimes^L_{O_S}\mathcal F_{S,an}^{GM}(L\mathbb D_SN)
\ar[d]^{T(\mathcal F_{S,an}^{GM},\otimes)(L\mathbb D_SM,L\mathbb D_SN)}
\ar[rrrr]^{T(\mathcal F_{S,an}^{GM},\mathcal F_{S,an}^{FDR})(M)\otimes T(\mathcal F_{S,an}^{GM},\mathcal F_{S,an}^{FDR})(N)} 
& \, & \, & \, & 
\mathcal F_{S,an}^{FDR}(M)\otimes^{Hdg}_{O_S}\mathcal F_{S,an}^{FDR}(N)\ar[d]^{T(\mathcal F_{S,an}^{FDR},\otimes)(M,N)} \\
\mathcal F_{S,an}^{GM}(\mathbb D_SL(M\otimes N))\ar[rrrr]^{T(\mathcal F_{S,an}^{GM},\mathcal F_{S,an}^{FDR})(M\otimes N)} 
& \, &  \, & \, & \mathcal F_{S,an}^{FDR}(M\otimes N)}
\end{equation*}
\end{prop}

\begin{proof}
Immediate from definition.
\end{proof}

\subsection{The transformation map between the analytic De Rahm functor and the analytification of the algebraic De Rahm functor}

\subsubsection{The transformation map between the analytic Gauss Manin realization functor 
and the analytification of the algebraic Gauss Manin realization functor}

Recall from section 2 that, for $f:T\to S$ a morphism with $T,S\in\Var(\mathbb C)$,
we have the following commutative diagram of sites (\ref{AnVar})
\begin{equation*}
\xymatrix{
\AnSp(\mathbb C)/T^{an}\ar[rr]^{\An_T}\ar[dd]^{P(f)}\ar[rd]^{\rho_T} & \, & 
\Var(\mathbb C)/T\ar[dd]^{P(f)}\ar[rd]^{\rho_T} & \, \\  
 \, & \AnSp(\mathbb C)^{sm}/T^{an}\ar[rr]^{\An_T}\ar[dd]^{P(f)} & \, & \Var(\mathbb C)^{sm}/T\ar[dd]^{P(f)} \\  
\AnSp(\mathbb C)/S^{an}\ar[rr]^{\An_S}\ar[rd]^{\rho_S} & \, & \Var(\mathbb C)/S\ar[rd]^{\rho_S} & \, \\  
\, &  \AnSp(\mathbb C)^{sm}/S^{an}\ar[rr]^{\An_S} & \, & \Var(\mathbb C)^{sm}/S}.  
\end{equation*}
This gives for $s:\mathcal I\to\mathcal J$ a functor with $\mathcal I,\mathcal J\in\Cat$ and
$f:T_{\bullet}\to S_{s(\bullet)}$ a morphism of diagram of algebraic varieties with 
$T_{\bullet}\in\Fun(\mathcal I,\Var(\mathbb C))$, $S_{\bullet}\in\Fun(\mathcal J,\Var(\mathbb C))$ 
the commutative diagram of sites (\ref{AnVarIJ})
\begin{equation*}
Dia(S):=\xymatrix{
\AnSp(\mathbb C)/T_{\bullet}^{an}\ar[rr]^{\An_{T_{\bullet}}}\ar[dd]^{P(f)}\ar[rd]^{\rho_{T_{\bullet}}} & \, & 
\Var(\mathbb C)/T_{\bullet}\ar[dd]^{P(f_{\bullet})}\ar[rd]^{\rho_{T_{\bullet}}} & \, \\  
 \, & \AnSp(\mathbb C)^{sm}/T_{\bullet}^{an}\ar[rr]^{\An_{T_{\bullet}}}\ar[dd]^{P(f_{\bullet})} & \, & 
\Var(\mathbb C)^{sm}/T_{\bullet}\ar[dd]^{P(f_{\bullet})} \\  
\AnSp(\mathbb C)/S_{\bullet}^{an}\ar[rr]^{\An_{S_{\bullet}}}\ar[rd]^{\rho_{S_{\bullet}}} & \, & 
\Var(\mathbb C)/S_{\bullet}\ar[rd]^{\rho_{S_{\bullet}}} & \, \\  
\, &  \AnSp(\mathbb C)^{sm}/S_{\bullet}^{an}\ar[rr]^{\An_{S_{\bullet}}} & \, & 
\Var(\mathbb C)^{sm}/S_{\bullet}}.  
\end{equation*}

We have the following canonical transformation map given by the pullback of (relative) differential forms:
Let $S\in\Var(\mathbb C)$.
Consider the following commutative diagram in $\RCat$ :
\begin{equation*}
D(an,e):\xymatrix{
(\AnSp(\mathbb C)^{sm}/S^{an},O_{\AnSp(\mathbb C)^{sm}/T})\ar[rr]^{\An_S}\ar[d]^{e(T)} & \, & 
(\Var(\mathbb C)^{sm}/S,O_{\Var(\mathbb C)^{sm}/S})\ar[d]^{e(S)} \\
(S^{an},O_{S^{an}})\ar[rr]^{\an_S} & \, & (S,O_S)}
\end{equation*}
It gives (see section 2) the canonical morphism in $C_{\an_S^*O_S}(\AnSp(\mathbb C)^{sm}/S^{an})$ 
\begin{eqnarray*}
\Omega_{/(S^{an}/S)}:=\Omega_{(O_{\AnSp(\mathbb C)^{sm}/S^{an}}/\An_S^*O_{\Var(\mathbb C)^{sm}/S})/(O_{S^{an}}/\an_S^*O_S}): \\
\An_S^*(\Omega^{\bullet}_{/S},F_b)=(\Omega^{\bullet}_{\An_S^*O_{\Var(\mathbb C)^{sm}/S}/\An_S^*e(S)^*O_S},F_b)\to
(\Omega^{\bullet}_{/S^{an}},F_b)=(\Omega^{\bullet}_{O_{\AnSp(\mathbb C)^{sm}/S^{an}}/e(S^{an})^*O_{S^{an}}},F_b)
\end{eqnarray*}
which is by definition given by the analytification on differential forms : for $(V/S^{an})=(V,h)\in\AnSp(\mathbb C)^{sm}/S^{an}$,
\begin{eqnarray*}
\Omega_{/(S^{an}/S)}(V/S^{an}):\hat\omega\in\An_S^*(\Omega^{r}_{/S})(V/S^{an}):=
\lim_{(h':U\to S \mbox{sm},g':V\to U^{an},h,g)}\Omega^r_{/S}(U/S) \\
\mapsto\Omega_{(V/U)/(S^{an}/S)}(V/S^{an})(\omega):=\widehat{\an_S^*(\omega)}\in\Omega^r_{S^{an}}(V/S^{an}) ;
\end{eqnarray*}
with $\omega\in\Gamma(U,\Omega^r_U)$ is such that $q(\omega)=\hat\omega$. 
If $S\in\SmVar(\mathbb C)$, the map $\Omega_{/(T/S)}:\An_S^*\Omega^{\bullet}_{/S}\to\Omega^{\bullet}_{/S^{an}}$
is a map in $C_{O_Sfil,\mathcal D}(\AnSp(\mathbb C)^{sm}/S^{an})$.
It induces the canonical morphism in $C_{O_Sfil,\mathcal D}(\AnSp(\mathbb C)^{sm}/S^{an})$:
\begin{eqnarray*}
E\Omega_{/(S^{an}/S)}:\An_S^*E_{et}(\Omega^{\bullet}_{/S},F_b)\xrightarrow{T(\An_S,E)(\Omega^{\bullet}_{/S},F_b)} 
E_{et}(\An_S^*(\Omega^{\bullet}_{/S},F_b))\xrightarrow{E(\Omega_{/(S^{an}/S)})} E_{et}(\Omega^{\bullet}_{/S^{an}},F_b)
\end{eqnarray*}

We have the following canonical transformation map given by the analytical functor:
\begin{defi}\label{TomegaAnGM}
Let $S\in\SmVar(\mathbb C)$. 
\begin{itemize}
\item[(i)]For $F\in C(\Var(\mathbb C)^{sm}/S)$, 
we have the canonical transformation map in $C_{Ofil,\mathcal D}(S^{an})$
\begin{eqnarray*}
T(an,\Omega_{/\cdot})(F): \\
((e(S)_*\mathcal Hom^{\bullet}(F,E_{et}(\Omega^{\bullet}_{/S},F)))^{an}):= 
O_{S^{an}}\otimes_{\an_S^*O_S}\an_S^*(e(S)_*\mathcal Hom^{\bullet}(F,E_{et}(\Omega^{\bullet}_{/S},F_b))) \\ 
\xrightarrow{T(an,e)(-)}
O_{S^{an}}\otimes_{\an_S^*O_S}(e(S^{an})_*\An_S^*\mathcal Hom^{\bullet}(F,E_{et}(\Omega^{\bullet}_{/S},F_b))) \\ 
\xrightarrow{T(\An,hom)(F,E_{et}(\Omega_{/S},F_b))}  
O_{S^{an}}\otimes_{\an_S^*O_S}(e(S^{an})_*\mathcal Hom^{\bullet}(\An_S^*F,\An_S^*E_{et}(\Omega^{\bullet}_{/S},F_b))) \\
\xrightarrow{\mathcal Hom(\An_S^*F,E\Omega_{/(S^{an}/S)}\otimes m)}
e(S^{an})_*\mathcal Hom^{\bullet}(\An_S^*F,E_{et}(\Omega^{\bullet}_{/S^{an}},F_b))
\end{eqnarray*}
\item[(ii)]We get from (i), for $F\in C(\Var(\mathbb C)^{sm}/S)$, 
the canonical transformation map in $\PSh_{\mathcal D^{\infty}}(S^{an})$
\begin{eqnarray*}
T^n(an,\Omega_{/\cdot})(F):
J_SH^n((e(S)_*\mathcal Hom^{\bullet}(F,E_{et}(\Omega^{\bullet}_{/S},F_b)))^{an}) \\
\xrightarrow{J_S(H^nT(an,\Omega_{/\cdot})(F))}
J_SH^n(e(S^{an})_*\mathcal Hom^{\bullet}(\An_S^*F,E_{et}(\Omega^{\bullet}_{/S^{an}},F_b))) \\
\xrightarrow{\mathcal J_S(-)}
H^ne(S^{an})_*\mathcal Hom^{\bullet}(\An_S^*F,E_{et}(\Omega^{\bullet}_{/S^{an}},F_b))
\end{eqnarray*}
\end{itemize}
\end{defi}

\begin{lem}\label{anTlGM}
Let $S\in\SmVar(\mathbb C)$. 
\begin{itemize}
\item[(i)]For $h:U\to S$ a smooth morphism with $U\in\SmVar(\mathbb C)$, the following diagram commutes
\begin{equation*}
\xymatrix{e(S)_*\mathcal Hom^{\bullet}(\mathbb Z(U/S),E_{et}(\Omega^{\bullet}_{/S},F_b))^{an}
\ar[rrr]^{T(\Omega_{/\cdot},an)(\mathbb Z(U/S))} & \, & \, &
e(S^{an})_*\mathcal Hom^{\bullet}(\mathbb Z(U^{an}/S^{an}),E_{et}(\Omega^{\bullet}_{/S^{an}},F_b)) \\
(h_*E_{zar}(\Omega_{U/S},F_b))^{an}\ar[u]^{k}\ar[rrr]^{T^O_{\omega}(an,h)} & \, & \, & 
h_{an*}E_{usu}(\Omega_{U^{an}/S^{an}},F_b)\ar[u]^{=}}.
\end{equation*}
\item[(ii)]For $h:U\to S$ a smooth morphism with $U\in\SmVar(\mathbb C)$, the following diagram commutes
\begin{equation*}
\xymatrix{J_SH^n((e(S)_*\mathcal Hom^{\bullet}(\mathbb Z(U/S),E_{et}(\Omega^{\bullet}_{/S},F)))^{an})
\ar[rrr]^{T^n(\Omega_{/\cdot},an)(\mathbb Z(U/S))} & \, & \, &
H^ne(S^{an})_*\mathcal Hom^{\bullet}(\mathbb Z(U^{an}/S^{an}),E_{et}(\Omega^{\bullet}_{/S^{an}},F_b)) \\
J_SH^n((h_*E_{zar}(\Omega_{U/S},F_b))^{an})\ar[u]^{k}\ar[rrr]^{\mathcal J_S(-)\circ J_ST^O_{\omega}(an,h)(O_U,F)} & \, & \, & 
H^nh_{an*}E_{usu}(\Omega_{U^{an}/S^{an}},F_b)\ar[u]^{=}}.
\end{equation*}
\end{itemize}
\end{lem}

\begin{proof}
Follows from Yoneda lemma.
\end{proof}

By definition of the algebraic an analytic De Rahm realization functor, we have a natural transformation between them :

\begin{defi}\label{TanDRdefGM}
Let $S\in\SmVar(\mathbb C)$.
Let $M\in\DA_c(S)$ and $Q\in C(\Var(\mathbb C)^{sm}/S)$ projectively cofibrant such that $M=D(\mathbb A^1_S,et)(Q)$. 
We have the canonical transformation in $D_{Ofil,\mathcal D}(S^{an})$
\begin{eqnarray*}
T(An,\mathcal F_{an}^{GM})(M):(\mathcal F^{S}_{GM}(M))^{an}:=
(e(S)_*\mathcal Hom^{\bullet}(Q,E_{et}(\Omega^{\bullet}_{/S},F_b)))^{an}[-d_S] \\
\xrightarrow{T(an,\Omega_{/\cdot})(Q)}
e(S)_*\mathcal Hom^{\bullet}(\An_S^*Q,E_{et}(\Omega^{\bullet}_{/S^{an}},F_b))[-d_S] \\
\xrightarrow{=}
e(S)_*\mathcal Hom^{\bullet}(\An_S^*Q,E_{usu}(\Omega^{\bullet}_{/S^{an}},F_b))[-d_S]
=:\mathcal F_{S,an}^{GM}(M) 
\end{eqnarray*} 
\end{defi}

We give now the definition in the non smooth case :
Let $S\in\Var(\mathbb C)$. 
Let $S=\cup_{i=1}^lS_i$ be an open cover such that there exist closed embeddings
$i_i:S_i\hookrightarrow\tilde S_i$ with $\tilde S_i\in\SmVar(\mathbb C)$.
For $I\subset J$, denote by $p_{IJ}:\tilde S_J\to\tilde S_I$ the projection.
Consider, for $I\subset J\subset\left[1,\ldots,l\right]$,
resp. for each $I\subset\left[1,\ldots,l\right]$, the following commutative diagrams in $\Var(\mathbb C)$
\begin{equation*}
D_{IJ}=\xymatrix{ S_I\ar[r]^{i_I} & \tilde S_I \\
S_J\ar[u]^{j_{IJ}}\ar[r]^{i_J} & \tilde S_J\ar[u]^{p_{IJ}}}.  
\end{equation*}

We then have the following lemma
\begin{lem}\label{TanlemGM}
The maps $T(an,\Omega_{\cdot})(L(i_{I*}j_I^*F))$ induce a morphism in $C_{Ofil,\mathcal D}(S/(\tilde S_I))$
\begin{eqnarray*}
(T(an,\Omega_{/\cdot})(L(i_{I*}j_I^*F))): 
(e(\tilde S_I)_*\mathcal Hom^{\bullet}(L(i_{I*}j_I^*F),E_{et}(\Omega^{\bullet}_{/\tilde S_I},F_b)))^{an}[-d_{\tilde S_I}],
(u^q_{IJ}(F))^{an}) \\ 
\to (e(\tilde S_I)_*\mathcal Hom(\An(\tilde S_I)^*L(i_{I*}j_I^*F),E_{et}(\Omega^{\bullet}_{/\tilde S_I},F_b))[-d_{\tilde S_I}],
u^q_{IJ}(F)) 
\end{eqnarray*}
\end{lem}

\begin{proof}
Obvious.
\end{proof}

\begin{defi}\label{TanDRdefsingGM}
Let $S\in\Var(\mathbb C)$. Let $S=\cup_{i=1}^lS_i$ be an open cover
such that there exist closed embeddings $i_i:S_i\hookrightarrow\tilde S_i$ with $\tilde S_i\in\SmVar(\mathbb C)$.  
For $I\subset\left[1,\ldots,l\right]$, denote $S_I=\cap_{i\in I} S_i$.
We have then closed embeddings $i_I:S_I\hookrightarrow\tilde S_I=\Pi_{i\in I}\tilde S_i$.
Let $M\in\DA_c(S)$ and $F\in C(\Var(\mathbb C)^{sm}/S)$ such that $M=D(\mathbb A^1_S,et)(F)$. 
We have, by lemma \ref{TanlemGM}, the canonical transformation in $D_{Ofil,\mathcal D,\infty}(S^{an})$
\begin{eqnarray*}
T(An,\mathcal F^{GM})(M):(\mathcal F_S^{GM}(M))^{an}:=
(e(\tilde S_I)_*\mathcal Hom^{\bullet}(L(i_{I*}j_I^*F),
E_{et}(\Omega^{\bullet}_{/\tilde S_I},F_b)))^{an}[-d_{\tilde S_I}],(u^q_{IJ}(F))^{an}) \\
\xrightarrow{T(an,\Omega_{/\cdot})(L(i_{I*}j_I^*F)))}
(e(\tilde S_I)_*\mathcal Hom(\An(\tilde S_I)^*L(i_{I*}j_I^*F),E_{et}(\Omega^{\bullet}_{/\tilde S_I},F))[-d_{\tilde S_I}],
u^q_{IJ}(F)) \\ 
\xrightarrow{=}
(e(\tilde S_I)_*\mathcal Hom(\An(\tilde S_I)^*L(i_{I*}j_I^*F),E_{usu}(\Omega^{\bullet}_{/\tilde S_I},F))[-d_{\tilde S_I}],
u^q_{IJ}(F))=:\mathcal F_{S,an}^{GM}(M)
\end{eqnarray*} 
\end{defi}

The following proposition says this transformation map 
between $\mathcal F^{S,an}$ and $(\mathcal F_S^{FDR})^{an}$ is functorial in $S\in\Var(\mathbb C)$,
hence define a commutative diagram of morphism of 2-functor :

\begin{prop}\label{2functdiaGM}
\begin{itemize}
\item[(i)] Let $g:T\to S$ a morphism with $T,S\in\Var(\mathbb C)$.
Assume there exist a factorization $g:T\xrightarrow{l} Y\times S\xrightarrow{p_S}S$
with, $Y\in\SmVar(\mathbb C)$, $l$ a closed embedding and $p_S$ the projection.
Let $S=\cup_{i=1}^lS_i$ be an open cover such that there exist closed embeddings 
$i_i:S_i\hookrightarrow\tilde S_i$ with $\tilde S_i\in\SmVar(\mathbb C)$.
We then have closed embedding $i_i\circ l:T_i\hookrightarrow Y\times\tilde S_i$
and $\tilde g_I:=p_{\tilde S_I}:Y\times\tilde S_I\to\tilde S_I$ is a lift of $g_I:=g_{|T_I}:T_I\to S_I$.  
Then, for $M\in\DA_c(S)$, the following diagram in $D_{Ofil,\mathcal D,\infty}(T^{an}/(Y^{an}\times\tilde S^{an}_I))$ commutes
\begin{equation*}
\xymatrix{
Rg^{*mod[-],\Gamma}(\mathcal F_S^{GM}(M))^{an}\ar[d]_{(T(g,\mathcal F^{GM})(M))^{an}}
\ar[rrr]^{Lg^{*mod[-],\Gamma}(T(An,\mathcal F_S^{FDR})(M))} & \, & \, &
Rg^{*mod[-],\Gamma}(\mathcal F_{S,an}^{GM}(M))\ar[d]^{T(g,\mathcal F_{an}^{GM})(M)} \\ 
(\mathcal F_T^{GM}(g^*M))^{an}\ar[rrr]^{(T(An,\mathcal F_T^{GM})(g^*M))} & \, & \, &
(\mathcal F_{T,an}^{GM}(g^*M))}
\end{equation*}
\item[(ii)] Let $S\in\Var(\mathbb C)$.
Let $S=\cup_{i=1}^lS_i$ be an open cover such that there exist closed embeddings 
$i_i:S_i\hookrightarrow\tilde S_i$ with $\tilde S_i\in\SmVar(\mathbb C)$.
Then, for $M,N\in\DA_c(S)$, the following diagram in $D_{Ofil,\mathcal D,\infty}(S^{an}/(\tilde S_I^{an}))$ commutes
\begin{equation*}
\xymatrix{(\mathcal F_S^{GM}(M)\otimes_{O_S}\mathcal F_S^{GM}(N))^{an} 
\ar[d]_{=}\ar[rr]^{(T(\otimes,\mathcal F^{GM})(M,N))^{an}} & \, & 
(\mathcal F_S^{GM}(M\otimes N))^{an}\ar[dd]^{T(An,\mathcal F_S^{GM})(M\otimes N)} \\
(\mathcal F_S^{GM}(M))^{an}\otimes O_{S^{an}}\mathcal F_S^{GM}(N)^{an}
\ar[d]_{T(An,\mathcal F_S^{GM})(M)\otimes_{O_{S^{an}}}T(An,\mathcal F_S^{GM})(M)} & \, & \, \\
\mathcal F_{S,an}^{GM}(M)\otimes_{O_{S^{an}}}\mathcal F_{S,an}^{GM}(N)
\ar[rr]^{(T(\otimes,\mathcal F_{an}^{GM})(M,N))} & \, &  
\mathcal F_{S,an}^{GM}(M\otimes N)}
\end{equation*}
\end{itemize}
\end{prop}

\begin{proof}
Immediate from definition.
\end{proof}

\begin{prop}\label{mainthmpropanTGM}
Let $f:X\to S$ a morphism with $S,X\in\Var(\mathbb C)$. Assume there exist a factorization
\begin{equation*}
f:X\xrightarrow{l}Y\times S\xrightarrow{p}S
\end{equation*}
with $Y\in\SmVar(\mathbb C)$, $l$ a closed embedding and $p$ the projection.
Let $S=\cup_iS_i$ an open cover such that there exists closed embeddings
$i_i:S_i\hookrightarrow\tilde S_i$ with $\tilde S_i\in\SmVar(\mathbb C)$.
\begin{itemize}
\item[(i)]We have then the following commutative diagram in $D_{Ofil,\mathcal D,\infty}(S^{an}/(\tilde S_I^{an}))$,
\begin{equation*}
\xymatrix{(\mathcal F_S^{GM}(M(X/S)))^{an}
\ar[rrr]^{T(\An,\mathcal F_{an}^{FDR})(M(X/S))}\ar[d]_{I^{GM}(X/S)^{an}} & \, & \, & 
\mathcal F_{S,an}^{GM}(M(X/S))\ar[d]^{I^{GM}(X/S)} \\
((p_{\tilde S_I*}\Gamma_{X_I}E_{zar}(\Omega^{\bullet}_{Y\times\tilde S_I/\tilde S_I}))^{an}[-d_{\tilde S_I}],w_{IJ}(X/S)^{an})
\ar[rrr]^{(T^O_{\omega}(an,p_{\tilde S_I})^{\gamma})}\ar[d]_{((p_*T^O_{\omega}(\otimes,\gamma)(-))^{an})} & \ ,& \, & 
(p_{\tilde S_I*}\Gamma_XE_{usu}(\Omega^{\bullet}_{(Y\times \tilde S_I)^{an}/\tilde S^{an}_I})[-d_{\tilde S_I}],w_{IJ}(X^{an}/S^{an}))
\ar[d]^{(p_*T^O_{\omega}(\otimes,\gamma)(-))}\\
(\int^{FDR}_{f}(\Gamma_{X_I}E_{zar}(O_{Y\times\tilde S_I},F_b)[-d_Y-d_{\tilde S_I}],x_{IJ}(X/S)))^{an}
\ar[rrr]^{(T(\an,\gamma)(-))\circ T^{\mathcal Dmod}(an,f)(-)} & \ ,& \, & 
\int^{FDR}_{f}(\Gamma_{X_I}E_{usu}(O_{(Y\times\tilde S_I)^{an}},F_b)[-d_Y-d_{\tilde S_I}],x_{IJ}(X/S))}.
\end{equation*} 
\item[(ii)]We have then the following commutative diagram in $\PSh_{\mathcal D^{\infty}fil}(S^{an}/(\tilde S_I^{an}))$,
\begin{equation*}
\xymatrix{J_SH^n(\mathcal F_S^{GM}(M^{BM}(X/S)))^{an}
\ar[rrr]^{\mathcal J_S(-)\circ H^nT(\An,\mathcal F_{an}^{FDR})(M(X/S))}\ar[d]_{H^n(I^{GM}(X/S)^{an})} & \, & \, & 
H^n\mathcal F_{S,an}^{GM}(M(X/S))\ar[d]^{H^nI^{GM}(X/S)} \\
J_SH^n((p_{\tilde S_I*}\Gamma_{X_I}E_{zar}(\Omega^{\bullet}_{Y\times\tilde S_I/\tilde S_I}))^{an}[-d_{\tilde S_I}],w_{IJ}(X/S)^{an})
\ar[rrr]^{H^n(T^O_{\omega}(an,p_{\tilde S_I})^{\gamma}))}\ar[d]_{H^n((p_*T^O_{\omega}(\otimes,\gamma)(-))^{an})} & \ ,& \, & 
H^n(p_{\tilde S_I*}\Gamma_XE_{usu}(\Omega^{\bullet}_{(Y\times \tilde S_I)^{an}/\tilde S^{an}_I})[-d_{\tilde S_I}],w_{IJ}(X/S))
\ar[d]^{H^n(p_*T^O_{\omega}(\otimes,\gamma)(-))}\\
J_SH^n(\int^{FDR}_{f}(\Gamma_{X_I}E_{zar}(O_{Y\times\tilde S_I},F_b)[-d_Y-d_{\tilde S_I}],x_{IJ}(X/S)))^{an}
\ar[rrr]^{H^n(T(\an,\gamma)(-))\circ H^nT^{\mathcal Dmod}(an,f)(-)} & \ ,& \, & 
\int^{FDR}_{f}(\Gamma_{X_I}E_{usu}(O_{(Y\times\tilde S_I)^{an}},F_b)[-d_Y-d_{\tilde S_I}],x_{IJ}(X^{an}/S^{an}))}.
\end{equation*} 
\end{itemize}
\end{prop}

\begin{proof}
\noindent(i):Immediate from definition.

\noindent(ii):Follows from (i).
\end{proof}

We deduce from proposition \ref{mainthmpropanTGM} and theorem \ref{GAGADmod} (GAGA for D-modules) the following :

\begin{thm}\label{mainthmGAGAGM}
\begin{itemize}
\item[(i)] Let $S\in\Var(\mathbb C)$. Then, for $M\in\DA_c(S)$
\begin{eqnarray*}
\mathcal J_S(-)\circ H^nT(\An,\mathcal F_{an}^{GM})(M):
J_S(H^n(\mathcal F_S^{GM}(M))^{an})\xrightarrow{\sim}H^n\mathcal F_{S,an}^{GM}(M)
\end{eqnarray*}
is an isomorphism in $\PSh_{\mathcal D}(S^{an}/(\tilde S_I^{an}))$.
\item[(ii)] A relative version of Grothendieck GAGA theorem for De Rham cohomology
Let $h:U\to S$ a smooth morphism with $S,U\in\SmVar(\mathbb C)$. Then,
\begin{equation*}
\mathcal J_S(-)\circ J_ST^O_{\omega}(an,h):
J_S((R^nh_*\Omega^{\bullet}_{U/S})^{an})\xrightarrow{\sim} R^nh_*\Omega^{\bullet}_{U^{an}/S^{an}}
\end{equation*}
is an isomorphism in $\PSh_{\mathcal D}(S^{an})$.
\end{itemize}
\end{thm}

\begin{proof}
\noindent(i):Follows from proposition \ref{mainthmpropanTGM}(ii) and theorem \ref{GAGADmod} using
a resolution by Corti-Hanamura motives.

\noindent(ii):Follows from (i) and lemma \ref{anTlGM}(ii).
\end{proof}

\subsubsection{The transformation map between the analytic filtered De Rham realization functor 
and the analytification of the filtered algebraic De Rham realization functor}

Recall from section 2 that, for $S\in\Var(\mathbb C)$ we have the following commutative diagrams of sites
\begin{equation*}
\xymatrix{\AnSp(\mathbb C)^2/S\ar[rr]^{\mu_S}\ar[dd]_{\An_S}\ar[rd]^{\rho_S} & \, & 
\AnSp(\mathbb C)^{2,pr}/S\ar[dd]^{\An_S}\ar[rd]^{\rho_S} & \, \\
\, & \AnSp(\mathbb C)^{2,sm}/S\ar[rr]^{\mu_S}\ar[dd]_{\An_S} & \, & \AnSp(\mathbb C)^{2,smpr}/S\ar[dd]^{\An_S} \\
\Var(\mathbb C)^2/S\ar[rr]^{\mu_S}\ar[rd]^{\rho_S} & \, &  \Var(\mathbb C)^{2,smpr}/S\ar[rd]^{\rho_S} & \, \\
\, & \Var(\mathbb C)^2/S\ar[rr]^{\mu_S} & \, &  \Var(\mathbb C)^{2,smpr}/S & \,}
\end{equation*}
and
\begin{equation}
\xymatrix{\AnSp(\mathbb C)^{2,pr}/S\ar[rr]^{\Gr_S^{12}}\ar[dd]_{\An_S}\ar[rd]^{\rho_S} & \, & 
\AnSp(\mathbb C)/S\ar[dd]^{\An_S}\ar[rd]^{\rho_S} & \, \\
\, & \AnSp(\mathbb C)^{2,smpr}/S\ar[rr]^{\Gr_S^{12}}\ar[dd]_{\An_S} & \, & \AnSp(\mathbb C)^{sm}/S\ar[dd]^{\An_S} \\
\Var(\mathbb C)^{2,pr}/S\ar[rr]^{\Gr_S^{12}}\ar[rd]^{\rho_S} & \, & \Var(\mathbb C)/S\ar[rd]^{\rho_S} & \, \\
\, & \Var(\mathbb C)^{2,sm}/S\ar[rr]^{\Gr_S^{12}} & \, & \Var(\mathbb C)^{sm}/S},
\end{equation}
and that for $f:T\to S$ a morphism with $T,S\in\Var(\mathbb C)$ we have the commutative diagram of site (\ref{AnVar12})
\begin{equation*}
\xymatrix{
\AnSp(\mathbb C)^2/T^{an}\ar[rr]^{\An_T}\ar[dd]^{P(f)}\ar[rd]^{\rho_T} & \, & 
\Var(\mathbb C)^2/T\ar[dd]^{P(f)}\ar[rd]^{\rho_T} & \, \\  
 \, & \AnSp(\mathbb C)^{2,sm}/T^{an}\ar[rr]^{\An_T}\ar[dd]^{P(f)} & \, & \Var(\mathbb C)^{2,sm}/T\ar[dd]^{P(f)} \\  
\AnSp(\mathbb C)^2/S^{an}\ar[rr]^{\An_S}\ar[rd]^{\rho_S} & \, & \Var(\mathbb C)^2/S\ar[rd]^{\rho_S} & \, \\  
\, &  \AnSp(\mathbb C)^{2,sm}/S^{an}\ar[rr]^{\An_S} & \, & \Var(\mathbb C)^{2,sm}/S}.  
\end{equation*}
This gives for $s:\mathcal I\to\mathcal J$ a functor with $\mathcal I,\mathcal J\in\Cat$ and
$f:T_{\bullet}\to S_{s(\bullet)}$ a morphism of diagram of algebraic varieties with 
$T_{\bullet}\in\Fun(\mathcal I,\Var(\mathbb C))$, $S_{\bullet}\in\Fun(\mathcal J,\Var(\mathbb C))$ 
the commutative diagram of sites (\ref{AnVar12IJ})
\begin{equation*}
Dia^{12}(S):=\xymatrix{
\AnSp(\mathbb C)^2/T_{\bullet}^{an}\ar[rr]^{\An_{T_{\bullet}}}\ar[dd]^{P(f_{\bullet})}\ar[rd]^{\rho_{T_{\bullet}}} & \, & 
\Var(\mathbb C)^2/T_{\bullet}\ar[dd]^{P(f_{\bullet})}\ar[rd]^{\rho_{T_{\bullet}}} & \, \\  
 \, & \AnSp(\mathbb C)^{2,sm}/T_{\bullet}^{an}\ar[rr]^{\An_{T_{\bullet}}}\ar[dd]^{P(f_{\bullet})} & \, & 
\Var(\mathbb C)^{2,sm}/T_{\bullet}\ar[dd]^{P(f_{\bullet})} \\  
\AnSp(\mathbb C)^2/S_{\bullet}^{an}\ar[rr]^{\An_{S_{\bullet}}}\ar[rd]^{\rho_{S_{\bullet}}} & \, & 
\Var(\mathbb C)^2/S_{\bullet}\ar[rd]^{\rho_{S_{\bullet}}} & \, \\  
\, &  \AnSp(\mathbb C)^{2,sm}/S_{\bullet}^{an}\ar[rr]^{\An_{S_{\bullet}}} & \, & 
\Var(\mathbb C)^{2,sm}/S_{\bullet}}.  
\end{equation*}

Let $S\in\SmVar(\mathbb C)$. We have the canonical map in $C(\Var(\mathbb C)^{2,smpr}/S)$ 
\begin{eqnarray*}
\Omega^{\Gamma,pr}_{/(S^{an}/S)}:(\Omega^{\bullet,\Gamma,pr}_{/S},F_{DR})\to
(\Omega^{\bullet,\Gamma,pr}_{/S^{an}},F_{DR})
\end{eqnarray*}
given by for $p:(Y\times S,Z)\to S$ the projection with $Y\in\SmVar(\mathbb C)$,
\begin{eqnarray*}
\Omega^{\Gamma,pr}_{/(S^{an}/S)}((Y\times S,Z)/S):
\Omega^{\bullet}_{Y\times S/S}\otimes_{O_{Y\times S}}\Gamma_Z^{\vee,Hdg}(O_{Y\times S},F_b) \\
\xrightarrow{\Omega_{((Y\times S)^{an}/Y\times S)/(S^{an}/S)}(-)} 
\Omega^{\bullet}_{(Y\times S)^{an}/S^{an}}\otimes_{O_{(Y\times S)^{an}}}(\Gamma_Z^{\vee,Hdg}(O_{Y\times S},F_b))^{an} 
\end{eqnarray*}

We have the following canonical transformation map given by the analytical functor:
\begin{defi}\label{TomegaAn}
Let $S\in\SmVar(\mathbb C)$. For $F\in C(\Var(\mathbb C)^{2,smpr}/S)$, 
we have the canonical transformation map in $C_{\mathcal D^{\infty}fil}(S^{an})$
\begin{eqnarray*}
T(an,\Omega^{\Gamma,pr}_{/S})(F): \\
(e(S)_*\mathcal Hom^{\bullet}(F,E_{et}(\Omega^{\bullet,\Gamma,pr}_{/S},F_{DR})))^{an}:= 
O_{S^{an}}\otimes_{\an_S^*O_S}\an_S^*(e(S)_*\mathcal Hom^{\bullet}(F,E_{et}(\Omega^{\bullet,\Gamma,pr}_{/S},F_{DR}))) \\ 
\xrightarrow{T(an,e)(-)}
O_{S^{an}}\otimes_{\an_S^*O_S}(e(S^{an})_*\An_S^*\mathcal Hom^{\bullet}(F,
E_{et}(\Omega^{\bullet,\Gamma,pr}_{/S},F_{DR}))) \\ 
\xrightarrow{T(\An,hom)(-,-)}  
O_{S^{an}}\otimes_{\an_S^*O_S}(e(S^{an})_*\mathcal Hom^{\bullet}(\An_S^*F,
\An_S^*E_{et}(\Omega^{\bullet,\Gamma,pr}_{/S},F_{DR}))) \\
\xrightarrow{\mathcal Hom(\An_S^*F,\An_S^*E_{et}(\Omega^{\Gamma,pr}_{/(S^{an}/S)})\otimes m)}
e(S^{an})_*\mathcal Hom^{\bullet}(\An_S^*F,E_{et}(\Omega^{\bullet,\Gamma,pr,an}_{/S^{an}},F_{DR}))
\end{eqnarray*}
\end{defi}

By definition of the algebraic an analytic De Rahm realization functor, we have a natural transformation between them :

\begin{defi}\label{TanDRdef}
Let $S\in\SmVar(\mathbb C)$.
Let $M\in\DA_c(S)$ and $(F,W)\in C_{fil}(\Var(\mathbb C)^{sm}/S)$ such that $(M,W)=D(\mathbb A^1_S,et)(F,W)$. 
We have the canonical transformation map $\pi_S(D(MHM(S)))\subset D_{\mathcal D(1,0)fil,\infty}(S^{an})$
\begin{eqnarray*}
T(An,\mathcal F_{an}^{FDR})(M):(\mathcal F^{S}_{FDR}(M))^{an}:=
(e(S)_*\mathcal Hom^{\bullet}(\hat R^{CH}(\rho_S^*L(F,W)), E_{zar}(\Omega^{\bullet,\Gamma,pr}_{/S},F_{DR})))^{an} \\
\xrightarrow{T(an,\Omega^{\Gamma,pr}_{/S})(-)}
e(S)_*\mathcal Hom^{\bullet}(\An_S^*\hat R^{CH}(\rho_S^*L(F,W)), E_{usu}(\Omega^{\bullet,\Gamma,pr,an}_{/S^{an}},F_{DR})) 
=:\mathcal F_{S,an}^{FDR}(M) 
\end{eqnarray*} 
\end{defi}

We give now the definition in the non smooth case :
Let $S\in\Var(\mathbb C)$. 
Let $S=\cup_{i=1}^lS_i$ be an open cover such that there exist closed embeddings
$i_i:S_i\hookrightarrow\tilde S_i$ with $\tilde S_i\in\SmVar(\mathbb C)$.
For $I\subset J$, denote by $p_{IJ}:\tilde S_J\to\tilde S_I$ the projection.
Consider, for $I\subset J\subset\left[1,\ldots,l\right]$,
resp. for each $I\subset\left[1,\ldots,l\right]$, the following commutative diagrams in $\Var(\mathbb C)$
\begin{equation*}
D_{IJ}=\xymatrix{ S_I\ar[r]^{i_I} & \tilde S_I \\
S_J\ar[u]^{j_{IJ}}\ar[r]^{i_J} & \tilde S_J\ar[u]^{p_{IJ}}}.  
\end{equation*}

We then have the following lemma
\begin{lem}\label{Tanlem}
The maps $T(an,\Omega_{\cdot})(-)$ induce a morphism in $C_{\mathcal D^{\infty}(fil)}(S/(\tilde S_I))$
\begin{eqnarray*}
(T(an,\Omega^{\Gamma,pr}_{/\tilde S_I})(\hat R^{CH}(\rho_{\tilde S_I}^*L(i_{I*}j_I^*F))): \\
((e'(\tilde S_I)_*\mathcal Hom^{\bullet}(\hat R^{CH}(\rho_{\tilde S_I}^*L(i_{I*}j_I^*F)),
E_{zar}(\Omega^{\bullet,\Gamma,pr}_{/\tilde S_I},F_{DR})))^{an},(u^q_{IJ}(F))^{an}) \\ 
\to (e'(\tilde S_I)_*\mathcal Hom^{\bullet}(\An(\tilde S_I)^*\hat R^{CH}(\rho_{\tilde S_I}^*L(i_{I*}j_I^*F)),
E_{usu}(\Omega^{\bullet,\Gamma,pr,an}_{/\tilde S_I},F_{DR})),u^q_{IJ}(F)) 
\end{eqnarray*}
\end{lem}

\begin{proof}
Obvious.
\end{proof}

\begin{defi}\label{TanDRdefsing}
Let $S\in\Var(\mathbb C)$. Let $S=\cup_{i=1}^lS_i$ be an open cover
such that there exist closed embeddings $i_i:S_i\hookrightarrow\tilde S_i$ with $\tilde S_i\in\SmVar(\mathbb C)$.  
For $I\subset\left[1,\ldots,l\right]$, denote $S_I=\cap_{i\in I} S_i$.
We have then closed embeddings $i_I:S_I\hookrightarrow\tilde S_I=\Pi_{i\in I}\tilde S_i$.
Let $M\in\DA_c(S)$ and $(F,W)\in C_{fil}(\Var(\mathbb C)^{sm}/S)$ such that $(M,W)=D(\mathbb A^1_S,et)(F,W)$. 
We have, by lemma \ref{Tanlem}, the canonical transformation map in 
$\pi_S(D(MHM(S)))\subset D_{\mathcal D(1,0)fil,\infty}(S^{an})$
\begin{eqnarray*}
T(An,\mathcal F_S^{FDR})(M):(\mathcal F_S^{FDR}(M))^{an}\xrightarrow{:=} \\
((e'(\tilde S_I)_*\mathcal Hom^{\bullet}(\hat R^{CH}(\rho_{\tilde S_I}^*L(i_{I*}j_I^*(F,W))), 
E_{zar}(\Omega^{\bullet,\Gamma,pr}_{/\tilde S_I},F_{DR})))^{an},(u^q_{IJ}(F,W))^{an}) \\ 
\xrightarrow{(T(an,\Omega^{\Gamma,pr}_{/\tilde S_I})(\hat R^{CH}(\rho_{\tilde S_I}^*L(i_{I*}j_I^*(F,W)))))} \\
(e'(\tilde S_I)_*\mathcal Hom^{\bullet}(\An(\tilde S_I)^*\hat R^{CH}(\rho_{\tilde S_I}^*L(i_{I*}j_I^*(F,W))), 
E_{usu}(\Omega^{\bullet,\Gamma,pr,an}_{/\tilde S_I},F_{DR})),u^q_{IJ}(F,W)) \\ 
\xrightarrow{=:}\mathcal F_{S,an}^{FDR}(M)
\end{eqnarray*} 
\end{defi}

The following proposition says this transformation map 
between $\mathcal F^{S,an}$ and $(\mathcal F_S^{FDR})^{an}$ is functorial in $S\in\Var(\mathbb C)$,
hence define a commutative diagram of morphism of 2-functor :

\begin{prop}\label{2functdia}
\begin{itemize}
\item[(i)]Let $g:T\to S$ a morphism with $T,S\in\Var(\mathbb C)$.
Assume there exist a factorization $g:T\xrightarrow{l} Y\times S\xrightarrow{p_S}S$
with, $Y\in\SmVar(\mathbb C)$, $l$ a closed embedding and $p_S$ the projection.
Let $S=\cup_{i=1}^lS_i$ be an open cover such that there exist closed embeddings 
$i_i:S_i\hookrightarrow\tilde S_i$ with $\tilde S_i\in\SmVar(\mathbb C)$.
We then have closed embedding $i_i\circ l:T_i\hookrightarrow Y\times\tilde S_i$
and $\tilde g_I:=p_{\tilde S_I}:Y\times\tilde S_I\to\tilde S_I$ is a lift of $g_I:=g_{|T_I}:T_I\to S_I$. 
\begin{itemize}
\item[(i0)] Then, for $M\in\DA_c(S)$, the following diagram in 
$\pi_T(D(MHM(T^{an}))\subset D_{\mathcal D(1,0)fil}(T^{an}/(Y^{an}\times\tilde S^{an}_I))$,
see definition \ref{TgDRdefsing} and definition \ref{TgDRdefsingan} commutes
\begin{equation*}
\begin{tikzcd}
g^{\hat*mod}_{Hdg}((\mathcal F_S^{FDR}(M))^{an})=(g^{\hat*mod}_{Hdg}(\mathcal F_S^{FDR}(M)))^{an}
\ar[rr,"g_{Hdg}^{\hat*mod}(T(An{,}\mathcal F_S^{FDR})(M))"]\ar[d,"'(T(g{,}\mathcal F^{FDR})(M))^{an}"'] & \, &
g^{\hat*mod}_{Hdg}(\mathcal F_{S,an}^{FDR}(M))\ar[d,"T(g{,}\mathcal F_{an}^{FDR})(M)"]  \\ 
(\mathcal F_T^{FDR}(g^*M))^{an}\ar[rr,"T(An{,}\mathcal F_T^{FDR})(g^*M)"] & \, &
\mathcal F_{T,an}^{FDR}(g^*M)
\end{tikzcd}
\end{equation*}
\item[(i1)] Then, for $M\in\DA_c(S)$, the following diagram in 
$\pi_T(D(MHM(T))\subset D_{\mathcal D(1,0)fil}(T^{an}/(Y^{an}\times\tilde S^{an}_I))$ commutes
\begin{equation*}
\begin{tikzcd}
g^{*mod}_{Hdg}((\mathcal F_S^{FDR}(M))^{an})=(g^{*mod}_{Hdg}(\mathcal F_S^{FDR}(M)))^{an}
\ar[rr,"g_{Hdg}^{\hat*mod}(T(An{,}\mathcal F_S^{FDR})(M))"] & \, &
g^{*mod}_{Hdg}(\mathcal F_{S,an}^{FDR}(M))  \\ 
(\mathcal F_T^{FDR}(g^*M))^{an}
\ar[rr,"T(An{,}\mathcal F_T^{FDR})(g^*M)"]\ar[u,"(T(g{,}\mathcal F^{FDR})(M))^{an}"] & \, &
\mathcal F_{T,an}^{FDR}(g^*M)\ar[u,"'T^!(g{,}\mathcal F_{an}^{FDR})(M)"']
\end{tikzcd}
\end{equation*}
\item[(i2)] Then, for $M\in\DA_c(T)$, the following diagram in 
$\pi_S(D(MHM(S^{an}))\subset D_{\mathcal D(1,0)fil}(S^{an}/(\tilde S^{an}_I))$ commutes
\begin{equation*}
\begin{tikzcd}
Rg^{Hdg}_*((\mathcal F_T^{FDR}(M))^{an})=(Rg^{Hdg}*(\mathcal F_T^{FDR}(M)))^{an}
\ar[rr,"Rg^{Hdg}_*(T(An{,}\mathcal F_S^{FDR})(M))"] & \, &
Rg^{Hdg}_*(\mathcal F_{S,an}^{FDR}(M))  \\ 
(\mathcal F_S^{FDR}(Rg_*M))^{an}
\ar[rr,"T(An{,}\mathcal F_S^{FDR})(Rg_*M)"]\ar[u,"(T_*(g{,}\mathcal F^{FDR})(M))^{an}"] & \, &
\mathcal F_{S,an}^{FDR}(Rg_*M)\ar[u,"'T_*(g{,}\mathcal F_{an}^{FDR})(M)"']
\end{tikzcd}
\end{equation*}
\item[(i3)] Then, for $M\in\DA_c(T)$, the following diagram in 
$\pi_S(D(MHM(S^{an}))\subset D_{\mathcal D(1,0)fil}(S^{an}/(\tilde S^{an}_I))$ commutes
\begin{equation*}
\begin{tikzcd}
Rg^{Hdg}_!((\mathcal F_T^{FDR}(M))^{an})=(Rg^{Hdg}!(\mathcal F_T^{FDR}(M)))^{an}
\ar[rr,"Rg^{Hdg}_*(T(An{,}\mathcal F_S^{FDR})(M))"]\ar[d,"'(T_!(g{,}\mathcal F^{FDR})(M))^{an}"'] & \, &
Rg^{Hdg}_!(\mathcal F_{S,an}^{FDR}(M))\ar[d,"T_!(g{,}\mathcal F_{an}^{FDR})(M)"]  \\ 
(\mathcal F_S^{FDR}(Rg_!M))^{an}\ar[rr,"T(An{,}\mathcal F_S^{FDR})(Rg_!M)"] & \, &
\mathcal F_{S,an}^{FDR}(Rg_!M)
\end{tikzcd}
\end{equation*}
\end{itemize}
\item[(ii)] Let $S\in\Var(\mathbb C)$.
Let $S=\cup_{i=1}^lS_i$ be an open cover such that there exist closed embeddings 
$i_i:S_i\hookrightarrow\tilde S_i$ with $\tilde S_i\in\SmVar(\mathbb C)$.
Then, for $M,N\in\DA_c(S)$, the following diagram in 
$\pi_S(D(MHM(S)))\subset D_{\mathcal D(1,0)fil}(S^{an}/(\tilde S_I^{an}))$ commutes
\begin{equation*}
\xymatrix{(\mathcal F_S^{FDR}(M)\otimes^{Hdg}_{O_S}\mathcal F_S^{FDR}(N))^{an} 
\ar[d]_{=}\ar[rr]^{(T(\otimes,\mathcal F^{FDR})(M,N))^{an}} & \, &
(\mathcal F_S^{FDR}(M\otimes N))^{an}\ar[dd]_{T(An,\mathcal F_S^{FDR})(M\otimes N)} \\
(\mathcal F_S^{FDR}(M))^{an})\otimes O_{S^{an}}(\mathcal F_S^{FDR}(N))^{an})
\ar[d]_{T(An,\mathcal F_S^{FDR})(M)\otimes^{Hdg}_{O_{S^{an}}}T(An,\mathcal F_S^{FDR})(M)} & \, & \, \\
\mathcal F_{S,an}^{FDR}(M)\otimes^{Hdg}_{O_{S^{an}}}\mathcal F_{S,an}^{FDR}(N)
\ar[rr]^{(T(\otimes,\mathcal F_{an}^{FDR})(M,N))} & \, &  
\mathcal F_{S,an}^{FDR}(M\otimes N)}
\end{equation*}
\end{itemize}
\end{prop}

\begin{proof}
Immediate from definition.
\end{proof}

\begin{prop}\label{mainthmpropanT}
Let $f:X\to S$ a morphism with $S,X\in\Var(\mathbb C)$. Assume there exist a factorization
\begin{equation*}
f:X\xrightarrow{l}Y\times S\xrightarrow{p}S
\end{equation*}
with $Y\in\SmVar(\mathbb C)$, $l$ a closed embedding and $p$ the projection.
Let $S=\cup_iS_i$ an open affine cover and $i_i:S_i\hookrightarrow\tilde S_i$ closed embeddings 
with $\tilde S_i\in\SmVar(\mathbb C)$. We have then the following commutative diagram 
in $D_{\mathcal D^{\infty}fil,\infty}(S^{an})$,
\begin{equation*}
\xymatrix{(\mathcal F_S^{FDR}(M(X/S)))^{an}
\ar[rr]^{T(\An,\mathcal F_S^{FDR})(M(X/S))}\ar[d]^{I(X/S)} & \, & 
\mathcal F_{S,an}^{FDR}(M(X/S))\ar[d]^{I(X/S)} \\
(Rf^{Hdg}_!(\Gamma_{X_I}^{\vee,Hdg}(O_{Y\times\tilde S_I},F_b)(d_Y)[2d_Y],x_{IJ}(X/S)))^{an}
\ar[rr]^{T^{\mathcal Dmod}(an,p_{\tilde S_I})(-)} & \, &  
Rf^{Hdg}_!((\Gamma^{\vee,Hdg}_{X_I}(O_{(Y\times\tilde S_I)^{an}},F_b)(d_Y)[2d_Y],x_{IJ}(X/S)))}.
\end{equation*} 
\end{prop}

\begin{proof}
Immediate from definition.
\end{proof}

We deduce from proposition \ref{mainthmpropanT} and theorem \ref{GAGADmod} (GAGA for D-modules) the following :

\begin{thm}\label{mainthmGAGA}
Let $S\in\Var(\mathbb C)$. For $M\in\DA_c(S)$, the map in 
$\pi_S(D(MHM(S^{an})))\subset D_{\mathcal D(1,0)fil,\infty}(S^{an})$
\begin{eqnarray*}
T(An,\mathcal F^{FDR})(M):(\mathcal F_S^{FDR}(M))^{an}\xrightarrow{\sim}\mathcal F_{S,an}^{FDR}(M)
\end{eqnarray*}
given in definition \ref{TanDRdefsing} is an isomorphism.
\end{thm}

\begin{proof}
Follows from proposition \ref{mainthmpropanT} and theorem \ref{GAGADmod}.
\end{proof}

We finish this subsection by the following easy proposition :

\begin{prop}\label{TanDRGMdia}
Let $S\in\Var(\mathbb C)$.
Let $S=\cup_{i=1}^lS_i$ be an open cover such that 
there exists closed embeddings $i_i:S_i\hookrightarrow\tilde S_i$ with $\tilde S_i\in\SmVar(\mathbb C)$
Let $M\in\DA_c(S)$ and $F\in C(\Var(\mathbb C)^{sm}/S)$ such that $M=D(\mathbb A^1_S,et)(F)$. 
Then the following diagram in $D_{Ofil,\mathcal D,\infty}(S/\tilde S_I)$ commutes 
\begin{equation*}
\xymatrix{  
\mathcal (F_S^{GM}(L\mathbb D_SM))^{an}
\ar[d]_{T(\An,\mathcal F_S^{GM})(\mathbb D_SM)}\ar[rrrr]^{\mathcal J\circ (T(\mathcal F_S^{GM},\mathcal F_S^{FDR})(M))^{an}}
& \, & \, & \, & J_S((\mathcal F_S^{FDR}(M))^{an})\ar[d]^{J_S(T(\An,\mathcal F_S^{FDR})(M))} \\ 
\mathcal F_{S,an}^{GM}(L\mathbb D_SM)\ar[rrrr]^{T(\mathcal F_{S,an}^{GM},\mathcal F_{S,an}^{FDR})(M)}
 & \, & \, & \, & J_S(\mathcal F_{S,an}^{FDR}(M))}
\end{equation*}
\end{prop}

\begin{proof}
Immediate from definition.
\end{proof}

\section{The Hodge realization functor for relative motives}

\subsection{The Betti realization functor}

We have two definition of the Betti realization functor which coincide at least for constructible motives,
one given by \cite{AyoubB} using the analytical functor and one given in \cite{B3} 
by composing the analytical functor with the forgetfull functor to the topological space of a complex analytic space
wich is a CW complex (see also \cite{LW} for the absolute case) .

\begin{defi}\label{bettidef}
Let $S\in\Var(\mathbb C)$.
\begin{itemize}
\item[(i)] The Ayoub's Betti realization functor is 
\begin{equation*}
\Bti_S^*:\DA(S)\to D(S^{an}) \; , \; 
M\in\DA(S)\mapsto\Bti_S^*M=Re(S^{an})_*\An_S^*M=e(S^{an})_*\underline{\sing}_{\mathbb D^*}\An_S^*F
\end{equation*}
where $F\in C(\Var(\mathbb C)^{sm}/S)$ is such that $M=D(\mathbb A^1,et)(F)$.
\item[(ii)] In \cite{B3}, we define the Betti realization functor as 
\begin{equation*}
\widetilde\Bti_S^*:\DA(S)\to D(S^{an})=D(S^{cw}) \; , \; 
M\mapsto\widetilde\Bti_S^*M=Re(S^{cw})_*\widetilde\Cw_S^*M=e(S^{cw})_*\underline{\sing}_{\mathbb I^*}\widetilde\Cw_S^*F
\end{equation*}
where $F\in C(\Var(\mathbb C)^{sm}/S)$ is such that $M=D(\mathbb A^1,et)(F)$.
\item[(iii)] For the Corti-Hanamura weight structure on $\DA^-(S)$, we have by functoriality of (i) the functor
\begin{equation*}
\Bti_S^*:\DA^-(S)\to D_{fil,\infty}(S^{an}) \; , \; 
M\mapsto\Bti_S^*M=e(S^{an})_*\underline{\sing}_{\mathbb D^*}\An_S^*(F,W)
\end{equation*}
where $(F,W)\in C_{fil}(\Var(\mathbb C)^{sm}/S)$ is such that $(M,W)=D(\mathbb A^1,et)(F,W)$.
\end{itemize}
Note that by \cite{B3}, $\An_S^*$ and $\widetilde\Cw_S^*$ derive trivially.
\end{defi}

Note that, by considering the explicit $\mathbb D_S^1$ local model for presheaves on $\AnSp(\mathbb C)^{sm}/S^{an}$,
$\Bti_S^*(\DA^-(S))\subset D^-(S^{an})$ ; 
by considering the explicit $\mathbb I_S^1$ local model for presheaves on $\CW^{sm}/S^{cw}$,
$\widetilde\Bti_S^*(\DA^-(S))\subset D^-(S^{an})$.

Let $f:T\to S$ a morphism, with $T,S\in\Var(\mathbb C)$. We have, 
for $M\in\DA(S)$, $(F,W)\in C_{fil}(\Var(\mathbb C)^{sm}/S)$ such that $(M,W)=D(\mathbb A^1,et)(F,W)$,
and an equivalence $(\mathbb A^1,et)$ local $e:f^*(F,W)\to (F',W)$ 
with $(F',W)\in C_{fil}(\Var(\mathbb C)^{sm}/S)$ such that $(f^*M,W)=D(\mathbb A^1,et)(F',W)$ 
the following canonical transformation map in $D_{fil}(T)$:
\begin{eqnarray}
T^0(f,\Bti)(M):f^*\Bti_S^*M:=f^*e(S^{an})_*\underline{\sing}_{\mathbb D^*}\An_S^*(F,W)
\xrightarrow{T(f,e)(-)} e(T^{an})_*f^*\underline{\sing}_{\mathbb D^*}\An_S^*(F,W) \\
\xrightarrow{e(T^{an})_*T(f,c)(F,W)}e(T^{an})_*\underline{\sing}_{\mathbb D^*}f^*\An_T^*(F,W)
\xrightarrow{=}e(T^{an})_*\underline{\sing}_{\mathbb D^*}\An_T^*f^*(F,W) \\
\xrightarrow{e(T^{an})_*\underline{\sing}_{\mathbb D^*}\An_T^*e}
e(T^{an})_*\underline{\sing}_{\mathbb D^*}\An_T^*(F',W)=:\Bti_T^*f^*M
\end{eqnarray}

\begin{defi}\label{TgBti}
Let $f:T\to S$ a morphism, with $T,S\in\Var(\mathbb C)$. 
Consider the graph factorization $f:T\xrightarrow{l}T\times S\xrightarrow{p}S$ of $f$
with $l$ the graph closed embedding and $p$ the projection.
We have, for $M\in\DA_c(S)$,  the following canonical transformation map in $D_{fil,c}(T^{an})$:
\begin{eqnarray*}
T(f,\Bti)(M,W):f^{*w}\Bti_S^*(M,W):=l^*\Gamma_T^{\vee,w}p^*\Bti_S^*(F,W) \\
\xrightarrow{T^0(p,\Bti)(-)}l^*\Gamma_T^{\vee,w}\Bti_{T\times S}^*p^*(F,W)
\xrightarrow{\gamma_T^{\vee}(p^*(F,W))}l^*\Gamma_T^{\vee,w}\Bti_{T\times S}^*\Gamma_T^{\vee}p^*(F,W) \\
\xrightarrow{=}l^*\Bti_{T\times S}^*\Gamma_T^{\vee}p^*(F,W)
\xrightarrow{T^0(l,\Bti)(-)}\Bti_T^*l^*\Gamma_T^{\vee}p^*(F,W)=\Bti_T^*f^*(M,W).
\end{eqnarray*}
where we use definition \ref{fw}.
\end{defi}

\begin{defi}\label{TBtiSix}
\begin{itemize}
\item Let $f:X\to S$ a morphism, with $X,S\in\Var(\mathbb C)$. 
We have, for $M\in\DA_c(X)$, the following transformation map in $D_{fil,c}(S^{an})$
\begin{eqnarray*}
T_*(f,\Bti)(M,W):\Bti_S^*(Rf_*(M,W))\xrightarrow{\ad(f^*,Rf_{*w})(\Bti_S^*(Rf_*(M,W)))}Rf_{*w}f^{*w}\Bti_S^*(Rf_*(M,W)) \\
\xrightarrow{T(f,\Bti)(Rf_*(M,W))}
Rf_{*w}\Bti_X^*(f^*Rf_*(M,W))\xrightarrow{\Bti_X^*(\ad(f^*,Rf_*)(M,W))}Rf_{*w}\Bti_X^*(M,W)
\end{eqnarray*}
Clearly if $l:Z\hookrightarrow S$ is a closed embedding, then $T_*(l,\Bti)(M,W)$ is an isomorphism
since $\ad(l^*,l_*)(-):l^*l_*(M,W)\to (M,W)$ is an isomorphism (see section 3).

\item Let $f:X\to S$ a morphism with $X,S\in\Var(\mathbb C)$. Assume there exist a factorization
$f:X\xrightarrow{l}Y\times S\xrightarrow{p_S}S$ with $Y\in\SmVar(\mathbb C)$, $l$ a closed embedding and $p_S$ the projection. 
We have then, for $M\in\DA_c(X)$, using theorem \ref{mainBti} for closed embeddings, 
the following transformation map in $D_{fil}((Y\times S)^{an})$
\begin{eqnarray*}
T_!(f,\Bti)(M):Rf_{!w}\Bti_X^*(M,W)=Rp_{S!w}l_*\Bti_X^*(M,W) \\
\xrightarrow{T_*(l,\Bti)(M,W)}^{-1} Rp_{S!w}\Bti(Y\times S)^*(l_*(M,W)) \\
\xrightarrow{\Bti(Y\times S)^*\ad(Lp_{S\sharp},p_S^*)(l_*(M,W))}
Rp_{S!w}\Bti(Y\times S)^*(p_S^*Lp_{S\sharp}l_*(M,W))\xrightarrow{T(p_S,\Bti)(p_{S\sharp}l_*(M,W))} \\
Rp_{S!w}p_S^*\Bti(Y\times S)^*(Lp_{S\sharp}l_*(M,W))=Rp_{S!w}p_S^{!w}\Bti(Y\times S)^*(Rf_!(M,W)) \\
\xrightarrow{\ad(Rp_{S!w},p_S^{!w})(-)}\Bti(Y\times S)^*(Rf_!(M,W))
\end{eqnarray*}
Clearly, for $f:X\to S$ a proper morphism, with $X,S\in\Var(\mathbb C)$ we have, for $M\in\DA_c(Y\times S)$,
$T_!(f,\Bti)(M,W)=T_*(f,\Bti)(M,W)$.

\item Let $f:X\to S$ a morphism with $X,S\in\Var(\mathbb C)$. 
We have, using the second point, for $M\in\DA(S)$, the following transformation map in $D_{fil}(X^{an})$
\begin{eqnarray*}
T^!(f,\Bti)(M,W):\Bti_X^*(f^!(M,W))\xrightarrow{\ad(f_!,Rf^!)(\Bti_X^*(f^!(M,W)))}f^{!w}Rf_{!w}\Bti_X^*(f^!(M,W)) \\
\xrightarrow{T_!(f,\Bti)((f^!(M,W)))}
f^{!w}\Bti_S^*(f_!f^!(M,W))\xrightarrow{\Bti_S^*(\ad(f_!,f^!)(M,W))}f^{!w}\Bti_S^*(M,W)
\end{eqnarray*}

\item Let $S\in\Var(\mathbb C)$. We have, for $M,N\in\DA(S)$ and $F,G\in C_(\Var(\mathbb C)^{sm}/S))$  
such that $M=D(\mathbb A^1,et)(F)$ and $N=D(\mathbb A^1,et)(G)$, 
the following transformation map in $D_{fil}(S^{an})$
\begin{eqnarray*}
\Bti_S^*(M,W)\otimes\Bti_S^*(N,W):=
(e(S)_*\underline{\sing}_{\mathbb D^*}\An_S^*(F,W))\otimes(e(S)_*\underline{\sing}_{\mathbb D^*}\An_S^*(G,F)) \\
\xrightarrow{T(\sing_{D^*},\otimes)(\An_S^*(F,W),\An_S^*(G,F))}
e(S)_*\underline{\sing}_{\mathbb D^*}\An_S^*((F,W)\otimes(G,W))=:\Bti_S^*((M,W)\otimes(N,W))
\end{eqnarray*}
\end{itemize}
\end{defi}

\begin{thm}\label{mainBti}
\begin{itemize}
\item[(i)]Let $f:X\to S$ a morphism, with $X,S\in\Var(\mathbb C)$. For $M\in\DA_c(S)$,
\begin{equation*}
T(f,\Bti)(M,W):f^{*w}\Bti_S^*(M,W)\xrightarrow{\sim}\Bti_X^*f^*(M,W)
\end{equation*}
is an isomorphism in $D_{fil}(X^{an})$.
\item[(ii)] Let $f:X\to S$ a morphism, with $X,S\in\Var(\mathbb C)$. For $M\in\DA_c(X)$,
\begin{equation*}
T_!(f,\Bti)(M,W):Rf_{!w}\Bti_X^*(M,W)\xrightarrow{\sim}\Bti_S^*Rf_!(M,W)
\end{equation*}
is an isomorphism.
\item[(iii)] Let $f:X\to S$ a morphism, with $X,S\in\Var(\mathbb C)$. For $M\in\DA_c(X)$,
\begin{equation*}
T_*(f,\Bti)(M,W):Rf_{*w}\Bti_X^*(M,W)\xrightarrow{\sim}\Bti_S^*Rf_*(M,W)
\end{equation*}
is an isomorphism.
\item[(iv)] Let $f:X\to S$ a morphism, with $X,S\in\Var(\mathbb C)$. For $M\in\DA_c(S)$,
\begin{equation*}
T^!(f,\Bti)(M,W):f^{!w}\Bti_S^*(M,W)\xrightarrow{\sim}\Bti_X^*f^!(M,W)
\end{equation*}
is an isomorphism.
\item[(v)] Let $S\in\Var(\mathbb C)$. For $M,N\in\DA_c(S)$,
\begin{equation*}
T(\otimes,\Bti)(M,W):\Bti_S^*(M,W)\otimes\Bti_S^*N\xrightarrow{\sim}\Bti_X^*((M,W)\otimes(N,W))
\end{equation*}
is an isomorphism.
\end{itemize}
\end{thm}

\begin{proof}
By functoriality it reduced to the case of Corti-Hanamura motives which then follows from \cite{AyoubB}.
\end{proof}

The main result on the Betti realization functor is the following
\begin{thm}\label{mainBtiCorCor}
\begin{itemize}
\item[(i)] We have $\Bti_S^*=\widetilde\Bti_S^*$ on $\DA^-(S)$
\item[(ii)] The canonical transformations $T(f,\Bti)$, for $f:T\to S$ a morphism in $\Var(\mathbb C)$,
define a morphism of 2 functor
\begin{equation*}
\Bti_{\cdot}^*:\DA(\cdot)\to D(\cdot^{an}), \; S\in\Var(\mathbb C) \mapsto\Bti_S^*:\DA(S)\to D(S^{an})
\end{equation*}
which is a morphism of homotopic 2 functor.
\item[(ii)'] The canonical transformations $T(f,\Bti)$, for $f:T\to S$ a morphism in $\Var(\mathbb C)$,
define a morphism of 2 functor
\begin{equation*}
\Bti_{\cdot}^*:\DA(\cdot)\to D_{fil}(\cdot^{an}), \; S\in\Var(\mathbb C) \mapsto\Bti_S^*:\DA(S)\to D_{fil}(S^{an})
\end{equation*}
which is a morphism of homotopic 2 functor.
\end{itemize}
\end{thm}

\begin{proof}
\noindent(i): See \cite{B3}

\noindent(ii) and (ii)':Follows from theorem \ref{mainBti}.
\end{proof}

\begin{rem}
For $X\in\Var(\mathbb C)$, the quasi-isomorphisms
\begin{equation*}
\mathbb Z\Hom(\bar{\mathbb D}_{et}^{\bullet},X)\xrightarrow{\An^*}\mathbb Z\Hom(\bar{\mathbb D}^n(0,1),X^{an})
\xrightarrow{\Hom(i,X^{cw})}\mathbb Z\Hom([0,1]^n,X^{cw}),
\end{equation*}
where, 
\begin{equation*}
\bar{\mathbb D}^n_{et}:=(e:U\to\mathbb A^n,\bar{\mathbb D}^n(0,1)\subset e(U))
\in\Fun(\mathcal V^{et}_{\mathbb A^n}(\bar{\mathbb D}^n(0,1)),\Var(\mathbb C)) 
\end{equation*}
is the system of etale neighborhood of the closed ball $\bar{\mathbb D}^n(0,1)\subset\mathbb A^n$,
and $i:[0,1]^n\hookrightarrow\bar{\mathbb D}^n(0,1)$ is the closed embedding,
shows that a closed singular chain $\alpha\in\mathbb Z\Hom^n([0,1]^n,X^{cw})$, is homologue to 
a closed singular chain 
\begin{equation*}
\beta=\alpha+\partial\gamma=\tilde{\beta}_{|[0,1]^n}\in\mathbb Z\Hom^n(\Delta^n,X^{cw}) 
\end{equation*}
which is the restriction by the closed embedding $[0,1]^n\hookrightarrow U^{cw}\xrightarrow{e}\mathbb A^n$,
where $e:U\to\mathbb A^n$ an etale morphism with $U\in\Var(\mathbb C)$,
of a complex algebraic morphism $\tilde\beta:U\to X$. 
Hence $\beta([0,1]^n)=\tilde{\beta}([0,1]^n)\subset X$ is the restriction
of a real algebraic subset of dimension $n$ in $Res_{\mathbb R}(X)$ 
(after restriction a scalar that is under the identification $\mathbb C\simeq\mathbb R^2$). 
\end{rem}

\begin{defi}
Let $S\in\Var(\mathbb C)$
The cohomological Betti realization functor is
\begin{eqnarray*}
\Bti_S^{\vee}:\DA(S)\to D(S^{cw}), \\
M\mapsto\Bti_S^{\vee}(M):=R\mathcal Hom(\Bti_S^*M,\mathbb Z_{S^{cw}})=R\mathcal Hom(M,\Bti_{S*}\mathbb Z_{S^{cw}})
\end{eqnarray*}
where for $\Bti_{S*}:K\in D(S^{cw})\mapsto R\An_{S*}e(S^{an})^*K\in\DA(S)$ is the right ajoint to $\Bti_S^*$.
\end{defi}

\subsection{The Hodge realization functor for relative motives} 

Let $S\in\Var(\mathbb C)$. Let $S=\cup_{i=1}^sS_i$ an open cover such that there exists
closed embedding $i_i:S\hookrightarrow\tilde S_i$ with $\tilde S_i\in\SmVar(\mathbb C)$.
Recall (see section 5.2) that $D_{\mathcal D(1,0)fil,rh}(S/(\tilde S_I))\times_I D_{fil}(S^{an})$ is the category  
\begin{itemize}
\item whose set of objects is the set of triples $\left\{(((M_I,F,W),u_{IJ}),(K,W),\alpha)\right\}$ with
\begin{eqnarray*} 
((M_I,F,W),u_{IJ})\in D_{\mathcal D(1,0)fil,rh}(S/(\tilde S_I)), \, (K,W)\in D_{fil}(S^{an}), \\ 
\alpha:T(S/(\tilde S_I))(K,W)\otimes\mathbb C_{S^{an}}\to DR(S)^{[-]}(((M_I,W),u_{IJ})^{an})
\end{eqnarray*}
where $\alpha$ is an morphism in $D_{fil}(S^{an}/(\tilde S_{I}^{an}))$,
\item and whose set of morphisms consists of 
\begin{equation*}
\phi=(\phi_D,\phi_C,[\theta]):(((M_{1I},F,W),u_{IJ}),(K_1,W),\alpha_1)\to(((M_{2I},F,W),u_{IJ}),(K_2,W),\alpha_2)
\end{equation*}
where $\phi_D:((M_1,F,W),u_{IJ})\to((M_2,F,W),u_{IJ})$ and $\phi_C:(K_1,W)\to (K_2,W)$ 
are morphisms and
\begin{eqnarray*}
\theta=(\theta^{\bullet},I(DR(S)(\phi^{an}_D))\circ I(\alpha_1),I(\alpha_2)\circ I(\phi_C\otimes I)): \\
I(T(S/(\tilde S_I))(K_1,W))\otimes\mathbb C_{S^{an}}[1]\to I(DR(S)(((M_{2I},W),u_{IJ})^{an}))  
\end{eqnarray*}
is an homotopy,  
$I:D_{fil}(S^{an}/(\tilde S^{an}_{I}))\to K_{fil}(S^{an}/(\tilde S^{an}_{I}))$
being the injective resolution functor, and for
\begin{itemize}
\item $\phi=(\phi_D,\phi_C,[\theta]):(((M_{1I},F,W),u_{IJ}),(K_1,W),\alpha_1)\to(((M_{2I},F,W),u_{IJ}),(K_2,W),\alpha_2)$
\item $\phi'=(\phi'_D,\phi'_C,[\theta']):(((M_{2I},F,W),u_{IJ}),(K_2,W),\alpha_2)\to(((M_{3I},F,W),u_{IJ}),(K_3,W),\alpha_3)$
\end{itemize}
the composition law is given by 
\begin{eqnarray*}
\phi'\circ\phi:=(\phi'_D\circ\phi_D,\phi'_C\circ\phi_C,
I(DR(S)(\phi^{'an}_D))\circ[\theta]+[\theta']\circ I(\phi_C\otimes I)[1]): \\
(((M_{1I},F,W),u_{IJ}),(K_1,W),\alpha_1)\to(((M_{3I},F,W),u_{IJ}),(K_3,W),\alpha_3),
\end{eqnarray*}
in particular for 
$(((M_I,F,W),u_{IJ}),(K,W),\alpha)\in D_{\mathcal D(1,0)fil,rh}(S/(\tilde S_I))\times_I D_{fil}(S^{an})$,
\begin{equation*}
I_{(((M_I,F,W),u_{IJ}),(K,W),\alpha)}=((I_{M_I}),I_K,0),
\end{equation*}
\end{itemize}
together with the localization functor
\begin{eqnarray*}
(D(zar),I):C_{\mathcal D(1,0)fil,rh}(S/(\tilde S_I))\times_I D_{fil}(S^{an})
\to D_{\mathcal D(1,0)fil,rh}(S/(\tilde S_I))\times_I D_{fil}(S^{an}) \\
\to D_{\mathcal D(1,0)fil,rh,\infty}(S/(\tilde S_I))\times_I D_{fil}(S^{an}).
\end{eqnarray*}
Moreover,
\begin{itemize}
\item For $(((M_{I},F,W),u_{IJ}),(K,W),\alpha)D_{\mathcal D(1,0)fil,rh}(S/(\tilde S_I))\times_I D_{fil}(S^{an})$, we set
\begin{equation*}
(((M_{I},F,W),u_{IJ}),(K,W),\alpha)[1]:=(((M_{I},F,W),u_{IJ})[1],(K,W)[1],\alpha[1]).
\end{equation*}
\item For 
\begin{equation*}
\phi=(\phi_D,\phi_C,[\theta]):(((M_{1I},F,W),u_{IJ}),(K_1,W),\alpha_1)\to(((M_{2I},F,W),u_{IJ}),(K_2,W),\alpha_2)
\end{equation*}
a morphism in $D_{\mathcal D(1,0)fil,rh}(S/(\tilde S_I))\times_I D_{fil}(S^{an})$, we set (see \cite{CG} definition 3.12)
\begin{eqnarray*}
\Cone(\phi):=(\Cone(\phi_D),\Cone(\phi_C),((\alpha_1,\theta),(\alpha_2,0)))
\in D_{\mathcal D(1,0)fil,rh}(S/(\tilde S_I))\times_I D_{fil}(S^{an}),
\end{eqnarray*}
$((\alpha_1,\theta),(\alpha_2,0))$ being the matrix given by the composition law, together with the canonical maps
\begin{itemize}
\item $c_1(-)=(c_1(\phi_D),c_1(\phi_C),0):(((M_{2I},F,W),u_{IJ}),(K_2,W),\alpha_2)\to\Cone(\phi)$
\item $c_2(-)=(c_2(\phi_D),c_2(\phi_C),0):\Cone(\phi)\to (((M_{1I},F,W),u_{IJ}),(K_1,W),\alpha_1)[1]$.
\end{itemize}
\end{itemize}

Consider the category 
\begin{equation*}
(D_{\mathcal D(1,0)fil}(\tilde S_I)\times_ID_{fil}(\tilde S_I^{an}))\subset\Fun(\Gamma(\tilde S_I),\TriCat)
\end{equation*}
\begin{itemize}
\item whose objects are $(((M_I,F,W),(K_I,W),\alpha_I),u_{IJ})\in\Fun(\Gamma(\tilde S_I),\TriCat)$ such that 
\begin{equation*}
((M_I,F,W),(K_I,W),\alpha_I)\in D_{\mathcal D(1,0)fil}(\tilde S_I)\times_ID_{fil}(\tilde S_I^{an})
=:\mathcal D(\tilde S_I) 
\end{equation*}
and for $I\subset J$, 
\begin{eqnarray*}
u_{IJ}:((M_I,F,W),(K_I,W),\alpha_I)\to \\
p_{IJ*}((M_J,F,W),(K_J,W),\alpha_J):=(p_{IJ*}(M_J,F,W),p_{IJ*}(K_J,W),p_{IJ*}\alpha_J)
\end{eqnarray*}
are morphisms in $\mathcal D(\tilde S_I)$,
\item whose morphisms $m=(m_I):(((M_I,F,W),(K_I,W),\alpha_I),u_{IJ})\to (((M'_I,F,W),(K'_I,W),\alpha'_I),v_{IJ})$
is a family of morphism such that $v_{IJ}\circ m_I=p_{IJ*}m_J\circ u_{IJ}$ in $\mathcal D(\tilde S_I)$
\end{itemize}
We have then the identity functor
\begin{eqnarray*}
I_S:D_{\mathcal D(1,0)fil}(S/(\tilde S_I))^0\times_ID_{fil}(S^{an})\to
(D_{\mathcal D(1,0)fil}(\tilde S_I)\times_ID_{fil}(\tilde S_I^{an})), \\
(((M_I,F,W),u_{IJ}),(K,W),\alpha)\mapsto (((M_I,F,W),i_{I*}j_I^*(K,W),j_I^*\alpha),(u_{IJ},I,0)), \\ 
m=(m_I,n)\mapsto m=(m_I,i_*j_I^*n)
\end{eqnarray*}
which is a full embedding since for $((M_I,F,W),u_{IJ})\in D_{\mathcal D(1,0)fil}(S/(\tilde S_I))^0$,
\begin{equation*}
u_{IJ}:(M_I,F,W)\to p_{IJ*}(M_J,F,W) 
\end{equation*}
are isomorphisms in $D_{\mathcal D(1,0)fil}(\tilde S_I)$, and hence for 
$((M_I,F,W),u_{IJ}),(K,W),\alpha)\in D_{\mathcal D(1,0)fil}(S/(\tilde S_I))^0\times_ID_{fil}(S^{an})$,
\begin{eqnarray*}
(u_{IJ},I,0):((M_I,F,W),i_{I*}j_I^*(K,W),j_I^*\alpha)\to \\
p_{IJ*}((M_J,F,W),i_{J*}j_J^*(K,W),j_J^*\alpha)=(p_{IJ*}(M_J,F,W),i_{I*}j_I^*(K,W),j_I^*\alpha)
\end{eqnarray*}
are isomorphisms in $\mathcal D(\tilde S_I)$.

\begin{defi}\label{IUSm}
For $h:U\to S$ a smooth morphism with $S,U\in\SmVar(\mathbb C)$ 
and $h:U\xrightarrow{n}X\xrightarrow{f}S$ a compactification of $h$ with $n$ an open embedding, $X\in\SmVar(\mathbb C)$
such that $D:=X\backslash U=\cup_{i=1}^sD_i\subset X$ is a normal crossing divisor, 
we denote by, using definition \ref{DHdgalpha} and definition \ref{wtildew}
\begin{eqnarray*}
I(U/S):h_{!Hdg}h^{!Hdg}\mathbb Z_S^{Hdg}\xrightarrow{:=} \\
(p_{S*}E_{zar}(\Omega^{\bullet}_{X\times S/S}\otimes_{O_{X\times S}}(n\times I)_{!Hdg}\Gamma_U^{\vee,Hdg}(O_{U\times S},F_b)),
\mathbb D_Sh_*E_{usu}\mathbb Z_{U^{an}},h_!\alpha(U)) \\
\xrightarrow{((DR(X\times S/S)(\ad((n\times I)_{!Hdg},(n\times I)^*)(-)),0),I,0)} \\
(\Cone((\Omega_{/S}^{\Gamma,pr}(i_{D_i}\times I))_{i\in[1,\ldots,s]}:
p_{S*}E_{zar}(\Omega^{\bullet}_{X\times S/S}\otimes_{O_{X\times S}}\Gamma_X^{\vee,Hdg}(O_{X\times S},F_b))\to \\
(\cdots\to(p_{S*}E_{zar}(\Omega^{\bullet}_{D_I\times S/S}\otimes_{O_{D_I\times S}}\Gamma_{D_I}^{\vee,Hdg}(O_{D_I\times S},F_b)))
\to\cdots)),\mathbb D_Sh_*E_{usu}\mathbb Z_{U^{an}},h_!\alpha(U)) \\
\xrightarrow{=:}
(\mathcal F_S^{FDR}(\mathbb Z(U/S)),\Bti_S^*\mathbb Z(U/S),\alpha(\mathbb Z(U/S)))
\end{eqnarray*}
the canonical isomorphism in $D_{\mathcal Dfil}(S)\times_ID(S^{an})$, with 
\begin{equation*}
h^{!Hdg}\mathbb Z_S^{Hdg}=(\Gamma_U^{\vee,Hdg}(O_{U\times S},F_b),\mathbb Z_{U^{an}},\alpha(U))\in HM(U), 
\end{equation*}
and $i_{D_i}:D_i\hookrightarrow X$ are the closed embeddings.
\end{defi}

\begin{lem}\label{thetam}
Let $S\in\SmVar(\mathbb C)$. Let $g:U'/S\to U/S$ o morphism with $U/S:=(U,h),U'/S:=(U',h)\in\Var(\mathbb C)^{sm}/S$.
Let $h:U\xrightarrow{n}X\xrightarrow{f}S$ a compactification of $h$ with $n$ an open embedding, $X\in\SmVar(\mathbb C)$
such that $D:=X\backslash U=\cup_{i=1}^sD_i\subset X$ is a normal crossing divisor, 
Let $h':U\xrightarrow{n'}X'\xrightarrow{f'}S$ a compactification of $h'$ with $n'$ an open embedding, $X'\in\SmVar(\mathbb C)$
such that $D':=X\backslash U=\cup_{i=1}^sD_i\subset X$ is a normal crossing divisor and such that
$g:U'\to U$ extend to $\bar g:X'\to X$, see definition-proposition \ref{RCHdef0}.
Then, using definition \ref{IUSm}, 
the following diagram in $D_{\mathcal Dfil}(S)\times_ID(S^{an})$ commutes
\begin{eqnarray*}
\xymatrix{
h'_{!Hdg}h^{'!Hdg}\mathbb Z_S^{Hdg}\ar[rrr]^{I(U'/S)}\ar[d]_{\ad(g_{!Hdg},g^{!Hdg})(h^{!Hdg}\mathbb Z_S^{Hdg})} & \, & \, &
(\mathcal F_S^{FDR}(\mathbb Z(U'/S)),\Bti_S^*\mathbb Z(U'/S),\alpha(\mathbb Z(U'/S)))
\ar[d]^{(\Omega_{/S}^{\Gamma,pr}(\hat R_S^{CH}(g)),\Bti_S^*(g),\theta(g))} \\ 
h_{!Hdg}h^{!Hdg}\mathbb Z_S^{Hdg}\ar[rrr]^{I(U/S)} & \, & \, &
(\mathcal F_S^{FDR}(\mathbb Z(U/S)),\Bti_S^*\mathbb Z(U/S),\alpha(\mathbb Z(U/S)))} 
\end{eqnarray*}
where 
\begin{eqnarray*}
\theta(g):=R_{\mathcal D}([\Gamma_g]):I(\Bti_S^*\mathbb Z(U'/S)\otimes\mathbb C)[1]
\to I(DR(S)(o_{fil}\mathcal F_S^{FDR}(\mathbb Z(U/S))^{an}))
\end{eqnarray*}
is the homotopy given by the third term of the Deligne homology class of the graph $\Gamma_g\subset U'\times_S U$
(see definition \ref{Delkdef}) and 
$o_{fil}:C_{\mathcal Dfil}(S)\to C_{\mathcal D}(S)$ 
is the forgetful functor and we recall (see section 5.2) that
$I:C(S_{\mathcal C}^{an}/\tilde S_I^{an})\to K(S_{\mathcal C}^{an}/\tilde S_I^{an})$
is the injective resolution functor.
\end{lem}

\begin{proof}
Immediate from definition.
\end{proof}

We now define the Hodge realization functor.

\begin{defi}\label{HodgeRealDAsing}
Let $k\subset\mathbb C$ a subfield. Let $S\in\Var(\mathbb C)$. Let $S=\cup_{i=1}^sS_i$ an open cover such that there exists
closed embedding $i_i:S\hookrightarrow\tilde S_i$ with $\tilde S_i\in\SmVar(\mathbb C)$.
We define the Hodge realization functor, using definition \ref{DRalgdefFunct}, 
definition \ref{bettidef}, and lemma \ref{thetam}
\begin{eqnarray*}
\mathcal F_S^{Hdg}:=(\mathcal F_S^{FDR},\Bti_S^*):
C(\Var(\mathbb C)^{sm}/S)\to D_{\mathcal D(1,0)fil}(S/(\tilde S_I))\times_I D_{fil}(S^{an}) 
\end{eqnarray*}
first on objects and then on morphisms :
\begin{itemize}
\item for $F\in C(\Var(\mathbb C)^{sm}/S)$, taking $(F,W)\in C_{fil}(\Var(\mathbb C)^{sm}/S)$
such that $D(\mathbb A^1,et)(F,W)$ gives the weight structure on $D(\mathbb A^1,et)(F)$,
\begin{eqnarray*}
\mathcal F_S^{Hdg}(F):=(\mathcal F_S^{FDR}(F),\Bti_S^*F,\alpha(F)):= \\
(e(S)_*\mathcal Hom((\hat R_{\tilde S_I}^{CH}(\rho_{\tilde S_I}^*Li_{I*}j_I^*(F,W)),
\hat R^{CH}(T^q(D_{IJ})(-))),(E_{zar}(\Omega^{\bullet,\Gamma,pr}_{/\tilde S_I},F_{DR}),T_{IJ})), \\
e(S)_*\underline{\sing}_{\mathbb D^*}\An_S^*L(F,W),\alpha(F)) 
\in D_{\mathcal D(1,0)fil}(S/(\tilde S_I))\times_I D_{fil}(S^{an})
\end{eqnarray*}
where $\alpha(F)$ is the map in $D_{\mathbb Dfil}(S^{an}/(\tilde S_I^{an}))$,
writing for short $DR(S):=DR(S)^{[-]}:=(DR(\tilde S_I)[-d_{\tilde S_I}])$
\begin{eqnarray*}
\alpha(F):T(S/(\tilde S_I))((\Bti_S^*(M,W))\otimes\mathbb C_S):=
(i_{I*}j_I^*((e(S)_*\underline{\sing}_{\mathbb D^*}\An_S^*L(F,W))\otimes\mathbb C_S),I) \\
\xrightarrow{=}
(e(\tilde S_I)_*\underline{\sing}_{\mathbb D^*}\An_{\tilde S_I}^*Li_{I*}j_I^*(F,W)\otimes\mathbb C_S,
T(p_{IJ},\An)(Li_{I*}j_I^*(F,W))) \\
\xrightarrow{=}
(((\cdot\to\oplus_{(U_{I\alpha},h_{I\alpha})\in V_I}h_{I\alpha !}h^!_{I\alpha}\mathbb Z_{\tilde S_I}
\xrightarrow{\ad(g^{\bullet !}_{I,\alpha,\beta},g^{\bullet}_{I,\alpha,\beta !})(-)} 
\oplus_{(U_{I\alpha},h_{I\alpha})\in V_I}h_{I\alpha !}h^!_{I\alpha}\mathbb Z_{\tilde S_I}\to\cdot),u_{IJ})) \\
\xrightarrow{(\alpha(\mathbb Z(U_{I\alpha}/\tilde S_I)),\theta(g^{\bullet}_{I,\alpha,\beta}))} \\
DR(S)(o_{fil}(e(S)_*\mathcal Hom((\hat R_{\tilde S_I}^{CH}(\rho_{\tilde S_I}^*Li_{I*}j_I^*(F,W)),
\hat R^{CH}(T^q(D_{IJ})(-))),(E_{zar}(\Omega^{\bullet,\Gamma,pr}_{/\tilde S_I},F_{DR}),T_{IJ})))^{an}) \\
\xrightarrow{=:}
DR(S)((o_{fil}\mathcal F^{FDR}_S(M,W))^{an})
\end{eqnarray*}
using lemma \ref{thetam}, $(\alpha(\mathbb Z(U_{I\alpha}/S)),\theta(g^{\bullet}_{I,\alpha,\beta}))$
being the matrix given inductively by the composition law in 
$D_{\mathcal D(1,0)fil}(\tilde S_I)\times_I D_{fil}(\tilde S_I^{an})$,
that is we have the following isomorphism in 
$(D_{\mathcal D(1,0)fil}(\tilde S_I)\times_I D_{fil}(\tilde S_I^{an}))$,
denoting for short $V_I:=\Var(\mathbb C)^{sm}/\tilde S_I$
\begin{eqnarray*}
(I^{\bullet}(U_{I\alpha}/\tilde S_I)): \\
((\cdot\to\oplus_{(U_{I\alpha},h_{I\alpha})\in V_I}h_{I\alpha !Hdg}h^{!Hdg}_{I\alpha}\mathbb Z^{Hdg}_{\tilde S_I}
\xrightarrow{\ad(g^{\bullet,!Hdg}_{I,\alpha,\beta},g^{\bullet}_{I,\alpha,\beta !Hdg})(-)} 
\oplus_{(U_{I\alpha},h_{I\alpha})\in V_I}h_{I\alpha !Hdg}h^{!Hdg}_{I\alpha}\mathbb Z^{Hdg}_{\tilde S_I}\to\cdot),u_{IJ}) \\
\xrightarrow{\sim}I_S(\mathcal F_S^{Hdg}(F):=(\mathcal F_S^{FDR}(F),\Bti_{\tilde S_I}^*Li_{I*}j_I^*(F,W),\alpha(F)))
\end{eqnarray*}
where we denote by $g^n_{I,\alpha,\beta}:U_{I\alpha}\to U_{I\beta}$ which satisfy 
$h_{I\beta}\circ g^n_{I,\alpha,\beta}=h_{I\alpha}$ the morphisms in the canonical projective resolution
\begin{eqnarray*}
q:Li_{I*}j_I^*F:=(\cdots\to\oplus_{(U_{I\alpha},h_{I\alpha})\in\Var(\mathbb C)^{sm}/\tilde S_I}\mathbb Z(U_{I\alpha}/\tilde S_I)
\xrightarrow{(\mathbb Z(g^{\bullet}_{\alpha,\beta}))} \\
\oplus_{(U_{I\alpha},h_{I\alpha})\in\Var(\mathbb C)^{sm}/\tilde S_I}\mathbb Z(U_{I\alpha}/\tilde S_I)\to\cdots)\to i_{I*}j_I^*F,
\end{eqnarray*}
\item for $m:F_1\to F_2$ a morphism in $C(\Var(\mathbb C)^{sm}/S)$, the morphism 
$\mathcal F_S^{Hdg}(m)$ in $D_{\mathcal D(1,0)fil}(S/(\tilde S_I))\times_I D_{fil}(S^{an})$ is given by
\begin{eqnarray*}
\mathcal F_S^{Hdg}(m):&=&{I_S^{-,-}}^{-1}((I^{\bullet}(U_{I\alpha}/(\tilde S_I)))\circ
(\ad(l_{I\alpha,\beta}^{\bullet !Hdg},l^{\bullet}_{I\alpha,\beta !Hdg})(\mathbb Z_{U_{I\alpha}}^{Hdg})) 
\circ (I^{\bullet}(U_{I\alpha}/(\tilde S_I)))^{-1}) \\
&=&(\mathcal F_S^{FDR}(m),\Bti_S^*(m),\theta(m):=(\theta(l_{I\alpha,\beta}))):
\mathcal F_S^{Hdg}(F_1)\to\mathcal F_S^{Hdg}(F_2)
\end{eqnarray*}
using lemma \ref{thetam}, that is we have the following commutative diagram in 
$(D_{\mathcal D(1,0)fil}(\tilde S_I)\times_I D_{fil}(\tilde S_I^{an}))$,
denoting for short $V_I:=\Var(\mathbb C)^{sm}/\tilde S_I$,
\begin{eqnarray*}
\begin{tikzcd}
((\cdot\to\oplus_{(U_{I\alpha},h_{I\alpha})\in V_I}h_{I\alpha !Hdg}h^{!Hdg}_{I\alpha}\mathbb Z^{Hdg}_{\tilde S_I}
\xrightarrow{A^{Hdg}_{g_{1I,\alpha,\beta}^{\bullet}}} 
\oplus_{(U_{I\alpha},h_{I\alpha})\in V_I}h_{I\alpha !Hdg}h^{!Hdg}_{I\alpha}\mathbb Z^{Hdg}_{\tilde S_I}\to\cdot),u_{IJ})
\ar[r,"(I^{\bullet}(U_{I\alpha}/\tilde S_I))"]
\ar[d,"\ad(l_{I{,}\alpha{,}\beta}^{\bullet !Hdg}{,}l^{I{,}\bullet}_{\alpha{,}\beta !Hdg})(-)"'] & 
\mathcal F_S^{Hdg}(F_1)\ar[d,"\mathcal F_S^{Hdg}(m)=(\mathcal F_S^{FDR}(m){,}\Bti_S^*(m){,}(\theta(l_{I\alpha{,}\beta})))"'] \\
((\cdot\to\oplus_{(U_{I\alpha},h_{I\alpha})\in V_I}h_{I\alpha !Hdg}h^{!Hdg}_{I\alpha}\mathbb Z^{Hdg}_{\tilde S_I}
\xrightarrow{A^{Hdg}_{g_{2I,\alpha,\beta}^{\bullet}}}
\oplus_{(U_{I\alpha},h_{I\alpha})\in V_I}h_{I\alpha !Hdg}h^{!Hdg}_{I\alpha}\mathbb Z^{Hdg}_{\tilde S_I}\to\cdot),u_{IJ})
\ar[r,"(I^{\bullet}(U_{\alpha}/\tilde S_I))"] & \mathcal F_S^{Hdg}(F_2)
\end{tikzcd}
\end{eqnarray*}
where
\begin{itemize}
\item we denoted for short
$A^{Hdg}_{g_{1I,\alpha,\beta}^{\bullet}}:=
\ad(g^{\bullet,!Hdg}_{1I,\alpha,\beta},g^{\bullet}_{1I,\alpha,\beta !Hdg})(h^{!Hdg}_{I\alpha}\mathbb Z^{Hdg}_{\tilde S_I})$
\item we denoted for short
$A^{Hdg}_{g_{2I,\alpha,\beta}^{\bullet}}:=
\ad(g^{\bullet,!Hdg}_{2I,\alpha,\beta},g^{\bullet}_{2I,\alpha,\beta !Hdg})(h^{!Hdg}_{I\alpha}\mathbb Z^{Hdg}_{\tilde S_I})$
\item we denote by $g^n_{1I,\alpha,\beta}:U_{I\alpha}\to U_{I\beta}$, which satisfy 
$h_{I\beta}\circ g^n_{1I,\alpha,\beta}=h_{I\alpha}$, the morphisms in the canonical projective resolution
\begin{eqnarray*}
q:Li_{I*}j_I^*F_1:=(\cdots\to\oplus_{(U_{I\alpha},h_{I\alpha})\in\Var(\mathbb C)^{sm}/\tilde S_I}\mathbb Z(U_{I\alpha}/\tilde S_I)
\xrightarrow{(\mathbb Z(g^{\bullet}_{1,\alpha,\beta}))} \\
\oplus_{(U_{I\alpha},h_{I\alpha})\in\Var(\mathbb C)^{sm}/\tilde S_I}\mathbb Z(U_{I\alpha}/\tilde S_I)\to\cdots)\to i_{I*}j_I^*F_1
\end{eqnarray*}
\item we denote by $g^n_{2I,\alpha,\beta}:U_{I\alpha}\to U_{I\beta}$, which satisfy 
$h_{I\beta}\circ g^n_{2I,\alpha,\beta}=h_{\alpha}$, the morphisms in the canonical projective resolution
\begin{eqnarray*}
q:Li_{I*}j_I^*F_2:=(\cdots\to\oplus_{(U_{I\alpha},h_{I\alpha})\in\Var(\mathbb C)^{sm}/\tilde S_I}\mathbb Z(U_{I\alpha}/\tilde S_I)
\xrightarrow{(\mathbb Z(g^{\bullet}_{2I,\alpha,\beta}))} \\
\oplus_{(U_{I\alpha},h_{I\alpha})\in\Var(\mathbb C)^{sm}/\tilde S_I}\mathbb Z(U_{I\alpha}/\tilde S_I)\to\cdots)\to i_{I*}j_I^*F_2
\end{eqnarray*}
\item we denote by $l_{I\alpha,\beta}^n:U_{I\alpha}\to U_{I\beta}$ which satisfy 
$h_{I\beta}\circ l^n_{I\alpha,\beta}=h_{I\alpha}$ and 
$l^{n+1}_{I\alpha,\beta}\circ g^n_{1I\alpha,\beta}=g^n_{2I\alpha,\beta}\circ l^n_{I\alpha,\beta}$ 
the morphisms in the morphism of canonical projective resolutions
\begin{eqnarray*}
Li_{I*}j_I^*(m):Li_I*j_I^*F_1:=
(\cdots\to\oplus_{(U_{I\alpha},h_{I\alpha})\in\Var(\mathbb C)^{sm}/\tilde S_I}\mathbb Z(U_{I\alpha}/\tilde S_I)\to\cdots) 
\xrightarrow{(\mathbb Z(l^{\bullet}_{I\alpha,\beta}))} \\
(\cdots\to\oplus_{(U_{I\alpha},h_{I\alpha})\in\Var(\mathbb C)^{sm}/\tilde S_I}\mathbb Z(U_{\alpha}/\tilde S_I)\to\cdots)
=:Li_{I*}j_I^*F_2,
\end{eqnarray*}
\item the maps $I^{\bullet}(U_{I\alpha})$ are given by definition \ref{IUSm} and lemma \ref{thetam}.
\end{itemize}
\end{itemize}
Obviously $\mathcal F_S^{Hdg}(F[1])=\mathcal F_S^{Hdg}(F)[1]$ and 
$\mathcal F_S^{Hdg}(\Cone(m))=\Cone(\mathcal F_S^{Hdg}(m))$. 
This functor induces by proposition \ref{projwach} the functor
\begin{eqnarray*}
\mathcal F_S^{Hdg}:=(\mathcal F_S^{FDR},\Bti_S^*):
\DA(S)\to D_{\mathcal D(1,0)fil}(S/(\tilde S_I))\times_I D_{fil}(S^{an}), \\ 
M=D(\mathbb A^1,et)(F)\mapsto
\mathcal F_S^{Hdg}(M):=\mathcal F_S^{Hdg}(F)=(\mathcal F_S^{FDR}(M),\Bti_S^*M,\alpha(M)),
\end{eqnarray*}
with $\alpha(M)=\alpha(F)$. 
\end{defi}

We now give the functoriality with respect to the five operation using the De Rahm realization case and the Betti realization case :

\begin{prop}\label{TbtiTFDR}
\begin{itemize}
\item[(i)]Let $g:T\to S$ a morphism with $T,S\in\Var(\mathbb C)$. 
Assume there exists a factorization $g:T\xrightarrow{l}Y\times S\xrightarrow{p}S$, with $Y\in\SmVar(\mathbb C)$,
$l$ a closed embedding and $p$ the projection.
Let $S=\cup_{i\in I}S_i$ an open cover and 
$i_i:S_i\hookrightarrow\tilde S_i$ closed embeddings with $\tilde S_i\in\SmVar(\mathbb C)$.
Then, $\tilde g_I:Y\times\tilde S_I\to\tilde S_I$ is a lift of $g_I=g_{|T_I}:T_I\to S_I$
and we have closed embeddings $i'_I:=i_I\circ l\circ j'_I:T_I\hookrightarrow Y\times\tilde S_I$.
Then, for $M\in DA_c(S)$, the following diagram commutes :
\begin{equation*}
\xymatrix{g^{*w}\Bti_S^*M\ar[rr]^{g^*(\alpha(M))}\ar[d]_{T(g,bti)(M)} & \, & 
DR(T)^{[-]}((g^{\hat{*}mod}_{Hdg}\mathcal F_S^{DR}(M))^{an})
\ar[d]^{DR(T)^{[-]}((T(g,\mathcal F^{FDR})(M))^{an})} \\
\Bti_T^*g^*M\ar[rr]^{\alpha(g^*M)} &  \, & DR(T)^{[-]}((\mathcal F_T^{DR}(g^*M))^{an})},
\end{equation*}
see section 5, definition \ref{TgDRdefsing} and definition \ref{TgBti}
\item[(ii)]Let $f:T\to S$ a morphism with $T,S\in\QPVar(\mathbb C)$. Then, for $M\in DA_c(T)$,the following diagram commutes :
\begin{equation*}
\xymatrix{Rf_{*w}\Bti_T^*M\ar[rr]^{f_*(\alpha(M)} & \, & DR(S)^{[-]}((Rf^{Hdg}_*\mathcal F_T^{DR}(M))^{an}) \\
\Bti_S^*Rf_*M\ar[u]^{T_*(f,bti)(M)}\ar[rr]^{\alpha(Rf_*M)} & \, &
DR(S)^{[-]}((\mathcal F_S^{DR}(Rf_*M))^{an})\ar[u]_{DR(S)^{[-]}((T_*(f,\mathcal F^{FDR})(M))^{an})}}
\end{equation*}
see section 5, definition \ref{SixTalg} and definition \ref{TBtiSix}
\item[(iii)]Let $f:T\to S$ a morphism with $T,S\in\QPVar(\mathbb C)$. Then, for $M\in DA_c(T)$,the following diagram commutes :
\begin{equation*}
\xymatrix{Rf_{!w}\Bti_T^*M\ar[d]_{T_!(f,bti)(M)}\ar[rr]^{f_!(\alpha(M))} & \, &  
DR(S)^{[-]}((Rf^{Hdg}_!\mathcal F^T_{DR}(M))^{an})\ar[d]^{DR(S)^{[-]}((T_!(f,\mathcal F_{DR})(M))^{an})} \\
\Bti_S^*f_!M\ar[rr]^{\alpha(Rf_!M)} & \, & DR(S)^{[-]}((\mathcal F^S_{DR}(Rf_!M))^{an})} 
\end{equation*}
see section 5, definition \ref{SixTalg} and definition \ref{TBtiSix}.
\item[(iv)]Let $f:T\to S$ a morphism with $T,S\in\QPVar(\mathbb C)$. Then, for $M\in DA_c(S)$,the following diagram commutes :
\begin{equation*}
\xymatrix{f^{!w}\Bti_S^*M\ar[rr]^{f^!(\alpha(M))} & \, & DR(T)^{[-]}((f^{*mod}_{Hdg}\mathcal F_S^{DR}(M))^{an}) \\
\Bti_T^*f^!M\ar[rr]^{\alpha(f^!M)}\ar[u]^{T^!(f,bti)(M)} &  \, &
DR(T)^{[-]}((\mathcal F_T^{DR}(f^!M))^{an})\ar[u]_{DR^{[-]}(T)((T^!(g,\mathcal F^{FDR})(M))^{an})}}
\end{equation*}
see section 5, definition \ref{SixTalg} and definition \ref{TBtiSix}.
\item[(v)] Let $S\in\Var(\mathbb C)$. Then, for $M,N\in DA_c(S)$,the following diagram commutes :
\begin{equation*}
\xymatrix{\Bti_S^*M\otimes\Bti_S^*N\ar[rrr]^{\alpha(M)\otimes\alpha(N)}\ar[d]_{T(\otimes,bti)(M,N)} & \, & \, &  
DR(S)((\mathcal F_S^{DR}(M)\otimes_{O_S}\mathcal F_S^{DR}(N))^{an})\ar[d]^{DR(S)((T(\otimes,\mathcal F^{DR})(M,N))^{an})} \\
\Bti_S^*(M\otimes N)\ar[rrr]^{(\alpha(M\otimes N))} & \, & \, &  DR(S)((\mathcal F^S_{DR}(M\otimes N))^{an})}
\end{equation*}
see definition \ref{SixTalg} and definition \ref{TBtiSix}.
\end{itemize}
\end{prop}

\begin{proof}

\noindent(i): Follows from the following commutative diagram in 
$(D_{\mathcal D(1,0)fil}(Y\times\tilde S_I)\times_I D_{fil}(Y\times\tilde S_I^{an}))$,
\begin{eqnarray*}
\begin{tikzcd}
((\to\oplus_{(U_{I\alpha},h_{I\alpha})\in V_I}\tilde g_I^{*Hdg}h_{I\alpha !Hdg}h^{!Hdg}_{I\alpha}\mathbb Z^{Hdg}_{\tilde S_I}
\xrightarrow{A^{Hdg}_{g_{I,\alpha,\beta}^{\bullet}}} 
\oplus_{(U_{I\alpha},h_{I\alpha})\in V_I}h_{I\alpha !Hdg}h^{!Hdg}_{I\alpha}\mathbb Z^{Hdg}_{\tilde S_I}\to),u_{IJ})
\ar[r,"(\tilde g_I^{*Hdg}I^{\bullet}(U_{I\alpha}/\tilde S_I))"]\ar[d,"T^{Hdg}(\tilde g_I{,}h_I)(-)"'] & 
\shortstack{$(g^{\hat*mod}_{Hdg}\mathcal F_T^{FDR}(F)$, \\ $g^{*w}\Bti_S^*F,g^*(\alpha(F)))$}
\ar[d,"(T(g{,}\mathcal F^{FDR})(M){,}T(g{,}\Bti)(M){,}0)"'] \\
((\to\oplus_{(U'_{I\alpha},h_{I\alpha})\in W_I}h'_{I\alpha !Hdg}h^{'!Hdg}_{I\alpha}\mathbb Z^{Hdg}_{Y\times\tilde S_I}
\xrightarrow{A^{Hdg}_{g_{I,\alpha,\beta}^{'\bullet}}}
\oplus_{(U'_{I\alpha},h'_{I\alpha})\in W_I}h'_{I\alpha !Hdg}h^{'!Hdg}_{I\alpha}\mathbb Z^{Hdg}_{Y\times\tilde S_I}\to),u_{IJ})
\ar[r,"(I^{\bullet}(U'_{\alpha}/Y\times\tilde S_I))"] & 
\shortstack{$(\mathcal F_T^{FDR}(g^*F)$, \\ $\Bti_T^*g^*F,\alpha(g^*F))$}
\end{tikzcd}
\end{eqnarray*}
where, we have denoted for short $V_I:=\Var(\mathbb C)^{sm}/\tilde S_I$ and 
$W_I:=\Var(\mathbb C)^{sm}/Y\times\tilde S_I$,
\begin{itemize}
\item we denoted for short
$A^{Hdg}_{g_{I,\alpha,\beta}^{\bullet}}:=
\ad(g^{\bullet,!Hdg}_{I,\alpha,\beta},g^{\bullet}_{I,\alpha,\beta !Hdg})(h^{!Hdg}_{I\alpha}\mathbb Z^{Hdg}_{\tilde S_I})$
\item we denoted for short
$A^{Hdg}_{g_{I,\alpha,\beta}^{'\bullet}}:=
\ad(g^{'\bullet,!Hdg}_{I,\alpha,\beta},g^{'\bullet}_{I,\alpha,\beta !Hdg})(h^{'!Hdg}_{I\alpha}\mathbb Z^{Hdg}_{Y\times\tilde S_I})$
\item we denote by $g^n_{I,\alpha,\beta}:U_{I\alpha}\to U_{I\beta}$, which satisfy 
$h_{I\beta}\circ g^n_{I,\alpha,\beta}=h_{I\alpha}$, the morphisms in the canonical projective resolution
\begin{eqnarray*}
q:Li_{I*}j_I^*F:=(\cdots\to\oplus_{(U_{I\alpha},h_{I\alpha})\in\Var(\mathbb C)^{sm}/\tilde S_I}\mathbb Z(U_{I\alpha}/\tilde S_I)
\xrightarrow{(\mathbb Z(g^{\bullet}_{\alpha,\beta}))} \\
\oplus_{(U_{I\alpha},h_{I\alpha})\in\Var(\mathbb C)^{sm}/\tilde S_I}\mathbb Z(U_{I\alpha}/\tilde S_I)\to\cdots)\to i_{I*}j_I^*F
\end{eqnarray*}
\item we denote by $g^{'n}_{I,\alpha,\beta}:U'_{I\alpha}\to U'_{I\beta}$, which satisfy 
$h'_{I\beta}\circ g^{'n}_{I,\alpha,\beta}=h'_{\alpha}$, the morphisms in the canonical projective resolution
\begin{eqnarray*}
q:Li'_{I*}j_I^{'*}g^*F:=(\cdots\to\oplus_{(U'_{I\alpha},h'_{I\alpha})\in\Var(\mathbb C)^{sm}/Y\times\tilde S_I}
\mathbb Z(U'_{I\alpha}/Y\times\tilde S_I)
\xrightarrow{(\mathbb Z(g^{'\bullet}_{I,\alpha,\beta}))} \\
\oplus_{(U'_{I\alpha},h'_{I\alpha})\in\Var(\mathbb C)^{sm}/Y\times\tilde S_I}
\mathbb Z(U'_{I\alpha}/Y\times\tilde S_I)\to\cdots)\to i'_{I*}j_I^{'*}g^*F
\end{eqnarray*}
\end{itemize}

\noindent(ii): Follows from (i) by adjonction.

\noindent(iii): The closed embedding case is given by (ii) and the smooth projection case follows from (i) by adjonction.

\noindent(iv): Follows from (iii) by adjonction.

\noindent(v):Obvious

\end{proof}

We can now state the following key proposition and the main theorem:

\begin{prop}\label{keyHdg}
Let $k\subset\mathbb C$ a subfield.
\begin{itemize}
\item[(i)] Let $S\in\Var(\mathbb C)$.
Let $S=\cup_i S_i$ an open cover such that there exist closed embeddings
$i_i:S_i\hookrightarrow\tilde S_i$ with $\tilde S_i\in\SmVar(\mathbb C)$. Then
we have the isomorphism in $D_{\mathcal D(1,0)fil}(S/(\tilde S_I))\times_ID_{fil}(S^{an})$
\begin{eqnarray*}
\mathcal F_S^{Hdg}(\mathbb Z_S)\xrightarrow{:=}
(\mathcal F_S^{FDR}(\mathbb Z_S),\Bti_S^*\mathbb Z_S,\alpha(\mathbb Z_S)) \\
\xrightarrow{((\Omega^{\Gamma,pr}_{/\tilde S_I}(\hat R^{CH}(\ad(i_I^*,i_{I*})(\mathbb Z_{\tilde S_I})^{\gamma}))),I,0)} \\
(e(S)_*\mathcal Hom((\hat R^{CH}(\Gamma^{\vee}_{S_I}\mathbb Z_{\tilde S_I}),\hat R^{CH}(x_{IJ})),
(E_{zar}(\Omega^{\bullet,\Gamma,pr}_{/\tilde S_I},F_{DR}),T_{IJ})),\mathbb Z_{S^{an}},\alpha(S)) \\
\xrightarrow{=}
\iota_S((\Gamma^{\vee,Hdg}_{S_I}(O_{\tilde S_I},F_b),x_{IJ}),\mathbb Z_{S^{an}},\alpha(S))=:\iota_S(\mathbb Z_S^{Hdg})
\end{eqnarray*}
with
\begin{eqnarray*}
\alpha(S):T(S/(\tilde S_I))(\mathbb Z_{S^{an}})\otimes\mathbb C_S:=(i_{I*}j_I^*\mathbb C_{S^{an}},I)
\xrightarrow{(\ad(i_I^*,i_{I*})(\mathbb C_{\tilde S^{an}_I})^q)^{-1}} \\
(\Gamma_{S_I}^{\vee}\mathbb C_{\tilde S_I},x_{IJ})\xrightarrow{(\Gamma_{S_I}^{\vee}\alpha(\tilde S_I))}
DR(S)(o_{fil}(\Gamma^{\vee,Hdg}_{S_I}(O_{\tilde S_I},F_b),x_{IJ}))
\end{eqnarray*}
\item[(ii)]Let $f:X\to S$ a morphism with $X,S\in\Var(\mathbb C)$, $X$ quasi-projective.
Consider a factorization $f:X\xrightarrow{l}Y\times S\xrightarrow{p_S}S$ with $Y=\mathbb P^{N,o}\subset\mathbb P^N$ an open subset,
$l$ a closed embedding and $p_S$ the projection. Let $S=\cup_i S_i$ an open cover such that there exist closed embeddings
$i_i:S_i\hookrightarrow\tilde S_i$ with $\tilde S_i\in\SmVar(\mathbb C)$. 
Recall that $S_I:=\cap_{i\in I} S_i$, $X_I=f^{-1}(S_I)$, and $\tilde S_I:=\Pi_{i\in I}\tilde S_i$. Then,
using proposition \ref{TbtiTFDR}(iii), the maps of definition \ref{SixTalg} and definition \ref{TBtiSix} gives
an isomorphism in $D_{\mathcal D(1,0)fil}(S/(\tilde S_I))\times_ID_{fil}(S^{an})$
\begin{eqnarray*}
(T_!(f,\mathcal F^{FDR})(\mathbb Z_X),T_!(f,\Bti)(\mathbb Z_X),0): \\
\mathcal F_S^{Hdg}(M^{BM}(X/S)):=
(\mathcal F_S^{FDR}(Rf_!\mathbb Z_X),\Bti_S^*Rf_!\mathbb Z_X,\alpha(Rf_!\mathbb Z_X)) \\
\xrightarrow{\sim} 
(Rf_{Hdg!}(\Gamma^{\vee,Hdg}_{X_I}(O_{Y\times\tilde S_I},F_b),x_{IJ}(X/S)),Rf_{!w}\mathbb Z_{X^{an}},
f_!(\alpha(X)))=:\iota_S(Rf_{!Hdg}\mathbb Z_X^{Hdg}).
\end{eqnarray*}
with
\begin{equation*}
\mathbb Z_X^{Hdg}:=((\Gamma^{\vee,Hdg}_{X_I}(O_{Y\times\tilde S_I},F_b),x_{IJ}(X/Y\times S)),
\mathbb Z_{X^{an}},\alpha(X))\in C(MHM(X))
\end{equation*}
\end{itemize}
\end{prop}

\begin{proof}
\noindent(i):Follows from proposition \ref{projwach}.

\noindent(ii): Follows from (i) by proposition \ref{TbtiTFDR}(iii),theorem \ref{mainthm}(i) and theorem \ref{mainBti}(i).
\end{proof}

The main theorem of this section is the following :

\begin{thm}\label{main}
Let $k\subset\mathbb C$ a subfield.
\begin{itemize}
\item[(i)] For $S\in\Var(\mathbb C)$, we have $\mathcal F_S^{Hdg}(\DA_c(S))\subset D(MHM(S))$,  
\begin{equation*}
\iota_S:D(MHM(S))\hookrightarrow D_{\mathcal D(1,0)fil}(S/(\tilde S_I))\times_I D_{fil}(S^{an}) 
\end{equation*}
being a full embedding by theorem \ref{Be}.
\item[(ii)] The Hodge realization functor $\mathcal F_{Hdg}(-)$ define a morphism of 2-functor on $\Var(\mathbb C)$
\begin{equation*}
\mathcal F^{Hdg}_{-}:\Var(\mathbb C)\to(\DA_c(-)\to D(MHM(-)))
\end{equation*}
whose restriction to $\QPVar(\mathbb C)$ is an homotopic 2-functor in sense of Ayoub. More precisely,
\begin{itemize}
\item[(ii0)] for $g:T\to S$ a morphism, with $T,S\in\QPVar(\mathbb C)$, and $M\in\DA_c(S)$, the
the maps of definition \ref{TgDRdefsing} and of definition \ref{TgBti} induce an isomorphism in $D(MHM(T))$
\begin{eqnarray*}
T(g,\mathcal F^{Hdg})(M):=(T(g,\mathcal F^{FDR})(M),T(g,bti)(M),0): \\
g^{\hat*Hdg}\mathcal F_S^{Hdg}(M):=\iota_T^{-1}(g^{\hat{*}mod}_{Hdg}\mathcal F_S^{FDR}(M),g^{*w}\Bti_S(M),g^*(\alpha(M))) \\
\xrightarrow{\sim}\iota_T^{-1}(\mathcal F_T^{FDR}(g^*M),\Bti_T^*(g^*M),\alpha(g^*M))=:\mathcal F_T^{Hdg}(g^*M),
\end{eqnarray*} 
\item[(ii1)] for $f:T\to S$ a morphism, with $T,S\in\QPVar(\mathbb C)$, and $M\in\DA_c(T)$,  
the maps of definition \ref{SixTalg} and of definition \ref{TBtiSix} induce an isomorphism in $D(MHM(S))$
\begin{eqnarray*}
T_*(f,\mathcal F^{Hdg})(M):=(T_*(f,\mathcal F^{FDR})(M),T_*(f,bti)(M),0): \\
Rf_{Hdg*}\mathcal F_T^{Hdg}(M):=\iota_S^{-1}(Rf^{Hdg}_*\mathcal F_T^{FDR}(M),Rf_{*w}\Bti_T(M),f_*(\alpha(M))) \\ 
\xrightarrow{\sim}\iota_S^{-1}(\mathcal F_S^{FDR}(Rf_*M),\Bti_S^*(Rf_*M),\alpha(Rf_*M))=:\mathcal F_S^{Hdg}(Rf_*M),
\end{eqnarray*}  
\item[(ii2)] for $f:T\to S$ a morphism, with $T,S\in\QPVar(\mathbb C)$, and $M\in\DA_c(T)$, 
the maps of definition \ref{SixTalg} and of definition \ref{TBtiSix} induce an isomorphism in $D(MHM(S))$
\begin{eqnarray*}
T_!(f,\mathcal F^{Hdg})(M):=(T_!(f,\mathcal F^{FDR})(M),T_!(f,bti)(M),0): \\
Rf_{!Hdg}\mathcal F_T^{Hdg}(M):=\iota_S^{-1}(Rf^{Hdg}_!\mathcal F_T^{FDR}(M),Rf_{!w}\Bti_T^*(M),f_!(\alpha(M))) \\ 
\xrightarrow{\sim}\iota_S^{-1}(\mathcal F_S^{FDR}(Rf_!M),\Bti_S^*(Rf_!M),\alpha(f_!M))=:\mathcal F_S^{Hdg}(Rf_!M),
\end{eqnarray*} 
\item[(ii3)] for $f:T\to S$ a morphism, with $T,S\in\QPVar(\mathbb C)$, and $M\in\DA_c(S)$,
the maps of definition \ref{SixTalg} and of definition \ref{TBtiSix} induce an isomorphism in $D(MHM(T))$
\begin{eqnarray*}
T^!(f,\mathcal F^{Hdg})(M):=(T^!(f,\mathcal F^{FDR})(M),T^!(f,bti)(M),0): \\
f^{*Hdg}\mathcal F_S^{Hdg}(M):=\iota_T^{-1}(f^{*mod}_{Hdg}\mathcal F_S^{FDR}(M),f^{!w}\Bti_S(M),f^!(\alpha(M))) \\ 
\xrightarrow{\sim}\iota_T^{-1}(\mathcal F_T^{FDR}(f^!M),\Bti_T^*(f^!M),\alpha(f^!M))=:\mathcal F_T^{Hdg}(f^!M),
\end{eqnarray*}
\item[(ii4)] for $S\in\Var(\mathbb C)$, and $M,N\in\DA_c(S)$,
the maps of definition \ref{SixTalg} and of definition \ref{TBtiSix} 
induce an isomorphism in $D(MHM(S))$
\begin{eqnarray*}
T(\otimes,\mathcal F^{Hdg})(M,N):=(T(\otimes,\mathcal F_S^{FDR})(M,N),T(\otimes,bti)(M,N),0): \\
\iota_S^{-1}(\mathcal F_S^{FDR}(M)\otimes^{Hdg}_{O_S}\mathcal F_S^{FDR}(N),
\Bti_S(M)\otimes\Bti_S(N),\alpha(M)\otimes\alpha(N)) \\
\xrightarrow{\sim}\mathcal F_S^{Hdg}(M\otimes N):=
\iota_S^{-1}(\mathcal F_S^{FDR}(M\otimes N),\Bti_S(M\otimes N),\alpha(M\otimes N)).  
\end{eqnarray*}
\end{itemize}
\item[(iii)] For $S\in\Var(\mathbb C)$, the following diagram commutes :
\begin{equation*}  
\xymatrix{\Var(\mathbb C)/S\ar[rrr]^{MH(/S)}\ar[d]_{M(/S)} & \, & \, & D(MHM(S))\ar[d]^{\iota^S} \\
\DA(S)\ar[rrr]^{\mathcal F_S^{Hdg}} & \, & \, & D_{\mathcal D(1,0)fil}(S/(\tilde S_I))\times_I D_{fil}(S^{an})}
\end{equation*}
\end{itemize}
\end{thm}

\begin{proof} 

\noindent(i): Let $M\in\DA_c(S)$. There exist by definition of constructible motives an isomorphism in $\DA(S)$ 
\begin{equation*}
w(M):M\xrightarrow{\sim}\Cone(M(X_0/S)[d_0]\xrightarrow{m_1}\cdots\xrightarrow{m_m} M(X_m/S)[d_m]),
\end{equation*}
with $f_n:X_n\to S$ morphisms and $X_n\in\QPVar(\mathbb C)$.
This gives the isomorphism in $D_{\mathcal D(1,0)fil}(S/(\tilde S_I))\times_I D_{fil}(S^{an})$
\begin{equation*}
\mathcal F_S^{Hdg}(w(M)):\mathcal F_S^{Hdg}(M)\xrightarrow{\sim}
\Cone(\mathcal F_S^{Hdg}(M(X_0/S))[d_0]\xrightarrow{\mathcal F_S^{Hdg}(m_1)}\cdots
\xrightarrow{\mathcal F_S^{Hdg}(m_1)}\mathcal F_S^{Hdg}(M(X_m/S))[d_m]),
\end{equation*}
On the other hand, by proposition \ref{keyHdg}(i), we have 
\begin{equation*}
\mathcal F_S^{Hdg}(M(X_n/S))\xrightarrow{\sim}Rf_{*Hdg}\mathbb Z_X^{Hdg}\in D(MHM(S)).
\end{equation*}
This prove (i).

\noindent(ii0): Follows from theorem \ref{mainthm}(i), proposition \ref{TbtiTFDR}(i) and theorem \ref{mainBti}.

\noindent(ii1): Follows from theorem \ref{mainthm}(iii), proposition \ref{TbtiTFDR}(ii), 
and theorem \ref{mainBti}(iii).

\noindent(ii2):Follows from theorem \ref{mainthm}(ii), proposition \ref{TbtiTFDR}(iii), 
and theorem \ref{mainBti}(ii).

\noindent(ii3): Follows from theorem \ref{mainthm}(iv), proposition \ref{TbtiTFDR}(iv), 
and theorem \ref{mainBti}(iv).

\noindent(ii4):Follows from theorem \ref{mainthm}(v), proposition \ref{TbtiTFDR}(v) and
theorem \ref{mainBti}(v).

\noindent(iii): By (ii), for $g:X'/S\to X/S$ a morphism, with $X',X,S\in\Var(\mathbb C)$ 
and $X/S=(X,f)$, $X'/S=(X',f')$, we have by adjonction the following commutative diagram
\begin{equation*}
\xymatrix{\mathcal F_S^{Hdg}(M(X'/S)=f'_!f^{'!}\mathbb Z_S=f_!g_!g^!f^!\mathbb Z_S)
\ar[d]_{T_!(f',\mathcal F^{Hdg})(f^{'!}M(X'/S))\circ T^!(f',\mathcal F^{Hdg})(M(X'/S))}
\ar[rr]^{\mathcal F_S^{Hdg}(M(/S)(g)=f_!\ad(g_!,g^!)(f^!\mathbb Z_S))} & \, & 
\mathcal F_S^{Hdg}(M(X/S)=f_!f^!\mathbb Z_S)
\ar[d]_{T_!(f,\mathcal F^{Hdg})(f^!M(X/S))\circ T^!(f,\mathcal F^{Hdg})(M(X/S))} \\
MH(X'/S):=Rf'_{!Hdg}f^{'!Hdg}\mathbb Z^{Hdg}_S=f_{!Hdg}g_{!Hdg}g^{!Hdg}f^{!Hdg}\mathbb Z^{Hdg}_S
\ar[rr]^{f_{!Hdg}\ad(g_{!Hdg},g^{!Hdg})(f^{!Hdg}\mathbb Z^{Hdg}_S)} & \, &  
MH(X/S):=f_{!Hdg}f^{!Hdg}\mathbb Z^{Hdg}_S}.
\end{equation*}
where the left and right columns are isomorphisms by (ii). This proves (iii).
\end{proof}

The theorem \ref{main} gives immediately the following :

\begin{cor}
Let $f:U\to S$, $f':U'\to S$ morphisms, with $U,U',S\in\Var(\mathbb C)$ irreducible, $U'$ smooth.
Let $\bar{S}\in\PVar(\mathbb C)$ a compactification of $S$.
Let $\bar{X},\bar{X'}\in\PVar(\mathbb C)$ compactification of $U$ and $U'$ respectively,
such that $f$ (resp. $f'$) extend to a morphism $\bar f:\bar X\to\bar S$, resp. $\bar{f'}:\bar{X'}\to\bar S$.
Denote $\bar D=\bar X\backslash U$ and $\bar D'=\bar{X'}\backslash U'$ and 
$\bar E=(\bar D\times_{\bar S}\bar X')\cup(\bar X\times_{\bar S}\bar D')$.
Denote $i:\bar D\hookrightarrow\bar X$, $i':\bar D\hookrightarrow\bar X$ denote the closed embeddings
and $j:U\hookrightarrow\bar X$, $j':U'\hookrightarrow\bar X'$ the open embeddings.
Denote $d=\dim(U)$ and $d'=\dim(U')$.
We have the following commutative diagram in $D(\mathbb Z)$
\begin{equation*}
\xymatrix{RHom_{\DA(\bar S)}^{\bullet}(M(U'/\bar S),M((\bar X,\bar D)/\bar S))
\ar[d]^{RI(-,-)}\ar[rr]^{{\mathcal F_S^{FDR}}^{(-,-)}} & \, & 
RHom_{DMHM(\bar S)}^{\bullet}(f'_{!Hdg}\mathbb Z_{U'}^{Hdg},f_{*Hdg}\mathbb Z_U^{Hdg})
\ar[d]^{RI(-,-)} \\
RHom^{\bullet}(M(\pt),M(\bar X'\times_{\bar S}\bar X,\bar E)(d')[2d'])
\ar[d]^{l}\ar[rr]^{{\mathcal F^{\pt}_{FDR}}^{(-,-)}} & \, & 
RHom^{\bullet}(\mathbb Z_{\pt}^{Hdg},a_{U'\times_SU!}\mathbb Z_{U\times_SU'}^{Hdg}(d')[2d'])
\ar[d]^{l} \\
\mathcal Z_d(\bar X'\times_{\bar S}\bar X,E,\bullet)
\ar[rr]^{\mathcal R^d_{\bar X'\times_{\bar S}\bar X}} & \, & 
C^{\mathcal D}_{2d+\bullet}(\bar X'\times_{\bar S}\bar X,E,Z(d))}
\end{equation*}
where 
\begin{equation*}
M((\bar X,\bar D)/\bar S):=\Cone(\ad(i_*,i^!):M(\bar D/\bar S)\to M(\bar X/\bar S))
=\bar f_*j_*E_{et}(\mathbb Z(U/U))\in\DA(\bar S)
\end{equation*}
and $l$ the isomorphisms given by canonical embedding of complexes.
\end{cor}

\begin{proof}
The upper square of this diagram follows from theorem \ref{main}(ii).
On the other side, the lower square follows from the absolute case.
\end{proof}


LAGA UMR CNRS 7539 \\
Universit\'e Paris 13, Sorbonne Paris Cit\'e, 99 av Jean-Baptiste Clement, \\
93430 Villetaneuse, France, \\
bouali@math.univ-paris13.fr

\end{document}